\RequirePackage[l2tabu,orthodox,abort]{nag} 
\RequirePackage{fix-cm} 

\documentclass[twoside,openright,titlepage,numbers=noenddot,headinclude,
               footinclude=true,cleardoublepage=empty,abstractoff, 
               BCOR=5mm,paper=a4,fontsize=11pt,
               british,
               ]{scrreprt}

\PassOptionsToPackage{pdfspacing, floatperchapter, linedheaders, subfig, manychapters}{classicthesis}

\newcommand{\myTitle}{Cellular Automata on Group Sets\xspace}
\newcommand{\myName}{Simon Wacker\xspace}

\newcommand{\myFaculty}{Fakultät für Informatik\xspace}

\newcommand{\myUni}{Karlsruher Institut für Technologie\xspace}

\pdfglyphtounicode{A}{0041}
\pdfglyphtounicode{AE}{00C6}
\pdfglyphtounicode{AEacute}{01FC}
\pdfglyphtounicode{AEmacron}{01E2}
\pdfglyphtounicode{AEsmall}{00E6}
\pdfglyphtounicode{Aacute}{00C1}
\pdfglyphtounicode{Aacutesmall}{00E1}
\pdfglyphtounicode{Abreve}{0102}
\pdfglyphtounicode{Abreveacute}{1EAE}
\pdfglyphtounicode{Abrevecyrillic}{04D0}
\pdfglyphtounicode{Abrevedotbelow}{1EB6}
\pdfglyphtounicode{Abrevegrave}{1EB0}
\pdfglyphtounicode{Abrevehookabove}{1EB2}
\pdfglyphtounicode{Abrevetilde}{1EB4}
\pdfglyphtounicode{Acaron}{01CD}
\pdfglyphtounicode{Acircle}{24B6}
\pdfglyphtounicode{Acircumflex}{00C2}
\pdfglyphtounicode{Acircumflexacute}{1EA4}
\pdfglyphtounicode{Acircumflexdotbelow}{1EAC}
\pdfglyphtounicode{Acircumflexgrave}{1EA6}
\pdfglyphtounicode{Acircumflexhookabove}{1EA8}
\pdfglyphtounicode{Acircumflexsmall}{00E2}
\pdfglyphtounicode{Acircumflextilde}{1EAA}
\pdfglyphtounicode{Acute}{00B4}
\pdfglyphtounicode{Acutesmall}{00B4}
\pdfglyphtounicode{Acyrillic}{0410}
\pdfglyphtounicode{Adblgrave}{0200}
\pdfglyphtounicode{Adieresis}{00C4}
\pdfglyphtounicode{Adieresiscyrillic}{04D2}
\pdfglyphtounicode{Adieresismacron}{01DE}
\pdfglyphtounicode{Adieresissmall}{00E4}
\pdfglyphtounicode{Adotbelow}{1EA0}
\pdfglyphtounicode{Adotmacron}{01E0}
\pdfglyphtounicode{Agrave}{00C0}
\pdfglyphtounicode{Agravesmall}{00E0}
\pdfglyphtounicode{Ahookabove}{1EA2}
\pdfglyphtounicode{Aiecyrillic}{04D4}
\pdfglyphtounicode{Ainvertedbreve}{0202}
\pdfglyphtounicode{Alpha}{0391}
\pdfglyphtounicode{Alphatonos}{0386}
\pdfglyphtounicode{Amacron}{0100}
\pdfglyphtounicode{Amonospace}{FF21}
\pdfglyphtounicode{Aogonek}{0104}
\pdfglyphtounicode{Aring}{00C5}
\pdfglyphtounicode{Aringacute}{01FA}
\pdfglyphtounicode{Aringbelow}{1E00}
\pdfglyphtounicode{Aringsmall}{00E5}
\pdfglyphtounicode{Asmall}{0061}
\pdfglyphtounicode{Atilde}{00C3}
\pdfglyphtounicode{Atildesmall}{00E3}
\pdfglyphtounicode{Aybarmenian}{0531}
\pdfglyphtounicode{B}{0042}
\pdfglyphtounicode{Bcircle}{24B7}
\pdfglyphtounicode{Bdotaccent}{1E02}
\pdfglyphtounicode{Bdotbelow}{1E04}
\pdfglyphtounicode{Becyrillic}{0411}
\pdfglyphtounicode{Benarmenian}{0532}
\pdfglyphtounicode{Beta}{0392}
\pdfglyphtounicode{Bhook}{0181}
\pdfglyphtounicode{Blinebelow}{1E06}
\pdfglyphtounicode{Bmonospace}{FF22}
\pdfglyphtounicode{Brevesmall}{02D8}
\pdfglyphtounicode{Bsmall}{0062}
\pdfglyphtounicode{Btopbar}{0182}
\pdfglyphtounicode{C}{0043}
\pdfglyphtounicode{Caarmenian}{053E}
\pdfglyphtounicode{Cacute}{0106}
\pdfglyphtounicode{Caron}{02C7}
\pdfglyphtounicode{Caronsmall}{02C7}
\pdfglyphtounicode{Ccaron}{010C}
\pdfglyphtounicode{Ccedilla}{00C7}
\pdfglyphtounicode{Ccedillaacute}{1E08}
\pdfglyphtounicode{Ccedillasmall}{00E7}
\pdfglyphtounicode{Ccircle}{24B8}
\pdfglyphtounicode{Ccircumflex}{0108}
\pdfglyphtounicode{Cdot}{010A}
\pdfglyphtounicode{Cdotaccent}{010A}
\pdfglyphtounicode{Cedillasmall}{00B8}
\pdfglyphtounicode{Chaarmenian}{0549}
\pdfglyphtounicode{Cheabkhasiancyrillic}{04BC}
\pdfglyphtounicode{Checyrillic}{0427}
\pdfglyphtounicode{Chedescenderabkhasiancyrillic}{04BE}
\pdfglyphtounicode{Chedescendercyrillic}{04B6}
\pdfglyphtounicode{Chedieresiscyrillic}{04F4}
\pdfglyphtounicode{Cheharmenian}{0543}
\pdfglyphtounicode{Chekhakassiancyrillic}{04CB}
\pdfglyphtounicode{Cheverticalstrokecyrillic}{04B8}
\pdfglyphtounicode{Chi}{03A7}
\pdfglyphtounicode{Chook}{0187}
\pdfglyphtounicode{Circumflexsmall}{02C6}
\pdfglyphtounicode{Cmonospace}{FF23}
\pdfglyphtounicode{Coarmenian}{0551}
\pdfglyphtounicode{Csmall}{0063}
\pdfglyphtounicode{D}{0044}
\pdfglyphtounicode{DZ}{01F1}
\pdfglyphtounicode{DZcaron}{01C4}
\pdfglyphtounicode{Daarmenian}{0534}
\pdfglyphtounicode{Dafrican}{0189}
\pdfglyphtounicode{Dbar}{0110}
\pdfglyphtounicode{Dcaron}{010E}
\pdfglyphtounicode{Dcedilla}{1E10}
\pdfglyphtounicode{Dcircle}{24B9}
\pdfglyphtounicode{Dcircumflexbelow}{1E12}
\pdfglyphtounicode{Dcroat}{0110}
\pdfglyphtounicode{Ddotaccent}{1E0A}
\pdfglyphtounicode{Ddotbelow}{1E0C}
\pdfglyphtounicode{Decyrillic}{0414}
\pdfglyphtounicode{Deicoptic}{03EE}
\pdfglyphtounicode{Delta}{2206}
\pdfglyphtounicode{Deltagreek}{0394}
\pdfglyphtounicode{Dhook}{018A}
\pdfglyphtounicode{Dieresis}{00A8}
\pdfglyphtounicode{DieresisAcute}{F6CC}
\pdfglyphtounicode{DieresisGrave}{F6CD}
\pdfglyphtounicode{Dieresissmall}{00A8}
\pdfglyphtounicode{Digamma}{D875 DFCB}
\pdfglyphtounicode{Digammagreek}{03DC}
\pdfglyphtounicode{Djecyrillic}{0402}
\pdfglyphtounicode{Dlinebelow}{1E0E}
\pdfglyphtounicode{Dmonospace}{FF24}
\pdfglyphtounicode{Dotaccentsmall}{02D9}
\pdfglyphtounicode{Dslash}{0110}
\pdfglyphtounicode{Dsmall}{0064}
\pdfglyphtounicode{Dtopbar}{018B}
\pdfglyphtounicode{Dz}{01F2}
\pdfglyphtounicode{Dzcaron}{01C5}
\pdfglyphtounicode{Dzeabkhasiancyrillic}{04E0}
\pdfglyphtounicode{Dzecyrillic}{0405}
\pdfglyphtounicode{Dzhecyrillic}{040F}
\pdfglyphtounicode{E}{0045}
\pdfglyphtounicode{Eacute}{00C9}
\pdfglyphtounicode{Eacutesmall}{00E9}
\pdfglyphtounicode{Ebreve}{0114}
\pdfglyphtounicode{Ecaron}{011A}
\pdfglyphtounicode{Ecedillabreve}{1E1C}
\pdfglyphtounicode{Echarmenian}{0535}
\pdfglyphtounicode{Ecircle}{24BA}
\pdfglyphtounicode{Ecircumflex}{00CA}
\pdfglyphtounicode{Ecircumflexacute}{1EBE}
\pdfglyphtounicode{Ecircumflexbelow}{1E18}
\pdfglyphtounicode{Ecircumflexdotbelow}{1EC6}
\pdfglyphtounicode{Ecircumflexgrave}{1EC0}
\pdfglyphtounicode{Ecircumflexhookabove}{1EC2}
\pdfglyphtounicode{Ecircumflexsmall}{00EA}
\pdfglyphtounicode{Ecircumflextilde}{1EC4}
\pdfglyphtounicode{Ecyrillic}{0404}
\pdfglyphtounicode{Edblgrave}{0204}
\pdfglyphtounicode{Edieresis}{00CB}
\pdfglyphtounicode{Edieresissmall}{00EB}
\pdfglyphtounicode{Edot}{0116}
\pdfglyphtounicode{Edotaccent}{0116}
\pdfglyphtounicode{Edotbelow}{1EB8}
\pdfglyphtounicode{Efcyrillic}{0424}
\pdfglyphtounicode{Egrave}{00C8}
\pdfglyphtounicode{Egravesmall}{00E8}
\pdfglyphtounicode{Eharmenian}{0537}
\pdfglyphtounicode{Ehookabove}{1EBA}
\pdfglyphtounicode{Eightroman}{2167}
\pdfglyphtounicode{Einvertedbreve}{0206}
\pdfglyphtounicode{Eiotifiedcyrillic}{0464}
\pdfglyphtounicode{Elcyrillic}{041B}
\pdfglyphtounicode{Elevenroman}{216A}
\pdfglyphtounicode{Emacron}{0112}
\pdfglyphtounicode{Emacronacute}{1E16}
\pdfglyphtounicode{Emacrongrave}{1E14}
\pdfglyphtounicode{Emcyrillic}{041C}
\pdfglyphtounicode{Emonospace}{FF25}
\pdfglyphtounicode{Encyrillic}{041D}
\pdfglyphtounicode{Endescendercyrillic}{04A2}
\pdfglyphtounicode{Eng}{014A}
\pdfglyphtounicode{Enghecyrillic}{04A4}
\pdfglyphtounicode{Enhookcyrillic}{04C7}
\pdfglyphtounicode{Eogonek}{0118}
\pdfglyphtounicode{Eopen}{0190}
\pdfglyphtounicode{Epsilon}{0395}
\pdfglyphtounicode{Epsilontonos}{0388}
\pdfglyphtounicode{Ercyrillic}{0420}
\pdfglyphtounicode{Ereversed}{018E}
\pdfglyphtounicode{Ereversedcyrillic}{042D}
\pdfglyphtounicode{Escyrillic}{0421}
\pdfglyphtounicode{Esdescendercyrillic}{04AA}
\pdfglyphtounicode{Esh}{01A9}
\pdfglyphtounicode{Esmall}{0065}
\pdfglyphtounicode{Eta}{0397}
\pdfglyphtounicode{Etarmenian}{0538}
\pdfglyphtounicode{Etatonos}{0389}
\pdfglyphtounicode{Eth}{00D0}
\pdfglyphtounicode{Ethsmall}{00F0}
\pdfglyphtounicode{Etilde}{1EBC}
\pdfglyphtounicode{Etildebelow}{1E1A}
\pdfglyphtounicode{Euro}{20AC}
\pdfglyphtounicode{Ezh}{01B7}
\pdfglyphtounicode{Ezhcaron}{01EE}
\pdfglyphtounicode{Ezhreversed}{01B8}
\pdfglyphtounicode{F}{0046}
\pdfglyphtounicode{FFIsmall}{0066 0066 0069}
\pdfglyphtounicode{FFLsmall}{0066 0066 006C}
\pdfglyphtounicode{FFsmall}{0066 0066}
\pdfglyphtounicode{FIsmall}{0066 0069}
\pdfglyphtounicode{FLsmall}{0066 006C}
\pdfglyphtounicode{Fcircle}{24BB}
\pdfglyphtounicode{Fdotaccent}{1E1E}
\pdfglyphtounicode{Feharmenian}{0556}
\pdfglyphtounicode{Feicoptic}{03E4}
\pdfglyphtounicode{Fhook}{0191}
\pdfglyphtounicode{Finv}{2132}
\pdfglyphtounicode{Fitacyrillic}{0472}
\pdfglyphtounicode{Fiveroman}{2164}
\pdfglyphtounicode{Fmonospace}{FF26}
\pdfglyphtounicode{Fourroman}{2163}
\pdfglyphtounicode{Fsmall}{0066}
\pdfglyphtounicode{G}{0047}
\pdfglyphtounicode{GBsquare}{3387}
\pdfglyphtounicode{Gacute}{01F4}
\pdfglyphtounicode{Gamma}{0393}
\pdfglyphtounicode{Gammaafrican}{0194}
\pdfglyphtounicode{Gangiacoptic}{03EA}
\pdfglyphtounicode{Gbreve}{011E}
\pdfglyphtounicode{Gcaron}{01E6}
\pdfglyphtounicode{Gcedilla}{0122}
\pdfglyphtounicode{Gcircle}{24BC}
\pdfglyphtounicode{Gcircumflex}{011C}
\pdfglyphtounicode{Gcommaaccent}{0122}
\pdfglyphtounicode{Gdot}{0120}
\pdfglyphtounicode{Gdotaccent}{0120}
\pdfglyphtounicode{Gecyrillic}{0413}
\pdfglyphtounicode{Germandbls}{0053 0053}
\pdfglyphtounicode{Germandblssmall}{0073 0073}
\pdfglyphtounicode{Ghadarmenian}{0542}
\pdfglyphtounicode{Ghemiddlehookcyrillic}{0494}
\pdfglyphtounicode{Ghestrokecyrillic}{0492}
\pdfglyphtounicode{Gheupturncyrillic}{0490}
\pdfglyphtounicode{Ghook}{0193}
\pdfglyphtounicode{Gimarmenian}{0533}
\pdfglyphtounicode{Gjecyrillic}{0403}
\pdfglyphtounicode{Gmacron}{1E20}
\pdfglyphtounicode{Gmir}{2141}
\pdfglyphtounicode{Gmonospace}{FF27}
\pdfglyphtounicode{Grave}{0060}
\pdfglyphtounicode{Gravesmall}{0060}
\pdfglyphtounicode{Gsmall}{0067}
\pdfglyphtounicode{Gsmallhook}{029B}
\pdfglyphtounicode{Gstroke}{01E4}
\pdfglyphtounicode{H}{0048}
\pdfglyphtounicode{H18533}{25CF}
\pdfglyphtounicode{H18543}{25AA}
\pdfglyphtounicode{H18551}{25AB}
\pdfglyphtounicode{H22073}{25A1}
\pdfglyphtounicode{HPsquare}{33CB}
\pdfglyphtounicode{Haabkhasiancyrillic}{04A8}
\pdfglyphtounicode{Hadescendercyrillic}{04B2}
\pdfglyphtounicode{Hardsigncyrillic}{042A}
\pdfglyphtounicode{Hbar}{0126}
\pdfglyphtounicode{Hbrevebelow}{1E2A}
\pdfglyphtounicode{Hcedilla}{1E28}
\pdfglyphtounicode{Hcircle}{24BD}
\pdfglyphtounicode{Hcircumflex}{0124}
\pdfglyphtounicode{Hdieresis}{1E26}
\pdfglyphtounicode{Hdotaccent}{1E22}
\pdfglyphtounicode{Hdotbelow}{1E24}
\pdfglyphtounicode{Hmonospace}{FF28}
\pdfglyphtounicode{Hoarmenian}{0540}
\pdfglyphtounicode{Horicoptic}{03E8}
\pdfglyphtounicode{Hsmall}{0068}
\pdfglyphtounicode{Hungarumlaut}{02DD}
\pdfglyphtounicode{Hungarumlautsmall}{02DD}
\pdfglyphtounicode{Hzsquare}{3390}
\pdfglyphtounicode{I}{0049}
\pdfglyphtounicode{IAcyrillic}{042F}
\pdfglyphtounicode{IJ}{0132}
\pdfglyphtounicode{IUcyrillic}{042E}
\pdfglyphtounicode{Iacute}{00CD}
\pdfglyphtounicode{Iacutesmall}{00ED}
\pdfglyphtounicode{Ibreve}{012C}
\pdfglyphtounicode{Icaron}{01CF}
\pdfglyphtounicode{Icircle}{24BE}
\pdfglyphtounicode{Icircumflex}{00CE}
\pdfglyphtounicode{Icircumflexsmall}{00EE}
\pdfglyphtounicode{Icyrillic}{0406}
\pdfglyphtounicode{Idblgrave}{0208}
\pdfglyphtounicode{Idieresis}{00CF}
\pdfglyphtounicode{Idieresisacute}{1E2E}
\pdfglyphtounicode{Idieresiscyrillic}{04E4}
\pdfglyphtounicode{Idieresissmall}{00EF}
\pdfglyphtounicode{Idot}{0130}
\pdfglyphtounicode{Idotaccent}{0130}
\pdfglyphtounicode{Idotbelow}{1ECA}
\pdfglyphtounicode{Iebrevecyrillic}{04D6}
\pdfglyphtounicode{Iecyrillic}{0415}
\pdfglyphtounicode{Ifractur}{2111}
\pdfglyphtounicode{Ifraktur}{2111}
\pdfglyphtounicode{Igrave}{00CC}
\pdfglyphtounicode{Igravesmall}{00EC}
\pdfglyphtounicode{Ihookabove}{1EC8}
\pdfglyphtounicode{Iicyrillic}{0418}
\pdfglyphtounicode{Iinvertedbreve}{020A}
\pdfglyphtounicode{Iishortcyrillic}{0419}
\pdfglyphtounicode{Imacron}{012A}
\pdfglyphtounicode{Imacroncyrillic}{04E2}
\pdfglyphtounicode{Imonospace}{FF29}
\pdfglyphtounicode{Iniarmenian}{053B}
\pdfglyphtounicode{Iocyrillic}{0401}
\pdfglyphtounicode{Iogonek}{012E}
\pdfglyphtounicode{Iota}{0399}
\pdfglyphtounicode{Iotaafrican}{0196}
\pdfglyphtounicode{Iotadieresis}{03AA}
\pdfglyphtounicode{Iotatonos}{038A}
\pdfglyphtounicode{Ismall}{0069}
\pdfglyphtounicode{Istroke}{0197}
\pdfglyphtounicode{Itilde}{0128}
\pdfglyphtounicode{Itildebelow}{1E2C}
\pdfglyphtounicode{Izhitsacyrillic}{0474}
\pdfglyphtounicode{Izhitsadblgravecyrillic}{0476}
\pdfglyphtounicode{J}{004A}
\pdfglyphtounicode{Jaarmenian}{0541}
\pdfglyphtounicode{Jcircle}{24BF}
\pdfglyphtounicode{Jcircumflex}{0134}
\pdfglyphtounicode{Jecyrillic}{0408}
\pdfglyphtounicode{Jheharmenian}{054B}
\pdfglyphtounicode{Jmonospace}{FF2A}
\pdfglyphtounicode{Jsmall}{006A}
\pdfglyphtounicode{K}{004B}
\pdfglyphtounicode{KBsquare}{3385}
\pdfglyphtounicode{KKsquare}{33CD}
\pdfglyphtounicode{Kabashkircyrillic}{04A0}
\pdfglyphtounicode{Kacute}{1E30}
\pdfglyphtounicode{Kacyrillic}{041A}
\pdfglyphtounicode{Kadescendercyrillic}{049A}
\pdfglyphtounicode{Kahookcyrillic}{04C3}
\pdfglyphtounicode{Kappa}{039A}
\pdfglyphtounicode{Kastrokecyrillic}{049E}
\pdfglyphtounicode{Kaverticalstrokecyrillic}{049C}
\pdfglyphtounicode{Kcaron}{01E8}
\pdfglyphtounicode{Kcedilla}{0136}
\pdfglyphtounicode{Kcircle}{24C0}
\pdfglyphtounicode{Kcommaaccent}{0136}
\pdfglyphtounicode{Kdotbelow}{1E32}
\pdfglyphtounicode{Keharmenian}{0554}
\pdfglyphtounicode{Kenarmenian}{053F}
\pdfglyphtounicode{Khacyrillic}{0425}
\pdfglyphtounicode{Kheicoptic}{03E6}
\pdfglyphtounicode{Khook}{0198}
\pdfglyphtounicode{Kjecyrillic}{040C}
\pdfglyphtounicode{Klinebelow}{1E34}
\pdfglyphtounicode{Kmonospace}{FF2B}
\pdfglyphtounicode{Koppacyrillic}{0480}
\pdfglyphtounicode{Koppagreek}{03DE}
\pdfglyphtounicode{Ksicyrillic}{046E}
\pdfglyphtounicode{Ksmall}{006B}
\pdfglyphtounicode{L}{004C}
\pdfglyphtounicode{LJ}{01C7}
\pdfglyphtounicode{LL}{004C 004C}
\pdfglyphtounicode{Lacute}{0139}
\pdfglyphtounicode{Lambda}{039B}
\pdfglyphtounicode{Lcaron}{013D}
\pdfglyphtounicode{Lcedilla}{013B}
\pdfglyphtounicode{Lcircle}{24C1}
\pdfglyphtounicode{Lcircumflexbelow}{1E3C}
\pdfglyphtounicode{Lcommaaccent}{013B}
\pdfglyphtounicode{Ldot}{013F}
\pdfglyphtounicode{Ldotaccent}{013F}
\pdfglyphtounicode{Ldotbelow}{1E36}
\pdfglyphtounicode{Ldotbelowmacron}{1E38}
\pdfglyphtounicode{Liwnarmenian}{053C}
\pdfglyphtounicode{Lj}{01C8}
\pdfglyphtounicode{Ljecyrillic}{0409}
\pdfglyphtounicode{Llinebelow}{1E3A}
\pdfglyphtounicode{Lmonospace}{FF2C}
\pdfglyphtounicode{Lslash}{0141}
\pdfglyphtounicode{Lslashsmall}{0142}
\pdfglyphtounicode{Lsmall}{006C}
\pdfglyphtounicode{M}{004D}
\pdfglyphtounicode{MBsquare}{3386}
\pdfglyphtounicode{Macron}{00AF}
\pdfglyphtounicode{Macronsmall}{00AF}
\pdfglyphtounicode{Macute}{1E3E}
\pdfglyphtounicode{Mcircle}{24C2}
\pdfglyphtounicode{Mdotaccent}{1E40}
\pdfglyphtounicode{Mdotbelow}{1E42}
\pdfglyphtounicode{Menarmenian}{0544}
\pdfglyphtounicode{Mmonospace}{FF2D}
\pdfglyphtounicode{Msmall}{006D}
\pdfglyphtounicode{Mturned}{019C}
\pdfglyphtounicode{Mu}{039C}
\pdfglyphtounicode{N}{004E}
\pdfglyphtounicode{NJ}{01CA}
\pdfglyphtounicode{Nacute}{0143}
\pdfglyphtounicode{Ncaron}{0147}
\pdfglyphtounicode{Ncedilla}{0145}
\pdfglyphtounicode{Ncircle}{24C3}
\pdfglyphtounicode{Ncircumflexbelow}{1E4A}
\pdfglyphtounicode{Ncommaaccent}{0145}
\pdfglyphtounicode{Ndotaccent}{1E44}
\pdfglyphtounicode{Ndotbelow}{1E46}
\pdfglyphtounicode{Ng}{014A}
\pdfglyphtounicode{Nhookleft}{019D}
\pdfglyphtounicode{Nineroman}{2168}
\pdfglyphtounicode{Nj}{01CB}
\pdfglyphtounicode{Njecyrillic}{040A}
\pdfglyphtounicode{Nlinebelow}{1E48}
\pdfglyphtounicode{Nmonospace}{FF2E}
\pdfglyphtounicode{Nowarmenian}{0546}
\pdfglyphtounicode{Nsmall}{006E}
\pdfglyphtounicode{Ntilde}{00D1}
\pdfglyphtounicode{Ntildesmall}{00F1}
\pdfglyphtounicode{Nu}{039D}
\pdfglyphtounicode{O}{004F}
\pdfglyphtounicode{OE}{0152}
\pdfglyphtounicode{OEsmall}{0153}
\pdfglyphtounicode{Oacute}{00D3}
\pdfglyphtounicode{Oacutesmall}{00F3}
\pdfglyphtounicode{Obarredcyrillic}{04E8}
\pdfglyphtounicode{Obarreddieresiscyrillic}{04EA}
\pdfglyphtounicode{Obreve}{014E}
\pdfglyphtounicode{Ocaron}{01D1}
\pdfglyphtounicode{Ocenteredtilde}{019F}
\pdfglyphtounicode{Ocircle}{24C4}
\pdfglyphtounicode{Ocircumflex}{00D4}
\pdfglyphtounicode{Ocircumflexacute}{1ED0}
\pdfglyphtounicode{Ocircumflexdotbelow}{1ED8}
\pdfglyphtounicode{Ocircumflexgrave}{1ED2}
\pdfglyphtounicode{Ocircumflexhookabove}{1ED4}
\pdfglyphtounicode{Ocircumflexsmall}{00F4}
\pdfglyphtounicode{Ocircumflextilde}{1ED6}
\pdfglyphtounicode{Ocyrillic}{041E}
\pdfglyphtounicode{Odblacute}{0150}
\pdfglyphtounicode{Odblgrave}{020C}
\pdfglyphtounicode{Odieresis}{00D6}
\pdfglyphtounicode{Odieresiscyrillic}{04E6}
\pdfglyphtounicode{Odieresissmall}{00F6}
\pdfglyphtounicode{Odotbelow}{1ECC}
\pdfglyphtounicode{Ogoneksmall}{02DB}
\pdfglyphtounicode{Ograve}{00D2}
\pdfglyphtounicode{Ogravesmall}{00F2}
\pdfglyphtounicode{Oharmenian}{0555}
\pdfglyphtounicode{Ohm}{2126}
\pdfglyphtounicode{Ohookabove}{1ECE}
\pdfglyphtounicode{Ohorn}{01A0}
\pdfglyphtounicode{Ohornacute}{1EDA}
\pdfglyphtounicode{Ohorndotbelow}{1EE2}
\pdfglyphtounicode{Ohorngrave}{1EDC}
\pdfglyphtounicode{Ohornhookabove}{1EDE}
\pdfglyphtounicode{Ohorntilde}{1EE0}
\pdfglyphtounicode{Ohungarumlaut}{0150}
\pdfglyphtounicode{Oi}{01A2}
\pdfglyphtounicode{Oinvertedbreve}{020E}
\pdfglyphtounicode{Omacron}{014C}
\pdfglyphtounicode{Omacronacute}{1E52}
\pdfglyphtounicode{Omacrongrave}{1E50}
\pdfglyphtounicode{Omega}{2126}
\pdfglyphtounicode{Omegacyrillic}{0460}
\pdfglyphtounicode{Omegagreek}{03A9}
\pdfglyphtounicode{Omegainv}{2127}
\pdfglyphtounicode{Omegaroundcyrillic}{047A}
\pdfglyphtounicode{Omegatitlocyrillic}{047C}
\pdfglyphtounicode{Omegatonos}{038F}
\pdfglyphtounicode{Omicron}{039F}
\pdfglyphtounicode{Omicrontonos}{038C}
\pdfglyphtounicode{Omonospace}{FF2F}
\pdfglyphtounicode{Oneroman}{2160}
\pdfglyphtounicode{Oogonek}{01EA}
\pdfglyphtounicode{Oogonekmacron}{01EC}
\pdfglyphtounicode{Oopen}{0186}
\pdfglyphtounicode{Oslash}{00D8}
\pdfglyphtounicode{Oslashacute}{01FE}
\pdfglyphtounicode{Oslashsmall}{00F8}
\pdfglyphtounicode{Osmall}{006F}
\pdfglyphtounicode{Ostrokeacute}{01FE}
\pdfglyphtounicode{Otcyrillic}{047E}
\pdfglyphtounicode{Otilde}{00D5}
\pdfglyphtounicode{Otildeacute}{1E4C}
\pdfglyphtounicode{Otildedieresis}{1E4E}
\pdfglyphtounicode{Otildesmall}{00F5}
\pdfglyphtounicode{P}{0050}
\pdfglyphtounicode{Pacute}{1E54}
\pdfglyphtounicode{Pcircle}{24C5}
\pdfglyphtounicode{Pdotaccent}{1E56}
\pdfglyphtounicode{Pecyrillic}{041F}
\pdfglyphtounicode{Peharmenian}{054A}
\pdfglyphtounicode{Pemiddlehookcyrillic}{04A6}
\pdfglyphtounicode{Phi}{03A6}
\pdfglyphtounicode{Phook}{01A4}
\pdfglyphtounicode{Pi}{03A0}
\pdfglyphtounicode{Piwrarmenian}{0553}
\pdfglyphtounicode{Pmonospace}{FF30}
\pdfglyphtounicode{Psi}{03A8}
\pdfglyphtounicode{Psicyrillic}{0470}
\pdfglyphtounicode{Psmall}{0070}
\pdfglyphtounicode{Q}{0051}
\pdfglyphtounicode{Qcircle}{24C6}
\pdfglyphtounicode{Qmonospace}{FF31}
\pdfglyphtounicode{Qsmall}{0071}
\pdfglyphtounicode{R}{0052}
\pdfglyphtounicode{Raarmenian}{054C}
\pdfglyphtounicode{Racute}{0154}
\pdfglyphtounicode{Rcaron}{0158}
\pdfglyphtounicode{Rcedilla}{0156}
\pdfglyphtounicode{Rcircle}{24C7}
\pdfglyphtounicode{Rcommaaccent}{0156}
\pdfglyphtounicode{Rdblgrave}{0210}
\pdfglyphtounicode{Rdotaccent}{1E58}
\pdfglyphtounicode{Rdotbelow}{1E5A}
\pdfglyphtounicode{Rdotbelowmacron}{1E5C}
\pdfglyphtounicode{Reharmenian}{0550}
\pdfglyphtounicode{Rfractur}{211C}
\pdfglyphtounicode{Rfraktur}{211C}
\pdfglyphtounicode{Rho}{03A1}
\pdfglyphtounicode{Ringsmall}{02DA}
\pdfglyphtounicode{Rinvertedbreve}{0212}
\pdfglyphtounicode{Rlinebelow}{1E5E}
\pdfglyphtounicode{Rmonospace}{FF32}
\pdfglyphtounicode{Rsmall}{0072}
\pdfglyphtounicode{Rsmallinverted}{0281}
\pdfglyphtounicode{Rsmallinvertedsuperior}{02B6}
\pdfglyphtounicode{S}{0053}
\pdfglyphtounicode{SF010000}{250C}
\pdfglyphtounicode{SF020000}{2514}
\pdfglyphtounicode{SF030000}{2510}
\pdfglyphtounicode{SF040000}{2518}
\pdfglyphtounicode{SF050000}{253C}
\pdfglyphtounicode{SF060000}{252C}
\pdfglyphtounicode{SF070000}{2534}
\pdfglyphtounicode{SF080000}{251C}
\pdfglyphtounicode{SF090000}{2524}
\pdfglyphtounicode{SF100000}{2500}
\pdfglyphtounicode{SF110000}{2502}
\pdfglyphtounicode{SF190000}{2561}
\pdfglyphtounicode{SF200000}{2562}
\pdfglyphtounicode{SF210000}{2556}
\pdfglyphtounicode{SF220000}{2555}
\pdfglyphtounicode{SF230000}{2563}
\pdfglyphtounicode{SF240000}{2551}
\pdfglyphtounicode{SF250000}{2557}
\pdfglyphtounicode{SF260000}{255D}
\pdfglyphtounicode{SF270000}{255C}
\pdfglyphtounicode{SF280000}{255B}
\pdfglyphtounicode{SF360000}{255E}
\pdfglyphtounicode{SF370000}{255F}
\pdfglyphtounicode{SF380000}{255A}
\pdfglyphtounicode{SF390000}{2554}
\pdfglyphtounicode{SF400000}{2569}
\pdfglyphtounicode{SF410000}{2566}
\pdfglyphtounicode{SF420000}{2560}
\pdfglyphtounicode{SF430000}{2550}
\pdfglyphtounicode{SF440000}{256C}
\pdfglyphtounicode{SF450000}{2567}
\pdfglyphtounicode{SF460000}{2568}
\pdfglyphtounicode{SF470000}{2564}
\pdfglyphtounicode{SF480000}{2565}
\pdfglyphtounicode{SF490000}{2559}
\pdfglyphtounicode{SF500000}{2558}
\pdfglyphtounicode{SF510000}{2552}
\pdfglyphtounicode{SF520000}{2553}
\pdfglyphtounicode{SF530000}{256B}
\pdfglyphtounicode{SF540000}{256A}
\pdfglyphtounicode{SS}{0053 0053}
\pdfglyphtounicode{SSsmall}{0073 0073}
\pdfglyphtounicode{Sacute}{015A}
\pdfglyphtounicode{Sacutedotaccent}{1E64}
\pdfglyphtounicode{Sampigreek}{03E0}
\pdfglyphtounicode{Scaron}{0160}
\pdfglyphtounicode{Scarondotaccent}{1E66}
\pdfglyphtounicode{Scaronsmall}{0161}
\pdfglyphtounicode{Scedilla}{015E}
\pdfglyphtounicode{Schwa}{018F}
\pdfglyphtounicode{Schwacyrillic}{04D8}
\pdfglyphtounicode{Schwadieresiscyrillic}{04DA}
\pdfglyphtounicode{Scircle}{24C8}
\pdfglyphtounicode{Scircumflex}{015C}
\pdfglyphtounicode{Scommaaccent}{0218}
\pdfglyphtounicode{Sdotaccent}{1E60}
\pdfglyphtounicode{Sdotbelow}{1E62}
\pdfglyphtounicode{Sdotbelowdotaccent}{1E68}
\pdfglyphtounicode{Seharmenian}{054D}
\pdfglyphtounicode{Sevenroman}{2166}
\pdfglyphtounicode{Shaarmenian}{0547}
\pdfglyphtounicode{Shacyrillic}{0428}
\pdfglyphtounicode{Shchacyrillic}{0429}
\pdfglyphtounicode{Sheicoptic}{03E2}
\pdfglyphtounicode{Shhacyrillic}{04BA}
\pdfglyphtounicode{Shimacoptic}{03EC}
\pdfglyphtounicode{Sigma}{03A3}
\pdfglyphtounicode{Sixroman}{2165}
\pdfglyphtounicode{Smonospace}{FF33}
\pdfglyphtounicode{Softsigncyrillic}{042C}
\pdfglyphtounicode{Ssmall}{0073}
\pdfglyphtounicode{Stigmagreek}{03DA}
\pdfglyphtounicode{T}{0054}
\pdfglyphtounicode{Tau}{03A4}
\pdfglyphtounicode{Tbar}{0166}
\pdfglyphtounicode{Tcaron}{0164}
\pdfglyphtounicode{Tcedilla}{0162}
\pdfglyphtounicode{Tcircle}{24C9}
\pdfglyphtounicode{Tcircumflexbelow}{1E70}
\pdfglyphtounicode{Tcommaaccent}{0162}
\pdfglyphtounicode{Tdotaccent}{1E6A}
\pdfglyphtounicode{Tdotbelow}{1E6C}
\pdfglyphtounicode{Tecyrillic}{0422}
\pdfglyphtounicode{Tedescendercyrillic}{04AC}
\pdfglyphtounicode{Tenroman}{2169}
\pdfglyphtounicode{Tetsecyrillic}{04B4}
\pdfglyphtounicode{Theta}{0398}
\pdfglyphtounicode{Thook}{01AC}
\pdfglyphtounicode{Thorn}{00DE}
\pdfglyphtounicode{Thornsmall}{00FE}
\pdfglyphtounicode{Threeroman}{2162}
\pdfglyphtounicode{Tildesmall}{02DC}
\pdfglyphtounicode{Tiwnarmenian}{054F}
\pdfglyphtounicode{Tlinebelow}{1E6E}
\pdfglyphtounicode{Tmonospace}{FF34}
\pdfglyphtounicode{Toarmenian}{0539}
\pdfglyphtounicode{Tonefive}{01BC}
\pdfglyphtounicode{Tonesix}{0184}
\pdfglyphtounicode{Tonetwo}{01A7}
\pdfglyphtounicode{Tretroflexhook}{01AE}
\pdfglyphtounicode{Tsecyrillic}{0426}
\pdfglyphtounicode{Tshecyrillic}{040B}
\pdfglyphtounicode{Tsmall}{0074}
\pdfglyphtounicode{Twelveroman}{216B}
\pdfglyphtounicode{Tworoman}{2161}
\pdfglyphtounicode{U}{0055}
\pdfglyphtounicode{Uacute}{00DA}
\pdfglyphtounicode{Uacutesmall}{00FA}
\pdfglyphtounicode{Ubreve}{016C}
\pdfglyphtounicode{Ucaron}{01D3}
\pdfglyphtounicode{Ucircle}{24CA}
\pdfglyphtounicode{Ucircumflex}{00DB}
\pdfglyphtounicode{Ucircumflexbelow}{1E76}
\pdfglyphtounicode{Ucircumflexsmall}{00FB}
\pdfglyphtounicode{Ucyrillic}{0423}
\pdfglyphtounicode{Udblacute}{0170}
\pdfglyphtounicode{Udblgrave}{0214}
\pdfglyphtounicode{Udieresis}{00DC}
\pdfglyphtounicode{Udieresisacute}{01D7}
\pdfglyphtounicode{Udieresisbelow}{1E72}
\pdfglyphtounicode{Udieresiscaron}{01D9}
\pdfglyphtounicode{Udieresiscyrillic}{04F0}
\pdfglyphtounicode{Udieresisgrave}{01DB}
\pdfglyphtounicode{Udieresismacron}{01D5}
\pdfglyphtounicode{Udieresissmall}{00FC}
\pdfglyphtounicode{Udotbelow}{1EE4}
\pdfglyphtounicode{Ugrave}{00D9}
\pdfglyphtounicode{Ugravesmall}{00F9}
\pdfglyphtounicode{Uhookabove}{1EE6}
\pdfglyphtounicode{Uhorn}{01AF}
\pdfglyphtounicode{Uhornacute}{1EE8}
\pdfglyphtounicode{Uhorndotbelow}{1EF0}
\pdfglyphtounicode{Uhorngrave}{1EEA}
\pdfglyphtounicode{Uhornhookabove}{1EEC}
\pdfglyphtounicode{Uhorntilde}{1EEE}
\pdfglyphtounicode{Uhungarumlaut}{0170}
\pdfglyphtounicode{Uhungarumlautcyrillic}{04F2}
\pdfglyphtounicode{Uinvertedbreve}{0216}
\pdfglyphtounicode{Ukcyrillic}{0478}
\pdfglyphtounicode{Umacron}{016A}
\pdfglyphtounicode{Umacroncyrillic}{04EE}
\pdfglyphtounicode{Umacrondieresis}{1E7A}
\pdfglyphtounicode{Umonospace}{FF35}
\pdfglyphtounicode{Uogonek}{0172}
\pdfglyphtounicode{Upsilon}{03A5}
\pdfglyphtounicode{Upsilon1}{03D2}
\pdfglyphtounicode{Upsilonacutehooksymbolgreek}{03D3}
\pdfglyphtounicode{Upsilonafrican}{01B1}
\pdfglyphtounicode{Upsilondieresis}{03AB}
\pdfglyphtounicode{Upsilondieresishooksymbolgreek}{03D4}
\pdfglyphtounicode{Upsilonhooksymbol}{03D2}
\pdfglyphtounicode{Upsilontonos}{038E}
\pdfglyphtounicode{Uring}{016E}
\pdfglyphtounicode{Ushortcyrillic}{040E}
\pdfglyphtounicode{Usmall}{0075}
\pdfglyphtounicode{Ustraightcyrillic}{04AE}
\pdfglyphtounicode{Ustraightstrokecyrillic}{04B0}
\pdfglyphtounicode{Utilde}{0168}
\pdfglyphtounicode{Utildeacute}{1E78}
\pdfglyphtounicode{Utildebelow}{1E74}
\pdfglyphtounicode{V}{0056}
\pdfglyphtounicode{Vcircle}{24CB}
\pdfglyphtounicode{Vdotbelow}{1E7E}
\pdfglyphtounicode{Vecyrillic}{0412}
\pdfglyphtounicode{Vewarmenian}{054E}
\pdfglyphtounicode{Vhook}{01B2}
\pdfglyphtounicode{Vmonospace}{FF36}
\pdfglyphtounicode{Voarmenian}{0548}
\pdfglyphtounicode{Vsmall}{0076}
\pdfglyphtounicode{Vtilde}{1E7C}
\pdfglyphtounicode{W}{0057}
\pdfglyphtounicode{Wacute}{1E82}
\pdfglyphtounicode{Wcircle}{24CC}
\pdfglyphtounicode{Wcircumflex}{0174}
\pdfglyphtounicode{Wdieresis}{1E84}
\pdfglyphtounicode{Wdotaccent}{1E86}
\pdfglyphtounicode{Wdotbelow}{1E88}
\pdfglyphtounicode{Wgrave}{1E80}
\pdfglyphtounicode{Wmonospace}{FF37}
\pdfglyphtounicode{Wsmall}{0077}
\pdfglyphtounicode{X}{0058}
\pdfglyphtounicode{Xcircle}{24CD}
\pdfglyphtounicode{Xdieresis}{1E8C}
\pdfglyphtounicode{Xdotaccent}{1E8A}
\pdfglyphtounicode{Xeharmenian}{053D}
\pdfglyphtounicode{Xi}{039E}
\pdfglyphtounicode{Xmonospace}{FF38}
\pdfglyphtounicode{Xsmall}{0078}
\pdfglyphtounicode{Y}{0059}
\pdfglyphtounicode{Yacute}{00DD}
\pdfglyphtounicode{Yacutesmall}{00FD}
\pdfglyphtounicode{Yatcyrillic}{0462}
\pdfglyphtounicode{Ycircle}{24CE}
\pdfglyphtounicode{Ycircumflex}{0176}
\pdfglyphtounicode{Ydieresis}{0178}
\pdfglyphtounicode{Ydieresissmall}{00FF}
\pdfglyphtounicode{Ydotaccent}{1E8E}
\pdfglyphtounicode{Ydotbelow}{1EF4}
\pdfglyphtounicode{Yen}{00A5}
\pdfglyphtounicode{Yericyrillic}{042B}
\pdfglyphtounicode{Yerudieresiscyrillic}{04F8}
\pdfglyphtounicode{Ygrave}{1EF2}
\pdfglyphtounicode{Yhook}{01B3}
\pdfglyphtounicode{Yhookabove}{1EF6}
\pdfglyphtounicode{Yiarmenian}{0545}
\pdfglyphtounicode{Yicyrillic}{0407}
\pdfglyphtounicode{Yiwnarmenian}{0552}
\pdfglyphtounicode{Ymonospace}{FF39}
\pdfglyphtounicode{Ysmall}{0079}
\pdfglyphtounicode{Ytilde}{1EF8}
\pdfglyphtounicode{Yusbigcyrillic}{046A}
\pdfglyphtounicode{Yusbigiotifiedcyrillic}{046C}
\pdfglyphtounicode{Yuslittlecyrillic}{0466}
\pdfglyphtounicode{Yuslittleiotifiedcyrillic}{0468}
\pdfglyphtounicode{Z}{005A}
\pdfglyphtounicode{Zaarmenian}{0536}
\pdfglyphtounicode{Zacute}{0179}
\pdfglyphtounicode{Zcaron}{017D}
\pdfglyphtounicode{Zcaronsmall}{017E}
\pdfglyphtounicode{Zcircle}{24CF}
\pdfglyphtounicode{Zcircumflex}{1E90}
\pdfglyphtounicode{Zdot}{017B}
\pdfglyphtounicode{Zdotaccent}{017B}
\pdfglyphtounicode{Zdotbelow}{1E92}
\pdfglyphtounicode{Zecyrillic}{0417}
\pdfglyphtounicode{Zedescendercyrillic}{0498}
\pdfglyphtounicode{Zedieresiscyrillic}{04DE}
\pdfglyphtounicode{Zeta}{0396}
\pdfglyphtounicode{Zhearmenian}{053A}
\pdfglyphtounicode{Zhebrevecyrillic}{04C1}
\pdfglyphtounicode{Zhecyrillic}{0416}
\pdfglyphtounicode{Zhedescendercyrillic}{0496}
\pdfglyphtounicode{Zhedieresiscyrillic}{04DC}
\pdfglyphtounicode{Zlinebelow}{1E94}
\pdfglyphtounicode{Zmonospace}{FF3A}
\pdfglyphtounicode{Zsmall}{007A}
\pdfglyphtounicode{Zstroke}{01B5}
\pdfglyphtounicode{a}{0061}
\pdfglyphtounicode{aabengali}{0986}
\pdfglyphtounicode{aacute}{00E1}
\pdfglyphtounicode{aadeva}{0906}
\pdfglyphtounicode{aagujarati}{0A86}
\pdfglyphtounicode{aagurmukhi}{0A06}
\pdfglyphtounicode{aamatragurmukhi}{0A3E}
\pdfglyphtounicode{aarusquare}{3303}
\pdfglyphtounicode{aavowelsignbengali}{09BE}
\pdfglyphtounicode{aavowelsigndeva}{093E}
\pdfglyphtounicode{aavowelsigngujarati}{0ABE}
\pdfglyphtounicode{abbreviationmarkarmenian}{055F}
\pdfglyphtounicode{abbreviationsigndeva}{0970}
\pdfglyphtounicode{abengali}{0985}
\pdfglyphtounicode{abopomofo}{311A}
\pdfglyphtounicode{abreve}{0103}
\pdfglyphtounicode{abreveacute}{1EAF}
\pdfglyphtounicode{abrevecyrillic}{04D1}
\pdfglyphtounicode{abrevedotbelow}{1EB7}
\pdfglyphtounicode{abrevegrave}{1EB1}
\pdfglyphtounicode{abrevehookabove}{1EB3}
\pdfglyphtounicode{abrevetilde}{1EB5}
\pdfglyphtounicode{acaron}{01CE}
\pdfglyphtounicode{acircle}{24D0}
\pdfglyphtounicode{acircumflex}{00E2}
\pdfglyphtounicode{acircumflexacute}{1EA5}
\pdfglyphtounicode{acircumflexdotbelow}{1EAD}
\pdfglyphtounicode{acircumflexgrave}{1EA7}
\pdfglyphtounicode{acircumflexhookabove}{1EA9}
\pdfglyphtounicode{acircumflextilde}{1EAB}
\pdfglyphtounicode{acute}{00B4}
\pdfglyphtounicode{acutebelowcmb}{0317}
\pdfglyphtounicode{acutecmb}{0301}
\pdfglyphtounicode{acutecomb}{0301}
\pdfglyphtounicode{acutedeva}{0954}
\pdfglyphtounicode{acutelowmod}{02CF}
\pdfglyphtounicode{acutetonecmb}{0341}
\pdfglyphtounicode{acyrillic}{0430}
\pdfglyphtounicode{adblgrave}{0201}
\pdfglyphtounicode{addakgurmukhi}{0A71}
\pdfglyphtounicode{adeva}{0905}
\pdfglyphtounicode{adieresis}{00E4}
\pdfglyphtounicode{adieresiscyrillic}{04D3}
\pdfglyphtounicode{adieresismacron}{01DF}
\pdfglyphtounicode{adotbelow}{1EA1}
\pdfglyphtounicode{adotmacron}{01E1}
\pdfglyphtounicode{ae}{00E6}
\pdfglyphtounicode{aeacute}{01FD}
\pdfglyphtounicode{aekorean}{3150}
\pdfglyphtounicode{aemacron}{01E3}
\pdfglyphtounicode{afii00208}{2015}
\pdfglyphtounicode{afii08941}{20A4}
\pdfglyphtounicode{afii10017}{0410}
\pdfglyphtounicode{afii10018}{0411}
\pdfglyphtounicode{afii10019}{0412}
\pdfglyphtounicode{afii10020}{0413}
\pdfglyphtounicode{afii10021}{0414}
\pdfglyphtounicode{afii10022}{0415}
\pdfglyphtounicode{afii10023}{0401}
\pdfglyphtounicode{afii10024}{0416}
\pdfglyphtounicode{afii10025}{0417}
\pdfglyphtounicode{afii10026}{0418}
\pdfglyphtounicode{afii10027}{0419}
\pdfglyphtounicode{afii10028}{041A}
\pdfglyphtounicode{afii10029}{041B}
\pdfglyphtounicode{afii10030}{041C}
\pdfglyphtounicode{afii10031}{041D}
\pdfglyphtounicode{afii10032}{041E}
\pdfglyphtounicode{afii10033}{041F}
\pdfglyphtounicode{afii10034}{0420}
\pdfglyphtounicode{afii10035}{0421}
\pdfglyphtounicode{afii10036}{0422}
\pdfglyphtounicode{afii10037}{0423}
\pdfglyphtounicode{afii10038}{0424}
\pdfglyphtounicode{afii10039}{0425}
\pdfglyphtounicode{afii10040}{0426}
\pdfglyphtounicode{afii10041}{0427}
\pdfglyphtounicode{afii10042}{0428}
\pdfglyphtounicode{afii10043}{0429}
\pdfglyphtounicode{afii10044}{042A}
\pdfglyphtounicode{afii10045}{042B}
\pdfglyphtounicode{afii10046}{042C}
\pdfglyphtounicode{afii10047}{042D}
\pdfglyphtounicode{afii10048}{042E}
\pdfglyphtounicode{afii10049}{042F}
\pdfglyphtounicode{afii10050}{0490}
\pdfglyphtounicode{afii10051}{0402}
\pdfglyphtounicode{afii10052}{0403}
\pdfglyphtounicode{afii10053}{0404}
\pdfglyphtounicode{afii10054}{0405}
\pdfglyphtounicode{afii10055}{0406}
\pdfglyphtounicode{afii10056}{0407}
\pdfglyphtounicode{afii10057}{0408}
\pdfglyphtounicode{afii10058}{0409}
\pdfglyphtounicode{afii10059}{040A}
\pdfglyphtounicode{afii10060}{040B}
\pdfglyphtounicode{afii10061}{040C}
\pdfglyphtounicode{afii10062}{040E}
\pdfglyphtounicode{afii10063}{F6C4}
\pdfglyphtounicode{afii10064}{F6C5}
\pdfglyphtounicode{afii10065}{0430}
\pdfglyphtounicode{afii10066}{0431}
\pdfglyphtounicode{afii10067}{0432}
\pdfglyphtounicode{afii10068}{0433}
\pdfglyphtounicode{afii10069}{0434}
\pdfglyphtounicode{afii10070}{0435}
\pdfglyphtounicode{afii10071}{0451}
\pdfglyphtounicode{afii10072}{0436}
\pdfglyphtounicode{afii10073}{0437}
\pdfglyphtounicode{afii10074}{0438}
\pdfglyphtounicode{afii10075}{0439}
\pdfglyphtounicode{afii10076}{043A}
\pdfglyphtounicode{afii10077}{043B}
\pdfglyphtounicode{afii10078}{043C}
\pdfglyphtounicode{afii10079}{043D}
\pdfglyphtounicode{afii10080}{043E}
\pdfglyphtounicode{afii10081}{043F}
\pdfglyphtounicode{afii10082}{0440}
\pdfglyphtounicode{afii10083}{0441}
\pdfglyphtounicode{afii10084}{0442}
\pdfglyphtounicode{afii10085}{0443}
\pdfglyphtounicode{afii10086}{0444}
\pdfglyphtounicode{afii10087}{0445}
\pdfglyphtounicode{afii10088}{0446}
\pdfglyphtounicode{afii10089}{0447}
\pdfglyphtounicode{afii10090}{0448}
\pdfglyphtounicode{afii10091}{0449}
\pdfglyphtounicode{afii10092}{044A}
\pdfglyphtounicode{afii10093}{044B}
\pdfglyphtounicode{afii10094}{044C}
\pdfglyphtounicode{afii10095}{044D}
\pdfglyphtounicode{afii10096}{044E}
\pdfglyphtounicode{afii10097}{044F}
\pdfglyphtounicode{afii10098}{0491}
\pdfglyphtounicode{afii10099}{0452}
\pdfglyphtounicode{afii10100}{0453}
\pdfglyphtounicode{afii10101}{0454}
\pdfglyphtounicode{afii10102}{0455}
\pdfglyphtounicode{afii10103}{0456}
\pdfglyphtounicode{afii10104}{0457}
\pdfglyphtounicode{afii10105}{0458}
\pdfglyphtounicode{afii10106}{0459}
\pdfglyphtounicode{afii10107}{045A}
\pdfglyphtounicode{afii10108}{045B}
\pdfglyphtounicode{afii10109}{045C}
\pdfglyphtounicode{afii10110}{045E}
\pdfglyphtounicode{afii10145}{040F}
\pdfglyphtounicode{afii10146}{0462}
\pdfglyphtounicode{afii10147}{0472}
\pdfglyphtounicode{afii10148}{0474}
\pdfglyphtounicode{afii10192}{F6C6}
\pdfglyphtounicode{afii10193}{045F}
\pdfglyphtounicode{afii10194}{0463}
\pdfglyphtounicode{afii10195}{0473}
\pdfglyphtounicode{afii10196}{0475}
\pdfglyphtounicode{afii10831}{F6C7}
\pdfglyphtounicode{afii10832}{F6C8}
\pdfglyphtounicode{afii10846}{04D9}
\pdfglyphtounicode{afii299}{200E}
\pdfglyphtounicode{afii300}{200F}
\pdfglyphtounicode{afii301}{200D}
\pdfglyphtounicode{afii57381}{066A}
\pdfglyphtounicode{afii57388}{060C}
\pdfglyphtounicode{afii57392}{0660}
\pdfglyphtounicode{afii57393}{0661}
\pdfglyphtounicode{afii57394}{0662}
\pdfglyphtounicode{afii57395}{0663}
\pdfglyphtounicode{afii57396}{0664}
\pdfglyphtounicode{afii57397}{0665}
\pdfglyphtounicode{afii57398}{0666}
\pdfglyphtounicode{afii57399}{0667}
\pdfglyphtounicode{afii57400}{0668}
\pdfglyphtounicode{afii57401}{0669}
\pdfglyphtounicode{afii57403}{061B}
\pdfglyphtounicode{afii57407}{061F}
\pdfglyphtounicode{afii57409}{0621}
\pdfglyphtounicode{afii57410}{0622}
\pdfglyphtounicode{afii57411}{0623}
\pdfglyphtounicode{afii57412}{0624}
\pdfglyphtounicode{afii57413}{0625}
\pdfglyphtounicode{afii57414}{0626}
\pdfglyphtounicode{afii57415}{0627}
\pdfglyphtounicode{afii57416}{0628}
\pdfglyphtounicode{afii57417}{0629}
\pdfglyphtounicode{afii57418}{062A}
\pdfglyphtounicode{afii57419}{062B}
\pdfglyphtounicode{afii57420}{062C}
\pdfglyphtounicode{afii57421}{062D}
\pdfglyphtounicode{afii57422}{062E}
\pdfglyphtounicode{afii57423}{062F}
\pdfglyphtounicode{afii57424}{0630}
\pdfglyphtounicode{afii57425}{0631}
\pdfglyphtounicode{afii57426}{0632}
\pdfglyphtounicode{afii57427}{0633}
\pdfglyphtounicode{afii57428}{0634}
\pdfglyphtounicode{afii57429}{0635}
\pdfglyphtounicode{afii57430}{0636}
\pdfglyphtounicode{afii57431}{0637}
\pdfglyphtounicode{afii57432}{0638}
\pdfglyphtounicode{afii57433}{0639}
\pdfglyphtounicode{afii57434}{063A}
\pdfglyphtounicode{afii57440}{0640}
\pdfglyphtounicode{afii57441}{0641}
\pdfglyphtounicode{afii57442}{0642}
\pdfglyphtounicode{afii57443}{0643}
\pdfglyphtounicode{afii57444}{0644}
\pdfglyphtounicode{afii57445}{0645}
\pdfglyphtounicode{afii57446}{0646}
\pdfglyphtounicode{afii57448}{0648}
\pdfglyphtounicode{afii57449}{0649}
\pdfglyphtounicode{afii57450}{064A}
\pdfglyphtounicode{afii57451}{064B}
\pdfglyphtounicode{afii57452}{064C}
\pdfglyphtounicode{afii57453}{064D}
\pdfglyphtounicode{afii57454}{064E}
\pdfglyphtounicode{afii57455}{064F}
\pdfglyphtounicode{afii57456}{0650}
\pdfglyphtounicode{afii57457}{0651}
\pdfglyphtounicode{afii57458}{0652}
\pdfglyphtounicode{afii57470}{0647}
\pdfglyphtounicode{afii57505}{06A4}
\pdfglyphtounicode{afii57506}{067E}
\pdfglyphtounicode{afii57507}{0686}
\pdfglyphtounicode{afii57508}{0698}
\pdfglyphtounicode{afii57509}{06AF}
\pdfglyphtounicode{afii57511}{0679}
\pdfglyphtounicode{afii57512}{0688}
\pdfglyphtounicode{afii57513}{0691}
\pdfglyphtounicode{afii57514}{06BA}
\pdfglyphtounicode{afii57519}{06D2}
\pdfglyphtounicode{afii57534}{06D5}
\pdfglyphtounicode{afii57636}{20AA}
\pdfglyphtounicode{afii57645}{05BE}
\pdfglyphtounicode{afii57658}{05C3}
\pdfglyphtounicode{afii57664}{05D0}
\pdfglyphtounicode{afii57665}{05D1}
\pdfglyphtounicode{afii57666}{05D2}
\pdfglyphtounicode{afii57667}{05D3}
\pdfglyphtounicode{afii57668}{05D4}
\pdfglyphtounicode{afii57669}{05D5}
\pdfglyphtounicode{afii57670}{05D6}
\pdfglyphtounicode{afii57671}{05D7}
\pdfglyphtounicode{afii57672}{05D8}
\pdfglyphtounicode{afii57673}{05D9}
\pdfglyphtounicode{afii57674}{05DA}
\pdfglyphtounicode{afii57675}{05DB}
\pdfglyphtounicode{afii57676}{05DC}
\pdfglyphtounicode{afii57677}{05DD}
\pdfglyphtounicode{afii57678}{05DE}
\pdfglyphtounicode{afii57679}{05DF}
\pdfglyphtounicode{afii57680}{05E0}
\pdfglyphtounicode{afii57681}{05E1}
\pdfglyphtounicode{afii57682}{05E2}
\pdfglyphtounicode{afii57683}{05E3}
\pdfglyphtounicode{afii57684}{05E4}
\pdfglyphtounicode{afii57685}{05E5}
\pdfglyphtounicode{afii57686}{05E6}
\pdfglyphtounicode{afii57687}{05E7}
\pdfglyphtounicode{afii57688}{05E8}
\pdfglyphtounicode{afii57689}{05E9}
\pdfglyphtounicode{afii57690}{05EA}
\pdfglyphtounicode{afii57694}{FB2A}
\pdfglyphtounicode{afii57695}{FB2B}
\pdfglyphtounicode{afii57700}{FB4B}
\pdfglyphtounicode{afii57705}{FB1F}
\pdfglyphtounicode{afii57716}{05F0}
\pdfglyphtounicode{afii57717}{05F1}
\pdfglyphtounicode{afii57718}{05F2}
\pdfglyphtounicode{afii57723}{FB35}
\pdfglyphtounicode{afii57793}{05B4}
\pdfglyphtounicode{afii57794}{05B5}
\pdfglyphtounicode{afii57795}{05B6}
\pdfglyphtounicode{afii57796}{05BB}
\pdfglyphtounicode{afii57797}{05B8}
\pdfglyphtounicode{afii57798}{05B7}
\pdfglyphtounicode{afii57799}{05B0}
\pdfglyphtounicode{afii57800}{05B2}
\pdfglyphtounicode{afii57801}{05B1}
\pdfglyphtounicode{afii57802}{05B3}
\pdfglyphtounicode{afii57803}{05C2}
\pdfglyphtounicode{afii57804}{05C1}
\pdfglyphtounicode{afii57806}{05B9}
\pdfglyphtounicode{afii57807}{05BC}
\pdfglyphtounicode{afii57839}{05BD}
\pdfglyphtounicode{afii57841}{05BF}
\pdfglyphtounicode{afii57842}{05C0}
\pdfglyphtounicode{afii57929}{02BC}
\pdfglyphtounicode{afii61248}{2105}
\pdfglyphtounicode{afii61289}{2113}
\pdfglyphtounicode{afii61352}{2116}
\pdfglyphtounicode{afii61573}{202C}
\pdfglyphtounicode{afii61574}{202D}
\pdfglyphtounicode{afii61575}{202E}
\pdfglyphtounicode{afii61664}{200C}
\pdfglyphtounicode{afii63167}{066D}
\pdfglyphtounicode{afii64937}{02BD}
\pdfglyphtounicode{agrave}{00E0}
\pdfglyphtounicode{agujarati}{0A85}
\pdfglyphtounicode{agurmukhi}{0A05}
\pdfglyphtounicode{ahiragana}{3042}
\pdfglyphtounicode{ahookabove}{1EA3}
\pdfglyphtounicode{aibengali}{0990}
\pdfglyphtounicode{aibopomofo}{311E}
\pdfglyphtounicode{aideva}{0910}
\pdfglyphtounicode{aiecyrillic}{04D5}
\pdfglyphtounicode{aigujarati}{0A90}
\pdfglyphtounicode{aigurmukhi}{0A10}
\pdfglyphtounicode{aimatragurmukhi}{0A48}
\pdfglyphtounicode{ainarabic}{0639}
\pdfglyphtounicode{ainfinalarabic}{FECA}
\pdfglyphtounicode{aininitialarabic}{FECB}
\pdfglyphtounicode{ainmedialarabic}{FECC}
\pdfglyphtounicode{ainvertedbreve}{0203}
\pdfglyphtounicode{aivowelsignbengali}{09C8}
\pdfglyphtounicode{aivowelsigndeva}{0948}
\pdfglyphtounicode{aivowelsigngujarati}{0AC8}
\pdfglyphtounicode{akatakana}{30A2}
\pdfglyphtounicode{akatakanahalfwidth}{FF71}
\pdfglyphtounicode{akorean}{314F}
\pdfglyphtounicode{alef}{05D0}
\pdfglyphtounicode{alefarabic}{0627}
\pdfglyphtounicode{alefdageshhebrew}{FB30}
\pdfglyphtounicode{aleffinalarabic}{FE8E}
\pdfglyphtounicode{alefhamzaabovearabic}{0623}
\pdfglyphtounicode{alefhamzaabovefinalarabic}{FE84}
\pdfglyphtounicode{alefhamzabelowarabic}{0625}
\pdfglyphtounicode{alefhamzabelowfinalarabic}{FE88}
\pdfglyphtounicode{alefhebrew}{05D0}
\pdfglyphtounicode{aleflamedhebrew}{FB4F}
\pdfglyphtounicode{alefmaddaabovearabic}{0622}
\pdfglyphtounicode{alefmaddaabovefinalarabic}{FE82}
\pdfglyphtounicode{alefmaksuraarabic}{0649}
\pdfglyphtounicode{alefmaksurafinalarabic}{FEF0}
\pdfglyphtounicode{alefmaksurainitialarabic}{FEF3}
\pdfglyphtounicode{alefmaksuramedialarabic}{FEF4}
\pdfglyphtounicode{alefpatahhebrew}{FB2E}
\pdfglyphtounicode{alefqamatshebrew}{FB2F}
\pdfglyphtounicode{aleph}{2135}
\pdfglyphtounicode{allequal}{224C}
\pdfglyphtounicode{alpha}{03B1}
\pdfglyphtounicode{alphatonos}{03AC}
\pdfglyphtounicode{amacron}{0101}
\pdfglyphtounicode{amonospace}{FF41}
\pdfglyphtounicode{ampersand}{0026}
\pdfglyphtounicode{ampersandmonospace}{FF06}
\pdfglyphtounicode{ampersandsmall}{0026}
\pdfglyphtounicode{amsquare}{33C2}
\pdfglyphtounicode{anbopomofo}{3122}
\pdfglyphtounicode{angbopomofo}{3124}
\pdfglyphtounicode{angbracketleft}{27E8}
\pdfglyphtounicode{angbracketright}{27E9}
\pdfglyphtounicode{angkhankhuthai}{0E5A}
\pdfglyphtounicode{angle}{2220}
\pdfglyphtounicode{anglebracketleft}{3008}
\pdfglyphtounicode{anglebracketleftvertical}{FE3F}
\pdfglyphtounicode{anglebracketright}{3009}
\pdfglyphtounicode{anglebracketrightvertical}{FE40}
\pdfglyphtounicode{angleleft}{2329}
\pdfglyphtounicode{angleright}{232A}
\pdfglyphtounicode{angstrom}{212B}
\pdfglyphtounicode{anoteleia}{0387}
\pdfglyphtounicode{anticlockwise}{27F2}
\pdfglyphtounicode{anudattadeva}{0952}
\pdfglyphtounicode{anusvarabengali}{0982}
\pdfglyphtounicode{anusvaradeva}{0902}
\pdfglyphtounicode{anusvaragujarati}{0A82}
\pdfglyphtounicode{aogonek}{0105}
\pdfglyphtounicode{apaatosquare}{3300}
\pdfglyphtounicode{aparen}{249C}
\pdfglyphtounicode{apostrophearmenian}{055A}
\pdfglyphtounicode{apostrophemod}{02BC}
\pdfglyphtounicode{apple}{F8FF}
\pdfglyphtounicode{approaches}{2250}
\pdfglyphtounicode{approxequal}{2248}
\pdfglyphtounicode{approxequalorimage}{2252}
\pdfglyphtounicode{approximatelyequal}{2245}
\pdfglyphtounicode{approxorequal}{224A}
\pdfglyphtounicode{araeaekorean}{318E}
\pdfglyphtounicode{araeakorean}{318D}
\pdfglyphtounicode{arc}{2312}
\pdfglyphtounicode{archleftdown}{21B6}
\pdfglyphtounicode{archrightdown}{21B7}
\pdfglyphtounicode{arighthalfring}{1E9A}
\pdfglyphtounicode{aring}{00E5}
\pdfglyphtounicode{aringacute}{01FB}
\pdfglyphtounicode{aringbelow}{1E01}
\pdfglyphtounicode{arrowboth}{2194}
\pdfglyphtounicode{arrowbothv}{2195}
\pdfglyphtounicode{arrowdashdown}{21E3}
\pdfglyphtounicode{arrowdashleft}{21E0}
\pdfglyphtounicode{arrowdashright}{21E2}
\pdfglyphtounicode{arrowdashup}{21E1}
\pdfglyphtounicode{arrowdblboth}{21D4}
\pdfglyphtounicode{arrowdblbothv}{21D5}
\pdfglyphtounicode{arrowdbldown}{21D3}
\pdfglyphtounicode{arrowdblleft}{21D0}
\pdfglyphtounicode{arrowdblright}{21D2}
\pdfglyphtounicode{arrowdblup}{21D1}
\pdfglyphtounicode{arrowdown}{2193}
\pdfglyphtounicode{arrowdownleft}{2199}
\pdfglyphtounicode{arrowdownright}{2198}
\pdfglyphtounicode{arrowdownwhite}{21E9}
\pdfglyphtounicode{arrowheaddownmod}{02C5}
\pdfglyphtounicode{arrowheadleftmod}{02C2}
\pdfglyphtounicode{arrowheadrightmod}{02C3}
\pdfglyphtounicode{arrowheadupmod}{02C4}
\pdfglyphtounicode{arrowhorizex}{F8E7}
\pdfglyphtounicode{arrowleft}{2190}
\pdfglyphtounicode{arrowleftbothalf}{21BD}
\pdfglyphtounicode{arrowleftdbl}{21D0}
\pdfglyphtounicode{arrowleftdblstroke}{21CD}
\pdfglyphtounicode{arrowleftoverright}{21C6}
\pdfglyphtounicode{arrowlefttophalf}{21BC}
\pdfglyphtounicode{arrowleftwhite}{21E6}
\pdfglyphtounicode{arrownortheast}{2197}
\pdfglyphtounicode{arrownorthwest}{2196}
\pdfglyphtounicode{arrowparrleftright}{21C6}
\pdfglyphtounicode{arrowparrrightleft}{21C4}
\pdfglyphtounicode{arrowright}{2192}
\pdfglyphtounicode{arrowrightbothalf}{21C1}
\pdfglyphtounicode{arrowrightdblstroke}{21CF}
\pdfglyphtounicode{arrowrightheavy}{279E}
\pdfglyphtounicode{arrowrightoverleft}{21C4}
\pdfglyphtounicode{arrowrighttophalf}{21C0}
\pdfglyphtounicode{arrowrightwhite}{21E8}
\pdfglyphtounicode{arrowsoutheast}{2198}
\pdfglyphtounicode{arrowsouthwest}{2199}
\pdfglyphtounicode{arrowtableft}{21E4}
\pdfglyphtounicode{arrowtabright}{21E5}
\pdfglyphtounicode{arrowtailleft}{21A2}
\pdfglyphtounicode{arrowtailright}{21A3}
\pdfglyphtounicode{arrowtripleleft}{21DA}
\pdfglyphtounicode{arrowtripleright}{21DB}
\pdfglyphtounicode{arrowup}{2191}
\pdfglyphtounicode{arrowupdn}{2195}
\pdfglyphtounicode{arrowupdnbse}{21A8}
\pdfglyphtounicode{arrowupdownbase}{21A8}
\pdfglyphtounicode{arrowupleft}{2196}
\pdfglyphtounicode{arrowupleftofdown}{21C5}
\pdfglyphtounicode{arrowupright}{2197}
\pdfglyphtounicode{arrowupwhite}{21E7}
\pdfglyphtounicode{arrowvertex}{F8E6}
\pdfglyphtounicode{asciicircum}{005E}
\pdfglyphtounicode{asciicircummonospace}{FF3E}
\pdfglyphtounicode{asciitilde}{007E}
\pdfglyphtounicode{asciitildemonospace}{FF5E}
\pdfglyphtounicode{ascript}{0251}
\pdfglyphtounicode{ascriptturned}{0252}
\pdfglyphtounicode{asmallhiragana}{3041}
\pdfglyphtounicode{asmallkatakana}{30A1}
\pdfglyphtounicode{asmallkatakanahalfwidth}{FF67}
\pdfglyphtounicode{asterisk}{002A}
\pdfglyphtounicode{asteriskaltonearabic}{066D}
\pdfglyphtounicode{asteriskarabic}{066D}
\pdfglyphtounicode{asteriskcentered}{2217}
\pdfglyphtounicode{asteriskmath}{2217}
\pdfglyphtounicode{asteriskmonospace}{FF0A}
\pdfglyphtounicode{asterisksmall}{FE61}
\pdfglyphtounicode{asterism}{2042}
\pdfglyphtounicode{asuperior}{0061}
\pdfglyphtounicode{asymptoticallyequal}{2243}
\pdfglyphtounicode{at}{0040}
\pdfglyphtounicode{atilde}{00E3}
\pdfglyphtounicode{atmonospace}{FF20}
\pdfglyphtounicode{atsmall}{FE6B}
\pdfglyphtounicode{aturned}{0250}
\pdfglyphtounicode{aubengali}{0994}
\pdfglyphtounicode{aubopomofo}{3120}
\pdfglyphtounicode{audeva}{0914}
\pdfglyphtounicode{augujarati}{0A94}
\pdfglyphtounicode{augurmukhi}{0A14}
\pdfglyphtounicode{aulengthmarkbengali}{09D7}
\pdfglyphtounicode{aumatragurmukhi}{0A4C}
\pdfglyphtounicode{auvowelsignbengali}{09CC}
\pdfglyphtounicode{auvowelsigndeva}{094C}
\pdfglyphtounicode{auvowelsigngujarati}{0ACC}
\pdfglyphtounicode{avagrahadeva}{093D}
\pdfglyphtounicode{aybarmenian}{0561}
\pdfglyphtounicode{ayin}{05E2}
\pdfglyphtounicode{ayinaltonehebrew}{FB20}
\pdfglyphtounicode{ayinhebrew}{05E2}
\pdfglyphtounicode{b}{0062}
\pdfglyphtounicode{babengali}{09AC}
\pdfglyphtounicode{backslash}{005C}
\pdfglyphtounicode{backslashmonospace}{FF3C}
\pdfglyphtounicode{badeva}{092C}
\pdfglyphtounicode{bagujarati}{0AAC}
\pdfglyphtounicode{bagurmukhi}{0A2C}
\pdfglyphtounicode{bahiragana}{3070}
\pdfglyphtounicode{bahtthai}{0E3F}
\pdfglyphtounicode{bakatakana}{30D0}
\pdfglyphtounicode{bar}{007C}
\pdfglyphtounicode{bardbl}{2225}
\pdfglyphtounicode{barmonospace}{FF5C}
\pdfglyphtounicode{bbopomofo}{3105}
\pdfglyphtounicode{bcircle}{24D1}
\pdfglyphtounicode{bdotaccent}{1E03}
\pdfglyphtounicode{bdotbelow}{1E05}
\pdfglyphtounicode{beamedsixteenthnotes}{266C}
\pdfglyphtounicode{because}{2235}
\pdfglyphtounicode{becyrillic}{0431}
\pdfglyphtounicode{beharabic}{0628}
\pdfglyphtounicode{behfinalarabic}{FE90}
\pdfglyphtounicode{behinitialarabic}{FE91}
\pdfglyphtounicode{behiragana}{3079}
\pdfglyphtounicode{behmedialarabic}{FE92}
\pdfglyphtounicode{behmeeminitialarabic}{FC9F}
\pdfglyphtounicode{behmeemisolatedarabic}{FC08}
\pdfglyphtounicode{behnoonfinalarabic}{FC6D}
\pdfglyphtounicode{bekatakana}{30D9}
\pdfglyphtounicode{benarmenian}{0562}
\pdfglyphtounicode{bet}{05D1}
\pdfglyphtounicode{beta}{03B2}
\pdfglyphtounicode{betasymbolgreek}{03D0}
\pdfglyphtounicode{betdagesh}{FB31}
\pdfglyphtounicode{betdageshhebrew}{FB31}
\pdfglyphtounicode{beth}{2136}
\pdfglyphtounicode{bethebrew}{05D1}
\pdfglyphtounicode{betrafehebrew}{FB4C}
\pdfglyphtounicode{between}{226C}
\pdfglyphtounicode{bhabengali}{09AD}
\pdfglyphtounicode{bhadeva}{092D}
\pdfglyphtounicode{bhagujarati}{0AAD}
\pdfglyphtounicode{bhagurmukhi}{0A2D}
\pdfglyphtounicode{bhook}{0253}
\pdfglyphtounicode{bihiragana}{3073}
\pdfglyphtounicode{bikatakana}{30D3}
\pdfglyphtounicode{bilabialclick}{0298}
\pdfglyphtounicode{bindigurmukhi}{0A02}
\pdfglyphtounicode{birusquare}{3331}
\pdfglyphtounicode{blackcircle}{25CF}
\pdfglyphtounicode{blackdiamond}{25C6}
\pdfglyphtounicode{blackdownpointingtriangle}{25BC}
\pdfglyphtounicode{blackleftpointingpointer}{25C4}
\pdfglyphtounicode{blackleftpointingtriangle}{25C0}
\pdfglyphtounicode{blacklenticularbracketleft}{3010}
\pdfglyphtounicode{blacklenticularbracketleftvertical}{FE3B}
\pdfglyphtounicode{blacklenticularbracketright}{3011}
\pdfglyphtounicode{blacklenticularbracketrightvertical}{FE3C}
\pdfglyphtounicode{blacklowerlefttriangle}{25E3}
\pdfglyphtounicode{blacklowerrighttriangle}{25E2}
\pdfglyphtounicode{blackrectangle}{25AC}
\pdfglyphtounicode{blackrightpointingpointer}{25BA}
\pdfglyphtounicode{blackrightpointingtriangle}{25B6}
\pdfglyphtounicode{blacksmallsquare}{25AA}
\pdfglyphtounicode{blacksmilingface}{263B}
\pdfglyphtounicode{blacksquare}{25A0}
\pdfglyphtounicode{blackstar}{2605}
\pdfglyphtounicode{blackupperlefttriangle}{25E4}
\pdfglyphtounicode{blackupperrighttriangle}{25E5}
\pdfglyphtounicode{blackuppointingsmalltriangle}{25B4}
\pdfglyphtounicode{blackuppointingtriangle}{25B2}
\pdfglyphtounicode{blank}{2423}
\pdfglyphtounicode{blinebelow}{1E07}
\pdfglyphtounicode{block}{2588}
\pdfglyphtounicode{bmonospace}{FF42}
\pdfglyphtounicode{bobaimaithai}{0E1A}
\pdfglyphtounicode{bohiragana}{307C}
\pdfglyphtounicode{bokatakana}{30DC}
\pdfglyphtounicode{bparen}{249D}
\pdfglyphtounicode{bqsquare}{33C3}
\pdfglyphtounicode{braceex}{F8F4}
\pdfglyphtounicode{braceleft}{007B}
\pdfglyphtounicode{braceleftbt}{F8F3}
\pdfglyphtounicode{braceleftmid}{F8F2}
\pdfglyphtounicode{braceleftmonospace}{FF5B}
\pdfglyphtounicode{braceleftsmall}{FE5B}
\pdfglyphtounicode{bracelefttp}{F8F1}
\pdfglyphtounicode{braceleftvertical}{FE37}
\pdfglyphtounicode{braceright}{007D}
\pdfglyphtounicode{bracerightbt}{F8FE}
\pdfglyphtounicode{bracerightmid}{F8FD}
\pdfglyphtounicode{bracerightmonospace}{FF5D}
\pdfglyphtounicode{bracerightsmall}{FE5C}
\pdfglyphtounicode{bracerighttp}{F8FC}
\pdfglyphtounicode{bracerightvertical}{FE38}
\pdfglyphtounicode{bracketleft}{005B}
\pdfglyphtounicode{bracketleftbt}{F8F0}
\pdfglyphtounicode{bracketleftex}{F8EF}
\pdfglyphtounicode{bracketleftmonospace}{FF3B}
\pdfglyphtounicode{bracketlefttp}{F8EE}
\pdfglyphtounicode{bracketright}{005D}
\pdfglyphtounicode{bracketrightbt}{F8FB}
\pdfglyphtounicode{bracketrightex}{F8FA}
\pdfglyphtounicode{bracketrightmonospace}{FF3D}
\pdfglyphtounicode{bracketrighttp}{F8F9}
\pdfglyphtounicode{breve}{02D8}
\pdfglyphtounicode{brevebelowcmb}{032E}
\pdfglyphtounicode{brevecmb}{0306}
\pdfglyphtounicode{breveinvertedbelowcmb}{032F}
\pdfglyphtounicode{breveinvertedcmb}{0311}
\pdfglyphtounicode{breveinverteddoublecmb}{0361}
\pdfglyphtounicode{bridgebelowcmb}{032A}
\pdfglyphtounicode{bridgeinvertedbelowcmb}{033A}
\pdfglyphtounicode{brokenbar}{00A6}
\pdfglyphtounicode{bstroke}{0180}
\pdfglyphtounicode{bsuperior}{0062}
\pdfglyphtounicode{btopbar}{0183}
\pdfglyphtounicode{buhiragana}{3076}
\pdfglyphtounicode{bukatakana}{30D6}
\pdfglyphtounicode{bullet}{2022}
\pdfglyphtounicode{bulletinverse}{25D8}
\pdfglyphtounicode{bulletoperator}{2219}
\pdfglyphtounicode{bullseye}{25CE}
\pdfglyphtounicode{c}{0063}
\pdfglyphtounicode{caarmenian}{056E}
\pdfglyphtounicode{cabengali}{099A}
\pdfglyphtounicode{cacute}{0107}
\pdfglyphtounicode{cadeva}{091A}
\pdfglyphtounicode{cagujarati}{0A9A}
\pdfglyphtounicode{cagurmukhi}{0A1A}
\pdfglyphtounicode{calsquare}{3388}
\pdfglyphtounicode{candrabindubengali}{0981}
\pdfglyphtounicode{candrabinducmb}{0310}
\pdfglyphtounicode{candrabindudeva}{0901}
\pdfglyphtounicode{candrabindugujarati}{0A81}
\pdfglyphtounicode{capslock}{21EA}
\pdfglyphtounicode{careof}{2105}
\pdfglyphtounicode{caron}{02C7}
\pdfglyphtounicode{caronbelowcmb}{032C}
\pdfglyphtounicode{caroncmb}{030C}
\pdfglyphtounicode{carriagereturn}{21B5}
\pdfglyphtounicode{cbopomofo}{3118}
\pdfglyphtounicode{ccaron}{010D}
\pdfglyphtounicode{ccedilla}{00E7}
\pdfglyphtounicode{ccedillaacute}{1E09}
\pdfglyphtounicode{ccircle}{24D2}
\pdfglyphtounicode{ccircumflex}{0109}
\pdfglyphtounicode{ccurl}{0255}
\pdfglyphtounicode{cdot}{010B}
\pdfglyphtounicode{cdotaccent}{010B}
\pdfglyphtounicode{cdsquare}{33C5}
\pdfglyphtounicode{cedilla}{00B8}
\pdfglyphtounicode{cedillacmb}{0327}
\pdfglyphtounicode{ceilingleft}{2308}
\pdfglyphtounicode{ceilingright}{2309}
\pdfglyphtounicode{cent}{00A2}
\pdfglyphtounicode{centigrade}{2103}
\pdfglyphtounicode{centinferior}{00A2}
\pdfglyphtounicode{centmonospace}{FFE0}
\pdfglyphtounicode{centoldstyle}{00A2}
\pdfglyphtounicode{centsuperior}{00A2}
\pdfglyphtounicode{chaarmenian}{0579}
\pdfglyphtounicode{chabengali}{099B}
\pdfglyphtounicode{chadeva}{091B}
\pdfglyphtounicode{chagujarati}{0A9B}
\pdfglyphtounicode{chagurmukhi}{0A1B}
\pdfglyphtounicode{chbopomofo}{3114}
\pdfglyphtounicode{cheabkhasiancyrillic}{04BD}
\pdfglyphtounicode{check}{2713}
\pdfglyphtounicode{checkmark}{2713}
\pdfglyphtounicode{checyrillic}{0447}
\pdfglyphtounicode{chedescenderabkhasiancyrillic}{04BF}
\pdfglyphtounicode{chedescendercyrillic}{04B7}
\pdfglyphtounicode{chedieresiscyrillic}{04F5}
\pdfglyphtounicode{cheharmenian}{0573}
\pdfglyphtounicode{chekhakassiancyrillic}{04CC}
\pdfglyphtounicode{cheverticalstrokecyrillic}{04B9}
\pdfglyphtounicode{chi}{03C7}
\pdfglyphtounicode{chieuchacirclekorean}{3277}
\pdfglyphtounicode{chieuchaparenkorean}{3217}
\pdfglyphtounicode{chieuchcirclekorean}{3269}
\pdfglyphtounicode{chieuchkorean}{314A}
\pdfglyphtounicode{chieuchparenkorean}{3209}
\pdfglyphtounicode{chochangthai}{0E0A}
\pdfglyphtounicode{chochanthai}{0E08}
\pdfglyphtounicode{chochingthai}{0E09}
\pdfglyphtounicode{chochoethai}{0E0C}
\pdfglyphtounicode{chook}{0188}
\pdfglyphtounicode{cieucacirclekorean}{3276}
\pdfglyphtounicode{cieucaparenkorean}{3216}
\pdfglyphtounicode{cieuccirclekorean}{3268}
\pdfglyphtounicode{cieuckorean}{3148}
\pdfglyphtounicode{cieucparenkorean}{3208}
\pdfglyphtounicode{cieucuparenkorean}{321C}
\pdfglyphtounicode{circle}{25CB}
\pdfglyphtounicode{circleR}{00AE}
\pdfglyphtounicode{circleS}{24C8}
\pdfglyphtounicode{circleasterisk}{229B}
\pdfglyphtounicode{circlecopyrt}{20DD}
\pdfglyphtounicode{circledivide}{2298}
\pdfglyphtounicode{circledot}{2299}
\pdfglyphtounicode{circleequal}{229C}
\pdfglyphtounicode{circleminus}{2296}
\pdfglyphtounicode{circlemultiply}{2297}
\pdfglyphtounicode{circleot}{2299}
\pdfglyphtounicode{circleplus}{2295}
\pdfglyphtounicode{circlepostalmark}{3036}
\pdfglyphtounicode{circlering}{229A}
\pdfglyphtounicode{circlewithlefthalfblack}{25D0}
\pdfglyphtounicode{circlewithrighthalfblack}{25D1}
\pdfglyphtounicode{circumflex}{02C6}
\pdfglyphtounicode{circumflexbelowcmb}{032D}
\pdfglyphtounicode{circumflexcmb}{0302}
\pdfglyphtounicode{clear}{2327}
\pdfglyphtounicode{clickalveolar}{01C2}
\pdfglyphtounicode{clickdental}{01C0}
\pdfglyphtounicode{clicklateral}{01C1}
\pdfglyphtounicode{clickretroflex}{01C3}
\pdfglyphtounicode{clockwise}{27F3}
\pdfglyphtounicode{club}{2663}
\pdfglyphtounicode{clubsuitblack}{2663}
\pdfglyphtounicode{clubsuitwhite}{2667}
\pdfglyphtounicode{cmcubedsquare}{33A4}
\pdfglyphtounicode{cmonospace}{FF43}
\pdfglyphtounicode{cmsquaredsquare}{33A0}
\pdfglyphtounicode{coarmenian}{0581}
\pdfglyphtounicode{colon}{003A}
\pdfglyphtounicode{colonmonetary}{20A1}
\pdfglyphtounicode{colonmonospace}{FF1A}
\pdfglyphtounicode{colonsign}{20A1}
\pdfglyphtounicode{colonsmall}{FE55}
\pdfglyphtounicode{colontriangularhalfmod}{02D1}
\pdfglyphtounicode{colontriangularmod}{02D0}
\pdfglyphtounicode{comma}{002C}
\pdfglyphtounicode{commaabovecmb}{0313}
\pdfglyphtounicode{commaaboverightcmb}{0315}
\pdfglyphtounicode{commaaccent}{F6C3}
\pdfglyphtounicode{commaarabic}{060C}
\pdfglyphtounicode{commaarmenian}{055D}
\pdfglyphtounicode{commainferior}{002C}
\pdfglyphtounicode{commamonospace}{FF0C}
\pdfglyphtounicode{commareversedabovecmb}{0314}
\pdfglyphtounicode{commareversedmod}{02BD}
\pdfglyphtounicode{commasmall}{FE50}
\pdfglyphtounicode{commasuperior}{002C}
\pdfglyphtounicode{commaturnedabovecmb}{0312}
\pdfglyphtounicode{commaturnedmod}{02BB}
\pdfglyphtounicode{compass}{263C}
\pdfglyphtounicode{complement}{2201}
\pdfglyphtounicode{compwordmark}{200C}
\pdfglyphtounicode{congruent}{2245}
\pdfglyphtounicode{contourintegral}{222E}
\pdfglyphtounicode{control}{2303}
\pdfglyphtounicode{controlACK}{0006}
\pdfglyphtounicode{controlBEL}{0007}
\pdfglyphtounicode{controlBS}{0008}
\pdfglyphtounicode{controlCAN}{0018}
\pdfglyphtounicode{controlCR}{000D}
\pdfglyphtounicode{controlDC1}{0011}
\pdfglyphtounicode{controlDC2}{0012}
\pdfglyphtounicode{controlDC3}{0013}
\pdfglyphtounicode{controlDC4}{0014}
\pdfglyphtounicode{controlDEL}{007F}
\pdfglyphtounicode{controlDLE}{0010}
\pdfglyphtounicode{controlEM}{0019}
\pdfglyphtounicode{controlENQ}{0005}
\pdfglyphtounicode{controlEOT}{0004}
\pdfglyphtounicode{controlESC}{001B}
\pdfglyphtounicode{controlETB}{0017}
\pdfglyphtounicode{controlETX}{0003}
\pdfglyphtounicode{controlFF}{000C}
\pdfglyphtounicode{controlFS}{001C}
\pdfglyphtounicode{controlGS}{001D}
\pdfglyphtounicode{controlHT}{0009}
\pdfglyphtounicode{controlLF}{000A}
\pdfglyphtounicode{controlNAK}{0015}
\pdfglyphtounicode{controlRS}{001E}
\pdfglyphtounicode{controlSI}{000F}
\pdfglyphtounicode{controlSO}{000E}
\pdfglyphtounicode{controlSOT}{0002}
\pdfglyphtounicode{controlSTX}{0001}
\pdfglyphtounicode{controlSUB}{001A}
\pdfglyphtounicode{controlSYN}{0016}
\pdfglyphtounicode{controlUS}{001F}
\pdfglyphtounicode{controlVT}{000B}
\pdfglyphtounicode{coproduct}{2A3F}
\pdfglyphtounicode{copyright}{00A9}
\pdfglyphtounicode{copyrightsans}{00A9}
\pdfglyphtounicode{copyrightserif}{00A9}
\pdfglyphtounicode{cornerbracketleft}{300C}
\pdfglyphtounicode{cornerbracketlefthalfwidth}{FF62}
\pdfglyphtounicode{cornerbracketleftvertical}{FE41}
\pdfglyphtounicode{cornerbracketright}{300D}
\pdfglyphtounicode{cornerbracketrighthalfwidth}{FF63}
\pdfglyphtounicode{cornerbracketrightvertical}{FE42}
\pdfglyphtounicode{corporationsquare}{337F}
\pdfglyphtounicode{cosquare}{33C7}
\pdfglyphtounicode{coverkgsquare}{33C6}
\pdfglyphtounicode{cparen}{249E}
\pdfglyphtounicode{cruzeiro}{20A2}
\pdfglyphtounicode{cstretched}{0297}
\pdfglyphtounicode{ct}{0063 0074}
\pdfglyphtounicode{curlyand}{22CF}
\pdfglyphtounicode{curlyleft}{21AB}
\pdfglyphtounicode{curlyor}{22CE}
\pdfglyphtounicode{curlyright}{21AC}
\pdfglyphtounicode{currency}{00A4}
\pdfglyphtounicode{cwm}{200C}
\pdfglyphtounicode{cyrBreve}{02D8}
\pdfglyphtounicode{cyrFlex}{00A0 0311}
\pdfglyphtounicode{cyrbreve}{02D8}
\pdfglyphtounicode{cyrflex}{00A0 0311}
\pdfglyphtounicode{d}{0064}
\pdfglyphtounicode{daarmenian}{0564}
\pdfglyphtounicode{dabengali}{09A6}
\pdfglyphtounicode{dadarabic}{0636}
\pdfglyphtounicode{dadeva}{0926}
\pdfglyphtounicode{dadfinalarabic}{FEBE}
\pdfglyphtounicode{dadinitialarabic}{FEBF}
\pdfglyphtounicode{dadmedialarabic}{FEC0}
\pdfglyphtounicode{dagesh}{05BC}
\pdfglyphtounicode{dageshhebrew}{05BC}
\pdfglyphtounicode{dagger}{2020}
\pdfglyphtounicode{daggerdbl}{2021}
\pdfglyphtounicode{dagujarati}{0AA6}
\pdfglyphtounicode{dagurmukhi}{0A26}
\pdfglyphtounicode{dahiragana}{3060}
\pdfglyphtounicode{dakatakana}{30C0}
\pdfglyphtounicode{dalarabic}{062F}
\pdfglyphtounicode{dalet}{05D3}
\pdfglyphtounicode{daletdagesh}{FB33}
\pdfglyphtounicode{daletdageshhebrew}{FB33}
\pdfglyphtounicode{daleth}{2138}
\pdfglyphtounicode{dalethatafpatah}{05D3 05B2}
\pdfglyphtounicode{dalethatafpatahhebrew}{05D3 05B2}
\pdfglyphtounicode{dalethatafsegol}{05D3 05B1}
\pdfglyphtounicode{dalethatafsegolhebrew}{05D3 05B1}
\pdfglyphtounicode{dalethebrew}{05D3}
\pdfglyphtounicode{dalethiriq}{05D3 05B4}
\pdfglyphtounicode{dalethiriqhebrew}{05D3 05B4}
\pdfglyphtounicode{daletholam}{05D3 05B9}
\pdfglyphtounicode{daletholamhebrew}{05D3 05B9}
\pdfglyphtounicode{daletpatah}{05D3 05B7}
\pdfglyphtounicode{daletpatahhebrew}{05D3 05B7}
\pdfglyphtounicode{daletqamats}{05D3 05B8}
\pdfglyphtounicode{daletqamatshebrew}{05D3 05B8}
\pdfglyphtounicode{daletqubuts}{05D3 05BB}
\pdfglyphtounicode{daletqubutshebrew}{05D3 05BB}
\pdfglyphtounicode{daletsegol}{05D3 05B6}
\pdfglyphtounicode{daletsegolhebrew}{05D3 05B6}
\pdfglyphtounicode{daletsheva}{05D3 05B0}
\pdfglyphtounicode{daletshevahebrew}{05D3 05B0}
\pdfglyphtounicode{dalettsere}{05D3 05B5}
\pdfglyphtounicode{dalettserehebrew}{05D3 05B5}
\pdfglyphtounicode{dalfinalarabic}{FEAA}
\pdfglyphtounicode{dammaarabic}{064F}
\pdfglyphtounicode{dammalowarabic}{064F}
\pdfglyphtounicode{dammatanaltonearabic}{064C}
\pdfglyphtounicode{dammatanarabic}{064C}
\pdfglyphtounicode{danda}{0964}
\pdfglyphtounicode{dargahebrew}{05A7}
\pdfglyphtounicode{dargalefthebrew}{05A7}
\pdfglyphtounicode{dasiapneumatacyrilliccmb}{0485}
\pdfglyphtounicode{dbar}{0111}
\pdfglyphtounicode{dblGrave}{00A0 030F}
\pdfglyphtounicode{dblanglebracketleft}{300A}
\pdfglyphtounicode{dblanglebracketleftvertical}{FE3D}
\pdfglyphtounicode{dblanglebracketright}{300B}
\pdfglyphtounicode{dblanglebracketrightvertical}{FE3E}
\pdfglyphtounicode{dblarchinvertedbelowcmb}{032B}
\pdfglyphtounicode{dblarrowdwn}{21CA}
\pdfglyphtounicode{dblarrowheadleft}{219E}
\pdfglyphtounicode{dblarrowheadright}{21A0}
\pdfglyphtounicode{dblarrowleft}{21D4}
\pdfglyphtounicode{dblarrowright}{21D2}
\pdfglyphtounicode{dblarrowup}{21C8}
\pdfglyphtounicode{dblbracketleft}{27E6}
\pdfglyphtounicode{dblbracketright}{27E7}
\pdfglyphtounicode{dbldanda}{0965}
\pdfglyphtounicode{dblgrave}{00A0 030F}
\pdfglyphtounicode{dblgravecmb}{030F}
\pdfglyphtounicode{dblintegral}{222C}
\pdfglyphtounicode{dbllowline}{2017}
\pdfglyphtounicode{dbllowlinecmb}{0333}
\pdfglyphtounicode{dbloverlinecmb}{033F}
\pdfglyphtounicode{dblprimemod}{02BA}
\pdfglyphtounicode{dblverticalbar}{2016}
\pdfglyphtounicode{dblverticallineabovecmb}{030E}
\pdfglyphtounicode{dbopomofo}{3109}
\pdfglyphtounicode{dbsquare}{33C8}
\pdfglyphtounicode{dcaron}{010F}
\pdfglyphtounicode{dcedilla}{1E11}
\pdfglyphtounicode{dcircle}{24D3}
\pdfglyphtounicode{dcircumflexbelow}{1E13}
\pdfglyphtounicode{dcroat}{0111}
\pdfglyphtounicode{ddabengali}{09A1}
\pdfglyphtounicode{ddadeva}{0921}
\pdfglyphtounicode{ddagujarati}{0AA1}
\pdfglyphtounicode{ddagurmukhi}{0A21}
\pdfglyphtounicode{ddalarabic}{0688}
\pdfglyphtounicode{ddalfinalarabic}{FB89}
\pdfglyphtounicode{dddhadeva}{095C}
\pdfglyphtounicode{ddhabengali}{09A2}
\pdfglyphtounicode{ddhadeva}{0922}
\pdfglyphtounicode{ddhagujarati}{0AA2}
\pdfglyphtounicode{ddhagurmukhi}{0A22}
\pdfglyphtounicode{ddotaccent}{1E0B}
\pdfglyphtounicode{ddotbelow}{1E0D}
\pdfglyphtounicode{decimalseparatorarabic}{066B}
\pdfglyphtounicode{decimalseparatorpersian}{066B}
\pdfglyphtounicode{decyrillic}{0434}
\pdfglyphtounicode{defines}{225C}
\pdfglyphtounicode{degree}{00B0}
\pdfglyphtounicode{dehihebrew}{05AD}
\pdfglyphtounicode{dehiragana}{3067}
\pdfglyphtounicode{deicoptic}{03EF}
\pdfglyphtounicode{dekatakana}{30C7}
\pdfglyphtounicode{deleteleft}{232B}
\pdfglyphtounicode{deleteright}{2326}
\pdfglyphtounicode{delta}{03B4}
\pdfglyphtounicode{deltaturned}{018D}
\pdfglyphtounicode{denominatorminusonenumeratorbengali}{09F8}
\pdfglyphtounicode{dezh}{02A4}
\pdfglyphtounicode{dhabengali}{09A7}
\pdfglyphtounicode{dhadeva}{0927}
\pdfglyphtounicode{dhagujarati}{0AA7}
\pdfglyphtounicode{dhagurmukhi}{0A27}
\pdfglyphtounicode{dhook}{0257}
\pdfglyphtounicode{dialytikatonos}{0385}
\pdfglyphtounicode{dialytikatonoscmb}{0344}
\pdfglyphtounicode{diamond}{2666}
\pdfglyphtounicode{diamondmath}{22C4}
\pdfglyphtounicode{diamondsolid}{2666}
\pdfglyphtounicode{diamondsuitwhite}{2662}
\pdfglyphtounicode{dieresis}{00A8}
\pdfglyphtounicode{dieresisacute}{00A0 0308 0301}
\pdfglyphtounicode{dieresisbelowcmb}{0324}
\pdfglyphtounicode{dieresiscmb}{0308}
\pdfglyphtounicode{dieresisgrave}{00A0 0308 0300}
\pdfglyphtounicode{dieresistonos}{0385}
\pdfglyphtounicode{difference}{224F}
\pdfglyphtounicode{dihiragana}{3062}
\pdfglyphtounicode{dikatakana}{30C2}
\pdfglyphtounicode{dittomark}{3003}
\pdfglyphtounicode{divide}{00F7}
\pdfglyphtounicode{dividemultiply}{22C7}
\pdfglyphtounicode{divides}{2223}
\pdfglyphtounicode{divisionslash}{2215}
\pdfglyphtounicode{djecyrillic}{0452}
\pdfglyphtounicode{dkshade}{2593}
\pdfglyphtounicode{dlinebelow}{1E0F}
\pdfglyphtounicode{dlsquare}{3397}
\pdfglyphtounicode{dmacron}{0111}
\pdfglyphtounicode{dmonospace}{FF44}
\pdfglyphtounicode{dnblock}{2584}
\pdfglyphtounicode{dochadathai}{0E0E}
\pdfglyphtounicode{dodekthai}{0E14}
\pdfglyphtounicode{dohiragana}{3069}
\pdfglyphtounicode{dokatakana}{30C9}
\pdfglyphtounicode{dollar}{0024}
\pdfglyphtounicode{dollarinferior}{0024}
\pdfglyphtounicode{dollarmonospace}{FF04}
\pdfglyphtounicode{dollaroldstyle}{0024}
\pdfglyphtounicode{dollarsmall}{FE69}
\pdfglyphtounicode{dollarsuperior}{0024}
\pdfglyphtounicode{dong}{20AB}
\pdfglyphtounicode{dorusquare}{3326}
\pdfglyphtounicode{dotaccent}{02D9}
\pdfglyphtounicode{dotaccentcmb}{0307}
\pdfglyphtounicode{dotbelowcmb}{0323}
\pdfglyphtounicode{dotbelowcomb}{0323}
\pdfglyphtounicode{dotkatakana}{30FB}
\pdfglyphtounicode{dotlessi}{0131}
\pdfglyphtounicode{dotlessj}{0237}
\pdfglyphtounicode{dotlessjstrokehook}{0284}
\pdfglyphtounicode{dotmath}{22C5}
\pdfglyphtounicode{dotplus}{2214}
\pdfglyphtounicode{dottedcircle}{25CC}
\pdfglyphtounicode{doubleyodpatah}{FB1F}
\pdfglyphtounicode{doubleyodpatahhebrew}{FB1F}
\pdfglyphtounicode{downfall}{22CE}
\pdfglyphtounicode{downslope}{29F9}
\pdfglyphtounicode{downtackbelowcmb}{031E}
\pdfglyphtounicode{downtackmod}{02D5}
\pdfglyphtounicode{dparen}{249F}
\pdfglyphtounicode{dsuperior}{0064}
\pdfglyphtounicode{dtail}{0256}
\pdfglyphtounicode{dtopbar}{018C}
\pdfglyphtounicode{duhiragana}{3065}
\pdfglyphtounicode{dukatakana}{30C5}
\pdfglyphtounicode{dz}{01F3}
\pdfglyphtounicode{dzaltone}{02A3}
\pdfglyphtounicode{dzcaron}{01C6}
\pdfglyphtounicode{dzcurl}{02A5}
\pdfglyphtounicode{dzeabkhasiancyrillic}{04E1}
\pdfglyphtounicode{dzecyrillic}{0455}
\pdfglyphtounicode{dzhecyrillic}{045F}
\pdfglyphtounicode{e}{0065}
\pdfglyphtounicode{eacute}{00E9}
\pdfglyphtounicode{earth}{2641}
\pdfglyphtounicode{ebengali}{098F}
\pdfglyphtounicode{ebopomofo}{311C}
\pdfglyphtounicode{ebreve}{0115}
\pdfglyphtounicode{ecandradeva}{090D}
\pdfglyphtounicode{ecandragujarati}{0A8D}
\pdfglyphtounicode{ecandravowelsigndeva}{0945}
\pdfglyphtounicode{ecandravowelsigngujarati}{0AC5}
\pdfglyphtounicode{ecaron}{011B}
\pdfglyphtounicode{ecedillabreve}{1E1D}
\pdfglyphtounicode{echarmenian}{0565}
\pdfglyphtounicode{echyiwnarmenian}{0587}
\pdfglyphtounicode{ecircle}{24D4}
\pdfglyphtounicode{ecircumflex}{00EA}
\pdfglyphtounicode{ecircumflexacute}{1EBF}
\pdfglyphtounicode{ecircumflexbelow}{1E19}
\pdfglyphtounicode{ecircumflexdotbelow}{1EC7}
\pdfglyphtounicode{ecircumflexgrave}{1EC1}
\pdfglyphtounicode{ecircumflexhookabove}{1EC3}
\pdfglyphtounicode{ecircumflextilde}{1EC5}
\pdfglyphtounicode{ecyrillic}{0454}
\pdfglyphtounicode{edblgrave}{0205}
\pdfglyphtounicode{edeva}{090F}
\pdfglyphtounicode{edieresis}{00EB}
\pdfglyphtounicode{edot}{0117}
\pdfglyphtounicode{edotaccent}{0117}
\pdfglyphtounicode{edotbelow}{1EB9}
\pdfglyphtounicode{eegurmukhi}{0A0F}
\pdfglyphtounicode{eematragurmukhi}{0A47}
\pdfglyphtounicode{efcyrillic}{0444}
\pdfglyphtounicode{egrave}{00E8}
\pdfglyphtounicode{egujarati}{0A8F}
\pdfglyphtounicode{eharmenian}{0567}
\pdfglyphtounicode{ehbopomofo}{311D}
\pdfglyphtounicode{ehiragana}{3048}
\pdfglyphtounicode{ehookabove}{1EBB}
\pdfglyphtounicode{eibopomofo}{311F}
\pdfglyphtounicode{eight}{0038}
\pdfglyphtounicode{eightarabic}{0668}
\pdfglyphtounicode{eightbengali}{09EE}
\pdfglyphtounicode{eightcircle}{2467}
\pdfglyphtounicode{eightcircleinversesansserif}{2791}
\pdfglyphtounicode{eightdeva}{096E}
\pdfglyphtounicode{eighteencircle}{2471}
\pdfglyphtounicode{eighteenparen}{2485}
\pdfglyphtounicode{eighteenperiod}{2499}
\pdfglyphtounicode{eightgujarati}{0AEE}
\pdfglyphtounicode{eightgurmukhi}{0A6E}
\pdfglyphtounicode{eighthackarabic}{0668}
\pdfglyphtounicode{eighthangzhou}{3028}
\pdfglyphtounicode{eighthnotebeamed}{266B}
\pdfglyphtounicode{eightideographicparen}{3227}
\pdfglyphtounicode{eightinferior}{2088}
\pdfglyphtounicode{eightmonospace}{FF18}
\pdfglyphtounicode{eightoldstyle}{0038}
\pdfglyphtounicode{eightparen}{247B}
\pdfglyphtounicode{eightperiod}{248F}
\pdfglyphtounicode{eightpersian}{06F8}
\pdfglyphtounicode{eightroman}{2177}
\pdfglyphtounicode{eightsuperior}{2078}
\pdfglyphtounicode{eightthai}{0E58}
\pdfglyphtounicode{einvertedbreve}{0207}
\pdfglyphtounicode{eiotifiedcyrillic}{0465}
\pdfglyphtounicode{ekatakana}{30A8}
\pdfglyphtounicode{ekatakanahalfwidth}{FF74}
\pdfglyphtounicode{ekonkargurmukhi}{0A74}
\pdfglyphtounicode{ekorean}{3154}
\pdfglyphtounicode{elcyrillic}{043B}
\pdfglyphtounicode{element}{2208}
\pdfglyphtounicode{elevencircle}{246A}
\pdfglyphtounicode{elevenparen}{247E}
\pdfglyphtounicode{elevenperiod}{2492}
\pdfglyphtounicode{elevenroman}{217A}
\pdfglyphtounicode{ellipsis}{2026}
\pdfglyphtounicode{ellipsisvertical}{22EE}
\pdfglyphtounicode{emacron}{0113}
\pdfglyphtounicode{emacronacute}{1E17}
\pdfglyphtounicode{emacrongrave}{1E15}
\pdfglyphtounicode{emcyrillic}{043C}
\pdfglyphtounicode{emdash}{2014}
\pdfglyphtounicode{emdashvertical}{FE31}
\pdfglyphtounicode{emonospace}{FF45}
\pdfglyphtounicode{emphasismarkarmenian}{055B}
\pdfglyphtounicode{emptyset}{2205}
\pdfglyphtounicode{enbopomofo}{3123}
\pdfglyphtounicode{encyrillic}{043D}
\pdfglyphtounicode{endash}{2013}
\pdfglyphtounicode{endashvertical}{FE32}
\pdfglyphtounicode{endescendercyrillic}{04A3}
\pdfglyphtounicode{eng}{014B}
\pdfglyphtounicode{engbopomofo}{3125}
\pdfglyphtounicode{enghecyrillic}{04A5}
\pdfglyphtounicode{enhookcyrillic}{04C8}
\pdfglyphtounicode{enspace}{2002}
\pdfglyphtounicode{eogonek}{0119}
\pdfglyphtounicode{eokorean}{3153}
\pdfglyphtounicode{eopen}{025B}
\pdfglyphtounicode{eopenclosed}{029A}
\pdfglyphtounicode{eopenreversed}{025C}
\pdfglyphtounicode{eopenreversedclosed}{025E}
\pdfglyphtounicode{eopenreversedhook}{025D}
\pdfglyphtounicode{eparen}{24A0}
\pdfglyphtounicode{epsilon}{03B5}
\pdfglyphtounicode{epsilon1}{03F5}
\pdfglyphtounicode{epsiloninv}{03F6}
\pdfglyphtounicode{epsilontonos}{03AD}
\pdfglyphtounicode{equal}{003D}
\pdfglyphtounicode{equaldotleftright}{2252}
\pdfglyphtounicode{equaldotrightleft}{2253}
\pdfglyphtounicode{equalmonospace}{FF1D}
\pdfglyphtounicode{equalorfollows}{22DF}
\pdfglyphtounicode{equalorgreater}{2A96}
\pdfglyphtounicode{equalorless}{2A95}
\pdfglyphtounicode{equalorprecedes}{22DE}
\pdfglyphtounicode{equalorsimilar}{2242}
\pdfglyphtounicode{equalsdots}{2251}
\pdfglyphtounicode{equalsmall}{FE66}
\pdfglyphtounicode{equalsuperior}{207C}
\pdfglyphtounicode{equivalence}{2261}
\pdfglyphtounicode{equivasymptotic}{224D}
\pdfglyphtounicode{erbopomofo}{3126}
\pdfglyphtounicode{ercyrillic}{0440}
\pdfglyphtounicode{ereversed}{0258}
\pdfglyphtounicode{ereversedcyrillic}{044D}
\pdfglyphtounicode{escyrillic}{0441}
\pdfglyphtounicode{esdescendercyrillic}{04AB}
\pdfglyphtounicode{esh}{0283}
\pdfglyphtounicode{eshcurl}{0286}
\pdfglyphtounicode{eshortdeva}{090E}
\pdfglyphtounicode{eshortvowelsigndeva}{0946}
\pdfglyphtounicode{eshreversedloop}{01AA}
\pdfglyphtounicode{eshsquatreversed}{0285}
\pdfglyphtounicode{esmallhiragana}{3047}
\pdfglyphtounicode{esmallkatakana}{30A7}
\pdfglyphtounicode{esmallkatakanahalfwidth}{FF6A}
\pdfglyphtounicode{estimated}{212E}
\pdfglyphtounicode{esuperior}{0065}
\pdfglyphtounicode{eta}{03B7}
\pdfglyphtounicode{etarmenian}{0568}
\pdfglyphtounicode{etatonos}{03AE}
\pdfglyphtounicode{eth}{00F0}
\pdfglyphtounicode{etilde}{1EBD}
\pdfglyphtounicode{etildebelow}{1E1B}
\pdfglyphtounicode{etnahtafoukhhebrew}{0591}
\pdfglyphtounicode{etnahtafoukhlefthebrew}{0591}
\pdfglyphtounicode{etnahtahebrew}{0591}
\pdfglyphtounicode{etnahtalefthebrew}{0591}
\pdfglyphtounicode{eturned}{01DD}
\pdfglyphtounicode{eukorean}{3161}
\pdfglyphtounicode{euro}{20AC}
\pdfglyphtounicode{evowelsignbengali}{09C7}
\pdfglyphtounicode{evowelsigndeva}{0947}
\pdfglyphtounicode{evowelsigngujarati}{0AC7}
\pdfglyphtounicode{exclam}{0021}
\pdfglyphtounicode{exclamarmenian}{055C}
\pdfglyphtounicode{exclamdbl}{203C}
\pdfglyphtounicode{exclamdown}{00A1}
\pdfglyphtounicode{exclamdownsmall}{00A1}
\pdfglyphtounicode{exclammonospace}{FF01}
\pdfglyphtounicode{exclamsmall}{0021}
\pdfglyphtounicode{existential}{2203}
\pdfglyphtounicode{ezh}{0292}
\pdfglyphtounicode{ezhcaron}{01EF}
\pdfglyphtounicode{ezhcurl}{0293}
\pdfglyphtounicode{ezhreversed}{01B9}
\pdfglyphtounicode{ezhtail}{01BA}
\pdfglyphtounicode{f}{0066}
\pdfglyphtounicode{fadeva}{095E}
\pdfglyphtounicode{fagurmukhi}{0A5E}
\pdfglyphtounicode{fahrenheit}{2109}
\pdfglyphtounicode{fathaarabic}{064E}
\pdfglyphtounicode{fathalowarabic}{064E}
\pdfglyphtounicode{fathatanarabic}{064B}
\pdfglyphtounicode{fbopomofo}{3108}
\pdfglyphtounicode{fcircle}{24D5}
\pdfglyphtounicode{fdotaccent}{1E1F}
\pdfglyphtounicode{feharabic}{0641}
\pdfglyphtounicode{feharmenian}{0586}
\pdfglyphtounicode{fehfinalarabic}{FED2}
\pdfglyphtounicode{fehinitialarabic}{FED3}
\pdfglyphtounicode{fehmedialarabic}{FED4}
\pdfglyphtounicode{feicoptic}{03E5}
\pdfglyphtounicode{female}{2640}
\pdfglyphtounicode{ff}{0066 0066}
\pdfglyphtounicode{ffi}{0066 0066 0069}
\pdfglyphtounicode{ffl}{0066 0066 006C}
\pdfglyphtounicode{fi}{0066 0069}
\pdfglyphtounicode{fifteencircle}{246E}
\pdfglyphtounicode{fifteenparen}{2482}
\pdfglyphtounicode{fifteenperiod}{2496}
\pdfglyphtounicode{figuredash}{2012}
\pdfglyphtounicode{filledbox}{25A0}
\pdfglyphtounicode{filledrect}{25AC}
\pdfglyphtounicode{finalkaf}{05DA}
\pdfglyphtounicode{finalkafdagesh}{FB3A}
\pdfglyphtounicode{finalkafdageshhebrew}{FB3A}
\pdfglyphtounicode{finalkafhebrew}{05DA}
\pdfglyphtounicode{finalkafqamats}{05DA 05B8}
\pdfglyphtounicode{finalkafqamatshebrew}{05DA 05B8}
\pdfglyphtounicode{finalkafsheva}{05DA 05B0}
\pdfglyphtounicode{finalkafshevahebrew}{05DA 05B0}
\pdfglyphtounicode{finalmem}{05DD}
\pdfglyphtounicode{finalmemhebrew}{05DD}
\pdfglyphtounicode{finalnun}{05DF}
\pdfglyphtounicode{finalnunhebrew}{05DF}
\pdfglyphtounicode{finalpe}{05E3}
\pdfglyphtounicode{finalpehebrew}{05E3}
\pdfglyphtounicode{finaltsadi}{05E5}
\pdfglyphtounicode{finaltsadihebrew}{05E5}
\pdfglyphtounicode{firsttonechinese}{02C9}
\pdfglyphtounicode{fisheye}{25C9}
\pdfglyphtounicode{fitacyrillic}{0473}
\pdfglyphtounicode{five}{0035}
\pdfglyphtounicode{fivearabic}{0665}
\pdfglyphtounicode{fivebengali}{09EB}
\pdfglyphtounicode{fivecircle}{2464}
\pdfglyphtounicode{fivecircleinversesansserif}{278E}
\pdfglyphtounicode{fivedeva}{096B}
\pdfglyphtounicode{fiveeighths}{215D}
\pdfglyphtounicode{fivegujarati}{0AEB}
\pdfglyphtounicode{fivegurmukhi}{0A6B}
\pdfglyphtounicode{fivehackarabic}{0665}
\pdfglyphtounicode{fivehangzhou}{3025}
\pdfglyphtounicode{fiveideographicparen}{3224}
\pdfglyphtounicode{fiveinferior}{2085}
\pdfglyphtounicode{fivemonospace}{FF15}
\pdfglyphtounicode{fiveoldstyle}{0035}
\pdfglyphtounicode{fiveparen}{2478}
\pdfglyphtounicode{fiveperiod}{248C}
\pdfglyphtounicode{fivepersian}{06F5}
\pdfglyphtounicode{fiveroman}{2174}
\pdfglyphtounicode{fivesuperior}{2075}
\pdfglyphtounicode{fivethai}{0E55}
\pdfglyphtounicode{fl}{0066 006C}
\pdfglyphtounicode{flat}{266D}
\pdfglyphtounicode{floorleft}{230A}
\pdfglyphtounicode{floorright}{230B}
\pdfglyphtounicode{florin}{0192}
\pdfglyphtounicode{fmonospace}{FF46}
\pdfglyphtounicode{fmsquare}{3399}
\pdfglyphtounicode{fofanthai}{0E1F}
\pdfglyphtounicode{fofathai}{0E1D}
\pdfglyphtounicode{follownotdbleqv}{2ABA}
\pdfglyphtounicode{follownotslnteql}{2AB6}
\pdfglyphtounicode{followornoteqvlnt}{22E9}
\pdfglyphtounicode{follows}{227B}
\pdfglyphtounicode{followsequal}{2AB0}
\pdfglyphtounicode{followsorcurly}{227D}
\pdfglyphtounicode{followsorequal}{227F}
\pdfglyphtounicode{fongmanthai}{0E4F}
\pdfglyphtounicode{forall}{2200}
\pdfglyphtounicode{forces}{22A9}
\pdfglyphtounicode{forcesbar}{22AA}
\pdfglyphtounicode{fork}{22D4}
\pdfglyphtounicode{four}{0034}
\pdfglyphtounicode{fourarabic}{0664}
\pdfglyphtounicode{fourbengali}{09EA}
\pdfglyphtounicode{fourcircle}{2463}
\pdfglyphtounicode{fourcircleinversesansserif}{278D}
\pdfglyphtounicode{fourdeva}{096A}
\pdfglyphtounicode{fourgujarati}{0AEA}
\pdfglyphtounicode{fourgurmukhi}{0A6A}
\pdfglyphtounicode{fourhackarabic}{0664}
\pdfglyphtounicode{fourhangzhou}{3024}
\pdfglyphtounicode{fourideographicparen}{3223}
\pdfglyphtounicode{fourinferior}{2084}
\pdfglyphtounicode{fourmonospace}{FF14}
\pdfglyphtounicode{fournumeratorbengali}{09F7}
\pdfglyphtounicode{fouroldstyle}{0034}
\pdfglyphtounicode{fourparen}{2477}
\pdfglyphtounicode{fourperiod}{248B}
\pdfglyphtounicode{fourpersian}{06F4}
\pdfglyphtounicode{fourroman}{2173}
\pdfglyphtounicode{foursuperior}{2074}
\pdfglyphtounicode{fourteencircle}{246D}
\pdfglyphtounicode{fourteenparen}{2481}
\pdfglyphtounicode{fourteenperiod}{2495}
\pdfglyphtounicode{fourthai}{0E54}
\pdfglyphtounicode{fourthtonechinese}{02CB}
\pdfglyphtounicode{fparen}{24A1}
\pdfglyphtounicode{fraction}{2044}
\pdfglyphtounicode{franc}{20A3}
\pdfglyphtounicode{frown}{2322}
\pdfglyphtounicode{g}{0067}
\pdfglyphtounicode{gabengali}{0997}
\pdfglyphtounicode{gacute}{01F5}
\pdfglyphtounicode{gadeva}{0917}
\pdfglyphtounicode{gafarabic}{06AF}
\pdfglyphtounicode{gaffinalarabic}{FB93}
\pdfglyphtounicode{gafinitialarabic}{FB94}
\pdfglyphtounicode{gafmedialarabic}{FB95}
\pdfglyphtounicode{gagujarati}{0A97}
\pdfglyphtounicode{gagurmukhi}{0A17}
\pdfglyphtounicode{gahiragana}{304C}
\pdfglyphtounicode{gakatakana}{30AC}
\pdfglyphtounicode{gamma}{03B3}
\pdfglyphtounicode{gammalatinsmall}{0263}
\pdfglyphtounicode{gammasuperior}{02E0}
\pdfglyphtounicode{gangiacoptic}{03EB}
\pdfglyphtounicode{gbopomofo}{310D}
\pdfglyphtounicode{gbreve}{011F}
\pdfglyphtounicode{gcaron}{01E7}
\pdfglyphtounicode{gcedilla}{0123}
\pdfglyphtounicode{gcircle}{24D6}
\pdfglyphtounicode{gcircumflex}{011D}
\pdfglyphtounicode{gcommaaccent}{0123}
\pdfglyphtounicode{gdot}{0121}
\pdfglyphtounicode{gdotaccent}{0121}
\pdfglyphtounicode{gecyrillic}{0433}
\pdfglyphtounicode{gehiragana}{3052}
\pdfglyphtounicode{gekatakana}{30B2}
\pdfglyphtounicode{geomequivalent}{224E}
\pdfglyphtounicode{geometricallyequal}{2251}
\pdfglyphtounicode{gereshaccenthebrew}{059C}
\pdfglyphtounicode{gereshhebrew}{05F3}
\pdfglyphtounicode{gereshmuqdamhebrew}{059D}
\pdfglyphtounicode{germandbls}{00DF}
\pdfglyphtounicode{gershayimaccenthebrew}{059E}
\pdfglyphtounicode{gershayimhebrew}{05F4}
\pdfglyphtounicode{getamark}{3013}
\pdfglyphtounicode{ghabengali}{0998}
\pdfglyphtounicode{ghadarmenian}{0572}
\pdfglyphtounicode{ghadeva}{0918}
\pdfglyphtounicode{ghagujarati}{0A98}
\pdfglyphtounicode{ghagurmukhi}{0A18}
\pdfglyphtounicode{ghainarabic}{063A}
\pdfglyphtounicode{ghainfinalarabic}{FECE}
\pdfglyphtounicode{ghaininitialarabic}{FECF}
\pdfglyphtounicode{ghainmedialarabic}{FED0}
\pdfglyphtounicode{ghemiddlehookcyrillic}{0495}
\pdfglyphtounicode{ghestrokecyrillic}{0493}
\pdfglyphtounicode{gheupturncyrillic}{0491}
\pdfglyphtounicode{ghhadeva}{095A}
\pdfglyphtounicode{ghhagurmukhi}{0A5A}
\pdfglyphtounicode{ghook}{0260}
\pdfglyphtounicode{ghzsquare}{3393}
\pdfglyphtounicode{gihiragana}{304E}
\pdfglyphtounicode{gikatakana}{30AE}
\pdfglyphtounicode{gimarmenian}{0563}
\pdfglyphtounicode{gimel}{05D2}
\pdfglyphtounicode{gimeldagesh}{FB32}
\pdfglyphtounicode{gimeldageshhebrew}{FB32}
\pdfglyphtounicode{gimelhebrew}{05D2}
\pdfglyphtounicode{gjecyrillic}{0453}
\pdfglyphtounicode{glottalinvertedstroke}{01BE}
\pdfglyphtounicode{glottalstop}{0294}
\pdfglyphtounicode{glottalstopinverted}{0296}
\pdfglyphtounicode{glottalstopmod}{02C0}
\pdfglyphtounicode{glottalstopreversed}{0295}
\pdfglyphtounicode{glottalstopreversedmod}{02C1}
\pdfglyphtounicode{glottalstopreversedsuperior}{02E4}
\pdfglyphtounicode{glottalstopstroke}{02A1}
\pdfglyphtounicode{glottalstopstrokereversed}{02A2}
\pdfglyphtounicode{gmacron}{1E21}
\pdfglyphtounicode{gmonospace}{FF47}
\pdfglyphtounicode{gohiragana}{3054}
\pdfglyphtounicode{gokatakana}{30B4}
\pdfglyphtounicode{gparen}{24A2}
\pdfglyphtounicode{gpasquare}{33AC}
\pdfglyphtounicode{gradient}{2207}
\pdfglyphtounicode{grave}{0060}
\pdfglyphtounicode{gravebelowcmb}{0316}
\pdfglyphtounicode{gravecmb}{0300}
\pdfglyphtounicode{gravecomb}{0300}
\pdfglyphtounicode{gravedeva}{0953}
\pdfglyphtounicode{gravelowmod}{02CE}
\pdfglyphtounicode{gravemonospace}{FF40}
\pdfglyphtounicode{gravetonecmb}{0340}
\pdfglyphtounicode{greater}{003E}
\pdfglyphtounicode{greaterdbleqlless}{2A8C}
\pdfglyphtounicode{greaterdblequal}{2267}
\pdfglyphtounicode{greaterdot}{22D7}
\pdfglyphtounicode{greaterequal}{2265}
\pdfglyphtounicode{greaterequalorless}{22DB}
\pdfglyphtounicode{greaterlessequal}{22DB}
\pdfglyphtounicode{greatermonospace}{FF1E}
\pdfglyphtounicode{greatermuch}{226B}
\pdfglyphtounicode{greaternotdblequal}{2A8A}
\pdfglyphtounicode{greaternotequal}{2A88}
\pdfglyphtounicode{greaterorapproxeql}{2A86}
\pdfglyphtounicode{greaterorequalslant}{2A7E}
\pdfglyphtounicode{greaterorequivalent}{2273}
\pdfglyphtounicode{greaterorless}{2277}
\pdfglyphtounicode{greaterornotdbleql}{2269}
\pdfglyphtounicode{greaterornotequal}{2269}
\pdfglyphtounicode{greaterorsimilar}{2273}
\pdfglyphtounicode{greateroverequal}{2267}
\pdfglyphtounicode{greatersmall}{FE65}
\pdfglyphtounicode{gscript}{0261}
\pdfglyphtounicode{gstroke}{01E5}
\pdfglyphtounicode{guhiragana}{3050}
\pdfglyphtounicode{guillemotleft}{00AB}
\pdfglyphtounicode{guillemotright}{00BB}
\pdfglyphtounicode{guilsinglleft}{2039}
\pdfglyphtounicode{guilsinglright}{203A}
\pdfglyphtounicode{gukatakana}{30B0}
\pdfglyphtounicode{guramusquare}{3318}
\pdfglyphtounicode{gysquare}{33C9}
\pdfglyphtounicode{h}{0068}
\pdfglyphtounicode{haabkhasiancyrillic}{04A9}
\pdfglyphtounicode{haaltonearabic}{06C1}
\pdfglyphtounicode{habengali}{09B9}
\pdfglyphtounicode{hadescendercyrillic}{04B3}
\pdfglyphtounicode{hadeva}{0939}
\pdfglyphtounicode{hagujarati}{0AB9}
\pdfglyphtounicode{hagurmukhi}{0A39}
\pdfglyphtounicode{haharabic}{062D}
\pdfglyphtounicode{hahfinalarabic}{FEA2}
\pdfglyphtounicode{hahinitialarabic}{FEA3}
\pdfglyphtounicode{hahiragana}{306F}
\pdfglyphtounicode{hahmedialarabic}{FEA4}
\pdfglyphtounicode{haitusquare}{332A}
\pdfglyphtounicode{hakatakana}{30CF}
\pdfglyphtounicode{hakatakanahalfwidth}{FF8A}
\pdfglyphtounicode{halantgurmukhi}{0A4D}
\pdfglyphtounicode{hamzaarabic}{0621}
\pdfglyphtounicode{hamzadammaarabic}{0621 064F}
\pdfglyphtounicode{hamzadammatanarabic}{0621 064C}
\pdfglyphtounicode{hamzafathaarabic}{0621 064E}
\pdfglyphtounicode{hamzafathatanarabic}{0621 064B}
\pdfglyphtounicode{hamzalowarabic}{0621}
\pdfglyphtounicode{hamzalowkasraarabic}{0621 0650}
\pdfglyphtounicode{hamzalowkasratanarabic}{0621 064D}
\pdfglyphtounicode{hamzasukunarabic}{0621 0652}
\pdfglyphtounicode{hangulfiller}{3164}
\pdfglyphtounicode{hardsigncyrillic}{044A}
\pdfglyphtounicode{harpoondownleft}{21C3}
\pdfglyphtounicode{harpoondownright}{21C2}
\pdfglyphtounicode{harpoonleftbarbup}{21BC}
\pdfglyphtounicode{harpoonleftright}{21CC}
\pdfglyphtounicode{harpoonrightbarbup}{21C0}
\pdfglyphtounicode{harpoonrightleft}{21CB}
\pdfglyphtounicode{harpoonupleft}{21BF}
\pdfglyphtounicode{harpoonupright}{21BE}
\pdfglyphtounicode{hasquare}{33CA}
\pdfglyphtounicode{hatafpatah}{05B2}
\pdfglyphtounicode{hatafpatah16}{05B2}
\pdfglyphtounicode{hatafpatah23}{05B2}
\pdfglyphtounicode{hatafpatah2f}{05B2}
\pdfglyphtounicode{hatafpatahhebrew}{05B2}
\pdfglyphtounicode{hatafpatahnarrowhebrew}{05B2}
\pdfglyphtounicode{hatafpatahquarterhebrew}{05B2}
\pdfglyphtounicode{hatafpatahwidehebrew}{05B2}
\pdfglyphtounicode{hatafqamats}{05B3}
\pdfglyphtounicode{hatafqamats1b}{05B3}
\pdfglyphtounicode{hatafqamats28}{05B3}
\pdfglyphtounicode{hatafqamats34}{05B3}
\pdfglyphtounicode{hatafqamatshebrew}{05B3}
\pdfglyphtounicode{hatafqamatsnarrowhebrew}{05B3}
\pdfglyphtounicode{hatafqamatsquarterhebrew}{05B3}
\pdfglyphtounicode{hatafqamatswidehebrew}{05B3}
\pdfglyphtounicode{hatafsegol}{05B1}
\pdfglyphtounicode{hatafsegol17}{05B1}
\pdfglyphtounicode{hatafsegol24}{05B1}
\pdfglyphtounicode{hatafsegol30}{05B1}
\pdfglyphtounicode{hatafsegolhebrew}{05B1}
\pdfglyphtounicode{hatafsegolnarrowhebrew}{05B1}
\pdfglyphtounicode{hatafsegolquarterhebrew}{05B1}
\pdfglyphtounicode{hatafsegolwidehebrew}{05B1}
\pdfglyphtounicode{hbar}{0127}
\pdfglyphtounicode{hbopomofo}{310F}
\pdfglyphtounicode{hbrevebelow}{1E2B}
\pdfglyphtounicode{hcedilla}{1E29}
\pdfglyphtounicode{hcircle}{24D7}
\pdfglyphtounicode{hcircumflex}{0125}
\pdfglyphtounicode{hdieresis}{1E27}
\pdfglyphtounicode{hdotaccent}{1E23}
\pdfglyphtounicode{hdotbelow}{1E25}
\pdfglyphtounicode{he}{05D4}
\pdfglyphtounicode{heart}{2665}
\pdfglyphtounicode{heartsuitblack}{2665}
\pdfglyphtounicode{heartsuitwhite}{2661}
\pdfglyphtounicode{hedagesh}{FB34}
\pdfglyphtounicode{hedageshhebrew}{FB34}
\pdfglyphtounicode{hehaltonearabic}{06C1}
\pdfglyphtounicode{heharabic}{0647}
\pdfglyphtounicode{hehebrew}{05D4}
\pdfglyphtounicode{hehfinalaltonearabic}{FBA7}
\pdfglyphtounicode{hehfinalalttwoarabic}{FEEA}
\pdfglyphtounicode{hehfinalarabic}{FEEA}
\pdfglyphtounicode{hehhamzaabovefinalarabic}{FBA5}
\pdfglyphtounicode{hehhamzaaboveisolatedarabic}{FBA4}
\pdfglyphtounicode{hehinitialaltonearabic}{FBA8}
\pdfglyphtounicode{hehinitialarabic}{FEEB}
\pdfglyphtounicode{hehiragana}{3078}
\pdfglyphtounicode{hehmedialaltonearabic}{FBA9}
\pdfglyphtounicode{hehmedialarabic}{FEEC}
\pdfglyphtounicode{heiseierasquare}{337B}
\pdfglyphtounicode{hekatakana}{30D8}
\pdfglyphtounicode{hekatakanahalfwidth}{FF8D}
\pdfglyphtounicode{hekutaarusquare}{3336}
\pdfglyphtounicode{henghook}{0267}
\pdfglyphtounicode{herutusquare}{3339}
\pdfglyphtounicode{het}{05D7}
\pdfglyphtounicode{hethebrew}{05D7}
\pdfglyphtounicode{hhook}{0266}
\pdfglyphtounicode{hhooksuperior}{02B1}
\pdfglyphtounicode{hieuhacirclekorean}{327B}
\pdfglyphtounicode{hieuhaparenkorean}{321B}
\pdfglyphtounicode{hieuhcirclekorean}{326D}
\pdfglyphtounicode{hieuhkorean}{314E}
\pdfglyphtounicode{hieuhparenkorean}{320D}
\pdfglyphtounicode{hihiragana}{3072}
\pdfglyphtounicode{hikatakana}{30D2}
\pdfglyphtounicode{hikatakanahalfwidth}{FF8B}
\pdfglyphtounicode{hiriq}{05B4}
\pdfglyphtounicode{hiriq14}{05B4}
\pdfglyphtounicode{hiriq21}{05B4}
\pdfglyphtounicode{hiriq2d}{05B4}
\pdfglyphtounicode{hiriqhebrew}{05B4}
\pdfglyphtounicode{hiriqnarrowhebrew}{05B4}
\pdfglyphtounicode{hiriqquarterhebrew}{05B4}
\pdfglyphtounicode{hiriqwidehebrew}{05B4}
\pdfglyphtounicode{hlinebelow}{1E96}
\pdfglyphtounicode{hmonospace}{FF48}
\pdfglyphtounicode{hoarmenian}{0570}
\pdfglyphtounicode{hohipthai}{0E2B}
\pdfglyphtounicode{hohiragana}{307B}
\pdfglyphtounicode{hokatakana}{30DB}
\pdfglyphtounicode{hokatakanahalfwidth}{FF8E}
\pdfglyphtounicode{holam}{05B9}
\pdfglyphtounicode{holam19}{05B9}
\pdfglyphtounicode{holam26}{05B9}
\pdfglyphtounicode{holam32}{05B9}
\pdfglyphtounicode{holamhebrew}{05B9}
\pdfglyphtounicode{holamnarrowhebrew}{05B9}
\pdfglyphtounicode{holamquarterhebrew}{05B9}
\pdfglyphtounicode{holamwidehebrew}{05B9}
\pdfglyphtounicode{honokhukthai}{0E2E}
\pdfglyphtounicode{hookabovecomb}{0309}
\pdfglyphtounicode{hookcmb}{0309}
\pdfglyphtounicode{hookpalatalizedbelowcmb}{0321}
\pdfglyphtounicode{hookretroflexbelowcmb}{0322}
\pdfglyphtounicode{hoonsquare}{3342}
\pdfglyphtounicode{horicoptic}{03E9}
\pdfglyphtounicode{horizontalbar}{2015}
\pdfglyphtounicode{horncmb}{031B}
\pdfglyphtounicode{hotsprings}{2668}
\pdfglyphtounicode{house}{2302}
\pdfglyphtounicode{hparen}{24A3}
\pdfglyphtounicode{hsuperior}{02B0}
\pdfglyphtounicode{hturned}{0265}
\pdfglyphtounicode{huhiragana}{3075}
\pdfglyphtounicode{huiitosquare}{3333}
\pdfglyphtounicode{hukatakana}{30D5}
\pdfglyphtounicode{hukatakanahalfwidth}{FF8C}
\pdfglyphtounicode{hungarumlaut}{02DD}
\pdfglyphtounicode{hungarumlautcmb}{030B}
\pdfglyphtounicode{hv}{0195}
\pdfglyphtounicode{hyphen}{002D}
\pdfglyphtounicode{hyphenchar}{002D}
\pdfglyphtounicode{hypheninferior}{002D}
\pdfglyphtounicode{hyphenmonospace}{FF0D}
\pdfglyphtounicode{hyphensmall}{FE63}
\pdfglyphtounicode{hyphensuperior}{002D}
\pdfglyphtounicode{hyphentwo}{2010}
\pdfglyphtounicode{i}{0069}
\pdfglyphtounicode{iacute}{00ED}
\pdfglyphtounicode{iacyrillic}{044F}
\pdfglyphtounicode{ibengali}{0987}
\pdfglyphtounicode{ibopomofo}{3127}
\pdfglyphtounicode{ibreve}{012D}
\pdfglyphtounicode{icaron}{01D0}
\pdfglyphtounicode{icircle}{24D8}
\pdfglyphtounicode{icircumflex}{00EE}
\pdfglyphtounicode{icyrillic}{0456}
\pdfglyphtounicode{idblgrave}{0209}
\pdfglyphtounicode{ideographearthcircle}{328F}
\pdfglyphtounicode{ideographfirecircle}{328B}
\pdfglyphtounicode{ideographicallianceparen}{323F}
\pdfglyphtounicode{ideographiccallparen}{323A}
\pdfglyphtounicode{ideographiccentrecircle}{32A5}
\pdfglyphtounicode{ideographicclose}{3006}
\pdfglyphtounicode{ideographiccomma}{3001}
\pdfglyphtounicode{ideographiccommaleft}{FF64}
\pdfglyphtounicode{ideographiccongratulationparen}{3237}
\pdfglyphtounicode{ideographiccorrectcircle}{32A3}
\pdfglyphtounicode{ideographicearthparen}{322F}
\pdfglyphtounicode{ideographicenterpriseparen}{323D}
\pdfglyphtounicode{ideographicexcellentcircle}{329D}
\pdfglyphtounicode{ideographicfestivalparen}{3240}
\pdfglyphtounicode{ideographicfinancialcircle}{3296}
\pdfglyphtounicode{ideographicfinancialparen}{3236}
\pdfglyphtounicode{ideographicfireparen}{322B}
\pdfglyphtounicode{ideographichaveparen}{3232}
\pdfglyphtounicode{ideographichighcircle}{32A4}
\pdfglyphtounicode{ideographiciterationmark}{3005}
\pdfglyphtounicode{ideographiclaborcircle}{3298}
\pdfglyphtounicode{ideographiclaborparen}{3238}
\pdfglyphtounicode{ideographicleftcircle}{32A7}
\pdfglyphtounicode{ideographiclowcircle}{32A6}
\pdfglyphtounicode{ideographicmedicinecircle}{32A9}
\pdfglyphtounicode{ideographicmetalparen}{322E}
\pdfglyphtounicode{ideographicmoonparen}{322A}
\pdfglyphtounicode{ideographicnameparen}{3234}
\pdfglyphtounicode{ideographicperiod}{3002}
\pdfglyphtounicode{ideographicprintcircle}{329E}
\pdfglyphtounicode{ideographicreachparen}{3243}
\pdfglyphtounicode{ideographicrepresentparen}{3239}
\pdfglyphtounicode{ideographicresourceparen}{323E}
\pdfglyphtounicode{ideographicrightcircle}{32A8}
\pdfglyphtounicode{ideographicsecretcircle}{3299}
\pdfglyphtounicode{ideographicselfparen}{3242}
\pdfglyphtounicode{ideographicsocietyparen}{3233}
\pdfglyphtounicode{ideographicspace}{3000}
\pdfglyphtounicode{ideographicspecialparen}{3235}
\pdfglyphtounicode{ideographicstockparen}{3231}
\pdfglyphtounicode{ideographicstudyparen}{323B}
\pdfglyphtounicode{ideographicsunparen}{3230}
\pdfglyphtounicode{ideographicsuperviseparen}{323C}
\pdfglyphtounicode{ideographicwaterparen}{322C}
\pdfglyphtounicode{ideographicwoodparen}{322D}
\pdfglyphtounicode{ideographiczero}{3007}
\pdfglyphtounicode{ideographmetalcircle}{328E}
\pdfglyphtounicode{ideographmooncircle}{328A}
\pdfglyphtounicode{ideographnamecircle}{3294}
\pdfglyphtounicode{ideographsuncircle}{3290}
\pdfglyphtounicode{ideographwatercircle}{328C}
\pdfglyphtounicode{ideographwoodcircle}{328D}
\pdfglyphtounicode{ideva}{0907}
\pdfglyphtounicode{idieresis}{00EF}
\pdfglyphtounicode{idieresisacute}{1E2F}
\pdfglyphtounicode{idieresiscyrillic}{04E5}
\pdfglyphtounicode{idotbelow}{1ECB}
\pdfglyphtounicode{iebrevecyrillic}{04D7}
\pdfglyphtounicode{iecyrillic}{0435}
\pdfglyphtounicode{ieungacirclekorean}{3275}
\pdfglyphtounicode{ieungaparenkorean}{3215}
\pdfglyphtounicode{ieungcirclekorean}{3267}
\pdfglyphtounicode{ieungkorean}{3147}
\pdfglyphtounicode{ieungparenkorean}{3207}
\pdfglyphtounicode{igrave}{00EC}
\pdfglyphtounicode{igujarati}{0A87}
\pdfglyphtounicode{igurmukhi}{0A07}
\pdfglyphtounicode{ihiragana}{3044}
\pdfglyphtounicode{ihookabove}{1EC9}
\pdfglyphtounicode{iibengali}{0988}
\pdfglyphtounicode{iicyrillic}{0438}
\pdfglyphtounicode{iideva}{0908}
\pdfglyphtounicode{iigujarati}{0A88}
\pdfglyphtounicode{iigurmukhi}{0A08}
\pdfglyphtounicode{iimatragurmukhi}{0A40}
\pdfglyphtounicode{iinvertedbreve}{020B}
\pdfglyphtounicode{iishortcyrillic}{0439}
\pdfglyphtounicode{iivowelsignbengali}{09C0}
\pdfglyphtounicode{iivowelsigndeva}{0940}
\pdfglyphtounicode{iivowelsigngujarati}{0AC0}
\pdfglyphtounicode{ij}{0133}
\pdfglyphtounicode{ikatakana}{30A4}
\pdfglyphtounicode{ikatakanahalfwidth}{FF72}
\pdfglyphtounicode{ikorean}{3163}
\pdfglyphtounicode{ilde}{02DC}
\pdfglyphtounicode{iluyhebrew}{05AC}
\pdfglyphtounicode{imacron}{012B}
\pdfglyphtounicode{imacroncyrillic}{04E3}
\pdfglyphtounicode{imageorapproximatelyequal}{2253}
\pdfglyphtounicode{imatragurmukhi}{0A3F}
\pdfglyphtounicode{imonospace}{FF49}
\pdfglyphtounicode{increment}{2206}
\pdfglyphtounicode{infinity}{221E}
\pdfglyphtounicode{iniarmenian}{056B}
\pdfglyphtounicode{integerdivide}{2216}
\pdfglyphtounicode{integral}{222B}
\pdfglyphtounicode{integralbottom}{2321}
\pdfglyphtounicode{integralbt}{2321}
\pdfglyphtounicode{integralex}{F8F5}
\pdfglyphtounicode{integraltop}{2320}
\pdfglyphtounicode{integraltp}{2320}
\pdfglyphtounicode{intercal}{22BA}
\pdfglyphtounicode{interrobang}{203D}
\pdfglyphtounicode{interrobangdown}{2E18}
\pdfglyphtounicode{intersection}{2229}
\pdfglyphtounicode{intersectiondbl}{22D2}
\pdfglyphtounicode{intersectionsq}{2293}
\pdfglyphtounicode{intisquare}{3305}
\pdfglyphtounicode{invbullet}{25D8}
\pdfglyphtounicode{invcircle}{25D9}
\pdfglyphtounicode{invsmileface}{263B}
\pdfglyphtounicode{iocyrillic}{0451}
\pdfglyphtounicode{iogonek}{012F}
\pdfglyphtounicode{iota}{03B9}
\pdfglyphtounicode{iotadieresis}{03CA}
\pdfglyphtounicode{iotadieresistonos}{0390}
\pdfglyphtounicode{iotalatin}{0269}
\pdfglyphtounicode{iotatonos}{03AF}
\pdfglyphtounicode{iparen}{24A4}
\pdfglyphtounicode{irigurmukhi}{0A72}
\pdfglyphtounicode{ismallhiragana}{3043}
\pdfglyphtounicode{ismallkatakana}{30A3}
\pdfglyphtounicode{ismallkatakanahalfwidth}{FF68}
\pdfglyphtounicode{issharbengali}{09FA}
\pdfglyphtounicode{istroke}{0268}
\pdfglyphtounicode{isuperior}{0069}
\pdfglyphtounicode{iterationhiragana}{309D}
\pdfglyphtounicode{iterationkatakana}{30FD}
\pdfglyphtounicode{itilde}{0129}
\pdfglyphtounicode{itildebelow}{1E2D}
\pdfglyphtounicode{iubopomofo}{3129}
\pdfglyphtounicode{iucyrillic}{044E}
\pdfglyphtounicode{ivowelsignbengali}{09BF}
\pdfglyphtounicode{ivowelsigndeva}{093F}
\pdfglyphtounicode{ivowelsigngujarati}{0ABF}
\pdfglyphtounicode{izhitsacyrillic}{0475}
\pdfglyphtounicode{izhitsadblgravecyrillic}{0477}
\pdfglyphtounicode{j}{006A}
\pdfglyphtounicode{jaarmenian}{0571}
\pdfglyphtounicode{jabengali}{099C}
\pdfglyphtounicode{jadeva}{091C}
\pdfglyphtounicode{jagujarati}{0A9C}
\pdfglyphtounicode{jagurmukhi}{0A1C}
\pdfglyphtounicode{jbopomofo}{3110}
\pdfglyphtounicode{jcaron}{01F0}
\pdfglyphtounicode{jcircle}{24D9}
\pdfglyphtounicode{jcircumflex}{0135}
\pdfglyphtounicode{jcrossedtail}{029D}
\pdfglyphtounicode{jdotlessstroke}{025F}
\pdfglyphtounicode{jecyrillic}{0458}
\pdfglyphtounicode{jeemarabic}{062C}
\pdfglyphtounicode{jeemfinalarabic}{FE9E}
\pdfglyphtounicode{jeeminitialarabic}{FE9F}
\pdfglyphtounicode{jeemmedialarabic}{FEA0}
\pdfglyphtounicode{jeharabic}{0698}
\pdfglyphtounicode{jehfinalarabic}{FB8B}
\pdfglyphtounicode{jhabengali}{099D}
\pdfglyphtounicode{jhadeva}{091D}
\pdfglyphtounicode{jhagujarati}{0A9D}
\pdfglyphtounicode{jhagurmukhi}{0A1D}
\pdfglyphtounicode{jheharmenian}{057B}
\pdfglyphtounicode{jis}{3004}
\pdfglyphtounicode{jmonospace}{FF4A}
\pdfglyphtounicode{jparen}{24A5}
\pdfglyphtounicode{jsuperior}{02B2}
\pdfglyphtounicode{k}{006B}
\pdfglyphtounicode{kabashkircyrillic}{04A1}
\pdfglyphtounicode{kabengali}{0995}
\pdfglyphtounicode{kacute}{1E31}
\pdfglyphtounicode{kacyrillic}{043A}
\pdfglyphtounicode{kadescendercyrillic}{049B}
\pdfglyphtounicode{kadeva}{0915}
\pdfglyphtounicode{kaf}{05DB}
\pdfglyphtounicode{kafarabic}{0643}
\pdfglyphtounicode{kafdagesh}{FB3B}
\pdfglyphtounicode{kafdageshhebrew}{FB3B}
\pdfglyphtounicode{kaffinalarabic}{FEDA}
\pdfglyphtounicode{kafhebrew}{05DB}
\pdfglyphtounicode{kafinitialarabic}{FEDB}
\pdfglyphtounicode{kafmedialarabic}{FEDC}
\pdfglyphtounicode{kafrafehebrew}{FB4D}
\pdfglyphtounicode{kagujarati}{0A95}
\pdfglyphtounicode{kagurmukhi}{0A15}
\pdfglyphtounicode{kahiragana}{304B}
\pdfglyphtounicode{kahookcyrillic}{04C4}
\pdfglyphtounicode{kakatakana}{30AB}
\pdfglyphtounicode{kakatakanahalfwidth}{FF76}
\pdfglyphtounicode{kappa}{03BA}
\pdfglyphtounicode{kappasymbolgreek}{03F0}
\pdfglyphtounicode{kapyeounmieumkorean}{3171}
\pdfglyphtounicode{kapyeounphieuphkorean}{3184}
\pdfglyphtounicode{kapyeounpieupkorean}{3178}
\pdfglyphtounicode{kapyeounssangpieupkorean}{3179}
\pdfglyphtounicode{karoriisquare}{330D}
\pdfglyphtounicode{kashidaautoarabic}{0640}
\pdfglyphtounicode{kashidaautonosidebearingarabic}{0640}
\pdfglyphtounicode{kasmallkatakana}{30F5}
\pdfglyphtounicode{kasquare}{3384}
\pdfglyphtounicode{kasraarabic}{0650}
\pdfglyphtounicode{kasratanarabic}{064D}
\pdfglyphtounicode{kastrokecyrillic}{049F}
\pdfglyphtounicode{katahiraprolongmarkhalfwidth}{FF70}
\pdfglyphtounicode{kaverticalstrokecyrillic}{049D}
\pdfglyphtounicode{kbopomofo}{310E}
\pdfglyphtounicode{kcalsquare}{3389}
\pdfglyphtounicode{kcaron}{01E9}
\pdfglyphtounicode{kcedilla}{0137}
\pdfglyphtounicode{kcircle}{24DA}
\pdfglyphtounicode{kcommaaccent}{0137}
\pdfglyphtounicode{kdotbelow}{1E33}
\pdfglyphtounicode{keharmenian}{0584}
\pdfglyphtounicode{kehiragana}{3051}
\pdfglyphtounicode{kekatakana}{30B1}
\pdfglyphtounicode{kekatakanahalfwidth}{FF79}
\pdfglyphtounicode{kenarmenian}{056F}
\pdfglyphtounicode{kesmallkatakana}{30F6}
\pdfglyphtounicode{kgreenlandic}{0138}
\pdfglyphtounicode{khabengali}{0996}
\pdfglyphtounicode{khacyrillic}{0445}
\pdfglyphtounicode{khadeva}{0916}
\pdfglyphtounicode{khagujarati}{0A96}
\pdfglyphtounicode{khagurmukhi}{0A16}
\pdfglyphtounicode{khaharabic}{062E}
\pdfglyphtounicode{khahfinalarabic}{FEA6}
\pdfglyphtounicode{khahinitialarabic}{FEA7}
\pdfglyphtounicode{khahmedialarabic}{FEA8}
\pdfglyphtounicode{kheicoptic}{03E7}
\pdfglyphtounicode{khhadeva}{0959}
\pdfglyphtounicode{khhagurmukhi}{0A59}
\pdfglyphtounicode{khieukhacirclekorean}{3278}
\pdfglyphtounicode{khieukhaparenkorean}{3218}
\pdfglyphtounicode{khieukhcirclekorean}{326A}
\pdfglyphtounicode{khieukhkorean}{314B}
\pdfglyphtounicode{khieukhparenkorean}{320A}
\pdfglyphtounicode{khokhaithai}{0E02}
\pdfglyphtounicode{khokhonthai}{0E05}
\pdfglyphtounicode{khokhuatthai}{0E03}
\pdfglyphtounicode{khokhwaithai}{0E04}
\pdfglyphtounicode{khomutthai}{0E5B}
\pdfglyphtounicode{khook}{0199}
\pdfglyphtounicode{khorakhangthai}{0E06}
\pdfglyphtounicode{khzsquare}{3391}
\pdfglyphtounicode{kihiragana}{304D}
\pdfglyphtounicode{kikatakana}{30AD}
\pdfglyphtounicode{kikatakanahalfwidth}{FF77}
\pdfglyphtounicode{kiroguramusquare}{3315}
\pdfglyphtounicode{kiromeetorusquare}{3316}
\pdfglyphtounicode{kirosquare}{3314}
\pdfglyphtounicode{kiyeokacirclekorean}{326E}
\pdfglyphtounicode{kiyeokaparenkorean}{320E}
\pdfglyphtounicode{kiyeokcirclekorean}{3260}
\pdfglyphtounicode{kiyeokkorean}{3131}
\pdfglyphtounicode{kiyeokparenkorean}{3200}
\pdfglyphtounicode{kiyeoksioskorean}{3133}
\pdfglyphtounicode{kjecyrillic}{045C}
\pdfglyphtounicode{klinebelow}{1E35}
\pdfglyphtounicode{klsquare}{3398}
\pdfglyphtounicode{kmcubedsquare}{33A6}
\pdfglyphtounicode{kmonospace}{FF4B}
\pdfglyphtounicode{kmsquaredsquare}{33A2}
\pdfglyphtounicode{kohiragana}{3053}
\pdfglyphtounicode{kohmsquare}{33C0}
\pdfglyphtounicode{kokaithai}{0E01}
\pdfglyphtounicode{kokatakana}{30B3}
\pdfglyphtounicode{kokatakanahalfwidth}{FF7A}
\pdfglyphtounicode{kooposquare}{331E}
\pdfglyphtounicode{koppacyrillic}{0481}
\pdfglyphtounicode{koreanstandardsymbol}{327F}
\pdfglyphtounicode{koroniscmb}{0343}
\pdfglyphtounicode{kparen}{24A6}
\pdfglyphtounicode{kpasquare}{33AA}
\pdfglyphtounicode{ksicyrillic}{046F}
\pdfglyphtounicode{ktsquare}{33CF}
\pdfglyphtounicode{kturned}{029E}
\pdfglyphtounicode{kuhiragana}{304F}
\pdfglyphtounicode{kukatakana}{30AF}
\pdfglyphtounicode{kukatakanahalfwidth}{FF78}
\pdfglyphtounicode{kvsquare}{33B8}
\pdfglyphtounicode{kwsquare}{33BE}
\pdfglyphtounicode{l}{006C}
\pdfglyphtounicode{labengali}{09B2}
\pdfglyphtounicode{lacute}{013A}
\pdfglyphtounicode{ladeva}{0932}
\pdfglyphtounicode{lagujarati}{0AB2}
\pdfglyphtounicode{lagurmukhi}{0A32}
\pdfglyphtounicode{lakkhangyaothai}{0E45}
\pdfglyphtounicode{lamaleffinalarabic}{FEFC}
\pdfglyphtounicode{lamalefhamzaabovefinalarabic}{FEF8}
\pdfglyphtounicode{lamalefhamzaaboveisolatedarabic}{FEF7}
\pdfglyphtounicode{lamalefhamzabelowfinalarabic}{FEFA}
\pdfglyphtounicode{lamalefhamzabelowisolatedarabic}{FEF9}
\pdfglyphtounicode{lamalefisolatedarabic}{FEFB}
\pdfglyphtounicode{lamalefmaddaabovefinalarabic}{FEF6}
\pdfglyphtounicode{lamalefmaddaaboveisolatedarabic}{FEF5}
\pdfglyphtounicode{lamarabic}{0644}
\pdfglyphtounicode{lambda}{03BB}
\pdfglyphtounicode{lambdastroke}{019B}
\pdfglyphtounicode{lamed}{05DC}
\pdfglyphtounicode{lameddagesh}{FB3C}
\pdfglyphtounicode{lameddageshhebrew}{FB3C}
\pdfglyphtounicode{lamedhebrew}{05DC}
\pdfglyphtounicode{lamedholam}{05DC 05B9}
\pdfglyphtounicode{lamedholamdagesh}{05DC 05B9 05BC}
\pdfglyphtounicode{lamedholamdageshhebrew}{05DC 05B9 05BC}
\pdfglyphtounicode{lamedholamhebrew}{05DC 05B9}
\pdfglyphtounicode{lamfinalarabic}{FEDE}
\pdfglyphtounicode{lamhahinitialarabic}{FCCA}
\pdfglyphtounicode{laminitialarabic}{FEDF}
\pdfglyphtounicode{lamjeeminitialarabic}{FCC9}
\pdfglyphtounicode{lamkhahinitialarabic}{FCCB}
\pdfglyphtounicode{lamlamhehisolatedarabic}{FDF2}
\pdfglyphtounicode{lammedialarabic}{FEE0}
\pdfglyphtounicode{lammeemhahinitialarabic}{FD88}
\pdfglyphtounicode{lammeeminitialarabic}{FCCC}
\pdfglyphtounicode{lammeemjeeminitialarabic}{FEDF FEE4 FEA0}
\pdfglyphtounicode{lammeemkhahinitialarabic}{FEDF FEE4 FEA8}
\pdfglyphtounicode{largecircle}{25EF}
\pdfglyphtounicode{latticetop}{22A4}
\pdfglyphtounicode{lbar}{019A}
\pdfglyphtounicode{lbelt}{026C}
\pdfglyphtounicode{lbopomofo}{310C}
\pdfglyphtounicode{lcaron}{013E}
\pdfglyphtounicode{lcedilla}{013C}
\pdfglyphtounicode{lcircle}{24DB}
\pdfglyphtounicode{lcircumflexbelow}{1E3D}
\pdfglyphtounicode{lcommaaccent}{013C}
\pdfglyphtounicode{ldot}{0140}
\pdfglyphtounicode{ldotaccent}{0140}
\pdfglyphtounicode{ldotbelow}{1E37}
\pdfglyphtounicode{ldotbelowmacron}{1E39}
\pdfglyphtounicode{leftangleabovecmb}{031A}
\pdfglyphtounicode{lefttackbelowcmb}{0318}
\pdfglyphtounicode{less}{003C}
\pdfglyphtounicode{lessdbleqlgreater}{2A8B}
\pdfglyphtounicode{lessdblequal}{2266}
\pdfglyphtounicode{lessdot}{22D6}
\pdfglyphtounicode{lessequal}{2264}
\pdfglyphtounicode{lessequalgreater}{22DA}
\pdfglyphtounicode{lessequalorgreater}{22DA}
\pdfglyphtounicode{lessmonospace}{FF1C}
\pdfglyphtounicode{lessmuch}{226A}
\pdfglyphtounicode{lessnotdblequal}{2A89}
\pdfglyphtounicode{lessnotequal}{2A87}
\pdfglyphtounicode{lessorapproxeql}{2A85}
\pdfglyphtounicode{lessorequalslant}{2A7D}
\pdfglyphtounicode{lessorequivalent}{2272}
\pdfglyphtounicode{lessorgreater}{2276}
\pdfglyphtounicode{lessornotdbleql}{2268}
\pdfglyphtounicode{lessornotequal}{2268}
\pdfglyphtounicode{lessorsimilar}{2272}
\pdfglyphtounicode{lessoverequal}{2266}
\pdfglyphtounicode{lesssmall}{FE64}
\pdfglyphtounicode{lezh}{026E}
\pdfglyphtounicode{lfblock}{258C}
\pdfglyphtounicode{lhookretroflex}{026D}
\pdfglyphtounicode{lira}{20A4}
\pdfglyphtounicode{liwnarmenian}{056C}
\pdfglyphtounicode{lj}{01C9}
\pdfglyphtounicode{ljecyrillic}{0459}
\pdfglyphtounicode{ll}{006C 006C}
\pdfglyphtounicode{lladeva}{0933}
\pdfglyphtounicode{llagujarati}{0AB3}
\pdfglyphtounicode{llinebelow}{1E3B}
\pdfglyphtounicode{llladeva}{0934}
\pdfglyphtounicode{llvocalicbengali}{09E1}
\pdfglyphtounicode{llvocalicdeva}{0961}
\pdfglyphtounicode{llvocalicvowelsignbengali}{09E3}
\pdfglyphtounicode{llvocalicvowelsigndeva}{0963}
\pdfglyphtounicode{lmiddletilde}{026B}
\pdfglyphtounicode{lmonospace}{FF4C}
\pdfglyphtounicode{lmsquare}{33D0}
\pdfglyphtounicode{lochulathai}{0E2C}
\pdfglyphtounicode{logicaland}{2227}
\pdfglyphtounicode{logicalnot}{00AC}
\pdfglyphtounicode{logicalnotreversed}{2310}
\pdfglyphtounicode{logicalor}{2228}
\pdfglyphtounicode{lolingthai}{0E25}
\pdfglyphtounicode{longdbls}{017F 017F}
\pdfglyphtounicode{longs}{017F}
\pdfglyphtounicode{longsh}{017F 0068}
\pdfglyphtounicode{longsi}{017F 0069}
\pdfglyphtounicode{longsl}{017F 006C}
\pdfglyphtounicode{longst}{017F 0074}
\pdfglyphtounicode{lowlinecenterline}{FE4E}
\pdfglyphtounicode{lowlinecmb}{0332}
\pdfglyphtounicode{lowlinedashed}{FE4D}
\pdfglyphtounicode{lozenge}{25CA}
\pdfglyphtounicode{lparen}{24A7}
\pdfglyphtounicode{lscript}{2113}
\pdfglyphtounicode{lslash}{0142}
\pdfglyphtounicode{lsquare}{2113}
\pdfglyphtounicode{lsuperior}{006C}
\pdfglyphtounicode{ltshade}{2591}
\pdfglyphtounicode{luthai}{0E26}
\pdfglyphtounicode{lvocalicbengali}{098C}
\pdfglyphtounicode{lvocalicdeva}{090C}
\pdfglyphtounicode{lvocalicvowelsignbengali}{09E2}
\pdfglyphtounicode{lvocalicvowelsigndeva}{0962}
\pdfglyphtounicode{lxsquare}{33D3}
\pdfglyphtounicode{m}{006D}
\pdfglyphtounicode{mabengali}{09AE}
\pdfglyphtounicode{macron}{00AF}
\pdfglyphtounicode{macronbelowcmb}{0331}
\pdfglyphtounicode{macroncmb}{0304}
\pdfglyphtounicode{macronlowmod}{02CD}
\pdfglyphtounicode{macronmonospace}{FFE3}
\pdfglyphtounicode{macute}{1E3F}
\pdfglyphtounicode{madeva}{092E}
\pdfglyphtounicode{magujarati}{0AAE}
\pdfglyphtounicode{magurmukhi}{0A2E}
\pdfglyphtounicode{mahapakhhebrew}{05A4}
\pdfglyphtounicode{mahapakhlefthebrew}{05A4}
\pdfglyphtounicode{mahiragana}{307E}
\pdfglyphtounicode{maichattawalowleftthai}{F895}
\pdfglyphtounicode{maichattawalowrightthai}{F894}
\pdfglyphtounicode{maichattawathai}{0E4B}
\pdfglyphtounicode{maichattawaupperleftthai}{F893}
\pdfglyphtounicode{maieklowleftthai}{F88C}
\pdfglyphtounicode{maieklowrightthai}{F88B}
\pdfglyphtounicode{maiekthai}{0E48}
\pdfglyphtounicode{maiekupperleftthai}{F88A}
\pdfglyphtounicode{maihanakatleftthai}{F884}
\pdfglyphtounicode{maihanakatthai}{0E31}
\pdfglyphtounicode{maitaikhuleftthai}{F889}
\pdfglyphtounicode{maitaikhuthai}{0E47}
\pdfglyphtounicode{maitholowleftthai}{F88F}
\pdfglyphtounicode{maitholowrightthai}{F88E}
\pdfglyphtounicode{maithothai}{0E49}
\pdfglyphtounicode{maithoupperleftthai}{F88D}
\pdfglyphtounicode{maitrilowleftthai}{F892}
\pdfglyphtounicode{maitrilowrightthai}{F891}
\pdfglyphtounicode{maitrithai}{0E4A}
\pdfglyphtounicode{maitriupperleftthai}{F890}
\pdfglyphtounicode{maiyamokthai}{0E46}
\pdfglyphtounicode{makatakana}{30DE}
\pdfglyphtounicode{makatakanahalfwidth}{FF8F}
\pdfglyphtounicode{male}{2642}
\pdfglyphtounicode{maltesecross}{2720}
\pdfglyphtounicode{mansyonsquare}{3347}
\pdfglyphtounicode{maqafhebrew}{05BE}
\pdfglyphtounicode{mars}{2642}
\pdfglyphtounicode{masoracirclehebrew}{05AF}
\pdfglyphtounicode{masquare}{3383}
\pdfglyphtounicode{mbopomofo}{3107}
\pdfglyphtounicode{mbsquare}{33D4}
\pdfglyphtounicode{mcircle}{24DC}
\pdfglyphtounicode{mcubedsquare}{33A5}
\pdfglyphtounicode{mdotaccent}{1E41}
\pdfglyphtounicode{mdotbelow}{1E43}
\pdfglyphtounicode{measuredangle}{2221}
\pdfglyphtounicode{meemarabic}{0645}
\pdfglyphtounicode{meemfinalarabic}{FEE2}
\pdfglyphtounicode{meeminitialarabic}{FEE3}
\pdfglyphtounicode{meemmedialarabic}{FEE4}
\pdfglyphtounicode{meemmeeminitialarabic}{FCD1}
\pdfglyphtounicode{meemmeemisolatedarabic}{FC48}
\pdfglyphtounicode{meetorusquare}{334D}
\pdfglyphtounicode{mehiragana}{3081}
\pdfglyphtounicode{meizierasquare}{337E}
\pdfglyphtounicode{mekatakana}{30E1}
\pdfglyphtounicode{mekatakanahalfwidth}{FF92}
\pdfglyphtounicode{mem}{05DE}
\pdfglyphtounicode{memdagesh}{FB3E}
\pdfglyphtounicode{memdageshhebrew}{FB3E}
\pdfglyphtounicode{memhebrew}{05DE}
\pdfglyphtounicode{menarmenian}{0574}
\pdfglyphtounicode{merkhahebrew}{05A5}
\pdfglyphtounicode{merkhakefulahebrew}{05A6}
\pdfglyphtounicode{merkhakefulalefthebrew}{05A6}
\pdfglyphtounicode{merkhalefthebrew}{05A5}
\pdfglyphtounicode{mhook}{0271}
\pdfglyphtounicode{mhzsquare}{3392}
\pdfglyphtounicode{middledotkatakanahalfwidth}{FF65}
\pdfglyphtounicode{middot}{00B7}
\pdfglyphtounicode{mieumacirclekorean}{3272}
\pdfglyphtounicode{mieumaparenkorean}{3212}
\pdfglyphtounicode{mieumcirclekorean}{3264}
\pdfglyphtounicode{mieumkorean}{3141}
\pdfglyphtounicode{mieumpansioskorean}{3170}
\pdfglyphtounicode{mieumparenkorean}{3204}
\pdfglyphtounicode{mieumpieupkorean}{316E}
\pdfglyphtounicode{mieumsioskorean}{316F}
\pdfglyphtounicode{mihiragana}{307F}
\pdfglyphtounicode{mikatakana}{30DF}
\pdfglyphtounicode{mikatakanahalfwidth}{FF90}
\pdfglyphtounicode{minus}{2212}
\pdfglyphtounicode{minusbelowcmb}{0320}
\pdfglyphtounicode{minuscircle}{2296}
\pdfglyphtounicode{minusmod}{02D7}
\pdfglyphtounicode{minusplus}{2213}
\pdfglyphtounicode{minute}{2032}
\pdfglyphtounicode{miribaarusquare}{334A}
\pdfglyphtounicode{mirisquare}{3349}
\pdfglyphtounicode{mlonglegturned}{0270}
\pdfglyphtounicode{mlsquare}{3396}
\pdfglyphtounicode{mmcubedsquare}{33A3}
\pdfglyphtounicode{mmonospace}{FF4D}
\pdfglyphtounicode{mmsquaredsquare}{339F}
\pdfglyphtounicode{mohiragana}{3082}
\pdfglyphtounicode{mohmsquare}{33C1}
\pdfglyphtounicode{mokatakana}{30E2}
\pdfglyphtounicode{mokatakanahalfwidth}{FF93}
\pdfglyphtounicode{molsquare}{33D6}
\pdfglyphtounicode{momathai}{0E21}
\pdfglyphtounicode{moverssquare}{33A7}
\pdfglyphtounicode{moverssquaredsquare}{33A8}
\pdfglyphtounicode{mparen}{24A8}
\pdfglyphtounicode{mpasquare}{33AB}
\pdfglyphtounicode{mssquare}{33B3}
\pdfglyphtounicode{msuperior}{006D}
\pdfglyphtounicode{mturned}{026F}
\pdfglyphtounicode{mu}{00B5}
\pdfglyphtounicode{mu1}{00B5}
\pdfglyphtounicode{muasquare}{3382}
\pdfglyphtounicode{muchgreater}{226B}
\pdfglyphtounicode{muchless}{226A}
\pdfglyphtounicode{mufsquare}{338C}
\pdfglyphtounicode{mugreek}{03BC}
\pdfglyphtounicode{mugsquare}{338D}
\pdfglyphtounicode{muhiragana}{3080}
\pdfglyphtounicode{mukatakana}{30E0}
\pdfglyphtounicode{mukatakanahalfwidth}{FF91}
\pdfglyphtounicode{mulsquare}{3395}
\pdfglyphtounicode{multicloseleft}{22C9}
\pdfglyphtounicode{multicloseright}{22CA}
\pdfglyphtounicode{multimap}{22B8}
\pdfglyphtounicode{multiopenleft}{22CB}
\pdfglyphtounicode{multiopenright}{22CC}
\pdfglyphtounicode{multiply}{00D7}
\pdfglyphtounicode{mumsquare}{339B}
\pdfglyphtounicode{munahhebrew}{05A3}
\pdfglyphtounicode{munahlefthebrew}{05A3}
\pdfglyphtounicode{musicalnote}{266A}
\pdfglyphtounicode{musicalnotedbl}{266B}
\pdfglyphtounicode{musicflatsign}{266D}
\pdfglyphtounicode{musicsharpsign}{266F}
\pdfglyphtounicode{mussquare}{33B2}
\pdfglyphtounicode{muvsquare}{33B6}
\pdfglyphtounicode{muwsquare}{33BC}
\pdfglyphtounicode{mvmegasquare}{33B9}
\pdfglyphtounicode{mvsquare}{33B7}
\pdfglyphtounicode{mwmegasquare}{33BF}
\pdfglyphtounicode{mwsquare}{33BD}
\pdfglyphtounicode{n}{006E}
\pdfglyphtounicode{nabengali}{09A8}
\pdfglyphtounicode{nabla}{2207}
\pdfglyphtounicode{nacute}{0144}
\pdfglyphtounicode{nadeva}{0928}
\pdfglyphtounicode{nagujarati}{0AA8}
\pdfglyphtounicode{nagurmukhi}{0A28}
\pdfglyphtounicode{nahiragana}{306A}
\pdfglyphtounicode{nakatakana}{30CA}
\pdfglyphtounicode{nakatakanahalfwidth}{FF85}
\pdfglyphtounicode{nand}{22BC}
\pdfglyphtounicode{napostrophe}{0149}
\pdfglyphtounicode{nasquare}{3381}
\pdfglyphtounicode{natural}{266E}
\pdfglyphtounicode{nbopomofo}{310B}
\pdfglyphtounicode{nbspace}{00A0}
\pdfglyphtounicode{ncaron}{0148}
\pdfglyphtounicode{ncedilla}{0146}
\pdfglyphtounicode{ncircle}{24DD}
\pdfglyphtounicode{ncircumflexbelow}{1E4B}
\pdfglyphtounicode{ncommaaccent}{0146}
\pdfglyphtounicode{ndotaccent}{1E45}
\pdfglyphtounicode{ndotbelow}{1E47}
\pdfglyphtounicode{negationslash}{0338}
\pdfglyphtounicode{nehiragana}{306D}
\pdfglyphtounicode{nekatakana}{30CD}
\pdfglyphtounicode{nekatakanahalfwidth}{FF88}
\pdfglyphtounicode{newsheqelsign}{20AA}
\pdfglyphtounicode{nfsquare}{338B}
\pdfglyphtounicode{ng}{014B}
\pdfglyphtounicode{ngabengali}{0999}
\pdfglyphtounicode{ngadeva}{0919}
\pdfglyphtounicode{ngagujarati}{0A99}
\pdfglyphtounicode{ngagurmukhi}{0A19}
\pdfglyphtounicode{ngonguthai}{0E07}
\pdfglyphtounicode{nhiragana}{3093}
\pdfglyphtounicode{nhookleft}{0272}
\pdfglyphtounicode{nhookretroflex}{0273}
\pdfglyphtounicode{nieunacirclekorean}{326F}
\pdfglyphtounicode{nieunaparenkorean}{320F}
\pdfglyphtounicode{nieuncieuckorean}{3135}
\pdfglyphtounicode{nieuncirclekorean}{3261}
\pdfglyphtounicode{nieunhieuhkorean}{3136}
\pdfglyphtounicode{nieunkorean}{3134}
\pdfglyphtounicode{nieunpansioskorean}{3168}
\pdfglyphtounicode{nieunparenkorean}{3201}
\pdfglyphtounicode{nieunsioskorean}{3167}
\pdfglyphtounicode{nieuntikeutkorean}{3166}
\pdfglyphtounicode{nihiragana}{306B}
\pdfglyphtounicode{nikatakana}{30CB}
\pdfglyphtounicode{nikatakanahalfwidth}{FF86}
\pdfglyphtounicode{nikhahitleftthai}{F899}
\pdfglyphtounicode{nikhahitthai}{0E4D}
\pdfglyphtounicode{nine}{0039}
\pdfglyphtounicode{ninearabic}{0669}
\pdfglyphtounicode{ninebengali}{09EF}
\pdfglyphtounicode{ninecircle}{2468}
\pdfglyphtounicode{ninecircleinversesansserif}{2792}
\pdfglyphtounicode{ninedeva}{096F}
\pdfglyphtounicode{ninegujarati}{0AEF}
\pdfglyphtounicode{ninegurmukhi}{0A6F}
\pdfglyphtounicode{ninehackarabic}{0669}
\pdfglyphtounicode{ninehangzhou}{3029}
\pdfglyphtounicode{nineideographicparen}{3228}
\pdfglyphtounicode{nineinferior}{2089}
\pdfglyphtounicode{ninemonospace}{FF19}
\pdfglyphtounicode{nineoldstyle}{0039}
\pdfglyphtounicode{nineparen}{247C}
\pdfglyphtounicode{nineperiod}{2490}
\pdfglyphtounicode{ninepersian}{06F9}
\pdfglyphtounicode{nineroman}{2178}
\pdfglyphtounicode{ninesuperior}{2079}
\pdfglyphtounicode{nineteencircle}{2472}
\pdfglyphtounicode{nineteenparen}{2486}
\pdfglyphtounicode{nineteenperiod}{249A}
\pdfglyphtounicode{ninethai}{0E59}
\pdfglyphtounicode{nj}{01CC}
\pdfglyphtounicode{njecyrillic}{045A}
\pdfglyphtounicode{nkatakana}{30F3}
\pdfglyphtounicode{nkatakanahalfwidth}{FF9D}
\pdfglyphtounicode{nlegrightlong}{019E}
\pdfglyphtounicode{nlinebelow}{1E49}
\pdfglyphtounicode{nmonospace}{FF4E}
\pdfglyphtounicode{nmsquare}{339A}
\pdfglyphtounicode{nnabengali}{09A3}
\pdfglyphtounicode{nnadeva}{0923}
\pdfglyphtounicode{nnagujarati}{0AA3}
\pdfglyphtounicode{nnagurmukhi}{0A23}
\pdfglyphtounicode{nnnadeva}{0929}
\pdfglyphtounicode{nohiragana}{306E}
\pdfglyphtounicode{nokatakana}{30CE}
\pdfglyphtounicode{nokatakanahalfwidth}{FF89}
\pdfglyphtounicode{nonbreakingspace}{00A0}
\pdfglyphtounicode{nonenthai}{0E13}
\pdfglyphtounicode{nonuthai}{0E19}
\pdfglyphtounicode{noonarabic}{0646}
\pdfglyphtounicode{noonfinalarabic}{FEE6}
\pdfglyphtounicode{noonghunnaarabic}{06BA}
\pdfglyphtounicode{noonghunnafinalarabic}{FB9F}
\pdfglyphtounicode{noonhehinitialarabic}{FEE7 FEEC}
\pdfglyphtounicode{nooninitialarabic}{FEE7}
\pdfglyphtounicode{noonjeeminitialarabic}{FCD2}
\pdfglyphtounicode{noonjeemisolatedarabic}{FC4B}
\pdfglyphtounicode{noonmedialarabic}{FEE8}
\pdfglyphtounicode{noonmeeminitialarabic}{FCD5}
\pdfglyphtounicode{noonmeemisolatedarabic}{FC4E}
\pdfglyphtounicode{noonnoonfinalarabic}{FC8D}
\pdfglyphtounicode{notapproxequal}{2247}
\pdfglyphtounicode{notarrowboth}{21AE}
\pdfglyphtounicode{notarrowleft}{219A}
\pdfglyphtounicode{notarrowright}{219B}
\pdfglyphtounicode{notbar}{2224}
\pdfglyphtounicode{notcontains}{220C}
\pdfglyphtounicode{notdblarrowboth}{21CE}
\pdfglyphtounicode{notdblarrowleft}{21CD}
\pdfglyphtounicode{notdblarrowright}{21CF}
\pdfglyphtounicode{notelement}{2209}
\pdfglyphtounicode{notelementof}{2209}
\pdfglyphtounicode{notequal}{2260}
\pdfglyphtounicode{notexistential}{2204}
\pdfglyphtounicode{notfollows}{2281}
\pdfglyphtounicode{notfollowsoreql}{2AB0 0338}
\pdfglyphtounicode{notforces}{22AE}
\pdfglyphtounicode{notforcesextra}{22AF}
\pdfglyphtounicode{notgreater}{226F}
\pdfglyphtounicode{notgreaterdblequal}{2267 0338}
\pdfglyphtounicode{notgreaterequal}{2271}
\pdfglyphtounicode{notgreaternorequal}{2271}
\pdfglyphtounicode{notgreaternorless}{2279}
\pdfglyphtounicode{notgreaterorslnteql}{2A7E 0338}
\pdfglyphtounicode{notidentical}{2262}
\pdfglyphtounicode{notless}{226E}
\pdfglyphtounicode{notlessdblequal}{2266 0338}
\pdfglyphtounicode{notlessequal}{2270}
\pdfglyphtounicode{notlessnorequal}{2270}
\pdfglyphtounicode{notlessorslnteql}{2A7D 0338}
\pdfglyphtounicode{notparallel}{2226}
\pdfglyphtounicode{notprecedes}{2280}
\pdfglyphtounicode{notprecedesoreql}{2AAF 0338}
\pdfglyphtounicode{notsatisfies}{22AD}
\pdfglyphtounicode{notsimilar}{2241}
\pdfglyphtounicode{notsubset}{2284}
\pdfglyphtounicode{notsubseteql}{2288}
\pdfglyphtounicode{notsubsetordbleql}{2AC5 0338}
\pdfglyphtounicode{notsubsetoreql}{228A}
\pdfglyphtounicode{notsucceeds}{2281}
\pdfglyphtounicode{notsuperset}{2285}
\pdfglyphtounicode{notsuperseteql}{2289}
\pdfglyphtounicode{notsupersetordbleql}{2AC6 0338}
\pdfglyphtounicode{notsupersetoreql}{228B}
\pdfglyphtounicode{nottriangeqlleft}{22EC}
\pdfglyphtounicode{nottriangeqlright}{22ED}
\pdfglyphtounicode{nottriangleleft}{22EA}
\pdfglyphtounicode{nottriangleright}{22EB}
\pdfglyphtounicode{notturnstile}{22AC}
\pdfglyphtounicode{nowarmenian}{0576}
\pdfglyphtounicode{nparen}{24A9}
\pdfglyphtounicode{nssquare}{33B1}
\pdfglyphtounicode{nsuperior}{207F}
\pdfglyphtounicode{ntilde}{00F1}
\pdfglyphtounicode{nu}{03BD}
\pdfglyphtounicode{nuhiragana}{306C}
\pdfglyphtounicode{nukatakana}{30CC}
\pdfglyphtounicode{nukatakanahalfwidth}{FF87}
\pdfglyphtounicode{nuktabengali}{09BC}
\pdfglyphtounicode{nuktadeva}{093C}
\pdfglyphtounicode{nuktagujarati}{0ABC}
\pdfglyphtounicode{nuktagurmukhi}{0A3C}
\pdfglyphtounicode{numbersign}{0023}
\pdfglyphtounicode{numbersignmonospace}{FF03}
\pdfglyphtounicode{numbersignsmall}{FE5F}
\pdfglyphtounicode{numeralsigngreek}{0374}
\pdfglyphtounicode{numeralsignlowergreek}{0375}
\pdfglyphtounicode{numero}{2116}
\pdfglyphtounicode{nun}{05E0}
\pdfglyphtounicode{nundagesh}{FB40}
\pdfglyphtounicode{nundageshhebrew}{FB40}
\pdfglyphtounicode{nunhebrew}{05E0}
\pdfglyphtounicode{nvsquare}{33B5}
\pdfglyphtounicode{nwsquare}{33BB}
\pdfglyphtounicode{nyabengali}{099E}
\pdfglyphtounicode{nyadeva}{091E}
\pdfglyphtounicode{nyagujarati}{0A9E}
\pdfglyphtounicode{nyagurmukhi}{0A1E}
\pdfglyphtounicode{o}{006F}
\pdfglyphtounicode{oacute}{00F3}
\pdfglyphtounicode{oangthai}{0E2D}
\pdfglyphtounicode{obarred}{0275}
\pdfglyphtounicode{obarredcyrillic}{04E9}
\pdfglyphtounicode{obarreddieresiscyrillic}{04EB}
\pdfglyphtounicode{obengali}{0993}
\pdfglyphtounicode{obopomofo}{311B}
\pdfglyphtounicode{obreve}{014F}
\pdfglyphtounicode{ocandradeva}{0911}
\pdfglyphtounicode{ocandragujarati}{0A91}
\pdfglyphtounicode{ocandravowelsigndeva}{0949}
\pdfglyphtounicode{ocandravowelsigngujarati}{0AC9}
\pdfglyphtounicode{ocaron}{01D2}
\pdfglyphtounicode{ocircle}{24DE}
\pdfglyphtounicode{ocircumflex}{00F4}
\pdfglyphtounicode{ocircumflexacute}{1ED1}
\pdfglyphtounicode{ocircumflexdotbelow}{1ED9}
\pdfglyphtounicode{ocircumflexgrave}{1ED3}
\pdfglyphtounicode{ocircumflexhookabove}{1ED5}
\pdfglyphtounicode{ocircumflextilde}{1ED7}
\pdfglyphtounicode{ocyrillic}{043E}
\pdfglyphtounicode{odblacute}{0151}
\pdfglyphtounicode{odblgrave}{020D}
\pdfglyphtounicode{odeva}{0913}
\pdfglyphtounicode{odieresis}{00F6}
\pdfglyphtounicode{odieresiscyrillic}{04E7}
\pdfglyphtounicode{odotbelow}{1ECD}
\pdfglyphtounicode{oe}{0153}
\pdfglyphtounicode{oekorean}{315A}
\pdfglyphtounicode{ogonek}{02DB}
\pdfglyphtounicode{ogonekcmb}{0328}
\pdfglyphtounicode{ograve}{00F2}
\pdfglyphtounicode{ogujarati}{0A93}
\pdfglyphtounicode{oharmenian}{0585}
\pdfglyphtounicode{ohiragana}{304A}
\pdfglyphtounicode{ohookabove}{1ECF}
\pdfglyphtounicode{ohorn}{01A1}
\pdfglyphtounicode{ohornacute}{1EDB}
\pdfglyphtounicode{ohorndotbelow}{1EE3}
\pdfglyphtounicode{ohorngrave}{1EDD}
\pdfglyphtounicode{ohornhookabove}{1EDF}
\pdfglyphtounicode{ohorntilde}{1EE1}
\pdfglyphtounicode{ohungarumlaut}{0151}
\pdfglyphtounicode{oi}{01A3}
\pdfglyphtounicode{oinvertedbreve}{020F}
\pdfglyphtounicode{okatakana}{30AA}
\pdfglyphtounicode{okatakanahalfwidth}{FF75}
\pdfglyphtounicode{okorean}{3157}
\pdfglyphtounicode{olehebrew}{05AB}
\pdfglyphtounicode{omacron}{014D}
\pdfglyphtounicode{omacronacute}{1E53}
\pdfglyphtounicode{omacrongrave}{1E51}
\pdfglyphtounicode{omdeva}{0950}
\pdfglyphtounicode{omega}{03C9}
\pdfglyphtounicode{omega1}{03D6}
\pdfglyphtounicode{omegacyrillic}{0461}
\pdfglyphtounicode{omegalatinclosed}{0277}
\pdfglyphtounicode{omegaroundcyrillic}{047B}
\pdfglyphtounicode{omegatitlocyrillic}{047D}
\pdfglyphtounicode{omegatonos}{03CE}
\pdfglyphtounicode{omgujarati}{0AD0}
\pdfglyphtounicode{omicron}{03BF}
\pdfglyphtounicode{omicrontonos}{03CC}
\pdfglyphtounicode{omonospace}{FF4F}
\pdfglyphtounicode{one}{0031}
\pdfglyphtounicode{onearabic}{0661}
\pdfglyphtounicode{onebengali}{09E7}
\pdfglyphtounicode{onecircle}{2460}
\pdfglyphtounicode{onecircleinversesansserif}{278A}
\pdfglyphtounicode{onedeva}{0967}
\pdfglyphtounicode{onedotenleader}{2024}
\pdfglyphtounicode{oneeighth}{215B}
\pdfglyphtounicode{onefitted}{0031}
\pdfglyphtounicode{onegujarati}{0AE7}
\pdfglyphtounicode{onegurmukhi}{0A67}
\pdfglyphtounicode{onehackarabic}{0661}
\pdfglyphtounicode{onehalf}{00BD}
\pdfglyphtounicode{onehangzhou}{3021}
\pdfglyphtounicode{oneideographicparen}{3220}
\pdfglyphtounicode{oneinferior}{2081}
\pdfglyphtounicode{onemonospace}{FF11}
\pdfglyphtounicode{onenumeratorbengali}{09F4}
\pdfglyphtounicode{oneoldstyle}{0031}
\pdfglyphtounicode{oneparen}{2474}
\pdfglyphtounicode{oneperiod}{2488}
\pdfglyphtounicode{onepersian}{06F1}
\pdfglyphtounicode{onequarter}{00BC}
\pdfglyphtounicode{oneroman}{2170}
\pdfglyphtounicode{onesuperior}{00B9}
\pdfglyphtounicode{onethai}{0E51}
\pdfglyphtounicode{onethird}{2153}
\pdfglyphtounicode{oogonek}{01EB}
\pdfglyphtounicode{oogonekmacron}{01ED}
\pdfglyphtounicode{oogurmukhi}{0A13}
\pdfglyphtounicode{oomatragurmukhi}{0A4B}
\pdfglyphtounicode{oopen}{0254}
\pdfglyphtounicode{oparen}{24AA}
\pdfglyphtounicode{openbullet}{25E6}
\pdfglyphtounicode{option}{2325}
\pdfglyphtounicode{ordfeminine}{00AA}
\pdfglyphtounicode{ordmasculine}{00BA}
\pdfglyphtounicode{orthogonal}{221F}
\pdfglyphtounicode{orunderscore}{22BB}
\pdfglyphtounicode{oshortdeva}{0912}
\pdfglyphtounicode{oshortvowelsigndeva}{094A}
\pdfglyphtounicode{oslash}{00F8}
\pdfglyphtounicode{oslashacute}{01FF}
\pdfglyphtounicode{osmallhiragana}{3049}
\pdfglyphtounicode{osmallkatakana}{30A9}
\pdfglyphtounicode{osmallkatakanahalfwidth}{FF6B}
\pdfglyphtounicode{ostrokeacute}{01FF}
\pdfglyphtounicode{osuperior}{006F}
\pdfglyphtounicode{otcyrillic}{047F}
\pdfglyphtounicode{otilde}{00F5}
\pdfglyphtounicode{otildeacute}{1E4D}
\pdfglyphtounicode{otildedieresis}{1E4F}
\pdfglyphtounicode{oubopomofo}{3121}
\pdfglyphtounicode{overline}{203E}
\pdfglyphtounicode{overlinecenterline}{FE4A}
\pdfglyphtounicode{overlinecmb}{0305}
\pdfglyphtounicode{overlinedashed}{FE49}
\pdfglyphtounicode{overlinedblwavy}{FE4C}
\pdfglyphtounicode{overlinewavy}{FE4B}
\pdfglyphtounicode{overscore}{00AF}
\pdfglyphtounicode{ovowelsignbengali}{09CB}
\pdfglyphtounicode{ovowelsigndeva}{094B}
\pdfglyphtounicode{ovowelsigngujarati}{0ACB}
\pdfglyphtounicode{owner}{220B}
\pdfglyphtounicode{p}{0070}
\pdfglyphtounicode{paampssquare}{3380}
\pdfglyphtounicode{paasentosquare}{332B}
\pdfglyphtounicode{pabengali}{09AA}
\pdfglyphtounicode{pacute}{1E55}
\pdfglyphtounicode{padeva}{092A}
\pdfglyphtounicode{pagedown}{21DF}
\pdfglyphtounicode{pageup}{21DE}
\pdfglyphtounicode{pagujarati}{0AAA}
\pdfglyphtounicode{pagurmukhi}{0A2A}
\pdfglyphtounicode{pahiragana}{3071}
\pdfglyphtounicode{paiyannoithai}{0E2F}
\pdfglyphtounicode{pakatakana}{30D1}
\pdfglyphtounicode{palatalizationcyrilliccmb}{0484}
\pdfglyphtounicode{palochkacyrillic}{04C0}
\pdfglyphtounicode{pansioskorean}{317F}
\pdfglyphtounicode{paragraph}{00B6}
\pdfglyphtounicode{parallel}{2225}
\pdfglyphtounicode{parenleft}{0028}
\pdfglyphtounicode{parenleftaltonearabic}{FD3E}
\pdfglyphtounicode{parenleftbt}{F8ED}
\pdfglyphtounicode{parenleftex}{F8EC}
\pdfglyphtounicode{parenleftinferior}{208D}
\pdfglyphtounicode{parenleftmonospace}{FF08}
\pdfglyphtounicode{parenleftsmall}{FE59}
\pdfglyphtounicode{parenleftsuperior}{207D}
\pdfglyphtounicode{parenlefttp}{F8EB}
\pdfglyphtounicode{parenleftvertical}{FE35}
\pdfglyphtounicode{parenright}{0029}
\pdfglyphtounicode{parenrightaltonearabic}{FD3F}
\pdfglyphtounicode{parenrightbt}{F8F8}
\pdfglyphtounicode{parenrightex}{F8F7}
\pdfglyphtounicode{parenrightinferior}{208E}
\pdfglyphtounicode{parenrightmonospace}{FF09}
\pdfglyphtounicode{parenrightsmall}{FE5A}
\pdfglyphtounicode{parenrightsuperior}{207E}
\pdfglyphtounicode{parenrighttp}{F8F6}
\pdfglyphtounicode{parenrightvertical}{FE36}
\pdfglyphtounicode{partialdiff}{2202}
\pdfglyphtounicode{paseqhebrew}{05C0}
\pdfglyphtounicode{pashtahebrew}{0599}
\pdfglyphtounicode{pasquare}{33A9}
\pdfglyphtounicode{patah}{05B7}
\pdfglyphtounicode{patah11}{05B7}
\pdfglyphtounicode{patah1d}{05B7}
\pdfglyphtounicode{patah2a}{05B7}
\pdfglyphtounicode{patahhebrew}{05B7}
\pdfglyphtounicode{patahnarrowhebrew}{05B7}
\pdfglyphtounicode{patahquarterhebrew}{05B7}
\pdfglyphtounicode{patahwidehebrew}{05B7}
\pdfglyphtounicode{pazerhebrew}{05A1}
\pdfglyphtounicode{pbopomofo}{3106}
\pdfglyphtounicode{pcircle}{24DF}
\pdfglyphtounicode{pdotaccent}{1E57}
\pdfglyphtounicode{pe}{05E4}
\pdfglyphtounicode{pecyrillic}{043F}
\pdfglyphtounicode{pedagesh}{FB44}
\pdfglyphtounicode{pedageshhebrew}{FB44}
\pdfglyphtounicode{peezisquare}{333B}
\pdfglyphtounicode{pefinaldageshhebrew}{FB43}
\pdfglyphtounicode{peharabic}{067E}
\pdfglyphtounicode{peharmenian}{057A}
\pdfglyphtounicode{pehebrew}{05E4}
\pdfglyphtounicode{pehfinalarabic}{FB57}
\pdfglyphtounicode{pehinitialarabic}{FB58}
\pdfglyphtounicode{pehiragana}{307A}
\pdfglyphtounicode{pehmedialarabic}{FB59}
\pdfglyphtounicode{pekatakana}{30DA}
\pdfglyphtounicode{pemiddlehookcyrillic}{04A7}
\pdfglyphtounicode{perafehebrew}{FB4E}
\pdfglyphtounicode{percent}{0025}
\pdfglyphtounicode{percentarabic}{066A}
\pdfglyphtounicode{percentmonospace}{FF05}
\pdfglyphtounicode{percentsmall}{FE6A}
\pdfglyphtounicode{period}{002E}
\pdfglyphtounicode{periodarmenian}{0589}
\pdfglyphtounicode{periodcentered}{00B7}
\pdfglyphtounicode{periodhalfwidth}{FF61}
\pdfglyphtounicode{periodinferior}{002E}
\pdfglyphtounicode{periodmonospace}{FF0E}
\pdfglyphtounicode{periodsmall}{FE52}
\pdfglyphtounicode{periodsuperior}{002E}
\pdfglyphtounicode{perispomenigreekcmb}{0342}
\pdfglyphtounicode{perpcorrespond}{2A5E}
\pdfglyphtounicode{perpendicular}{22A5}
\pdfglyphtounicode{pertenthousand}{2031}
\pdfglyphtounicode{perthousand}{2030}
\pdfglyphtounicode{peseta}{20A7}
\pdfglyphtounicode{pfsquare}{338A}
\pdfglyphtounicode{phabengali}{09AB}
\pdfglyphtounicode{phadeva}{092B}
\pdfglyphtounicode{phagujarati}{0AAB}
\pdfglyphtounicode{phagurmukhi}{0A2B}
\pdfglyphtounicode{phi}{03C6}
\pdfglyphtounicode{phi1}{03D5}
\pdfglyphtounicode{phieuphacirclekorean}{327A}
\pdfglyphtounicode{phieuphaparenkorean}{321A}
\pdfglyphtounicode{phieuphcirclekorean}{326C}
\pdfglyphtounicode{phieuphkorean}{314D}
\pdfglyphtounicode{phieuphparenkorean}{320C}
\pdfglyphtounicode{philatin}{0278}
\pdfglyphtounicode{phinthuthai}{0E3A}
\pdfglyphtounicode{phisymbolgreek}{03D5}
\pdfglyphtounicode{phook}{01A5}
\pdfglyphtounicode{phophanthai}{0E1E}
\pdfglyphtounicode{phophungthai}{0E1C}
\pdfglyphtounicode{phosamphaothai}{0E20}
\pdfglyphtounicode{pi}{03C0}
\pdfglyphtounicode{pi1}{03D6}
\pdfglyphtounicode{pieupacirclekorean}{3273}
\pdfglyphtounicode{pieupaparenkorean}{3213}
\pdfglyphtounicode{pieupcieuckorean}{3176}
\pdfglyphtounicode{pieupcirclekorean}{3265}
\pdfglyphtounicode{pieupkiyeokkorean}{3172}
\pdfglyphtounicode{pieupkorean}{3142}
\pdfglyphtounicode{pieupparenkorean}{3205}
\pdfglyphtounicode{pieupsioskiyeokkorean}{3174}
\pdfglyphtounicode{pieupsioskorean}{3144}
\pdfglyphtounicode{pieupsiostikeutkorean}{3175}
\pdfglyphtounicode{pieupthieuthkorean}{3177}
\pdfglyphtounicode{pieuptikeutkorean}{3173}
\pdfglyphtounicode{pihiragana}{3074}
\pdfglyphtounicode{pikatakana}{30D4}
\pdfglyphtounicode{pisymbolgreek}{03D6}
\pdfglyphtounicode{piwrarmenian}{0583}
\pdfglyphtounicode{planckover2pi}{210F}
\pdfglyphtounicode{planckover2pi1}{210F}
\pdfglyphtounicode{plus}{002B}
\pdfglyphtounicode{plusbelowcmb}{031F}
\pdfglyphtounicode{pluscircle}{2295}
\pdfglyphtounicode{plusminus}{00B1}
\pdfglyphtounicode{plusmod}{02D6}
\pdfglyphtounicode{plusmonospace}{FF0B}
\pdfglyphtounicode{plussmall}{FE62}
\pdfglyphtounicode{plussuperior}{207A}
\pdfglyphtounicode{pmonospace}{FF50}
\pdfglyphtounicode{pmsquare}{33D8}
\pdfglyphtounicode{pohiragana}{307D}
\pdfglyphtounicode{pointingindexdownwhite}{261F}
\pdfglyphtounicode{pointingindexleftwhite}{261C}
\pdfglyphtounicode{pointingindexrightwhite}{261E}
\pdfglyphtounicode{pointingindexupwhite}{261D}
\pdfglyphtounicode{pokatakana}{30DD}
\pdfglyphtounicode{poplathai}{0E1B}
\pdfglyphtounicode{postalmark}{3012}
\pdfglyphtounicode{postalmarkface}{3020}
\pdfglyphtounicode{pparen}{24AB}
\pdfglyphtounicode{precedenotdbleqv}{2AB9}
\pdfglyphtounicode{precedenotslnteql}{2AB5}
\pdfglyphtounicode{precedeornoteqvlnt}{22E8}
\pdfglyphtounicode{precedes}{227A}
\pdfglyphtounicode{precedesequal}{2AAF}
\pdfglyphtounicode{precedesorcurly}{227C}
\pdfglyphtounicode{precedesorequal}{227E}
\pdfglyphtounicode{prescription}{211E}
\pdfglyphtounicode{prime}{2032}
\pdfglyphtounicode{primemod}{02B9}
\pdfglyphtounicode{primereverse}{2035}
\pdfglyphtounicode{primereversed}{2035}
\pdfglyphtounicode{product}{220F}
\pdfglyphtounicode{projective}{2305}
\pdfglyphtounicode{prolongedkana}{30FC}
\pdfglyphtounicode{propellor}{2318}
\pdfglyphtounicode{propersubset}{2282}
\pdfglyphtounicode{propersuperset}{2283}
\pdfglyphtounicode{proportion}{2237}
\pdfglyphtounicode{proportional}{221D}
\pdfglyphtounicode{psi}{03C8}
\pdfglyphtounicode{psicyrillic}{0471}
\pdfglyphtounicode{psilipneumatacyrilliccmb}{0486}
\pdfglyphtounicode{pssquare}{33B0}
\pdfglyphtounicode{puhiragana}{3077}
\pdfglyphtounicode{pukatakana}{30D7}
\pdfglyphtounicode{punctdash}{2014}
\pdfglyphtounicode{pvsquare}{33B4}
\pdfglyphtounicode{pwsquare}{33BA}
\pdfglyphtounicode{q}{0071}
\pdfglyphtounicode{qadeva}{0958}
\pdfglyphtounicode{qadmahebrew}{05A8}
\pdfglyphtounicode{qafarabic}{0642}
\pdfglyphtounicode{qaffinalarabic}{FED6}
\pdfglyphtounicode{qafinitialarabic}{FED7}
\pdfglyphtounicode{qafmedialarabic}{FED8}
\pdfglyphtounicode{qamats}{05B8}
\pdfglyphtounicode{qamats10}{05B8}
\pdfglyphtounicode{qamats1a}{05B8}
\pdfglyphtounicode{qamats1c}{05B8}
\pdfglyphtounicode{qamats27}{05B8}
\pdfglyphtounicode{qamats29}{05B8}
\pdfglyphtounicode{qamats33}{05B8}
\pdfglyphtounicode{qamatsde}{05B8}
\pdfglyphtounicode{qamatshebrew}{05B8}
\pdfglyphtounicode{qamatsnarrowhebrew}{05B8}
\pdfglyphtounicode{qamatsqatanhebrew}{05B8}
\pdfglyphtounicode{qamatsqatannarrowhebrew}{05B8}
\pdfglyphtounicode{qamatsqatanquarterhebrew}{05B8}
\pdfglyphtounicode{qamatsqatanwidehebrew}{05B8}
\pdfglyphtounicode{qamatsquarterhebrew}{05B8}
\pdfglyphtounicode{qamatswidehebrew}{05B8}
\pdfglyphtounicode{qarneyparahebrew}{059F}
\pdfglyphtounicode{qbopomofo}{3111}
\pdfglyphtounicode{qcircle}{24E0}
\pdfglyphtounicode{qhook}{02A0}
\pdfglyphtounicode{qmonospace}{FF51}
\pdfglyphtounicode{qof}{05E7}
\pdfglyphtounicode{qofdagesh}{FB47}
\pdfglyphtounicode{qofdageshhebrew}{FB47}
\pdfglyphtounicode{qofhatafpatah}{05E7 05B2}
\pdfglyphtounicode{qofhatafpatahhebrew}{05E7 05B2}
\pdfglyphtounicode{qofhatafsegol}{05E7 05B1}
\pdfglyphtounicode{qofhatafsegolhebrew}{05E7 05B1}
\pdfglyphtounicode{qofhebrew}{05E7}
\pdfglyphtounicode{qofhiriq}{05E7 05B4}
\pdfglyphtounicode{qofhiriqhebrew}{05E7 05B4}
\pdfglyphtounicode{qofholam}{05E7 05B9}
\pdfglyphtounicode{qofholamhebrew}{05E7 05B9}
\pdfglyphtounicode{qofpatah}{05E7 05B7}
\pdfglyphtounicode{qofpatahhebrew}{05E7 05B7}
\pdfglyphtounicode{qofqamats}{05E7 05B8}
\pdfglyphtounicode{qofqamatshebrew}{05E7 05B8}
\pdfglyphtounicode{qofqubuts}{05E7 05BB}
\pdfglyphtounicode{qofqubutshebrew}{05E7 05BB}
\pdfglyphtounicode{qofsegol}{05E7 05B6}
\pdfglyphtounicode{qofsegolhebrew}{05E7 05B6}
\pdfglyphtounicode{qofsheva}{05E7 05B0}
\pdfglyphtounicode{qofshevahebrew}{05E7 05B0}
\pdfglyphtounicode{qoftsere}{05E7 05B5}
\pdfglyphtounicode{qoftserehebrew}{05E7 05B5}
\pdfglyphtounicode{qparen}{24AC}
\pdfglyphtounicode{quarternote}{2669}
\pdfglyphtounicode{qubuts}{05BB}
\pdfglyphtounicode{qubuts18}{05BB}
\pdfglyphtounicode{qubuts25}{05BB}
\pdfglyphtounicode{qubuts31}{05BB}
\pdfglyphtounicode{qubutshebrew}{05BB}
\pdfglyphtounicode{qubutsnarrowhebrew}{05BB}
\pdfglyphtounicode{qubutsquarterhebrew}{05BB}
\pdfglyphtounicode{qubutswidehebrew}{05BB}
\pdfglyphtounicode{question}{003F}
\pdfglyphtounicode{questionarabic}{061F}
\pdfglyphtounicode{questionarmenian}{055E}
\pdfglyphtounicode{questiondown}{00BF}
\pdfglyphtounicode{questiondownsmall}{00BF}
\pdfglyphtounicode{questiongreek}{037E}
\pdfglyphtounicode{questionmonospace}{FF1F}
\pdfglyphtounicode{questionsmall}{003F}
\pdfglyphtounicode{quotedbl}{0022}
\pdfglyphtounicode{quotedblbase}{201E}
\pdfglyphtounicode{quotedblleft}{201C}
\pdfglyphtounicode{quotedblmonospace}{FF02}
\pdfglyphtounicode{quotedblprime}{301E}
\pdfglyphtounicode{quotedblprimereversed}{301D}
\pdfglyphtounicode{quotedblright}{201D}
\pdfglyphtounicode{quoteleft}{2018}
\pdfglyphtounicode{quoteleftreversed}{201B}
\pdfglyphtounicode{quotereversed}{201B}
\pdfglyphtounicode{quoteright}{2019}
\pdfglyphtounicode{quoterightn}{0149}
\pdfglyphtounicode{quotesinglbase}{201A}
\pdfglyphtounicode{quotesingle}{0027}
\pdfglyphtounicode{quotesinglemonospace}{FF07}
\pdfglyphtounicode{r}{0072}
\pdfglyphtounicode{raarmenian}{057C}
\pdfglyphtounicode{rabengali}{09B0}
\pdfglyphtounicode{racute}{0155}
\pdfglyphtounicode{radeva}{0930}
\pdfglyphtounicode{radical}{221A}
\pdfglyphtounicode{radicalex}{F8E5}
\pdfglyphtounicode{radoverssquare}{33AE}
\pdfglyphtounicode{radoverssquaredsquare}{33AF}
\pdfglyphtounicode{radsquare}{33AD}
\pdfglyphtounicode{rafe}{05BF}
\pdfglyphtounicode{rafehebrew}{05BF}
\pdfglyphtounicode{ragujarati}{0AB0}
\pdfglyphtounicode{ragurmukhi}{0A30}
\pdfglyphtounicode{rahiragana}{3089}
\pdfglyphtounicode{rakatakana}{30E9}
\pdfglyphtounicode{rakatakanahalfwidth}{FF97}
\pdfglyphtounicode{ralowerdiagonalbengali}{09F1}
\pdfglyphtounicode{ramiddlediagonalbengali}{09F0}
\pdfglyphtounicode{ramshorn}{0264}
\pdfglyphtounicode{rangedash}{2013}
\pdfglyphtounicode{ratio}{2236}
\pdfglyphtounicode{rbopomofo}{3116}
\pdfglyphtounicode{rcaron}{0159}
\pdfglyphtounicode{rcedilla}{0157}
\pdfglyphtounicode{rcircle}{24E1}
\pdfglyphtounicode{rcommaaccent}{0157}
\pdfglyphtounicode{rdblgrave}{0211}
\pdfglyphtounicode{rdotaccent}{1E59}
\pdfglyphtounicode{rdotbelow}{1E5B}
\pdfglyphtounicode{rdotbelowmacron}{1E5D}
\pdfglyphtounicode{referencemark}{203B}
\pdfglyphtounicode{reflexsubset}{2286}
\pdfglyphtounicode{reflexsuperset}{2287}
\pdfglyphtounicode{registered}{00AE}
\pdfglyphtounicode{registersans}{00AE}
\pdfglyphtounicode{registerserif}{00AE}
\pdfglyphtounicode{reharabic}{0631}
\pdfglyphtounicode{reharmenian}{0580}
\pdfglyphtounicode{rehfinalarabic}{FEAE}
\pdfglyphtounicode{rehiragana}{308C}
\pdfglyphtounicode{rehyehaleflamarabic}{0631 FEF3 FE8E 0644}
\pdfglyphtounicode{rekatakana}{30EC}
\pdfglyphtounicode{rekatakanahalfwidth}{FF9A}
\pdfglyphtounicode{resh}{05E8}
\pdfglyphtounicode{reshdageshhebrew}{FB48}
\pdfglyphtounicode{reshhatafpatah}{05E8 05B2}
\pdfglyphtounicode{reshhatafpatahhebrew}{05E8 05B2}
\pdfglyphtounicode{reshhatafsegol}{05E8 05B1}
\pdfglyphtounicode{reshhatafsegolhebrew}{05E8 05B1}
\pdfglyphtounicode{reshhebrew}{05E8}
\pdfglyphtounicode{reshhiriq}{05E8 05B4}
\pdfglyphtounicode{reshhiriqhebrew}{05E8 05B4}
\pdfglyphtounicode{reshholam}{05E8 05B9}
\pdfglyphtounicode{reshholamhebrew}{05E8 05B9}
\pdfglyphtounicode{reshpatah}{05E8 05B7}
\pdfglyphtounicode{reshpatahhebrew}{05E8 05B7}
\pdfglyphtounicode{reshqamats}{05E8 05B8}
\pdfglyphtounicode{reshqamatshebrew}{05E8 05B8}
\pdfglyphtounicode{reshqubuts}{05E8 05BB}
\pdfglyphtounicode{reshqubutshebrew}{05E8 05BB}
\pdfglyphtounicode{reshsegol}{05E8 05B6}
\pdfglyphtounicode{reshsegolhebrew}{05E8 05B6}
\pdfglyphtounicode{reshsheva}{05E8 05B0}
\pdfglyphtounicode{reshshevahebrew}{05E8 05B0}
\pdfglyphtounicode{reshtsere}{05E8 05B5}
\pdfglyphtounicode{reshtserehebrew}{05E8 05B5}
\pdfglyphtounicode{revasymptequal}{22CD}
\pdfglyphtounicode{reversedtilde}{223D}
\pdfglyphtounicode{reviahebrew}{0597}
\pdfglyphtounicode{reviamugrashhebrew}{0597}
\pdfglyphtounicode{revlogicalnot}{2310}
\pdfglyphtounicode{revsimilar}{223D}
\pdfglyphtounicode{rfishhook}{027E}
\pdfglyphtounicode{rfishhookreversed}{027F}
\pdfglyphtounicode{rhabengali}{09DD}
\pdfglyphtounicode{rhadeva}{095D}
\pdfglyphtounicode{rho}{03C1}
\pdfglyphtounicode{rho1}{03F1}
\pdfglyphtounicode{rhook}{027D}
\pdfglyphtounicode{rhookturned}{027B}
\pdfglyphtounicode{rhookturnedsuperior}{02B5}
\pdfglyphtounicode{rhosymbolgreek}{03F1}
\pdfglyphtounicode{rhotichookmod}{02DE}
\pdfglyphtounicode{rieulacirclekorean}{3271}
\pdfglyphtounicode{rieulaparenkorean}{3211}
\pdfglyphtounicode{rieulcirclekorean}{3263}
\pdfglyphtounicode{rieulhieuhkorean}{3140}
\pdfglyphtounicode{rieulkiyeokkorean}{313A}
\pdfglyphtounicode{rieulkiyeoksioskorean}{3169}
\pdfglyphtounicode{rieulkorean}{3139}
\pdfglyphtounicode{rieulmieumkorean}{313B}
\pdfglyphtounicode{rieulpansioskorean}{316C}
\pdfglyphtounicode{rieulparenkorean}{3203}
\pdfglyphtounicode{rieulphieuphkorean}{313F}
\pdfglyphtounicode{rieulpieupkorean}{313C}
\pdfglyphtounicode{rieulpieupsioskorean}{316B}
\pdfglyphtounicode{rieulsioskorean}{313D}
\pdfglyphtounicode{rieulthieuthkorean}{313E}
\pdfglyphtounicode{rieultikeutkorean}{316A}
\pdfglyphtounicode{rieulyeorinhieuhkorean}{316D}
\pdfglyphtounicode{rightangle}{221F}
\pdfglyphtounicode{rightanglene}{231D}
\pdfglyphtounicode{rightanglenw}{231C}
\pdfglyphtounicode{rightanglese}{231F}
\pdfglyphtounicode{rightanglesw}{231E}
\pdfglyphtounicode{righttackbelowcmb}{0319}
\pdfglyphtounicode{righttriangle}{22BF}
\pdfglyphtounicode{rihiragana}{308A}
\pdfglyphtounicode{rikatakana}{30EA}
\pdfglyphtounicode{rikatakanahalfwidth}{FF98}
\pdfglyphtounicode{ring}{02DA}
\pdfglyphtounicode{ringbelowcmb}{0325}
\pdfglyphtounicode{ringcmb}{030A}
\pdfglyphtounicode{ringhalfleft}{02BF}
\pdfglyphtounicode{ringhalfleftarmenian}{0559}
\pdfglyphtounicode{ringhalfleftbelowcmb}{031C}
\pdfglyphtounicode{ringhalfleftcentered}{02D3}
\pdfglyphtounicode{ringhalfright}{02BE}
\pdfglyphtounicode{ringhalfrightbelowcmb}{0339}
\pdfglyphtounicode{ringhalfrightcentered}{02D2}
\pdfglyphtounicode{ringinequal}{2256}
\pdfglyphtounicode{rinvertedbreve}{0213}
\pdfglyphtounicode{rittorusquare}{3351}
\pdfglyphtounicode{rlinebelow}{1E5F}
\pdfglyphtounicode{rlongleg}{027C}
\pdfglyphtounicode{rlonglegturned}{027A}
\pdfglyphtounicode{rmonospace}{FF52}
\pdfglyphtounicode{rohiragana}{308D}
\pdfglyphtounicode{rokatakana}{30ED}
\pdfglyphtounicode{rokatakanahalfwidth}{FF9B}
\pdfglyphtounicode{roruathai}{0E23}
\pdfglyphtounicode{rparen}{24AD}
\pdfglyphtounicode{rrabengali}{09DC}
\pdfglyphtounicode{rradeva}{0931}
\pdfglyphtounicode{rragurmukhi}{0A5C}
\pdfglyphtounicode{rreharabic}{0691}
\pdfglyphtounicode{rrehfinalarabic}{FB8D}
\pdfglyphtounicode{rrvocalicbengali}{09E0}
\pdfglyphtounicode{rrvocalicdeva}{0960}
\pdfglyphtounicode{rrvocalicgujarati}{0AE0}
\pdfglyphtounicode{rrvocalicvowelsignbengali}{09C4}
\pdfglyphtounicode{rrvocalicvowelsigndeva}{0944}
\pdfglyphtounicode{rrvocalicvowelsigngujarati}{0AC4}
\pdfglyphtounicode{rsuperior}{0072}
\pdfglyphtounicode{rtblock}{2590}
\pdfglyphtounicode{rturned}{0279}
\pdfglyphtounicode{rturnedsuperior}{02B4}
\pdfglyphtounicode{ruhiragana}{308B}
\pdfglyphtounicode{rukatakana}{30EB}
\pdfglyphtounicode{rukatakanahalfwidth}{FF99}
\pdfglyphtounicode{rupeemarkbengali}{09F2}
\pdfglyphtounicode{rupeesignbengali}{09F3}
\pdfglyphtounicode{rupiah}{20A8}
\pdfglyphtounicode{ruthai}{0E24}
\pdfglyphtounicode{rvocalicbengali}{098B}
\pdfglyphtounicode{rvocalicdeva}{090B}
\pdfglyphtounicode{rvocalicgujarati}{0A8B}
\pdfglyphtounicode{rvocalicvowelsignbengali}{09C3}
\pdfglyphtounicode{rvocalicvowelsigndeva}{0943}
\pdfglyphtounicode{rvocalicvowelsigngujarati}{0AC3}
\pdfglyphtounicode{s}{0073}
\pdfglyphtounicode{sabengali}{09B8}
\pdfglyphtounicode{sacute}{015B}
\pdfglyphtounicode{sacutedotaccent}{1E65}
\pdfglyphtounicode{sadarabic}{0635}
\pdfglyphtounicode{sadeva}{0938}
\pdfglyphtounicode{sadfinalarabic}{FEBA}
\pdfglyphtounicode{sadinitialarabic}{FEBB}
\pdfglyphtounicode{sadmedialarabic}{FEBC}
\pdfglyphtounicode{sagujarati}{0AB8}
\pdfglyphtounicode{sagurmukhi}{0A38}
\pdfglyphtounicode{sahiragana}{3055}
\pdfglyphtounicode{sakatakana}{30B5}
\pdfglyphtounicode{sakatakanahalfwidth}{FF7B}
\pdfglyphtounicode{sallallahoualayhewasallamarabic}{FDFA}
\pdfglyphtounicode{samekh}{05E1}
\pdfglyphtounicode{samekhdagesh}{FB41}
\pdfglyphtounicode{samekhdageshhebrew}{FB41}
\pdfglyphtounicode{samekhhebrew}{05E1}
\pdfglyphtounicode{saraaathai}{0E32}
\pdfglyphtounicode{saraaethai}{0E41}
\pdfglyphtounicode{saraaimaimalaithai}{0E44}
\pdfglyphtounicode{saraaimaimuanthai}{0E43}
\pdfglyphtounicode{saraamthai}{0E33}
\pdfglyphtounicode{saraathai}{0E30}
\pdfglyphtounicode{saraethai}{0E40}
\pdfglyphtounicode{saraiileftthai}{F886}
\pdfglyphtounicode{saraiithai}{0E35}
\pdfglyphtounicode{saraileftthai}{F885}
\pdfglyphtounicode{saraithai}{0E34}
\pdfglyphtounicode{saraothai}{0E42}
\pdfglyphtounicode{saraueeleftthai}{F888}
\pdfglyphtounicode{saraueethai}{0E37}
\pdfglyphtounicode{saraueleftthai}{F887}
\pdfglyphtounicode{sarauethai}{0E36}
\pdfglyphtounicode{sarauthai}{0E38}
\pdfglyphtounicode{sarauuthai}{0E39}
\pdfglyphtounicode{satisfies}{22A8}
\pdfglyphtounicode{sbopomofo}{3119}
\pdfglyphtounicode{scaron}{0161}
\pdfglyphtounicode{scarondotaccent}{1E67}
\pdfglyphtounicode{scedilla}{015F}
\pdfglyphtounicode{schwa}{0259}
\pdfglyphtounicode{schwacyrillic}{04D9}
\pdfglyphtounicode{schwadieresiscyrillic}{04DB}
\pdfglyphtounicode{schwahook}{025A}
\pdfglyphtounicode{scircle}{24E2}
\pdfglyphtounicode{scircumflex}{015D}
\pdfglyphtounicode{scommaaccent}{0219}
\pdfglyphtounicode{sdotaccent}{1E61}
\pdfglyphtounicode{sdotbelow}{1E63}
\pdfglyphtounicode{sdotbelowdotaccent}{1E69}
\pdfglyphtounicode{seagullbelowcmb}{033C}
\pdfglyphtounicode{second}{2033}
\pdfglyphtounicode{secondtonechinese}{02CA}
\pdfglyphtounicode{section}{00A7}
\pdfglyphtounicode{seenarabic}{0633}
\pdfglyphtounicode{seenfinalarabic}{FEB2}
\pdfglyphtounicode{seeninitialarabic}{FEB3}
\pdfglyphtounicode{seenmedialarabic}{FEB4}
\pdfglyphtounicode{segol}{05B6}
\pdfglyphtounicode{segol13}{05B6}
\pdfglyphtounicode{segol1f}{05B6}
\pdfglyphtounicode{segol2c}{05B6}
\pdfglyphtounicode{segolhebrew}{05B6}
\pdfglyphtounicode{segolnarrowhebrew}{05B6}
\pdfglyphtounicode{segolquarterhebrew}{05B6}
\pdfglyphtounicode{segoltahebrew}{0592}
\pdfglyphtounicode{segolwidehebrew}{05B6}
\pdfglyphtounicode{seharmenian}{057D}
\pdfglyphtounicode{sehiragana}{305B}
\pdfglyphtounicode{sekatakana}{30BB}
\pdfglyphtounicode{sekatakanahalfwidth}{FF7E}
\pdfglyphtounicode{semicolon}{003B}
\pdfglyphtounicode{semicolonarabic}{061B}
\pdfglyphtounicode{semicolonmonospace}{FF1B}
\pdfglyphtounicode{semicolonsmall}{FE54}
\pdfglyphtounicode{semivoicedmarkkana}{309C}
\pdfglyphtounicode{semivoicedmarkkanahalfwidth}{FF9F}
\pdfglyphtounicode{sentisquare}{3322}
\pdfglyphtounicode{sentosquare}{3323}
\pdfglyphtounicode{seven}{0037}
\pdfglyphtounicode{sevenarabic}{0667}
\pdfglyphtounicode{sevenbengali}{09ED}
\pdfglyphtounicode{sevencircle}{2466}
\pdfglyphtounicode{sevencircleinversesansserif}{2790}
\pdfglyphtounicode{sevendeva}{096D}
\pdfglyphtounicode{seveneighths}{215E}
\pdfglyphtounicode{sevengujarati}{0AED}
\pdfglyphtounicode{sevengurmukhi}{0A6D}
\pdfglyphtounicode{sevenhackarabic}{0667}
\pdfglyphtounicode{sevenhangzhou}{3027}
\pdfglyphtounicode{sevenideographicparen}{3226}
\pdfglyphtounicode{seveninferior}{2087}
\pdfglyphtounicode{sevenmonospace}{FF17}
\pdfglyphtounicode{sevenoldstyle}{0037}
\pdfglyphtounicode{sevenparen}{247A}
\pdfglyphtounicode{sevenperiod}{248E}
\pdfglyphtounicode{sevenpersian}{06F7}
\pdfglyphtounicode{sevenroman}{2176}
\pdfglyphtounicode{sevensuperior}{2077}
\pdfglyphtounicode{seventeencircle}{2470}
\pdfglyphtounicode{seventeenparen}{2484}
\pdfglyphtounicode{seventeenperiod}{2498}
\pdfglyphtounicode{seventhai}{0E57}
\pdfglyphtounicode{sfthyphen}{00AD}
\pdfglyphtounicode{shaarmenian}{0577}
\pdfglyphtounicode{shabengali}{09B6}
\pdfglyphtounicode{shacyrillic}{0448}
\pdfglyphtounicode{shaddaarabic}{0651}
\pdfglyphtounicode{shaddadammaarabic}{FC61}
\pdfglyphtounicode{shaddadammatanarabic}{FC5E}
\pdfglyphtounicode{shaddafathaarabic}{FC60}
\pdfglyphtounicode{shaddafathatanarabic}{0651 064B}
\pdfglyphtounicode{shaddakasraarabic}{FC62}
\pdfglyphtounicode{shaddakasratanarabic}{FC5F}
\pdfglyphtounicode{shade}{2592}
\pdfglyphtounicode{shadedark}{2593}
\pdfglyphtounicode{shadelight}{2591}
\pdfglyphtounicode{shademedium}{2592}
\pdfglyphtounicode{shadeva}{0936}
\pdfglyphtounicode{shagujarati}{0AB6}
\pdfglyphtounicode{shagurmukhi}{0A36}
\pdfglyphtounicode{shalshelethebrew}{0593}
\pdfglyphtounicode{sharp}{266F}
\pdfglyphtounicode{shbopomofo}{3115}
\pdfglyphtounicode{shchacyrillic}{0449}
\pdfglyphtounicode{sheenarabic}{0634}
\pdfglyphtounicode{sheenfinalarabic}{FEB6}
\pdfglyphtounicode{sheeninitialarabic}{FEB7}
\pdfglyphtounicode{sheenmedialarabic}{FEB8}
\pdfglyphtounicode{sheicoptic}{03E3}
\pdfglyphtounicode{sheqel}{20AA}
\pdfglyphtounicode{sheqelhebrew}{20AA}
\pdfglyphtounicode{sheva}{05B0}
\pdfglyphtounicode{sheva115}{05B0}
\pdfglyphtounicode{sheva15}{05B0}
\pdfglyphtounicode{sheva22}{05B0}
\pdfglyphtounicode{sheva2e}{05B0}
\pdfglyphtounicode{shevahebrew}{05B0}
\pdfglyphtounicode{shevanarrowhebrew}{05B0}
\pdfglyphtounicode{shevaquarterhebrew}{05B0}
\pdfglyphtounicode{shevawidehebrew}{05B0}
\pdfglyphtounicode{shhacyrillic}{04BB}
\pdfglyphtounicode{shiftleft}{21B0}
\pdfglyphtounicode{shiftright}{21B1}
\pdfglyphtounicode{shimacoptic}{03ED}
\pdfglyphtounicode{shin}{05E9}
\pdfglyphtounicode{shindagesh}{FB49}
\pdfglyphtounicode{shindageshhebrew}{FB49}
\pdfglyphtounicode{shindageshshindot}{FB2C}
\pdfglyphtounicode{shindageshshindothebrew}{FB2C}
\pdfglyphtounicode{shindageshsindot}{FB2D}
\pdfglyphtounicode{shindageshsindothebrew}{FB2D}
\pdfglyphtounicode{shindothebrew}{05C1}
\pdfglyphtounicode{shinhebrew}{05E9}
\pdfglyphtounicode{shinshindot}{FB2A}
\pdfglyphtounicode{shinshindothebrew}{FB2A}
\pdfglyphtounicode{shinsindot}{FB2B}
\pdfglyphtounicode{shinsindothebrew}{FB2B}
\pdfglyphtounicode{shook}{0282}
\pdfglyphtounicode{sigma}{03C3}
\pdfglyphtounicode{sigma1}{03C2}
\pdfglyphtounicode{sigmafinal}{03C2}
\pdfglyphtounicode{sigmalunatesymbolgreek}{03F2}
\pdfglyphtounicode{sihiragana}{3057}
\pdfglyphtounicode{sikatakana}{30B7}
\pdfglyphtounicode{sikatakanahalfwidth}{FF7C}
\pdfglyphtounicode{siluqhebrew}{05BD}
\pdfglyphtounicode{siluqlefthebrew}{05BD}
\pdfglyphtounicode{similar}{223C}
\pdfglyphtounicode{similarequal}{2243}
\pdfglyphtounicode{sindothebrew}{05C2}
\pdfglyphtounicode{siosacirclekorean}{3274}
\pdfglyphtounicode{siosaparenkorean}{3214}
\pdfglyphtounicode{sioscieuckorean}{317E}
\pdfglyphtounicode{sioscirclekorean}{3266}
\pdfglyphtounicode{sioskiyeokkorean}{317A}
\pdfglyphtounicode{sioskorean}{3145}
\pdfglyphtounicode{siosnieunkorean}{317B}
\pdfglyphtounicode{siosparenkorean}{3206}
\pdfglyphtounicode{siospieupkorean}{317D}
\pdfglyphtounicode{siostikeutkorean}{317C}
\pdfglyphtounicode{six}{0036}
\pdfglyphtounicode{sixarabic}{0666}
\pdfglyphtounicode{sixbengali}{09EC}
\pdfglyphtounicode{sixcircle}{2465}
\pdfglyphtounicode{sixcircleinversesansserif}{278F}
\pdfglyphtounicode{sixdeva}{096C}
\pdfglyphtounicode{sixgujarati}{0AEC}
\pdfglyphtounicode{sixgurmukhi}{0A6C}
\pdfglyphtounicode{sixhackarabic}{0666}
\pdfglyphtounicode{sixhangzhou}{3026}
\pdfglyphtounicode{sixideographicparen}{3225}
\pdfglyphtounicode{sixinferior}{2086}
\pdfglyphtounicode{sixmonospace}{FF16}
\pdfglyphtounicode{sixoldstyle}{0036}
\pdfglyphtounicode{sixparen}{2479}
\pdfglyphtounicode{sixperiod}{248D}
\pdfglyphtounicode{sixpersian}{06F6}
\pdfglyphtounicode{sixroman}{2175}
\pdfglyphtounicode{sixsuperior}{2076}
\pdfglyphtounicode{sixteencircle}{246F}
\pdfglyphtounicode{sixteencurrencydenominatorbengali}{09F9}
\pdfglyphtounicode{sixteenparen}{2483}
\pdfglyphtounicode{sixteenperiod}{2497}
\pdfglyphtounicode{sixthai}{0E56}
\pdfglyphtounicode{slash}{002F}
\pdfglyphtounicode{slashmonospace}{FF0F}
\pdfglyphtounicode{slong}{017F}
\pdfglyphtounicode{slongdotaccent}{1E9B}
\pdfglyphtounicode{slurabove}{2322}
\pdfglyphtounicode{slurbelow}{2323}
\pdfglyphtounicode{smile}{2323}
\pdfglyphtounicode{smileface}{263A}
\pdfglyphtounicode{smonospace}{FF53}
\pdfglyphtounicode{sofpasuqhebrew}{05C3}
\pdfglyphtounicode{softhyphen}{00AD}
\pdfglyphtounicode{softsigncyrillic}{044C}
\pdfglyphtounicode{sohiragana}{305D}
\pdfglyphtounicode{sokatakana}{30BD}
\pdfglyphtounicode{sokatakanahalfwidth}{FF7F}
\pdfglyphtounicode{soliduslongoverlaycmb}{0338}
\pdfglyphtounicode{solidusshortoverlaycmb}{0337}
\pdfglyphtounicode{sorusithai}{0E29}
\pdfglyphtounicode{sosalathai}{0E28}
\pdfglyphtounicode{sosothai}{0E0B}
\pdfglyphtounicode{sosuathai}{0E2A}
\pdfglyphtounicode{space}{0020}
\pdfglyphtounicode{spacehackarabic}{0020}
\pdfglyphtounicode{spade}{2660}
\pdfglyphtounicode{spadesuitblack}{2660}
\pdfglyphtounicode{spadesuitwhite}{2664}
\pdfglyphtounicode{sparen}{24AE}
\pdfglyphtounicode{sphericalangle}{2222}
\pdfglyphtounicode{square}{25A1}
\pdfglyphtounicode{squarebelowcmb}{033B}
\pdfglyphtounicode{squarecc}{33C4}
\pdfglyphtounicode{squarecm}{339D}
\pdfglyphtounicode{squarediagonalcrosshatchfill}{25A9}
\pdfglyphtounicode{squaredot}{22A1}
\pdfglyphtounicode{squarehorizontalfill}{25A4}
\pdfglyphtounicode{squareimage}{228F}
\pdfglyphtounicode{squarekg}{338F}
\pdfglyphtounicode{squarekm}{339E}
\pdfglyphtounicode{squarekmcapital}{33CE}
\pdfglyphtounicode{squareln}{33D1}
\pdfglyphtounicode{squarelog}{33D2}
\pdfglyphtounicode{squaremg}{338E}
\pdfglyphtounicode{squaremil}{33D5}
\pdfglyphtounicode{squareminus}{229F}
\pdfglyphtounicode{squaremm}{339C}
\pdfglyphtounicode{squaremsquared}{33A1}
\pdfglyphtounicode{squaremultiply}{22A0}
\pdfglyphtounicode{squareoriginal}{2290}
\pdfglyphtounicode{squareorthogonalcrosshatchfill}{25A6}
\pdfglyphtounicode{squareplus}{229E}
\pdfglyphtounicode{squaresolid}{25A0}
\pdfglyphtounicode{squareupperlefttolowerrightfill}{25A7}
\pdfglyphtounicode{squareupperrighttolowerleftfill}{25A8}
\pdfglyphtounicode{squareverticalfill}{25A5}
\pdfglyphtounicode{squarewhitewithsmallblack}{25A3}
\pdfglyphtounicode{squiggleleftright}{21AD}
\pdfglyphtounicode{squiggleright}{21DD}
\pdfglyphtounicode{srsquare}{33DB}
\pdfglyphtounicode{ssabengali}{09B7}
\pdfglyphtounicode{ssadeva}{0937}
\pdfglyphtounicode{ssagujarati}{0AB7}
\pdfglyphtounicode{ssangcieuckorean}{3149}
\pdfglyphtounicode{ssanghieuhkorean}{3185}
\pdfglyphtounicode{ssangieungkorean}{3180}
\pdfglyphtounicode{ssangkiyeokkorean}{3132}
\pdfglyphtounicode{ssangnieunkorean}{3165}
\pdfglyphtounicode{ssangpieupkorean}{3143}
\pdfglyphtounicode{ssangsioskorean}{3146}
\pdfglyphtounicode{ssangtikeutkorean}{3138}
\pdfglyphtounicode{ssuperior}{0073}
\pdfglyphtounicode{st}{0073 0074}
\pdfglyphtounicode{star}{22C6}
\pdfglyphtounicode{sterling}{00A3}
\pdfglyphtounicode{sterlingmonospace}{FFE1}
\pdfglyphtounicode{strokelongoverlaycmb}{0336}
\pdfglyphtounicode{strokeshortoverlaycmb}{0335}
\pdfglyphtounicode{subset}{2282}
\pdfglyphtounicode{subsetdbl}{22D0}
\pdfglyphtounicode{subsetdblequal}{2AC5}
\pdfglyphtounicode{subsetnoteql}{228A}
\pdfglyphtounicode{subsetnotequal}{228A}
\pdfglyphtounicode{subsetorequal}{2286}
\pdfglyphtounicode{subsetornotdbleql}{2ACB}
\pdfglyphtounicode{subsetsqequal}{2291}
\pdfglyphtounicode{succeeds}{227B}
\pdfglyphtounicode{suchthat}{220B}
\pdfglyphtounicode{suhiragana}{3059}
\pdfglyphtounicode{sukatakana}{30B9}
\pdfglyphtounicode{sukatakanahalfwidth}{FF7D}
\pdfglyphtounicode{sukunarabic}{0652}
\pdfglyphtounicode{summation}{2211}
\pdfglyphtounicode{sun}{263C}
\pdfglyphtounicode{superset}{2283}
\pdfglyphtounicode{supersetdbl}{22D1}
\pdfglyphtounicode{supersetdblequal}{2AC6}
\pdfglyphtounicode{supersetnoteql}{228B}
\pdfglyphtounicode{supersetnotequal}{228B}
\pdfglyphtounicode{supersetorequal}{2287}
\pdfglyphtounicode{supersetornotdbleql}{2ACC}
\pdfglyphtounicode{supersetsqequal}{2292}
\pdfglyphtounicode{svsquare}{33DC}
\pdfglyphtounicode{syouwaerasquare}{337C}
\pdfglyphtounicode{t}{0074}
\pdfglyphtounicode{tabengali}{09A4}
\pdfglyphtounicode{tackdown}{22A4}
\pdfglyphtounicode{tackleft}{22A3}
\pdfglyphtounicode{tadeva}{0924}
\pdfglyphtounicode{tagujarati}{0AA4}
\pdfglyphtounicode{tagurmukhi}{0A24}
\pdfglyphtounicode{taharabic}{0637}
\pdfglyphtounicode{tahfinalarabic}{FEC2}
\pdfglyphtounicode{tahinitialarabic}{FEC3}
\pdfglyphtounicode{tahiragana}{305F}
\pdfglyphtounicode{tahmedialarabic}{FEC4}
\pdfglyphtounicode{taisyouerasquare}{337D}
\pdfglyphtounicode{takatakana}{30BF}
\pdfglyphtounicode{takatakanahalfwidth}{FF80}
\pdfglyphtounicode{tatweelarabic}{0640}
\pdfglyphtounicode{tau}{03C4}
\pdfglyphtounicode{tav}{05EA}
\pdfglyphtounicode{tavdages}{FB4A}
\pdfglyphtounicode{tavdagesh}{FB4A}
\pdfglyphtounicode{tavdageshhebrew}{FB4A}
\pdfglyphtounicode{tavhebrew}{05EA}
\pdfglyphtounicode{tbar}{0167}
\pdfglyphtounicode{tbopomofo}{310A}
\pdfglyphtounicode{tcaron}{0165}
\pdfglyphtounicode{tccurl}{02A8}
\pdfglyphtounicode{tcedilla}{0163}
\pdfglyphtounicode{tcheharabic}{0686}
\pdfglyphtounicode{tchehfinalarabic}{FB7B}
\pdfglyphtounicode{tchehinitialarabic}{FB7C}
\pdfglyphtounicode{tchehmedialarabic}{FB7D}
\pdfglyphtounicode{tchehmeeminitialarabic}{FB7C FEE4}
\pdfglyphtounicode{tcircle}{24E3}
\pdfglyphtounicode{tcircumflexbelow}{1E71}
\pdfglyphtounicode{tcommaaccent}{0163}
\pdfglyphtounicode{tdieresis}{1E97}
\pdfglyphtounicode{tdotaccent}{1E6B}
\pdfglyphtounicode{tdotbelow}{1E6D}
\pdfglyphtounicode{tecyrillic}{0442}
\pdfglyphtounicode{tedescendercyrillic}{04AD}
\pdfglyphtounicode{teharabic}{062A}
\pdfglyphtounicode{tehfinalarabic}{FE96}
\pdfglyphtounicode{tehhahinitialarabic}{FCA2}
\pdfglyphtounicode{tehhahisolatedarabic}{FC0C}
\pdfglyphtounicode{tehinitialarabic}{FE97}
\pdfglyphtounicode{tehiragana}{3066}
\pdfglyphtounicode{tehjeeminitialarabic}{FCA1}
\pdfglyphtounicode{tehjeemisolatedarabic}{FC0B}
\pdfglyphtounicode{tehmarbutaarabic}{0629}
\pdfglyphtounicode{tehmarbutafinalarabic}{FE94}
\pdfglyphtounicode{tehmedialarabic}{FE98}
\pdfglyphtounicode{tehmeeminitialarabic}{FCA4}
\pdfglyphtounicode{tehmeemisolatedarabic}{FC0E}
\pdfglyphtounicode{tehnoonfinalarabic}{FC73}
\pdfglyphtounicode{tekatakana}{30C6}
\pdfglyphtounicode{tekatakanahalfwidth}{FF83}
\pdfglyphtounicode{telephone}{2121}
\pdfglyphtounicode{telephoneblack}{260E}
\pdfglyphtounicode{telishagedolahebrew}{05A0}
\pdfglyphtounicode{telishaqetanahebrew}{05A9}
\pdfglyphtounicode{tencircle}{2469}
\pdfglyphtounicode{tenideographicparen}{3229}
\pdfglyphtounicode{tenparen}{247D}
\pdfglyphtounicode{tenperiod}{2491}
\pdfglyphtounicode{tenroman}{2179}
\pdfglyphtounicode{tesh}{02A7}
\pdfglyphtounicode{tet}{05D8}
\pdfglyphtounicode{tetdagesh}{FB38}
\pdfglyphtounicode{tetdageshhebrew}{FB38}
\pdfglyphtounicode{tethebrew}{05D8}
\pdfglyphtounicode{tetsecyrillic}{04B5}
\pdfglyphtounicode{tevirhebrew}{059B}
\pdfglyphtounicode{tevirlefthebrew}{059B}
\pdfglyphtounicode{tfm:cmbsy10/diamond}{2662}
\pdfglyphtounicode{tfm:cmbsy10/heart}{2661}
\pdfglyphtounicode{tfm:cmbsy5/diamond}{2662}
\pdfglyphtounicode{tfm:cmbsy5/heart}{2661}
\pdfglyphtounicode{tfm:cmbsy6/diamond}{2662}
\pdfglyphtounicode{tfm:cmbsy6/heart}{2661}
\pdfglyphtounicode{tfm:cmbsy7/diamond}{2662}
\pdfglyphtounicode{tfm:cmbsy7/heart}{2661}
\pdfglyphtounicode{tfm:cmbsy8/diamond}{2662}
\pdfglyphtounicode{tfm:cmbsy8/heart}{2661}
\pdfglyphtounicode{tfm:cmbsy9/diamond}{2662}
\pdfglyphtounicode{tfm:cmbsy9/heart}{2661}
\pdfglyphtounicode{tfm:cmmi10/phi}{03D5}
\pdfglyphtounicode{tfm:cmmi10/phi1}{03C6}
\pdfglyphtounicode{tfm:cmmi12/phi}{03D5}
\pdfglyphtounicode{tfm:cmmi12/phi1}{03C6}
\pdfglyphtounicode{tfm:cmmi5/phi}{03D5}
\pdfglyphtounicode{tfm:cmmi5/phi1}{03C6}
\pdfglyphtounicode{tfm:cmmi6/phi}{03D5}
\pdfglyphtounicode{tfm:cmmi6/phi1}{03C6}
\pdfglyphtounicode{tfm:cmmi7/phi}{03D5}
\pdfglyphtounicode{tfm:cmmi7/phi1}{03C6}
\pdfglyphtounicode{tfm:cmmi8/phi}{03D5}
\pdfglyphtounicode{tfm:cmmi8/phi1}{03C6}
\pdfglyphtounicode{tfm:cmmi9/phi}{03D5}
\pdfglyphtounicode{tfm:cmmi9/phi1}{03C6}
\pdfglyphtounicode{tfm:cmmib10/phi}{03D5}
\pdfglyphtounicode{tfm:cmmib10/phi1}{03C6}
\pdfglyphtounicode{tfm:cmmib5/phi}{03D5}
\pdfglyphtounicode{tfm:cmmib5/phi1}{03C6}
\pdfglyphtounicode{tfm:cmmib6/phi}{03D5}
\pdfglyphtounicode{tfm:cmmib6/phi1}{03C6}
\pdfglyphtounicode{tfm:cmmib7/phi}{03D5}
\pdfglyphtounicode{tfm:cmmib7/phi1}{03C6}
\pdfglyphtounicode{tfm:cmmib8/phi}{03D5}
\pdfglyphtounicode{tfm:cmmib8/phi1}{03C6}
\pdfglyphtounicode{tfm:cmmib9/phi}{03D5}
\pdfglyphtounicode{tfm:cmmib9/phi1}{03C6}
\pdfglyphtounicode{tfm:cmsy10/diamond}{2662}
\pdfglyphtounicode{tfm:cmsy10/heart}{2661}
\pdfglyphtounicode{tfm:cmsy5/heart}{2661}
\pdfglyphtounicode{tfm:cmsy6/diamond}{2662}
\pdfglyphtounicode{tfm:cmsy6/heart}{2661}
\pdfglyphtounicode{tfm:cmsy7/diamond}{2662}
\pdfglyphtounicode{tfm:cmsy7/heart}{2661}
\pdfglyphtounicode{tfm:cmsy8/diamond}{2662}
\pdfglyphtounicode{tfm:cmsy8/heart}{2661}
\pdfglyphtounicode{tfm:cmsy9/diamond}{2662}
\pdfglyphtounicode{tfm:cmsy9/heart}{2661}
\pdfglyphtounicode{tfm:eurb10/phi}{03D5}
\pdfglyphtounicode{tfm:eurb10/phi1}{03C6}
\pdfglyphtounicode{tfm:eurb5/phi}{03D5}
\pdfglyphtounicode{tfm:eurb5/phi1}{03C6}
\pdfglyphtounicode{tfm:eurb6/phi}{03D5}
\pdfglyphtounicode{tfm:eurb6/phi1}{03C6}
\pdfglyphtounicode{tfm:eurb7/phi}{03D5}
\pdfglyphtounicode{tfm:eurb7/phi1}{03C6}
\pdfglyphtounicode{tfm:eurb8/phi}{03D5}
\pdfglyphtounicode{tfm:eurb8/phi1}{03C6}
\pdfglyphtounicode{tfm:eurb9/phi}{03D5}
\pdfglyphtounicode{tfm:eurb9/phi1}{03C6}
\pdfglyphtounicode{tfm:eurm10/phi}{03D5}
\pdfglyphtounicode{tfm:eurm10/phi1}{03C6}
\pdfglyphtounicode{tfm:eurm5/phi}{03D5}
\pdfglyphtounicode{tfm:eurm5/phi1}{03C6}
\pdfglyphtounicode{tfm:eurm6/phi}{03D5}
\pdfglyphtounicode{tfm:eurm6/phi1}{03C6}
\pdfglyphtounicode{tfm:eurm7/phi}{03D5}
\pdfglyphtounicode{tfm:eurm7/phi1}{03C6}
\pdfglyphtounicode{tfm:eurm8/phi}{03D5}
\pdfglyphtounicode{tfm:eurm8/phi1}{03C6}
\pdfglyphtounicode{tfm:eurm9/phi}{03D5}
\pdfglyphtounicode{tfm:eurm9/phi1}{03C6}
\pdfglyphtounicode{tfm:fplmbi/phi}{03D5}
\pdfglyphtounicode{tfm:fplmbi/phi1}{03C6}
\pdfglyphtounicode{tfm:fplmri/phi}{03D5}
\pdfglyphtounicode{tfm:fplmri/phi1}{03C6}
\pdfglyphtounicode{tfm:lmbsy10/diamond}{2662}
\pdfglyphtounicode{tfm:lmbsy10/heart}{2661}
\pdfglyphtounicode{tfm:lmbsy5/diamond}{2662}
\pdfglyphtounicode{tfm:lmbsy5/heart}{2661}
\pdfglyphtounicode{tfm:lmbsy7/diamond}{2662}
\pdfglyphtounicode{tfm:lmbsy7/heart}{2661}
\pdfglyphtounicode{tfm:lmmi10/phi}{03D5}
\pdfglyphtounicode{tfm:lmmi10/phi1}{03C6}
\pdfglyphtounicode{tfm:lmmi12/phi}{03D5}
\pdfglyphtounicode{tfm:lmmi12/phi1}{03C6}
\pdfglyphtounicode{tfm:lmmi5/phi}{03D5}
\pdfglyphtounicode{tfm:lmmi5/phi1}{03C6}
\pdfglyphtounicode{tfm:lmmi6/phi}{03D5}
\pdfglyphtounicode{tfm:lmmi6/phi1}{03C6}
\pdfglyphtounicode{tfm:lmmi7/phi}{03D5}
\pdfglyphtounicode{tfm:lmmi7/phi1}{03C6}
\pdfglyphtounicode{tfm:lmmi8/phi}{03D5}
\pdfglyphtounicode{tfm:lmmi8/phi1}{03C6}
\pdfglyphtounicode{tfm:lmmi9/phi}{03D5}
\pdfglyphtounicode{tfm:lmmi9/phi1}{03C6}
\pdfglyphtounicode{tfm:lmmib10/phi}{03D5}
\pdfglyphtounicode{tfm:lmmib10/phi1}{03C6}
\pdfglyphtounicode{tfm:lmmib5/phi}{03D5}
\pdfglyphtounicode{tfm:lmmib5/phi1}{03C6}
\pdfglyphtounicode{tfm:lmmib7/phi}{03D5}
\pdfglyphtounicode{tfm:lmmib7/phi1}{03C6}
\pdfglyphtounicode{tfm:lmsy10/diamond}{2662}
\pdfglyphtounicode{tfm:lmsy10/heart}{2661}
\pdfglyphtounicode{tfm:lmsy5/diamond}{2662}
\pdfglyphtounicode{tfm:lmsy5/heart}{2661}
\pdfglyphtounicode{tfm:lmsy6/diamond}{2662}
\pdfglyphtounicode{tfm:lmsy6/heart}{2661}
\pdfglyphtounicode{tfm:lmsy7/diamond}{2662}
\pdfglyphtounicode{tfm:lmsy7/heart}{2661}
\pdfglyphtounicode{tfm:lmsy8/diamond}{2662}
\pdfglyphtounicode{tfm:lmsy8/heart}{2661}
\pdfglyphtounicode{tfm:lmsy9/diamond}{2662}
\pdfglyphtounicode{tfm:lmsy9/heart}{2661}
\pdfglyphtounicode{tfm:msam10/diamond}{2662}
\pdfglyphtounicode{tfm:msam5/diamond}{2662}
\pdfglyphtounicode{tfm:msam6/diamond}{2662}
\pdfglyphtounicode{tfm:msam7/diamond}{2662}
\pdfglyphtounicode{tfm:msam8/diamond}{2662}
\pdfglyphtounicode{tfm:msam9/diamond}{2662}
\pdfglyphtounicode{tfm:pxbmia/phi}{03D5}
\pdfglyphtounicode{tfm:pxbmia/phi1}{03C6}
\pdfglyphtounicode{tfm:pxbsy/diamond}{2662}
\pdfglyphtounicode{tfm:pxbsy/heart}{2661}
\pdfglyphtounicode{tfm:pxbsya/diamond}{2662}
\pdfglyphtounicode{tfm:pxmia/phi}{03D5}
\pdfglyphtounicode{tfm:pxmia/phi1}{03C6}
\pdfglyphtounicode{tfm:pxsy/diamond}{2662}
\pdfglyphtounicode{tfm:pxsy/heart}{2661}
\pdfglyphtounicode{tfm:pxsya/diamond}{2662}
\pdfglyphtounicode{tfm:pzdr/a1}{2701}
\pdfglyphtounicode{tfm:pzdr/a10}{2721}
\pdfglyphtounicode{tfm:pzdr/a100}{275E}
\pdfglyphtounicode{tfm:pzdr/a101}{2761}
\pdfglyphtounicode{tfm:pzdr/a102}{2762}
\pdfglyphtounicode{tfm:pzdr/a103}{2763}
\pdfglyphtounicode{tfm:pzdr/a104}{2764}
\pdfglyphtounicode{tfm:pzdr/a105}{2710}
\pdfglyphtounicode{tfm:pzdr/a106}{2765}
\pdfglyphtounicode{tfm:pzdr/a107}{2766}
\pdfglyphtounicode{tfm:pzdr/a108}{2767}
\pdfglyphtounicode{tfm:pzdr/a109}{2660}
\pdfglyphtounicode{tfm:pzdr/a11}{261B}
\pdfglyphtounicode{tfm:pzdr/a110}{2665}
\pdfglyphtounicode{tfm:pzdr/a111}{2666}
\pdfglyphtounicode{tfm:pzdr/a112}{2663}
\pdfglyphtounicode{tfm:pzdr/a117}{2709}
\pdfglyphtounicode{tfm:pzdr/a118}{2708}
\pdfglyphtounicode{tfm:pzdr/a119}{2707}
\pdfglyphtounicode{tfm:pzdr/a12}{261E}
\pdfglyphtounicode{tfm:pzdr/a120}{2460}
\pdfglyphtounicode{tfm:pzdr/a121}{2461}
\pdfglyphtounicode{tfm:pzdr/a122}{2462}
\pdfglyphtounicode{tfm:pzdr/a123}{2463}
\pdfglyphtounicode{tfm:pzdr/a124}{2464}
\pdfglyphtounicode{tfm:pzdr/a125}{2465}
\pdfglyphtounicode{tfm:pzdr/a126}{2466}
\pdfglyphtounicode{tfm:pzdr/a127}{2467}
\pdfglyphtounicode{tfm:pzdr/a128}{2468}
\pdfglyphtounicode{tfm:pzdr/a129}{2469}
\pdfglyphtounicode{tfm:pzdr/a13}{270C}
\pdfglyphtounicode{tfm:pzdr/a130}{2776}
\pdfglyphtounicode{tfm:pzdr/a131}{2777}
\pdfglyphtounicode{tfm:pzdr/a132}{2778}
\pdfglyphtounicode{tfm:pzdr/a133}{2779}
\pdfglyphtounicode{tfm:pzdr/a134}{277A}
\pdfglyphtounicode{tfm:pzdr/a135}{277B}
\pdfglyphtounicode{tfm:pzdr/a136}{277C}
\pdfglyphtounicode{tfm:pzdr/a137}{277D}
\pdfglyphtounicode{tfm:pzdr/a138}{277E}
\pdfglyphtounicode{tfm:pzdr/a139}{277F}
\pdfglyphtounicode{tfm:pzdr/a14}{270D}
\pdfglyphtounicode{tfm:pzdr/a140}{2780}
\pdfglyphtounicode{tfm:pzdr/a141}{2781}
\pdfglyphtounicode{tfm:pzdr/a142}{2782}
\pdfglyphtounicode{tfm:pzdr/a143}{2783}
\pdfglyphtounicode{tfm:pzdr/a144}{2784}
\pdfglyphtounicode{tfm:pzdr/a145}{2785}
\pdfglyphtounicode{tfm:pzdr/a146}{2786}
\pdfglyphtounicode{tfm:pzdr/a147}{2787}
\pdfglyphtounicode{tfm:pzdr/a148}{2788}
\pdfglyphtounicode{tfm:pzdr/a149}{2789}
\pdfglyphtounicode{tfm:pzdr/a15}{270E}
\pdfglyphtounicode{tfm:pzdr/a150}{278A}
\pdfglyphtounicode{tfm:pzdr/a151}{278B}
\pdfglyphtounicode{tfm:pzdr/a152}{278C}
\pdfglyphtounicode{tfm:pzdr/a153}{278D}
\pdfglyphtounicode{tfm:pzdr/a154}{278E}
\pdfglyphtounicode{tfm:pzdr/a155}{278F}
\pdfglyphtounicode{tfm:pzdr/a156}{2790}
\pdfglyphtounicode{tfm:pzdr/a157}{2791}
\pdfglyphtounicode{tfm:pzdr/a158}{2792}
\pdfglyphtounicode{tfm:pzdr/a159}{2793}
\pdfglyphtounicode{tfm:pzdr/a16}{270F}
\pdfglyphtounicode{tfm:pzdr/a160}{2794}
\pdfglyphtounicode{tfm:pzdr/a161}{2192}
\pdfglyphtounicode{tfm:pzdr/a162}{27A3}
\pdfglyphtounicode{tfm:pzdr/a163}{2194}
\pdfglyphtounicode{tfm:pzdr/a164}{2195}
\pdfglyphtounicode{tfm:pzdr/a165}{2799}
\pdfglyphtounicode{tfm:pzdr/a166}{279B}
\pdfglyphtounicode{tfm:pzdr/a167}{279C}
\pdfglyphtounicode{tfm:pzdr/a168}{279D}
\pdfglyphtounicode{tfm:pzdr/a169}{279E}
\pdfglyphtounicode{tfm:pzdr/a17}{2711}
\pdfglyphtounicode{tfm:pzdr/a170}{279F}
\pdfglyphtounicode{tfm:pzdr/a171}{27A0}
\pdfglyphtounicode{tfm:pzdr/a172}{27A1}
\pdfglyphtounicode{tfm:pzdr/a173}{27A2}
\pdfglyphtounicode{tfm:pzdr/a174}{27A4}
\pdfglyphtounicode{tfm:pzdr/a175}{27A5}
\pdfglyphtounicode{tfm:pzdr/a176}{27A6}
\pdfglyphtounicode{tfm:pzdr/a177}{27A7}
\pdfglyphtounicode{tfm:pzdr/a178}{27A8}
\pdfglyphtounicode{tfm:pzdr/a179}{27A9}
\pdfglyphtounicode{tfm:pzdr/a18}{2712}
\pdfglyphtounicode{tfm:pzdr/a180}{27AB}
\pdfglyphtounicode{tfm:pzdr/a181}{27AD}
\pdfglyphtounicode{tfm:pzdr/a182}{27AF}
\pdfglyphtounicode{tfm:pzdr/a183}{27B2}
\pdfglyphtounicode{tfm:pzdr/a184}{27B3}
\pdfglyphtounicode{tfm:pzdr/a185}{27B5}
\pdfglyphtounicode{tfm:pzdr/a186}{27B8}
\pdfglyphtounicode{tfm:pzdr/a187}{27BA}
\pdfglyphtounicode{tfm:pzdr/a188}{27BB}
\pdfglyphtounicode{tfm:pzdr/a189}{27BC}
\pdfglyphtounicode{tfm:pzdr/a19}{2713}
\pdfglyphtounicode{tfm:pzdr/a190}{27BD}
\pdfglyphtounicode{tfm:pzdr/a191}{27BE}
\pdfglyphtounicode{tfm:pzdr/a192}{279A}
\pdfglyphtounicode{tfm:pzdr/a193}{27AA}
\pdfglyphtounicode{tfm:pzdr/a194}{27B6}
\pdfglyphtounicode{tfm:pzdr/a195}{27B9}
\pdfglyphtounicode{tfm:pzdr/a196}{2798}
\pdfglyphtounicode{tfm:pzdr/a197}{27B4}
\pdfglyphtounicode{tfm:pzdr/a198}{27B7}
\pdfglyphtounicode{tfm:pzdr/a199}{27AC}
\pdfglyphtounicode{tfm:pzdr/a2}{2702}
\pdfglyphtounicode{tfm:pzdr/a20}{2714}
\pdfglyphtounicode{tfm:pzdr/a200}{27AE}
\pdfglyphtounicode{tfm:pzdr/a201}{27B1}
\pdfglyphtounicode{tfm:pzdr/a202}{2703}
\pdfglyphtounicode{tfm:pzdr/a203}{2750}
\pdfglyphtounicode{tfm:pzdr/a204}{2752}
\pdfglyphtounicode{tfm:pzdr/a205}{276E}
\pdfglyphtounicode{tfm:pzdr/a206}{2770}
\pdfglyphtounicode{tfm:pzdr/a21}{2715}
\pdfglyphtounicode{tfm:pzdr/a22}{2716}
\pdfglyphtounicode{tfm:pzdr/a23}{2717}
\pdfglyphtounicode{tfm:pzdr/a24}{2718}
\pdfglyphtounicode{tfm:pzdr/a25}{2719}
\pdfglyphtounicode{tfm:pzdr/a26}{271A}
\pdfglyphtounicode{tfm:pzdr/a27}{271B}
\pdfglyphtounicode{tfm:pzdr/a28}{271C}
\pdfglyphtounicode{tfm:pzdr/a29}{2722}
\pdfglyphtounicode{tfm:pzdr/a3}{2704}
\pdfglyphtounicode{tfm:pzdr/a30}{2723}
\pdfglyphtounicode{tfm:pzdr/a31}{2724}
\pdfglyphtounicode{tfm:pzdr/a32}{2725}
\pdfglyphtounicode{tfm:pzdr/a33}{2726}
\pdfglyphtounicode{tfm:pzdr/a34}{2727}
\pdfglyphtounicode{tfm:pzdr/a35}{2605}
\pdfglyphtounicode{tfm:pzdr/a36}{2729}
\pdfglyphtounicode{tfm:pzdr/a37}{272A}
\pdfglyphtounicode{tfm:pzdr/a38}{272B}
\pdfglyphtounicode{tfm:pzdr/a39}{272C}
\pdfglyphtounicode{tfm:pzdr/a4}{260E}
\pdfglyphtounicode{tfm:pzdr/a40}{272D}
\pdfglyphtounicode{tfm:pzdr/a41}{272E}
\pdfglyphtounicode{tfm:pzdr/a42}{272F}
\pdfglyphtounicode{tfm:pzdr/a43}{2730}
\pdfglyphtounicode{tfm:pzdr/a44}{2731}
\pdfglyphtounicode{tfm:pzdr/a45}{2732}
\pdfglyphtounicode{tfm:pzdr/a46}{2733}
\pdfglyphtounicode{tfm:pzdr/a47}{2734}
\pdfglyphtounicode{tfm:pzdr/a48}{2735}
\pdfglyphtounicode{tfm:pzdr/a49}{2736}
\pdfglyphtounicode{tfm:pzdr/a5}{2706}
\pdfglyphtounicode{tfm:pzdr/a50}{2737}
\pdfglyphtounicode{tfm:pzdr/a51}{2738}
\pdfglyphtounicode{tfm:pzdr/a52}{2739}
\pdfglyphtounicode{tfm:pzdr/a53}{273A}
\pdfglyphtounicode{tfm:pzdr/a54}{273B}
\pdfglyphtounicode{tfm:pzdr/a55}{273C}
\pdfglyphtounicode{tfm:pzdr/a56}{273D}
\pdfglyphtounicode{tfm:pzdr/a57}{273E}
\pdfglyphtounicode{tfm:pzdr/a58}{273F}
\pdfglyphtounicode{tfm:pzdr/a59}{2740}
\pdfglyphtounicode{tfm:pzdr/a6}{271D}
\pdfglyphtounicode{tfm:pzdr/a60}{2741}
\pdfglyphtounicode{tfm:pzdr/a61}{2742}
\pdfglyphtounicode{tfm:pzdr/a62}{2743}
\pdfglyphtounicode{tfm:pzdr/a63}{2744}
\pdfglyphtounicode{tfm:pzdr/a64}{2745}
\pdfglyphtounicode{tfm:pzdr/a65}{2746}
\pdfglyphtounicode{tfm:pzdr/a66}{2747}
\pdfglyphtounicode{tfm:pzdr/a67}{2748}
\pdfglyphtounicode{tfm:pzdr/a68}{2749}
\pdfglyphtounicode{tfm:pzdr/a69}{274A}
\pdfglyphtounicode{tfm:pzdr/a7}{271E}
\pdfglyphtounicode{tfm:pzdr/a70}{274B}
\pdfglyphtounicode{tfm:pzdr/a71}{25CF}
\pdfglyphtounicode{tfm:pzdr/a72}{274D}
\pdfglyphtounicode{tfm:pzdr/a73}{25A0}
\pdfglyphtounicode{tfm:pzdr/a74}{274F}
\pdfglyphtounicode{tfm:pzdr/a75}{2751}
\pdfglyphtounicode{tfm:pzdr/a76}{25B2}
\pdfglyphtounicode{tfm:pzdr/a77}{25BC}
\pdfglyphtounicode{tfm:pzdr/a78}{25C6}
\pdfglyphtounicode{tfm:pzdr/a79}{2756}
\pdfglyphtounicode{tfm:pzdr/a8}{271F}
\pdfglyphtounicode{tfm:pzdr/a81}{25D7}
\pdfglyphtounicode{tfm:pzdr/a82}{2758}
\pdfglyphtounicode{tfm:pzdr/a83}{2759}
\pdfglyphtounicode{tfm:pzdr/a84}{275A}
\pdfglyphtounicode{tfm:pzdr/a85}{276F}
\pdfglyphtounicode{tfm:pzdr/a86}{2771}
\pdfglyphtounicode{tfm:pzdr/a87}{2772}
\pdfglyphtounicode{tfm:pzdr/a88}{2773}
\pdfglyphtounicode{tfm:pzdr/a89}{2768}
\pdfglyphtounicode{tfm:pzdr/a9}{2720}
\pdfglyphtounicode{tfm:pzdr/a90}{2769}
\pdfglyphtounicode{tfm:pzdr/a91}{276C}
\pdfglyphtounicode{tfm:pzdr/a92}{276D}
\pdfglyphtounicode{tfm:pzdr/a93}{276A}
\pdfglyphtounicode{tfm:pzdr/a94}{276B}
\pdfglyphtounicode{tfm:pzdr/a95}{2774}
\pdfglyphtounicode{tfm:pzdr/a96}{2775}
\pdfglyphtounicode{tfm:pzdr/a97}{275B}
\pdfglyphtounicode{tfm:pzdr/a98}{275C}
\pdfglyphtounicode{tfm:pzdr/a99}{275D}
\pdfglyphtounicode{tfm:rpxbmi/phi}{03D5}
\pdfglyphtounicode{tfm:rpxbmi/phi1}{03C6}
\pdfglyphtounicode{tfm:rpxmi/phi}{03D5}
\pdfglyphtounicode{tfm:rpxmi/phi1}{03C6}
\pdfglyphtounicode{tfm:rpzdr/a1}{2701}
\pdfglyphtounicode{tfm:rpzdr/a10}{2721}
\pdfglyphtounicode{tfm:rpzdr/a100}{275E}
\pdfglyphtounicode{tfm:rpzdr/a101}{2761}
\pdfglyphtounicode{tfm:rpzdr/a102}{2762}
\pdfglyphtounicode{tfm:rpzdr/a103}{2763}
\pdfglyphtounicode{tfm:rpzdr/a104}{2764}
\pdfglyphtounicode{tfm:rpzdr/a105}{2710}
\pdfglyphtounicode{tfm:rpzdr/a106}{2765}
\pdfglyphtounicode{tfm:rpzdr/a107}{2766}
\pdfglyphtounicode{tfm:rpzdr/a108}{2767}
\pdfglyphtounicode{tfm:rpzdr/a109}{2660}
\pdfglyphtounicode{tfm:rpzdr/a11}{261B}
\pdfglyphtounicode{tfm:rpzdr/a110}{2665}
\pdfglyphtounicode{tfm:rpzdr/a111}{2666}
\pdfglyphtounicode{tfm:rpzdr/a112}{2663}
\pdfglyphtounicode{tfm:rpzdr/a117}{2709}
\pdfglyphtounicode{tfm:rpzdr/a118}{2708}
\pdfglyphtounicode{tfm:rpzdr/a119}{2707}
\pdfglyphtounicode{tfm:rpzdr/a12}{261E}
\pdfglyphtounicode{tfm:rpzdr/a120}{2460}
\pdfglyphtounicode{tfm:rpzdr/a121}{2461}
\pdfglyphtounicode{tfm:rpzdr/a122}{2462}
\pdfglyphtounicode{tfm:rpzdr/a123}{2463}
\pdfglyphtounicode{tfm:rpzdr/a124}{2464}
\pdfglyphtounicode{tfm:rpzdr/a125}{2465}
\pdfglyphtounicode{tfm:rpzdr/a126}{2466}
\pdfglyphtounicode{tfm:rpzdr/a127}{2467}
\pdfglyphtounicode{tfm:rpzdr/a128}{2468}
\pdfglyphtounicode{tfm:rpzdr/a129}{2469}
\pdfglyphtounicode{tfm:rpzdr/a13}{270C}
\pdfglyphtounicode{tfm:rpzdr/a130}{2776}
\pdfglyphtounicode{tfm:rpzdr/a131}{2777}
\pdfglyphtounicode{tfm:rpzdr/a132}{2778}
\pdfglyphtounicode{tfm:rpzdr/a133}{2779}
\pdfglyphtounicode{tfm:rpzdr/a134}{277A}
\pdfglyphtounicode{tfm:rpzdr/a135}{277B}
\pdfglyphtounicode{tfm:rpzdr/a136}{277C}
\pdfglyphtounicode{tfm:rpzdr/a137}{277D}
\pdfglyphtounicode{tfm:rpzdr/a138}{277E}
\pdfglyphtounicode{tfm:rpzdr/a139}{277F}
\pdfglyphtounicode{tfm:rpzdr/a14}{270D}
\pdfglyphtounicode{tfm:rpzdr/a140}{2780}
\pdfglyphtounicode{tfm:rpzdr/a141}{2781}
\pdfglyphtounicode{tfm:rpzdr/a142}{2782}
\pdfglyphtounicode{tfm:rpzdr/a143}{2783}
\pdfglyphtounicode{tfm:rpzdr/a144}{2784}
\pdfglyphtounicode{tfm:rpzdr/a145}{2785}
\pdfglyphtounicode{tfm:rpzdr/a146}{2786}
\pdfglyphtounicode{tfm:rpzdr/a147}{2787}
\pdfglyphtounicode{tfm:rpzdr/a148}{2788}
\pdfglyphtounicode{tfm:rpzdr/a149}{2789}
\pdfglyphtounicode{tfm:rpzdr/a15}{270E}
\pdfglyphtounicode{tfm:rpzdr/a150}{278A}
\pdfglyphtounicode{tfm:rpzdr/a151}{278B}
\pdfglyphtounicode{tfm:rpzdr/a152}{278C}
\pdfglyphtounicode{tfm:rpzdr/a153}{278D}
\pdfglyphtounicode{tfm:rpzdr/a154}{278E}
\pdfglyphtounicode{tfm:rpzdr/a155}{278F}
\pdfglyphtounicode{tfm:rpzdr/a156}{2790}
\pdfglyphtounicode{tfm:rpzdr/a157}{2791}
\pdfglyphtounicode{tfm:rpzdr/a158}{2792}
\pdfglyphtounicode{tfm:rpzdr/a159}{2793}
\pdfglyphtounicode{tfm:rpzdr/a16}{270F}
\pdfglyphtounicode{tfm:rpzdr/a160}{2794}
\pdfglyphtounicode{tfm:rpzdr/a161}{2192}
\pdfglyphtounicode{tfm:rpzdr/a162}{27A3}
\pdfglyphtounicode{tfm:rpzdr/a163}{2194}
\pdfglyphtounicode{tfm:rpzdr/a164}{2195}
\pdfglyphtounicode{tfm:rpzdr/a165}{2799}
\pdfglyphtounicode{tfm:rpzdr/a166}{279B}
\pdfglyphtounicode{tfm:rpzdr/a167}{279C}
\pdfglyphtounicode{tfm:rpzdr/a168}{279D}
\pdfglyphtounicode{tfm:rpzdr/a169}{279E}
\pdfglyphtounicode{tfm:rpzdr/a17}{2711}
\pdfglyphtounicode{tfm:rpzdr/a170}{279F}
\pdfglyphtounicode{tfm:rpzdr/a171}{27A0}
\pdfglyphtounicode{tfm:rpzdr/a172}{27A1}
\pdfglyphtounicode{tfm:rpzdr/a173}{27A2}
\pdfglyphtounicode{tfm:rpzdr/a174}{27A4}
\pdfglyphtounicode{tfm:rpzdr/a175}{27A5}
\pdfglyphtounicode{tfm:rpzdr/a176}{27A6}
\pdfglyphtounicode{tfm:rpzdr/a177}{27A7}
\pdfglyphtounicode{tfm:rpzdr/a178}{27A8}
\pdfglyphtounicode{tfm:rpzdr/a179}{27A9}
\pdfglyphtounicode{tfm:rpzdr/a18}{2712}
\pdfglyphtounicode{tfm:rpzdr/a180}{27AB}
\pdfglyphtounicode{tfm:rpzdr/a181}{27AD}
\pdfglyphtounicode{tfm:rpzdr/a182}{27AF}
\pdfglyphtounicode{tfm:rpzdr/a183}{27B2}
\pdfglyphtounicode{tfm:rpzdr/a184}{27B3}
\pdfglyphtounicode{tfm:rpzdr/a185}{27B5}
\pdfglyphtounicode{tfm:rpzdr/a186}{27B8}
\pdfglyphtounicode{tfm:rpzdr/a187}{27BA}
\pdfglyphtounicode{tfm:rpzdr/a188}{27BB}
\pdfglyphtounicode{tfm:rpzdr/a189}{27BC}
\pdfglyphtounicode{tfm:rpzdr/a19}{2713}
\pdfglyphtounicode{tfm:rpzdr/a190}{27BD}
\pdfglyphtounicode{tfm:rpzdr/a191}{27BE}
\pdfglyphtounicode{tfm:rpzdr/a192}{279A}
\pdfglyphtounicode{tfm:rpzdr/a193}{27AA}
\pdfglyphtounicode{tfm:rpzdr/a194}{27B6}
\pdfglyphtounicode{tfm:rpzdr/a195}{27B9}
\pdfglyphtounicode{tfm:rpzdr/a196}{2798}
\pdfglyphtounicode{tfm:rpzdr/a197}{27B4}
\pdfglyphtounicode{tfm:rpzdr/a198}{27B7}
\pdfglyphtounicode{tfm:rpzdr/a199}{27AC}
\pdfglyphtounicode{tfm:rpzdr/a2}{2702}
\pdfglyphtounicode{tfm:rpzdr/a20}{2714}
\pdfglyphtounicode{tfm:rpzdr/a200}{27AE}
\pdfglyphtounicode{tfm:rpzdr/a201}{27B1}
\pdfglyphtounicode{tfm:rpzdr/a202}{2703}
\pdfglyphtounicode{tfm:rpzdr/a203}{2750}
\pdfglyphtounicode{tfm:rpzdr/a204}{2752}
\pdfglyphtounicode{tfm:rpzdr/a205}{276E}
\pdfglyphtounicode{tfm:rpzdr/a206}{2770}
\pdfglyphtounicode{tfm:rpzdr/a21}{2715}
\pdfglyphtounicode{tfm:rpzdr/a22}{2716}
\pdfglyphtounicode{tfm:rpzdr/a23}{2717}
\pdfglyphtounicode{tfm:rpzdr/a24}{2718}
\pdfglyphtounicode{tfm:rpzdr/a25}{2719}
\pdfglyphtounicode{tfm:rpzdr/a26}{271A}
\pdfglyphtounicode{tfm:rpzdr/a27}{271B}
\pdfglyphtounicode{tfm:rpzdr/a28}{271C}
\pdfglyphtounicode{tfm:rpzdr/a29}{2722}
\pdfglyphtounicode{tfm:rpzdr/a3}{2704}
\pdfglyphtounicode{tfm:rpzdr/a30}{2723}
\pdfglyphtounicode{tfm:rpzdr/a31}{2724}
\pdfglyphtounicode{tfm:rpzdr/a32}{2725}
\pdfglyphtounicode{tfm:rpzdr/a33}{2726}
\pdfglyphtounicode{tfm:rpzdr/a34}{2727}
\pdfglyphtounicode{tfm:rpzdr/a35}{2605}
\pdfglyphtounicode{tfm:rpzdr/a36}{2729}
\pdfglyphtounicode{tfm:rpzdr/a37}{272A}
\pdfglyphtounicode{tfm:rpzdr/a38}{272B}
\pdfglyphtounicode{tfm:rpzdr/a39}{272C}
\pdfglyphtounicode{tfm:rpzdr/a4}{260E}
\pdfglyphtounicode{tfm:rpzdr/a40}{272D}
\pdfglyphtounicode{tfm:rpzdr/a41}{272E}
\pdfglyphtounicode{tfm:rpzdr/a42}{272F}
\pdfglyphtounicode{tfm:rpzdr/a43}{2730}
\pdfglyphtounicode{tfm:rpzdr/a44}{2731}
\pdfglyphtounicode{tfm:rpzdr/a45}{2732}
\pdfglyphtounicode{tfm:rpzdr/a46}{2733}
\pdfglyphtounicode{tfm:rpzdr/a47}{2734}
\pdfglyphtounicode{tfm:rpzdr/a48}{2735}
\pdfglyphtounicode{tfm:rpzdr/a49}{2736}
\pdfglyphtounicode{tfm:rpzdr/a5}{2706}
\pdfglyphtounicode{tfm:rpzdr/a50}{2737}
\pdfglyphtounicode{tfm:rpzdr/a51}{2738}
\pdfglyphtounicode{tfm:rpzdr/a52}{2739}
\pdfglyphtounicode{tfm:rpzdr/a53}{273A}
\pdfglyphtounicode{tfm:rpzdr/a54}{273B}
\pdfglyphtounicode{tfm:rpzdr/a55}{273C}
\pdfglyphtounicode{tfm:rpzdr/a56}{273D}
\pdfglyphtounicode{tfm:rpzdr/a57}{273E}
\pdfglyphtounicode{tfm:rpzdr/a58}{273F}
\pdfglyphtounicode{tfm:rpzdr/a59}{2740}
\pdfglyphtounicode{tfm:rpzdr/a6}{271D}
\pdfglyphtounicode{tfm:rpzdr/a60}{2741}
\pdfglyphtounicode{tfm:rpzdr/a61}{2742}
\pdfglyphtounicode{tfm:rpzdr/a62}{2743}
\pdfglyphtounicode{tfm:rpzdr/a63}{2744}
\pdfglyphtounicode{tfm:rpzdr/a64}{2745}
\pdfglyphtounicode{tfm:rpzdr/a65}{2746}
\pdfglyphtounicode{tfm:rpzdr/a66}{2747}
\pdfglyphtounicode{tfm:rpzdr/a67}{2748}
\pdfglyphtounicode{tfm:rpzdr/a68}{2749}
\pdfglyphtounicode{tfm:rpzdr/a69}{274A}
\pdfglyphtounicode{tfm:rpzdr/a7}{271E}
\pdfglyphtounicode{tfm:rpzdr/a70}{274B}
\pdfglyphtounicode{tfm:rpzdr/a71}{25CF}
\pdfglyphtounicode{tfm:rpzdr/a72}{274D}
\pdfglyphtounicode{tfm:rpzdr/a73}{25A0}
\pdfglyphtounicode{tfm:rpzdr/a74}{274F}
\pdfglyphtounicode{tfm:rpzdr/a75}{2751}
\pdfglyphtounicode{tfm:rpzdr/a76}{25B2}
\pdfglyphtounicode{tfm:rpzdr/a77}{25BC}
\pdfglyphtounicode{tfm:rpzdr/a78}{25C6}
\pdfglyphtounicode{tfm:rpzdr/a79}{2756}
\pdfglyphtounicode{tfm:rpzdr/a8}{271F}
\pdfglyphtounicode{tfm:rpzdr/a81}{25D7}
\pdfglyphtounicode{tfm:rpzdr/a82}{2758}
\pdfglyphtounicode{tfm:rpzdr/a83}{2759}
\pdfglyphtounicode{tfm:rpzdr/a84}{275A}
\pdfglyphtounicode{tfm:rpzdr/a85}{276F}
\pdfglyphtounicode{tfm:rpzdr/a86}{2771}
\pdfglyphtounicode{tfm:rpzdr/a87}{2772}
\pdfglyphtounicode{tfm:rpzdr/a88}{2773}
\pdfglyphtounicode{tfm:rpzdr/a89}{2768}
\pdfglyphtounicode{tfm:rpzdr/a9}{2720}
\pdfglyphtounicode{tfm:rpzdr/a90}{2769}
\pdfglyphtounicode{tfm:rpzdr/a91}{276C}
\pdfglyphtounicode{tfm:rpzdr/a92}{276D}
\pdfglyphtounicode{tfm:rpzdr/a93}{276A}
\pdfglyphtounicode{tfm:rpzdr/a94}{276B}
\pdfglyphtounicode{tfm:rpzdr/a95}{2774}
\pdfglyphtounicode{tfm:rpzdr/a96}{2775}
\pdfglyphtounicode{tfm:rpzdr/a97}{275B}
\pdfglyphtounicode{tfm:rpzdr/a98}{275C}
\pdfglyphtounicode{tfm:rpzdr/a99}{275D}
\pdfglyphtounicode{tfm:rtxbmi/phi}{03D5}
\pdfglyphtounicode{tfm:rtxbmi/phi1}{03C6}
\pdfglyphtounicode{tfm:rtxmi/phi}{03D5}
\pdfglyphtounicode{tfm:rtxmi/phi1}{03C6}
\pdfglyphtounicode{tfm:txbmia/phi}{03D5}
\pdfglyphtounicode{tfm:txbmia/phi1}{03C6}
\pdfglyphtounicode{tfm:txbsy/diamond}{2662}
\pdfglyphtounicode{tfm:txbsy/heart}{2661}
\pdfglyphtounicode{tfm:txbsya/diamond}{2662}
\pdfglyphtounicode{tfm:txmia/phi}{03D5}
\pdfglyphtounicode{tfm:txmia/phi1}{03C6}
\pdfglyphtounicode{tfm:txsy/diamond}{2662}
\pdfglyphtounicode{tfm:txsy/heart}{2661}
\pdfglyphtounicode{tfm:txsya/diamond}{2662}
\pdfglyphtounicode{tfm:zd/a1}{2701}
\pdfglyphtounicode{tfm:zd/a10}{2721}
\pdfglyphtounicode{tfm:zd/a100}{275E}
\pdfglyphtounicode{tfm:zd/a101}{2761}
\pdfglyphtounicode{tfm:zd/a102}{2762}
\pdfglyphtounicode{tfm:zd/a103}{2763}
\pdfglyphtounicode{tfm:zd/a104}{2764}
\pdfglyphtounicode{tfm:zd/a105}{2710}
\pdfglyphtounicode{tfm:zd/a106}{2765}
\pdfglyphtounicode{tfm:zd/a107}{2766}
\pdfglyphtounicode{tfm:zd/a108}{2767}
\pdfglyphtounicode{tfm:zd/a109}{2660}
\pdfglyphtounicode{tfm:zd/a11}{261B}
\pdfglyphtounicode{tfm:zd/a110}{2665}
\pdfglyphtounicode{tfm:zd/a111}{2666}
\pdfglyphtounicode{tfm:zd/a112}{2663}
\pdfglyphtounicode{tfm:zd/a117}{2709}
\pdfglyphtounicode{tfm:zd/a118}{2708}
\pdfglyphtounicode{tfm:zd/a119}{2707}
\pdfglyphtounicode{tfm:zd/a12}{261E}
\pdfglyphtounicode{tfm:zd/a120}{2460}
\pdfglyphtounicode{tfm:zd/a121}{2461}
\pdfglyphtounicode{tfm:zd/a122}{2462}
\pdfglyphtounicode{tfm:zd/a123}{2463}
\pdfglyphtounicode{tfm:zd/a124}{2464}
\pdfglyphtounicode{tfm:zd/a125}{2465}
\pdfglyphtounicode{tfm:zd/a126}{2466}
\pdfglyphtounicode{tfm:zd/a127}{2467}
\pdfglyphtounicode{tfm:zd/a128}{2468}
\pdfglyphtounicode{tfm:zd/a129}{2469}
\pdfglyphtounicode{tfm:zd/a13}{270C}
\pdfglyphtounicode{tfm:zd/a130}{2776}
\pdfglyphtounicode{tfm:zd/a131}{2777}
\pdfglyphtounicode{tfm:zd/a132}{2778}
\pdfglyphtounicode{tfm:zd/a133}{2779}
\pdfglyphtounicode{tfm:zd/a134}{277A}
\pdfglyphtounicode{tfm:zd/a135}{277B}
\pdfglyphtounicode{tfm:zd/a136}{277C}
\pdfglyphtounicode{tfm:zd/a137}{277D}
\pdfglyphtounicode{tfm:zd/a138}{277E}
\pdfglyphtounicode{tfm:zd/a139}{277F}
\pdfglyphtounicode{tfm:zd/a14}{270D}
\pdfglyphtounicode{tfm:zd/a140}{2780}
\pdfglyphtounicode{tfm:zd/a141}{2781}
\pdfglyphtounicode{tfm:zd/a142}{2782}
\pdfglyphtounicode{tfm:zd/a143}{2783}
\pdfglyphtounicode{tfm:zd/a144}{2784}
\pdfglyphtounicode{tfm:zd/a145}{2785}
\pdfglyphtounicode{tfm:zd/a146}{2786}
\pdfglyphtounicode{tfm:zd/a147}{2787}
\pdfglyphtounicode{tfm:zd/a148}{2788}
\pdfglyphtounicode{tfm:zd/a149}{2789}
\pdfglyphtounicode{tfm:zd/a15}{270E}
\pdfglyphtounicode{tfm:zd/a150}{278A}
\pdfglyphtounicode{tfm:zd/a151}{278B}
\pdfglyphtounicode{tfm:zd/a152}{278C}
\pdfglyphtounicode{tfm:zd/a153}{278D}
\pdfglyphtounicode{tfm:zd/a154}{278E}
\pdfglyphtounicode{tfm:zd/a155}{278F}
\pdfglyphtounicode{tfm:zd/a156}{2790}
\pdfglyphtounicode{tfm:zd/a157}{2791}
\pdfglyphtounicode{tfm:zd/a158}{2792}
\pdfglyphtounicode{tfm:zd/a159}{2793}
\pdfglyphtounicode{tfm:zd/a16}{270F}
\pdfglyphtounicode{tfm:zd/a160}{2794}
\pdfglyphtounicode{tfm:zd/a161}{2192}
\pdfglyphtounicode{tfm:zd/a162}{27A3}
\pdfglyphtounicode{tfm:zd/a163}{2194}
\pdfglyphtounicode{tfm:zd/a164}{2195}
\pdfglyphtounicode{tfm:zd/a165}{2799}
\pdfglyphtounicode{tfm:zd/a166}{279B}
\pdfglyphtounicode{tfm:zd/a167}{279C}
\pdfglyphtounicode{tfm:zd/a168}{279D}
\pdfglyphtounicode{tfm:zd/a169}{279E}
\pdfglyphtounicode{tfm:zd/a17}{2711}
\pdfglyphtounicode{tfm:zd/a170}{279F}
\pdfglyphtounicode{tfm:zd/a171}{27A0}
\pdfglyphtounicode{tfm:zd/a172}{27A1}
\pdfglyphtounicode{tfm:zd/a173}{27A2}
\pdfglyphtounicode{tfm:zd/a174}{27A4}
\pdfglyphtounicode{tfm:zd/a175}{27A5}
\pdfglyphtounicode{tfm:zd/a176}{27A6}
\pdfglyphtounicode{tfm:zd/a177}{27A7}
\pdfglyphtounicode{tfm:zd/a178}{27A8}
\pdfglyphtounicode{tfm:zd/a179}{27A9}
\pdfglyphtounicode{tfm:zd/a18}{2712}
\pdfglyphtounicode{tfm:zd/a180}{27AB}
\pdfglyphtounicode{tfm:zd/a181}{27AD}
\pdfglyphtounicode{tfm:zd/a182}{27AF}
\pdfglyphtounicode{tfm:zd/a183}{27B2}
\pdfglyphtounicode{tfm:zd/a184}{27B3}
\pdfglyphtounicode{tfm:zd/a185}{27B5}
\pdfglyphtounicode{tfm:zd/a186}{27B8}
\pdfglyphtounicode{tfm:zd/a187}{27BA}
\pdfglyphtounicode{tfm:zd/a188}{27BB}
\pdfglyphtounicode{tfm:zd/a189}{27BC}
\pdfglyphtounicode{tfm:zd/a19}{2713}
\pdfglyphtounicode{tfm:zd/a190}{27BD}
\pdfglyphtounicode{tfm:zd/a191}{27BE}
\pdfglyphtounicode{tfm:zd/a192}{279A}
\pdfglyphtounicode{tfm:zd/a193}{27AA}
\pdfglyphtounicode{tfm:zd/a194}{27B6}
\pdfglyphtounicode{tfm:zd/a195}{27B9}
\pdfglyphtounicode{tfm:zd/a196}{2798}
\pdfglyphtounicode{tfm:zd/a197}{27B4}
\pdfglyphtounicode{tfm:zd/a198}{27B7}
\pdfglyphtounicode{tfm:zd/a199}{27AC}
\pdfglyphtounicode{tfm:zd/a2}{2702}
\pdfglyphtounicode{tfm:zd/a20}{2714}
\pdfglyphtounicode{tfm:zd/a200}{27AE}
\pdfglyphtounicode{tfm:zd/a201}{27B1}
\pdfglyphtounicode{tfm:zd/a202}{2703}
\pdfglyphtounicode{tfm:zd/a203}{2750}
\pdfglyphtounicode{tfm:zd/a204}{2752}
\pdfglyphtounicode{tfm:zd/a205}{276E}
\pdfglyphtounicode{tfm:zd/a206}{2770}
\pdfglyphtounicode{tfm:zd/a21}{2715}
\pdfglyphtounicode{tfm:zd/a22}{2716}
\pdfglyphtounicode{tfm:zd/a23}{2717}
\pdfglyphtounicode{tfm:zd/a24}{2718}
\pdfglyphtounicode{tfm:zd/a25}{2719}
\pdfglyphtounicode{tfm:zd/a26}{271A}
\pdfglyphtounicode{tfm:zd/a27}{271B}
\pdfglyphtounicode{tfm:zd/a28}{271C}
\pdfglyphtounicode{tfm:zd/a29}{2722}
\pdfglyphtounicode{tfm:zd/a3}{2704}
\pdfglyphtounicode{tfm:zd/a30}{2723}
\pdfglyphtounicode{tfm:zd/a31}{2724}
\pdfglyphtounicode{tfm:zd/a32}{2725}
\pdfglyphtounicode{tfm:zd/a33}{2726}
\pdfglyphtounicode{tfm:zd/a34}{2727}
\pdfglyphtounicode{tfm:zd/a35}{2605}
\pdfglyphtounicode{tfm:zd/a36}{2729}
\pdfglyphtounicode{tfm:zd/a37}{272A}
\pdfglyphtounicode{tfm:zd/a38}{272B}
\pdfglyphtounicode{tfm:zd/a39}{272C}
\pdfglyphtounicode{tfm:zd/a4}{260E}
\pdfglyphtounicode{tfm:zd/a40}{272D}
\pdfglyphtounicode{tfm:zd/a41}{272E}
\pdfglyphtounicode{tfm:zd/a42}{272F}
\pdfglyphtounicode{tfm:zd/a43}{2730}
\pdfglyphtounicode{tfm:zd/a44}{2731}
\pdfglyphtounicode{tfm:zd/a45}{2732}
\pdfglyphtounicode{tfm:zd/a46}{2733}
\pdfglyphtounicode{tfm:zd/a47}{2734}
\pdfglyphtounicode{tfm:zd/a48}{2735}
\pdfglyphtounicode{tfm:zd/a49}{2736}
\pdfglyphtounicode{tfm:zd/a5}{2706}
\pdfglyphtounicode{tfm:zd/a50}{2737}
\pdfglyphtounicode{tfm:zd/a51}{2738}
\pdfglyphtounicode{tfm:zd/a52}{2739}
\pdfglyphtounicode{tfm:zd/a53}{273A}
\pdfglyphtounicode{tfm:zd/a54}{273B}
\pdfglyphtounicode{tfm:zd/a55}{273C}
\pdfglyphtounicode{tfm:zd/a56}{273D}
\pdfglyphtounicode{tfm:zd/a57}{273E}
\pdfglyphtounicode{tfm:zd/a58}{273F}
\pdfglyphtounicode{tfm:zd/a59}{2740}
\pdfglyphtounicode{tfm:zd/a6}{271D}
\pdfglyphtounicode{tfm:zd/a60}{2741}
\pdfglyphtounicode{tfm:zd/a61}{2742}
\pdfglyphtounicode{tfm:zd/a62}{2743}
\pdfglyphtounicode{tfm:zd/a63}{2744}
\pdfglyphtounicode{tfm:zd/a64}{2745}
\pdfglyphtounicode{tfm:zd/a65}{2746}
\pdfglyphtounicode{tfm:zd/a66}{2747}
\pdfglyphtounicode{tfm:zd/a67}{2748}
\pdfglyphtounicode{tfm:zd/a68}{2749}
\pdfglyphtounicode{tfm:zd/a69}{274A}
\pdfglyphtounicode{tfm:zd/a7}{271E}
\pdfglyphtounicode{tfm:zd/a70}{274B}
\pdfglyphtounicode{tfm:zd/a71}{25CF}
\pdfglyphtounicode{tfm:zd/a72}{274D}
\pdfglyphtounicode{tfm:zd/a73}{25A0}
\pdfglyphtounicode{tfm:zd/a74}{274F}
\pdfglyphtounicode{tfm:zd/a75}{2751}
\pdfglyphtounicode{tfm:zd/a76}{25B2}
\pdfglyphtounicode{tfm:zd/a77}{25BC}
\pdfglyphtounicode{tfm:zd/a78}{25C6}
\pdfglyphtounicode{tfm:zd/a79}{2756}
\pdfglyphtounicode{tfm:zd/a8}{271F}
\pdfglyphtounicode{tfm:zd/a81}{25D7}
\pdfglyphtounicode{tfm:zd/a82}{2758}
\pdfglyphtounicode{tfm:zd/a83}{2759}
\pdfglyphtounicode{tfm:zd/a84}{275A}
\pdfglyphtounicode{tfm:zd/a85}{276F}
\pdfglyphtounicode{tfm:zd/a86}{2771}
\pdfglyphtounicode{tfm:zd/a87}{2772}
\pdfglyphtounicode{tfm:zd/a88}{2773}
\pdfglyphtounicode{tfm:zd/a89}{2768}
\pdfglyphtounicode{tfm:zd/a9}{2720}
\pdfglyphtounicode{tfm:zd/a90}{2769}
\pdfglyphtounicode{tfm:zd/a91}{276C}
\pdfglyphtounicode{tfm:zd/a92}{276D}
\pdfglyphtounicode{tfm:zd/a93}{276A}
\pdfglyphtounicode{tfm:zd/a94}{276B}
\pdfglyphtounicode{tfm:zd/a95}{2774}
\pdfglyphtounicode{tfm:zd/a96}{2775}
\pdfglyphtounicode{tfm:zd/a97}{275B}
\pdfglyphtounicode{tfm:zd/a98}{275C}
\pdfglyphtounicode{tfm:zd/a99}{275D}
\pdfglyphtounicode{tfm:zpzdr-reversed/a1}{2701}
\pdfglyphtounicode{tfm:zpzdr-reversed/a10}{2721}
\pdfglyphtounicode{tfm:zpzdr-reversed/a100}{275E}
\pdfglyphtounicode{tfm:zpzdr-reversed/a101}{2761}
\pdfglyphtounicode{tfm:zpzdr-reversed/a102}{2762}
\pdfglyphtounicode{tfm:zpzdr-reversed/a103}{2763}
\pdfglyphtounicode{tfm:zpzdr-reversed/a104}{2764}
\pdfglyphtounicode{tfm:zpzdr-reversed/a105}{2710}
\pdfglyphtounicode{tfm:zpzdr-reversed/a106}{2765}
\pdfglyphtounicode{tfm:zpzdr-reversed/a107}{2766}
\pdfglyphtounicode{tfm:zpzdr-reversed/a108}{2767}
\pdfglyphtounicode{tfm:zpzdr-reversed/a109}{2660}
\pdfglyphtounicode{tfm:zpzdr-reversed/a11}{261B}
\pdfglyphtounicode{tfm:zpzdr-reversed/a110}{2665}
\pdfglyphtounicode{tfm:zpzdr-reversed/a111}{2666}
\pdfglyphtounicode{tfm:zpzdr-reversed/a112}{2663}
\pdfglyphtounicode{tfm:zpzdr-reversed/a117}{2709}
\pdfglyphtounicode{tfm:zpzdr-reversed/a118}{2708}
\pdfglyphtounicode{tfm:zpzdr-reversed/a119}{2707}
\pdfglyphtounicode{tfm:zpzdr-reversed/a12}{261E}
\pdfglyphtounicode{tfm:zpzdr-reversed/a120}{2460}
\pdfglyphtounicode{tfm:zpzdr-reversed/a121}{2461}
\pdfglyphtounicode{tfm:zpzdr-reversed/a122}{2462}
\pdfglyphtounicode{tfm:zpzdr-reversed/a123}{2463}
\pdfglyphtounicode{tfm:zpzdr-reversed/a124}{2464}
\pdfglyphtounicode{tfm:zpzdr-reversed/a125}{2465}
\pdfglyphtounicode{tfm:zpzdr-reversed/a126}{2466}
\pdfglyphtounicode{tfm:zpzdr-reversed/a127}{2467}
\pdfglyphtounicode{tfm:zpzdr-reversed/a128}{2468}
\pdfglyphtounicode{tfm:zpzdr-reversed/a129}{2469}
\pdfglyphtounicode{tfm:zpzdr-reversed/a13}{270C}
\pdfglyphtounicode{tfm:zpzdr-reversed/a130}{2776}
\pdfglyphtounicode{tfm:zpzdr-reversed/a131}{2777}
\pdfglyphtounicode{tfm:zpzdr-reversed/a132}{2778}
\pdfglyphtounicode{tfm:zpzdr-reversed/a133}{2779}
\pdfglyphtounicode{tfm:zpzdr-reversed/a134}{277A}
\pdfglyphtounicode{tfm:zpzdr-reversed/a135}{277B}
\pdfglyphtounicode{tfm:zpzdr-reversed/a136}{277C}
\pdfglyphtounicode{tfm:zpzdr-reversed/a137}{277D}
\pdfglyphtounicode{tfm:zpzdr-reversed/a138}{277E}
\pdfglyphtounicode{tfm:zpzdr-reversed/a139}{277F}
\pdfglyphtounicode{tfm:zpzdr-reversed/a14}{270D}
\pdfglyphtounicode{tfm:zpzdr-reversed/a140}{2780}
\pdfglyphtounicode{tfm:zpzdr-reversed/a141}{2781}
\pdfglyphtounicode{tfm:zpzdr-reversed/a142}{2782}
\pdfglyphtounicode{tfm:zpzdr-reversed/a143}{2783}
\pdfglyphtounicode{tfm:zpzdr-reversed/a144}{2784}
\pdfglyphtounicode{tfm:zpzdr-reversed/a145}{2785}
\pdfglyphtounicode{tfm:zpzdr-reversed/a146}{2786}
\pdfglyphtounicode{tfm:zpzdr-reversed/a147}{2787}
\pdfglyphtounicode{tfm:zpzdr-reversed/a148}{2788}
\pdfglyphtounicode{tfm:zpzdr-reversed/a149}{2789}
\pdfglyphtounicode{tfm:zpzdr-reversed/a15}{270E}
\pdfglyphtounicode{tfm:zpzdr-reversed/a150}{278A}
\pdfglyphtounicode{tfm:zpzdr-reversed/a151}{278B}
\pdfglyphtounicode{tfm:zpzdr-reversed/a152}{278C}
\pdfglyphtounicode{tfm:zpzdr-reversed/a153}{278D}
\pdfglyphtounicode{tfm:zpzdr-reversed/a154}{278E}
\pdfglyphtounicode{tfm:zpzdr-reversed/a155}{278F}
\pdfglyphtounicode{tfm:zpzdr-reversed/a156}{2790}
\pdfglyphtounicode{tfm:zpzdr-reversed/a157}{2791}
\pdfglyphtounicode{tfm:zpzdr-reversed/a158}{2792}
\pdfglyphtounicode{tfm:zpzdr-reversed/a159}{2793}
\pdfglyphtounicode{tfm:zpzdr-reversed/a16}{270F}
\pdfglyphtounicode{tfm:zpzdr-reversed/a160}{2794}
\pdfglyphtounicode{tfm:zpzdr-reversed/a161}{2192}
\pdfglyphtounicode{tfm:zpzdr-reversed/a162}{27A3}
\pdfglyphtounicode{tfm:zpzdr-reversed/a163}{2194}
\pdfglyphtounicode{tfm:zpzdr-reversed/a164}{2195}
\pdfglyphtounicode{tfm:zpzdr-reversed/a165}{2799}
\pdfglyphtounicode{tfm:zpzdr-reversed/a166}{279B}
\pdfglyphtounicode{tfm:zpzdr-reversed/a167}{279C}
\pdfglyphtounicode{tfm:zpzdr-reversed/a168}{279D}
\pdfglyphtounicode{tfm:zpzdr-reversed/a169}{279E}
\pdfglyphtounicode{tfm:zpzdr-reversed/a17}{2711}
\pdfglyphtounicode{tfm:zpzdr-reversed/a170}{279F}
\pdfglyphtounicode{tfm:zpzdr-reversed/a171}{27A0}
\pdfglyphtounicode{tfm:zpzdr-reversed/a172}{27A1}
\pdfglyphtounicode{tfm:zpzdr-reversed/a173}{27A2}
\pdfglyphtounicode{tfm:zpzdr-reversed/a174}{27A4}
\pdfglyphtounicode{tfm:zpzdr-reversed/a175}{27A5}
\pdfglyphtounicode{tfm:zpzdr-reversed/a176}{27A6}
\pdfglyphtounicode{tfm:zpzdr-reversed/a177}{27A7}
\pdfglyphtounicode{tfm:zpzdr-reversed/a178}{27A8}
\pdfglyphtounicode{tfm:zpzdr-reversed/a179}{27A9}
\pdfglyphtounicode{tfm:zpzdr-reversed/a18}{2712}
\pdfglyphtounicode{tfm:zpzdr-reversed/a180}{27AB}
\pdfglyphtounicode{tfm:zpzdr-reversed/a181}{27AD}
\pdfglyphtounicode{tfm:zpzdr-reversed/a182}{27AF}
\pdfglyphtounicode{tfm:zpzdr-reversed/a183}{27B2}
\pdfglyphtounicode{tfm:zpzdr-reversed/a184}{27B3}
\pdfglyphtounicode{tfm:zpzdr-reversed/a185}{27B5}
\pdfglyphtounicode{tfm:zpzdr-reversed/a186}{27B8}
\pdfglyphtounicode{tfm:zpzdr-reversed/a187}{27BA}
\pdfglyphtounicode{tfm:zpzdr-reversed/a188}{27BB}
\pdfglyphtounicode{tfm:zpzdr-reversed/a189}{27BC}
\pdfglyphtounicode{tfm:zpzdr-reversed/a19}{2713}
\pdfglyphtounicode{tfm:zpzdr-reversed/a190}{27BD}
\pdfglyphtounicode{tfm:zpzdr-reversed/a191}{27BE}
\pdfglyphtounicode{tfm:zpzdr-reversed/a192}{279A}
\pdfglyphtounicode{tfm:zpzdr-reversed/a193}{27AA}
\pdfglyphtounicode{tfm:zpzdr-reversed/a194}{27B6}
\pdfglyphtounicode{tfm:zpzdr-reversed/a195}{27B9}
\pdfglyphtounicode{tfm:zpzdr-reversed/a196}{2798}
\pdfglyphtounicode{tfm:zpzdr-reversed/a197}{27B4}
\pdfglyphtounicode{tfm:zpzdr-reversed/a198}{27B7}
\pdfglyphtounicode{tfm:zpzdr-reversed/a199}{27AC}
\pdfglyphtounicode{tfm:zpzdr-reversed/a2}{2702}
\pdfglyphtounicode{tfm:zpzdr-reversed/a20}{2714}
\pdfglyphtounicode{tfm:zpzdr-reversed/a200}{27AE}
\pdfglyphtounicode{tfm:zpzdr-reversed/a201}{27B1}
\pdfglyphtounicode{tfm:zpzdr-reversed/a202}{2703}
\pdfglyphtounicode{tfm:zpzdr-reversed/a203}{2750}
\pdfglyphtounicode{tfm:zpzdr-reversed/a204}{2752}
\pdfglyphtounicode{tfm:zpzdr-reversed/a205}{276E}
\pdfglyphtounicode{tfm:zpzdr-reversed/a206}{2770}
\pdfglyphtounicode{tfm:zpzdr-reversed/a21}{2715}
\pdfglyphtounicode{tfm:zpzdr-reversed/a22}{2716}
\pdfglyphtounicode{tfm:zpzdr-reversed/a23}{2717}
\pdfglyphtounicode{tfm:zpzdr-reversed/a24}{2718}
\pdfglyphtounicode{tfm:zpzdr-reversed/a25}{2719}
\pdfglyphtounicode{tfm:zpzdr-reversed/a26}{271A}
\pdfglyphtounicode{tfm:zpzdr-reversed/a27}{271B}
\pdfglyphtounicode{tfm:zpzdr-reversed/a28}{271C}
\pdfglyphtounicode{tfm:zpzdr-reversed/a29}{2722}
\pdfglyphtounicode{tfm:zpzdr-reversed/a3}{2704}
\pdfglyphtounicode{tfm:zpzdr-reversed/a30}{2723}
\pdfglyphtounicode{tfm:zpzdr-reversed/a31}{2724}
\pdfglyphtounicode{tfm:zpzdr-reversed/a32}{2725}
\pdfglyphtounicode{tfm:zpzdr-reversed/a33}{2726}
\pdfglyphtounicode{tfm:zpzdr-reversed/a34}{2727}
\pdfglyphtounicode{tfm:zpzdr-reversed/a35}{2605}
\pdfglyphtounicode{tfm:zpzdr-reversed/a36}{2729}
\pdfglyphtounicode{tfm:zpzdr-reversed/a37}{272A}
\pdfglyphtounicode{tfm:zpzdr-reversed/a38}{272B}
\pdfglyphtounicode{tfm:zpzdr-reversed/a39}{272C}
\pdfglyphtounicode{tfm:zpzdr-reversed/a4}{260E}
\pdfglyphtounicode{tfm:zpzdr-reversed/a40}{272D}
\pdfglyphtounicode{tfm:zpzdr-reversed/a41}{272E}
\pdfglyphtounicode{tfm:zpzdr-reversed/a42}{272F}
\pdfglyphtounicode{tfm:zpzdr-reversed/a43}{2730}
\pdfglyphtounicode{tfm:zpzdr-reversed/a44}{2731}
\pdfglyphtounicode{tfm:zpzdr-reversed/a45}{2732}
\pdfglyphtounicode{tfm:zpzdr-reversed/a46}{2733}
\pdfglyphtounicode{tfm:zpzdr-reversed/a47}{2734}
\pdfglyphtounicode{tfm:zpzdr-reversed/a48}{2735}
\pdfglyphtounicode{tfm:zpzdr-reversed/a49}{2736}
\pdfglyphtounicode{tfm:zpzdr-reversed/a5}{2706}
\pdfglyphtounicode{tfm:zpzdr-reversed/a50}{2737}
\pdfglyphtounicode{tfm:zpzdr-reversed/a51}{2738}
\pdfglyphtounicode{tfm:zpzdr-reversed/a52}{2739}
\pdfglyphtounicode{tfm:zpzdr-reversed/a53}{273A}
\pdfglyphtounicode{tfm:zpzdr-reversed/a54}{273B}
\pdfglyphtounicode{tfm:zpzdr-reversed/a55}{273C}
\pdfglyphtounicode{tfm:zpzdr-reversed/a56}{273D}
\pdfglyphtounicode{tfm:zpzdr-reversed/a57}{273E}
\pdfglyphtounicode{tfm:zpzdr-reversed/a58}{273F}
\pdfglyphtounicode{tfm:zpzdr-reversed/a59}{2740}
\pdfglyphtounicode{tfm:zpzdr-reversed/a6}{271D}
\pdfglyphtounicode{tfm:zpzdr-reversed/a60}{2741}
\pdfglyphtounicode{tfm:zpzdr-reversed/a61}{2742}
\pdfglyphtounicode{tfm:zpzdr-reversed/a62}{2743}
\pdfglyphtounicode{tfm:zpzdr-reversed/a63}{2744}
\pdfglyphtounicode{tfm:zpzdr-reversed/a64}{2745}
\pdfglyphtounicode{tfm:zpzdr-reversed/a65}{2746}
\pdfglyphtounicode{tfm:zpzdr-reversed/a66}{2747}
\pdfglyphtounicode{tfm:zpzdr-reversed/a67}{2748}
\pdfglyphtounicode{tfm:zpzdr-reversed/a68}{2749}
\pdfglyphtounicode{tfm:zpzdr-reversed/a69}{274A}
\pdfglyphtounicode{tfm:zpzdr-reversed/a7}{271E}
\pdfglyphtounicode{tfm:zpzdr-reversed/a70}{274B}
\pdfglyphtounicode{tfm:zpzdr-reversed/a71}{25CF}
\pdfglyphtounicode{tfm:zpzdr-reversed/a72}{274D}
\pdfglyphtounicode{tfm:zpzdr-reversed/a73}{25A0}
\pdfglyphtounicode{tfm:zpzdr-reversed/a74}{274F}
\pdfglyphtounicode{tfm:zpzdr-reversed/a75}{2751}
\pdfglyphtounicode{tfm:zpzdr-reversed/a76}{25B2}
\pdfglyphtounicode{tfm:zpzdr-reversed/a77}{25BC}
\pdfglyphtounicode{tfm:zpzdr-reversed/a78}{25C6}
\pdfglyphtounicode{tfm:zpzdr-reversed/a79}{2756}
\pdfglyphtounicode{tfm:zpzdr-reversed/a8}{271F}
\pdfglyphtounicode{tfm:zpzdr-reversed/a81}{25D7}
\pdfglyphtounicode{tfm:zpzdr-reversed/a82}{2758}
\pdfglyphtounicode{tfm:zpzdr-reversed/a83}{2759}
\pdfglyphtounicode{tfm:zpzdr-reversed/a84}{275A}
\pdfglyphtounicode{tfm:zpzdr-reversed/a85}{276F}
\pdfglyphtounicode{tfm:zpzdr-reversed/a86}{2771}
\pdfglyphtounicode{tfm:zpzdr-reversed/a87}{2772}
\pdfglyphtounicode{tfm:zpzdr-reversed/a88}{2773}
\pdfglyphtounicode{tfm:zpzdr-reversed/a89}{2768}
\pdfglyphtounicode{tfm:zpzdr-reversed/a9}{2720}
\pdfglyphtounicode{tfm:zpzdr-reversed/a90}{2769}
\pdfglyphtounicode{tfm:zpzdr-reversed/a91}{276C}
\pdfglyphtounicode{tfm:zpzdr-reversed/a92}{276D}
\pdfglyphtounicode{tfm:zpzdr-reversed/a93}{276A}
\pdfglyphtounicode{tfm:zpzdr-reversed/a94}{276B}
\pdfglyphtounicode{tfm:zpzdr-reversed/a95}{2774}
\pdfglyphtounicode{tfm:zpzdr-reversed/a96}{2775}
\pdfglyphtounicode{tfm:zpzdr-reversed/a97}{275B}
\pdfglyphtounicode{tfm:zpzdr-reversed/a98}{275C}
\pdfglyphtounicode{tfm:zpzdr-reversed/a99}{275D}
\pdfglyphtounicode{thabengali}{09A5}
\pdfglyphtounicode{thadeva}{0925}
\pdfglyphtounicode{thagujarati}{0AA5}
\pdfglyphtounicode{thagurmukhi}{0A25}
\pdfglyphtounicode{thalarabic}{0630}
\pdfglyphtounicode{thalfinalarabic}{FEAC}
\pdfglyphtounicode{thanthakhatlowleftthai}{F898}
\pdfglyphtounicode{thanthakhatlowrightthai}{F897}
\pdfglyphtounicode{thanthakhatthai}{0E4C}
\pdfglyphtounicode{thanthakhatupperleftthai}{F896}
\pdfglyphtounicode{theharabic}{062B}
\pdfglyphtounicode{thehfinalarabic}{FE9A}
\pdfglyphtounicode{thehinitialarabic}{FE9B}
\pdfglyphtounicode{thehmedialarabic}{FE9C}
\pdfglyphtounicode{thereexists}{2203}
\pdfglyphtounicode{therefore}{2234}
\pdfglyphtounicode{theta}{03B8}
\pdfglyphtounicode{theta1}{03D1}
\pdfglyphtounicode{thetasymbolgreek}{03D1}
\pdfglyphtounicode{thieuthacirclekorean}{3279}
\pdfglyphtounicode{thieuthaparenkorean}{3219}
\pdfglyphtounicode{thieuthcirclekorean}{326B}
\pdfglyphtounicode{thieuthkorean}{314C}
\pdfglyphtounicode{thieuthparenkorean}{320B}
\pdfglyphtounicode{thirteencircle}{246C}
\pdfglyphtounicode{thirteenparen}{2480}
\pdfglyphtounicode{thirteenperiod}{2494}
\pdfglyphtounicode{thonangmonthothai}{0E11}
\pdfglyphtounicode{thook}{01AD}
\pdfglyphtounicode{thophuthaothai}{0E12}
\pdfglyphtounicode{thorn}{00FE}
\pdfglyphtounicode{thothahanthai}{0E17}
\pdfglyphtounicode{thothanthai}{0E10}
\pdfglyphtounicode{thothongthai}{0E18}
\pdfglyphtounicode{thothungthai}{0E16}
\pdfglyphtounicode{thousandcyrillic}{0482}
\pdfglyphtounicode{thousandsseparatorarabic}{066C}
\pdfglyphtounicode{thousandsseparatorpersian}{066C}
\pdfglyphtounicode{three}{0033}
\pdfglyphtounicode{threearabic}{0663}
\pdfglyphtounicode{threebengali}{09E9}
\pdfglyphtounicode{threecircle}{2462}
\pdfglyphtounicode{threecircleinversesansserif}{278C}
\pdfglyphtounicode{threedeva}{0969}
\pdfglyphtounicode{threeeighths}{215C}
\pdfglyphtounicode{threegujarati}{0AE9}
\pdfglyphtounicode{threegurmukhi}{0A69}
\pdfglyphtounicode{threehackarabic}{0663}
\pdfglyphtounicode{threehangzhou}{3023}
\pdfglyphtounicode{threeideographicparen}{3222}
\pdfglyphtounicode{threeinferior}{2083}
\pdfglyphtounicode{threemonospace}{FF13}
\pdfglyphtounicode{threenumeratorbengali}{09F6}
\pdfglyphtounicode{threeoldstyle}{0033}
\pdfglyphtounicode{threeparen}{2476}
\pdfglyphtounicode{threeperiod}{248A}
\pdfglyphtounicode{threepersian}{06F3}
\pdfglyphtounicode{threequarters}{00BE}
\pdfglyphtounicode{threequartersemdash}{F6DE}
\pdfglyphtounicode{threeroman}{2172}
\pdfglyphtounicode{threesuperior}{00B3}
\pdfglyphtounicode{threethai}{0E53}
\pdfglyphtounicode{thzsquare}{3394}
\pdfglyphtounicode{tihiragana}{3061}
\pdfglyphtounicode{tikatakana}{30C1}
\pdfglyphtounicode{tikatakanahalfwidth}{FF81}
\pdfglyphtounicode{tikeutacirclekorean}{3270}
\pdfglyphtounicode{tikeutaparenkorean}{3210}
\pdfglyphtounicode{tikeutcirclekorean}{3262}
\pdfglyphtounicode{tikeutkorean}{3137}
\pdfglyphtounicode{tikeutparenkorean}{3202}
\pdfglyphtounicode{tilde}{02DC}
\pdfglyphtounicode{tildebelowcmb}{0330}
\pdfglyphtounicode{tildecmb}{0303}
\pdfglyphtounicode{tildecomb}{0303}
\pdfglyphtounicode{tildedoublecmb}{0360}
\pdfglyphtounicode{tildeoperator}{223C}
\pdfglyphtounicode{tildeoverlaycmb}{0334}
\pdfglyphtounicode{tildeverticalcmb}{033E}
\pdfglyphtounicode{timescircle}{2297}
\pdfglyphtounicode{tipehahebrew}{0596}
\pdfglyphtounicode{tipehalefthebrew}{0596}
\pdfglyphtounicode{tippigurmukhi}{0A70}
\pdfglyphtounicode{titlocyrilliccmb}{0483}
\pdfglyphtounicode{tiwnarmenian}{057F}
\pdfglyphtounicode{tlinebelow}{1E6F}
\pdfglyphtounicode{tmonospace}{FF54}
\pdfglyphtounicode{toarmenian}{0569}
\pdfglyphtounicode{tohiragana}{3068}
\pdfglyphtounicode{tokatakana}{30C8}
\pdfglyphtounicode{tokatakanahalfwidth}{FF84}
\pdfglyphtounicode{tonebarextrahighmod}{02E5}
\pdfglyphtounicode{tonebarextralowmod}{02E9}
\pdfglyphtounicode{tonebarhighmod}{02E6}
\pdfglyphtounicode{tonebarlowmod}{02E8}
\pdfglyphtounicode{tonebarmidmod}{02E7}
\pdfglyphtounicode{tonefive}{01BD}
\pdfglyphtounicode{tonesix}{0185}
\pdfglyphtounicode{tonetwo}{01A8}
\pdfglyphtounicode{tonos}{0384}
\pdfglyphtounicode{tonsquare}{3327}
\pdfglyphtounicode{topatakthai}{0E0F}
\pdfglyphtounicode{tortoiseshellbracketleft}{3014}
\pdfglyphtounicode{tortoiseshellbracketleftsmall}{FE5D}
\pdfglyphtounicode{tortoiseshellbracketleftvertical}{FE39}
\pdfglyphtounicode{tortoiseshellbracketright}{3015}
\pdfglyphtounicode{tortoiseshellbracketrightsmall}{FE5E}
\pdfglyphtounicode{tortoiseshellbracketrightvertical}{FE3A}
\pdfglyphtounicode{totaothai}{0E15}
\pdfglyphtounicode{tpalatalhook}{01AB}
\pdfglyphtounicode{tparen}{24AF}
\pdfglyphtounicode{trademark}{2122}
\pdfglyphtounicode{trademarksans}{2122}
\pdfglyphtounicode{trademarkserif}{2122}
\pdfglyphtounicode{tretroflexhook}{0288}
\pdfglyphtounicode{triagdn}{25BC}
\pdfglyphtounicode{triaglf}{25C4}
\pdfglyphtounicode{triagrt}{25BA}
\pdfglyphtounicode{triagup}{25B2}
\pdfglyphtounicode{triangle}{25B3}
\pdfglyphtounicode{triangledownsld}{25BC}
\pdfglyphtounicode{triangleinv}{25BD}
\pdfglyphtounicode{triangleleft}{25C1}
\pdfglyphtounicode{triangleleftequal}{22B4}
\pdfglyphtounicode{triangleleftsld}{25C0}
\pdfglyphtounicode{triangleright}{25B7}
\pdfglyphtounicode{trianglerightequal}{22B5}
\pdfglyphtounicode{trianglerightsld}{25B6}
\pdfglyphtounicode{trianglesolid}{25B2}
\pdfglyphtounicode{ts}{02A6}
\pdfglyphtounicode{tsadi}{05E6}
\pdfglyphtounicode{tsadidagesh}{FB46}
\pdfglyphtounicode{tsadidageshhebrew}{FB46}
\pdfglyphtounicode{tsadihebrew}{05E6}
\pdfglyphtounicode{tsecyrillic}{0446}
\pdfglyphtounicode{tsere}{05B5}
\pdfglyphtounicode{tsere12}{05B5}
\pdfglyphtounicode{tsere1e}{05B5}
\pdfglyphtounicode{tsere2b}{05B5}
\pdfglyphtounicode{tserehebrew}{05B5}
\pdfglyphtounicode{tserenarrowhebrew}{05B5}
\pdfglyphtounicode{tserequarterhebrew}{05B5}
\pdfglyphtounicode{tserewidehebrew}{05B5}
\pdfglyphtounicode{tshecyrillic}{045B}
\pdfglyphtounicode{tsuperior}{0074}
\pdfglyphtounicode{ttabengali}{099F}
\pdfglyphtounicode{ttadeva}{091F}
\pdfglyphtounicode{ttagujarati}{0A9F}
\pdfglyphtounicode{ttagurmukhi}{0A1F}
\pdfglyphtounicode{tteharabic}{0679}
\pdfglyphtounicode{ttehfinalarabic}{FB67}
\pdfglyphtounicode{ttehinitialarabic}{FB68}
\pdfglyphtounicode{ttehmedialarabic}{FB69}
\pdfglyphtounicode{tthabengali}{09A0}
\pdfglyphtounicode{tthadeva}{0920}
\pdfglyphtounicode{tthagujarati}{0AA0}
\pdfglyphtounicode{tthagurmukhi}{0A20}
\pdfglyphtounicode{tturned}{0287}
\pdfglyphtounicode{tuhiragana}{3064}
\pdfglyphtounicode{tukatakana}{30C4}
\pdfglyphtounicode{tukatakanahalfwidth}{FF82}
\pdfglyphtounicode{turnstileleft}{22A2}
\pdfglyphtounicode{turnstileright}{22A3}
\pdfglyphtounicode{tusmallhiragana}{3063}
\pdfglyphtounicode{tusmallkatakana}{30C3}
\pdfglyphtounicode{tusmallkatakanahalfwidth}{FF6F}
\pdfglyphtounicode{twelvecircle}{246B}
\pdfglyphtounicode{twelveparen}{247F}
\pdfglyphtounicode{twelveperiod}{2493}
\pdfglyphtounicode{twelveroman}{217B}
\pdfglyphtounicode{twentycircle}{2473}
\pdfglyphtounicode{twentyhangzhou}{5344}
\pdfglyphtounicode{twentyparen}{2487}
\pdfglyphtounicode{twentyperiod}{249B}
\pdfglyphtounicode{two}{0032}
\pdfglyphtounicode{twoarabic}{0662}
\pdfglyphtounicode{twobengali}{09E8}
\pdfglyphtounicode{twocircle}{2461}
\pdfglyphtounicode{twocircleinversesansserif}{278B}
\pdfglyphtounicode{twodeva}{0968}
\pdfglyphtounicode{twodotenleader}{2025}
\pdfglyphtounicode{twodotleader}{2025}
\pdfglyphtounicode{twodotleadervertical}{FE30}
\pdfglyphtounicode{twogujarati}{0AE8}
\pdfglyphtounicode{twogurmukhi}{0A68}
\pdfglyphtounicode{twohackarabic}{0662}
\pdfglyphtounicode{twohangzhou}{3022}
\pdfglyphtounicode{twoideographicparen}{3221}
\pdfglyphtounicode{twoinferior}{2082}
\pdfglyphtounicode{twomonospace}{FF12}
\pdfglyphtounicode{twonumeratorbengali}{09F5}
\pdfglyphtounicode{twooldstyle}{0032}
\pdfglyphtounicode{twoparen}{2475}
\pdfglyphtounicode{twoperiod}{2489}
\pdfglyphtounicode{twopersian}{06F2}
\pdfglyphtounicode{tworoman}{2171}
\pdfglyphtounicode{twostroke}{01BB}
\pdfglyphtounicode{twosuperior}{00B2}
\pdfglyphtounicode{twothai}{0E52}
\pdfglyphtounicode{twothirds}{2154}
\pdfglyphtounicode{u}{0075}
\pdfglyphtounicode{uacute}{00FA}
\pdfglyphtounicode{ubar}{0289}
\pdfglyphtounicode{ubengali}{0989}
\pdfglyphtounicode{ubopomofo}{3128}
\pdfglyphtounicode{ubreve}{016D}
\pdfglyphtounicode{ucaron}{01D4}
\pdfglyphtounicode{ucircle}{24E4}
\pdfglyphtounicode{ucircumflex}{00FB}
\pdfglyphtounicode{ucircumflexbelow}{1E77}
\pdfglyphtounicode{ucyrillic}{0443}
\pdfglyphtounicode{udattadeva}{0951}
\pdfglyphtounicode{udblacute}{0171}
\pdfglyphtounicode{udblgrave}{0215}
\pdfglyphtounicode{udeva}{0909}
\pdfglyphtounicode{udieresis}{00FC}
\pdfglyphtounicode{udieresisacute}{01D8}
\pdfglyphtounicode{udieresisbelow}{1E73}
\pdfglyphtounicode{udieresiscaron}{01DA}
\pdfglyphtounicode{udieresiscyrillic}{04F1}
\pdfglyphtounicode{udieresisgrave}{01DC}
\pdfglyphtounicode{udieresismacron}{01D6}
\pdfglyphtounicode{udotbelow}{1EE5}
\pdfglyphtounicode{ugrave}{00F9}
\pdfglyphtounicode{ugujarati}{0A89}
\pdfglyphtounicode{ugurmukhi}{0A09}
\pdfglyphtounicode{uhiragana}{3046}
\pdfglyphtounicode{uhookabove}{1EE7}
\pdfglyphtounicode{uhorn}{01B0}
\pdfglyphtounicode{uhornacute}{1EE9}
\pdfglyphtounicode{uhorndotbelow}{1EF1}
\pdfglyphtounicode{uhorngrave}{1EEB}
\pdfglyphtounicode{uhornhookabove}{1EED}
\pdfglyphtounicode{uhorntilde}{1EEF}
\pdfglyphtounicode{uhungarumlaut}{0171}
\pdfglyphtounicode{uhungarumlautcyrillic}{04F3}
\pdfglyphtounicode{uinvertedbreve}{0217}
\pdfglyphtounicode{ukatakana}{30A6}
\pdfglyphtounicode{ukatakanahalfwidth}{FF73}
\pdfglyphtounicode{ukcyrillic}{0479}
\pdfglyphtounicode{ukorean}{315C}
\pdfglyphtounicode{umacron}{016B}
\pdfglyphtounicode{umacroncyrillic}{04EF}
\pdfglyphtounicode{umacrondieresis}{1E7B}
\pdfglyphtounicode{umatragurmukhi}{0A41}
\pdfglyphtounicode{umonospace}{FF55}
\pdfglyphtounicode{underscore}{005F}
\pdfglyphtounicode{underscoredbl}{2017}
\pdfglyphtounicode{underscoremonospace}{FF3F}
\pdfglyphtounicode{underscorevertical}{FE33}
\pdfglyphtounicode{underscorewavy}{FE4F}
\pdfglyphtounicode{union}{222A}
\pdfglyphtounicode{uniondbl}{22D3}
\pdfglyphtounicode{unionmulti}{228E}
\pdfglyphtounicode{unionsq}{2294}
\pdfglyphtounicode{universal}{2200}
\pdfglyphtounicode{uogonek}{0173}
\pdfglyphtounicode{uparen}{24B0}
\pdfglyphtounicode{upblock}{2580}
\pdfglyphtounicode{upperdothebrew}{05C4}
\pdfglyphtounicode{uprise}{22CF}
\pdfglyphtounicode{upsilon}{03C5}
\pdfglyphtounicode{upsilondieresis}{03CB}
\pdfglyphtounicode{upsilondieresistonos}{03B0}
\pdfglyphtounicode{upsilonlatin}{028A}
\pdfglyphtounicode{upsilontonos}{03CD}
\pdfglyphtounicode{upslope}{29F8}
\pdfglyphtounicode{uptackbelowcmb}{031D}
\pdfglyphtounicode{uptackmod}{02D4}
\pdfglyphtounicode{uragurmukhi}{0A73}
\pdfglyphtounicode{uring}{016F}
\pdfglyphtounicode{ushortcyrillic}{045E}
\pdfglyphtounicode{usmallhiragana}{3045}
\pdfglyphtounicode{usmallkatakana}{30A5}
\pdfglyphtounicode{usmallkatakanahalfwidth}{FF69}
\pdfglyphtounicode{ustraightcyrillic}{04AF}
\pdfglyphtounicode{ustraightstrokecyrillic}{04B1}
\pdfglyphtounicode{utilde}{0169}
\pdfglyphtounicode{utildeacute}{1E79}
\pdfglyphtounicode{utildebelow}{1E75}
\pdfglyphtounicode{uubengali}{098A}
\pdfglyphtounicode{uudeva}{090A}
\pdfglyphtounicode{uugujarati}{0A8A}
\pdfglyphtounicode{uugurmukhi}{0A0A}
\pdfglyphtounicode{uumatragurmukhi}{0A42}
\pdfglyphtounicode{uuvowelsignbengali}{09C2}
\pdfglyphtounicode{uuvowelsigndeva}{0942}
\pdfglyphtounicode{uuvowelsigngujarati}{0AC2}
\pdfglyphtounicode{uvowelsignbengali}{09C1}
\pdfglyphtounicode{uvowelsigndeva}{0941}
\pdfglyphtounicode{uvowelsigngujarati}{0AC1}
\pdfglyphtounicode{v}{0076}
\pdfglyphtounicode{vadeva}{0935}
\pdfglyphtounicode{vagujarati}{0AB5}
\pdfglyphtounicode{vagurmukhi}{0A35}
\pdfglyphtounicode{vakatakana}{30F7}
\pdfglyphtounicode{vav}{05D5}
\pdfglyphtounicode{vavdagesh}{FB35}
\pdfglyphtounicode{vavdagesh65}{FB35}
\pdfglyphtounicode{vavdageshhebrew}{FB35}
\pdfglyphtounicode{vavhebrew}{05D5}
\pdfglyphtounicode{vavholam}{FB4B}
\pdfglyphtounicode{vavholamhebrew}{FB4B}
\pdfglyphtounicode{vavvavhebrew}{05F0}
\pdfglyphtounicode{vavyodhebrew}{05F1}
\pdfglyphtounicode{vcircle}{24E5}
\pdfglyphtounicode{vdotbelow}{1E7F}
\pdfglyphtounicode{vector}{20D7}
\pdfglyphtounicode{vecyrillic}{0432}
\pdfglyphtounicode{veharabic}{06A4}
\pdfglyphtounicode{vehfinalarabic}{FB6B}
\pdfglyphtounicode{vehinitialarabic}{FB6C}
\pdfglyphtounicode{vehmedialarabic}{FB6D}
\pdfglyphtounicode{vekatakana}{30F9}
\pdfglyphtounicode{venus}{2640}
\pdfglyphtounicode{verticalbar}{007C}
\pdfglyphtounicode{verticallineabovecmb}{030D}
\pdfglyphtounicode{verticallinebelowcmb}{0329}
\pdfglyphtounicode{verticallinelowmod}{02CC}
\pdfglyphtounicode{verticallinemod}{02C8}
\pdfglyphtounicode{vewarmenian}{057E}
\pdfglyphtounicode{vhook}{028B}
\pdfglyphtounicode{vikatakana}{30F8}
\pdfglyphtounicode{viramabengali}{09CD}
\pdfglyphtounicode{viramadeva}{094D}
\pdfglyphtounicode{viramagujarati}{0ACD}
\pdfglyphtounicode{visargabengali}{0983}
\pdfglyphtounicode{visargadeva}{0903}
\pdfglyphtounicode{visargagujarati}{0A83}
\pdfglyphtounicode{visiblespace}{2423}
\pdfglyphtounicode{visualspace}{2423}
\pdfglyphtounicode{vmonospace}{FF56}
\pdfglyphtounicode{voarmenian}{0578}
\pdfglyphtounicode{voicediterationhiragana}{309E}
\pdfglyphtounicode{voicediterationkatakana}{30FE}
\pdfglyphtounicode{voicedmarkkana}{309B}
\pdfglyphtounicode{voicedmarkkanahalfwidth}{FF9E}
\pdfglyphtounicode{vokatakana}{30FA}
\pdfglyphtounicode{vparen}{24B1}
\pdfglyphtounicode{vtilde}{1E7D}
\pdfglyphtounicode{vturned}{028C}
\pdfglyphtounicode{vuhiragana}{3094}
\pdfglyphtounicode{vukatakana}{30F4}
\pdfglyphtounicode{w}{0077}
\pdfglyphtounicode{wacute}{1E83}
\pdfglyphtounicode{waekorean}{3159}
\pdfglyphtounicode{wahiragana}{308F}
\pdfglyphtounicode{wakatakana}{30EF}
\pdfglyphtounicode{wakatakanahalfwidth}{FF9C}
\pdfglyphtounicode{wakorean}{3158}
\pdfglyphtounicode{wasmallhiragana}{308E}
\pdfglyphtounicode{wasmallkatakana}{30EE}
\pdfglyphtounicode{wattosquare}{3357}
\pdfglyphtounicode{wavedash}{301C}
\pdfglyphtounicode{wavyunderscorevertical}{FE34}
\pdfglyphtounicode{wawarabic}{0648}
\pdfglyphtounicode{wawfinalarabic}{FEEE}
\pdfglyphtounicode{wawhamzaabovearabic}{0624}
\pdfglyphtounicode{wawhamzaabovefinalarabic}{FE86}
\pdfglyphtounicode{wbsquare}{33DD}
\pdfglyphtounicode{wcircle}{24E6}
\pdfglyphtounicode{wcircumflex}{0175}
\pdfglyphtounicode{wdieresis}{1E85}
\pdfglyphtounicode{wdotaccent}{1E87}
\pdfglyphtounicode{wdotbelow}{1E89}
\pdfglyphtounicode{wehiragana}{3091}
\pdfglyphtounicode{weierstrass}{2118}
\pdfglyphtounicode{wekatakana}{30F1}
\pdfglyphtounicode{wekorean}{315E}
\pdfglyphtounicode{weokorean}{315D}
\pdfglyphtounicode{wgrave}{1E81}
\pdfglyphtounicode{whitebullet}{25E6}
\pdfglyphtounicode{whitecircle}{25CB}
\pdfglyphtounicode{whitecircleinverse}{25D9}
\pdfglyphtounicode{whitecornerbracketleft}{300E}
\pdfglyphtounicode{whitecornerbracketleftvertical}{FE43}
\pdfglyphtounicode{whitecornerbracketright}{300F}
\pdfglyphtounicode{whitecornerbracketrightvertical}{FE44}
\pdfglyphtounicode{whitediamond}{25C7}
\pdfglyphtounicode{whitediamondcontainingblacksmalldiamond}{25C8}
\pdfglyphtounicode{whitedownpointingsmalltriangle}{25BF}
\pdfglyphtounicode{whitedownpointingtriangle}{25BD}
\pdfglyphtounicode{whiteleftpointingsmalltriangle}{25C3}
\pdfglyphtounicode{whiteleftpointingtriangle}{25C1}
\pdfglyphtounicode{whitelenticularbracketleft}{3016}
\pdfglyphtounicode{whitelenticularbracketright}{3017}
\pdfglyphtounicode{whiterightpointingsmalltriangle}{25B9}
\pdfglyphtounicode{whiterightpointingtriangle}{25B7}
\pdfglyphtounicode{whitesmallsquare}{25AB}
\pdfglyphtounicode{whitesmilingface}{263A}
\pdfglyphtounicode{whitesquare}{25A1}
\pdfglyphtounicode{whitestar}{2606}
\pdfglyphtounicode{whitetelephone}{260F}
\pdfglyphtounicode{whitetortoiseshellbracketleft}{3018}
\pdfglyphtounicode{whitetortoiseshellbracketright}{3019}
\pdfglyphtounicode{whiteuppointingsmalltriangle}{25B5}
\pdfglyphtounicode{whiteuppointingtriangle}{25B3}
\pdfglyphtounicode{wihiragana}{3090}
\pdfglyphtounicode{wikatakana}{30F0}
\pdfglyphtounicode{wikorean}{315F}
\pdfglyphtounicode{wmonospace}{FF57}
\pdfglyphtounicode{wohiragana}{3092}
\pdfglyphtounicode{wokatakana}{30F2}
\pdfglyphtounicode{wokatakanahalfwidth}{FF66}
\pdfglyphtounicode{won}{20A9}
\pdfglyphtounicode{wonmonospace}{FFE6}
\pdfglyphtounicode{wowaenthai}{0E27}
\pdfglyphtounicode{wparen}{24B2}
\pdfglyphtounicode{wreathproduct}{2240}
\pdfglyphtounicode{wring}{1E98}
\pdfglyphtounicode{wsuperior}{02B7}
\pdfglyphtounicode{wturned}{028D}
\pdfglyphtounicode{wynn}{01BF}
\pdfglyphtounicode{x}{0078}
\pdfglyphtounicode{xabovecmb}{033D}
\pdfglyphtounicode{xbopomofo}{3112}
\pdfglyphtounicode{xcircle}{24E7}
\pdfglyphtounicode{xdieresis}{1E8D}
\pdfglyphtounicode{xdotaccent}{1E8B}
\pdfglyphtounicode{xeharmenian}{056D}
\pdfglyphtounicode{xi}{03BE}
\pdfglyphtounicode{xmonospace}{FF58}
\pdfglyphtounicode{xparen}{24B3}
\pdfglyphtounicode{xsuperior}{02E3}
\pdfglyphtounicode{y}{0079}
\pdfglyphtounicode{yaadosquare}{334E}
\pdfglyphtounicode{yabengali}{09AF}
\pdfglyphtounicode{yacute}{00FD}
\pdfglyphtounicode{yadeva}{092F}
\pdfglyphtounicode{yaekorean}{3152}
\pdfglyphtounicode{yagujarati}{0AAF}
\pdfglyphtounicode{yagurmukhi}{0A2F}
\pdfglyphtounicode{yahiragana}{3084}
\pdfglyphtounicode{yakatakana}{30E4}
\pdfglyphtounicode{yakatakanahalfwidth}{FF94}
\pdfglyphtounicode{yakorean}{3151}
\pdfglyphtounicode{yamakkanthai}{0E4E}
\pdfglyphtounicode{yasmallhiragana}{3083}
\pdfglyphtounicode{yasmallkatakana}{30E3}
\pdfglyphtounicode{yasmallkatakanahalfwidth}{FF6C}
\pdfglyphtounicode{yatcyrillic}{0463}
\pdfglyphtounicode{ycircle}{24E8}
\pdfglyphtounicode{ycircumflex}{0177}
\pdfglyphtounicode{ydieresis}{00FF}
\pdfglyphtounicode{ydotaccent}{1E8F}
\pdfglyphtounicode{ydotbelow}{1EF5}
\pdfglyphtounicode{yeharabic}{064A}
\pdfglyphtounicode{yehbarreearabic}{06D2}
\pdfglyphtounicode{yehbarreefinalarabic}{FBAF}
\pdfglyphtounicode{yehfinalarabic}{FEF2}
\pdfglyphtounicode{yehhamzaabovearabic}{0626}
\pdfglyphtounicode{yehhamzaabovefinalarabic}{FE8A}
\pdfglyphtounicode{yehhamzaaboveinitialarabic}{FE8B}
\pdfglyphtounicode{yehhamzaabovemedialarabic}{FE8C}
\pdfglyphtounicode{yehinitialarabic}{FEF3}
\pdfglyphtounicode{yehmedialarabic}{FEF4}
\pdfglyphtounicode{yehmeeminitialarabic}{FCDD}
\pdfglyphtounicode{yehmeemisolatedarabic}{FC58}
\pdfglyphtounicode{yehnoonfinalarabic}{FC94}
\pdfglyphtounicode{yehthreedotsbelowarabic}{06D1}
\pdfglyphtounicode{yekorean}{3156}
\pdfglyphtounicode{yen}{00A5}
\pdfglyphtounicode{yenmonospace}{FFE5}
\pdfglyphtounicode{yeokorean}{3155}
\pdfglyphtounicode{yeorinhieuhkorean}{3186}
\pdfglyphtounicode{yerahbenyomohebrew}{05AA}
\pdfglyphtounicode{yerahbenyomolefthebrew}{05AA}
\pdfglyphtounicode{yericyrillic}{044B}
\pdfglyphtounicode{yerudieresiscyrillic}{04F9}
\pdfglyphtounicode{yesieungkorean}{3181}
\pdfglyphtounicode{yesieungpansioskorean}{3183}
\pdfglyphtounicode{yesieungsioskorean}{3182}
\pdfglyphtounicode{yetivhebrew}{059A}
\pdfglyphtounicode{ygrave}{1EF3}
\pdfglyphtounicode{yhook}{01B4}
\pdfglyphtounicode{yhookabove}{1EF7}
\pdfglyphtounicode{yiarmenian}{0575}
\pdfglyphtounicode{yicyrillic}{0457}
\pdfglyphtounicode{yikorean}{3162}
\pdfglyphtounicode{yinyang}{262F}
\pdfglyphtounicode{yiwnarmenian}{0582}
\pdfglyphtounicode{ymonospace}{FF59}
\pdfglyphtounicode{yod}{05D9}
\pdfglyphtounicode{yoddagesh}{FB39}
\pdfglyphtounicode{yoddageshhebrew}{FB39}
\pdfglyphtounicode{yodhebrew}{05D9}
\pdfglyphtounicode{yodyodhebrew}{05F2}
\pdfglyphtounicode{yodyodpatahhebrew}{FB1F}
\pdfglyphtounicode{yohiragana}{3088}
\pdfglyphtounicode{yoikorean}{3189}
\pdfglyphtounicode{yokatakana}{30E8}
\pdfglyphtounicode{yokatakanahalfwidth}{FF96}
\pdfglyphtounicode{yokorean}{315B}
\pdfglyphtounicode{yosmallhiragana}{3087}
\pdfglyphtounicode{yosmallkatakana}{30E7}
\pdfglyphtounicode{yosmallkatakanahalfwidth}{FF6E}
\pdfglyphtounicode{yotgreek}{03F3}
\pdfglyphtounicode{yoyaekorean}{3188}
\pdfglyphtounicode{yoyakorean}{3187}
\pdfglyphtounicode{yoyakthai}{0E22}
\pdfglyphtounicode{yoyingthai}{0E0D}
\pdfglyphtounicode{yparen}{24B4}
\pdfglyphtounicode{ypogegrammeni}{037A}
\pdfglyphtounicode{ypogegrammenigreekcmb}{0345}
\pdfglyphtounicode{yr}{01A6}
\pdfglyphtounicode{yring}{1E99}
\pdfglyphtounicode{ysuperior}{02B8}
\pdfglyphtounicode{ytilde}{1EF9}
\pdfglyphtounicode{yturned}{028E}
\pdfglyphtounicode{yuhiragana}{3086}
\pdfglyphtounicode{yuikorean}{318C}
\pdfglyphtounicode{yukatakana}{30E6}
\pdfglyphtounicode{yukatakanahalfwidth}{FF95}
\pdfglyphtounicode{yukorean}{3160}
\pdfglyphtounicode{yusbigcyrillic}{046B}
\pdfglyphtounicode{yusbigiotifiedcyrillic}{046D}
\pdfglyphtounicode{yuslittlecyrillic}{0467}
\pdfglyphtounicode{yuslittleiotifiedcyrillic}{0469}
\pdfglyphtounicode{yusmallhiragana}{3085}
\pdfglyphtounicode{yusmallkatakana}{30E5}
\pdfglyphtounicode{yusmallkatakanahalfwidth}{FF6D}
\pdfglyphtounicode{yuyekorean}{318B}
\pdfglyphtounicode{yuyeokorean}{318A}
\pdfglyphtounicode{yyabengali}{09DF}
\pdfglyphtounicode{yyadeva}{095F}
\pdfglyphtounicode{z}{007A}
\pdfglyphtounicode{zaarmenian}{0566}
\pdfglyphtounicode{zacute}{017A}
\pdfglyphtounicode{zadeva}{095B}
\pdfglyphtounicode{zagurmukhi}{0A5B}
\pdfglyphtounicode{zaharabic}{0638}
\pdfglyphtounicode{zahfinalarabic}{FEC6}
\pdfglyphtounicode{zahinitialarabic}{FEC7}
\pdfglyphtounicode{zahiragana}{3056}
\pdfglyphtounicode{zahmedialarabic}{FEC8}
\pdfglyphtounicode{zainarabic}{0632}
\pdfglyphtounicode{zainfinalarabic}{FEB0}
\pdfglyphtounicode{zakatakana}{30B6}
\pdfglyphtounicode{zaqefgadolhebrew}{0595}
\pdfglyphtounicode{zaqefqatanhebrew}{0594}
\pdfglyphtounicode{zarqahebrew}{0598}
\pdfglyphtounicode{zayin}{05D6}
\pdfglyphtounicode{zayindagesh}{FB36}
\pdfglyphtounicode{zayindageshhebrew}{FB36}
\pdfglyphtounicode{zayinhebrew}{05D6}
\pdfglyphtounicode{zbopomofo}{3117}
\pdfglyphtounicode{zcaron}{017E}
\pdfglyphtounicode{zcircle}{24E9}
\pdfglyphtounicode{zcircumflex}{1E91}
\pdfglyphtounicode{zcurl}{0291}
\pdfglyphtounicode{zdot}{017C}
\pdfglyphtounicode{zdotaccent}{017C}
\pdfglyphtounicode{zdotbelow}{1E93}
\pdfglyphtounicode{zecyrillic}{0437}
\pdfglyphtounicode{zedescendercyrillic}{0499}
\pdfglyphtounicode{zedieresiscyrillic}{04DF}
\pdfglyphtounicode{zehiragana}{305C}
\pdfglyphtounicode{zekatakana}{30BC}
\pdfglyphtounicode{zero}{0030}
\pdfglyphtounicode{zeroarabic}{0660}
\pdfglyphtounicode{zerobengali}{09E6}
\pdfglyphtounicode{zerodeva}{0966}
\pdfglyphtounicode{zerogujarati}{0AE6}
\pdfglyphtounicode{zerogurmukhi}{0A66}
\pdfglyphtounicode{zerohackarabic}{0660}
\pdfglyphtounicode{zeroinferior}{2080}
\pdfglyphtounicode{zeromonospace}{FF10}
\pdfglyphtounicode{zerooldstyle}{0030}
\pdfglyphtounicode{zeropersian}{06F0}
\pdfglyphtounicode{zerosuperior}{2070}
\pdfglyphtounicode{zerothai}{0E50}
\pdfglyphtounicode{zerowidthjoiner}{FEFF}
\pdfglyphtounicode{zerowidthnonjoiner}{200C}
\pdfglyphtounicode{zerowidthspace}{200B}
\pdfglyphtounicode{zeta}{03B6}
\pdfglyphtounicode{zhbopomofo}{3113}
\pdfglyphtounicode{zhearmenian}{056A}
\pdfglyphtounicode{zhebrevecyrillic}{04C2}
\pdfglyphtounicode{zhecyrillic}{0436}
\pdfglyphtounicode{zhedescendercyrillic}{0497}
\pdfglyphtounicode{zhedieresiscyrillic}{04DD}
\pdfglyphtounicode{zihiragana}{3058}
\pdfglyphtounicode{zikatakana}{30B8}
\pdfglyphtounicode{zinorhebrew}{05AE}
\pdfglyphtounicode{zlinebelow}{1E95}
\pdfglyphtounicode{zmonospace}{FF5A}
\pdfglyphtounicode{zohiragana}{305E}
\pdfglyphtounicode{zokatakana}{30BE}
\pdfglyphtounicode{zparen}{24B5}
\pdfglyphtounicode{zretroflexhook}{0290}
\pdfglyphtounicode{zstroke}{01B6}
\pdfglyphtounicode{zuhiragana}{305A}
\pdfglyphtounicode{zukatakana}{30BA}
\usepackage[utf8]{inputenc}
\usepackage[T1]{fontenc}
\usepackage[british, ngerman]{babel}
\usepackage[english = british]{csquotes} 

\usepackage{textcomp} 
\usepackage{scrhack} 
\usepackage{xspace} 
\usepackage{mparhack} 
\usepackage{fixltx2e} 
\usepackage[latest]{latexrelease}
\usepackage[printonlyused, smaller]{acronym} 
    
\usepackage{ragged2e} 
\newcounter{dummy} 

\usepackage{ifoddpage}

\usepackage[fleqn]{amsmath}
\usepackage[warning, all]{onlyamsmath}\AtBeginDocument{\catcode`\$=3} 
\usepackage{amssymb}
\usepackage{mathtools}
\usepackage{mathdots} 
\usepackage{gensymb} 
\usepackage{dsfont}
\usepackage{stmaryrd}
\usepackage{tensor} 

\allowdisplaybreaks

\usepackage{amsthm}
\usepackage{thmtools}

\usepackage{tabularx} 
\usepackage{caption}
\usepackage{subfig}

\PassOptionsToPackage{pdftex}{graphicx}
\PassOptionsToPackage{svgnames}{xcolor}
\usepackage{standalone}

\usepackage[pdftex]{graphicx}
    \graphicspath{{.}}
\usepackage{trimclip, adjustbox} 

\usepackage[svgnames]{xcolor} 
\usepackage{tikz}
\usetikzlibrary{positioning, patterns, arrows, calc}
\usepackage{tikzscale}

\usepackage{calc}
\usepackage[strict]{changepage} 

\usepackage{scalerel} 

\usepackage[xindy]{imakeidx}
\makeindex[program = texindy]
\makeindex[program = texindy, name = symbols, title = Symbols] 

\usepackage[pdftex, hyperfootnotes = false, pdfpagelabels]{hyperref} 
\hypersetup{%
    colorlinks = true, linktocpage = true, pdfstartpage = 3, pdfstartview = FitV,%
    breaklinks = true, pdfpagemode = UseNone, pageanchor = true, pdfpagemode = UseOutlines,%
    plainpages = false, bookmarksnumbered, bookmarksopen = true, bookmarksopenlevel = 1,%
    hypertexnames = true, pdfhighlight = /O,
    urlcolor = webbrown, linkcolor = RoyalBlue, citecolor = webgreen, 
    pdftitle = {\myTitle},%
    pdfauthor = {\textcopyright\ \myName, \myUni, \myFaculty},%
    pdfsubject = {},%
    pdfkeywords = {},%
    pdfcreator = {pdfLaTeX},%
    pdfproducer = {LaTeX with hyperref and classicthesis}%
}
\usepackage[noabbrev]{cleveref} 

\usepackage[
  backend = biber,
  style = alphabetic, 
  sortlocale = en_GB,
  hyperref = true,
  backref = true,
  abbreviate = false,
  url = true,
  doi = true,
  eprint = true,
  ibidtracker = false, 
  citetracker = false, 
  autocite = footnote
]{biblatex}
\DeclareFieldFormat*{citetitle}{#1} 

\usepackage{snotez}
\setsidenotes{
  note-mark-format = #1.,
  text-format = \itshape\footnotesize,
  perpage = true,
  footnote = true
}
\makeatletter 
\def\snotez@text#1#2{%
      \snotez@marginpar{%
        \checkoddpage
        \ifoddpage
          \RaggedRight
          \snotez@format%
          \snotez@write@mark{\snotez@note@mark@format{\@the@snotez@mark}}%
          \snotez@note@mark@sep#2%
        \else
          \RaggedLeft
          \snotez@format%
          \snotez@write@mark{\snotez@note@mark@format{\@the@snotez@mark}}%
          \snotez@note@mark@sep#2%
        \fi%
      }
}
\makeatother

\usepackage{classicthesis}

\usepackage{cfr-lm} 

\newcommand\clearToOddPage{\clearpage\ifodd\value{page}\else\null\thispagestyle{empty}\clearpage\fi} 
\renewcommand{\slash}{/\penalty\exhyphenpenalty\hspace{0pt}} 

\newcommand*{\define}[1]{\emph{#1}\index{#1}} 
\newcommand*{\defineX}[2]{\emph{#1}\index{#2}} 

\newcommand{\mathnote}[1]{
  \checkoddpage%
  \ifoddpage%
    \tag*{\rlap{\hspace{\marginparsep}\smash{\parbox[t]{\marginparwidth}{\itshape\footnotesize\leavevmode\color{Black}\raggedright\hspace{0pt}{#1}}}}}%
  \else%
    \tag*{\llap{\smash{\parbox[t]{\marginparwidth}{\itshape\footnotesize\leavevmode\color{Black}\raggedleft\hspace{0pt}{#1}}}\hspace{\marginparsep}\hspace{\textwidth}}}%
  \fi%
}

\providecommand\ifAndOnlyIf{\DOTSB\;\Longleftrightarrow\;} 
\providecommand\implies{\DOTSB\;\Longrightarrow\;} 

\DeclarePairedDelimiter\parens{\lparen}{\rparen} 
\DeclarePairedDelimiter\brackets{\lbrack}{\rbrack} 

\DeclarePairedDelimiter\absoluteValueOf{\lvert}{\rvert} 
\DeclarePairedDelimiter\cardinalityOf{\lvert}{\rvert} 
\DeclarePairedDelimiter\lengthOf{\lvert}{\rvert} 
\DeclarePairedDelimiter\lengthOfPath{\lvert}{\rvert} 
\DeclarePairedDelimiter\sizeOf{\lvert}{\rvert} 

\DeclarePairedDelimiter\normOf{\lVert}{\rVert} 
\DeclarePairedDelimiter\radiusOf{\lVert}{\rVert} 
\DeclarePairedDelimiterX\innerProductOf[2]{\langle}{\rangle}{#1,#2} 

\DeclarePairedDelimiter\setOf{\{}{\}} 
\newcommand*\suchThat{\mid} 

\DeclarePairedDelimiter\family{\{}{\}} 
\DeclarePairedDelimiter\net{\{}{\}} 
\DeclarePairedDelimiter\sequence{\lparen}{\rparen} 

\DeclarePairedDelimiter\equivalenceClassOf{\lbrack}{\rbrack} 

\DeclarePairedDelimiter\groupGeneratedBy{\langle}{\rangle} 
\DeclarePairedDelimiter\subshiftGeneratedBy{\langle}{\rangle} 

\DeclarePairedDelimiter\closedInterval{\lbrack}{\rbrack}
\DeclarePairedDelimiter\openInterval{\rbrack}{\lbrack}
\DeclarePairedDelimiter\leftOpenAndRightClosedInterval{\rbrack}{\rbrack}
\DeclarePairedDelimiter\leftClosedAndRightOpenInterval{\lbrack}{\lbrack}

\DeclarePairedDelimiter\ceil{\lceil}{\rceil} 
\DeclarePairedDelimiter\floor{\lfloor}{\rfloor} 

\mathchardef\breakingcomma\mathcode`\, 
{\catcode`,=\active
  \gdef,{\breakingcomma\discretionary{}{}{}}
}
\DeclareRobustCommand\ntuple[1]{\lparen\mathcode`\,=\string"8000 #1\rparen} 

\newcommand*{\N}{\mathbb{N}} 
\newcommand*{\Z}{{\mathbb{Z}}} 
\newcommand*{\R}{\mathbb{R}} 
\newcommand*{\C}{\mathbb{C}} 
\newcommand*{\K}{\mathbb{K}} 
\newcommand*{\F}{\mathbb{F}} 

\DeclareMathOperator{\fractionalPart}{frac} 

\newcommand*{\ball}{\mathbb{B}} 
\newcommand*{\sphere}{\mathbb{S}} 
\newcommand*{\hyperbolicSpace}{\mathbb{H}} 
\newcommand*{\upperHalfPlane}{\mathbb{H}^+} 

\DeclareMathOperator{\boundedFunctionsOn}{\ell^\infty} 
\DeclareMathOperator{\simpleFunctionsOn}{\mathcal{E}} 
\DeclareMathOperator{\continuousMaps}{C} 

\DeclareMathOperator{\automorphismsOf}{Aut} 

\DeclareMathOperator{\probabilityMeasuresOn}{\mathcal{PM}} 
\DeclareMathOperator{\meansOn}{\mathcal{MS}} 

\newcommand*{\imaginaryUnit}{\imath} 
\newcommand*{\EulersNumber}{\mathrm{e}} 

\newcommand*{\indicatorFunction}{\mathds{1}} 
\newcommand*{\functionThatIsIdenticalToOne}{\mathds{1}} 
\newcommand*{\functionThatIsIdenticalToZero}{0} 

\DeclareMathOperator{\identityMap}{id}

\DeclareMathOperator{\Exists}{\exists} 
\DeclareMathOperator{\notExists}{\nexists} 
\DeclareMathOperator{\ForEach}{\forall} 
\newcommand*\Holds{:} 
\newcommand*\SuchThat{:} 

\newcommand*\actsOnPoint{\triangleright} 
\newcommand*\cosetActsOnPoint{\mathbin{\protect\scalerel*{\trianglerighteqslant}{\rhd}}} 
\newcommand*\actsOnMap{\mathbin{\protect\scalerel*{\blacktriangleright}{\triangleright}}} 
\newcommand*\isActedUponBy{\triangleleft} 
\newcommand*\isSemiActedUponBy{\mathbin{\protect\scalerel*{\trianglelefteqslant}{\rhd}}} 
\newcommand*\actsByItsCoordinateOn{\mathbin{\protect\scalerel*{\ooalign{$\trianglelefteqslant$\cr\hidewidth\raise.225ex\hbox{$\blacktriangleleft$}\cr}}{\triangleleft}}} 
\newcommand*\conjugates{\odot} 

\makeatletter
\providecommand*{\Dashv}{%
  \mathrel{%
    \mathpalette\@Dashv\vDash
  }%
}
\newcommand*{\@Dashv}[2]{%
  \reflectbox{$\m@th#1#2$}%
}
\makeatother

\newcommand*\actsOnMeasure{\mathbin{\protect\scalerel*{\vDash}{\rhd}}} 
\newcommand*\measureIsKindOfSemiActedUponBy{\mathbin{\protect\scalerel*{\Dashv}{\rhd}}} 

\newcommand{\VDash}{%
  \mathrel{
    \text{\clipbox{0pt 0pt {.8\width} 0pt}{$\Vdash$}}
    \mkern.9mu
    \text{\adjustbox{width=.87\width,height=\height}{$\vDash$}}
  }
}

\makeatletter
\providecommand*{\DashV}{%
  \mathrel{%
    \mathpalette\@DashV\VDash
  }%
}
\newcommand*{\@DashV}[2]{%
  \reflectbox{$\m@th#1#2$}%
}
\makeatother
\makeatletter
\providecommand*{\dashV}{%
  \mathrel{%
    \mathpalette\@dashV\Vdash
  }%
}
\newcommand*{\@dashV}[2]{%
  \reflectbox{$\m@th#1#2$}%
}
\makeatother

\newcommand*\actsOnBoundedFunction{\mathbin{\protect\scalerel*{\Vdash}{\rhd}}} 
\newcommand*\boundedFunctionIsKindOfSemiActedUponBy{\mathbin{\protect\scalerel*{\dashV}{\rhd}}} 
\newcommand*\actsOnMean{\mathbin{\protect\scalerel*{\VDash}{\rhd}}} 
\newcommand*\meanIsKindOfSemiActedUponBy{\mathbin{\protect\scalerel*{\DashV}{\rhd}}} 

\DeclareMathOperator{\symmetricGroupOf}{Sym} 
\DeclareMathOperator{\dihedralGroup}{D} 
\DeclareMathOperator{\centreOf}{Z} 

\newcommand*\modulo{\slash}
\newcommand*\reverseModulo{\backslash} 

\newcommand*\isIsomorphicTo{\cong}

\DeclareMathOperator{\generalLinearGroup}{GL} 
\DeclareMathOperator{\orthogonalGroup}{O} 
\DeclareMathOperator{\specialOrthogonalGroup}{SO} 
\DeclareMathOperator{\specialLinearGroup}{SL} 
\DeclareMathOperator{\unitaryGroup}{U} 
\DeclareMathOperator{\EuclideanGroup}{E} 
\DeclareMathOperator{\antiDeSitterSpace}{AdS} 
\DeclareMathOperator{\isometriesOf}{I} 

\newcommand*{\transposed}{\intercal} 
\DeclareMathOperator{\traceOf}{trace} 
\DeclareMathOperator{\linearSpanOf}{span} 

\newcommand*\from{\colon} 
\newcommand*\givenBy{\colon} 
\newcommand*\after{\circ} 
\DeclareMathOperator{\domainOf}{dom} 
\DeclareMathOperator{\imageOf}{im} 
\newcommand*{\blank}{\mathord{\_}} 
\newcommand*{\restrictedTo}{\mathord{\upharpoonright}} 

\DeclareMathOperator*{\supportOf}{supp} 
\DeclareMathOperator{\kernelOf}{ker} 

\newcommand*\headOf{\sigma} 
\newcommand*\tailOf{\tau} 
\DeclareMathOperator{\distanceOf}{d} 
\newcommand*\graphProductOf{\mathbin{\square}} 
\DeclareMathOperator{\degreeOf}{deg} 

\DeclareMathOperator{\periodic}{Per} 
\newcommand*\boundaryOf{\partial} 
\DeclareMathOperator{\entropyOf}{ent} 
\DeclareMathOperator{\differenceOf}{diff} 

\newcommand*\occursIn{\sqsubseteq} 
\newcommand*\semiOccursIn{\mathrel{\ooalign{$\occursIn$\cr\hidewidth\raise.225ex\hbox{$\circ\mkern.5mu$}\cr}}} 

\DeclareMathOperator{\powerSetOf}{\mathcal{P}} 
\newcommand*\symmetricSetDifferenceOf{\mathbin{\triangle}} 
\newcommand{\closureOf}[1]{\mkern 1.5mu\overline{\mkern-1.5mu#1\mkern-1.5mu}\mkern 1.5mu} 

\makeatletter
\def\moverlay{\mathpalette\mov@rlay}
\def\mov@rlay#1#2{\leavevmode\vtop{%
   \baselineskip\z@skip \lineskiplimit-\maxdimen
   \ialign{\hfil$\m@th#1##$\hfil\cr#2\crcr}}}
\newcommand{\charfusion}[3][\mathord]{
    #1{\ifx#1\mathop\vphantom{#2}\fi
        \mathpalette\mov@rlay{#2\cr#3}
      }
    \ifx#1\mathop\expandafter\displaylimits\fi}
\makeatother

\newcommand{\disjointUnionWith}{\charfusion[\mathbin]{\cup}{\cdot}} 
\newcommand{\bigDisjointUnionOf}{\charfusion[\mathop]{\bigcup}{\cdot}} 

\newcommand*\dominates{\succcurlyeq} 
\newcommand*\isDominatedBy{\preccurlyeq} 
\newcommand*\isEquivalentTo{\sim} 
\newcommand*\isNotEquivalentTo{\nsim} 

\DeclareMathOperator{\evaluationMap}{ev} 
\DeclareMathOperator{\cylinder}{Cyl} 

\newcommand{\myfloatalign}{\centering} 

\newenvironment{wide}{%
  \begin{adjustwidth*}{0pt}{-\marginparsep-\marginparwidth}
    \myfloatalign
}{%
  \end{adjustwidth*}%
}

\newcommand*{\Graph}{\mathcal{G}}

\newcommand*{\Tree}{\mathcal{T}}
\newcommand*{\continuumGraph}{\mathfrak{G}}
\DeclareMathOperator{\Signals}{Sgnl}
\DeclareMathOperator{\Kinds}{Knd}
\DeclareMathOperator{\kindOf}{knd}
\DeclareMathOperator{\Data}{Dt}
\DeclareMathOperator{\NoData}{None}
\DeclareMathOperator{\Vertices}{V}
\newcommand*{\continuumVertices}{\mathfrak{V}}
\DeclareMathOperator{\Edges}{E}
\newcommand*{\continuumEdges}{\mathfrak{E}}
\DeclareMathOperator{\Paths}{P}
\newcommand*{\continuumPaths}{\mathfrak{P}}
\newcommand*{\directionPreserving}{\shortrightarrow}
\newcommand*{\directed}{\rightarrow}
\DeclareMathOperator{\Directions}{Dir}
\DeclareMathOperator{\Arrows}{Arr}
\DeclareMathOperator{\every}{vry}
\DeclareMathOperator{\Times}{T}
\newcommand*{\extendedTimes}{\closureOf{\Times}}
\DeclareMathOperator{\States}{Q}
\DeclareMathOperator{\Configurations}{Cnf}
\DeclareMathOperator{\ReachOf}{R}
\DeclareMathOperator{\directionOf}{dir}
\newcommand*{\start}{\sourceOf}
\newcommand*{\xend}{\targetOf}
\DeclareMathOperator{\speedOf}{spd}
\DeclareMathOperator{\velocityOf}{vel}
\DeclareMathOperator{\datumOf}{dt}
\DeclareMathOperator{\signOf}{sgn}
\DeclareMathOperator{\subtreesOf}{sbtrs}
\DeclareMathOperator{\maxSubtreesOf}{max\,sbtrs}
\newcommand*{\secondMaxSubtreesOf}{\closureOf{\subtreesOf}}
\DeclareMathOperator{\maxPaths}{max\,pths}
\DeclareMathOperator{\maxVertices}{max\,vrtcs}
\DeclareMathOperator{\maxTree}{\hat{\Tree}} 

\newcommand*{\localTransitionFunction}{\delta}
\newcommand*{\globalTransitionFunction}{\mathbin{\boxdot}} 

\makeatletter
\newcommand{\dotDelta}{{\vphantom{\Delta}\mathpalette\d@tD@lta\relax}}
\newcommand{\d@tD@lta}[2]{%
  \ooalign{\hidewidth$\m@th#1\mkern-1mu\cdot$\hidewidth\cr$\m@th#1\Delta$\cr}%
}
\makeatother

\newcommand*{\downTo}{\downarrow} 

\newcommand*{\gtfJump}{\mathbin{\dotDelta}} 
\newcommand*{\gtfIgnoreCollisions}{\mathbin{\boxplus}} 
\newcommand*{\gtfNonNegativeSingularities}{\mathbin{\boxminus}} 
\newcommand*{\gtfNegativeSingularityOfOrderMinusOne}{\mathbin{\boxast}} 
\newcommand*{\concat}{\bullet}
\newcommand*{\reverse}{-}
\newcommand*{\multipli}{\cdot}
\newcommand{\continuumRepresentationOf}[1]{\mkern 1.5mu\overline{\mkern-1.5mu#1\mkern-1.5mu}\mkern 1.5mu} 

\makeatletter
\DeclareRobustCommand{\cev}[1]{%
  \mathpalette\do@cev{#1}%
}
\newcommand{\do@cev}[2]{%
  \fix@cev{#1}{+}%
  \reflectbox{$\m@th#1\vec{\reflectbox{$\fix@cev{#1}{-}\m@th#1#2\fix@cev{#1}{+}$}}$}%
  \fix@cev{#1}{-}%
}
\newcommand{\fix@cev}[2]{%
  \ifx#1\displaystyle
    \mkern#23mu
  \else
    \ifx#1\textstyle
      \mkern#23mu
    \else
      \ifx#1\scriptstyle
        \mkern#22mu
      \else
        \mkern#22mu
      \fi
    \fi
  \fi
}
\makeatother

\newcommand{\character}[1]{\mathtt{#1}}
\newcommand*{\initiateKind}{\character{I}}
\newcommand*{\leafKind}{\character{L}}
\newcommand*{\countKind}{\character{C}}
\newcommand*{\midpointKind}{\character{M}}
\newcommand*{\findMidpointKind}{\character{U}}
\newcommand*{\reflectedFindMidpointKind}{\cev{\findMidpointKind}}
\newcommand*{\slowedDownFindMidpointKind}{\character{V}}
\newcommand*{\freezeKind}{\character{F}}
\newcommand*{\thawKind}{\character{T}}
\newcommand*{\divideKind}{\character{D}}
\newcommand*{\frozenDivideKind}{\freezeKind\divideKind}
\newcommand*{\reflectedDivideKind}{\cev{\divideKind}}
\newcommand*{\boundaryKind}{\character{B}}
\newcommand*{\fireKind}{\character{X}}
\newcommand*{\frozenFireKind}{\freezeKind\fireKind}

\newcommand*{\initiateSignal}[1]{\initiateKind_{#1}}
\newcommand*{\leafSignal}[1]{\leafKind_{#1}}
\newcommand*{\midpointSignal}[2]{\midpointKind_{\setOf{#1, #2}}}
\newcommand*{\countSignal}[1]{\countKind_{#1}}
\newcommand*{\findMidpointSignal}[2]{\findMidpointKind_{#1, #2}} 
\newcommand*{\reflectedFindMidpointSignal}[4]{\reflectedFindMidpointKind_{#1, #2, #3}^{#4}} 
\newcommand*{\slowedDownFindMidpointSignal}[4]{\slowedDownFindMidpointKind_{#1, #2, #3}^{#4}} 
\newcommand*{\freezeSignal}[1]{\freezeKind_{#1}} 
\newcommand*{\thawSignal}[3]{\thawKind_{#1, #2}^{#3}} 
\newcommand*{\divideSignal}[2]{\divideKind_{#1, #2}} 
\newcommand*{\frozenDivideSignal}[2]{\frozenDivideKind_{#1, #2}} 
\newcommand*{\reflectedDivideSignal}[1]{\reflectedDivideKind_{#1}}
\newcommand*{\boundarySignal}{\boundaryKind}
\newcommand*{\fireSignal}{\fireKind}
\newcommand*{\frozenFireSignal}{\frozenFireKind}

\newcommand*{\midpoint}{\mathfrak{m}}
\newcommand*{\general}{\mathfrak{g}}

\newcommand*{\sourceOf}{\sigma}
\newcommand*{\targetOf}{\tau}
\newcommand*{\bedOf}{\beta}
\newcommand*{\eendsOf}{\varepsilon}
\DeclareMathOperator{\invert}{inv}

\newcommand*{\weightOf}{\omega}
\newcommand*{\length}{\weightOf} 

\newcommand*{\direct}[1]{\vec{#1}}

\newcommand*{\emptyWord}{\lambda}

\newcommand*{\booleans}{\mathbb{B}}
\newcommand*{\no}{\character{no}}
\newcommand*{\yes}{\character{yes}}

\DeclareMathOperator{\tree}{tree}
\DeclareMathOperator{\virtualTree}{virtualTree}
\DeclareMathOperator{\boundaryCases}{bndryCases}

\DeclareMathOperator{\areEmpty}{areEmpty}
\DeclareMathOperator{\isInLeaf}{inLeaf}
\DeclareMathOperator{\isPenultimate}{penultimate}

\DeclareMathOperator*{\argmax}{arg\,max}

\newcommand*{\discreteInterval}[2]{[#1 : #2]}

\newcommand*{\deadState}{\mathtt{\#}}
\newcommand*{\generalState}{\mathtt{g}}
\newcommand*{\soldierState}{\mathtt{s}}
\newcommand*{\fireState}{\mathtt{f}}

\declaretheoremstyle[
  headfont = \scshape, headpunct = ., postheadspace = 0.5em,
  bodyfont = \itshape,
  spaceabove = \baselineskip, spacebelow = \baselineskip,
  qed = \ensuremath{_\square} 
]{theorem}

\declaretheorem[numberwithin=section, style=theorem, refname={main theorem, main theorems}, Refname={Main Theorem, Main Theorems}, name=Main Theorem]{main-theorem}
\declaretheorem[sibling=main-theorem, style=theorem, refname={theorem,theorems}, Refname={Theorem,Theorems}]{theorem}
\declaretheorem[sibling=main-theorem, style=theorem, refname={lemma,lemmata}, Refname={Lemma,Lemmata}]{lemma}
\declaretheorem[sibling=main-theorem, style=theorem, refname={corollary,corollaries}, Refname={Corollary,Corollaries}]{corollary}

\declaretheoremstyle[
  headfont = \scshape, headpunct = ., postheadspace = 0.5em,
  bodyfont = \upshape,
  spaceabove = \baselineskip, spacebelow = \baselineskip,
  qed = \ensuremath{_\square} 
]{definition}

\declaretheorem[sibling=main-theorem, style=definition, refname={definition,definitions}, Refname={Definition,Definitions}]{definition}
\declaretheorem[sibling=main-theorem, style=definition, refname={example,examples}, Refname={Example,Examples}]{example}
\declaretheorem[sibling=main-theorem, style=definition, refname={counterexample,counterexamples}, Refname={Counterexample,Counterexamples}]{counterexample}
\declaretheorem[sibling=main-theorem, style=definition, refname={remark,remarks}, Refname={Remark,Remarks}]{remark}
\declaretheorem[sibling=main-theorem, style=definition, refname={open problem,open problems}, Refname={Open Problem,Open Problems}, name=Open Problem]{open-problem}

\declaretheoremstyle[
  headfont = \scshape, headpunct = ., postheadspace = 0.5em,
  bodyfont = \normalfont,
  spaceabove = \baselineskip, spacebelow = \baselineskip,
  qed = \ensuremath{_\blacksquare}
]{proof}

\declaretheorem[unnumbered, style=proof, refname={proof,proofs}, Refname={Proof,Proofs}]{proof}
\declaretheorem[unnumbered, style=proof, refname={proof sektch,proof sketches}, Refname={Proof Sketch,Proof Sketches}, name=Proof Sketch]{proof-sketch}
\declaretheorem[unnumbered, style=proof, refname={usage note,usage notes}, Refname={Usage Note,Usage Notes}, name=Usage Note]{usage-note}

\makeatletter
\newcommand{\proofPart}[1]{%
  \par%
  \addvspace{\medskipamount}%
  \noindent\emph{#1.}
  \@afterheading%
}
\makeatother

\def\vonNeumannNeighbourhood(#1,#2,#3,#4,#5,(#6,#7),#8){ 
  \begin{scope}[shift={(#6,#7)},rotate=#8]
    \draw[fill=#4] (-3,-1) rectangle (-1,+1); 
    \draw[fill=#3] (-1,+1) rectangle (+1,+3); 
    \draw[fill=#1] (-1,-1) rectangle (+1,+1); 
    \draw[fill=#5] (-1,-3) rectangle (+1,-1); 
    \draw[fill=#2] (+1,-1) rectangle (+3,+1); 
  \end{scope}
}

\def\vonNeumannNeighbourhoodTwoRotationsHorizontal(#1,#2,#3,#4,#5,(#6,#7)){
  \begin{scope}[shift={(#6,#7)}]
    \foreach \x/\angle in {0/0,8/90} {
      \vonNeumannNeighbourhood(#1,#2,#3,#4,#5,(\x,0),\angle)
    }
  \end{scope}
}

\def\vonNeumannNeighbourhoodTwoRotationsVertical(#1,#2,#3,#4,#5,(#6,#7)){
  \begin{scope}[shift={(#6,#7)}]
    \foreach \y/\angle in {0/0,8/90} {
      \vonNeumannNeighbourhood(#1,#2,#3,#4,#5,(0,-\y),\angle)
    }
  \end{scope}
}

\def\vonNeumannNeighbourhoodFourRotationsHorizontal(#1,#2,#3,#4,#5,(#6,#7)){
  \begin{scope}[shift={(#6,#7)}]
    \foreach \x/\angle in {0/0,8/90,16/180,24/270} {
      \vonNeumannNeighbourhood(#1,#2,#3,#4,#5,(\x,0),\angle);
    }
  \end{scope}
}

\def\vonNeumannNeighbourhoodFourRotationsVertical(#1,#2,#3,#4,#5,(#6,#7)){
  \begin{scope}[shift={(#6,#7)}]
    \foreach \y/\angle in {0/0,8/90,16/180,24/270} {
      \vonNeumannNeighbourhood(#1,#2,#3,#4,#5,(0,-\y),\angle);
    }
  \end{scope}
}

\def\vonNeumannNeighbourhoodFourRotationsVerticalBlocksOfTwo(#1,#2,#3,#4,#5,(#6,#7)){
  \begin{scope}[shift={(#6,#7)}]
    \foreach \y/\angle in {0/0,8/180} {
      \vonNeumannNeighbourhood(#1,#2,#3,#4,#5,(0,-\y),\angle);
    }
    \foreach \y/\angle in {0/90,8/270} {
      \vonNeumannNeighbourhood(#1,#2,#3,#4,#5,(8,-\y),\angle);
    }
  \end{scope}
}

\def\vonNeumannNeighbourhoodCenter(#1,(#2,#3)){
  \begin{scope}[shift={(#2,#3)}]
    \vonNeumannNeighbourhood(#1,white,white,white,white,(0,0),0)
    \vonNeumannNeighbourhoodFourRotationsHorizontal(#1,black,white,white,white,(8,0))
    \vonNeumannNeighbourhoodFourRotationsHorizontal(#1,black,black,white,white,(40,0))
    \vonNeumannNeighbourhoodTwoRotationsHorizontal(#1,black,white,black,white,(72,0))
    \vonNeumannNeighbourhoodFourRotationsHorizontal(#1,black,black,black,white,(88,0))
    \vonNeumannNeighbourhood(#1,black,black,black,black,(120,0),0)
  \end{scope}
}

\def\vonNeumannNeighbourhoodCenterRotations(#1,(#2,#3)){
  \begin{scope}[shift={(#2,#3)}]
    \vonNeumannNeighbourhood(#1,white,white,white,white,(0,0),0)
    \vonNeumannNeighbourhoodFourRotationsVertical(#1,black,white,white,white,(8,0))
    \vonNeumannNeighbourhoodFourRotationsVertical(#1,black,black,white,white,(16,0))
    \vonNeumannNeighbourhoodTwoRotationsVertical(#1,black,white,black,white,(24,0))
    \vonNeumannNeighbourhoodFourRotationsVertical(#1,black,black,black,white,(32,0))
    \vonNeumannNeighbourhood(#1,black,black,black,black,(40,0),0)
  \end{scope}
}

\def\vonNeumannNeighbourhoodCenterHorizontalAndVerticalReflections(#1,(#2,#3)){
  \begin{scope}[shift={(#2,#3)}]
    \vonNeumannNeighbourhood(#1,white,white,white,white,(0,0),0)
    \vonNeumannNeighbourhoodFourRotationsVerticalBlocksOfTwo(#1,black,white,white,white,(8,0))
    \vonNeumannNeighbourhoodFourRotationsVertical(#1,black,black,white,white,(24,0))
    \vonNeumannNeighbourhoodTwoRotationsHorizontal(#1,black,white,black,white,(32,0))
    \vonNeumannNeighbourhoodFourRotationsVerticalBlocksOfTwo(#1,black,black,black,white,(48,0))
    \vonNeumannNeighbourhood(#1,black,black,black,black,(64,0),0)
  \end{scope}
}

\def\vonNeumannNeighbourhoodCenterDiagonalReflections(#1,(#2,#3)){
  \begin{scope}[shift={(#2,#3)}]
    \vonNeumannNeighbourhood(#1,white,white,white,white,(0,0),0)
    \vonNeumannNeighbourhoodFourRotationsVertical(#1,black,white,white,white,(8,0))
    \vonNeumannNeighbourhoodFourRotationsVerticalBlocksOfTwo(#1,black,black,white,white,(16,0))
    \vonNeumannNeighbourhoodTwoRotationsVertical(#1,black,white,black,white,(32,0))
    \vonNeumannNeighbourhoodFourRotationsVertical(#1,black,black,black,white,(40,0))
    \vonNeumannNeighbourhood(#1,black,black,black,black,(48,0),0)
  \end{scope}
}

\def\vonNeumannNeighbourhoodCenterPointReflections(#1,(#2,#3)){
  \begin{scope}[shift={(#2,#3)}]
    \vonNeumannNeighbourhood(#1,white,white,white,white,(0,0),0)
    \vonNeumannNeighbourhoodFourRotationsVerticalBlocksOfTwo(#1,black,white,white,white,(8,0))
    \vonNeumannNeighbourhoodFourRotationsVerticalBlocksOfTwo(#1,black,black,white,white,(24,0))
    \vonNeumannNeighbourhoodTwoRotationsVertical(#1,black,white,black,white,(40,0))
    \vonNeumannNeighbourhoodFourRotationsVerticalBlocksOfTwo(#1,black,black,black,white,(48,0))
    \vonNeumannNeighbourhood(#1,black,black,black,black,(64,0),0)
  \end{scope}
}

\def\TShapedCross(#1,(#2,#3),#4){%
  \begin{scope}[shift = {(#2, #3)}, rotate = #4]
    \draw (0, 0) -- (#1, 0);
    \draw (0, 0) -- (0, #1);
    \draw (0, 0) -- (0, -#1);
  \end{scope}%
}

\def\AMinusContour((#1,#2),#3){
  \begin{scope}[shift = {(#1, #2)}, scale = #3]
    \draw[loosely dashed, tension = 0.5] plot[smooth cycle] coordinates { 
        (-0 + \t, -\t)
        (-0 + \k + \l + \t, -\j)
        (-0, -\j - \k - \l - \t)
        (-0 - \j - \k - \l - \t, 0)
        (-0, \j + \k + \l + \t)
        (-0 + \k + \l + \t, \j)
        (-0 + \t, \t)
    };
  \end{scope}
}
\def\APlusContour((#1,#2),#3){
  \begin{scope}[shift = {(#1, #2)}, scale = #3]
    \draw[densely dotted, tension = 0.5] plot[smooth cycle] coordinates { 
        (0 - 2*\t, -2*\t)
        (0 - \k - \l - \t, -\j)
        (0, -\j - \k - \l - \t)
        (0 + \k + \l + \t, -\j)
        (0 + 2*\t, -2*\t)
        (0 + \j, -\k - \l - \t)
        (0 + \j + \k + \l + \t, 0)
        (0 + \j, \k + \l + \t)
        (0 + 2*\t, 2*\t)
        (0 + \k + \l + \t, \j)
        (0, \j + \k + \l + \t)
        (0 - \k - \l - \t, \j)
        (0 - 2*\t, 2*\t)
    };
  \end{scope}
}
\def\BPlusContour((#1,#2),#3,#4,#5,#6){
  \begin{scope}[shift = {(#1, #2)}, scale = #3, rotate = 90] 
    \draw[densely dashed, tension = 0.5] plot[smooth cycle] coordinates { 
        (#4 - \t, 0)
        (#4, -#6)
        (0 - #5, -\t)
        (0 - \k - \l - \t, -\j)
        (0, -\j - \k - \l - \t)
        (0 + \k + \l + \t, -\j)
        (0 + \t, -\t)
        (0 + \j, -\k - \l - \t)
        (0 + \j + \k + \l + \t, 0)
        (0 + \j, \k + \l + \t)
        (0 + \t, \t)
        (0 + \k + \l + \t, \j)
        (0, \j + \k + \l + \t)
        (0 - \k - \l - \t, \j)
        (0 - #5, +\t)
        (#4, #6)
    };
  \end{scope}
}
\def\BMinusContourTwoBlobs((#1,#2),#3){
  \begin{scope}[shift = {(#1, #2)}, scale = #3]
    \foreach \r in {0, 180} {
      \begin{scope}[rotate = \r]
        \draw[loosely dashdotted, tension = 0.5] plot[smooth cycle] coordinates { 
            (-\j + \k + \l + \t, 0)
            (-\j, \k + \l + \t)
            (-\j - \k - \l - \t, 0)
            (-\j, -\k - \l - \t)
        };
      \end{scope}
    };
  \end{scope}
}
\def\BMinusContour((#1,#2),#3){
  \BMinusContourTwoBlobs((#1,#2),#3);
  \begin{scope}[shift = {(#1, #2)}, scale = #3]
    \foreach \r in {0, 180} {
      \begin{scope}[shift = {(0, -\j)}, rotate = \r]
        \draw[dashdotted, tension = 0.5] plot[smooth cycle] coordinates { 
            (-\k + \l + \t/3, 0)
            (-\k, \l + \t/3)
            (-\k - \l - \t/3, 0)
            (-\k, -\l - \t/3)
        };
      \end{scope}
      \begin{scope}[shift = {(0, -\j - \k)}, rotate = \r]
        \draw[densely dashdotted, tension = 0.5] plot[smooth cycle] coordinates { 
            (-\l + \t/3, 0)
            (-\l, \t/3)
            (-\l - \t/3, 0)
            (-\l, -\t/3)
        };
      \end{scope}
    };
  \end{scope}
}

\def\GridWithCoordinates((#1,#2),(#3,#4)){%
  \draw[step = 1cm, gray, very thin] (#1, #2) grid (#3, #4);
  \draw[->] (#1, 0) -- (#3, 0) node[right] {$x$};
  \draw[->] (0, #2) -- (0, #4) node[above] {$y$};
}
\def\TaxicabBallOutline((#1,#2),#3,#4){%
  \begin{scope}[shift = {(#1, #2)}]
    \draw[#4] (0, #3) -- (#3, 0) -- (0, -#3) -- (-#3, 0) -- cycle;
  \end{scope}%
}
\def\TwoTaxicabBallOutlines((#1,#2),#3,#4){%
  \TaxicabBallOutline((#1,#2),#3,solid);
  \begin{scope}[shift = {(#1, #2)}]
    \draw[loosely dashed] (0, #4) -- (#4, 0) -- (0, -#4);
  \end{scope}%
}
\def\TaxicabBallDots((#1,#2),#3){%
  \begin{scope}[shift = {(#1, #2)}]
    \foreach \x in {-#3, ..., #3} {
      \pgfmathtruncatemacro{\f}{abs(\x) - #3}
      \pgfmathtruncatemacro{\t}{#3 - abs(\x)}%
      \foreach \y in {\f, ..., \t} {
        \fill (\x, \y) circle (0.1cm);
      }
    }
  \end{scope}%
}
\def\TaxicabSphereDots((#1,#2),#3){%
  \begin{scope}[shift = {(#1, #2)}]
    \foreach \x in {-#3, ..., #3} {
      \pgfmathtruncatemacro{\f}{abs(\x) - #3}
      \pgfmathtruncatemacro{\t}{#3 - abs(\x)}%
      \foreach \y in {\f, \t} {
        \fill (\x, \y) circle (0.1cm);
      }
    }
  \end{scope}%
}
\def\TaxicabSphereDotsLeft((#1,#2),#3){%
  \begin{scope}[shift = {(#1, #2)}]
    \foreach \x in {-#3, ..., 0} {
      \pgfmathtruncatemacro{\f}{abs(\x) - #3}
      \pgfmathtruncatemacro{\t}{#3 - abs(\x)}%
      \foreach \y in {\f, \t} {
        \fill (\x, \y) circle (0.1cm);
      }
    }
  \end{scope}%
}
\def\TaxicabSphereDotsRight((#1,#2),#3){%
  \begin{scope}[shift = {(#1, #2)}]
    \foreach \x in {0, ..., #3} {
      \pgfmathtruncatemacro{\f}{abs(\x) - #3}
      \pgfmathtruncatemacro{\t}{#3 - abs(\x)}%
      \foreach \y in {\f, \t} {
        \fill (\x, \y) circle (0.1cm);
      }
    }
  \end{scope}%
}

\def\TShapedCross(#1,(#2,#3),#4){%
  \begin{scope}[shift = {(#2, #3)}, rotate = #4]
    \draw (0, 0) -- (#1, 0);
    \draw (0, 0) -- (0, #1);
    \draw (0, 0) -- (0, -#1);
  \end{scope}%
}

\def\AaInverseDots{%
  \pgfmathsetmacro\www{0.166}
  \draw (0 - \i - \www, 0 - \www) rectangle (0 - \i + \www, 0 + \www);
  \draw (0 - \www, 0 - \www) rectangle (0 + \www, 0 + \www);
  \draw (0 - \j - \www, 0 + \i - \www) rectangle (0 - \j + \www, 0 + \i + \www);
  \draw (0 - \i - \j - \www, 0 - \www) rectangle (0 - \i - \j + \www, 0 + \www);
  \draw (0 - \j - \www, 0 - \i - \www) rectangle (0 - \j + \www, 0 - \i + \www);
  \draw (0 + \i - \www, 0 - \www) rectangle (0 + \i + \www, 0 + \www);
  \draw (0 + \i - \k - \www, 0 + \j - \www) rectangle (0 + \i - \k + \www, 0 + \j + \www);
  \draw (0 + \i - \k - \www, 0 - \j - \www) rectangle (0 + \i - \k + \www, 0 - \j + \www);
  \draw (0 - \www, 0 + \i - \www) rectangle (0 + \www, 0 + \i + \www);
  \draw (0 - \k - \www, 0 + \i + \j - \www) rectangle (0 - \k + \www, 0 + \i + \j + \www);
  \draw (0 - \j - \k - \www, 0 + \i - \www) rectangle (0 - \j - \k + \www, 0 + \i + \www);
  \draw (0 - \i - \k - \www, 0 + \j - \www) rectangle (0 - \i - \k + \www, 0 + \j + \www);
  \draw (0 - \i - \j - \k - \www, 0 - \www) rectangle (0 - \i - \j - \k + \www, 0 + \www);
  \draw (0 - \i - \k - \www, 0 - \j - \www) rectangle (0 - \i - \k + \www, 0 - \j + \www);
  \draw (0 - \www, 0 - \i - \www) rectangle (0 + \www, 0 - \i + \www);
  \draw (0 - \j - \k - \www, 0 - \i - \www) rectangle (0 - \j - \k + \www, 0 - \i + \www);
  \draw (0 - \k - \www, 0 - \i - \j - \www) rectangle (0 - \k + \www, 0 - \i - \j + \www);
}
\def\AbInverseDots{%
  \pgfmathsetmacro\www{0.166}
  \draw (0, 0 - \i) circle (\www);
  \draw (0 + \i, 0 - \j) circle (\www);
  \draw (0, 0) circle (\www);
  \draw (0 - \i, 0 - \j) circle (\www);
  \draw (0, 0 - \i - \j) circle (\www);
  \draw (0 + \i + \j, 0 - \k) circle (\www);
  \draw (0 + \i, 0) circle (\www);
  \draw (0 + \i, 0 - \j - \k) circle (\www);
  \draw (0 + \j, 0 + \i - \k) circle (\www);
  \draw (0, 0 + \i) circle (\www);
  \draw (0 - \j, 0 + \i - \k) circle (\www);
  \draw (0 - \i, 0) circle (\www);
  \draw (0 - \i - \j, 0 - \k) circle (\www);
  \draw (0 - \i, 0 - \j - \k) circle (\www);
  \draw (0 + \j, 0 - \i - \k) circle (\www);
  \draw (0 - \j, 0 - \i - \k) circle (\www);
  \draw (0, 0 - \i - \j - \k) circle (\www);
}
\def\freeGroupCayleyGraph(#1){%
  \draw (0, 0) -- node[above] {$a$} (\i, 0); 
    \TShapedCross(\j,(\i,0),0);
    \begin{scope}[every path/.append style = {#1}]
      \TShapedCross(\k,(\i+\j,0),0);
        \TShapedCross(\l,(\i+\j+\k,0),0);
        \TShapedCross(\l,(\i+\j,\k),90);
        \TShapedCross(\l,(\i+\j,-\k),-90);
      \TShapedCross(\k,(\i,\j),90);
        \TShapedCross(\l,(\i,\j+\k),90);
        \TShapedCross(\l,(\i+\k,\j),0);
        \TShapedCross(\l,(\i-\k,\j),180);
      \TShapedCross(\k,(\i,-\j),-90);
        \TShapedCross(\l,(\i,-\j-\k),-90);
        \TShapedCross(\l,(\i+\k,-\j),0);
        \TShapedCross(\l,(\i-\k,-\j),180);
    \end{scope}

  \draw (0, 0) -- node[below] {$a^{-1}$} (-\i, 0); 
    \TShapedCross(\j,(-\i,0),180);
    \begin{scope}[every path/.append style = {#1}]
      \TShapedCross(\k,(-\i-\j,0),180);
        \TShapedCross(\l,(-\i-\j-\k,0),180);
        \TShapedCross(\l,(-\i-\j,\k),90);
        \TShapedCross(\l,(-\i-\j,-\k),-90);
      \TShapedCross(\k,(-\i,\j),90);
        \TShapedCross(\l,(-\i,\j+\k),90);
        \TShapedCross(\l,(-\i+\k,\j),0);
        \TShapedCross(\l,(-\i-\k,\j),180);
      \TShapedCross(\k,(-\i,-\j),-90);
        \TShapedCross(\l,(-\i,-\j-\k),-90);
        \TShapedCross(\l,(-\i+\k,-\j),0);
        \TShapedCross(\l,(-\i-\k,-\j),180);
    \end{scope}

  \draw (0, 0) -- node[left] {$b$} (0, \i); 
    \TShapedCross(\j,(0,\i),90);
    \begin{scope}[every path/.append style = {#1}]
      \TShapedCross(\k,(0,\i+\j),90);
        \TShapedCross(\l,(0,\i+\j+\k),90);
        \TShapedCross(\l,(\k,\i+\j),0);
        \TShapedCross(\l,(-\k,\i+\j),180);
      \TShapedCross(\k,(\j,\i),0);
        \TShapedCross(\l,(\j+\k,\i),0);
        \TShapedCross(\l,(\j,\i+\k),90);
        \TShapedCross(\l,(\j,\i-\k),-90);
      \TShapedCross(\k,(-\j,\i),180);
        \TShapedCross(\l,(-\j-\k,\i),180);
        \TShapedCross(\l,(-\j,\i+\k),90);
        \TShapedCross(\l,(-\j,\i-\k),-90);
    \end{scope}

  \draw (0, 0) -- node[right] {$b^{-1}$} (0, -\i); 
    \TShapedCross(\j,(0,-\i),-90);
    \begin{scope}[every path/.append style = {#1}]
      \TShapedCross(\k,(0,-\i-\j),-90);
        \TShapedCross(\l,(0,-\i-\j-\k),-90);
        \TShapedCross(\l,(\k,-\i-\j),0);
        \TShapedCross(\l,(-\k,-\i-\j),180);
      \TShapedCross(\k,(\j,-\i),0);
        \TShapedCross(\l,(\j+\k,-\i),0);
        \TShapedCross(\l,(\j,-\i+\k),90);
        \TShapedCross(\l,(\j,-\i-\k),-90);
      \TShapedCross(\k,(-\j,-\i),180);
        \TShapedCross(\l,(-\j-\k,-\i),180);
        \TShapedCross(\l,(-\j,-\i+\k),90);
        \TShapedCross(\l,(-\j,-\i-\k),-90);
    \end{scope}
}
\def\freeGroupCayleyGraphWithoutInnerCross{%
  \begin{scope}[every path/.append style = {dash pattern = on 1pt off 1pt}]
      \TShapedCross(\j,(\i,0),0);
        \TShapedCross(\k,(\i+\j,0),0);
          \TShapedCross(\l,(\i+\j+\k,0),0);
          \TShapedCross(\l,(\i+\j,\k),90);
          \TShapedCross(\l,(\i+\j,-\k),-90);
        \TShapedCross(\k,(\i,\j),90);
          \TShapedCross(\l,(\i,\j+\k),90);
          \TShapedCross(\l,(\i+\k,\j),0);
          \TShapedCross(\l,(\i-\k,\j),180);
        \TShapedCross(\k,(\i,-\j),-90);
          \TShapedCross(\l,(\i,-\j-\k),-90);
          \TShapedCross(\l,(\i+\k,-\j),0);
          \TShapedCross(\l,(\i-\k,-\j),180);
  \end{scope}

  \begin{scope}[every path/.append style = {dash pattern = on 1pt off 1pt on \the\pgflinewidth off 1pt}]
      \TShapedCross(\j,(-\i,0),180);
        \TShapedCross(\k,(-\i-\j,0),180);
          \TShapedCross(\l,(-\i-\j-\k,0),180);
          \TShapedCross(\l,(-\i-\j,\k),90);
          \TShapedCross(\l,(-\i-\j,-\k),-90);
        \TShapedCross(\k,(-\i,\j),90);
          \TShapedCross(\l,(-\i,\j+\k),90);
          \TShapedCross(\l,(-\i+\k,\j),0);
          \TShapedCross(\l,(-\i-\k,\j),180);
        \TShapedCross(\k,(-\i,-\j),-90);
          \TShapedCross(\l,(-\i,-\j-\k),-90);
          \TShapedCross(\l,(-\i+\k,-\j),0);
          \TShapedCross(\l,(-\i-\k,-\j),180);
  \end{scope}

  \begin{scope}[every path/.append style = {densely dotted}]
      \TShapedCross(\j,(0,\i),90);
        \TShapedCross(\k,(0,\i+\j),90);
          \TShapedCross(\l,(0,\i+\j+\k),90);
          \TShapedCross(\l,(\k,\i+\j),0);
          \TShapedCross(\l,(-\k,\i+\j),180);
        \TShapedCross(\k,(\j,\i),0);
          \TShapedCross(\l,(\j+\k,\i),0);
          \TShapedCross(\l,(\j,\i+\k),90);
          \TShapedCross(\l,(\j,\i-\k),-90);
        \TShapedCross(\k,(-\j,\i),180);
          \TShapedCross(\l,(-\j-\k,\i),180);
          \TShapedCross(\l,(-\j,\i+\k),90);
          \TShapedCross(\l,(-\j,\i-\k),-90);
  \end{scope}

    \TShapedCross(\j,(0,-\i),-90);
      \TShapedCross(\k,(0,-\i-\j),-90);
        \TShapedCross(\l,(0,-\i-\j-\k),-90);
        \TShapedCross(\l,(\k,-\i-\j),0);
        \TShapedCross(\l,(-\k,-\i-\j),180);
      \TShapedCross(\k,(\j,-\i),0);
        \TShapedCross(\l,(\j+\k,-\i),0);
        \TShapedCross(\l,(\j,-\i+\k),90);
        \TShapedCross(\l,(\j,-\i-\k),-90);
      \TShapedCross(\k,(-\j,-\i),180);
        \TShapedCross(\l,(-\j-\k,-\i),180);
        \TShapedCross(\l,(-\j,-\i+\k),90);
        \TShapedCross(\l,(-\j,-\i-\k),-90);
}

\def\freeGroupBallRadiusOne{%
  \draw (0, 0) -- node[above] {$a$} (\i, 0); 

  \begin{scope}[every path/.append style = {dotted}]
    \TShapedCross(\j,(\i,0),0);

      \TShapedCross(\k,(\i+\j,0),0);
      \TShapedCross(\k,(\i,\j),90);
      \TShapedCross(\k,(\i,-\j),-90);

        \TShapedCross(\l,(\i+\j+\k,0),0);
        \TShapedCross(\l,(\i+\j,\k),90);
        \TShapedCross(\l,(\i+\j,-\k),-90);

        \TShapedCross(\l,(\i,\j+\k),90);
        \TShapedCross(\l,(\i+\k,\j),0);
        \TShapedCross(\l,(\i-\k,\j),180);

        \TShapedCross(\l,(\i,-\j-\k),-90);
        \TShapedCross(\l,(\i+\k,-\j),0);
        \TShapedCross(\l,(\i-\k,-\j),180);
  \end{scope}

  \draw (0, 0) -- node[below] {$a^{-1}$} (-\i, 0); 

  \begin{scope}[every path/.append style = {dotted}]
    \TShapedCross(\j,(-\i,0),180);
      \TShapedCross(\k,(-\i-\j,0),180);
      \TShapedCross(\k,(-\i,\j),90);
      \TShapedCross(\k,(-\i,-\j),-90);

        \TShapedCross(\l,(-\i-\j-\k,0),180);
        \TShapedCross(\l,(-\i-\j,\k),90);
        \TShapedCross(\l,(-\i-\j,-\k),-90);

        \TShapedCross(\l,(-\i,\j+\k),90);
        \TShapedCross(\l,(-\i+\k,\j),0);
        \TShapedCross(\l,(-\i-\k,\j),180);

        \TShapedCross(\l,(-\i,-\j-\k),-90);
        \TShapedCross(\l,(-\i+\k,-\j),0);
        \TShapedCross(\l,(-\i-\k,-\j),180);
  \end{scope}

  \draw (0, 0) -- node[left] {$b$} (0, \i); 

  \begin{scope}[every path/.append style = {dotted}]
    \TShapedCross(\j,(0,\i),90);
      \TShapedCross(\k,(0,\i+\j),90);
      \TShapedCross(\k,(\j,\i),0);
      \TShapedCross(\k,(-\j,\i),180);

        \TShapedCross(\l,(0,\i+\j+\k),90);
        \TShapedCross(\l,(\k,\i+\j),0);
        \TShapedCross(\l,(-\k,\i+\j),180);

        \TShapedCross(\l,(\j+\k,\i),0);
        \TShapedCross(\l,(\j,\i+\k),90);
        \TShapedCross(\l,(\j,\i-\k),-90);

        \TShapedCross(\l,(-\j-\k,\i),180);
        \TShapedCross(\l,(-\j,\i+\k),90);
        \TShapedCross(\l,(-\j,\i-\k),-90);
  \end{scope}

  \draw (0, 0) -- node[right] {$b^{-1}$} (0, -\i); 

  \begin{scope}[every path/.append style = {dotted}]
    \TShapedCross(\j,(0,-\i),-90);
      \TShapedCross(\k,(0,-\i-\j),-90);
      \TShapedCross(\k,(\j,-\i),0);
      \TShapedCross(\k,(-\j,-\i),180);

        \TShapedCross(\l,(0,-\i-\j-\k),-90);
        \TShapedCross(\l,(\k,-\i-\j),0);
        \TShapedCross(\l,(-\k,-\i-\j),180);

        \TShapedCross(\l,(\j+\k,-\i),0);
        \TShapedCross(\l,(\j,-\i+\k),90);
        \TShapedCross(\l,(\j,-\i-\k),-90);

        \TShapedCross(\l,(-\j-\k,-\i),180);
        \TShapedCross(\l,(-\j,-\i+\k),90);
        \TShapedCross(\l,(-\j,-\i-\k),-90);
  \end{scope}

  \begin{scope}[shift = {(0, 0)}, scale = 1]
    \draw[dashed, tension = 0.5] plot[smooth cycle, tension = 0] coordinates { 
        (0 - \t, -\t)
        (0 + \i - \t, -\t)
        (0 + \i - \t, -\j - \t)
        (0 + \i + \t, -\j - \t)
        (0 + \i + \t, -\t)
        (0 + \i + \j + \t, -\t)
        (0 + \i + \j + \t, \t)
        (0 + \i + \t, \t)
        (0 + \i + \t, \j + \t)
        (0 + \i - \t, \j + \t)
        (0 + \i - \t, \t)
        (0 - \t, \t)
    };
  \end{scope}

  \fill (0, \i) circle (0.1cm);
  \fill (0, -\i) circle (0.1cm);
  \fill (-\i, 0) circle (0.1cm);
}

\def\freeGroupBallRadiusTwo{%
  \draw (0, 0) -- node[above] {$a$} (\i, 0); 
    \TShapedCross(\j,(\i,0),0);

    \begin{scope}[every path/.append style = {densely dotted}]
      \TShapedCross(\k,(\i+\j,0),0);
      \TShapedCross(\k,(\i,\j),90);
      \TShapedCross(\k,(\i,-\j),-90);

        \TShapedCross(\l,(\i+\j+\k,0),0);
        \TShapedCross(\l,(\i+\j,\k),90);
        \TShapedCross(\l,(\i+\j,-\k),-90);

        \TShapedCross(\l,(\i,\j+\k),90);
        \TShapedCross(\l,(\i+\k,\j),0);
        \TShapedCross(\l,(\i-\k,\j),180);

        \TShapedCross(\l,(\i,-\j-\k),-90);
        \TShapedCross(\l,(\i+\k,-\j),0);
        \TShapedCross(\l,(\i-\k,-\j),180);
    \end{scope}

  \draw (0, 0) -- node[below] {$a^{-1}$} (-\i, 0); 
    \TShapedCross(\j,(-\i,0),180);

    \begin{scope}[every path/.append style = {densely dotted}]
      \TShapedCross(\k,(-\i-\j,0),180);
      \TShapedCross(\k,(-\i,\j),90);
      \TShapedCross(\k,(-\i,-\j),-90);

        \TShapedCross(\l,(-\i-\j-\k,0),180);
        \TShapedCross(\l,(-\i-\j,\k),90);
        \TShapedCross(\l,(-\i-\j,-\k),-90);

        \TShapedCross(\l,(-\i,\j+\k),90);
        \TShapedCross(\l,(-\i+\k,\j),0);
        \TShapedCross(\l,(-\i-\k,\j),180);

        \TShapedCross(\l,(-\i,-\j-\k),-90);
        \TShapedCross(\l,(-\i+\k,-\j),0);
        \TShapedCross(\l,(-\i-\k,-\j),180);
    \end{scope}

  \draw (0, 0) -- node[left] {$b$} (0, \i); 
    \TShapedCross(\j,(0,\i),90);

    \begin{scope}[every path/.append style = {densely dotted}]
      \TShapedCross(\k,(0,\i+\j),90);
      \TShapedCross(\k,(\j,\i),0);
      \TShapedCross(\k,(-\j,\i),180);

        \TShapedCross(\l,(0,\i+\j+\k),90);
        \TShapedCross(\l,(\k,\i+\j),0);
        \TShapedCross(\l,(-\k,\i+\j),180);

        \TShapedCross(\l,(\j+\k,\i),0);
        \TShapedCross(\l,(\j,\i+\k),90);
        \TShapedCross(\l,(\j,\i-\k),-90);

        \TShapedCross(\l,(-\j-\k,\i),180);
        \TShapedCross(\l,(-\j,\i+\k),90);
        \TShapedCross(\l,(-\j,\i-\k),-90);
    \end{scope}

  \draw (0, 0) -- node[right] {$b^{-1}$} (0, -\i); 
    \TShapedCross(\j,(0,-\i),-90);

    \begin{scope}[every path/.append style = {densely dotted}]
      \TShapedCross(\k,(0,-\i-\j),-90);
      \TShapedCross(\k,(\j,-\i),0);
      \TShapedCross(\k,(-\j,-\i),180);

        \TShapedCross(\l,(0,-\i-\j-\k),-90);
        \TShapedCross(\l,(\k,-\i-\j),0);
        \TShapedCross(\l,(-\k,-\i-\j),180);

        \TShapedCross(\l,(\j+\k,-\i),0);
        \TShapedCross(\l,(\j,-\i+\k),90);
        \TShapedCross(\l,(\j,-\i-\k),-90);

        \TShapedCross(\l,(-\j-\k,-\i),180);
        \TShapedCross(\l,(-\j,-\i+\k),90);
        \TShapedCross(\l,(-\j,-\i-\k),-90);
    \end{scope}

  \begin{scope}[shift = {(0, 0)}, scale = 1]
    \draw[dashed, tension = 0.5] plot[smooth cycle, tension = 0] coordinates { 
        (0 - \i - \t, -\t)
        (0 - \t, -\t)
        (0 - \t, -\i - \t)
        (0 + \t, -\i - \t)
        (0 + \t, -\t)
        (0 + \i - \t, -\t)
        (0 + \i - \t, -\j + \t)
        (0 + \i - \k - \t, -\j + \t)
        (0 + \i - \k - \t, -\j - \t)
        (0 + \i - \t, -\j - \t)
        (0 + \i - \t, -\j - \k - \t)
        (0 + \i + \t, -\j - \k - \t)
        (0 + \i + \t, -\j - \t)
        (0 + \i + \k + \t, -\j - \t)
        (0 + \i + \k + \t, -\j + \t)
        (0 + \i + \t, -\j + \t)
        (0 + \i + \t, -\t)
        (0 + \i + \j - \t, -\t)
        (0 + \i + \j - \t, -\k - \t)
        (0 + \i + \j + \t, -\k - \t)
        (0 + \i + \j + \t, -\t)
        (0 + \i + \j + \k + \t, -\t)
        (0 + \i + \j + \k + \t, \t)
        (0 + \i + \j + \t, \t)
        (0 + \i + \j + \t, \k + \t)
        (0 + \i + \j - \t, \k + \t)
        (0 + \i + \j - \t, \t)
        (0 + \i + \t, \t)
        (0 + \i + \t, \j - \t)
        (0 + \i + \k + \t, \j - \t)
        (0 + \i + \k + \t, \j + \t)
        (0 + \i + \t, \j + \t)
        (0 + \i + \t, \j + \k + \t)
        (0 + \i - \t, \j + \k + \t)
        (0 + \i - \t, \j + \t)
        (0 + \i - \k - \t, \j + \t)
        (0 + \i - \k - \t, \j - \t)
        (0 + \i - \t, \j - \t)
        (0 + \i - \t, \t)
        (0 + \t, \t)
        (0 + \t, \i + \t)
        (0 - \t, \i + \t)
        (0 - \t, \t)
        (0 - \i - \t, \t)
    };
  \end{scope}

  \fill (0, \i + \j) circle (0.1cm);
    \fill (\j, \i) circle (0.1cm);
    \fill (-\j, \i) circle (0.1cm);
  \fill (0, -\i - \j) circle (0.1cm);
    \fill (\j, -\i) circle (0.1cm);
    \fill (-\j, -\i) circle (0.1cm);
  \fill (-\i - \j, 0) circle (0.1cm);
    \fill (-\i, \j) circle (0.1cm);
    \fill (-\i, -\j) circle (0.1cm);
}

\def\freeGroupBallRadiusThree{%
  \draw (0, 0) -- node[above] {$a$} (\i, 0); 
    \TShapedCross(\j,(\i,0),0);

      \TShapedCross(\k,(\i+\j,0),0);
      \TShapedCross(\k,(\i,\j),90);
      \TShapedCross(\k,(\i,-\j),-90);

      \begin{scope}[every path/.append style = {densely dotted}]
        \TShapedCross(\l,(\i+\j+\k,0),0);
        \TShapedCross(\l,(\i+\j,\k),90);
        \TShapedCross(\l,(\i+\j,-\k),-90);

        \TShapedCross(\l,(\i,\j+\k),90);
        \TShapedCross(\l,(\i+\k,\j),0);
        \TShapedCross(\l,(\i-\k,\j),180);

        \TShapedCross(\l,(\i,-\j-\k),-90);
        \TShapedCross(\l,(\i+\k,-\j),0);
        \TShapedCross(\l,(\i-\k,-\j),180);
      \end{scope}

  \draw (0, 0) -- node[below] {$a^{-1}$} (-\i, 0); 
    \TShapedCross(\j,(-\i,0),180);

      \TShapedCross(\k,(-\i-\j,0),180);
      \TShapedCross(\k,(-\i,\j),90);
      \TShapedCross(\k,(-\i,-\j),-90);

      \begin{scope}[every path/.append style = {densely dotted}]
        \TShapedCross(\l,(-\i-\j-\k,0),180);
        \TShapedCross(\l,(-\i-\j,\k),90);
        \TShapedCross(\l,(-\i-\j,-\k),-90);

        \TShapedCross(\l,(-\i,\j+\k),90);
        \TShapedCross(\l,(-\i+\k,\j),0);
        \TShapedCross(\l,(-\i-\k,\j),180);

        \TShapedCross(\l,(-\i,-\j-\k),-90);
        \TShapedCross(\l,(-\i+\k,-\j),0);
        \TShapedCross(\l,(-\i-\k,-\j),180);
      \end{scope}

  \draw (0, 0) -- node[left] {$b$} (0, \i); 
    \TShapedCross(\j,(0,\i),90);

      \TShapedCross(\k,(0,\i+\j),90);
      \TShapedCross(\k,(\j,\i),0);
      \TShapedCross(\k,(-\j,\i),180);
      \begin{scope}[every path/.append style = {densely dotted}]
        \TShapedCross(\l,(0,\i+\j+\k),90);
        \TShapedCross(\l,(\k,\i+\j),0);
        \TShapedCross(\l,(-\k,\i+\j),180);

        \TShapedCross(\l,(\j+\k,\i),0);
        \TShapedCross(\l,(\j,\i+\k),90);
        \TShapedCross(\l,(\j,\i-\k),-90);

        \TShapedCross(\l,(-\j-\k,\i),180);
        \TShapedCross(\l,(-\j,\i+\k),90);
        \TShapedCross(\l,(-\j,\i-\k),-90);
      \end{scope}

  \draw (0, 0) -- node[right] {$b^{-1}$} (0, -\i); 
    \TShapedCross(\j,(0,-\i),-90);

      \TShapedCross(\k,(0,-\i-\j),-90);
      \TShapedCross(\k,(\j,-\i),0);
      \TShapedCross(\k,(-\j,-\i),180);

      \begin{scope}[every path/.append style = {densely dotted}]
        \TShapedCross(\l,(0,-\i-\j-\k),-90);
        \TShapedCross(\l,(\k,-\i-\j),0);
        \TShapedCross(\l,(-\k,-\i-\j),180);

        \TShapedCross(\l,(\j+\k,-\i),0);
        \TShapedCross(\l,(\j,-\i+\k),90);
        \TShapedCross(\l,(\j,-\i-\k),-90);

        \TShapedCross(\l,(-\j-\k,-\i),180);
        \TShapedCross(\l,(-\j,-\i+\k),90);
        \TShapedCross(\l,(-\j,-\i-\k),-90);
      \end{scope}

  \foreach \r in {0, 90, 180} {
    \begin{scope}[rotate = \r]
      \fill (0, \i + \j + \k) circle (0.1cm);
      \fill (\k, \i + \j) circle (0.1cm);
      \fill (-\k, \i + \j) circle (0.1cm);

      \fill (\j, \i + \k) circle (0.1cm);
      \fill (\j, \i - \k) circle (0.1cm);
      \fill (\j + \k, \i) circle (0.1cm);

      \fill (-\j, \i + \k) circle (0.1cm);
      \fill (-\j, \i - \k) circle (0.1cm);
      \fill (-\j - \k, \i) circle (0.1cm);
    \end{scope}
  }
}

\def\PlusShapedCross(#1,(#2,#3)){%
  \begin{scope}[shift = {(#2, #3)}]
    \fill (0, 0) circle (1.333 * 0.1cm);
    \draw (0, 0) edge node[auto] {$a$} (#1, 0);
    \draw (0, 0) edge node[auto] {$a^{-1}$} (-#1, 0);
    \draw (0, 0) edge node[auto] {$b$} (0, #1);
    \draw (0, 0) edge node[auto] {$b^{-1}$} (0, -#1);
  \end{scope}%
}
\def\PlusShapedCrossWithOneLongerArm(#1,#2,(#3,#4),#5){%
  \begin{scope}[shift = {(#3, #4)}, rotate = #5]
    \fill (0, 0) circle (1.333 * #1 * 0.1cm);
    \draw (0, 0) -- (#2, 0);
    \draw (0, 0) -- (-#1, 0);
    \draw (0, 0) -- (0, #2);
    \draw (0, 0) -- (0, -#2);
  \end{scope}%
}
\def\HShapedCross(#1,#2,(#3,#4),#5){%
  \begin{scope}[shift = {(#3, #4)}, rotate = #5, every path/.append style = {densely dotted}]
    \draw (0, 0) -- (#1, 0)
          (0, 0) -- (-#1, 0);
    \draw (#1, 0) -- (#1, #2)
          (#1, 0) -- (#1, -#2);
    \draw (-#1, 0) -- (-#1, #2)
          (-#1, 0) -- (-#1, -#2);
  \end{scope}%
}
\def\TShapedCrossDotted(#1,(#2,#3),#4){%
  \begin{scope}[shift = {(#2, #3)}, rotate = #4, every path/.append style = {densely dotted}]
    \draw (0, 0) -- (#1, 0)
          (0, 0) -- (0, #1)
          (0, 0) -- (0, -#1);
  \end{scope}%
}
\def\TShapedCrossHalfDotted(#1,(#2,#3),#4){%
  \begin{scope}[shift = {(#2, #3)}, rotate = #4]
    \fill (#1, 0) circle (1.333 * #1 * 0.1cm);
    \draw (0, 0) -- (#1, 0);
    \draw[densely dotted]
          (0, 0) -- (0, #1)
          (0, 0) -- (0, -#1);
  \end{scope}%
}
\def\FreeGroupTilingOneArm((#1,#2),#3){%
  \begin{scope}[shift = {(#1, #2)}, rotate = #3]
    \PlusShapedCrossWithOneLongerArm(\k,\l,(0+\j+\k,0),0);

    \PlusShapedCrossWithOneLongerArm(\l,\m,(0+\j,\k+\l),90);
    \PlusShapedCrossWithOneLongerArm(\l,\m,(0+\j,-\k-\l),-90);

    \begin{scope}[shift = {(0, 0)}]
      \draw[densely dotted]
          (0, 0) -- (\j, 0);
    \end{scope}
    \HShapedCross(\k,\l,(0+\j,0),90);

    \TShapedCrossDotted(\m,(0+\j+\k+\l,0),0);
    \TShapedCrossDotted(\m,(0+\j+\k,\l),90);
    \TShapedCrossDotted(\m,(0+\j+\k,-\l),-90);

    \TShapedCrossHalfDotted(\m,(0+\j+\l,\k),0);
    \TShapedCrossHalfDotted(\m,(0+\j-\l,\k),180);

    \TShapedCrossHalfDotted(\m,(0+\j+\l,-\k),0);
    \TShapedCrossHalfDotted(\m,(0+\j-\l,-\k),180);
  \end{scope}
}

\def\lineSolidDottedOfVaryingSizesWithArms(#1,#2,(#3,#4),#5){%
  \begin{scope}[shift = {(#3,#4)}, rotate = #5]
    \fill (0, 0) circle (0.666 * #1 * 0.1cm);
    \draw (0, 0) -- (#1, 0);
    \draw[densely dotted]
          (#1, 0) -- (#1 + #2, 0)
          (#1, 0) -- (#1, #2)
          (#1, 0) -- (#1, -#2);
    \fill (#1, #2) circle (0.666 * #2 * 0.1cm);
    \fill (#1, -#2) circle (0.666 * #2 * 0.1cm);
  \end{scope}
}
\def\lineSolidDottedSolidOfVaryingSizesWithArms(#1,#2,#3,(#4,#5),#6){%
  \lineSolidDottedOfVaryingSizesWithArms(#1,#2,(#4,#5),#6)
  \begin{scope}[shift = {(#4,#5)}, rotate = #6]
    \fill (#1 + #2 + #3, 0) circle (0.666 * #3 * 0.1cm);
    \draw (#1 + #2, 0) -- (#1 + #2 + #3, 0);
    \draw[densely dotted]
          (#1 + #2, 0) -- (#1 + #2, #3)
          (#1 + #2, 0) -- (#1 + #2, -#3);
    \fill (#1 + #2, #3) circle (0.666 * #3 * 0.1cm);
    \fill (#1 + #2, -#3) circle (0.666 * #3 * 0.1cm);
  \end{scope}
}

\def\figureOr{%
  \begin{tikzpicture}[scale = 0.55, y = -1cm, arrow/.style = {->, shorten >= 0.2cm, shorten <= 0.2cm}, composed/.style = {densely dotted}]
    \draw (-3.5, 0) grid (4.5, 2);
    \foreach \x in {-3, ..., 3} {
      \node at (\x + 0.5, -0.5) {$\x$};
      \node at (\x + 0.5, 1.5) {$\scriptstyle\lor$};
    }
    \foreach \x in {-4, ..., 3} {
      \draw[arrow] (\x + 1.5, 0.5) -- (\x + 0.5, 1.5);
      \draw[arrow] (\x + 0.5, 0.5) -- (\x + 1.5, 1.5);
    }
    \draw (-4, 0.5) edge[|->, bend right = 90] node[left] {$\Delta$} (-4, 1.5);
  \end{tikzpicture}
}

\def\figureShifts{%
  \begin{minipage}{.5\textwidth}
    \myfloatalign
    \subfloat[]{%
      \begin{tikzpicture}[scale = 0.55, y = -1cm, arrow/.style = {->, shorten >= 0.2cm, shorten <= 0.2cm}, composed/.style = {densely dotted}]
        \draw (-3.5, 0) grid (4.5, 2);
        \foreach \x in {-3, ..., 3} {
          \node at (\x + 0.5, -0.5) {$\x$};
        }
        \foreach \x in {-3, ..., 4} {
          \draw[arrow] (\x + 0.5, 0.5) -- (\x - 0.5, 1.5);
        }
        \draw (-4, 0.5) edge[|->, bend right = 90] node[left] {$\Delta$} (-4, 1.5);
      \end{tikzpicture}
      \label{figure:left-shift}
    }
  \end{minipage}%
  \begin{minipage}{.5\textwidth}
    \myfloatalign
    \subfloat[]{%
      \begin{tikzpicture}[scale = 0.55, y = -1cm, arrow/.style = {->, shorten >= 0.2cm, shorten <= 0.2cm}, composed/.style = {densely dotted}]
        \draw (-3.5, 0) grid (4.5, 2);
        \foreach \x in {-3, ..., 3} {
          \node at (\x + 0.5, -0.5) {$\x$};
        }
        \foreach \x in {-3, ..., 4} {
          \draw[arrow] (\x - 0.5, 0.5) -- (\x + 0.5, 1.5);
        }
        \draw (-4, 0.5) edge[|->, bend right = 90] node[left] {$\Delta'$} (-4, 1.5);
      \end{tikzpicture}
      \label{figure:right-shift}
    }
  \end{minipage}
  \begin{minipage}{.5\textwidth}
    \myfloatalign
    \subfloat[]{%
      \begin{tikzpicture}[scale = 0.55, y = -1cm, arrow/.style = {->, shorten >= 0.2cm, shorten <= 0.2cm}, composed/.style = {densely dotted}]
        \draw (-3.5, 0) grid (4.5, 2);
        \foreach \x in {-3, ..., 3} {
          \node at (\x + 0.5, -0.5) {$\x$};
        }
        \foreach \x in {-3, -1, ..., 3} {
          \draw[arrow] (\x - 0.5, 0.5) -- (\x + 0.5, 1.5);
        }
        \foreach \x in {-3, -1, ..., 4} {
          \draw[arrow] (\x + 0.5, 0.5) -- (\x - 0.5, 1.5);
        }
        \draw (-4, 0.5) edge[|->, bend right = 90] node[left] {$\Delta''$} (-4, 1.5);
      \end{tikzpicture}
      \label{figure:odd-even-shift}
    }
  \end{minipage}%
  \begin{minipage}{.5\textwidth}
    \myfloatalign
    \subfloat[]{%
      \begin{tikzpicture}[scale = 0.55, y = -1cm, arrow/.style = {->, shorten >= 0.2cm, shorten <= 0.2cm}, composed/.style = {densely dotted}]
        \draw (-3.5, 0) grid (4.5, 2);
        \foreach \x in {-3, ..., 3} {
          \node at (\x + 0.5, -0.5) {$\x$};
        }
        \foreach \x in {-4, -3, -2, -1, 0, 2, 3} {
          \draw[arrow] (\x + 1.5, 0.5) -- (\x + 0.5, 1.5);
        }
        \draw[arrow] (1 - 0.5, 0.5) -- (1 + 0.5, 1.5);
        \draw (-4, 0.5) edge[|->, bend right = 90] node[left] {$\Delta'''$} (-4, 1.5);
      \end{tikzpicture}
      \label{figure:strange-shift}
    }
  \end{minipage}
}

\def\figureCellularAutomatonTwoDimensionalIntegerLatticeLocalConfiguraitonOfOrigin{%
  \begin{wide}
    \begin{tikzpicture}[scale=0.06]
      \vonNeumannNeighbourhoodCenter(white,(0,0))
      \vonNeumannNeighbourhoodCenter(black,(128,0))
    \end{tikzpicture}
  \end{wide}
}

\def\figureExampleLocalConfigurationConstraintsHorizontalAndVerticalReflections{%
  \begin{tikzpicture}[scale=0.06]
    \vonNeumannNeighbourhoodCenterHorizontalAndVerticalReflections(white,(0,0))
    \vonNeumannNeighbourhoodCenterHorizontalAndVerticalReflections(black,(72,0))
  \end{tikzpicture}
}

\def\figureExampleLocalConfigurationConstraintsDiagonalReflections{%
  \begin{tikzpicture}[scale=0.06]
    \vonNeumannNeighbourhoodCenterDiagonalReflections(white,(0,0))
    \vonNeumannNeighbourhoodCenterDiagonalReflections(black,(56,0))
  \end{tikzpicture}
}

\def\figureExampleLocalConfigurationConstraintsRotations{%
  \begin{tikzpicture}[scale=0.06]
    \vonNeumannNeighbourhoodCenterRotations(white,(0,0))
    \vonNeumannNeighbourhoodCenterRotations(black,(48,0))
  \end{tikzpicture}
}

\def\figureExampleLocalConfigurationConstraintsPointReflection{%
  \begin{tikzpicture}[scale=0.06]
    \vonNeumannNeighbourhoodCenterPointReflections(white,(0,0))
    \vonNeumannNeighbourhoodCenterPointReflections(black,(72,0))
  \end{tikzpicture}
}

\def\figureStrangeShiftSquared{%
  \begin{tikzpicture}[scale = 0.55, y = -1cm, arrow/.style = {->, shorten >= 0.2cm, shorten <= 0.2cm}, composed/.style = {densely dotted}]
    \draw (-3.5, 0) grid (4.5, 3);
    \foreach \x in {-3, ..., 3} {
      \node at (\x + 0.5, -0.5) {$\x$};
    }
    \foreach \x in {-4, -3, -2, -1, 0, 2, 3} {
      \foreach \y in {0, 1} {
        \draw[arrow] (\x + 1.5, \y + 0.5) -- (\x + 0.5, \y + 1.5);
      }
    }
    \foreach \x in {-4, -3, -2, -1, 2} {
      \draw[arrow, composed] (\x + 2.5, 0 + 0.5) -- (\x + 0.5, 1 + 1.5);
    }
    \foreach \y in {0, 1} {
      \draw[arrow] (1 - 0.5, \y + 0.5) -- (1 + 0.5, \y + 1.5);
      \draw[arrow, composed, shorten >= 0.175cm, shorten <= 0.175cm] (\y + 0.5, 0.5) -- (\y + 0.5, 2.5);
    }
    \draw (-4, 0.5) edge[|->, bend right = 90] node[left] {$\Delta'''$} (-4, 1.25);
    \draw (-4, 1.75) edge[|->, bend right = 90] node[left] {$\Delta'''$} (-4, 2.5);
  \end{tikzpicture}
}

\def\figureCompositionsOfGlobalTransitionFunctions{%
  \begin{minipage}{.5\textwidth}
    \myfloatalign
    \subfloat[]{%
      \begin{tikzpicture}[scale = 0.55, y = -1cm, arrow/.style = {->, shorten >= 0.2cm, shorten <= 0.2cm}, composed/.style = {densely dotted}]
        \draw (-3.5, 0) grid (4.5, 3);
        \foreach \x in {-3, ..., 3} {
          \node at (\x + 0.5, -0.5) {$\x$};
        }
        \foreach \x in {-4, ..., 3} {
          \draw[arrow] (\x + 1.5, 0.5) -- (\x + 0.5, 1.5);
        }
        \foreach \x in {-4, ..., 3} {
          \draw[arrow] (\x + 1.5, 1 + 0.5) -- (\x + 0.5, 1 + 1.5);
        }
        \foreach \x in {-4, ..., 2} {
          \draw[arrow, composed] (\x + 2.5, 0 + 0.5) -- (\x + 0.5, 1 + 1.5);
        }
        \draw (-4, 0.5) edge[|->, bend right = 90] node[left] {$\Delta$} (-4, 1.25);
        \draw (-4, 1.75) edge[|->, bend right = 90] node[left] {$\Delta$} (-4, 2.5);
      \end{tikzpicture}
      \label{figure:left-shift-followed-by-left-shift}
    }
  \end{minipage}%
  \hfill
  \begin{minipage}{.5\textwidth}
    \myfloatalign
    \subfloat[]{%
      \begin{tikzpicture}[scale = 0.55, y = -1cm, arrow/.style = {->, shorten >= 0.2cm, shorten <= 0.2cm}, composed/.style = {densely dotted}]
        \draw (-3.5, 0) grid (4.5, 3);
        \foreach \x in {-3, ..., 3} {
          \node at (\x + 0.5, -0.5) {$\x$};
          \node at (\x + 0.5, 1.5) {$\scriptstyle\lor$};
          \node at (\x + 0.5, 2.5) {$\scriptstyle\lor$};
        }
        \foreach \y in {0, 1} {
          \foreach \x in {-4, ..., 3} {
            \draw[arrow] (\x + 1.5, \y + 0.5) -- (\x + 0.5, \y + 1.5);
            \draw[arrow] (\x + 0.5, \y + 0.5) -- (\x + 1.5, \y + 1.5);
          }
        }
        \foreach \x in {-4, ..., 2} {
          \draw[arrow, composed] (\x + 2.5, 0.5) -- (\x + 0.5, 2.5);
          \draw[arrow, composed] (\x + 0.5, 0.5) -- (\x + 2.5, 2.5);
          \draw[arrow, composed, shorten >= 0.125cm, shorten <= 0.125cm] (\x + 1.5, 0.5) -- (\x + 1.5, 2.5);
        }
        \draw (-4, 0.5) edge[|->, bend right = 90] node[left] {$\Delta''$} (-4, 1.25);
        \draw (-4, 1.75) edge[|->, bend right = 90] node[left] {$\Delta''$} (-4, 2.5);
      \end{tikzpicture}
      \label{figure:or-followed-by-or}
    }
  \end{minipage}
  \begin{minipage}{.5\textwidth}
    \myfloatalign
    \subfloat[]{%
      \begin{tikzpicture}[scale = 0.55, y = -1cm, arrow/.style = {->, shorten >= 0.2cm, shorten <= 0.2cm}, composed/.style = {densely dotted}]
        \draw (-3.5, 0) grid (4.5, 3);
        \foreach \x in {-3, ..., 3} {
          \node at (\x + 0.5, -0.5) {$\x$};
          \node at (\x + 0.5, 2.5) {$\scriptstyle\lor$};
        }
        \foreach \x in {-4, ..., 3} {
          \draw[arrow] (\x + 1.5, 0.5) -- (\x + 0.5, 1.5);
        }
        \foreach \x in {-4, ..., 3} {
          \draw[arrow] (\x + 1.5, 1.5) -- (\x + 0.5, 2.5);
          \draw[arrow] (\x + 0.5, 1.5) -- (\x + 1.5, 2.5);
        }
        \foreach \x in {-4, ..., 2} {
          \draw[arrow, composed] (\x + 2.5, 0 + 0.5) -- (\x + 0.5, 1 + 1.5);
        }
        \foreach \x in {-3, ..., 3} {
          \draw[arrow, composed, shorten >= 0.125cm, shorten <= 0.125cm] (\x + 0.5, 0 + 0.5) -- (\x + 0.5, 1 + 1.5);
        }
        \draw (-4, 0.5) edge[|->, bend right = 90] node[left] {$\Delta$} (-4, 1.25);
        \draw (-4, 1.75) edge[|->, bend right = 90] node[left] {$\Delta''$} (-4, 2.5);
      \end{tikzpicture}
      \label{figure:left-shift-followed-by-or}
    }
  \end{minipage}%
  \hfill
  \begin{minipage}{.5\textwidth}
    \myfloatalign
    \subfloat[]{%
      \begin{tikzpicture}[scale = 0.55, y = -1cm, arrow/.style = {->, shorten >= 0.2cm, shorten <= 0.2cm}, composed/.style = {densely dotted}]
        \draw (-3.5, 0) grid (4.5, 3);
        \foreach \x in {-3, ..., 3} {
          \node at (\x + 0.5, -0.5) {$\x$};
          \node at (\x + 0.5, 1.5) {$\scriptstyle\lor$};
          \node at (\x + 0.5, 2.5) {$\scriptstyle\lor$};
        }
        \foreach \x in {-4, ..., 3} {
          \draw[arrow] (\x + 1.5, 1.5) -- (\x + 0.5, 2.5);
        }
        \foreach \x in {-4, ..., 3} {
          \draw[arrow] (\x + 1.5, 0.5) -- (\x + 0.5, 1.5);
          \draw[arrow] (\x + 0.5, 0.5) -- (\x + 1.5, 1.5);
        }
        \foreach \x in {-4, ..., 2} {
          \draw[arrow, composed] (\x + 2.5, 0 + 0.5) -- (\x + 0.5, 1 + 1.5);
        }
        \foreach \x in {-3, ..., 3} {
          \draw[arrow, composed, shorten >= 0.125cm, shorten <= 0.125cm] (\x + 0.5, 0 + 0.5) -- (\x + 0.5, 1 + 1.5);
        }
        \draw (-4, 0.5) edge[|->, bend right = 90] node[left] {$\Delta''$} (-4, 1.25);
        \draw (-4, 1.75) edge[|->, bend right = 90] node[left] {$\Delta$} (-4, 2.5);
      \end{tikzpicture}
      \label{figure:or-followed-by-left-shift}
    }
  \end{minipage}
  \begin{minipage}{.5\textwidth}
    \myfloatalign
    \subfloat[]{%
      \begin{tikzpicture}[scale = 0.55, y = -1cm, arrow/.style = {->, shorten >= 0.2cm, shorten <= 0.2cm}, composed/.style = {densely dotted}]
        \draw (-3.5, 0) grid (4.5, 3);
        \foreach \x in {-3, ..., 3} {
          \node at (\x + 0.5, -0.5) {$\x$};
        }
        \foreach \x in {-4, ..., 3} {
          \draw[arrow] (\x + 1.5, 0.5) -- (\x + 0.5, 1.5);
        }
        \foreach \x in {-4, -3, -2, -1, 0, 2, 3} {
          \draw[arrow] (\x + 1.5, 1 + 0.5) -- (\x + 0.5, 1 + 1.5);
        }
        \foreach \x in {-4, -3, -2, -1, 0, 2} {
          \draw[arrow, composed] (\x + 2.5, 0 + 0.5) -- (\x + 0.5, 1 + 1.5);
        }
        \draw[arrow] (1 - 0.5, 1.5) -- (1 + 0.5, 2.5);
        \draw[arrow, composed, shorten >= 0.125cm, shorten <= 0.125cm] (1.5, 0.5) -- (1.5, 2.5);
        \draw (-4, 0.5) edge[|->, bend right = 90] node[left] {$\Delta$} (-4, 1.25);
        \draw (-4, 1.75) edge[|->, bend right = 90] node[left] {$\Delta'$} (-4, 2.5);
      \end{tikzpicture}
      \label{figure:left-shift-followed-by-strange-shift}
    }
  \end{minipage}%
  \hfill
  \begin{minipage}{.5\textwidth}
    \myfloatalign
    \subfloat[]{%
      \begin{tikzpicture}[scale = 0.55, y = -1cm, arrow/.style = {->, shorten >= 0.2cm, shorten <= 0.2cm}, composed/.style = {densely dotted}]
        \draw (-3.5, 0) grid (4.5, 3);
        \foreach \x in {-3, ..., 3} {
          \node at (\x + 0.5, -0.5) {$\x$};
        }
        \foreach \x in {-4, ..., 3} {
          \draw[arrow] (\x + 1.5, 1 + 0.5) -- (\x + 0.5, 1 + 1.5);
        }
        \foreach \x in {-4, -3, -2, -1, 0, 2, 3} {
          \draw[arrow] (\x + 1.5, 0.5) -- (\x + 0.5, 1.5);
        }
        \foreach \x in {-4, -3, -2, -1, 1, 2} {
          \draw[arrow, composed] (\x + 2.5, 0 + 0.5) -- (\x + 0.5, 1 + 1.5);
        }
        \draw[arrow] (1 - 0.5, 0.5) -- (1 + 0.5, 1.5);
        \draw[arrow, composed, shorten >= 0.125cm, shorten <= 0.125cm] (0.5, 0.5) -- (0.5, 2.5);
        \draw (-4, 0.5) edge[|->, bend right = 90] node[left] {$\Delta'$} (-4, 1.25);
        \draw (-4, 1.75) edge[|->, bend right = 90] node[left] {$\Delta$} (-4, 2.5);
      \end{tikzpicture}
      \label{figure:strange-shift-followed-by-left-shift}
    }
  \end{minipage}
  \begin{minipage}{.5\textwidth}
    \myfloatalign
    \subfloat[]{%
      \begin{tikzpicture}[scale = 0.55, y = -1cm, arrow/.style = {->, shorten >= 0.2cm, shorten <= 0.2cm}, composed/.style = {densely dotted}]
        \draw (-3.5, 0) grid (4.5, 3);
        \foreach \x in {-3, ..., 3} {
          \node at (\x + 0.5, -0.5) {$\x$};
          \node at (\x + 0.5, 2.5) {$\scriptstyle\lor$};
        }
        \foreach \x in {-4, -3, -2, -1, 0, 2, 3} {
          \draw[arrow] (\x + 1.5, 0.5) -- (\x + 0.5, 1.5);
        }
        \foreach \x in {-4, ..., 3} {
          \draw[arrow] (\x + 1.5, 1.5) -- (\x + 0.5, 2.5);
          \draw[arrow] (\x + 0.5, 1.5) -- (\x + 1.5, 2.5);
        }
        \foreach \x in {-4, -3, -2, -1, 1, 2} {
          \draw[arrow, composed] (\x + 2.5, 0 + 0.5) -- (\x + 0.5, 1 + 1.5);
        }
        \foreach \x in {-3, -2, -1, 1, 3} {
          \draw[arrow, composed, shorten >= 0.125cm, shorten <= 0.125cm] (\x + 0.5, 0 + 0.5) -- (\x + 0.5, 1 + 1.5);
        }
        \draw[arrow, composed] (0.5, 0.5) -- (2.5, 2.5);
        \draw[arrow] (1 - 0.5, 0.5) -- (1 + 0.5, 1.5);
        \draw[arrow, composed, shorten >= 0.125cm, shorten <= 0.125cm] (0.4, 0.5) -- (0.4, 2.5);
        \draw[arrow, composed, shorten >= 0.125cm, shorten <= 0.125cm] (0.6, 0.5) -- (0.6, 2.5);
        \draw (-4, 0.5) edge[|->, bend right = 90] node[left] {$\Delta'$} (-4, 1.25);
        \draw (-4, 1.75) edge[|->, bend right = 90] node[left] {$\Delta''$} (-4, 2.5);
      \end{tikzpicture}
      \label{figure:strange-shift-followed-by-or}
    }
  \end{minipage}%
  \hfill
  \begin{minipage}{.5\textwidth}
    \myfloatalign
    \subfloat[]{%
      \begin{tikzpicture}[scale = 0.55, y = -1cm, arrow/.style = {->, shorten >= 0.2cm, shorten <= 0.2cm}, composed/.style = {densely dotted}]
        \draw (-3.5, 0) grid (4.5, 3);
        \foreach \x in {-3, ..., 3} {
          \node at (\x + 0.5, -0.5) {$\x$};
          \node at (\x + 0.5, 1.5) {$\scriptstyle\lor$};
          \node at (\x + 0.5, 2.5) {$\scriptstyle\lor$};
        }
        \foreach \x in {-4, -3, -2, -1, 0, 2, 3} {
          \draw[arrow] (\x + 1.5, 1.5) -- (\x + 0.5, 2.5);
        }
        \draw[arrow] (1 - 0.5, 1.5) -- (1 + 0.5, 2.5);
        \foreach \x in {-4, ..., 3} {
          \draw[arrow] (\x + 1.5, 0.5) -- (\x + 0.5, 1.5);
          \draw[arrow] (\x + 0.5, 0.5) -- (\x + 1.5, 1.5);
        }
        \foreach \x in {-4, -3, -2, -1, 0, 2} {
          \draw[arrow, composed] (\x + 2.5, 0 + 0.5) -- (\x + 0.5, 1 + 1.5);
        }
        \foreach \x in {-3, -2, -1, 0, 1, 2, 3} {
          \draw[arrow, composed, shorten >= 0.125cm, shorten <= 0.125cm] (\x + 0.5, 0 + 0.5) -- (\x + 0.5, 1 + 1.5);
        }
        \draw[arrow, composed] (-0.5, 0.5) -- (1.5, 2.5);
        \draw (-4, 0.5) edge[|->, bend right = 90] node[left] {$\Delta''$} (-4, 1.25);
        \draw (-4, 1.75) edge[|->, bend right = 90] node[left] {$\Delta'$} (-4, 2.5);
      \end{tikzpicture}
      \label{figure:or-followed-by-strange-shift}
    }
  \end{minipage}
}

\def\figureToriPointsArcsAndRibbons{%
  \begin{minipage}{.5\linewidth}
    \myfloatalign
    \subfloat[]{%
      \begin{tikzpicture}[scale = 0.55]
        \begin{scope}[shift = {(0, 0)}]
          \draw (0, 0) -- (4, 0);
          \draw (0, 4) -- (4, 4);
          \draw[dashed] (0, 2) -- (0, 4);
          \draw[dashed] (4, 0) -- (4, 4);

          \node[circle, fill, inner sep = 0pt, minimum size = 2pt] (m) at (0.5, 0) [label = -90:$m$] {};
          \node[circle, fill, inner sep = 0pt, minimum size = 2pt] (m') at (2.5, 0) [label = -90:$m'$] {};

          \draw[very thick] (m) -- (m');
          \draw[very thick, dashed] (0, 0) -- (0, 2);

          \node[circle, fill, inner sep = 0pt, minimum size = 2pt] at (0.5, 2) {};
          \node[circle, fill, inner sep = 0pt, minimum size = 2pt] (a) at (0, 2) [label = 180:$a$] {};
        \end{scope}
        \begin{scope}[shift = {(5.5, 0)}]
          \draw (0, 0) -- (4, 0);
          \draw (0, 4) -- (4, 4);
          \draw[dashed] (0, 0) -- (0, 3);
          \draw[dashed] (4, 0) -- (4, 4);

          \pgfmathsetmacro{\r}{{4 / (2 * pi)}}
          \begin{scope}[shift = {(2, 4.5 + \r)}]
            \draw [dashdotted, domain = 90:270] plot ({\r * cos(\x)}, {\r * sin(\x)});
            \draw [dotted, domain = -90:90] plot ({\r * cos(\x)}, {\r * sin(\x)});

            \draw[dashdotted] (-2, -\r) -- (0, -\r);
            \draw[dotted] (0, -\r) -- (2, -\r);
          \end{scope}

          \node[circle, fill, inner sep = 0pt, minimum size = 2pt] (m) at (0.5, 0) [label = -90:$m$] {};
          \node[circle, fill, inner sep = 0pt, minimum size = 2pt] (m') at (2.5, 0) [label = -90:$m'$] {};
          \node[circle, fill, inner sep = 0pt, minimum size = 2pt] (mr) at (3.5, 0) [label = 90:$\varrho_{-1}(m)$] {};

          \draw[very thick] (mr) -- (m');
          \draw[very thick, dashed] (0, 4) -- (0, 3);

          \node[circle, fill, inner sep = 0pt, minimum size = 2pt] at (0.5, 3) {};
          \node[circle, fill, inner sep = 0pt, minimum size = 2pt] (a') at (0, 3) [label = 180:$a'$] {};
        \end{scope}
      \end{tikzpicture}
      \label{figure:tori:points}
    }
  \end{minipage}%
  \begin{minipage}{.5\linewidth}
    \myfloatalign
    \subfloat[]{%
      \begin{tikzpicture}[scale = 0.55]
        \begin{scope}[shift = {(0, 0)}]
          \draw (0, 0) -- (4, 0);
          \draw (0, 4) -- (4, 4);
          \draw[dashed] (0, 0) -- (0, 4);
          \draw[dashed] (4, 0) -- (4, 4);

          \node[circle, fill, inner sep = 0pt, minimum size = 2pt] (m) at (0.5, 0) [label = -90:$m$] {};
          \node[yscale = 0.6] (Yl) at (2.5, 0) {$\{$};
          \node[yscale = 0.6] (Yr) at (3.75, 0) {$\}$};
          \node (Y) at (3.125, 0) [label = -90:$Y$] {};

          \node[rotate = 90, yscale = 0.6] (Bl) at (0, 2) {$\{$};
          \node[rotate = 90, yscale = 0.6] (Br) at (0, 3.25) {$\}$};
          \node (B) at (0, 2.625) [label = 180:$B$] {};
          \draw (0.5, 2) -- (0.5, 3.25);
        \end{scope}
        \begin{scope}[shift = {(5.5, 0)}]
          \draw (0, 0) -- (4, 0);
          \draw (0, 4) -- (4, 4);
          \draw[dashed] (0, 0) -- (0, 4);
          \draw[dashed] (4, 0) -- (4, 4);

          \pgfmathsetmacro{\r}{{4 / (2 * pi)}}
          \begin{scope}[shift = {(2, 4.5 + \r)}]
            \draw [dashdotted, domain = 90:270] plot ({\r * cos(\x)}, {\r * sin(\x)});
            \draw [dotted, domain = -90:90] plot ({\r * cos(\x)}, {\r * sin(\x)});

            \draw[dashdotted] (-2, -\r) -- (0, -\r);
            \draw[dotted] (0, -\r) -- (2, -\r);
          \end{scope}

          \node[circle, fill, inner sep = 0pt, minimum size = 2pt] (m) at (0.5, 0) [label = -90:$m$] {};
          \node[yscale = 0.6] (Yl) at (2.5, 0) {$\{$};
          \node[yscale = 0.6] (Yr) at (3.75, 0) {$\}$};
          \node (Y) at (3.125, 0) [label = -90:$Y$] {};

          \node[circle, fill, inner sep = 0pt, minimum size = 2pt] (mr) at (3.5, 0) [label = $\varrho_{-1}(m)$] {};

          \node[rotate = 90, yscale = 0.6] (Bl) at (0, 3) {$\{$};
          \node[rotate = 90, yscale = 0.6] (Br) at (0, 0.25) {$\}$};
          \node (B') at (0, 3.5) [label = 180:$B'$] {};
          \draw (0.5, 3) -- (0.5, 4);
          \draw (0.5, 0) -- (0.5, 0.25);
        \end{scope}
      \end{tikzpicture}
      \label{figure:tori:arc-Y}
    }
  \end{minipage}%

  \begin{minipage}{.5\linewidth}
    \myfloatalign
    \subfloat[]{%
      \begin{tikzpicture}[scale = 0.55]
        \begin{scope}[shift = {(0, 0)}]
          \draw (0, 0) -- (4, 0);
          \draw (0, 4) -- (4, 4);
          \draw[dashed] (0, 0) -- (0, 4);
          \draw[dashed] (4, 0) -- (4, 4);

          \node[yscale = 0.6] (Xl) at (0.5, 0) {$\{$};
          \node[yscale = 0.6] (Xr) at (1.75, 0) {$\}$};
          \node (X) at (1.125, 0) [label = -90:$X$] {};
          \node[circle, fill, inner sep = 0pt, minimum size = 2pt] (m') at (2.5, 0) [label = -90:$m'$] {};

          \draw (0.5, 2) -- (1.75, 0.75);
        \end{scope}
        \begin{scope}[shift = {(5, 0)}]
          \draw (0, 0) -- (4, 0);
          \draw (0, 4) -- (4, 4);
          \draw[dashed] (0, 0) -- (0, 4);
          \draw[dashed] (4, 0) -- (4, 4);

          \pgfmathsetmacro{\r}{{4 / (2 * pi)}}
          \begin{scope}[shift = {(2, 4.5 + \r)}]
            \draw [dashdotted, domain = 90:270] plot ({\r * cos(\x)}, {\r * sin(\x)});
            \draw [dotted, domain = -90:90] plot ({\r * cos(\x)}, {\r * sin(\x)});

            \draw[dashdotted] (-2, -\r) -- (0, -\r);
            \draw[dotted] (0, -\r) -- (2, -\r);
          \end{scope}

          \node[yscale = 0.6] (Xl) at (0.5, 0) {$\{$};
          \node[yscale = 0.6] (Xr) at (1.75, 0) {$\}$};
          \node (X) at (1.125, 0) [label = -90:$X$] {};
          \node[circle, fill, inner sep = 0pt, minimum size = 2pt] (m') at (2.5, 0) [label = -90:$m'$] {};

          \node[yscale = 0.6] (Xl') at (3.5, 0) {$\}$};
          \node[yscale = 0.6] (Xr') at (2.25, 0) {$\{$};
          \node (X') at (2.875, 0) [label = $\varrho_{-1}(X)$] {};

          \draw (0.5, 3) -- (1.5, 4);
          \draw (1.5, 0) -- (1.75, 0.25);
        \end{scope}
      \end{tikzpicture}
      \label{figure:tori:arc-X}
    }
  \end{minipage}%
  \begin{minipage}{.5\linewidth}
    \myfloatalign
    \subfloat[]{%
      \begin{tikzpicture}[scale = 0.55]
        \begin{scope}[shift = {(0, 0)}]
          \draw (0, 0) -- (4, 0);
          \draw (0, 4) -- (4, 4);
          \draw[dashed] (0, 0) -- (0, 4);
          \draw[dashed] (4, 0) -- (4, 4);

          \node[yscale = 0.6] (Xl) at (0.5, 0) {$\{$};
          \node[yscale = 0.6] (Xr) at (1.75, 0) {$\}$};
          \node (X) at (1.125, 0) [label = -90:$X$] {};

          \node[yscale = 0.8] (Yl) at (3.25, 0) {$\{$};
          \node[yscale = 0.8] (Yr) at (3.75, 0) {$\}$};
          \node (Y) at (3.5, 0) [label = -90:$Y$] {};

          \draw[fill] (0.5, 2.75) -- (0.5, 3.25) -- (1.75, 2) -- (1.75, 1.5) -- cycle;
        \end{scope}
        \begin{scope}[shift = {(5, 0)}]
          \draw (0, 0) -- (4, 0);
          \draw (0, 4) -- (4, 4);
          \draw[dashed] (0, 0) -- (0, 4);
          \draw[dashed] (4, 0) -- (4, 4);

          \pgfmathsetmacro{\r}{{4 / (2 * pi)}}
          \begin{scope}[shift = {(2, 4.5 + \r)}]
            \draw [dashdotted, domain = 90:270] plot ({\r * cos(\x)}, {\r * sin(\x)});
            \draw [dotted, domain = -90:90] plot ({\r * cos(\x)}, {\r * sin(\x)});

            \draw[dashdotted] (-2, -\r) -- (0, -\r);
            \draw[dotted] (0, -\r) -- (2, -\r);
          \end{scope}

          \node[yscale = 0.6] (Xl) at (0.5, 0) {$\{$};
          \node[yscale = 0.6] (Xr) at (1.75, 0) {$\}$};
          \node (X) at (1.125, 0) [label = -90:$X$] {};

          \node[yscale = 0.8] (Yl) at (3.25, 0) {$\{$};
          \node[yscale = 0.8] (Yr) at (3.75, 0) {$\}$};
          \node (Y) at (3.5, 0) [label = -90:$Y$] {};

          \node[yscale = 0.6] (Xl') at (3.5, 0) {$\}$};
          \node[yscale = 0.6] (Xr') at (2.25, 0) {$\{$};
          \node (X') at (2.875, 0) [label = $\varrho_{-1}(X)$] {};

          \draw[fill] (0.5, 0.25) -- (1.75, 1.5) -- (1.75, 1) -- (0.75, 0) -- (0.5, 0) -- cycle;
          \draw[fill] (0.5, 4) -- (0.5, 3.75) -- (0.75, 4);
        \end{scope}
      \end{tikzpicture}
      \label{figure:tori:ribbons}
    }
  \end{minipage}
}

\def\figureLatticeFolnerOriginal{%
      \subfloat[$\rho = 1$]{
        \begin{tikzpicture}[scale = 0.4]
          \GridWithCoordinates((-5,-4),(5,4));
          \TaxicabBallOutline((0,0),1,solid); 
          \TaxicabBallOutline((1,0),1,dashed);
          \TaxicabSphereDotsLeft((0,0),1);
        \end{tikzpicture}
      }
      \hfill
      \subfloat[$\rho = 2$]{
        \begin{tikzpicture}[scale = 0.4]
          \GridWithCoordinates((-5,-4),(5,4));
          \TaxicabBallOutline((0,0),2,solid); 
          \TaxicabBallOutline((1,0),2,dashed);
          \TaxicabSphereDotsLeft((0,0),2);
        \end{tikzpicture}
      }
      \hfill
      \subfloat[$\rho = 3$]{
        \begin{tikzpicture}[scale = 0.4]
          \GridWithCoordinates((-5,-4),(5,4));
          \TaxicabBallOutline((0,0),3,solid); 
          \TaxicabBallOutline((1,0),3,dashed);
          \TaxicabSphereDotsLeft((0,0),3);
        \end{tikzpicture}
      }
}

\def\figureTreeTranslationsRotations{%
  \begin{wide}
    \subfloat{
      \begin{tikzpicture}[scale = 0.4]
        \pgfmathsetmacro\i{3}
        \pgfmathsetmacro\j{1.5}
        \pgfmathsetmacro\k{0.75}
        \pgfmathsetmacro\l{0.375}

        \draw (4, -3) edge[|->, bend right = 60] node[below] {rotation about the midpoint of the edge labelled $a$} (20, -3);
        \draw[dash pattern = on 1pt off 1pt] (0, 0) -- node[auto] {$a$} (\i, 0); 
        \draw[dash pattern = on 1pt off 1pt on \the\pgflinewidth off 1pt] (0, 0) -- node[auto] {$a^{-1}$} (-\i, 0); 
        \draw[densely dotted] (0, 0) -- node[auto] {$b$} (0, \i); 
        \draw (0, 0) -- node[auto] {$b^{-1}$} (0, -\i); 
        \freeGroupCayleyGraphWithoutInnerCross
        \draw (4, 3) edge[|->, bend left] node[above] {$\varphi(180)$} (8, 3); 
        \begin{scope}[shift = {(12, 0)}]
          \draw[dash pattern = on 1pt off 1pt on \the\pgflinewidth off 1pt] (0, 0) -- node[auto] {$a$} (\i, 0); 
          \draw[dash pattern = on 1pt off 1pt] (0, 0) -- node[auto] {$a^{-1}$} (-\i, 0); 
          \draw (0, 0) -- node[auto] {$b$} (0, \i); 
          \draw[densely dotted] (0, 0) -- node[auto] {$b^{-1}$} (0, -\i); 
          \begin{scope}[rotate = 180]
            \freeGroupCayleyGraphWithoutInnerCross
          \end{scope}
        \end{scope}
        \draw (16, 3) edge[|->, bend left] node[above] {$a \cdot \blank$} (20, 3);
        \begin{scope}[shift = {(24, 0)}]
            \draw[dash pattern = on 1pt off 1pt] (0, 0) -- node[auto] {$a$} (\i, 0); 
              \begin{scope}[shift = {(\i, 0)}]
                \draw[dash pattern = on 1pt off 1pt on \the\pgflinewidth off 1pt] (0, 0) -- (\j, 0);
                \draw (0, 0) -- (0, \j);
                \draw[densely dotted] (0, 0) -- (0, -\j);
              \end{scope}
              \begin{scope}[every path/.append style = {dash pattern = on 1pt off 1pt on \the\pgflinewidth off 1pt}]
                \TShapedCross(\k,(\i+\j,0),0);
                  \TShapedCross(\l,(\i+\j+\k,0),0);
                  \TShapedCross(\l,(\i+\j,\k),90);
                  \TShapedCross(\l,(\i+\j,-\k),-90);
              \end{scope}
                \TShapedCross(\k,(\i,\j),90);
                  \TShapedCross(\l,(\i,\j+\k),90);
                  \TShapedCross(\l,(\i+\k,\j),0);
                  \TShapedCross(\l,(\i-\k,\j),180);
              \begin{scope}[every path/.append style = {densely dotted}]
                \TShapedCross(\k,(\i,-\j),-90);
                  \TShapedCross(\l,(\i,-\j-\k),-90);
                  \TShapedCross(\l,(\i+\k,-\j),0);
                  \TShapedCross(\l,(\i-\k,-\j),180);
              \end{scope}

          \begin{scope}[every path/.append style = {dash pattern = on 1pt off 1pt}]
            \draw (0, 0) -- node[auto] {$a^{-1}$} (-\i, 0); 
              \TShapedCross(\j,(-\i,0),180);
                \TShapedCross(\k,(-\i-\j,0),180);
                  \TShapedCross(\l,(-\i-\j-\k,0),180);
                  \TShapedCross(\l,(-\i-\j,\k),90);
                  \TShapedCross(\l,(-\i-\j,-\k),-90);
                \TShapedCross(\k,(-\i,\j),90);
                  \TShapedCross(\l,(-\i,\j+\k),90);
                  \TShapedCross(\l,(-\i+\k,\j),0);
                  \TShapedCross(\l,(-\i-\k,\j),180);
                \TShapedCross(\k,(-\i,-\j),-90);
                  \TShapedCross(\l,(-\i,-\j-\k),-90);
                  \TShapedCross(\l,(-\i+\k,-\j),0);
                  \TShapedCross(\l,(-\i-\k,-\j),180);

            \draw (0, 0) -- node[auto] {$b$} (0, \i); 
              \TShapedCross(\j,(0,\i),90);
                \TShapedCross(\k,(0,\i+\j),90);
                  \TShapedCross(\l,(0,\i+\j+\k),90);
                  \TShapedCross(\l,(\k,\i+\j),0);
                  \TShapedCross(\l,(-\k,\i+\j),180);
                \TShapedCross(\k,(\j,\i),0);
                  \TShapedCross(\l,(\j+\k,\i),0);
                  \TShapedCross(\l,(\j,\i+\k),90);
                  \TShapedCross(\l,(\j,\i-\k),-90);
                \TShapedCross(\k,(-\j,\i),180);
                  \TShapedCross(\l,(-\j-\k,\i),180);
                  \TShapedCross(\l,(-\j,\i+\k),90);
                  \TShapedCross(\l,(-\j,\i-\k),-90);

            \draw (0, 0) -- node[auto] {$b^{-1}$} (0, -\i); 
              \TShapedCross(\j,(0,-\i),-90);
                \TShapedCross(\k,(0,-\i-\j),-90);
                  \TShapedCross(\l,(0,-\i-\j-\k),-90);
                  \TShapedCross(\l,(\k,-\i-\j),0);
                  \TShapedCross(\l,(-\k,-\i-\j),180);
                \TShapedCross(\k,(\j,-\i),0);
                  \TShapedCross(\l,(\j+\k,-\i),0);
                  \TShapedCross(\l,(\j,-\i+\k),90);
                  \TShapedCross(\l,(\j,-\i-\k),-90);
                \TShapedCross(\k,(-\j,-\i),180);
                  \TShapedCross(\l,(-\j-\k,-\i),180);
                  \TShapedCross(\l,(-\j,-\i+\k),90);
                  \TShapedCross(\l,(-\j,-\i-\k),-90);
          \end{scope}
        \end{scope}
      \end{tikzpicture}
    }
  \end{wide}

  \vspace{2em}

  \begin{wide}
    \subfloat{
      \begin{tikzpicture}[scale = 0.4]
        \pgfmathsetmacro\i{3}
        \pgfmathsetmacro\j{1.5}
        \pgfmathsetmacro\k{0.75}
        \pgfmathsetmacro\l{0.375}

        \draw (4, -3) edge[|->, bend right = 60] node[below] {$a b^{-1} \cdot \varphi(180)(\blank)$} (20, -3);
        \draw[dash pattern = on 1pt off 1pt] (0, 0) -- node[auto] {$a$} (\i, 0); 
        \draw[dash pattern = on 1pt off 1pt on \the\pgflinewidth off 1pt] (0, 0) -- node[auto] {$a^{-1}$} (-\i, 0); 
        \draw[densely dotted] (0, 0) -- node[auto] {$b$} (0, \i); 
        \draw (0, 0) -- node[auto] {$b^{-1}$} (0, -\i); 
        \freeGroupCayleyGraphWithoutInnerCross
        \draw (4, 3) edge[|->, bend left] node[above] {$\varphi(90)$} (8, 3); 
        \begin{scope}[shift = {(12, 0)}]
          \draw (0, 0) -- node[auto] {$a$} (\i, 0); 
          \draw[densely dotted] (0, 0) -- node[auto] {$a^{-1}$} (-\i, 0); 
          \draw[dash pattern = on 1pt off 1pt] (0, 0) -- node[auto] {$b$} (0, \i); 
          \draw[dash pattern = on 1pt off 1pt on \the\pgflinewidth off 1pt] (0, 0) -- node[auto] {$b^{-1}$} (0, -\i); 
          \begin{scope}[rotate = 90]
            \freeGroupCayleyGraphWithoutInnerCross
          \end{scope}
        \end{scope}
        \draw (16, 3) edge[|->, bend left] node[above] {$a \cdot \varphi(90)(a^{-1} \cdot \blank)$} (20, 3);
        \begin{scope}[shift = {(24, 0)}]
          \draw (0, 0) -- node[auto] {$a$} (\i, 0); 
            \TShapedCross(\j,(\i,0),0);
              \TShapedCross(\k,(\i+\j,0),0);
                \TShapedCross(\l,(\i+\j+\k,0),0);
                \TShapedCross(\l,(\i+\j,\k),90);
                \TShapedCross(\l,(\i+\j,-\k),-90);
              \TShapedCross(\k,(\i,\j),90);
                \TShapedCross(\l,(\i,\j+\k),90);
                \TShapedCross(\l,(\i+\k,\j),0);
                \TShapedCross(\l,(\i-\k,\j),180);
              \begin{scope}[shift = {(\i, -\j)}]
                \draw[densely dotted]
                      (0, 0) -- (0, -\k)
                      (0, -\k) -- (0, -\k - \l)
                      (0, -\k) -- (\l, -\k)
                      (0, -\k) -- (-\l, -\k);
                \draw[dash pattern=on 1pt off 1pt on \the\pgflinewidth off 1pt] 
                      (0, 0) -- (\k, 0)
                      (\k, 0) -- (\k + \l, 0)
                      (\k, 0) -- (\k, \l)
                      (\k, 0) -- (\k, -\l);
                \draw[dash pattern = on 1pt off 1pt]
                      (0, 0) -- (-\k, 0)
                      (-\k, 0) -- (-\k - \l, 0)
                      (-\k, 0) -- (-\k, \l)
                      (-\k, 0) -- (-\k, -\l);
              \end{scope}

          \draw (0, 0) -- node[auto] {$a^{-1}$} (-\i, 0); 
            \TShapedCross(\j,(-\i,0),180);
              \TShapedCross(\k,(-\i-\j,0),180);
                \TShapedCross(\l,(-\i-\j-\k,0),180);
                \TShapedCross(\l,(-\i-\j,\k),90);
                \TShapedCross(\l,(-\i-\j,-\k),-90);
              \TShapedCross(\k,(-\i,\j),90);
                \TShapedCross(\l,(-\i,\j+\k),90);
                \TShapedCross(\l,(-\i+\k,\j),0);
                \TShapedCross(\l,(-\i-\k,\j),180);
              \TShapedCross(\k,(-\i,-\j),-90);
                \TShapedCross(\l,(-\i,-\j-\k),-90);
                \TShapedCross(\l,(-\i+\k,-\j),0);
                \TShapedCross(\l,(-\i-\k,-\j),180);

          \draw (0, 0) -- node[auto] {$b$} (0, \i); 
            \TShapedCross(\j,(0,\i),90);
              \TShapedCross(\k,(0,\i+\j),90);
                \TShapedCross(\l,(0,\i+\j+\k),90);
                \TShapedCross(\l,(\k,\i+\j),0);
                \TShapedCross(\l,(-\k,\i+\j),180);
              \TShapedCross(\k,(\j,\i),0);
                \TShapedCross(\l,(\j+\k,\i),0);
                \TShapedCross(\l,(\j,\i+\k),90);
                \TShapedCross(\l,(\j,\i-\k),-90);
              \TShapedCross(\k,(-\j,\i),180);
                \TShapedCross(\l,(-\j-\k,\i),180);
                \TShapedCross(\l,(-\j,\i+\k),90);
                \TShapedCross(\l,(-\j,\i-\k),-90);

          \draw (0, 0) -- node[auto] {$b^{-1}$} (0, -\i); 
            \TShapedCross(\j,(0,-\i),-90);
              \TShapedCross(\k,(0,-\i-\j),-90);
                \TShapedCross(\l,(0,-\i-\j-\k),-90);
                \TShapedCross(\l,(\k,-\i-\j),0);
                \TShapedCross(\l,(-\k,-\i-\j),180);
              \TShapedCross(\k,(\j,-\i),0);
                \TShapedCross(\l,(\j+\k,-\i),0);
                \TShapedCross(\l,(\j,-\i+\k),90);
                \TShapedCross(\l,(\j,-\i-\k),-90);
              \TShapedCross(\k,(-\j,-\i),180);
                \TShapedCross(\l,(-\j-\k,-\i),180);
                \TShapedCross(\l,(-\j,-\i+\k),90);
                \TShapedCross(\l,(-\j,-\i-\k),-90);
        \end{scope}
      \end{tikzpicture}
    }
  \end{wide}
}

\def\figureTreeFolnerOriginal{%
    \subfloat[$\rho = 1$]{
      \begin{tikzpicture}[scale = 0.4]
        \pgfmathsetmacro\i{3}
        \pgfmathsetmacro\j{1.5}
        \pgfmathsetmacro\k{0.75}
        \pgfmathsetmacro\l{0.375}
        \pgfmathsetmacro\t{\l} 

        \freeGroupBallRadiusOne;
      \end{tikzpicture}
    }
    \hfill
    \subfloat[$\rho = 2$]{
      \begin{tikzpicture}[scale = 0.4]
        \pgfmathsetmacro\i{3}
        \pgfmathsetmacro\j{1.5}
        \pgfmathsetmacro\k{0.75}
        \pgfmathsetmacro\l{0.375}
        \pgfmathsetmacro\t{0.1875} 

        \freeGroupBallRadiusTwo;
      \end{tikzpicture}
    }
    \hfill
    \subfloat[$\rho = 3$]{
      \begin{tikzpicture}[scale = 0.4]
        \pgfmathsetmacro\i{3}
        \pgfmathsetmacro\j{1.5}
        \pgfmathsetmacro\k{0.75}
        \pgfmathsetmacro\l{0.375}
        \pgfmathsetmacro\t{0.09375} 

        \freeGroupBallRadiusThree; 
      \end{tikzpicture}
    }
}

\def\figureTreeParadoxicalDecomposition{%
    \subfloat[]{
      \begin{tikzpicture}[scale = 0.65]
        \pgfmathsetmacro\i{3}
        \pgfmathsetmacro\j{1}
        \pgfmathsetmacro\k{0.333}
        \pgfmathsetmacro\l{0.166}
        \pgfmathsetmacro\t{\l} 
          \freeGroupCayleyGraph(solid);
          \AMinusContour((-\i,0),1);
            \node[anchor = north east] at (-\i - \j, -\j) {$A_{e_{F_2}}$}; 
          \APlusContour((0,0),3.25);
            \node[anchor = base] at (\i, -\i) {$A_a$}; 
      \end{tikzpicture}
      \label{figure:tree:paradoxical-decomposition:one}
    }
    \hfill
    \subfloat[]{
      \begin{tikzpicture}[scale = 0.65]
        \pgfmathsetmacro\i{3}
        \pgfmathsetmacro\j{1}
        \pgfmathsetmacro\k{0.333}
        \pgfmathsetmacro\l{0.166}
        \pgfmathsetmacro\t{\l} 
        \begin{scope}[every path/.style={draw=white}]
          \APlusContour((0,0),3.25);
        \end{scope}
          \freeGroupCayleyGraph(solid);
          \BMinusContour((0,-\i),1);
            \node[anchor = north east] at (-\j, -\i - \j) {$B_{e_{F_2}}$}; 
          \BPlusContour((0,0),3,-1.4,\t,\t/6);
            \node[anchor = base] at (\i, \i) {$B_b$}; 
      \end{tikzpicture}
      \label{figure:tree:paradoxical-decomposition:two}
    }

    \vspace{2em}

    \subfloat[]{
      \begin{tikzpicture}[scale = 0.65]
        \pgfmathsetmacro\i{3}
        \pgfmathsetmacro\j{1}
        \pgfmathsetmacro\k{0.333}
        \pgfmathsetmacro\l{0.166}
        \pgfmathsetmacro\t{\l} 
          \freeGroupCayleyGraph(solid);
          \AMinusContour((-\i,0),1);
            \node[anchor = north east] at (-\i - \j, -\j) {$A_{e_{F_2}} \cdot e_{F_2}$}; 
          \APlusContour((\i,0),1);
            \node[anchor = north west] at (\i + \j, -\j) {$A_a \cdot a$}; 
          \BMinusContour((0,-\i),1);
            \node[anchor = north east] at (-\j, -\i - \j) {$B_{e_{F_2}} \cdot e_{F_2}$}; 
          \BPlusContour((0,\i),1,-\i-\i-\j-\k-\l,4*\t,\t/3);
            \node[anchor = south west] at (\j, \i + \j) {$B_b \cdot b$}; 
      \end{tikzpicture}
      \label{figure:tree:paradoxical-decomposition:three}
    }
}

\def\figureFiniteNImpliesTwoToOneSurjectiveMap{%
    \begin{tikzpicture} 
      \pgfmathsetmacro\w{5} 
      \pgfmathsetmacro\cox{8} 

      \draw (0, 0) rectangle (\w, 8);
      \draw (\cox, 0) rectangle (\cox + \w, 8);

      \node (M) at (\w/2, 8.5) {$M$};
      \node (M) at (\cox + \w/2, 8.5) {$M$};

      \node[fill, circle, inner sep = 1pt, label = left: $m$] (m) at (\w/2 + 0.625, 5.75) {}; 
      \node[draw, circle, inner sep = 15pt, label = above: $F$] (F) at (\w/2, 5.75) {};
      \node[draw, circle, inner sep = 15pt, label = above: $F$] (Fx) at (\cox + \w/2, 5.75) {};

      \node[draw, circle, inner sep = 27pt, label = {above: $\setOf{m \in M \suchThat (m \isSemiActedUponBy E) \cap F \neq \emptyset}$}] (Nl) at (\w/2, 2) {}; 
      \node[inner sep = 51pt, label = {above: $=$}] (Nl) at (\w/2, 2) {}; 
      \node[inner sep = 60pt, label = {above: $\mathcal{N}_l(F')$}] (Nl) at (\w/2, 2) {}; 

      \node[fill, circle, inner sep = 0.25pt] (mxx) at (\cox + \w/2 + 0.625, 5.75) {};
      \node[draw, circle, inner sep = 12pt, dashed] (m lib N) at (\cox + \w/2 + 0.625, 5.75) {}; 
      \node[fill, circle, inner sep = 1pt] (m lib n) at (\cox + \w/2 + 1, 5.5) {}; 
      \node[draw, circle, inner sep = 27pt, label = {above: $\mathcal{N}_r(F) = F \isSemiActedUponBy E$}] (Nr) at (\cox + \w/2, 5.75) {};
      \node[fill, circle, inner sep = 1pt, label = right: $m'$] (m') at (\cox + \w/2 - 0.5, 2) {};
      \node[draw, circle, inner sep = 15pt, label = above: $F'$] (F') at (\cox + \w/2, 2) {};
      \node[draw, circle, inner sep = 15pt, label = above: $F'$] (F'x) at (\w/2, 2) {};
      \node[fill, circle, inner sep = 0.25pt] (m'xx) at (\w/2 - 0.5, 2) {};
      \node[draw, circle, inner sep = 12pt, densely dotted] (m lib N) at (\w/2 - 0.5, 2) {}; 
      \node[fill, circle, inner sep = 1pt] (lib n inv m') at (\w/2 - 0.875, 1.75) {};
      \node[draw, circle, inner sep = 12pt, dashed] (lib n inv m' lib N) at (\w/2 - 0.875, 1.75) {};

      \draw[->] (m) edge node[above] {$\blank \isSemiActedUponBy e$} (m lib n)
                (lib n inv m') edge node[above] {$\blank \isSemiActedUponBy e'$} (m');
    \end{tikzpicture}
}

\def\figureTwoToOneSurjectiveMapImpliesTarskiFolnerNoParadoxicalDecomposition{%
    \begin{tikzpicture}
      \pgfmathsetmacro\w{4} 
      \pgfmathsetmacro\cox{7} 

      \foreach \x in {0} {
        \draw (\x, 0) rectangle (\x + \w, 8);
        \foreach \y in {1, 2, ..., 7} {
          \draw (\x, \y) -- (\x + \w, \y);
        }
      }
      \foreach \x in {\cox} {
        \draw (\x, 0) rectangle (\x + \w, 8);
        \foreach \y in {2, 4, 6} {
          \draw (\x, \y) -- (\x + \w, \y);
        }
      }
      \foreach \x in {\cox + 0.125} {
        \draw[gray, dashed] (\x, -0.125) rectangle (\x + \w, 8 - 0.125);
        \foreach \y in {1.875, 3.875, 5.875} {
          \draw[gray, dashed] (\x, \y) -- (\x + \w, \y);
        }
      }
      \draw[very thick] (0, 4) -- (\w, 4);

      \node (M) at (\w / 2, 8.5) {$M$};
      \node (M) at (\cox + \w / 2, 8.5) {$M$};
      \node (psi M) at (-1, 6) {$\psi(M)$};
      \node (psi' M) at (-1, 2) {$\psi'(M)$};

      \node[fill, circle, inner sep = 1pt] (psi m) at (\w - 0.5, 6.5) {};
      \node (psi A n) at (1.45, 6.5) {$\psi(A_e) = A_e \isSemiActedUponBy e$};
      \node[fill, circle, inner sep = 1pt] (psi' m) at (\w - 0.5, 2.5) {};
      \node (psi' B n') at (1.65, 2.5) {$\psi'(B_{e'}) = B_{e'} \isSemiActedUponBy e'$};
      \node[fill, circle, inner sep = 1pt, label = right: $m$] (m) at (\cox + 0.5, 5) {}; 
      \node (A n) at (\cox + \w / 2, 5) {$A_e$};
      \node[gray] (B n') at (\cox + \w /2 + 1.25, 4.5) {$B_{e'}$};

      \draw[->] (psi m) edge[out = -60, in = -180] node[above] {$\phi$} (m)
                (m) edge[out = 130, in = 10] node[above] {$\psi$} (psi m)
                (m) edge[out = 70, in = 70] node[above] {$\blank \isSemiActedUponBy e$} (psi m)
                (psi' m) edge[out = 60, in = -180] node[below] {$\phi$} (m)
                (m) edge[out = -130, in = -10] node[below] {$\psi'$} (psi' m)
                (m) edge[out = -70, in = -70] node[below] {$\blank \isSemiActedUponBy e'$} (psi' m);
    \end{tikzpicture}
}

\def\figureTreeTarskiFolnerTwoToOneMapPhi{%
    \begin{tikzpicture}[scale = 0.7]
      \pgfmathsetmacro\i{3}
      \pgfmathsetmacro\j{1}
      \pgfmathsetmacro\k{0.333}
      \pgfmathsetmacro\l{0.166}

      \pgfmathsetmacro\t{\l} 

      \freeGroupCayleyGraph(solid); 

      \AMinusContour((-\i,0),1);
          \node[anchor = north east] at (-\i - \j, -\j) {$A^{-}$}; 
      \APlusContour((\i,0),1);
          \node[anchor = north west] at (\i + \j, -\j) {$A^{+}$}; 
      \BMinusContour((0,-\i),1);
          \node[anchor = north east] at (-\j, -\i - \j) {$B^{-}$}; 
      \BPlusContour((0,\i),1,-\i-\i-\j-\k-\l,4*\t,\t/3);
          \node[anchor = south west] at (\j, \i + \j) {$B^{+}$}; 

      \draw (4, 3) edge[|->, bend left] node[above] {$\phi$} (8, 3);

      \begin{scope}[shift = {(12, 0)}]
        \freeGroupCayleyGraph(solid);

        \AMinusContour((-\i,0),1);
        \APlusContour((0,0),3.333);
            \node[anchor = north west] at (\i + \j + \k + \l, -\j - \k - \l) {$A^{+} \cdot a^{-1}$}; 
        \BMinusContour((0,-\i),1);
        \BPlusContour((0,0),3,-1.4,\t,\t/6);
            \node[anchor = south east] at (-\i - \j, \j + \k + \l) {$B^{+} \cdot b^{-1}$}; 
      \end{scope}
    \end{tikzpicture}
}

\def\figureTreeTarskiFolnerTwoToOneMapPsiAndPsiPrime{%
    \subfloat[]{
      \begin{tikzpicture}[scale = 0.35]
        \pgfmathsetmacro\i{3}
        \pgfmathsetmacro\j{1}
        \pgfmathsetmacro\k{0.333}
        \pgfmathsetmacro\l{0.166}

        \pgfmathsetmacro\t{\l} 

        \freeGroupCayleyGraph(solid); 

        \AMinusContour((-\i,0),1);
        \APlusContour((\i,0),1);

        \draw (8, 3) edge[|->, bend right] node[above] {$\psi$} (4, 3);
        \draw (4, -3) edge[|->, bend right] node[below] {$\phi$} (8, -3);

        \begin{scope}[shift = {(12, 0)}]
          \freeGroupCayleyGraph(solid);

          \AMinusContour((-\i,0),1);
          \APlusContour((0,0),3.25);
        \end{scope}
      \end{tikzpicture}
    }
    \subfloat[]{
      \begin{tikzpicture}[scale = 0.35]
        \pgfmathsetmacro\i{3}
        \pgfmathsetmacro\j{1}
        \pgfmathsetmacro\k{0.333}
        \pgfmathsetmacro\l{0.166}

        \pgfmathsetmacro\t{\l} 

        \freeGroupCayleyGraph(solid); 

        \BMinusContour((0,-\i),1);
        \BPlusContour((0,\i),1,-\i-\i-\j-\k-\l,4*\t,\t/3);

        \draw (8, 3) edge[|->, bend right] node[above] {$\psi'$} (4, 3);
        \draw (4, -3) edge[|->, bend right] node[below] {$\phi$} (8, -3);

        \begin{scope}[shift = {(12, 0)}]
          \freeGroupCayleyGraph(solid);

          \BMinusContour((0,-\i),1);
          \BPlusContour((0,0),3,-1.4,\t,\t/6);
        \end{scope}
      \end{tikzpicture}
    }
}

\def\figureTreeTarskiFolnerFiniteSubsetE{%
    \begin{tikzpicture}[scale = 0.7]
      \pgfmathsetmacro\i{3}
      \pgfmathsetmacro\j{1}
      \pgfmathsetmacro\k{0.333}
      \pgfmathsetmacro\l{0.166}

      \pgfmathsetmacro\t{\l} 

      \freeGroupCayleyGraph(solid);

      \AMinusContour((-\i,0),1); 
          \node[anchor = north east] at (-\i - \j, -\j) {$A^{-} \cdot e_{F_2}$}; 
      \AMinusContour((-\j,\i),0.333); 
          \node[anchor = east] at (-\j - \k - \l - \l, \i) {$A^{-} \cdot b$}; 

      \APlusContour((\i+\j,0),0.333); 
          \node[anchor = north west] at (\i + \j + \k, -\k) {$A^{+} \cdot a$}; 
      \APlusContour((\j,\i),0.333); 
          \node[anchor = west] at (\j + \k + \l + \l, \i) {$A^{+} \cdot b$}; 

      \BMinusContour((0,-\i),1); 
          \node[anchor = north east] at (-\j, -\i - \j) {$B^{-} \cdot e_{F_2}$}; 
      \BMinusContour((\i,-\j),0.333); 
          \node[anchor = north west] at (\i + \k, -\j - \k) {$B^{-} \cdot a$}; 

      \BPlusContour((\i,\j),0.333,-\i-\i-\j-\k-\l,4*\t,\t/3); 
          \node[anchor = south west] at (\i + \k, \j + \k) {$B^{+} \cdot a$}; 
      \BPlusContour((0,\i+\j),0.333,-\i-\j-\k-\l-\i-\j-\k-\l-\i-\j-\k-\l-\i-\j-\k-\l-\i-\j-\k-\l-\i-\j,14*\t,\t/3); 
          \node[anchor = south west] at (\k, \i + \j + \k) {$B^{+} \cdot b$}; 

      \draw (8, 3) edge[|->, bend right] node[above] {$\blank \cdot e_{F_2}$, $\blank \cdot a$, or $\blank \cdot b$} (4, 3);
      \draw (4, -3) edge[|->, bend right] node[below] {$\phi$} (8, -3);

      \begin{scope}[shift = {(12, 0)}]
        \freeGroupCayleyGraph(solid); 

        \AMinusContour((-\i,0),1);
            \node[anchor = north east] at (-\i - \j, -\j) {$A^{-}$}; 
        \APlusContour((\i,0),1);
            \node[anchor = north west] at (\i + \j, -\j) {$A^{+}$}; 
        \BMinusContour((0,-\i),1);
            \node[anchor = north east] at (-\j, -\i - \j) {$B^{-}$}; 
        \BPlusContour((0,\i),1,-\i-\i-\j-\k-\l,4*\t,\t/3);
            \node[anchor = south west] at (\j, \i + \j) {$B^{+}$}; 
      \end{scope}
    \end{tikzpicture} 
}

\def\figureLatticeInteriorClosureBoundary{%
  \begin{wide}
      \subfloat[$A^{-\setOf{(1,0)}}$]{
        \begin{tikzpicture}[scale = 0.4]
          \GridWithCoordinates((-4,-4),(4,4));
          \TaxicabBallOutline((0,0),2,solid); 
          \TaxicabBallOutline((-1,0),2,dashed);
          \TaxicabBallDots((-1,0),2);
        \end{tikzpicture}
      }
      \subfloat[$A^{-\setOf{(0,0), (1,0)}}$]{
        \begin{tikzpicture}[scale = 0.4]
          \GridWithCoordinates((-4,-4),(4,4));
          \TaxicabBallOutline((0,0),2,solid); 
          \TaxicabBallOutline((-1,0),2,dashed);
          \TaxicabBallDots((-1/2,0),1.5);
        \end{tikzpicture}
      }
      \subfloat[$A^{-\setOf{(-1,0), (0,0), (1,0)}}$]{
        \begin{tikzpicture}[scale = 0.4]
          \GridWithCoordinates((-4,-4),(4,4));
          \TaxicabBallOutline((0,0),2,solid); 
          \TaxicabBallOutline((-1,0),2,dashed);
          \TaxicabBallOutline((1,0),2,dashed);
          \TaxicabBallDots((0,0),1);
        \end{tikzpicture}
      }
      \subfloat[$A^{-E}$]{
        \begin{tikzpicture}[scale = 0.4]
          \GridWithCoordinates((-4,-4),(4,4));
          \TaxicabBallOutline((0,0),2,solid); 
          \TaxicabBallOutline((-1,0),2,dashdotted);
          \TaxicabBallOutline((1,0),2,dashdotted);
          \TaxicabBallOutline((0,1),2,densely dotted);
          \TaxicabBallOutline((0,-1),2,densely dotted);
          \TaxicabBallDots((0,0),1);
        \end{tikzpicture}
      }
  \end{wide}
  \begin{wide}
      \subfloat[$A^{+\setOf{(1,0)}}$]{
        \begin{tikzpicture}[scale = 0.4]
          \GridWithCoordinates((-4,-4),(4,4));
          \TaxicabBallOutline((0,0),2,solid); 
          \TaxicabBallOutline((-1,0),2,dashed);
          \TaxicabBallDots((-1,0),2);
        \end{tikzpicture}
      }
      \subfloat[$A^{+\setOf{(0,0), (1,0)}}$]{
        \begin{tikzpicture}[scale = 0.4]
          \GridWithCoordinates((-4,-4),(4,4));
          \TaxicabBallOutline((0,0),2,solid); 
          \TaxicabBallOutline((-1,0),2,dashed);
          \TaxicabBallDots((0,0),2);
          \TaxicabBallDots((-1,0),2);
        \end{tikzpicture}
      }
      \subfloat[$A^{+\setOf{(-1,0), (0,0), (1,0)}}$]{
        \begin{tikzpicture}[scale = 0.4]
          \GridWithCoordinates((-4,-4),(4,4));
          \TaxicabBallOutline((0,0),2,solid); 
          \TaxicabBallOutline((-1,0),2,dashed);
          \TaxicabBallOutline((1,0),2,dashed);
          \TaxicabBallDots((0,0),2); 
          \TaxicabBallDots((-1,0),2);
          \TaxicabBallDots((1,0),2);
        \end{tikzpicture}
      }
      \subfloat[$A^{+E}$]{
        \begin{tikzpicture}[scale = 0.4]
          \GridWithCoordinates((-4,-4),(4,4));
          \TaxicabBallOutline((0,0),2,solid); 
          \TaxicabBallOutline((-1,0),2,dashdotted);
          \TaxicabBallOutline((1,0),2,dashdotted);
          \TaxicabBallOutline((0,1),2,densely dotted);
          \TaxicabBallOutline((0,-1),2,densely dotted);
          \TaxicabBallDots((-1,0),2);
          \TaxicabBallDots((1,0),2);
          \TaxicabBallDots((0,1),2);
          \TaxicabBallDots((0,-1),2);
        \end{tikzpicture}
      }
  \end{wide}
  \begin{wide}
      \subfloat[$\boundaryOf_{\setOf{(1,0)}} A$]{
        \begin{tikzpicture}[scale = 0.4]
          \GridWithCoordinates((-4,-4),(4,4));
          \TaxicabBallOutline((0,0),2,solid); 
        \end{tikzpicture}
      }
      \subfloat[$\boundaryOf_{\setOf{(0,0), (1,0)}} A$]{
        \begin{tikzpicture}[scale = 0.4]
          \GridWithCoordinates((-4,-4),(4,4));
          \TaxicabBallOutline((0,0),2,solid); 
          \TaxicabSphereDots((-1/2,0),2.5);
        \end{tikzpicture}
      }
      \subfloat[$\boundaryOf_{\setOf{(-1,0), (0,0), (1,0)}} A$]{
        \begin{tikzpicture}[scale = 0.4]
          \GridWithCoordinates((-4,-4),(4,4));
          \TaxicabBallOutline((0,0),2,solid); 
          \TaxicabSphereDots((0,0),2); 
          \TaxicabSphereDotsLeft((-1,0),2);
          \TaxicabSphereDotsRight((1,0),2);
        \end{tikzpicture}
      }
      \subfloat[$\boundaryOf_{E} A$]{
        \begin{tikzpicture}[scale = 0.4]
          \GridWithCoordinates((-4,-4),(4,4));
          \TaxicabBallOutline((0,0),2,solid); 
          \TaxicabSphereDots((0,0),2);
          \TaxicabSphereDots((0,0),3);
        \end{tikzpicture}
      }
  \end{wide}
}

\def\figureTreeInteriorClosureBoundary{%
  \begin{wide}
    \subfloat[$A^{-\setOf{a}}$]{
      \begin{tikzpicture}[scale = 0.4]
        \pgfmathsetmacro\i{3}
        \pgfmathsetmacro\j{1.5}
        \pgfmathsetmacro\k{0.75}
        \pgfmathsetmacro\l{0.375}

        \freeGroupCayleyGraph(densely dotted);

        \fill (0 - \i, 0) circle (0.1cm);
        \fill (0, 0) circle (0.1cm);
        \fill (0 - \j, 0 + \i) circle (0.1cm);
        \fill (0 - \i - \j, 0) circle (0.1cm);
        \fill (0 - \j, 0 - \i) circle (0.1cm);
        \fill (0 + \i, 0) circle (0.1cm);
        \fill (0 + \i - \k, 0 + \j) circle (0.1cm);
        \fill (0 + \i - \k, 0 - \j) circle (0.1cm);
        \fill (0, 0 + \i) circle (0.1cm);
        \fill (0 - \k, 0 + \i + \j) circle (0.1cm);
        \fill (0 - \j - \k, 0 + \i) circle (0.1cm);
        \fill (0 - \i - \k, 0 + \j) circle (0.1cm);
        \fill (0 - \i - \j - \k, 0) circle (0.1cm);
        \fill (0 - \i - \k, 0 - \j) circle (0.1cm);
        \fill (0, 0 - \i) circle (0.1cm);
        \fill (0 - \j - \k, 0 - \i) circle (0.1cm);
        \fill (0 - \k, 0 - \i - \j) circle (0.1cm);
      \end{tikzpicture}
    }
    \hfill
    \subfloat[$A^{-\setOf{a, b}}$]{
      \begin{tikzpicture}[scale = 0.4]
        \pgfmathsetmacro\i{3}
        \pgfmathsetmacro\j{1.5}
        \pgfmathsetmacro\k{0.75}
        \pgfmathsetmacro\l{0.375}

        \freeGroupCayleyGraph(densely dotted);

        \AaInverseDots
        \AbInverseDots

        \fill (0 - \i, 0) circle (0.1cm);
        \fill (0, 0) circle (0.1cm);
        \fill (0 + \i, 0) circle (0.1cm);
        \fill (0, 0 + \i) circle (0.1cm);
        \fill (0, 0 - \i) circle (0.1cm);
      \end{tikzpicture}
    }
    \hfill
    \subfloat[$A^{-\setOf{a, b, a^{-1}, b^{-1}}}$]{
      \begin{tikzpicture}[scale = 0.4]
        \pgfmathsetmacro\i{3}
        \pgfmathsetmacro\j{1.5}
        \pgfmathsetmacro\k{0.75}
        \pgfmathsetmacro\l{0.375}

        \freeGroupCayleyGraph(densely dotted);

        \fill (0 - \i, 0) circle (0.1cm);
        \fill (0, 0) circle (0.1cm);
        \fill (0 + \i, 0) circle (0.1cm);
        \fill (0, 0 + \i) circle (0.1cm);
        \fill (0, 0 - \i) circle (0.1cm);
      \end{tikzpicture}
    }
  \end{wide}
  \begin{wide}
    \subfloat[$A^{+\setOf{a}}$]{
      \begin{tikzpicture}[scale = 0.4]
        \pgfmathsetmacro\i{3}
        \pgfmathsetmacro\j{1.5}
        \pgfmathsetmacro\k{0.75}
        \pgfmathsetmacro\l{0.375}

        \freeGroupCayleyGraph(densely dotted);

        \fill (0 - \i, 0) circle (0.1cm);
        \fill (0, 0) circle (0.1cm);
        \fill (0 - \j, 0 + \i) circle (0.1cm);
        \fill (0 - \i - \j, 0) circle (0.1cm);
        \fill (0 - \j, 0 - \i) circle (0.1cm);
        \fill (0 + \i, 0) circle (0.1cm);
        \fill (0 + \i - \k, 0 + \j) circle (0.1cm);
        \fill (0 + \i - \k, 0 - \j) circle (0.1cm);
        \fill (0, 0 + \i) circle (0.1cm);
        \fill (0 - \k, 0 + \i + \j) circle (0.1cm);
        \fill (0 - \j - \k, 0 + \i) circle (0.1cm);
        \fill (0 - \i - \k, 0 + \j) circle (0.1cm);
        \fill (0 - \i - \j - \k, 0) circle (0.1cm);
        \fill (0 - \i - \k, 0 - \j) circle (0.1cm);
        \fill (0, 0 - \i) circle (0.1cm);
        \fill (0 - \j - \k, 0 - \i) circle (0.1cm);
        \fill (0 - \k, 0 - \i - \j) circle (0.1cm);
      \end{tikzpicture}
    }
    \hfill
    \subfloat[$A^{+\setOf{a, b}}$]{
      \begin{tikzpicture}[scale = 0.4]
        \pgfmathsetmacro\i{3}
        \pgfmathsetmacro\j{1.5}
        \pgfmathsetmacro\k{0.75}
        \pgfmathsetmacro\l{0.375}

        \freeGroupCayleyGraph(densely dotted);

        \AaInverseDots
        \AbInverseDots

        \fill (0 - \i, 0) circle (0.1cm);
        \fill (0, 0) circle (0.1cm);
        \fill (0 - \j, 0 + \i) circle (0.1cm);
        \fill (0 - \i - \j, 0) circle (0.1cm);
        \fill (0 - \j, 0 - \i) circle (0.1cm);
        \fill (0 + \i, 0) circle (0.1cm);
        \fill (0 + \i - \k, 0 + \j) circle (0.1cm);
        \fill (0 + \i - \k, 0 - \j) circle (0.1cm);
        \fill (0, 0 + \i) circle (0.1cm);
        \fill (0 - \k, 0 + \i + \j) circle (0.1cm);
        \fill (0 - \j - \k, 0 + \i) circle (0.1cm);
        \fill (0 - \i - \k, 0 + \j) circle (0.1cm);
        \fill (0 - \i - \j - \k, 0) circle (0.1cm);
        \fill (0 - \i - \k, 0 - \j) circle (0.1cm);
        \fill (0, 0 - \i) circle (0.1cm);
        \fill (0 - \j - \k, 0 - \i) circle (0.1cm);
        \fill (0 - \k, 0 - \i - \j) circle (0.1cm);
        \fill (0 + \i, 0 - \j) circle (0.1cm);
        \fill (0 - \i, 0 - \j) circle (0.1cm);
        \fill (0, 0 - \i - \j) circle (0.1cm);
        \fill (0 + \i + \j, 0 - \k) circle (0.1cm);
        \fill (0 + \i, 0 - \j - \k) circle (0.1cm);
        \fill (0 + \j, 0 + \i - \k) circle (0.1cm);
        \fill (0 - \j, 0 + \i - \k) circle (0.1cm);
        \fill (0 - \i - \j, 0 - \k) circle (0.1cm);
        \fill (0 - \i, 0 - \j - \k) circle (0.1cm);
        \fill (0 + \j, 0 - \i - \k) circle (0.1cm);
        \fill (0 - \j, 0 - \i - \k) circle (0.1cm);
        \fill (0, 0 - \i - \j - \k) circle (0.1cm);
      \end{tikzpicture}
    }
    \hfill
    \subfloat[$A^{+\setOf{a, b, a^{-1}, b^{-1}}}$]{
      \begin{tikzpicture}[scale = 0.4]
        \pgfmathsetmacro\i{3}
        \pgfmathsetmacro\j{1.5}
        \pgfmathsetmacro\k{0.75}
        \pgfmathsetmacro\l{0.375}

        \freeGroupCayleyGraph(densely dotted);

        \fill (0 - \i, 0) circle (0.1cm);
        \fill (0, 0) circle (0.1cm);
        \fill (0 - \j, 0 + \i) circle (0.1cm);
        \fill (0 - \i - \j, 0) circle (0.1cm);
        \fill (0 - \j, 0 - \i) circle (0.1cm);
        \fill (0 + \i, 0) circle (0.1cm);
        \fill (0 + \i - \k, 0 + \j) circle (0.1cm);
        \fill (0 + \i - \k, 0 - \j) circle (0.1cm);
        \fill (0, 0 + \i) circle (0.1cm);
        \fill (0 - \k, 0 + \i + \j) circle (0.1cm);
        \fill (0 - \j - \k, 0 + \i) circle (0.1cm);
        \fill (0 - \i - \k, 0 + \j) circle (0.1cm);
        \fill (0 - \i - \j - \k, 0) circle (0.1cm);
        \fill (0 - \i - \k, 0 - \j) circle (0.1cm);
        \fill (0, 0 - \i) circle (0.1cm);
        \fill (0 - \j - \k, 0 - \i) circle (0.1cm);
        \fill (0 - \k, 0 - \i - \j) circle (0.1cm);
        \fill (0 + \i + \j, 0) circle (0.1cm);
        \fill (0 + \j, 0 + \i) circle (0.1cm);
        \fill (0 + \j, 0 - \i) circle (0.1cm);
        \fill (0 + \i + \j + \k, 0) circle (0.1cm);
        \fill (0 + \i + \k, 0 + \j) circle (0.1cm);
        \fill (0 + \i + \k, 0 - \j) circle (0.1cm);
        \fill (0 + \j + \k, 0 + \i) circle (0.1cm);
        \fill (0 + \k, 0 + \i + \j) circle (0.1cm);
        \fill (0 - \i + \k, 0 + \j) circle (0.1cm);
        \fill (0 - \i + \k, 0 - \j) circle (0.1cm);
        \fill (0 + \j + \k, 0 - \i) circle (0.1cm);
        \fill (0 + \k, 0 - \i - \j) circle (0.1cm);
        \fill (0 + \i, 0 - \j) circle (0.1cm);
        \fill (0 - \i, 0 - \j) circle (0.1cm);
        \fill (0, 0 - \i - \j) circle (0.1cm);
        \fill (0 + \i + \j, 0 - \k) circle (0.1cm);
        \fill (0 + \i, 0 - \j - \k) circle (0.1cm);
        \fill (0 + \j, 0 + \i - \k) circle (0.1cm);
        \fill (0 - \j, 0 + \i - \k) circle (0.1cm);
        \fill (0 - \i - \j, 0 - \k) circle (0.1cm);
        \fill (0 - \i, 0 - \j - \k) circle (0.1cm);
        \fill (0 + \j, 0 - \i - \k) circle (0.1cm);
        \fill (0 - \j, 0 - \i - \k) circle (0.1cm);
        \fill (0, 0 - \i - \j - \k) circle (0.1cm);
        \fill (0 + \i, 0 + \j) circle (0.1cm);
        \fill (0, 0 + \i + \j) circle (0.1cm);
        \fill (0 - \i, 0 + \j) circle (0.1cm);
        \fill (0 + \i + \j, 0 + \k) circle (0.1cm);
        \fill (0 + \i, 0 + \j + \k) circle (0.1cm);
        \fill (0 + \j, 0 + \i + \k) circle (0.1cm);
        \fill (0, 0 + \i + \j + \k) circle (0.1cm);
        \fill (0 - \j, 0 + \i + \k) circle (0.1cm);
        \fill (0 - \i, 0 + \j + \k) circle (0.1cm);
        \fill (0 - \i - \j, 0 + \k) circle (0.1cm);
        \fill (0 + \j, 0 - \i + \k) circle (0.1cm);
        \fill (0 - \j, 0 - \i + \k) circle (0.1cm);
      \end{tikzpicture}
    }
  \end{wide}
  \begin{wide}
    \subfloat[$\boundaryOf_{\setOf{a}} A$]{
      \begin{tikzpicture}[scale = 0.4]
        \pgfmathsetmacro\i{3}
        \pgfmathsetmacro\j{1.5}
        \pgfmathsetmacro\k{0.75}
        \pgfmathsetmacro\l{0.375}

        \freeGroupCayleyGraph(densely dotted);
      \end{tikzpicture}
    }
    \hfill
    \subfloat[$\boundaryOf_{\setOf{a, b}} A$]{
      \begin{tikzpicture}[scale = 0.4]
        \pgfmathsetmacro\i{3}
        \pgfmathsetmacro\j{1.5}
        \pgfmathsetmacro\k{0.75}
        \pgfmathsetmacro\l{0.375}

        \freeGroupCayleyGraph(densely dotted);

        \fill (0 - \j, 0 + \i) circle (0.1cm);
        \fill (0 - \i - \j, 0) circle (0.1cm);
        \fill (0 - \j, 0 - \i) circle (0.1cm);
        \fill (0 + \i - \k, 0 + \j) circle (0.1cm);
        \fill (0 + \i - \k, 0 - \j) circle (0.1cm);
        \fill (0 - \k, 0 + \i + \j) circle (0.1cm);
        \fill (0 - \j - \k, 0 + \i) circle (0.1cm);
        \fill (0 - \i - \k, 0 + \j) circle (0.1cm);
        \fill (0 - \i - \j - \k, 0) circle (0.1cm);
        \fill (0 - \i - \k, 0 - \j) circle (0.1cm);
        \fill (0 - \j - \k, 0 - \i) circle (0.1cm);
        \fill (0 - \k, 0 - \i - \j) circle (0.1cm);
        \fill (0 + \i, 0 - \j) circle (0.1cm);
        \fill (0 - \i, 0 - \j) circle (0.1cm);
        \fill (0, 0 - \i - \j) circle (0.1cm);
        \fill (0 + \i + \j, 0 - \k) circle (0.1cm);
        \fill (0 + \i, 0 - \j - \k) circle (0.1cm);
        \fill (0 + \j, 0 + \i - \k) circle (0.1cm);
        \fill (0 - \j, 0 + \i - \k) circle (0.1cm);
        \fill (0 - \i - \j, 0 - \k) circle (0.1cm);
        \fill (0 - \i, 0 - \j - \k) circle (0.1cm);
        \fill (0 + \j, 0 - \i - \k) circle (0.1cm);
        \fill (0 - \j, 0 - \i - \k) circle (0.1cm);
        \fill (0, 0 - \i - \j - \k) circle (0.1cm);
      \end{tikzpicture}
    }
    \hfill
    \subfloat[$\boundaryOf_{\setOf{a, b, a^{-1}, b^{-1}}} A$]{
      \begin{tikzpicture}[scale = 0.4]
        \pgfmathsetmacro\i{3}
        \pgfmathsetmacro\j{1.5}
        \pgfmathsetmacro\k{0.75}
        \pgfmathsetmacro\l{0.375}

        \freeGroupCayleyGraph(densely dotted);

        \fill (0 - \j, 0 + \i) circle (0.1cm);
        \fill (0 - \i - \j, 0) circle (0.1cm);
        \fill (0 - \j, 0 - \i) circle (0.1cm);
        \fill (0 + \i - \k, 0 + \j) circle (0.1cm);
        \fill (0 + \i - \k, 0 - \j) circle (0.1cm);
        \fill (0 - \k, 0 + \i + \j) circle (0.1cm);
        \fill (0 - \j - \k, 0 + \i) circle (0.1cm);
        \fill (0 - \i - \k, 0 + \j) circle (0.1cm);
        \fill (0 - \i - \j - \k, 0) circle (0.1cm);
        \fill (0 - \i - \k, 0 - \j) circle (0.1cm);
        \fill (0 - \j - \k, 0 - \i) circle (0.1cm);
        \fill (0 - \k, 0 - \i - \j) circle (0.1cm);
        \fill (0 + \i + \j, 0) circle (0.1cm);
        \fill (0 + \j, 0 + \i) circle (0.1cm);
        \fill (0 + \j, 0 - \i) circle (0.1cm);
        \fill (0 + \i + \j + \k, 0) circle (0.1cm);
        \fill (0 + \i + \k, 0 + \j) circle (0.1cm);
        \fill (0 + \i + \k, 0 - \j) circle (0.1cm);
        \fill (0 + \j + \k, 0 + \i) circle (0.1cm);
        \fill (0 + \k, 0 + \i + \j) circle (0.1cm);
        \fill (0 - \i + \k, 0 + \j) circle (0.1cm);
        \fill (0 - \i + \k, 0 - \j) circle (0.1cm);
        \fill (0 + \j + \k, 0 - \i) circle (0.1cm);
        \fill (0 + \k, 0 - \i - \j) circle (0.1cm);
        \fill (0 + \i, 0 - \j) circle (0.1cm);
        \fill (0 - \i, 0 - \j) circle (0.1cm);
        \fill (0, 0 - \i - \j) circle (0.1cm);
        \fill (0 + \i + \j, 0 - \k) circle (0.1cm);
        \fill (0 + \i, 0 - \j - \k) circle (0.1cm);
        \fill (0 + \j, 0 + \i - \k) circle (0.1cm);
        \fill (0 - \j, 0 + \i - \k) circle (0.1cm);
        \fill (0 - \i - \j, 0 - \k) circle (0.1cm);
        \fill (0 - \i, 0 - \j - \k) circle (0.1cm);
        \fill (0 + \j, 0 - \i - \k) circle (0.1cm);
        \fill (0 - \j, 0 - \i - \k) circle (0.1cm);
        \fill (0, 0 - \i - \j - \k) circle (0.1cm);
        \fill (0 + \i, 0 + \j) circle (0.1cm);
        \fill (0, 0 + \i + \j) circle (0.1cm);
        \fill (0 - \i, 0 + \j) circle (0.1cm);
        \fill (0 + \i + \j, 0 + \k) circle (0.1cm);
        \fill (0 + \i, 0 + \j + \k) circle (0.1cm);
        \fill (0 + \j, 0 + \i + \k) circle (0.1cm);
        \fill (0, 0 + \i + \j + \k) circle (0.1cm);
        \fill (0 - \j, 0 + \i + \k) circle (0.1cm);
        \fill (0 - \i, 0 + \j + \k) circle (0.1cm);
        \fill (0 - \i - \j, 0 + \k) circle (0.1cm);
        \fill (0 + \j, 0 - \i + \k) circle (0.1cm);
        \fill (0 - \j, 0 - \i + \k) circle (0.1cm);
      \end{tikzpicture}
    }
  \end{wide}
}

\def\figureLatticeTiling{%
  \begin{tikzpicture}[scale = 0.4]
    \GridWithCoordinates((-6,-6),(6,6));
    \begin{scope}
      \clip (-6, -6) rectangle (6, 6);
      \foreach \x in {-4, 0, 4} {
        \foreach \y in {-4, 0, 4} {
          \fill (\x, \y) circle (0.1cm);
          \TwoTaxicabBallOutlines((\x,\y),1,2);
        }
      }
      \foreach \x in {-6, -2, 2, 6} {
        \foreach \y in {-6, -2, 2, 6} {
          \fill (\x, \y) circle (0.1cm);
          \TwoTaxicabBallOutlines((\x,\y),1,2);
        }
      }
    \end{scope}
  \end{tikzpicture}
}

\def\figureTreeTilingBalls{%
    \begin{tikzpicture}
      \pgfmathsetmacro\i{3}
      \pgfmathsetmacro\j{1.5}
      \pgfmathsetmacro\k{0.75}
      \pgfmathsetmacro\l{0.375}
      \pgfmathsetmacro\m{0.1875}

      \PlusShapedCross(\i,(0,0));
      \FreeGroupTilingOneArm((\i,0),0);
      \FreeGroupTilingOneArm((\i,0),90);
      \FreeGroupTilingOneArm((\i,0),-90);

      \node[label = 60:$e_{F_2} a b^2$] at (\i, \j + \k) {};
    \end{tikzpicture} 
}

\def\figureTreeTilingLines{%
  \begin{tikzpicture}
    \pgfmathsetmacro\i{3}
    \pgfmathsetmacro\j{1.5}
    \pgfmathsetmacro\k{0.75}
    \pgfmathsetmacro\l{0.375}
    \pgfmathsetmacro\m{0.1875}
    \pgfmathsetmacro\n{0.09375}

    \path (0, 0) edge[draw = none] node[auto] {$a$} (\i, 0);
    \path (0, 0) edge[draw = none] node[auto] {$b$} (0, \i);
    \path (0, 0) edge[draw = none] node[auto] {$a^{-1}$} (-\i, 0);
    \path (0, 0) edge[draw = none] node[auto] {$b^{-1}$} (0, -\i);

    \lineSolidDottedSolidOfVaryingSizesWithArms(\i,\j,\k,(0,0),0);
      \lineSolidDottedSolidOfVaryingSizesWithArms(\l,\m,\n,(\i+\j+\k,0),0);
    \lineSolidDottedSolidOfVaryingSizesWithArms(\i,\j,\k,(0,0),180);
      \lineSolidDottedSolidOfVaryingSizesWithArms(\l,\m,\n,(-\i-\j-\k,0),180);
    \begin{scope}[shift = {(\i + \j + \k, 0)}]
      \draw[densely dotted]
            (0, 0) -- (0, \l)
            (0, 0) -- (0, -\l);
      \fill (0, \l) circle (0.666 * \l * 0.1cm);
      \fill (0, -\l) circle (0.666 * \l * 0.1cm);
    \end{scope}
    \begin{scope}[shift = {(-\i - \j - \k, 0)}]
      \draw[densely dotted]
            (0, 0) -- (0, \l)
            (0, 0) -- (0, -\l);
      \fill (0, \l) circle (0.666 * \l * 0.1cm);
      \fill (0, -\l) circle (0.666 * \l * 0.1cm);
    \end{scope}

    \draw[densely dotted] 
          (0, 0) -- (0, \i)
          (0, 0) -- (0, -\i);
    \fill (0, -\i) circle (0.666 * \j * 0.1cm);

    \lineSolidDottedSolidOfVaryingSizesWithArms(\k,\l,\m,(\i,\j),0);
    \lineSolidDottedSolidOfVaryingSizesWithArms(\k,\l,\m,(\i,\j),180);
    \draw[densely dotted] (\i, \j) -- (\i, \j + \k);
    \fill (\i, \j + \k) circle (0.666 * \k * 0.1cm);

    \lineSolidDottedSolidOfVaryingSizesWithArms(\k,\l,\m,(-\i,\j),0);
    \lineSolidDottedSolidOfVaryingSizesWithArms(\k,\l,\m,(-\i,\j),180);
    \draw[densely dotted] (-\i, \j) -- (-\i, \j + \k);
    \fill (-\i, \j + \k) circle (0.666 * \k * 0.1cm);

    \lineSolidDottedSolidOfVaryingSizesWithArms(\j,\k,\l,(0,\i),0);
    \lineSolidDottedSolidOfVaryingSizesWithArms(\j,\k,\l,(0,\i),180);
    \draw[densely dotted] (0, \i) -- (0, \i + \j);
    \fill (0, \i + \j) circle (0.666 * \k * 0.1cm);
  \end{tikzpicture}
}

\def\figureTilingOfSpherePiEighth{%
  \begin{wide}
      \subfloat[The circles $E_m$, drawn solid, about $m \in \setOf{m_0, S} \cup \setOf{m_1, m_2, m_3, m_4} \cup \setOf{c_1, c_2, \dotsc, c_8}$.]{
        \includegraphics[trim = 65px 68px 35px 85px, clip]{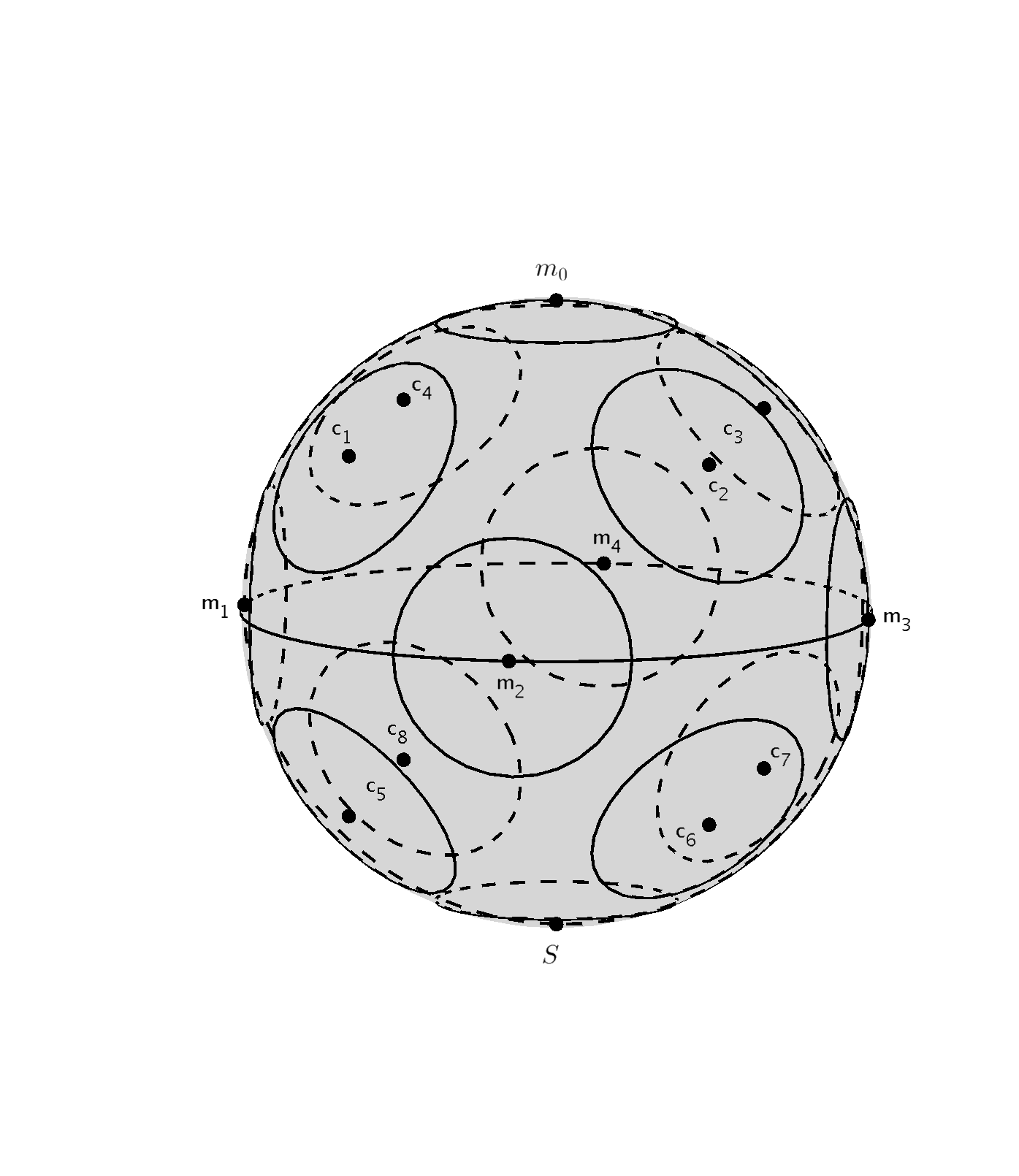}
      }
      \subfloat[The boundaries of the curved circular disks $E'_m$, drawn dotted, about $m \in \setOf{m_0, S} \cup \setOf{m_1, m_2, m_3, m_4}$.]{
        \includegraphics[trim = 65px 68px 35px 85px, clip]{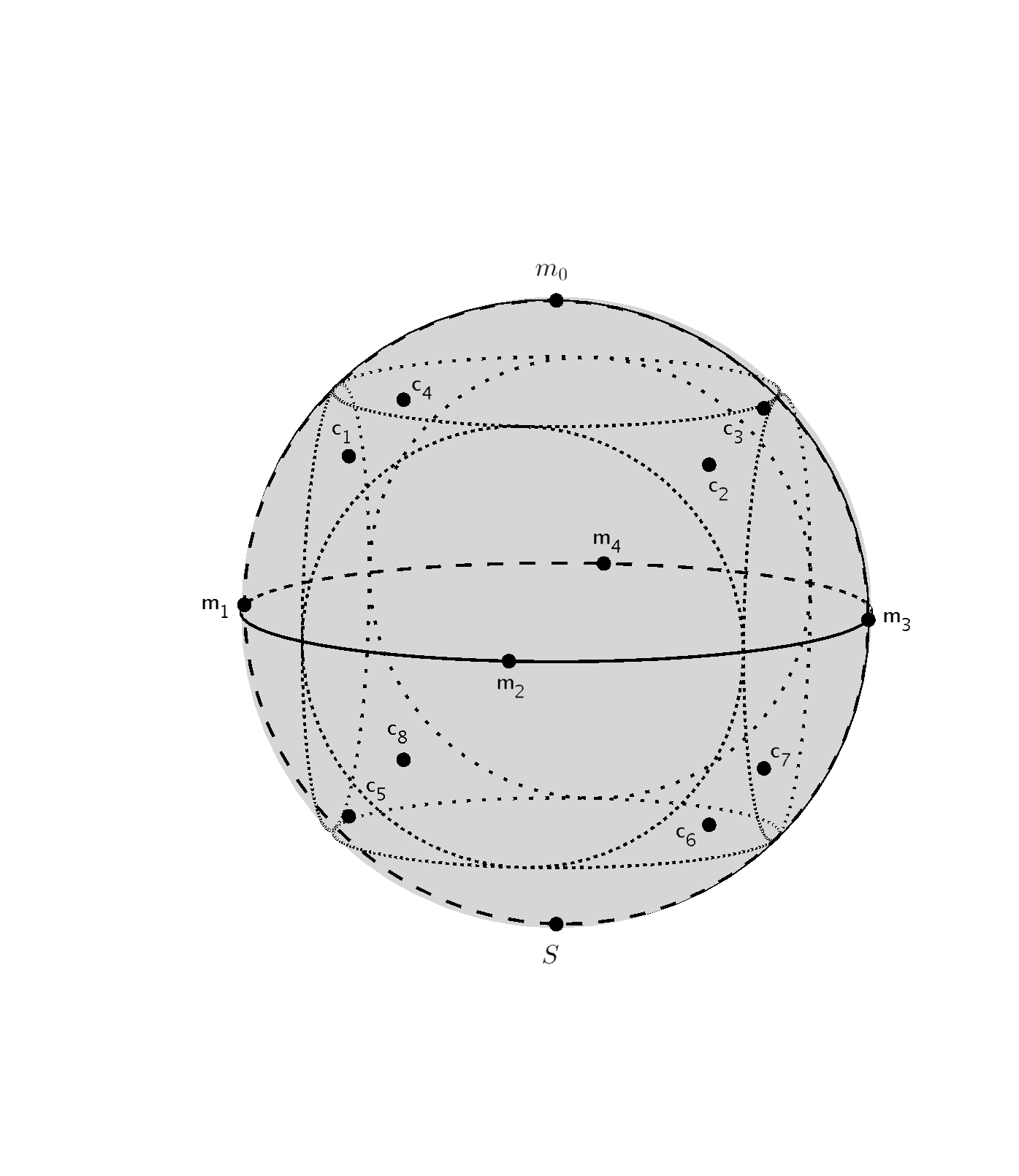}
      }
  \end{wide}
  \begin{wide}
      \subfloat[The boundaries of the curved circular disks $E'_m$, drawn dotted, about $m \in \setOf{m_0, S} \cup \setOf{m_1, m_2, m_3, m_4} \cup \setOf{c_1, c_2, \dotsc, c_8}$.]{
        \includegraphics[trim = 65px 68px 35px 85px, clip]{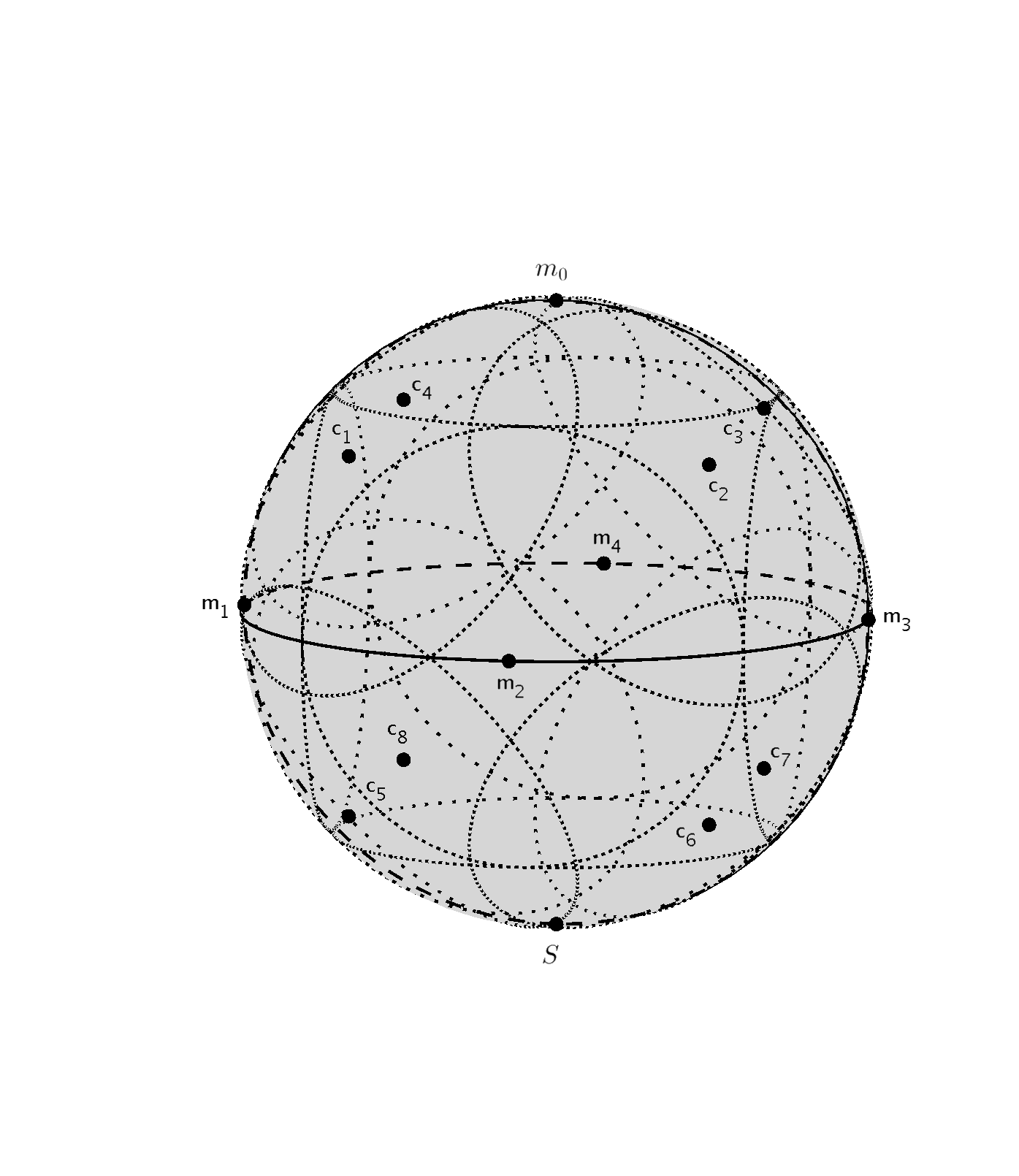}
      }
      \subfloat[The circles $E_m$, drawn solid, and the boundaries of the curved circular disks $E'_m$, drawn dotted, about $m \in \setOf{m_0, S} \cup \setOf{m_1, m_2, m_3, m_4} \cup \setOf{c_1, c_2, \dotsc, c_8}$.]{
        \includegraphics[trim = 65px 68px 35px 85px, clip]{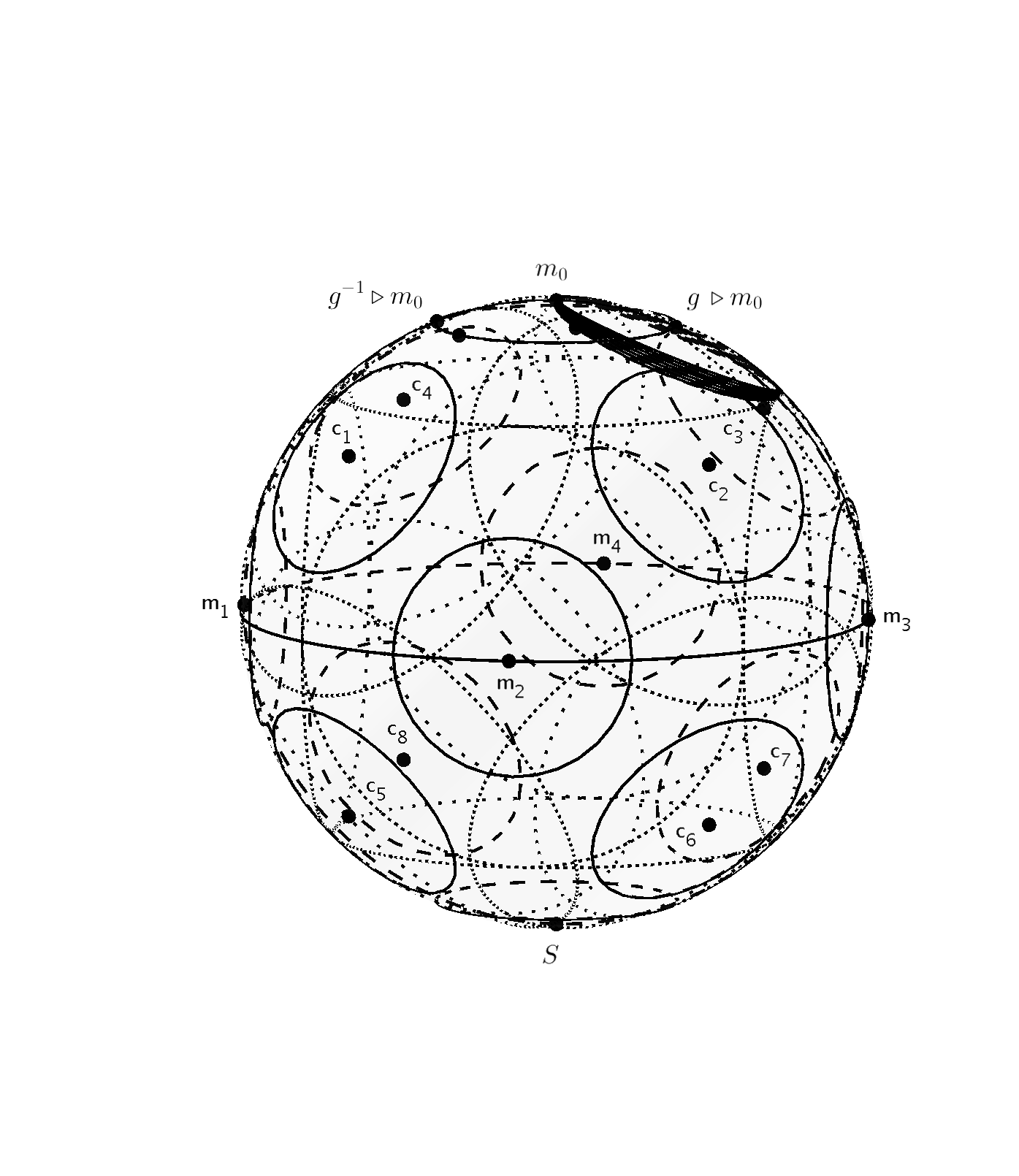}
      }
  \end{wide}
}

\def\figureExistenceOfTiling{%
  \begin{tikzpicture}[> = To]
    \node (m) at (0,0) {$m$};
    \node[below right = of m] (x) {$\bullet$};
    \node[above right = of x] (t) {$t$};

    \draw[->] (m) edge[bend left] node[right] {$e'$} (x);
    \draw[->] (x) edge[bend left] node[left] {$(g')^{-1} G_0$} (m); 
    \draw[->] (t) edge[bend left] node[right] {$e$} (x);
    \draw[->] (t) edge[bend right] node[above] {$g (g')^{-1} G_0$} (m); 
  \end{tikzpicture}
}

\def\figureUpperBoundOfTilingCapFolnerNet{%
  \begin{wide}
    \begin{minipage}[c]{0.2\textwidth} 
      \begin{tikzpicture}[circle dotted/.style = {dash pattern = on .05mm off 1.5pt, line cap = round}] 
        \pgfmathsetmacro\E{0.4} 
        \pgfmathsetmacro\Ep{0.65} 
        \pgfmathsetmacro\xFi{1 - 0.5} 
        \pgfmathsetmacro\yFi{1 - 0.5} 
        \pgfmathsetmacro\wFi{2.7} 
        \pgfmathsetmacro\hFi{3} 

        \begin{scope} 
          \clip (\xFi, \yFi) rectangle (\xFi + \wFi, \yFi + \E + \Ep) 
                (\xFi, \yFi + \hFi) rectangle (\xFi + \wFi, \yFi + \hFi - \E - \Ep) 
                (\xFi, \yFi + \E + \Ep) rectangle (\xFi + \E + \Ep, \yFi + \hFi - \E - \Ep) 
                (\xFi + \wFi, \yFi + \E + \Ep) rectangle (\xFi + \wFi - \E - \Ep, \yFi + \hFi - \E - \Ep); 
          \draw[line width = 0.5pt, dashed, pattern = north east lines, pattern color = gray] (\xFi, \yFi) rectangle (\xFi + \wFi, \yFi + \hFi); 
        \end{scope}
        \foreach \x in {1, ..., 2} { 
          \foreach \y in {1, ..., 3} {
            \draw[fill = white] (\x, \y) circle (1.25pt);
          }
        }
        \foreach \x in {0, 3, 4} { 
          \foreach \y in {0, ..., 4} {
            \fill (\x, \y) circle (1.25pt);
          }
        }
        \foreach \x in {1, 2} { 
          \foreach \y in {0, 4} {
            \fill (\x, \y) circle (1.25pt);
          }
        }
        \foreach \x in {0, ..., 4} {
          \foreach \y in {0, ..., 4} {
            \draw (\x - \E, \y - \E) rectangle (\x + \E, \y + \E); 
            \draw[gray, dashdotted] (\x - \Ep, \y - \Ep) rectangle (\x + \Ep, \y + \Ep); 
          }
        }
        \draw[line width = 0.5pt, dashed] (\xFi, \yFi) rectangle (\xFi + \wFi, \yFi + \hFi); 
        \draw[line width = 0.75pt, circle dotted] (\xFi + \E, \yFi + \E) rectangle (\xFi + \wFi - \E, \yFi + \hFi - \E); 
        \draw[line width = 0.75pt, circle dotted] (\xFi + \E + \Ep, \yFi + \E + \Ep) rectangle (\xFi + \wFi - \E - \Ep, \yFi + \hFi - \E - \Ep); 
      \end{tikzpicture}
    \end{minipage}%
    \hfill
    \begin{minipage}[c]{0.8\textwidth}
      The whole space is $M$; the dots and circles are the elements of the tiling $T$; for each element $t \in T$, the region enclosed by the rectangle with solid border centred at $t$ is the set $t \isSemiActedUponBy E$ and the region enclosed by the rectangle with dash-dotted border centred at $t$ is the set $t \isSemiActedUponBy E'$; the region enclosed by the rectangle with dashed border is $F_i$; the region enclosed by the largest rectangle with dotted border is $F_i^{-E}$; the region enclosed by the smallest rectangle with dotted border is $(F_i^{-E})^{-E'}$, which, in this depiction, is equal to $F_i^{-E''}$, where $E'' = E' \cdot (G_0 \cdot E)$; the circles are the elements of $T_i = T \cap F_i^{-E}$, which, in this depiction, is equal to $T_i' = T \cap ((F_i^{-E})^{-E'})^{+E'}$; the hatched region is $\boundaryOf_{E''}^- F_i = F_i \smallsetminus F_i^{-E''}$. 
    \end{minipage}
  \end{wide}
}

\def\figureTreeEntropyDoesIncreaseGlobalTransition{%
  \begin{tikzpicture}[scale = 0.5]
    \pgfmathsetmacro\i{3}
    \pgfmathsetmacro\j{1.5}
    \pgfmathsetmacro\k{0.75}
    \pgfmathsetmacro\l{0.375}

    \freeGroupCayleyGraph(solid);

    \draw (0 + \i + \j, 0) circle (0.1666cm);
    \draw (0 + \i, 0 + \j) circle (0.1666cm);
    \draw (0 + \i, 0 - \j) circle (0.1666cm);
    \fill (0 + \j, 0 + \i) circle (0.1666cm);
    \fill (0, 0 + \i + \j) circle (0.1666cm);
    \fill (0 - \j, 0 + \i) circle (0.1666cm);
    \fill (0 - \i, 0 + \j) circle (0.1666cm);
    \fill (0 - \i - \j, 0) circle (0.1666cm);
    \fill (0 - \i, 0 - \j) circle (0.1666cm);
    \fill (0 + \j, 0 - \i) circle (0.1666cm);
    \fill (0 - \j, 0 - \i) circle (0.1666cm);
    \fill (0, 0 - \i - \j) circle (0.1666cm);

    \draw (4, 3) edge[|->, bend left] node[above] {$\Delta$} (8, 3);

    \begin{scope}[shift = {(12, 0)}]
      \freeGroupCayleyGraph(solid);

      \draw (0 + \i, 0) circle (0.1666cm);
      \fill (0, 0 + \i) circle (0.1666cm);
      \fill (0 - \i, 0) circle (0.1666cm);
      \fill (0, 0 - \i) circle (0.1666cm);
    \end{scope}
  \end{tikzpicture}
}

\def\figureEntropyBoundedAboveIfSpecialTilingExists{%
  \begin{minipage}[c]{0.5\textwidth} 
    \begin{tikzpicture}[circle dotted/.style = {dash pattern = on .05mm off 1.5pt, line cap = round}] 
      \pgfmathsetmacro\E{0.4} 
      \pgfmathsetmacro\Ep{0.7} 
      \pgfmathsetmacro\xFi{1 - 0.5} 
      \pgfmathsetmacro\yFi{1 - 0.5} 
      \pgfmathsetmacro\wFi{2.7} 
      \pgfmathsetmacro\hFi{3} 

      \foreach \x in {0, ..., 4} {
        \foreach \y in {0, ..., 4} {
          \path[fill] (\x, \y) circle (1pt); 
          \draw (\x - \E, \y - \E) rectangle (\x + \E, \y + \E); 
        }
      }
      \draw[dashed, pattern = north east lines] (\xFi, \yFi) rectangle (\xFi + \wFi, \yFi + \hFi); 
      \foreach \x in {1, ..., 2} {
        \foreach \y in {1, ..., 3} {
          \draw[fill = white] (\x - \E, \y - \E) rectangle (\x + \E, \y + \E); 
          \draw (\x, \y) circle (1pt); 
        }
      }
      \draw[line width = 0.75pt, circle dotted] (\xFi + \E, \yFi + \E) rectangle (\xFi + \wFi - \E, \yFi + \hFi - \E); 
    \end{tikzpicture}
  \end{minipage}%
  \hfill
  \begin{minipage}[c]{0.5\textwidth}
    The whole space is $M$; the dots and circles are the elements of the tiling $T$; for each element $t \in T$, the region enclosed by the rectangle with solid border centred at $t$ is the set $t \isSemiActedUponBy E$; the region enclosed by the rectangle with dashed border is $F_i$; the region enclosed by the rectangle with dotted border is $F_i^{-E}$; the circles are the elements of $T_i = T \cap F_i^{-E}$; the hatched region is $F_i^* = F_i \smallsetminus (\bigcup_{t \in T_i} t \isSemiActedUponBy E)$. 
  \end{minipage}
}

\def\figureNonMaximalEntropyButPreInjective{%
  \begin{tikzpicture}[scale = 0.5]
    \pgfmathsetmacro\i{3}
    \pgfmathsetmacro\j{1.5}
    \pgfmathsetmacro\k{0.75}
    \pgfmathsetmacro\l{0.375}
    \pgfmathsetmacro\m{0.1875}
    \pgfmathsetmacro\n{0.09375}
    \pgfmathsetmacro\a{40}

    \foreach \r in {0, 120, 240} {
      \begin{scope}[rotate = \r]
        \draw (0, 0) -- (\i, 0);
        \foreach \r in {\a, -\a} {
          \begin{scope}[shift = {(\i, 0)}, rotate = \r]
            \draw (0, 0) -- (\j, 0);
            \foreach \r in {\a, -\a} {
              \begin{scope}[shift = {(\j, 0)}, rotate = \r]
                \draw (0, 0) -- (\k, 0);
                \foreach \r in {\a, -\a} {
                  \begin{scope}[shift = {(\k, 0)}, rotate = \r]
                    \draw (0, 0) -- (\l, 0);
                    \foreach \r in {\a, -\a} {
                      \begin{scope}[shift = {(\l, 0)}, rotate = \r]
                        \draw (0, 0) -- (\m, 0);
                        \foreach \r in {\a, -\a} {
                          \begin{scope}[shift = {(\m, 0)}, rotate = \r]
                            \draw (0, 0) -- (\n, 0);
                          \end{scope}
                        }
                      \end{scope}
                    }
                  \end{scope}
                }
              \end{scope}
            }
          \end{scope}
        }
      \end{scope}
    }
    \fill (0, 0) circle (0.125);
    \node[below right = 0 and -0.1] at (0, 0) {$e_G$};
    \fill (\i, 0) circle (0.1);
    \node[above] at (\i, 0) {$g$};
    \begin{scope}[shift = {(\i, 0)}, rotate = \a]
      \fill (\j, 0) circle (0.075);
      \node[below right = 0 and -0.2] at (\j, 0) {$g s$};
    \end{scope}
    \begin{scope}[shift = {(\i, 0)}, rotate = -\a]
      \fill (\j, 0) circle (0.075);
      \node[above right = 0 and -0.2] at (\j, 0) {$g s'$};
    \end{scope}
  \end{tikzpicture}
}

\def\figureMLiberationAIsEitherOutsideOfAOrInsideOfClosureOfA{%
  \begin{tikzpicture}[> = To]
    \node (m) at (0,0) {$m$};
    \node[right = of m] (x) {$\bullet$};
    \node[below = of m] (y) {$\bullet$};
    \node[draw, shape = circle, minimum size = 2cm, above right = -0.7cm of x] (A) {$A$};

    \draw[->] (m) edge node[above] {$e'$} (x);
    \draw[->] (m) edge node[left] {$e$} (y);
    \draw[->] (y) edge[bend right] node[right] {$g^{-1} \cdot e'$} (x);
  \end{tikzpicture}
}

\def\figureThereIsARhoKappaThetaTiling{%
  \begin{wide}
    \begin{tikzpicture}[> = To, x=1cm, y=1cm, circle dotted/.style = {dash pattern = on .05mm off 1.5pt, line cap = round}]
      \pgfmathsetmacro\t{0.5} 
      \pgfmathsetmacro\k{1} 
      \pgfmathsetmacro\c{5} 
      \pgfmathsetmacro\r{0.0625} 

      \draw[fill = black] (0,0) circle[radius = \r] node {}; 
      \draw (0,0) circle[radius = {1 * ((2 * \t) + \k + 1)}] node {}; 
      \draw (0,0) circle[radius = {2 * ((2 * \t) + \k + 1)}] node {}; 

      \foreach \a in {1,2,...,\c} {
        \draw[fill = black] (\a * 360 / \c: {1 * ((2 * \t) + \k + 1)}) circle[radius = \r] node[below] {$m_{i, \a}$};
        \draw[circle dotted] (\a * 360 / \c: {1 * ((2 * \t) + \k + 1)}) circle[radius = \t];
        \draw[dashdotted] (\a * 360 / \c: {1 * ((2 * \t) + \k + 1)}) circle[radius = {\t + \k}];
        \draw[dashed] (\a * 360 / \c: {1 * ((2 * \t) + \k + 1)}) circle[radius = {(2 * \t) + \k}];
      }
      \draw[fill = black] (1.75 * 360 / \c: {1.75 * ((2 * \t) + \k + 1)}) circle[radius = \r] node[above] {$m$}
                          -- 
                          (1.75 * 360 / \c: {1 * ((2 * \t) + \k + 1)}) circle[radius = \r] node[below] {$m'$}
                          -- 
                          (2 * 360 / \c: {1 * ((2 * \t) + \k + 1)}) node[below] {$\phantom{m_{i, 2}}\ \qquad = m''$};
    \end{tikzpicture}
  \end{wide}
}

\def\figureTechnicalInequalityForTheoremEntropyLessIfRealSubsetsForRhoKappaThetaTiling{%
  \begin{minipage}[c]{0.5\textwidth} 
    \begin{tikzpicture}[circle dotted/.style = {dash pattern = on .05mm off 1.5pt, line cap = round}] 
      \pgfmathsetmacro\t{0.3} 
      \pgfmathsetmacro\k{0.1} 
      \pgfmathsetmacro\xFi{1 - 0.5} 
      \pgfmathsetmacro\yFi{1 - 0.5} 
      \pgfmathsetmacro\wFi{2.7} 
      \pgfmathsetmacro\hFi{3} 

      \foreach \x in {0, ..., 4} {
        \foreach \y in {0, ..., 4} {
          \path[fill] (\x, \y) circle (1pt); 
          \draw (\x - \t, \y - \t) rectangle (\x + \t, \y + \t); 
          \draw[dashdotted] (\x - \t - \k, \y - \t - \k) rectangle (\x + \t + \k, \y + \t + \k); 
        }
      }
      \draw[dashed, pattern = north east lines] (\xFi, \yFi) rectangle (\xFi + \wFi, \yFi + \hFi); 
      \foreach \x in {1, ..., 2} {
        \foreach \y in {1, ..., 3} {
          \draw[fill = white] (\x - \t, \y - \t) rectangle (\x + \t, \y + \t); 
          \draw (\x, \y) circle (1pt); 
        }
      }
      \draw[line width = 0.75pt, circle dotted] (\xFi + \t + \k, \yFi + \t + \k) rectangle (\xFi + \wFi - \t - \k, \yFi + \hFi - \t - \k); 
    \end{tikzpicture}
  \end{minipage}%
  \hfill
  \begin{minipage}[c]{0.5\textwidth} 
    The whole space is $M$; the dots and circles are the elements of the set $T$; for each element $t \in T$, the region enclosed by the rectangle with solid border about $t$ is the set $\ball(t, \theta)$ and the region enclosed by the rectangle with dash-dotted border about $t$ is the set $\ball(t, \theta)^{+\kappa}$; the region enclosed by the rectangle with dashed border is $F$; the region enclosed by the rectangle with dotted border is $F^{-(\theta + \kappa)}$; the circles are the elements of $S = T \cap F^{-(\theta + \kappa)}$; the hatched region is the set $F \smallsetminus (\bigcup_{s \in S} \ball(s, \theta))$. 
  \end{minipage}
}

\def\figureSetContainedInUnionOfBallsAndOfInteriorOfSet{%
  \begin{minipage}[c]{\textwidth / 2} 
    \begin{tikzpicture}[circle dotted/.style = {dash pattern = on .05mm off 1.5pt, line cap = round}] 
      \pgfmathsetmacro\t{0.4} 
      \pgfmathsetmacro\tp{0.7} 
      \pgfmathsetmacro\xF{1 - 0.55} 
      \pgfmathsetmacro\yF{1 - 0.55} 
      \pgfmathsetmacro\wF{2.8} 
      \pgfmathsetmacro\hF{3.1} 

      \foreach \x in {1, ..., 2} {
        \foreach \y in {1, ..., 3} {
          \draw (\x, \y) circle (1pt);
        }
      }
      \foreach \x in {0, 3, 4} {
        \foreach \y in {0, ..., 4} {
          \fill (\x, \y) circle (1pt);
        }
      }
      \foreach \x in {1, 2} {
        \foreach \y in {0, 4} {
          \fill (\x, \y) circle (1pt);
        }
      }
      \foreach \x in {0, ..., 4} {
        \foreach \y in {0, ..., 4} {
          \draw (\x - \t, \y - \t) rectangle (\x + \t, \y + \t); 
          \draw[gray, dashdotted] (\x - \tp, \y - \tp) rectangle (\x + \tp, \y + \tp); 
        } 
      }
      \draw[line width = 0.5pt, dashed] (\xF, \yF) rectangle (\xF + \wF, \yF + \hF); 
      \draw[line width = 0.75pt, circle dotted] (\xF + \t, \yF + \t) rectangle (\xF + \wF - \t, \yF + \hF - \t); 
      \draw[line width = 0.75pt, circle dotted] (\xF + \t + \tp, \yF + \t + \tp) rectangle (\xF + \wF - \t - \tp, \yF + \hF - \t - \tp); 
    \end{tikzpicture}
  \end{minipage}%
  \begin{minipage}[c]{\textwidth / 2}
    The whole space is $M$; the dots and circles are the elements of the set $T$; for each element $t \in T$, the region enclosed by the rectangle with solid border about $t$ is the set $\ball(t, \theta)^{+\kappa}$ and the region enclosed by the rectangle with dash-dotted border about $t$ is the set $\ball(t, \theta')$; the region enclosed by the rectangle with dashed border is $F$; the region enclosed by the smallest rectangle with dotted border is $F^{-(\theta + \kappa + \theta')}$ and the region enclosed by the largest rectangle with dotted border is $F^{-(\theta + \kappa)}$; the circles are the elements of $S = T \cap F^{-(\theta + \kappa)}$.
  \end{minipage}
}

\definecolor{myred}{HTML}{e41a1c}
\definecolor{myblue}{HTML}{377eb8}
\definecolor{mygreen}{HTML}{4daf4a}
\definecolor{myviolet}{HTML}{984ea3}

\tikzset{%
  edge/.style = {draw = black},
  vertex/.style = {draw = black, thick},
  midpoint/.style = {draw = black, thick},
  boundary/.style = {draw = black, thick},
  find/.style = {draw = black, dotted, very thick}, 
  reflected/.style = {draw = black, densely dotted, thick},
  syncLeft/.style = {pattern = horizontal lines, pattern color = black, draw = none}, 
  syncMiddle/.style = {pattern = fivepointed stars, pattern color = black, draw = none},
  syncRight/.style = {pattern = dots, pattern color = black, draw = none}, 
  slowed/.style = {draw = black}, 
  freeze/.style = {draw = black, dashed}, 
  traversingThaw/.style = {draw = black, loosely dotted, thick},
  thawingThaw/.style = {draw = black, dashed},
  divide0/.style = {draw = black, solid},
  divideN/.style = {draw = black, densely dotted},
  reflectedDivide/.style = {draw = black},
  leftColour/.style = {draw = mygreen, pattern color = mygreen},
  leftFill/.style = {color = mygreen},
  rightColour/.style = {draw = myblue, pattern color = myblue},
  rightFill/.style = {color = myblue},
  overlayLeftColour/.style = {draw = myred, pattern color = myred},
  overlayLeftFill/.style = {color = myred},
  overlayRightColour/.style = {draw = myviolet, pattern color = myviolet},
  overlayRightFill/.style = {color = myviolet}
}

\def\figureSingularities{%
  \subfloat[A singularity of order $1$.]{%
    \resizebox{(\linewidth-1em)/2}{!}{%
      \begin{tikzpicture}[> = To, x = 1cm, y = {-1cm}]
        \pgfmathsetmacro{\l}{3} 
        \pgfmathsetmacro{\h}{\l + \l} 

        \coordinate (gg) at (0, 0); 
        \coordinate (dr) at (\l, \l); 
        \coordinate (rdl) at (0, 2 * \l); 

        \draw[divide0] (gg) -- (dr);
        \draw[reflectedDivide] (dr) -- (rdl);

        \foreach \type in {1, 2, ..., 15} {
          \pgfmathsetmacro{\v}{((2/3)^\type) * \l};
          \draw[divideN] (gg) -- (\v, {2 * \l - \v});
        }

        \draw[vertex] (0, 0) -- (0, 6)
                      (3, 0) -- (3, 6);
        \draw[edge] (0, 0) -- (3, 0)
                    (0, 6) -- (3, 6);
      \end{tikzpicture}
    }
    \label{figure:singularity-of-order-one}
  }
  \hfill
  \subfloat[A singularity of order $-1$.]{%
    \resizebox{(\linewidth-1em)/2}{!}{%
      \begin{tikzpicture}[> = To, x = 1cm, y = {-1cm}]
        \pgfmathsetmacro{\l}{3} 
        \pgfmathsetmacro{\h}{\l + \l} 

        \coordinate (gg) at (0, 0); 
        \coordinate (dr) at (\l, \l); 
        \coordinate (rdl) at (0, 2 * \l); 

        \draw[divide0] (0, 0) -- (-1, 1);
        \draw[reflectedDivide] (-1, 1) -- (2, 4);

        \begin{scope}
          \clip (0, 0) -- (0, 2) -- (2, 4) -- (2, 0) -- cycle;
          \foreach \type in {1, 2, ..., 15} {
            \pgfmathsetmacro{\v}{((2/3)^\type) * \l};
            \draw[divideN] (gg) -- (\v, {2 * \l - \v});
          }
        \end{scope}

        \draw[vertex] (-1, 0) -- (-1, 4)
                      (0, 0) -- (0, 4)
                      (2, 0) -- (2, 4);
        \draw[edge] (-1, 0) -- (2, 0)
                    (-1, 4) -- (2, 4);
      \end{tikzpicture}
    }
    \label{figure:singularity-of-order-minus-one}
  }
}

\def\figureMidpointOneEdge{%
  \begin{tikzpicture}[> = To, x = 1cm, y = {-1cm}]
    \pgfmathsetmacro{\l}{3} 
    \pgfmathsetmacro{\h}{2 * \l} 

    \draw[find] (0, 0) -- (\l, \l);
    \draw[reflected] (\l, \l) -- (0, \h);
    \draw[slowed] (0, 0) -- ({(1/3) * \h}, \h);
    \draw[midpoint] (0.5 * \l, 1.5 * \l) -- (0.5 * \l, \h);

    \draw[vertex] (0, 0) -- (0, \h); 
    \draw[vertex] (\l, 0) -- (\l, \h); 

    \draw[edge] (0, 0) -- (\l, 0); 
    \draw[edge] (0, \h) -- (\l, \h); 
  \end{tikzpicture}
}

\def\figureMidpointTwoEdges{%
  \begin{tikzpicture}[> = To, x = 1cm, y = {-1cm}] 
    \pgfmathsetmacro{\l}{3} 
    \pgfmathsetmacro{\r}{2} 
    \pgfmathsetmacro{\w}{\l + \r} 
    \pgfmathsetmacro{\h}{2 * \l + \r} 

    \begin{scope}[shift = {(\l, -1.5)}, scale = 0.5]
      \draw (0, 0) -- (-2.1213, 2.1213);
      \draw (0, 0) -- (1.414, 1.414);
    \end{scope}

    \coordinate (g) at (\l, 0); 

    \draw[find] (g) -- (0, \l);
    \draw[reflected] (0, \l) -- (\l, 2 * \l);
    \draw[slowed] (\l, 2 * \l) -- ({\l + (1/3) * \r}, \h);

    \draw[find] (g) -- (\w, \r);
    \draw[reflected] (\w, \r) -- (\l, 2 * \r);
    \draw[slowed] (\l, 2 * \r) -- ({0.5 * \w - (1/3) * 0.5 * \w}, \h);

    \draw[midpoint] (0.5 * \w, \h - 0.5 * \w) -- (0.5 * \w, \h);

    \draw[vertex] (0, 0) -- (0, \h);
    \draw[vertex] (g) -- (\l, \h);
    \draw[vertex] (\w, 0) -- (\w, \h);

    \draw[edge] (0, 0) -- (g); 
    \draw[edge] (g) -- (\w, 0); 
    \draw[edge] (0, \h) -- (\l, \h); 
    \draw[edge] (\l, \h) -- (\w, \h); 
  \end{tikzpicture}
}

\def\mazoyerOneLevel(#1,#2){
  \pgfmathsetmacro{\l}{#1} 
  \pgfmathsetmacro{\h}{\l + \l} 

  \coordinate (gg) at (0, 0); 
  \coordinate (dr) at (\l, \l); 
  \coordinate (rdl) at (0, 2 * \l); 

  \draw[divide0] (gg) -- (dr);
  \draw[reflectedDivide] (dr) -- (rdl);

  \foreach \type in {1, 2, ..., 15} {
    \pgfmathsetmacro{\v}{((2/3)^\type) * \l};
    \coordinate (b #2 \type) at (\v, {2 * \l - \v}); 
    \coordinate (be #2 \type) at (\v, \h); 
    \draw[divideN] (gg) -- (b #2 \type);
    \draw[boundary] (b #2 \type) -- (be #2 \type);
  }
}

\def\mazoyerTwoLevels(#1){%
  \pgfmathsetmacro{\le}{#1} 
  \mazoyerOneLevel(\le,0); 
  \foreach \type in {1, 2, ..., 15} {
    \begin{scope}[shift = (b 0 \type)]
      \pgfmathsetmacro{\typeMinusOne}{\type - 1};
      \mazoyerOneLevel({(1/3) * (2/3)^\typeMinusOne * \le},\type);
    \end{scope}
  }
}

\def\figureMazoyer{%
  \begin{tikzpicture}[> = To, x = 1cm, y = {-1cm}]
    \mazoyerTwoLevels(3);
    \draw[vertex] (0, 0) -- (0, 6)
                  (3, 0) -- (3, 6);
    \draw[edge] (0, 0) -- (3, 0)
                (0, 6) -- (3, 6);
  \end{tikzpicture}
}

\def\figureFSSPMazoyerWithTwoEdges{%
  \begin{tikzpicture}[> = To, x = 1cm, y = {-1cm}]
    \pgfmathsetmacro{\l}{3} 
    \pgfmathsetmacro{\r}{2} 
    \pgfmathsetmacro{\w}{\l + \r} 
    \pgfmathsetmacro{\h}{\l + \l + \r} 

    \begin{scope}[shift = {(\l, -1.5)}, scale = 0.5]
      \draw (0, 0) -- (-2.1213, 2.1213);
      \draw (0, 0) -- (1.414, 1.414);
    \end{scope}

    \coordinate (g) at (\l, 0); 
    \coordinate (lb) at (0, \l); 
    \coordinate (rb) at (\w, \r); 
    \coordinate (lg) at (\l, 2 * \l); 
    \coordinate (rg) at (\l, 2 * \r); 
    \coordinate (m) at (0.5 * \w, \l + 0.5 * \w); 
    \coordinate (mlf) at (0.5 * \l, \l + 0.5 * \l); 
    \coordinate (mlt) at (0.5 * \l, \h - 0.5 * \l); 
    \coordinate (mrf) at (\l + 0.5 * \r, \r + 0.5 * \r); 
    \coordinate (mrt) at (\l + 0.5 * \r, \h - 0.5 * \r); 
    \coordinate (lfb) at ($(mlf) + (- 0.5 * \l, 0.5 * \l)$); 
    \coordinate (rfb) at ($(mrf) + (0.5 * \r, 0.5 * \r)$); 
    \coordinate (lfg) at ($(mlf) + (0.5 * \l, 0.5 * \l)$); 
    \coordinate (rfg) at ($(mrf) + (- 0.5 * \r, 0.5 * \r)$); 

    \begin{scope} 
      \clip (0, 0) -- (0, 2 * \l) -- (0.5 * \l, 1.5 * \l) -- (\l, 2 * \l) -- (\l, 0) -- cycle;
      \begin{scope}[xscale = -1, shift = {(-\l, 0)}]
        \mazoyerTwoLevels(\l);
      \end{scope}
    \end{scope}
    \draw[freeze] (0.5 * \l, 1.5 * \l) -- (0, 2 * \l)
                  (0.5 * \l, 1.5 * \l) -- (\l, 2 * \l); 
    \begin{scope}[shift = {(0, \h - 2 * \l)}] 
      \begin{scope}
        \clip (0, 2 * \l) -- (0.5 * \l, 1.5 * \l) -- (\l, 2 * \l) -- cycle;
        \begin{scope}[xscale = -1, shift = {(-\l, 0)}]
          \mazoyerTwoLevels(\l);
        \end{scope}
      \end{scope}
      \draw[thawingThaw] (0.5 * \l, 1.5 * \l) -- (0, 2 * \l)
                  (0.5 * \l, 1.5 * \l) -- (\l, 2 * \l); 
    \end{scope}

    \begin{scope}[shift = {(\l, 0)}]
      \begin{scope} 
        \clip (0, 0) -- (0, 2 * \r) -- (0.5 * \r, 1.5 * \r) -- (\r, 2 * \r) -- (\r, 0) -- cycle;
        \mazoyerTwoLevels(\r);
      \end{scope}
      \draw[freeze] (0.5 * \r, 1.5 * \r) -- (0, 2 * \r)
                    (0.5 * \r, 1.5 * \r) -- (\r, 2 * \r); 
      \begin{scope}[shift = {(0, \h - 2 * \r)}] 
        \begin{scope} 
          \clip (0, 2 * \r) -- (0.5 * \r, 1.5 * \r) -- (\r, 2 * \r) -- cycle;
          \mazoyerTwoLevels(\r);
        \end{scope}
        \draw[thawingThaw] (0.5 * \r, 1.5 * \r) -- (0, 2 * \r)
                    (0.5 * \r, 1.5 * \r) -- (\r, 2 * \r); 
      \end{scope}
    \end{scope}

    \draw[slowed] (g) -- (mlf);
    \draw[slowed] (g) -- (mrf);
    \draw[slowed] (\l, \r + \r) -- (0.5 * \w, {\r + \r + 1.5 * (\l - \r)});

    \draw[white, semithick] (mlf) -- (lfg)
                            (mrf) -- (rfg);
    \draw[freeze] (mlf) -- (lfg)
                  (mlf) -- (lfb)
                  (mrf) -- (rfb)
                  (mrf) -- (rfg);

    \draw[traversingThaw] (m) -- (mlt)
                          (m) -- (mrt);

    \draw[white, semithick] (0, \h) -- (0.5 * \l, \h - 0.5 * \l) -- (\l, \h);
    \draw[white, semithick] (\w, \h) -- (\l + 0.5 * \r, \h - 0.5 * \r) -- (\w - \r, \h);
    \draw[thawingThaw] (0.5 * \l, \h - 0.5 * \l) -- (0, \h)
                       (0.5 * \l, \h - 0.5 * \l) -- (\l, \h);
    \draw[thawingThaw] (\l + 0.5 * \r, \h - 0.5 * \r) -- (\w, \h)
                       (\l + 0.5 * \r, \h - 0.5 * \r) -- (\w - \r, \h);

    \draw[midpoint] (mlf) -- (mlt)
                    (mrf) -- (mrt);

    \draw[vertex] (0, 0) -- (0, \h)
                  (g) -- (\l, \h)
                  (\w, 0) -- (\w, \h);

    \draw[edge] (0, 0) -- (g) -- (\w, 0)
                (0, \h) -- (\l, \h) -- (\w, \h);
  \end{tikzpicture}
}

\def\figureFSSPWithTwoEdgesAndGeneralInBetween{%
  \begin{tikzpicture}[> = To, x = 1cm, y = {-1cm}] 
    \pgfmathsetmacro{\l}{3} 
    \pgfmathsetmacro{\r}{2} 
    \pgfmathsetmacro{\o}{1.333} 
    \pgfmathsetmacro{\x}{max(\l, \r)} 
    \pgfmathsetmacro{\y}{min(\l, \r)} 
    \pgfmathsetmacro{\w}{\l + \r} 
    \pgfmathsetmacro{\h}{\x + \x + \y} 

    \begin{scope}[shift = {(\l, -1.5)}, scale = 0.5]
      \draw (0, 0) -- (-2.1213, 2.1213);
      \draw (0, 0) -- (1.414, 1.414);
    \end{scope}

    \coordinate (g) at (\l, 0); 
    \coordinate (lb) at (0, \l); 
    \coordinate (rb) at (\w, \r); 
    \coordinate (lg) at (\l, 2 * \l); 
    \coordinate (rg) at (\l, 2 * \r); 
    \coordinate (m) at (0.5 * \w, \x + 0.5 * \w); 
    \coordinate (mlf) at (0.5 * \l, \l + 0.5 * \l); 
    \coordinate (mlt) at (0.5 * \l, \h - 0.5 * \l); 
    \coordinate (mrf) at (\l + 0.5 * \r, \r + 0.5 * \r); 
    \coordinate (mrt) at (\l + 0.5 * \r, \h - 0.5 * \r); 
    \coordinate (lfb) at ($(mlf) + (- 0.5 * \l, 0.5 * \l)$); 
    \coordinate (rfb) at ($(mrf) + (0.5 * \r, 0.5 * \r)$); 
    \coordinate (lfg) at ($(mlf) + (0.5 * \l, 0.5 * \l)$); 
    \coordinate (rfg) at ($(mrf) + (- 0.5 * \r, 0.5 * \r)$); 

    \draw[syncLeft] (g) -- (lb) -- (lfb) -- (mlf) -- (lg) -- cycle;
    \draw[syncRight] (g) -- (rb) -- (rfb) -- (mrf) -- (rg) -- cycle;

    \draw[syncLeft, thawingThaw]
        (0, \h) -- (0.5 * \l, \h - 0.5 * \l)
                -- (\l, \h);
    \draw[syncRight, thawingThaw]
        (\w, \h) -- (\l + 0.5 * \r, \h - 0.5 * \r)
                 -- (\w - \r, \h);

    \draw[find] (g) -- (lb);
    \draw[reflected] (lb) -- (mlf);
    \draw[find] (g) -- (rb);
    \draw[reflected] (rb) -- (\x, \y + \y);
    \draw[slowed] (\x, \y + \y) --  (0.5 * \w, {\y + \y + 1.5 * (\x - \y)});

    \draw[slowed] (g) -- (mlf);
    \draw[slowed] (g) -- (mrf);

    \draw[midpoint] (mlf) -- (mlt)
                    (mrf) -- (mrt);

    \draw[traversingThaw] (m) -- (mlt)
                          (m) -- (mrt);

    \draw[freeze] (mlf) -- (lfg)
                  (mlf) -- (lfb)
                  (mrf) -- (rfb)
                  (mrf) -- (rfg);

    \draw[vertex] (0, 0) -- (0, \h)
                  (g) -- (\l, \h)
                  (\w, 0) -- (\w, \h);

    \draw[edge] (0, 0) -- (g) -- (\w, 0)
                (0, \h) -- (\l, \h) -- (\w, \h);
  \end{tikzpicture}
}

\def\figureFSSPWithTwoEdgesAndGeneralAtTheLeft{%
  \begin{tikzpicture}[> = To, x = 1cm, y = {-1cm}] 
    \pgfmathsetmacro{\l}{3} 
    \pgfmathsetmacro{\r}{2} 
    \pgfmathsetmacro{\o}{1.333} 
    \pgfmathsetmacro{\w}{\l + \r} 
    \pgfmathsetmacro{\h}{2 * \w} 

    \coordinate (g) at (0, 0); 
    \coordinate (v) at (\l, 0); 
    \coordinate (iv) at (\l, \l); 
    \coordinate (lb) at (\l, \l); 
    \coordinate (rb) at (\w, \l + \r); 
    \coordinate (lg) at (0, 2 * \l); 
    \coordinate (rv) at (\l, \l + 2 * \r); 
    \coordinate (m) at (0.5 * \w, \w + 0.5 * \w); 
    \coordinate (mlf) at (0.5 * \l, \l + 0.5 * \l); 
    \coordinate (mlt) at (0.5 * \l, \h - 0.5 * \l); 
    \coordinate (mrf) at (\l + 0.5 * \r, \l + \r + 0.5 * \r); 
    \coordinate (mrt) at (\l + 0.5 * \r, \h - 0.5 * \r); 
    \coordinate (lfv) at ($(mlf) + (0.5 * \l, 0.5 * \l)$); 
    \coordinate (rfb) at ($(mrf) + (0.5 * \r, 0.5 * \r)$); 
    \coordinate (lfg) at ($(mlf) + (- 0.5 * \l, 0.5 * \l)$); 
    \coordinate (rfg) at ($(mrf) + (- 0.5 * \r, 0.5 * \r)$); 

    \draw[syncLeft] (g) -- (lb) -- (lfv) -- (mlf) -- (lg) -- cycle;
    \draw[syncRight] (iv) -- (rb) -- (rfb) -- (mrf) -- (rv) -- cycle;

    \draw[syncLeft, thawingThaw]
        (0, \h) -- (0.5 * \l, \h - 0.5 * \l)
                -- (\l, \h);
    \draw[syncRight, thawingThaw]
        (\w, \h) -- (\l + 0.5 * \r, \h - 0.5 * \r)
                 -- (\w - \r, \h);

    \draw[find] (g) -- (lb);
    \draw[reflected] (lb) -- (mlf); 
    \draw[find] (iv) -- (rb);
    \draw[reflected] (rb) -- (rv) -- (m);

    \draw[traversingThaw] (m) -- (mlt);
    \draw[traversingThaw] (m) -- (mrt);

    \draw[midpoint] (mlf) -- (mlt);
    \draw[midpoint] (mrf) -- (mrt);

    \draw[slowed] (g) -- (mlf)
                      -- (\w/2, \h - \w/2);
    \draw[slowed] (iv) -- (mrf);

    \draw[freeze] (mlf) -- (lfg);
    \draw[freeze] (mlf) -- (lfv);
    \draw[freeze] (mrf) -- (rfb);
    \draw[freeze] (mrf) -- (rfg);

    \draw[vertex] (v) -- (\l, \h);
    \draw[vertex] (0, 0) -- (0, \h);
    \draw[vertex] (\w, 0) -- (\w, \h);

    \draw[edge] (0, 0) -- (v); 
    \draw[edge] (v) -- (\w, 0); 
    \draw[edge] (0, \h) -- (\l, \h); 
    \draw[edge] (\l, \h) -- (\w, \h); 
  \end{tikzpicture}
}

\def\drawMidpointAtOf(#1,#2,#3){%
  \begin{scope}[shift = #1]
    \begin{scope}
      \clip (0, 3pt) rectangle (-3pt, -3pt);
      \fill[#2, draw = black] (0, 0) circle (2pt);
    \end{scope}
    \begin{scope}
      \clip (0, -3pt) rectangle (3pt, 3pt);
      \fill[#3, draw = black] (0, 0) circle (2pt);
    \end{scope}
  \end{scope}
}

\def\drawMidpointAtOfThreeEdges(#1,#2,#3,#4){%
  \begin{scope}[shift = #1]
    \fill[#3, draw = none] (-1pt, -2pt) rectangle (1pt, 2pt);
    \draw[black] (-1pt, -2pt) -- (1pt, -2pt)
                 (-1pt, 2pt) -- (1pt, 2pt);
    \begin{scope}[shift = {(-1pt,0)}]
      \clip (0, 3pt) rectangle (-3pt, -3pt);
      \fill[#2, draw = black] (0, 0) circle (2pt);
    \end{scope}
    \begin{scope}[shift = {(1pt,0)}]
      \clip (0, -3pt) rectangle (3pt, 3pt);
      \fill[#4, draw = black] (0, 0) circle (2pt);
    \end{scope}
  \end{scope}
}

\def\figureFSSPWithThreeEdgesInARowAndGeneralAtTheSecondVertexFromTheLeft{%
  \begin{tikzpicture}[> = To, x = 1cm, y = {-1cm}] 
    \pgfmathsetmacro{\l}{1} 
    \pgfmathsetmacro{\b}{3} 
    \pgfmathsetmacro{\r}{2} 
    \pgfmathsetmacro{\w}{\l + \b + \r} 
    \pgfmathsetmacro{\h}{\w + \b + \r} 

    \coordinate (g) at (\l, 0); 
    \coordinate (m) at (\w/2, \h - \w/2); 

    \coordinate (mf2) at ({(\l + \b)/2}, {\b + (\l + \b)/2}); 
    \coordinate (ml2) at ({\l + (\b + \r)/2}, {1.5 * (\b + \r)}); 
    \coordinate (ml) at (\l/2, 1.5 * \l); 
    \coordinate (mr) at (\l + \b + \r/2, \b + 1.5 * \r); 
    \coordinate (mb) at (\l + \b/2, 1.5 * \b); 

    \coordinate (mf2e) at ({(\l + \b)/2}, {\h - (\l + \b)/2}); 
    \coordinate (ml2e) at ({\l + (\b + \r)/2}, {\h - (\b + \r)/2}); 
    \coordinate (mle) at (\l/2, \h - \l/2); 
    \coordinate (mre) at (\l + \b + \r/2, \h - \r/2); 
    \coordinate (mbe) at (\l + \b/2, \h - \b/2); 

    \draw[syncLeft] (g) -- ($(ml) + (\l/2, \l/2)$) -- (ml) -- ($(ml) + (-\l/2, \l/2)$) -- (0, \l) -- cycle;
    \draw[syncMiddle] (g) -- (\l + \b, \b) -- ($(mb) + (\b/2, \b/2)$) -- (mb) -- ($(mb) + (-\b/2, \b/2)$) -- cycle;
    \draw[syncRight] (\l + \b, \b) -- (\w, \b + \r) -- ($(mr) + (\r/2, \r/2)$) -- (mr) -- ($(mr) + (-\r/2, \r/2)$) -- cycle;

    \draw[find] (g) -- (0, \l);
    \draw[reflected] (0, \l) -- (\l, 2 * \l);
    \draw[slowed] (\l, 2 * \l) -- (m); 
    \draw[midpoint] (mf2) -- (mf2e); 

    \draw[find] (g) -- (\l + \b, \b);
    \draw[reflected] (\l + \b, \b) -- (mf2); 
    \draw[find] (\l + \b, \b) -- (\w, \b + \r);
    \draw[reflected] (\w, \b + \r) -- (m);
    \draw[slowed] (\l + \b, \b) -- (mr);
    \draw[midpoint] (mr) -- (mre);

    \draw[slowed] (g) -- (\l/2, 1.5 * \l);
    \draw[midpoint] (ml) -- (mle);

    \draw[slowed] (g) -- (mb) -- (ml2);
    \draw[midpoint] (mb) -- (mbe);
    \draw[midpoint] (ml2) -- (ml2e);

    \draw[freeze] (ml) -- ($(ml) + (\l/2, \l/2)$)
                  (ml) -- ($(ml) + (-\l/2, \l/2)$)
                  (mb) -- ($(mb) + (\b/2, \b/2)$)
                  (mb) -- ($(mb) + (-\b/2, \b/2)$)
                  (mr) -- ($(mr) + (\r/2, \r/2)$)
                  (mr) -- ($(mr) + (-\r/2, \r/2)$);

    \draw[traversingThaw] (m) -- (mf2e)
                          (m) -- (ml2e)
                          (mf2e) -- (mle)
                          (mf2e) -- (mbe)
                          (ml2e) -- (mbe)
                          (ml2e) -- (mre);

    \draw[syncLeft] (mle) -- ($(mle) + (\l/2, \l/2)$) -- ($(mle) + (-\l/2, \l/2)$);
    \draw[thawingThaw] (mle) -- ($(mle) + (\l/2, \l/2)$)
                       (mle) -- ($(mle) + (-\l/2, \l/2)$);
    \draw[syncMiddle] (mbe) -- ($(mbe) + (\b/2, \b/2)$) -- ($(mbe) + (-\b/2, \b/2)$);
    \draw[thawingThaw] (mbe) -- ($(mbe) + (\b/2, \b/2)$)
                       (mbe) -- ($(mbe) + (-\b/2, \b/2)$);
    \draw[syncRight] (mre) -- ($(mre) + (\r/2, \r/2)$) -- ($(mre) + (-\r/2, \r/2)$);
    \draw[thawingThaw] (mre) -- ($(mre) + (\r/2, \r/2)$)
                       (mre) -- ($(mre) + (-\r/2, \r/2)$);

    \drawMidpointAtOfThreeEdges((m),leftFill,overlayLeftFill,rightFill);
    \drawMidpointAtOf((mf2),leftFill,overlayLeftFill);
    \drawMidpointAtOf((mf2e),leftFill,overlayLeftFill);
    \drawMidpointAtOf((mr),rightFill,rightFill);
    \drawMidpointAtOf((mre),rightFill,rightFill);
    \drawMidpointAtOf((ml),leftFill,leftFill);
    \drawMidpointAtOf((mle),leftFill,leftFill);
    \drawMidpointAtOf((mb),overlayLeftFill,overlayLeftFill);
    \drawMidpointAtOf((mbe),overlayLeftFill,overlayLeftFill);
    \drawMidpointAtOf((ml2),overlayLeftFill,rightFill);
    \drawMidpointAtOf((ml2e),overlayLeftFill,rightFill);

    \draw[vertex] (0, 0) -- (0, \h)
                  (\l, 0) -- (\l, \h)
                  (\l + \b, 0) -- (\l + \b, \h)
                  (\w, 0) -- (\w, \h);
    \draw[edge, leftColour] (0, 0) -- (\l, 0)
                            (0, \h) -- (\l, \h);
    \draw[edge, overlayLeftColour] (\l, 0) -- (\l + \b, 0)
                              (\l, \h) -- (\l + \b, \h);
    \draw[edge, rightColour] (\l + \b, 0) -- (\w, 0)
                             (\l + \b, \h) -- (\w, \h);
  \end{tikzpicture}
}

\def\figureTheTreeThatIsSynchronised{%
  \begin{tikzpicture}[> = To, x = 1cm, y = {-1cm}]
    \draw[leftColour] (0, 0) -- (-0.707, 0.707);
    \draw[overlayLeftColour] (0, 0) -- (2.1213, 2.1213);
    \draw[rightColour] (2.1213, 2.1213) -- (2.1213, 4.1213);
  \end{tikzpicture}
}

\def\lengthsForFSSPWithThreeEdgesIncidentToOneVertex{%
  \pgfmathsetmacro{\l}{3} 
  \pgfmathsetmacro{\r}{2} 
  \pgfmathsetmacro{\o}{1.333} 
  \pgfmathsetmacro{\w}{\l + \r} 
}

\def\verticesAndEdgesForFSSPWithThreeEdgesIncidentToOneVertex{%
  \draw[vertex] (v) -- (\l, \h);
  \draw[vertex, leftColour] (0, 0) -- (0, \h);
  \draw[vertex, rightColour] (\w, 0) -- (\w, \h);
  \draw[vertex, overlayLeftColour] (\l - \o, 0) -- (\l - \o, \h);
  \draw[vertex, overlayRightColour] (\l + \o, 0) -- (\l + \o, \h);

  \draw[edge, leftColour] (0, 0) -- (v) 
                          (0, \h) -- (\l, \h); 
  \draw[edge, rightColour] (v) -- (\w, 0) 
                           (\l, \h) -- (\w, \h); 
  \draw[edge, overlayLeftColour] ($(v) - (0, 0.4pt)$) -- ($(v) + (-\o, 0) - (0, 0.4pt)$)
                                 ($(v) + (0, \h) - (0, 0.4pt)$) -- ($(v) + (-\o, \h) - (0, 0.4pt)$);
  \draw[edge, overlayRightColour] ($(v) - (0, 0.4pt)$) -- ($(v) + (\o, 0) - (0, 0.4pt)$)
                                  ($(v) + (0, \h) - (0, 0.4pt)$) -- ($(v) + (\o, \h) - (0, 0.4pt)$);
}

\def\figureFSSPWithThreeEdgesIncidentToTheGeneralVertex{%
  \begin{tikzpicture}[> = To, x = 1cm, y = {-1cm}] 
    \lengthsForFSSPWithThreeEdgesIncidentToOneVertex;
    \pgfmathsetmacro{\h}{\l + \l + \r} 

    \begin{scope}[shift = {(0.8, -0.8)}, scale = 0.2]
      \draw[leftColour] (0, 0) -- (-2.1213, 2.1213);
      \draw[overlayLeftColour] (0, 0) -- (0, 1.333);
      \draw[overlayRightColour] ($(0, 0) + (0.4pt, 0)$) -- ($(0, 1.333) + (0.4pt, 0)$);
      \draw[rightColour] (0, 0) -- (1.414, 1.414);
    \end{scope}
    \begin{scope}[shift = {(4.3, -0.8)}, scale = 0.2]
      \draw[leftColour] (0, 0) -- (-2.1213, 2.1213);
      \draw[overlayLeftColour] (0, 0) -- (0, -1.333);
      \draw[overlayRightColour] ($(0, 0) + (0.4pt, 0)$) -- ($(0, -1.333) + (0.4pt, 0)$);
      \draw[rightColour] (0, 0) -- (1.414, 1.414);
    \end{scope}

    \coordinate (v) at (\l, 0); 
    \coordinate (lb) at (0, \l); 
    \coordinate (rb) at (\w, \r); 
    \coordinate (obl) at (\l - \o, \o); 
    \coordinate (obr) at (\l + \o, \o); 
    \coordinate (lg) at (\l, 2 * \l); 
    \coordinate (rg) at (\l, 2 * \r); 
    \coordinate (og) at (\l, 2 * \o); 
    \coordinate (m) at (0.5 * \w, \l + 0.5 * \w); 
    \coordinate (m1) at ({0.5 * (\l + \o)}, {\l + 0.5 * (\l + \o)}); 
    \coordinate (m2) at ({\l - \o + 0.5 * (\r + \o)}, {\r + 0.5 * (\r + \o)}); 
    \coordinate (m1t) at ({0.5 * (\l + \o)}, {\h - 0.5 * (\l + \o)}); 
    \coordinate (m2t) at ({\l - \o + 0.5 * (\r + \o)}, {\h - 0.5 * (\r + \o)}); 
    \coordinate (mlf) at (0.5 * \l, \l + 0.5 * \l); 
    \coordinate (mlt) at (0.5 * \l, \h - 0.5 * \l); 
    \coordinate (mrf) at (\l + 0.5 * \r, \r + 0.5 * \r); 
    \coordinate (mrt) at (\l + 0.5 * \r, \h - 0.5 * \r); 
    \coordinate (mofl) at (\l - 0.5 * \o, \o + 0.5 * \o); 
    \coordinate (motl) at (\l - 0.5 * \o, \h - 0.5 * \o); 
    \coordinate (mofr) at (\l + 0.5 * \o, \o + 0.5 * \o); 
    \coordinate (motr) at (\l + 0.5 * \o, \h - 0.5 * \o); 
    \coordinate (lfb) at ($(mlf) + (- 0.5 * \l, 0.5 * \l)$); 
    \coordinate (rfb) at ($(mrf) + (0.5 * \r, 0.5 * \r)$); 
    \coordinate (ofbl) at ($(mofl) + (- 0.5 * \o, 0.5 * \o)$); 
    \coordinate (ofbr) at ($(mofr) + (0.5 * \o, 0.5 * \o)$); 
    \coordinate (lfg) at ($(mlf) + (0.5 * \l, 0.5 * \l)$); 
    \coordinate (rfg) at ($(mrf) + (- 0.5 * \r, 0.5 * \r)$); 
    \coordinate (ofg) at ($(mofr) + (- 0.5 * \o, 0.5 * \o)$); 

    \fill[syncLeft, leftColour, draw = none] (v) -- (lb) -- (lfb) -- (mlf) -- (lg) -- cycle;
    \fill[syncRight, rightColour, draw = none] (v) -- (rb) -- (rfb) -- (mrf) -- (rg) -- cycle;
    \fill[syncMiddle, overlayLeftColour, draw = none] (v) -- (obl) -- (ofbl) -- (mofl) -- (og) -- cycle;
    \fill[syncMiddle, overlayRightColour, draw = none] (v) -- (obr) -- (ofbr) -- (mofr) -- (og) -- cycle;

    \draw[find, leftColour] (v) -- (lb);
      \draw[reflected, leftColour] (lb) -- (m); 
    \draw[find, rightColour] (v) -- (rb);
      \draw[reflected, rightColour] (rb) -- (rg);
    \draw[find, overlayLeftColour] (v) -- (obl);
      \draw[reflected, overlayLeftColour] (obl) -- (og);
      \draw[slowed, rightColour] (og) -- (m2);
    \draw[find, overlayRightColour] (v) -- (obr);
      \draw[reflected, overlayRightColour] (obr) -- (og);
      \draw[slowed, leftColour] (og) -- (m1);

    \draw[thawingThaw, syncLeft, leftColour]
        (0, \h) -- (0.5 * \l, \h - 0.5 * \l)
                -- (\l, \h);
    \draw[thawingThaw, syncRight, rightColour]
        (\w, \h) -- (\l + 0.5 * \r, \h - 0.5 * \r)
                 -- (\w - \r, \h);
    \draw[thawingThaw, syncMiddle, overlayLeftColour]
        (\l, \h) -- (\l - 0.5 * \o, \h - 0.5 * \o)
                 -- (\l - \o, \h);
    \draw[thawingThaw, syncMiddle, overlayRightColour]
        (\l, \h) -- (\l + 0.5 * \o, \h - 0.5 * \o)
                 -- (\l + \o, \h);

    \draw[traversingThaw, leftColour] (m) -- (m1t);
      \draw[traversingThaw, leftColour] (m1t) -- (mlt);
      \draw[traversingThaw, leftColour] (m1t) -- ($(m1t) + ({\l - 0.5 * (\l + \o)}, {\l - 0.5 * (\l + \o)})$); 
          \draw[traversingThaw, overlayRightColour] ($(m1t) + ({\l - 0.5 * (\l + \o)}, {\l - 0.5 * (\l + \o)})$) -- (motr); 

    \draw[traversingThaw, leftColour] (m) -- ($(m) + ({\l - 0.5 * (\l + \r)}, {\l - 0.5 * (\l + \r)})$); 
      \draw[traversingThaw, rightColour] ($(m) + ({\l - 0.5 * (\l + \r)}, {\l - 0.5 * (\l + \r)})$) -- (m2t); 
          \draw[traversingThaw, rightColour] (m2t) -- (mrt); 
          \draw[traversingThaw, rightColour] (m2t) -- ($(m2t) + ({-(\l - \o + 0.5 * (\r + \o) - \l)}, {\l - \o + 0.5 * (\r + \o) - \l})$); 
              \draw[traversingThaw, overlayLeftColour] ($(m2t) + ({-(\l - \o + 0.5 * (\r + \o) - \l)}, {\l - \o + 0.5 * (\r + \o) - \l})$) -- (motl); 

    \draw[slowed, leftColour] (v) -- (mlf); 
    \draw[slowed, rightColour] (v) -- (mrf); 
    \draw[slowed, overlayLeftColour] ($(v) + (0.4pt, 0)$) -- ($(mofl) + (0.4pt, 0)$); 
    \draw[slowed, overlayRightColour] ($(v) - (0.4pt, 0)$) -- ($(mofr) - (0.4pt, 0)$); 
    \draw[slowed, leftColour] (\l, \r + \r) -- (0.5 * \w, {\r + \r + 1.5 * (\l - \r)});

    \draw[freeze, leftColour] (mlf) -- (lfg)
                              (mlf) -- (lfb);
    \draw[freeze, rightColour] (mrf) -- (rfb)
                               (mrf) -- (rfg);
    \draw[freeze, overlayLeftColour] (mofl) -- (ofbl)
                                     (mofl) -- (ofg);
    \draw[freeze, overlayRightColour] (mofr) -- (ofbr)
                                      (mofr) -- (ofg);

    \drawMidpointAtOf((m),leftFill,rightFill);
    \draw[midpoint, leftColour] (mlf) -- (mlt); 
        \drawMidpointAtOf((mlf),leftFill,leftFill);
        \drawMidpointAtOf((mlt),leftFill,leftFill);
    \draw[midpoint, rightColour] (mrf) -- (mrt); 
        \drawMidpointAtOf((mrf),rightFill,rightFill);
        \drawMidpointAtOf((mrt),rightFill,rightFill);
    \draw[midpoint, overlayLeftColour] (mofl) -- (motl); 
        \drawMidpointAtOf((mofl),overlayLeftFill,overlayLeftFill);
        \drawMidpointAtOf((motl),overlayLeftFill,overlayLeftFill);
    \draw[midpoint, overlayRightColour] (mofr) -- (motr); 
        \drawMidpointAtOf((mofr),overlayRightFill,overlayRightFill);
        \drawMidpointAtOf((motr),overlayRightFill,overlayRightFill);
    \draw[midpoint, leftColour] (m1) -- (m1t); 
        \drawMidpointAtOf((m1),leftFill,overlayRightFill);
        \drawMidpointAtOf((m1t),leftFill,overlayRightFill);
    \draw[midpoint, rightColour] (m2) -- (m2t); 
        \drawMidpointAtOf((m2),overlayLeftFill,rightFill);
        \drawMidpointAtOf((m2t),overlayLeftFill,rightFill);

    \verticesAndEdgesForFSSPWithThreeEdgesIncidentToOneVertex;
  \end{tikzpicture}
}

\def\figureFSSPWithThreeEdgesIncidentToTheSameVertexAndTheGeneralIsNotOnTheLongestPath{%
  \begin{tikzpicture}[> = To, x = 1cm, y = {-1cm}] 
    \lengthsForFSSPWithThreeEdgesIncidentToOneVertex
    \pgfmathsetmacro{\h}{\o + \l + \l + \r} 

    \coordinate (v) at (\l, 0); 

    \coordinate (gl) at (\l - \o, 0); 
    \coordinate (gr) at (\l + \o, 0); 
    \coordinate (lgv) at (\l, \o); 
    \coordinate (fov) at (\l, 2 * \o); 
    \coordinate (folg) at (\l - \o, 2 * \o); 
    \coordinate (forg) at (\l + \o, 2 * \o); 
    \coordinate (mol) at (\l - \o/2, 1.5 * \o); 
    \coordinate (mor) at (\l + \o/2, 1.5 * \o); 
    \path let \p1 = (mol) in coordinate (stol) at (\x1, \h - \o/2); 
    \path let \p1 = (mor) in coordinate (stor) at (\x1, \h - \o/2); 
    \coordinate (sv) at ($(mol) + (\o/2, 3 * \o/2)$); 

    \coordinate (lvl) at ($(lgv) + (-\l, \l)$); 
    \coordinate (rlv) at ($(lvl) + (\l, \l)$); 
    \coordinate (ml) at ($(lvl) + (\l/2, \l/2)$); 
    \coordinate (fll) at ($(ml) + (-\l/2, \l/2)$); 
    \coordinate (flv) at ($(ml) + (\l/2, \l/2)$); 
    \path let \p1 = (ml) in coordinate (stl) at (\x1, \h - \l/2); 

    \coordinate (lvr) at ($(lgv) + (\r, \r)$); 
    \coordinate (rfv) at ($(lvr) + (-\r, \r)$); 
    \coordinate (mr) at ($(lvr) + (-\r/2, \r/2)$); 
    \coordinate (frr) at ($(mr) + (\r/2, \r/2)$); 
    \coordinate (frv) at ($(mr) + (-\r/2, \r/2)$); 
    \path let \p1 = (mr) in coordinate (str) at (\x1, \h - \r/2); 

    \coordinate (mlo) at ({(\l + \o)/2}, {1.5 * (\l + \o)}); 
    \coordinate (mro) at ({\l - \o + (\o + \r)/2}, {1.5 * (\o + \r)}); 
    \coordinate (m) at ({(\l + \r)/2}, {\o + \l + (\l + \r)/2}); 
    \coordinate (stlo) at ({(\l + \o)/2}, {\h - (\l + \o)/2}); 
    \coordinate (stro) at ({\l - \o + (\o + \r)/2}, {\h - (\o + \r)/2}); 

    \coordinate (tv) at (\l, \h); 
    \coordinate (tr) at (\w, \h); 
    \coordinate (tl) at (0, \h); 
    \coordinate (tol) at (\l - \o, \h); 
    \coordinate (tor) at (\l + \o, \h); 
    \path let \p1 = (m) in coordinate (tmhv) at (\l, \y1 + \x1 - \l cm); 
    \path let \p1 = (stlo) in coordinate (tmlohv) at (\l, \y1 + \x1 - \l cm); 
    \path let \p1 = (stro) in coordinate (tmrohv) at (\l, \y1 - \x1 + \l cm); 

    \fill[syncLeft, leftColour, draw = none] (lgv) -- (lvl) -- (fll) -- (ml) -- (flv) -- cycle;
    \fill[syncRight, rightColour, draw = none] (lgv) -- (lvr) -- (frr) -- (mr) -- (frv) -- cycle;
    \fill[syncMiddle, overlayLeftColour, draw = none] (gl) -- (lgv) -- (fov) -- (mol) -- (folg) -- cycle;
    \fill[syncMiddle, overlayRightColour, draw = none] (gr) -- (lgv) -- (fov) -- (mor) -- (forg) -- cycle;

    \fill[syncLeft, leftColour, draw = none] (tl) -- (stl) -- (tv) -- cycle;
    \fill[syncRight, rightColour, draw = none] (tr) -- (str) -- (tv) -- cycle;
    \fill[syncMiddle, overlayLeftColour, draw = none] (tol) -- (stol) -- (tv) -- cycle;
    \fill[syncMiddle, overlayRightColour, draw = none] (tor) -- (stor) -- (tv) -- cycle;

    \draw[find, overlayLeftColour] (gl) -- (lgv); 
    \draw[find, rightColour] (lgv) -- (lvr); 
    \draw[reflected, overlayLeftColour] (lgv) -- (mol); 
    \draw[slowed, rightColour] (lgv) -- (mr); 
    \draw[reflected, rightColour] (lvr) -- ($(lvr) + (-\r, \r)$); 
        \draw[slowed, leftColour] ($(lvr) + (-\r, \r)$) -- (m); 

    \draw[find, overlayRightColour] (gr) -- (lgv); 
    \draw[find, leftColour] (lgv) -- (lvl); 
    \draw[reflected, overlayRightColour] (lgv) -- (mor); 
    \draw[slowed, leftColour] (lgv) -- (ml); 
    \draw[reflected, leftColour] (lvl) -- (m); 

    \draw[slowed, overlayLeftColour] (gl) -- (mol) -- (sv);
        \draw[slowed, rightColour] (sv) -- (mro);
    \draw[slowed, overlayRightColour] (gr) -- (mor) -- (sv);
        \draw[slowed, leftColour] (sv) -- (mlo);

    \draw[midpoint, leftColour] (ml) -- (stl)
                                (mlo) -- (stlo);
    \draw[midpoint, rightColour] (mr) -- (str)
                                 (mro) -- (stro);
    \draw[midpoint, overlayLeftColour] (mol) -- (stol);
    \draw[midpoint, overlayRightColour] (mor) -- (stor);

    \draw[freeze, leftColour] (ml) -- (fll)
                              (ml) -- (flv);
    \draw[freeze, rightColour] (mr) -- (frr)
                               (mr) -- (frv);
    \draw[freeze, overlayLeftColour] (mol) -- (fov)
                                     (mol) -- (folg);
    \draw[freeze, overlayRightColour] (mor) -- (fov)
                                      (mor) -- (forg);

    \draw[traversingThaw, leftColour] (m) -- (stlo) -- (stl);
    \draw[traversingThaw, leftColour] (m) -- (tmhv);
        \draw[traversingThaw, rightColour] (tmhv) -- (stro) -- (str);
    \draw[traversingThaw, leftColour] (stlo) -- (tmlohv);
        \draw[traversingThaw, overlayRightColour] (tmlohv) -- (stor);
    \draw[traversingThaw, rightColour] (stro) -- (tmrohv);
        \draw[traversingThaw, overlayLeftColour] (tmrohv) -- (stol);

    \draw[thawingThaw, leftColour] (stl) -- (tl)
                                   (stl) -- (tv);
    \draw[thawingThaw, rightColour] (str) -- (tr)
                                    (str) -- (tv);
    \draw[thawingThaw, overlayLeftColour] (stol) -- (tol)
                                          (stol) -- (tv);
    \draw[thawingThaw, overlayRightColour] (stor) -- (tor)
                                           (stor) -- (tv);

    \drawMidpointAtOf((ml),leftFill,leftFill);
        \drawMidpointAtOf((stl),leftFill,leftFill);
    \drawMidpointAtOf((mr),rightFill,rightFill);
        \drawMidpointAtOf((str),rightFill,rightFill);
    \drawMidpointAtOf((mol),overlayLeftFill,overlayLeftFill);
        \drawMidpointAtOf((stol),overlayLeftFill,overlayLeftFill);
    \drawMidpointAtOf((mor),overlayRightFill,overlayRightFill);
        \drawMidpointAtOf((stor),overlayRightFill,overlayRightFill);
    \drawMidpointAtOf((mlo),leftFill,overlayRightFill);
        \drawMidpointAtOf((stlo),leftFill,overlayRightFill);
    \drawMidpointAtOf((mro),overlayLeftFill,rightFill);
        \drawMidpointAtOf((stro),overlayLeftFill,rightFill);
    \drawMidpointAtOf((m),leftFill,rightFill);

    \verticesAndEdgesForFSSPWithThreeEdgesIncidentToOneVertex
  \end{tikzpicture}
}

\def\figureTheTimeAtWhichTheMidpointsOfLongestPathsAreFound{%
  \begin{tikzpicture}
    \fill (0, 1) circle (2pt) node[right] {$\general$};
    \fill (0, -1) circle (2pt) node[below] {$v''$};
    \fill (-3, -1) circle (2pt) node[below] {$\sourceOf(\hat{p})$};
    \fill (4, -1) circle (2pt) node[below] {$\targetOf(\hat{p})$};
    \draw (0, 1) -- (0, -1)
          (-3, -1) -- (0, -1) -- (4, -1);

    \fill (-2, 0.5) circle (2pt) node[left] {$v$};
    \fill (0, 0.5) circle (2pt) node[right] {$v'$};
    \draw[loosely dashed] (-2, 0.5) -- (0, 0.5);

    \fill (-2, -2) circle (2pt) node[below] {$v$};
    \fill (-1, -1) circle (2pt) node[above] {$v'$};
    \draw[thick, loosely dotted] (-2, -2) -- (-1, -1);

    \fill (3, -3) circle (2pt) node[below] {$v$};
    \fill (2, -1) circle (2pt) node[above] {$v'$};
    \draw[dashdotted] (3, -3) -- (2, -1);
  \end{tikzpicture}
}

\def\figureRemarkMidpointsOfMaximumWeightPathsAreRecognised{%
  \begin{tikzpicture}[tree/.style = {dashed, gray}]
    \coordinate (g) at (0, 0);
    \coordinate (u1) at (1, 0);
    \coordinate (u+) at (2, 0);
    \coordinate (uk) at (3, 0);
    \coordinate (hatv) at (4, 0);
    \coordinate (w1) at (4.5, -0.25);
    \coordinate (wl) at (5.25, -0.6125);
    \coordinate (hatv1) at (4 + 5, 5/2); 
    \coordinate (hatv2) at (4 + 8, -8/2); 
    \coordinate (m) at (4 + 1.47987, -1.47987/2); 

    \draw[tree] (g) -- ++(-0.8, -0.4) -- ++(0, 0.8) -- cycle;
        \draw[tree] (g) -- ++(-1, -0.7) -- ++(0, 0.9) -- cycle;
    \draw[tree] (u1) -- ++(-0.2, -1.5) -- ++(0.4, 0) -- cycle;
        \draw[tree] (u1) -- ++(-0.4, -1.7) -- ++(0.4, 0) -- cycle;
    \draw[tree] (uk) -- ++(-0.2, -2) -- ++(0.4, 0) -- cycle;
        \draw[tree] (uk) -- ++(-0.4, -2.2) -- ++(0.4, 0) -- cycle;
    \draw[tree] (hatv) -- ++(-0.2, -1.8) -- ++(0.4, 0) -- cycle;
        \draw[tree] (hatv) -- ++(-0.4, -2) -- ++(0.4, 0) -- cycle;
    \draw[tree] (w1) -- ++(-0.2, -1) -- ++(0.4, 0) -- cycle;
        \draw[tree] (w1) -- ++(-0.1, -1.1) -- ++(0.4, 0) -- cycle;
    \draw[tree] (wl) -- ++(-0.2, -2) -- ++(0.4, 0) -- cycle;
        \draw[tree] (wl) -- ++(-0.1, -2.1) -- ++(0.4, 0) -- cycle;
    \draw[tree] (hatv) -- ($(hatv1) + (-1, 1)$) -- ($(hatv1) - (-1, 1)$) -- cycle;
    \draw[tree] (hatv) -- ($(hatv2) + (-1, -1)$) -- ($(hatv2) - (-1, -1)$) -- cycle;

    \draw[dotted, thick] (g) -- (hatv);
    \draw (hatv) -- (hatv1)
          (hatv) -- (hatv2);

    \fill (g) circle (2pt) node[above] {$\general$};
    \fill (u1) circle (2pt) node[above] {$u_1$};
    \fill (u+) circle (0pt) node[above] {$\dotsb$};
    \fill (uk) circle (2pt) node[above] {$u_k$};
    \fill (hatv) circle (2pt) node[above] {$\hat{v}$};
    \fill (w1) circle (2pt) node[above right] {$w_1$};
    \fill (wl) circle (2pt) node[above right] {$w_\ell$}; 
    \fill (hatv1) circle (2pt) node[above right] {$\hat{v}_1$};
    \fill (hatv2) circle (2pt) node[below right] {$\hat{v}_2$};
    \fill[draw = black, fill = white] (m) circle (2pt) node[below] {$\hat{\midpoint}$};
  \end{tikzpicture}
}

\def\figureWhenDoesCountSignalMemoriseThePenultimateDirectionBaseCase{%
  \begin{tikzpicture}[tree/.style = {dashed, gray}]
    \coordinate (g) at (0, 0);
    \draw[tree] (g) -- ++(-0.8, -0.4) -- ++(0, 0.8) -- cycle;
        \draw[tree] (g) -- ++(-1, -0.7) -- ++(0, 0.9) -- cycle;
    \fill (g) circle (2pt) node[above] {$\general$};

    \coordinate (w) at (1.5, 0);
    \draw[tree] (w) -- ++(-0.2, -1.5) -- ++(0.4, 0) -- cycle;
        \draw[tree] (w) -- ++(-0.4, -1.7) -- ++(0.4, 0) -- cycle;
    \fill (w) circle (2pt) node[above] {$w$};

    \coordinate (v) at (5, 0);
    \fill (v) circle (2pt) node[below left] {$v$};
    \coordinate (v') at ($(v) + (60: 0.5)$);
    \fill (v') circle (2pt) node[left] {$v'$};
    \coordinate (v'') at ($(v) + (0: 1)$);
    \fill (v'') circle (2pt) node[below] {$v''$};
    \coordinate (v''') at ($(v) + (-60: 1.5)$);
    \fill (v''') circle (2pt) node[left] {$v'''$};

    \draw (g) -- (w)
          (w) -- (v)
          (v) -- (v')
          (v) -- (v'')
          (v) -- (v''');
  \end{tikzpicture}
}

\def\figureWhenDoesCountSignalMemoriseThePenultimateDirectionInductiveStep{%
  \begin{tikzpicture}[tree/.style = {dashed, gray}]
    \coordinate (g) at (0, 0);
    \draw[tree] (g) -- ++(-0.8, -0.4) -- ++(0, 0.8) -- cycle;
        \draw[tree] (g) -- ++(-1, -0.7) -- ++(0, 0.9) -- cycle;
    \fill (g) circle (2pt) node[above] {$\general$};

    \coordinate (w) at (1.5, 0);
    \draw[tree] (w) -- ++(-0.2, -1.5) -- ++(0.4, 0) -- cycle;
        \draw[tree] (w) -- ++(-0.4, -1.7) -- ++(0.4, 0) -- cycle;
    \fill (w) circle (2pt) node[above] {$w$};

    \coordinate (v) at (5, 0);
    \fill (v) circle (2pt) node[below left] {$v$};

    \coordinate (v') at ($(v) + (60: 0.5)$);
    \draw[tree] (v') -- ++(0.6, 1) -- ++(0.4, 0) -- cycle;
        \draw[tree] (v') -- ++(0.4, 1.1) -- ++(0.5, 0) -- cycle;
    \fill (v') circle (2pt) node[left] {$v'$};

    \coordinate (v'') at ($(v) + (0: 1)$);
    \draw[tree] (v'') -- ++(1, -0.3) -- ++(0, 0.8) -- cycle;
        \draw[tree] (v'') -- ++(1.2, -0.5) -- ++(0, 0.9) -- cycle;
    \fill (v'') circle (2pt) node[below] {$v''$};

    \coordinate (v''') at ($(v) + (-60: 1.5)$);
    \draw[tree] (v''') -- ++(0.2, -1) -- ++(0.4, 0) -- cycle;
        \draw[tree] (v''') -- ++(0.4, -1.1) -- ++(0.5, 0) -- cycle;
    \fill (v''') circle (2pt) node[left] {$v'''$};

    \draw (g) -- (w)
          (w) -- (v)
          (v) -- (v')
          (v) -- (v'')
          (v) -- (v''');
  \end{tikzpicture}
}

\def\figureWGeneralV{%
  \begin{tikzpicture}[tree/.style = {dashed, gray}]
    \coordinate (g) at (0, 0);
    \draw[tree] (g) -- ++(-0.8, -0.4) -- ++(0, 0.8) -- cycle;
        \draw[tree] (g) -- ++(-1, -0.7) -- ++(0, 0.9) -- cycle;
    \fill (g) circle (2pt) node[above] {$\general$};

    \coordinate (v) at (3, 0);
    \draw[tree] (v) -- ++(7, -1) -- ++(0, 1.7) -- cycle;
    \draw[tree] (v) -- ++(-0.2, -1.5) -- ++(0.4, 0) -- cycle;
        \draw[tree] (v) -- ++(-0.4, -1.7) -- ++(0.4, 0) -- cycle;
    \fill (v) circle (2pt) node[above] {$v$};

    \coordinate (w) at (2.5, 0);
    \draw[tree] (w) -- ++(-0.3, -0.7) -- ++(0.5, 0) -- cycle;
    \fill (w) circle (2pt) node[above] {$w$};

    \coordinate (w') at (2, 0);
    \draw[tree] (w') -- ++(-0.6, -1) -- ++(0.4, 0) -- cycle;
        \draw[tree] (w') -- ++(-0.4, -1.1) -- ++(0.5, 0) -- cycle;
    \fill (w') circle (2pt) node[above] {$w'$};

    \coordinate (w'') at (0.5, 0);
    \draw[tree] (w'') -- ++(-1, -2.3) -- ++(1.5, 0) -- cycle;
        \draw[tree] (w'') -- ++(-0.8, -1.1) -- ++(0.8, 0) -- cycle;
    \fill (w'') circle (2pt) node[above] {$w''$};

    \draw (g) -- (v);
  \end{tikzpicture}
}

\def\figureMarkedSignalsInnerVerticesInducitveStep{%
  \begin{tikzpicture}[tree/.style = {dashed, gray}]
    \coordinate (g) at (0, 0);
    \draw[tree] (g) -- ++(-0.8, -0.4) -- ++(0, 0.8) -- cycle;
        \draw[tree] (g) -- ++(-1, -0.7) -- ++(0, 0.9) -- cycle;
    \fill (g) circle (2pt) node[above] {$\general$};

    \coordinate (v) at (3, 0);
    \draw[tree] (v) -- ++(7, -1) -- ++(0, 1.7) -- cycle;
    \draw[tree] (v) -- ++(-0.2, -1.5) -- ++(0.4, 0) -- cycle;
        \draw[tree] (v) -- ++(-0.4, -1.7) -- ++(0.4, 0) -- cycle;
    \fill (v) circle (2pt) node[above] {$v$};

    \coordinate (v') at (2, 0);
    \draw[tree] (v') -- ++(-0.6, -1) -- ++(0.4, 0) -- cycle;
        \draw[tree] (v') -- ++(-0.4, -1.1) -- ++(0.5, 0) -- cycle;
    \fill (v') circle (2pt) node[above] {$v'$};

    \draw (g) -- (v')
          (v') -- (v);
  \end{tikzpicture}
}

\def\figureMidpointsOfMaximumWeightPathsAreRecognised{%
  \subfloat[$g = \hat{v} = \hat{\midpoint}$]{%
    \begin{tikzpicture}[tree/.style = {dashed, gray}]
      \coordinate (g) at (0, 0);
      \coordinate (v1) at (-0.5, 0);
      \coordinate (v2) at (1, 0);
      \coordinate (hatv1) at (-4, 0);
      \coordinate (hatv2) at (4, 0);

      \draw[tree] (v1) -- ($(hatv2) + (0, -1)$) -- ($(hatv2) + (0, 1)$) --node[above] {$\maxTree_{v_1}$} cycle;
      \draw[tree] (v2) -- ($(hatv1) + (0, -1.5)$) -- ($(hatv1) + (0, 1.5)$) --node[above] {$\maxTree_{v_2}$} cycle;

      \draw[dotted, thick]
          (hatv1) -- (v1)
          (hatv2) -- (v2);
      \draw (v1) -- (g) -- (v2);

      \fill (g) circle (2pt) node[above] {$\general$}; 
          \path (g) circle (2pt) node[below] {$\hat{\midpoint}$};
      \fill (v1) circle (2pt) node[above] {$v_1$};
      \fill (v2) circle (2pt) node[above] {$v_2$};
      \fill (hatv1) circle (2pt) node[left] {$\hat{v}_1$};
      \fill (hatv2) circle (2pt) node[right] {$\hat{v}_2$};
    \end{tikzpicture}
    \label{figure:midpoints-of-maximum-weight-paths-are-recognised:general-lies-on-path-and-midpoint-is-nearest-to-general}
  }

  \subfloat[$g = \hat{v} \neq \hat{\midpoint}$]{%
    \begin{tikzpicture}[tree/.style = {dashed, gray}]
      \coordinate (g) at (-2, 0);
      \coordinate (m) at (1.5, 0);
      \coordinate (v1) at (1, 0);
      \coordinate (v2) at (2.5, 0);
      \coordinate (hatv1) at (-4, 0);
      \coordinate (hatv2) at (7, 0);

      \draw[tree] (v1) -- ($(hatv2) + (0, -1)$) -- ($(hatv2) + (0, 1)$) --node[above] {$\maxTree_{v_1}$} cycle;
      \draw[tree] (v2) -- ($(hatv1) + (0, -1.5)$) -- ($(hatv1) + (0, 1.5)$) --node[above] {$\maxTree_{v_2}$} cycle;

      \draw[dotted, thick]
          (g) -- (v1)
          (hatv1) -- (g)
          (hatv2) -- (v2);
      \draw (v1) -- (v2);

      \fill (g) circle (2pt) node[above] {$\general$}; 
      \fill[draw = black, fill = white] (m) circle (2pt) node[below] {$\hat{\midpoint}$};
      \fill (v1) circle (2pt) node[above] {$v_1$};
      \fill (v2) circle (2pt) node[above] {$v_2$};
      \fill (hatv1) circle (2pt) node[left] {$\hat{v}_1$};
      \fill (hatv2) circle (2pt) node[right] {$\hat{v}_2$};
    \end{tikzpicture}
    \label{figure:midpoints-of-maximum-weight-paths-are-recognised:general-lies-on-path-and-midpoint-is-not-nearest-to-general}
  }

  \subfloat[$g \neq \hat{v} = \hat{\midpoint}$]{%
    \begin{tikzpicture}[tree/.style = {dashed, gray}]
      \coordinate (g) at (-1, 0);
      \coordinate (v1) at (0.5, 0.25);
      \coordinate (v2) at (1, -0.5);
      \coordinate (hatv) at (0, 0);
      \coordinate (hatv1) at (4, 2);
      \coordinate (hatv2) at (4, -2);

      \draw[tree] (g) -- ++(-0.8, -0.4) -- ++(0, 0.8) -- cycle;
          \draw[tree] (g) -- ++(-1, -0.7) -- ++(0, 0.9) -- cycle;

      \draw[dotted, thick]
          (g) -- (hatv)
          (hatv1) -- (v1)
          (hatv2) -- (v2);
      \draw (hatv) -- (v1)
            (hatv) -- (v2);

      \fill (g) circle (2pt) node[above] {$\general$}; 
      \fill (v1) circle (2pt) node[above] {$v_1$};
      \fill (v2) circle (2pt) node[above] {$v_2$};
      \fill (hatv) circle (2pt) node[above] {$\hat{v}$};
          \path (hatv) circle (2pt) node[below] {$\hat{\midpoint}$};
      \fill (hatv1) circle (2pt) node[right] {$\hat{v}_1$};
      \fill (hatv2) circle (2pt) node[right] {$\hat{v}_2$};
    \end{tikzpicture}
    \label{figure:midpoints-of-maximum-weight-paths-are-recognised:midpoint-is-nearest-to-general}
  }

  \subfloat[$g \neq \hat{v} \neq \hat{\midpoint}$]{%
    \begin{tikzpicture}[tree/.style = {dashed, gray}]
      \coordinate (g) at (0, 0);
      \coordinate (hatv) at (4, 0);
      \coordinate (v1) at (5.25, -0.6125);
      \coordinate (v2) at (6.5, -1.25);
      \coordinate (hatv1) at (4 + 5, 5/2); 
      \coordinate (hatv2) at (4 + 8, -8/2); 
      \coordinate (m) at (4 + 1.47987, -1.47987/2); 

      \draw[tree] (g) -- ++(-0.8, -0.4) -- ++(0, 0.8) -- cycle;
          \draw[tree] (g) -- ++(-1, -0.7) -- ++(0, 0.9) -- cycle;

      \draw[dotted, thick]
          (g) -- (hatv)
          (hatv) -- (v1)
          (hatv) -- (hatv1)
          (v2) -- (hatv2);
      \draw
          (v1) -- (m)
          (m) -- (v2);

      \fill (g) circle (2pt) node[above] {$\general$};
      \fill (hatv) circle (2pt) node[above] {$\hat{v}$};
      \fill (v1) circle (2pt) node[above] {$v_1$};
      \fill (v2) circle (2pt) node[above] {$v_2$};
      \fill (hatv1) circle (2pt) node[above right] {$\hat{v}_1$};
      \fill (hatv2) circle (2pt) node[below right] {$\hat{v}_2$};
      \fill[draw = black, fill = white] (m) circle (2pt) node[below] {$\hat{\midpoint}$};
    \end{tikzpicture}
    \label{figure:midpoints-of-maximum-weight-paths-are-recognised:midpoint-is-not-nearest-to-general}
  }
}

\begin{document}

  \frenchspacing
  \raggedbottom
  \selectlanguage{british} 

  \pagenumbering{roman}

  \pagestyle{plain}
  \begin{titlepage}
    \begin{addmargin}[-1cm]{-3cm}
      \begin{center}
        \hfill\vfill
        \begingroup
          {\huge \myTitle}\\
        \endgroup
        \bigskip\bigskip\bigskip
        {\large zur Erlangung des akademischen Grades eines}\\
        \bigskip
        {\Large Doktors der Naturwissenschaften}\\
        \bigskip\bigskip\bigskip
        der KIT-Fakultät für Informatik\\
        des Karlsruher Instituts für Technologie (KIT)\\
        \bigskip\bigskip\bigskip
        genehmigte\\
        \bigskip
        {\Large Dissertation}\\
        \bigskip\bigskip\bigskip
        von\\
        \medskip
        {\Large\myName}\\
        \medskip
        aus Heilbronn\\
        \bigskip\bigskip\bigskip\bigskip\bigskip\bigskip
        \begin{tabular}{l}
          Datum der mündlichen Prüfung: 8. Juni 2017\\
          \hfill\\
          Erster Gutachter: Prof. Dr. Jörn Müller-Quade\\
          Zweiter Gutachter: MCF Dr. Michel Coornaert
        \end{tabular}
        \vfill
      \end{center}%
    \end{addmargin}
  \end{titlepage}

  \thispagestyle{empty}

  \hfill

  \vfill

  \clearToOddPage
  \pdfbookmark[1]{Abstract}{Abstract}
  \begingroup
  \let\clearpage\relax
  \let\cleardoublepage\relax
  \let\cleardoublepage\relax

  \chapter*{Abstract}
  The present dissertation studies cellular automata whose cell spaces are not as usual the integers or groups but merely sets that are acted transitively upon by groups. These cell spaces can be continuous like the real numbers acted upon by translations or the sphere acted upon by rotations, or discrete like vertex-transitive graphs acted upon by graph automorphisms or uniform tilings of the hyperbolic plane acted upon by tiling-respecting bijections. As usual all cells have a state which they change synchronously in discrete time steps depending on the states of neighbouring cells as described by a local transition function. This global behaviour is represented by the global transition function, which maps global states of the automaton, so-called global configurations, to such.

  Among others the following questions are investigated: Which properties of the local transition function are necessary and sufficient for the global transition function to be equivariant under translations (commute with the induced action on global configurations)? Is the composition of global transition functions itself a global transition function and if it is, of which cellular automaton (see \cref{chapter:automata})? Can global transition functions be pulled or pushed onto quotients, products, restrictions, and extensions of their cell spaces and if they can, how do the corresponding cellular automata change (see \cref{chapter:quotients-and-periodicity})? Are global transition functions for a well chosen topology (or uniformity) on the set of global configurations characterised by equivariance under translations and (uniform) continuity? Is the inverse of a bijective global transition function itself a global transition function (see \cref{chapter:Curtis-Hedlund-Lyndon})? On which cell spaces are global transition functions surjective if and only if they are pre-injective (see \cref{chapter:garden})? How can such cell spaces be characterised (see \cref{chapter:amenability})? Can these questions be answered for restrictions of global transition functions to translation invariant and compact subsets of the set of global configurations (see \cref{chapter:growth,chapter:Moore})? Is there an optimal-time algorithm for the firing squad synchronisation problem on (continuous) graph-shaped cell spaces (see \cref{chapter:fssp})?

  \section*{Equivariance and Composition}
  On groups global transition functions are equivariant under translations and the composition of global transition functions is a global transition function. On general cell spaces this is only the case if the local transition function has certain symmetries. Broadly speaking, these symmetries swallow changes of orientation that can be caused by translations and they prevent that composition produces asymmetries that contradict the local and uniform inner workings of cellular automata.

  \section*{Quotients, Products, Restrictions, and Extensions}
  If the neighbourhood of a cellular automaton is included in some dimensions, then the automaton and its global transition function can be restricted to these dimensions and, analogously, both can be pulled onto the quotient of the cell space modulo the other dimensions. Conversely, a cellular automaton can be extended to a superset of its cell space and, if the cell space is a quotient, it can be pushed onto the set that underlies the quotient.

  \section*{Characterisation by Equivariance and Continuity}
  Because a global step of a cellular automaton is local and uniform in each cell, global transition functions of cellular automata on groups are characterised by equivariance under translations and (uniform) continuity with respect to a well chosen topology (or uniformity) on the set of global configurations. On general cell spaces such a characterisation is again only possible if the local transition function has certain symmetries. It follows from such a characterisation and topological theorems that the inverse of a bijective global transition function is itself a global transition function.

  \section*{Equivalence of Surjectivity and Pre-Injectivity}
  Maps on finite sets are surjective if and only if they are injective. The set of global configurations is in general infinite but, as compensation, global transition functions are local and uniform. And indeed, global transition functions of cellular automata over the integers are surjective if and only if they are pre-injective, where pre-injectivity is essentially injectivity on global configurations with finite support. On groups this equivalence is only satisfied if the group is amenable, where amenability means in a sense that the group behaves like a finite group. On general cell spaces one can also find a notion of amenability such that the before-mentioned equivalence holds.

  \section*{Restrictions to Shift Spaces}
  Sometimes it is reasonable to restrict the domain of a global transition function to certain global configurations. For example in the case that the input and the intermediate global configurations that occur in calculations performed by the automaton are of a certain form. Such sets of global configurations should for technical reasons be invariant under translations and compact and they are called shift spaces. Are such restrictions of global transition functions surjective if and only if they are pre-injective? Yes, at least if the cell space is amenable and the shift space is of finite type and strongly irreducible, which more or less means that it is generated by finitely many forbidden finite patterns and that two allowed finite patterns that are far enough apart can be embedded in a global configuration of the shift space.

  \section*{The Firing Squad Synchronisation Problem}
  A cellular automaton solves the firing squad synchronisation problem if, for each connected and finite region on which one cell is distinguished and whose other cells are in a quiescent state, after finitely many steps all cells of the region transit at the same time into a so-called fire state and this state does not occur before. In dimension one, on rectangles in dimension two, and for other special cases, there are optimal-time algorithms for this problem but not in general. In the present thesis, an optimal-time algorithm on (continuous) graph-shaped cell spaces is presented, which needs unbounded many states but uses these solely for simple geometrical constructions.


  \vfill

  \begin{otherlanguage}{ngerman}
    \pdfbookmark[1]{Kurzfassung}{Kurzfassung}
    \chapter*{Kurzfassung}
    In der vorliegenden Dissertation werden Zellularautomaten untersucht, deren Zellraum nicht wie üblich die ganzen Zahlen oder eine Gruppe ist sondern lediglich eine Menge auf der eine Gruppe transitiv operiert. Diese Zellräume können kontinuierlich sein, wie beispielsweise die reellen Zahlen mit der Operation durch Verschiebungen oder die Sphäre mit der Operation durch Drehungen, oder diskret, wie beispielsweise knoten-transitive Graphen mit der Operation durch Graphautomorphismen oder gleichförmige Parkettierungen der hyperbolischen Ebene mit der Operation durch parkettierungserhaltende Bijektionen. Wie üblich haben alle Zellen einen Zustand, den sie synchron, in diskreten Zeitschritten, in Abhängigkeit der Zustände benachbarter Zellen und wie von einer lokalen Überführungsfunktion vorgeschrieben, wechseln. Dieses globale Verhalten wird von der globalen Überführungsfunktion beschrieben, welche Gesamtzustände des Automaten, sogenannte globale Konfigurationen, auf ebensolche abbildet.

    Unter anderen werden die folgenden Fragestellungen behandelt: Unter welchen Bedingungen an die lokale Überführungsfunktion ist die globale Überführungsfunktion verschiebungsäquivariant (verträgt sich mit der induzierten Gruppenoperation auf globalen Konfigurationen)? Ist die Komposition zweier globaler Überführungsfunktionen wieder eine globale Überführungsfunktion und wenn ja, von welchem Zellularautomaten (siehe Kapitel~\ref{chapter:automata})? Können globale Überführungsfunktionen auf Quotienten, Produkte, Einschränkungen und Erweiterungen ihrer Zellräume gezogen werden und wenn ja, wie wandeln sich dabei die zugehörigen Zellularautomaten (siehe Kapitel~\ref{chapter:quotients-and-periodicity})? Sind globale Überführungsfunktionen nach geeigneter Wahl einer Topologie (oder Uniformität) auf der Menge der globalen Konfigurationen durch Verschiebungsäquivarianz und (gleichmäßige) Stetigkeit charakterisiert? Ist die Umkehrfunktion einer bijektiven globalen Überführungsfunktion wieder eine globale Überführungsfunktion (siehe Kapitel~\ref{chapter:Curtis-Hedlund-Lyndon})? Auf welchen Zellräumen sind globale Überführungsfunktionen genau dann surjektiv, wenn sie prä-injektiv sind (siehe Kapitel~\ref{chapter:garden})? Auf welche Weisen können derartige Zellräume charakterisiert werden (siehe Kapitel~\ref{chapter:amenability})? Können diese Fragen auch für Einschränkungen globaler Überführungsfunktionen auf verschiebungsinvariante und kompakte Teilmengen der Menge der globalen Konfigurationen beantwortet werden (siehe Kapitel~\ref{chapter:growth} und~\ref{chapter:Moore})? Gibt es einen Optimalzeitalgorithmus für das \enquote{Firing Squad Synchronisation Problem} auf (kontinuierlichen) graph-förmigen Zellräumen (siehe Kapitel~\ref{chapter:fssp})?

    \section*{Verschiebungsäquivarianz und Komposition}
    Auf Gruppen sind globale Überführungsfunktionen verschiebungsäquivariant und die Komposition zweier globaler Überführungsfunktionen ist eine globale Überführungsfunktion. Auf allgemeineren Zellräumen ist dies nur dann der Fall, wenn die lokalen Überführungsfunktionen gewisse Symmetrien aufweisen. Salopp ausgedrückt schlucken diese Symmetrien Orientierungsänderungen, die bei Verschiebungen auftreten können, und sie verhindern, dass bei der Komposition Asymmetrien entstehen, die der lokalen und uniformen Arbeitsweise von Zellularautomaten widersprechen.

    \section*{Quotienten, Produkte, Einschränkungen und Erweiterungen}
    Erstreckt sich die Nachbarschaft eines Zellularautomaten nur in einige Dimensionen, so ist dieser und seine globale Überführungsfunktion auf eben diese Dimensionen einschränkbar und analog können beide auf den Quotienten des Zellraums mit den anderen Dimensionen gezogen werden. Umgekehrt lässt sich ein Zellularautomat auf eine Obermenge seines Zellraums erweitern und, falls dieser ein Quotient ist, auf die dem Quotienten zugrundeliegende Menge.

    \section*{Charakterisierung durch Verschiebungsäquivarianz und Stetigkeit}
    Da ein globaler Schritt eines Zellularautomaten lokal und uniform in jeder Zelle ist, sind globale Überführungsfunktionen von Zellularautomaten auf Gruppen nach geeigneter Wahl einer Topologie (oder Uniformität) auf der Menge der globalen Konfigurationen durch Verschiebungsäquivarianz und (gleichmäßige) Stetigkeit charakterisiert. Auf allgemeineren Zellräumen ist eine solche Charakterisierung abermals nur unter zusätzlichen Bedingungen an die lokale Überführungsfunktion möglich. Insbesondere folgt aus einer solchen Charakterisierung und Sätzen der Topologie, dass die Umkehrfunktion einer bijektiven globalen Überführungsfunktion wieder eine globale Überführungsfunktion ist.

    \section*{Äquivalenz von Surjektivität und Prä-Injektivität}
    Abbildungen auf endlichen Mengen sind genau dann surjektiv, wenn sie injektiv sind. Die Menge der globalen Konfigurationen ist im Allgemeinen unendlich, aber dafür sind globale Überführungsfunktionen lokal und uniform. Tatsächlich sind globale Überführungsfunktionen von Zellularautomaten auf den ganzen Zahlen genau dann surjektiv, wenn sie prä-injektiv sind, wobei Prä-Injektivität im Wesentlichen Injektivität auf globalen Konfigurationen mit endlichem Träger ist. Auf Gruppen gilt diese Äquivalenz nur dann, wenn die Gruppe mittelbar ist, wobei Mittelbarkeit in gewissem Sinne sicherstellt, dass sich die Gruppe wie eine endliche Gruppe verhält. Auf allgemeineren Zellräumen kann man ebenfalls einen derartigen Mittelbarkeitsbegriff so finden, dass die vorhin erwähnte Äquivalenz gilt.

    \section*{Einschränkungen auf Verschiebungsräume}
    Manchmal ist es sinnvoll den Definitionsbereich einer globalen Überführungsfunktion auf bestimmte globale Konfigurationen einzuschränken, beispielsweise dann, wenn die Eingaben und die bei der Berechnung auftretenden Zwischenergebnisse nur von einer gewissen Form sind. Solche Mengen globaler Konfigurationen sollten aus technischen Gründen verschiebungsinvariant und kompakt sein und werden Verschiebungsräume genannt. Auch für solche globale Überführungsfunktionen stellt sich die Frage ob sie genau dann surjektiv sind, wenn sie prä-injektiv sind. Dies ist zumindest dann der Fall, wenn der zugrundeliegende Zellraum mittelbar ist und der Verschiebungsraum streng irreduzibel ist, was in etwa bedeutet, dass man zwei erlaubte Muster, sofern man sie weit genug voneinander entfernt, in eine globale Konfiguration des Verschiebungsraums einbetten kann.

    \section*{The Firing Squad Synchronisation Problem}
    Ein Zellularautomat löst das \enquote{Firing Squad Synchronisation Problem}, wenn für jedes zusammenhängende und endliche Gebiet auf dem eine Zelle ausgezeichnet ist und deren andere Zellen sich im Ruhezustand befinden, nach endlich vielen Schritten alle Zellen des Gebiets gleichzeitig in einen sogenannten Feuerzustand übergehen und dieser vorher nicht vorkommt. Im Eindimensionalen, auf Rechtecken im Zweidimensionalen, und für andere Spezialfälle gibt es Optimalzeitalgorithmen für dieses Problem, jedoch nicht im Allgemeinen. In dieser Arbeit wird erstmals ein Optimalzeitalgorithmus auf (kontinuierlichen) graph-förmigen Zellräumen vorgestellt, der zwar unbeschränkt viele Zustände benötigt, diese aber nur für einfache geometrische Konstruktionen verwendet.
  \end{otherlanguage}

  \endgroup

  \vfill

  \clearToOddPage
  \pdfbookmark[1]{Publications}{publications}
  \chapter*{Publications}
  Some ideas and figures have appeared previously in the following publications:

  \begin{refsection}[my]
      \small
      \nocite{*} 
      \printbibliography[heading=none]
  \end{refsection}

  \clearToOddPage
  \pdfbookmark[1]{Acknowledgement}{acknowledgement}

  \bigskip

  \begingroup
  \let\clearpage\relax
  \let\cleardoublepage\relax
  \let\cleardoublepage\relax
  \chapter*{Acknowledgement}
  Special thanks go to my unofficial supervisor Dr. Thomas Worsch, who always took the time to answer any --- simple or hard, silly or smart, dull or exciting, pleasing or provocative --- questions I had, for his excellent support and for keeping me motivated.

  \vfill
  \endgroup

  \pagestyle{scrheadings}
  \clearToOddPage
  \refstepcounter{dummy}
  \pdfbookmark[1]{\contentsname}{tableofcontents}
  \setcounter{tocdepth}{2} 
  \setcounter{secnumdepth}{3} 
  \manualmark
  \markboth{\spacedlowsmallcaps{\contentsname}}{\spacedlowsmallcaps{\contentsname}}
  \tableofcontents
  \automark[section]{chapter}
  \renewcommand{\chaptermark}[1]{\markboth{\spacedlowsmallcaps{#1}}{\spacedlowsmallcaps{#1}}}
  \renewcommand{\sectionmark}[1]{\markright{\thesection\enspace\spacedlowsmallcaps{#1}}}


  \pagenumbering{arabic}


  \clearToOddPage
  \chapter{Cellular Automata}
  \label{chapter:automata}

  \paragraph{Abstract.} We introduce cellular automata whose cell spaces are left-ho\-mo\-ge\-neous spaces; show that their global transition functions are equivariant, determined by their behaviour at one point, and closed under composition; and construct automata from others by taking quotients and products, restrictions and extensions. Examples of left-ho\-mo\-ge\-neous spaces are spheres, Euclidean spaces, as well as hyperbolic spaces acted on by isometries; uniform tilings acted on by symmetries; vertex-transitive graphs, in particular, Cayley graphs, acted on by automorphisms; groups acting on themselves by multiplication; and integer lattices acted on by translations.

  \paragraph{Remark.} Some parts of this chapter appeared in the paper \enquote{\citetitle*{wacker:automata:2016}}\cite{wacker:automata:2016} and they generalise parts of sections~1.1 to~1.5 of the monograph \enquote{\citetitle*{ceccherini-silberstein:coornaert:2010}}\cite{ceccherini-silberstein:coornaert:2010} and parts of the paper \enquote{\citetitle*{moriceau:2011}}\cite{moriceau:2011}.

  \paragraph{Motivation.} In the first chapter of the monograph \enquote{Cellular Automata and Groups}\cite{ceccherini-silberstein:coornaert:2010}, Tullio Ceccherini-Silberstein and Michel Coornaert develop the theory of cellular automata whose cell spaces are groups. Examples of groups are abound: The integer lattices and Euclidean spaces with addition (translation), the one-dimensional unit sphere embedded in the complex plane with complex multiplication (rotation), and the vertices of a Cayley graph with the group structure it encodes (graph automorphisms).

  Yet, there are many structured sets that do not admit a structure-preserving group structure. For example: Each Euclidean $n$-sphere, except for the zero-, one-, and three-dimensional, does not admit a topological group structure; and the Petersen graph does not admit an adjacency-preserving group structure on its vertices. However, these structured sets can be acted on by subgroups of their automorphism group by function application. For example Euclidean $n$-spheres can be acted on by rotations about their centres and graphs can be acted on by adjacency-preserving permutations of their vertices. 

  Moreover, there are structured groups that have more symmetries than can be expressed by the group structure. The integer lattices and the Euclidean spaces under addition, for example, are groups, but addition expresses only their translational symmetries but not their rotational and reflectional ones. Though, they can be acted on by arbitrary subgroups of their symmetry groups, like the ones generated by translations and rotations.

  The general notion that encompasses these structure-preserving actions is that of a group set, that is, a set that is acted on by a group. A group set $M$ acted on by $G$ such that for each tuple $(m, m') \in M \times M$ there is a symmetry $g \in G$ that transports $m$ to $m'$ is called \emph{left-ho\-mo\-ge\-neous space} and the action of $G$ on $M$ is said to be \emph{transitive}. In particular, groups are left-ho\-mo\-ge\-neous spaces --- they act on themselves on the left by multiplication.

  In this chapter, we develop the theory of cellular automata whose cell spaces are left-ho\-mo\-ge\-neous spaces. These cellular automata are defined so that their global transition functions are equivariant under the induced group action on global configurations. Depending on the choice of the cell space, these actions may be plain translations but also rotations and reflections. Exemplary for the first case are integer lattices that are acted on by translations; and for the second case Euclidean $n$-spheres that are acted on by rotations, but also the two-dimensional integer lattice that is acted on by the group generated by translations and the rotation by 90\degree.

  Sébastien Moriceau defines and studies a more restricted notion of cellular automata over group sets in his paper \enquote{\citetitle*{moriceau:2011}}\cite{moriceau:2011}. He requires sets of states and neighbourhoods to be finite. His automata are the global transition functions of what we call semi-cellular automata with finite set of states and finite sufficient neighbourhood. For these he proves many results that are analogous to those in the present chapter, though he uses different techniques.

  His automata obtain the next state of a cell by translating the global configuration such that the cell is moved to the origin, restricting that configuration to the neighbourhood of the origin, and applying the local transition function to that local configuration. Our automata obtain the next state of a cell by determining the neighbours of the cell, observing the states of that neighbours, and applying the local transition function to that observed local configuration. The obtained states are the same, but the viewpoints are different, which manifests itself in proofs and constructions.

  To determine the neighbourhood of a cell we let the relative neighbourhood semi-act on the right on the cell. That right semi-action is to the left group action what right multiplication is to the corresponding left group multiplication. Many properties of cellular automata are a consequence of the interplay between properties of that semi-action, translations of global configurations, and rotations of local configurations. That semi-action allows us to define the notion of right amenability for left-ho\-mo\-ge\-neous spaces (see \cref{chapter:amenability}) and to prove the Garden of Eden theorem for automata over such spaces (see \cref{chapter:garden}), which states that each such automaton is surjective if and only if it is pre-injective. For example finitely right-gen\-er\-at\-ed left-ho\-mo\-ge\-neous spaces of sub-exponential growth are right amenable, in particular, quotients of finitely generated groups of sub-exponential growth by finite subgroups acted upon by left multiplication.

  \paragraph{Contents.} In \cref{section:introduction} we motivate the definition of semi-cellular and cellular automata on left-ho\-mo\-ge\-neous spaces by a geometrical interpretation of traditional cellular automata on the two-dimensional integer lattice. In \cref{section:actions} we introduce left group actions and our prime examples, which illustrate phenomena that cannot be encountered in groups acting on themselves on the left by multiplication. In \cref{section:semi-action} we introduce coordinate systems, cell spaces, and right quotient set semi-actions that are induced by transitive left group actions and coordinate systems. In \cref{section:automata} we introduce semi-cellular and cellular automata. In \cref{section:invariances} we show that a global transition function does not depend on the choice of coordinate system, is equivariant under the induced left group action on global configurations, is determined by its behaviour in the origin, and that the composition of two global transition functions is a global transition function.

  \section{Introduction}
  \label{section:introduction}

  Informally, a traditional two-dimensional cellular automaton is a regular grid of similar finite-state machines working in synchrony whose inputs are the states of neighbouring machines. Formally, it is a quadruple made up of the set $M = \Z^2$ of cells, a finite set $Q$ of states, a finite subset $N$ of $\Z^2$ --- the (relative) neighbourhood (think of vectors) ---, and a local transition function $\delta$ from $Q^N$ to $Q$. The maps in $Q^N$ are local configurations and the maps in $Q^M$ are global configurations. The neighbourhood of a cell $m$ is $m + N$ and the local configuration that is observed by a cell $m$ in a global configuration $c$ is the map
  \begin{align*} 
    c(m + \blank)\restrictedTo_N \from N &\to Q,\\
    n &\mapsto c(m + n).
  \end{align*}
  The global transition function is the map
  \begin{align*}
    \Delta \from Q^M &\to Q^M,\\
    c &\mapsto \big[m \mapsto \delta\big(c(m + \blank)\restrictedTo_N\big)\big].
  \end{align*} 
  The state $\Delta(c)(m)$ is determined by applying the local transition function to the local configuration observed by $m$ in $c$. Because the local transition function is the same for all cells, the map $\Delta$ is uniform. And, because the local configuration that is observed by a cell is determined by the states of its finitely many neighbours, the map $\Delta$ is local.

  The translations of global configurations are the maps $t \actsOnMap \blank$, for $t \in \Z^2$, where
  \begin{align*}
    \actsOnMap \from \Z^2 \times Q^M &\to Q^M,\\
    (t, c) &\mapsto \big[m \mapsto c\big((-t) + m\big)\big].
  \end{align*}
  The global transition function is equivariant under translations (because it is uniform), which means that
  \begin{equation*}
    \ForEach t \in \Z^2 \ForEach c \in Q^M \Holds \Delta(t \actsOnMap c) = t \actsOnMap \Delta(c),
  \end{equation*}
  and it is continuous with respect to the prodiscrete topology on $Q^M$ (because it is local), which means that
  \begin{equation*}
    \ForEach O \subseteq Q^M \text{ open} \Holds \Delta^{-1}(O) \text{ is open}.
  \end{equation*}
  Conversely, each map $\Delta$ on $Q^M$ with these two properties is the global transition function of a cellular automaton. This characterisation of global transition functions is known as Curtis-Hedlund-Lyndon theorem. It follows from this theorem and basic topology that the inverse of a bijective global transition function is itself a global transition function. Moreover, because global transition functions are uniform, they are determined by their behaviour in the origin; and, because the composition of two global transition functions is uniform and local, it is itself a global transition function.

  The cells of the cellular automaton can be restricted to the group that is generated by the neighbourhood. For example, if $N$ is the set $\setOf{(-1, 0), (0, 0), (1, 0)}$, then the quadruple made up of the set $\Z \times \setOf{0}$ of cells, the set $Q$ of states, the neighbourhood $N$, and the local transition function $\delta$ is a cellular automaton over $\Z \times \setOf{0}$. Its global transition function is the map $\Delta' \from Q^{\Z \times \setOf{0}} \to Q^{\Z \times \setOf{0}}$, $c \mapsto [(m_1, 0) \mapsto \delta(c((m_1, 0) + \blank)\restrictedTo_N)]$ and the global transition function $\Delta$ of the original automaton is essentially the product $\prod_{z \in \Z} \Delta'$. Similarly, because $N$ has only the origin in common with the normal subgroup $\setOf{0} \times \Z$ of $\Z^2$, the original automaton induces one whose set of cells is the quotient group $\Z^2 \modulo (\setOf{0} \times \Z)$.

  Each subset $X$ of the phase space $Q^M$ has an entropy, which is a number that measures the complexity of $X$. More precisely, it is the asymptotic growth rate of the number of square-shaped patterns that occur in $X$ for squares about the origin of increasing sizes. Because the boundaries of large enough squares in $M$ are negligible with respect to the squares themselves, it follows from the locality of the global transition function $\Delta$ that the entropy of $\Delta(Q^M)$ is not greater than the one of $Q^M$. Using this one can show that $\Delta$ is surjective if and only if $\Delta(Q^M)$ has maximal entropy, and that $\Delta(Q^M)$ has maximal entropy if and only if $\Delta$ is pre-injective. This establishes the Garden of Eden theorem, which states that $\Delta$ is surjective if and only if it is pre-injective.

  Cellular automata are capable of synchronising all cells of a nice enough region, in the sense that a synchronisation process can be started in one cell and after some time all cells of the region go into the same predetermined state and this state did not occur before. This is achieved with a divide and conquer strategy: The region is divided in a regular fashion into smaller and smaller regions until no further division is possible, which happens for all regions at the same time and at which point all cells go into the predetermined state.

  The above definition of automata can be interpreted algebraically or geometrically. Algebraically, the cells $M$ form a group under addition and the (relative) neighbourhood $N$ is given such that the neighbourhood of the neutral element $0$ is $N$ itself. If above we replace $\Z^2$ by any group $G$, the element $0$ by its neutral element, and the addition $+$ by its operation, then we get the definition of cellular automata over $G$. Their global transition functions are still characterised by equivariance and continuity, invertible if and only if they are bijective, determined by their behaviour in the origin, closed under composition, and, if the group is amenable, then they are surjective if and only if they are pre-injective, where amenability broadly speaking means that there are subsets of $G$ whose boundaries are negligible with respect to the subsets themselves.

  Geometrically, the cells $M$ form a grid. This grid has the translational symmetries $g + \blank$, for $g \in \Z^2$. They form a group under composition. This group is isomorphic to the group $G = \Z^2$ under addition. In other words, the translation vectors in $G$ encode the translations of $M$. Hence, the group $G$ acts on $M$ on the left by translations, more formally, by $(g, m) \mapsto g + m$. Dedicate the cell $m_0 = 0$ as the origin of $M$. The (relative) neighbourhood $N$ is given such that the neighbourhood of the origin $m_0$ is $N$ itself. For each cell $m$, there is a translation vector $g_{m_0, m}$ such that $g_{m_0, m} + m_0 = m$, namely $m$; this property of the action of $G$ on $M$ is known as transitivity. The neighbourhood of $m$ is the translation of $N$ by $g_{m_0, m}$, namely $g_{m_0, m} + N$. A cell $m$ uses these translations to observe the local configuration $c(m + \blank)\restrictedTo_N$ in a global configuration $c$.
  Moreover, the action of $G$ on $M$ induces the action $\actsOnMap$ of $G$ on $Q^M$ (translations of global configurations). That interpretation suggests the following generalisation of cellular automata. 

  Let $M$ be a set, let $G$ be a group, and let $\actsOnPoint$ be a transitive left group action of $G$ on $M$. The action $\actsOnPoint$ induces a left group action $\actsOnMap$ of $G$ on $Q^M$. Moreover, let $m_0$ be an element of $M$, let $g_{m_0, m_0}$ be the neutral element of $G$, and, for each element $m \in M \smallsetminus \setOf{m_0}$, let $g_{m_0, m}$ be an element of $G$ such that $g_{m_0, m} \actsOnPoint m_0 = m$. If in the geometrical interpretation we use $M$ as the grid, the group $G$ as the translation vectors, the action $\actsOnPoint$ as the action of $G$ on $M$ by translations, namely $+$, the element $m_0$ as the origin, and, for each cell $m$, the element $g_{m_0, m}$ as the dedicated translation vector of $m$, then we get a definition of semi-cellular automata over $M$. They have the prefix \enquote{semi} because their global transition functions are in general not equivariant under $\actsOnMap$.

  For example, the right shift map on the one-dimensional grid, acted upon by the group of translations and reflections, is equivariant under translations but not under reflections. The reason is that translations leave the meanings of right and left, whereas reflections reverse them. That the right shift map is not equivariant under reflections can already be seen by the fact that its local transition function depends on the distinction between right and left, in other words, by the fact that it is not invariant under reflections. In general, to get equivariance of a global transition function under all symmetries, its local transition function must be invariant under certain symmetries. 

  The stabiliser $G_0$ of $m_0$ is the set of all group elements that fix $m_0$ (think of rotations about $m_0$). Let us assume that $N$ is invariant under $G_0$. Then, the restriction of $\actsOnPoint$ to $G_0$ induces a left group action $\bullet$ of $G_0$ on $Q^N$ (think of rotations of local configurations). A global transition function is equivariant under $\actsOnMap$ if and only if its local transition function is invariant under $\bullet$. A semi-cellular automaton with the latter property is a cellular automaton. For these automata, the global transition function does not depend on the choice of $\family{g_{m_0, m}}_{m \in M}$, the Curtis-Hedlund-Lyndon theorem holds, and the other statements made above hold also. Actually, for many properties it is sufficient that the local transition function is invariant under the restriction of $\bullet$ to $G_0 \cap H$, where $H$ is the subgroup of $G$ generated by $\setOf{g_{m_0, m} \suchThat m \in M}$.

  Note that we fix an origin $m_0$, because we need a reference cell for the (relative) neighbourhood $N$; we fix the group elements $\family{g_{m_0, m}}_{m \in M}$, because we need them to define the neighbourhood of each cell $m$ as $g_{m_0, m} \actsOnPoint N$; we choose $g_{m_0, m_0}$ as the neutral element, because we want the neighbour of the origin that corresponds to the relative neighbour $n$ to be $n$ itself; we require the left group action $\actsOnPoint$ to be transitive, because otherwise there would be a cell $m$ for which there would be no group element $g$ such that $g \actsOnPoint m_0 = m$ and hence we could not fix a group element $g_{m_0, m}$.
  Fixing $m_0$ and $\family{g_{m_0, m}}_{m \in M}$ is necessary for the definition of global transition functions by local transition functions. However, it is not necessary for the characterisation of global transition functions by equivariance and continuity. 

  We could regard $N$ as the neighbourhood of the origin (a set of points) and $g_{m_0, m} \actsOnPoint N$ as translating the neighbourhood of the origin to $m$. However, we regard $N$ as a relative neighbourhood (a set of vectors) and $m \isSemiActedUponBy N = g_{m_0, m} \actsOnPoint N$ as adding the relative neighbourhood on the right to $m$ (think of point-vector addition). Because we have chosen an origin $m_0$, there is a canonical bijection between $M$ and the quotient set $G \modulo G_0$. Under the identification of $M$ with $G \modulo G_0$, the relative neighbourhood $N$ is a subset of $G \modulo G_0$. And, the map $\isSemiActedUponBy \from M \times G \modulo G_0 \to M$, $(m, \mathfrak{g}) \mapsto g_{m_0, m} \actsOnPoint \mathfrak{g}$, is a right semi-action that semi-commutes with $\actsOnPoint$. We could have defined $\isSemiActedUponBy$ as a map with domain $M \times M$, but to state its properties it is convenient to have the domain $M \times G \modulo G_0$. In the definition of $\isSemiActedUponBy$ we identified $M$ with $G \modulo G_0$ to regard $\mathfrak{g} \in G \modulo G_0$ as an element of $M$. There is an equivalent definition that works without this identification. 

  The transporter $G_{m, m'}$ is the set of all group elements that transport $m$ to $m'$. The quotient set $G \modulo G_0$ is the set of all transporters from $m_0$. Under the identification of $M$ with $G \modulo G_0$, a (relative) neighbour $n$ is the transporter $G_{m_0, n}$. And, under the identification of singleton sets with their only element, we have $m \isSemiActedUponBy n = g_{m_0, m} \actsOnPoint n = G_{m, g_{m_0, m} \actsOnPoint n} \actsOnPoint m = g_{m_0, m} n g_{m_0, m}^{-1} \actsOnPoint m$. In particular, because $g_{m_0, m_0}$ is the neutral element, we have $m_0 \isSemiActedUponBy n = n \actsOnPoint m_0$. So, we can think of $n$ as a localised vector with initial cell $m_0$, of conjugating $n$ with $g_{m_0, m}$ as changing the initial cell to $m$, of $\actsOnPoint$ as a means to turn localised vectors into a cell by adding it to its initial cell, and of $\isSemiActedUponBy$ as a means to add a localised vector with initial cell $m_0$ to any cell. And, we can define $\isSemiActedUponBy$ by $m \isSemiActedUponBy n = g_{m_0, m} n g_{m_0, m}^{-1} \actsOnPoint m$.

  The right semi-action $\isSemiActedUponBy$ can be used to define the notion of right amenability for the left action $\actsOnPoint$. Right amenability is characterised by the existence of invariant finitely additive probability measures, the existence of invariant means, the existence of right Følner nets, and the non-existence of right paradoxical decompositions. This characterisation is known as Tarski-Følner theorem. In the case that the left action $\actsOnPoint$ is right amenable, the Garden of Eden theorem holds.

  Some properties of cellular automata still hold if the set of states and the neighbourhood are infinite. And, on continuous spaces, infinite and compact neighbourhoods are more natural than finite ones. Therefore, unless stated otherwise, we drop the usual finiteness requirements. However, in the uniform variant of the Curtis-Hedlund-Lyndon theorem, we require the neighbourhood to be compact (a generalisation of finiteness for continuous spaces); and in the topological variant, we require the set of states and the neighbourhood to be finite.

  \section{Left Group Actions}
  \label{section:actions}

  \paragraph{Introduction.} The symmetries of a geometric object are distance-preserving bijections on the object. The identity map is the trivial symmetry that maps each point to itself. The composition of two symmetries is again a symmetry. The symmetries under composition form a group. One says that this group (or some subgroup) \emph{acts} on the geometric object by mapping points to points. 
  On each point the identity map acts trivially and the composition of two symmetries acts in the same way as the right symmetry does followed by the left. 

  The symmetries of a circle are the reflections about lines through its centre and the rotations about it. Note that the reflection about the centre of the circle is not missing, because it is identical to the rotation by $180\degree$. For each pair of points on the circle, there is a symmetry, even a rotation, that maps one point to the other. One says that the symmetries, even only the rotations, act \emph{transitively} on the circle. 
  Moreover, for each pair of points on the circle, there is at most, even exactly, one rotation that maps one point to the other. One says that the rotations act \emph{freely} on the circle. 
  However, the symmetries do not act freely on the circle, because for each pair of points, in addition to the rotation that maps one point to the other, the reflection about the line through the centre of the circle and the centres of the circular arcs connecting the two points does so too.

  For each point on the circle, the points it can be mapped to by a symmetry, called \emph{orbit of the point}, is the circle; the symmetries that map the point to itself, called \emph{stabiliser of the point}, are the identity map and the reflection about the line through the centre of the circle and the point; and the symmetries that map a point to another point, called \emph{transporter from the point to the other one}, are the reflection about the line through the centre of the circle and the centres of the circular arcs connecting the points, and the rotation by the central angle between the points with the correct sign. 

  The symmetries of a square are the four reflections (about the horizontal line through the centre of the square, the vertical line through the centre of the square, and the two diagonals) and the four rotations (by $0\degree$, $90\degree$, $180\degree$, and $270\degree$). For each pair of vertices of the square, there is a symmetry, even a rotation, that maps one vertex to the other. However, there is no symmetry that maps a vertex to a point on an edge and vice versa. Thus, the symmetries, even only the rotations, act transitively on the vertices; whereas the symmetries do not act transitively on the square itself. 
  Therefore, it is often appropriate to regard \enquote{pointy} geometrical objects as graphs, that is, as their vertices equipped with structural information induced by their edges. Their symmetries are graph automorphisms, that is, bijections on the vertices that preserve adjacency. 

  \paragraph{Contents.} In \cref{definition:left-group-set} we introduce left group actions and left group sets. In \cref{example:group:action,example:sphere:action,example:plane:action} we introduce three examples of left group sets that we use to illustrate new notions throughout the present chapter. In \cref{remark:group-actions-as-symmetries} we show in which sense left group actions act by symmetries. In \cref{definition:left-subgroup-action} we introduce restrictions of left group actions to subgroups. In \cref{definition:transitive-free} we introduce faithfulness, freeness, transitivity, and regularity of left group actions. In \cref{definition:homogeneous-space,definition:principal-homogeneous-space} we introduce left-ho\-mo\-ge\-neous spaces and principal ones. In \cref{definition:orbit-stabiliser-transporter} we introduce orbits, stabilisers, and transporters. In \cref{lemma:stabiliser-versus-transporter} we show how stabilisers and transporters of two elements with the same orbit relate to each other. In \cref{definition:quotient-set,remark:quotient-set-and-transporters} we introduce quotient sets and note that the quotient set by the stabiliser of a point is the set of transporters from that point. In \cref{definition:orbit-space,remark:orbit-space-is-partition} we introduce orbit spaces and note that they partition the set of points. In \cref{definition:invariant-map,definition:invariant-set} we introduce invariance of maps and sets under left group actions. In \cref{definition:equivariant-tuple,definition:equivariant-map} we introduce equivariance of maps from one group set to another; such maps are homomorphisms. In \cref{lemma:inverse-of-tuple-is-equivariant,corollary:inverse-is-equivariant} we show that the inverse of an equivariant and bijective map is again equivariant. In \cref{lemma:induced-left-group-action-on-quotient-set} we show that a group acts transitively on each of its quotient sets on the left by multiplication. In \cref{lemma:iota-is-bijective-and-equivariant} we show that each transitive left group action is isomorphic to an action as in \cref{lemma:induced-left-group-action-on-quotient-set}. In \cref{corollary:iota-is-bijective-and-equivariant} we show that each free and transitive left group action is isomorphic to a group multiplication, which we elaborate on in \cref{example:principal-homogeneous-space}. In \cref{example:normal:homogeneous-space} we note that each left action of an abelian group induces a faithful, free, and transitive left group action. And in \cref{lemma:induced-left-action-on-functions} we show that left group actions induce left group actions on maps over invariant subsets of points.

  \begin{definition}
  \label{definition:left-group-set}
    Let $M$ be a set, let $G$ be a group, let $\actsOnPoint$ be a map from $G \times M$ to $M$, and let $e_G$ be the neutral element of $G$. The map $\actsOnPoint$ is called \define{left group action of $G$ on $M$}\graffito{left group action $\actsOnPoint$ of $G$ on $M$}\index{group action of $G$ on $M$!left}\index[symbols]{arrow right@$\actsOnPoint$}, the group $G$ is said to \define{act on $M$ on the left by $\actsOnPoint$}\graffito{$G$ acts on $M$ on the left by $\actsOnPoint$}, and the triple $\ntuple{M, G, \actsOnPoint}$ is called \define{left group set}\graffito{left group set $\ntuple{M, G, \actsOnPoint}$}\index{group set!left}\index[symbols]{Mcalligraphic@$\mathcal{M}$} if and only if
    \begin{equation*}
      \ForEach m \in M \Holds e_G \actsOnPoint m = m,
    \end{equation*}
    and
    \begin{equation*}
      \ForEach m \in M \ForEach g \in G \ForEach g' \in G \Holds g g' \actsOnPoint m = g \actsOnPoint (g' \actsOnPoint m). \qedhere
    \end{equation*}
  \end{definition}

  \begin{example}[Group]
  \label{example:group:action}
    Let $G$ be a group. It acts on itself on the left by multiplication.
  \end{example}

  \begin{example}[Plane]
  \label{example:plane:action}
    Let $M$ be the Euclidean plane $\R^2$ and let $G$ be the special Euclidean group $\EuclideanGroup^+(2)$, that is, the group generated by translations and rotations of $M$. The group $G$ acts on $M$ on the left by function application. This action is denoted by $\actsOnPoint$.
  \end{example}

  \begin{example}[Sphere]
  \label{example:sphere:action}
    Let $M$ be the Euclidean unit $2$-sphere, that is, the surface of the ball of radius $1$ in three-dimensional Euclidean space, and let $G$ be the rotation group. The group $G$ acts on $M$ on the left by function application, that is, by rotation about the centre of the sphere. This action is denoted by $\actsOnPoint$.
  \end{example}

  \begin{definition}
    Let $M$, $M'$, and $M''$ be three sets, let $f$ be a map from $M \times M'$ to $M''$, and let $m_0$ and $m_0'$ be two elements of $M$ and $M'$ respectively. The maps
    \begin{align*}
      M &\to M'', \mathnote{partially applied map $f(\blank, m_0')$}\index[symbols]{funderscorem0prime@$f(\blank, m_0')$}\\
      m &\mapsto f(m, m_0'),
    \end{align*}
    and
    \begin{align*}
      M' &\to M'', \mathnote{partially applied map $f(m_0, \blank)$}\index[symbols]{fm0underscore@$f(m_0, \blank)$}\\
      m' &\mapsto f(m_0, m'),
    \end{align*}
    are denoted by $f(\blank, m_0')$ and $f(m_0, \blank)$ respectively and called \define{partially applied maps}\index[symbols]{underscore@$\blank$}. We will also use analogous notations for maps with more than two arguments and for maps that are written in infix notation.
  \end{definition}

  \begin{remark}
  \label{remark:group-actions-as-symmetries}
    Let $M$ be a set, let $G$ be a group, and let $\symmetricGroupOf(M)$ be the symmetric group of $M$. For each left group action $\actsOnPoint$ of $G$ on $M$, the map
    \begin{align*}
      f \from G &\to \symmetricGroupOf(M),\\
      g &\mapsto g \actsOnPoint \blank,
    \end{align*}
    is a group homomorphism. And, for each group homomorphism $f$ from $G$ to $\symmetricGroupOf(M)$, the map
    \begin{align*}
      \actsOnPoint \from G \times M &\to M,\\
      (g, m) &\mapsto f(g)(m),
    \end{align*}
    is a left group action.
  \end{remark}


  \begin{definition}
    Let $M$ and $M'$ be two sets, let $N$ and $N'$ be two subsets of $M$ and $M'$ respectively, and let $f$ be a map from $M$ to $M'$ such that $f(N) \subseteq N'$.
    \begin{aenumerate}
      \item The map
            \begin{align*}
              f\restrictedTo_{N \to N'} \from N &\to N', \mathnote{restriction $f\restrictedTo_{N \to N'}$ of $f$ to $N$ and $N'$}\index[symbols]{harpoon up right N to N prime@$\restrictedTo_{N \to N'}$}\\
              n &\mapsto f(n),
            \end{align*}
            is called \define{restriction of $f$ to $N$ and $N'$}.
      \item The map
            \begin{align*}
              f\restrictedTo_N \from N &\to M', \mathnote{domain restriction $f\restrictedTo_N$ of $f$ to $N$}\index[symbols]{harpoon up right N@$\restrictedTo_N$}\\
              n &\mapsto f(n),
            \end{align*}
            is called \define{domain restriction of $f$ to $N$}. \qedhere
    \end{aenumerate}%
  \end{definition}

  \begin{definition}
  \label{definition:left-subgroup-action}
    Let $\actsOnPoint$ be a left group action of $G$ on $M$ and let $H$ be a subgroup of $G$. The left group action $\actsOnPoint\restrictedTo_{H \times M}$ of $H$ on $M$ is denoted by $\actsOnPoint_H$.
  \end{definition}

  \begin{definition}
  \label{definition:transitive-free}
    Let $\actsOnPoint$ be a left group action of $G$ on $M$. It is called
    \begin{aenumerate}
      \item \define{faithful}\graffito{faithful} if and only if
            \begin{equation*}
              \ForEach g \in G \ForEach g' \in G \smallsetminus \setOf{g} \Exists m \in M \SuchThat g \actsOnPoint m \neq g' \actsOnPoint m;
            \end{equation*}
      \item\label{item:transitive-free:free}
            \define{free}\graffito{free} if and only if
            \begin{equation*}
              \ForEach g \in G \ForEach g' \in G \Holds (\Exists m \in M \SuchThat g \actsOnPoint m = g' \actsOnPoint m) \implies g = g';
            \end{equation*}
      \item\label{item:transitive-free:transitive}
            \define{transitive}\graffito{transitive} if and only if the set $M$ is non-empty and
            \begin{equation*}
              \ForEach m \in M \ForEach m' \in M \Exists g \in G \SuchThat g \actsOnPoint m = m';
            \end{equation*}
      \item \define{regular}\graffito{regular} if and only if it is both free and transitive. \qedhere
    \end{aenumerate}
  \end{definition}

  \begin{example}[Group]
  \label{example:group:transitive-and-free}
    In the situation of \cref{example:group:action}, the left group action is faithful, transitive, and free.
  \end{example}

  \begin{example}[Plane]
  \label{example:plane:transitive}
    In the situation of \cref{example:plane:action}, the left group action is faithful and transitive but not free.
  \end{example}

  \begin{example}[Sphere]
  \label{example:sphere:transitive}
    In the situation of \cref{example:sphere:action}, the left group action is faithful and transitive but not free.
  \end{example}

  \begin{definition}
    Let $\actsOnPoint$ be a left group action of $G$ on $M$ and let $P$ be an adjective. The group $G$ is said to \define{act $P$ly on $M$ on the left by $\actsOnPoint$}\graffito{act $P$ly on $M$ on the left by $\actsOnPoint$} if and only if the action $\actsOnPoint$ is $P$.
  \end{definition}

  %

  \begin{definition}
  \label{definition:homogeneous-space}
    Let $\ntuple{M, G, \actsOnPoint}$ be a left group set. It is called \graffito{left-ho\-mo\-ge\-neous space}\define{left-ho\-mo\-ge\-neous space}\index{homogeneous space!left} if and only if the action $\actsOnPoint$ is transitive.
  \end{definition}

  \begin{definition}
  \label{definition:principal-homogeneous-space}
    Let $\ntuple{M, G, \actsOnPoint}$ be a left-ho\-mo\-ge\-neous space. It is called \define{principal}\graffito{principal}\index{left-ho\-mo\-ge\-neous space!principal}\index{homogeneous space!left!principal} if and only if the action $\actsOnPoint$ is free.
  \end{definition}

  \begin{definition}
  \label{definition:orbit-stabiliser-transporter}
    Let $\actsOnPoint$ be a left group action of $G$ on $M$, and let $m$ and $m'$ be two elements of $M$.
    \begin{aenumerate}
      \item The set
            \begin{equation*}
              G \actsOnPoint m = \setOf{g \actsOnPoint m \suchThat g \in G}
              \mathnote{orbit $G \actsOnPoint m$ of $m$ under $\actsOnPoint$}
              \index[symbols]{Garrowrightm@$G \actsOnPoint m$}
            \end{equation*}
            is called \define{orbit of $m$ under $\actsOnPoint$}.
      \item The set
            \begin{equation*}
              G_m = \setOf{g \in G \suchThat g \actsOnPoint m = m}
              \mathnote{stabiliser $G_m$ of $m$ under $\actsOnPoint$}
              \index[symbols]{Gm@$G_m$}
            \end{equation*}
            is called \define{stabiliser of $m$ under $\actsOnPoint$}.
      \item The set
            \begin{equation*}
              G_{m, m'} = \setOf{g \in G \suchThat g \actsOnPoint m = m'}
              \mathnote{transporter $G_{m, m'}$ of $m$ to $m'$ under $\actsOnPoint$}
              \index[symbols]{Gmmprime@$G_{m, m'}$}
            \end{equation*}
            is called \define{transporter of $m$ to $m'$ under $\actsOnPoint$}. \qedhere
    \end{aenumerate}
  \end{definition}

  \begin{example}[Group]
  \label{example:group:orbit-stabiliser}
    In the situation of \cref{example:group:transitive-and-free}, each orbit is $G$ and each stabiliser is $\setOf{e_G}$. 
  \end{example}

  \begin{example}[Plane]
  \label{example:plane:orbit-stabiliser}
    In the situation of \cref{example:plane:transitive}, for each point $m \in M$, its orbit is $M$ and its stabiliser is the group of rotations about itself.
  \end{example}

  \begin{example}[Sphere]
  \label{example:sphere:orbit-stabiliser}
    In the situation of \cref{example:sphere:transitive}, for each point $m \in M$, its orbit is $M$ and its stabiliser is the group of rotations about the line through the centre and itself.
  \end{example}

  \begin{example}[Riemannian Symmetric Space]
  \label{example:riemannian-symmetric-space}
    Let $(M, \innerProductOf{\blank}{\blank})$ be a Riemannian manifold, let $\isometriesOf(M)$ be the isometry group of $M$, and let the geodesic reflection at any point of $M$ be an isometry. Then, $M$ is geodesically complete and $\isometriesOf(M)$ acts transitively on $M$ by function application (see theorems~1.4 and~1.3 of \cite{nikonorov:2007}). Examples of such Riemannian manifolds include:
    \begin{aenumerate}
      \item The Euclidean space $\R^d$ of dimension $d \in \N_+$ with the Euclidean metric. Its isometry group $\isometriesOf(M)$ is the Euclidean group $\EuclideanGroup(d)$ generated by translations and orthogonal linear maps. The stabiliser of the origin $0$ is the orthogonal group $\orthogonalGroup(d)$.
      \item The sphere $\sphere^d = \setOf{v \in \R^{d+1} \suchThat \normOf{v}_2 = 1}$ of dimension $d \in \N_+$ with the Riemannian metric induced by the dot product $\cdot$. Its isometry group $\isometriesOf(M)$ is the orthogonal group $\orthogonalGroup(d + 1)$. The stabiliser of the north pole $\ntuple{0, 0, \dotsc, 0, 1}^\transposed$ is $\orthogonalGroup(d) \subseteq \orthogonalGroup(d + 1)$.
      \item The real hyperbolic space $\hyperbolicSpace^d = \setOf{v \in \R^{d+1} \suchThat \innerProductOf{v}{v} = -1, v_{d+1} > 0}$ of dimension $d \in \N_+$ with the Riemannian metric induced by the Lorentzian indefinite inner product
            \begin{align*}
              \innerProductOf{\blank}{\blank} \from \R^{d+1} \times \R^{d+1} &\to \R,\\
              (v,v') &\mapsto \sum_{i = 1}^d v_i v'_i - v_{i+1} v'_{i+1}.
            \end{align*}
            Its isometry group $\isometriesOf(M)$ is the group of \enquote{future preserving} Lorentz transformations $\orthogonalGroup(d,1)^+$. The stabiliser of the north pole $(0,0,\dotsc,0,1)^\transposed$ is $\orthogonalGroup(d) \subseteq \orthogonalGroup(d, 1)^+$.
      \item The orthogonal group $\orthogonalGroup(d) = \setOf{A \in \R^{d \times d} \suchThat A^\transposed A = I}$ in dimension $d \in \N_+$ with the Riemannian metric induced by the trace inner product
            \begin{align*}
              \innerProductOf{\blank}{\blank} \from \R^{d \times d} \times \R^{d \times d} &\to \R,\\
              (A, A') &\mapsto \traceOf A^\transposed A'.
            \end{align*} 
      \item Each compact Lie group $M$ with bi-invariant Riemannian metric $\innerProductOf{\blank}{\blank}$. \qedhere 
    \end{aenumerate}
  \end{example}

  \begin{lemma}
  \label{lemma:stabiliser-versus-transporter}
    Let $\actsOnPoint$ be a left group action of $G$ on $M$, let $m$ and $m'$ be two elements of $M$ that have the same orbit under $\actsOnPoint$, and let $g$ be an element of $G_{m, m'}$. Then, $G_{m'} = g G_m g^{-1}$ and $g G_m = G_{m, m'} = G_{m'} g$.
  \end{lemma}

  \begin{proof}
    First, let $g' \in G_m$. Then,
    \begin{equation*}
      g g' g^{-1} \actsOnPoint m'
      = g g' \actsOnPoint m
      = g \actsOnPoint m
      = m'.
    \end{equation*}
    Hence, $g g' g^{-1} \in G_{m'}$. In conclusion, $g G_m g^{-1} \subseteq G_{m'}$. Secondly, let $g' \in G_{m'}$. Then, as above, $g'' = g^{-1} g' g \in G_m$. Hence, $g' = g g'' g^{-1} \in g G_m g^{-1}$. In conclusion, $G_{m'} \subseteq g G_m g^{-1}$. To sum up, $G_{m'} = g G_m g^{-1}$. Moreover, $g G_m = G_{m'} g$.

    Thirdly, let $g' \in G_m$. Then, $g g' \in G_{m, m'}$. In conclusion, $g G_m \subseteq G_{m, m'}$. Lastly, let $g' \in G_{m, m'}$. Then, $g'' = g^{-1} g' \in G_m$. Hence, $g' = g g'' \in g G_m$. In conclusion, $G_{m, m'} \subseteq g G_m$. To sum up, $g G_m = G_{m, m'}$.
  \end{proof}

  \begin{definition}
  \label{definition:quotient-set}
    Let $G$ be a group and let $H$ be a subgroup of $G$. The set
    \begin{equation*}
      G \modulo H = \setOf{g H \suchThat g \in G}
      \mathnote{quotient set $G \modulo H$ of $G$ by $H$}
      \index[symbols]{GmoduloH@$G \modulo H$}
    \end{equation*}
    is called \define{quotient set of $G$ by $H$}.
  \end{definition}

  \begin{remark}
  \label{remark:quotient-set-and-transporters}
    Let $\actsOnPoint$ be a left group action of $G$ on $M$ and let $m$ be an element of $M$. The quotient set $G \modulo G_m$ is equal to $\setOf{G_{m, m'} \suchThat m' \in G \actsOnPoint m}$.
  \end{remark}

  \begin{definition}
  \label{definition:orbit-space}
    Let $\actsOnPoint$ be a left group action of $G$ on $M$. The set\index[symbols]{G modulo-reverse M@$G \reverseModulo M$} 
    \begin{equation*}
      G \reverseModulo M = \setOf{G \actsOnPoint m \suchThat m \in M}
      \mathnote{orbit space $G \reverseModulo M$ of $\actsOnPoint$}
      \index[symbols]{GmoduloMback@$G \reverseModulo M$}
    \end{equation*}
    is called \define{orbit space of $\actsOnPoint$}.
  \end{definition}

  \begin{definition}
    Let $M$ be a set and let $\mathfrak{N}$ be a subset of the power set of $M$. The set $\mathfrak{N}$ is called
    \begin{aenumerate}
      \item \defineX{pairwise disjoint}{pairwise disjoint!set}\graffito{pairwise disjoint} if and only if
            \begin{equation*}
              \ForEach N \in \mathfrak{N} \ForEach N' \in \mathfrak{N} \Holds (N \neq N' \implies N \cap N' = \emptyset);
            \end{equation*}
      \item \defineX{cover of $M$}{cover of $M$!set}\graffito{cover of $M$} if and only if $\bigcup_{N \in \mathfrak{N}} N = M$;
      \item \defineX{partition of $M$}{partition of $M$!set}\graffito{partition of $M$} if and only if it is pairwise disjoint and a cover of $M$. \qedhere
    \end{aenumerate}
  \end{definition}

  \begin{remark}
  \label{remark:orbit-space-is-partition}
    The orbit space of $\actsOnPoint$ is a partition of $M$.
  \end{remark}

  \begin{definition}
    Let $M$ be a set, let $N$ be a subset of $M$, and let $f$ be a map from $M$ to $M$. The set $N$ is called \define{invariant under $f$}\graffito{set invariant under $f$} if and only if $f(N) \subseteq N$.
  \end{definition}

  \begin{definition}
  \label{definition:invariant-map}
    Let $M$ and $M'$ be two sets, let $f$ be a map from $M$ to $M'$, and let $\actsOnPoint$ be a left group action of $G$ on $M$. The map $f$ is called \defineX{$\actsOnPoint$-invariant}{invariant map@$\actsOnPoint$-invariant map}\graffito{$\actsOnPoint$-invariant map} if and only if
    \begin{equation*}
      \ForEach g \in G \ForEach m \in M \Holds f(g \actsOnPoint m) = f(m). \qedhere
    \end{equation*}
  \end{definition}

  \begin{definition}
  \label{definition:invariant-set}
    Let $\actsOnPoint$ be a left group action of $G$ on $M$ and let $N$ be a subset of $M$. The set $N$ is called \defineX{$\actsOnPoint$-invariant}{invariant set@$\actsOnPoint$-invariant set}\graffito{$\actsOnPoint$-invariant set} if and only if $G \actsOnPoint N \subseteq N$.
  \end{definition}

  \begin{definition}
  \label{definition:equivariant-tuple}
    Let $M$ and $M'$ be two sets, let $f$ be a map from $M$ to $M'$, let $G$ and $G'$ be two groups, let $\varphi$ be a group homomorphism from $G$ to $G'$, and let $\actsOnPoint$ and $\actsOnPoint'$ be two left group actions of $G$ on $M$ and of $G'$ on $M'$ respectively. The tuple $(f, \varphi)$ is called
    \begin{aenumerate} 
      \item \defineX{$(\actsOnPoint, \actsOnPoint')$-e\-qui\-var\-i\-ant}{equivariant tuple 10@$(\actsOnPoint, \actsOnPoint')$-e\-qui\-var\-i\-ant tuple}\graffito{$(\actsOnPoint, \actsOnPoint')$-e\-qui\-var\-i\-ant tuple} if and only if 
            \begin{equation*}
              \ForEach g \in G \ForEach m \in M \Holds f(g \actsOnPoint m) = \varphi(g) \actsOnPoint' f(m);
            \end{equation*}
      \item \defineX{$\actsOnPoint$-e\-qui\-var\-i\-ant}{equivariant tuple 20@$\actsOnPoint$-e\-qui\-var\-i\-ant tuple}\graffito{$\actsOnPoint$-e\-qui\-var\-i\-ant tuple} if and only if it is $(\actsOnPoint, \actsOnPoint')$-e\-qui\-var\-i\-ant, $M = M'$, and $\actsOnPoint = \actsOnPoint'$. \qedhere
    \end{aenumerate}
  \end{definition}

  %
  %

  \begin{definition} 
  \label{definition:equivariant-map}
    Let $M$ and $M'$ be two sets, let $f$ be a map from $M$ to $M'$, let $G$ be a group, and let $\actsOnPoint$ and $\actsOnPoint'$ be two left group actions of $G$ on $M$ and $M'$ respectively. The map $f$ is called
    \begin{aenumerate}
      \item \defineX{$(\actsOnPoint, \actsOnPoint')$-e\-qui\-var\-i\-ant}{equivariant map 1@$(\actsOnPoint, \actsOnPoint')$-e\-qui\-var\-i\-ant map}\graffito{$(\actsOnPoint, \actsOnPoint')$-e\-qui\-var\-i\-ant map} if and only if the tuple $(f, \identityMap_G)$ is $(\actsOnPoint, \actsOnPoint')$-e\-qui\-var\-i\-ant;
      \item \defineX{$\actsOnPoint$-e\-qui\-var\-i\-ant}{equivariant map 2@$\actsOnPoint$-e\-qui\-var\-i\-ant map}\graffito{$\actsOnPoint$-e\-qui\-var\-i\-ant map} if and only if it is $(\actsOnPoint, \actsOnPoint')$-e\-qui\-var\-i\-ant, $M = M'$, and $\actsOnPoint = \actsOnPoint'$. \qedhere
    \end{aenumerate}
  \end{definition}

  %

  \begin{lemma}
  \label{lemma:inverse-of-tuple-is-equivariant}
    In the situation of \cref{definition:equivariant-tuple}, let $f$ and $\varphi$ be bijective, and let $(f, \varphi)$ be $(\actsOnPoint, \actsOnPoint')$-e\-qui\-var\-i\-ant. The tuple $(f^{-1}, \varphi^{-1})$ is $(\actsOnPoint', \actsOnPoint)$-e\-qui\-var\-i\-ant.
  \end{lemma}

  \begin{proof}
    For each $g' \in G'$ and each $m' \in M'$,
    \begin{align*}
      f^{-1}(g' \actsOnPoint' m')
      &= f^{-1}(\varphi(\varphi^{-1}(g')) \actsOnPoint' f(f^{-1}(m')))\\
      &= f^{-1}(f(\varphi^{-1}(g') \actsOnPoint f^{-1}(m')))\\
      &= \varphi^{-1}(g') \actsOnPoint f^{-1}(m'). \qedhere
    \end{align*}
  \end{proof}

  \begin{corollary}
  \label{corollary:inverse-is-equivariant}
    In the situation of \cref{definition:equivariant-map}, let $f$ be bijective and $(\actsOnPoint, \actsOnPoint')$-e\-qui\-var\-i\-ant. The inverse $f^{-1}$ is $(\actsOnPoint', \actsOnPoint)$-e\-qui\-var\-i\-ant.
  \end{corollary}

  \begin{proof}
    This is a direct consequence of \cref{lemma:inverse-of-tuple-is-equivariant}.
  \end{proof}


  \begin{lemma}
  \label{lemma:induced-left-group-action-on-quotient-set}
    Let $G$ be a group and let $H$ be a subgroup of $G$. The group $G$ acts transitively on the quotient set $G \modulo H$ on the left by
    \begin{align*}
      \cdot \from G \times G \modulo H &\to G \modulo H, \mathnote{left group action $\cdot$ of $G$ on $G \modulo H$}\index[symbols]{dotcentre@$\cdot$}\index[symbols]{centredot@$\cdot$}\\
      (g, g' H) &\mapsto g g' H. \qedhere
    \end{align*}
  \end{lemma}

  \begin{proof}
    The map $\cdot$ is well-defined, because, for each $g \in G$, each $g_1' \in G$, and each $g_2' \in G$,
    \begin{equation*}
      g_1' H = g_2' H \ifAndOnlyIf g g_1' H = g g_2' H.
    \end{equation*}
    It is a left group action, because, for each $g' H \in G \modulo H$, 
    \begin{equation*}
      e_G \cdot g' H = g' H,
    \end{equation*}
    and, for each $g_1 \in G$, each $g_2 \in G$, and each $g' H \in G \modulo H$,
    \begin{align*}
      g_1 g_2 \cdot g' H
      &= g_1 g_2 g' H\\
      &= g_1 \cdot g_2 g' H\\
      &= g_1 \cdot (g_2 \cdot g' H).
    \end{align*}
    It is transitive, because, for each $g_1' H \in G \modulo H$ and each $g_2' H \in G \modulo H$,
    \begin{equation*}
      g_2' {g_1'}^{-1} \cdot g_1' H = g_2' H. \qedhere
    \end{equation*}
  \end{proof}

  After choosing an element $m_0$, a left-ho\-mo\-ge\-neous space is isomorphic to $\ntuple{G \modulo G_{m_0}, G, \cdot}$, which is shown in

  \begin{lemma}
  \label{lemma:iota-is-bijective-and-equivariant}
    Let $\actsOnPoint$ be a transitive left group action of $G$ on $M$, let $m_0$ be an element of $M$, and let $G_0$ be the stabiliser of $m_0$ under $\actsOnPoint$. The map
    \begin{align*}
      \iota \from M &\to G \modulo G_0, \mathnote{$(\actsOnPoint, \cdot)$-e\-qui\-var\-i\-ant bijection $\iota$ from $M$ to $G \modulo G_0$}\index[symbols]{iota@$\iota$}\\
      m &\mapsto G_{m_0, m},
    \end{align*}
    is $(\actsOnPoint, \cdot)$-e\-qui\-var\-i\-ant and bijective. 
  \end{lemma}

  \begin{proof}
    For each $g \in G$ and each $m \in M$,
    \begin{equation*}
      \iota(g \actsOnPoint m)
      = G_{m_0, g \actsOnPoint m}
      = g \cdot G_{m_0, m}
      = g \cdot \iota(m).
    \end{equation*}
    Hence, $\iota$ is $(\actsOnPoint, \cdot)$-e\-qui\-var\-i\-ant. Moreover, for each $(m, m') \in M \times M$ with $m \neq m'$, we have $G_{m_0, m} \neq G_{m_0, m'}$. Thus, $\iota$ is injective. Furthermore, for each $g G_0 \in G \modulo G_0$, we have $G_{m_0, g \actsOnPoint m_0} = g G_0$. Therefore, $\iota$ is surjective. 
  \end{proof}

  After choosing an element $m_0$, a principal left-ho\-mo\-ge\-neous space is isomorphic to $G$, which is shown in

  \begin{corollary}
  \label{corollary:iota-is-bijective-and-equivariant}
    Let $\actsOnPoint$ be a free and transitive left group action of $G$ on $M$ and let $m_0$ be an element of $M$. The map
    \begin{align*}
      \iota \from M &\to G,\\
      m &\mapsto g_{m_0, m}, \text{ where $g_{m_0, m} \in G_{m_0, m}$},
    \end{align*}
    is $(\actsOnPoint, \cdot)$-e\-qui\-var\-i\-ant and bijective.
  \end{corollary}

  \begin{proof}
    Because $\cardinalityOf{G_{m_0, m}} = 1$ and $G_{m_0} = \setOf{e_G}$, the map $\iota$ is well-defined, the quotient set $G \modulo G_{m_0}$ can be identified with $G$, and the statement is a direct consequence of \cref{lemma:iota-is-bijective-and-equivariant}.
  \end{proof}

  \begin{example}[Principal]
  \label{example:principal-homogeneous-space}
    Let $\mathcal{M} = \ntuple{M, G, \actsOnPoint}$ be a principal left-ho\-mo\-ge\-neous space. Then, for each element $m \in M$ and each element $m' \in M$, there is one and only one element $g_{m, m'} \in G$ such that $g_{m, m'} \actsOnPoint m = m'$, in particular, because $e_G \actsOnPoint m = m$, we have $g_{m, m} = e_G$, and, because $g_{m, m'}^{-1} \actsOnPoint m' = m$, we have $g_{m', m} = g_{m, m'}^{-1}$.

    Let $m_0$ be an element of $M$ and let $M$ be equipped with the group multiplication
    \begin{align*}
      M \times M &\to M,\\
      (m, m') &\mapsto g_{m_0, m} g_{m_0, m'} g_{m_0, m}^{-1} \actsOnPoint m \quad (= g_{m_0, m} g_{m_0, m'} \actsOnPoint m_0).
    \end{align*}
    Then, the element $m_0$ is the neutral element and, for each element $m \in M$, the element $g_{m, m_0} \actsOnPoint m_0$ is the inverse element of $m$. Moreover, the maps
    \begin{align*}
      \left\{
        \begin{aligned}
          \iota \from M &\to G,\\
          m &\mapsto g_{m_0, m},
        \end{aligned}
      \right\}
      \text{ and }
      \left\{
        \begin{aligned}
          \iota^{-1} \from G &\to M,\\
          g &\mapsto g \actsOnPoint m_0,
        \end{aligned}
      \right\}
    \end{align*}
    are group isomorphisms that are inverse to each other. Under the identification of $M$ with $G$ by either isomorphism, which depends on the arbitrary choice of $m_0$, the left group action $\actsOnPoint$ is the group multiplication of $M$ and of $G$. In the words of John Baez: \enquote{A torsor [principal left-ho\-mo\-ge\-neous space] is like a group that has forgotten its identity}\cite{baez:2009}.
  \end{example}

  \begin{example}[Normal]
  \label{example:normal:homogeneous-space}
    Let $\mathcal{M} = \ntuple{M, G, \actsOnPoint}$ be a left-ho\-mo\-ge\-neous space whose stabilisers are normal subgroups of $G$, which is for example the case if the group $G$ is abelian. Then, because the stabilisers are conjugate to each other, they are all equal and we denote them by $G_0$. The map
    \begin{align*}
      \cosetActsOnPoint \from G \modulo G_0 \times M &\to M,\\
      (g G_0, m) &\mapsto g \actsOnPoint m,
    \end{align*}
    is a faithful, free, and transitive left group action of $G \modulo G_0$ on $M$ and the triple $\mathcal{M}' = \ntuple{M, G \modulo G_0, \cosetActsOnPoint}$ is a principal left-ho\-mo\-ge\-neous space. Moreover, the quotient group $G \modulo G_0$ is isomorphic to $G$ if and only if the stabiliser $G_0$ is trivial, which is the case if and only if the action $\actsOnPoint$ is free.
  \end{example}

  \begin{lemma}
  \label{lemma:induced-left-action-on-functions}
    Let $\actsOnPoint$ be a left group action of $G$ on $M$, let $H$ be a subgroup of $G$, let $N$ be a subset of $M$ such that $H \actsOnPoint N \subseteq N$, and let $Q$ be a set. The group $H$ acts on $Q^N$ on the left by
    \begin{align*}
      \actsOnMap \from H \times Q^N &\to Q^N,\\
      (h, f) &\mapsto [n \mapsto f(h^{-1} \actsOnPoint n)]. \qedhere
    \end{align*}
  \end{lemma}

  \begin{proof}
    The map $\actsOnMap$ is well-defined, because $H \actsOnPoint N \subseteq N$. It is a left group action, because, for each $f \in Q^N$ and each $n \in N$,
    \begin{equation*}
      (e_H \actsOnMap f)(n) = f(e_H \actsOnPoint n) = f(n),
    \end{equation*}
    and, for each $h \in H$, each $h' \in H$, each $f \in Q^N$, and each $n \in N$,
    \begin{align*}
      (h h' \actsOnMap f)(n)
      &= f((h')^{-1} h^{-1} \actsOnPoint n)\\
      &= f((h')^{-1} \actsOnPoint (h^{-1} \actsOnPoint n))\\
      &= (h' \actsOnMap f)(h^{-1} \actsOnPoint n)\\
      &= (h \actsOnMap (h' \actsOnMap f))(n). \qedhere
    \end{align*}
  \end{proof}

  \section{Right Quotient Set Semi-Actions}
  \label{section:semi-action}

  \paragraph{Introduction.} The rotations of a circle act on it on the left by function application, because a rotation by $x$ degrees, $\rho_x$, followed by a rotation by $y$ degrees, $\rho_y$, is the same as the rotation by $y + x$ degrees, $\rho_y \after \rho_x$, symbolically, $\rho_y \actsOnPoint (\rho_x \actsOnPoint \blank) = (\rho_y \after \rho_x) \actsOnPoint \blank$. They also act on the circle on the right by function application, because a rotation by $x$ degrees, $\rho_x$, followed by a rotation by $y$ degrees, $\rho_y$, is the same as the rotation by $x + y$ degrees, $\rho_x \after \rho_y$, symbolically, $(\blank \isActedUponBy \rho_x) \isActedUponBy \rho_y = \blank \isActedUponBy (\rho_x \after \rho_y)$. These actions commute with each other, because a rotation by $x$ degrees followed by a rotation by $y$ degrees is the same as a rotation by $y$ degrees followed by a rotation by $x$ degrees, symbolically, $(\rho_x \actsOnPoint \blank) \isActedUponBy \rho_y = \rho_x \actsOnPoint (\blank \isActedUponBy \rho_y)$. More succinctly, rotations act on the right by function application, and the left and right actions commute, because rotations commute under composition. 

  The symmetries of a circle act on it on the left by function application. However, they do not act on it on the right by function application, because, for example, the rotation by $90\degree$ followed by the reflection about the vertical line $v$ through the centre of the circle is not the same as the reflection about $v$ followed by the rotation by $90\degree$, symbolically, $(\blank \isActedUponBy \rho_{90\degree}) \isActedUponBy \varrho_v \neq \blank \isActedUponBy (\rho_{90\degree} \after \varrho_v)$. In a sense, the problem is that reflections treat different points differently: Some points stay put, others are reflected to points close by, and still others to points far away.

  To solve this, let us fix a point on the circle and call it \emph{origin} (beware, do not mistake this point for the origin of the space the circle may be embedded in). We want to define the right group semi-action such that a symmetry acts on each point as it does on the origin. For example, under the right semi-action, if a symmetry stabilises the origin, then it shall stabilise each point; and, if a symmetry throws the origin to its opposite point, then so it shall do with each point. So, a symmetry semi-acts on the right on a point by first rotating the point to the origin, secondly acting with the symmetry on the left, and lastly undoing the first rotation, symbolically, $m \isActedUponBy \sigma = (\rho_m \after \sigma \after \rho_m^{-1}) \actsOnPoint m$, where $\rho_m$ denotes the rotation that rotates the origin to $m$. Note that $m_0 \isActedUponBy \sigma = \sigma \actsOnPoint m_0$. And, that this semi-action agrees with the right group action of the rotations on the circle. In particular, it is transitive.

  The identity map semi-acts trivially on each point on the right, symbolically, $m \isActedUponBy \identityMap = m$. However, in general, the composition of two symmetries semi-acts in a different way on the right than the first symmetry does followed by the second, symbolically, $m \isActedUponBy (\sigma \after \varsigma) \neq (m \isActedUponBy \sigma) \isActedUponBy \varsigma$. Yet, it can be seen that the difference is little in the sense that there is a symmetry $\varsigma_0$ that stabilises the origin and may depend on $m$ and $\sigma$ such that $m \isActedUponBy (\sigma \after \varsigma) = (m \isActedUponBy \sigma) \isActedUponBy (\varsigma_0 \after \varsigma)$. Because of this property, the map $\isActedUponBy$ is a semi-action. Note that only the identity map and the reflection about the line $\ell$ through the centre of the circle and the origin stabilise the origin; and that $\varsigma_0$ is the identity map, if $\sigma$ and $\varsigma$ are both rotations or both reflections, and the reflection about $\ell$, otherwise.

  The right semi-action semi-commutes with the left action, which means that first acting on the left and then semi-acting on the right is almost the same as first semi-acting on the right and then acting on the left, where the defect is again a symmetry that stabilises the origin. Symbolically, $(\sigma \actsOnPoint m) \isActedUponBy \varsigma = \sigma \actsOnPoint (m \isActedUponBy (\varsigma_0 \after \varsigma))$, where $\varsigma_0$ stabilises the origin and may depend on $m$ and $\sigma$.

  For a point $m$ on the circle, there are two symmetries that map the origin to $m$, the rotation $\rho_m$ and the (roto-)reflection $\rho_m \after \varrho_{\ell}$, where $\varrho_{\ell}$ is the reflection that stabilises the origin. Because these two symmetries semi-act the same way on each point on the right, the semi-action $\isActedUponBy$ is not free. The elements of the quotient set of the symmetries of the circle by the stabiliser of the origin, namely $\setOf{\identityMap, \varrho_{\ell}}$, are the sets $\setOf{\rho_m, \rho_m \after \varrho_{\ell}}$ for points $m$. This quotient set semi-acts on the circle by $m \isSemiActedUponBy \setOf{\rho_m, \rho_m \after \varrho_{\ell}} = m \isActedUponBy \rho_m$. So, it acts in the same way as $\isActedUponBy$ but is free, which means that, if $m \isSemiActedUponBy \Sigma = m \isSemiActedUponBy T$, then $\Sigma = T$.

  Under the identification of the quotient set, which is even a quotient group, with the rotations, the right quotient set semi-action $\isSemiActedUponBy$ is identical to the right group action $\isActedUponBy$ of the rotations on the circle we considered at the beginning of this introduction. However, while the symmetry group of the circle has this nice subgroup, namely the rotation group, that acts freely and transitively on it on the right, the symmetry groups of other geometrical objects do not have such subgroups. Nevertheless, we can construct a free and transitive right quotient set semi-action on these geometrical objects as we did for the circle.

  \paragraph{Contents.} In \cref{definition:coordinate-system} we introduce coordinate systems for left-ho\-mo\-ge\-neous spaces as tuples made up of an origin and, for each point, a group element (think of a coordinate) that transports the origin to that point. In \cref{definition:cell-space} we introduce cell spaces as left-ho\-mo\-ge\-neous spaces equipped with coordinate systems. In \cref{definition:big-subgroup} we introduce bigness of subgroups with respect to a coordinate system as containing all coordinates. In \cref{lemma:liberation} we introduce right quotient set semi-actions induced by cell spaces of the quotient set of the group by the stabiliser of the origin on the points, which is to the left group action what right multiplication is to the corresponding left group multiplication. In \cref{lemma:liberation:transitive-and-free} we show that semi-actions are free and transitive. In \cref{lemma:semi-commutativity-of-liberation} we show that semi-actions semi-commute with their corresponding left group action and exhaust their defect with respect to this semi-commutativity in the origin. And in \cref{lemma:identification-of-G-quotient-Gzero-with-M-by-right-semiaction} we show that under the identification of the quotient set with the points, left group actions on quotient sets by multiplication and right quotient set semi-actions can be expressed in terms of left group actions on points.

  \begin{definition}
    Let $M$ be a set, let $I$ be a set, and let $f$ be a map from $I$ to $M$. The map $f$ is called \graffito{family $\family{m_i}_{i \in I}$ of elements in $M$ indexed by $I$}\define{family of elements in $M$ indexed by $I$} and denoted by $\family{m_i}_{i \in I}$\index[symbols]{miiinI@$\family{m_i}_{i \in I}$}, where $m_i = f(i)$, for $i \in I$.
  \end{definition}

  \begin{definition} 
  \label{definition:coordinate-system}
    Let $\mathcal{M} = \ntuple{M, G, \actsOnPoint}$ be a left-ho\-mo\-ge\-neous space, let $m_0$ be an element of $M$, let $g_{m_0, m_0}$ be the neutral element of $G$, and, for each element $m \in M \smallsetminus \setOf{m_0}$, let $g_{m_0, m}$ be an element of $G$ such that $g_{m_0, m} \actsOnPoint m_0 = m$. The tuple $\mathcal{K} = \ntuple{m_0, \family{g_{m_0, m}}_{m \in M}}$ is called \define{coordinate system for $\mathcal{M}$}\graffito{coordinate system $\mathcal{K}$ for $\mathcal{M}$}\index[symbols]{Kcalligraphic@$\mathcal{K}$}; the element $m_0$ is called \define{origin}\graffito{origin $m_0$}\index[symbols]{m0@$m_0$}; for each element $m \in M$, the element $g_{m_0, m}$ is called \define{coordinate of $m$}\graffito{coordinate $g_{m_0, m}$ of $m$}\index[symbols]{gm0m@$g_{m_0, m}$}; for each subgroup $H$ of $G$, the stabiliser of the origin $m_0$ under $\actsOnPoint_H$, which is $G_{m_0} \cap H$, is denoted by $H_0$\graffito{stabiliser $H_0$ of $m_0$ under $\actsOnPoint_H$}\index[symbols]{H0@$H_0$}, in particular, the stabiliser $G_{m_0}$ is denoted by $G_0$\graffito{stabiliser $G_0$ of $m_0$ under $\actsOnPoint$}\index[symbols]{G0@$G_0$}.
  \end{definition}

  \begin{definition}
  \label{definition:cell-space}
    Let $\mathcal{M} = \ntuple{M, G, \actsOnPoint}$ be a left-ho\-mo\-ge\-neous space and let $\mathcal{K} = \ntuple{m_0, \family{g_{m_0, m}}_{m \in M}}$ be a coordinate system for $\mathcal{M}$. The tuple $\mathcal{R} = \ntuple{\mathcal{M}, \mathcal{K}}$ is called \define{cell space}\graffito{cell space $\mathcal{R}$}\index[symbols]{Rcalligraphic@$\mathcal{R}$}, each element $m \in M$ is called \define{cell}\graffito{cell $m$}\index[symbols]{m@$m$}, and each element $g \in G$ is called \define{symmetry}\graffito{symmetry $g$}\index[symbols]{g@$g$}.
  \end{definition}

  \begin{example}[Group]
  \label{example:group:cell-space}
    In the situation of \cref{example:group:orbit-stabiliser}, let $m_0$ be the neutral element $e_G$ of $G$ and, for each element $m \in G$, let $g_{m_0, m}$ be the only element in $G$ such that $g_{m_0, m} m_0 = m$, namely $m$. The tuple $\mathcal{K} = \ntuple{m_0, \family{g_{m_0, m}}_{m \in G}}$ is a coordinate system for $\mathcal{M} = \ntuple{G, G, \cdot}$ and the tuple $\mathcal{R} = \ntuple{\mathcal{M}, \mathcal{K}}$ is a cell space.
  \end{example}

  \begin{example}[Plane]
  \label{example:plane:cell-space}
    In the situation of \cref{example:plane:orbit-stabiliser}, let $m_0$ be the origin $(0,0)^\transposed$ of $M$ and, for each point $m \in M$, let $g_{m_0, m}$ be the translation $\blank + m$ that translates $m_0$ to $m$. Note that $g_{m_0, m_0}$ is the identity map. The tuple $\mathcal{K} = \ntuple{m_0, \family{g_{m_0, m}}_{m \in M}}$ is a coordinate system for $\mathcal{M} = \ntuple{M, G, \actsOnPoint}$ and the tuple $\mathcal{R} = \ntuple{\mathcal{M}, \mathcal{K}}$ is a cell space.
  \end{example}

  \begin{example}[Sphere]
  \label{example:sphere:cell-space}
    In the situation of \cref{example:sphere:orbit-stabiliser}, let $m_0$ be the north pole $(0,0,1)^\transposed$ of $M$ and, for each point $m \in M$, let $g_{m_0, m}$ be a rotation about an axis in the $(x, y)$-plane that rotates $m_0$ to $m$, which is unique unless $m$ is the south pole. Note that $g_{m_0, m_0}$ is the identity map. The tuple $\mathcal{K} = \ntuple{m_0, \family{g_{m_0, m}}_{m \in M}}$ is a coordinate system for $\mathcal{M} = \ntuple{M, G, \actsOnPoint}$ and the tuple $\mathcal{R} = \ntuple{\mathcal{M}, \mathcal{K}}$ is a cell space.
  \end{example}

  \begin{definition}
  \label{definition:big-subgroup}
    Let $\ntuple{\mathcal{M}, \mathcal{K}} = \ntuple{\ntuple{M, G, \actsOnPoint}, \ntuple{m_0, \family{g_{m_0, m}}_{m \in M}}}$ be a cell space and let $H$ be a subgroup of $G$. The group $H$ is called \defineX{$\mathcal{K}$-big}{big@$\mathcal{K}$-big}\graffito{$\mathcal{K}$-big} if and only if
    \begin{equation*}
      \ForEach m \in M \Holds g_{m_0, m} \in H. \qedhere
    \end{equation*}
  \end{definition}

  \begin{example}[Group]
  \label{example:group:big}
    In the situation of \cref{example:group:cell-space}, because the set $\setOf{g_{m_0, m} \suchThat m \in G}$ is $G$, the only $\mathcal{K}$-big subgroup of $G$ is the group $G$.
  \end{example}

  \begin{example}[Plane]
  \label{example:plane:big}
    In the situation of \cref{example:plane:cell-space}, because the set $\setOf{g_{m_0, m} \suchThat m \in M}$ is the set $T$ of translations, the subgroup $T$ of $G$ is $\mathcal{K}$-big; and, because $G$ is the inner semi-direct product of the rotations $R_0$ about $m_0$ acting on $T$, each inner semi-direct product of a subgroup of $R_0$ acting on $T$ is a $\mathcal{K}$-big subgroup of $G$. 
  \end{example}

  \begin{example}[Sphere]
  \label{example:sphere:big}
    In the situation of \cref{example:sphere:cell-space}, because the set $\setOf{g_{m_0, m} \suchThat m \in M}$ generates $G$, the only $\mathcal{K}$-big subgroup of $G$ is the group $G$.
  \end{example}

  \begin{remark}
    The terms \enquote{coordinate system} and \enquote{big} are due to \cite{moriceau:2011}.
  \end{remark}

  \begin{remark}
  \label{remark:smallest-big-subgroup}
    The subgroup of $G$ that is generated by the set $\setOf{g_{m_0, m} \suchThat m \in M}$ of coordinates, is the smallest $\mathcal{K}$-big subgroup of $G$, where \emph{smallest} means that it is included in each $\mathcal{K}$-big subgroup of $G$.
  \end{remark}

  In the remainder of this section, let $\mathcal{R} = \ntuple{\mathcal{M}, \mathcal{K}} = \ntuple{\ntuple{M, G, \actsOnPoint}, \ntuple{m_0, \family{g_{m_0, m}}_{m \in M}}}$ be a cell space.

  \begin{lemma}
  \label{lemma:liberation}
    The map
    \begin{align*}
      \isSemiActedUponBy \from M \times G \modulo G_0 &\to M, \mathnote{right quotient set semi-action $\isSemiActedUponBy$ of $G \modulo G_0$ on $M$ with defect $G_0$}\index[symbols]{arrowleftunderscore@$\isSemiActedUponBy$}\\
      (m, g G_0) &\mapsto g_{m_0, m} g g_{m_0, m}^{-1} \actsOnPoint m \quad (= g_{m_0, m} g \actsOnPoint m_0),
    \end{align*}
    is a \define{right quotient set semi-action of $G \modulo G_0$ on $M$ with defect $G_0$}, which means that, for each $\mathcal{K}$-big subgroup $H$ of $G$,
    \begin{equation*}
      \ForEach m \in M \Holds m \isSemiActedUponBy G_0 = m,
    \end{equation*}
    and
    \begin{multline*}
      \ForEach m \in M \ForEach h \in H \Exists h_0 \in H_0 \SuchThat \ForEach \mathfrak{g} \in G \modulo G_0 \Holds\\
            m \isSemiActedUponBy h \cdot \mathfrak{g} = (m \isSemiActedUponBy h G_0) \isSemiActedUponBy h_0 \cdot \mathfrak{g}. \qedhere
    \end{multline*} 
  \end{lemma}

  \begin{proof} 
    For each $m \in M$,
    \begin{equation*}
      m \isSemiActedUponBy G_0
      = m \isSemiActedUponBy e_G G_0
      = g_{m_0, m} e_G \actsOnPoint m_0
      = g_{m_0, m} \actsOnPoint m_0
      = m.
    \end{equation*}
    Let $H$ be a $\mathcal{K}$-big subgroup of $G$. Furthermore, let $m \in M$ and let $h \in H$. Put $h_0 = g_{m_0, g_{m_0, m} h \actsOnPoint m_0}^{-1} g_{m_0, m} h$. Because $g_{m_0, g_{m_0, m} h \actsOnPoint m_0}^{-1} \in H$, $g_{m_0, m} \in H$, and
    \begin{align*}
      h_0 \actsOnPoint m_0
      &= g_{m_0, g_{m_0, m} h \actsOnPoint m_0}^{-1} \actsOnPoint (g_{m_0, m} h \actsOnPoint m_0)\\
      &= m_0,
    \end{align*}
    we have $h_0 \in H_0$. Moreover, for each $g G_0 \in G \modulo G_0$,
    \begin{align*}
      m \isSemiActedUponBy h \cdot g G_0
      &= m \isSemiActedUponBy h g G_0\\
      &= g_{m_0, m} h g \actsOnPoint m_0\\
      &= g_{m_0, g_{m_0, m} h \actsOnPoint m_0} h_0 g \actsOnPoint m_0\\
      &= (g_{m_0, m} h \actsOnPoint m_0) \isSemiActedUponBy h_0 g G_0\\
      &= (m \isSemiActedUponBy h G_0) \isSemiActedUponBy h_0 \cdot g G_0. \qedhere
    \end{align*}
  \end{proof}

  \begin{remark}
    The second property of the right quotient set semi-action $\isSemiActedUponBy$ is equivalent to the following: For each $\mathcal{K}$-big subgroup $H$ of $G$,
    \begin{multline*}
       \ForEach m \in M \ForEach h \in H \Exists h_0 \in H_0 \SuchThat \ForEach \mathfrak{g} \in G \modulo G_0 \Holds\\
             (m \isSemiActedUponBy h G_0) \isSemiActedUponBy \mathfrak{g} = m \isSemiActedUponBy h h_0 \cdot \mathfrak{g}. \qedhere
     \end{multline*}
  \end{remark}

  \begin{example}[Principal]
  \label{example:principal-homogeneous-space:liberation}
    In the situation of \cref{example:principal-homogeneous-space}, the tuple $\mathcal{K} = \ntuple{m_0, \family{g_{m_0, m}}_{m \in M}}$ is a coordinate system for $\mathcal{M}$, the stabiliser of $m_0$ under $\actsOnPoint$ and defect of $\isSemiActedUponBy$ is the trivial subgroup of $G$, and, under the identification of $G \modulo G_0$ with $G$ and of $G$ with $M$ as in \cref{example:principal-homogeneous-space}, the induced semi-action $\isSemiActedUponBy$ is the group multiplication on $M$ from \cref{example:principal-homogeneous-space} (note the similarity of their definitions).
  \end{example}

  \begin{example}[Group]
  \label{example:group:liberation}
    In the situation of \cref{example:group:big}, the stabiliser $G_0$ of the neutral element $m_0$ under $\cdot$ is the trivial subgroup $\setOf{e_G}$ of $G$ and, for each element $m \in G$ and each element $g \in G$, we have $m \isSemiActedUponBy g G_0 = g_{m_0, m} g g_{m_0, m}^{-1} m = m g m^{-1} m = m g$. Under the natural identification of $G \modulo G_0$ with $G$, the induced semi-action $\isSemiActedUponBy$ is the right group action of $G$ on itself by right multiplication.
  \end{example}

  \begin{example}[Plane]
  \label{example:plane:liberation} 
    In the situation of \cref{example:plane:big}, the stabiliser $G_0$ of the origin $m_0$ under $\actsOnPoint$ is the group of rotations about $m_0$. The special Euclidean group $G$ is the inner semi-direct product of $G_0$ acting on the abelian group $T$ of translations. Under the identification of $G \modulo G_0$ with $T$ by $t G_0 \mapsto t$, the induced semi-action $\isSemiActedUponBy$ is the right group action of $T$ on $M$ by function application. 
  \end{example}

  \begin{example}[Sphere]
  \label{example:sphere:liberation}
    In the situation of \cref{example:sphere:big}, the stabiliser $G_0$ of the north pole $m_0$ under $\actsOnPoint$ is the group of rotations about the $z$-axis. An element $g G_0 \in G \modulo G_0$ semi-acts on a point $m$ on the right by the induced semi-action $\isSemiActedUponBy$ by first rotating $m$ to $m_0$, $g_{m_0, m}^{-1} \actsOnPoint m = m_0$, secondly rotating $m_0$ as prescribed by $g$, $g g_{m_0, m}^{-1} \actsOnPoint m = g \actsOnPoint m_0$, and thirdly undoing the first rotation, $g_{m_0, m} g g_{m_0, m}^{-1} \actsOnPoint m = g_{m_0, m} \actsOnPoint (g \actsOnPoint m_0)$, in other words, by first changing the rotation axis of $g$ such that the new axis stands to the line through the centre and $m$ as the old one stood to the line through the centre and $m_0$, $g_{m_0, m} g g_{m_0, m}^{-1}$, and secondly rotating $m$ as prescribed by this new rotation.

    Let $N_0$ be a subset of the sphere $M$, which we think of as a geometrical object on the sphere that has its centre at $m_0$, for example, a circle of latitude. The set $N = \setOf{g G_0 \in G \modulo G_0 \suchThat g \actsOnPoint m_0 \in N_0} = \setOf{G_{m_0, m} \suchThat m \in N_0} = \setOf{g_{m_0, m} G_0 \suchThat m \in N_0}$ can be thought of as a realisation of $N_0$ in $G \modulo G_0$, because $m_0 \isSemiActedUponBy N = g_{m_0, m_0} \setOf{g_{m_0, m} \suchThat m \in N_0} \actsOnPoint m_0 = N_0$. Furthermore, for each point $m \in M$, the set $m \isSemiActedUponBy N = g_{m_0, m} \actsOnPoint N_0$ has the same shape and size as $N_0$ but its centre at $m$. Note that $N = \iota(N_0)$, where $\iota$ is the bijection from \cref{lemma:iota-is-bijective-and-equivariant}. 
  \end{example}


  \begin{lemma}
  \label{lemma:identification-of-G-quotient-Gzero-with-M-by-right-semiaction}
    The maps
    \begin{equation*}
      \left\{
      \begin{aligned}
        \iota \from M &\to G \modulo G_0,\\
        m &\mapsto G_{m_0, m},
      \end{aligned}
      \right\}
      \text{ and }
      \left\{
      \begin{aligned}
        m_0 \isSemiActedUponBy \blank \from G \modulo G_0 &\to M,\\
        \mathfrak{g} &\mapsto m_0 \isSemiActedUponBy \mathfrak{g},
      \end{aligned}
      \right\}
    \end{equation*}
    are inverse to each other and, under the identification of $G \modulo G_0$ with $M$ by either of these maps,
    \begin{equation*}
      \ForEach g \in G \ForEach \mathfrak{g} \in G \modulo G_0 \simeq M \Holds g \cdot \mathfrak{g} = g \actsOnPoint \mathfrak{g},
    \end{equation*}
    and
    \begin{equation*}
      \ForEach m \in M \ForEach \mathfrak{g} \in G \modulo G_0 \simeq M \Holds m \isSemiActedUponBy \mathfrak{g} = g_{m_0, m} \actsOnPoint \mathfrak{g}. \qedhere
    \end{equation*}
  \end{lemma}

  \begin{proof}
    According to \cref{lemma:iota-is-bijective-and-equivariant}, the map $\iota$ is bijective and, for each $m \in M$,
    \begin{align*}
      m_0 \isSemiActedUponBy \iota(m)
      &= m_0 \isSemiActedUponBy G_{m_0, m}\\
      &= m_0 \isSemiActedUponBy g_{m_0, m} G_0\\
      &= g_{m_0, m_0} g_{m_0, m} \actsOnPoint m_0\\
      &= e_G \actsOnPoint m\\
      &= m.
    \end{align*}
    Therefore, $m_0 \isSemiActedUponBy \blank = \iota^{-1}$. Moreover, according to \cref{lemma:iota-is-bijective-and-equivariant}, the map $\iota$ is $(\actsOnPoint, \cdot)$-e\-qui\-var\-i\-ant. Hence, for each $g \in G$ and each $\mathfrak{g} \in G \modulo G_0$,
    \begin{align*}
      g \cdot \mathfrak{g}
      = g \cdot \iota(\iota^{-1}(\mathfrak{g}))
      = \iota(g \actsOnPoint \iota^{-1}(\mathfrak{g})).
    \end{align*}
    And, for each $m \in M$ and each $g G_0 \in G \modulo G_0$,
    \begin{align*}
      m \isSemiActedUponBy g G_0
      &= g_{m_0, m} g \actsOnPoint m_0\\
      &= g_{m_0, m} \actsOnPoint (g \actsOnPoint m_0)\\
      &= g_{m_0, m} \actsOnPoint \iota^{-1}(G_{m_0, g \actsOnPoint m_0})\\
      &= g_{m_0, m} \actsOnPoint \iota^{-1}(g G_0). \qedhere
    \end{align*}
  \end{proof}

  \begin{example}[Group]
  \label{example:group:identification}
    In the situation of \cref{example:group:liberation}, the maps $m_0 \isSemiActedUponBy \blank$ and $\iota$ are the identity maps on $G \simeq G \modulo G_0$.
  \end{example}

  \begin{example}[Plane]
  \label{example:plane:identification}
    In the situation of \cref{example:plane:liberation}, the maps $m_0 \isSemiActedUponBy \blank$ and $\iota$ map translations to points, which encode translation vectors, and vice versa.
  \end{example}

  \begin{example}[Sphere]
  \label{example:sphere:identification}
    In the situation of \cref{example:sphere:liberation}, under the identification of $G \modulo G_0$ with $\setOf{g_{m_0, m} \suchThat m \in M}$ by $g_{m_0, m} G_0 \mapsto g_{m_0, m}$, the maps $m_0 \isSemiActedUponBy \blank$ and $\iota$ map rotations to points, which encode rotation angles and axes in the $(x, y)$-plane, and vice versa.
  \end{example}

  \begin{lemma}
  \label{lemma:right-semi-action-similar-to-leftaction-in-m0}
    The semi-action $\isSemiActedUponBy$ is \define{similar to $\actsOnPoint$ in $m_0$}, which means that
    \begin{equation*}
      \ForEach g \in G \Holds m_0 \isSemiActedUponBy g G_0 = g \actsOnPoint m_0.
    \end{equation*}
    In particular,
    \begin{equation*}
      \ForEach g \in G \ForEach \mathfrak{g} \in G \modulo G_0 \Holds m_0 \isSemiActedUponBy g \cdot \mathfrak{g} = g \actsOnPoint (m_0 \isSemiActedUponBy \mathfrak{g}),
    \end{equation*}
    and
    \begin{equation*}
      \ForEach m \in M \ForEach \mathfrak{g} \in G \modulo G_0 \Holds m \isSemiActedUponBy \mathfrak{g} = g_{m_0, m} \actsOnPoint (m_0 \isSemiActedUponBy \mathfrak{g}). \qedhere
    \end{equation*}
  \end{lemma}

  \begin{proof}
    The similarity follows from the fact that $g_{m_0, m_0} = e_G$. And the other two properties follow from the fact that $\actsOnPoint$ is a left group action.
  \end{proof}

  \begin{lemma}
  \label{lemma:liberation:transitive-and-free}
    The semi-action $\isSemiActedUponBy$ is
    \begin{aenumerate}
      \item \label{item:liberation:transitive-and-free:free}
            \define{free}\graffito{free}, which means that
            \begin{multline*}
              \ForEach \mathfrak{g} \in G \modulo G_0 \ForEach \mathfrak{g}' \in G \modulo G_0 \Holds\\
                    (\Exists m \in M \SuchThat m \isSemiActedUponBy \mathfrak{g} = m \isSemiActedUponBy \mathfrak{g}') \implies \mathfrak{g} = \mathfrak{g}';
            \end{multline*}
      \item \label{item:liberation:transitive-and-free:transitive}
            \define{transitive}\graffito{transitive}, which means that the set $M$ is non-empty and
            \begin{equation*}
              \ForEach m \in M \ForEach m' \in M \Exists \mathfrak{g} \in G \modulo G_0 \SuchThat m \isSemiActedUponBy \mathfrak{g} = m'. \qedhere
            \end{equation*}
    \end{aenumerate}
  \end{lemma}

  \begin{proof}
    \begin{aenumerate}
      \item Let $m \in M$ and let $m' \in M$. Put $m'' = g_{m_0, m}^{-1} \actsOnPoint m'$. Because $\actsOnPoint$ is transitive, there is a $g \in G$ such that $g \actsOnPoint m_0 = m''$. Hence,
            \begin{align*}
              m \isSemiActedUponBy g G_0
              &= g_{m_0, m} g \actsOnPoint m_0\\
              &= g_{m_0, m} \actsOnPoint (g \actsOnPoint m_0)\\
              &= g_{m_0, m} \actsOnPoint m''\\
              &= g_{m_0, m} \actsOnPoint (g_{m_0, m}^{-1} \actsOnPoint m')\\
              &= e_G \actsOnPoint m'\\
              &= m'.
            \end{align*}
      \item Let $g G_0$ and $g' G_0$ be two elements of $G \modulo G_0$, and let $m$ be an element of $M$ such that $m \isSemiActedUponBy g G_0 = m \isSemiActedUponBy g' G_0$. Then, $g_{m_0, m} g \actsOnPoint m_0 = g_{m_0, m} g' \actsOnPoint m_0$. Hence, $g \actsOnPoint m_0 = g' \actsOnPoint m_0$. Therefore, $g^{-1} g' \actsOnPoint m_0 = m_0$. Thus, $g^{-1} g' \in G_0$. In conclusion, $g G_0 = g' G_0$. \qedhere
    \end{aenumerate}
  \end{proof}

  \begin{example}[Group]
  \label{example:group:liberation:transitive-and-free}
    In the situation of \cref{example:group:identification}, the semi-action $\isSemiActedUponBy$, being but right multiplication, is free and transitive.
  \end{example}

  \begin{example}[Plane]
  \label{example:plane:liberation:transitive-and-free}
    In the situation of \cref{example:plane:identification}, the semi-action $\isSemiActedUponBy$, being but translation, is free and transitive.
  \end{example}

  \begin{example}[Sphere]
  \label{example:sphere:liberation:transitive-and-free}
    In the situation of \cref{example:sphere:identification}, under the identification of $G \modulo G_0$ with $M$ by $\iota$, according to \cref{lemma:identification-of-G-quotient-Gzero-with-M-by-right-semiaction}, for each point $m \in M$, the map $m \isSemiActedUponBy \blank$ is the rotation by $g_{m_0, m}$, which is injective and surjective, and therefore the semi-action $\isSemiActedUponBy$ is free and transitive.
  \end{example}

  \begin{lemma}
  \label{lemma:semi-commutativity-of-liberation} 
    The semi-action $\isSemiActedUponBy$
    \begin{aenumerate}
      \item \label{item:semi-commutativity-of-liberation:semi-commutes}
            \define{semi-commutes with $\actsOnPoint$}\graffito{semi-commutes with $\actsOnPoint$}\index{commutes with@semi-commutes with $\actsOnPoint$}, which means that, for each $\mathcal{K}$-big subgroup $H$ of $G$,
            \begin{multline*}
              \ForEach m \in M \ForEach h \in H \Exists h_0 \in H_0 \SuchThat \ForEach \mathfrak{g} \in G \modulo G_0 \Holds\\
                    (h \actsOnPoint m) \isSemiActedUponBy \mathfrak{g} = h \actsOnPoint (m \isSemiActedUponBy h_0 \cdot \mathfrak{g});
            \end{multline*}
      \item \label{item:semi-commutativity-of-liberation:exhausts}
            \define{exhausts its defect with respect to its semi-commutativity with $\actsOnPoint$ in $m_0$}\graffito{exhaust ones defect with respect to ones semi-commutativity with $\actsOnPoint$ in $m_0$}, which means that 
            \begin{equation*}
              \ForEach g_0 \in G_0 \ForEach \mathfrak{g} \in G \modulo G_0 \Holds
                    (g_0^{-1} \actsOnPoint m_0) \isSemiActedUponBy \mathfrak{g} = g_0^{-1} \actsOnPoint (m_0 \isSemiActedUponBy g_0 \cdot \mathfrak{g}). \qedhere
            \end{equation*}
    \end{aenumerate}
  \end{lemma}

  \begin{proof}
    Let $H$ be a $\mathcal{K}$-big subgroup of $G$.
    \begin{aenumerate}
      \item Let $h \in H$ and let $m \in M$. Put $h_0 = g_{m_0, m}^{-1} h^{-1} g_{m_0, h \actsOnPoint m}$. Because $g_{m_0, m}^{-1} \in H$, $g_{m_0, h \actsOnPoint m}^{-1} \in H$, and
            \begin{align*}
              g_{m_0, m}^{-1} h^{-1} g_{m_0, h \actsOnPoint m} \actsOnPoint m_0
              &= g_{m_0, m}^{-1} h^{-1} \actsOnPoint (h \actsOnPoint m)\\
              &= g_{m_0, m}^{-1} \actsOnPoint m\\
              &= m_0,
            \end{align*}
            we have $h_0 \in H_0$. Moreover, for each $g G_0 \in G \modulo G_0$,
            \begin{align*}
              (h \actsOnPoint m) \isSemiActedUponBy g G_0
              &= g_{m_0, h \actsOnPoint m} g \actsOnPoint m_0\\
              &= h g_{m_0, m} h_0 g \actsOnPoint m_0\\
              &= h \actsOnPoint (g_{m_0, m} h_0 g \actsOnPoint m_0)\\
              &= h \actsOnPoint (m \isSemiActedUponBy h_0 \cdot g G_0).
            \end{align*}
      \item For each $g_0 \in G_0$ and each $g G_0 \in G \modulo G_0$, because $g_{m_0, m_0} = e_G$,
            \begin{align*}
              (g_0^{-1} \actsOnPoint m_0) \isSemiActedUponBy g G_0
              &= m_0 \isSemiActedUponBy g G_0\\
              &= g_{m_0, m_0} g \actsOnPoint m_0\\
              &= g_0^{-1} g_{m_0, m_0} g_0 g \actsOnPoint m_0\\
              &= g_0^{-1} \actsOnPoint (g_{m_0, m_0} g_0 g \actsOnPoint m_0)\\
              &= g_0^{-1} \actsOnPoint (m_0 \isSemiActedUponBy g_0 g G_0)\\
              &= g_0^{-1} \actsOnPoint (m_0 \isSemiActedUponBy g_0 \cdot g G_0). \qedhere
            \end{align*}
    \end{aenumerate}
  \end{proof}

  \begin{example}[Group]
  \label{example:group:liberation:commute}
    In the situation of \cref{example:group:liberation:transitive-and-free}, because group multiplication is associative, the right multiplication $\isSemiActedUponBy$ commutes with the left multiplication $\actsOnPoint$.
  \end{example}

  \begin{example}[Plane]
  \label{example:plane:liberation:commute}
    In the situation of \cref{example:plane:liberation:transitive-and-free}, because the group $T$ of translations is abelian, the right semi-action $\isSemiActedUponBy$ commutes with the left action $\actsOnPoint_T$. Indeed,
    \begin{equation*}
      \ForEach t \in T \ForEach \mathfrak{t} \in T \Holds (t \actsOnPoint_T \blank) \isSemiActedUponBy \mathfrak{t} = \mathfrak{t}(t(\blank)) = t(\mathfrak{t}(\blank)) = t \actsOnPoint_T (\blank \isSemiActedUponBy \mathfrak{t}). \qedhere
    \end{equation*}
  \end{example}

  \begin{example}[Cayley Graph]
  \label{example:cell-space-cayley-graph}
    Let $M$ be a group, let $S$ be a generating set of $M$, and let $\mathcal{G} = \ntuple{M, E, \lambda}$ be the coloured $S$-Cayley graph of $M$. The graph $\mathcal{G}$ is edge-labelled and directed, its set of vertices $M$ is the set underlying $M$, its set of edges $E$ is the subset $\setOf{(m, m s) \suchThat m \in M, s \in S}$ of $M \times M$, and its edge-labelling $\lambda$ is the map $E \to S$, $(m, m s) \mapsto s$.

    The automorphism group of $\mathcal{G}$, namely
    \begin{align*}
      G = \setOf{
        &g \from M \to M \text{ bijective} \suchThat
        \ForEach m \in M \ForEach m' \in M \ForEach s \in S \Holds\\
        &(m, m') \in E \land \lambda(m, m') = s\\
        &\ifAndOnlyIf (g(m), g(m')) \in E \land \lambda(g(m), g(m')) = s
      },
    \end{align*}
    acts freely and transitively on $M$ on the left by function application, which we denote by $\actsOnPoint$. For each element $m \in M$, left multiplication by $m$, which we denote by $m \cdot \blank$, is the unique graph automorphism that maps the neutral element $e_M$ to $m$ (note that all graph automorphisms are of this form).

    The triple $\mathcal{M} = \ntuple{M, G, \actsOnPoint}$ is a principal left-ho\-mo\-ge\-neous space and the tuple $\mathcal{K} = \ntuple{e_M, \family{m \cdot \blank}_{m \in M}}$ is a coordinate system for $\mathcal{M}$. Under the identification of $G \modulo G_0$ with $G$, the induced semi-action $\isSemiActedUponBy$ is the right group action of $G$ on $M$ given by $(m, g) \mapsto m \cdot g(e_M)$.
  \end{example}

  \begin{example}[Normal]
  \label{example:normal:homogeneous-space:liberation}
    In the situation of \cref{example:normal:homogeneous-space}, let $\mathcal{K} = \ntuple{m_0, \family{g_{m_0, m}}_{m \in M}}$ be a coordinate system for $\mathcal{M}$. The tuple $\mathcal{K}' = \ntuple{m_0, \family{g_{m_0, m} G_0}_{m \in M}}$ is the one and only coordinate system for $\mathcal{M}'$ with origin $m_0$. Both coordinate systems induce the right quotient group action $\isSemiActedUponBy$ of $G \modulo G_0$ on $M$ given by $(m, g G_0) \mapsto g_{m_0, m} g \actsOnPoint m_0$. In the case that the group $G$ is abelian, this action is given by $(m, g G_0) \mapsto g \actsOnPoint m$, which is, apart from the order of the arguments, identical to $\cosetActsOnPoint$. In any case, under the identification of $M$ with $G \modulo G_0$ by $\iota$, the left and right quotient group actions $\cosetActsOnPoint$ and $\isSemiActedUponBy$ are identical to the group multiplication of $G \modulo G_0$. And, the right quotient group action $\isSemiActedUponBy$ commutes with the left group action $\actsOnPoint$.
  \end{example}

  \section{Semi-Cellular, Big-Cellular, and Cellular Automata}
  \label{section:automata}

  \paragraph{Introduction.} Let us consider the following discrete-time dynamical system on a circle whose points can be coloured black and white. In one time step a point becomes white if there is a white point nearby in the clockwise direction, where two points are said to be near each other if their arc distance is not greater than, say, one-hundredth of the circle's circumference; and otherwise retains its colour. The time evolution of that system is uniform, in the sense that each point determines its next colour by the same rule, and it is local, in the sense that each point determines its next colour by means of the colours of nearby points. Moreover, it is equivariant under rotations of the circle, but, despite its uniformity, it is not equivariant under reflections.

  For example, if the left side of the circle, without the top and bottom points, is white and the other points are black, then in one time step the top point stays black; but, if we first reflect the circle about the vertical line through the centre of the circle, secondly do one time step, and lastly reflect again, then the top point is white. The reason is that the local rule uses the direction of rotation, namely clockwise, which stays the same under rotations but changes under reflections. If the rule considered all nearby points regardless of whether they lie in the clockwise or anticlockwise direction, then time evolution would be equivariant under all symmetries of the circle.

  The time evolution of this dynamical system is the global transition function of a semi-cellular automaton whose cells are the points of the circle, whose states are the colours black and white, whose neighbourhood is the set of all nearby points (in either direction) of a designated point, and whose local transition function is the local rule for the designated point described above. Actually, the nearby points in the clockwise direction would be sufficient as neighbourhood, but those are not invariant under the reflection that stabilises the designated point and for technical reasons we want this invariance.

  Without designating a point, the neighbourhood can be described by all rotations that map a point to a point nearby and the local transition function by a map that maps a local configuration to the colour white if there is a white clockwise neighbour, where we call a neighbour clockwise if it maps a point to a point nearby in the clockwise direction; and otherwise to the colour of the identity map, which plays the role of the designated point. The neighbourhood of a point can be recovered by applying the (relative) neighbours to the point, in other words, by acting with the (relative) neighbours on the point on the right. As we have seen in the introduction of \cref{section:semi-action}, this right action is equivalent to the right semi-action induced by the left action of the symmetries on the circle, where the rotations are identified with the quotient set of the symmetries of the circle by the stabiliser of a designated point.

  \paragraph{Contents.} A semi-cellular automaton is made up of a cell space, a set of states, a neighbourhood (think of a disk about the origin as points or vectors), and a local transition function (see \cref{definition:semi-cellular-automaton}). The stabiliser of the origin acts on local configurations (think of rotations of disk-shaped patterns; see \cref{definition:local-configuration}). The group acts on global configurations (think of rotations and translations of unbounded patterns; see \cref{definition:global-configuration}). A cell observes a local configuration by first semi-acting on the right on itself with the (relative) neighbourhood to determine its neighbours (think of point-vector additions) and secondly reading the states of these neighbours (see \cref{definition:observed-local-configuration}); or, alternatively, by first translating the global configuration such that the cell is moved to the origin and secondly restricting this translated configuration to the neighbourhood (see \cref{remark:identification:observed-local-configuration}). The global transition function applies the local transition function synchronously to the observed local configuration of each cell to determine its new state (see \cref{definition:global-transition-function} and \cref{remark:identification:global-transition-function}). This function is equivariant under the action on global configurations (traditionally known as shift-invariance) if and only if the local transition function is invariant under the action on local configurations (see \cref{theorem:local-invariance-versus-global-equivariance}). If the latter is the case, then the automaton is a cellular automaton (see \cref{definition:cellular-automaton}). And, if the local transition function is only invariant under a restriction of the action on local configurations to a big subgroup, then the automaton is a big-cellular automaton (see \cref{definition:big-cellular-automaton}).

  \paragraph{Body.} In this section, let $\mathcal{R} = \ntuple{\mathcal{M}, \mathcal{K}} = \ntuple{\ntuple{M, G, \actsOnPoint}, \ntuple{m_0, \family{g_{m_0, m}}_{m \in M}}}$ be a cell space.

  \begin{definition}
  \label{definition:semi-cellular-automaton}
    Let $Q$ be a set, let $N$ be a subset of $G \modulo G_0$ such that $G_0 \cdot N \subseteq N$, and let $\delta$ be a map from $Q^N$ to $Q$. The quadruple $\mathcal{C} = \ntuple{\mathcal{R}, Q, N, \delta}$ is called \define{semi-cellular automaton}\graffito{semi-cellular automaton $\mathcal{C}$}\index{automaton!semi-cellular}\index[symbols]{Ccalligraphic@$\mathcal{C}$}, each element $q \in Q$ is called \define{state}\graffito{state $q$}\index[symbols]{q@$q$}, the set $N$ is called \define{neighbourhood}\graffito{neighbourhood $N$}\index[symbols]{N@$N$}, each element $n \in N$ is called \define{neighbour}\graffito{neighbour $n$}\index[symbols]{n@$n$}, and the map $\delta$ is called \define{local transition function}\graffito{local transition function $\delta$}\index{transition function!local}\index[symbols]{delta@$\delta$}.
  \end{definition}

  \begin{remark}
    Under the identification of $G \modulo G_0$ with $M$ by $\iota$, the neighbourhood $N$ is a subset of $M$ such that $G_0 \actsOnPoint N \subseteq N$
  \end{remark}

  \begin{example}[Group]
  \label{example:group:semi-cellular-automaton}
    In the situation of \cref{example:group:liberation:commute}, the semi-cellular automata over $\mathcal{R}$ are the usual cellular automata over the group $G$.
  \end{example}

  \begin{example}[Plane]
  \label{example:plane:semi-cellular-automaton}
    In the situation of \cref{example:plane:liberation:commute}, let $Q$ be the set of real numbers, let $M$ be identified with $G \modulo G_0$ by $\iota$, let $\varepsilon$ be a positive real number, let $N$ be the open disk of radius $\varepsilon$ about $m_0$, and let
    \begin{align*}
      \delta \from Q^N &\to Q,\\
      \ell &\mapsto \begin{dcases*}
                      \frac{\partial^2 \ell}{\partial x^2}(m_0) + \frac{\partial^2 \ell}{\partial y^2}(m_0), &\parbox[t]{.41\textwidth}{if $\ell$ is twice continuously differentiable on an open disk about $m_0$,}\\ 
                      0, &otherwise,
                    \end{dcases*}
    \end{align*}
    where $\partial^2 \ell / \partial x^2(m_0)$ and $\partial^2 \ell / \partial y^2(m_0)$ are the second-order partial derivatives of $\ell$ by $x$ and $y$ at $m_0$. The quadruple $\mathcal{C} = \ntuple{\mathcal{R}, Q, N, \delta}$ is a semi-cellular automaton.
  \end{example}

  \begin{example}[Sphere]
  \label{example:sphere:semi-cellular-automaton}
    In the situation of \cref{example:sphere:liberation:transitive-and-free}, let $Q$ be the set $\setOf{0,1}$, let $N_0$ be the union of all circles of latitude between $45\degree$ and $90\degree$ north, which is a curved circular disk of radius $\pi/4$ with the north pole $m_0$ at its centre, let $N$ be the set $\iota(N_0) = \setOf{G_{m_0, m} \suchThat m \in N_0}$, and let 
    \begin{align*}
      \delta \from Q^N &\to Q,\\
      \ell &\mapsto \begin{dcases*}
                      0, &if $\forall n \in N \Holds \ell(n) = 0$,\\
                      1, &if $\exists n \in N \SuchThat \ell(n) = 1$.
                    \end{dcases*}
    \end{align*}
    The quadruple $\mathcal{C} = \ntuple{\mathcal{R}, Q, N, \delta}$ is a semi-cellular automaton.
  \end{example}

  In the remainder of this section, let $\mathcal{C} = \ntuple{\mathcal{R}, Q, N, \delta}$ be a semi-cellular automaton.

  %


  \begin{definition}
  \label{definition:local-configuration}
    Each map $\ell \in Q^N$ is called \define{local configuration}\graffito{local configuration $\ell$}\index{configuration!local}\index[symbols]{l@$\ell$}. The stabiliser $G_0$ acts on $Q^N$ on the left by
    \begin{align*}
      \bullet \from G_0 \times Q^N &\to Q^N, \mathnote{left group action $\bullet$ of $G_0$ on $Q^N$}\index[symbols]{bullet@$\bullet$}\\
      (g_0, \ell) &\mapsto [n \mapsto \ell(g_0^{-1} \cdot n)]. \qedhere
    \end{align*}
  \end{definition}

  \begin{remark}
    Under the identification of $G \modulo G_0$ with $M$ by $\iota$, we have $\bullet \from (g_0, \ell) \mapsto [n \mapsto \ell(g_0^{-1} \actsOnPoint n)]$.
  \end{remark}



  \begin{definition}
  \label{definition:big-cellular-automaton}
    The semi-cellular automaton $\mathcal{C}$ is called \graffito{big-cellular automaton $\mathcal{C}$}\define{big-cellular automaton}\index{automaton!big-cellular}\index{cellular automaton!big}\index[symbols]{Ccalligraphic@$\mathcal{C}$} if and only if there is a $\mathcal{K}$-big subgroup $H$ of $G$ such that the local transition function $\delta$ is $\bullet_{H_0}$-invariant.
  \end{definition}

  \begin{remark}
  \label{remark:induced-local-transition-function-on-quotient}
    Let $H$ be a $\mathcal{K}$-big subgroup of $G$ and let $\pi$ be the canonical projection of $Q^N$ onto $H_0 \reverseModulo Q^N$. The local transition function $\delta$ is $\bullet_{H_0}$-invariant if and only if there is a map $\mathfrak{d} \from H_0 \reverseModulo Q^N \to Q$ such that $\delta = \mathfrak{d} \after \pi$, in other words, if and only if the map $\mathfrak{d} \from H_0 \reverseModulo Q^N \to Q$, $H_0 \bullet \ell \mapsto \delta(\ell)$, is well-defined. 
  \end{remark}

  \begin{remark} 
    Because each $\mathcal{K}$-big subgroup of $G$ includes the $\mathcal{K}$-big subgroup of $G$ that is generated by the set $\setOf{g_{m_0, m} \suchThat m \in M}$ of coordinates, the semi-cellular automaton $\mathcal{C}$ is a big-cellular automaton if and only if its local transition function $\delta$ is $\bullet_{G_0 \cap \groupGeneratedBy{g_{m_0, m} \suchThat m \in M}}$-invariant.
  \end{remark}

  \begin{definition}
  \label{definition:cellular-automaton}
    The semi-cellular automaton $\mathcal{C}$ is called \graffito{cellular automaton $\mathcal{C}$}\define{cellular automaton}\index{automaton!cellular}\index[symbols]{Ccalligraphic@$\mathcal{C}$} if and only if its local transition function $\delta$ is $\bullet$-invariant.
  \end{definition}

  \begin{example}[Group]
  \label{example:group:cellular-automaton}
    In the situation of \cref{example:group:semi-cellular-automaton}, the stabiliser $G_0$ of the neutral element $m_0$ is the trivial subgroup $\setOf{e_G}$ of $G$. Therefore, each semi-cellular automaton over $\mathcal{R}$ has a $\bullet$-invariant local transition function and is hence a cellular automaton.
  \end{example}

  \begin{example}[Plane]
  \label{example:plane:cellular-automaton}
    In the situation of \cref{example:plane:semi-cellular-automaton}, think of the neighbourhood $N$ as being embedded in the $(x,y)$-plane of the three-dimensional Euclidean space and of local configurations as graphs extending into the $z$-direction. The rotations $G_0$ about the $z$-axis act on these graphs by $\bullet$ by rotating them. For each rotation $g_0 \in G_0$ and each local configuration $\ell \in Q^N$, the map $\ell$ is twice continuously differentiable on an open disk about $m_0$ if and only if $g_0 \bullet \ell$ is; if they both are, then $\delta(\ell) = \delta(g_0 \bullet \ell)$ by elementary calculus; otherwise, $\delta(\ell) = 0 = \delta(g_0 \bullet \ell)$. Hence, the local transition function $\delta$ is $\bullet$-invariant. Therefore, the quadruple $\mathcal{C}$ is a cellular automaton.
  \end{example} 

  \begin{example}[Sphere]
  \label{example:sphere:cellular-automaton}
    In the situation of \cref{example:sphere:semi-cellular-automaton}, think of $0$ as black, $1$ as white, and of local configurations as black-and-white patterns on $N_0 = m_0 \isSemiActedUponBy N$. The rotations $G_0$ about the $z$-axis act on these patterns by $\bullet$ by rotating them. The local transition function $\delta$ maps the black pattern to $0$ and all others to $1$, which is invariant under rotations. Therefore, the quadruple $\mathcal{C}$ is a cellular automaton.
  \end{example}

  \begin{definition}
  \label{definition:global-configuration}
    The set $Q^M$ is called \define{phase space}\graffito{phase space $Q^M$}\index[symbols]{QM@$Q^M$} and each map $c \in Q^M$ is called \define{global configuration}\graffito{global configuration $c$}\index{configuration!global}\index[symbols]{c@$c$}. The group $G$ acts on $Q^M$ on the left by
    \begin{align*}
      \actsOnMap \from G \times Q^M &\to Q^M, \mathnote{left group action $\actsOnMap$ of $G$ on $Q^M$}\index[symbols]{arrowrightblack@$\actsOnMap$}\\
      (g, c) &\mapsto [m \mapsto c(g^{-1} \actsOnPoint m)]. \qedhere
    \end{align*}
  \end{definition} 


  \begin{definition} 
  \label{definition:observed-local-configuration}
    For each cell $m \in M$, the set $m \isSemiActedUponBy N$ is called \graffito{neighbourhood $m \isSemiActedUponBy N$ of $m$}\defineX{neighbourhood of $m$}{neighbourhood!of $m$}\index[symbols]{marrowleftunderscoreN@$m \isSemiActedUponBy N$}. And, for each global configuration $c \in Q^M$ and each cell $m \in M$, the local configuration
    \begin{align*}
      c(m \isSemiActedUponBy \blank)\restrictedTo_N \from N &\to Q, \mathnote{local configuration $c(m \isSemiActedUponBy \blank)\restrictedTo_N$ observed by $m$ in $c$}\\
      n &\mapsto c(m \isSemiActedUponBy n),
    \end{align*}
    is called \define{observed by $m$ in $c$}\index{local configuration!observed by $m$ in $c$}\index{configuration!local!observed by $m$ in $c$}.
  \end{definition}

  \begin{remark}
  \label{remark:identification:observed-local-configuration}
    Under the identification of $G \modulo G_0$ with $M$ by $\iota$, according to the proof of the forthcoming \cref{lemma:global-transition-function-without-liberation}, for each cell $m \in M$, we have $c(m \isSemiActedUponBy \blank)\restrictedTo_N = (g_{m_0, m}^{-1} \actsOnMap c)\restrictedTo_N$. 
  \end{remark}

  \begin{remark}
  \label{remark:local-configuration-induces-global-configuration} 
    Because the semi-action $\isSemiActedUponBy$ is free (see \cref{item:liberation:transitive-and-free:free} of lemma \ref{lemma:liberation:transitive-and-free}), for each local configuration $\ell \in Q^N$ and each cell $m \in M$, there is a global configuration $c \in Q^M$ such that the local configuration observed by $m$ in $c$ is $\ell$. 
  \end{remark}

  \begin{definition}
  \label{definition:global-transition-function}
    The map
    \begin{align*}
      \Delta \from Q^M &\to Q^M, \mathnote{global transition function $\Delta$}\index{transition function!global}\index[symbols]{Delta@$\Delta$}\\
      c &\mapsto [m \mapsto \delta(n \mapsto c(m \isSemiActedUponBy n))], 
    \end{align*}
    is called \define{global transition function}.
  \end{definition}

  \begin{remark}
  \label{remark:identification:global-transition-function}
    Under the identification of $G \modulo G_0$ with $M$ by $\iota$, according to \cref{lemma:global-transition-function-without-liberation}, we have $\Delta \from c \mapsto [m \mapsto \delta((g_{m_0, m}^{-1} \actsOnMap c)\restrictedTo_N)]$.
  \end{remark}


  \begin{example}[Plane]
  \label{example:plane:global-transition-function}
    In the situation of \cref{example:plane:cellular-automaton}, the restriction of the global transition function $\Delta$ to the twice continuously differentiable maps is known as \define{Laplace operator}\graffito{Laplace operator}. Recall that the neighbourhood $N$ is an open disk of radius $\varepsilon$ about $m_0$, where $\varepsilon$ is a positive real number. The map $\Delta$ does not depend on the radius $\varepsilon$ --- it can be chosen arbitrarily small without affecting $\Delta$. In other words, there is no smallest neighbourhood for cellular automata whose global transition functions are $\Delta$. 
  \end{example}

  \begin{example}[Sphere]
  \label{example:sphere:global-transition-function}
    In the situation of \cref{example:sphere:cellular-automaton}, repeated applications of the global transition function of $\mathcal{C}$ grows white regions on $M$.
  \end{example}

  \begin{example}[Hyperbolic Game of Life]
    Sébastien Moriceau presents an adaptation of Conway's Game of Life cellular automaton on a tessellation of the hyperbolic plane in example~3.3 (a) in \cite{moriceau:2011}.
  \end{example}

  \begin{remark}
  \label{remark:dependence-of-new-states-of-cells-on-other-cells}
    For each subset $A$ of $M$ and each global configuration $c \in Q^M$, the states of the cells $A$ in $\Delta(c)$ \define{depend at most on}\graffito{depend at most on} the states of the cells $A \isSemiActedUponBy N$ in $c$, which means that
    \begin{multline*}
      \ForEach A \subseteq M \ForEach c \in Q^M \ForEach c' \in Q^M \Holds\\
          c\restrictedTo_{A \isSemiActedUponBy N} = c'\restrictedTo_{A \isSemiActedUponBy N} \implies \Delta(c)\restrictedTo_A = \Delta(c')\restrictedTo_A. \qedhere
    \end{multline*}
  \end{remark}

  \begin{remark}
  \label{remark:determined-by-behaviour-at-origin-without-assumptions}
    The global transition function $\Delta$ is \graffito{determined by its behaviour at the origin $m_0$}\define{determined by its behaviour at the origin $m_0$}, which means that
    \begin{equation*}
      \ForEach c \in Q^M \ForEach m \in M \Holds \Delta(c)(m) = \Delta(g_{m_0, m}^{-1} \actsOnMap c)(m_0).
    \end{equation*}
    Indeed, for each global configuration $c \in Q^M$ and each cell $m \in M$, according to \cref{lemma:right-semi-action-similar-to-leftaction-in-m0},
    \begin{align*}
      \Delta(c)(m)
      &= \delta\parens{n \mapsto c(m \isSemiActedUponBy n)}\\
      &= \delta\parens{n \mapsto c(g_{m_0, m} \actsOnPoint (m_0 \isSemiActedUponBy n))}\\
      &= \delta\parens{n \mapsto (g_{m_0, m}^{-1} \actsOnMap c)(m_0 \isSemiActedUponBy n)}\\
      &= \Delta(g_{m_0, m}^{-1} \actsOnMap c)(m_0). \qedhere
    \end{align*}
  \end{remark}

  \begin{definition}
    Let $E$ be a subset of $N$. It is called \graffito{sufficient neighbourhood $E$}\define{sufficient neighbourhood}\index{neighbourhood!sufficient}\index[symbols]{E@$E$} if and only if
    \begin{equation*}
      \ForEach \ell \in Q^N \ForEach \ell' \in Q^N \Holds \parens[\big]{\ell\restrictedTo_E = \ell'\restrictedTo_E \implies \delta(\ell) = \delta(\ell')},
    \end{equation*}
    in which case the map
    \begin{align*}
      \eta \from Q^E &\to Q,\\
      \ell_E &\mapsto \delta(\ell), \text{ where $\ell \in Q^N$ such that $\ell\restrictedTo_E = \ell_E$},
    \end{align*}
    is called \define{sufficient local transition function}\graffito{sufficient local transition function $\eta$}\index{local transition function!sufficient}\index[symbols]{eta@$\eta$}.
  \end{definition}

  \begin{remark}
    The neighbourhood itself is a sufficient neighbourhood. And, if the local transition function depends on all neighbours, then the neighbourhood is the only one. In general, it is impossible to choose the neighbourhood such that the local transition function depends on all neighbours, because it may depend only on an arbitrarily small disk about the origin (as in \cref{example:plane:global-transition-function}) or it may depend on a neighbour $n$ but not on $g_0 \cdot n$ for some stabiliser $g_0$ of the origin (as in the example in the introduction of \cref{section:automata}). Note that, because $G_0 \cdot N \subseteq N$, for each $n \in N$, we have $g_0 \cdot n \in N$. So, the assumption $G_0 \cdot N \subseteq N$ may force the neighbourhood to be larger than is actually necessary to determine the next state of a cell. 
  \end{remark}

  Recall that the right quotient set semi-action $\isSemiActedUponBy$ semi-commutes with the left group action $\actsOnPoint$, symbolically, $(g^{-1} \actsOnPoint m) \isSemiActedUponBy n = g^{-1} \actsOnPoint (m \isSemiActedUponBy g_0 \cdot n)$, where $g_0$ does not depend on $n$. Thus, the local configuration that is observed by a cell $g^{-1} \actsOnPoint m$ in a global configuration $c$ is a rotation by $g_0$ of the local configuration that is observed by $m$ in $g \actsOnMap c$, symbolically, $c((g^{-1} \actsOnPoint m) \isSemiActedUponBy \blank) = g_0^{-1} \bullet ((g \actsOnMap c)(m \isSemiActedUponBy \blank))$ (see \cref{lemma:semi-commutativity-induced-left-action-and-bullet}). Hence, because local transition functions of cellular automata are $\bullet$-invariant, their global transition functions are $\actsOnMap$-e\-qui\-var\-i\-ant (see \cref{theorem:local-invariance-versus-global-equivariance}). Moreover, this property of observed local configurations is also essential in the proofs of other theorems of \cref{section:invariances}.

  \begin{lemma}
  \label{lemma:semi-commutativity-induced-left-action-and-bullet}
    Let $m$ be an element of $M$, let $g$ be an element of $G$, and let $g_0$ be an element of $G_0$ such that
    \begin{equation*}
      \ForEach n \in N \Holds (g^{-1} \actsOnPoint m) \isSemiActedUponBy n = g^{-1} \actsOnPoint (m \isSemiActedUponBy g_0 \cdot n).
    \end{equation*}
    For each global configuration $c \in Q^M$,
    \begin{equation*}
      [n \mapsto c((g^{-1} \actsOnPoint m) \isSemiActedUponBy n)] = g_0^{-1} \bullet [n \mapsto (g \actsOnMap c)(m \isSemiActedUponBy n)]. \qedhere
    \end{equation*}
  \end{lemma}

  \begin{proof}
    For each global configuration $c \in Q^M$,
    \begin{align*}
      [n \mapsto c((g^{-1} \actsOnPoint m) \isSemiActedUponBy n)]
      &= [n \mapsto c(g^{-1} \actsOnPoint (m \isSemiActedUponBy g_0 \cdot n))]\\
      &= g_0^{-1} \bullet [n \mapsto c(g^{-1} \actsOnPoint (m \isSemiActedUponBy n))]\\
      &= g_0^{-1} \bullet [n \mapsto (g \actsOnMap c)(m \isSemiActedUponBy n)]. \qedhere
    \end{align*}
  \end{proof}

  In the definition of a semi-cellular automaton, instead of a (relative) neighbourhood $N$ and a local transition function $\delta$, we could have used a neighbourhood $N_0$ of $m_0$ (see \cref{definition:neighbourhood-of-origin}) and a local transition function $\delta_0$ of $m_0$ (see \cref{definition:local-trans-fnc-of-origin}). Then, in the definition of the global transition function, to determine the next state of a cell, we could have translated the global configuration such that the cell is translated to $m_0$, restricted this translation to the neighbourhood of $m_0$, and applied the local transition function of $m_0$ (see \cref{lemma:global-transition-function-without-liberation}). Note that under the identification of $G \modulo G_0$ with $M$ by $\iota$, we have $N = N_0$ and $\delta = \delta_0$. 

  \begin{definition}
  \label{definition:neighbourhood-of-origin}
    The set $N_0 = m_0 \isSemiActedUponBy N$ is called \graffito{neighbourhood $N_0$ of $m_0$}\defineX{neighbourhood of $m_0$}{neighbourhood!of $m_0$}\index[symbols]{N0@$N_0$}. 
  \end{definition}

  \begin{definition}
  \label{definition:local-trans-fnc-of-origin}
    The map
    \begin{align*}
      \delta_0 \from Q^{N_0} &\to Q, \mathnote{local transition function $\delta_0$ of $m_0$}\index[symbols]{delta0@$\delta_0$}\\
      \ell_0 &\mapsto \delta(n \mapsto \ell_0(m_0 \isSemiActedUponBy n)),
    \end{align*}
    is called \defineX{local transition function of $m_0$}{local transition function!of $m_0$}.
  \end{definition}

  \begin{lemma}
  \label{lemma:global-transition-function-without-liberation}
    The global transition function $\Delta$ of $\mathcal{C}$ is identical to the map
    \begin{align*} 
      \Delta_0 \from Q^M &\to Q^M,\\
      c &\mapsto [m \mapsto \delta_0((g_{m_0, m}^{-1} \actsOnMap c)\restrictedTo_{N_0})]. \qedhere
    \end{align*}
  \end{lemma}

  \begin{proof}
    Let $c \in Q^M$ and let $m \in M$. For each $n = g G_0 \in N$,
    \begin{align*} 
      m \isSemiActedUponBy n
      &= g_{m_0, m} g \actsOnPoint m_0\\
      &= g_{m_0, m} g_{m_0, m_0} g \actsOnPoint m_0\\
      &= g_{m_0, m} \actsOnPoint (g_{m_0, m_0} g \actsOnPoint m_0)\\
      &= g_{m_0, m} \actsOnPoint (m_0 \isSemiActedUponBy n)
    \end{align*}
    and thus
    \begin{equation*}
      c(m \isSemiActedUponBy n) = c(g_{m_0, m} \actsOnPoint (m_0 \isSemiActedUponBy n)) = (g_{m_0, m}^{-1} \actsOnMap c)(m_0 \isSemiActedUponBy n).
    \end{equation*}
    Therefore,
    \begin{align*}
      \Delta(c)(m)
      &= \delta(n \mapsto c(m \isSemiActedUponBy n))\\
      &= \delta(n \mapsto (g_{m_0, m}^{-1} \actsOnMap c)(m_0 \isSemiActedUponBy n))\\
      &= \delta_0(n_0 \mapsto (g_{m_0, m}^{-1} \actsOnMap c)(n_0))\\
      &= \delta_0((g_{m_0, m}^{-1} \actsOnMap c)\restrictedTo_{N_0})\\
      &= \Delta_0(c)(m).
    \end{align*}
    In conclusion, $\Delta = \Delta_0$.
  \end{proof}

  \begin{example}[Left Shift Map]
  \label{example:lattice:the-left-shift}
    Let $M$ be the one-dimensional integer lattice $\Z$, let $G$ be the group $\setOf{\tau_t \from m \mapsto t + m \suchThat t \in \Z}$ of translations of $M$, and let $\actsOnPoint$ be the left group action of $G$ on $M$ by function application. Moreover, let $m_0$ be the origin $0$ and, for each point $m \in M$, let $g_{m_0, m}$ be the translation $\tau_m$. Furthermore, let $Q$ be the set $\setOf{0, 1}$, let $N$ be the set $\setOf{\tau_1}$, let $\delta$ be the map $Q^N \to Q$, $\ell \mapsto \ell(\tau_1)$, and let $G$ be identified with $G \modulo G_0$ by $\tau_t \mapsto \tau_t G_0$, where $G_0$ is the stabiliser $\setOf{\tau_0}$ of the origin $m_0$ under $\actsOnPoint$. The triple $\mathcal{M} = \ntuple{M, G, \actsOnPoint}$ is a principal left-ho\-mo\-ge\-neous space, the tuple $\mathcal{K} = \ntuple{m_0, \family{g_{m_0, m}}_{m \in M}}$ is a coordinate system for $\mathcal{M}$, the quadruple $\mathcal{C} = \ntuple{\ntuple{\mathcal{M}, \mathcal{K}}, Q, N, \delta}$ is a cellular automaton whose global transition function $\Delta$ is the \graffito{left shift map}\define{left shift map}\index{shift map!left} $Q^M \to Q^M$, $c \mapsto c(\blank + 1)$. Under the additional identification of $\Z$ with $G$ by $t \mapsto \tau_t$, the quadruple $\mathcal{C}$ is a traditional cellular automaton, more precisely, the action $\actsOnPoint$ is the addition $+$ on $\Z$, the neighbourhood $N$ is the set $\setOf{1}$, and the local transition function $\delta$ is the map $\ell \mapsto \ell(1)$.
  \end{example}

  \section{Invariance, Equivariance, Determination, and Composition of Global Transition Functions}
  \sectionmark{Invariance, Equivariance, Determination, and Composition}
  \label{section:invariances}

  \paragraph{Contents.} In \cref{theorem:independence-of-coordinate-system} we show how to turn a big-cellular automaton over one coordinate system into an automaton over another system that has the same global transition function. Conversely, in \cref{theorem:global-transition-function-determines-local-transition-function} we show that two big-cellular automata with the same global transition function are related as in \cref{theorem:independence-of-coordinate-system} except for superfluous neighbours. In \cref{corollary:independence-of-coordinate-system} we show that a global transition function does not depend on the choice of coordinates. In \cref{theorem:local-invariance-versus-global-equivariance} we show that a global transition function is $\actsOnMap_H$-e\-qui\-var\-i\-ant if and only if its local transition function is $\bullet_{H_0}$-invariant. In \cref{theorem:determination-of-cellular-automata-by-behaviour-at-origin} we show that a global transition function is determined by its behaviour in the origin. And in \cref{theorem:composition-of-cellular-automata} we show that the composition of two global transition functions is a global transition function.

  \subsection{Invariance Under Change of Coordinates of Global Transition Functions}

  \paragraph{Summary.} Conjugation by a group element $g$ is a bijection from a quotient set $G \modulo G_m$ to $G \modulo G_{g \actsOnPoint m}$ (see \cref{lemma:left-group-action-on-union-of-quotients-by-stabilisers}). Under the identifications of such quotient sets with $M$, all these conjugations together are in a sense the left group action $\actsOnPoint$ (see \cref{remark:action-on-union-of-quotients-is-in-a-sense-leftaction}). They are used to relate the right quotient set semi-action induced by one coordinate system to the semi-action induced by another system (see \cref{lemma:liberation-and-coordinate-systems}) and to turn a semi-cellular automaton over one coordinate system into an automaton in another system that has the same global transition function (see \cref{theorem:independence-of-coordinate-system}). It follows from the specifics of the latter construction that a global transition function does not depend on the choice of coordinates (see \cref{corollary:independence-of-coordinate-system}). And it follows that the set of global transition functions of big-cellular automata over one coordinate system is the same as the one over another system (see \cref{corollary:existence-in-other-coordinate-system}). 
  Conversely, if two big-cellular automata have identical global transition functions, then they are related as in \cref{theorem:independence-of-coordinate-system} except for superfluous neighbours (see \cref{theorem:global-transition-function-determines-local-transition-function}). It follows that two such automata over coordinate systems with the same origin are the same except for superfluous neighbours.


  \begin{lemma} 
  \label{lemma:left-group-action-on-union-of-quotients-by-stabilisers}
    Let $\actsOnPoint$ be a left group action of $G$ on $M$. The group $G$ acts on $\bigcup_{m \in M} G \modulo G_m$ on the left by
    \begin{align*}
      \conjugates \from G \times \bigcup_{m \in M} G \modulo G_m &\to \bigcup_{m \in M} G \modulo G_m, \mathnote{left group action $\conjugates$ of $G$ on $\bigcup_{m \in M} G \modulo G_m$}\index[symbols]{circle@$\conjugates$}\\
      (g, g' G_m) &\mapsto g g' G_m g^{-1} \quad (= g g' g^{-1} G_{g \actsOnPoint m}),
    \end{align*}
    such that, for each element $g \in G$ and each element $m \in M$, the map
    \begin{align*}
      (g \conjugates \blank)\restrictedTo_{G \modulo G_m \to G \modulo G_{g \actsOnPoint m}} \from G \modulo G_m &\to G \modulo G_{g \actsOnPoint m},\\
      g' G_m &\mapsto g \conjugates g' G_m,
    \end{align*}
    is bijective.
  \end{lemma}

  \begin{proof}
    First, let $g \in G$, let $m \in M$, let $g' G_m \in G \modulo G_m$, and let $m' = g \actsOnPoint m$. Then, according to \cref{lemma:stabiliser-versus-transporter}, we have $G_{m'} = g G_m g^{-1}$. Hence,
    \begin{align*}
      g g' G_m g^{-1}
      &= g g' (g^{-1} g) G_m g^{-1}\\
      &= (g g' g^{-1}) (g G_m g^{-1})\\
      &= g g' g^{-1} G_{m'}\\
      &\in G \modulo G_{m'}.
    \end{align*}
    In conclusion, the maps $\conjugates$ and $(g \conjugates \blank)\restrictedTo_{G \modulo G_m \to G \modulo G_{g \actsOnPoint m}}$ are well-defined.

    Secondly, let $g G_m \in \bigcup_{m \in M} G \modulo G_m$. Then, $e_G \conjugates g G_m = g G_m$. Moreover, for each $g' \in G$ and each $g'' \in G$,
    \begin{align*}
      g' g'' \conjugates g G_m
      &= g' g'' g G_m (g'')^{-1} (g')^{-1}\\
      &= g' \conjugates g'' g G_m (g'')^{-1}\\
      &= g' \conjugates (g'' \conjugates g G_m).
    \end{align*}
    In conclusion, the map $\conjugates$ is a left group action.

    Thirdly, for each $g \in G$ and each $m \in M$, the map $(g \conjugates \blank)\restrictedTo_{G \modulo G_m \to G \modulo G_{g \actsOnPoint m}}$ is bijective, because its inverse is $(g^{-1} \conjugates \blank)\restrictedTo_{G \modulo G_{g \actsOnPoint m} \to G \modulo G_m}$.
  \end{proof}

  \begin{remark}
  \label{remark:action-on-union-of-quotients-is-in-a-sense-leftaction}
    For each element $m_0 \in M$, let $\iota_{m_0}$ be the $(\actsOnPoint, \cdot)$-e\-qui\-var\-i\-ant bijection from \cref{lemma:iota-is-bijective-and-equivariant}. For each element $g \in G$, each element $m \in M$, and each element $m' \in M$,
    \begin{equation*}
      g \conjugates \iota_m(m')
      = g \conjugates G_{m, m'}
      = G_{g \actsOnPoint m, g \actsOnPoint m'}
      = \iota_{g \actsOnPoint m}(g \actsOnPoint m').
    \end{equation*}
    In this sense, the map $\conjugates$ is the left group action $\actsOnPoint$.
  \end{remark}

  \begin{remark}
  \label{remark:circ-of-cdot-invariant-set-is-invariant}
    Let $g$ be an element of $G$, let $m$ be an element of $M$, and let $N$ be a subset of $G \modulo G_m$ such that $G_m \cdot N \subseteq N$. Then, $G_{g \actsOnPoint m} \cdot (g \conjugates N) \subseteq g \conjugates N$. Indeed,
    \begin{align*}
      G_{g \actsOnPoint m} \cdot (g \conjugates N)
      &= (g G_m g^{-1}) \cdot g N g^{-1}\\
      &= g (G_m \cdot N) g^{-1}\\
      &\subseteq g N g^{-1}\\
      &= g \conjugates N. \qedhere
    \end{align*}
  \end{remark}

  The group actions $\conjugates$ and $\cdot$ commute in the sense given in

  \begin{remark}
  \label{remark:actions-circ-and-cdot-commute}
    Let $\actsOnPoint$ be a left group action of $G$ on $M$, let $m$ be an element of $M$, let $g$ and $g'$ be two elements of $G$, and let $g'' G_m$ be an element of $G \modulo G_m$. Then, $g \conjugates (g' \cdot g'' G_m) = g g' g^{-1} \cdot (g \conjugates g'' G_m)$. Indeed,
    \begin{align*}
      g \conjugates (g' \cdot g'' G_m) 
                                 &= g g' g'' g^{-1} G_{g \actsOnPoint m}\\
                                 &= g g' g^{-1} \cdot (g g'' g^{-1} G_{g \actsOnPoint m})\\
                                 &= g g' g^{-1} \cdot (g \conjugates g'' G_m). \qedhere
    \end{align*}
  \end{remark}

  \begin{lemma} 
  \label{lemma:liberation-and-coordinate-systems}
    Let $\mathcal{M} = \ntuple{M, G, \actsOnPoint}$ be a left-ho\-mo\-ge\-neous space, and let $\mathcal{K} = \ntuple{m_0, \family{g_{m_0, m}}_{m \in M}}$ and $\mathcal{K}' = \ntuple{m_0', \family{g_{m_0', m}'}_{m \in M}}$ be two coordinate systems for $\mathcal{M}$. The right quotient set semi-actions $\isSemiActedUponBy$ and $\isSemiActedUponBy'$ of $G \modulo G_0$ and $G \modulo G_0'$ on $M$ are \define{similar}\graffito{similar}, which means that, for each $\mathcal{K}$- and~$\mathcal{K}'$-big subgroup $H$ of $G$ and each element $h \in H$ such that $h \actsOnPoint m_0 = m_0'$,
    \begin{equation*}
      \ForEach m \in M \Exists h_0 \in H_0 \SuchThat \ForEach \mathfrak{g}' \in G \modulo G_0' \Holds
          m \isSemiActedUponBy' \mathfrak{g}' = m \isSemiActedUponBy h_0 \cdot (h^{-1} \conjugates \mathfrak{g}'). \qedhere
    \end{equation*}
  \end{lemma}

  \begin{proof} 
    Let $H$ be a $\mathcal{K}$- and~$\mathcal{K}'$-big subgroup of $G$, let $h \in H$ such that $h \actsOnPoint m_0 = m_0'$, and let $m \in M$. Put $h_0 = g_{m_0, m}^{-1} g_{m_0', m}' h$. Then, $h_0 \in H_0$ and $g_{m_0', m}' = g_{m_0, m} h_0 h^{-1}$. Furthermore, let $g G_0' \in G \modulo G_0'$. Then,
    \begin{align*}
      m \isSemiActedUponBy' g G_0'
      &= g_{m_0', m}' g (g_{m_0', m}')^{-1} \actsOnPoint m\\
      &= g_{m_0, m} h_0 h^{-1} g h h_0^{-1} g_{m_0, m}^{-1} \actsOnPoint m\\
      &= m \isSemiActedUponBy h_0 h^{-1} g h h_0^{-1} G_0.
    \end{align*}
    Thus, because $h_0^{-1} G_0 = G_0$ and $h G_0 h^{-1} = h \conjugates G_0 = G_0'$,
    \begin{align*}
      m \isSemiActedUponBy' g G_0'
      &= m \isSemiActedUponBy h_0 \cdot h^{-1} g h G_0\\
      &= m \isSemiActedUponBy h_0 \cdot h^{-1} g h G_0 h^{-1} h\\
      &= m \isSemiActedUponBy h_0 \cdot h^{-1} g G_0' h\\
      &= m \isSemiActedUponBy h_0 \cdot (h^{-1} \conjugates g G_0'). \qedhere
    \end{align*}
  \end{proof}

  \begin{lemma}
  \label{definition:left-group-action-on-local-transition-functions}
    Let $\actsOnPoint$ be a left group action of $G$ on $M$, let $Q$ be a set, and let
    \begin{equation*}
      D = \bigcup_{m \in M} \setOf{\delta \from Q^N \to Q \suchThat N \subseteq G \modulo G_m \text{ with } G_m \cdot N \subseteq N}.
    \end{equation*}
    The group $G$ acts on $D$ on the left by
    \begin{align*}
      \otimes \from G \times D &\to D, \mathnote{left group action $\otimes$ of $G$ on $D$}\\
      (g, \delta \from Q^N \to Q)
      &\mapsto \left[
                 \begin{aligned}
                   Q^{g \conjugates N} &\to Q,\\
                   \ell' &\mapsto \delta(n \mapsto \ell'(g \conjugates n)). 
                 \end{aligned}
               \right] \qedhere
    \end{align*}
  \end{lemma}

  \begin{proof}
    According to \cref{remark:circ-of-cdot-invariant-set-is-invariant}, the map $\otimes$ is well-defined. Let $(\delta \from Q^N \to Q) \in D$. Then, $e_G \otimes \delta = \delta$. And, for each $g \in G$, each $g' \in G$, and each $\ell'' \in Q^{g g' \conjugates N}$,
    \begin{align*}
      (g g' \otimes \delta)(\ell'')
      &= \delta(n \mapsto \ell''(g g' \conjugates n))\\
      &= \delta(n \mapsto \ell''(g \conjugates (g' \conjugates n)))\\
      &= \delta(n \mapsto [n' \mapsto \ell''(g \conjugates n')](g' \conjugates n))\\
      &= (g' \otimes \delta)(n' \mapsto \ell''(g \conjugates n'))\\
      &= (g \otimes (g' \otimes \delta))(\ell'').
    \end{align*}
    In conclusion, $\otimes$ is a left group action.
  \end{proof}


  \begin{theorem}
  \label{theorem:independence-of-coordinate-system}
    In the situation of \cref{lemma:liberation-and-coordinate-systems}, let $H$ be a $\mathcal{K}$- and~$\mathcal{K}'$-big subgroup of $G$, let $\mathcal{C} = \ntuple{\ntuple{\mathcal{M}, \mathcal{K}}, Q, N, \delta}$ be a semi-cellular automaton such that $\delta$ is $\bullet_{H_0}$-invariant, let $N'$ be the set $h \conjugates N$, and let $\delta'$ be the map $h \otimes \delta$. The quadruple $\mathcal{C}' = \ntuple{\ntuple{\mathcal{M}, \mathcal{K}'}, Q, N', \delta'}$ is a semi-cellular automaton whose local transition function is $\bullet_{H_0'}$-invariant and whose global transition function is identical to the one of $\mathcal{C}$.
  \end{theorem}

  \begin{proof}
    We have $N' \subseteq G \modulo G_0'$ and, according to \cref{remark:circ-of-cdot-invariant-set-is-invariant}, we have $G_0' \cdot N' \subseteq N'$. In conclusion, $\mathcal{C}'$ is a semi-cellular automaton.

    Moreover, for each $h_0' \in H_0'$ and each $\ell' \in Q^{N'}$, according to \cref{remark:actions-circ-and-cdot-commute} and because $\delta$ is $\bullet_{H_0}$-invariant,
    \begin{align*}
      \delta'(h_0' \bullet \ell')
      &= \delta(n \mapsto (h_0' \bullet \ell')(h \conjugates n))\\
      &= \delta(n \mapsto \ell'((h_0')^{-1} \cdot (h \conjugates n)))\\
      &= \delta(n \mapsto \ell'(h \conjugates (h^{-1} (h_0')^{-1} h \cdot n)))\\
      &= \delta((h^{-1} h_0' h) \bullet [n \mapsto \ell'(h \conjugates n)])\\
      &= \delta(n \mapsto \ell'(h \conjugates n))\\
      &= \delta'(\ell').
    \end{align*}
    In conclusion, $\delta'$ is $\bullet_{H_0'}$-invariant.

    Furthermore, let $c \in Q^M$ and let $m \in M$. According to \cref{lemma:liberation-and-coordinate-systems}, there is an $h_0 \in H_0$ such that
    \begin{equation*}
      \ForEach n' \in N' \Holds m \isSemiActedUponBy' n' = m \isSemiActedUponBy h_0 \cdot (h^{-1} \conjugates n').
    \end{equation*}
    Therefore, because $\delta$ is $\bullet_{H_0}$-invariant,
    \begin{align*}
      \Delta'(c)(m) &= \delta'(n' \mapsto c(m \isSemiActedUponBy' n'))\\
                    &= \delta'(n' \mapsto c(m \isSemiActedUponBy h_0 \cdot (h^{-1} \conjugates n')))\\
                    &= \delta(n \mapsto c(m \isSemiActedUponBy h_0 \cdot (h^{-1} \conjugates (h \conjugates n))))\\
                    &= \delta(n \mapsto c(m \isSemiActedUponBy h_0 \cdot n))\\
                    &= \delta(h_0^{-1} \bullet [n \mapsto c(m \isSemiActedUponBy n)])\\ 
                    &= \delta(n \mapsto c(m \isSemiActedUponBy n))\\
                    &= \Delta(c)(m).
    \end{align*}
    In conclusion, $\Delta' = \Delta$.
  \end{proof}

  \begin{corollary}
  \label{corollary:independence-of-coordinate-system}
    In the situation of \cref{lemma:liberation-and-coordinate-systems}, let $m_0 = m_0'$, let $H$ be a $\mathcal{K}$- and~$\mathcal{K}'$-big subgroup of $G$, and let $\mathcal{C} = \ntuple{\ntuple{\mathcal{M}, \mathcal{K}}, Q, N, \delta}$ be a semi-cellular automaton such that $\delta$ is $\bullet_{H_0}$-invariant. The global transition function of the semi-cellular automaton $\ntuple{\ntuple{\mathcal{M}, \mathcal{K}'}, Q, N, \delta}$ is identical to the one of $\mathcal{C}$.
  \end{corollary}

  \begin{proof}
    This is a direct consequence of \cref{theorem:independence-of-coordinate-system} with $h = e_G$.
  \end{proof}

  \begin{corollary}
  \label{corollary:existence-in-other-coordinate-system}
    Let $\ntuple{\mathcal{M}, \mathcal{K}} = \ntuple{\ntuple{M, G, \actsOnPoint}, \ntuple{m_0, \family{g_{m_0, m}}_{m \in M}}}$ be a cell space and let $\mathcal{C}$ be a cellular automaton over $\ntuple{\mathcal{M}, \mathcal{K}}$. For each coordinate system $\mathcal{K}' = \ntuple{m_0', \family{g_{m_0', m}'}_{m \in M}}$ for $\mathcal{M}$, there is a cellular automaton over $\ntuple{\mathcal{M}, \mathcal{K}'}$ whose global transition function is identical to the one of $\mathcal{C}$.
  \end{corollary}

  \begin{proof}
    This is a direct consequence of \cref{theorem:independence-of-coordinate-system} with $H = G$ and $h = g_{m_0, m_0'}$.
  \end{proof}

  \begin{example}[Lattice]
  \label{example:lattice:homogeneous-space}
    Let $M$ be the one-dimensional integer lattice $\Z$, let $T$ be the group $\setOf{\tau_t \from m \mapsto t + m \suchThat t \in \Z}$ of translations of $M$, let $R$ be the set $\setOf{\varrho_t \from m \mapsto t - m \suchThat t \in \Z}$ of reflections of $M$, let $G$ be the group $T \cup R$ of translations and reflections of $M$, and let $\actsOnPoint$ be the left group action of $G$ on $M$ by function application. The triple $\mathcal{M} = \ntuple{M, G, \actsOnPoint}$ is a left-ho\-mo\-ge\-neous space and the stabiliser $G_0$ of the origin $0$ under $\actsOnPoint$ is the group $\setOf{\tau_0, \varrho_0}$.
  \end{example}

  \begin{example}[Logical Or]
  \label{example:lattice:automaton:or}
    In the situation of \cref{example:lattice:homogeneous-space}, let $\mathcal{K}$ be a coordinate system for $\mathcal{M}$, let $Q$ be the set $\setOf{0, 1}$, let $N$ be the set $\setOf{-1, 1}$, let $\delta$ be the map $Q^N \to Q$, $\ell \mapsto \ell(-1) \lor \ell(1)$, and let $M$ be identified with $G \modulo G_0$ by $\iota \givenBy m \mapsto G_{0, m}$. The quadruple $\mathcal{C} = \ntuple{\ntuple{\mathcal{M}, \mathcal{K}}, Q, N, \delta}$ is a cellular automaton whose global transition function $\Delta$ is the map $Q^M \to Q^M$, $c \mapsto c(\blank - 1) \lor c(\blank + 1)$ (see \cref{figure:or}). Note that, because the binary operator $\lor$ is commutative, the local transition function $\delta$ is $\bullet$-invariant and the global transition function $\Delta$ does not depend on the coordinate system.
    \begin{figure}
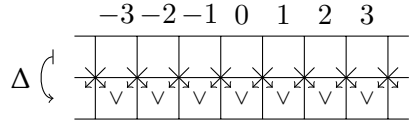

      \myfloatalign
      \figureOr
      \caption{The first row depicts the part of a global configuration that corresponds to the cells $\setOf{-3, -2, \dotsc, 3}$ and the second row the same part of the image of that global configuration under $\Delta$. And, the arrows from cells in the first to cells in the second row together with the disjunction symbol depict which states are combined disjunctively to yield the new states.} 
      \label{figure:or}
    \end{figure}
  \end{example}

  \begin{counterexample}[Peculiar Shift Maps]
  \label{example:lattice:automata:shifts}
    In \cref{theorem:independence-of-coordinate-system}, if the assumption that the subgroup $H$ of $G$ is $\mathcal{K}$-big or the one that the local transition function $\delta$ is $\bullet_{H_0}$-invariant does not hold, then the global transition function of $\mathcal{C}'$ may be different from the one of $\mathcal{C}$, which is illustrated by the following examples.

    In the situation of \cref{example:lattice:homogeneous-space}, let $Q$ be the set $\setOf{0, 1}$, let $N$ be the set $\setOf{-1, 1}$, let $\delta$ be the map $Q^N \to Q$, $\ell \mapsto \ell(1)$, and let $M$ be identified with $G \modulo G_0$ by $\iota \givenBy m \mapsto G_{0, m}$. Because $\delta$ distinguishes between left and right, more precisely, because $\delta(\varrho_0 \bullet \blank) \neq \delta$, the map $\delta$ is not $\bullet$-invariant. And, the global transition functions of semi-cellular automata over $\mathcal{M}$ with local transition function $\delta$ depend on the choice of coordinates, which is illustrated by the following examples.
    \begin{aenumerate}
      \item\label{item:lattice:automata:shifts:left-and-right}
            The tuples $\mathcal{K} = \ntuple{0, \family{\tau_m}_{m \in M}}$ and $\mathcal{K}' = \ntuple{0, \family{\varrho_m}_{m \in M}}$ are two coordinate systems for $\mathcal{M}$. For each cell $m \in M$ and each element $\mathfrak{g} \in G \modulo G_0 \simeq M$, we have $m \isSemiActedUponBy \mathfrak{g} = m + \mathfrak{g}$ and $m \isSemiActedUponBy' \mathfrak{g} = m - \mathfrak{g}$, in particular, the (actual) neighbour of $m$ that corresponds to the (relative) right neighbour $1$ in $\mathcal{K}$ is the cell $m + 1$, which lies to the right of $m$, and in $\mathcal{K}'$ it is the cell $m - 1$, which lies to the left. The reason is that the coordinates $\family{\tau_m}_{m \in M}$ maintain the meanings of left and right, whereas the coordinates $\family{\varrho_m}_{m \in M}$ reverse them. 

            The quadruples $\mathcal{C} = \ntuple{\ntuple{\mathcal{M}, \mathcal{K}}, Q, N, \delta}$ and $\mathcal{C}' = \ntuple{\ntuple{\mathcal{M}, \mathcal{K}'}, Q, N, \delta}$ are two semi-cellular automata whose global transition functions $\Delta$ and $\Delta'$ are the left shift map $c \mapsto c(\blank + 1)$ (see \cref{figure:left-shift}) and the right shift map $c \mapsto c(\blank - 1)$ (see \cref{figure:right-shift}). The reason is that $\delta$ depends on the meanings of left and right, which are maintained by $\isSemiActedUponBy$ but reversed by $\isSemiActedUponBy'$.

            Note that $\mathcal{K}'$ is actually not a coordinate system, because, by definition, the coordinate of the origin must be the identity map $\tau_0$. That requirement though could be discarded with minor changes to some statements. It was merely made for convenience.
      \item\label{item:lattice:automata:shifts:odd-even} 
            The tuple $\mathcal{K}'' = \ntuple{0, \family{\tau_m}_{m \in 2 M} \times \family{\varrho_m}_{m \in 2 M + 1}}$ is a coordinate system for $\mathcal{M}$. The right quotient set semi-action of $G \modulo G_0 \simeq M$ on $M$ induced by $\ntuple{\mathcal{M}, \mathcal{K}''}$ is the map
            \begin{align*}
              \isSemiActedUponBy'' \from M \times G \modulo G_0 &\to M,\\
              (m, \mathfrak{g}) &\mapsto \begin{dcases*}
                m + \mathfrak{g}, &if $m$ is even,\\
                m - \mathfrak{g}, &if $m$ is odd.
              \end{dcases*}
            \end{align*} 
            In particular, the (actual) neighbour of a cell that corresponds to the (relative) neighbour $1$ is the cell to its right, if it is even, and the one to its left, if it is odd. The quadruple $\mathcal{C}'' = \ntuple{\ntuple{\mathcal{M}, \mathcal{K}''}, Q, N, \delta}$ is a semi-cellular automaton whose global transition function is the map 
            \begin{align*}
              \Delta'' \from Q^M &\to Q^M,\\
              c &\mapsto [m \mapsto \left\{
                                      \begin{aligned}
                                        &c(m + 1), &&\text{if $m$ is even,}\\
                                        &c(m - 1), &&\text{if $m$ is odd,}
                                      \end{aligned}
                                    \right\}].
            \end{align*}
            In one step, each even cell exchanges states with the odd cell to its right (see \cref{figure:odd-even-shift}).
      \item\label{item:lattice:automata:shifts:strange}
            The tuple $\mathcal{K}''' = \ntuple{0, \family{\tau_m}_{m \in M \smallsetminus \setOf{1}} \times \family{\varrho_m}_{m \in \setOf{1}}}$ is a coordinate system for $\mathcal{M}$. The right quotient set semi-action of $G \modulo G_0 \simeq M$ on $M$ induced by $\ntuple{\mathcal{M}, \mathcal{K}'''}$ is the map
            \begin{align*}
              \isSemiActedUponBy''' \from M \times G \modulo G_0 &\to M,\\
              (m, \mathfrak{g}) &\mapsto \begin{dcases*}
                m + \mathfrak{g}, &if $m \neq 1$,\\
                m - \mathfrak{g}, &if $m = 1$.
              \end{dcases*}
            \end{align*}
            The quadruple $\mathcal{C}''' = \ntuple{\ntuple{\mathcal{M}, \mathcal{K}'''}, Q, N, \delta}$ is a semi-cellular automaton whose global transition function is the map
            \begin{align*}
              \Delta''' \from Q^M &\to Q^M,\\
              c &\mapsto [m \mapsto \left\{
                                      \begin{aligned}
                                        &c(m + 1), &&\text{if $m \neq 1$,}\\
                                        &c(m - 1), &&\text{if $m = 1$,}
                                      \end{aligned}
                                    \right\}].
            \end{align*}
            In one step, the cells $0$ and $1$ exchange states and each other cell takes the state from the cell to its right (see \cref{figure:strange-shift}).
    \end{aenumerate} 
    The semi-cellular automata $\mathcal{C}$, $\mathcal{C}'$, $\mathcal{C}''$, and $\mathcal{C}'''$ have different global transition functions, although they are equal except for their coordinates, in particular, they are related as in the construction in \cref{theorem:independence-of-coordinate-system} with $h = e_G$. However, the group $G$ is the only subgroup of $G$ that is big with respect to each pair of the coordinate systems $\mathcal{K}$, $\mathcal{K}'$, $\mathcal{K}''$, and $\mathcal{K}'''$, but the local transition function $\delta$ is not $\bullet_{G_0}$-invariant. And, the subgroups $H$ of $T$ are the only subgroups of $G$ such that the local transition function $\delta$ is $\bullet_{H_0}$-invariant, but those subgroups of $G$ are only $\mathcal{K}$-big and neither big with respect to $\mathcal{K}'$, $\mathcal{K}''$, nor $\mathcal{K}'''$.
    \begin{figure}
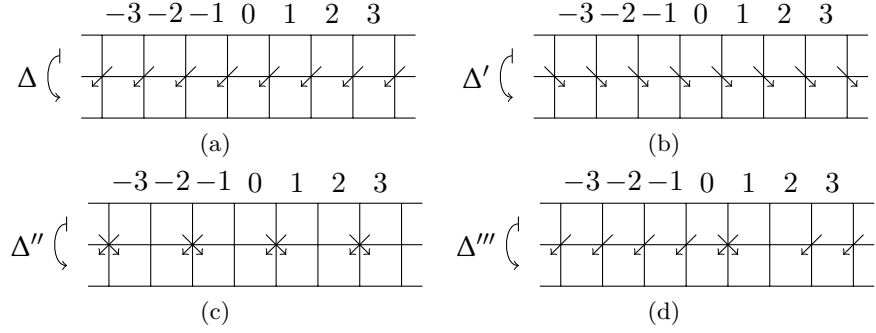

      \myfloatalign
      \figureShifts
      \caption{In each subfigure, the first row depicts the part of a global configuration that corresponds to the cells $\setOf{-3, -2, \dotsc, 3}$ and the second row the same part of the image of that global configuration under $\Delta$, $\Delta'$, $\Delta''$, or $\Delta'''$. And, the arrows from cells in the first to cells in the second row depict the flow of states.}
    \end{figure}
  \end{counterexample}

  \begin{theorem} 
  \label{theorem:global-transition-function-determines-local-transition-function} 
    In the situation of \cref{lemma:liberation-and-coordinate-systems}, let $H$ be a $\mathcal{K}$- and~$\mathcal{K}'$-big subgroup of $G$, let $\mathcal{C} = \ntuple{\ntuple{\mathcal{M}, \mathcal{K}}, Q, N, \delta}$ and $\mathcal{C}' = \ntuple{\ntuple{\mathcal{M}, \mathcal{K}'}, Q, N', \delta'}$ be two semi-cellular automata such that $\delta$ and $\delta'$ are $\bullet_{H_0}$-invariant and $\Delta$ and $\Delta'$ are identical, let $N_* = N \cap (h^{-1} \conjugates N')$, let $N'_* = (h \conjugates N) \cap N'$, let
    \begin{align*}
      \delta_* \from Q^{N_*} &\to Q,\\
      \ell_* &\mapsto \delta(\ell), \text{ where $\ell \in Q^N$ such that $\ell\restrictedTo_{N_*} = \ell_*$},
    \end{align*}
    and let
    \begin{align*}
      \delta'_* \from Q^{N'_*} &\to Q,\\
      \ell'_* &\mapsto \delta'(\ell'), \text{ where $\ell' \in Q^{N'}$ such that $\ell'\restrictedTo_{N'_*} = \ell'_*$}.
    \end{align*}
    The quadruples $\mathcal{C}_* = \ntuple{\ntuple{\mathcal{M}, \mathcal{K}}, Q, N_*, \delta_*}$ and $\mathcal{C}'_* = \ntuple{\ntuple{\mathcal{M}, \mathcal{K}'}, Q, N'_*, \delta'_*}$ are two semi-cellular automata such that $\delta_*$ and $\delta'_*$ are $\bullet_{H_0}$-invariant, and $\Delta_* = \Delta = \Delta' = \Delta'_*$. Moreover, $N'_* = h \conjugates N_*$, $\delta'_* = h \otimes \delta_*$, and 
    \begin{multline*}
      \ForEach \ell_* \in Q^{N_*} \ForEach \ell \in Q^N \ForEach \ell' \in Q^{N'} \Holds\\
      \begin{aligned}
        \big(
          &\ell\restrictedTo_{N_*} = \ell_* \land \ell'\restrictedTo_{N'_*} = \ell_*(h^{-1} \conjugates \blank)\\
          &{}\implies \delta(\ell) = \delta_*(\ell_*) = \delta'_*(\ell_*(h^{-1} \conjugates \blank)) = \delta'(\ell')
        \big). \qedhere 
      \end{aligned}
    \end{multline*}
  \end{theorem}

  \begin{proof}
    First, we have $h^{-1} \conjugates N' \subseteq G \modulo G_{h^{-1} \actsOnPoint m_0'} = G \modulo G_0$ and, according to \cref{remark:circ-of-cdot-invariant-set-is-invariant}, we have $G_0 \cdot (h^{-1} \conjugates N') \subseteq h^{-1} \conjugates N'$. Therefore, $G_0 \cdot N_* = (G_0 \cdot N) \cap (G_0 \cdot (h^{-1} \conjugates N')) \subseteq N \cap (h^{-1} \conjugates N') = N_*$. Analogously, $G_0 \cdot N'_* \subseteq N'_*$. Moreover, $h \conjugates N_* = (h \conjugates N) \cap (h \conjugates (h^{-1} \conjugates N')) = (h \conjugates N) \cap N' = N'_*$.

    Secondly, let $\ell_* \in Q^{N_*}$, let $\ell \in Q^N$, and let $\ell' \in Q^{N'}$ such that $\ell\restrictedTo_{N_*} = \ell_*$ and $\ell'\restrictedTo_{N'_*} = \ell_*(h^{-1} \conjugates \blank)$. Then, because $m_0 \isSemiActedUponBy \blank$ and $m_0 \isSemiActedUponBy' \blank$ are injective, there are $c$, $c' \in Q^M$ such that $c(m_0 \isSemiActedUponBy \blank)\restrictedTo_N = \ell$ and $c'(m_0 \isSemiActedUponBy' \blank)\restrictedTo_{N'} = \ell'$ (see \cref{remark:local-configuration-induces-global-configuration}).
    And, according to \cref{lemma:liberation-and-coordinate-systems}, there is an $h_0 \in H_0$ such that
    \begin{equation*}
      \ForEach \mathfrak{g}' \in G \modulo G_0' \Holds m_0 \isSemiActedUponBy' \mathfrak{g}' = m_0 \isSemiActedUponBy h_0 \cdot (h^{-1} \conjugates \mathfrak{g}').
    \end{equation*}
    Thus, because $h_0 \cdot (h^{-1} \conjugates N') = h^{-1} \conjugates N'$,
    \begin{equation*}
      m_0 \isSemiActedUponBy' N' = m_0 \isSemiActedUponBy h_0 \cdot (h^{-1} \conjugates N') = m_0 \isSemiActedUponBy h^{-1} \conjugates N'.
    \end{equation*}
    And, because $h_0 \cdot N_* = N_*$,
    \begin{equation*}
      m_0 \isSemiActedUponBy' N'_* = m_0 \isSemiActedUponBy h_0 \cdot (h^{-1} \conjugates N'_*) = m_0 \isSemiActedUponBy h_0 \cdot N_* = m_0 \isSemiActedUponBy N_*.
    \end{equation*}
    Hence, because $m_0 \isSemiActedUponBy \blank$ is injective,
    \begin{align*}
      (m_0 \isSemiActedUponBy N) \cap (m_0 \isSemiActedUponBy' N')
      &= (m_0 \isSemiActedUponBy N) \cap (m_0 \isSemiActedUponBy h^{-1} \conjugates N')\\
      &= m_0 \isSemiActedUponBy (N \cap (h^{-1} \conjugates N'))\\
      &= m_0 \isSemiActedUponBy N_*\\
      &= m_0 \isSemiActedUponBy' N'_*.
    \end{align*}
    Moreover, for each $n_* \in N_*$, because $m_0 \isSemiActedUponBy h_0 \cdot n_* = h_0 \actsOnPoint (m_0 \isSemiActedUponBy n_*)$ (see \cref{lemma:right-semi-action-similar-to-leftaction-in-m0}),
    \begin{align*}
      c(m_0 \isSemiActedUponBy n_*)
      &= \ell(n_*)\\
      &= \ell_*(n_*)\\
      &= \ell_*(h^{-1} \conjugates (h \conjugates n_*))\\
      &= \ell'(h \conjugates n_*)\\
      &= c'(m_0 \isSemiActedUponBy' h \conjugates n_*)\\
      &= c'(m_0 \isSemiActedUponBy h_0 \cdot n_*)\\
      &= (h_0^{-1} \actsOnMap c')(m_0 \isSemiActedUponBy n_*).
    \end{align*}
    Hence, because $m_0 \isSemiActedUponBy N_* = m_0 \isSemiActedUponBy' N'_*$, we have $c\restrictedTo_{m_0 \isSemiActedUponBy N_*} = (h_0^{-1} \actsOnMap c')\restrictedTo_{m_0 \isSemiActedUponBy' N'_*}$. Therefore, because $(m_0 \isSemiActedUponBy N) \cap (m_0 \isSemiActedUponBy' N') = m_0 \isSemiActedUponBy' N'_*$, there is a $c'' \in Q^M$ such that $c''\restrictedTo_{m_0 \isSemiActedUponBy N} = c\restrictedTo_{m_0 \isSemiActedUponBy N}$ and $c''\restrictedTo_{m_0 \isSemiActedUponBy' N'} = (h_0^{-1} \actsOnMap c')\restrictedTo_{m_0 \isSemiActedUponBy' N'}$. Thus, because $\Delta = \Delta'$, and $(h_0^{-1} \actsOnMap c')(m_0 \isSemiActedUponBy \blank)\restrictedTo_N = h_0^{-1} \bullet \ell'$, and $\delta'$ is $\bullet_{H_0}$-invariant, 
    \begin{align*}
      \delta(\ell)
      &= \Delta(c)(m_0)\\
      &= \Delta(c'')(m_0)\\
      &= \Delta'(c'')(m_0)\\
      &= \Delta'(h_0^{-1} \actsOnMap c')(m_0)\\
      &= \delta'(h_0^{-1} \bullet \ell')\\
      &= \delta'(\ell').
    \end{align*}
    In conclusion, $\delta(\ell) = \delta(\ell')$.

    Thirdly, it follows that, for each $\ell_* \in Q^{N_*}$, each $\ell_1 \in Q^N$, and each $\ell_2 \in Q^N$ such that $\ell_1\restrictedTo_{N_*} = \ell_*$ and $\ell_2\restrictedTo_{N_*} = \ell_*$, we have $\delta(\ell_1) = \delta'(\ell') = \delta(\ell_2)$, where $\ell' \in Q^{N'}$ such that $\ell'\restrictedTo_{N'_*} = \ell_*(h^{-1} \conjugates \blank)$. In conclusion, $\delta_*$ is well-defined and, analogously, $\delta'_*$ is well-defined.

    Fourthly, let $\ell_* \in Q^{N_*}$. Then, there are $\ell \in Q^N$ and $\ell' \in Q^{N'}$ such that $\ell\restrictedTo_{N_*} = \ell_*$ and $\ell'\restrictedTo_{N'_*} = \ell_*(h^{-1} \conjugates \blank)$. Hence, because $\delta'(\ell') = \delta(\ell)$,
    \begin{equation*}
      (h^{-1} \otimes \delta'_*)(\ell_*)
      = \delta'_*(\ell_*(h^{-1} \conjugates \blank))
      = \delta'(\ell')
      = \delta(\ell)
      = \delta_*(\ell_*).
    \end{equation*}
    In conclusion, $h^{-1} \otimes \delta'_* = \delta_*$ and hence $\delta'_* = h \otimes \delta_*$.

    Fifthly, let $h_0 \in H_0$ and let $\ell_* \in N_*$. Then, there is an $\ell \in Q^N$ such that $\ell\restrictedTo_{N_*} = \ell_*$, in particular, $(h_0 \bullet \ell)\restrictedTo_{N_*} = h_0 \bullet \ell_*$. Hence, because $\delta$ is $\bullet_{H_0}$-invariant,
    \begin{equation*}
      \delta_*(h_0 \bullet \ell_*)
      = \delta(h_0 \bullet \ell)
      = \delta(\ell)
      = \delta_*(\ell_*).
    \end{equation*}
    In conclusion, $\delta_*$ is $\bullet$-invariant and, analogously, $\delta'_*$ is $\bullet$-invariant.

    Lastly, for each $c \in Q^M$ and each $m \in M$, by the definition of $\delta_*$,
    \begin{align*}
      \Delta_*(c)(m) &= \delta_*(n_* \mapsto c(m \isSemiActedUponBy n_*))\\
                     &= \delta(n \mapsto c(m \isSemiActedUponBy n))\\
                     &= \Delta(c)(m).
    \end{align*}
    Therefore, $\Delta_* = \Delta$ and, analogously, $\Delta'_* = \Delta'$. In conclusion, $\Delta_* = \Delta = \Delta' = \Delta'_*$.
  \end{proof}

  \begin{corollary}
  \label{corollary:global-transition-function-determines-local-transition-function}
    In the situation of \cref{lemma:liberation-and-coordinate-systems}, let $m_0 = m_0'$, let $H$ be a $\mathcal{K}$- and~$\mathcal{K}'$-big subgroup of $G$, let $\mathcal{C} = \ntuple{\ntuple{\mathcal{M}, \mathcal{K}}, Q, N, \delta}$ and $\mathcal{C}' = \ntuple{\ntuple{\mathcal{M}, \mathcal{K}'}, Q, N', \delta'}$ be two semi-cellular automata such that $\delta$ and $\delta'$ are $\bullet_{H_0}$-invariant and $\Delta$ and $\Delta'$ are identical, let $N_* = N \cap N'$, and let
    \begin{align*}
      \delta_* \from Q^{N_*} &\to Q,\\
      \ell_* &\mapsto \delta(\ell), \text{ where $\ell \in Q^N$ such that $\ell\restrictedTo_{N_*} = \ell_*$},\\
      (\ell_* &\mapsto \delta'(\ell'), \text{ where $\ell' \in Q^{N'}$ such that $\ell'\restrictedTo_{N_*} = \ell_*$}).
    \end{align*}
    The quadruples $\mathcal{C}_* = \ntuple{\ntuple{\mathcal{M}, \mathcal{K}}, Q, N_*, \delta_*}$ and $\mathcal{C}'_* = \ntuple{\ntuple{\mathcal{M}, \mathcal{K}'}, Q, N_*, \delta_*}$ are two semi-cellular automata such that $\delta_*$ is $\bullet_{H_0}$-invariant and $\Delta_* = \Delta = \Delta' = \Delta'_*$. Moreover,
    \begin{multline*}
      \ForEach \ell_* \in Q^{N_*} \ForEach \ell \in Q^N \ForEach \ell' \in Q^{N'} \Holds\\
          \big(\ell\restrictedTo_{N_*} = \ell_* = \ell'\restrictedTo_{N_*}
              \implies \delta(\ell) = \delta_*(\ell_*) = \delta'(\ell')\big). \qedhere 
    \end{multline*}
  \end{corollary}

  \begin{proof}
    This is a direct consequence of \cref{theorem:global-transition-function-determines-local-transition-function} with $h = e_G$.
  \end{proof}

  \begin{counterexample}[Peculiar Shift Maps]
    In \cref{theorem:global-transition-function-determines-local-transition-function}, if the assumption that the subgroup $H$ of $G$ is $\mathcal{K}$- and~$\mathcal{K}'$-big does not hold, then the construction may not work, which is illustrated by the following examples.

    In the situation of \cref{item:lattice:automata:shifts:left-and-right} of \cref{example:lattice:automata:shifts}, let $\delta'$ be the map $Q^N \to Q$, $\ell \mapsto \ell(-1)$ and let $\mathcal{C}''$ be the semi-cellular automaton $\ntuple{\ntuple{\mathcal{M}, \mathcal{K}'}, Q, N, \delta'}$. The local transition functions $\delta$ and $\delta'$ are $\bullet_{T_0}$-invariant and the global transition functions $\Delta$ and $\Delta''$ are both the left shift map, but the subgroup $T$ of $G$ is only $\mathcal{K}$-big but not~$\mathcal{K}'$-big. The construction in \cref{theorem:global-transition-function-determines-local-transition-function} applied to $\mathcal{C}$ and $\mathcal{C}''$ with $h = e_G$ yields $\delta_* = \delta$ and $\delta'_* = \delta'$. However, contrary to the statement in the theorem, we have $\delta_* \neq e_G \otimes \delta'_*$. This could be the case because $\mathcal{K}'$ is actually not a coordinate system, which is why we also give the following counterexample.

    In the situation of \cref{item:lattice:automata:shifts:odd-even} of \cref{example:lattice:automata:shifts}, let $N'$ be the set $\setOf{0, 2}$, let $\delta'$ be the map $Q^{N'} \to Q$, $\ell \mapsto \ell(0)$, and let $M$ be identified with $G \modulo G_1$ by $m \mapsto G_{1, m}$ (we have identified $M$ twice now, but it will be clear from the context which identification applies). The tuple $\mathcal{K}''' = \ntuple{1, \family{\tau_{m - 1}}_{m \in 2 M + 1} \times \family{\varrho_{m + 1}}_{m \in 2 M}}$ is a coordinate system for $\mathcal{M}$. The right quotient set semi-action of $G \modulo G_1 \simeq M$ on $M$ induced by $\ntuple{\mathcal{M}, \mathcal{K}'''}$ is the map
    \begin{align*}
      \isSemiActedUponBy''' \from M \times G \modulo G_1 &\to M,\\
      (m, \mathfrak{g}) &\mapsto
          \begin{dcases*}
            m - 1 + \mathfrak{g}, &if $m$ is odd,\\
            m + 1 - \mathfrak{g}, &if $m$ is even.
          \end{dcases*}
    \end{align*}
    The quadruple $\mathcal{C}''' = \ntuple{\ntuple{\mathcal{M}, \mathcal{K}'''}, Q, N, \delta'}$ is a semi-cellular automaton whose global transition function is the map
    \begin{align*}
      \Delta''' \from Q^M &\to Q^M,\\
      c &\mapsto [m \mapsto \left\{
                              \begin{aligned}
                                &c(m - 1), &&\text{if $m$ is odd,}\\
                                &c(m + 1), &&\text{if $m$ is even,}
                              \end{aligned}
                            \right\}].
    \end{align*}
    The local transition functions $\delta$ and $\delta'$ are $\bullet_{T_0}$-invariant and the global transition functions $\Delta''$ and $\Delta'''$ are identical, but the subgroup $T$ of $G$ is neither $\mathcal{K}''$- nor~$\mathcal{K}'''$-big. The construction in \cref{theorem:global-transition-function-determines-local-transition-function} applied to $\mathcal{C}''$ and $\mathcal{C}'''$ with $h = \tau_1$ yields $\delta_* = \delta$ and $\delta'_* = \delta'$. However, because $(\tau_1 \otimes \delta)(\ell') = \delta(n \mapsto \ell'(\tau_1(n))) = \ell'(\tau_1(1)) = \ell'(2)$, for $\ell' \in Q^{N'}$ (where we used \cref{remark:action-on-union-of-quotients-is-in-a-sense-leftaction}), we have $\tau_1 \otimes \delta_* \neq \delta'_*$, contrary to the statement in the theorem.
  \end{counterexample}

  \subsection{Equivariance, Determination, and Composition of Global Transition Functions}

  \paragraph{Summary.} The local configuration that is observed by a translated cell is identical to a rotation of the one observed by the original cell in the reversely translated global configuration, symbolically, $c((g \actsOnPoint m) \isSemiActedUponBy \blank) = g_0 \bullet ((g^{-1} \actsOnMap c)(m \isSemiActedUponBy \blank))$ (see \cref{lemma:semi-commutativity-induced-left-action-and-bullet}). Hence, if the local transition function is invariant under rotations, then the global transition function is equivariant under translations (see \cref{theorem:local-invariance-versus-global-equivariance}).

  Local configurations can be embedded in global configurations and, through this embedding, the local transition function is identical to the global transition function evaluated at the origin, symbolically, $\delta(\ell) = \Delta(\bar{\ell})(m_0)$. And, rotations of local configurations translate to rotations of global configurations about the origin, symbolically, $\delta(h_0 \bullet \ell) = \Delta(h_0 \actsOnMap \bar{\ell})(m_0)$. Hence, if the global transition function is equivariant under rotations, then the local transition function is invariant under rotations (see \cref{theorem:local-invariance-versus-global-equivariance}). 

  \begin{lemma}
  \label{lemma:local-invariance-versus-global-equivariance:global-to-local}
    Let $\mathcal{R} = \ntuple{\mathcal{M}, \mathcal{K}} = \ntuple{\ntuple{M, G, \actsOnPoint}, \ntuple{m_0, \family{g_{m_0, m}}_{m \in M}}}$ be a cell space, let $\mathcal{C} = \ntuple{\mathcal{R}, Q, N, \delta}$ be a semi-cellular automaton, let $H$ be a subgroup of $G$, and let $\Delta_0$ be an $\actsOnMap_H$-e\-qui\-var\-i\-ant map from $Q^M$ to $Q^M$ such that
    \begin{align}
    \label{equation:local-invariance-versus-global-equivariance:global-transition-function-at-origins}
      \ForEach c \in Q^M \Holds \Delta_0(c)(m_0) = \delta(n \mapsto c(m_0 \isSemiActedUponBy n)).
    \end{align}
    The local transition function $\delta$ is $\bullet_{H_0}$-invariant.
  \end{lemma}

  \begin{proof}
    Let $\ell \in Q^N$ and let $h_0 \in H_0$. Then, there is a $c \in Q^M$ such that $\ell$ is observed by $m_0$ in $c$ (see \cref{remark:local-configuration-induces-global-configuration}), that is to say, that 
    \begin{equation*}
      \ForEach n \in N \Holds \ell(n) = c(m_0 \isSemiActedUponBy n).
    \end{equation*}
    And, because $\isSemiActedUponBy$ exhausts its defect with respect to its semi-com\-mu\-ta\-tiv\-i\-ty with $\actsOnPoint$ in $m_0$ (see \cref{item:semi-commutativity-of-liberation:exhausts} of \cref{lemma:semi-commutativity-of-liberation}),
    \begin{equation*}
      \ForEach n \in N \Holds (h_0^{-1} \actsOnPoint m_0) \isSemiActedUponBy n = h_0^{-1} \actsOnPoint (m_0 \isSemiActedUponBy h_0 \cdot n).
    \end{equation*}
    Hence, according to \cref{lemma:semi-commutativity-induced-left-action-and-bullet},
    \begin{align*}
      \delta(h_0 \bullet \ell)
      &= \delta(h_0 \bullet [n \mapsto c(m_0 \isSemiActedUponBy n)])\\
      &= \delta(h_0 \bullet [n \mapsto c((h_0^{-1} \actsOnPoint m_0) \isSemiActedUponBy n)])\\
      &= \delta(n \mapsto (h_0 \actsOnMap c)(m_0 \isSemiActedUponBy n))\\
      &= \Delta_0(h_0 \actsOnMap c)(m_0).
    \end{align*}
    And, because $\Delta_0$ is $\actsOnMap_H$-e\-qui\-var\-i\-ant,
    \begin{align*}
      \Delta_0(h_0 \actsOnMap c)(m_0)
      &= (h_0 \actsOnMap \Delta_0(c))(m_0)\\
      &= \Delta_0(c)(h_0^{-1} \actsOnPoint m_0)\\
      &= \Delta_0(c)(m_0)\\
      &= \delta(\ell).
    \end{align*}
    Put the last two chains of equalities together to see that $\delta(h_0 \bullet \ell) = \delta(\ell)$. In conclusion, $\delta$ is $\bullet_{H_0}$-invariant.
  \end{proof}

  \begin{theorem}
  \label{theorem:local-invariance-versus-global-equivariance}
    Let $\mathcal{R} = \ntuple{\mathcal{M}, \mathcal{K}} = \ntuple{\ntuple{M, G, \actsOnPoint}, \ntuple{m_0, \family{g_{m_0, m}}_{m \in M}}}$ be a cell space, let $\mathcal{C} = \ntuple{\mathcal{R}, Q, N, \delta}$ be a semi-cellular automaton, and let $H$ be a $\mathcal{K}$-big subgroup of $G$. The local transition function $\delta$ is $\bullet_{H_0}$-invariant if and only if the global transition function $\Delta$ is $\actsOnMap_H$-e\-qui\-var\-i\-ant.
  \end{theorem}

  \begin{proof}
    First, let $\delta$ be $\bullet_{H_0}$-invariant. Furthermore, let $h \in H$, let $c \in Q^M$, and let $m \in M$. Then, because $\isSemiActedUponBy$ semi-commutes with $\actsOnPoint$ (see \cref{item:semi-commutativity-of-liberation:semi-commutes} of \cref{lemma:semi-commutativity-of-liberation}), there is an $h_0 \in H_0$ such that
    \begin{equation*}
      \ForEach n \in N \Holds (h^{-1} \actsOnPoint m) \isSemiActedUponBy n = h^{-1} \actsOnPoint (m \isSemiActedUponBy h_0 \cdot n).
    \end{equation*}
    Hence, according to \cref{lemma:semi-commutativity-induced-left-action-and-bullet},
    \begin{align*}
      (h \actsOnMap \Delta(c))(m)
      &= \Delta(c)(h^{-1} \actsOnPoint m)\\
      &= \delta(n \mapsto c((h^{-1} \actsOnPoint m) \isSemiActedUponBy n))\\
      &= \delta(h_0^{-1} \bullet [n \mapsto (h \actsOnMap c)(m \isSemiActedUponBy n)])\\
      &= \delta(n \mapsto (h \actsOnMap c)(m \isSemiActedUponBy n))\\
      &= \Delta(h \actsOnMap c)(m).
    \end{align*}
    In conclusion, $\Delta$ is $\actsOnMap_H$-e\-qui\-var\-i\-ant.

    Secondly, let $\Delta$ be $\actsOnMap_H$-e\-qui\-var\-i\-ant. Then, according to \cref{definition:global-transition-function} and \cref{lemma:local-invariance-versus-global-equivariance:global-to-local}, the local transition function $\delta$ is $\bullet_{H_0}$-invariant.
  \end{proof}

  \begin{corollary}
  \label{corollary:semi-is-non-semi-if-and-only-if-global-transition-function-equivariant}
    Let $\mathcal{C}$ be a semi-cellular automaton. It is a cellular automaton if and only if its global transition function is $\actsOnMap$-e\-qui\-var\-i\-ant.
  \end{corollary}

  \begin{proof}
    This is a direct consequence of \cref{theorem:local-invariance-versus-global-equivariance} with $H = G$.
  \end{proof}

  \begin{example}[Lattice]
  \label{example:cell-space-integer-lattice}
    Let $M$ be the integer lattice $\Z^d$ of dimension $d \in \N_+$, let $G$ be the symmetry group of $M$, which is generated by translations, rotations, and reflections, let $\actsOnPoint$ be the transitive left group action of $G$ on $M$ by function application, let $m_0$ be the origin $0$, and, for each point $m \in M$, let $g_{m_0, m}$ be the translation that maps $m_0$ to $m$, namely $m + \blank$.

    The tuple $\mathcal{R} = \ntuple{\ntuple{M, G, \actsOnPoint}, \ntuple{m_0, \family{g_{m_0, m}}_{m \in M}}}$ is a cell space. The stabiliser $G_0$ of $m_0$ under $\actsOnPoint$ is generated by the rotations about and reflections through the origin. Under the identification of $G \modulo G_0$ with $\Z^d$ by $\iota$, the induced semi-action $\isSemiActedUponBy$ is the map $M \times \Z^d \to M$, $(m, t) \mapsto m + t$ (the elements of $\Z^d$ are translation vectors).

    In the remainder of this example, let $d = 2$. Then, for each cell $m \in M$, there are four and only four rotations about $m$, namely those by $0\degree$, $90\degree$, $180\degree$, and $270\degree$; and there are five and only five reflections through $m$, namely those about the vertical, horizontal, descending diagonal, and ascending diagonal line through $m$, and the point reflection.

    Let $Q$ be the binary set $\setOf{0, 1}$, let $N$ be the subset $\setOf{v \in \Z^2 \suchThat \normOf{v} \leq 1}$ of $\Z^2$ (note that $G_0 \actsOnPoint N \subseteq N$), and let $\delta$ be a map from $Q^N \to Q$. The quadruple $\mathcal{C} = \ntuple{\mathcal{R}, Q, N, \delta}$ is a semi-cellular automaton. The neighbourhood $N$ is called \define{von Neumann neighbourhood}\graffito{von Neumann neighbourhood}\index{neighbourhood!von Neumann}\index{Neumann neighbourhood von@von Neumann neighbourhood}. The local configurations of $Q^N$ are depicted in \cref{figure:cellular-automaton-two-dimensional-integer-lattice:local-configuration-of-origin}. 
    \begin{figure}
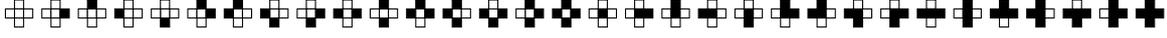

      \myfloatalign
      \figureCellularAutomatonTwoDimensionalIntegerLatticeLocalConfiguraitonOfOrigin
      \caption{The local configurations of $Q^N$.}
      \label{figure:cellular-automaton-two-dimensional-integer-lattice:local-configuration-of-origin}
    \end{figure}
    \begin{aenumerate}
      \item Let $H$ be the $\mathcal{K}$-big subgroup of $G$ that is generated by the translations. Then, the stabiliser $H_0$ of $m_0$ under $\actsOnPoint_H$ is the trivial subgroup of $H$, the local transition function $\delta$ is $\bullet_{H_0}$-invariant, and the global transition function $\Delta$ of $\mathcal{C}$ is $\actsOnMap_H$-e\-qui\-var\-i\-ant.

            Under the identification of $(H, \after)$ with $(\Z^2, +)$ by $t + \blank \mapsto t$, that equivariance follows directly from the associativity of $+$, more precisely, from $-t + (m + n) = (-t + m) + n$, for $t \in \Z^2$, $m \in M$, and $n \in N$. Indeed, for each translation vector $t \in \Z^2$, each global configuration $c \in Q^M$, and each cell $m \in M$,
            \begin{align*}
              \Delta(t \actsOnMap c)(m)
              &= \delta(n \mapsto c(-t + (m + n)))\\
              &= \delta(n \mapsto c((-t + m) + n))\\
              &= (t \actsOnMap \Delta(c))(m).
            \end{align*}
      \item Let $H$ be the $\mathcal{K}$-big subgroup of $G$ that is generated by the translations and reflections about the horizontal and vertical axis. Then, the stabiliser $H_0$ is generated by the latter, the local transition function $\delta$ is $\bullet_{H_0}$-invariant if and only if it maps the local configurations that are in the same column in \cref{figure:example-local-configuration-constraints-horizontal-and-vertical-reflections} to the same state, which are exactly those that are in the same orbit under $\bullet_{H_0}$.
            \begin{figure}
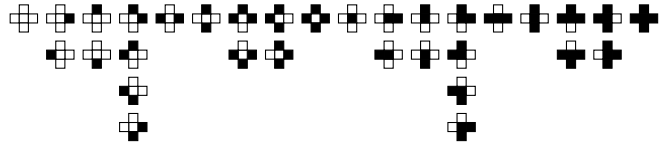

              \myfloatalign
              \figureExampleLocalConfigurationConstraintsHorizontalAndVerticalReflections
              \caption{The local configurations of $Q^N$ arranged such that the ones in the same column can be transformed into each other by reflections about the horizontal or vertical axis or combinations of these.}
              \label{figure:example-local-configuration-constraints-horizontal-and-vertical-reflections}
            \end{figure}
      \item Let $H$ be the $\mathcal{K}$-big subgroup of $G$ that is generated by the translations and reflections about the ascending and descending diagonal line through the origin. Then, the stabiliser $H_0$ is generated by the latter, the local transition function $\delta$ is $\bullet_{H_0}$-invariant if and only if it maps the local configurations that are in the same column in \cref{figure:example-local-configuration-constraints-diagonal-reflections} to the same state.
            \begin{figure}
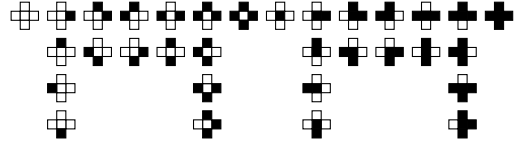

              \myfloatalign
              \figureExampleLocalConfigurationConstraintsDiagonalReflections
              \caption{The local configurations of $Q^N$ arranged such that the ones in the same column can be transformed into each other by reflections about the ascending or descending diagonal through the origin or combinations of these.}
              \label{figure:example-local-configuration-constraints-diagonal-reflections}
            \end{figure}
      \item Let $H$ be the $\mathcal{K}$-big subgroup of $G$ that is generated by the translations and rotations about the origin. Then, the stabiliser $H_0$ consists of the latter, the local transition function $\delta$ is $\bullet_{H_0}$-invariant if and only if it maps the local configurations that are in the same column in \cref{figure:example-local-configuration-constraints-rotations} to the same state.
            \begin{figure}
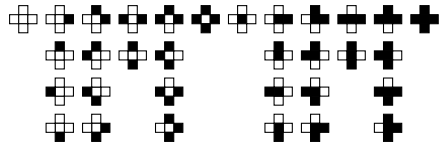

              \myfloatalign
              \figureExampleLocalConfigurationConstraintsRotations
              \caption{The local configurations of $Q^M$ arranged such that the ones in the same column can be transformed into each other by rotations about the origin.}
              \label{figure:example-local-configuration-constraints-rotations}
            \end{figure}
      \item Let $H$ be the $\mathcal{K}$-big subgroup of $G$ that is generated by the translations and the point reflection through the origin. Then, the stabiliser $H_0$ consists of the latter, the local transition function $\delta$ is $\bullet_{H_0}$-invariant if and only if it maps the local configurations that are in same column in \cref{figure:example-local-configuration-constraints-point-reflection} to the same state.
            \begin{figure}
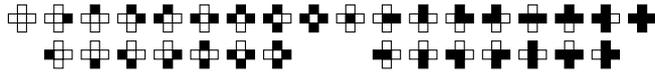

              \myfloatalign
              \figureExampleLocalConfigurationConstraintsPointReflection
              \caption{The local configurations of $Q^M$ arranged such that the ones in the same column can be transformed into each other by the point reflection through the origin.}
              \label{figure:example-local-configuration-constraints-point-reflection}
            \end{figure}
      \item Let $H$ be the $\mathcal{K}$-big subgroup of $G$ that is generated by the translations, rotations and reflections. Then, $H$ is the group $G$ and the stabiliser $H_0$ is generated by the rotations about and reflections through the origin. However, because there are only two states, the orbits of $\bullet_{H_0}$ are the same as if $H_0$ were generated only by the rotations about the origin. Hence, the local transition function $\delta$ is $\bullet_{H_0}$-invariant if and only if it maps the local configurations that are in the same column in \cref{figure:example-local-configuration-constraints-rotations} to the same state.
    \end{aenumerate}
    In each case, the global transition function $\Delta$ is equivariant under the respective symmetries of $H$. Figuratively speaking, being equivariant under translations means being oblivious of global positions, being equivariant under rotations means being oblivious of orientation, and being equivariant under reflections means being oblivious of orientation-reversal.
  \end{example}

  \begin{example}[Plane]
  \label{example:plane:equivariance}
    In the situation of \cref{example:plane:liberation}, let $\mathcal{C} = \ntuple{\mathcal{R}, Q, N, \delta}$ be a semi-cellular automaton, let $R$ be a subgroup of the group $G_0$ of rotations about $m_0$, and let $T \rtimes R$ be the inner semi-direct product of $R$ acting on $T$. The group $T \rtimes R$ is $\mathcal{K}$-big and the stabiliser $(T \rtimes R)_0$ of $m_0$ under $\actsOnMap_{T \rtimes R}$ is the group $R$. According to \cref{theorem:local-invariance-versus-global-equivariance}, the local transition function $\delta$ is $\bullet_R$-invariant if and only if the global transition function $\Delta$ is $\actsOnMap_{T \rtimes R}$-e\-qui\-var\-i\-ant. In particular, because $\delta$ is $\bullet_{\setOf{e_G}}$-invariant, the map $\Delta$ is $\actsOnMap_T$-e\-qui\-var\-i\-ant. Broadly speaking, the more invariant $\delta$ is, the more equivariant is $\Delta$, and vice versa. 
  \end{example}

  \begin{counterexample}[Peculiar Shift Maps]
    In \cref{theorem:local-invariance-versus-global-equivariance}, if the assumption that the subgroup $H$ of $G$ is $\mathcal{K}$-big does not hold, then the stated equivalence may not hold, which is illustrated by the following example.

    In the situation of \cref{item:lattice:automata:shifts:odd-even} of \cref{example:lattice:automata:shifts}, the subgroup $T$ of $G$ is not $\mathcal{K}''$-big and, although the local transition function $\delta$ is $\bullet_{T_0}$-invariant, the global transition function $\Delta''$ is not $\actsOnMap_T$-e\-qui\-var\-i\-ant. However, according to \cref{lemma:local-invariance-versus-global-equivariance:global-to-local}, even for a subgroup $H$ of $G$ that is not $\mathcal{K}$-big, there is no semi-cellular automaton whose local transition function is not $\bullet_{H_0}$-invariant, although its global transition function is $\actsOnMap_H$-e\-qui\-var\-i\-ant.
  \end{counterexample}

  If a map on global configurations is equivariant under translations and is determined by a local transition function at the origin, then, using translations from the origin to all other cells, we see that the map is determined by the local transition function at all cells, in other words, the map is the induced global transition function (see \cref{theorem:determination-of-cellular-automata-by-behaviour-at-origin}).

  \begin{theorem}
  \label{theorem:determination-of-cellular-automata-by-behaviour-at-origin}
    Let $\mathcal{R} = \ntuple{\mathcal{M}, \mathcal{K}} = \ntuple{\ntuple{M, G, \actsOnPoint}, \ntuple{m_0, \family{g_{m_0, m}}_{m \in M}}}$ be a cell space, let $\mathcal{C} = \ntuple{\mathcal{R}, Q, N, \delta}$ be a semi-cellular automaton, let $\Delta_0$ be a map from $Q^M$ to $Q^M$, and let $H$ be a $\mathcal{K}$-big subgroup of $G$. The following two statements are equivalent: 
    \begin{aenumerate}
      \item\label{item:determination-of-cellular-automata-by-behaviour-at-origin:invariance}
            The local transition function $\delta$ is $\bullet_{H_0}$-invariant and the global transition function of $\mathcal{C}$ is $\Delta_0$.
      \item\label{item:determination-of-cellular-automata-by-behaviour-at-origin:equivariance}
            The map $\Delta_0$ is $\actsOnMap_H$-e\-qui\-var\-i\-ant and
            \begin{equation}
            \label{equation:global-transition-function-at-origins}
              \ForEach c \in Q^M \Holds \Delta_0(c)(m_0) = \delta(n \mapsto c(m_0 \isSemiActedUponBy n)). \qedhere
            \end{equation}
    \end{aenumerate}
  \end{theorem}

  \begin{proof}
    First, let $\delta$ be $\bullet_{H_0}$-invariant and let $\Delta_0$ be the global transition function of $\mathcal{C}$. According to \cref{theorem:local-invariance-versus-global-equivariance}, the map $\Delta_0$ is $\actsOnMap_H$-e\-qui\-var\-i\-ant and, according to \cref{definition:global-transition-function}, \cref{equation:global-transition-function-at-origins} holds.

    Secondly, let $\Delta_0$ be $\actsOnMap_H$-e\-qui\-var\-i\-ant and let \cref{equation:global-transition-function-at-origins} hold. According to \cref{lemma:local-invariance-versus-global-equivariance:global-to-local}, the local transition function $\delta$ is $\bullet_{H_0}$-invariant. Furthermore, let $c \in Q^M$ and let $m \in M$. Put $h = g_{m_0, m}^{-1} \in H$. Then,
    \begin{equation*}
      \Delta_0(c)(m) = \Delta_0(c)(h^{-1} \actsOnPoint m_0) = (h \actsOnMap \Delta_0(c))(m_0).
    \end{equation*}
    And, because $\Delta_0$ is $\actsOnMap_H$-e\-qui\-var\-i\-ant,
    \begin{equation*}
      (h \actsOnMap \Delta_0(c))(m_0)
      = \Delta_0(h \actsOnMap c)(m_0)
      = \delta(n \mapsto (h \actsOnMap c)(m_0 \isSemiActedUponBy n)). 
    \end{equation*}
    And, because $\isSemiActedUponBy$ semi-commutes with $\actsOnPoint$ (see \cref{item:semi-commutativity-of-liberation:semi-commutes} of \cref{lemma:semi-commutativity-of-liberation}), there is an $h_0 \in H_0$ such that, for each $n \in N$, we have $(h^{-1} \actsOnPoint m_0) \isSemiActedUponBy n = h^{-1} \actsOnPoint (m_0 \isSemiActedUponBy h_0 \cdot n)$; and therefore, according to \cref{lemma:semi-commutativity-induced-left-action-and-bullet},
    \begin{align*}
      \delta(n \mapsto (h \actsOnMap c)(m_0 \isSemiActedUponBy n))
      &= \delta(h_0 \bullet [n \mapsto c((h^{-1} \actsOnPoint m_0) \isSemiActedUponBy n)])\\
      &= \delta(h_0 \bullet [n \mapsto c(m \isSemiActedUponBy n)]).
    \end{align*}
    And, because $\delta$ is $\bullet_{H_0}$-invariant,
    \begin{equation*}
      \delta(h_0 \bullet [n \mapsto c(m \isSemiActedUponBy n)]) = \delta(n \mapsto c(m \isSemiActedUponBy n)).
    \end{equation*}
    Put the last four chains of equalities together to see that $\Delta_0(c)(m) = \delta(n \mapsto c(m \isSemiActedUponBy n))$. In conclusion, $\Delta_0$ is the global transition function of $\mathcal{C}$.
  \end{proof}

  \begin{corollary}
    Let $\mathcal{C}$ be a cellular automaton with set of cells $M$ and set of states $Q$ and let $\Delta_0$ be a map from $Q^M$ to $Q^M$. The global transition function of $\mathcal{C}$ is $\Delta_0$ if and only if 
    \begin{equation*}
      \ForEach c \in Q^M \Holds \Delta_0(c)(m_0) = \delta(n \mapsto c(m_0 \isSemiActedUponBy n)). \qedhere
    \end{equation*}
  \end{corollary}

  \begin{proof}
    This is a direct consequence of \cref{corollary:semi-is-non-semi-if-and-only-if-global-transition-function-equivariant} and theorem \ref{theorem:determination-of-cellular-automata-by-behaviour-at-origin} with $H = G$. 
  \end{proof}

  \begin{counterexample}[Peculiar Complement and Shift Maps]
    In \cref{theorem:determination-of-cellular-automata-by-behaviour-at-origin}, if the assumption that the subgroup $H$ of $G$ is $\mathcal{K}$-big does not hold, then \cref{item:determination-of-cellular-automata-by-behaviour-at-origin:invariance} may not be equivalent to \cref{item:determination-of-cellular-automata-by-behaviour-at-origin:equivariance}, which is illustrated by the following examples. 

    First, let $\mathcal{R}$ be the cell space $\ntuple{\ntuple{\Z, \Z, +}, \ntuple{0, \family{m}_{m \in \Z}}}$, let $\mathcal{C}$ be the cellular automaton $\ntuple{\mathcal{R}, \Z \modulo 2\Z, \setOf{0}, \delta \from \ell \mapsto \ell(0)}$, let $\Delta_0$ be the map $(\Z \modulo 2\Z)^\Z \to (\Z \modulo 2\Z)^\Z$, $c \mapsto c(m) + (m + 2\Z)$, and let $H$ be the subgroup $2 \Z$ of $\Z$, which is not $\ntuple{0, \family{m}_{m \in \Z}}$-big. The map $\Delta_0$ is the identity map on $2\Z$ and the bitwise complement map on $2\Z + 1$, where we call $0 + 2\Z$ and $1 + 2\Z$ \emph{complements of each other}. The map $\Delta_0$ is $\actsOnMap_H$-e\-qui\-var\-i\-ant and $\Delta_0(\blank)(0) = \blank(0) = \delta(n \mapsto \blank(0 \isSemiActedUponBy n))$. However, the global transition function of $\mathcal{C}$ is not $\Delta_0$ but the identity map on $(\Z \modulo 2\Z)^\Z$. Actually, the map $\Delta_0$ is not the global transition function of any semi-cellular automaton. 

    Secondly, in the situation of \cref{item:lattice:automata:shifts:strange} of \cref{example:lattice:automata:shifts}, the subgroup $T$ of $G$ is not $\mathcal{K}'''$-big and, although the local transition function $\delta$ is $\bullet_{T_0}$-invariant and the global transition function of $\mathcal{C}'''$ is $\Delta'''$, the map $\Delta'''$ is not $\actsOnMap_T$-e\-qui\-var\-i\-ant.
  \end{counterexample}

  Multiplying the neighbourhoods of two cellular automata and chaining their local transition functions yields an automaton whose global transition function is the composition of the ones of the other two automata (see \cref{theorem:composition-of-cellular-automata}).

  \begin{theorem}
  \label{theorem:composition-of-cellular-automata}
    Let $\mathcal{R} = \ntuple{\mathcal{M}, \mathcal{K}} = \ntuple{\ntuple{M, G, \actsOnPoint}, \ntuple{m_0, \family{g_{m_0, m}}_{m \in M}}}$ be a cell space, let $\mathcal{C} = \ntuple{\mathcal{R}, Q, N, \delta}$ and $\mathcal{C}' = \ntuple{\mathcal{R}, Q, N', \delta'}$ be two semi-cellular automata, and let $H$ be a $\mathcal{K}$-big subgroup of $G$ such that $\delta$ and $\delta'$ are $\bullet_{H_0}$-invariant. Furthermore, let 
    \begin{equation*}
      N'' = \setOf{g \cdot n' \suchThat n \in N, n' \in N', g \in n} 
    \end{equation*}
    and let
    \begin{align*}
      \delta'' \from Q^{N''} &\to Q,\\
      \ell'' &\mapsto \delta(n \mapsto \delta'(n' \mapsto \ell''(g_{m_0, m_0 \isSemiActedUponBy n} \cdot n'))). 
    \end{align*}
    The quadruple $\mathcal{C}'' = \ntuple{\mathcal{R}, Q, N'', \delta''}$ is a semi-cellular automaton whose local transition function is $\bullet_{H_0}$-invariant and whose global transition function is $\Delta \after \Delta'$.
  \end{theorem}

  \begin{proof}
    Because $G_0 \cdot N \subseteq N$, we have $G_0 \cdot N'' \subseteq N''$. And, because $g_{m_0, m_0 \isSemiActedUponBy n} \in n$, the map $\delta''$ is well-defined. Therefore, the quadruple $\mathcal{C}'' = \ntuple{\mathcal{R}, Q, N'', \delta''}$ is a semi-cellular automaton.

    Moreover, because $\delta$ and $\delta'$ are $\bullet_{H_0}$-invariant, according to theorem \ref{theorem:local-invariance-versus-global-equivariance}, the maps $\Delta$ and $\Delta'$ are $\actsOnMap_H$-e\-qui\-var\-i\-ant and thus $\Delta \after \Delta'$ also. And, because $g_{m_0, m_0} = e_G$, 
    \begin{equation*}
      \ForEach n \in N \ForEach n' \in N' \Holds (m_0 \isSemiActedUponBy n) \isSemiActedUponBy n' = m_0 \isSemiActedUponBy g_{m_0, m_0 \isSemiActedUponBy n} \cdot n'.
    \end{equation*}
    Hence, for each $c \in Q^M$,
    \begin{align*}
      (\Delta \after \Delta')(c)(m_0) &= \delta(n \mapsto \Delta'(c)(m_0 \isSemiActedUponBy n))\\
                                      &= \delta(n \mapsto \delta'(n' \mapsto c((m_0 \isSemiActedUponBy n) \isSemiActedUponBy n')))\\
                                      &= \delta(n \mapsto \delta'(n' \mapsto c(m_0 \isSemiActedUponBy g_{m_0, m_0 \isSemiActedUponBy n} \cdot n')))\\
                                      &= \delta''(n'' \mapsto c(m_0 \isSemiActedUponBy n'')).
    \end{align*}
    Therefore, according to \cref{theorem:determination-of-cellular-automata-by-behaviour-at-origin}, the local transition function $\delta''$ is $\bullet_{H_0}$-invariant and the global transition function of $\mathcal{C}''$ is $\Delta \after \Delta'$.
  \end{proof}

  \begin{corollary}
  \label{corollary:composition-of-cellular-automata}
    Let $\mathcal{C}$ and $\mathcal{C}'$ be two cellular automata over $\mathcal{R}$ and $\mathcal{R}'$ respectively. There is a cellular automaton whose global transition function is $\Delta \after \Delta'$.
  \end{corollary}

  \begin{proof}
    This is a direct consequence of \cref{corollary:existence-in-other-coordinate-system} and theorem \ref{theorem:composition-of-cellular-automata} with $H = G$. 
  \end{proof}

  \begin{remark}
  \label{remark:identification:composition-of-cellular-automata}
    In \cref{theorem:composition-of-cellular-automata}, under the identification of $G \modulo G_0$ with $M$ by $\iota$, we have $N'' = \setOf{g \actsOnPoint n' \suchThat n \in N, n' \in N', g \in G_{m_0, n}}$ and, for each neighbour $n \in N$ and each neighbour $n' \in N'$, we have $g_{m_0, m_0 \isSemiActedUponBy n} \cdot n' = g_{m_0, n} \actsOnPoint n'$.
  \end{remark}

  \begin{remark}
    That the local transition function $\delta''$ in \cref{theorem:composition-of-cellular-automata} is $\bullet_{H_0}$-invariant can be shown directly as follows: Let $h_0 \in H_0$. Then, for each $n \in N$,
    \begin{equation*}
      h_0^{-1} g_{m_0, m_0 \isSemiActedUponBy n} G_0
      = h_0^{-1} \cdot n
      = g_{m_0, m_0 \isSemiActedUponBy h_0^{-1} \cdot n} G_0, 
    \end{equation*}
    and hence, because $H$ is $\mathcal{K}$-big, the element $h_{n,0} = g_{m_0, m_0 \isSemiActedUponBy h_0^{-1} \cdot n}^{-1} h_0^{-1} g_{m_0, m_0 \isSemiActedUponBy n} \in G_0 \cap H = H_0$ satisfies $h_0^{-1} g_{m_0, m_0 \isSemiActedUponBy n} = g_{m_0, m_0 \isSemiActedUponBy h_0^{-1} \cdot n} h_{n,0}$. Therefore, because $\delta$ and $\delta'$ are $\bullet_{H_0}$-invariant, for each $\ell'' \in Q^{N''}$,
    \begin{align*}
      \delta''(h_0 \bullet \ell'')
      &= \delta(n \mapsto \delta'(n' \mapsto (h_0 \bullet \ell'')(g_{m_0, m_0 \isSemiActedUponBy n} \cdot n')))\\
      &= \delta(n \mapsto \delta'(n' \mapsto \ell''(h_0^{-1} g_{m_0, m_0 \isSemiActedUponBy n} \cdot n')))\\
      &= \delta(n \mapsto \delta'(n' \mapsto \ell''(g_{m_0, m_0 \isSemiActedUponBy h_0^{-1} \cdot n} \cdot (h_{n,0} \cdot n'))))\\
      &= \delta(n \mapsto \delta'(h_{n,0}^{-1} \bullet [n' \mapsto \ell''(g_{m_0, m_0 \isSemiActedUponBy h_0^{-1} \cdot n} \cdot n')]))\\
      &= \delta(n \mapsto \delta'(n' \mapsto \ell''(g_{m_0, m_0 \isSemiActedUponBy h_0^{-1} \cdot n} \cdot n')))\\
      &= \delta(h_0 \bullet [n \mapsto \delta'(n' \mapsto \ell''(g_{m_0, m_0 \isSemiActedUponBy n} \cdot n'))])\\
      &= \delta(n \mapsto \delta'(n' \mapsto \ell''(g_{m_0, m_0 \isSemiActedUponBy n} \cdot n')))\\
      &= \delta''(\ell'').
    \end{align*}
    In conclusion, $\delta''$ is $\bullet_{H_0}$-invariant.
  \end{remark}

  \begin{counterexample}[Peculiar Shift Maps]
  \label{example:necessity:composition-of-cellular-automata}
    In \cref{theorem:composition-of-cellular-automata}, if the assumption that the subgroup $H$ of $G$ is $\mathcal{K}$-big or the one that the local transition functions $\delta$ and $\delta'$ are $\bullet_{H_0}$-invariant does not hold, then the global transition function of $\mathcal{C}''$ may not be $\Delta \after \Delta'$, which is illustrated by the following example.

    Let $\mathcal{C}''' = \ntuple{\ntuple{\mathcal{M}, \mathcal{K}'''}, Q, N, \delta}$ be the semi-cellular automata from \cref{item:lattice:automata:shifts:strange} of \cref{example:lattice:automata:shifts}. The subgroup $T$ of $G$ is not $\mathcal{K}'''$-big (and the local transition function $\delta$ is not $\bullet_{G_0}$-invariant) and although the local transition function $\delta$ is $\bullet_{T_0}$-invariant (and although the subgroup $G$ of $G$ is $\mathcal{K}'''$-big), the map $(\Delta''')^2$, which is depicted in \cref{figure:strange-shift-squared}, is not a global transition function of a semi-cellular automaton, in particular, not of the one the construction in \cref{theorem:composition-of-cellular-automata} applied to $\mathcal{C}'''$ yields.

    Broadly speaking, for each global configuration $c \in Q^M$ and each cell $m \in M$, the state
    \begin{equation*}
      (\Delta''')^2(c)(m) = \begin{dcases*}
        c(m + 2), &if $m \notin \setOf{0, 1}$,\\
        c(m), &if $m \in \setOf{0, 1}$,
      \end{dcases*}
    \end{equation*}
    depends in an asymmetric way on the cell that cannot be induced uniformly by a local transition function and a coordinate system. And, the construction in \cref{theorem:composition-of-cellular-automata} applied to $\mathcal{C}'''$ yields the semi-cellular automaton $\mathcal{C}'' = \ntuple{\ntuple{\mathcal{M}, \mathcal{K}'''}, Q, \setOf{-2, -1, 0, 1, 2}, \delta'' \from \ell'' \mapsto \ell''(0)}$ whose global transition function is the identity map on $Q^M$. Clearly, the square $(\Delta''')^2$ is not the identity map on $Q^M$.
    \begin{figure}
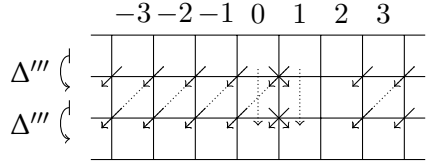

      \myfloatalign
      \figureStrangeShiftSquared
      \caption{The first row depicts the part of a global configuration that corresponds to the cells $\setOf{-3, -2, \dotsc, 3}$, the second row the same part of the image of that global configuration under $\Delta'''$, and the third row the part of the image under $(\Delta''')^2$. And, the arrows from cells in one row to another depict the flow of states.}
      \label{figure:strange-shift-squared}
    \end{figure}
  \end{counterexample}

  \begin{proof}
    First, suppose that there is a coordinate system $\mathcal{K} = \ntuple{m_0, \family{g_{m_0, m}}_{m \in M}}$ for $\mathcal{M}$ and there is a semi-cellular automaton $\mathcal{C}'' = \ntuple{\ntuple{\mathcal{M}, \mathcal{K}}, Q, N'', \delta''}$ whose global transition function is $(\Delta''')^2$. Let $m \in M \smallsetminus \setOf{0, 1}$ and let $c \in Q^M$ such that $c(m + 2) \neq c(m)$, for example, $m = 4$ and $c \in Q^M$ with $c \equiv 0$ on $M \smallsetminus \setOf{6}$ and $c \equiv 1$ on $\setOf{6}$. Then, letting $c' = g_{m_0, 0} g_{m_0, m}^{-1} \actsOnMap c$, according to \cref{remark:identification:global-transition-function},
    \begin{align*}
      c(m + 2)
      &= (\Delta''')^2(c)(m)\\
      &= \delta''((g_{m_0, m}^{-1} \actsOnMap c)\restrictedTo_{N''})\\
      &= \delta''((g_{m_0, 0}^{-1} g_{m_0, 0} g_{m_0, m}^{-1} \actsOnMap c)\restrictedTo_{N''})\\
      &= \delta''((g_{m_0, 0}^{-1} \actsOnMap c')\restrictedTo_{N''})\\
      &= (\Delta''')^2(c')(0)\\
      &= c'(0)\\
      &= c(g_{m_0, m} g_{m_0, 0}^{-1} \actsOnPoint 0)\\
      &= c(m),
    \end{align*}
    which contradicts that $c(m + 2) \neq c(m)$. In conclusion, there is no such coordinate system and semi-cellular automaton.

    Secondly, the construction in \cref{theorem:composition-of-cellular-automata} applied to $\mathcal{C}'''$ yields the semi-cellular automaton $\mathcal{C}''$ with the neighbourhood $N'' = N + N = \setOf{-2, -1, 0, 1, 2}$ and the local transition function $\delta''$ such that, for each local configuration $\ell'' \in Q^{N''}$,
    \begin{multline*}
      \delta''(\ell'')
       = \delta(n \mapsto \delta(n' \mapsto \ell''(g_{0, n}''' \actsOnPoint n')))\\
       = \ell''(g_{0, 1}''' \actsOnPoint 1)
       = \ell''(\varrho_1(1))
       = \ell''(0),
     \end{multline*}
     where $\ntuple{0, \family{g_{0, m}'''}_{m \in M}}$ is the coordinate system $\mathcal{K}'''$. In conclusion, the global transition function of $\mathcal{C}''$ is the identity map on $Q^M$.
  \end{proof}

  \begin{example}[Peculiar Shift Maps and Logical Or]
    In the situation of \cref{example:lattice:homogeneous-space}, let $\mathcal{K}$ and $\mathcal{K}'$ be the two coordinate systems $\ntuple{0, \family{\tau_m}_{m \in M}}$ and $\ntuple{0, \family{\tau_m}_{m \in M \smallsetminus \setOf{1}} \times \family{\varrho_m}_{m \in \setOf{1}}}$ for $\mathcal{M}$, let $\mathcal{K}''$ be a coordinate system for $\mathcal{M}$, let $Q$ be the set $\setOf{0, 1}$, let $N$ be the set $\setOf{-1, 1}$, let $\delta$ be the map $Q^N \to Q$, $\ell \mapsto \ell(1)$, let $\delta''$ be the map $Q^N \to Q$, $\ell \mapsto \ell(-1) \lor \ell(1)$, and let $M$ be identified with $G \modulo G_0$ by $\iota \givenBy m \mapsto G_{0, m}$. The semi-cellular automaton $\mathcal{C} = \ntuple{\ntuple{\mathcal{M}, \mathcal{K}}, Q, N, \delta}$ is the one from \cref{item:lattice:automata:shifts:left-and-right} of \cref{example:lattice:automata:shifts} whose global transition function $\Delta$ is the left shift map, the semi-cellular automaton $\mathcal{C}' = \ntuple{\ntuple{\mathcal{M}, \mathcal{K}'}, Q, N, \delta}$ is the one from \cref{item:lattice:automata:shifts:strange} of \cref{example:lattice:automata:shifts} whose global transition function $\Delta'$ is the left shift map with one defect in cell $1$, and the cellular automaton $\mathcal{C}'' = \ntuple{\ntuple{\mathcal{M}, \mathcal{K}''}, Q, N, \delta''}$ is the one from \cref{example:lattice:automaton:or} whose global transition function $\Delta''$ is a pairwise logical or map (it does not depend on the choice of $\mathcal{K}''$).

    The group $T$ of translations is $\mathcal{K}$-big and, if we choose $\mathcal{K}$ for $\mathcal{K}''$, it is also $\mathcal{K}''$-big; and the local transition functions $\delta$ and $\delta''$ are $\bullet_{T_0}$-invariant. The construction in \cref{theorem:composition-of-cellular-automata} yields the following: The square $\Delta^2$ is the global transition function of the semi-cellular automaton $\ntuple{\ntuple{\mathcal{M}, \mathcal{K}}, Q, \setOf{-2, -1, 0, 1, 2}, \ell \mapsto \ell(2)}$ (see \cref{figure:left-shift-followed-by-left-shift}), the square $(\Delta'')^2$ is the one of the cellular automaton $\ntuple{\ntuple{\mathcal{M}, \blank}, Q, \setOf{-2, -1, 0, 1, 2}, \ell \mapsto \ell(-2) \lor \ell(0) \lor \ell(2)}$ (see \cref{figure:or-followed-by-or}), and the compositions $\Delta'' \after \Delta$ and $\Delta \after \Delta''$ are the one of the semi-cellular automaton $\ntuple{\ntuple{\mathcal{M}, \mathcal{K}}, Q, \setOf{-2, -1, 0, 1, 2}, \ell \mapsto \ell(0) \lor \ell(2)}$ (see \cref{figure:left-shift-followed-by-or,figure:or-followed-by-left-shift}). 

    Explicitly, for each global configuration $c \in Q^M$ and each cell $m \in M$, we have $\Delta^2(c)(m) = c(m + 2)$, and $(\Delta'')^2(c)(m) = c(m - 2) \lor c(m) \lor c(m + 2)$, and $(\Delta'' \after \Delta)(c)(m) = (\Delta \after \Delta'')(c)(m) = c(m) \lor c(m + 2)$.

    There is no $\mathcal{K}'$-big subgroup $H$ of $G$ such that the local transition function $\delta'$ is $\bullet_{H_0}$-invariant. And indeed, the square $(\Delta')^2$ (see \cref{figure:strange-shift-squared}), the compositions $\Delta' \after \Delta$ and $\Delta \after \Delta'$ (see \cref{figure:left-shift-followed-by-strange-shift,figure:strange-shift-followed-by-left-shift}), and the composition $\Delta'' \after \Delta'$ (see \cref{figure:strange-shift-followed-by-or}) are not global transition functions of semi-cellular automata, which is shown for $(\Delta')^2$ in \cref{example:necessity:composition-of-cellular-automata}. Nevertheless, the composition $\Delta' \after \Delta''$ (see \cref{figure:or-followed-by-strange-shift}) is the global transition function of the semi-cellular automaton $\ntuple{\ntuple{\mathcal{M}, \mathcal{K}'}, Q, \setOf{-2, -1, 0, 1, 2}, \ell \mapsto \ell(0) \lor \ell(2)}$; the reason is that, broadly speaking, the global transition function $\Delta''$, which is applied first to global configurations, does not have a defect that could be propagated by $\Delta'$. 

    Explicitly, for each global configuration $c \in Q^M$ and each cell $m \in M$,
    \begin{equation*}
      (\Delta')^2(c)(m) = \begin{dcases*}
        c(m + 2), &if $m \notin \setOf{0, 1}$,\\
        c(m), &if $m \in \setOf{0, 1}$,
      \end{dcases*}
    \end{equation*}
    \begin{equation*}
      (\Delta' \after \Delta)(c)(m) = \begin{dcases*}
        c(m + 2), &if $m \neq 1$,\\
        c(m), &if $m = 1$,
      \end{dcases*}
    \end{equation*}
    \begin{equation*}
      (\Delta \after \Delta')(c)(m) = \begin{dcases*}
        c(m + 2), &if $m \neq 0$,\\
        c(m), &if $m = 0$,
      \end{dcases*}
    \end{equation*}
    \begin{equation*}
      (\Delta'' \after \Delta')(c)(m) = \begin{dcases*}
        c(m) \lor c(m + 2), &if $m \notin \setOf{0, 2}$,\\
        c(m), &if $m = 0$,\\
        c(m - 2) \lor c(m + 2), &if $m = 2$,
      \end{dcases*}
    \end{equation*}
    and
    \begin{equation*}
      (\Delta' \after \Delta'')(c)(m) = \begin{dcases*}
        c(m) \lor c(m + 2), &if $m \neq 1$,\\
        c(m - 2) \lor c(m), &if $m = 1$.
      \end{dcases*}
    \end{equation*}
    Broadly speaking, the states of the first four compositions depend in an asymmetric way on the cell that cannot be induced uniformly by a local transition function and a coordinate system; whereas the state of the last composition depends on the state of the cell itself and on the state of its second to the right or to the left neighbour and this reversal of orientation can be achieved by using translations in the one case and reflections in the other.

    It is notable that $\Delta$ and $\Delta'$ are induced by the same local transition function but by different coordinates (the ones for $\Delta$ are even the translations); and, if we choose $\mathcal{K}$ for $\mathcal{K}''$, that $\Delta'$ and $\Delta''$ are induced by the same coordinates but different local transition functions (the one for $\Delta''$ is even $\bullet$-invariant).
    \begin{figure}
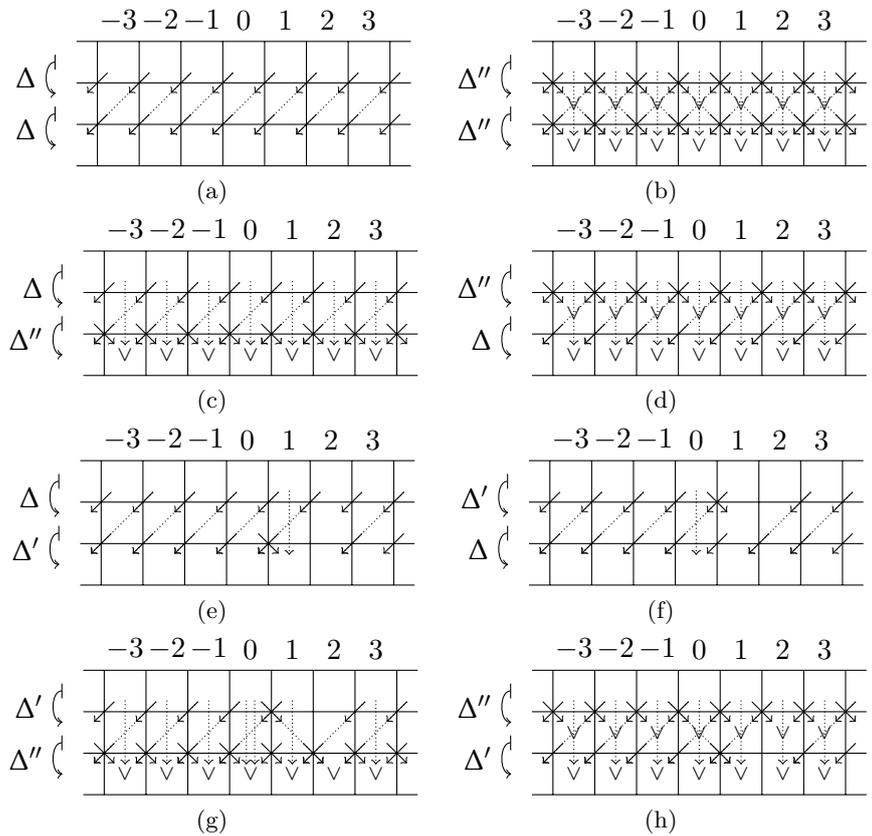

      \myfloatalign
      \figureCompositionsOfGlobalTransitionFunctions
      \caption{In each subfigure, the first row depicts the part of a global configuration that corresponds to the cells $\setOf{-3, -2, \dotsc, 3}$, the second row the same part of the image of that global configuration under a global transition function, and the third row the part of the image under the composition of two global transition functions. And, the arrows from cells in one row to another depict the flow of states and the disjunction symbols in cells depict that incoming states from arrows of the same length are combined disjunctively.} 
      \label{figure:compositions-of-global-transition-functions}
    \end{figure}
  \end{example}

  %
  %
  %

  \clearToOddPage
  \chapter{Quotients, Products, Restrictions, Extensions, and Decompositions of Spaces, Automata and Global Transition Functions}
  \chaptermark{Quotients, Products, Restrictions, Extensions, ...}
  \label{chapter:quotients-and-periodicity}

  \paragraph{Abstract.} We introduce and study periodicity of global configurations, and quotients, products, restrictions, extensions, and decompositions of left-ho\-mo\-ge\-neous spaces, coordinate systems, semi-cellular and cellular automata, and global transition functions of cellular automata.

  \paragraph{Remark.} This chapter generalises parts of sections~1.6 and~1.7 of the monograph \enquote{\citetitle*{ceccherini-silberstein:coornaert:2010}}\cite{ceccherini-silberstein:coornaert:2010}.

  \paragraph{Introduction.} The orientation preserving symmetries of an infinite circular cylinder are the rotations about its axis and the translations along its axis. These symmetries and its subgroups act on the cylinder itself and also on patterns over the cylinder.

  The patterns over the cylinder that are invariant under rotations are precisely those that are unicoloured on intersections of planes that are perpendicular to the axis, which are circles. Such a pattern can squeezed to a pattern over the axis, which is essentially the cylinder modulo rotations. Conversely, such a squeezed pattern can be stretched to the original pattern over the cylinder. So, there is a one-to-one correspondence between patterns that are invariant under rotations and patterns over the axis. When we squeeze the cylinder to its axis, we squeeze its symmetries to the translations, which are essentially the symmetries modulo rotations.

  Analogously, the patterns that are invariant under translations are precisely those that are unicoloured along the axis. Such a pattern can be squashed to a pattern over the circle that has the same radius as the cylinder, which is essentially the cylinder modulo translations. Conversely, such a squashed pattern can can be stretched to the original pattern over the cylinder. So, there is a one-to-one correspondence between patterns that are invariant under translations and patterns over the circle. When we squash the cylinder to the circle, we squash its symmetries to the rotations, which are essentially the symmetries modulo translations.

  A map on patterns over the cylinder that is equivariant under rotations (or translations), maps patterns that are invariant under rotations (or translations) to such patterns. It can be squeezed (or squashed) to a map on patterns over the axis (or the circle). Conversely, such a squeezed (or squashed) map can be stretched to a map on patterns over the cylinder, which however may not be the original map. These constructions can for example be performed with global transition functions of cellular automata over the cylinder.

  If the neighbourhood of a semi-cellular automaton over the cylinder is such that the neighbours of a cell lie on a circle (or on the line through the cell that is parallel to the axis), then its global transition function can be restricted to the circle (or to the line) and this restriction can in turn be extended to the cylinder by repeating it along the axis (or in a circle around the axis). The restriction itself is a global transition function and its extension is the original function. As the original function consists of multiple copies of the restriction, it is injective, surjective, or bijective if and only if the restriction has the respective property.


  \paragraph{Contents.} In \cref{section:periodicity} we introduce and study the notion of periodicity for global configurations: Those of the same period are in a bijective relationship with the global configurations over the orbit under the period and periodicity is preserved by global transition functions of cellular automata. In \cref{section:quotients} we introduce quotients of left group sets, coordinate systems, cell spaces, semi-cellular and cellular automata, and global transition functions. In \cref{section:naive-products} we introduce naïve products of semi-cellular automata; they are called naïve because naïve products of cellular automata are in general not cellular automata and because their global transition function depends on many arbitrary choices that have to be made. In \cref{section:products} we introduce products of cellular automata and products of global transition functions by subgroups that relate in a certain way to the stabiliser. In \cref{section:restrictions} we introduce restrictions of left group sets, coordinate systems, cell spaces, semi-cellular and cellular automata, and global transition functions. In \cref{section:extensions} we introduce extensions of semi-cellular automata and global transition functions. In \cref{section:decompositions} we show that global transition functions of some cellular automata are products of restrictions of themselves to subgroups that relate in a certain way to the stabiliser.

  \paragraph{Preliminary Notions.} A \emph{left group set} is a triple $\ntuple{M, G, \actsOnPoint}$, where $M$ is a set, $G$ is a group, and $\actsOnPoint$ is a map from $G \times M$ to $M$, called \emph{left group action of $G$ on $M$}, such that $G \to \symmetricGroupOf(M)$, $g \mapsto [g \actsOnPoint \blank]$, is a group homomorphism. The action $\actsOnPoint$ is \emph{transitive} if $M$ is non-empty and for each $m \in M$ the map $\blank \actsOnPoint m$ is surjective; and \emph{free} if for each $m \in M$ the map $\blank \actsOnPoint m$ is injective. For each $m \in M$, the set $G \actsOnPoint m$ is the \emph{orbit of $m$}, the set $G_m = (\blank \actsOnPoint m)^{-1}(m)$ is the \emph{stabiliser of $m$}, and, for each $m' \in M$, the set $G_{m, m'} = (\blank \actsOnPoint m)^{-1}(m')$ is the \emph{transporter of $m$ to $m'$}.

  A \emph{left-ho\-mo\-ge\-neous space} is a left group set $\mathcal{M} = \ntuple{M, G, \actsOnPoint}$ such that $\actsOnPoint$ is transitive. A \emph{coordinate system for $\mathcal{M}$} is a tuple $\mathcal{K} = \ntuple{m_0, \family{g_{m_0, m}}_{m \in M}}$, where $m_0 \in M$ and, for each $m \in M$, we have $g_{m_0, m} \actsOnPoint m_0 = m$. The stabiliser $G_{m_0}$ is denoted by $G_0$. The tuple $\mathcal{R} = \ntuple{\mathcal{M}, \mathcal{K}}$ is a \emph{cell space}. The set $\setOf{g G_0 \suchThat g \in G}$ of left cosets of $G_0$ in $G$ is denoted by $G \modulo G_0$. The map $\isSemiActedUponBy \from M \times G \modulo G_0 \to M$, $(m, g G_0) \mapsto g_{m_0, m} g \actsOnPoint m_0$ is a \emph{right semi-action of $G \modulo G_0$ on $M$ with defect $G_0$}, which means that
  \begin{equation*}
    \ForEach m \in M \Holds m \isSemiActedUponBy G_0 = m,
  \end{equation*}
  and
  \begin{multline*}
    \ForEach m \in M \ForEach g \in G \Exists g_0 \in G_0 \SuchThat \ForEach \mathfrak{g}' \in G \modulo G_0 \Holds\\
          m \isSemiActedUponBy g \cdot \mathfrak{g}' = (m \isSemiActedUponBy g G_0) \isSemiActedUponBy g_0 \cdot \mathfrak{g}'.
  \end{multline*}
  It is \emph{transitive}, which means that the set $M$ is non-empty and for each $m \in M$ the map $m \isSemiActedUponBy \blank$ is surjective; and \emph{free}, which means that for each $m \in M$ the map $m \isSemiActedUponBy \blank$ is injective; and \emph{semi-commutes with $\actsOnPoint$}, which means that
  \begin{multline*}
    \ForEach m \in M \ForEach g \in G \Exists g_0 \in G_0 \SuchThat \ForEach \mathfrak{g}' \in G \modulo G_0 \Holds\\
          (g \actsOnPoint m) \isSemiActedUponBy \mathfrak{g}' = g \actsOnPoint (m \isSemiActedUponBy g_0 \cdot \mathfrak{g}').
  \end{multline*}

  A \emph{semi-cellular automaton} is a quadruple $\mathcal{C} = \ntuple{\mathcal{R}, Q, N, \delta}$, where $\mathcal{R}$ is a cell space; $Q$, called \emph{set of states}, is a set; $N$, called \emph{neighbourhood}, is a subset of $G \modulo G_0$ such that $G_0 \cdot N \subseteq N$; and $\delta$, called \emph{local transition function}, is a map from $Q^N$ to $Q$. A \emph{local configuration} is a map $\ell \in Q^N$ and a \emph{global configuration} is a map $c \in Q^M$. The stabiliser $G_0$ acts on $Q^N$ on the left by $\bullet \from G_0 \times Q^N \to Q^N$, $(g_0, \ell) \mapsto [n \mapsto \ell(g_0^{-1} \cdot n)]$, and the group $G$ acts on $Q^M$ on the left by $\actsOnMap \from G \times Q^M \to Q^M$, $(g, c) \mapsto [m \mapsto c(g^{-1} \actsOnPoint m)]$. The \emph{global transition function of $\mathcal{C}$} is the map $\Delta \from Q^M \to Q^M$, $c \mapsto [m \mapsto \delta(n \mapsto c(m \isSemiActedUponBy n))]$. A \emph{sufficient neighbourhood of $\mathcal{C}$} is a subset $E$ of $N$ such that, for each $\ell \in Q^N$ and each $\ell' \in Q^N$ with $\ell\restrictedTo_E = \ell'\restrictedTo_E$, we have $\delta(\ell) = \delta(\ell')$.

  A \define{cellular automaton} is a semi-cellular automaton $\mathcal{C} = \ntuple{\mathcal{R}, Q, N, \delta}$ such that the local transition function $\delta$ is \define{$\bullet$-invariant}, which means that, for each $g_0 \in G_0$, we have $\delta(g_0 \bullet \blank) = \delta(\blank)$. Its global transition function is $\actsOnMap$-e\-qui\-var\-i\-ant, which means that, for each $g \in G$, we have $\Delta(g \actsOnMap \blank) = g \actsOnMap \Delta(\blank)$. (See \cref{chapter:automata}.)

  \paragraph{Remark.} Because the notation in the present chapter is rather heavy, I only present the statements and proofs for semi-cellular and cellular but not for big-cellular automata, which would require the introduction of a further subgroup. Moreover, for lack of a better option, I overloaded some notation, but it should be clear from the context what is meant. 

  \section{Periodicity}
  \label{section:periodicity}

  \begin{definition}
    Let $\ntuple{M, G, \actsOnPoint}$ be a left group set, let $H$ be a subgroup of $G$, let $Q$ be a set, and let $c$ be a map from $M$ to $Q$. The map $c$ is called \defineX{$H$-periodic}{periodic@$H$-periodic}\graffito{$H$-periodic} if and only if
    \begin{equation*}
      \ForEach h \in H \Holds h \actsOnMap c = c. \qedhere
    \end{equation*}
  \end{definition}

  \begin{definition}
    Let $\ntuple{M, G, \actsOnPoint}$ be a left group set, let $H$ be a subgroup of $G$, and let $Q$ be a set. The \graffito{set $\periodic_H(Q^M)$ of $H$-periodic maps in $Q^M$}set of all $H$-periodic maps in $Q^M$ is denoted by $\periodic_H(Q^M)$\index[symbols]{PerHQM@$\periodic_H(Q^M)$}.
  \end{definition}

  \begin{remark}
    The set $\periodic_H(Q^M)$ is the set of fixed points of the left group action $\actsOnMap\restrictedTo_{H \times Q^M}$. 
  \end{remark}

  \begin{example}[Extreme Cases]
    The set $\periodic_{\setOf{e_G}}(Q^M)$ is equal to $Q^M$ and the set $\periodic_G(Q^M)$ is equal to the set of constant maps from $M$ to $Q$.
  \end{example} 

  \begin{example}[Cylinder]
  \label{example:42-lattice:periodic-configurations}
    Let $\Z_{42}$ be the additive cyclic group $\Z \modulo 42\Z$ of order $42$, let $M$ be the set $\Z_{42} \times \Z$, let $G$ be the additive group $\Z^2$, let $\actsOnPoint$ be the left group action of $G$ on $M$ by $((r_1, t_2), (m_1 + 42\Z, m_2)) \mapsto ((r_1 + m_1) + 42\Z, t_2 + m_2)$, let $H$ be the normal subgroup $\setOf{0} \times \Z$ of $G$, and let $Q$ be the binary set $\setOf{0, 1}$. The triple $\mathcal{M} = \ntuple{M, G, \actsOnPoint}$ is a left-ho\-mo\-ge\-neous space.

    Geometrically, the set $M$ is a discrete cylinder with axis $a$. Its rotations about $a$ can be encoded by $\Z_{42}$ and its translations along $a$ can be encoded by $\Z$. The group $G$ is a cover of these encodings, where its first component $\Z$ covers $\Z_{42}$. The rotational and translational symmetries themselves are $\setOf{(r_1, t_2) \actsOnPoint \blank \suchThat (r_1, t_2) \in G}$. Pictorially, maps from $M$ to $Q$ are black-and-white patterns over the cylinder $M$.

    For each black-and-white pattern $c \from M \to Q$, it is $H$-periodic if and only if 
    \begin{multline*}
      \ForEach t_2 \in \Z \ForEach (m_1 + 42\Z, m_2) \in M \Holds\\
          c(m_1 + 42\Z, m_2 + t_2) = c(m_1 + 42\Z, m_2),
    \end{multline*}
    in other words, if and only if the pattern $c$ is unicoloured along the axis $a$. 
    Therefore, the set $\periodic_H(Q^M)$ is equal to
    \begin{equation*}
      \setOf{c \in Q^M \suchThat \ForEach m_1 + 42\Z \in \Z_{42} \Holds c(m_1 + 42\Z, \blank) \text{ is constant}}. \qedhere
    \end{equation*}
  \end{example}

  The limit map of a convergent sequence of $H$-periodic maps with respect to the prodiscrete topology is $H$-periodic, which is shown in

  \begin{lemma} 
  \label{lemma:periodic-configurations-closed-in-phase-space}
    Let $\ntuple{M, G, \actsOnPoint}$ be a left group set, let $H$ be a subgroup of $G$, and let $Q$ be a set. Equip the set $Q^M$ with the prodiscrete topology. The set $\periodic_H(Q^M)$ is closed in $Q^M$.
  \end{lemma}

  \begin{proof}
    According to \cref{lemma:phase-space-is-Hausdorff}, the set $Q^M$ is Hausdorff, and, according to \cref{lemma:action-on-configurations-is-continuous}, the group action $\actsOnMap$ is continuous. Let $g \in G$. Then, the map $g \actsOnMap \blank$ is continuous. Thus, the map
    \begin{align*}
      \phi \from Q^M &\to Q^M \times Q^M,\\
      c &\mapsto (g \actsOnMap \blank, c),
    \end{align*}
    is continuous. Hence, because the diagonal $D = \setOf{(c,c) \suchThat c \in Q^M}$ is closed in $Q^M \times Q^M$, its preimage under $\phi$, namely
    \begin{equation*} 
      \phi^{-1}(D) = \setOf{c \in Q^M \suchThat g \actsOnMap c = c},
    \end{equation*}
    is closed in $Q^M$. In conclusion, because the intersection of closed sets is closed, the set
    \begin{equation*}
      \periodic_H(Q^M) = \bigcap_{h \in H} \setOf{c \in Q^M \suchThat h \actsOnMap c = c}
    \end{equation*}
    is closed in $Q^M$.
  \end{proof}

  \begin{example}[Cylinder]
  \label{example:42-lattice:periodic-configurations-closed}
    In the situation of \cref{example:42-lattice:periodic-configurations}, the limit patterns of convergent sequences of black-and-white patterns over the cylinder $M$ that are unicoloured along its axis $a$ are unicoloured along $a$, in other words, the set $\periodic_H(Q^M)$ is closed.
  \end{example}

  If we first project points to their orbits under $H$ and secondly map those orbits to states, the resulting map is $H$-periodic, because all points in the same orbit under $H$ are mapped to the same state, which is shown in

  \begin{lemma}
  \label{lemma:quotient-configuration-after-projection-is-periodic}
    Let $\ntuple{M, G, \actsOnPoint}$ be a left group set, let $H$ be a subgroup of $G$, let $H \reverseModulo M$ be the orbit space of $\actsOnPoint\restrictedTo_{H \times M}$, let $Q$ be a set, let $\pi$ be the canonical projection from $M$ onto $H \reverseModulo M$, and let $c_{\modulo H}$ be a map from $H \reverseModulo M$ to $Q$. The map $c_{\modulo H} \after \pi$ is $H$-periodic.
  \end{lemma}

  \begin{proof}
    For each $h \in H$ and each $m \in M$,
    \begin{align*}
      \parens[\big]{h \actsOnMap (c_{\modulo H} \after \pi)}(m)
      &= (c_{\modulo H} \after \pi)(h^{-1} \actsOnPoint m)\\
      &= c_{\modulo H}\parens[\big]{H \actsOnPoint (h^{-1} \actsOnPoint m)}\\
      &= c_{\modulo H}\parens[\big]{(H h^{-1}) \actsOnPoint m}\\
      &= c_{\modulo H}(H \actsOnPoint m)\\
      &= (c_{\modulo H} \after \pi)(m).
    \end{align*}
    In conclusion, the map $c_{\modulo H} \after \pi$ is $H$-periodic.
  \end{proof}

  All $H$-periodic maps can be constructed as in \cref{lemma:quotient-configuration-after-projection-is-periodic}, which is shown in

  \begin{lemma}
  \label{lemma:quotient-configurations-versus-periodic-ones}
    Let $\ntuple{M, G, \actsOnPoint}$ be a left group set, let $H$ be a subgroup of $G$, let $H \reverseModulo M$ be the orbit space of $\actsOnPoint\restrictedTo_{H \times M}$, let $Q$ be a set, and let $\pi$ be the canonical projection from $M$ onto $H \reverseModulo M$. The map
    \begin{align*}
      \pi_* \from Q^{H \reverseModulo M} &\to \periodic_H(Q^M), \mathnote{bijection $\pi_*$ from $Q^{H \reverseModulo M}$ to $\periodic_H(Q^M)$}\index[symbols]{pistar@$\pi_*$}\\
      c_{\modulo H} &\mapsto c_{\modulo H} \after \pi,
    \end{align*}
    is bijective.
  \end{lemma}

  \begin{proof}
    Note that, according to \cref{lemma:quotient-configuration-after-projection-is-periodic}, the map $\pi_*$ is well-defined.

    First, let $c_{\modulo H}$, $c_{\modulo H}' \in Q^{H \reverseModulo M}$ such that $c_{\modulo H} \neq c'_{\modulo H}$. Then, there is an $H \actsOnPoint m \in H \reverseModulo M$ such that $c_{\modulo H}(H \actsOnPoint m) \neq c'_{\modulo H}(H \actsOnPoint m)$. Hence,
    \begin{equation*}
      \pi_*(c_{\modulo H})(m) = c_{\modulo H}(\pi(m)) \neq c'_{\modulo H}(\pi(m)) = \pi_*(c'_{\modulo H})(m).
    \end{equation*}
    Therefore, $\pi_*(c_{\modulo H}) \neq \pi_*(c'_{\modulo H})$. In conclusion, the map $\pi_*$ is injective.

    Secondly, let $c \in \periodic_H(Q^M)$. Then, for each $h \in H$ and each $m \in M$,
    \begin{equation*}
      c(h \actsOnPoint m) = (h^{-1} \actsOnMap c)(m) = c(m).
    \end{equation*}
    Hence, the map
    \begin{align*}
      c_{\modulo H} \from H \reverseModulo M &\to Q,\\
      H \actsOnPoint m &\mapsto c(m),
    \end{align*}
    is well-defined. Moreover, $\pi_*(c_{\modulo H}) = c$. In conclusion, the map $\pi_*$ is surjective.
  \end{proof}

  It follows that the number of $H$-periodic maps is equal to the number of maps from the set of orbits under $H$, as stated in

  \begin{corollary}
    Let $\ntuple{M, G, \actsOnPoint}$ be a left group set, let $H$ be a subgroup of $G$ such that the orbit space $H \reverseModulo M$ of $\actsOnPoint\restrictedTo_{H \times M}$ is finite, and let $Q$ be a finite set. The set $\periodic_H(Q^M)$ is finite and $\cardinalityOf{\periodic_H(Q^M)} = \cardinalityOf{Q}^{\cardinalityOf{H \reverseModulo M}}$.
  \end{corollary}

  \begin{proof}
    This is a direct consequence of \cref{lemma:quotient-configurations-versus-periodic-ones}.
  \end{proof}

  \begin{example}[Cylinder]
  \label{example:42-lattice:after-projection-is-periodic}
    In the situation of \cref{example:42-lattice:periodic-configurations}, identify the orbit space $H \reverseModulo M$ with the discrete circle $\Z_{42}$ by $H + \ntuple{m_1 + 42\Z, m_2}\allowbreak \mapsto m_1 + 42\Z$ (note that $H + (m_1 + 42\Z, m_2) = \setOf{m_1 + 42\Z} \times \Z$). Then, the canonical projection $\pi$ is the map $M \to \Z_{42}$, $(m_1 + 42\Z, m_2) \mapsto m_1 + 42\Z$, which is the projection onto the first component. Hence, for each black-and-white pattern $c_{\modulo H}$ over $\Z_{42}$, the pattern $c_{\modulo H} \after \pi$ over $M$ is unicoloured along the axis $a$. And, the map $\pi_*$ extends black-and-white patterns over $\Z_{42}$ unicolouredly along $a$ to patterns over $M$. And, there are precisely $2^{42}$ black-and-white patterns that are unicoloured along $a$.
  \end{example}

  Equivariant maps preserve periodicity, as stated in

  \begin{lemma} 
  \label{lemma:fixed-points-invariance}
    Let $\ntuple{M, G, \actsOnPoint}$ and $\ntuple{M', G', \actsOnPoint'}$ be two left group sets, let $H$ be a subgroup of $G$, let $Q$ be a set, let $\Delta$ be a map from $Q^M$ to $Q^{M'}$, and let $\varphi$ be a group homomorphism from $G$ to $G'$ such that the tuple $(\Delta, \varphi)$ is $(\actsOnMap, \actsOnMap')$-e\-qui\-var\-i\-ant. Then, $\Delta(\periodic_H(Q^M)) \subseteq \periodic_{\varphi(H)}(Q^{M'})$.
  \end{lemma}

  \begin{proof}
    Let $c \in \periodic_H(Q^M)$. Then, for each $h \in H$,
    \begin{equation*}
      \varphi(h) \actsOnMap' \Delta(c)
      = \Delta(h \actsOnMap c)
      = \Delta(c).
    \end{equation*}
    Thus, $\Delta(c)$ is $\varphi(H)$-periodic.
  \end{proof}

  Global transition functions of cellular automata preserve periodicity, as stated in

  \begin{corollary} 
  \label{corollary:fixed-points-invariance}
    Let $\Delta$ be the global transition function of a cellular automaton over $\ntuple{M, G, \actsOnPoint}$ and let $H$ be a subgroup of $G$. The set $\periodic_H(Q^M)$ is invariant under $\Delta$.
  \end{corollary}

  \begin{proof}
    According to \cref{theorem:local-invariance-versus-global-equivariance}, the map $\Delta$ is $\actsOnMap$-e\-qui\-var\-i\-ant, in other words, $(\Delta, \identityMap_G)$ is $(\actsOnMap, \actsOnMap)$-e\-qui\-var\-i\-ant. Hence, according to lemma \ref{lemma:fixed-points-invariance}, we have $\Delta(\periodic_H(Q^M)) \subseteq \periodic_H(Q^M)$. 
  \end{proof}

  \begin{example}[Cylinder]
  \label{example:42-lattice:fixed-points-invariance}
    In the situation of \cref{example:42-lattice:periodic-configurations}, let $r + 42\Z$ be an element of $\Z_{42}$ and let $\Delta$ be the $\actsOnMap$-e\-qui\-var\-i\-ant map $Q^M \to Q^M$, $c \mapsto c(\blank + (r + 42\Z, 0))$, which rotates black-and-white patterns over the cylinder $M$ about its axis $a$. The map $\Delta$ maps patterns that are unicoloured along $a$ to such patterns.
  \end{example}


  Shifts preserve periodicity under normal groups, as stated in

  \begin{lemma} 
  \label{lemma:periodics-invariant-under-action-on-phase-space}
    Let $\ntuple{M, G, \actsOnPoint}$ be a left group set, let $H$ be a normal subgroup of $G$, and let $Q$ be a set. The set $\periodic_H(Q^M)$ is $\actsOnMap$-invariant.
  \end{lemma}

  \begin{proof}
    Let $g \in G$ and let $c \in \periodic_H(Q^M)$. Moreover, let $h \in H$. Then, because $g H = H g$, there is an $h' \in H$ such that $h g = g h'$. Hence,
    \begin{equation*}
      h \actsOnMap (g \actsOnMap c)
      = h g \actsOnMap c
      = g h' \actsOnMap c
      = g \actsOnMap (h' \actsOnMap c)
      = g \actsOnMap c.
    \end{equation*}
    Therefore, $g \actsOnMap c \in \periodic_H(Q^M)$. In conclusion, $\periodic_H(Q^M)$ is $\actsOnMap$-invariant.
  \end{proof}

  Quotient groups act on periodic maps on the left, as stated in

  \begin{lemma}
    Let $\ntuple{M, G, \actsOnPoint}$ be a left group set, let $H$ be a normal subgroup of $G$, and let $Q$ be a set. The map
    \begin{align*}
      \actsOnMap_{\periodic_H} \from G \modulo H \times \periodic_H(Q^M) &\to \periodic_H(Q^M), \mathnote{left group action $\actsOnMap_{\periodic_H}$ of $G \modulo H$ on $\periodic_H(Q^M)$}\index[symbols]{arrow right black periodic H@$\actsOnMap_{\periodic_H}$}\\
      (g H, c) &\mapsto g \actsOnMap c,
    \end{align*}
    is a left group action of $G \modulo H$ on $\periodic_H(Q^M)$.
  \end{lemma}

  \begin{proof}
    The map $\actsOnMap_{\periodic_H}$ is well-defined, because, according to lemma \ref{lemma:periodics-invariant-under-action-on-phase-space}, for each $g \in G$, each $h \in H$, and each $c \in \periodic_H(Q^M)$, 
    \begin{equation*}
      g h \actsOnMap c
      = g \actsOnMap (h \actsOnMap c)
      = g \actsOnMap c
      \in \periodic_H(Q^M).
    \end{equation*}
    And it is a left group action, because, for each $g H \in G \modulo H$ and each $g' H \in G \modulo H$, we have $g H g' H = g g' H$, and $\actsOnMap$ is a left group action.
  \end{proof}

  \begin{example}[Cylinder]
  \label{example:42-lattice:periodics-invariant-under-action-on-phase-space}
    In the situation of \cref{example:42-lattice:periodic-configurations}, the subgroup $H$ of $G$ is normal and we identify the quotient group $G \modulo H$ with $\Z$ by $(r_1, t_2) + H \mapsto r_1$ (note that $(r_1, t_2) + H = \setOf{r_1} \times \Z$). Then, the left group action $\actsOnMap$ is the map $G \times Q^M \to Q^M$, $((r_1, t_2), c) \mapsto c(\blank - (r_1 + 42\Z), \blank - t_2)$. It rotates and translates black-and-white patterns over the cylinder $M$ about and along its axis $a$. In particular, it leaves the set of patterns that are unicoloured along $a$ invariant. Moreover, the left quotient group action $\actsOnMap_{\periodic_H}$ is the map $\Z \times \periodic_H(Q^M) \to \periodic_H(Q^M)$, $(r_1, c) \mapsto c(\blank - (r_1 + 42\Z), \blank)$. It rotates black-and-white patterns over the cylinder $M$ that are unicoloured along its axis about its axis.
  \end{example}

  Restrictions of equivariant maps to periodic patterns are equivariant under induced quotient group actions, as stated in

  \begin{lemma} 
  \label{lemma:restriction-to-periodic-is-equivariant}
    Let $\ntuple{M, G, \actsOnPoint}$ and $\ntuple{M', G', \actsOnPoint'}$ be two left group sets, let $H$ be a normal subgroup of $G$, let $Q$ be a set, let $\Delta$ be a map from $Q^M$ to $Q^{M'}$, and let $\varphi$ be a surjective group homomorphism from $G$ to $G'$ such that the tuple $(\Delta, \varphi)$ is $(\actsOnMap, \actsOnMap')$-e\-qui\-var\-i\-ant. The restriction $\Delta\restrictedTo_{\periodic_H(Q^M) \to \periodic_H(Q^{M'})}$ is $(\actsOnMap_{\periodic_H}, \actsOnMap'_{\periodic_{\varphi(H)}})$-e\-qui\-var\-i\-ant.
  \end{lemma}

  \begin{proof}
    According to \cref{lemma:fixed-points-invariance}, the restriction $\Delta_{\periodic_H} = \Delta\restrictedTo_{\periodic_H(Q^M) \to \periodic_H(Q^{M'})}$ is well-defined. And, because $H$ is normal in $G$ and $\varphi$ is a surjective group homomorphism, $\varphi(H)$ is a normal subgroup of $G'$ and thus $\actsOnMap'_{\periodic_{\varphi(H)}}$ is well-defined. Moreover, for each $g H \in G \modulo H$ and each $c \in \periodic_H(Q^M)$,
    \begin{align*}
      \Delta_{\periodic_H}(g H \actsOnMap_{\periodic_H} c)
      &= \Delta(g \actsOnMap c)\\
      &= \varphi(g) \actsOnMap' \Delta(c)\\
      &= \varphi(g) \varphi(H) \actsOnMap'_{\periodic_{\varphi(H)}} \Delta_{\periodic_H}(c).
    \end{align*}
    In conclusion, $\Delta_{\periodic_H}$ is $(\actsOnMap_{\periodic_H}, \actsOnMap'_{\periodic_{\varphi(H)}})$-e\-qui\-var\-i\-ant.
  \end{proof}

  This holds in particular for global transition functions of cellular automata, as stated in

  \begin{corollary}
  \label{corollary:restriction-to-periodic-is-equivariant}
    Let $\Delta$ be the global transition function of a cellular automaton over $\ntuple{M, G, \actsOnPoint}$, let $H$ be a normal subgroup of $G$, and let $Q$ be a set. The restriction $\Delta\restrictedTo_{\periodic_H(Q^M) \to \periodic_H(Q^M)}$ is $\actsOnMap_{\periodic_H}$-e\-qui\-var\-i\-ant.
  \end{corollary}

  \begin{proof}
    According to \cref{theorem:local-invariance-versus-global-equivariance}, the map $\Delta$ is $\actsOnMap$-e\-qui\-var\-i\-ant, in other words, $(\Delta, \identityMap_G)$ is $(\actsOnMap, \actsOnMap)$-e\-qui\-var\-i\-ant. Hence, according to lemma \ref{lemma:restriction-to-periodic-is-equivariant}, the tuple $(\Delta\restrictedTo_{\periodic_H(Q^M) \to \periodic_H(Q^M)}, \identityMap_G)$ is $(\actsOnMap_{\periodic_H}, \actsOnMap_{\periodic_H})$-e\-qui\-var\-i\-ant, in other words, $\Delta\restrictedTo_{\periodic_H(Q^M) \to \periodic_H(Q^M)}$ is $\actsOnMap_{\periodic_H}$-e\-qui\-var\-i\-ant. 
  \end{proof}

  \begin{example}[Cylinder]
    In the situation of \cref{example:42-lattice:fixed-points-invariance}, the map $\Delta$ is $\actsOnMap$-e\-qui\-var\-i\-ant, which means, according to \cref{example:42-lattice:periodics-invariant-under-action-on-phase-space}, that it is equivariant under rotations about the axis $a$ of the cylinder $M$ and under translations along $a$. So, its restriction to black-and-white patterns that are unicoloured along $a$ is still equivariant under rotations about $a$, which is but $\actsOnMap_{\periodic_H}$-equivariance.
  \end{example}

  \section{Quotients}
  \label{section:quotients}

  Quotient groups act on orbit spaces on the left, as stated in

  \begin{lemma}
    Let $\mathcal{M} = \ntuple{M, G, \actsOnPoint}$ be a left group set, let $H$ be a normal subgroup of $G$, let $H \reverseModulo M$ be the orbit space of $\actsOnPoint\restrictedTo_{H \times M}$, and let
    \begin{align*}
      \actsOnPoint_{\modulo H} \from G \modulo H \times H \reverseModulo M &\to H \reverseModulo M, \mathnote{left group action $\actsOnPoint_{\modulo H}$ of $G \modulo H$ on $H \reverseModulo M$}\index[symbols]{arrow right modulo H@$\actsOnPoint_{\modulo H}$}\\
      (g H, H \actsOnPoint m) &\mapsto H \actsOnPoint (g \actsOnPoint m).
    \end{align*}
    The triple \graffito{quotient $\mathcal{M} \modulo H$ of $\mathcal{M}$ by $H$}$\mathcal{M} \modulo H = \ntuple{H \reverseModulo M, G \modulo H, \actsOnPoint_{\modulo H}}$ is a left group set and is called \define{quotient of $\mathcal{M}$ by $H$}\index[symbols]{MmoduloH@$\mathcal{M} \modulo H$}.
  \end{lemma}

  \begin{proof}
    First, let $g H$, $g' H \in G \modulo H$ such that $g H = g' H$ and let $H \actsOnPoint m$, $H \actsOnPoint m' \in H \reverseModulo M$ such that $H \actsOnPoint m = H \actsOnPoint m'$. Then, there is an $h \in H$ such that $m = h \actsOnPoint m'$. And, because $g H = g' H = H g'$, there is an $h' \in H$ such that $g h = h' g'$. Hence,
    \begin{align*}
      H \actsOnPoint (g \actsOnPoint m)
      &= H \actsOnPoint \parens[\big]{g \actsOnPoint (h \actsOnPoint m')}\\
      &= H \actsOnPoint (g h \actsOnPoint m')\\
      &= H \actsOnPoint (h' g' \actsOnPoint m')\\
      &= H \actsOnPoint \parens[\big]{h' \actsOnPoint (g' \actsOnPoint m')}\\
      &= H h' \actsOnPoint (g' \actsOnPoint m')\\
      &= H \actsOnPoint (g' \actsOnPoint m').
    \end{align*}
    In conclusion, the map $\actsOnPoint_{\modulo H}$ is well-defined.

    Secondly, for each $H \actsOnPoint m \in H \reverseModulo M$,
    \begin{align*}
      e_{G \modulo H} \actsOnPoint_{\modulo H} (H \actsOnPoint m)
      &= e_G H \actsOnPoint_{\modulo H} (H \actsOnPoint m)\\
      &= H \actsOnPoint (e_G \actsOnPoint m)\\
      &= H \actsOnPoint m.
    \end{align*}
    And, for each $g H \in G \modulo H$, each $g' H \in G \modulo H$, and each $H \actsOnPoint m \in H \reverseModulo M$,
    \begin{align*}
      g H \cdot g' H \actsOnPoint_{\modulo H} (H \actsOnPoint m)
      &= g g' H \actsOnPoint_{\modulo H} (H \actsOnPoint m)\\
      &= H \actsOnPoint (g g' \actsOnPoint m)\\
      &= H \actsOnPoint \parens[\big]{g \actsOnPoint (g' \actsOnPoint m)}\\
      &= g H \actsOnPoint_{\modulo H} \parens[\big]{H \actsOnPoint (g' \actsOnPoint m)}\\
      &= g H \actsOnPoint_{\modulo H} \parens[\big]{g' H \actsOnPoint (H \actsOnPoint m)}.
    \end{align*}
    In conclusion, the map $\actsOnPoint_{\modulo H}$ is a left group action.
  \end{proof}

  \begin{example}[Cylinder]
  \label{example:42-lattice:quotient:of-group-set} 
    In the situation of \cref{example:42-lattice:periodic-configurations}, the subgroup $H$ is normal in $G$, we identify the orbit space $H \reverseModulo M$ with the discrete circle $\Z_{42}$ by $H + (m_1 + 42\Z, m_2) \mapsto m_1 + 42\Z$ (note that $H + (m_1 + 42\Z, m_2) = \setOf{m_1 + 42\Z} \times \Z$), and we identify the quotient group $G \modulo H$ with $\Z$ by $(r_1, t_2) + H \mapsto r_1$ (note that $(r_1, t_2) + H = \setOf{r_1} \times \Z$). The left quotient group action $\actsOnPoint_{\modulo H}$ is the map $\Z \times \Z_{42} \to \Z_{42}$, $(r_1, m_1 + 42\Z) \mapsto (r_1 + m_1) + 42\Z$ and the quotient of $\mathcal{M}$ by $H$ is the triple $\ntuple{\Z_{42}, \Z, \actsOnPoint_{\modulo H}}$.
  \end{example} 

  Quotients of left-ho\-mo\-ge\-neous spaces are left-ho\-mo\-ge\-neous spaces, which is shown in

  \begin{lemma}
    Let $\mathcal{M} = \ntuple{M, G, \actsOnPoint}$ be a left-ho\-mo\-ge\-neous space and let $H$ be a normal subgroup of $G$. The quotient of $\mathcal{M}$ by $H$ is a left-ho\-mo\-ge\-neous space.
  \end{lemma}

  \begin{proof}
    For each $H \actsOnPoint m \in H \reverseModulo M$ and each $H \actsOnPoint m' \in H \reverseModulo M$, there is a $g \in G$ such that $g \actsOnPoint m = m'$, and hence $g H \actsOnPoint_{\modulo H} (H \actsOnPoint m) = H \actsOnPoint m'$. In conclusion, the group action $\actsOnPoint_{\modulo H}$ is transitive.
  \end{proof}

  \begin{example}[Cylinder]
  \label{example:42-lattice:quotient:homogeneous-space}
    In the situation of \cref{example:42-lattice:quotient:of-group-set}, the group set $\mathcal{M}$ is a left-ho\-mo\-ge\-neous space and so is the quotient $\mathcal{M} \modulo H$.
  \end{example}

  \begin{definition} 
    Let $G$ be a group and let $G_0$ and $H$ be two subgroups of $G$. The set $\setOf{g_0 H \suchThat g_0 \in G_0} \subseteq G \modulo H$ is denoted by \graffito{$G_0 \modulo H$}$G_0 \modulo H$\index[symbols]{G0moduloH@$G_0 \modulo H$}.
  \end{definition}

  \begin{remark}
    If $H$ is normal in $G$, then $G_0 \modulo H$ is a subgroup of $G \modulo H$.
  \end{remark}

  Stabilisers in quotients are quotients of stabilisers, as stated in

  \begin{lemma}
  \label{lemma:stabilisers-and-quotients}
    Let $\ntuple{H \reverseModulo M, G \modulo H, \actsOnPoint_{\modulo H}}$ be a quotient of $\ntuple{M, G, \actsOnPoint}$ by $H$ and let $m$ be an element of $M$. The stabiliser $(G \modulo H)_{H \actsOnPoint m}$ of $H \actsOnPoint m$ under $\actsOnPoint_{\modulo H}$ is $G_m \modulo H$.
  \end{lemma}

  \begin{proof}
    First, let $g_m H \in G_m \modulo H$. Then,
    \begin{equation*}
      g_m H \actsOnPoint_{\modulo H} (H \actsOnPoint m)
      = H \actsOnPoint (g_m \actsOnPoint m)
      = H \actsOnPoint m.
    \end{equation*}
    Hence, $g_m H \in (G \modulo H)_{H \actsOnPoint m}$. In conclusion, $G_m \modulo H \subseteq (G \modulo H)_{H \actsOnPoint m}$.

    Secondly, let $g H \in (G \modulo H)_{H \actsOnPoint m}$. Then,
    \begin{align*}
      H \actsOnPoint m
      = g H \actsOnPoint_{\modulo H} (H \actsOnPoint m)
      = H \actsOnPoint (g \actsOnPoint m).
    \end{align*}
    Hence, there is an $h \in H$ such that $h \actsOnPoint m = g \actsOnPoint m$. Thus, $h^{-1} g \actsOnPoint m = m$. Therefore, $g_m = h^{-1} g \in G_m$ and $g = h g_m$. Because $H$ is normal in $G$,
    \begin{equation*}
      g H = h g_m H = H h g_m = H g_m = g_m H \in G_m \modulo H.
    \end{equation*}
    In conclusion, $(G \modulo H)_{H \actsOnPoint m} \subseteq G_m \modulo H$. Altogether, $(G \modulo H)_{H \actsOnPoint m} = G_m \modulo H$.
  \end{proof}

  \begin{example}[Cylinder]
  \label{example:42-lattice:quotient:stabilisers}
    In the situation of \cref{example:42-lattice:quotient:homogeneous-space}, for each element $(m_1 + 42\Z, m_2) \in M$, the stabiliser of $(m_1, m_2)$ under $\actsOnPoint$ is the subgroup $42\Z \times \setOf{0}$ of $G$ and the stabiliser of $m_1 \simeq H \actsOnPoint (m_1, m_2)$ under $\actsOnPoint_{\modulo H}$ is the subgroup $42\Z$ of $\Z \simeq G \modulo H$.
  \end{example}

  Quotient actions on quotient patterns are essentially the original actions on periodic patterns, as stated in

  \begin{lemma}
    Let $\ntuple{H \reverseModulo M, G \modulo H, \actsOnPoint_{\modulo H}}$ be a quotient of $\ntuple{M, G, \actsOnPoint}$ by $H$, let $Q$ be a set, and let $\pi$ be the canonical projection from $M$ onto $H \reverseModulo M$. For each symmetry $g \in G$, each map $c_{\modulo H} \in Q^{H \reverseModulo M}$, and each point $m \in M$,
    \begin{equation*}
      (g H \actsOnMap_{\modulo H} c_{\modulo H})(H \actsOnPoint m) = \parens[\big]{g \actsOnMap \pi_*(c_{\modulo H})}(m). \qedhere
    \end{equation*}
  \end{lemma}

  \begin{proof}
    For each $g \in G$, each $c_{\modulo H} \in Q^{H \reverseModulo M}$, and each $m \in M$,
    \begin{align*}
      (g H \actsOnMap_{\modulo H} c_{\modulo H})(H \actsOnPoint m)
      &= c_{\modulo H}\parens[\big]{g^{-1} H \actsOnPoint_{\modulo H} (H \actsOnPoint m)}\\
      &= c_{\modulo H}\parens[\big]{H \actsOnPoint (g^{-1} \actsOnPoint m)}\\
      &= (c_{\modulo H} \after \pi)(g^{-1} \actsOnPoint m)\\
      &= \pi_*(c_{\modulo H})(g^{-1} \actsOnPoint m)\\
      &= \parens[\big]{g \actsOnMap \pi_*(c_{\modulo H})}(m). \qedhere
    \end{align*}
  \end{proof}


  The projection $\pi_*$ is equivariant, as stated in

  \begin{lemma}
  \label{lemma:rho-star-is-equivariant}
    Let $\ntuple{H \reverseModulo M, G \modulo H, \actsOnPoint_{\modulo H}}$ be a quotient of $\ntuple{M, G, \actsOnPoint}$ by $H$, let $Q$ be a set, and let $\pi$ be the canonical projection from $M$ onto $H \reverseModulo M$. The map $\pi_*$ is $(\actsOnMap_{\modulo H}, \actsOnMap_{\periodic_H})$-e\-qui\-var\-i\-ant.
  \end{lemma}

  \begin{proof}
    For each $g H \in G \modulo H$, each $c_{\modulo H} \in Q^{H \reverseModulo M}$, and each $m \in M$,
    \begin{align*}
      \pi_*(g H \actsOnMap_{\modulo H} c_{\modulo H})(m)
      &= (g H \actsOnMap_{\modulo H} c_{\modulo H})\parens[\big]{\pi(m)}\\
      &= (g H \actsOnMap_{\modulo H} c_{\modulo H})(H \actsOnPoint m)\\
      &= c_{\modulo H}\parens[\big]{g^{-1} H \actsOnPoint_{\modulo H} (H \actsOnPoint m)}\\
      &= c_{\modulo H}\parens[\big]{H \actsOnPoint (g^{-1} \actsOnPoint m)}\\
      &= c_{\modulo H}\parens[\big]{\pi(g^{-1} \actsOnPoint m)}\\
      &= \pi_*(c_{\modulo H})(g^{-1} \actsOnPoint m)\\
      &= \parens[\big]{g \actsOnMap \pi_*(c_{\modulo H})}(m)\\
      &= \parens[\big]{g H \actsOnMap_{\periodic_H} \pi_*(c_{\modulo H})}(m). \qedhere
    \end{align*}
  \end{proof}

  Given a right inverse of the canonical projection onto the orbit space, a coordinate system induces one on the quotient space, as stated in

  \begin{lemma}
    Let $\mathcal{M} = \ntuple{M, G, \actsOnPoint}$ be a left-ho\-mo\-ge\-neous space, let $\mathcal{K} = \ntuple{m_0, \family{g_{m_0, m}}_{m \in M}}$ be a coordinate system for $\mathcal{M}$, let $\mathcal{R}$ be the cell space $\ntuple{\mathcal{M}, \mathcal{K}}$, let $H$ be a normal subgroup of $G$, and let $\rho$ be a right inverse of the canonical projection $\pi \from M \to H \reverseModulo M$ such that $\rho(H \actsOnPoint m_0) = m_0$. The tuple \graffito{quotient $\mathcal{K} \modulo (H, \rho)$ of $\mathcal{K}$ by $H$ and $\rho$}$\mathcal{K} \modulo (H, \rho) = \ntuple{H \actsOnPoint m_0, \family{g_{m_0, \rho(H \actsOnPoint m)}}_{H \actsOnPoint m \in H \reverseModulo M}}$ is a coordinate system for $\mathcal{M} \modulo H$ and is called \define{quotient of $\mathcal{K}$ by $H$ and $\rho$}\index[symbols]{KmoduloHrcalligraphic@$\mathcal{K} \modulo (H, \rho)$}. And, the tuple \graffito{quotient $\mathcal{R} \modulo (H, \rho)$ of $\mathcal{R}$ by $H$ and $\rho$}$\mathcal{R} \modulo (H, \rho) = \ntuple{\mathcal{M} \modulo H, \mathcal{K} \modulo (H, \rho)}$ is a cell space and is called \define{quotient of $\mathcal{R}$ by $H$ and $\rho$}\index[symbols]{RmoduloHrcalligraphic@$\mathcal{R} \modulo (H, \rho)$}.
  \end{lemma}

  \begin{proof}
    Because $\rho(H \actsOnPoint m_0) = m_0$ and $g_{m_0, m_0} = e_G$, we have $g_{m_0, \rho(H \actsOnPoint m_0)} H = e_G H = e_{G \modulo H}$. And, for each $H \actsOnPoint m \in H \reverseModulo M$,
    \begin{align*}
      g_{m_0, \rho(H \actsOnPoint m)} H \actsOnPoint_{\modulo H} (H \actsOnPoint m_0)
      &= H \actsOnPoint (g_{m_0, \rho(H \actsOnPoint m)} \actsOnPoint m_0)\\
      &= H \actsOnPoint \rho(H \actsOnPoint m)\\
      &= H \actsOnPoint m.
    \end{align*}
    In conclusion, $\mathcal{K} \modulo (H, \rho)$ is a coordinate system for $\mathcal{M} \modulo H$ and $\mathcal{R} \modulo (H, \rho)$ is a cell space.
  \end{proof}

  \begin{example}[Cylinder]
  \label{example:42-lattice:quotient:coordinate-system}
    In the situation of \cref{example:42-lattice:quotient:stabilisers}, the tuple $\mathcal{K} = \ntuple{(0 + 42\Z, 0), \family{(m_1 \bmod 42, m_2)}_{(m_1 + 42\Z, m_2) \in M}}$ is a coordinate system for $\ntuple{M, G, \actsOnPoint}$, where $\blank \bmod 42$ denotes the remainder of the Euclidean division by $42$; the tuple $\mathcal{R} = \ntuple{\mathcal{M}, \mathcal{K}}$ is a cell space; the canonical projection $\pi$ from $M$ onto $H \reverseModulo M \simeq \Z_{42}$ is given by $(m_1 + 42\Z, m_2) \mapsto m_1 + 42\Z$; a right inverse $\rho$ of $\pi$ is given by $m_1 + 42\Z \mapsto (m_1 + 42\Z, 0)$; the quotient of $\mathcal{K}$ by $H$ and $\rho$ is the tuple $\ntuple{0 + 42\Z, \family{m_1 \bmod 42}_{m_1 \in \Z}}$, which is a coordinate system for $\mathcal{M} \modulo H$.
  \end{example}

  The right semi-actions induced by quotients of cell spaces can be expressed in terms of the right semi-actions induced by the cell spaces themselves, as stated in

  \begin{lemma} 
    Let $\mathcal{R} \modulo (H, \rho)$ be a quotient of $\ntuple{\ntuple{M, G, \actsOnPoint}, \ntuple{m_0, \family{g_{m_0, m}}_{m \in M}}}$ by $H$ and $\rho$. The right quotient set semi-action of $(G \modulo H) \modulo (G_0 \modulo H)$ on $H \reverseModulo M$ with defect $G_0 \modulo H$ induced by $\mathcal{R} \modulo (H, \rho)$ is
    \begin{align*}
      \isSemiActedUponBy_{\modulo (H, \rho)} \from H \reverseModulo M \times (G \modulo H) \modulo (G_0 \modulo H) &\to H \reverseModulo M, \mathnote{right quotient set semi-action $\isSemiActedUponBy_{\modulo (H, \rho)}$ induced by $\mathcal{R} \modulo (H, \rho)$}\index[symbols]{arrow left scored modulo H rho@$\isSemiActedUponBy_{\modulo (H, \rho)}$}\\
      \parens[\big]{H \actsOnPoint m, g H (G_0 \modulo H)} &\mapsto H \actsOnPoint \parens[\big]{\rho(H \actsOnPoint m) \isSemiActedUponBy g G_0}. \qedhere
    \end{align*}
  \end{lemma}

  \begin{proof}
    According to \cref{lemma:stabilisers-and-quotients}, we have $(G \modulo H)_{H \actsOnPoint m_0} = G_0 \modulo H$. Moreover, for each $H \actsOnPoint m \in H \reverseModulo M$ and each $g H (G_0 \modulo H) \in (G \modulo H) \modulo (G_0 \modulo H)$,
    \begin{align*}
      &(H \actsOnPoint m) \isSemiActedUponBy_{\modulo (H, \rho)} g H (G_0 \modulo H)\\
      &= H \actsOnPoint \parens[\big]{\rho(H \actsOnPoint m) \isSemiActedUponBy g G_0}\\
      &= H \actsOnPoint \parens[\big]{g_{m_0, \rho(H \actsOnPoint m)} g g_{m_0, \rho(H \actsOnPoint m)}^{-1} \actsOnPoint \rho(H \actsOnPoint m)}\\
      &= g_{m_0, \rho(H \actsOnPoint m)} g g_{m_0, \rho(H \actsOnPoint m)}^{-1} H \actsOnPoint_{\modulo H} \parens[\big]{H \actsOnPoint \parens[\big]{\rho(H \actsOnPoint m)}}\\
      &= g_{m_0, \rho(H \actsOnPoint m)} H \cdot g H \cdot g_{m_0, \rho(H \actsOnPoint m)}^{-1} H \actsOnPoint_{\modulo H} (H \actsOnPoint m).
    \end{align*}
    Hence, $\isSemiActedUponBy_{\modulo (H, \rho)}$ is well-defined and the right quotient set semi-action induced by $\mathcal{R} \modulo (H, \rho)$.
  \end{proof}

  \begin{example}[Cylinder]
  \label{example:42-lattice:quotient:right-semi-action}
    In the situation of \cref{example:42-lattice:quotient:coordinate-system}, identify the quotient set $G \modulo G_0$ with $\Z_{42} \times \Z$ by $(r_1, t_2) + G_0 \mapsto (r_1 + 42\Z, t_2)$. Then, the right quotient set semi-action $\isSemiActedUponBy$ induced by $\mathcal{R}$ is the map $M \times (\Z_{42} \times \Z) \to M$, $((m_1 + 42\Z, m_2), (r_1 + 42\Z, t_2)) \mapsto ((m_1 + r_1) + 42\Z, m_2 + t_2)$, which is addition in the direct product $\Z_{42} \times \Z$.

    Recall that the quotient group $G \modulo H$ is identified with $\Z$. Hence, its subgroup $G_0 \modulo H$ is identified with $42\Z$ and the quotient set $(G \modulo H) \modulo (G_0 \modulo H)$ is identified with $\Z_{42}$. Also recall that the orbit space $H \reverseModulo M$ is identified with $\Z_{42}$. Therefore, the right quotient set semi-action $\isSemiActedUponBy_{\modulo (H, \rho)}$ induced by $\mathcal{R} \modulo (H, \rho)$ is the map $\Z_{42} \times \Z_{42} \to \Z_{42}$, $(m_1 + 42\Z, r_1 + 42\Z) \mapsto (m_1 + r_1) + 42\Z$, which is addition in $\Z_{42}$.
  \end{example}

  The quotient of a cell space by a group that includes the stabiliser of the origin does not depend on the chosen right inverse of the projection, as stated in

  \begin{lemma} 
  \label{lemma:independence-of-quotient-loc-lib-from-r}
    Let $\mathcal{R} \modulo (H, \rho)$ be a quotient of $\ntuple{\mathcal{M}, \mathcal{K}} = \ntuple{\ntuple{M, G, \actsOnPoint}, \ntuple{m_0, \family{g_{m_0, m}}_{m \in M}}}$ by $H$ and $\rho$ such that the stabiliser $G_0$ is included in $H$. The quotient $\mathcal{R} \modulo (H, \rho)$ and the right quotient set semi-action $\isSemiActedUponBy_{\modulo (H, \rho)}$ do not depend on $\rho$.
  \end{lemma}

  \begin{proof}
    Let $m$, $m' \in M$ such that $H \actsOnPoint m = H \actsOnPoint m'$. Then, there is an $h \in H$ such that $h \actsOnPoint m = m'$. Thus, $h g_{m_0, m} \in G_{m_0, m'}$. Hence, because $G_{m_0, m'} = G_0 g_{m_0, m'}$, there is a $g_0 \in G_0$ such that $h g_{m_0, m} = g_0 g_{m_0, m'}$. Therefore, because $H$ is normal, $h H = H$, and $H g_0 = H$,
    \begin{multline*}
      g_{m_0, m} H
      = H g_{m_0, m}
      = H h g_{m_0, m}\\
      = H g_0 g_{m_0, m'}
      = H g_{m_0, m'}
      = g_{m_0, m'} H.
    \end{multline*}
    Hence, for each $m \in M$, we have $g_{m_0, m} H = g_{\rho(H \actsOnPoint m),m_0} H$. Therefore, the coordinates $g_{m_0, \rho(H \actsOnPoint m)} H$, for $H \actsOnPoint m \in H \reverseModulo M$, do not depend on $\rho$. In conclusion, $\mathcal{R} \modulo (H, \rho)$ and $\isSemiActedUponBy_{\modulo (H, \rho)}$ do not depend on $\rho$.
  \end{proof}

  \begin{remark}
    Because $G_0$ is included in $H$, the stabiliser $G_0 \modulo H$ of $H \actsOnPoint m_0$ under $\actsOnPoint_{\modulo H}$ is trivial, hence the left group action $\actsOnPoint_{\modulo H}$ is free, and therefore $\mathcal{M} \modulo H$ is a principal left-ho\-mo\-ge\-neous space and $\mathcal{K} \modulo (H, \rho)$ is its unique coordinate system for the origin $m_0$.
  \end{remark}

  The neighbourhood of a semi-cellular automaton over $M$ is a subset of $G \modulo G_0$ and the neighbourhood of one over $H \reverseModulo M$ is a subset of $(G \modulo H) \modulo (G_0 \modulo H)$. There is a canonical projection from $G \modulo G_0$ onto $(G \modulo H) \modulo (G_0 \modulo H)$, which under suitable identifications of $G \modulo G_0$ with $M$ and of $(G \modulo H) \modulo (G_0 \modulo H)$ with $H \reverseModulo M$ is the canonical projection from $M$ onto $H \reverseModulo M$, as stated in 

  \begin{lemma}
    Let $G$ be a group, let $G_0$ be a subgroup of $G$, and let $H$ be a normal subgroup of $G$. The map
    \begin{align*}
      \varpi \from G \modulo G_0 &\to (G \modulo H) \modulo (G_0 \modulo H), \mathnote{canonical projection $\varpi$ of $G \modulo G_0$ onto $(G \modulo H) \modulo (G_0 \modulo H)$}\index[symbols]{pivar@$\varpi$}\index[symbols]{pivar@$\varpi$}\\
      g G_0 &\mapsto g H (G_0 \modulo H),
    \end{align*}
    is well-defined, is surjective, and is called \define{canonical projection from $G \modulo G_0$ onto $(G \modulo H) \modulo (G_0 \modulo H)$}.
  \end{lemma}

  \begin{proof}
    Let $g G_0$, $g' G_0 \in G \modulo G_0$ such that $g G_0 = g' G_0$. Then, $g^{-1} g' \in G_0$. Hence, $g^{-1} H \cdot g' H = g^{-1} g' H \in G_0 \modulo H$. Therefore, $g H (G_0 \modulo H) = g' H (G_0 \modulo H)$. In conclusion, $\varpi$ is well-defined. Moreover, because, for each $g H (G_0 \modulo H) \in (G \modulo H) \modulo (G_0 \modulo H)$, we have $\varpi(g G_0) = g H (G_0 \modulo H)$, the map $\varpi$ is surjective.
  \end{proof}

  \begin{example}[Cylinder]
  \label{example:42-lattice:quotient:projection-of-quotient-set}
    In the situation of \cref{example:42-lattice:quotient:right-semi-action}, the canonical projection $\varpi$ is the map $\Z_{42} \times \Z \to \Z_{42}$, $(r_1 + 42\Z, t_2) \mapsto r_1 + 42\Z$, which is the projection to the first component.
  \end{example}


  Semi-cellular automata can be projected onto ones over orbit spaces, as shown in

  \begin{lemma}
    Let $\mathcal{R} \modulo (H, \rho)$ be a quotient of $\ntuple{\ntuple{M, G, \actsOnPoint}, \ntuple{m_0, \family{g_{m_0, m}}_{m \in M}}}$ by $H$ and $\rho$, let $\mathcal{C} = \ntuple{\mathcal{R}, Q, N, \delta}$ be a semi-cellular or cellular automaton, let $\varpi$ be the canonical projection from $G \modulo G_0$ onto $(G \modulo H) \modulo (G_0 \modulo H)$, let
    \begin{equation*}
      N_{\modulo H} = \varpi(N) = \setOf{g H (G_0 \modulo H) \suchThat g G_0 \in N},
    \end{equation*}
    and let
    \begin{align*}
      \delta_{\modulo H} \from Q^{N_{\modulo H}} &\to Q,\\
      \ell_{\modulo H} &\mapsto \delta\parens[\big]{n \mapsto \ell_{\modulo H}\parens[\big]{\varpi(n)}} \quad \parens{= \delta\parens[\big]{g G_0 \mapsto \ell_{\modulo H}\parens[\big]{g H (G_0 \modulo H)}}}.
    \end{align*}
    The quadruple \graffito{quotient $\mathcal{C} \modulo (H, \rho)$ of $\mathcal{C}$ by $H$ and $\rho$}$\mathcal{C} \modulo (H, \rho) = \ntuple{\mathcal{R} \modulo (H, \rho), Q, N_{\modulo H}, \delta_{\modulo H}}$ is a semi-cellular or cellular automaton respectively and is called \define{quotient of $\mathcal{C}$ by $H$ and $\rho$}\index[symbols]{CmoduloHrcalligraphic@$\mathcal{C} \modulo (H, \rho)$}.
  \end{lemma}

  \begin{proof}
    For each $g_0 \in G_0$ and each $g G_0 \in N$, we have $g_0 g G_0 \in N$. Hence,
    \begin{align*}
      (G_0 \modulo H) \cdot N_{\modulo H}
      &= \setOf{g_0 H \suchThat g_0 \in G_0} \cdot \setOf{g H (G_0 \modulo H) \suchThat g G_0 \in N}\\
      &= \setOf{g_0 H \cdot g H (G_0 \modulo H) \suchThat g_0 \in G_0, g G_0 \in N}\\
      &= \setOf{g_0 g H (G_0 \modulo H) \suchThat g_0 \in G_0, g G_0 \in N}\\
      &\subseteq \setOf{g H (G_0 \modulo H) \suchThat g G_0 \in N}\\
      &= N_{\modulo H}.
    \end{align*}
    Therefore, the quadruple $\mathcal{C} \modulo (H, \rho)$ is a semi-cellular automaton. From now on, let $\mathcal{C}$ be a cellular automaton. Furthermore, let $g_0 H \in G_0 \modulo H$ and let $\ell_{\modulo H} \in Q^{N_{\modulo H}}$. Then,
    \begin{align*}
      \delta_{\modulo H}(g_0 H \bullet_{\modulo H} \ell_{\modulo H})
      &= \delta\parens[\big]{g G_0 \mapsto (g_0 H \bullet_{\modulo H} \ell_{\modulo H})(\varpi(g G_0))}\\
      &= \delta\parens[\big]{g G_0 \mapsto (g_0 H \bullet_{\modulo H} \ell_{\modulo H})(g H (G_0 \modulo H))}\\
      &= \delta\parens[\big]{g G_0 \mapsto \ell_{\modulo H}(g_0^{-1} H \cdot g H (G_0 \modulo H))}\\
      &= \delta\parens[\big]{g G_0 \mapsto \ell_{\modulo H}(g_0^{-1} g H (G_0 \modulo H))}\\
      &= \delta\parens[\big]{g G_0 \mapsto \ell_{\modulo H}(\varpi(g_0^{-1} g G_0))}\\
      &= \delta\parens[\big]{g_0 \bullet \brackets[\big]{g G_0 \mapsto \ell_{\modulo H}(\varpi(g G_0))}}\\
      &= \delta\parens[\big]{g G_0 \mapsto \ell_{\modulo H}(\varpi(g G_0))}\\
      &= \delta_{\modulo H}(\ell_{\modulo H}).
    \end{align*}
    Hence, $\delta_{\modulo H}$ is $\bullet_{\modulo H}$-invariant. In conclusion, $\mathcal{C} \modulo (H, \rho)$ is a cellular automaton.
  \end{proof}

  \begin{example}[Cylinder]
  \label{example:42-lattice:quotient:shift-automaton}
    In the situation of \cref{example:42-lattice:quotient:projection-of-quotient-set}, let $Q$ be the binary set $\setOf{0, 1}$, let $z$ be an integer, let $N$ be the singleton set $\setOf{(-1 + 42\Z, z)}$, let $\delta$ be the $\bullet$-invariant map $Q^N \to Q$, $\ell \mapsto \ell(-1 + 42\Z, z)$, and let $\mathcal{C}$ be the cellular automaton $\ntuple{\mathcal{R}, Q, N, \delta}$. The global transition function $\Delta$ of $\mathcal{C}$ is a shift over $M$. For example, if $z = -1$, then it is the diagonal shift from the bottom-left to the top-right.

    The quotient $\mathcal{C} \modulo (H, \rho)$ is the quadruple $\ntuple{\mathcal{R} \modulo (H, \rho), Q, N_{\modulo H}, \delta_{\modulo H}}$, where the neighbourhood $N_{\modulo H}$ is the singleton set $\setOf{-1 + 42\Z}$ and the local transition function $\delta_{\modulo H}$ is the map $Q^{N_{\modulo H}} \to Q$, $\ell_{\modulo H} \mapsto \ell_{\modulo H}(-1 + 42\Z)$. The global transition function $\Delta_{\modulo H}$ of $\mathcal{C} \modulo (H, \rho)$ is the shift from left to right over $\Z_{42}$.
  \end{example}

  The global transition function of the quotient of a semi-cellular automaton can be expressed in terms of the global transition function of the automaton, as shown in

  \begin{lemma}
  \label{lemma:global-transition-function-of-quotient}
    Let $\mathcal{R} \modulo (H, \rho)$ be a quotient of $\ntuple{\ntuple{M, G, \actsOnPoint}, \ntuple{m_0, \family{g_{m_0, m}}_{m \in M}}}$ by $H$ and $\rho$, let $\mathcal{C} = \ntuple{\mathcal{R}, Q, N, \delta}$ be a semi-cellular automaton, and let $\pi$ be the canonical projection from $M$ onto $H \reverseModulo M$. The global transition function $\Delta_{\modulo H}$ of $\mathcal{C} \modulo (H, \rho)$ is identical to $\pi_*^{-1} \after \Delta\restrictedTo_{\periodic_H(Q^M) \to \periodic_H(Q^M)} \after \pi_*$, in particular, it does not depend on $\rho$ and it is uniquely determined by $\ntuple{M, G, \actsOnPoint}$, $H$, and $\Delta$.
  \end{lemma}

  \begin{proof}
    Let $\Delta_{\periodic_H} = \Delta\restrictedTo_{\periodic_H(Q^M) \to \periodic_H(Q^M)}$. Then, for each $c_{\modulo H} \in Q^{H \reverseModulo M}$ and each $H \actsOnPoint m \in H \reverseModulo M$,
    \begin{align*}
      &\Delta_{\modulo H}(c_{\modulo H})(H \actsOnPoint m)\\
      &= \delta_{\modulo H}(n_{\modulo H} \mapsto c_{\modulo H}((H \actsOnPoint m) \isSemiActedUponBy_{\modulo (H, \rho)} n_{\modulo H}))\\
      &= \delta_{\modulo H}(g H (G_0 \modulo H) \mapsto c_{\modulo H}((H \actsOnPoint m) \isSemiActedUponBy_{\modulo (H, \rho)} g H (G_0 \modulo H)))\\
      &= \delta(g G_0 \mapsto c_{\modulo H}((H \actsOnPoint m) \isSemiActedUponBy_{\modulo (H, \rho)} g H (G_0 \modulo H)))\\
      &= \delta(g G_0 \mapsto c_{\modulo H}(H \actsOnPoint (\rho(H \actsOnPoint m) \isSemiActedUponBy g G_0)))\\
      &= \delta(n \mapsto c_{\modulo H}(H \actsOnPoint (\rho(H \actsOnPoint m) \isSemiActedUponBy n)))\\
      &= \delta(n \mapsto (c_{\modulo H} \after \pi)((\rho(H \actsOnPoint m) \isSemiActedUponBy n)))\\
      &= \Delta(c_{\modulo H} \after \pi)(\rho(H \actsOnPoint m))\\
      &= \Delta(\pi_*(c_{\modulo H}))(\rho(H \actsOnPoint m))\\
      &= \Delta_{\periodic_H}(\pi_*(c_{\modulo H}))(\rho(H \actsOnPoint m))\\
      &= \pi_*^{-1}(\Delta_{\periodic_H}(\pi_*(c_{\modulo H})))(H \actsOnPoint \rho(H \actsOnPoint m))\\
      &= \pi_*^{-1}(\Delta_{\periodic_H}(\pi_*(c_{\modulo H})))(H \actsOnPoint m)\\
      &= (\pi_*^{-1} \after \Delta_{\periodic_H} \after \pi_*)(c_{\modulo H})(H \actsOnPoint m).
    \end{align*}
    Therefore, $\Delta_{\modulo H} = \pi_*^{-1} \after \Delta_{\periodic_H} \after \pi_*$, which does not depend on $\rho$.
  \end{proof}

  \begin{remark} 
    Let $\mathcal{C}$ be a cellular automaton. Then, according to \cref{corollary:restriction-to-periodic-is-equivariant}, the restriction $\Delta\restrictedTo_{\periodic_H(Q^M) \to \periodic_H(Q^M)}$ is $\actsOnMap_{\periodic_H}$-e\-qui\-var\-i\-ant, and, according to \cref{lemma:rho-star-is-equivariant}, the map $\pi_*$ is $(\actsOnMap_{\modulo H}, \actsOnMap_{\periodic_H})$-e\-qui\-var\-i\-ant. Hence, because $\Delta_{\modulo H}$ is identical to $\pi_*^{-1} \after \Delta\restrictedTo_{\periodic_H(Q^M) \to \periodic_H(Q^M)} \after \pi_*$, the global transition function $\Delta_{\modulo H}$ is $\actsOnMap_{\modulo H}$-e\-qui\-var\-i\-ant.
  \end{remark}


  \begin{definition}
  \label{definition:quotient-of-cellular-automaton-function}
    Let $\Delta$ be the global transition function of a semi-cellular or cellular automaton over $\ntuple{M, G, \actsOnPoint}$ and let $H$ be a normal subgroup of $G$. The map \graffito{quotient $\Delta_{\modulo H}$ of $\Delta$ by $H$}$\Delta_{\modulo H} = \pi_*^{-1} \after \Delta\restrictedTo_{\periodic_H(Q^M) \to \periodic_H(Q^M)} \after \pi_*$ is the global transition function of a semi-cellular or cellular automaton over $\ntuple{H \reverseModulo M, G \modulo H, \actsOnPoint_{\modulo H}}$ and is called \define{quotient of $\Delta$ by $H$}\index[symbols]{DeltamoduloH@$\Delta_{\modulo H}$}.
  \end{definition}


  \begin{example}[Cylinder]
  \label{example:42-lattice:quotient:shift-automaton-nice-representation}
    In the situation of \cref{example:42-lattice:quotient:shift-automaton}, according to \cref{example:42-lattice:after-projection-is-periodic}, the canonical projection $\pi$ is the map $M \to \Z_{42}$, $(m_1 + 42\Z, m_2) \mapsto m_1 + 42\Z$ and the map $\pi_*$ is the bijection $Q^{\Z_{42}} \to \periodic_H(Q^M)$, $c_{\modulo H} \mapsto c_{\modulo H} \after \pi$. Hence, the inverse $\pi_*^{-1}$ is the bijection $\periodic_H(Q^M) \to Q^{\Z_{42}}$, $c \mapsto c(\blank, 0)$. Therefore, global transition function $\Delta_{\modulo H}$ is the map $Q^{\Z_{42}} \to Q^{\Z_{42}}$, $c_{\modulo H} \mapsto \Delta(c_{\modulo H} \after \pi)(\blank, 0)$.
  \end{example}


  \section{Naïve Products}
  \label{section:naive-products}

  Given a right inverse of the canonical projection from $G \modulo G_0$ onto $(G \modulo H) \modulo (G_0 \modulo H)$, a semi-cellular automaton over a quotient space can be extended to the original space, as shown in

  \begin{definition}
  \label{definition:naive-product-of-cellular-automaton}
    Let $\mathcal{R} \modulo (H, \rho)$ be a quotient of $\ntuple{\ntuple{M, G, \actsOnPoint}, \ntuple{m_0, \family{g_{m_0, m}}_{m \in M}}}$ by $H$ and $\rho$, let $\mathcal{C}_{\modulo H} = \ntuple{\mathcal{R} \modulo (H, \rho), Q, N_{\modulo H}, \delta_{\modulo H}}$ be a semi-cellular automaton, let $\kappa$ be a right inverse of the canonical projection $\varpi \from G \modulo G_0 \to (G \modulo H) \modulo (G_0 \modulo H)$, let $N = G_0 \cdot \kappa(N_{\modulo H})$, let $\varkappa = \kappa\restrictedTo_{N_{\modulo H} \to N}$, and let
    \begin{align*}
      \delta \from Q^N &\to Q,\\
      \ell &\mapsto \delta_{\modulo H}\parens[\big]{n_{\modulo H} \mapsto \ell\parens[\big]{\kappa(n_{\modulo H})}} \quad \parens[\big]{= \delta_{\modulo H}(\ell \after \varkappa)}.
    \end{align*}
    The quadruple \graffito{naïve product $\mathcal{C}$ of $\mathcal{C}_{\modulo H}$ by $H$ and $\kappa$}$\mathcal{C} = \ntuple{\mathcal{R}, Q, N, \delta}$ is a semi-cellular automaton and is called \define{naïve product of $\mathcal{C}_{\modulo H}$ by $H$ and $\kappa$}\index{product of $\mathcal{C}_{\modulo H}$ by $H$ and $\kappa$!naïve}. 
  \end{definition}

  \begin{remark}
    The subscript $\modulo H$ of $\mathcal{C}_{\modulo H}$, $N_{\modulo H}$, and $\delta_{\modulo H}$ shall suggest that $\mathcal{C}_{\modulo H}$ is a semi-cellular automaton over the quotient $\mathcal{R} \modulo (H, \rho)$. It shall \emph{not} mean that $\mathcal{C}_{\modulo H}$ is the quotient of a semi-cellular automaton by $H$.
  \end{remark}

  \begin{example}[Cylinder]
  \label{example:42-lattice:product:shift-automaton}
    In the situation of \cref{example:42-lattice:quotient:shift-automaton}, let $\kappa$ be a right inverse of the canonical projection $\varpi \from \Z_{42} \times \Z \to \Z_{42}$, let $\mathcal{C}_\kappa$ be the naïve product of $\mathcal{C}_{\modulo H}$ by $H$ and $\kappa$, and let $\Delta_\kappa$ be the global transition function of $\mathcal{C}$.
    \begin{aenumerate}
      \item If $\kappa$ is given by $r_1 + 42\Z \mapsto (r_1 + 42\Z, z)$, then $\mathcal{C}_\kappa = \mathcal{C}$ and $\Delta_\kappa = \Delta$.
      \item If $\kappa$ is given by $r_1 + 42\Z \mapsto (r_1 + 42\Z, 0)$, then $\Delta_\kappa$ is the shift from left to right over $M$.
      \item If $\kappa$ is given by $r_1 + 42\Z \mapsto (r_1 + 42\Z, r_1 \bmod 42)$, then $\Delta_\kappa$ is the diagonal shift from the bottom-left to the top-right over $M$.
      \item If $\kappa$ is given by $r_1 + 42\Z \mapsto (r_1 + 42\Z, -(r_1 \bmod 42))$, then $\Delta_\kappa$ is the diagonal shift from the top-left to the bottom-right over $M$.
    \end{aenumerate}
    So, it depends on the choice of $\kappa$, whether the naïve product of the quotient of $\mathcal{C}$ is again $\mathcal{C}$ or not.
  \end{example}

  If the right inverse $\kappa$ is in a certain sense equivariant, then the product of a cellular automaton is a cellular automaton, as stated in

  \begin{lemma}
  \label{lemma:in-which-case-is-naive-product-cellular-automaton}
    In the situation of \cref{definition:naive-product-of-cellular-automaton}, let $\mathcal{C}_{\modulo H}$ be a cellular automaton and let $\kappa$ be such that
    \begin{equation}
    \label{equation:in-which-case-is-naive-product-cellular-automaton}
      \ForEach g_0 \in G_0 \ForEach \mathfrak{g}_{\modulo H} \in (G \modulo H) \modulo (G_0 \modulo H) \Holds
          g_0 \cdot \kappa(\mathfrak{g}_{\modulo H}) = \kappa(g_0 H \cdot_{\modulo H} \mathfrak{g}_{\modulo H}).
    \end{equation}
    The naïve product $\mathcal{C}$ of $\mathcal{C}_{\modulo H}$ by $H$ and $\kappa$ is a cellular automaton.
  \end{lemma}

  \begin{proof}
    Let $g_0 \in G_0$ and let $\ell \in Q^N$. Then, for each $n_{\modulo H} \in N_{\modulo H}$, because \cref{equation:in-which-case-is-naive-product-cellular-automaton} holds,
    \begin{align*}
      ((g_0 \bullet \ell) \after \varkappa)(n_{\modulo H})
      &= (g_0 \bullet \ell)(\kappa(n_{\modulo H}))\\
      &= \ell(g_0^{-1} \cdot \kappa(n_{\modulo H}))\\
      &= \ell(\kappa(g_0^{-1} H \cdot_{\modulo H} n_{\modulo H}))\\
      &= (\ell \after \varkappa)(g_0^{-1} H \cdot_{\modulo H} n_{\modulo H})\\
      &= (g_0 H \bullet_{\modulo H} (\ell \after \varkappa))(n_{\modulo H}).
    \end{align*}
    Thus, $(g_0 \bullet \ell) \after \varkappa = g_0 H \bullet_{\modulo H} (\ell \after \varkappa)$. Hence, because $\delta_{\modulo H}$ is $\bullet_{\modulo H}$-invariant,
    \begin{align*}
      \delta(g_0 \bullet \ell)
      &= \delta_{\modulo H}((g_0 \bullet \ell) \after \varkappa)\\
      &= \delta_{\modulo H}(g_0 H \bullet_{\modulo H} (\ell \after \varkappa))\\
      &= \delta_{\modulo H}(\ell \after \varkappa)\\
      &= \delta(\ell).
    \end{align*}
    Therefore, $\delta$ is $\bullet$-invariant. In conclusion, $\mathcal{C}$ is a cellular automaton.
  \end{proof}

  \begin{remark}
  \label{remark:in-which-case-is-naive-product-cellular-automaton}
    In the situation of \cref{definition:naive-product-of-cellular-automaton}, if, for each $g_0 \in G_0$, we have $g_0 \cdot \kappa(\blank) = \kappa(\blank)$, then $\mathcal{C}$ is a cellular automaton, regardless of whether $\mathcal{C}_{\modulo H}$ is one or not.
  \end{remark}

  \begin{example}[Cylinder] 
    In the situation of \cref{example:42-lattice:product:shift-automaton}, because the group $G$ is abelian, for each element $(t_0 + 42\Z, 0) \in G_0$, we have $(t_0 + 42\Z, 0) \cdot \kappa(\blank) = \kappa(\blank)$, and hence, according to \cref{remark:in-which-case-is-naive-product-cellular-automaton}, the naïve product $\mathcal{C}_\kappa$ is a cellular automaton. However, \cref{equation:in-which-case-is-naive-product-cellular-automaton} of \cref{lemma:in-which-case-is-naive-product-cellular-automaton}, which in this example is equivalent to
    \begin{multline*}
      \ForEach (t_0 + 42\Z, 0) \in G_0 \ForEach r_1 + 42\Z \in \Z_{42} \simeq (G \modulo H) \modulo (G_0 \modulo H) \Holds\\
          (t_0 + 42\Z, 0) \cdot \kappa(r_1 + 42\Z) = \kappa((t_0 + r_1) + 42\Z),
    \end{multline*}
    holds if and only if the projection of $\kappa$ to the second component is constant, which need not be the case. Therefore, although $\mathcal{C}_\kappa$ is a cellular automaton, the map $\kappa$ need not satisfy \cref{equation:in-which-case-is-naive-product-cellular-automaton}.
  \end{example}

  If \cref{equation:in-which-case-is-naive-product-cellular-automaton} of \cref{lemma:in-which-case-is-naive-product-cellular-automaton} does not hold, then the naïve product may not be a cellular automaton. For example, if we factor by a normal subgroup $H$ of $G$ that includes the stabiliser $G_0$, then all local transition functions $\delta_{\modulo H}$ are $\bullet_{\modulo H}$-invariant but their products are in general not $\bullet$-invariant, which is illustrated in

  \begin{example}[Tree] 
  \label{example:tree:product-is-not-a-cellular-automaton}
    Let $G$ be the free group over $\setOf{a, b}$, where $a \neq b$, let $G_0$ be the subgroup of $G$ that is generated by $a$, let $M$ be the quotient set $G \modulo G_0$, let $\actsOnPoint$ be the transitive left group action of $G$ on $M$ by left multiplication, let $\mathcal{K} = \ntuple{m_0, \family{g_{m_0, m}}_{m \in M}}$ be a coordinate system for $\mathcal{M} = \ntuple{M, G, \actsOnPoint}$, and let $\mathcal{R}$ be the cell space $\ntuple{\mathcal{M}, \mathcal{K}}$.

    Moreover, let $H$ be the kernel of the group homomorphism $\varphi \from G \to \Z \modulo 2\Z$ given by $a \mapsto 2\Z$, $b \mapsto 1 + 2\Z$, let $\rho$ be a right inverse of the canonical projection $\pi \from M \to H \reverseModulo M$ such that $\rho(H \actsOnPoint m_0) = m_0$, and let $\kappa$ be a right inverse of the canonical projection $\varpi \from G \modulo G_0 \to (G \modulo H) \modulo (G_0 \modulo H)$.

    Furthermore, let $Q$ be the binary set $\setOf{0, 1}$, let $N_{\modulo H}$ be the singleton set $\setOf{b H (G_0 \modulo H)}$, let $\delta_{\modulo H}$ be the map $Q^{N_{\modulo H}} \to Q$, $\ell_{\modulo H} \mapsto \ell_{\modulo H}(b H (G_0 \modulo H))$, let $\mathcal{C}_{\modulo H}$ be the semi-cellular automaton $\ntuple{\mathcal{R} \modulo (H, \rho), Q, N_{\modulo H}, \delta_{\modulo H}}$, and let $\mathcal{C} = \ntuple{\mathcal{R}, Q, N, \delta}$ be the naïve product of $\mathcal{C}_{\modulo H}$ by $H$ and $\kappa$.

    The group $H$ is a normal subgroup of $G$ of index $2$; the orbit space $H \reverseModulo M$ is the double quotient set $H \reverseModulo (G \modulo G_0)$, which we identify with $G \modulo H$ by $H (g G_0) \mapsto g H$; the quotient group $G \modulo H$ is equal to $\setOf{H, b H}$; the left group action $\actsOnPoint_{\modulo H}$ and the right quotient set semi-action $\isSemiActedUponBy_{\modulo (H, \rho)}$ it induces are identical to the group multiplication of $G \modulo H$; and the global transition function of $\mathcal{C}_{\modulo H}$ is the shift operator on $Q^{G \modulo H}$. 

    The stabiliser $G_0$ of $m_0$ under $\actsOnPoint$ is a subgroup of $H$, the stabiliser $G_0 \modulo H$ of $m_0 H$ under $\actsOnPoint_{\modulo H}$ is trivial, we identify the quotient set $(G \modulo H) \modulo (G_0 \modulo H)$ with $G \modulo H$ by $g H (G_0 \modulo H) \mapsto g H$, the local transition function $\delta_{\modulo H}$ is trivially $\bullet_{\modulo H}$-invariant, and the automaton $\mathcal{C}_{\modulo H}$ is a cellular automaton.

    Because $\varpi(\kappa(b H)) = b H$, there is an element $g \in G$ that is represented by a reduced word in which the symbol $b$ occurs such that $\kappa(b H) = g G_0$. Thus, for each element $g_0 \in G_0 \smallsetminus \setOf{e_G}$, we have $g_0 \cdot \kappa(b H) \neq \kappa(b H)$, in particular, \cref{equation:in-which-case-is-naive-product-cellular-automaton} of \cref{lemma:in-which-case-is-naive-product-cellular-automaton} does not hold. Hence, for each element $g_0 \in G_0 \smallsetminus \setOf{e_G}$, there is a local configuration $\ell \in Q^N$ such that $\ell(g_0^{-1} \cdot \kappa(b H)) \neq \ell(\kappa(b H))$ and therefore
    \begin{multline*}
      \delta(g_0 \bullet \ell)
      =    \delta_{\modulo H}(n_{\modulo H} \mapsto \ell(g_0^{-1} \cdot \kappa(n_{\modulo H})))
      =    \ell(g_0^{-1} \cdot \kappa(b H))\\
      \neq \ell(\kappa(b H))
      =    \delta_{\modulo H}(n_{\modulo H} \mapsto \ell(\kappa(n_{\modulo H})))
      =    \delta(\ell).
    \end{multline*}
    More specifically, if $\kappa$ is the map $G \modulo H \to G \modulo G_0$, $H \mapsto G_0$, $b H \mapsto b G_0$, and $\ell$ is the local configuration $N \to Q$, $b G_0 \mapsto 0$, $a^k b G_0 \mapsto 1$, for $k \in \Z \smallsetminus \setOf{0}$, then $\delta(a \bullet \ell) = 1 \neq 0 = \delta(\ell)$. In conclusion, the local transition function $\delta$ is not $\bullet$-invariant and hence the naïve product $\mathcal{C}$ is not a cellular automaton.
  \end{example}

  The image of a global configuration and a cell under the global transition function of the naïve product of a semi-cellular automaton can be expressed in terms of the global transition function of the original automaton, as stated in

  \begin{lemma} 
  \label{lemma:global-transition-function-of-product-in-terms-of-original-one}
    Let $\mathcal{R} \modulo (H, \rho)$ be a quotient of $\ntuple{\ntuple{M, G, \actsOnPoint}, \ntuple{m_0, \family{g_{m_0, m}}_{m \in M}}}$ by $H$ and $\rho$, let $\mathcal{C} = \ntuple{\mathcal{R}, Q, N, \delta}$ be a naïve product of $\mathcal{C}_{\modulo H}$ by $H$ and $\kappa$, let $c$ be a global configuration of $Q^M$, let $m$ be an element of $M$, let $f_m$ be the map $((H \actsOnPoint m) \isSemiActedUponBy_{\modulo (H, \rho)} \blank)^{-1} \from H \reverseModulo M \to (G \modulo H) \modulo (G_0 \modulo H)$, and let
    \begin{align*}
      c_{\modulo H, m} \from H \reverseModulo M &\to Q,\\
      H \actsOnPoint m' &\mapsto c(m \isSemiActedUponBy \kappa(f_m(H \actsOnPoint m'))).
    \end{align*}
    Then, $\Delta(c)(m) = \Delta_{\modulo H}(c_{\modulo H, m})(H \actsOnPoint m)$.
  \end{lemma}

  \begin{proof}
    Because $\isSemiActedUponBy_{\modulo (H, \rho)}$ is free and transitive, $(H \actsOnPoint m) \isSemiActedUponBy_{\modulo (H, \rho)} \blank$ is injective and surjective. Hence, $f_m$ is well-defined and therefore $c_{\modulo H, m}$ is also well-defined. Moreover,
    \begin{align*}
      \Delta(c)(m)
      &= \delta(n \mapsto c(m \isSemiActedUponBy n))\\
      &= \delta_{\modulo H}(n_{\modulo H} \mapsto c(m \isSemiActedUponBy \kappa(n_{\modulo H})))\\
      &= \delta_{\modulo H}(n_{\modulo H} \mapsto c_{\modulo H, m}((H \actsOnPoint m) \isSemiActedUponBy_{\modulo (H, \rho)} n_{\modulo H}))\\
      &= \Delta_{\modulo H}(c_{\modulo H, m})(H \actsOnPoint m). \qedhere
    \end{align*}
  \end{proof}


  \begin{example}[Cylinder]
  \label{example:42-lattice:product:global-transition-function-of-product-in-terms-of-original-one}
    In the situation of \cref{example:42-lattice:product:shift-automaton}, let $\kappa$ be given by $r_1 + 42\Z \mapsto (r_1 + 42\Z, 0)$, let $c$ be a global configuration of $Q^M$, let $m = (m_1 + 42\Z, m_2)$ be an element of $M$, and let $f_m$ and $c_{\modulo H, m}$ be as in \cref{lemma:global-transition-function-of-product-in-terms-of-original-one}. Recall that, according to \cref{example:42-lattice:quotient:right-semi-action}, the semi-action $\isSemiActedUponBy_{\modulo (H, \rho)}$ is addition in $\Z_{42}$ and the semi-action $\isSemiActedUponBy$ is addition in $\Z_{42} \times \Z$. Hence, the map $f_m$ is the map $\Z_{42} \to \Z_{42}$, $m_1' + 42\Z \mapsto (m_1' - m_1) + 42\Z$, and the global configuration $c_{\modulo H, m}$ is the map $\Z_{42} \to Q$, $m_1' + 42\Z \mapsto c(m_1' + 42\Z, m_2)$, which is essentially the restriction of $c$ to $\Z_{42} \times \setOf{m_2}$.
  \end{example}

  The quotient of the naïve product of a semi-cellular automaton is the original automaton, as shown in

  \begin{lemma}
  \label{lemma:product-then-quotient-is-identity}
    Let $\mathcal{R} \modulo (H, \rho)$ be a quotient of $\ntuple{\ntuple{M, G, \actsOnPoint}, \ntuple{m_0, \family{g_{m_0, m}}_{m \in M}}}$ by $H$ and $\rho$, let $\mathcal{C}_{\modulo H} = \ntuple{\mathcal{R} \modulo (H, \rho), Q, N_{\modulo H}, \delta_{\modulo H}}$ be a semi-cellular automaton, let $\mathcal{C}$ be a naïve product of $\mathcal{C}_{\modulo H}$ by $H$ and $\kappa$, and let $\mathcal{C} \modulo (H, \rho) = \ntuple{\mathcal{R} \modulo (H, \rho), Q, N_{\modulo H}', \delta_{\modulo H}'}$ be a quotient of $\mathcal{C} = \ntuple{\mathcal{R}, Q, N, \delta}$ by $H$ and $\rho$. Then, $\mathcal{C} \modulo (H, \rho) = \mathcal{C}_{\modulo H}$.
  \end{lemma}

  \begin{proof}
    For each $g_0 \in G_0$ and each $g G_0 \in G \modulo G_0$, $\varpi(g_0 \cdot g G_0) = g_0 g H (G_0 \modulo H) = g_0 H \cdot_{\modulo H} g H (G_0 \modulo H) = g_0 H \cdot_{\modulo H} \varpi(g G_0)$. Hence, 
    \begin{align*}
      N_{\modulo H}'
      &= \varpi(N)\\
      &= \varpi(G_0 \cdot \kappa(N_{\modulo H}))\\
      &= \setOf{g_0 H \suchThat g_0 \in G_0} \cdot_{\modulo H} \varpi(\kappa(N_{\modulo H}))\\
      &= (G_0 \modulo H) \cdot_{\modulo H} N_{\modulo H}\\
      &= N_{\modulo H}.
    \end{align*}
    Therefore, for each $\ell_{\modulo H} \in Q^{N_{\modulo H}}$,
    \begin{align*}
      \delta_{\modulo H}'(\ell_{\modulo H})
      &= \delta(n \mapsto \ell_{\modulo H}(\varpi(n)))\\
      &= \delta_{\modulo H}(n_{\modulo H} \mapsto \ell_{\modulo H}(\varpi(\kappa(n_{\modulo H}))))\\
      &= \delta_{\modulo H}(n_{\modulo H} \mapsto \ell_{\modulo H}(n_{\modulo H}))\\
      &= \delta_{\modulo H}(\ell_{\modulo H}).
    \end{align*}
    In conclusion, $\mathcal{C} \modulo (H, \rho) = \mathcal{C}_{\modulo H}$.
  \end{proof}

  \begin{remark}
    According to \cref{example:42-lattice:product:shift-automaton}, the naïve product of the quotient of a semi-cellular automaton $\mathcal{C}$ may not be $\mathcal{C}$.
  \end{remark}

  \section{Products}
  \label{section:products}

  The naïve product by $H$ and $\kappa$ of a cellular automaton is in general not a cellular automaton (see \cref{example:tree:product-is-not-a-cellular-automaton}). However, as we have seen in \cref{lemma:in-which-case-is-naive-product-cellular-automaton}, if the map $\kappa$ is in a certain sense equivariant (see \cref{definition:funny-equivariance-for-product}), then the naïve product by $H$ and $\kappa$ of a cellular automaton is again a cellular automaton and is simply called \emph{product by $H$ and $\kappa$} (see \cref{lemma:non-naive-product-of-cellular-automaton}). If the group $H$ is $G_0$-producible (see \cref{definition:producible}), then such maps $\kappa$ can be constructed as in \cref{lemma:equivariant-right-inverse-of-iota}. The global transition function of the product of a cellular automaton does not depend on the coordinate system, in particular, the right inverse $\kappa$ can in a sense be chosen independently of the origin (see \cref{lemma:non-naive-product-of-cellular-automaton}). 

  \begin{definition}
  \label{definition:producible}
    Let $G$ be a group, and let $G_0$ and $H$ be two subgroups of $G$. The group $H$ is called \defineX{$G_0$-producible}{producible@$G_0$-producible}\graffito{$G_0$-producible} if and only if
    \begin{equation*}
      \ForEach g \in G \Holds (g G_0 g^{-1}) \cap (G_0 H) \subseteq G_0. \qedhere
    \end{equation*}
  \end{definition}


  \begin{remark}
  \label{remark:transitive-action-of-abelian-group-is-producible}
    If the group $G$ is abelian or the group $H$ is included in $G_0$, then the group $H$ is $G_0$-producible. 
  \end{remark}

  \begin{lemma}
    Let $\actsOnPoint$ be a transitive left group action of $G$ on $M$, let $m_0$ be an element of $M$, and let $H$ be a subgroup of $G$. The group $H$ is $G_{m_0}$-producible if and only if
    \begin{equation*}
      \ForEach m \in M \Holds G_m \cap (G_{m_0} H) \subseteq G_{m_0}. \qedhere 
    \end{equation*}
  \end{lemma}

  \begin{proof}
    This is a direct consequence of \cref{lemma:stabiliser-versus-transporter}.
  \end{proof}

  \begin{lemma}
    Let $G$ be a group, let $G_0$ be a subgroup of $G$, let $H$ be a normal subgroup of $G$, and let $g$ be an element of $G$. The group $H$ is $G_0$-producible if and only if it is $(g G_0 g^{-1})$-producible.
  \end{lemma}

  \begin{proof}
    For reasons of symmetry, it suffices to prove one implication. To this end, let $H$ be $G_0$-producible. Furthermore, let $g' \in G$. Then, because $H$ is normal in $G$, we have $g G_0 g^{-1} H = g G_0 H g^{-1}$. And, because $H$ is $G_0$-producible, we have $(g^{-1} g' g G_0 g^{-1} (g')^{-1} g) \cap (G_0 H) \subseteq G_0$. Hence,
    \begin{align*}
      &(g' g G_0 g^{-1} (g')^{-1}) \cap (g G_0 g^{-1} H)\\
      &= g \parens[\big]{(g^{-1} g' g G_0 g^{-1} (g')^{-1} g) \cap (G_0 H)} g^{-1}\\
      &\subseteq g G_0 g^{-1}.
    \end{align*}
    In conclusion, $H$ is $g G_0 g^{-1}$-producible.
  \end{proof}

  \begin{definition}
    Let $G$ be a group and let $H$ be a subgroup of $G$. The subgroup $H$ is called \define{malnormal}\graffito{malnormal} if and only if
    \begin{equation*}
      \ForEach g \in G \smallsetminus H \Holds (g H g^{-1}) \cap H = \setOf{e_G}. \qedhere
    \end{equation*}
  \end{definition}

  \begin{remark}[{\cite[Item~i of Proposition~2]{delaHarpe:weber:2014}}]
    The subgroups $\setOf{e_G}$ and $G$ of $G$ are the only ones that are both normal and malnormal.
  \end{remark}

  We can use semi-direct products to construct groups with producible subgroups, as shown in

  \begin{lemma} 
  \label{lemma:semidirect-product-is-producible}
    Let $F$ and $H$ be two groups, let $G_0$ be a malnormal subgroup of $F$, and let $\varphi$ be a group homomorphism from $F$ to $\automorphismsOf(H)$ such that $G_0 \subseteq \kernelOf(\varphi)$, and let $G = H \rtimes_\varphi F$ be the outer semi-direct product of $F$ acting on $H$ by $\varphi$. The normal subgroup $H \times \setOf{e_F}$ of $G$ is $(\setOf{e_H} \times G_0)$-producible and, if there is a tuple $(g_0, f) \in G_0 \times (F \smallsetminus G_0)$ such that $f g_0 f^{-1} \notin G_0$, then $\setOf{e_H} \times G_0$ is not normal in $G$.
  \end{lemma}

  \begin{proof}
    Let $(h, f) \in G$. Then,
    \begin{align*}
      &(h, f) \cdot (\setOf{e_H} \times G_0) \cdot (h, f)^{-1}\\
      &= (h, f) \cdot (\setOf{e_H} \times G_0) \cdot (\varphi(f^{-1})(h^{-1}), f^{-1})\\
      &= \setOf{(h \varphi(f g_0 f^{-1})(h^{-1}), f g_0 f^{-1}) \suchThat g_0 \in G_0}.
    \end{align*}
    And, because the codomain of $\varphi$ is $\automorphismsOf(H)$,
    \begin{align*}
      (\setOf{e_H} \times G_0) \cdot (H \times \setOf{e_F})
      &= \setOf{(\varphi(g_0)(h), g_0) \suchThat h \in H, g_0 \in G_0}\\
      &= H \times G_0.
    \end{align*}

    First, let it be the case that $f \in G_0$. Then, $f g_0 f^{-1} \in G_0$. Hence, because $G_0 \subseteq \kernelOf(\varphi)$, we have $\varphi(f g_0 f^{-1}) = \identityMap$. Thus, $h \varphi(f g_0 f^{-1})(h^{-1}) = e_H$. Therefore, $(h, f) \cdot (\setOf{e_H} \times G_0) \cdot (h, f)^{-1} \subseteq \setOf{e_H} \times G_0$. 

    Secondly, let it be the case that $f \notin G_0$. Then, because $G_0$ is malnormal in $F$, we have $f g_0 f^{-1} = e_F$ or $f g_0 f^{-1} \notin G_0$. Hence, we have $(h \varphi(f g_0 f^{-1})(h^{-1}), f g_0 f^{-1}) = (e_H, e_F)$ or $f g_0 f^{-1} \notin G_0$. Therefore, $((h, f) \cdot (\setOf{e_H} \times G_0) \cdot (h, f)^{-1}) \cap (\setOf{e_H} \times G_0) \cdot (H \times \setOf{e_F}) \subseteq \setOf{(e_H, e_F)}$. 

    In either case, $((h, f) \cdot (\setOf{e_H} \times G_0) \cdot (h, f)^{-1}) \cap (\setOf{e_H} \times G_0) \cdot (H \times \setOf{e_F}) \subseteq \setOf{e_H} \times G_0$. In conclusion, the group $H \times \setOf{e_F}$ is $(\setOf{e_H} \times G_0)$-producible.

    Moreover, if there is a tuple $(g_0, f) \in G_0 \times (F \smallsetminus G_0)$ such that $f g_0 f^{-1} \notin G_0$, then, according to the second case above, $(e_H, f g_0 f^{-1}) \in (e_H, f) \cdot (\setOf{e_H} \times G_0) \cdot (e_H, f)^{-1}$, and therefore $(e_H, f) \cdot (\setOf{e_H} \times G_0) \cdot (e_H, f)^{-1} \nsubseteq \setOf{e_H} \times G_0$. In conclusion, $\setOf{e_H} \times G_0$ is not normal in $G$.
  \end{proof}

  We can use direct products to construct groups with producible subgroups, as shown in

  \begin{lemma}
  \label{lemma:direct-product-is-producible}
    Let $F$ and $H$ be two groups, let $G_0$ be a subgroup of $F$, and let $G$ be the direct product $H \times F$. The normal subgroup $H \times \setOf{e_F}$ of $G$ is $(\setOf{e_H} \times G_0)$-producible and, if $G_0$ is not normal in $F$, then $\setOf{e_H} \times G_0$ is not normal in $G$.
  \end{lemma}

  \begin{proof}
    The proof is omitted here. It is similar to the one of \cref{lemma:semidirect-product-is-producible}.
  \end{proof}

  \begin{example}[Thick Tree] 
  \label{example:thick-tree:producible}
    Let $F$ be the free group over $\setOf{a, b}$, where $a \neq b$, let $H$ be the additive cyclic group $\Z \modulo 7\Z$ of order $7$, let $G_0$ be the non-normal subgroup of $F$ that is generated by $a$, let $G$ be the direct product $H \times F$, let $\actsOnPoint$ be the transitive left group action of $G$ on $G \modulo (\setOf{e_H} \times G_0)$ by left multiplication, and identify $G \modulo (\setOf{e_H} \times G_0)$ with $H \times (F \modulo G_0)$, $H \times \setOf{e_F}$ with $H$, and $(\setOf{e_H} \times G_0)$ with $G_0$. The normal subgroup $H$ of $G$ is $G_0$-producible and the subgroup $G_0$ of $G$ is not normal.
  \end{example}

  As we have seen in \cref{lemma:in-which-case-is-naive-product-cellular-automaton}, if the product of a cellular automaton is made by a right inverse $\kappa$ that is in a certain sense equivariant, which we give a name below, then that product is itself a cellular automaton.

  \begin{definition} 
  \label{definition:funny-equivariance-for-product}
    Let $G$ be a group, let $G_0$ be a subgroup of $G$, let $H$ be a normal subgroup of $G$, let $\varpi$ be the canonical projection from $G \modulo G_0$ onto $(G \modulo H) \modulo (G_0 \modulo H)$, and let $\kappa$ be a right inverse of $\varpi$. The map $\kappa$ is called \defineX{$G_0$-e\-qui\-var\-i\-ant}{equivariantG0@$G_0$-e\-qui\-var\-i\-ant}\graffito{$G_0$-e\-qui\-var\-i\-ant} if and only if 
    \begin{multline*}
      \ForEach g_0 \in G_0 \ForEach \mathfrak{g}_{\modulo H} \in (G \modulo H) \modulo (G_0 \modulo H) \Holds\\
          g_0 \cdot \kappa(\mathfrak{g}_{\modulo H}) = \kappa(g_0 H \cdot_{\modulo H} \mathfrak{g}_{\modulo H}). \qedhere
    \end{multline*}
  \end{definition}


  \begin{example}[Thick Tree]
  \label{example:thick-tree:funny-equivariance-for-product}
    In the situation of \cref{example:thick-tree:producible}, recall that $G \modulo G_0$ is identified with $H \times (F \modulo G_0)$, and identify $G \modulo H$ with $F$ and $G_0 \modulo H$ with $G_0$. Then, the canonical projection $\varpi$ from $G \modulo G_0$ onto $(G \modulo H) \modulo (G_0 \modulo H)$ is the projection to the second component, namely $(h, f G_0) \mapsto f G_0$. Hence, for each right inverse $\kappa$ of $\varpi$, it is $G_0$-e\-qui\-var\-i\-ant if and only if
    \begin{equation*}
      \ForEach z \in \Z \ForEach f G_0 \in F \modulo G_0 \Holds a^z \cdot \kappa(f G_0) = \kappa(a^z \cdot f G_0),
    \end{equation*}
    which is the case if and only if, for each element $f G_0 \in F \modulo G_0$, the projection to the first component of the restriction $\kappa\restrictedTo_{G_0 \cdot f G_0}$ is constant. In particular, for each element $h \in H$, the map $\kappa \from F \modulo G_0 \to H \times (F \modulo G_0)$, $f G_0 \mapsto (h, f G_0)$, is a $G_0$-e\-qui\-var\-i\-ant right inverse of $\varpi$.
  \end{example}

  If $H$ is a $G_0$-producible and normal subgroup of $G$, then we can construct a $G_0$-e\-qui\-var\-i\-ant right inverse of the canonical projection from $G \modulo G_0$ onto $(G \modulo H) \modulo (G_0 \modulo H)$, as shown in

  \begin{lemma} 
  \label{lemma:equivariant-right-inverse-of-iota}
    Let $G$ be a group, let $G_0$ be a subgroup of $G$, let $H$ be a $G_0$-producible and normal subgroup of $G$, and let $\varpi$ be the canonical projection from $G \modulo G_0$ onto $(G \modulo H) \modulo (G_0 \modulo H)$. Furthermore, let $X$ be the quotient set $(G \modulo H) \modulo (G_0 \modulo H)$, let $\sim$ be the equivalence relation on $X$ given by
    \begin{equation*}
      \ForEach x \in X \ForEach x' \in X \Holds
          \parens*{x \sim x' \ifAndOnlyIf \Exists g_0 \in G_0 \SuchThat g_0 H \cdot_{\modulo H} x = x'}, 
    \end{equation*}
    let $\setOf{x_i \suchThat i \in I}$ be a transversal of $X \modulo {\sim}$, and let $\xi$ be a right inverse of $\varpi$. The map
    \begin{align*} 
      \kappa \from X &\to G \modulo G_0,\\
      g_0 H \cdot_{\modulo H} x_i &\mapsto g_0 \cdot \xi(x_i),
    \end{align*}
    is a $G_0$-e\-qui\-var\-i\-ant right inverse of $\varpi$.
  \end{lemma}

  \begin{proof}
    \proofPart{Subproof of well-definedness}
      The equivalence classes $[x_i]_\sim$, for $i \in I$, partition $X$ and, for each $i \in I$ and each $x \in [x_i]_\sim$, there is a $g_0 \in G_0$ such that $x = g_0 H \cdot_{\modulo H} x_i$. Hence, for each $x \in X$, there is a unique $i \in I$ and a $g_0 \in G_0$ such that $x = g_0 H \cdot_{\modulo H} x_i$. Therefore, if, for each $i \in I$, each $g_0 \in G_0$, and each $g_0' \in G_0$ such that $g_0 H \cdot_{\modulo H} x_i = g_0' H \cdot_{\modulo H} x_i$, we have $g_0 \cdot \xi(x_i) = g_0' \cdot \xi(x_i)$, then $\kappa$ is well-defined.

      Let $i \in I$, let $g_0 \in G_0$, and let $g_0' \in G_0$ such that $g_0 H \cdot_{\modulo H} x_i = g_0' H \cdot_{\modulo H} x_i$. Then, there is a $g \in G$ such that $\xi(x_i) = g G_0$. Thus, $x_i = \varpi(\xi(x_i)) = g H (G_0 \modulo H)$. Hence, with $g_0'' = g_0^{-1} g_0'$,
      \begin{align*}
        g_0 H \cdot_{\modulo H} x_i = g_0' H \cdot_{\modulo H} x_i
        &\ifAndOnlyIf g_0'' H \cdot_{\modulo H} x_i = x_i\\
        &\ifAndOnlyIf g_0'' g H (G_0 \modulo H) = g H (G_0 \modulo H)\\
        &\ifAndOnlyIf g^{-1} g_0'' g H (G_0 \modulo H) = G_0 \modulo H\\
        &\ifAndOnlyIf \Exists g_0 \in G_0 \SuchThat g^{-1} g_0'' g H = g_0 H\\
        &\ifAndOnlyIf \Exists g_0 \in G_0 \SuchThat g^{-1} g_0'' g \in g_0 H\\
        &\ifAndOnlyIf g^{-1} g_0'' g \in G_0 H.
      \end{align*}
      And, because $g^{-1} g_0'' g \in G_{g^{-1} \actsOnPoint m}$ and $H$ is $G_0$-producible,
      \begin{equation*}
        g^{-1} g_0'' g \in G_0 H
        \implies g^{-1} g_0'' g \in G_0.
      \end{equation*}
      And,
      \begin{align*}
        g^{-1} g_0'' g \in G_0
        &\ifAndOnlyIf g_0'' g G_0 = g G_0\\
        &\ifAndOnlyIf g_0'' \cdot \xi(x_i) = \xi(x_i)\\
        &\ifAndOnlyIf g_0 \cdot \xi(x_i) = g_0' \cdot \xi(x_i).
      \end{align*}
      In conclusion, because $g_0 H \cdot_{\modulo H} x_i = g_0' H \cdot_{\modulo H} x_i$, we have $g_0 \cdot \xi(x_i) = g_0' \cdot \xi(x_i)$.

    \proofPart{Subproof of right inverseness}
      Let $x \in X$. Then, there is a $g_0 \in G_0$ and there is an $i \in I$ such that $x = g_0 H \cdot_{\modulo H} x_i$. And, there is a $g \in G$ such that $\xi(x_i) = g G_0$. Hence,
      \begin{align*}
        (\varpi \after \kappa)(x)
        &= \varpi(\kappa(x))\\
        &= \varpi(g_0 \cdot \xi(x_i))\\
        &= \varpi(g_0 g G_0)\\
        &= g_0 g H (G_0 \modulo H)\\
        &= g_0 H \cdot_{\modulo H} g H (G_0 \modulo H)\\
        &= g_0 H \cdot_{\modulo H} \varpi(g G_0)\\
        &= g_0 H \cdot_{\modulo H} \varpi(\xi(x_i))\\
        &= g_0 H \cdot_{\modulo H} x_i\\
        &= x.
      \end{align*}
      In conclusion, $\kappa$ is a right inverse of $\varpi$.

    \proofPart{Subproof of $G_0$-equivariance}
      Let $x \in X$ and let $g_0 \in G_0$. Then, there is a $g_0' \in G_0$ and there is an $i \in I$ such that $x = g_0' H \cdot_{\modulo H} x_i$. Hence,
      \begin{align*}
        g_0 \cdot \kappa(x) &= g_0 \cdot (g_0' \cdot \xi(x_i))\\
        &= g_0 g_0' \cdot \xi(x_i)\\
        &= \kappa(g_0 g_0' H \cdot_{\modulo H} x_i)\\
        &= \kappa(g_0 H \cdot_{\modulo H} (g_0' H \cdot_{\modulo H} x_i))\\
        &= \kappa(g_0 H \cdot_{\modulo H} x).
      \end{align*}
      In conclusion, $\kappa$ is $G_0$-e\-qui\-var\-i\-ant.
  \end{proof}

  \begin{corollary}
    Let $G$ be a group, let $G_0$ be a subgroup of $G$, let $H$ be a $G_0$-producible and normal subgroup of $G$, and let $\varpi$ be the canonical projection from $G \modulo G_0$ onto $(G \modulo H) \modulo (G_0 \modulo H)$. There is a $G_0$-e\-qui\-var\-i\-ant right inverse of $\varpi$.
  \end{corollary}

  \begin{proof}
    This is a direct consequence of \cref{lemma:equivariant-right-inverse-of-iota}, because there is a right inverse $\xi$ of $\varpi$.
  \end{proof}

  \begin{remark} 
  \label{remark:every-equivariant-right-inverse-has-our-form}
    Each $G_0$-e\-qui\-var\-i\-ant right inverse $\kappa$ of $\varpi$ is of the form as in \cref{lemma:equivariant-right-inverse-of-iota}. Indeed, let $\kappa$ be such a right inverse, let $\setOf{x_i \suchThat i \in I}$ be a transversal of $X \modulo {\sim}$, and let $\xi$ be identical to $\kappa$. Then, for each $g_0 \in G_0$ and each $i \in I$, we have $\kappa(g_0 H \cdot_{\modulo H} x_i) = g_0 \cdot \kappa(x_i) = g_0 \cdot \xi(x_i)$.
  \end{remark}

  \begin{example}[Thick Tree]
  \label{example:thick-tree:equivariant-inverse}
    In the situation of \cref{example:thick-tree:funny-equivariance-for-product}, broadly speaking, given a right inverse $\xi$ of $\varpi$ and representatives $f_i G_0$, for $i \in I$, of the orbits $G_0 \cdot f G_0$, for $f G_0 \in F \modulo G_0$, forcing the images of the elements of $G_0 \cdot f_i G_0$ under $\xi$ to have the same first component as the image of $f_i G_0$ under $\xi$, yields a $G_0$-e\-qui\-var\-i\-ant right inverse $\kappa$ of $\varpi$.
  \end{example}

  The naïve product by a $G_0$-e\-qui\-var\-i\-ant right inverse is a cellular automaton, which we give a name in

  \begin{lemma}
    Let $\mathcal{R} \modulo (H, \rho)$ be a quotient of $\ntuple{\ntuple{M, G, \actsOnPoint}, \ntuple{m_0, \family{g_{m_0, m}}_{m \in M}}}$ by $H$ and $\rho$, let $\kappa$ be a $G_0$-e\-qui\-var\-i\-ant right inverse of the canonical projection $\varpi \from G \modulo G_0 \to (G \modulo H) \modulo (G_0 \modulo H)$, and let $\mathcal{C}_{\modulo H} = \ntuple{\mathcal{R} \modulo (H, \rho), Q, N_{\modulo H}, \delta_{\modulo H}}$ be a cellular automaton. The naïve product \graffito{product $\mathcal{C}$ of $\mathcal{C}_{\modulo H}$ by $H$ and $\kappa$}$\mathcal{C}$ of $\mathcal{C}_{\modulo H}$ by $H$ and $\kappa$ is a cellular automaton and is called \define{product of $\mathcal{C}_{\modulo H}$ by $H$ and $\kappa$}.
  \end{lemma}

  \begin{proof}
    According to \cref{lemma:in-which-case-is-naive-product-cellular-automaton}, the semi-cellular automaton $\mathcal{C}$ is a cellular automaton.
  \end{proof}

  \begin{example}[Thick Tree]
  \label{example:thick-tree:product-cellular-automaton}
    In the situation of \cref{example:thick-tree:funny-equivariance-for-product}, let $\mathcal{K} = \ntuple{e_G G_0, \family{g_{e_G G_0, g G_0}}_{g G_0 \in G \modulo G_0}}$ be a coordinate system for $\mathcal{M} = \ntuple{G \modulo G_0, G, \actsOnPoint}$, let $\rho$ be the right inverse of the canonical projection $\pi \from H \times (F \modulo G_0) \simeq G \modulo G_0 \to H \reverseModulo (G \modulo G_0) \simeq F \modulo G_0$ given by $\rho(f G_0) = (0 + 7\Z, f G_0)$, and let $\kappa$ be the right inverse of the canonical projection $\varpi \from H \times (F \modulo G_0) \simeq G \modulo G_0 \to (G \modulo H) \modulo (G_0 \modulo H) \simeq F \modulo G_0$ given by $f G_0 \mapsto (0 + 7\Z, f G_0)$ (note that the right inverses $\rho$ and $\kappa$ are identical).

    Furthermore, let $Q$ be the binary set $\setOf{0, 1}$; let $N_{\modulo H}$ be the $(G_0 \simeq G_0 \modulo H)$-invariant subset $\setOf{a^z \cdot b G_0 \suchThat z \in \Z}$ of $F \modulo G_0 \simeq (G \modulo H) \modulo (G_0 \modulo H)$; let $\delta_{\modulo H}$ be the $\bullet_{\modulo H}$-invariant map
    \begin{align*}
      Q^{N_{\modulo H}} &\to Q,\\
      \ell_{\modulo H} &\mapsto \begin{dcases*}
        0, &if $\sum_{n_{\modulo H} \in N_{\modulo H}} \ell_{\modulo H}(n_{\modulo H}) < \infty$,\\
        1, &otherwise;
      \end{dcases*}
    \end{align*}
    let $\mathcal{C}_{\modulo H}$ be the cellular automaton $\ntuple{\ntuple{\mathcal{M}, \mathcal{K}} \modulo (H, \rho), Q, N_{\modulo H}, \delta_{\modulo H}}$; and let $\Delta_{\modulo H}$ be the global transition function of $\mathcal{C}_{\modulo H}$.

    The product of $\mathcal{C}_{\modulo H}$ by $H$ and $\kappa$ has the neighbourhood $N = \setOf{0 + 7\Z} \times N_{\modulo H}$; the local transition function $\delta \from Q^N \to Q$, $\ell \mapsto 0$, if $\sum_{n \in N} \ell(n) < \infty$, and $\ell \mapsto 1$, otherwise; and the product $\prod_{h \in H} \Delta_{\modulo H}$ for global transition function $\Delta$.
  \end{example}

  \begin{definition}
    Let $\ntuple{M, G, \actsOnPoint}$ be a left group set, let $H$ be a normal subgroup of $G$, and, for each point $m \in M$, let $\varpi_m$ be the canonical projection from $G \modulo G_m$ onto $(G \modulo H) \modulo (G_m \modulo H)$.
    \begin{aenumerate}
      \item For each point $m \in M$, the set of all $G_m$-e\-qui\-var\-i\-ant right inverses of $\varpi_m$ is denoted by $\mathfrak{k}_m$\index[symbols]{Km fraktur@$\mathfrak{k}_m$}\graffito{$\mathfrak{k}_m$}.
      \item The union of the sets $\mathfrak{k}_m$, for $m \in M$, is denoted by $\mathfrak{k}$\index[symbols]{K fraktur@$\mathfrak{k}$}\graffito{$\mathfrak{k}$}. \qedhere
    \end{aenumerate}
  \end{definition}

  The group acts on the set of all equivariant right inverses, as shown in

  \begin{lemma}
  \label{lemma:action-on-right-inverses-of-varpi}
    Let $\ntuple{M, G, \actsOnPoint}$ be a left group set and let $H$ be a normal subgroup of $G$. The group $G$ acts on $\mathfrak{k}$ on the left by
    \begin{align*} 
      \boxdot \from G \times \mathfrak{k} &\to \mathfrak{k},\\
      (g, \kappa_m) &\mapsto \left[
        \begin{aligned} 
          (G \modulo H) \modulo (G_{g \actsOnPoint m} \modulo H) &\to G \modulo G_{g \actsOnPoint m},\\
          x &\mapsto g \conjugates \kappa_m(g^{-1} H \conjugates_{\modulo H} x),
        \end{aligned}
      \right] 
    \end{align*}
    such that, for each symmetry $g \in G$ and each point $m \in M$, the map
    \begin{align*}
      (g \boxdot \blank)\restrictedTo_{\mathfrak{k}_m \to \mathfrak{k}_{g \actsOnPoint m}}
      \from \mathfrak{k}_m &\to \mathfrak{k}_{g \actsOnPoint m},\\
      \kappa_m &\mapsto g \boxdot \kappa_m,
    \end{align*}
    is bijective. 
  \end{lemma}

  \begin{proof}
    For each $m \in M$, let $X_m = (G \modulo H) \modulo (G_m \modulo H) = (G \modulo H) \modulo (G \modulo H)_{H \actsOnPoint m}$. Then, for each $m \in M$, we have $X_{g \actsOnPoint m} = (G \modulo H) \modulo (G \modulo H)_{g H \cdot_H (H \actsOnPoint m)}$. Moreover, for each $m \in M$, let $\varpi_m$ be the canonical projection from $G \modulo G_m$ onto $(G \modulo H) \modulo (G_m \modulo H)$.

    First, let $g \in G$, let $m \in M$, and let $\kappa_m \in \mathfrak{k}_m$. Then, for each $x \in X_{g \actsOnPoint m}$, according to \cref{lemma:left-group-action-on-union-of-quotients-by-stabilisers}, we have $g^{-1} H \conjugates_{\modulo H} x \in X_m$ and hence $g \conjugates \kappa_m(g^{-1} H \conjugates_{\modulo H} x) \in G \modulo G_{g \actsOnPoint m}$. Therefore, the map $g \boxdot \kappa_m$ is well-defined.

    Moreover, for each $x \in X_{g \actsOnPoint m}$, there is a $g' \in G$ such that $\kappa_m(g^{-1} H \conjugates_{\modulo H} x) = g' G_m$ and hence
    \begin{align*}
      \varpi_{g \actsOnPoint m}\parens[\big]{(g \boxdot \kappa_m)(x)}
      &= \varpi_{g \actsOnPoint m}\parens[\big]{g \conjugates \kappa_m(g^{-1} H \conjugates_{\modulo H} x)}\\
      &= \varpi_{g \actsOnPoint m}(g \conjugates g' G_m)\\
      &= \varpi_{g \actsOnPoint m}(g g' g^{-1} G_{g \actsOnPoint m})\\
      &= g g' g^{-1} H (G \modulo H)_{H \actsOnPoint (g \actsOnPoint m)}\\
      &= g H \conjugates_{\modulo H} g' H (G \modulo H)_{H \actsOnPoint m}\\
      &= g H \conjugates_{\modulo H} \varpi_m(g' G_m)\\
      &= g H \conjugates_{\modulo H} \varpi_m\parens[\big]{\kappa_m(g^{-1} H \conjugates_{\modulo H} x)}\\
      &= g H \conjugates_{\modulo H} (g^{-1} H \conjugates_{\modulo H} x)\\
      &= x.
    \end{align*}
    Therefore, the map $g \boxdot \kappa_m$ is a right inverse of $\varpi_{g \actsOnPoint m}$.

    Furthermore, for each $g_{g \actsOnPoint m} \in G_{g \actsOnPoint m}$ and each $x \in X_{g \actsOnPoint m}$, according to \cref{remark:actions-circ-and-cdot-commute} and because $g^{-1} g_{g \actsOnPoint m} g \in G_m$ and $\kappa_m$ is $G_m$-e\-qui\-var\-i\-ant,
    \begin{align*}
      (g \boxdot \kappa_m)(g_{g \actsOnPoint m} H \cdot_{\modulo H} x)
      &= g \conjugates \kappa_m\parens[\big]{g^{-1} H \conjugates_{\modulo H} (g_{g \actsOnPoint m} H \cdot_{\modulo H} x)}\\
      &= g \conjugates \kappa_m\parens[\big]{g^{-1} g_{g \actsOnPoint m} g H \cdot_{\modulo H} (g^{-1} H \conjugates_{\modulo H} x)}\\
      &= g \conjugates \parens[\big]{g^{-1} g_{g \actsOnPoint m} g \cdot \kappa_m(g^{-1} H \conjugates_{\modulo H} x)}\\
      &= g_{g \actsOnPoint m} \cdot \parens[\big]{g \conjugates \kappa_m(g^{-1} H \conjugates_{\modulo H} x)}\\
      &= g_{g \actsOnPoint m} \cdot (g \boxdot \kappa_m)(x).
    \end{align*}
    Therefore, the right inverse $g \boxdot \kappa_m$ is $G_{g \actsOnPoint m}$-e\-qui\-var\-i\-ant. Altogether, the right inverse $g \boxdot \kappa_m$ is an element of $\mathfrak{k}_{g \actsOnPoint m}$. In conclusion, the maps $\boxdot$ and $(g \boxdot \blank)\restrictedTo_{\mathfrak{k}_m \to \mathfrak{k}_{g \actsOnPoint m}}$ are well-defined.

    Secondly, for each $\kappa \in \mathfrak{k}$, we have $e_G \boxdot \kappa = \kappa$. And, for each $g \in G$, each $g' \in G$, each $\kappa \in \mathfrak{k}$, and each $x \in X_m$,
    \begin{align*}
      (g g' \boxdot \kappa)(x)
      &= g g' \conjugates \kappa\parens[\big]{(g')^{-1} g^{-1} H \conjugates_{\modulo H} x}\\
      &= g \conjugates \parens[\big]{g' \conjugates \kappa((g')^{-1} H \conjugates_{\modulo H} (g^{-1} H \conjugates_{\modulo H} x))}\\
      &= g \conjugates (g' \boxdot \kappa)(g^{-1} H \conjugates_{\modulo H} x)\\
      &= (g \boxdot (g' \boxdot \kappa))(x),
    \end{align*}
    and hence $g g' \boxdot \kappa = g \boxdot (g' \boxdot \kappa)$. In conclusion, the map $\boxdot$ is a left group action.
  \end{proof}

  \begin{example}[Thick Tree]
  \label{example:thick-tree:action}
    In the situation of \cref{example:thick-tree:product-cellular-automaton}, let $(h, f)$ be an element of $G$. Then, the subgroup $G_0 \simeq \setOf{e_H} \times G_0$ of $G$ is the stabiliser of $e_G G_0 \simeq \setOf{e_H} \times G_0$ under $\actsOnPoint$ and its conjugate $f G_0 f^{-1} \simeq \setOf{e_H} \times (f G_0 f^{-1})$ is the stabiliser of $(h, f) G_0 \simeq \setOf{h} \times (f G_0)$ under $\actsOnPoint$. Hence, the right inverse $(h, f) \boxdot \kappa$ of the canonical projection $\varpi_{(h, f) G_0} \from H \times (F \modulo (f G_0 f^{-1})) \simeq G \modulo (f G_0 f^{-1}) \to (G \modulo H) \modulo (f G_0 f^{-1} \modulo H) \simeq F \modulo (f G_0 f^{-1})$ is given by $f' f G_0 f^{-1} \mapsto (0 + 7\Z, f' f G_0 f^{-1})$. 
  \end{example}

  The orbit of a right inverse contains one for each stabiliser, which is shown in

  \begin{lemma}
    Let $\ntuple{M, G, \actsOnPoint}$ be a left-ho\-mo\-ge\-neous space and let $H$ be a normal subgroup of $G$. Then,
    \begin{equation*}
      \ForEach \kappa \in \mathfrak{k} \ForEach m \in M \Exists g \in G \SuchThat g \boxdot \kappa \in \mathfrak{k}_m. \qedhere
    \end{equation*} 
  \end{lemma}

  \begin{proof}
    Let $\kappa \in \mathfrak{k}$ and let $m \in M$. Then, there is an $m' \in M$ such that $\kappa \in \mathfrak{k}_{m'}$. And, because $\actsOnPoint$ is transitive, there is a $g \in G$ such that $g \actsOnPoint m' = m$. Hence, according to \cref{lemma:action-on-right-inverses-of-varpi}, we have $g \boxdot \kappa \in \mathfrak{k}_{g \actsOnPoint m'} = \mathfrak{k}_m$.
  \end{proof}

  The global transition function of the product of a cellular automaton does not depend on the coordinate system, which is shown in

  \begin{lemma} 
  \label{lemma:global-transition-function-of-product-does-not-depend-on-a-lot-of-things}
    Let $\mathcal{M} = \ntuple{M, G, \actsOnPoint}$ be a left-ho\-mo\-ge\-neous space, let $\mathcal{K} = \ntuple{m_0, \family{g_{m_0, m}}_{m \in M}}$ and $\mathcal{K}' = \ntuple{m_0', \family{g_{m_0', m}'}_{m \in M}}$ be two coordinate systems for $\mathcal{M}$, let $H$ be a normal subgroup of $G$, and let $\rho$ and $\rho'$ be two right inverses of the canonical projection $\pi \from M \to H \reverseModulo M$ such that $\rho(H \actsOnPoint m_0) = m_0$ and $\rho'(H \actsOnPoint m_0') = m_0'$

    Furthermore, let $\mathcal{C}_{\modulo H} = \ntuple{\ntuple{\mathcal{M}, \mathcal{K}} \modulo (H, \rho), Q, N_{\modulo H}, \delta_{\modulo H}}$ and $\mathcal{C}_{\modulo H}' = \ntuple{\ntuple{\mathcal{M}, \mathcal{K}'} \modulo (H, \rho'), Q, N_{\modulo H}', \delta_{\modulo H}'}$ be two cellular automata such that the global transition function $\Delta_{\modulo H}$ of $\mathcal{C}_{\modulo H}$ is identical to the one, namely $\Delta_{\modulo H}'$, of $\mathcal{C}_{\modulo H}'$.

    Moreover, let $\kappa$ be an element of $\mathfrak{k}_{m_0}$ and let $\kappa'$ be an element of $\mathfrak{k}_{m_0'}$ such that there is an element $g \in G_{m_0, m_0'}$ that satisfies $g \boxdot \kappa = \kappa'$, let $\mathcal{C}$ be the product of $\mathcal{C}_{\modulo H}$ by $H$ and $\kappa$, and let $\mathcal{C}'$ be the product of $\mathcal{C}_{\modulo H}'$ by $H$ and $\kappa'$.

    The global transition function $\Delta$ of $\mathcal{C}$ is identical to the one, namely $\Delta'$, of $\mathcal{C}'$.
  \end{lemma}

  \begin{proof}
    Because $\isSemiActedUponBy$ and $\isSemiActedUponBy'$ are similar (see \cref{lemma:liberation-and-coordinate-systems}), for each $m \in M$, there is a $g_{m,0}' \in G_0'$ such that
    \begin{equation}
    \label{equation:global-transition-function-of-product-does-not-depend-on-a-lot-of-things:zero}
      \ForEach \mathfrak{g} \in G \modulo G_0 \Holds m \isSemiActedUponBy \mathfrak{g} = m \isSemiActedUponBy' g_{m,0}' \cdot (g \conjugates \mathfrak{g}).
    \end{equation}
    And, because $\Delta_{\modulo H} = \Delta_{\modulo H}'$, according to \cref{theorem:global-transition-function-determines-local-transition-function}, there are cellular automata $\mathcal{C}_{\modulo H, *} = \ntuple{\ntuple{\mathcal{M}, \mathcal{K}} \modulo (H, \rho), Q, N_{\modulo H, *}, \delta_{\modulo H, *}}$ and $\mathcal{C}_{\modulo H, *}' = \ntuple{\ntuple{\mathcal{M}, \mathcal{K}'} \modulo (H, \rho'), Q, N_{\modulo H, *}', \delta_{\modulo H, *}'}$ such that $\delta_{\modulo H, *} = g^{-1} H \otimes_{\modulo H} \delta_{\modulo H, *}'$ and
    \begin{gather}
    \label{equation:global-transition-function-of-product-does-not-depend-on-a-lot-of-things:one}
      \ForEach \ell_{\modulo H} \in Q^{N_{\modulo H}} \Holds \delta_{\modulo H}(\ell_{\modulo H}) = \delta_{\modulo H, *}(\ell_{\modulo H}\restrictedTo_{N_{\modulo H, *}}),\\
    \label{equation:global-transition-function-of-product-does-not-depend-on-a-lot-of-things:two}
      \ForEach \ell_{\modulo H}' \in Q^{N_{\modulo H}'} \Holds \delta_{\modulo H}'(\ell_{\modulo H}') = \delta_{\modulo H, *}'(\ell_{\modulo H}'\restrictedTo_{N_{\modulo H, *}'}).
    \end{gather}

    Let $c \in Q^M$ and let $m \in M$. Then, according to \cref{equation:global-transition-function-of-product-does-not-depend-on-a-lot-of-things:zero},
    \begin{align*}
      \Delta(c)(m)
      &= \delta(n \mapsto c(m \isSemiActedUponBy n))\\
      &= \delta_{\modulo H}(n_{\modulo H} \mapsto c(m \isSemiActedUponBy \kappa(n_{\modulo H})))\\
      &= \delta_{\modulo H}(n_{\modulo H} \mapsto c(m \isSemiActedUponBy' g_{m,0}' \cdot (g \conjugates \kappa(n_{\modulo H})))).
    \end{align*}
    Thus, according to \cref{remark:actions-circ-and-cdot-commute},
    \begin{equation*}
      \Delta(c)(m) = \delta_{\modulo H}(n_{\modulo H} \mapsto c(m \isSemiActedUponBy' g \conjugates (g^{-1} g_{m,0}' g \cdot \kappa(n_{\modulo H})))).
    \end{equation*}
    Thus, because $\kappa$ is $G_0$-e\-qui\-var\-i\-ant,
    \begin{align*}
      \Delta(c)(m)
      &= \delta_{\modulo H}(n_{\modulo H} \mapsto c(m \isSemiActedUponBy' g \conjugates \kappa(g^{-1} g_{m,0}' g \cdot n_{\modulo H})))\\
      &= \delta_{\modulo H}((g^{-1} g_{m,0}' g)^{-1} \bullet_{\modulo H} [n_{\modulo H} \mapsto c(m \isSemiActedUponBy' g \conjugates \kappa(n_{\modulo H}))]).
    \end{align*}
    Thus, because $\delta_{\modulo H}$ is $\bullet_{\modulo H}$-invariant,
    \begin{equation*}
      \Delta(c)(m) = \delta_{\modulo H}(n_{\modulo H} \mapsto c(m \isSemiActedUponBy' g \conjugates \kappa(n_{\modulo H}))).
    \end{equation*}
    Thus, according to \cref{equation:global-transition-function-of-product-does-not-depend-on-a-lot-of-things:one},
    \begin{equation*}
      \Delta(c)(m) = \delta_{\modulo H, *}(n_{\modulo H, *} \mapsto c(m \isSemiActedUponBy' g \conjugates \kappa(n_{\modulo H, *}))).
    \end{equation*}
    Thus, because $\delta_{\modulo H, *} = g^{-1} H \otimes_{\modulo H} \delta_{\modulo H, *}'$,
    \begin{align*}
      \Delta(c)(m)
      &= (g^{-1} H \otimes_{\modulo H} \delta_{\modulo H, *}')(n_{\modulo H, *} \mapsto c(m \isSemiActedUponBy' g \conjugates \kappa(n_{\modulo H, *})))\\
      &= \delta_{\modulo H, *}'(n_{\modulo H, *}' \mapsto c(m \isSemiActedUponBy' g \conjugates \kappa(g^{-1} H \conjugates_{\modulo H} n_{\modulo H, *}')))\\
      &= \delta_{\modulo H, *}'(n_{\modulo H, *}' \mapsto c(m \isSemiActedUponBy' (g \boxdot \kappa)(n_{\modulo H, *}'))).
    \end{align*}
    Thus, because $g \boxdot \kappa = \kappa'$,
    \begin{equation*}
      \Delta(c)(m) = \delta_{\modulo H, *}'(n_{\modulo H, *}' \mapsto c(m \isSemiActedUponBy' \kappa'(n_{\modulo H, *}'))).
    \end{equation*}
    Thus, according to \cref{equation:global-transition-function-of-product-does-not-depend-on-a-lot-of-things:two},
    \begin{align*}
      \Delta(c)(m)
      &= \delta_{\modulo H}(n_{\modulo H}' \mapsto c(m \isSemiActedUponBy' \kappa'(n_{\modulo H}')))\\
      &= \delta'(n \mapsto c(m \isSemiActedUponBy' n))\\
      &= \Delta'(c)(m).
    \end{align*}
    In conclusion, $\Delta = \Delta'$.
  \end{proof}

  The notion of product of a global transition function is well-defined and we give it a name in

  \begin{lemma} 
  \label{lemma:non-naive-product-of-cellular-automaton}
    Let $\mathcal{M} = \ntuple{M, G, \actsOnPoint}$ be a left-ho\-mo\-ge\-neous space, let $H$ be a normal subgroup of $G$, let $\Delta_{\modulo H}$ be the global transition function of a cellular automaton over $\mathcal{M} \modulo H$, and let $K$ be an element of the orbit space $G \reverseModulo \mathfrak{k}$ under $\boxdot$. The global transition function $\Delta$ of a product of $\mathcal{C}_{\modulo H}$ by $H$ and $\kappa$ -- where $\mathcal{K} = \ntuple{m_0, \family{g_{m_0, m}}_{m \in M}}$ is a coordinate system for $\mathcal{M}$, $\rho$ is a right inverse of the canonical projection $\pi \from M \to H \reverseModulo M$ such that $\rho(H \actsOnPoint m_0) = m_0$, $\mathcal{C}_{\modulo H}$ is a cellular automaton over $\ntuple{\mathcal{M}, \mathcal{K}} \modulo (H, \rho)$ whose global transition function is $\Delta_{\modulo H}$, and $\kappa$ is an element of $K \cap \mathfrak{k}_{m_0}$ -- is the global transition function of a cellular automaton over $\mathcal{M}$ and is called \define{product of $\Delta_{\modulo H}$ by $H$ and $K$}\graffito{product $\Delta$ of $\Delta_{\modulo H}$ by $H$ and $K$}.
  \end{lemma}

  \begin{proof}
    According to \cref{lemma:global-transition-function-of-product-does-not-depend-on-a-lot-of-things}, the global transition function $\Delta$ does not depend on the choice of $\mathcal{K}$, $\rho$, $\mathcal{C}_{\modulo H}$, and $\kappa$; and it is hence well-defined.
  \end{proof}

  \begin{example}[Thick Tree]
  \label{example:thick-tree:product}
    In the situation of \cref{example:thick-tree:product-cellular-automaton}, the product $\Delta$ of $\Delta_{\modulo H}$ by $H$ and $G \boxdot \kappa$ is the product $\prod_{h \in H} \Delta_{\modulo H}$, which does neither depend on $\mathcal{K}$, nor on $\rho$, nor on $\kappa$.
  \end{example}

  The quotient of the product of a global transition function is identical to the original global transition function, as shown in

  \begin{lemma}
  \label{lemma:product-then-quotient-is-identity-ca-function}
    Let $\mathcal{M} = \ntuple{M, G, \actsOnPoint}$ be a left-ho\-mo\-ge\-neous space, let $H$ be a normal subgroup of $G$, let $\Delta_{\modulo H}$ be the global transition function of a cellular automaton over $\mathcal{M} \modulo H$, let $\Delta$ be a product of $\Delta_{\modulo H}$ by $H$ and $K$, and let $\Delta_{\modulo H}'$ be the quotient of $\Delta$ by $H$. Then, $\Delta_{\modulo H}' = \Delta_{\modulo H}$.
  \end{lemma}

  \begin{proof}
    This is a direct consequence of \cref{lemma:product-then-quotient-is-identity}.
  \end{proof}

  The product of the quotient of a global transition function may not be the original global transition function, as illustrated in

  \begin{example}[Thick Tree]
  \label{example:thick-tree:product-of-quotient}
    In the situation of \cref{example:thick-tree:funny-equivariance-for-product}, let $\mathcal{K} = \ntuple{e_G G_0, \family{g_{e_G G_0, g G_0}}_{g G_0 \in G \modulo G_0}}$ be a coordinate system for $\mathcal{M} = \ntuple{G \modulo G_0, G, \actsOnPoint}$; let $Q$ be the binary set $\setOf{0, 1}$; let $N$ be the subset $\setOf{-1 + 7\Z, 0 + 7\Z, 1 + 7\Z} \times \setOf{e_F G_0}$ of $G \modulo G_0$ (note that $G_0 \cdot N \subseteq N$); let $\delta$ be the $\bullet$-invariant map $Q^N \to Q$, $\ell \mapsto 0$, if $\sum_{n \in N} \ell(n) < 3/2$, and $\ell \mapsto 1$, otherwise, which is known as \define{majority rule}\graffito{majority rule}\index{rule!majority}; let $\mathcal{C}$ be the cellular automaton $\ntuple{\ntuple{\mathcal{M}, \mathcal{K}}, Q, N, \delta}$; and let $\Delta$ be the global transition function of $\mathcal{C}$, which realises the majority rule on each of the discrete circles $H \times \setOf{f G_0}$, for $f G_0 \in F \modulo G_0$. Each quotient of $\mathcal{C}$ by $H$ has the neighbourhood $N_{\modulo H} = \setOf{e_F G_0}$, the local transition function $\delta_{\modulo H} \from Q^{\setOf{e_F G_0}} \to Q$, $\ell_{\modulo H} \mapsto \ell_{\modulo H}(e_F G_0)$, and the identity map on $Q^{F \modulo G_0}$ for global transition function $\Delta_{\modulo H}$. And, each product of $\Delta_{\modulo H}$ by $H$ is the identity map on $Q^{G \modulo G_0}$, which is not $\Delta$.
  \end{example}

  \section{Restrictions}
  \label{section:restrictions}


  There is a canonical bijection from $H \modulo G_0$ to $H \modulo H_0$, which is given in

  \begin{lemma}
  \label{lemma:canonical-bijection-of-H-quotient-G-zero-to-H-quotient-H-zero}
    Let $G$ be a group, let $G_0$ and $H$ be two subgroups of $G$, and let $H_0$ be the subgroup $G_0 \cap H$ of $H$. The map
    \begin{align*}
      \zeta \from H \modulo G_0 &\to H \modulo H_0, \mathnote{canonical bijection $\zeta$ from $H \modulo G_0$ to $H \modulo H_0$}\index[symbols]{zeta@$\zeta$}\\
      h G_0 &\mapsto h H_0,
    \end{align*}
    is well-defined, is bijective, is $\cdot$-e\-qui\-var\-i\-ant, and is called \define{canonical bijection from $H \modulo G_0$ to $H \modulo H_0$}.
  \end{lemma}

  \begin{proof}
    For each $h \in H$ and each $h' \in H$,
    \begin{align*}
      h G_0 = h' G_0
      &\ifAndOnlyIf h^{-1} h' \in G_0\\
      &\ifAndOnlyIf h^{-1} h' \in H_0\\
      &\ifAndOnlyIf h H_0 = h' H_0.
    \end{align*}
    Hence, $\zeta$ is well-defined and bijective. Moreover, for each $h \in H$ and each $h' G_0 \in H \modulo G_0$, we have $h \cdot \zeta(h' G_0) = h h' H_0 = \zeta(h \cdot h' G_0)$. Therefore, $\zeta$ is $\cdot$-e\-qui\-var\-i\-ant.
  \end{proof}

  \begin{lemma}
  \label{lemma:normal-subgroup-left-and-right-actions}
    Let $\mathcal{R} = \ntuple{\ntuple{M, G, \actsOnPoint}, \ntuple{m_0, \family{g_{m_0, m}}_{m \in M}}}$ be a cell space and let $H$ be a normal subgroup of $G$. Then,
    \begin{equation*}
      \ForEach m \in M \Holds H \actsOnPoint m = m \isSemiActedUponBy (H \modulo G_0). \qedhere 
    \end{equation*}
  \end{lemma}

  \begin{proof}
    For each $m \in M$, because $g_{m_0, m}^{-1} H = H g_{m_0, m}^{-1}$,
    \begin{align*}
      H \actsOnPoint m
      &= g_{m_0, m} g_{m_0, m}^{-1} H \actsOnPoint m\\
      &= g_{m_0, m} H g_{m_0, m}^{-1} \actsOnPoint m\\
      &= m \isSemiActedUponBy (H \modulo G_0). \qedhere
    \end{align*}
  \end{proof}

  Restrictions of left-ho\-mo\-ge\-neous spaces, coordinate systems, and cell spaces are introduced in the two forthcoming definitions.

  \begin{definition}
    Let $\mathcal{M} = \ntuple{M, G, \actsOnPoint}$ be a left group set, let $m_0$ be an element of $M$, and let $H$ be a subgroup of $G$. The triple $\mathcal{M}\restrictedTo_{m_0, H} = \ntuple{H \actsOnPoint m_0, H, \actsOnPoint\restrictedTo_{H \times (H \actsOnPoint m_0) \to H \actsOnPoint m_0}}$\graffito{restriction $\mathcal{M}\restrictedTo_{m_0, H}$ of $\mathcal{M}$ at $m_0$ to $H$} is a left-ho\-mo\-ge\-neous space and is called \define{restriction of $\mathcal{M}$ at $m_0$ to $H$}\index[symbols]{Mm0Hcalligraphicharpoon@$\mathcal{M}\restrictedTo_{m_0, H}$}.
  \end{definition}

  \begin{definition}
    Let $\mathcal{M} = \ntuple{M, G, \actsOnPoint}$ be a left-ho\-mo\-ge\-neous space, let $\mathcal{K} = \ntuple{m_0, \family{g_{m_0, m}}_{m \in M}}$ be a coordinate system for $\mathcal{M}$, let $\mathcal{R}$ be the cell space $\ntuple{\mathcal{M}, \mathcal{K}}$, and let $H$ be a subgroup of $G$ such that $\setOf{g_{m_0, h \actsOnPoint m_0} \suchThat h \in H} \subseteq H$. The tuple $\mathcal{K}\restrictedTo_H = \ntuple{m_0, \family{g_{m_0, h \actsOnPoint m_0}}_{h \actsOnPoint m_0 \in H \actsOnPoint m_0}}$\graffito{restriction $\mathcal{K}\restrictedTo_H$ of $\mathcal{K}$ to $H$} is a coordinate system for $\mathcal{M}\restrictedTo_{m_0, H}$ and is called \define{restriction of $\mathcal{K}$ to $H$}\index[symbols]{KHcalligraphicharpoon@$\mathcal{K}\restrictedTo_H$}. And, the tuple $\mathcal{R}\restrictedTo_H = \ntuple{\mathcal{M}\restrictedTo_{m_0, H}, \mathcal{K}\restrictedTo_H}$\graffito{restriction $\mathcal{R}\restrictedTo_H$ of $\mathcal{R}$ to $H$} is a cell space and is called \define{restriction of $\mathcal{R}$ to $H$}\index[symbols]{RHcalligraphicharpoon@$\mathcal{R}\restrictedTo_H$}.
  \end{definition}

  \begin{example}[Sphere]
  \label{example:sphere:restriction-of-cell-space}
    Let $M$ be the Euclidean unit $2$-sphere, let $G$ be the rotation group, let $\actsOnPoint$ be the left group action of $G$ on $M$ by function application, let $m_0$ be the north pole $(0,0,1)^\transposed$ of $M$ and, for each point $m \in M$, let $g_{m_0, m}$ be a rotation about an axis in the $(x, y)$-plane that rotates $m_0$ to $m$, which is unique unless $m$ is the south pole, $\mathcal{R} = \ntuple{\mathcal{M}, \mathcal{K}}$ be the cell space $\ntuple{\ntuple{M, G, \actsOnPoint}, \ntuple{m_0, \family{g_{m_0, m}}_{m \in M}}}$ (this is the situation of \cref{example:sphere:cell-space}).

    Furthermore, let $a$ be the axis of the rotation $g_{m_0, S}$, where $S$ is the south pole $(0, 0, -1)^{\transposed}$, and let $H$ be the subgroup of $G$ that consists of the rotations about $a$. The set $H \actsOnPoint m_0$ is a great circle through the north and the south pole, the abelian group $H$ is the group of orientation-preserving symmetries of this circle, the restriction $\actsOnPoint\restrictedTo_{H \times (H \actsOnPoint m_0) \to H \actsOnPoint m_0}$ is the free left group action of $H$ on $H \actsOnPoint m_0$ by function application, and we have $\setOf{g_{m_0, h \actsOnPoint m_0} \suchThat h \actsOnPoint m_0 \in H \actsOnPoint m_0} \subseteq H$. Hence, the restriction of $\mathcal{M}$ at $m_0$ to $H$ is isomorphic to the abelian group $\R \modulo \Z$ and the restriction of $\mathcal{K}$ to $H$ is the unique coordinate system for $\mathcal{M}\restrictedTo_{m_0, H}$ whose origin is $m_0$.
  \end{example}

  In the remainder of this section, let $\mathcal{R} = \ntuple{\ntuple{M, G, \actsOnPoint}, \ntuple{m_0, \family{g_{m_0, m}}_{m \in M}}}$ be a cell space, let $H$ be a subgroup of $G$, and let $\mathcal{R}\restrictedTo_H = \ntuple{\ntuple{M_{\restrictedTo_H}, H, \actsOnPoint_{\restrictedTo_H}}, \ntuple{m_0, \family{g_{m_0, m_{\restrictedTo_H}}}_{m_{\restrictedTo_H} \in M_{\restrictedTo_H}}}}$ be the restriction of $\mathcal{R}$ to $H$.

  The right quotient set semi-actions of a cell space and its restriction relate as shown in

  \begin{lemma}
  \label{lemma:liberation-of-restriction}
    The stabiliser of $m_0$ under $\actsOnPoint_{\restrictedTo_H}$ is $H_0 = G_0 \cap H$ and the right quotient set semi-action induced by $\mathcal{R}\restrictedTo_H$ is
    \begin{align*}
      \isSemiActedUponBy_{\restrictedTo_H} \from M_{\restrictedTo_H} \times H \modulo H_0 &\to M_{\restrictedTo_H},\\
      (m_{\restrictedTo_H}, h H_0) &\mapsto m_{\restrictedTo_H} \isSemiActedUponBy h G_0. \qedhere
    \end{align*}
  \end{lemma}

  \begin{proof}
    The subgroup $H_0$ of $H$ is the stabiliser of $m_0$ under $\actsOnPoint_{\restrictedTo_H}$. Moreover, according to \cref{lemma:canonical-bijection-of-H-quotient-G-zero-to-H-quotient-H-zero}, for each $h \in H$, the coset $h G_0$ is uniquely determined by $h H_0$. And, for each $m_{\restrictedTo_H} \in M_{\restrictedTo_H}$ and each $h H_0 \in H \modulo H_0$,
    \begin{align*}
      m_{\restrictedTo_H} \isSemiActedUponBy_{\restrictedTo_H} h H_0
      &= m_{\restrictedTo_H} \isSemiActedUponBy h G_0\\
      &= g_{m_0, m_{\restrictedTo_H}} h g_{m_0, m_{\restrictedTo_H}}^{-1} \actsOnPoint m_{\restrictedTo_H}\\
      &= g_{m_0, m_{\restrictedTo_H}} h g_{m_0, m_{\restrictedTo_H}}^{-1} \actsOnPoint_{\restrictedTo_H} m_{\restrictedTo_H}\\
      &\in M_{\restrictedTo_H}.
    \end{align*}
    In conclusion, the map $\isSemiActedUponBy_{\restrictedTo_H}$ is well-defined and it is the right quotient set semi-action induced by $\mathcal{R}\restrictedTo_H$.
  \end{proof}

  Restrictions of (semi-)cellular automata are introduced in

  \begin{lemma}
    Let $\mathcal{C} = \ntuple{\mathcal{R}, Q, N, \delta}$ be a semi-cellular or cellular automaton with a sufficient neighbourhood that is included in $H \modulo G_0$, let $E$ be the biggest sufficient neighbourhood of $\mathcal{C}$ such that $H_0 \cdot E \subseteq E$ and $E \subseteq H \modulo G_0$, let $\eta$ be the sufficient local transition function from $Q^E$ to $Q$, let $\zeta$ be the canonical bijection from $H \modulo G_0$ to $H \modulo H_0$, let
    \begin{equation*}
      N_{\restrictedTo_H} = \zeta(E) \quad (= \setOf{h H_0 \suchThat h G_0 \in E}),
    \end{equation*}
    and let
    \begin{align*}
      \delta_{\restrictedTo_H} \from Q^{N_{\restrictedTo_H}} &\to Q,\\
      \ell_{\restrictedTo_H} &\mapsto \eta\parens{e \mapsto \ell_{\restrictedTo_H}(\zeta(e))} \quad (= \eta\parens{h G_0 \mapsto \ell_{\restrictedTo_H}(h H_0)}).
    \end{align*}
    The quadruple $\mathcal{C}\restrictedTo_H = \ntuple{\mathcal{R}\restrictedTo_H, Q, N_{\restrictedTo_H}, \delta_{\restrictedTo_H}}$\graffito{restriction $\mathcal{C}\restrictedTo_H$ of $\mathcal{C}$ to $H$} is a semi-cellular or cellular automaton respectively, and is called \define{restriction of $\mathcal{C}$ to $H$}\index[symbols]{CHcalligraphicharpoon@$\mathcal{C}\restrictedTo_H$}.
  \end{lemma}

  \begin{proof} 
    First, because $H_0 \cdot \zeta(E) = \zeta(H_0 \cdot E)$ and $H_0 \cdot E \subseteq E$, we have $H_0 \cdot N_{\restrictedTo_H} \subseteq N_{\restrictedTo_H}$. In conclusion, the quadruple $\mathcal{C}\restrictedTo_H$ is a semi-cellular automaton.

    Secondly, let $\mathcal{C}$ be a cellular automaton. Furthermore, let $h_0 \in H_0$. Then, because $\zeta$ is $\cdot$-e\-qui\-var\-i\-ant and $\eta$ is $\bullet_{H_0}$-invariant, for each $\ell_{\restrictedTo_H} \in Q^{N_{\restrictedTo_H}}$,
    \begin{align*} 
      \delta_{\restrictedTo_H}(h_0 \bullet_{\restrictedTo_H} \ell_{\restrictedTo_H})
      &= \eta\parens[\big]{e \mapsto (h_0 \bullet_{\restrictedTo_H} \ell_{\restrictedTo_H})(\zeta(e)}\\
      &= \eta\parens[\big]{e \mapsto \ell_{\restrictedTo_H}(h_0^{-1} \cdot \zeta(e)}\\
      &= \eta\parens[\big]{e \mapsto \ell_{\restrictedTo_H}(\zeta(h_0^{-1} \cdot e)}\\
      &= \eta\parens[\big]{h_0 \bullet \brackets[\big]{e \mapsto \ell_{\restrictedTo_H}(\zeta(e))}}\\
      &= \eta\parens[\big]{e \mapsto \ell_{\restrictedTo_H}(\zeta(e))}\\
      &= \delta_{\restrictedTo_H}(\ell_{\restrictedTo_H}).
    \end{align*}
    In conclusion, $\delta_{\restrictedTo_H}$ is $\bullet_{\restrictedTo_H}$-invariant and hence $\mathcal{C}\restrictedTo_H$ is a cellular automaton.
  \end{proof}

  \begin{example}[Sphere]
  \label{example:sphere:restriction-of-automaton}
    In the situation of \cref{example:sphere:restriction-of-cell-space}, let $h$ be one of the two rotations about $a$ by $1\degree$, let $Q$ be the set $\setOf{0, 1}$, let $E$ be the singleton set $\setOf{h G_0}$, let $N$ be the set $G_0 \actsOnPoint E$, let $\eta$ be the map $Q^E \to Q$, $\ell_E \mapsto \ell_E(h G_0)$, and let $\delta$ be the map $Q^N \to Q$, $\ell \mapsto \eta(\ell\restrictedTo_E)$. The restriction of the semi-cellular automaton $\mathcal{C} = \ntuple{\mathcal{R}, Q, N, \delta}$ to $H$ is the cellular automaton $\mathcal{C}\restrictedTo_H = \ntuple{\mathcal{R}\restrictedTo_H, Q, E, \eta}$ whose global transition function is the rotation map $h^{-1} \actsOnMap_{\restrictedTo_H} \blank$.
  \end{example}

  The global transition functions of a semi-cellular automaton and its restriction relate as shown in

  \begin{lemma}
  \label{lemma:global-transition-function-of-restriction-versus-the-unrestricted-one}
    Let $\mathcal{C}\restrictedTo_H$ be a restriction of $\mathcal{C} = \ntuple{\mathcal{R}, Q, N, \delta}$ to $H$. Then,
    \begin{equation*}
      \ForEach c \in Q^M \Holds \Delta_{\restrictedTo_H}(c\restrictedTo_{M_{\restrictedTo_H}}) = \Delta(c)\restrictedTo_{M_{\restrictedTo_H}}. \qedhere
    \end{equation*}
  \end{lemma}

  \begin{proof}
    Let $\eta \from Q^E \to Q$ be a sufficient local transition function of $\mathcal{C}$ such that $H_0 \cdot E \subseteq E$ and $E \subseteq H \modulo G_0$, and let $\zeta$ be the canonical bijection from $H \modulo G_0$ to $H \modulo H_0$. Furthermore, let $c \in Q^M$. Then, for each $m_{\restrictedTo_H} \in M_{\restrictedTo_H}$, according to \cref{lemma:liberation-of-restriction},
    \begin{align*}
      \Delta_{\restrictedTo_H}(c\restrictedTo_{M_{\restrictedTo_H}})(m_{\restrictedTo_H})
      &= \delta_{\restrictedTo_H}\parens[\big]{n_{\restrictedTo_H} \mapsto c\restrictedTo_{M_{\restrictedTo_H}}(m_{\restrictedTo_H} \isSemiActedUponBy_{\restrictedTo_H} n_{\restrictedTo_H})}\\
      &= \eta\parens[\big]{e \mapsto c\restrictedTo_{M_{\restrictedTo_H}}(m_{\restrictedTo_H} \isSemiActedUponBy_{\restrictedTo_H} \zeta(e))}\\
      &= \eta\parens[\big]{e \mapsto c\restrictedTo_{M_{\restrictedTo_H}}(m_{\restrictedTo_H} \isSemiActedUponBy e)}\\
      &= \eta\parens[\big]{e \mapsto c(m_{\restrictedTo_H} \isSemiActedUponBy e)}\\
      &= \delta\parens[\big]{n \mapsto c(m_{\restrictedTo_H} \isSemiActedUponBy n)}\\
      &= \Delta(c)(m_{\restrictedTo_H}).
    \end{align*}
    In conclusion, $\Delta_{\restrictedTo_H}(c\restrictedTo_{M_{\restrictedTo_H}}) = \Delta(c)\restrictedTo_{M_{\restrictedTo_H}}$.
  \end{proof}

  Restrictions of global transition functions of cellular automata are introduced in

  \begin{lemma}
  \label{lemma:restriction-of-global-transition-function}
    Let $\mathcal{M} = \ntuple{M, G, \actsOnPoint}$ be a left-ho\-mo\-ge\-neous space, let $H$ be a subgroup of $G$, let $m_0$ be an element of $M$, and let $\Delta \from Q^M \to Q^M$ be the global transition function of a cellular automaton over $\mathcal{M}$ with sufficient neighbourhood $E$ such that $E \subseteq H \modulo G_0$. Furthermore, let $\pi$ be the canonical projection from $Q^M$ onto $Q^{H \actsOnPoint m_0}$ and let $\rho$ be a right inverse of $\pi$. The map $\Delta_{\restrictedTo_{m_0,H}} = \pi \after \Delta \after \rho$\graffito{restriction $\Delta_{\restrictedTo_{m_0, H}}$ of $\Delta$ at $m_0$ to $H$} is the global transition function of a cellular automaton over $\ntuple{H \actsOnPoint m_0, H, \actsOnPoint\restrictedTo_{H \times (H \actsOnPoint m_0) \to H \actsOnPoint m_0}}$, does not depend on $\rho$, and is called \define{restriction of $\Delta$ at $m_0$ to $H$}\index[symbols]{Deltam0Hharpoonm@$\Delta_{\restrictedTo_{m_0,H}}$}.
  \end{lemma}

  \begin{proof}
    There is a coordinate system $\mathcal{K}$ with origin $m_0$ and a cellular automaton $\mathcal{C} = \ntuple{\ntuple{\mathcal{M}, \mathcal{K}}, Q, N, \delta}$ with sufficient neighbourhood $E$ whose global transition function is $\Delta$. Moreover, there is a coordinate system $\mathcal{K}' = \ntuple{m_0, \family{g_{m_0, m}'}_{m \in M}}$ for $\mathcal{M}$ such that, for each $m \in H \actsOnPoint m_0$, we have $g_{m_0, m}' \in H$. And, according to \cref{corollary:independence-of-coordinate-system}, the global transition function of the cellular automaton $\mathcal{C}' = \ntuple{\ntuple{\mathcal{M}, \mathcal{K}'}, Q, N, \delta}$ is $\Delta$. And, the set $E' = H_0 \cdot E$ is a sufficient neighbourhood of $\mathcal{C}'$ such that $H_0 \cdot E' \subseteq E'$ and $E' \subseteq H \modulo G_0$. Hence, according to \cref{lemma:global-transition-function-of-restriction-versus-the-unrestricted-one}, the map $\Delta_{\restrictedTo_{m_0,H}}$ is the global transition function of the restriction of $\mathcal{C}'$ to $H$, which does not depend on $\rho$.
  \end{proof}

  %

  \section{Extensions}
  \label{section:extensions}

  Extensions of semi-cellular automata are introduced in

  \begin{definition}
    Let $\mathcal{C}_{\restrictedTo_H} = \ntuple{\mathcal{R}\restrictedTo_H, Q, N_{\restrictedTo_H}, \delta_{\restrictedTo_H}}$ be a semi-cellular automaton, let $\zeta$ be the canonical bijection from $H \modulo G_0$ to $H \modulo H_0$, let
    \begin{equation*}
      N = G_0 \cdot \zeta^{-1}(N_{\restrictedTo_H}) = \setOf{g_0 h G_0 \suchThat g_0 \in G_0, h H_0 \in N_{\restrictedTo_H}},
    \end{equation*}
    and let
    \begin{align*}
      \delta \from Q^N &\to Q,\\
      \ell &\mapsto \delta_{\restrictedTo_H}\parens[\big]{n_{\restrictedTo_H} \mapsto \ell(\zeta^{-1}(n_{\restrictedTo_H}))} \quad (= \delta_{\restrictedTo_H}\parens[\big]{h H_0 \mapsto \ell(h G_0)}).
    \end{align*}
    The quadruple $\mathcal{C} = \ntuple{\mathcal{R}, Q, N, \delta}$\graffito{extension $\mathcal{C}$ of $\mathcal{C}_{\restrictedTo_H}$ to $G$} is a semi-cellular automaton and is called \define{extension of $\mathcal{C}_{\restrictedTo_H}$ to $G$}\index[symbols]{CHcalligraphicharpoon@$\mathcal{C}_{\restrictedTo_H}$}.
  \end{definition}

  \begin{example}[Sphere]
  \label{example:sphere:extension-of-automaton}
    In the situation of \cref{example:sphere:restriction-of-automaton}, the extension of the cellular automaton $\mathcal{C}\restrictedTo_H$ to $G$ is the semi-cellular automaton $\mathcal{C}$.
  \end{example}

  The global transition functions of a semi-cellular automaton and its extension relate as shown in

  \begin{lemma}
  \label{lemma:global-transition-function-of-extension}
    Let $\mathcal{C}$ be an extension of $\mathcal{C}_{\restrictedTo_H} = \ntuple{\mathcal{R}\restrictedTo_H, Q, N_{\restrictedTo_H}, \delta_{\restrictedTo_H}}$ to $G$. Then,
    \begin{equation*}
      \ForEach c \in Q^M \Holds \Delta(c)\restrictedTo_{M_{\restrictedTo_H}} = \Delta_{\restrictedTo_H}(c\restrictedTo_{M_{\restrictedTo_H}}),
    \end{equation*}
    and
    \begin{equation*}
      \ForEach c \in Q^M \ForEach m \in M \Holds \Delta(c)(m) = \Delta_{\restrictedTo_H}((g_{m_0, m}^{-1} \actsOnMap c)\restrictedTo_{M_{\restrictedTo_H}})(m_0). \qedhere
    \end{equation*}
  \end{lemma}

  \begin{proof}
    Let $\zeta$ be the canonical bijection from $H \modulo G_0$ to $H \modulo H_0$.

    First, let $c \in Q^M$. Then, for each $m_{\restrictedTo_H} \in M_{\restrictedTo_H}$, according to \cref{lemma:liberation-of-restriction},
    \begin{align*}
      \Delta(c)(m_{\restrictedTo_H})
      &= \delta\parens[\big]{n \mapsto c(m_{\restrictedTo_H} \isSemiActedUponBy n)}\\
      &= \delta_{\restrictedTo_H}\parens[\big]{n_{\restrictedTo_H} \mapsto c(m_{\restrictedTo_H} \isSemiActedUponBy \zeta^{-1}(n_{\restrictedTo_H}))}\\
      &= \delta_{\restrictedTo_H}\parens[\big]{n_{\restrictedTo_H} \mapsto c(m_{\restrictedTo_H} \isSemiActedUponBy_{\restrictedTo_H} n_{\restrictedTo_H})}\\
      &= \delta_{\restrictedTo_H}\parens[\big]{n_{\restrictedTo_H} \mapsto c\restrictedTo_{M_{\restrictedTo_H}}(m_{\restrictedTo_H} \isSemiActedUponBy_{\restrictedTo_H} n_{\restrictedTo_H})}\\
      &= \Delta_{\restrictedTo_H}(c\restrictedTo_{M_{\restrictedTo_H}})(m_{\restrictedTo_H}).
    \end{align*}
    In conclusion, $\Delta(c)\restrictedTo_{M_{\restrictedTo_H}} = \Delta_{\restrictedTo_H}(c\restrictedTo_{M_{\restrictedTo_H}})$.

    Secondly, for each $c \in Q^M$ and each $m \in M$, according to \cref{remark:determined-by-behaviour-at-origin-without-assumptions},
    \begin{equation*}
      \Delta(c)(m)
      = \Delta(g_{m_0, m}^{-1} \actsOnMap c)(m_0)
      = \Delta_{\restrictedTo_H}((g_{m_0, m}^{-1} \actsOnMap c)\restrictedTo_{M_{\restrictedTo_H}})(m_0). \qedhere
    \end{equation*}
  \end{proof}

  The global transition functions of a semi-cellular automaton and the restriction of its extension are identical as shown in

  \begin{corollary}
    Let $\mathcal{C}$ be an extension of $\mathcal{C}_{\restrictedTo_H} = \ntuple{\mathcal{R}\restrictedTo_H, Q, N_{\restrictedTo_H}, \delta_{\restrictedTo_H}}$ to $G$ and let $\mathcal{C}\restrictedTo_H = \ntuple{\mathcal{R}\restrictedTo_H, Q, N_{\restrictedTo_H}', \delta_{\restrictedTo_H}'}$ be the restriction of $\mathcal{C}$ to $H$. The global transition functions of $\mathcal{C}_{\restrictedTo_H}$ and $\mathcal{C}\restrictedTo_H$ are identical.
  \end{corollary}

  \begin{proof}
    This is a direct consequence of \cref{lemma:global-transition-function-of-restriction-versus-the-unrestricted-one,lemma:global-transition-function-of-extension}.
  \end{proof}

  The global transition functions of a semi-cellular automaton and the extension of its restriction are identical as shown in

  \begin{corollary}
    Let $\mathcal{C}\restrictedTo_H$ be a restriction of $\mathcal{C} = \ntuple{\mathcal{R}, Q, N, \delta}$ to $H$ and let $\mathcal{C}' = \ntuple{\mathcal{R}, Q, N', \delta'}$ be the extension of $\mathcal{C}\restrictedTo_H$ to $G$. The global transition functions of $\mathcal{C}$ and $\mathcal{C}'$ are identical.
  \end{corollary}

  \begin{proof}
    For each $c \in Q^M$ and each $m \in M$, according to \cref{lemma:global-transition-function-of-extension}, \cref{lemma:global-transition-function-of-restriction-versus-the-unrestricted-one}, and \cref{remark:determined-by-behaviour-at-origin-without-assumptions},
    \begin{align*}
      \Delta'(c)(m)
      &= \Delta_{\restrictedTo_H}((g_{m_0, m}^{-1} \actsOnMap c)\restrictedTo_{M_{\restrictedTo_H}})(m_0)\\
      &= \Delta(g_{m_0, m}^{-1} \actsOnMap c)(m_0)\\
      &= \Delta(c)(m).
    \end{align*}
    In conclusion, $\Delta' = \Delta$.
  \end{proof}

  \begin{remark}
    The extension of the restriction of a semi-cellular automaton and the restriction of the extension of a semi-cellular automaton is in general not the automaton we started out with, because in the first case superfluous neighbours may be removed by the restriction that are not re-added by the extension and in the second case superfluous neighbours may be added by the extension that are not removed by the restriction.
  \end{remark}

  The extension of a cellular automaton is in general not a cellular automaton. However, if the local transition function is in a certain sense invariant under $\bullet$, then the extension is a cellular automaton, which is shown in

  \begin{lemma}
    Let $\mathcal{C}_{\restrictedTo_H} = \ntuple{\mathcal{R}\restrictedTo_H, Q, N_{\restrictedTo_H}, \delta_{\restrictedTo_H}}$ be a semi-cellular automaton, let $\zeta$ be the canonical bijection from $H \modulo G_0$ to $H \modulo H_0$, let $N$ be the set $G_0 \cdot \zeta^{-1}(N_{\restrictedTo_H})$, let $\xi$ be the restriction of $\zeta^{-1}$ to $N_{\restrictedTo_H} \to N$, such that
    \begin{align} 
    \label{equation:bullet-invariance-on-steroids}
      \ForEach g_0 \in G_0 \ForEach \ell \in Q^N \Holds
          \delta_{\restrictedTo_H}\parens[\big]{(g_0 \bullet \ell) \after \xi} = \delta_{\restrictedTo_H}(\ell \after \xi).
    \end{align} 
    The semi-cellular automaton $\mathcal{C}_{\restrictedTo_H}$ and its extension $\mathcal{C}$ to $G$ are cellular automata.
  \end{lemma}

  \begin{proof}
    Let $h_0 \in H_0$ and let $\ell_{\restrictedTo_H} \in Q^{N_{\restrictedTo_H}}$. Then, because $\xi$ is $\cdot$-e\-qui\-var\-i\-ant, for each $n_{\restrictedTo_H} \in N_{\restrictedTo_H}$,
    \begin{align*}
      (h_0 \bullet_{\restrictedTo_H} \ell_{\restrictedTo_H})(n_{\restrictedTo_H})
      &= \ell_{\restrictedTo_H}(h_0^{-1} \cdot n_{\restrictedTo_H})\\
      &= (\ell_{\restrictedTo_H} \after \xi^{-1})\parens[\big]{\xi(h_0^{-1} \cdot n_{\restrictedTo_H})}\\
      &= (\ell_{\restrictedTo_H} \after \xi^{-1})\parens[\big]{h_0^{-1} \cdot \xi(n_{\restrictedTo_H})}\\
      &= \parens[\big]{h_0 \bullet (\ell_{\restrictedTo_H} \after \xi^{-1})}\parens[\big]{\xi(n_{\restrictedTo_H})}\\
      &= \parens[\Big]{\parens[\big]{h_0 \bullet (\ell_{\restrictedTo_H} \after \xi^{-1})} \after \xi}(n_{\restrictedTo_H}).
    \end{align*}
    Thus, $h_0 \bullet_{\restrictedTo_H} \ell_{\restrictedTo_H} = (h_0 \bullet (\ell_{\restrictedTo_H} \after \xi^{-1})) \after \xi$. Hence, according to \cref{equation:bullet-invariance-on-steroids},
    \begin{align*}
      \delta_{\restrictedTo_H}(h_0 \bullet_{\restrictedTo_H} \ell_{\restrictedTo_H})
      &= \delta_{\restrictedTo_H}\parens[\big]{(h_0 \bullet (\ell_{\restrictedTo_H} \after \xi^{-1})) \after \xi}\\
      &= \delta_{\restrictedTo_H}\parens[\big]{(\ell_{\restrictedTo_H} \after \xi^{-1}) \after \xi}\\
      &= \delta_{\restrictedTo_H}(\ell_{\restrictedTo_H}).
    \end{align*}
    Therefore, $\delta_{\restrictedTo_H}$ is $\bullet_{\restrictedTo_H}$-invariant. In conclusion, $\mathcal{C}_{\restrictedTo_H}$ is a cellular automaton.

    Moreover, according to \cref{equation:bullet-invariance-on-steroids}, for each $g_0 \in G_0$ and each $\ell \in Q^N$,
    \begin{align*}
      \delta(g_0 \bullet \ell)
      &= \delta_{\restrictedTo_H}\parens[\big]{n_{\restrictedTo_H} \mapsto (g_0 \bullet \ell)(\zeta^{-1}(n_{\restrictedTo_H}))}\\
      &= \delta_{\restrictedTo_H}\parens[\big]{(g_0 \bullet \ell) \after \xi}\\
      &= \delta_{\restrictedTo_H}(\ell \after \xi)\\
      &= \delta(\ell).
    \end{align*}
    Hence, $\delta$ is $\bullet$-invariant. In conclusion, $\mathcal{C}$ is a cellular automaton.
  \end{proof}

  \begin{remark}
  \label{remark:extension-trivial-if-stabiliser-included-subgroup}
    If $G_0$ is included in $H$, then $G_0 = H_0$, $N = N_{\restrictedTo_H}$, $\delta = \delta_{\restrictedTo_H}$, $\xi = \identityMap_N$, and \cref{equation:bullet-invariance-on-steroids} just states that $\delta$ is $\bullet$-invariant.
  \end{remark}


  Extensions of global transition functions of cellular automata are introduced in

  \begin{lemma}
    Let $\mathcal{M} = \ntuple{M, G, \actsOnPoint}$ be a left-ho\-mo\-ge\-neous space, let $m_0$ be an element of $M$, let $H$ be a subgroup of $G$ that includes $G_0$, and let $\Delta_{\restrictedTo_H} \from Q^{H \actsOnPoint m_0} \to Q^{H \actsOnPoint m_0}$ be the global transition function of a cellular automaton over $\mathcal{M}\restrictedTo_{m_0, H}$. Furthermore, let $\mathcal{K} = \ntuple{m_0, \family{g_{m_0, m}}_{m \in M}}$ be a coordinate system for $\mathcal{M}$ such that $\setOf{g_{m_0, h \actsOnPoint m_0} \suchThat h \actsOnPoint m_0 \in H \actsOnPoint m_0} \subseteq H$. The map
    \begin{align*}
      \Delta \from Q^M &\to Q^M, \mathnote{extension $\Delta$ of $\Delta_{\restrictedTo_H}$ at $m_0$ to $G$}\\
      c &\mapsto [m \mapsto \Delta_{\restrictedTo_H}((g_{m_0, m}^{-1} \actsOnMap c)\restrictedTo_{H \actsOnPoint m_0})(m_0)],
    \end{align*}
    is the global transition function of a cellular automaton over $\mathcal{M}$, does not depend on the coordinates $\family{g_{m_0, m}}_{m \in M}$, and is called \define{extension of $\Delta_{\restrictedTo_H}$ at $m_0$ to $G$}.
  \end{lemma}

  \begin{proof}
    According to \cref{corollary:existence-in-other-coordinate-system}, there is a cellular automaton $\mathcal{C}_{\restrictedTo_H} = \ntuple{\ntuple{\mathcal{M}\restrictedTo_{m_0, H}, \mathcal{K}\restrictedTo_H}, Q, N, \delta}$ whose global transition function is $\Delta_{\restrictedTo_H}$. According to \cref{remark:extension-trivial-if-stabiliser-included-subgroup}, its extension to $G$ is the cellular automaton $\mathcal{C} = \ntuple{\ntuple{\mathcal{M}, \mathcal{K}}, Q, N, \delta}$. According to \cref{lemma:global-transition-function-of-extension}, the map $\Delta$ is the global transition function of $\mathcal{C}$. According to \cref{corollary:independence-of-coordinate-system}, this global transition function does not depend on the coordinates $\family{g_{m_0, m}}_{m \in M}$ and hence neither does $\Delta$.
  \end{proof}

  %

  \section{Decompositions}
  \label{section:decompositions}

  Given a subgroup $H$ of $G$, the orbit space of $\isSemiActedUponBy\restrictedTo_{M \times H \modulo G_0}$ does in general not partition $M$ (see \cref{example:sphere:partition-of-cell-space}). However, if $H$ satisfies the property $G_0 \cdot (H \modulo G_0) \subseteq H \modulo G_0$, then the aforementioned orbit space partitions $M$ (see \cref{definition:G0-G0-invariant,lemma:partition-of-M-induced-by-semi-action-restriction}). For such a group $H$ and a semi-cellular automaton whose neighbourhood is included in $H \modulo G_0$, the phase space is the product of configurations on orbits and the global transition function of the automaton is the product of its restrictions, in a certain sense, to orbits (see \cref{lemma:partition-of-global-transition-function}). The latter restrictions are conjugations of the restriction of the global transition function to $H$ as introduced above (compare \cref{lemma:restrictions-of-global-transition-function-are-conjugate-to-each-other}) and hence the global transition function is injective, surjective, or bijective if and only if its restriction to $H$ has the respective property (see \cref{theorem:cellular-automaton-injective-if-and-only-if-restriction-injective}).

  \begin{definition}
  \label{definition:G0-G0-invariant}
    Let $G$ be a group, and let $G_0$ and $H$ be two subgroups of $G$. The group $H$ is called \defineX{$(G_0, G_0)$-invariant}{invariantG0G0@$(G_0, G_0)$-invariant}\graffito{$(G_0, G_0)$-invariant}\index{G0G0invariant@$(G_0, G_0)$-invariant} if and only if
    \begin{equation*}
      G_0 \cdot (H \modulo G_0) \subseteq H \modulo G_0. \qedhere
    \end{equation*}
  \end{definition}

  \begin{lemma}
    Let $G$ be a group, let $G_0$ be a subgroup of $G$, and let $H$ be a normal subgroup of $G$. The group $H$ is $(G_0, G_0)$-invariant.
  \end{lemma}

  \begin{proof}
    Let $g_0 \in G_0$ and let $h \in H$. Then, because $H$ is normal in $G$, we have $g_0 H = H g_0$. Hence, there is an $h' \in H$ such that $g_0 h = h' g_0$. Therefore, $g_0 h G_0 = h' g_0 G_0 = h' G_0$. In conclusion, $G_0 \cdot (H \modulo G_0) \subseteq H \modulo G_0$.
  \end{proof}

  \begin{lemma}
    Let $G$ be a group, and let $G_0$ and $H$ be two subgroups of $G$ such that $G_0 \subseteq H$. The group $H$ is $(G_0, G_0)$-invariant.
  \end{lemma}

  \begin{proof}
    Let $g_0 \in G_0$ and let $h \in H$. Then, because $G_0 \subseteq H$, we have $g_0 h \in H$. Hence, $g_0 h G_0 \in H \modulo G_0$. In conclusion, $G_0 \cdot (H \modulo G_0) \subseteq H \modulo G_0$.
  \end{proof}

  %

  For a $(G_0, G_0)$-invariant subgroup $H$ of $G$, the orbit space of $\isSemiActedUponBy\restrictedTo_{M \times H \modulo G_0}$ partitions $M$, which is shown in

  \begin{lemma}
  \label{lemma:partition-of-M-induced-by-semi-action-restriction}
    Let $\mathcal{R} = \ntuple{\ntuple{M, G, \actsOnPoint}, \ntuple{m_0, \family{g_{m_0, m}}_{m \in M}}}$ be a cell space, let $H$ be a $(G_0, G_0)$-invariant subgroup of $G$, and let $\mathfrak{H}$ be the quotient set $H \modulo G_0$. The orbit space $\setOf{m \isSemiActedUponBy \mathfrak{H} \suchThat m \in M}$ partitions $M$.
  \end{lemma}

  \begin{proof}
    Let $\sim$ be the binary relation on $M$ given by
    \begin{equation*}
      \ForEach m \in M \ForEach m' \in M \Holds m \sim m' \ifAndOnlyIf m \in m' \isSemiActedUponBy \mathfrak{H}.
    \end{equation*}

    First, let $m \in M$. Then, $m = m \isSemiActedUponBy G_0$. Thus, $m \in m \isSemiActedUponBy \mathfrak{H}$. Hence, $m \sim m$. In conclusion, $\sim$ is reflexive.

    Secondly, let $m$ and $m' \in M$ such that $m \sim m'$. Then, there is an $h' \in H$ such that $m = m' \isSemiActedUponBy h' G_0$. And, there is a $g_0' \in G_0$ such that
    \begin{equation*}
      \ForEach \mathfrak{g} \in G \modulo G_0 \Holds m' \isSemiActedUponBy h' \cdot \mathfrak{g} = (m' \isSemiActedUponBy h' G_0) \isSemiActedUponBy g_0' \cdot \mathfrak{g}.
    \end{equation*}
    Put $h = (h')^{-1}$. Then,
    \begin{equation*}
      m' = m' \isSemiActedUponBy G_0
         = m' \isSemiActedUponBy h' h G_0
         = (m' \isSemiActedUponBy h' G_0) \isSemiActedUponBy g_0' \cdot h G_0
         = m \isSemiActedUponBy g_0' \cdot h G_0.
    \end{equation*}
    And, because $H$ is $(G_0, G_0)$-invariant, there is an $h'' G_0 \in \mathfrak{H}$ such that $g_0' \cdot h G_0 = h'' G_0$. Hence, $m' = m \isSemiActedUponBy h'' G_0 \in m \isSemiActedUponBy \mathfrak{H}$. Therefore, $m' \sim m$. In conclusion, $\sim$ is symmetric.

    Thirdly, let $m$, $m'$, $m'' \in M$ such that $m \sim m'$ and $m' \sim m''$. Then, there are $h'$, $h'' \in H$ such that $m = m' \isSemiActedUponBy h' G_0$ and $m' = m'' \isSemiActedUponBy h'' G_0$. Hence, $m = (m'' \isSemiActedUponBy h'' G_0) \isSemiActedUponBy h' G_0$. And, there is a $g_0'' \in G_0$ such that
    \begin{equation*}
      \ForEach \mathfrak{g} \in G \modulo G_0 \Holds (m'' \isSemiActedUponBy h'' G_0) \isSemiActedUponBy g_0'' \cdot \mathfrak{g} = m'' \isSemiActedUponBy h'' \cdot \mathfrak{g}.
    \end{equation*}
    Put $g = (g_0'')^{-1} h'$. Then,
    \begin{align*}
      m &= (m'' \isSemiActedUponBy h'' G_0) \isSemiActedUponBy h' G_0\\
        &= (m'' \isSemiActedUponBy h'' G_0) \isSemiActedUponBy g_0'' \cdot g G_0\\
        &= m'' \isSemiActedUponBy h'' \cdot g G_0.
    \end{align*}
    And, because $H$ is $(G_0, G_0)$-invariant, there is an $h G_0 \in \mathfrak{H}$ such that $g G_0 = (g_0'')^{-1} h' G_0 = h G_0$. Therefore, $m = m'' \isSemiActedUponBy h'' \cdot h G_0 \in m'' \isSemiActedUponBy \mathfrak{H}$. Hence, $m \sim m''$. In conclusion, $\sim$ is transitive.

    Altogether, $\sim$ is an equivalence relation. In conclusion, $M \modulo {\sim} = \setOf{m \isSemiActedUponBy \mathfrak{H} \suchThat m \in M}$ partitions $M$.
  \end{proof}

  %

  %

  \begin{counterexample}[Sphere]
  \label{example:sphere:partition-of-cell-space}
    In the situation of \cref{example:sphere:restriction-of-cell-space}, under the identification of $M$ with $G \modulo G_0$ by $\iota$, the set $H \modulo G_0$ is the great circle $H \actsOnPoint m_0$ and the set $G_0 \cdot (H \modulo G_0)$ is the union $\bigcup_{g_0 \in G_0} g_0 \actsOnPoint (H \actsOnPoint m_0)$ of the rotations of the great circle $H \actsOnPoint m_0$ about the vertical axis, which is the whole sphere $M$. It follows that the subgroup $H$ of $G$ is not $(G_0, G_0)$-invariant.

    Moreover, for each point $m \in M$, the set $m \isSemiActedUponBy (H \modulo G_0)$ is the great circle $g_{m_0, m} \actsOnPoint (H \actsOnPoint m_0)$ through $m$, which is equal to $H \actsOnPoint m_0$ if and only if the rotation axis of $g_{m_0, m}$ is $a$, which in turn holds if and only if $m$ lies on $H \actsOnPoint m_0$. As each two different great circles intersect in precisely two points, the orbit space $\setOf{m \isSemiActedUponBy (H \modulo G_0) \suchThat m \in M}$ does not partition the sphere $M$ although it covers it.
  \end{counterexample}

  \begin{lemma}
  \label{lemma:liberation-after-liberation-of-invariant-subgroup-stays-in-liberation}
    Let $\mathcal{R} = \ntuple{\ntuple{M, G, \actsOnPoint}, \ntuple{m_0, \family{g_{m_0, m}}_{m \in M}}}$ be a cell space, let $H$ be a $(G_0, G_0)$-invariant subgroup of $G$, and let $\mathfrak{H}$ be the quotient set $H \modulo G_0$. Then,
    \begin{equation*}
      \ForEach m \in M \Holds (m \isSemiActedUponBy \mathfrak{H}) \isSemiActedUponBy \mathfrak{H}' \subseteq m \isSemiActedUponBy \mathfrak{H}. \qedhere
    \end{equation*}
  \end{lemma}

  \begin{proof}
    Let $m \in M$, let $h G_0 \in \mathfrak{H}$, and let $h' G_0 \in \mathfrak{H}$. Then, because $\isSemiActedUponBy$ has defect $G_0$, there is a $g_0 \in G_0$ such that $(m \isSemiActedUponBy h G_0) \isSemiActedUponBy h' G_0 = m \isSemiActedUponBy h g_0 h' G_0$. And, because $H$ is $(G_0, G_0)$-invariant, there is an $h'' \in H$ such that $g_0 h' G_0 = h'' G_0$. Therefore, $(m \isSemiActedUponBy h G_0) \isSemiActedUponBy h' G_0 = m \isSemiActedUponBy h h'' G_0 \in m \isSemiActedUponBy \mathfrak{H}$. In conclusion, $(m \isSemiActedUponBy \mathfrak{H}) \isSemiActedUponBy \mathfrak{H}' \subseteq m \isSemiActedUponBy \mathfrak{H}$.
  \end{proof} 

  For a $(G_0, G_0)$-invariant subgroup $H$ of $G$ and a semi-cellular automaton $\mathcal{C}$ whose neighbourhood is included in $H \modulo G_0$, the partition of $M$ by the orbit space of $\isSemiActedUponBy\restrictedTo_{M \times H \modulo G_0}$ induces a partition of the phase space $Q^M$ and a partition of the global transition function of $\mathcal{C}$ as shown in

  \begin{lemma} 
  \label{lemma:partition-of-global-transition-function}
    Let $\mathcal{R} = \ntuple{\ntuple{M, G, \actsOnPoint}, \ntuple{m_0, \family{g_{m_0, m}}_{m \in M}}}$ be a cell space, let $\mathcal{C} = \ntuple{\mathcal{R}, Q, N, \delta}$ be a semi-cellular automaton, and let $\Delta$ be the global transition function of $\mathcal{C}$. Furthermore, let $H$ be a $(G_0, G_0)$-invariant subgroup of $G$ such that $N \subseteq H \modulo G_0$, let $\mathfrak{H}$ be the quotient set $H \modulo G_0$, and let $\setOf{m_i \suchThat i \in I}$ be a transversal of $\setOf{m \isSemiActedUponBy \mathfrak{H} \suchThat m \in M}$. Moreover, for each index $i \in I$, let
    \begin{align*}
      \Delta_{m_i \isSemiActedUponBy \mathfrak{H}} \from Q^{m_i \isSemiActedUponBy \mathfrak{H}} &\to Q^{m_i \isSemiActedUponBy \mathfrak{H}},\\
      c_{m_i \isSemiActedUponBy \mathfrak{H}} &\mapsto \brackets[\big]{m_i \isSemiActedUponBy \mathfrak{h} \mapsto \delta\parens[\big]{n \mapsto c_{m_i \isSemiActedUponBy \mathfrak{H}}((m_i \isSemiActedUponBy \mathfrak{h}) \isSemiActedUponBy n)}}.
    \end{align*}
    Then, $M = \bigDisjointUnionOf_{i \in I} m_i \isSemiActedUponBy \mathfrak{H}$, $Q^M = \prod_{i \in I} Q^{m_i \isSemiActedUponBy \mathfrak{H}}$, and $\Delta = \prod_{i \in I} \Delta_{m_i \isSemiActedUponBy \mathfrak{H}}$.
  \end{lemma}

  \begin{proof}
    Because $\isSemiActedUponBy$ is a right group semi-action of $G \modulo G_0$ on $M$ with defect $G_0$ and $H$ is a $(G_0, G_0)$-invariant subgroup of $G$, according to \cref{lemma:partition-of-M-induced-by-semi-action-restriction}, the family $\family{m \isSemiActedUponBy \mathfrak{H}}_{m \in M}$ partitions $M$. Therefore, the transversal $\setOf{m_i \suchThat i \in I}$ is well-defined, $M = \bigDisjointUnionOf_{i \in I} m_i \isSemiActedUponBy \mathfrak{H}$, and $Q^M = \prod_{i \in I} Q^{m_i \isSemiActedUponBy \mathfrak{H}}$. 

    Let $i \in I$. Then, for each $m_i \isSemiActedUponBy \mathfrak{h} \in m_i \isSemiActedUponBy \mathfrak{H}$ and each $n \in N$, because $H$ is $(G_0, G_0)$-invariant and $N \subseteq \mathfrak{H}$, according to \cref{lemma:liberation-after-liberation-of-invariant-subgroup-stays-in-liberation}, we have $(m_i \isSemiActedUponBy \mathfrak{h}) \isSemiActedUponBy n \in m_i \isSemiActedUponBy \mathfrak{H}$. In conclusion, $\Delta_{m_i \isSemiActedUponBy \mathfrak{H}}$ is well-defined.

    Moreover, for each $c_{m_i \isSemiActedUponBy \mathfrak{H}} \in Q^{m_i \isSemiActedUponBy \mathfrak{H}}$ and each extension $c \in Q^M$ of $c_{m_i \isSemiActedUponBy \mathfrak{H}}$, we have $\Delta_{m_i \isSemiActedUponBy \mathfrak{H}}(c_{m_i \isSemiActedUponBy \mathfrak{H}})(m) = \Delta(c)(m)$. In conclusion, $\Delta = \prod_{i \in I} \Delta_{m_i \isSemiActedUponBy \mathfrak{H}}$.
  \end{proof}

  The maps $\Delta_{m_i \isSemiActedUponBy \mathfrak{H}}$, for $i \in I$, are conjugate to each other as shown in

  \begin{lemma}
  \label{lemma:restrictions-of-global-transition-function-are-conjugate-to-each-other}
    In the situation of \cref{lemma:partition-of-global-transition-function}, let $\mathcal{C}$ be a cellular automaton, let $m$ and $m'$ be two elements of $M$, let
    \begin{align*}
      \phi_{m,m'} \from m \isSemiActedUponBy \mathfrak{H} &\to m' \isSemiActedUponBy \mathfrak{H},\\
      m \isSemiActedUponBy \mathfrak{h} &\mapsto m' \isSemiActedUponBy \mathfrak{h},
    \end{align*}
    and let
    \begin{align*}
      \phi_{m,m'}^* \from Q^{m' \isSemiActedUponBy \mathfrak{H}} &\to Q^{m \isSemiActedUponBy \mathfrak{H}},\\
      c_{m' \isSemiActedUponBy \mathfrak{H}} &\mapsto c_{m' \isSemiActedUponBy \mathfrak{H}} \after \phi_{m,m'}.
    \end{align*}
    Then, $\Delta_{m' \isSemiActedUponBy \mathfrak{H}} = (\phi_{m,m'}^*)^{-1} \after \Delta_{m \isSemiActedUponBy \mathfrak{H}} \after \phi_{m,m'}^*$.
  \end{lemma}

  \begin{proof}
    Because $\isSemiActedUponBy$ is free, the map $\phi_{m,m'}$ is well-defined and bijective. Hence, the map $\phi_{m,m'}^*$ is also well-defined and bijective.

    Let $c \in Q^M$ and let $m \isSemiActedUponBy \mathfrak{h} \in m \isSemiActedUponBy \mathfrak{H}$. Then
    \begin{align*}
      (\phi_{m,m'}^* \after \Delta_{m' \isSemiActedUponBy \mathfrak{H}})(c\restrictedTo_{m' \isSemiActedUponBy \mathfrak{H}})(m \isSemiActedUponBy \mathfrak{h})
      &= \phi_{m,m'}^*\parens[\big]{\Delta_{m' \isSemiActedUponBy \mathfrak{H}}(c\restrictedTo_{m' \isSemiActedUponBy \mathfrak{H}})}(m \isSemiActedUponBy \mathfrak{h})\\
      &= \parens[\big]{\Delta_{m' \isSemiActedUponBy \mathfrak{H}}(c\restrictedTo_{m' \isSemiActedUponBy \mathfrak{H}}) \after \phi_{m,m'}}(m \isSemiActedUponBy \mathfrak{h})\\
      &= \Delta_{m' \isSemiActedUponBy \mathfrak{H}}(c\restrictedTo_{m' \isSemiActedUponBy \mathfrak{H}})(m' \isSemiActedUponBy \mathfrak{h})\\
      &= \Delta(c)(m' \isSemiActedUponBy \mathfrak{h}).
    \end{align*}
    And, because $m' \isSemiActedUponBy \mathfrak{h} = g_{m_0, m'} g_{m_0, m}^{-1} \actsOnPoint (m \isSemiActedUponBy \mathfrak{h})$ and $\Delta$ is $\actsOnMap$-e\-qui\-var\-i\-ant, 
    \begin{align*}
      \Delta(c)(m' \isSemiActedUponBy \mathfrak{h})
      &= \Delta(c)\parens[\big]{(g_{m_0, m} g_{m_0, m'}^{-1})^{-1} \actsOnPoint (m \isSemiActedUponBy \mathfrak{h})}\\
      &= \parens[\big]{g_{m_0, m} g_{m_0, m'}^{-1} \actsOnMap \Delta(c)}(m \isSemiActedUponBy \mathfrak{h})\\
      &= \Delta(g_{m_0, m} g_{m_0, m'}^{-1} \actsOnMap c)(m \isSemiActedUponBy \mathfrak{h}).
    \end{align*}
    And, because, for each $m \isSemiActedUponBy \mathfrak{h}' \in m \isSemiActedUponBy \mathfrak{H}$,
    \begin{align*}
      (g_{m_0, m} g_{m_0, m'}^{-1} \actsOnMap c)(m \isSemiActedUponBy \mathfrak{h}')
      &= c\parens[\big]{g_{m_0, m'} g_{m_0, m}^{-1} \actsOnPoint (m \isSemiActedUponBy \mathfrak{h}')}\\
      &= c(m' \isSemiActedUponBy \mathfrak{h}')\\
      &= c\restrictedTo_{m' \isSemiActedUponBy \mathfrak{H}}(m' \isSemiActedUponBy \mathfrak{h}')\\
      &= c\restrictedTo_{m' \isSemiActedUponBy \mathfrak{H}}\parens[\big]{\phi_{m,m'}(m \isSemiActedUponBy \mathfrak{h}')}\\
      &= (c\restrictedTo_{m' \isSemiActedUponBy \mathfrak{H}} \after \phi_{m,m'})(m \isSemiActedUponBy \mathfrak{h}'),
    \end{align*}
    we have $(g_{m_0, m} g_{m_0, m'}^{-1} \actsOnMap c)\restrictedTo_{m \isSemiActedUponBy \mathfrak{H}} = (c \after \phi_{m,m'})\restrictedTo_{m \isSemiActedUponBy \mathfrak{H}}$, and hence
    \begin{align*}
      \Delta(g_{m_0, m} g_{m_0, m'}^{-1} \actsOnMap c)(m \isSemiActedUponBy \mathfrak{h})
      &= \Delta_{m \isSemiActedUponBy \mathfrak{H}}\parens[\big]{(g_{m_0, m} g_{m_0, m'}^{-1} \actsOnMap c)\restrictedTo_{m \isSemiActedUponBy \mathfrak{H}}}(m \isSemiActedUponBy \mathfrak{h})\\
      &= \Delta_{m \isSemiActedUponBy \mathfrak{H}}(c\restrictedTo_{m \isSemiActedUponBy \mathfrak{H}} \after \phi_{m,m'})(m \isSemiActedUponBy \mathfrak{h})\\
      &= \Delta_{m \isSemiActedUponBy \mathfrak{H}}\parens[\big]{\phi_{m,m'}^*(c\restrictedTo_{m \isSemiActedUponBy \mathfrak{H}})}(m \isSemiActedUponBy \mathfrak{h})\\
      &= (\Delta_{m \isSemiActedUponBy \mathfrak{H}} \after \phi_{m,m'}^*)(c\restrictedTo_{m \isSemiActedUponBy \mathfrak{H}})(m \isSemiActedUponBy \mathfrak{h}).
    \end{align*}
    Altogether,
    \begin{equation*}
      (\phi_{m,m'}^* \after \Delta_{m' \isSemiActedUponBy \mathfrak{H}})(c\restrictedTo_{m \isSemiActedUponBy \mathfrak{H}})(m \isSemiActedUponBy \mathfrak{h}) = (\Delta_{m \isSemiActedUponBy \mathfrak{H}} \after \phi_{m,m'}^*)(c\restrictedTo_{m \isSemiActedUponBy \mathfrak{H}})(m \isSemiActedUponBy \mathfrak{h}).
    \end{equation*}
    Therefore, $\phi_{m,m'}^* \after \Delta_{m' \isSemiActedUponBy \mathfrak{H}} = \Delta_{m \isSemiActedUponBy \mathfrak{H}} \after \phi_{m,m'}^*$. In conclusion, $\Delta_{m' \isSemiActedUponBy \mathfrak{H}} = (\phi_{m,m'}^*)^{-1} \after \Delta_{m \isSemiActedUponBy \mathfrak{H}} \after \phi_{m,m'}^*$.
  \end{proof}

  Because of the above decomposition of $\Delta$, the conjugacy of the maps $\Delta_{m_i \isSemiActedUponBy \mathfrak{H}}$, for $i \in I$, and the equality of $\Delta_{m_0 \isSemiActedUponBy \mathfrak{H}}$ to the restriction of $\Delta$ to $H$, properties of the latter restriction translate to properties of $\Delta$ as stated in

  \begin{theorem} 
  \label{theorem:cellular-automaton-injective-if-and-only-if-restriction-injective}
    In the situation of \cref{lemma:partition-of-global-transition-function}, let $\mathcal{C}$ be a cellular automaton, let $\mathcal{C}\restrictedTo_H$ be the restriction of $\mathcal{C}$ to $H$, and let $\Delta_{\restrictedTo_H}$ be the global transition function of $\mathcal{C}\restrictedTo_H$. The global transition function $\Delta$ is injective, surjective, or bijective if and only if the global transition function $\Delta_{\restrictedTo_H}$ is injective, surjective, or bijective respectively.
  \end{theorem}

  \begin{proof}
    According to \cref{lemma:right-semi-action-similar-to-leftaction-in-m0}, we have $m_0 \isSemiActedUponBy \mathfrak{H} = H \actsOnPoint m_0$.
    Thus, $Q^{m_0 \isSemiActedUponBy \mathfrak{H}} = Q^{H \actsOnPoint m_0}$. Hence, for each $c_H \in Q^{m_0 \isSemiActedUponBy \mathfrak{H}} = Q^{H \actsOnPoint m_0}$ and each extension $c \in Q^M$ of $c_H$, according to \cref{lemma:partition-of-global-transition-function} and \cref{lemma:global-transition-function-of-restriction-versus-the-unrestricted-one},
    \begin{equation*}
      \Delta_{m_0 \isSemiActedUponBy \mathfrak{H}}(c_H)
      = \Delta(c)\restrictedTo_{m_0 \isSemiActedUponBy \mathfrak{H}}
      = \Delta(c)\restrictedTo_{H \actsOnPoint m_0}
      = \Delta_{\restrictedTo_H}(c_H).
    \end{equation*}
    Therefore, $\Delta_{m_0 \isSemiActedUponBy \mathfrak{H}} = \Delta_{\restrictedTo_H}$.

    Moreover, according to \cref{lemma:partition-of-global-transition-function}, the map $\Delta$ is the product of $\Delta_{m_i \isSemiActedUponBy \mathfrak{H}}$, for $i \in I$, and thus injective, surjective, or bijective if and only if all $\Delta_{m_i \isSemiActedUponBy \mathfrak{H}}$, for $i \in I$, have the respective property. And, for each $i \in I$, according to \cref{lemma:restrictions-of-global-transition-function-are-conjugate-to-each-other}, the map $\Delta_{m_i \isSemiActedUponBy \mathfrak{H}}$ is conjugate to $\Delta_{m_0 \isSemiActedUponBy \mathfrak{H}}$ and hence injective, surjective, or bijective if and only if $\Delta_{m_0 \isSemiActedUponBy \mathfrak{H}}$ has the respective property. Therefore, $\Delta$ is injective, surjective, or bijective if and only if $\Delta_{m_0 \isSemiActedUponBy \mathfrak{H}}\ (= \Delta_{\restrictedTo_H})$ has the respective property.
  \end{proof}

  \begin{corollary}
  \label{corollary:cellular-automaton-function-injective-if-and-only-if-restriction-injective}
    Let $\mathcal{M} = \ntuple{M, G, \actsOnPoint}$ be a left-ho\-mo\-ge\-neous space, let $m_0$ be an element of $M$, let $H$ be a $(G_0, G_0)$-invariant subgroup of $G$, let $\Delta$ be the global transition function of a cellular automaton over $\mathcal{M}$ with a sufficient neighbourhood $E$ such that $E \subseteq H \modulo G_0$, and let $\Delta_{\restrictedTo_H}$ be the restriction of $\Delta$ at $m_0$ to $H$. The global transition function $\Delta$ is injective, surjective, or bijective if and only if the global transition function $\Delta_{\restrictedTo_H}$ is injective, surjective, or bijective respectively.
  \end{corollary}

  \begin{proof}
    This is a direct consequence of \cref{theorem:cellular-automaton-injective-if-and-only-if-restriction-injective} with the proof of \cref{lemma:restriction-of-global-transition-function}.
  \end{proof}

  \begin{example}[Cylinder] 
    Let $M$ be the infinite circular cylinder $(\R \modulo \Z) \times \R$, let $G$ be the additive group $\R^2$, let $\actsOnPoint$ be the left group action of $G$ on $M$ by $((t_1, t_2), (m_1 + \Z, m_2)) \mapsto ((t_1 + m_1) + \Z, t_2 + m_2)$, let $\mathcal{M}$ be the left-ho\-mo\-ge\-neous space $\ntuple{M, G, \actsOnPoint}$, let $\mathcal{K}$ be the coordinate system $\ntuple{(0 + \Z, 0), \family{(\fractionalPart(m_1), m_2)}_{(m_1 + \Z, m_2) \in M}}$ for $\mathcal{M}$, where $\fractionalPart(m_1)$ denotes the fractional part of $m_1$, let $Q$ be the binary set $\setOf{0, 1}$, let $N$ be the singleton set $\setOf{(0 + \Z, -1)}$, let $\delta$ be the $\bullet$-invariant map $Q^N \to Q$, $\ell \mapsto \ell(0 + \Z, -1)$, let $M$ be identified with $G \modulo G_0$ by $\iota$, let $\mathcal{C}$ be the cellular automaton $\ntuple{\ntuple{\mathcal{M}, \mathcal{K}}, Q, N, \delta}$, and let $H$ be the normal and hence $(G_0, G_0)$-invariant subgroup $\setOf{0} \times \R$ of $G$.

    The global transition function $\Delta$ of $\mathcal{C}$ is the shift map $(0 + \Z, 1) \actsOnMap \blank$ along the axis of $M$, that is, the map $Q^M \to Q^M$, $c \mapsto [(m_1 + \Z, m_2) \mapsto c((m_1 + \Z, m_2 - 1))]$. Under the canonical identification of $\setOf{0 + \Z} \times \R$ with $\R$, its restriction $\Delta_{\restrictedTo_H}$ at $(0 + \Z, 0)$ to $H$ is the shift map $Q^{\R} \to Q^{\R}$, $c_{\restrictedTo_H} \mapsto [m_2 \mapsto c_{\restrictedTo_H}(m_2 - 1)]$. Under the canonical identifications of $\R$ with $\setOf{m_1 + \Z} \times \R$, for $m_1 + \Z \in \R \modulo \Z$, the map $\Delta$ is the product $\prod_{m_1 + \Z \in \R \modulo \Z} \Delta_{\restrictedTo_H}$. As the map $\Delta_{\restrictedTo_H}$ is bijective, so is $\Delta$.
  \end{example}

  \begin{remark}
    Everything that has been done with and said about cellular automata in this chapter could have been done with and said about big-cellular automata under suitable assumptions. We chose not to do so, because the presentation would have been even more cumbersome.
  \end{remark}

  \clearToOddPage
  \chapter{Curtis-Hedlund-Lyndon Theorems}
  \label{chapter:Curtis-Hedlund-Lyndon}

  \paragraph{Abstract.} We prove a topological as well as a uniform variant of the Curtis-Hedlund-Lyndon theorem for big-cellular automata with compact sufficient neighbourhoods over tame left-ho\-mo\-ge\-neous spaces. The latter states that a map on the phase space is the global transition function of such an automaton if and only if it is equivariant and uniformly continuous. It follows from topological theorems that such an automaton is invertible if and only if its global transition function is a uniform isomorphism, which, in a more special setting, is equivalent to being bijective.

  \paragraph{Remark.} Some parts of this chapter appeared in the paper \enquote{\citetitle*{wacker:automata:2016}}\cite{wacker:automata:2016} and they generalise parts of sections~1.2, 1.8, 1.9, and~1.10 of the monograph \enquote{\citetitle*{ceccherini-silberstein:coornaert:2010}}\cite{ceccherini-silberstein:coornaert:2010}.

  \paragraph{Summary.} The \emph{prodiscrete topology} on $Q^M$ has for a base the \emph{cylinders} $\setOf{c \in Q^M \suchThat c\restrictedTo_F = b}$, for $b \in Q^F$ and $F \subseteq M$ finite; and the \emph{prodiscrete uniformity} on $Q^M$ has for a base the \emph{cylinders} $\setOf{(c, c') \in Q^M \times Q^M \suchThat c\restrictedTo_F = c'\restrictedTo_F}$, for $F \subseteq M$ finite. In the case that $Q$ is finite, a topological variant of the Curtis-Hedlund-Lyndon theorem states that a map from $Q^M$ to $Q^M$ is the global transition function of a cellular automaton with finite sufficient neighbourhood if and only if it is equivariant and continuous; and in any case, a uniform variant of the Curtis-Hedlund-Lyndon theorem states that a map from $Q^M$ to $Q^M$ is the global transition function of a cellular automaton with compact sufficient neighbourhood if and only if it is equivariant and uniformly continuous. The finiteness of the sufficient neighbourhood stems from the finiteness of $F$ in the cylinders and the qualifier sufficient is needed because the neighbourhood itself may have to be infinite due to the requirement to be invariant under left multiplication by $G_0$.

  In more detail, let $M$ be equipped with a topology and equip $G \modulo G_0$ with the topology induced by $m_0 \isSemiActedUponBy \blank$. The \emph{uniformity of discrete convergence on compacta} on $Q^M$ has for a base the \emph{cylinders} $\setOf{(c, c') \in Q^M \times Q^M \suchThat c\restrictedTo_K = c'\restrictedTo_K}$, for $K \subseteq M$ compact. If the right semi-action $\isSemiActedUponBy$ maps compacta to sets included in compacta, which is called \emph{semi-tameness}, then the uniform variant of the Curtis-Hedlund-Lyndon theorem holds. It follows that cellular automata with compact sufficient neighbourhoods are invertible if and only if their global transition functions are equivariant uniform isomorphisms, which is, in the case that $Q$ is finite and $M$ carries the discrete topology, equivalent to being equivariant, continuous, and bijective. All these statements also hold for big-cellular automata if we choose corresponding notions of equivariance.

  The Curtis-Hedlund-Lyndon theorem is a famous theorem by Morton Landers Curtis, Gustav Arnold Hedlund, and Roger Conant Lyndon from 1969 and was published in the paper \enquote{\citetitle*{hedlund:1969}}\cite{hedlund:1969}.

  \paragraph{Contents.} In \cref{section:prodiscrete-topology} we prove a topological variant of the Curtis-Hedlund-Lyndon theorem, which characterises global transition functions of big-cellular automata by equivariance and continuity. In \cref{section:tameness} we introduce tameness and semi-tameness of left-ho\-mo\-ge\-neous spaces, which are essential in the proof of the uniform variant of the Curtis-Hedlund-Lyndon theorem. In \cref{section:properness} we introduce properness and semi-properness of left-ho\-mo\-ge\-neous spaces, which are sufficient conditions for semi-tameness. In \cref{section:characterisation} we prove a uniform variant of the Curtis-Hedlund-Lyndon theorem, which characterises global transition functions of big-cellular automata by equivariance and uniform continuity. And in \cref{section:invertibility} we characterise invertibility of big-cellular automata.

  \paragraph{Preliminary Notions.} A \emph{left group set} is a triple $\ntuple{M, G, \actsOnPoint}$, where $M$ is a set, $G$ is a group, and $\actsOnPoint$ is a map from $G \times M$ to $M$, called \emph{left group action of $G$ on $M$}, such that $G \to \symmetricGroupOf(M)$, $g \mapsto [g \actsOnPoint \blank]$, is a group homomorphism. The action $\actsOnPoint$ is \emph{transitive} if $M$ is non-empty and for each $m \in M$ the map $\blank \actsOnPoint m$ is surjective; and \emph{free} if for each $m \in M$ the map $\blank \actsOnPoint m$ is injective. For each $m \in M$, the set $G \actsOnPoint m$ is the \emph{orbit of $m$}, the set $G_m = (\blank \actsOnPoint m)^{-1}(m)$ is the \emph{stabiliser of $m$}, and, for each $m' \in M$, the set $G_{m, m'} = (\blank \actsOnPoint m)^{-1}(m')$ is the \emph{transporter of $m$ to $m'$}.

  A \emph{left-ho\-mo\-ge\-neous space} is a left group set $\mathcal{M} = \ntuple{M, G, \actsOnPoint}$ such that $\actsOnPoint$ is transitive. A \emph{coordinate system for $\mathcal{M}$} is a tuple $\mathcal{K} = \ntuple{m_0, \family{g_{m_0, m}}_{m \in M}}$, where $m_0 \in M$ and, for each $m \in M$, we have $g_{m_0, m} \actsOnPoint m_0 = m$. The stabiliser $G_{m_0}$ is denoted by $G_0$. The tuple $\mathcal{R} = \ntuple{\mathcal{M}, \mathcal{K}}$ is a \emph{cell space}. The set $\setOf{g G_0 \suchThat g \in G}$ of left cosets of $G_0$ in $G$ is denoted by $G \modulo G_0$. The map $\isSemiActedUponBy \from M \times G \modulo G_0 \to M$, $(m, g G_0) \mapsto g_{m_0, m} g \actsOnPoint m_0$ is a \emph{right semi-action of $G \modulo G_0$ on $M$ with defect $G_0$}, which means that
  \begin{equation*}
    \ForEach m \in M \Holds m \isSemiActedUponBy G_0 = m,
  \end{equation*}
  and
  \begin{multline*}
    \ForEach m \in M \ForEach g \in G \Exists g_0 \in G_0 \SuchThat \ForEach \mathfrak{g}' \in G \modulo G_0 \Holds\\
          m \isSemiActedUponBy g \cdot \mathfrak{g}' = (m \isSemiActedUponBy g G_0) \isSemiActedUponBy g_0 \cdot \mathfrak{g}'.
  \end{multline*}
  It is \emph{transitive}, which means that the set $M$ is non-empty and for each $m \in M$ the map $m \isSemiActedUponBy \blank$ is surjective; and \emph{free}, which means that for each $m \in M$ the map $m \isSemiActedUponBy \blank$ is injective; and \emph{semi-commutes with $\actsOnPoint$}, which means that
  \begin{multline*}
    \ForEach m \in M \ForEach g \in G \Exists g_0 \in G_0 \SuchThat \ForEach \mathfrak{g}' \in G \modulo G_0 \Holds\\
          (g \actsOnPoint m) \isSemiActedUponBy \mathfrak{g}' = g \actsOnPoint (m \isSemiActedUponBy g_0 \cdot \mathfrak{g}').
  \end{multline*}

  A \emph{semi-cellular automaton} is a quadruple $\mathcal{C} = \ntuple{\mathcal{R}, Q, N, \delta}$, where $\mathcal{R}$ is a cell space; $Q$, called \emph{set of states}, is a set; $N$, called \emph{neighbourhood}, is a subset of $G \modulo G_0$ such that $G_0 \cdot N \subseteq N$; and $\delta$, called \emph{local transition function}, is a map from $Q^N$ to $Q$. A \emph{local configuration} is a map $\ell \in Q^N$ and a \emph{global configuration} is a map $c \in Q^M$. The stabiliser $G_0$ acts on $Q^N$ on the left by $\bullet \from G_0 \times Q^N \to Q^N$, $(g_0, \ell) \mapsto [n \mapsto \ell(g_0^{-1} \cdot n)]$, and the group $G$ acts on $Q^M$ on the left by $\actsOnMap \from G \times Q^M \to Q^M$, $(g, c) \mapsto [m \mapsto c(g^{-1} \actsOnPoint m)]$. The \emph{global transition function of $\mathcal{C}$} is the map $\Delta \from Q^M \to Q^M$, $c \mapsto [m \mapsto \delta(n \mapsto c(m \isSemiActedUponBy n))]$. A \emph{sufficient neighbourhood of $\mathcal{C}$} is a subset $E$ of $N$ such that, for each $\ell \in Q^N$ and each $\ell' \in Q^N$ with $\ell\restrictedTo_E = \ell'\restrictedTo_E$, we have $\delta(\ell) = \delta(\ell')$.

  A subgroup $H$ of $G$ is \emph{$\mathcal{K}$-big} if the set $\setOf{g_{m_0, m} \suchThat m \in M}$ is included in $H$. A \emph{big-cellular automaton} is a semi-cellular automaton $\mathcal{C} = \ntuple{\mathcal{R}, Q, N, \delta}$ such that, for some $\mathcal{K}$-big subgroup $H$ of $G$, the local transition function $\delta$ is \emph{$\bullet_{G_0 \cap H}$-invariant}, which means that, for each $h_0 \in G_0 \cap H$, we have $\delta(h_0 \bullet \blank) = \delta(\blank)$. Its global transition function is $\actsOnMap_H$-e\-qui\-var\-i\-ant, which means that, for each $h \in H$, we have $\Delta(h \actsOnMap \blank) = h \actsOnMap \Delta(\blank)$. A \emph{cellular automaton} is a big-cellular automaton whose local transition function is $\bullet$-invariant. (See \cref{chapter:automata}.)

  In the present chapter, we assume that the reader is familiar with the basics of the theories of topological and uniform spaces. A recapitulation of the required basics is given in \cref{chapter:topologies,chapter:uniformities}.

  \section{Topological Curtis-Hedlund-Lyndon Theorem}
  \label{section:prodiscrete-topology}

  \paragraph{Contents.} In \cref{definition:prodiscrete-topology,definition:generalised-prodiscrete-topology} we introduce the prodiscrete topology and a generalisation. In \cref{lemma:phase-space-is-Hausdorff,lemma:phase-space-is-compact,lemma:action-on-configurations-is-continuous} we show that the phase space is Hausdorff and compact, and that the left group action on it is continuous. In \cref{lemma:cellular-automata-with-finite-neighbourhoods-are-continuous,lemma:image-of-phase-space-is-closed} we show that the global transition function of each big-cellular automaton with a finite sufficient neighbourhood is continuous and hence has a closed image. And in \cref{theorem:topological-Curtis-Hedlund-Lyndon-direct-proof} we prove a generalised topological variant of the Curtis-Hedlund-Lyndon theorem.

  \paragraph{Body.} The prodiscrete topology is introduced in

  \begin{definition}
  \label{definition:prodiscrete-topology}
    Let $M$ be a set and let $Q$ be a set. Equip $Q$ with the discrete topology and $Q^M = \prod_{m \in M} Q$ with the product topology. This topology on $Q^M$ is called \defineX{prodiscrete}{prodiscrete!topology}\graffito{prodiscrete topology}\index{topology!prodiscrete}.
  \end{definition}

  \begin{remark}
  \label{remark:prodiscrete-topology}
    The prodiscrete topology on $Q^M$ is the coarsest topology such that, for each element $m \in M$, the projection
    \begin{align*}
      \pi_m \from Q^M &\to Q, \mathnote{$\pi_m$, for $m \in M$}\index[symbols]{pim@$\pi_m$}\\
      c &\mapsto c(m),
    \end{align*}
    is continuous. Thus, it has for a subbase the sets
    \begin{equation*}
      \pi_m^{-1}(q) = \setOf{c \in Q^M \suchThat c(m) = q},
      \text{ for $q \in Q$ and $m \in M$}.
    \end{equation*}
    Hence, it has for a base the sets
    \begin{multline*}
      \bigcap_{m \in F} \pi_m^{-1}(b(m)) = \setOf{c \in Q^M \suchThat c\restrictedTo_F = b},\\
      \text{ for $b \in Q^F$ and $F \subseteq M$ finite}.
    \end{multline*}
    Therefore, for each map $c \in Q^M$, the sets, called \define{cylinders},
    \begin{equation*}
      \cylinder(c, F) = \setOf{c' \in Q^M \suchThat c'\restrictedTo_F = c\restrictedTo_F},
      \text{ for $F \subseteq M$ finite},
      \mathnote{cylinders $\cylinder(c, F)$, for $c \in Q^M$ and $F \subseteq M$ finite}
      \index[symbols]{CylcF@$\cylinder(c, F)$}
    \end{equation*}
    constitute a neighbourhood base of $c$.
  \end{remark}

  A generalisation of the prodiscrete topology is introduced in

  \begin{definition} 
  \label{definition:generalised-prodiscrete-topology}
    Let $\actsOnPoint$ be a left group action of $G$ on $M$, let $L$ be a subgroup of $G$, and let $Q$ be a set. The group $L$ acts on $M$ on the left by $\actsOnPoint\restrictedTo_{L \times M}$. Let $\setOf{m_i \suchThat i \in I}$\graffito{$m_i$, for $i \in I$}\index[symbols]{mi@$m_i$} be a transversal of the orbit space of $\actsOnPoint\restrictedTo_{L \times M}$. Then, $M = \bigDisjointUnionOf_{i \in I} L \actsOnPoint m_i$ and $Q^M = \prod_{i \in I} Q^{L \actsOnPoint m_i}$. For each index $i \in I$, equip $Q^{L \actsOnPoint m_i}$ with the discrete topology, and equip $Q^M$ with the product topology. This topology on $Q^M$ is called \defineX{$(\actsOnPoint, L)$-prodiscrete}{prodiscrete topology@$(\actsOnPoint, L)$-prodiscrete topology}\graffito{$(\actsOnPoint, L)$-prodiscrete topology}\index{topology!$(\actsOnPoint, L)$-prodiscrete}.
  \end{definition}

  \begin{remark}
    The $(\actsOnPoint, \setOf{e_G})$-prodiscrete topology is but the prodiscrete topology. 
  \end{remark}

  \begin{remark}
  \label{remark:funny-prodiscrete-topology}
    The $(\actsOnPoint, L)$-prodiscrete topology on $Q^M$ is the coarsest topology such that, 
    for each element $m \in M$, the projection
    \begin{align*}
      \pi_{L \actsOnPoint m} \from Q^M &\to Q^{L \actsOnPoint m}, \mathnote{$\pi_{L \actsOnPoint m}$, for $m \in M$}\index[symbols]{piLarrowrightm@$\pi_{L \actsOnPoint m}$}\\
      c &\mapsto c\restrictedTo_{L \actsOnPoint m},
    \end{align*}
    is continuous. Thus, it has for a subbase the sets
    \begin{equation*}
      \pi_{L \actsOnPoint m}^{-1}(b) = \setOf{c \in Q^M \suchThat c\restrictedTo_{L \actsOnPoint m} = b},
      \text{ for $b \in Q^{L \actsOnPoint m}$ and $m \in M$}.
    \end{equation*}
    Hence, it has for a base the sets
    \begin{multline*}
      \bigcap_{m \in F} \pi_{L \actsOnPoint m}^{-1}(b\restrictedTo_{L \actsOnPoint m})
          = \setOf{c \in Q^M \suchThat c\restrictedTo_{L \actsOnPoint F} = b},\\
      \text{ for $b \in Q^{L \actsOnPoint F}$ and $F \subset M$ finite}.
    \end{multline*}
    Therefore, for each map $c \in Q^M$, the sets, called \define{cylinders},
    \begin{equation*}
      \cylinder(c, L \actsOnPoint F) = \setOf{c' \in Q^M \suchThat c'\restrictedTo_{L \actsOnPoint F} = c\restrictedTo_{L \actsOnPoint F}},
      \text{ for $F \subseteq M$ finite},
      \mathnote{cylinders $\cylinder(c, L \actsOnPoint F)$, for $c \in Q^M$ and $F \subseteq M$ finite}
      \index[symbols]{CylcLarrowrightM0@$\cylinder(c, L \actsOnPoint F)$}
    \end{equation*}
    constitute a neighbourhood base of $c$.
  \end{remark}

  The phase space is Hausdorff and totally disconnected, which is shown in

  \begin{lemma} 
  \label{lemma:phase-space-is-Hausdorff}
    Let $\actsOnPoint$ be a left group action of $G$ on $M$, let $L$ be a subgroup of $G$, and let $Q$ be a set. The set $Q^M$, equipped with the $(\actsOnPoint, L)$-prodiscrete topology, is Hausdorff and totally disconnected.
  \end{lemma}

  \begin{proof}
    Let $\setOf{m_i \suchThat i \in I}$ be a transversal of the orbit space of $\actsOnPoint\restrictedTo_{L \times M}$. Then, for each index $i \in I$, the set $Q^{L \actsOnPoint m_i}$, equipped with the discrete topology, is Hausdorff and totally disconnected. Therefore, according to \cref{lemma:product-of-Hausdorff-is-Hausdorff,lemma:product-of-totally-disconnected-is-disconnected}, the set $Q^M = \prod_{i \in I} Q^{L \actsOnPoint m_i}$, equipped with the product topology, is Hausdorff and totally disconnected.
  \end{proof}

  If the set of states is finite, then the phase space is compact, which is shown in

  \begin{lemma} 
  \label{lemma:phase-space-is-compact}
    Let $\actsOnPoint$ be a left group action of $G$ on $M$, let $L$ be a finite subgroup of $G$, and let $Q$ be a finite set. The set $Q^M$, equipped with the $(\actsOnPoint, L)$-prodiscrete topology, is compact.
  \end{lemma}

  \begin{proof}
    Let $\setOf{m_i \suchThat i \in I}$ be a transversal of the orbit space of $\actsOnPoint\restrictedTo_{L \times M}$. Then, for each index $i \in I$, because $\cardinalityOf{Q^{L \actsOnPoint m_i}} \leq \cardinalityOf{Q}^{\cardinalityOf{L}} < \infty$, the set $Q^{L \actsOnPoint m_i}$ is finite and hence, equipped with the discrete topology, it is compact. Therefore, according to Tychonoff's \cref{theorem:Tychonoff}, the set $Q^M = \prod_{i \in I} Q^{L \actsOnPoint m_i}$, equipped with the product topology, is compact.
  \end{proof}

  The left group action on the phase space is continuous, which is shown in

  \begin{lemma} 
  \label{lemma:action-on-configurations-is-continuous}
    Let $\actsOnPoint$ be a left group action of $G$ on $M$, let $Q$ be a set, and let $Q^M$ be equipped with the prodiscrete topology. The left group action $\actsOnMap$ is \define{continuous}\graffito{continuous left group action}, which means that, for each element $g \in G$, the map $g \actsOnMap \blank$ is continuous.
  \end{lemma}

  \begin{proof}
    Let $g \in G$ and let
    \begin{align*}
      \phi_g \from Q^M &\to Q^M,\\
      c &\mapsto g \actsOnMap c.
    \end{align*}
    Furthermore, let $m \in M$. Then, for each $c \in Q^M$,
    \begin{align*}
      (\pi_m \after \phi_g)(c)
      &= (g \actsOnMap c)(m)\\
      &= c(g^{-1} \actsOnPoint m)\\
      &= \pi_{g^{-1} \actsOnPoint m}(c).
    \end{align*}
    Thus, the map $\pi_m \after \phi_g$ is identical to $\pi_{g^{-1} \actsOnPoint m}$ and is hence continuous. Therefore, according to \cref{lemma:map-to-initial-topology-continuous-if-and-only-if-gens-after-map-continuous}, the map $\phi_g = g \actsOnMap \blank$ is continuous. In conclusion, the action $\actsOnMap$ is continuous.
  \end{proof}

  The global transition function of a big-cellular automaton is continuous, which is shown in

  \begin{lemma}
  \label{lemma:cellular-automata-with-finite-neighbourhoods-are-continuous}
    Let $\mathcal{R} = \ntuple{\mathcal{M}, \mathcal{K}} = \ntuple{\ntuple{M, G, \actsOnPoint}, \ntuple{m_0, \family{g_{m_0, m}}_{m \in M}}}$ be a cell space, let $H$ be a $\mathcal{K}$-big subgroup of $G$, let $\mathcal{C} = \ntuple{\mathcal{R}, Q, N, \delta}$ be a semi-cellular automaton with $\bullet_{H_0}$-invariant local transition function $\delta$ and finite sufficient neighbourhood $E$, let $L$ be a subgroup of $H$, and let $Q^M$ be equipped with the $(\actsOnPoint, L)$-prodiscrete topology. The global transition function $\Delta$ of $\mathcal{C}$ is continuous.
  \end{lemma}

  \begin{proof}
    Let $c \in Q^M$. Furthermore, let $O$ be an open neighbourhood of $\Delta(c)$. Then, there is a finite subset $F$ of $M$ such that $\cylinder(\Delta(c), L \actsOnPoint F) \subseteq O$. And, because $E$ is finite, the set $F \isSemiActedUponBy E$ is finite.

    Let $c' \in \cylinder(c, L \actsOnPoint (F \isSemiActedUponBy E))$, let $f \in F$, and let $l \in L$. Then, because $L \subseteq H$ and $\Delta$ is $\actsOnMap_H$-e\-qui\-var\-i\-ant by \cref{theorem:local-invariance-versus-global-equivariance},
    \begin{align*}
      \Delta(c')(l \actsOnPoint f)
      &= (l^{-1} \actsOnMap \Delta(c'))(f)\\
      &= \Delta(l^{-1} \actsOnMap c')(f)\\
      &= \delta(n \mapsto (l^{-1} \actsOnMap c')(f \isSemiActedUponBy n))\\
      &= \delta(n \mapsto c'(l \actsOnPoint (f \isSemiActedUponBy n)))
    \end{align*}
    and analogously
    \begin{equation*}
      \Delta(c)(l \actsOnPoint f) = \delta(n \mapsto c(l \actsOnPoint (f \isSemiActedUponBy n))).
    \end{equation*}
    And, because $c'$ is in $\cylinder(c, L \actsOnPoint (F \isSemiActedUponBy E))$,
    \begin{equation*}
      \ForEach e \in E \Holds c'(l \actsOnPoint (f \isSemiActedUponBy e)) = c(l \actsOnPoint (f \isSemiActedUponBy e)).
    \end{equation*}
    Hence, because $E$ is a sufficient neighbourhood, we have $\Delta(c')(l \actsOnPoint f) = \Delta(c)(l \actsOnPoint f)$. Therefore, $\Delta(c') \in \cylinder(\Delta(c), L \actsOnPoint F)$. Thus,
    \begin{equation*}
      \Delta(\cylinder(c, L \actsOnPoint (F \isSemiActedUponBy N))) \subseteq \cylinder(\Delta(c), L \actsOnPoint F) \subseteq O.
    \end{equation*}
    In conclusion, the global transition function $\Delta$ is continuous. 
  \end{proof}

  \begin{remark}
  \label{remark:cellular-automata-with-finite-neighbourhoods-are-continuous}
    Let $\mathcal{R} = \ntuple{\ntuple{M, G, \actsOnPoint}, \ntuple{m_0, \family{g_{m_0, m}}_{m \in M}}}$ be a cell space, let $\mathcal{C} = \ntuple{\mathcal{R}, Q, N, \delta}$ be a semi-cellular automaton with finite sufficient neighbourhood, and let $Q^M$ be equipped with the prodiscrete topology. As in the proof of \cref{lemma:cellular-automata-with-finite-neighbourhoods-are-continuous}, one can show that the global transition function $\Delta$ of $\mathcal{C}$ is continuous.
  \end{remark}

  If the set of states is finite, then the image of the global transition function of a big-cellular automaton is closed, which is shown in

  \begin{lemma}
  \label{lemma:image-of-phase-space-is-closed}
    In the situation of \cref{lemma:cellular-automata-with-finite-neighbourhoods-are-continuous}, let $Q$ and $L$ be finite. The set $\Delta(Q^M)$ is closed in $Q^M$.
  \end{lemma}

  \begin{proof}
    According to \cref{lemma:cellular-automata-with-finite-neighbourhoods-are-continuous}, the global transition function $\Delta$ is continuous. And, according to \cref{lemma:phase-space-is-compact}, the phase space $Q^M$ is compact. Hence, the set $\Delta(Q^M)$ is a compact subset of $Q^M$. Moreover, according to \cref{lemma:phase-space-is-Hausdorff}, the phase space $Q^M$ is Hausdorff. Therefore, the compact set $\Delta(Q^M)$ is in particular closed.
  \end{proof}

  \begin{remark}
    In the situation of \cref{remark:cellular-automata-with-finite-neighbourhoods-are-continuous}, as in the proof of \cref{lemma:image-of-phase-space-is-closed}, one can show that $\Delta(Q^M)$ is closed in $Q^M$.
  \end{remark}

  If the set of states is finite, then a map on the phase space is the global transition function of a big-cellular automaton with finite neighbourhood if and only if the map is equivariant and continuous, which is shown in

  \begin{theorem}[Generalised Topological Variant; Morton Landers Curtis, Gustav Arnold Hedlund, and Roger Conant Lyndon, 1969]
  \label{theorem:topological-Curtis-Hedlund-Lyndon-direct-proof}
    Let $\mathcal{R} = \ntuple{\mathcal{M}, \mathcal{K}} = \ntuple{\ntuple{M, G, \actsOnPoint}, \ntuple{m_0, \family{g_{m_0, m}}_{m \in M}}}$ be a cell space, let $Q$ be a finite set, let $\Delta$ be a map from $Q^M$ to $Q^M$, let $H$ be a $\mathcal{K}$-big subgroup of $G$, let $L$ be a finite subgroup of $H$, and let $Q^M$ be equipped with the $(\actsOnPoint, L)$-prodiscrete topology. The following two statements are equivalent:
    \begin{aenumerate}
      \item \label{item:topological-Curtis-Hedlund-Lyndon-direct-proof:global-transition-function}
            The map $\Delta$ is the global transition function of a semi-cellular automaton over $\mathcal{R}$ with $\bullet_{H_0}$-invariant local transition function and finite sufficient neighbourhood.
      \item \label{item:topological-Curtis-Hedlund-Lyndon-direct-proof:equivariant-and-continuous}
            The map $\Delta$ is $\actsOnMap_H$-e\-qui\-var\-i\-ant and continuous. \qedhere
    \end{aenumerate}
  \end{theorem}

  \begin{proof}
    First, let $\Delta$ be the global transition function of a semi-cellular automaton $\mathcal{C} = \ntuple{\mathcal{R}, Q, N, \delta}$ with $\bullet_{H_0}$-invariant local transition function $\delta$ and finite sufficient neighbourhood $E$. Then, according to \cref{theorem:local-invariance-versus-global-equivariance}, the map $\Delta$ is $\actsOnMap_H$-e\-qui\-var\-i\-ant, and, according to \cref{lemma:cellular-automata-with-finite-neighbourhoods-are-continuous}, it is continuous.

    Secondly, let $\Delta$ be $\actsOnMap_H$-e\-qui\-var\-i\-ant and continuous. Then, the map
    \begin{align*}
      \Lambda = \pi_{L \actsOnPoint m_0} \after \Delta \from Q^M &\to Q^{L \actsOnPoint m_0},\\
      c &\mapsto \Delta(c)\restrictedTo_{L \actsOnPoint m_0},
    \end{align*}
    is continuous.

    For a moment, let $c \in Q^M$. Then, because $Q^{L \actsOnPoint m_0}$ is equipped with the discrete topology, the preimage $\Lambda^{-1}(\Lambda(c))$ is open. Hence, there is a finite subset $F_c$ of $M$ such that $\cylinder(c, L \actsOnPoint F_c) \subseteq \Lambda^{-1}(\Lambda(c))$. In other words,
    \begin{equation}
    \label{equation:topological-Curtis-Hedlund-Lyndon-direct-proof:finite}
      \ForEach c' \in \cylinder(c, L \actsOnPoint F_c) \Holds \Lambda(c') = \Lambda(c).
    \end{equation}

    Moreover, because $c \in \cylinder(c, L \actsOnPoint F_c)$, for $c \in Q^M$, the sets $\cylinder(c, L \actsOnPoint F_c)$, for $c \in Q^M$, constitute an open cover of $Q^M$. Hence, because the space $Q^M$ is compact by \cref{lemma:phase-space-is-compact}, there is a finite subset $C$ of $Q^M$ such that
    \begin{equation}
    \label{equation:topological-Curtis-Hedlund-Lyndon-direct-proof:cover}
      Q^M = \bigcup_{c \in C} \cylinder(c, L \actsOnPoint F_c).
    \end{equation}

    Furthermore, because $C$ and $L$ are finite, the set $E_0 = \bigcup_{c \in C} L \actsOnPoint F_c$ is a finite subset of $M$ and the set $E = (m_0 \isSemiActedUponBy \blank)^{-1}(E_0)$ is a finite subset of $G \modulo G_0$. Let $N = G_0 \cdot E$. Then, $G_0 \cdot N \subseteq N$.

    For a while, let $c$ and $c'$ be two global configurations of $Q^M$ such that $c\restrictedTo_{m_0 \isSemiActedUponBy E} = c'\restrictedTo_{m_0 \isSemiActedUponBy E}$. Then, according to \cref{equation:topological-Curtis-Hedlund-Lyndon-direct-proof:cover}, there is an element $c_0 \in C$ such that $c \in \cylinder(c_0, L \actsOnPoint F_{c_0})$. Hence, because $L \actsOnPoint F_{c_0} \subseteq E_0 = m_0 \isSemiActedUponBy E$,
    \begin{equation*}
      c'\restrictedTo_{L \actsOnPoint F_{c_0}} = c\restrictedTo_{L \actsOnPoint F_{c_0}} = c_0\restrictedTo_{L \actsOnPoint F_{c_0}}.
    \end{equation*}
    Thus, $c' \in \cylinder(c_0, L \actsOnPoint F_{c_0})$. Therefore, according to \cref{equation:topological-Curtis-Hedlund-Lyndon-direct-proof:finite}, we have $\Lambda(c) = \Lambda(c_0) = \Lambda(c')$. Hence, because $m_0 \in L \actsOnPoint m_0$,
    \begin{equation*}
      \Delta(c)(m_0) = \Lambda(c)(m_0) = \Lambda(c')(m_0) = \Delta(c')(m_0).
    \end{equation*}
    Therefore, because $E \subseteq N$, there is a map $\delta \from Q^N \to Q$ such that
    \begin{equation*}
      \ForEach c \in Q^M \Holds \Delta(c)(m_0) = \delta(n \mapsto c(m_0 \isSemiActedUponBy n)).
    \end{equation*}
    The quadruple $\mathcal{C} = \ntuple{\mathcal{M}, Q, N, \delta}$ is a semi-cellular automaton with finite sufficient neighbourhood $E$. Conclude with \cref{theorem:determination-of-cellular-automata-by-behaviour-at-origin} that $\delta$ is $\bullet_{H_0}$-invariant and that $\Delta$ is the global transition function of $\mathcal{C}$.
  \end{proof}

  The Curtis-Hedlund-Lyndon theorem presented above is a generalisation of a known theorem for cellular automata over groups.

  \begin{remark}
  \label{remark:generalised-Curtis-Hedlund-Lyndon-on-groups}
    In the case that $L = \setOf{e_G}$, \cref{theorem:topological-Curtis-Hedlund-Lyndon-direct-proof} is essentially theorem~6.7 in \cite{moriceau:2011} and, if additionally $M = G$ and $\actsOnPoint$ is the group multiplication of $G$, then it is theorem~1.8.1 in \cite{ceccherini-silberstein:coornaert:2010}.
  \end{remark}

  \begin{counterexample}[Group {\cite[Example~1.8.2]{ceccherini-silberstein:coornaert:2010}}]
    In theorem \ref{theorem:topological-Curtis-Hedlund-Lyndon-direct-proof}, if the assumption that the set $Q$ of states is finite does not hold, then \cref{item:topological-Curtis-Hedlund-Lyndon-direct-proof:global-transition-function} does not follow from \cref{item:topological-Curtis-Hedlund-Lyndon-direct-proof:equivariant-and-continuous}, which is illustrated by the following example. 

    Let $G$ be an infinite group, let $Q$ be the set $G$, and let $\Delta$ be the map $Q^G \to Q^G$, $c \mapsto [g \mapsto c(g \cdot c(g))]$. According to example~1.8.2 in \cite{ceccherini-silberstein:coornaert:2010}, the map $\Delta$ is $\actsOnMap$-e\-qui\-var\-i\-ant and continuous but not the global transition function of a cellular automaton over $G$ with finite neighbourhood.

    Note that, even if the set $Q$ is infinite, \cref{item:topological-Curtis-Hedlund-Lyndon-direct-proof:equivariant-and-continuous} follows from \cref{item:topological-Curtis-Hedlund-Lyndon-direct-proof:global-transition-function}.
  \end{counterexample}

  \section{Tameness and Semi-Tameness of Spaces}
  \label{section:tameness}

  \paragraph{Introduction.} The symmetries of a circle, namely rotations and (roto-)reflections, act on it on the left by function application. A rotation is uniquely determined by its angle and a (roto-)reflection can be uniquely identified by an angle with respect to a designated line. Hence, the symmetries of the circle can be identified with the disjoint union of the angles, say, $A \times \setOf{1, -1} = \leftClosedAndRightOpenInterval{0, 360} \times \setOf{1, -1}$. We do the identification by $(a, r) \mapsto \rho_a \after \varrho_r$, where $\rho_a$ is the rotation by $a$ degrees, $\varrho_1$ is the identity map, and $\varrho_{-1}$ is the reflection about the vertical line through the centre of the circle. This identification induces the group structure $(a, r) + (a', r') \mapsto (a + r \cdot a' \bmod 360, r \cdot r')$ on $A \times \setOf{1, -1}$ and the left group action $(a, r) \actsOnPoint m = \rho_a(\varrho_r(m))$ of $A \times \setOf{1, -1}$ on the circle. 

  Geometrically, the group $A \times \setOf{1, -1}$ is made up of two circles of circumference $360$. The geometry on each circle induces a topology on it and both topologies together induce the product topology on $A \times \setOf{1, -1}$. It can be shown that addition and inversion in $A \times \setOf{1, -1}$ are continuous, and that the left group action $\actsOnPoint$ is continuous also. After all, the action just rotates and reflects points, so points that are close stay close. 

  As we have seen in the introduction of \cref{section:semi-action}, the right quotient set semi-action induced by $\actsOnPoint$ can be identified with the right group action of the rotations on the circle. Under the identification of rotations with angles, this action is $m \isSemiActedUponBy a = \rho_a(m)$, which is continuous. Hence, if we act with a compact subset of angles on a compact subset of points by that action, then we get a compact subset of points. This is not the case for all right semi-actions, but it is for those induced by so called \emph{tame} left group actions. This property is of interest for cellular automata, because it implies that, if an automaton has a compact (relative) neighbourhood, then the actual neighbourhood of each cell is compact, even the union of all actual neighbourhoods of a compact subset of points is compact. For our purposes, being included in a compact set is sufficient, which is already implied by \emph{semi-tameness}. 

  \paragraph{Contents.} In \cref{definition:topology-on-quotient-by-stabiliser} we equip $G \modulo G_0$ with the topology of $M$. In \cref{definition:tame-and-semi-tame,definition:tame-and-semi-tame-group-set} we introduce tameness and semi-tameness. In \cref{remark:discrete-tame-cell-space} we note that the action of a group on a discrete space is tame and semi-tame. In \cref{example:affine-reals} we illustrate that (semi-)tameness of cell spaces depends on the coordinate system. In \cref{definition:topological-group,definition:topological-group-set} we introduce topological groups and group sets. In \cref{lemma:continuous-coordinate-map-implies-tameness} we show that a cell space with a continuous coordinate map is tame. And in \cref{example:affine-reals:coordinate-maps} we illustrate that continuity of coordinate maps is not necessary for semi-tameness.

  \paragraph{Body.} We equip $G \modulo G_0$ with the topology of $M$ in

  \begin{definition}
  \label{definition:topology-on-quotient-by-stabiliser}
    Let $\mathcal{R} = \ntuple{\ntuple{M, G, \actsOnPoint}, \ntuple{m_0, \family{g_{m_0, m}}_{m \in M}}}$ be a cell space and let $M$ be equipped with a topology. Equip $G \modulo G_0$ with the topology induced by $m_0 \isSemiActedUponBy \blank$\graffito{topology on $G \modulo G_0$}. 
  \end{definition}

  A cell space is (semi-)tame if its right semi-action maps compacta to (sets included in) compacta. This property is essential in the proof of \cref{theorem:uniform-Curtis-Hedlund-Lyndon}; there it is used to show that global transition functions of big-cellular automata are uniformly continuous.

  \begin{definition}
  \label{definition:tame-and-semi-tame}
    Let $M$ be a topological space and let $\mathcal{R} = \ntuple{\ntuple{M, G, \actsOnPoint}, \ntuple{m_0, \family{g_{m_0, m}}_{m \in M}}}$ be a cell space. The cell space $\mathcal{R}$ is called
    \begin{aenumerate}
      \item \define{tame}\graffito{tame cell space}\index{cell space!tame}\index{tame!cell space} if and only if, for each compact subset $K$ of $M$ and each compact subset $E$ of $G \modulo G_0$, the set $K \isSemiActedUponBy E$ is a compact subset of $M$;
      \item \define{semi-tame}\graffito{semi-tame cell space}\index{tame!semi-}\index{cell space!semi-tame}\index{semi-tame!cell space} if and only if, for each compact subset $K$ of $M$ and each compact subset $E$ of $G \modulo G_0$, the set $K \isSemiActedUponBy E$ is included in a compact subset of $M$. \qedhere
    \end{aenumerate}
  \end{definition}

  A left-ho\-mo\-ge\-neous space is (semi-)tame if all its right semi-actions map compacta to (sets included in) compacta.

  \begin{definition}
  \label{definition:tame-and-semi-tame-group-set}
    Let $\mathcal{M}$ be a left-ho\-mo\-ge\-neous space. It is called \define{tame}\graffito{tame/semi-tame left-ho\-mo\-ge\-neous space}\index{left-ho\-mo\-ge\-neous space!tame}\index{homogeneous space!tame}\index{tame!left-ho\-mo\-ge\-neous space} or \define{semi-tame}\index{tame!semi-}\index{left-ho\-mo\-ge\-neous space!semi-tame}\index{homogeneous space!semi-tame}\index{semi-tame!left-ho\-mo\-ge\-neous space} if and only if, for each coordinate system $\mathcal{K}$ for $\mathcal{M}$, the cell space $\ntuple{\mathcal{M}, \mathcal{K}}$ is tame or semi-tame respectively.
  \end{definition}

  \begin{remark}
    Each tame cell space and each tame left-ho\-mo\-ge\-neous space is semi-tame.
  \end{remark}

  \begin{remark}
  \label{remark:discrete-tame-cell-space}
    If $M$ is equipped with the discrete topology, then the left-ho\-mo\-ge\-neous space $\mathcal{M}$ is tame.
  \end{remark}

  Tameness and semi-tameness of cell spaces depend on the coordinate system, which is illustrated by

  \begin{example}[Affine Space]
  \label{example:affine-reals}
    Let $M$ be the one-dimensional Euclidean space $\R$, let $G$ be the group $\setOf{\alpha_{t, d} \from \R \to \R, r \mapsto t + d \cdot r \suchThat t \in \R \text{ and } d \in \R \smallsetminus \setOf{0}}$ of invertible affine transformations of $M$, which is called \define{affine group of $M$}\graffito{affine group of $M$}, and let $\actsOnPoint$ be the left group action of $G$ on $M$ by function application. The triple $\mathcal{M} = \ntuple{M, G, \actsOnPoint}$ is a left-ho\-mo\-ge\-neous space and the stabiliser $G_0$ of the origin $0$ under $\actsOnPoint$ is the group $\setOf{\alpha_{0, d} \suchThat d \in \R \smallsetminus \setOf{0}}$ of invertible linear transformations of $M$. Identify $G \modulo G_0$ with $M$ by $\iota \givenBy m \mapsto G_{0, m}$. 
    \begin{aenumerate}
      \item The tuple $\mathcal{K} = \ntuple{0, \family{\alpha_{m, 1}}_{m \in M}}$ is a coordinate system for $\mathcal{M}$ such that $\isSemiActedUponBy = +$. Hence, because the map $+$ is continuous, for each compact subset $K$ of $M$ and each compact subset $E$ of $G \modulo G_0 \simeq M$, the set $K \isSemiActedUponBy E = K + E$ is compact. Therefore, the cell space $\ntuple{\mathcal{M}, \mathcal{K}}$ is tame. 
      \item The tuple $\mathcal{K}' = \ntuple{0, \family{\alpha_{m, 1}}_{m \in M \smallsetminus \setOf{-1, 1}} \times \family{\alpha_{m, 2}}_{m \in \setOf{-1, 1}}}$ is a coordinate system for $\mathcal{M}$ such that
            \begin{equation*}
              \isSemiActedUponBy \from (m, \mathfrak{g}) \mapsto \begin{dcases*}
                m + \mathfrak{g}, &if $m \notin \setOf{-1, 1}$,\\
                m + 2 \mathfrak{g}, &if $m \in \setOf{-1, 1}$.
              \end{dcases*}
            \end{equation*}
            Hence, for each compact subset $K$ of $M$ and each compact subset $E$ of $G \modulo G_0 \simeq M$, the set $K \isSemiActedUponBy E$ is included in the compact set $K + 2 E$ but it need not be compact itself. For example, although the bounded and closed interval $\closedInterval{0, 1}$ and the singleton set $\setOf{1}$ are compact, the set $\closedInterval{0, 1} \isSemiActedUponBy \setOf{1} = \leftClosedAndRightOpenInterval{1, 2} \cup \setOf{3}$ is bounded but not compact. Therefore, the cell space $\ntuple{\mathcal{M}, \mathcal{K}'}$ is semi-tame but not tame. 
      \item\label{item:affine-reals:unbounded}
            The tuple $\mathcal{K}'' = \ntuple{0, \family{\alpha_{m, 1/m}}_{m \in M \smallsetminus \setOf{0}} \times \family{\alpha_{m, 1}}_{m \in \setOf{0}}}$ is a coordinate system for $\mathcal{M}$ such that
            \begin{equation*}
              \isSemiActedUponBy \from (m, \mathfrak{g}) \mapsto \begin{dcases*}
                m + \frac{1}{m} \cdot \mathfrak{g}, &if $m \neq 0$,\\
                m + \mathfrak{g}, &if $m = 0$.
              \end{dcases*}
            \end{equation*}
            The map $f \from M \smallsetminus \setOf{0} \to M$, $m \mapsto m + 1/m$, is continuous and thus it maps intervals to intervals, in particular, because $f(1) = 2$ and $\lim_{m \downarrow 0} = \infty$, we have $f(\leftOpenAndRightClosedInterval{0, 1}) = \leftClosedAndRightOpenInterval{2, \infty}$.
            Hence, although the sets $\closedInterval{0, 1}$ and $\setOf{1}$ are compact, the set $\closedInterval{0, 1} \isSemiActedUponBy \setOf{1} = \setOf{1} \cup f(\leftOpenAndRightClosedInterval{0, 1}) = \setOf{1} \cup \leftClosedAndRightOpenInterval{2, \infty}$ is unbounded, in particular, it is not included in a compact set. Therefore, the cell space $\ntuple{\mathcal{M}, \mathcal{K}''}$ is not semi-tame.
    \end{aenumerate}
    There are similar examples for the $d$-dimensional Euclidean space $\R^d$ acted on by the affine group $\R^d \rtimes \generalLinearGroup(d, \R)$ and for the $d$-dimensional integer lattice $\Z^d$ acted on by the symmetric group $\symmetricGroupOf(\Z^d)$ of bijective maps on $\Z^d$, for $d \in \N_+$.
  \end{example}

  A topological group is a group equipped with a topology such that multiplication and inversion are continuous.

  \begin{definition}
  \label{definition:topological-group}
    Let $G$ be a group equipped with a topology. It is called \define{topological}\graffito{topological group}\index{topological!group}\index{group topological@topological group} if and only if the maps
    \begin{equation*}
      \left\{
      \begin{aligned}
        G \times G &\to G,\\
        (g, g') &\mapsto g g',
      \end{aligned}
      \right\}
      \text{ and }
      \left\{
      \begin{aligned}
        G &\to G,\\
        g &\mapsto g^{-1},
      \end{aligned}
      \right\}
    \end{equation*}
    are continuous, where $G \times G$ is equipped with the product topology.
  \end{definition}

  A group set is topological if its group action is continuous.

  \begin{definition}
  \label{definition:topological-group-set}
    Let $M$ be a topological space, let $G$ be a topological group, let $\mathcal{M} = \ntuple{M, G, \actsOnPoint}$ be a left group set, and let $G \times M$ be equipped with the product topology. The group set $\mathcal{M}$ is called \define{topological}\graffito{topological left group set}\index{topological!left group set}\index{left group set!topological}\index{group set!topological} if and only if the map $\actsOnPoint$ is continuous.
  \end{definition}

  Equipping a bare left group set with the discrete topology yields a topological group set.

  \begin{lemma}
  \label{lemma:discrete-topological-cell-space}
    Let $\mathcal{M} = \ntuple{M, G, \actsOnPoint}$ be a left group set. Equip $M$ and $G$ with their respective discrete topology. The group set $\mathcal{M}$ is topological.
  \end{lemma}

  \begin{proof}
    The product topology on $G \times M$ and the one on $M \times M$ are discrete. Hence, each subset of $G \times M$ and each of $M \times M$ is open. Thus, the map $\actsOnPoint$ is continuous. In conclusion, the group set $\mathcal{M}$ is topological.
  \end{proof}

  A topological cell space is tame if its coordinate map is continuous. 

  \begin{lemma}
  \label{lemma:continuous-coordinate-map-implies-tameness}
    Let $\mathcal{R} = \ntuple{\ntuple{M, G, \actsOnPoint}, \ntuple{m_0, \family{g_{m_0, m}}_{m \in M}}}$ be a topological cell space such that the \define{coordinate map}\graffito{coordinate map $\chi$}\index[symbols]{chi@$\chi$} $\chi \from M \to G$, $m \mapsto g_{m_0, m}$, is continuous. The map $\isSemiActedUponBy$ is continuous and the space $\mathcal{R}$ is tame.
  \end{lemma}

  \begin{proof}
    According to the definition of the product topology on $M \times G \modulo G_0$, the projections $\pi_1 \from M \times G \modulo G_0 \to M$, $(m, \mathfrak{g}) \mapsto m$ and $\pi_2 \from M \times G \modulo G_0 \to M$, $(m, \mathfrak{g}) \mapsto \mathfrak{g}$, are continuous. And, according to the universal property of the product topology on $M \times M$, because the maps $\chi \after \pi_1$ and $\pi_2$ are continuous, the map $\psi \from M \times G \modulo G_0 \to M \times M$, $(m, \mathfrak{g}) \mapsto (g_{m_0, m}, \mathfrak{g}) = ((\chi \after \pi_1)(m, \mathfrak{g}), \pi_2(m, \mathfrak{g}))$ is continuous. Hence, under the identification of $G \modulo G_0$ with $M$ by $\iota$, because $\actsOnPoint$ is continuous, the map $\isSemiActedUponBy = \actsOnPoint \after \psi$ is continuous. Therefore, the map $\isSemiActedUponBy$ maps compacta to compacta, in particular, for each compact subset $K$ of $M$ and each compact subset $E$ of $G \modulo G_0$, the set $K \isSemiActedUponBy E$ is compact. In conclusion, the space $\mathcal{R}$ is tame. 
  \end{proof}

  Continuity of coordinate maps is not necessary for semi-tameness, which is illustrated by

  \begin{example}[Affine Space]
  \label{example:affine-reals:coordinate-maps}
    In the situation of \cref{example:affine-reals}, the set $\R \smallsetminus \setOf{0}$ is a group under multiplication and the set $\R$ is a group under addition. Let $\varphi$ be the group homomorphism $\R \smallsetminus \setOf{0} \to \automorphismsOf(\R)$, $d \mapsto [t \mapsto d \cdot t]$, let $\R \rtimes_\varphi (\R \smallsetminus \setOf{0})$ be the outer semi-direct product of $\R \smallsetminus \setOf{0}$ acting on $\R$ by $\varphi$, and let $\R \rtimes_\varphi (\R \smallsetminus \setOf{0})$ be equipped with the subspace metric and subspace topology induced by $\R \times \R$.

    The group multiplication of $\R \rtimes_\varphi (\R \smallsetminus \setOf{0})$ is the continuous map $((t, d), (t', d')) \mapsto (t + d \cdot t', d \cdot d')$ and the group inversion of $\R \rtimes_\varphi (\R \smallsetminus \setOf{0})$ is the continuous map $(t, d) \mapsto (-t/d, 1/d)$. Hence, the group $\R \rtimes_\varphi (\R \smallsetminus \setOf{0})$ is topological. Moreover, the map $\R \rtimes_\varphi (\R \smallsetminus \setOf{0}) \to G$, $(t, d) \mapsto \alpha_{t, d}$, is a group isomorphism and, under the identification of $\R \rtimes_\varphi (\R \smallsetminus \setOf{0})$ with the group $G$ by that isomorphism, which we shall henceforth do, the left group action $\actsOnPoint$ is the continuous map $((t, d), r) \mapsto t + d \cdot r$. Therefore, the left group set $\mathcal{M}$ is topological.
    \begin{aenumerate}
      \item The coordinate map $\chi \from m \mapsto (m, 1)$ induced by $\mathcal{K}$ is continuous. 
      \item The coordinate map
            \begin{equation*}
              \chi' \from m \mapsto \begin{dcases*}
                (m, 1), &if $m \notin \setOf{-1, 1}$,\\
                (m, 2), &if $m \in \setOf{-1, 1}$,
              \end{dcases*}
            \end{equation*}
            induced by $\mathcal{K}'$ is not continuous. However, as we have already seen, the cell space $\ntuple{\mathcal{M}, \mathcal{K}'}$ is semi-tame. Hence, continuity of the coordinate map is sufficient but not necessary for semi-tameness. 
      \item The coordinate map
            \begin{equation*}
              \chi'' \from m \mapsto \begin{dcases*}
                (m, \frac{1}{m}), &if $m \neq 0$,\\
                (m, 1), &if $m = 0$,
              \end{dcases*}
            \end{equation*}
            induced by $\mathcal{K}''$ is not continuous. \qedhere 
    \end{aenumerate}
  \end{example}

  \begin{open-problem}
    Is there a tame cell space whose coordinate map is not continuous?
  \end{open-problem}

  \begin{open-problem}
    Is there a semi-tame topological left-ho\-mo\-ge\-neous space for which there is no coordinate system such that the coordinate map is continuous? I conjecture that the sphere acted upon by rotations is such a space.
  \end{open-problem} 

  \section{Properness and Semi-Properness of Left Homogeneous Spaces}
  \label{section:properness}

  \paragraph{Introduction.} In the situation of the introduction of \cref{section:tameness}, the action map of $\actsOnPoint$ is $\alpha \from (A \times \setOf{1, -1}) \times M \to M \times M$, $((a, r), m) \mapsto ((a, r) \actsOnPoint m, m)$, where $M$ denotes the circle. Its domain is, topologically, the disjoint union of two tori equipped with the product topology; and its codomain is, topologically, a torus. The preimage of $(m', m)$ under $\alpha$ consists of two tuples, one in each torus, namely $((a, 1), m)$ and $((a', -1), m)$, where $a$ is the angle such that $\rho_a(m) = m'$ and $a'$ is the angle such that $\rho_{a'}(\varrho_{-1}(m)) = m'$ (see \cref{figure:tori:points}); the preimage of $Y \times \setOf{m}$, where $Y$ is an arc of the circle, consists of two \enquote{similar} paths in the poloidal or toroidal direction, one in each torus, of the form $(B \times \setOf{1}) \times \setOf{m}$ and $(B' \times \setOf{-1}) \times \setOf{m}$, where $B$ and $B'$ are arcs of $A$ of the same length (see \cref{figure:tori:arc-Y}); the preimage of $\setOf{m'} \times X$, where $X$ is an arc of the circle, consists of two \enquote{similar} paths in the \enquote{diagonal} direction, one in each torus, of the form $\setOf{((a_m, 1), m) \suchThat m \in X}$ and $\setOf{((a_m', -1), m) \suchThat m \in X}$, where $a_m$ is the angle such that $\rho_{a_m}(m) = m'$ and $a_m'$ is the angle such that $\rho_{a_m'}(\varrho_{-1}(m)) = m'$ (see \cref{figure:tori:arc-X}); and, the preimage of $Y \times X$, where $X$ and $Y$ are arcs of the circle, consist of two \enquote{similar} ribbons in the \enquote{diagonal} direction, one in each torus (see \cref{figure:tori:ribbons}). It can be shown that the preimages of compact subsets of the torus under the action map are compact. 

  A left group action with such an action map is called \emph{proper}. And, if each transversal of the preimages of the elements of compact subsets are only included in compact subsets of $M$, it is called \emph{semi-proper}. For each coordinate system, properness implies \emph{almost} tameness and semi-properness implies semi-tameness. 
  \begin{figure}
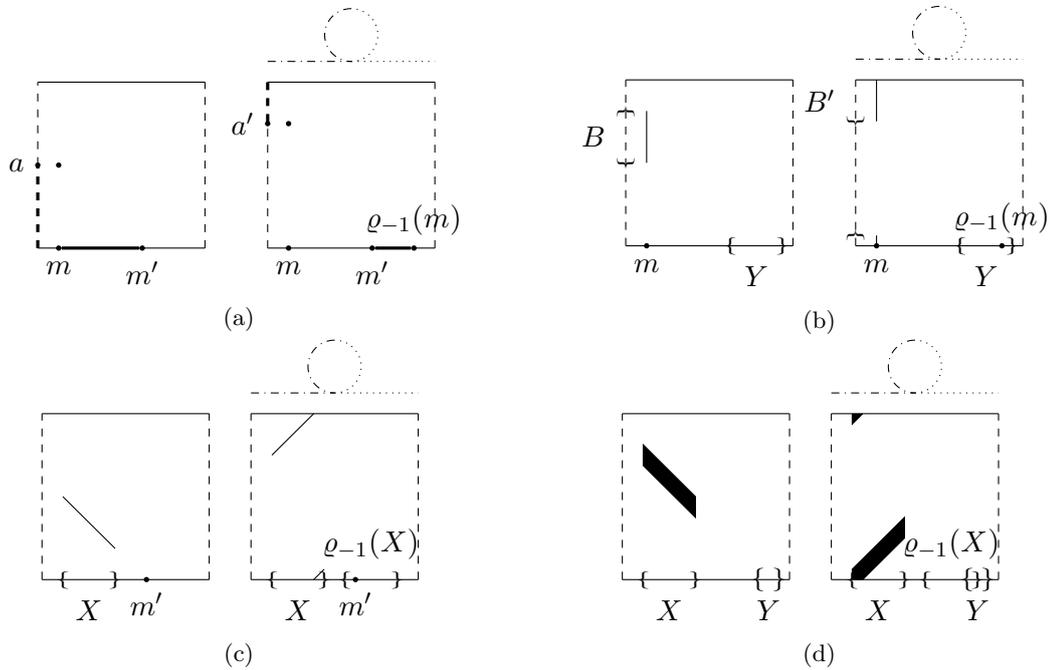

    \myfloatalign
    \begin{wide}
      \figureToriPointsArcsAndRibbons
    \end{wide}
    \caption{%
      In each subfigure, the left square depicts the torus $A \times \setOf{1} \times M$ and the right square the torus $A \times \setOf{-1} \times M$, where, in each square, the solid edges are identified and the dashed edges are identified. The angles $A$ are plotted on the vertical axis from bottom to top beginning with $0$ and proceeding in the positive direction; and the points $M$ are plotted on the horizontal axis from left to right beginning with the top most point of $M$ and proceeding anticlockwise. The reflection of $M$ about the vertical axis is in the depictions the reflections about the midpoints of the solid horizontal lines. Elements are depicted as dots, and sets as lines or ribbons that are enclosed in curly braces if necessary.
      In \cref{figure:tori:points} the two points in the squares depict the preimage of $(m', m)$ under the action map $\alpha$ of $\actsOnPoint$. In \cref{figure:tori:arc-Y} the vertical line segments in the squares depict the preimage of $Y \times \setOf{m}$ under $\alpha$. In \cref{figure:tori:arc-X} the diagonal line segments in the squares depict the preimage of $\setOf{m'} \times X$ under $\alpha$. And in \cref{figure:tori:ribbons} the diagonal ribbons in the squares depict the preimage of $Y \times X$ under $\alpha$.
    }
    \label{figure:tori:points-arcs-and-ribbons}
  \end{figure}


  \paragraph{Contents.} In \cref{definition:transversal} we introduce transversals of sets of sets. In \cref{definition:proper-map,definition:proper-action} we introduce properness and semi-properness of maps and actions. In \cref{exp:proper,exp:semi-proper} we present examples of proper and semi-proper actions. In \cref{lemma:proper-compact-subset-of-group} we show that for semi-proper actions each transversal of the transporters from a compact set to a compact set is included in a compact set. In \cref{lemma:discrete-semi-proper-cell-space} we show that the action of a discrete group on a discrete space is semi-proper. In \cref{lemma:semi-proper-implies-semi-tame} we show that each transitive semi-proper action is semi-tame for each coordinate system. In \cref{example:not-semi-proper-but-tame-for-each-coordinate-system} we present an example that demonstrates that for the previous property semi-properness is not necessary. And in \cref{example:not-semi-proper-but-tame-for-some-coordinate-systems} we present an action that is not semi-proper but semi-tame for some but not all coordinate systems.

  \paragraph{Body.} A transversal of a set of sets is a set that contains one element from each set.

  \begin{definition} 
  \label{definition:transversal}
    Let $M$ be a set, let $\mathfrak{L}$ be a subset of the power set of $M$, and let $T$ be a subset of $M$. The set $T$ is called \define{transversal of $\mathfrak{L}$}\graffito{transversal $T$ of $\mathfrak{L}$} if and only if there is a surjective map $f \from \mathfrak{L} \to T$ such that for each set $A \in \mathfrak{L}$ we have $f(A) \in A$.
  \end{definition}

  The image of a compact set under a continuous map is compact, but the preimage is in general not. A proper map is a continuous map whose preimages of compact sets are compact. And, a semi-proper map is one whose transversals of its fibres are included in compact sets.

  \begin{definition}
  \label{definition:proper-map}
    Let $M$ and $M'$ be two topological spaces and let $f$ be a continuous map from $M$ to $M'$. The map $f$ is called
    \begin{aenumerate}
      \item \define{proper}\graffito{proper map}\index{proper!map}\index{map!proper} if and only if, for each compact subset $K$ of $M'$, its preimage $f^{-1}(K)$ is a compact subset of $M$; 
      \item \define{semi-proper}\graffito{semi-proper map}\index{proper!semi-}\index{semi-proper!map}\index{map!semi-proper} if and only if, for each compact subset $K$ of $M'$, each transversal of $\setOf{f^{-1}(k) \suchThat k \in K}$ is included in a compact subset of $M$. \qedhere
    \end{aenumerate}
  \end{definition}

  A topological group set is proper if its action map is proper, and semi-proper if its action map is semi-proper. If a group action is proper as a map, then its action map is proper; the converse though is false (see problem~21-1 in \cite{lee:2013}).

  \begin{definition}
  \label{definition:proper-action}
    Let $\mathcal{M} = \ntuple{M, G, \actsOnPoint}$ be a topological left group set, and let $G \times M$ and $M \times M$ be equipped with their respective product topology. The group set $\mathcal{M}$ is called
    \begin{aenumerate}
      \item \define{proper}\graffito{proper left group set}\index{proper!left group set}\index{left group set!proper}\index{group set!proper} if and only if its so-called \define{action map}
            \begin{align*}
              \alpha \from G \times M &\to M \times M, \mathnote{action map $\alpha$}\index{map!action}\index[symbols]{alpha@$\alpha$}\\
              (g, m) &\mapsto (g \actsOnPoint m, m),
            \end{align*}
            is proper;
      \item \define{semi-proper}\graffito{semi-proper left group set}\index{proper!semi-}\index{semi-proper!left group set}\index{left group set!semi-proper}\index{group set!semi-proper} if and only if its action map is semi-proper. \qedhere
    \end{aenumerate}
  \end{definition}

  \begin{remark}
  \label{remark:proper-implies-semi-proper}
    Each proper map is semi-proper and each proper left group set is semi-proper.
  \end{remark}

  There are numerous examples of proper and of semi-proper group sets. The Curtis-Hedlund-Lyndon theorem though is mostly interesting for non-compact ones. Therefore, we are especially interested in non-compact and proper or semi-proper homogeneous spaces.

  \begin{example}[Proper] 
  \label{exp:proper}
    \begin{aenumerate}
      \item Each left group set $\ntuple{M, G, \actsOnPoint}$ with finite stabilisers, where $M$ and $G$ are equipped with their respective discrete topology, is proper.
      \item\label{item:proper:proper-if-and-only-if-stabs-compact} 
            Let $M$ be a Hausdorff topological space, let $G$ be a Lie group, and let $\actsOnPoint$ be a transitive and continuous left group action of $G$ on $M$. According to example~2 in paragraph~1.1 in chapter~3 in \cite{gorbatsevich:2013}, the left-ho\-mo\-ge\-neous space $\ntuple{M, G, \actsOnPoint}$ is proper if and only if its stabilisers are compact.
      \item Let $M$ be a Riemannian manifold, let $G$ be the isometry group of $M$, and let $\actsOnPoint$ be the left group action of $G$ on $M$ by function application. According to example~3 in paragraph~1.1 in chapter~3 in \cite{gorbatsevich:2013}, the left group set $\ntuple{M, G, \actsOnPoint}$ is proper. 
      \item Let $M$ be a manifold, let $G$ be a compact Lie group, and let $\actsOnPoint$ be a continuous left group action of $G$ on $M$. According to corollary~21.6 in \cite{lee:2013}, the left group set $\ntuple{M, G, \actsOnPoint}$ is proper.
      \item
            Let $G$ be a Lie group, let $H$ be a compact subgroup of $G$, and let $\actsOnPoint$ be the transitive left group action of $G$ on $G \modulo H$ by left multiplication. According to theorem~21.17 in \cite{lee:2013}, the quotient space $G \modulo H$ is a smooth manifold and the left group action $\actsOnPoint$ is smooth. Because $H$ is compact, all stabilisers of $\actsOnPoint$ are compact and hence, according to \cref{item:proper:proper-if-and-only-if-stabs-compact}, the left-ho\-mo\-ge\-neous space $\ntuple{G \modulo H, G, \actsOnPoint}$ is proper. For example: 
            \begin{aenumerate}
              \item Let $n$ be a positive integer. The general linear group $\generalLinearGroup(n, \R)$ is a non-compact Lie group, the orthogonal group $\orthogonalGroup(n, \R)$ is a maximal compact subgroup of $\generalLinearGroup(n, \R)$, and the special orthogonal group $\specialOrthogonalGroup(n, \R)$ is a compact subgroup of $\generalLinearGroup(n, \R)$. Hence, the left-ho\-mo\-ge\-neous spaces $\ntuple{\generalLinearGroup(n, \R) \modulo \orthogonalGroup(n, \R), \generalLinearGroup(n, \R), \cdot}$ and $\ntuple{\generalLinearGroup(n, \R) \modulo \specialOrthogonalGroup(n, \R), \generalLinearGroup(n, \R), \cdot}$ are proper. 
              \item 
                    The special linear group $\specialLinearGroup(2, \R)$ is a non-compact Lie group and the special orthogonal group $\specialOrthogonalGroup(2, \R)$ is a maximal compact subgroup of $\specialLinearGroup(2, \R)$. Hence, the left-ho\-mo\-ge\-neous space $\ntuple{\specialLinearGroup(2, \R) \modulo \specialOrthogonalGroup(2, \R), \specialLinearGroup(2, \R), \cdot}$ is proper.

                    The upper half-plane $\upperHalfPlane = \setOf{x + \imaginaryUnit y \in \C \suchThat x \in \R, y \in \R_{> 0}}$ equipped with the Poincaré metric is the Poincaré half-plane, a model of the hyperbolic plane. The isometries of the Poincaré half-plane are the Möbius transformations
                    \begin{align*}
                      \mu_{a, b, c, d} \from \upperHalfPlane &\to \upperHalfPlane,\\
                      z &\mapsto \frac{a z + b}{c z + d},
                    \end{align*} 
                    where $a$, $b$, $c$, and $d$ are four real numbers satisfying $a d - b c \in \setOf{-1, 1}$. Note that because $\mu_{a, b, c, d} = \mu_{-a, -b, -c, -d}$, this naming of Möbius transformations is not unique.

                    The group $\specialLinearGroup(2, \R)$ acts transitively and smoothly, but not faithfully, on the Poincaré half-plane by
                    \begin{align*} 
                      \actsOnPoint \from \specialLinearGroup(2, \R) \times \upperHalfPlane &\to \upperHalfPlane,\\
                      (\left(\begin{matrix}
                        a & b\\
                        c & d
                      \end{matrix}\right), z) &\mapsto \mu_{a, b, c, d}(z).
                    \end{align*}
                    Note that because the entries $a$, $b$, $c$, and $d$ of a matrix in $\specialLinearGroup(2, \R)$ satisfy $a d - b c = 1$, this is an action by orientation-preserving isometries. The stabiliser of $\imaginaryUnit$ under $\actsOnPoint$ is $\specialOrthogonalGroup(2, \R)$. According to theorem~21.18 in \cite{lee:2013}, there is an $(\cdot, \actsOnPoint)$-e\-qui\-var\-i\-ant diffeomorphism from $\specialLinearGroup(2, \R) \modulo \specialOrthogonalGroup(2, \R)$ to $\upperHalfPlane$. In particular, the left-ho\-mo\-ge\-neous space $\ntuple{\upperHalfPlane, \specialLinearGroup(2, \R), \actsOnPoint}$ is proper. 
              \item Let $n$ be a positive integer. The Euclidean group $\EuclideanGroup(n)$ is a non-compact Lie group and the orthogonal group $\orthogonalGroup(n, \R)$ is a compact subgroup of $\EuclideanGroup(n)$. Hence, the left-ho\-mo\-ge\-neous space $\ntuple{\EuclideanGroup(n) \modulo \orthogonalGroup(n, \R), \EuclideanGroup(n), \cdot}$ is proper.

                    The group $\EuclideanGroup(n)$ acts transitively and smoothly on the Euclidean space $\R^n$ by function application denoted by $\actsOnPoint$. The stabiliser of $0$ under $\actsOnPoint$ is $\orthogonalGroup(n, \R)$. According to theorem~21.18 in \cite{lee:2013}, there is an $(\cdot, \actsOnPoint)$-e\-qui\-var\-i\-ant diffeomorphism from $\EuclideanGroup(n) \modulo \orthogonalGroup(n, \R)$ to $\R^n$. In particular, the left-ho\-mo\-ge\-neous space $\ntuple{\R^n, \EuclideanGroup(n), \actsOnPoint}$ is proper. 
            \end{aenumerate}
      \item According to the first paragraph of section~6 in \cite{kassel:2012}, the indefinite unitary group $\unitaryGroup(n, 1)$ acts properly and transitively on the $(2n + 1)$-dimensional anti-de Sitter space $\antiDeSitterSpace^{2n + 1}$. 
      \item Let $M$ be the vertices of a uniform tiling of the Euclidean or hyperbolic plane, let $G$ be the symmetry group of the tiling, let $\actsOnPoint$ be the transitive left group action of $G$ on $M$ by function application, and let $M$ and $G$ be equipped with their respective discrete topology. The left-ho\-mo\-ge\-neous space $\ntuple{M, G, \actsOnPoint}$ is proper. \qedhere 
    \end{aenumerate}
  \end{example} 

  \begin{example}[Semi-Proper]
  \label{exp:semi-proper}
    \begin{aenumerate} 
      \item According to \cref{remark:proper-implies-semi-proper}, each left group set from \cref{exp:proper} is semi-proper.
      \item According to the forthcoming \cref{lemma:discrete-semi-proper-cell-space}, each left group set $\ntuple{M, G, \actsOnPoint}$, where $M$ and $G$ are equipped with their respective discrete topology, is semi-proper.
      \item\label{item:semi-proper:but-not-proper}
            Let $\mathcal{M} = \ntuple{M, G, \actsOnPoint}$ be a left group set, let $M$ and $G$ be equipped with their respective discrete topology, let $m$ be an element of $M$ such that the stabiliser of $m$ under $\actsOnPoint$ is infinite. According to the forthcoming \cref{lemma:discrete-semi-proper-cell-space}, the group set $\mathcal{M}$ is semi-proper. However, because the singleton set $\setOf{(m, m)}$ is compact (that is, finite) and its preimage under the action map of $\mathcal{M}$ is $G_m \times \setOf{m}$ is not compact (that is, not finite), the group set $\mathcal{M}$ is not proper. For example:
            \begin{aenumerate}
              \item The group $\Z^2$ acts transitively on $\Z$ on the left by
                    \begin{align*}
                      \oplus \from \Z^2 \times \Z &\to \Z,\\
                      ((x, y), z) \mapsto z + x.
                    \end{align*}
                    Each element of $\Z$ has the infinite stabiliser $\setOf{0} \times \Z$. Hence, the left-ho\-mo\-ge\-neous space $\ntuple{\Z^2, \Z, \oplus}$ is semi-proper but not proper.
              \item Let $F_2$ be the free group over $\setOf{a, b}$, where $a \neq b$, let $\varphi$ be the homomorphism from $F_2$ to $\Z$ that is uniquely determined by $\varphi(a) = 1$ and $\varphi(b) = 0$, and let $\actsOnPoint$ be the transitive left group action of $F_2$ on $\Z$ given by $(g, z) \mapsto \varphi(g) + z$. For each integer $z \in \Z$, the stabiliser of $z$ is $\varphi^{-1}(0)$, which is the infinite set of all elements of $F_2$ that can be written as products of $a$, $b$, $a^{-1}$, and $b^{-1}$ with the same number of occurrences of $a$ and $a^{-1}$. Hence, the left-ho\-mo\-ge\-neous space $\ntuple{F_2, \Z, \actsOnPoint}$ is semi-proper but not proper. 
              \item Let $S$ be an infinite set, let $F$ be the free group over $S$, let $M$ be the uncoloured $S$-Cayley graph of $F$, let $G$ be the graph automorphisms of $M$, and let $\actsOnPoint$ be the transitive left group action of $G$ on $M$ by function application. For each element $m \in M$, because there are $\cardinalityOf{S}$-many outgoing edges to isomorphic subgraphs, the stabiliser of $m$ is infinite. Hence, the left-ho\-mo\-ge\-neous space $\ntuple{G, M, \actsOnPoint}$ is semi-proper but not proper. \qedhere 
            \end{aenumerate}
    \end{aenumerate}
  \end{example}

  The union of the transporters from a compact set to a compact set under a proper group action is compact, and each transversal of the transporters from a compact set to a compact set under a semi-proper group action is included in a compact set.

  \begin{lemma}
  \label{lemma:proper-compact-subset-of-group}
    Let $\mathcal{M} = \ntuple{M, G, \actsOnPoint}$ be a left group set, and let $K$ and $K'$ be two compact subsets of $M$.
    \begin{aenumerate}
      \item\label{item:proper-compact-subset-of-group:proper}
            If $\mathcal{M}$ is proper, then the set
            \begin{equation*}
              \bigcup_{(k, k') \in K \times K'} G_{k', k}
            \end{equation*}
            is a compact subset of $G$.
      \item\label{item:proper-compact-subset-of-group:semi-proper}
            If $\mathcal{M}$ is semi-proper, then each transversal of the set
            \begin{equation*}
              \setOf{G_{k', k} \suchThat (k, k') \in K \times K'}
            \end{equation*}
            is included in a compact subset of $G$. \qedhere
    \end{aenumerate}
  \end{lemma}

  \begin{proof} 
    According to Tychonoff's \cref{theorem:Tychonoff}, and because product and subspace topologies behave well with each other, the set $K \times K'$ is a compact subset of $M \times M$.
    \begin{aenumerate}
      \item Let $\mathcal{M}$ be proper. Then, because the action map $\alpha$ of $\mathcal{M}$ is proper, the preimage
            \begin{equation*}
              \alpha^{-1}(K \times K') = \setOf{(g, k') \in G \times K' \suchThat g \actsOnPoint k' \in K}
            \end{equation*}
            is compact. Hence, because the canonical projection $\pi$ of $G \times M$ onto $G$ is continuous, the image
            \begin{equation*}
              \pi(\alpha^{-1}(K \times K')) = \setOf{g \in G \suchThat \Exists k' \in K' \SuchThat g \actsOnPoint k' \in K}
            \end{equation*}
            is compact. This set is just
            \begin{equation*}
              \bigcup_{k \in K} \bigcup_{k' \in K'} \setOf{g \in G \suchThat g \actsOnPoint k' = k} = \bigcup_{(k, k') \in K \times K'} G_{k', k}.
            \end{equation*}
      \item Let $\mathcal{M}$ be semi-proper. Then, because the action map $\alpha$ of $\mathcal{M}$ is semi-proper and the preimage of each tuple $(k, k') \in K \times K'$ under $\alpha$ is $G_{k', k} \times \setOf{k'}$, each transversal of
            \begin{equation*}
              \setOf{G_{k', k} \times \setOf{k'} \suchThat (k, k') \in K \times K'}
            \end{equation*}
            is included in a compact subset of $G \times M$. Hence, because the canonical projection $\pi$ of $G \times M$ onto $G$ is continuous, each transversal of
            \begin{equation*} 
              \setOf{G_{k', k} \suchThat (k, k') \in K \times K'}
            \end{equation*}
            is included in a compact subset of $G$. \qedhere
    \end{aenumerate}
  \end{proof}


  \begin{corollary}
  \label{corollary:stabiliser-of-proper-compact}
    Let $\ntuple{M, G, \actsOnPoint}$ be a proper left group set and let $m$ be an element of $M$. The stabiliser $G_m$ of $m$ under $\actsOnPoint$ is a compact subset of $G$.
  \end{corollary}

  \begin{proof}
    This is a direct consequence of \cref{item:proper-compact-subset-of-group:proper} of \cref{lemma:proper-compact-subset-of-group} with $K = K' = \setOf{m}$. 
  \end{proof}

  \begin{remark}
  \label{remark:proper-implies-compact-stabilisers-semi-proper-does-not}
    The stabilisers of semi-proper group sets need not be compact, as \cref{item:semi-proper:but-not-proper} of \cref{exp:semi-proper} illustrates.
  \end{remark}

  Equipping a bare left group set with the discrete topology yields a semi-proper group set. If its stabilisers are finite, then it is even proper; but if they are infinite, then it is not proper (see \cref{corollary:stabiliser-of-proper-compact}).

  \begin{lemma}
  \label{lemma:discrete-semi-proper-cell-space}
    Let $\mathcal{M} = \ntuple{M, G, \actsOnPoint}$ be a left group set. Equip $M$ and $G$ with their respective discrete topology. The group set $\mathcal{M}$ is semi-proper. And, if its stabilisers are finite, then it is even proper.
  \end{lemma}

  \begin{proof}
    According to \cref{lemma:discrete-topological-cell-space}, the group set $\mathcal{M}$ is topological. And, because the product topology on $G \times M$ and the one on $M \times M$ are discrete, each subset of $G \times M$ and each of $M \times M$ is compact if and only if it is finite.

    First, let $K$ be a finite subset of $M \times M$. Then, because the set $\setOf{\alpha^{-1}(k) \suchThat k \in K}$, where $\alpha$ is the action map of $\mathcal{M}$, is finite, so is each of its transversals. Thus, the map $\alpha$ is semi-proper. In conclusion, the group set $\mathcal{M}$ is semi-proper.

    Secondly, let the stabilisers of $\mathcal{M}$ be finite. Furthermore, let $K$ be a finite subset of $M \times M$. Then, for each tuple $(m', m) \in K$, either the transporter $G_{m, m'}$ is empty, and hence finite, or there is an element $g \in G_{m, m'}$ and thus $G_{m, m'} = g G_m$, and therefore $G_{m, m'}$ is finite; thus, in either case, the preimage $\alpha^{-1}((m', m)) = G_{m, m'} \times \setOf{m}$ is finite. Hence, the preimage $\alpha^{-1}(K) = \bigcup_{(m', m) \in K} \alpha^{-1}((m', m))$ is finite. Therefore, the action map $\alpha$ is proper. In conclusion, the group set $\mathcal{M}$ is proper.
  \end{proof}


  Each proper left-ho\-mo\-ge\-neous space is \emph{almost} tame and each semi-proper left-ho\-mo\-ge\-neous space is semi-tame (see \cref{lemma:semi-proper-implies-semi-tame}). However, (semi-)properness is not necessary for (semi-)tameness (see \cref{example:not-semi-proper-but-tame-for-each-coordinate-system}).

  \begin{lemma}
  \label{lemma:semi-proper-implies-semi-tame}
    Let $\mathcal{M} = \ntuple{M, G, \actsOnPoint}$ be a left-ho\-mo\-ge\-neous space. If it is semi-proper, then it is semi-tame. And, if it is proper, then it is \emph{almost} tame, in the sense that, for each compact subset $K$ of $M$ and each compact subset $E$ of $G \modulo G_0$ such that $G_0 \cdot E \subseteq E$, the set $K \isSemiActedUponBy E$ is a compact subset of $M$.
  \end{lemma}

  \begin{proof}
    First, let $\mathcal{M}$ be semi-proper. Furthermore, let $\mathcal{K} = \ntuple{m_0, \family{g_{m_0, m}}_{m \in M}}$ be a coordinate system for $\mathcal{M}$, let $K$ be a compact subset of $M$, and let $E$ be a compact subset of $G \modulo G_0$. Then, because $k \isSemiActedUponBy e = g_{m_0, k} \actsOnPoint (m_0 \isSemiActedUponBy e)$, for $k \in K$ and $e \in E$,
    \begin{equation*}
      K \isSemiActedUponBy E
      = \bigcup_{k \in K} g_{m_0, k} \actsOnPoint (m_0 \isSemiActedUponBy E)
      = \setOf{g_{m_0, k} \suchThat k \in K} \actsOnPoint (m_0 \isSemiActedUponBy E).
    \end{equation*}
    Both the singleton set $\setOf{m_0}$ and the set $K$ are compact. Thus, according to \cref{item:proper-compact-subset-of-group:semi-proper} of \cref{lemma:proper-compact-subset-of-group}, the transversal $\setOf{g_{m_0, k} \suchThat k \in K}$ of $\setOf{G_{m_0, k} \suchThat (m_0, k) \in \setOf{m_0} \times K}$ is included in a compact subset of $G$. Hence, the product $\setOf{g_{m_0, k} \suchThat k \in K} \times (m_0 \isSemiActedUponBy E)$ is included in a compact subset of $G \times M$. Therefore, because $\actsOnPoint$ is continuous, the image $\setOf{g_{m_0, k} \suchThat k \in K} \actsOnPoint (m_0 \isSemiActedUponBy E)$ is included in a compact subset of $M$. Thus, the set $K \isSemiActedUponBy E$ is included in a compact subset of $M$. Hence, the cell space $\ntuple{\mathcal{M}, \mathcal{K}}$ is semi-tame. In conclusion, the homogeneous space $\mathcal{M}$ is semi-tame.

    Secondly, let $\mathcal{M}$ be proper. Furthermore, let $\mathcal{K} = \ntuple{m_0, \family{g_{m_0, m}}_{m \in M}}$ be a coordinate system for $\mathcal{M}$, let $K$ be a compact subset of $M$, and let $E$ be a compact subset of $G \modulo G_0$ such that $G_0 \cdot E \subseteq E$. Then, because $k \isSemiActedUponBy g_0 \cdot e = g_{m_0, k} g_0 \actsOnPoint (m_0 \isSemiActedUponBy e)$, for $k \in K$, $g_0 \in G_0$, and $e \in E$,
    \begin{align*}
      K \isSemiActedUponBy E
      &= \bigcup_{k \in K} k \isSemiActedUponBy G_0 \cdot E\\
      &= \bigcup_{k \in K} g_{m_0, k} G_0 \actsOnPoint (m_0 \isSemiActedUponBy E)\\
      &= (\bigcup_{k \in K} G_{m_0, k}) \actsOnPoint (m_0 \isSemiActedUponBy E).
    \end{align*}
    The set $K$ is compact and so is the singleton set $\setOf{m_0}$. Thus, according to \cref{item:proper-compact-subset-of-group:proper} of \cref{lemma:proper-compact-subset-of-group}, the set $\bigcup_{k \in K} G_{m_0,k}$ is compact. And, because $E$ is compact, the set $m_0 \isSemiActedUponBy E$ is compact. Hence, the product $(\bigcup_{k \in K} G_{m_0,k}) \times (H \actsOnPoint m_0)$ is compact. Therefore, because $\actsOnPoint$ is continuous, the image $(\bigcup_{k \in K} G_{m_0,k}) \actsOnPoint (H \actsOnPoint m_0)$ is compact. In conclusion, the set $K \isSemiActedUponBy E$ is compact.
  \end{proof}

  We take a closer look at \emph{almost} tameness in

  \begin{remark}
    Let $\mathcal{M} = \ntuple{M, G, \actsOnPoint}$ be a left-ho\-mo\-ge\-neous space and let it be called \define{almost tame}\graffito{almost tame}\index{tame!almost} if and only if it is topological and, for each coordinate system $\mathcal{K} = \ntuple{m_0, \family{g_{m_0, m}}_{m \in M}}$ for $\mathcal{M}$, the stabiliser $G_0$ is compact, and, for each compact subset $K$ of $M$ and each compact subset $E$ of $G \modulo G_0$ such that $G_0 \cdot E \subseteq E$, the set $K \isSemiActedUponBy E$ is a compact subset of $M$.

    If $\mathcal{M}$ is proper, then it is almost tame (this follows from \cref{corollary:stabiliser-of-proper-compact} and \cref{lemma:semi-proper-implies-semi-tame}). And, if it is almost tame, then it is semi-tame (this follows from the fact that, if $E$ is a compact subset of $G \modulo G_0$, then $G_0 \cdot E$ is a compact subset of $G \modulo G_0$, where we use that $\mathcal{M}$ is topological and that $G_0$ is compact). However, for tameness and semi-tameness we do not need $\mathcal{M}$ to be topological, and there is no simple relationship between tameness and almost tameness.
  \end{remark}

  There are left-ho\-mo\-ge\-neous spaces that are not semi-proper but nevertheless tame, which is illustrated by

  \begin{example}[Euclidean Space]
  \label{example:not-semi-proper-but-tame-for-each-coordinate-system}
    For each $d \in \N_+$, equip $\R^d$ with the Euclidean topology. Note that, for each $d \in \N_+$ and each $d' \in \N_+$, the product topology on $\R^d \times \R^{d'}$ is the Euclidean topology on $\R^{d + d'}$.

    The groups $\R^2$ and $\R$ under addition are topological. And, the group $\R^2$ acts continuously and transitively on $\R$ on the left by $\actsOnPoint \from ((x, y), z) \mapsto x + z$. Each stabiliser under $\actsOnPoint$ is $\setOf{0} \times \R$. And, each right quotient set semi-action of $\R^2 \modulo (\setOf{0} \times \R)$ on $\R$ induced by $\actsOnPoint$ is $\isSemiActedUponBy \from (z, (x, y) + \setOf{0} \times \R) \mapsto z + x$.

    Under the identification of $\R$ with $\R^2 \modulo (\setOf{0} \times \R)$ by the homeomorphic map $\iota \from x \mapsto \setOf{x} \times \R$, the semi-action $\isSemiActedUponBy$ is the continuous map $(z, x) \mapsto z + x$ from $\R^2$ to $\R$. Recall that images of compact subsets under continuous maps are compact. Hence, for each compact subset $K$ of $\R$ and each compact subset $E$ of $\R$, because $K \times E$ is a compact subset of $\R^2$, the set $K \isSemiActedUponBy E$ is a compact subset of $\R$. Therefore, the left-ho\-mo\-ge\-neous space $\setOf{\R, \R^2, \actsOnPoint}$ is tame.

    The action map of $\actsOnPoint$ is $\alpha \from ((x, y), z) \mapsto (x + z, z)$. The diagonal $K$ of the unit square in $\R^2$, namely $\setOf{(z, z) \suchThat z \in \closedInterval{0, 1}}$, is compact. Moreover, for each element $(z, z) \in K$, its preimage under $\alpha$ is $(\setOf{0} \times \R) \times \setOf{z}$. Hence, the set $T = \setOf{((0, 0), 0)} \cup \setOf{((0, 1/z), z) \suchThat z \in \leftOpenAndRightClosedInterval{0, 1}}$ is a transversal of $\setOf{\alpha^{-1}((z, z)) \suchThat (z, z) \in K}$. The transversal $T$ is not included in a compact subset of $\R^2 \times \R$. Therefore, the left-ho\-mo\-ge\-neous space $\setOf{\R, \R^2, \actsOnPoint}$ is not semi-proper.

    Analogously, the group $\R^2$ acts on the \emph{compact} topological circle $\R \modulo \Z$ on the left by $\actsOnPoint \from ((x, y), z + \Z) \mapsto (x + z) + \Z$ and the left-ho\-mo\-ge\-neous space $\setOf{\R \modulo \Z, \R^2, \actsOnPoint}$ is tame but not semi-proper.
  \end{example} 

  There are left-ho\-mo\-ge\-neous spaces that are not semi-proper but, if equipped with certain coordinate systems, nevertheless tame or semi-tame, which is illustrated by

  \begin{example}[Affine Space]
  \label{example:not-semi-proper-but-tame-for-some-coordinate-systems}
    In the situation of \cref{example:affine-reals:coordinate-maps}, the left-ho\-mo\-ge\-neous space $\mathcal{M}$ is not semi-proper, but, for some coordinate systems $\mathcal{K}$, the cell space $\ntuple{\mathcal{M}, \mathcal{K}}$ is tame or semi-tame. Indeed, the action map $\alpha$ is $((t, d), m) \mapsto (t + d \cdot m, m)$, the preimage $\alpha^{-1}(\closedInterval{0, 1} \times \setOf{1})$ is $\setOf{\setOf{((m' - d, d), 1) \suchThat d \in \R \smallsetminus \setOf{0}} \suchThat m' \in \closedInterval{0, 1}}$, and an unbounded transversal of that preimage is $\setOf{((-1, 1), 1)} \cup \setOf{((m' - 1/m', 1/m'), 1) \suchThat m' \in \leftOpenAndRightClosedInterval{0, 1}}$, which is not included in a compact set. 
  \end{example}

  \section{Uniform Curtis-Hedlund-Lyndon Theorem} 
  \label{section:characterisation}

  \paragraph{Introduction.} A metric space is a uniform space and a uniform space is a topological space. The uniformity induced by a metric is generated by the set of entourages that contains for each positive size the entourage of open balls of that size about each point (think of a family of neighbourhoods, one for each point, that are comparable in size). The topology induced by a uniformity is the set that contains the neighbourhoods of all points that can be extracted from the entourages (by throwing all neighbourhoods in one pot we lose the ability to compare the sizes of neighbourhoods of different points).

  A map on a metric space is continuous if, by bringing two points sufficiently close together, their images can be brought arbitrarily close together; and it is uniformly continuous if the sufficient closeness of two points that is needed to get a certain closeness of their images is the same for all pairs of points. In topological terms, continuity means that preimages of open sets are open; and in uniform terms, uniform continuity means that preimages of entourages are entourages. The concept of uniform continuity cannot be expressed in topological terms, because the sizes of neighbourhoods of different points are not comparable.

  The global transition function of a traditional cellular automaton is uniform and local (see \cref{section:introduction}). It is uniform because all cells use the same local transition function, and local because the neighbourhood is finite. Uniformity is equivalent to equivariance under translations of global configurations and locality to continuity with respect to the prodiscrete topology on global configurations.

  If the set of cells is discrete, then locality is naturally expressed by requiring neighbourhoods to be finite and characterised by uniform continuity, if the set of states is infinite, and continuity, otherwise, where the set of global configurations is equipped with the prodiscrete uniformity or topology, which is generated by cylinders with finite bases. And, if the set of cells is continuous, then locality is naturally expressed by requiring neighbourhoods to be compact and characterised by uniform continuity, where the set of global configurations is equipped with the uniformity that is generated by cylinders with compact bases. In the first case, if the set of states is finite, then continuity suffices in the characterisation; in the second case though, if the neighbourhood is infinite, then even for finite sets of states, it seems that continuity is not sufficient in the characterisation. Moreover, it may seem more natural to choose a bounded and closed neighbourhood. However, boundedness is only defined on metric spaces. And, according to the Heine-Borel theorem, on Euclidean spaces, compactness is characterised by boundedness and closedness. So, compactness seems and turns out to be a good alternative. 

  For the characterisation, if the set of states is infinite, then we require the left group action to be semi-tame. The reason is that for a global transition function to be continuous, for each cylinder $\mathfrak{C}$ with a compact base $K$, there must be a cylinder $\mathfrak{C}'$ with a compact base whose image is included in $\mathfrak{C}$. Because the new states of the cells in $K$ are uniquely determined by the current states of cells in $K \isSemiActedUponBy N$, the base of the cylinder $\mathfrak{C}'$ must include $K \isSemiActedUponBy N$. Hence, for there to be such a cylinder $\mathfrak{C}'$, the set $K \isSemiActedUponBy N$ must be included in a compact set, which is the case if the left group action is semi-tame. 

  \paragraph{Contents.} Given a semi-tame topological structure on the set of cells, we equip the set of global configurations with a uniform and a topological structure (see \cref{definition:topology-and-continuity-on-configurations}) that generalise the prodiscrete topology and uniformity (see \cref{remark:topology-and-uniformity-of-discrete-convergence-on-compacta:prodiscrete-uniformity}), and prove a uniform and a topological variant of the Curtis-Hedlund-Lyndon theorem. In \cref{theorem:uniform-Curtis-Hedlund-Lyndon}, the uniform variant, we show that global transition functions are characterised by $\actsOnMap$-equivariance and uniform continuity. And in its \cref{corollary:Curtis-Hedlund-Lyndon}, the topological variant, that under the assumption that the set of states is finite, they are characterised by $\actsOnMap$-equivariance and continuity. Note that, if the set of cells is compact or finite respectively, then uniform continuity or continuity are insubstantial (see \cref{remark:if-compact-then-Curtis-Hedlund-Lyndon-degenerates}).


  \paragraph{Body.} We equip the set of global configurations with the topology and the uniformity generated by cylinders with compact bases in

  \begin{definition} 
  \label{definition:topology-and-continuity-on-configurations}
    Let $M$ be a topological space and let $Q$ be a set.
    \begin{aenumerate}
      \item The topology on $Q^M$ that has for a subbase the sets\graffito{$\mathfrak{E}(K, b)$, for $b \in Q^K$ and $K \subseteq M$ compact} 
            \begin{multline*}
              \mathfrak{E}(K, b) = \setOf{c \in Q^M \suchThat c\restrictedTo_K = b}, \index[symbols]{EKbfraktur@$\mathfrak{E}(K, b)$}\\
              \text{ for } b \in Q^K \text{ and } K \subseteq M \text{ compact},
            \end{multline*}
            is called \define{topology of discrete convergence on compacta}\graffito{topology of discrete convergence on compacta}\index{of discrete convergence on compacta!topology}\index{discrete convergence on compacta@of discrete convergence on compacta}\index{discrete convergence on compacta!topology@topology}.
      \item The uniformity on $Q^M$ that has for a subbase the sets\graffito{$\mathfrak{E}(K)$, for $K \subseteq M$ compact} 
            \begin{multline*}
              \mathfrak{E}(K) = \setOf{(c, c') \in Q^M \times Q^M \suchThat c\restrictedTo_K = c'\restrictedTo_K}, \index[symbols]{EKfraktur@$\mathfrak{E}(K)$}\\
              \text{ for } K \subseteq M \text{ compact},
            \end{multline*}
            is called \define{uniformity of discrete convergence on compacta}\graffito{uniformity of discrete convergence on compacta}\index{of discrete convergence on compacta!uniformity}\index{discrete convergence on compacta!uniformity@uniformity}. \qedhere
    \end{aenumerate}
  \end{definition}

  \begin{remark}
    Because finite intersections of compact sets are compact, the topology and uniformity of discrete convergence on compacta have for a base the sets $\mathfrak{E}(K, b)$, for $b \in Q^K$ and $K \subseteq M$ compact, and $\mathfrak{E}(K)$, for $K \subseteq M$ compact, respectively.
  \end{remark}

  \begin{remark} 
  \label{remark:topology-induced-by-uniformity-of-discrete-convergence-on-compacta}
    The topology induced by the uniformity of discrete convergence on compacta on $Q^M$ is the topology of discrete convergence on compacta on $Q^M$.
  \end{remark}

  \begin{remark}
  \label{remark:compact-cell-space-uniformity-of-discrete-is-discrete}
    If the topological space $M$ is compact, then the topology and uniformity of discrete convergence on compacta on $Q^M$ are the discrete topology and uniformity on $Q^M$.
  \end{remark}

  \begin{remark}
    Because each finite subset of $M$ is compact, the topology and uniformity of discrete convergence on compacta on $Q^M$ are finer than the prodiscrete topology and uniformity on $Q^M$. The prodiscrete uniformity is introduced in \cref{definition:prodiscrete-uniformity} of \cref{chapter:uniformities} and a base for it is given in \cref{remark:base-of-prodiscrete-uniformity}.
  \end{remark}

  \begin{remark}
    Let $M$ be a topological space and let $Q$ be a set. For each compact subset $K$ of $M$, equip $Q^K$ with the discrete topology. The topology of discrete convergence on compacta on $Q^M$ is the coarsest topology on $Q^M$ such that, for each compact subset $K$ of $M$, the projection
    \begin{align*}
      \pi_K \from Q^M &\to Q^K, \mathnote{$\pi_K$, $K \subseteq M$ compact}\index[symbols]{piK@$\pi_K$}\\
      c &\mapsto c\restrictedTo_K,
    \end{align*}
    is continuous. Note that $\pi_K^{-1}(b) = \mathfrak{E}(K, b)$, for $b \in Q^K$ and $K \subseteq M$ compact.
  \end{remark}

  If the set of cells carries the discrete topology, then the topology and uniformity on global configurations reduce to the prodiscrete topology and uniformity.

  \begin{remark}
  \label{remark:topology-and-uniformity-of-discrete-convergence-on-compacta:prodiscrete-uniformity}
    Let $M$ be equipped with the discrete topology. Then, for each subset $A$ of $M$, the set $A$ is compact if and only if it is finite. Therefore, the topology and uniformity of discrete convergence on compacta on $Q^M$ are the prodiscrete topology and uniformity on $Q^M$.
  \end{remark}

  If the set of cells and the group that acts on it carry the discrete topology, and the set of states is finite, then uniform notions on global configurations reduce to topological ones.

  \begin{remark} 
  \label{remark:discrete-topological-cell-space}
    Let $\mathcal{R} = \ntuple{\ntuple{M, G, \actsOnPoint}, \ntuple{m_0, \family{g_{m_0, m}}_{m \in M}}}$ be a cell space and let $Q$ be a finite set. Equip $M$ with the discrete topology, and equip $Q^M$ with the prodiscrete topology. According to \cref{remark:discrete-tame-cell-space}, the cell space $\mathcal{R}$ is tame. Because the topology on $G \modulo G_0$ is discrete, each subset $E$ of $G \modulo G_0$ is compact if and only if it is finite. Because $Q^M$ is compact (see \cref{lemma:phase-space-is-compact}), each map $\Delta \from Q^M \to Q^M$ is uniformly continuous if and only if it is continuous. And, because $Q^M$ is Hausdorff (see \cref{lemma:phase-space-is-Hausdorff}), each map $\Delta \from Q^M \to Q^M$, is a uniform isomorphism if and only if it is continuous and bijective.
  \end{remark}

  If the set of states carries the discrete topology, then the well-known uniformity of uniform convergence on compacta is the uniformity of discrete convergence on compacta.

  \begin{definition}
    Let $M$ be a topological space, let $Q$ be a uniform space, and let $\mathcal{U}$ be the uniformity of $Q$. The uniformity on $Q^M$ that has for a subbase the sets\graffito{$\mathfrak{E}(K, E)$, for $E \in \mathcal{U}$ and $K \subseteq M$ compact} 
    \begin{multline*}
      \mathfrak{E}(K, E) = \setOf{(c,c') \in Q^M \times Q^M \suchThat \ForEach k \in K \Holds (c(k),c'(k)) \in E}, \index[symbols]{EKEfraktur@$\mathfrak{E}(K, E)$}\\
      \text{ for } E \in \mathcal{U} \text{ and } K \subseteq M \text{ compact},
    \end{multline*}
    is called \define{uniformity of uniform convergence on compacta}\graffito{uniformity of uniform convergence on compacta}.
  \end{definition}

  \begin{remark}
    If the uniform space $Q$ is (uniformly) discrete, then the uniformity of uniform convergence on compacta on $Q^M$ is the uniformity of discrete convergence on compacta on $Q^M$. Indeed, because the diagonal $D = \setOf{(q, q) \suchThat q \in Q}$ is an entourage of the uniform space $Q$, for each compact subset $K$ of $M$, the set
    \begin{equation*}
      \mathfrak{E}(K, D) = \setOf{(c, c') \in Q^M \times Q^M \suchThat \ForEach k \in K \Holds c(k) = c'(k)}
    \end{equation*}
    is an entourage of the uniform space $Q^M$, which is included in each of the entourages 
    \begin{multline*}
      \mathfrak{E}(K, E) = \setOf{(c,c') \in Q^M \times Q^M \suchThat \ForEach k \in K \Holds (c(k), c'(k)) \in E},\\
      \text{ for } E \in \mathcal{U}.
    \end{multline*}
    Therefore, the uniformity on $Q^M$ has for a subbase the sets $\mathfrak{E}(K, D) = \mathfrak{E}(K)$, for $K \subseteq M$ compact.
  \end{remark}

  \begin{remark}
    Let $Q^M$ be equipped with the topology of uniform convergence on compacta and let
    \begin{equation*}
      \continuousMaps(M, Q) = \setOf{c \from M \to Q \suchThat c \text{ is continuous}} \subseteq Q^M
    \end{equation*}
    be equipped with the subspace topology. According to theorem~43.7 in \cite{willard:2012}, this topology on $\continuousMaps(M, Q)$ is the well-known compact-open topology. 
  \end{remark}

  Global transition functions of cellular automata over semi-tame cell spaces with compact sufficient neighbourhoods are characterised by equivariance and continuity.

  \begin{main-theorem}[Uniform Variant; Morton Landers Curtis, Gustav Arnold Hedlund, and Roger Conant Lyndon, 1969]
  \label{theorem:uniform-Curtis-Hedlund-Lyndon}
    Let $\mathcal{R} = \ntuple{\mathcal{M}, \mathcal{K}} = \ntuple{\ntuple{M, G, \actsOnPoint}, \ntuple{m_0, \family{g_{m_0, m}}_{m \in M}}}$ be a semi-tame cell space, let $Q$ be a set, let $\Delta$ be a map from $Q^M$ to $Q^M$, let $Q^M$ be equipped with the uniformity of discrete convergence on compacta, and let $H$ be a $\mathcal{K}$-big subgroup of $G$. The following two statements are equivalent:
    \begin{aenumerate}
      \item \label{item:uniform-Curtis-Hedlund-Lyndon:global-transition-function}
            The map $\Delta$ is the global transition function of a semi-cellular automaton over $\mathcal{R}$ with $\bullet_{H_0}$-invariant local transition function and compact sufficient neighbourhood.
      \item \label{item:uniform-Curtis-Hedlund-Lyndon:equivariant-and-continuous}
            The map $\Delta$ is $\actsOnMap_H$-e\-qui\-var\-i\-ant and uniformly continuous. \qedhere
    \end{aenumerate}
  \end{main-theorem}

  \begin{proof} 
    First, let $\Delta$ be the global transition function of a semi-cellular automaton $\mathcal{C} = \ntuple{\mathcal{R}, Q, N, \delta}$ with $\bullet_{H_0}$-invariant local transition function $\delta$ and compact sufficient neighbourhood $E$. Then, according to \cref{theorem:local-invariance-versus-global-equivariance}, the map $\Delta$ is $\actsOnMap_H$-e\-qui\-var\-i\-ant. Moreover, let $K$ be a compact subset of $M$. Because $\mathcal{R}$ is semi-tame, the set $K \isSemiActedUponBy E$ is included in a compact subset $L$ of $M$. For each $c \in Q^M$ and each $c' \in Q^M$, if $c\restrictedTo_{K \isSemiActedUponBy E} = c'\restrictedTo_{K \isSemiActedUponBy E}$, then $\Delta(c)\restrictedTo_K = \Delta(c')\restrictedTo_K$, in particular, if $c\restrictedTo_L = c'\restrictedTo_L$, then $\Delta(c)\restrictedTo_K = \Delta(c')\restrictedTo_K$. Thus,
    \begin{equation*}
      (\Delta \times \Delta)(\mathfrak{E}(L)) \subseteq \mathfrak{E}(K).
    \end{equation*}
    Because the sets $\mathfrak{E}(K)$, for $K \subseteq M$ compact, constitute a base of the uniformity on $Q^M$, the global transition function $\Delta$ is uniformly continuous.

    Secondly, let $\Delta$ be $\actsOnMap_H$-e\-qui\-var\-i\-ant and uniformly continuous. Then, because $\Delta$ is uniformly continuous, there is a compact subset $E_0$ of $M$ such that
    \begin{equation*}
      (\Delta \times \Delta)(\mathfrak{E}(E_0)) \subseteq \mathfrak{E}(\setOf{m_0}).
    \end{equation*}
    Therefore, for each $c \in Q^M$, the state $\Delta(c)(m_0)$ depends at most on $c\restrictedTo_{E_0}$. The subset $E = (m_0 \isSemiActedUponBy \blank)^{-1}(E_0)$ of $G \modulo G_0$ is compact. Let $N$ be the set $G_0 \cdot E$. Then, $G_0 \cdot N \subseteq N$. And, because $E_0 = m_0 \isSemiActedUponBy E \subseteq m_0 \isSemiActedUponBy N$, for each $c \in Q^M$, the state $\Delta(c)(m_0)$ depends at most on $c\restrictedTo_{m_0 \isSemiActedUponBy E}$, in particular, it depends at most on $c\restrictedTo_{m_0 \isSemiActedUponBy N}$. Hence, there is a map $\delta \from Q^N \to Q$ such that
    \begin{equation*}
      \ForEach c \in Q^M \Holds \Delta(c)(m_0) = \delta(n \mapsto c(m_0 \isSemiActedUponBy n)).
    \end{equation*}
    The quadruple $\mathcal{C} = \ntuple{\mathcal{R}, Q, N, \delta}$ is a semi-cellular automaton with compact sufficient neighbourhood $E$. Conclude with \cref{theorem:determination-of-cellular-automata-by-behaviour-at-origin} that $\delta$ is $\bullet_{H_0}$-invariant and that $\Delta$ is the global transition function of $\mathcal{C}$.
  \end{proof}

  \begin{remark}
  \label{remark:tameness-not-used-in-the-backward-direction-of-proof-of-Curtis-Hedlund}
    Note that in the proof semi-tameness of $\mathcal{R}$ is used to deduce \cref{item:uniform-Curtis-Hedlund-Lyndon:equivariant-and-continuous} from \cref{item:uniform-Curtis-Hedlund-Lyndon:global-transition-function} but not to deduce \cref{item:uniform-Curtis-Hedlund-Lyndon:global-transition-function} from \cref{item:uniform-Curtis-Hedlund-Lyndon:equivariant-and-continuous}.
  \end{remark}

  Global transition functions of cellular automata with finite sets of states and finite sufficient neighbourhoods are characterised by equivariance and continuity. We already proved this in \cref{theorem:topological-Curtis-Hedlund-Lyndon-direct-proof}, but a slightly less general version also follows from \cref{theorem:uniform-Curtis-Hedlund-Lyndon}.

  \begin{corollary}[Topological Variant; Morton Landers Curtis, Gustav Arnold Hedlund, and Roger Conant Lyndon, 1969]
  \label{corollary:Curtis-Hedlund-Lyndon}
    Let $\mathcal{R} = \ntuple{\mathcal{M}, \mathcal{K}} = \ntuple{\ntuple{M, G, \actsOnPoint}, \ntuple{m_0, \family{g_{m_0, m}}_{m \in M}}}$ be a cell space, let $Q$ be a finite set, let $\Delta$ be a map from $Q^M$ to $Q^M$, let $Q^M$ be equipped with the prodiscrete topology, and let $H$ be a $\mathcal{K}$-big subgroup of $G$. The following two statements are equivalent:
    \begin{aenumerate}
      \item \label{item:Curtis-Hedlund-Lyndon:global-transition-function}
            The map $\Delta$ is the global transition function of a semi-cellular automaton over $\mathcal{R}$ with $\bullet_{H_0}$-invariant local transition function and finite sufficient neighbourhood.
      \item \label{item:Curtis-Hedlund-Lyndon:equivariant-and-continuous}
            The map $\Delta$ is $\actsOnMap_H$-e\-qui\-var\-i\-ant and continuous. \qedhere
    \end{aenumerate}
  \end{corollary}

  \begin{proof}
    With \cref{remark:topology-and-uniformity-of-discrete-convergence-on-compacta:prodiscrete-uniformity,remark:discrete-topological-cell-space} this follows directly from \cref{theorem:uniform-Curtis-Hedlund-Lyndon}.
  %
  %
  \end{proof}

  If the set of cells is too small, then the properties concerning locality in the Curtis-Hedlund-Lyndon theorems, namely having a compact or finite sufficient neighbourhood and being uniformly continuous or continuous, are always satisfied.

  \begin{remark}
  \label{remark:if-compact-then-Curtis-Hedlund-Lyndon-degenerates}
    In the situation of \cref{theorem:uniform-Curtis-Hedlund-Lyndon} or \cref{corollary:Curtis-Hedlund-Lyndon}, let $M$ be compact or finite respectively. Then, \cref{theorem:uniform-Curtis-Hedlund-Lyndon} and \cref{corollary:Curtis-Hedlund-Lyndon} reduce to the statement that the following two statements are equivalent:
    \begin{aenumerate}
      \item The map $\Delta$ is the global transition function of a semi-cellular automaton over $\mathcal{R}$ with $\bullet_{H_0}$-invariant local transition function.
      \item The map $\Delta$ is $\actsOnMap_H$-e\-qui\-var\-i\-ant.
    \end{aenumerate}
    Indeed, if $M$ is compact, then, for each semi-cellular automaton, there is one with the compact neighbourhood $G \modulo G_0 \simeq M$ that has the same global transition function; and, because, according to \cref{remark:compact-cell-space-uniformity-of-discrete-is-discrete}, the uniformity on $Q^M$ is the discrete uniformity, the map $\Delta$ is uniformly continuous. And, if $M$ is finite, then each semi-cellular automaton has a finite neighbourhood; and, because the topology on $Q^M$ is the discrete topology, the map $\Delta$ is continuous.
  \end{remark}

  The Curtis-Hedlund-Lyndon theorems presented above are generalisations of known theorems for cellular automata over groups.

  \begin{remark}
  \label{remark:Curtis-Hedlund-Lyndon-on-groups}
    In the case that $M = G$ and $\actsOnPoint$ is the group multiplication of $G$, \cref{theorem:uniform-Curtis-Hedlund-Lyndon} is theorem~1.9.1 in \cite{ceccherini-silberstein:coornaert:2010} and \cref{corollary:Curtis-Hedlund-Lyndon} is theorem~1.8.1 in \cite{ceccherini-silberstein:coornaert:2010}.
  \end{remark}

  The uniform variant of the Curtis-Hedlund-Lyndon theorem can be applied to some automata that we have already encountered and to generalisations of them.

  \begin{example}[Plane, Sphere, Left Shift Map]
    The global transition functions of \cref{example:plane:global-transition-function,example:sphere:global-transition-function,example:lattice:the-left-shift} are $\actsOnMap$-e\-qui\-var\-i\-ant and uniformly continuous.
  \end{example}

  \begin{example}[Riemannian Symmetric Space] 
    Extensions of the Laplace–Beltrami operators on Riemannian symmetric spaces of constant curvature, like Euclidean spaces, spheres, and hyperbolic spaces, by the constant $0$ like in \cref{example:plane:global-transition-function} are global transition functions of cellular automata and therefore $\actsOnMap$-e\-qui\-var\-i\-ant and uniformly continuous. Broadly speaking, this is true for all differential operators on such spaces that do not depend on global positions.
  \end{example} 

  On cell spaces that are not semi-tame and for subgroups that are not big, the uniform variant of the Curtis-Hedlund-Lyndon theorem does not hold, which is illustrated by

  \begin{counterexample}[Affine Space] 
    In the situation of \cref{item:affine-reals:unbounded} in \cref{example:affine-reals}, the cell space $\mathcal{R} = \ntuple{\mathcal{M}, \mathcal{K}''}$ is not semi-tame, in particular, neither is its coordinate map continuous nor is it semi-proper. Let $Q$ be the set $\setOf{0, 1}$, let $N$ be the set $\R$, and let $\delta$ be the map $Q^N \to Q$, $\ell \mapsto \ell(1)$. Recall that $M$ is identified with $G \modulo G_0$ by $\iota$.

    The quadruple $\mathcal{C} = \ntuple{\mathcal{R}, Q, N, \delta}$ is a semi-cellular automaton that has the compact sufficient neighbourhood $\setOf{1}$ and whose local transition function is $\bullet_{T_0}$-invariant, where $T$ is the subgroup $\setOf{\alpha_{t, 1} \suchThat t \in \R}$ of translations of $G$, which is not $\mathcal{K}''$-big and $T_0$ is the stabiliser of $m_0$ under $\actsOnPoint_T$. Note that, because $G_0 \actsOnPoint \setOf{1} = \R$ and we need $1 \in N$ and $G_0 \actsOnPoint N \subseteq N$, we had to choose $N = \R$.

    The global transition function $\Delta$ of $\mathcal{C}$ is the map $Q^M \to Q^M$, $c \mapsto [m \mapsto c(m \isSemiActedUponBy 1)]$. The state $\Delta(\blank)(0)$ depends on the state of cell $0 \isSemiActedUponBy 1 = 1$ and, for each cell $m \in M \smallsetminus \setOf{0}$, the state $\Delta(\blank)(m)$ depends on the state of $m \isSemiActedUponBy 1 = 1/m + m$, which diverges to $\pm \infty$ as $m$ tends to $0$ from above or below. In particular, the restriction $\Delta(\blank)\restrictedTo_{\closedInterval{0, 1}}$ depends on the states of all cells in $D = \setOf{1} \cup \leftClosedAndRightOpenInterval{2, \infty}$.

    The set $\mathfrak{D} = \setOf{(c, c') \in Q^M \times Q^M \suchThat c\restrictedTo_D = c'\restrictedTo_D}$ is the smallest set such that $(\Delta \times \Delta)(D) \subseteq \mathfrak{E}(\closedInterval{0, 1})$. However, because $D$ is unbounded, it is neither an entourage nor included in one, hence $(\Delta \times \Delta)^{-1}(\mathfrak{E}(\closedInterval{0, 1}))$ is not an entourage, and therefore $\Delta$ is not uniformly continuous.

    In conclusion, for cell spaces that are not semi-tame, \cref{item:uniform-Curtis-Hedlund-Lyndon:equivariant-and-continuous} of \cref{theorem:uniform-Curtis-Hedlund-Lyndon} does not follow from \cref{item:uniform-Curtis-Hedlund-Lyndon:global-transition-function}. However, according to \cref{remark:tameness-not-used-in-the-backward-direction-of-proof-of-Curtis-Hedlund}, even for cell spaces that are not semi-tame, \cref{item:uniform-Curtis-Hedlund-Lyndon:global-transition-function} follows from \cref{item:uniform-Curtis-Hedlund-Lyndon:equivariant-and-continuous}.
  \end{counterexample}

  Whether, under the assumption that the set of states is finite, continuity is sufficient in the uniform variant of the Curtis-Hedlund-Lyndon theorem is an

  \begin{open-problem}
    Is there a semi-tame cell space $\mathcal{R}$, a \emph{finite} set $Q$, a map $\Delta$ from $Q^M$ to $Q^M$, a $\mathcal{K}$-big subgroup $H$ of $G$ such that $\Delta$ is $\actsOnMap_H$-e\-qui\-var\-i\-ant and continuous with respect to the topology of discrete convergence on compacta on $Q^M$, but such that $\Delta$ is not the global transition function of a semi-cellular automaton over $\mathcal{R}$ with $\bullet_{H_0}$-invariant local transition function and compact sufficient neighbourhood?

    In example~1.8.2 in \cite{ceccherini-silberstein:coornaert:2010}, for each infinite group $G$, an $\actsOnMap$-e\-qui\-var\-i\-ant and continuous map $\Delta$ on $G^G$, equipped with the prodiscrete topology, is constructed that is not the global transition function of a cellular automaton over $G$ with finite neighbourhood. One may try to use the general idea of the construction to construct a $\actsOnMap$-e\-qui\-var\-i\-ant and continuous map $\Delta$ on $\setOf{0, 1}^\R$ that is not the global transition function of a cellular automaton over $(\R, +)$ with compact neighbourhood.
  \end{open-problem}

  \section{Invertibility of Big-Cellular Automata}
  \label{section:invertibility}

  \paragraph{Summary.} A cellular automaton is invertible if its computations can be made in reverse by another automaton (see \cref{definition:invertible}). If we consider only automata with compact sufficient neighbourhoods, it follows from \cref{theorem:uniform-Curtis-Hedlund-Lyndon} that a map is the global transition function of an invertible automaton if and only if it is an equivariant uniform isomorphism (see \cref{theorem:invertible-uniform-Curtis-Hedlund-Lyndon}). Hence, if we consider only those with finite set of states and finite sufficient neighbourhoods, a map is the global transition function of an invertible automaton if and only if it is equivariant, continuous, and bijective (see \cref{corollary:invertible-Curtis-Hedlund-Lyndon}). In particular, such an automaton is invertible if and only if its global transition function is bijective (see \cref{corollary:invertible-if-and-only-if-bijective}). 


  \begin{definition}
  \label{definition:invertible}
    Let $\mathcal{C} = \ntuple{\mathcal{R}, Q, N, \delta}$ be a semi-cellular automaton. It is called \define{invertible}\graffito{invertible} if and only if there is a semi-cellular automaton $\mathcal{C}'$, called \define{inverse to $\mathcal{C}$}\graffito{inverse to $\mathcal{C}$}, such that the global transition functions of $\mathcal{C}$ and $\mathcal{C}'$ are inverse to each other.
  \end{definition}

  \begin{theorem}
  \label{theorem:invertible-uniform-Curtis-Hedlund-Lyndon}
    Let $\mathcal{R} = \ntuple{\mathcal{M}, \mathcal{K}} = \ntuple{\ntuple{M, G, \actsOnPoint}, \ntuple{m_0, \family{g_{m_0, m}}_{m \in M}}}$ be a semi-tame cell space, let $Q$ be a set, let $\Delta$ be a map from $Q^M$ to $Q^M$, let $Q^M$ be equipped with the uniformity of discrete convergence on compacta, and let $H$ be a $\mathcal{K}$-big subgroup of $G$. The following two statements are equivalent:
    \begin{aenumerate}
      \item \label{item:invertible-Curtis-Hedlund-Lyndon:automaton}
            The map $\Delta$ is the global transition function of an invertible semi-cellular automaton $\mathcal{C}$ over $\mathcal{R}$ that has an inverse $\mathcal{C}'$ such that the local transition functions of $\mathcal{C}$ and $\mathcal{C}'$ are $\bullet_{H_0}$-invariant, and $\mathcal{C}$ and $\mathcal{C}'$ have compact sufficient neighbourhoods.
      \item \label{item:invertible-Curtis-Hedlund-Lyndon:equivariant}
            The map $\Delta$ is an $\actsOnMap_H$-e\-qui\-var\-i\-ant uniform isomorphism. \qedhere
    \end{aenumerate}
  \end{theorem}

  \begin{proof}
    With \cref{corollary:inverse-is-equivariant} this follows directly from main theorem \ref{theorem:uniform-Curtis-Hedlund-Lyndon}. 
  %
  \end{proof}

  \begin{corollary}
  \label{corollary:invertible-Curtis-Hedlund-Lyndon}
    Let $\mathcal{R} = \ntuple{\mathcal{M}, \mathcal{K}} = \ntuple{\ntuple{M, G, \actsOnPoint}, \ntuple{m_0, \family{g_{m_0, m}}_{m \in M}}}$ be a cell space, let $Q$ be a finite set, let $\Delta$ be a map from $Q^M$ to $Q^M$, let $Q^M$ be equipped with the prodiscrete topology, and let $H$ be a $\mathcal{K}$-big subgroup of $G$. The following two statements are equivalent:
    \begin{aenumerate}
      \item \label{item:cor-invertible-Curtis-Hedlund-Lyndon:automaton}
            The map $\Delta$ is the global transition function of an invertible semi-cellular automaton $\mathcal{C}$ over $\mathcal{R}$ that has an inverse $\mathcal{C}'$ such that the local transition functions of $\mathcal{C}$ and $\mathcal{C}'$ are $\bullet_{H_0}$-invariant, and $\mathcal{C}$ and $\mathcal{C}'$ have finite sufficient neighbourhoods.
      \item \label{item:cor-invertible-Curtis-Hedlund-Lyndon:equivariant}
            The map $\Delta$ is $\actsOnMap_H$-e\-qui\-var\-i\-ant, continuous, and bijective. \qedhere
    \end{aenumerate}
  \end{corollary}

  \begin{proof}
    With \cref{remark:discrete-topological-cell-space} this follows directly from \cref{theorem:invertible-uniform-Curtis-Hedlund-Lyndon}.
  \end{proof}


  \begin{corollary}
  \label{corollary:invertible-if-and-only-if-bijective}
    Let $\mathcal{R} = \ntuple{\mathcal{M}, \mathcal{K}} = \ntuple{\ntuple{M, G, \actsOnPoint}, \ntuple{m_0, \family{g_{m_0, m}}_{m \in M}}}$ be a cell space, let $H$ be a $\mathcal{K}$-big subgroup of $G$. Furthermore, let $\mathcal{C}$ be a semi-cellular automaton over $\mathcal{R}$ with finite set of states, finite sufficient neighbourhood, and $\bullet_{H_0}$-invariant local transition function. The automaton $\mathcal{C}$ is invertible if and only if its global transition function is bijective.
  \end{corollary}

  \begin{proof}
    With \cref{theorem:local-invariance-versus-global-equivariance} and \cref{corollary:Curtis-Hedlund-Lyndon} this follows directly from \cref{corollary:invertible-Curtis-Hedlund-Lyndon}.
  \end{proof}

  \begin{remark}
    It follows from \cref{theorem:topological-Curtis-Hedlund-Lyndon-direct-proof} that \cref{corollary:invertible-Curtis-Hedlund-Lyndon} also holds if $Q^M$ is equipped with the $(\actsOnPoint, L)$-prodiscrete topology, where $L$ is a finite subgroup of $H$.
  \end{remark}

  \begin{counterexample}
    In \cref{corollary:invertible-if-and-only-if-bijective}, if the assumption that the semi-cellular automaton $\mathcal{C}$ has a finite set of states does not hold, then it may not be invertible even if its global transition function is bijective, which is shown by the cellular automaton of example~1.10.3 in \cite{ceccherini-silberstein:coornaert:2010}.
  \end{counterexample}

  \clearToOddPage
  \chapter{Right Amenability and the Tarski-Følner Theorem}
  \label{chapter:amenability}

  \paragraph{Abstract.} We introduce right amenability, right Følner nets, and right paradoxical decompositions for left-ho\-mo\-ge\-neous spaces and prove the Tarski-Følner theorem for left-ho\-mo\-ge\-neous spaces with finite stabilisers. It states that right amenability, the existence of right Følner nets, and the non-existence of right paradoxical decompositions are equivalent.

  \paragraph{Remark.} Most parts of this chapter appeared in the paper \enquote{\citetitle*{wacker:amenable:2016}}\cite{wacker:amenable:2016} and they generalise parts of chapter~4 of the monograph \enquote{\citetitle*{ceccherini-silberstein:coornaert:2010}}\cite{ceccherini-silberstein:coornaert:2010}.

  \paragraph{Introduction.} The notion of amenability for groups was introduced by John von Neumann in 1929 in the paper \enquote{\citetitle*{von-neumann:1929}}\cite{von-neumann:1929}. It generalises the notion of finiteness. A group $G$ is \emph{left} or \emph{right amenable} if there is a finitely additive probability measure on the power set of $G$ that is invariant under left and right multiplication respectively. Groups are left amenable if and only if they are right amenable. A group is \emph{amenable} if it is left or right amenable.

  Examples of amenable groups are abound: Finite groups like the Fischer–Griess Monster group, abelian groups like the group of integers, nilpotent groups like the Heisenberg group, solvable groups like the group of invertible upper triangular matrices under multiplication, finitely generated groups of sub-exponential growth like the Grigorchuk group, and many more. Examples of non-amenable groups are the free groups whose rank is greater than $1$ and the groups that have a subgroup that is isomorphic to the free group of rank $2$. 

  The definitions of left and right amenability generalise to left and right group sets respectively. A left group set $\ntuple{M, G, \actsOnPoint}$ is \emph{left amenable} if there is a finitely additive probability measure on $\powerSetOf(M)$ that is invariant under $\actsOnPoint$. There is in general no natural action on the right that is to a left group action what right multiplication is to left group multiplication. Therefore, for a left group set there is no natural notion of right amenability.

  A transitive left group action $\actsOnPoint$ of $G$ on $M$ induces, for each element $m_0 \in M$ and each family $\family{g_{m_0, m}}_{m \in M}$ of elements in $G$ such that, for each point $m \in M$, we have $g_{m_0, m} \actsOnPoint m_0 = m$, a right quotient set semi-action $\isSemiActedUponBy$ of $G \modulo G_0$ on $M$ with defect $G_0$ given by $m \isSemiActedUponBy g G_0 = g_{m_0, m} g g_{m_0, m}^{-1} \actsOnPoint m$, where $G_0$ is the stabiliser of $m_0$ under $\actsOnPoint$. Each of these right semi-actions is to the left group action what right multiplication is to left group multiplication. They occur in the definition of global transition functions of semi-cellular automata over left-ho\-mo\-ge\-neous spaces. A \emph{coordinate system} is a choice of $m_0$ and $\family{g_{m_0, m}}_{m \in M}$. 

  A left-ho\-mo\-ge\-neous space is \emph{right amenable} if there is a coordinate system such that there is a finitely additive probability measure on $\powerSetOf(M)$ that is semi-invariant under $\isSemiActedUponBy$. For example finite left-ho\-mo\-ge\-neous spaces, abelian groups, and finitely right-gen\-er\-at\-ed left-ho\-mo\-ge\-neous spaces of sub-exponential growth are right amenable, in particular, quotients of finitely generated groups of sub-exponential growth by finite subgroups acted upon by left multiplication.

  A net of non-empty and finite subsets of $M$ is a \emph{right Følner net} if, broadly speaking, these subsets are asymptotically invariant under $\isSemiActedUponBy$. A finite subset $E$ of $G \modulo G_0$ and two partitions $\family{A_e}_{e \in E}$ and $\family{B_e}_{e \in E}$ of $M$ constitute a \emph{right paradoxical decomposition} if the map $\blank \isSemiActedUponBy e$ is injective on $A_e$ and $B_e$, and the family $\family{(A_e \isSemiActedUponBy e) \disjointUnionWith (B_e \isSemiActedUponBy e)}_{e \in E}$ is a partition of $M$. The Tarski-Følner theorem states that right amenability, the existence of right Følner nets, and the non-existence of right paradoxical decompositions are equivalent.

  The Tarski alternative theorem and the theorem of Følner, which constitute the Tarski-Følner theorem, are famous theorems by Alfred Tarski and Erling Følner from 1938 and 1955, see the papers \enquote{\citetitle*{tarski:1938}}\cite{tarski:1938} and \enquote{\citetitle*{folner:1955}}\cite{folner:1955}.


  \paragraph{Contents.} In \cref{section:measures-and-means} we introduce finitely additive probability measures and means, and kind of right semi-actions on them. In \cref{section:right-amenability} we introduce right amenability. In \cref{section:right-Folner-nets} we introduce right Følner nets and present examples that are used throughout \cref{chapter:amenability,chapter:garden}. In \cref{section:right-paradoxical-decomposition} we introduce right paradoxical decompositions. In \cref{section:Tarski-Folner-theorem} we prove the Tarski alternative theorem and the theorem of Følner. And in \cref{section:left-versus-right} we show under which assumptions left implies right amenability and give two examples of right-a\-me\-na\-ble left-ho\-mo\-ge\-neous spaces.

  \paragraph{Preliminary Notions.} A \emph{left group set} is a triple $\ntuple{M, G, \actsOnPoint}$, where $M$ is a set, $G$ is a group, and $\actsOnPoint$ is a map from $G \times M$ to $M$, called \emph{left group action of $G$ on $M$}, such that $G \to \symmetricGroupOf(M)$, $g \mapsto [g \actsOnPoint \blank]$, is a group homomorphism. The action $\actsOnPoint$ is \emph{transitive} if $M$ is non-empty and for each $m \in M$ the map $\blank \actsOnPoint m$ is surjective; and \emph{free} if for each $m \in M$ the map $\blank \actsOnPoint m$ is injective. For each $m \in M$, the set $G \actsOnPoint m$ is the \emph{orbit of $m$}, the set $G_m = (\blank \actsOnPoint m)^{-1}(m)$ is the \emph{stabiliser of $m$}, and, for each $m' \in M$, the set $G_{m, m'} = (\blank \actsOnPoint m)^{-1}(m')$ is the \emph{transporter of $m$ to $m'$}.

  A \emph{left-ho\-mo\-ge\-neous space} is a left group set $\mathcal{M} = \ntuple{M, G, \actsOnPoint}$ such that $\actsOnPoint$ is transitive. A \emph{coordinate system for $\mathcal{M}$} is a tuple $\mathcal{K} = \ntuple{m_0, \family{g_{m_0, m}}_{m \in M}}$, where $m_0 \in M$ and, for each $m \in M$, we have $g_{m_0, m} \actsOnPoint m_0 = m$. The stabiliser $G_{m_0}$ is denoted by $G_0$. The tuple $\mathcal{R} = \ntuple{\mathcal{M}, \mathcal{K}}$ is a \emph{cell space}. The set $\setOf{g G_0 \suchThat g \in G}$ of left cosets of $G_0$ in $G$ is denoted by $G \modulo G_0$. The map $\isSemiActedUponBy \from M \times G \modulo G_0 \to M$, $(m, g G_0) \mapsto g_{m_0, m} g \actsOnPoint m_0$ is a \emph{right semi-action of $G \modulo G_0$ on $M$ with defect $G_0$}, which means that
  \begin{equation*}
    \ForEach m \in M \Holds m \isSemiActedUponBy G_0 = m,
  \end{equation*}
  and
  \begin{multline*}
    \ForEach m \in M \ForEach g \in G \Exists g_0 \in G_0 \SuchThat \ForEach \mathfrak{g}' \in G \modulo G_0 \Holds\\
          m \isSemiActedUponBy g \cdot \mathfrak{g}' = (m \isSemiActedUponBy g G_0) \isSemiActedUponBy g_0 \cdot \mathfrak{g}'.
  \end{multline*}
  It is \emph{transitive}, which means that the set $M$ is non-empty and for each $m \in M$ the map $m \isSemiActedUponBy \blank$ is surjective; and \emph{free}, which means that for each $m \in M$ the map $m \isSemiActedUponBy \blank$ is injective; and \emph{semi-commutes with $\actsOnPoint$}, which means that
  \begin{multline*}
    \ForEach m \in M \ForEach g \in G \Exists g_0 \in G_0 \SuchThat \ForEach \mathfrak{g}' \in G \modulo G_0 \Holds\\
          (g \actsOnPoint m) \isSemiActedUponBy \mathfrak{g}' = g \actsOnPoint (m \isSemiActedUponBy g_0 \cdot \mathfrak{g}').
  \end{multline*}
  (See \cref{chapter:automata}.)

  A cell space $\mathcal{R}$ is \emph{finitely and symmetrically right generated} if there is a finite subset $S$ of $G \modulo G_0$ with $G_0 \cdot S \subseteq S$ and $S^{-1} \subseteq S$, where $S^{-1} = \setOf{g^{-1} G_0 \suchThat s \in S, g \in s}$, such that
  \begin{multline*} 
    \ForEach m \in M \Exists k \in \N_0 \Exists \sequence{s_i}_{i \in \setOf{1, 2, \dotsc, k}} \text{ in } S \cup S^{-1} \SuchThat\\
        \parens[\Big]{\parens[\big]{(m_0 \isSemiActedUponBy s_1) \isSemiActedUponBy s_2} \isSemiActedUponBy \dotsb} \isSemiActedUponBy s_k = m.
  \end{multline*}
  The \emph{uncoloured $S$-Cayley graph of $\mathcal{R}$} is the symmetric and $2 \cardinalityOf{S}$-regular directed graph $\mathcal{G} = \ntuple{M, \setOf{(m, m \isSemiActedUponBy s) \suchThat m \in M, s \in S}}$, the \emph{$S$-metric on $\mathcal{R}$} is the distance $\distanceOf$ on $\mathcal{G}$, and the \emph{$S$-length on $\mathcal{R}$} is the map $\lengthOf{\blank} = \distanceOf(m_0, \blank)$. For each $m \in M$ and each $\rho \in \Z$, the \emph{$S$-ball of radius $\rho$ centred at $m$} is the set $\ball(m, \rho) = \setOf{m' \in M \suchThat \distanceOf(m, m') \leq \rho}$, the \emph{$S$-sphere of radius $\rho$ centred at $m$} is the set $\sphere(m, \rho) = \setOf{m' \in M \suchThat \distanceOf(m, m') = \rho}$, the ball $\ball(m_0, \rho)$ is denoted by $\ball(\rho)$, and the sphere $\sphere(m_0, \rho)$ by $\sphere(\rho)$. (See \cref{chapter:growth}.)

  In the present chapter, we assume that the reader is familiar with the basics of the theory of dual spaces and with Hall's harem theorem. A recapitulation of the required basics and of the theory surrounding Hall's theorems is given in \cref{chapter:dual-spaces,chapter:halls-theorems}.

  \section{Finitely Additive Probability Measures and Means}
  \label{section:measures-and-means}

  In this section, let $\mathcal{R} = \ntuple{\ntuple{M, G, \actsOnPoint}, \ntuple{m_0, \family{g_{m_0, m}}_{m \in M}}}$ be a cell space.

  \paragraph{Summary.} First, we introduce finitely additive probability measures on $M$ (see \cref{definition:finitely-additive-probability-measure}); a left group action of $G$ and a kind of right quotient set semi-action of $G \modulo G_0$, both on the set of maps from $\powerSetOf(M)$ to $\closedInterval{0, 1}$, which contains all finitely additive probability measures on $M$ (see \cref{definition:left-group-action-on-measures,definition:kind-of-right-semi-action-on-measures}); and what it means for a measure to be semi-invariant under the kind of right semi-action (see \cref{definition:semi-invariance-of-measures}). Secondly, we introduce means on $M$ (see \cref{definition:mean}); a left group action of $G$ and a kind of right quotient set semi-action of $G \modulo G_0$, both on the topological dual space of $\boundedFunctionsOn(M)$, which contains all means on $M$ (see \cref{definition:left-group-action-on-means,definition:kind-of-right-semi-action-on-means}); and what it means for a mean to be invariant under the kind of right semi-action (see \cref{definition:invariance-of-means}). Lastly, we give a bijection from the set of means on $M$ to the set of finitely additive probability measures on $M$ (see \cref{theorem:means-versus-measures}). We show in \cref{section:right-amenability} that this bijection preserves (semi-)invariance (see in particular \cref{theorem:mean-characterisation-of-right-amenable}). 

  \begin{definition}
  \label{definition:finitely-additive-probability-measure}
    Let $\mu \from \powerSetOf(M) \to \closedInterval{0, 1}$ be a map. It is called
    \begin{aenumerate}
      \item \define{normalised}\graffito{normalised} if and only if $\mu(M) = 1$;
      \item \define{finitely additive}\graffito{finitely additive} if and only if, 
            \begin{multline*}
              \ForEach A \subseteq M \ForEach B \subseteq M \Holds\\
                  \parens[\big]{A \cap B = \emptyset \implies \mu(A \cup B) = \mu(A) + \mu(B)};
            \end{multline*}
      \item \define{finitely additive probability measure on $M$}\graffito{finitely additive probability measure $\mu$ on $M$} if and only if it is normalised and finitely additive.
    \end{aenumerate}
    The set of all finitely additive probability measures on $M$ is denoted by $\probabilityMeasuresOn(M)$\graffito{set $\probabilityMeasuresOn(M)$}\index[symbols]{PMM@$\probabilityMeasuresOn(M)$}.
  \end{definition}

  \begin{remark}
    Finitely additive measures, which do not need to be normalised, are known as \define{content}\graffito{content}. 
  \end{remark}

  \begin{definition}
  \label{definition:left-group-action-on-measures}
    The group $G$ acts on $\closedInterval{0, 1}^{\powerSetOf(M)}$ on the left by
    \begin{align*} 
      \actsOnMeasure \from G \times \closedInterval{0, 1}^{\powerSetOf(M)} &\to \closedInterval{0, 1}^{\powerSetOf(M)}, \mathnote{left group action $\actsOnMeasure$ of $G$ on $\closedInterval{0, 1}^{\powerSetOf(M)}$}\index[symbols]{vDash@$\actsOnMeasure$}\\
      (g, \varphi) &\mapsto [A \mapsto \varphi(g^{-1} \actsOnPoint A)],
    \end{align*}
    such that $G \actsOnMeasure \probabilityMeasuresOn(M) \subseteq \probabilityMeasuresOn(M)$. 
  \end{definition}

  \begin{definition} 
  \label{definition:kind-of-right-semi-action-on-measures}
    The quotient set $G \modulo G_0$ kind of semi-acts on $\closedInterval{0, 1}^{\powerSetOf(M)}$ on the right by\index[symbols]{Dashv@$\measureIsKindOfSemiActedUponBy$}
    \begin{align*}
      \measureIsKindOfSemiActedUponBy \from \closedInterval{0, 1}^{\powerSetOf(M)} \times G \modulo G_0 &\to \closedInterval{0, 1}^{\powerSetOf(M)}, \mathnote{kind of right quotient set semi-action $\measureIsKindOfSemiActedUponBy$ of $G \modulo G_0$ on $\closedInterval{0, 1}^{\powerSetOf(M)}$}\\
      (\varphi, \mathfrak{g}) &\mapsto [A \mapsto \varphi(A \isSemiActedUponBy \mathfrak{g})]. \qedhere
    \end{align*}
  \end{definition}


  \begin{remark}
    Let $\mu$ be a finitely additive probability measure on $M$ and let $\mathfrak{g}$ be an element of $G \modulo G_0$. Then, $M \isSemiActedUponBy \mathfrak{g}$ need not be equal to $M$ and, for each subset $A$ of $M$ and each subset $B$ of $M$ such that the sets $A$ and $B$ are disjoint, the sets $A \isSemiActedUponBy \mathfrak{g}$ and $B \isSemiActedUponBy \mathfrak{g}$ need not be disjoint. Therefore, $\mu \measureIsKindOfSemiActedUponBy \mathfrak{g}$ need neither be normalised nor finitely additive, in particular, it need not be a finitely additive probability measure on $M$.
  \end{remark}

  \begin{definition}
  \label{definition:semi-invariance-of-measures}
    Let $\varphi$ be an element of $\closedInterval{0, 1}^{\powerSetOf(M)}$. It is called \graffito{$\measureIsKindOfSemiActedUponBy$-semi-invariant}\defineX{$\measureIsKindOfSemiActedUponBy$-semi-invariant}{semi-invariant} if and only if, for each element $\mathfrak{g} \in G \modulo G_0$ and each subset $A$ of $M$ such that the map $\blank \isSemiActedUponBy \mathfrak{g}$ is injective on $A$, we have $(\varphi \measureIsKindOfSemiActedUponBy \mathfrak{g})(A) = \varphi(A)$.
  \end{definition}

  \begin{remark} 
  \label{remark:groups:measureIsKindOfSemiActedUponBy}
    Let $\mathcal{R}$ be the cell space $\ntuple{\ntuple{G, G, \cdot}, \ntuple{e_G, \family{g}_{g \in G}}}$. Then, $G_0 = \setOf{e_G}$ and $\isSemiActedUponBy = \cdot$. Hence, $\measureIsKindOfSemiActedUponBy \from (\varphi, g) \mapsto [A \mapsto \varphi(A \cdot g)]$. Except for $g$ not being inverted, this is the right group action of $G$ on $\probabilityMeasuresOn(M)$ as defined in paragraph~4 of section~4.3 in \cite{ceccherini-silberstein:coornaert:2010}. Moreover, for each element $g \in G$, the map $\blank \isSemiActedUponBy g$ is injective. Hence, being $\measureIsKindOfSemiActedUponBy$-semi-invariant is the same as being right-invariant as defined in paragraph~2 of section~4.4 in \cite{ceccherini-silberstein:coornaert:2010}. 
  \end{remark}

  \begin{remark}[Finite]
  \label{remark:finite:measure}
    If the set $M$ is finite, then the map
    \begin{align*}
      \mu \from \powerSetOf(M) &\to \closedInterval{0, 1},\\
      A &\mapsto \frac{\cardinalityOf{A}}{\cardinalityOf{M}},
    \end{align*}
    is a $\actsOnMeasure$-invariant and $\measureIsKindOfSemiActedUponBy$-semi-invariant finitely additive probability measure on $M$.
  \end{remark} 

  \begin{proof} 
    The map $\mu$ is normalised and finitely additive. Thus, it is a finitely additive probability measure on $M$. Moreover, for each element $g \in G$ and each subset $A$ of $M$, because the map $g^{-1} \actsOnPoint \blank$ is injective,
    \begin{equation*}
      (g \actsOnMeasure \mu)(A)
      = \mu(g^{-1} \actsOnPoint A)
      = \frac{\cardinalityOf{g^{-1} \actsOnPoint A}}{\cardinalityOf{M}}
      = \frac{\cardinalityOf{A}}{\cardinalityOf{M}}
      = \mu(A).
    \end{equation*}
    And, for each element $\mathfrak{g} \in G \modulo G_0$ and each subset $A$ of $M$ such that the map $\blank \isSemiActedUponBy \mathfrak{g}$ is injective on $A$,
    \begin{equation*}
      (\mu \measureIsKindOfSemiActedUponBy \mathfrak{g})(A)
      = \mu(A \isSemiActedUponBy \mathfrak{g})
      = \frac{\cardinalityOf{A \isSemiActedUponBy \mathfrak{g}}}{\cardinalityOf{M}}
      = \frac{\cardinalityOf{A}}{\cardinalityOf{M}}
      = \mu(A).
    \end{equation*}
    Hence, the finitely additive probability measure $\mu$ is $\actsOnMeasure$-invariant and $\measureIsKindOfSemiActedUponBy$-semi-invariant.
  \end{proof}

  \begin{remark}[Infinite]
  \label{remark:infinite:measure}
    If the set $M$ is infinite, then, for each finite subset $F$ of $M$, the map
    \begin{align*}
      \mu_F \from \powerSetOf(M) &\to \closedInterval{0, 1},\\
      A &\mapsto \frac{\cardinalityOf{A \cap F}}{\cardinalityOf{F}},
    \end{align*}
    is a finitely additive probability measure on $M$ that is neither $\actsOnMeasure$-invariant nor $\measureIsKindOfSemiActedUponBy$-semi-invariant. Nevertheless, for certain cell spaces over $M$, there is a net $\net{F_i}_{i \in I}$ of finite subsets of $M$ such that the net $\net{\mu_{F_i}}_{i \in I}$ converges in a weak sense to a $\actsOnMeasure$-invariant or $\measureIsKindOfSemiActedUponBy$-semi-invariant finitely additive probability measure on $M$. However, even on $\Z$, there is no explicit formula for that measure. 
  \end{remark}

  \begin{definition} 
    The vector space of bounded real-valued functions on $M$ with pointwise addition and scalar multiplication is denoted by $\boundedFunctionsOn(M)$\graffito{vector space $\boundedFunctionsOn(M)$ of bounded real-valued functions on $M$}\index[symbols]{linfinityM@$\boundedFunctionsOn(M)$},
    the supremum norm on $\boundedFunctionsOn(M)$ is denoted by $\normOf{\blank}_\infty$\graffito{supremum norm $\normOf{\blank}_\infty$ on $\boundedFunctionsOn(M)$}\index[symbols]{norminfinity@$\normOf{\blank}_\infty$},
    the topological dual space of $\boundedFunctionsOn(M)$ is denoted by $\boundedFunctionsOn(M)^*$\graffito{topological dual space $\boundedFunctionsOn(M)^*$ of $\boundedFunctionsOn(M)$}\index[symbols]{linfinityMstar@$\boundedFunctionsOn(M)^*$},
    the pointwise partial order on $\boundedFunctionsOn(M)$ is denoted by $\leq$\graffito{pointwise partial order $\leq$ on $\boundedFunctionsOn(M)$}\index[symbols]{lessthanorequalto@$\leq$},
    and the constant function $[m \mapsto 0]$ is denoted by $\functionThatIsIdenticalToZero$\graffito{zero function $\functionThatIsIdenticalToZero = [m \mapsto 0]$}\index[symbols]{0 nullity function@$\functionThatIsIdenticalToZero$}. 
  \end{definition}



  \begin{definition}
    Let $A$ be a subset of $M$. The map
    \begin{align*}
      \indicatorFunction_A \from M &\to \setOf{0,1},\index[symbols]{1doublestrikeA@$\indicatorFunction_A$}\\
      m &\mapsto \begin{dcases*}
                   1, & if $m \in A$,\\
                   0, & if $m \notin A$,
                 \end{dcases*}
    \end{align*}
    is called \define{indicator function of $A$ on $M$}\graffito{indicator function $\indicatorFunction_A$ of $A$ on $M$}.
  \end{definition}

  \begin{definition} 
  \label{definition:mean}
    Let $\nu \from \boundedFunctionsOn(M) \to \R$ be a map. It is called
    \begin{aenumerate}
      \item \define{normalised}\graffito{normalised} if and only if $\nu(\indicatorFunction_M) = 1$;
      \item \define{non-negativity preserving}\graffito{non-negativity preserving} if and only if
            \begin{equation*}
              \ForEach f \in \boundedFunctionsOn(M) \Holds (f \geq \functionThatIsIdenticalToZero \implies \nu(f) \geq 0); \mathnote{non-negativity preserving}
            \end{equation*}
      \item \define{mean on $M$}\graffito{mean $\nu$ on $M$} if and only if it is linear, normalised, and non-negativity preserving. 
    \end{aenumerate}
    The set of all means on $M$ is denoted by $\meansOn(M)$\graffito{set $\meansOn(M)$}\index[symbols]{MSM@$\meansOn(M)$}. 
  \end{definition}

  %

  \begin{definition} 
    Let $\Psi$ be a map from $\boundedFunctionsOn(M)$ to $\boundedFunctionsOn(M)$. It is called \define{non-negativity preserving}\graffito{non-negativity preserving} if and only if 
    \begin{equation*}
      \ForEach f \in \boundedFunctionsOn(M) \Holds (f \geq \functionThatIsIdenticalToZero \implies \Psi(f) \geq 0). \qedhere
    \end{equation*}
  \end{definition}

  \begin{lemma}
  \label{lemma:liberation-preimage}
    Let $G_0$ be finite, let $A$ be a finite subset of $M$, and let $\mathfrak{g}$ be an element of $G \modulo G_0$. Then, $\cardinalityOf{(\blank \isSemiActedUponBy \mathfrak{g})^{-1}(A)} \leq \cardinalityOf{G_0} \cdot \cardinalityOf{A}$.
  \end{lemma}

  \begin{proof}
    Let $a \in A$ such that $(\blank \isSemiActedUponBy \mathfrak{g})^{-1}(a) \neq \emptyset$. There are $m$ and $m' \in M$ such that $G_{m_0,m} = \mathfrak{g}$ and $m' \isSemiActedUponBy \mathfrak{g} = a$. For each $m'' \in M$, we have $m'' \isSemiActedUponBy \mathfrak{g} = g_{m_0, m''} \actsOnPoint m$ and hence 
    \begin{align*} 
      m'' \isSemiActedUponBy \mathfrak{g} = a
      &\ifAndOnlyIf m'' \isSemiActedUponBy \mathfrak{g} = m' \isSemiActedUponBy \mathfrak{g}\\
      &\ifAndOnlyIf g_{m_0, m'}^{-1} g_{m_0, m''} \actsOnPoint m = m\\
      &\ifAndOnlyIf g_{m_0, m'}^{-1} g_{m_0, m''} \in G_m\\
      &\ifAndOnlyIf g_{m_0, m''} \in g_{m_0, m'} G_m.
    \end{align*}
    Moreover, for each $m''$ and each $m''' \in M$ with $m'' \neq m'''$, we have $g_{m_0, m''} \neq g_{m_0, m'''}$. 
    Thus,
    \begin{align*}
      \cardinalityOf{(\blank \isSemiActedUponBy \mathfrak{g})^{-1}(a)}
      &=    \cardinalityOf{\setOf{m'' \in M \suchThat m'' \isSemiActedUponBy \mathfrak{g} = a}}\\
      &=    \cardinalityOf{\setOf{m'' \in M \suchThat g_{m_0, m''} \in g_{m_0, m'} G_m}}\\
      &\leq \cardinalityOf{g_{m_0, m'} G_m}\\
      &=    \cardinalityOf{G_m}\\
      &=    \cardinalityOf{G_0}.
    \end{align*}
    Therefore, because $(\blank \isSemiActedUponBy \mathfrak{g})^{-1}(A) = \bigcup_{a \in A} (\blank \isSemiActedUponBy \mathfrak{g})^{-1}(a)$, we have $\cardinalityOf{(\blank \isSemiActedUponBy \mathfrak{g})^{-1}(A)} \leq \cardinalityOf{G_0} \cdot \cardinalityOf{A}$.
  \end{proof}

  \begin{definition}
    The group $G$ acts on $\boundedFunctionsOn(M)$ on the left by
    \begin{align*}
      \actsOnBoundedFunction \from G \times \boundedFunctionsOn(M) &\to \boundedFunctionsOn(M), \mathnote{left group action $\actsOnBoundedFunction$ of $G$ on $\boundedFunctionsOn(M)$}\index[symbols]{Vdash@$\actsOnBoundedFunction$}\\
      (g, f) &\mapsto [m \mapsto f(g^{-1} \actsOnPoint m)]. \qedhere
    \end{align*}
  \end{definition}


  \begin{lemma}
  \label{lemma:kind-of-semi-action-of-quotient-on-boundeds}
    Let $G_0$ be finite. The quotient set $G \modulo G_0$ kind of semi-acts on $\boundedFunctionsOn(M)$ on the right by\index[symbols]{dashV@$\boundedFunctionIsKindOfSemiActedUponBy$}
    \begin{align*}
      \boundedFunctionIsKindOfSemiActedUponBy \from \boundedFunctionsOn(M) \times G \modulo G_0 &\to \boundedFunctionsOn(M), \mathnote{kind of right quotient set semi-action $\boundedFunctionIsKindOfSemiActedUponBy$ of $G \modulo G_0$ on $\boundedFunctionsOn(M)$}\\
                                         (f, \mathfrak{g}) &\mapsto [m \mapsto \sum_{m' \in (\blank \isSemiActedUponBy \mathfrak{g})^{-1}(m)} f(m')],
    \end{align*}
    such that, for each tuple $(f, \mathfrak{g}) \in \boundedFunctionsOn(M) \times G \modulo G_0$, we have $\normOf{f \boundedFunctionIsKindOfSemiActedUponBy \mathfrak{g}}_\infty \leq \cardinalityOf{G_0} \cdot \normOf{f}_\infty$.
  \end{lemma}

  \begin{proof}
    Let $\mathfrak{g} \in G \modulo G_0$. Furthermore, let $f \in \boundedFunctionsOn(M)$. Moreover, let $m \in M$. Because $G_0$ is finite, according to \cref{lemma:liberation-preimage}, we have $\cardinalityOf{(\blank \isSemiActedUponBy \mathfrak{g})^{-1}(m)} \leq \cardinalityOf{G_0} < \infty$. Hence, the sum in the definition of $\boundedFunctionIsKindOfSemiActedUponBy$ is finite. Furthermore,
    \begin{align*}
      \absoluteValueOf{(f \boundedFunctionIsKindOfSemiActedUponBy \mathfrak{g})(m)}
      &\leq \sum_{m' \in (\blank \isSemiActedUponBy \mathfrak{g})^{-1}(m)} \absoluteValueOf{f(m')}\\
      &\leq \parens*{\sum_{m' \in (\blank \isSemiActedUponBy \mathfrak{g})^{-1}(m)} 1} \cdot \normOf{f}_\infty\\
      &=    \cardinalityOf{(\blank \isSemiActedUponBy \mathfrak{g})^{-1}(m)} \cdot \normOf{f}_\infty\\
      &\leq \cardinalityOf{G_0} \cdot \normOf{f}_\infty.
    \end{align*} 
    Therefore, $f \boundedFunctionIsKindOfSemiActedUponBy \mathfrak{g} \in \boundedFunctionsOn(M)$, $\normOf{f \boundedFunctionIsKindOfSemiActedUponBy \mathfrak{g}}_\infty \leq \cardinalityOf{G_0} \cdot \normOf{f}_\infty$, and $\boundedFunctionIsKindOfSemiActedUponBy$ is well-defined.
  \end{proof}

  \begin{remark}
  \label{remark:groups:boundedFunctionIsKindOfSemiActedUponBy}
    In the situation of \cref{remark:groups:measureIsKindOfSemiActedUponBy}, we have $\boundedFunctionIsKindOfSemiActedUponBy \from (f, g) \mapsto [m \mapsto f(m \cdot g^{-1})]$. Hence, $\boundedFunctionIsKindOfSemiActedUponBy$ is the right group action of $G$ on $\R^G$ as defined in paragraph~5 of section~4.3 in \cite{ceccherini-silberstein:coornaert:2010}. 
  \end{remark}

  \begin{lemma}
  \label{lemma:boundedFunctionIsKindOfSemiActedUponBy-linear-continuous-non-negativity-preserving}
    Let $G_0$ be finite and let $\mathfrak{g}$ be an element of $G \modulo G_0$. The map $\blank \boundedFunctionIsKindOfSemiActedUponBy \mathfrak{g}$ is linear, continuous, and non-negativity preserving.
  \end{lemma}

  \begin{proof}
    Linearity follows from linearity of summation, continuity follows from linearity and $\normOf{\blank \boundedFunctionIsKindOfSemiActedUponBy \mathfrak{g}}_\infty \leq \cardinalityOf{G_0} \cdot \normOf{\blank}_\infty$, and non-negativity preservation follows from non-negativity preservation of summation.  
  \end{proof}

  \begin{lemma}
  \label{lemma:means-are-in-dual-of-boundeds-and-have-norm-one}
    Let $\nu$ be a mean on $M$. Then, $\nu \in \boundedFunctionsOn(M)^*$ and $\normOf{\nu}_{\boundedFunctionsOn(M)^*} = 1$. In particular, $\nu$ is continuous.
  \end{lemma}

  \begin{proof}
    Compare proposition~4.1.7 in \cite{ceccherini-silberstein:coornaert:2010}.
  \end{proof}

  \begin{definition}
  \label{definition:left-group-action-on-means}
    The group $G$ acts on $\boundedFunctionsOn(M)^*$ on the left by\index[symbols]{VDash@$\actsOnMean$}
    \begin{align*}
      \actsOnMean \from G \times \boundedFunctionsOn(M)^* &\to \boundedFunctionsOn(M)^*, \mathnote{left group action $\actsOnMean$ of $G$ on $\boundedFunctionsOn(M)^*$}\\
      (g, \psi) &\mapsto [f \mapsto \psi(g^{-1} \actsOnBoundedFunction f)],
    \end{align*}
    such that $G \actsOnMean \meansOn(M) \subseteq \meansOn(M)$.
  \end{definition}


  \begin{definition} 
  \label{definition:kind-of-right-semi-action-on-means}
    Let $G_0$ be finite. The quotient set $G \modulo G_0$ kind of semi-acts on $\boundedFunctionsOn(M)^*$ on the right by\index[symbols]{DashV@$\meanIsKindOfSemiActedUponBy$}
    \begin{align*}
      \meanIsKindOfSemiActedUponBy \from \boundedFunctionsOn(M)^* \times G \modulo G_0 &\to \boundedFunctionsOn(M)^*, \mathnote{kind of right quotient set semi-action $\meanIsKindOfSemiActedUponBy$ of $G \modulo G_0$ on $\boundedFunctionsOn(M)^*$}\\
      (\psi, \mathfrak{g}) &\mapsto [f \mapsto \psi(f \boundedFunctionIsKindOfSemiActedUponBy \mathfrak{g})]. \qedhere
    \end{align*}
  \end{definition}

  \begin{proof}
    Let $\psi \in \boundedFunctionsOn(M)^*$ and let $\mathfrak{g} \in G \modulo G_0$. Then, $\psi \meanIsKindOfSemiActedUponBy \mathfrak{g} = \psi \after (\blank \boundedFunctionIsKindOfSemiActedUponBy \mathfrak{g})$. Because $\psi$ and $\blank \boundedFunctionIsKindOfSemiActedUponBy \mathfrak{g}$ are linear and continuous, so is $\psi \meanIsKindOfSemiActedUponBy \mathfrak{g}$.
  \end{proof}

  \begin{definition}
  \label{definition:invariance-of-means}
    Let $G_0$ be finite and let $\psi$ be an element of $\boundedFunctionsOn(M)^*$. It is called \defineX{$\meanIsKindOfSemiActedUponBy$-invariant}{invariant}\graffito{$\meanIsKindOfSemiActedUponBy$-invariant} if and only if, for each element $\mathfrak{g} \in G \modulo G_0$ and each function $f \in \boundedFunctionsOn(M)$, we have $(\psi \meanIsKindOfSemiActedUponBy \mathfrak{g})(f) = \psi(f)$.
  \end{definition}



  \begin{remark}
  \label{remark:groups:meanIsKindOfSemiActedUponBy} 
    In the situation of \cref{remark:groups:boundedFunctionIsKindOfSemiActedUponBy}, we have $\meanIsKindOfSemiActedUponBy \from (\psi, g) \mapsto [f \mapsto \psi(f \boundedFunctionIsKindOfSemiActedUponBy g]$. Except for $g$ not being inverted, this is the right group action of $G$ on $\boundedFunctionsOn(G)^*$ as defined in paragraph~6 of section~4.3 in \cite{ceccherini-silberstein:coornaert:2010}. Hence, being $\meanIsKindOfSemiActedUponBy$-invariant is the same as being right-invariant as defined in paragraph~3 of section~4.4 in \cite{ceccherini-silberstein:coornaert:2010}. 
  \end{remark}

  \begin{remark}[Finite]
    If the set $M$ is finite, then the map 
    \begin{align*}
      \nu \from \boundedFunctionsOn(M) &\to \R,\\
      f &\mapsto \frac{1}{\cardinalityOf{M}} \sum_{m \in M} f(m) \quad (= \sum_{r \in \R} r \cdot \mu({f^{-1}(r)}), 
    \end{align*}
    is a $\actsOnMean$-invariant and $\meanIsKindOfSemiActedUponBy$-invariant mean on $M$, where $\mu$ is the finitely additive probability measure of \cref{remark:finite:measure}. Note that $\nu$ is the integral on $\boundedFunctionsOn(M)$ induced by $\mu$, that, for each function $f \in \boundedFunctionsOn(M)$, we have $f = \sum_{r \in \R} r \cdot \indicatorFunction_{f^{-1}(r)}$, and that, for each subset $A$ of $M$, we have $\mu(A) = \nu(\indicatorFunction_A)$. 
  \end{remark}

  \begin{proof} 
    The map $\nu$ is linear, normalised and non-negativity preserving. Thus, it is a mean on $M$. Moreover, for each element $g \in G$ and each function $f \in \boundedFunctionsOn(M)$, because the map $g \actsOnPoint \blank$ is bijective,
    \begin{align*}
      (g \actsOnMean \nu)(f)
      &= \nu(g^{-1} \actsOnBoundedFunction f)\\
      &= \frac{1}{\cardinalityOf{M}} \sum_{m \in M} (g^{-1} \actsOnBoundedFunction f)(m)\\
      &= \frac{1}{\cardinalityOf{M}} \sum_{m \in M} f(g \actsOnPoint m)\\
      &= \frac{1}{\cardinalityOf{M}} \sum_{m \in M} f(m)\\
      &= \nu(f).
    \end{align*}
    And, for each element $\mathfrak{g} \in G \modulo G_0$ and each function $f \in \boundedFunctionsOn(M)$, because $M = \bigDisjointUnionOf_{m \in M} (\blank \isSemiActedUponBy \mathfrak{g})^{-1}(m)$,
    \begin{align*} 
      (\nu \meanIsKindOfSemiActedUponBy \mathfrak{g})(f)
      &= \nu(f \boundedFunctionIsKindOfSemiActedUponBy \mathfrak{g})\\
      &= \frac{1}{\cardinalityOf{M}} \sum_{m \in M} (f \boundedFunctionIsKindOfSemiActedUponBy \mathfrak{g})(m)\\
      &= \frac{1}{\cardinalityOf{M}} \sum_{m \in M} \sum_{m' \in (\blank \isSemiActedUponBy \mathfrak{g})^{-1}(m)} f(m')\\
      &= \frac{1}{\cardinalityOf{M}} \sum_{m \in M} f(m)\\
      &= \nu(f).
    \end{align*}
    Hence, the mean $\nu$ is $\actsOnMean$-invariant and $\meanIsKindOfSemiActedUponBy$-invariant.

    Alternatively, we can deduce that $\nu$ is $\meanIsKindOfSemiActedUponBy$-invariant from the $\measureIsKindOfSemiActedUponBy$-semi-invariance of $\mu$ in a way that generalises to the non-finite case with the axiom of choice (compare \cref{lemma:mean-invariant-if-and-only-if-semi-invariant,theorem:mean-characterisation-of-right-amenable}). Let $\mathfrak{g}$ be an element of $G \modulo G_0$ and let $f$ be a function of $\boundedFunctionsOn(M)$. Because $\cardinalityOf{(\blank \isSemiActedUponBy \mathfrak{g})^{-1}(m)} \leq \cardinalityOf{G_0}$, for $m \in M$, (in other words, the map $\blank \isSemiActedUponBy \mathfrak{g}$ is $(\leq \cardinalityOf{G_0})$-to-$1$), there is a partition $\family{M_k}_{k \in \setOf{1, 2, \dotsc, \cardinalityOf{G_0}}}$ of $M$ such that, for each index $k \in \setOf{1, 2, \dotsc, \cardinalityOf{G_0}}$, the map $\blank \isSemiActedUponBy \mathfrak{g}$ is injective on $M_k$. We have
    \begin{equation*}
      f = \sum_{r \in \R} r \cdot \indicatorFunction_{f^{-1}(r)}
        = \sum_{r \in \R} \sum_{k \in \setOf{1, 2, \dotsc, \cardinalityOf{G_0}}} r \cdot \indicatorFunction_{f^{-1}(r) \cap M_k}.
    \end{equation*}
    Therefore, because $\nu \meanIsKindOfSemiActedUponBy \mathfrak{g}$ is linear,
    \begin{equation*}
      (\nu \meanIsKindOfSemiActedUponBy \mathfrak{g})(f)
      = \sum_{r \in \R} \sum_{k \in \setOf{1, 2, \dotsc, \cardinalityOf{G_0}}} r \cdot (\nu \meanIsKindOfSemiActedUponBy \mathfrak{g})(\indicatorFunction_{f^{-1}(r) \cap M_k}).
    \end{equation*}
    Let $r$ be a real number and let $k$ be an index of $\setOf{1, 2, \dotsc, \cardinalityOf{G_0}}$. Then, because $\blank \isSemiActedUponBy \mathfrak{g}$ is injective on $f^{-1}(r) \cap M_k$, we have $\indicatorFunction_{f^{-1}(r) \cap M_k} \boundedFunctionIsKindOfSemiActedUponBy \mathfrak{g} = \indicatorFunction_{f^{-1}(r) \cap M_k \isSemiActedUponBy \mathfrak{g}}$. Hence, by definition of $\nu$ and because $\mu$ is $\measureIsKindOfSemiActedUponBy$-semi-invariant,
    \begin{align*}
      (\nu \meanIsKindOfSemiActedUponBy \mathfrak{g})(\indicatorFunction_{f^{-1}(r) \cap M_k})
      &= \nu(\indicatorFunction_{f^{-1}(r) \cap M_k \isSemiActedUponBy \mathfrak{g}})\\
      &= \mu(f^{-1}(r) \cap M_k \isSemiActedUponBy \mathfrak{g})\\
      &= \mu(f^{-1}(r) \cap M_k)\\
      &= \nu(\indicatorFunction_{f^{-1}(r) \cap M_k}).
    \end{align*}
    Therefore, with the equation above, because $\nu$ is linear,
    \begin{equation*}
      (\nu \meanIsKindOfSemiActedUponBy \mathfrak{g})(f)
      = \sum_{r \in \R} \sum_{k \in \setOf{1, 2, \dotsc, \cardinalityOf{G_0}}} r \cdot \nu(\indicatorFunction_{f^{-1}(r) \cap M_k})
      = \nu(f). \qedhere
    \end{equation*}
  \end{proof}

  \begin{remark}[Infinite]
    If the set $M$ is infinite, then, for each finite subset $F$ of $M$, the map
    \begin{align*}
      \nu_F \from \boundedFunctionsOn(M) &\to \R,\\
      f &\mapsto \frac{1}{\cardinalityOf{F}} \sum_{m \in F} f(m) \quad (= \sum_{r \in \R} r \cdot \mu_F({f^{-1}(r)}), 
    \end{align*}
    is a mean on $M$ that is neither $\actsOnMean$-invariant nor $\meanIsKindOfSemiActedUponBy$-invariant, where $\mu_F$ is the finitely additive probability measure of \cref{remark:infinite:measure}. Nevertheless, for certain cell spaces over $M$, there is a net $\net{F_i}_{i \in I}$ of finite subsets of $M$ such that the net $\net{\nu_{F_i}}_{i \in I}$ converges in a weak sense to a $\actsOnMean$-invariant or $\meanIsKindOfSemiActedUponBy$-semi-invariant mean on $M$. However, even on $\Z$, there is no explicit formula for that mean. This construction is used in the subproof \textsc{\ref{item:Tarski-Folner:not-right-amenable} implies \ref{item:Tarski-Folner:no-Folner-net}} of \cref{theorem:Tarski-Folner}.
  \end{remark}

  %

  \begin{theorem}
  \label{theorem:means-versus-measures}
    The map
    \begin{align*}
      \Phi \from \meansOn(M) &\to \probabilityMeasuresOn(M), \mathnote{map $\Phi$ from $\meansOn(M)$ to $\probabilityMeasuresOn(M)$}\index[symbols]{Phi@$\Phi$}\\
      \nu &\mapsto [A \mapsto \nu(\indicatorFunction_A)],
    \end{align*}
    is bijective.
  \end{theorem} 

  \begin{proof}
    Compare theorem~4.1.8 in \cite{ceccherini-silberstein:coornaert:2010}.
  \end{proof}

  \begin{theorem}
  \label{theorem:means-convex-and-compact-wrt-weak-star-topology}
    The set $\meansOn(M)$ is a convex and compact subset of $\boundedFunctionsOn(M)^*$ equipped with the weak-$*$ topology.
  \end{theorem}

  \begin{proof}
    Compare theorem~4.2.1 in \cite{ceccherini-silberstein:coornaert:2010}.
  \end{proof}

  \section{Right Amenability}
  \label{section:right-amenability}

  \paragraph{Contents.} In \cref{definition:right-amenable} we introduce the notion of right amenability using finitely additive probability measures. And in \cref{theorem:mean-characterisation-of-right-amenable} we characterise right amenability of cell spaces with finite stabilisers using means.

  \begin{definition} 
    Let $\ntuple{M, G, \actsOnPoint}$ be a left group set. It is called \graffito{left amenable}\define{left amenable}\index{amenable!left} if and only if there is a $\actsOnMeasure$-invariant finitely additive probability measure on $M$.
  \end{definition}

  \begin{definition}
  \label{definition:right-amenable}
    Let $\mathcal{M} = \ntuple{M, G, \actsOnPoint}$ be a left-ho\-mo\-ge\-neous space. It is called \define{right amenable}\graffito{right-a\-me\-na\-ble left-ho\-mo\-ge\-neous space}\index{amenable!right} if and only if there is a coordinate system $\mathcal{K} = \ntuple{m_0, \family{g_{m_0, m}}_{m \in M}}$ for $\mathcal{M}$ such that there is a $\measureIsKindOfSemiActedUponBy$-semi-invariant finitely additive probability measure on $M$, in which case the cell space $\mathcal{R} = \ntuple{\mathcal{M}, \mathcal{K}}$ is called \define{right amenable}\graffito{right-a\-me\-na\-ble cell space}\index{amenable!right}. 
  \end{definition}

  \begin{remark} 
    In the situation of \cref{remark:groups:measureIsKindOfSemiActedUponBy}, being right amenable is the same as being amenable as defined in definition~4.4.5 in \cite{ceccherini-silberstein:coornaert:2010}. 
  \end{remark}

  In the remainder of this section, let $\mathcal{R} = \ntuple{\ntuple{M, G, \actsOnPoint}, \ntuple{m_0, \family{g_{m_0, m}}_{m \in M}}}$ be a cell space such that the stabiliser $G_0$ of $m_0$ under $\actsOnPoint$ is finite. 


  The (kind of) right semi-actions $\isSemiActedUponBy$ and $\boundedFunctionIsKindOfSemiActedUponBy$ are compatible in the sense given in

  \begin{lemma}
  \label{lemma:indicator-function-of-A-liberation-n-is-indicator-function-of-A-acted-upon}
    Let $\mathfrak{g}$ be an element of $G \modulo G_0$ and let $A$ be a subset of $M$ such that the map $\blank \isSemiActedUponBy \mathfrak{g}$ is injective on $A$. Then, $\indicatorFunction_{A \isSemiActedUponBy \mathfrak{g}} = \indicatorFunction_A \boundedFunctionIsKindOfSemiActedUponBy \mathfrak{g}$.
  \end{lemma}

  \begin{proof}
    For each $m \in M$, because $\blank \isSemiActedUponBy \mathfrak{g}$ is injective on $A$,
    \begin{align*}
      \indicatorFunction_{A \isSemiActedUponBy \mathfrak{g}}(m)
      &= \begin{dcases*}
           1, &if $m \in A \isSemiActedUponBy \mathfrak{g}$,\\
           0, &otherwise,
         \end{dcases*}\\
      &= \cardinalityOf{\setOf{m' \in A \suchThat m' \isSemiActedUponBy \mathfrak{g} = m}}\\ 
      &= \sum_{m' \in (\blank \isSemiActedUponBy \mathfrak{g})^{-1}(m)} \indicatorFunction_A(m')\\
      &= (\indicatorFunction_A \boundedFunctionIsKindOfSemiActedUponBy \mathfrak{g})(m).
    \end{align*}
    In conclusion, $\indicatorFunction_{A \isSemiActedUponBy \mathfrak{g}} = \indicatorFunction_A \boundedFunctionIsKindOfSemiActedUponBy \mathfrak{g}$.
  \end{proof}



  Simple functions approximate bounded functions arbitrarily good as stated in

  \begin{lemma} 
  \label{lemma:simples-dense-in-boundeds}
    The vector space
    \begin{equation*}
      \simpleFunctionsOn(M) = \setOf{f \from M \to \R \suchThat f(M) \text{ is finite}} \quad (= \linearSpanOf\setOf{\indicatorFunction_A \suchThat A \subseteq M}) \mathnote{$\simpleFunctionsOn(M)$}\index[symbols]{Ecalligraphic@$\simpleFunctionsOn(M)$}
    \end{equation*}
    is dense in the Banach space $\ntuple{\boundedFunctionsOn(M), \normOf{\blank}_\infty}$.
  \end{lemma}

  \begin{proof}
    Compare lemma~4.1.9 in \cite{ceccherini-silberstein:coornaert:2010}.
  \end{proof}

  For a mean, being invariant at certain indicator functions is sufficient for being invariant everywhere, which is shown in

  \begin{lemma} 
  \label{lemma:mean-invariant-if-and-only-if-semi-invariant}
    Let $\psi$ be an element of $\boundedFunctionsOn(M)^*$ such that, for each element $\mathfrak{g} \in G \modulo G_0$ and each subset $A$ of $M$ such that the map $\blank \isSemiActedUponBy \mathfrak{g}$ is injective on $A$, we have $(\psi \meanIsKindOfSemiActedUponBy \mathfrak{g})(\indicatorFunction_A) = \psi(\indicatorFunction_A)$. The map $\psi$ is $\meanIsKindOfSemiActedUponBy$-invariant.
  \end{lemma}

  \begin{proof}
    Let $\mathfrak{g} \in G \modulo G_0$.

    First, let $A \subseteq M$. Moreover, let $m \in M$. According to \cref{lemma:liberation-preimage}, we have $k_m = \cardinalityOf{(\blank \isSemiActedUponBy \mathfrak{g})^{-1}(m)} \leq \cardinalityOf{G_0}$. Hence, there are pairwise distinct $m_{m,1}$, $m_{m,2}$, $\dotsc$, $m_{m,k_m} \in M$ such that $(\blank \isSemiActedUponBy \mathfrak{g})^{-1}(m) = \setOf{m_{m,1}, m_{m,2}, \dotsc, m_{m,k_m}}$. For each $i \in \setOf{1, 2, \dotsc, \cardinalityOf{G_0}}$, put
    \begin{equation*} 
      A_i = \setOf{m_{m,i} \suchThat m \in M, k_m \geq i} \cap A.
    \end{equation*} 
    Because, for each $m \in M$ and each $m' \in M$ such that $m \neq m'$, we have $(\blank \isSemiActedUponBy \mathfrak{g})^{-1}(m) \cap (\blank \isSemiActedUponBy \mathfrak{g})^{-1}(m') = \emptyset$, the sets $A_1$, $A_2$, $\dotsc$, $A_{\cardinalityOf{G_0}}$ are pairwise disjoint and the map $\blank \isSemiActedUponBy \mathfrak{g}$ is injective on each of these sets. Moreover, because $\bigcup_{m \in M} (\blank \isSemiActedUponBy \mathfrak{g})^{-1}(m) = M$, we have $\bigcup_{i = 1}^{\cardinalityOf{G_0}} A_i = A$. Therefore, $\indicatorFunction_A = \sum_{i = 1}^{\cardinalityOf{G_0}} \indicatorFunction_{A_i}$. Thus, because $\psi \meanIsKindOfSemiActedUponBy \mathfrak{g}$ and $\psi$ are linear,
    \begin{align*}
      (\psi \meanIsKindOfSemiActedUponBy \mathfrak{g})(\indicatorFunction_A)
      &= (\psi \meanIsKindOfSemiActedUponBy \mathfrak{g})\parens*{\sum_{i = 1}^{\cardinalityOf{G_0}} \indicatorFunction_{A_i}}\\
      &= \sum_{i = 1}^{\cardinalityOf{G_0}} (\psi \meanIsKindOfSemiActedUponBy \mathfrak{g})(\indicatorFunction_{A_i})\\
      &= \sum_{i = 1}^{\cardinalityOf{G_0}} \psi(\indicatorFunction_{A_i})\\
      &= \psi(\indicatorFunction_A).
    \end{align*}
    Therefore, $\psi \meanIsKindOfSemiActedUponBy \mathfrak{g} = \psi$ on the set of indicator functions.
  %
    Thus, because the indicator functions span $\simpleFunctionsOn(M)$, and $\psi \meanIsKindOfSemiActedUponBy \mathfrak{g}$ and $\psi$ are linear, $\psi \meanIsKindOfSemiActedUponBy \mathfrak{g} = \psi$ on $\simpleFunctionsOn(M)$. Hence, because $\simpleFunctionsOn(M)$ is dense in $\boundedFunctionsOn(M)$, and $\psi \meanIsKindOfSemiActedUponBy \mathfrak{g}$ and $\psi$ are continuous, $\psi \meanIsKindOfSemiActedUponBy \mathfrak{g} = \psi$ on $\boundedFunctionsOn(M)$.
    In conclusion, $\psi$ is $\meanIsKindOfSemiActedUponBy$-invariant.
  \end{proof}

  A characterisation of right amenability by means is given in

  \begin{theorem}
  \label{theorem:mean-characterisation-of-right-amenable}
    The cell space $\mathcal{R}$ is right amenable if and only if there is a $\meanIsKindOfSemiActedUponBy$-invariant mean on $M$.
  \end{theorem}

  \begin{proof} 
    Let $\Phi$ be the map in \cref{theorem:means-versus-measures}.

    First, let $\mathcal{R}$ be right amenable. Then, there is $\measureIsKindOfSemiActedUponBy$-semi-invariant finitely additive probability measure $\mu$ on $M$. Put $\nu = \Phi^{-1}(\mu)$. Then, for each $\mathfrak{g} \in G \modulo G_0$ and each $A \subseteq M$ such that $\blank \isSemiActedUponBy \mathfrak{g}$ is injective on $A$, according to \cref{lemma:indicator-function-of-A-liberation-n-is-indicator-function-of-A-acted-upon},
    \begin{equation*}
      (\nu \meanIsKindOfSemiActedUponBy \mathfrak{g})(\indicatorFunction_A)
      = \nu(\indicatorFunction_A \boundedFunctionIsKindOfSemiActedUponBy \mathfrak{g})
      = \nu(\indicatorFunction_{A \isSemiActedUponBy \mathfrak{g}})
      = \mu(A \isSemiActedUponBy \mathfrak{g}) 
      = \mu(A)\\
      = \nu(\indicatorFunction_A).
    \end{equation*}
    Thus, according to \cref{lemma:mean-invariant-if-and-only-if-semi-invariant}, the mean $\nu$ is $\meanIsKindOfSemiActedUponBy$-invariant.

    Secondly, let there be a $\meanIsKindOfSemiActedUponBy$-invariant mean $\nu$ on $M$. Put $\mu = \Phi(\nu)$.
    Then, for each $\mathfrak{g} \in G \modulo G_0$ and each $A \subseteq M$ such that $\blank \isSemiActedUponBy \mathfrak{g}$ is injective on $A$, according to \cref{lemma:indicator-function-of-A-liberation-n-is-indicator-function-of-A-acted-upon},
    \begin{align*}
      (\mu \measureIsKindOfSemiActedUponBy \mathfrak{g})(A)
      = \mu(A \isSemiActedUponBy \mathfrak{g})
      = \nu(\indicatorFunction_{A \isSemiActedUponBy \mathfrak{g}})
      = \nu(\indicatorFunction_A \boundedFunctionIsKindOfSemiActedUponBy \mathfrak{g})
      = \nu(\indicatorFunction_A)
      = \mu(A).
    \end{align*}
    Hence, $\mu$ is $\measureIsKindOfSemiActedUponBy$-semi-invariant.
  \end{proof}

  \section{Right Følner Nets}
  \label{section:right-Folner-nets}

  In this section, let $\mathcal{R} = \ntuple{\mathcal{M}, \mathcal{K}} = \ntuple{\ntuple{M, G, \actsOnPoint}, \ntuple{m_0, \family{g_{m_0, m}}_{m \in M}}}$ be a cell space.

  \paragraph{Contents.} In \cref{definition:right-Folner-net} we introduce right Følner nets, and in \cref{theorem:epsilon-characterisation-of-Folner-net} we give a necessary and sufficient condition for the existence of such nets. In the case that the stabiliser $G_0$ is finite, in \cref{theorem:reversed-sets-characterisation-of-Folner-nets} we characterise right Følner nets, in \cref{corollary:epsilon-characterisation-of-Folner-net} we give another necessary and sufficient condition for the existence of such nets, and in \cref{theorem:Folner-net-independent-of-coordinate-system} we show that being a right Følner net does not depend on the choice of coordinate system.

  \begin{definition} 
  \label{definition:right-Folner-net}
    Let $\net{F_i}_{i \in I}$ be a net in $\setOf{F \subseteq M \suchThat F \neq \emptyset, F \text{ finite}}$ indexed by $(I, \leq)$. It is called \define{right Følner net in $\mathcal{R}$ indexed by $(I, \leq)$}\graffito{right Følner net $\net{F_i}_{i \in I}$ in $\mathcal{R}$ indexed by $(I, \leq)$}\index{Følner net in $\mathcal{R}$ indexed by $(I, \leq)$@right Følner net in $\mathcal{R}$ indexed by $(I, \leq)$}\index{net!Følner} if and only if 
    \begin{equation*}
      \ForEach \mathfrak{g} \in G \modulo G_0 \Holds \lim_{i \in I} \frac{\cardinalityOf{F_i \smallsetminus (\blank \isSemiActedUponBy \mathfrak{g})^{-1}(F_i)}}{\cardinalityOf{F_i}} = 0. \qedhere
    \end{equation*} 
  \end{definition}

  \begin{remark} 
    In the situation of \cref{remark:groups:measureIsKindOfSemiActedUponBy}, for each element $g \in G$ and each index $i \in I$, we have $(\blank \isSemiActedUponBy g)^{-1}(F_i) = F_i \cdot g^{-1}$. Hence, right Følner nets in $\mathcal{R}$ are exactly right Følner nets for $G$ as defined in the first paragraph after definition~4.7.2 in \cite{ceccherini-silberstein:coornaert:2010}. 
  \end{remark}

  \begin{remark}[Finite]
    If the set $M$ is finite, then the constant sequence $\net{M}_{k \in \N_0}$ is a right Følner net in $\mathcal{R}$.
  \end{remark}


  The necessary and sufficient conditions for the existence of right Følner nets that are given in \cref{theorem:epsilon-characterisation-of-Folner-net,corollary:epsilon-characterisation-of-Folner-net} follow directly from

  \begin{lemma} 
  \label{lemma:characterisation-of-zero-convergent-net-that-depends-on-a-parameter}
    Let $V$ be a set, let $W$ be a set, and let $\Psi$ be a map from $V \times W$ to $\R$. There is a net $\net{v_i}_{i \in I}$ in $V$ indexed by $(I, \leq)$ such that
    \begin{equation}
    \label{equation:characterisation-of-zero-convergent-net-that-depends-on-a-parameter:zero-convergence}
      \ForEach w \in W \Holds \lim_{i \in I} \Psi(v_i, w) = 0,
    \end{equation}
    if and only if, for each finite subset $Q$ of $W$ and each positive real number $\varepsilon \in \R_{> 0}$, there is an element $v \in V$ such that
    \begin{equation}
    \label{equation:characterisation-of-zero-convergent-net-that-depends-on-a-parameter:inequality}
      \ForEach q \in Q \Holds \Psi(v, q) < \varepsilon. \qedhere
    \end{equation}
  \end{lemma}

  \begin{proof}
    First, let there be a net $\net{v_i}_{i \in I}$ in $V$ indexed by $(I, \leq)$ such that \cref{equation:characterisation-of-zero-convergent-net-that-depends-on-a-parameter:zero-convergence} holds. Furthermore, let $Q$ be a finite subset of $W$ and let $\varepsilon \in \R_{> 0}$. Because \cref{equation:characterisation-of-zero-convergent-net-that-depends-on-a-parameter:zero-convergence} holds, for each $q \in Q$, there is an $i_q \in I$ such that,
    \begin{equation*}
      \ForEach i \in I \Holds (i \geq i_q \implies \Psi(v_i, q) < \varepsilon).
    \end{equation*}
    Because $(I, \leq)$ is a directed set and $Q$ is finite, there is an $i \in I$ such that, for each $q \in Q$, we have $i \geq i_q$. Put $v = v_i$. Then, \cref{equation:characterisation-of-zero-convergent-net-that-depends-on-a-parameter:inequality} holds.

    Secondly, for each finite $Q \subseteq W$ and each $\varepsilon \in \R_{> 0}$, let there be a $v \in V$ such that \cref{equation:characterisation-of-zero-convergent-net-that-depends-on-a-parameter:inequality} holds. Furthermore, let
    \begin{equation*}
      I = \setOf{Q \subseteq W \suchThat Q \text{ is finite}} \times \R_{> 0}
    \end{equation*}
    and let $\leq$ be the preorder on $I$ given by
    \begin{equation*}
      \ForEach (Q, \varepsilon) \in I \ForEach (Q', \varepsilon') \in I \Holds
          (Q, \varepsilon) \leq (Q', \varepsilon') \ifAndOnlyIf Q \subseteq Q' \land \varepsilon \geq \varepsilon'.
    \end{equation*}
    For each $(Q, \varepsilon) \in I$ and each $(Q', \varepsilon') \in I$, the element $(Q \cup Q', \min(\varepsilon, \varepsilon'))$ of $I$ is an upper bound of $(Q, \varepsilon)$ and of $(Q', \varepsilon')$. Hence, $(I, \leq)$ is a directed set.

    By precondition, for each $i = (Q, \varepsilon) \in I$, there is a $v_i \in V$ such that
    \begin{equation*}
      \ForEach q \in Q \Holds \Psi(v_i, q) < \varepsilon.
    \end{equation*}
    Let $w \in W$ and let $\varepsilon_0 \in \R_{> 0}$. Put $i_0 = (\setOf{w}, \varepsilon_0)$. For each $i = (Q, \varepsilon) \in I$ with $i \geq i_0$, we have $w \in Q$ and $\varepsilon \leq \varepsilon_0$. Hence,
    \begin{equation*}
      \ForEach i \in I \Holds (i \geq i_0 \implies \Psi(v_i, w) < \varepsilon_0).
    \end{equation*}
    Therefore, $\net{v_i}_{i \in I}$ is a net in $V$ indexed by $(I, \leq)$ such that \cref{equation:characterisation-of-zero-convergent-net-that-depends-on-a-parameter:zero-convergence} holds.
  \end{proof}

  \begin{theorem} 
  \label{theorem:epsilon-characterisation-of-Folner-net} 
    There is a right Følner net in $\mathcal{R}$ if and only if, for each finite subset $E$ of $G \modulo G_0$ and each positive real number $\varepsilon \in \R_{> 0}$, there is a non-empty and finite subset $F$ of $M$ such that
    \begin{equation*}
      \ForEach e \in E \Holds \frac{\cardinalityOf{F \smallsetminus (\blank \isSemiActedUponBy e)^{-1}(F)}}{\cardinalityOf{F}} < \varepsilon. \qedhere
    \end{equation*}
  \end{theorem}

  \begin{proof}
    This is a direct consequence of \cref{lemma:characterisation-of-zero-convergent-net-that-depends-on-a-parameter} with
    \begin{align*}
      \Psi \from \setOf{F \subseteq M \suchThat F \neq \emptyset, F \text{ finite}} \times G \modulo G_0 &\to \R,\\
      (F, \mathfrak{g}) &\mapsto \frac{\cardinalityOf{F \smallsetminus (\blank \isSemiActedUponBy \mathfrak{g})^{-1}(F)}}{\cardinalityOf{F}}. \qedhere
    \end{align*}
  \end{proof}

  In the proof of \cref{theorem:reversed-sets-characterisation-of-Folner-nets}, the upper bound given in \cref{lemma:cardinality-of-inverse-image-of-liberation-minus-the-same-less-than-or-equal-to-whatever} is essential, which itself follows from the upper bound given in lemma \ref{lemma:liberation-preimage} and the inclusion given in \cref{lemma:liberation-by-n-yields-element-in-set-setminus-bigcup-liberation-by-inverse-times-n-prime}, which in turn follows from the equality given in \cref{lemma:rightsemiaction-can-be-undone}. 

  %


  \begin{lemma}
  \label{lemma:rightsemiaction-can-be-undone}
    Let $m$ be an element of $M$, and let $\mathfrak{g}$ be an element of $G \modulo G_0$. There is an element $g \in \mathfrak{g}$ such that
    \begin{equation*}
      \ForEach \mathfrak{g}' \in G \modulo G_0 \Holds (m \isSemiActedUponBy \mathfrak{g}) \isSemiActedUponBy \mathfrak{g}' = m \isSemiActedUponBy g \cdot \mathfrak{g}',
    \end{equation*}
    in particular, for said $g \in \mathfrak{g}$, we have $(m \isSemiActedUponBy \mathfrak{g}) \isSemiActedUponBy g^{-1} G_0 = m$.
  \end{lemma}

  \begin{proof}
    There is a $g \in G$ such that $g G_0 = \mathfrak{g}$. Moreover, because $\isSemiActedUponBy$ is a semi-action with defect $G_0$, there is a $g_0 \in G_0$ such that
    \begin{equation*}
      \ForEach \mathfrak{g}' \in G \modulo G_0 \Holds (m \isSemiActedUponBy g G_0) \isSemiActedUponBy \mathfrak{g}' = m \isSemiActedUponBy g \cdot (g_0^{-1} \cdot \mathfrak{g}').
    \end{equation*}
    Because $g \cdot (g_0^{-1} \cdot \mathfrak{g}') = g g_0^{-1} \cdot \mathfrak{g}'$ and $g g_0^{-1} \in \mathfrak{g}$, the statement holds.
  \end{proof}

  \begin{lemma}
  \label{lemma:liberation-by-n-yields-element-in-set-setminus-bigcup-liberation-by-inverse-times-n-prime}
    Let $A$ and $A'$ be two subsets of $M$, and let $\mathfrak{g}$ and $\mathfrak{g}'$ be two elements of $G \modulo G_0$. Then, for each element $m \in (\blank \isSemiActedUponBy \mathfrak{g})^{-1}(A) \smallsetminus (\blank \isSemiActedUponBy \mathfrak{g}')^{-1}(A')$,
    \begin{equation*}
      m \isSemiActedUponBy \mathfrak{g} \in \bigcup_{g \in \mathfrak{g}} A \smallsetminus (\blank \isSemiActedUponBy g^{-1} \cdot \mathfrak{g}')^{-1}(A')
    \end{equation*}
    and
    \begin{equation*}
      m \isSemiActedUponBy \mathfrak{g}' \in \bigcup_{g' \in \mathfrak{g}'} (\blank \isSemiActedUponBy (g')^{-1} \cdot \mathfrak{g})^{-1}(A) \smallsetminus A'. \qedhere
    \end{equation*}
  \end{lemma}

  \begin{proof}
    Let $m \in (\blank \isSemiActedUponBy \mathfrak{g})^{-1}(A) \smallsetminus (\blank \isSemiActedUponBy \mathfrak{g}')^{-1}(A')$. Then, $m \isSemiActedUponBy \mathfrak{g} \in A$ and $m \isSemiActedUponBy \mathfrak{g}' \notin A'$. According to \cref{lemma:rightsemiaction-can-be-undone}, there is a $g \in \mathfrak{g}$ and a $g' \in \mathfrak{g}'$ such that $(m \isSemiActedUponBy \mathfrak{g}) \isSemiActedUponBy g^{-1} \cdot \mathfrak{g}' = m \isSemiActedUponBy \mathfrak{g}' \notin A'$ and $(m \isSemiActedUponBy \mathfrak{g}') \isSemiActedUponBy (g')^{-1} \cdot \mathfrak{g} = m \isSemiActedUponBy \mathfrak{g} \in A$. Hence, $m \isSemiActedUponBy \mathfrak{g} \notin (\blank \isSemiActedUponBy g^{-1} \cdot \mathfrak{g}')^{-1}(A')$ and $m \isSemiActedUponBy \mathfrak{g}' \in (\blank \isSemiActedUponBy (g')^{-1} \cdot \mathfrak{g})^{-1}(A)$. Therefore, $m \isSemiActedUponBy \mathfrak{g} \in A \smallsetminus (\blank \isSemiActedUponBy g^{-1} \cdot \mathfrak{g}')^{-1}(A')$ and $m \isSemiActedUponBy \mathfrak{g}' \in (\blank \isSemiActedUponBy (g')^{-1} \cdot \mathfrak{g})^{-1}(A) \smallsetminus A'$. In conclusion, $m \isSemiActedUponBy \mathfrak{g} \in \bigcup_{g \in \mathfrak{g}} A \smallsetminus (\blank \isSemiActedUponBy g^{-1} \cdot \mathfrak{g}')^{-1}(A')$ and $m \isSemiActedUponBy \mathfrak{g}' \in \bigcup_{g' \in \mathfrak{g}'} (\blank \isSemiActedUponBy (g')^{-1} \cdot \mathfrak{g})^{-1}(A) \smallsetminus A'$.
  \end{proof}

  \begin{lemma}
  \label{lemma:cardinality-of-inverse-image-of-liberation-minus-the-same-less-than-or-equal-to-whatever}
    Let $G_0$ be finite, let $F$ and $F'$ be two finite subsets of $M$, and let $\mathfrak{g}$ and $\mathfrak{g}'$ be two elements of $G \modulo G_0$. Then,
    \begin{multline*}
      \cardinalityOf{(\blank \isSemiActedUponBy \mathfrak{g})^{-1}(F) \smallsetminus (\blank \isSemiActedUponBy \mathfrak{g}')^{-1}(F')}\\
      \leq
      \begin{dcases*} 
        \cardinalityOf{G_0}^2 \cdot \max_{g \in \mathfrak{g}} \cardinalityOf{F \smallsetminus (\blank \isSemiActedUponBy g^{-1} \cdot \mathfrak{g}')^{-1}(F')},\\
        \cardinalityOf{G_0}^2 \cdot \max_{g' \in \mathfrak{g}'} \cardinalityOf{(\blank \isSemiActedUponBy (g')^{-1} \cdot \mathfrak{g})^{-1}(F) \smallsetminus F'}. \qedhere
      \end{dcases*}
    \end{multline*}
  \end{lemma}

  \begin{proof}
    Put $A = (\blank \isSemiActedUponBy \mathfrak{g})^{-1}(F) \smallsetminus (\blank \isSemiActedUponBy \mathfrak{g}')^{-1}(F')$. For each $g \in \mathfrak{g}$, put $B_g = F \smallsetminus (\blank \isSemiActedUponBy g^{-1} \cdot \mathfrak{g}')^{-1}(F')$. For each $g' \in \mathfrak{g}'$, put $B_{g'}' = (\blank \isSemiActedUponBy (g')^{-1} \cdot \mathfrak{g})^{-1}(F) \smallsetminus F'$.

    According to \cref{lemma:liberation-by-n-yields-element-in-set-setminus-bigcup-liberation-by-inverse-times-n-prime}, the restrictions $(\blank \isSemiActedUponBy \mathfrak{g})\restrictedTo_{A \to \bigcup_{g \in \mathfrak{g}} B_g}$ and $(\blank \isSemiActedUponBy \mathfrak{g}')\restrictedTo_{A \to \bigcup_{g' \in \mathfrak{g}'} B_{g'}'}$ are well-defined. Moreover, for each $m \in M$, according to \cref{lemma:liberation-preimage}, we have $\cardinalityOf{(\blank \isSemiActedUponBy \mathfrak{g})^{-1}(m)} \leq \cardinalityOf{G_0}$ and $\cardinalityOf{(\blank \isSemiActedUponBy \mathfrak{g}')^{-1}(m)} \leq \cardinalityOf{G_0}$. Therefore, because $\cardinalityOf{\mathfrak{g}} = \cardinalityOf{G_0}$, 
    \begin{equation*}
      \cardinalityOf{A}
      \leq \cardinalityOf{G_0} \cdot \cardinalityOf{\bigcup_{g \in \mathfrak{g}} B_g}
      \leq \cardinalityOf{G_0} \cdot \sum_{g \in \mathfrak{g}} \cardinalityOf{B_g}
      \leq \cardinalityOf{G_0}^2 \cdot \max_{g \in \mathfrak{g}} \cardinalityOf{B_g}
    \end{equation*}
    and analogously
    \begin{equation*}
      \cardinalityOf{A} \leq \cardinalityOf{G_0}^2 \cdot \max_{g' \in \mathfrak{g}'} \cardinalityOf{B_{g'}'}. \qedhere
    \end{equation*}
  \end{proof}

  A characterisation of right Følner nets is given in

  \begin{theorem} 
  \label{theorem:reversed-sets-characterisation-of-Folner-nets}
    Let $G_0$ be finite and let $\net{F_i}_{i \in I}$ be a net in $\setOf{F \subseteq M \suchThat F \neq \emptyset, F \text{ finite}}$ indexed by $(I, \leq)$. The net $\net{F_i}_{i \in I}$ is a right Følner net in $\mathcal{R}$ if and only if
    \begin{equation}
    \label{equation:reversed-sets-characterisation-of-Folner-nets:condition}
      \ForEach \mathfrak{g} \in G \modulo G_0 \Holds \lim_{i \in I} \frac{\cardinalityOf{(\blank \isSemiActedUponBy \mathfrak{g})^{-1}(F_i) \smallsetminus F_i}}{\cardinalityOf{F_i}} = 0. \qedhere
    \end{equation}
  \end{theorem}

  \begin{proof}
    Let $\mathfrak{g} \in G \modulo G_0$. Furthermore, let $i \in I$. Because $F_i = (\blank \isSemiActedUponBy G_0)^{-1}(F_i)$, according to \cref{lemma:cardinality-of-inverse-image-of-liberation-minus-the-same-less-than-or-equal-to-whatever},
    \begin{equation*}
      \cardinalityOf{(\blank \isSemiActedUponBy \mathfrak{g})^{-1}(F_i) \smallsetminus F_i}
      \leq \cardinalityOf{G_0}^2 \cdot \max_{g \in \mathfrak{g}} \cardinalityOf{F_i \smallsetminus (\blank \isSemiActedUponBy g^{-1} G_0)^{-1}(F_i)}
    \end{equation*}
    and
    \begin{equation*}
      \cardinalityOf{F_i \smallsetminus (\blank \isSemiActedUponBy \mathfrak{g})^{-1}(F_i)}
      \leq \cardinalityOf{G_0}^2 \cdot \max_{g \in \mathfrak{g}} \cardinalityOf{(\blank \isSemiActedUponBy g^{-1} G_0)^{-1}(F_i) \smallsetminus F_i}.
    \end{equation*}
    Moreover, $\cardinalityOf{\mathfrak{g}} = \cardinalityOf{G_0} < \infty$. Therefore, if $\net{F_i}_{i \in I}$ is a right Følner net in $\mathcal{R}$, then
    \begin{equation*}
      \lim_{i \in I} \frac{\cardinalityOf{(\blank \isSemiActedUponBy \mathfrak{g})^{-1}(F_i) \smallsetminus F_i}}{\cardinalityOf{F_i}} = 0;
    \end{equation*}
    and, if \cref{equation:reversed-sets-characterisation-of-Folner-nets:condition} holds, then
    \begin{equation*}
      \lim_{i \in I} \frac{\cardinalityOf{F_i \smallsetminus (\blank \isSemiActedUponBy \mathfrak{g})^{-1}(F_i)}}{\cardinalityOf{F_i}} = 0.
    \end{equation*}
    In conclusion, $\net{F_i}_{i \in I}$ is a right Følner net in $\mathcal{R}$ if and only if \cref{equation:reversed-sets-characterisation-of-Folner-nets:condition} holds.
  %
  \end{proof}

  Necessary and sufficient conditions for the existence of right Følner nets are given in

  \begin{corollary} 
  \label{corollary:epsilon-characterisation-of-Folner-net}
    Let $G_0$ be finite. There is a right Følner net in $\mathcal{R}$ if and only if, for each finite subset $E$ of $G \modulo G_0$ and each positive real number $\varepsilon \in \R_{> 0}$, there is a non-empty and finite subset $F$ of $M$ such that
    \begin{equation*}
      \ForEach e \in E \Holds \frac{\cardinalityOf{(\blank \isSemiActedUponBy e)^{-1}(F) \smallsetminus F}}{\cardinalityOf{F}} < \varepsilon. \qedhere
    \end{equation*}
  \end{corollary}

  \begin{proof}
    This is a direct consequence of \cref{theorem:reversed-sets-characterisation-of-Folner-nets} and \cref{lemma:characterisation-of-zero-convergent-net-that-depends-on-a-parameter} with
    \begin{align*}
      \Psi \from \setOf{F \subseteq M \suchThat F \neq \emptyset, F \text{ finite}} \times G \modulo G_0 &\to \R,\\
      (F, \mathfrak{g}) &\mapsto \frac{\cardinalityOf{(\blank \isSemiActedUponBy \mathfrak{g})^{-1}(F) \smallsetminus F}}{\cardinalityOf{F}}. \qedhere
    \end{align*}
  \end{proof}

  That being a right Følner net does not depend on the choice of coordinate system is shown in

  \begin{theorem} 
  \label{theorem:Folner-net-independent-of-coordinate-system}
    Let $\mathcal{M} = \ntuple{M, G, \actsOnPoint}$ be a left-ho\-mo\-ge\-neous space with finite stabilisers, let $\mathcal{K} = \ntuple{m_0, \family{g_{m_0, m}}_{m \in M}}$ and $\mathcal{K}' = \ntuple{m_0', \family{g_{m_0', m}'}_{m \in M}}$ be two coordinate systems for $\mathcal{M}$, and let $\mathcal{F} = \net{F_i}_{i \in I}$ be a net in $\setOf{F \subseteq M \suchThat F \neq \emptyset, F \text{ finite}}$ indexed by $(I, \leq)$. The net $\mathcal{F}$ is a right Følner net in $\ntuple{\mathcal{M}, \mathcal{K}}$ if and only if it is one in $\ntuple{\mathcal{M}, \mathcal{K}'}$.
  \end{theorem}

  \begin{proof} 
    First, let $\mathcal{F}$ be a right Følner net in $\ntuple{\mathcal{M}, \mathcal{K}}$. Furthermore, let $g$ be an element of $G$ such that $g \actsOnPoint m_0 = m_0'$, let $\mathfrak{g}'$ be an element of $G \modulo G_0'$ and let $i$ be an index of $I$. Then, according to \cref{lemma:liberation-and-coordinate-systems},
    \begin{equation*}
      \ForEach m \in M \Exists g_0 \in G_0 \SuchThat m \isSemiActedUponBy' \mathfrak{g}' = m \isSemiActedUponBy g_0 \cdot (g^{-1} \conjugates \mathfrak{g}').
    \end{equation*}
    Thus,
    \begin{equation*}
      (\blank \isSemiActedUponBy' \mathfrak{g}')^{-1}(F_i)
      \subseteq \bigcup_{g_0 \in G_0} (\blank \isSemiActedUponBy g_0 \cdot (g^{-1} \conjugates \mathfrak{g}'))^{-1}(F_i),
    \end{equation*}
    in particular,
    \begin{equation*}
      (\blank \isSemiActedUponBy' \mathfrak{g}')^{-1}(F_i) \smallsetminus F_i
      \subseteq \bigcup_{g_0 \in G_0} (\blank \isSemiActedUponBy g_0 \cdot (g^{-1} \conjugates \mathfrak{g}'))^{-1}(F_i) \smallsetminus F_i.
    \end{equation*}
    Hence, because the stabiliser $G_0$ is finite,
    \begin{equation*}
      \cardinalityOf{(\blank \isSemiActedUponBy' \mathfrak{g}')^{-1}(F_i) \smallsetminus F_i}
      \leq \sum_{g_0 \in G_0} \cardinalityOf{(\blank \isSemiActedUponBy g_0 \cdot (g^{-1} \conjugates \mathfrak{g}'))^{-1}(F_i) \smallsetminus F_i}.
    \end{equation*}
    Therefore, because $\mathcal{F}$ is a right Følner net in $\ntuple{\mathcal{M}, \mathcal{K}}$,
    \begin{equation*}
      \lim_{i \in I} \frac{\cardinalityOf{(\blank \isSemiActedUponBy' \mathfrak{g}')^{-1}(F_i)}}{\cardinalityOf{F_i}} = 0.
    \end{equation*}
    In conclusion, the net $\mathcal{F}$ is a right Følner net in $\ntuple{\mathcal{M}, \mathcal{K}'}$.

    Secondly, let $\mathcal{F}$ be a right Følner net in $\ntuple{\mathcal{M}, \mathcal{K}'}$. It follows as above that $\mathcal{F}$ is a right Følner net in $\ntuple{\mathcal{M}, \mathcal{K}}$. \qedhere
  \end{proof}

  \begin{example}[Lattice]
  \label{example:lattice}
    Let $M$ be the two-dimensional integer lattice $\Z^2$, let $S_2$ be the symmetric group on $\setOf{1, 2}$, let $V$ be the multiplicative group $\setOf{-1, 1}$, let $\varphi$ be the group homomorphism $S_2 \to \automorphismsOf(V^2)$, $\pi \mapsto [(v_1, v_2) \mapsto (v_{\pi(1)}, v_{\pi(2)})]$, let $V^2 \rtimes_\varphi S_2$ be the outer semi-direct product of $S_2$ acting on $V^2$ by $\varphi$, let $\psi$ be the group homomorphism $V^2 \rtimes_\varphi S_2 \to \automorphismsOf(\Z^2)$, $((v_1, v_2), \pi) \mapsto [(t_1, t_2) \mapsto (v_1 \cdot t_{\pi(1)}, v_2 \cdot t_{\pi(2)})]$, let $G = \Z^2 \rtimes_\psi (V^2 \rtimes_\varphi S_2)$ be the outer semi-direct product of $V^2 \rtimes_\varphi S_2$ acting on $\Z^2$ by $\psi$, and let $\actsOnPoint$ be the transitive left group action of $G$ on $M$ by $(((t_1, t_2), ((v_1, v_2), \pi)), (z_1, z_2)) \mapsto (t_1 + v_1 \cdot z_{\pi(1)}, t_2 + v_2 \cdot z_{\pi(2)})$. The triple $\mathcal{M} = \ntuple{M, G, \actsOnPoint}$ is a left-ho\-mo\-ge\-neous space. 

    The group $G$ encodes the symmetries of the lattice $M$, more precisely, the component $\Z^2$ encodes the translational symmetries and the component $V^2 \times S_2$ the reflectional and rotational ones that stabilise the origin. For example, the element $((4, 2), ((-1, 1), \identityMap))$ encodes the reflection about the $x$-axis, followed by a translation by $(4, 2)$; the element $((0, 0), ((1, 1), (1\;2)))$ encodes the reflection about the line through the origin of slope $1$; and the element $((0, 0), ((-1, 1), (1\;2)))$ encodes the anticlockwise rotation about the origin through $90\degree$, where the permutation $(1\;2)$, written in cycle notation, is the transposition that swaps $1$ and $2$. The symmetry group of $M$ is the group $\setOf{g \actsOnPoint \blank \suchThat g \in G}$ under composition, which is isomorphic to $G$. 

    Let $m_0$ be the origin $(0, 0)$ and, for each point $m \in M$, let $g_{m_0, m}$ be the translation $(m, ((1, 1), \identityMap))$. The tuple $\mathcal{K} = \ntuple{m_0, \family{g_{m_0, m}}_{m \in M}}$ is a coordinate system for $\mathcal{M}$. The stabiliser $G_0$ of $m_0$ under $\actsOnPoint$ is the set $\setOf{(0, 0)} \times (V^2 \times S_2)$, the quotient group $G \modulo G_0$ is isomorphic to the group $\Z^2$ by $(t, ((v_1, v_2), \pi)) G_0 \mapsto t$, and, under this isomorphism, the right quotient group semi-action $\isSemiActedUponBy$ of $G \modulo G_0$ on $M$ is the map $M \times \Z^2 \to M$, $(z, t) \mapsto z + t$.

    The cell space $\mathcal{R} = \ntuple{\mathcal{M}, \mathcal{K}}$ is finitely and symmetrically right generated by $S = \setOf{(-1, 0), (0, -1), (0, 1), (1, 0)}$. The $S$-metric $\distanceOf$, the $S$-length $\lengthOf{\blank} = \distanceOf(m_0, \blank)$, the $S$-balls $\ball(\dotso)$, and the $S$-spheres $\sphere(\dotso)$ are restrictions of the corresponding notions on $\R^2$ with respect to the taxicab metric on $\R^2$. The balls $\ball_{\R^2}(\dotso)$ and spheres $\sphere_{\R^2}(\dotso)$ induced by the taxicab metric on $\R^2$ are diamonds, that is, filled and unfilled squares with sides oriented at $45\degree$ to the coordinate axes. 

    The former notions will be rigorously introduced in \cref{chapter:growth}. In the present situation, the $S$-metric is the graph metric on the $S$-Cayley graph of $\Z^2$, the $S$-length is the distance to the origin, for each cell $m \in M$ and each integer $\rho$, the $S$-ball $\ball(m, \rho)$ and $S$-sphere $\sphere(m, \rho)$ of radius $\rho$ centred at $m$ are the sets of points whose distances to $m$ are not greater than and, respectively, equal to $\rho$, and we write $\ball(\rho)$ and $\sphere(\rho)$ for the balls and spheres centred at the origin. 

    Let $t$ be an element of $\Z^2$. For each non-negative integer $\rho$, the preimage $(\blank \isSemiActedUponBy t)^{-1}(\ball(\rho))$ is the ball $\ball(-t, \rho)$ and the complement $\ball(\rho) \smallsetminus (\blank \isSemiActedUponBy t)^{-1}(\ball(\rho))$ is included in the thickened sphere $\ball(\rho) \smallsetminus \ball(\rho - \lengthOf{t})$ (see \cref{figure:lattice:Folner-original} for the case that $t = (-1, 0)$). Because the size of the latter grows linearly in $\rho$, the sequence $\sequence{\ball(\rho) \smallsetminus (\blank \isSemiActedUponBy t)^{-1}(\ball(\rho))}_{\rho \in \N_0}$ grows linearly in size (in the case that $t = (-1, 0)$, it grows in size like $\sequence{2 \rho + 1}_{\rho \in \N_0}$). Moreover, the sequence $\sequence{\ball(\rho)}_{\rho \in \N_0}$ grows polynomially in size, more precisely, it grows in size like $\sequence{2\rho (\rho + 1) + 1}_{\rho \in \N_0}$. Hence, the quotient $\cardinalityOf{\ball(\rho) \smallsetminus (\blank \isSemiActedUponBy t)^{-1}(\ball(\rho))} / \cardinalityOf{\ball(\rho)}$ converges to $0$ as $\rho$ tends to $\infty$. Therefore, the sequence $\sequence{\ball(\rho)}_{\rho \in \N_0}$ is a right Følner net in $\mathcal{R}$.
    \begin{figure}
      \myfloatalign
      \begin{wide}
        \figureLatticeFolnerOriginal
        \caption{In each subfigure, the whole space is $\R^2$, the grid points are elements of $M = \Z^2$, the grid points in the region enclosed by the diamond with solid border are the elements of $\ball(\rho)$, the grid points in the region enclosed by the diamond with dashed border are the elements of $(\blank \isSemiActedUponBy (-1, 0))^{-1}(\ball(\rho))$, and the dots are the elements of $\ball(\rho) \smallsetminus (\blank \isSemiActedUponBy (-1, 0))^{-1}(\ball(\rho))$, for the respective $\rho \in \setOf{1, 2, 3}$.}
        \label{figure:lattice:Folner-original}
      \end{wide}
    \end{figure}
  \end{example}

  \begin{example}[Tree]
  \label{example:tree}
    Let $M$ be the vertices of the uncoloured $\setOf{a, b,\allowbreak a^{-1}, b^{-1}}$-Cayley graph of the free group $F_2$ over $\setOf{a, b}$, where $a \neq b$, let $\varsigma$ be the group automorphism from $F_2$ to $F_2$ determined by $a \mapsto b$ and $b \mapsto a^{-1}$, let $R$ be the cyclic group $\setOf{0, 90, 180, 270}$ under addition modulo $360$, let $\varphi$ be the group homomorphism $R \to \automorphismsOf(F_2)$, $r \mapsto \varsigma^{r / 90}$, let $G = F_2 \rtimes_\varphi R$ be the outer semi-direct product of $R$ acting on $F_2$ by $\varphi$, and let $\actsOnPoint$ be the transitive left group action of $G$ on $M$ by $((f, r), m) \mapsto f \cdot \varphi(r)(m)$. The triple $\mathcal{M} = \ntuple{M, G, \actsOnPoint}$ is a left-ho\-mo\-ge\-neous space.

    The group $G$ encodes some graph automorphisms of $M$, more precisely, the component $F_2$ encodes the translational automorphisms and the component $R$ the rotational ones that stabilise the origin. For example, the element $(a b^{-1}, 90)$ encodes the anticlockwise rotation about the origin through $90\degree$, followed by a translation by $a b^{-1}$, which is the anticlockwise rotation about $a$ through $90\degree$; see \cref{figure:tree:translations-rotations} for further examples. In general, for each vertex $m \in M$ and each angle $r \in R$, the anticlockwise rotation about $m$ through $r\degree$, is the graph automorphism $m \cdot \varphi(r)(m^{-1} \cdot \blank) = m \cdot (\varphi(r)(m))^{-1} \cdot \varphi(r)(\blank)$, which is encoded by $(m, 0) \cdot (e_{F_2}, r) \cdot (m^{-1}, 0) = (m \cdot (\varphi(r)(m))^{-1}, r)$. The map $g \mapsto g \actsOnPoint \blank$ embeds the group $G$ into the graph-automorphism group of $M$. 
    \begin{figure}
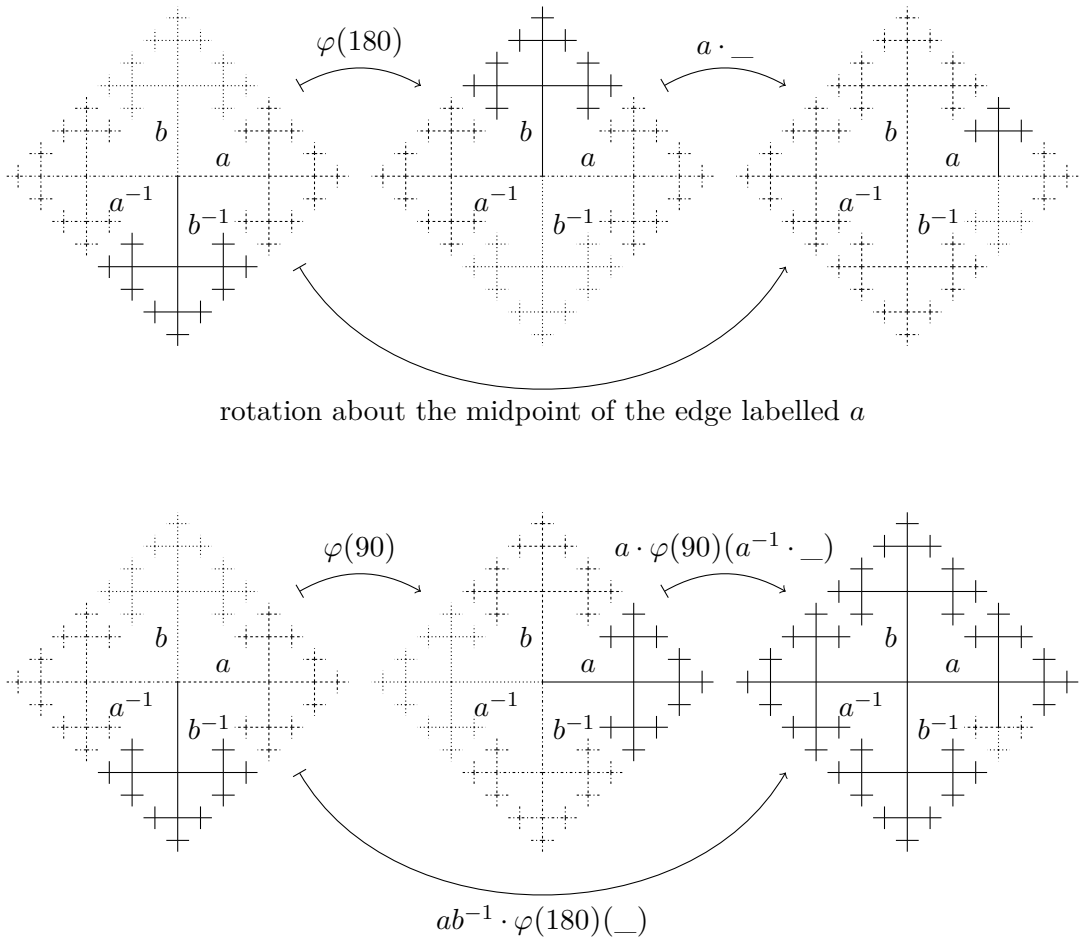

      \myfloatalign
      \figureTreeTranslationsRotations
      \begin{wide}
        \caption{The same part of the uncoloured $\setOf{a, b, a^{-1}, b^{-1}}$-Cayley graph of $F_2$ is depicted six times, each time the vertex in the centre is the neutral element $e_{F_2}$, the arrows of the form $\mapsto$ are translations, rotations, or their compositions, and the different line styles illustrate which vertices are mapped to which vertices.}
        \label{figure:tree:translations-rotations}
      \end{wide}
    \end{figure}

    Let $m_0$ be the neutral element $e_{F_2}$ of $F_2$ and, for each vertex $m \in M$, let $g_{m_0, m}$ be the translation $(m, 0)$. The tuple $\mathcal{K} = \ntuple{m_0, \family{g_{m_0, m}}_{m \in M}}$ is a coordinate system for $\mathcal{M}$. The stabiliser $G_0$ of $m_0$ under $\actsOnPoint$ is the set $\setOf{e_{F_2}} \times R$, the quotient group $G \modulo G_0$ is isomorphic to the group $F_2$ by $(f, r) G_0 \mapsto f$, and, under this isomorphism, the right quotient group semi-action $\isSemiActedUponBy$ of $G \modulo G_0$ on $M$ is the map $M \times F_2 \to M$, $(m, f) \mapsto m \cdot f$.

    The cell space $\mathcal{R} = \ntuple{\mathcal{M}, \mathcal{K}}$ is finitely and symmetrically right generated by $S = \setOf{a, b, a^{-1}, b^{-1}}$. The uncoloured $S$-Cayley graph of $\mathcal{R}$ is equal to the uncoloured $S$-Cayley graph of $F_2$. Hence, the $S$-metric $\distanceOf$ is identical to the $S$-word metric on $F_2$ and the $S$-length $\lengthOf{\blank} = \distanceOf(m_0, \blank)$ is identical to the $S$-word norm on $F_2$. And, for each element $m \in M$ and each integer $\rho$, the sets $\ball(m, \rho) = \setOf{m' \in M \suchThat \distanceOf(m, m') \leq \rho}$ and $\sphere(m, \rho) = \setOf{m' \in M \suchThat \distanceOf(m, m') = \rho}$ are the ball and sphere of radius $\rho$ centred at $m$, and we denote $\ball(m_0, \rho)$ by $\ball(\rho)$ and $\sphere(m_0, \rho)$ by $\sphere(\rho)$. The former notions will be rigorously introduced in \cref{chapter:growth}.

    For each non-negative integer $\rho$, the preimage $(\blank \isSemiActedUponBy a^{-1})^{-1}(\ball(\rho))$ is the set $\ball(\rho) \cdot a$, which is not the ball $\ball(a, \rho)$ in the left $\setOf{a, b, a^{-1}, b^{-1}}$-Cayley graph of $F_2$, but it is in the right one. The sequence $\sequence{\ball(\rho) \smallsetminus (\blank \isSemiActedUponBy a^{-1})^{-1}(\ball(\rho))}_{\rho \in \N_0}$ grows exponentially in size, more precisely, it grows in size like $\sequence{3^\rho}_{\rho \in \N_0}$ (see \cref{figure:tree:Folner-original}). Moreover, the sequence $\sequence{\ball(\rho)}_{\rho \in \N_0}$ also grows exponentially in size, more precisely, it grows in size like $\sequence{\sum_{\varrho = 0}^\rho \cardinalityOf{\sphere(\varrho)}}_{\rho \in \N_0} = \sequence{2 \cdot 3^\rho - 1}_{\rho \in \N_0}$. Hence, the quotient $\cardinalityOf{\ball(\rho) \smallsetminus (\blank \isSemiActedUponBy a^{-1})^{-1}(\ball(\rho))} / \cardinalityOf{\ball(\rho)}$ converges to $1/2$ as $\rho$ tends to $\infty$. Therefore, neither the sequence $\sequence{\ball(\rho)}_{\rho \in \N_0}$ nor any of its subsequences is a right Følner net in $\mathcal{R}$. Actually, as we will see, the cell space $\mathcal{R}$ is not right amenable and hence there is no right Følner net in $\mathcal{R}$.
    \begin{figure}
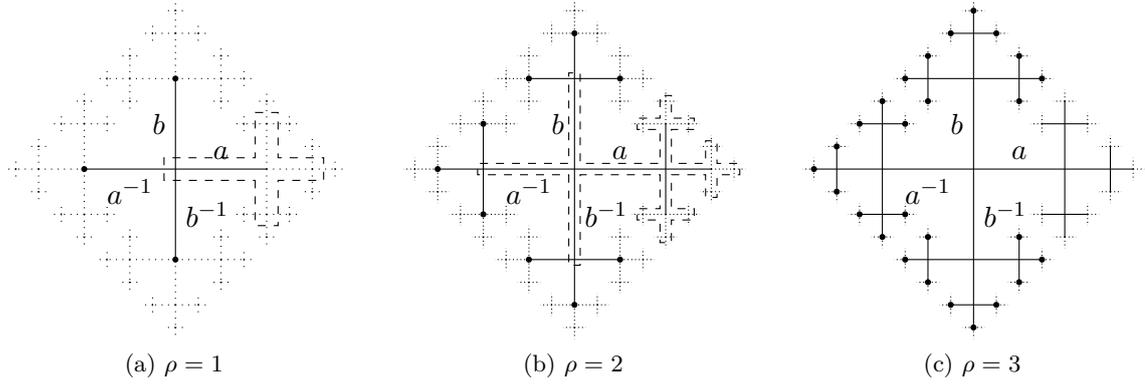

      \myfloatalign
      \begin{wide}
        \figureTreeFolnerOriginal
        \caption{In each subfigure, the same part of the uncoloured \emph{right} $\setOf{a, b, a^{-1}, b^{-1}}$-Cayley graph of $F_2$ is depicted, the vertex in the centre is the neutral element $e_{F_2}$, the vertices that are adjacent to solid edges are the elements of $\ball(\rho)$, the vertices in the region enclosed by the dashed lineation are the elements of $(\blank \isSemiActedUponBy a^{-1})^{-1}(\ball(\rho))$ (note that this lineation is not depicted in the case that $\rho = 3$), and the dots are the elements of $\ball(\rho) \smallsetminus (\blank \isSemiActedUponBy a^{-1})^{-1}(\ball(\rho))$, for the respective $\rho \in \setOf{1, 2, 3}$. Note that in a right Cayley graph, for each element $f \in F_2$ and each generator $s \in \setOf{a, b, a^{-1}, b^{-1}}$, there is an edge from $f$ to $s f$.}
        \label{figure:tree:Folner-original}
      \end{wide}
    \end{figure}
  \end{example}

  \section{Right Paradoxical Decompositions}
  \label{section:right-paradoxical-decomposition}

  In this section, let $\mathcal{R} = \ntuple{\ntuple{M, G, \actsOnPoint}, \ntuple{m_0, \family{g_{m_0, m}}_{m \in M}}}$ be a cell space.

  \begin{definition} 
            Let $A$ and $A'$ be two sets. The set $A \cup A'$ is denoted by $A \disjointUnionWith A'$\graffito{disjoint union $A \disjointUnionWith A'$}\index[symbols]{cupdot@$\disjointUnionWith$} if and only if the sets $A$ and $A'$ are disjoint.
  \end{definition}

  \begin{definition}
  \label{definition:right-paradoxical-decomposition}
    Let $E$ be a finite subset of $G \modulo G_0$, and let $\family{A_e}_{e \in E}$ and $\family{B_e}_{e \in E}$ be two families of subsets of $M$ indexed by $E$ such that, for each index $e \in E$, the map $\blank \isSemiActedUponBy e$ is injective on $A_e$ and on $B_e$, and
    \begin{equation*}
      M = \bigDisjointUnionOf_{e \in E} A_e
        = \bigDisjointUnionOf_{e \in E} B_e
        = \parens*{\bigDisjointUnionOf_{e \in E} A_e \isSemiActedUponBy e} \disjointUnionWith \parens*{\bigDisjointUnionOf_{e \in E} B_e \isSemiActedUponBy e}.
    \end{equation*}
    The triple $\ntuple{E, \family{A_e}_{e \in E}, \family{B_e}_{e \in E}}$ is called \graffito{right paradoxical decomposition $\ntuple{E, \family{A_e}_{e \in E}, \family{B_e}_{e \in E}}$ of $\mathcal{R}$}\index{paradoxical decomposition of $\mathcal{R}$@right paradoxical decomposition of $\mathcal{R}$}\define{right paradoxical decomposition of $\mathcal{R}$}.
  \end{definition}

  \begin{remark} 
    In the situation of \cref{remark:groups:measureIsKindOfSemiActedUponBy}, for each element $g \in G$, the map $\blank \isSemiActedUponBy g$ is injective. Hence, right paradoxical decompositions of $\mathcal{R}$ are the same as right paradoxical decompositions of $G$ as defined in definition~4.8.1 in \cite{ceccherini-silberstein:coornaert:2010}.
  \end{remark}

  \begin{example}[Tree] 
  \label{example:tree:paradoxical-decomposition}
    In the situation of \cref{example:tree}, let $A^{+}$ and $A^{-}$ be the sets of the elements of $M$ whose reduced form ends with $a$ and $a^{-1}$ respectively, let $B^{+}$ be the set of the elements of $M$ whose reduced form ends with $b$ or is $b^{-n}$ for some $n \in \N_0$, and let $B^{-}$ be the set of the elements of $M$ whose reduced form ends with $b^{-1}$ but is not $b^{-n}$ for any $n \in \N_+$. Furthermore, let $E = \setOf{e_{F_2}, a, b}$, let $A_{e_{F_2}} = A^{-}$, let $A_a = A^{+} \cdot a^{-1}$, let $A_b = \emptyset$, let $B_{e_{F_2}} = B^{-}$, let $B_a = \emptyset$, and let $B_b = B^{+} \cdot b^{-1}$. The set $A_a$ is the set of the elements of $M$ whose reduced form does not end with $a^{-1}$ and the set $B_b$ is the set of the elements of $M$ whose reduced form does not end with $b^{-1}$ or is $b^{-n}$ for some $n \in \N_0$. The triple $\ntuple{E, \family{A_e}_{e \in E}, \family{B_e}_{e \in E}}$ is a right paradoxical decomposition of $\mathcal{R}$ (see \cref{figure:tree:paradoxical-decomposition}; compare example~4.8.2 in \cite{ceccherini-silberstein:coornaert:2010}).
    \begin{figure}
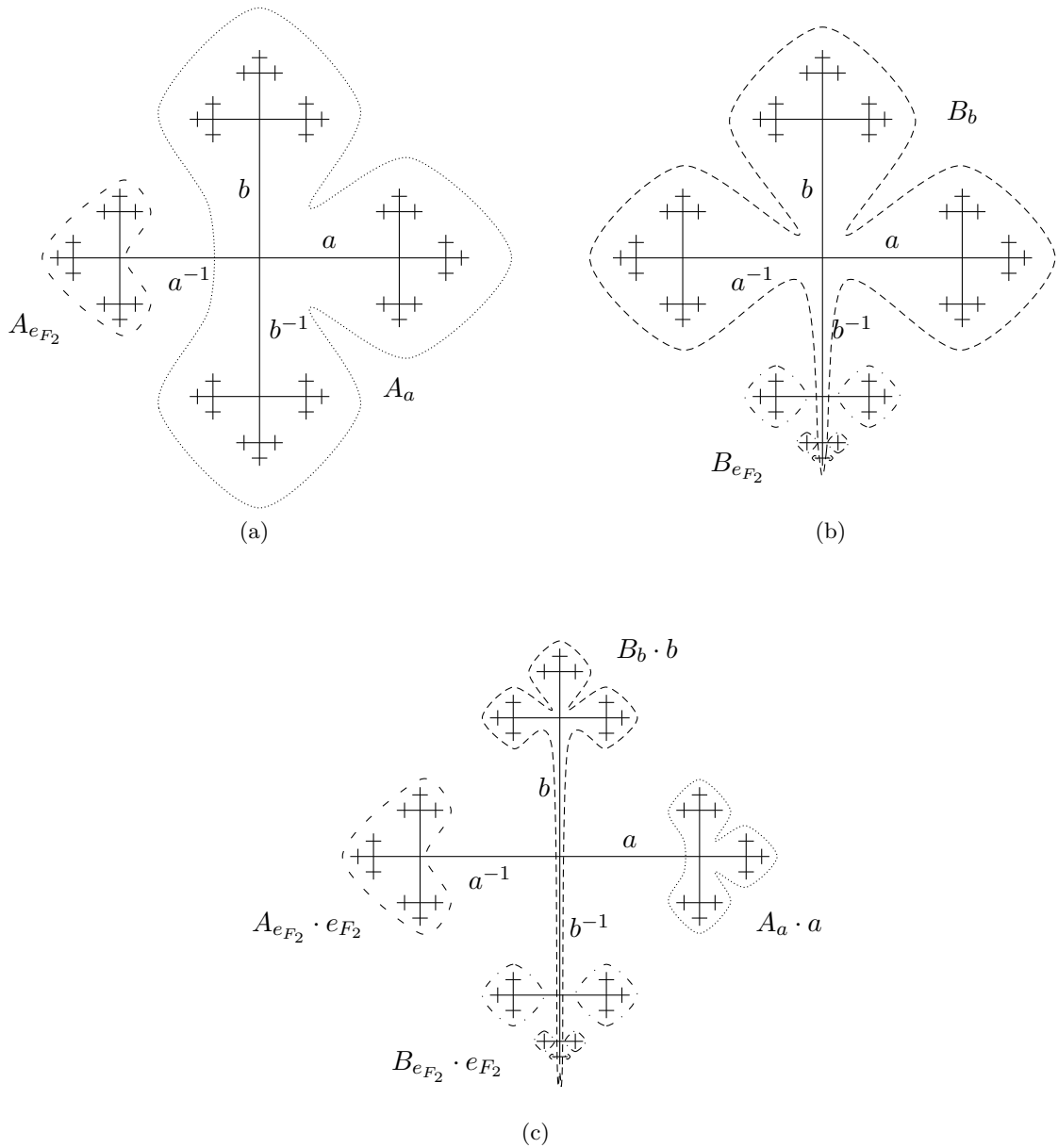

      \myfloatalign
      \begin{wide}
        \figureTreeParadoxicalDecomposition
        \caption{In each subfigure, the same part of the \emph{right} $\setOf{a, b, a^{-1}, b^{-1}}$-Cayley graph of $F_2$ is depicted and the vertex in the centre is the neutral element $e_{F_2}$. Note that in a right Cayley graph, for each element $f \in F_2$ and each generator $s \in \setOf{a, b, a^{-1}, b^{-1}}$, there is an edge from $f$ to $s f$. The graph in \cref{figure:tree:paradoxical-decomposition:one} illustrates the partition of $F_2$ by $A_{e_{F_2}}$, $A_a$, and $A_b = \emptyset$; the one in \cref{figure:tree:paradoxical-decomposition:two} the partition by $B_{e_{F_2}}$, $B_a = \emptyset$, and $B_b$; and the one in \cref{figure:tree:paradoxical-decomposition:three} the partition by $A_{e_{F_2}} \cdot e_{F_2}$, $A_a \cdot a$, $A_b \cdot b = \emptyset$, $B_{e_{F_2}} \cdot e_{F_2}$, $B_a \cdot a = \emptyset$, and $B_b \cdot b$.} 
        \label{figure:tree:paradoxical-decomposition}
      \end{wide}
    \end{figure}
  \end{example}

  A right paradoxical decomposition of $M$ is one of $\indicatorFunction_M$ in the sense given in

  \begin{lemma}
  \label{lemma:functional-right-paradoxical-decomposition}
    Let $G_0$ be finite and let $\ntuple{E, \family{A_e}_{e \in E}, \family{B_e}_{e \in E}}$ be a right paradoxical decomposition of $\mathcal{R}$. Then,
    \begin{equation*}
      \indicatorFunction_M = \sum_{e \in E} \indicatorFunction_{A_e}
                  = \sum_{e \in E} \indicatorFunction_{B_e}
                  = \sum_{e \in E} (\indicatorFunction_{A_e} \boundedFunctionIsKindOfSemiActedUponBy e) + \sum_{e \in E} (\indicatorFunction_{B_e} \boundedFunctionIsKindOfSemiActedUponBy e). \qedhere
    \end{equation*}
  \end{lemma}

  \begin{proof}
    This is a direct consequence of \cref{definition:right-paradoxical-decomposition} and \cref{lemma:indicator-function-of-A-liberation-n-is-indicator-function-of-A-acted-upon}.
  \end{proof}

  \section{Tarski's and Følner's Theorem}
  \label{section:Tarski-Folner-theorem}

  In this section, let $\mathcal{R} = \ntuple{\ntuple{M, G, \actsOnPoint}, \ntuple{m_0, \family{g_{m_0, m}}_{m \in M}}}$ be a cell space such that the stabiliser $G_0$ of $m_0$ under $\actsOnPoint$ is finite.

  \paragraph{Contents.} In \cref{lemma:meanIsKindOfSemiActedUponBy-is-continuous,lemma:convergent-net-sufficient-for-amenability,lemma:double-liberation-with-sets-to-one-liberation} we show some technical properties that are used in the proof of the Tarski-Følner theorem~\ref{theorem:Tarski-Folner}, which states for a cell space that right amenability, the existence of a right Følner net, and the non-existence of a right paradoxical decomposition are equivalent. The analogous statement for left-ho\-mo\-ge\-neous spaces is given in \cref{corollary:Tarski-Folner-for-left-homogeneous-spaces}.

  \begin{lemma} 
  \label{lemma:meanIsKindOfSemiActedUponBy-is-continuous}
    Let $\mathfrak{g}$ be an element of $G \modulo G_0$. The map $\blank \meanIsKindOfSemiActedUponBy \mathfrak{g}$ is continuous, where $\boundedFunctionsOn(M)^*$ is equipped with the weak-$*$ topology.
  \end{lemma}

  \begin{proof}
    For each $f \in \boundedFunctionsOn(M)$, let $\evaluationMap_f$ be the evaluation map $\boundedFunctionsOn(M)^* \to \R$, $\psi \mapsto \psi(f)$. Furthermore, let $f \in \boundedFunctionsOn(M)$. Then, for each $\psi \in \boundedFunctionsOn(M)^*$,
    \begin{equation*}
      (\evaluationMap_f \after (\blank \meanIsKindOfSemiActedUponBy \mathfrak{g}))(\psi)
      = \evaluationMap_f(\psi \meanIsKindOfSemiActedUponBy \mathfrak{g})
      = (\psi \meanIsKindOfSemiActedUponBy \mathfrak{g})(f)
      = \psi(f \boundedFunctionIsKindOfSemiActedUponBy \mathfrak{g})
      = \evaluationMap_{f \boundedFunctionIsKindOfSemiActedUponBy \mathfrak{g}}(\psi).
    \end{equation*}
    Thus, $\evaluationMap_f \after (\blank \meanIsKindOfSemiActedUponBy \mathfrak{g}) = \evaluationMap_{f \boundedFunctionIsKindOfSemiActedUponBy \mathfrak{g}}$. Hence, because $\evaluationMap_{f \boundedFunctionIsKindOfSemiActedUponBy \mathfrak{g}}$ is continuous, so is $\evaluationMap_f \after (\blank \meanIsKindOfSemiActedUponBy \mathfrak{g})$. Therefore, according to \cref{lemma:map-to-initial-topology-continuous-if-and-only-if-gens-after-map-continuous}, the map $\blank \meanIsKindOfSemiActedUponBy \mathfrak{g}$ is continuous.
  \end{proof}

  \begin{lemma} 
  \label{lemma:convergent-net-sufficient-for-amenability}
    Let $\net{\nu_i}_{i \in I}$ be a net in $\meansOn(M)$ such that, for each element $\mathfrak{g} \in G \modulo G_0$, the net $\net{(\nu_i \meanIsKindOfSemiActedUponBy \mathfrak{g}) - \nu_i}_{i \in I}$ converges to $\functionThatIsIdenticalToZero$ in $\boundedFunctionsOn(M)^*$ equipped with the weak-$*$ topology. The cell space $\mathcal{R}$ is right amenable.
  \end{lemma}

  \begin{proof} 
    Let $\mathfrak{g} \in G \modulo G_0$. According to \cref{theorem:means-convex-and-compact-wrt-weak-star-topology}, the set $\meansOn(M)$ is compact in $\boundedFunctionsOn(M)^*$ equipped with the weak-$*$ topology. Hence, there is a subnet $\net{\nu_{i_j}}_{j \in J}$ of $\net{\nu_i}_{i \in I}$ that converges to a $\nu \in \meansOn(M)$. Because, according to \cref{lemma:meanIsKindOfSemiActedUponBy-is-continuous}, the map $\blank \meanIsKindOfSemiActedUponBy \mathfrak{g}$ is continuous, the net $\net{(\nu_{i_j} \meanIsKindOfSemiActedUponBy \mathfrak{g}) - \nu_{i_j}}_{j \in J}$ converges to $(\nu \meanIsKindOfSemiActedUponBy \mathfrak{g}) - \nu$ in $\boundedFunctionsOn(M)^*$. Because it is a subnet of $\net{(\nu_i \meanIsKindOfSemiActedUponBy \mathfrak{g}) - \nu_i}_{i \in I}$, it also converges to $\functionThatIsIdenticalToZero$ in $\boundedFunctionsOn(M)^*$. Because, according to \cref{lemma:weak-star-topology-is-Hausdorff}, the space $\boundedFunctionsOn(M)^*$ is Hausdorff, we have $(\nu \meanIsKindOfSemiActedUponBy \mathfrak{g}) - \nu = \functionThatIsIdenticalToZero$ and hence $\nu \meanIsKindOfSemiActedUponBy \mathfrak{g} = \nu$. Altogether, $\nu$ is a $\meanIsKindOfSemiActedUponBy$-invariant mean. In conclusion, according to \cref{theorem:mean-characterisation-of-right-amenable}, the cell space $\mathcal{R}$ is right amenable.
  \end{proof}

  \begin{lemma}
  \label{lemma:double-liberation-with-sets-to-one-liberation}
    Let $m$ be an element of $M$, and let $E$ and $E'$ be two subsets of $G \modulo G_0$. There is a subset $E''$ of $G \modulo G_0$ such that $(m \isSemiActedUponBy E) \isSemiActedUponBy E' = m \isSemiActedUponBy E''$; if $G_0 \in E \cap E'$, then $G_0 \in E''$; if $E$ and $E'$ are finite, then $\cardinalityOf{E''} \leq \cardinalityOf{E} \cdot \cardinalityOf{E'}$; and if $G_0 \cdot E' \subseteq E'$, then $E'' = \setOf{g \cdot e' \suchThat e \in E, e' \in E', g \in e}$.
  \end{lemma}

  \begin{proof}
    For each $e \in E$, according to \cref{lemma:rightsemiaction-can-be-undone}, there is a $g_e \in e$ such that 
    \begin{equation*}
      \ForEach \mathfrak{g} \in G \modulo G_0 \Holds (m \isSemiActedUponBy e) \isSemiActedUponBy \mathfrak{g} = m \isSemiActedUponBy g_e \cdot \mathfrak{g}. 
    \end{equation*}
    Put $E'' = \setOf{g_e \cdot e' \suchThat e \in E, e' \in E'}$. Then, $(m \isSemiActedUponBy E) \isSemiActedUponBy E' = m \isSemiActedUponBy E''$. Moreover, if $G_0 \in E \cap E'$, then $G_0 = g_{G_0} \cdot G_0 \in E''$; if $E$ and $E'$ are finite, then $\cardinalityOf{E''} \leq \cardinalityOf{E} \cdot \cardinalityOf{E'}$; and if $G_0 \cdot E' \subseteq E'$, then $E''$ is as stated. 
  \end{proof}

  \begin{main-theorem} 
  \label{theorem:Tarski-Folner}
    Let $\mathcal{R} = \ntuple{\ntuple{M, G, \actsOnPoint}, \ntuple{m_0, \family{g_{m_0, m}}_{m \in M}}}$ be a cell space such that the stabiliser $G_0$ of $m_0$ under $\actsOnPoint$ is finite. The following statements are equivalent:
    \begin{aenumerate}
      \item \label{item:Tarski-Folner:not-right-amenable}
            The cell space $\mathcal{R}$ is not right amenable;
      \item \label{item:Tarski-Folner:no-Folner-net}
            There is no right Følner net in $\mathcal{R}$;
      \item \label{item:Tarski-Folner:finite-N}
            There is a finite subset $E$ of $G \modulo G_0$ such that $G_0 \in E$ and, for each finite subset $F$ of $M$, we have $\cardinalityOf{F \isSemiActedUponBy E} \geq 2 \cardinalityOf{F}$;
      \item \label{item:Tarski-Folner:two-to-one-surjective-map}
            There is a $2$-to-$1$ surjective map $\phi \from M \to M$ and there is a finite subset $E$ of $G \modulo G_0$ such that
            \begin{equation*}
              \ForEach m \in M \Exists e \in E \SuchThat \phi(m) \isSemiActedUponBy e = m;
            \end{equation*}
      \item \label{item:Tarski-Folner:no-paradoxical-decomposition}
            There is a right paradoxical decomposition of $\mathcal{R}$. \qedhere
    \end{aenumerate}
  \end{main-theorem}

  \begin{proof}
    \begin{description}
      \item[\ref{item:Tarski-Folner:not-right-amenable} implies \ref{item:Tarski-Folner:no-Folner-net}.] 
            Let there be a right Følner net $\net{F_i}_{i \in I}$ in $\mathcal{R}$. Furthermore, let $i \in I$. Put
            \begin{align*}
              \nu_i \from \boundedFunctionsOn(M) &\to \R,\\
              f &\mapsto \frac{1}{\cardinalityOf{F_i}} \sum_{m \in F_i} f(m).
            \end{align*}
            Then, $\nu_i \in \meansOn(M)$. 
            Moreover, let $\mathfrak{g} \in G \modulo G_0$ and let $f \in \boundedFunctionsOn(M)$. Then,
            \begin{align*}
              (\nu_i \meanIsKindOfSemiActedUponBy \mathfrak{g})(f)
              &= \nu_i(f \boundedFunctionIsKindOfSemiActedUponBy \mathfrak{g})\\
              &= \frac{1}{\cardinalityOf{F_i}} \sum_{m \in F_i} (f \boundedFunctionIsKindOfSemiActedUponBy \mathfrak{g})(m)\\
              &= \frac{1}{\cardinalityOf{F_i}} \sum_{m \in F_i} \sum_{m' \in (\blank \isSemiActedUponBy \mathfrak{g})^{-1}(m)} f(m')\\
              &= \frac{1}{\cardinalityOf{F_i}} \sum_{m \in (\blank \isSemiActedUponBy \mathfrak{g})^{-1}(F_i)} f(m).
            \end{align*}
            Hence,
            \begin{equation*}
              (\nu_i \meanIsKindOfSemiActedUponBy \mathfrak{g} - \nu_i)(f)
              = \frac{1}{\cardinalityOf{F_i}}
                \begin{aligned}[t]
                  \bigg(&\sum_{m \in (\blank \isSemiActedUponBy \mathfrak{g})^{-1}(F_i) \smallsetminus F_i} f(m)\\
                        &{}- \sum_{m \in F_i \smallsetminus (\blank \isSemiActedUponBy \mathfrak{g})^{-1}(F_i)} f(m)\bigg).
                \end{aligned}
            \end{equation*}
            Therefore,
            \begin{align*}
              \absoluteValueOf{(\nu_i \meanIsKindOfSemiActedUponBy \mathfrak{g} - \nu_i)(f)}
              &\leq \frac{1}{\cardinalityOf{F_i}}
                    \begin{aligned}[t]
                      \bigg(&\sum_{m \in (\blank \isSemiActedUponBy \mathfrak{g})^{-1}(F_i) \smallsetminus F_i} \absoluteValueOf{f(m)}\\
                            &{}+ \sum_{m \in F_i \smallsetminus (\blank \isSemiActedUponBy \mathfrak{g})^{-1}(F_i)} \absoluteValueOf{f(m)}\bigg)
                    \end{aligned}\\
              &\leq 
                    \begin{aligned}[t]
                      \bigg(&\frac{\cardinalityOf{(\blank \isSemiActedUponBy \mathfrak{g})^{-1}(F_i) \smallsetminus F_i}}{\cardinalityOf{F_i}}\\
                            &{}+ \frac{\cardinalityOf{F_i \smallsetminus (\blank \isSemiActedUponBy \mathfrak{g})^{-1}(F_i)}}{\cardinalityOf{F_i}}\bigg) \cdot \normOf{f}_\infty.
                    \end{aligned}
            \end{align*}
            According to \cref{definition:right-Folner-net} and \cref{theorem:reversed-sets-characterisation-of-Folner-nets}, the nets $\net{\cardinalityOf{(\blank \isSemiActedUponBy \mathfrak{g})^{-1}(F_i) \smallsetminus F_i} / \cardinalityOf{F_i}}_{i \in I}$ and $\net{\cardinalityOf{F_i \smallsetminus (\blank \isSemiActedUponBy \mathfrak{g})^{-1}(F_i)} / \cardinalityOf{F_i}}_{i \in I}$ converge to $0$. Hence, so does $\net{\absoluteValueOf{(\nu_i \meanIsKindOfSemiActedUponBy \mathfrak{g} - \nu_i)(f)}}_{i \in I}$. Thus, the net $\net{\nu_i \meanIsKindOfSemiActedUponBy \mathfrak{g} - \nu_i}_{i \in I}$ converges to $\functionThatIsIdenticalToZero$ in $\boundedFunctionsOn(M)^*$ equipped with the weak-$*$ topology. Hence, according to \cref{lemma:convergent-net-sufficient-for-amenability}, the cell space $\mathcal{R}$ is right amenable. In conclusion, by contraposition, if $\mathcal{R}$ is not right amenable, then there is no right Følner net in $\mathcal{R}$.
      \item[\ref{item:Tarski-Folner:no-Folner-net} implies \ref{item:Tarski-Folner:finite-N}.]
            Let there be no right Følner net in $\mathcal{R}$. According to \cref{theorem:epsilon-characterisation-of-Folner-net}, there is a finite $E_1 \subseteq G \modulo G_0$ and an $\varepsilon \in \R_{> 0}$ such that, for each non-empty and finite $F \subseteq M$, there is an $e_F \in E_1$ such that 
            \begin{equation*}
              \frac{\cardinalityOf{F \smallsetminus (\blank \isSemiActedUponBy e_F)^{-1}(F)}}{\cardinalityOf{F}} \geq \varepsilon.
            \end{equation*}
            Put $E_2 = \setOf{G_0} \cup E_1$.

              Let $F$ be a non-empty and finite subset of $M$. Then, $F \subseteq F \cup (F \isSemiActedUponBy E_1) = F \isSemiActedUponBy E_2$. Thus,
              \begin{align*}
                \cardinalityOf{F \isSemiActedUponBy E_2} - \cardinalityOf{F}
                &=    \cardinalityOf{(F \isSemiActedUponBy E_2) \smallsetminus F}\\
                &=    \cardinalityOf{(F \isSemiActedUponBy E_1) \smallsetminus F}\\
                &\geq \cardinalityOf{(F \isSemiActedUponBy e_F) \smallsetminus F}.
              \end{align*}
              Moreover, according to \cref{lemma:liberation-preimage}, we have $\cardinalityOf{(\blank \isSemiActedUponBy e_F)^{-1}((F \isSemiActedUponBy e_F) \smallsetminus F)} \leq \cardinalityOf{G_0} \cdot \cardinalityOf{(F \isSemiActedUponBy e_F) \smallsetminus F}$. Hence,
              \begin{equation*}
                \cardinalityOf{F \isSemiActedUponBy E_2} - \cardinalityOf{F}
                \geq \frac{\cardinalityOf{(\blank \isSemiActedUponBy e_F)^{-1}((F \isSemiActedUponBy e_F) \smallsetminus F)}}{\cardinalityOf{G_0}}.
              \end{equation*}
              Therefore, because $F \smallsetminus (\blank \isSemiActedUponBy e_F)^{-1}(F) \subseteq (\blank \isSemiActedUponBy e_F)^{-1}((F \isSemiActedUponBy e_F) \smallsetminus F)$, 
              \begin{align*}
                \cardinalityOf{F \isSemiActedUponBy E_2} - \cardinalityOf{F}
                &\geq \frac{\cardinalityOf{F \smallsetminus (\blank \isSemiActedUponBy e_F)^{-1}(F)}}{\cardinalityOf{G_0}}\\
                &\geq \frac{\varepsilon}{\cardinalityOf{G_0}} \cardinalityOf{F}.
              \end{align*}
              Put $\xi = 1 + \varepsilon / \cardinalityOf{G_0}$. Then, $\cardinalityOf{F \isSemiActedUponBy E_2} \geq \xi \cardinalityOf{F}$. Because $\varepsilon$ does not depend on $F$, neither does $\xi$. Therefore, for each non-empty and finite $F \subseteq M$, we have $\cardinalityOf{F \isSemiActedUponBy E_2} \geq \xi \cardinalityOf{F}$.

            Let $F$ be a non-empty and finite subset of $M$. Because $\xi > 1$, there is an $n \in \N_+$ such that $\xi^n \geq 2$. Hence,
            \begin{align*} 
              \cardinalityOf{\underbrace{(((F \isSemiActedUponBy E_2) \isSemiActedUponBy \dotsb) \isSemiActedUponBy E_2) \isSemiActedUponBy E_2}_{\text{$n$ times}}}
              &\geq \xi \cardinalityOf{((F \isSemiActedUponBy E_2) \isSemiActedUponBy \dotsb) \isSemiActedUponBy E_2}\\
              &\geq \dotsb\\
              &\geq \xi^n \cardinalityOf{F}\\
              &\geq 2 \cardinalityOf{F}.
            \end{align*}
            Moreover, according to \cref{lemma:double-liberation-with-sets-to-one-liberation}, there is an $E \subseteq G \modulo G_0$ such that $E$ is finite, $G_0 \in E$, and $F \isSemiActedUponBy E = (((F \isSemiActedUponBy E_2) \isSemiActedUponBy \dotsb) \isSemiActedUponBy E_2) \isSemiActedUponBy E_2$. In conclusion, $\cardinalityOf{F \isSemiActedUponBy E} \geq 2 \cardinalityOf{F}$.
      \item[\ref{item:Tarski-Folner:finite-N} implies \ref{item:Tarski-Folner:two-to-one-surjective-map} (see \cref{figure:finite-n-implies-two-to-one-surjective-map}).]
            Let there be a finite $E \subseteq G \modulo G_0$ such that, for each finite $F \subseteq M$, we have $\cardinalityOf{F \isSemiActedUponBy E} \geq 2 \cardinalityOf{F}$. Furthermore, let $\mathcal{G}$ be the bipartite graph 
            \begin{equation*}
              \ntuple{M, M, \setOf{(m, m') \in M \times M \suchThat \Exists e \in E \SuchThat m \isSemiActedUponBy e = m'}}.
            \end{equation*}
            Moreover, let $F$ be a finite subset of $M$. The right neighbourhood of $F$ in $\mathcal{G}$ is
            \begin{equation*}
              \mathcal{N}_r(F)
              = \setOf{m' \in M \suchThat \Exists e \in E \SuchThat F \isSemiActedUponBy e \ni m'}
              = F \isSemiActedUponBy E
            \end{equation*}
            and the left neighbourhood of $F$ in $\mathcal{G}$ is
            \begin{equation*}
              \mathcal{N}_l(F)
              = \setOf{m \in M \suchThat \Exists e \in E \SuchThat m \isSemiActedUponBy e \in F}
              = \bigcup_{e \in E} (\blank \isSemiActedUponBy e)^{-1}(F).
            \end{equation*}
            By precondition $\cardinalityOf{\mathcal{N}_r(F)} = \cardinalityOf{F \isSemiActedUponBy E} \geq 2 \cardinalityOf{F}$. Moreover, because $G_0 \in E$,
            we have $F = (\blank \isSemiActedUponBy G_0)^{-1}(F) \subseteq \mathcal{N}_l(F)$ and hence $\cardinalityOf{\mathcal{N}_l(F)} \geq \cardinalityOf{F} \geq 2^{-1} \cardinalityOf{F}$. Therefore, according to
            the Hall harem \cref{theorem:hall-harem}, there is a perfect $(1,2)$-matching for $\mathcal{G}$. In conclusion, there is a $2$-to-$1$ surjective map $\phi \from M \to M$ such that, for each $m \in M$, the tuple $(\phi(m), m)$ is an edge in $\mathcal{G}$, that is, there is an $e \in E$ such that $\phi(m) \isSemiActedUponBy e = m$.
            \begin{figure}
              \myfloatalign
              \begin{wide}
                \figureFiniteNImpliesTwoToOneSurjectiveMap
                \caption{
                  Schematic representation of the set-up of the proof of \textsc{\ref{item:Tarski-Folner:finite-N} implies \ref{item:Tarski-Folner:two-to-one-surjective-map}} of \cref{theorem:Tarski-Folner}: Each region enclosed by one of the two rectangles is $M$; the regions enclosed by the smaller circles with solid borders are subsets $F$ and $F'$ of $M$ respectively; the regions enclosed by the circles with dashed borders are $m \isSemiActedUponBy E$ and $m' \isSemiActedUponBy E$ respectively; the region enclosed by the circle with dotted border is $\bigcup_{e' \in E} (\blank \isSemiActedUponBy e')^{-1}(m')$; the regions enclosed by the larger circles are $\mathcal{N}_r(F)$ and $\mathcal{N}_l(F')$ respectively.
                }
                \label{figure:finite-n-implies-two-to-one-surjective-map}
              \end{wide}
            \end{figure}
      \item[\ref{item:Tarski-Folner:two-to-one-surjective-map} implies \ref{item:Tarski-Folner:no-paradoxical-decomposition} (see \cref{figure:two-to-one-surjective-map-implies-tarski-folner:no-paradoxical-decomposition}).]
            Let there be a $2$-to-$1$ surjective map $\phi \from M \to M$ and a finite subset $E$ of $G \modulo G_0$ such that
            \begin{equation*}
              \ForEach m \in M \Exists e \in E \SuchThat \phi(m) \isSemiActedUponBy e = m.
            \end{equation*}
            By the axiom of choice, there are two injective maps $\psi$ and $\psi' \from M \to M$ such that, for each $m \in M$, we have $\phi^{-1}(m) = \setOf{\psi(m), \psi'(m)}$. For each $e \in E$, let
            \begin{equation*}
              A_e = \setOf{m \in M \suchThat m \isSemiActedUponBy e = \psi(m)}
            \end{equation*}
            and let
            \begin{equation*}
              B_e = \setOf{m \in M \suchThat m \isSemiActedUponBy e = \psi'(m)}.
            \end{equation*}
            Let $m \in M$. There is an $e \in E$ such that $\phi(\psi(m)) \isSemiActedUponBy e = \psi(m)$. Because $\phi(\psi(m)) = m$, we have $m \in A_e$. And, because $\isSemiActedUponBy$ is free, for each $e' \in E \smallsetminus \setOf{e}$, we have $m \isSemiActedUponBy e' \neq m \isSemiActedUponBy e = \psi(m)$ and thus $m \notin A_{e'}$. Therefore,
            \begin{equation*}
              M = \bigDisjointUnionOf_{e \in E} A_e,
            \end{equation*}
            and analogously
            \begin{equation*}
              M = \bigDisjointUnionOf_{e \in E} B_e.
            \end{equation*}
            Moreover, $\psi(A_e) = A_e \isSemiActedUponBy e$ and $\psi'(B_e) = B_e \isSemiActedUponBy e$. Hence, because $M = \psi(M) \disjointUnionWith \psi'(M)$, and $\psi$ and $\psi'$ are injective,
            \begin{align*}
              M 
                &= \parens*{\bigDisjointUnionOf_{e \in E} \psi(A_e)} \disjointUnionWith \parens*{\bigDisjointUnionOf_{e \in E} \psi'(B_e)}\\
                &= \parens*{\bigDisjointUnionOf_{e \in E} A_e \isSemiActedUponBy e} \disjointUnionWith \parens*{\bigDisjointUnionOf_{e \in E} B_e \isSemiActedUponBy e}.
            \end{align*}
            Furthermore, because $\psi$ and $\psi'$ are injective, for each $e \in E$, the maps $(\blank \isSemiActedUponBy e)\restrictedTo_{A_e} = \psi\restrictedTo_{A_e}$ and $(\blank \isSemiActedUponBy e)\restrictedTo_{B_e} = \psi'\restrictedTo_{B_e}$ are injective. In conclusion, $\ntuple{E, \family{A_e}_{e \in E}, \family{B_e}_{e \in E}}$ is a right paradoxical decomposition of $\mathcal{R}$.
            \begin{figure}
              \myfloatalign
              \begin{wide}
                \figureTwoToOneSurjectiveMapImpliesTarskiFolnerNoParadoxicalDecomposition
                \caption{
                  Schematic representation of the set-up of the proof of \textsc{\ref{item:Tarski-Folner:two-to-one-surjective-map} implies \ref{item:Tarski-Folner:no-paradoxical-decomposition}} of \cref{theorem:Tarski-Folner}: Each region enclosed by one of the three columns, two with solid and one with dashed border, is $M$; the dot in the right column, called $m$, is an element of $M$, and the two dots in the left column are its preimages under $\phi$, which are its images under $\psi$ and $\psi'$; there are elements $e$ and $e'$ of $E$ such that $m \isSemiActedUponBy e = \psi(m)$ and $m \isSemiActedUponBy e' = \psi'(m)$, in other words, $m \in A_e$ and $m \in B_{e'}$; as depicted in the right columns with solid and dashed borders, the families $\family{A_e}_{e \in E}$ and $\family{B_{e'}}_{e' \in E}$ are partitions of $M$; as depicted in the left column, the set $\setOf{\psi(M), \psi'(M)}$ is a partition of $M$, the family $\family{\psi(A_e)}_{e \in E} = \family{A_e \isSemiActedUponBy e}_{e \in E}$ is a partition of $\psi(M)$, and the family $\family{\psi'(B_{e'})}_{e' \in E} = \family{B_{e'} \isSemiActedUponBy e'}_{e' \in E}$ is a partition of $\psi'(M)$.
                }
                \label{figure:two-to-one-surjective-map-implies-tarski-folner:no-paradoxical-decomposition}
              \end{wide}
            \end{figure}
      \item[\ref{item:Tarski-Folner:no-paradoxical-decomposition} implies \ref{item:Tarski-Folner:not-right-amenable}.]
            Let there be a right paradoxical decomposition $\ntuple{E, \family{A_e}_{e \in E}, \family{B_e}_{e \in E}}$ of $\mathcal{R}$. According to \cref{lemma:functional-right-paradoxical-decomposition},
            \begin{equation*}
              \indicatorFunction_M = \sum_{e \in E} \indicatorFunction_{A_e}
                          = \sum_{e \in E} \indicatorFunction_{B_e}
                          = \sum_{e \in E} (\indicatorFunction_{A_e} \boundedFunctionIsKindOfSemiActedUponBy e) + \sum_{e \in E} (\indicatorFunction_{B_e} \boundedFunctionIsKindOfSemiActedUponBy e).
            \end{equation*}
            Suppose that $\mathcal{R}$ is right amenable. Then, according to theorem \ref{theorem:mean-characterisation-of-right-amenable}, there is a $\meanIsKindOfSemiActedUponBy$-invariant mean $\nu$ on $M$. Because $\nu$ is linear and normalised, 
            \begin{align*}
              1 &= \nu(\indicatorFunction_M)\\
                &= \sum_{e \in E} \nu(\indicatorFunction_{A_e} \boundedFunctionIsKindOfSemiActedUponBy e) + \sum_{e \in E} \nu(\indicatorFunction_{B_e} \boundedFunctionIsKindOfSemiActedUponBy e)\\
                &= \sum_{e \in E} (\nu \meanIsKindOfSemiActedUponBy e)(\indicatorFunction_{A_e}) + \sum_{e \in E} (\nu \meanIsKindOfSemiActedUponBy e)(\indicatorFunction_{B_e})\\
                &= \sum_{e \in E} \nu(\indicatorFunction_{A_e}) + \sum_{e \in E} \nu(\indicatorFunction_{B_e})\\
                &= \nu(\indicatorFunction_M) + \nu(\indicatorFunction_M)\\
                &= 1 + 1\\
                &= 2,
            \end{align*}
            which contradicts that $1 \neq 2$. In conclusion, $\mathcal{R}$ is not right amenable. \qedhere  
    \end{description}
  \end{proof}

  \begin{corollary}[Tarski alternative theorem; Alfred Tarski, 1938] 
  \label{corollary:tarski}
    Let $\mathcal{M}$ be a left-ho\-mo\-ge\-neous space with finite stabilisers. It is right amenable if and only if there is a coordinate system $\mathcal{K}$ for $\mathcal{M}$ such that there is no right paradoxical decomposition of $\ntuple{\mathcal{M}, \mathcal{K}}$.
  \end{corollary}

  \begin{corollary}[Theorem of Følner; Erling Følner, 1955]
  \label{corollary:Folner}
    Let $\mathcal{M}$ be a left-ho\-mo\-ge\-neous space with finite stabilisers. It is right amenable if and only if there is a coordinate system $\mathcal{K}$ for $\mathcal{M}$ such that there is a right Følner net in $\ntuple{\mathcal{M}, \mathcal{K}}$.
  \end{corollary}

  \begin{corollary}
  \label{corollary:Tarski-Folner-for-left-homogeneous-spaces}
    Let $\mathcal{M}$ be a left-ho\-mo\-ge\-neous space with finite stabilisers. The following statements are equivalent:
    \begin{aenumerate}
      \item The space $\mathcal{M}$ is right amenable;
      \item For each coordinate system $\mathcal{K}$ for $\mathcal{M}$, there is a $\measureIsKindOfSemiActedUponBy$-semi-invariant finitely additive probability measure on $M$;
      \item For each coordinate system $\mathcal{K}$ for $\mathcal{M}$, there is a $\meanIsKindOfSemiActedUponBy$-invariant mean on $M$;
      \item For each coordinate system $\mathcal{K}$ for $\mathcal{M}$, there is a right Følner net in $\ntuple{\mathcal{M}, \mathcal{K}}$;
      \item For each coordinate system $\mathcal{K}$ for $\mathcal{M}$, there is no right paradoxical decomposition of $\ntuple{\mathcal{M}, \mathcal{K}}$. \qedhere
    \end{aenumerate}
  \end{corollary}

  \begin{proof} 
    This is a direct consequence of \cref{theorem:Folner-net-independent-of-coordinate-system}, \cref{theorem:mean-characterisation-of-right-amenable}, and \cref{theorem:Tarski-Folner}. \qedhere
  \end{proof}

  \begin{remark}
    In the situation of \cref{remark:groups:measureIsKindOfSemiActedUponBy}, \cref{corollary:tarski,corollary:Folner} constitute theorem~4.9.1 in \cite{ceccherini-silberstein:coornaert:2010}.
  \end{remark}

  \begin{example}[Tree]
  \label{example:tree:Tarski-Folner}
    In \cref{example:tree:paradoxical-decomposition} we constructed a right paradoxical decomposition of $\mathcal{R}$. Hence, according to \cref{theorem:Tarski-Folner}, the cell space $\mathcal{R}$ is not right amenable, there is a subset $E$ of $G \modulo G_0$ as in \cref{item:Tarski-Folner:finite-N} of \cref{theorem:Tarski-Folner}, there is a map $\phi$ from $M$ to $M$ as in \cref{item:Tarski-Folner:two-to-one-surjective-map} of \cref{theorem:Tarski-Folner}, and there are maps $\psi$ and $\psi'$ from $M$ to $M$ as in the subproof \textsc{\ref{item:Tarski-Folner:two-to-one-surjective-map} implies \ref{item:Tarski-Folner:no-paradoxical-decomposition}} of \cref{theorem:Tarski-Folner}. For the right paradoxical decomposition we constructed, these sets and maps can be given explicitly and we do so in the following.

    The map
    \begin{align*}
      \phi \from M &\to M,\\
      m &\mapsto \begin{dcases*}
                   m, &if $m \in A^{-}\ (= A_{e_{F_2}})$,\\
                   m \cdot a^{-1}, &if $m \in A^{+}\ (= A_a \cdot a)$,\\
                   m, &if $m \in B^{-}\ (= B_{e_{F_2}})$,\\
                   m \cdot b^{-1}, &if $m \in B^{+}\ (= B_b \cdot b)$,
                 \end{dcases*}
    \end{align*}
    is well-defined, because the family $\family{A^{-}, A^{+}, B^{-}, B^{+}}$ is a partition of $M$; it is $2$-to-$1$ surjective, because it is bijective from $A^{-} \disjointUnionWith A^{+}$ to $M$ as well as from $B^{-} \disjointUnionWith B^{+}$ to $M$ (see \cref{figure:tree:Tarski-Folner:2-to-1-map:phi}); and, for each cell $m \in M$, there is an element $e \in E$ such that $\phi(m) \cdot e = m$, because if $m \in A^{-} \cup B^{-}$, then $\phi(m) \cdot e_{F_2} = m$, if $m \in A^{+}$, then $\phi(m) \cdot a = m$, and if $m \in B^{+}$, then $\phi(m) \cdot b = m$. Hence, the map $\phi$ satisfies the properties of \cref{item:Tarski-Folner:two-to-one-surjective-map} of \cref{theorem:Tarski-Folner}.
    \begin{figure}
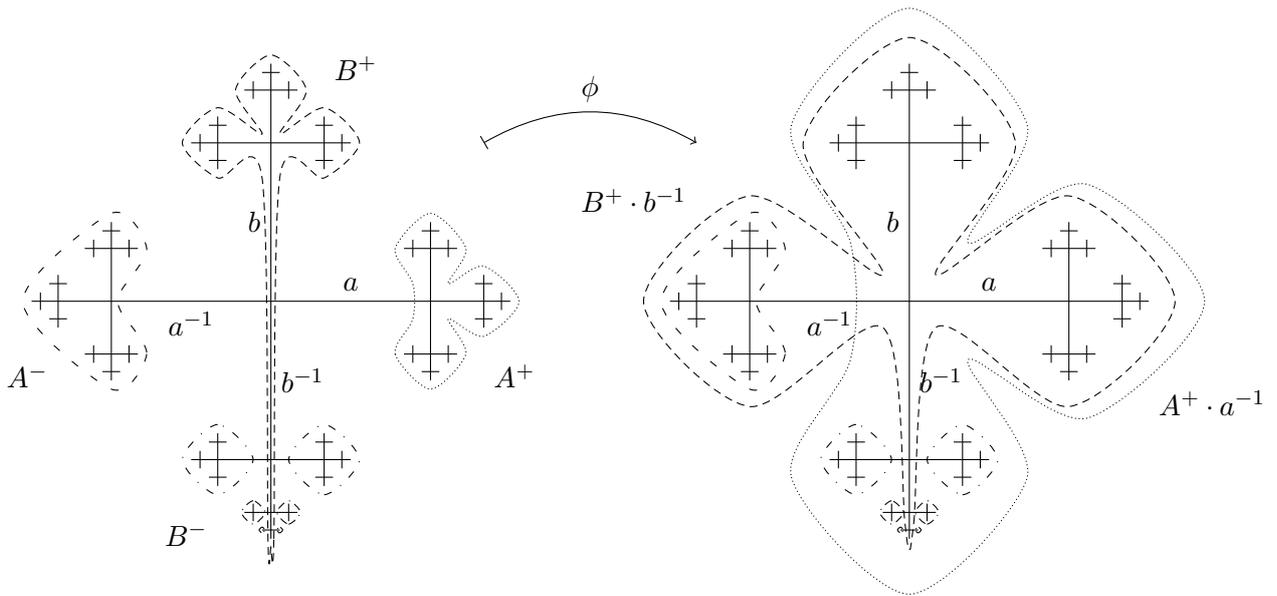

      \myfloatalign
      \begin{wide}
        \figureTreeTarskiFolnerTwoToOneMapPhi
        \caption{The same part of the right $\setOf{a, b, a^{-1}, b^{-1}}$-Cayley graph of $F_2$ is depicted two times, where each time the vertex in the centre is the neutral element $e_{F_2}$. The vertices of the graph on the left enclosed by a contour are mapped by $\phi$ to the vertices of the graph on the right enclosed by the contour with the same pattern and a similar shape. As can be seen, the images of $A^{-}$ and $A^{+}$ under $\phi$ are disjoint and cover $M$, and so are and do the images of $B^{-}$ and $B^{+}$. In other words, each vertex of the graph on the right is enclosed by precisely two contours, where one encloses the vertices of the image of $A^{-}$ or $A^{+}$ under $\phi$ and the other one the vertices of the image of $B^{-}$ or $B^{+}$.}
        \label{figure:tree:Tarski-Folner:2-to-1-map:phi}
      \end{wide}
    \end{figure}

    The maps
    \begin{align*}
      \psi \from M &\to M,\\
      m &\mapsto \begin{dcases*}
                  m, &if $m \in A^{-}\ (= A_{e_{F_2}})$,\\
                  m \cdot a, &if $m \in A^{+} \cdot a^{-1}\ (= A_a)$,\\
                \end{dcases*}
    \end{align*}
    and
    \begin{align*}
      \psi' \from M &\to M,\\
      m &\mapsto \begin{dcases*}
                  m, &if $m \in B^{-}\ (= B_{e_{F_2}})$,\\
                  m \cdot b, &if $m \in B^{+} \cdot b^{-1}\ (= B_b)$,\\
                \end{dcases*}
    \end{align*}
    are well-defined, because the families $\family{A^{-}, A^{+} \cdot a^{-1}}$ and $\family{B^{-}, B^{+} \cdot b^{-1}}$ are partitions of $M$; they are injective, because the sets $A^{-}$ and $A^{+}$ as well as the sets $B^{-}$ and $B^{+}$ are disjoint (see \cref{figure:tree:Tarski-Folner:2-to-1-map:psi-and-psi-prime}); the family $\family{\psi(M), \psi'(M)}$ is a partition of $M$, because the family $\family{A^{-}, A^{+}, B^{-}, B^{+}}$ is a partition of $M$; for each cell $m \in M$, we have $\phi^{-1}(m) = \setOf{\psi(m), \psi'(m)}$, as can be verified; and, for each element $e \in E$, we have $A_e = \setOf{m \in M \suchThat m \cdot e = \psi(m)}$ and $B_e = \setOf{m \in M \suchThat m \cdot e = \psi'(m)}$, as can be verified. Hence, the maps $\psi$ and $\psi'$ are the ones from the subproof \textsc{\ref{item:Tarski-Folner:two-to-one-surjective-map} implies \ref{item:Tarski-Folner:no-paradoxical-decomposition}} of \cref{theorem:Tarski-Folner} that correspond to $\phi$.
    \begin{figure}
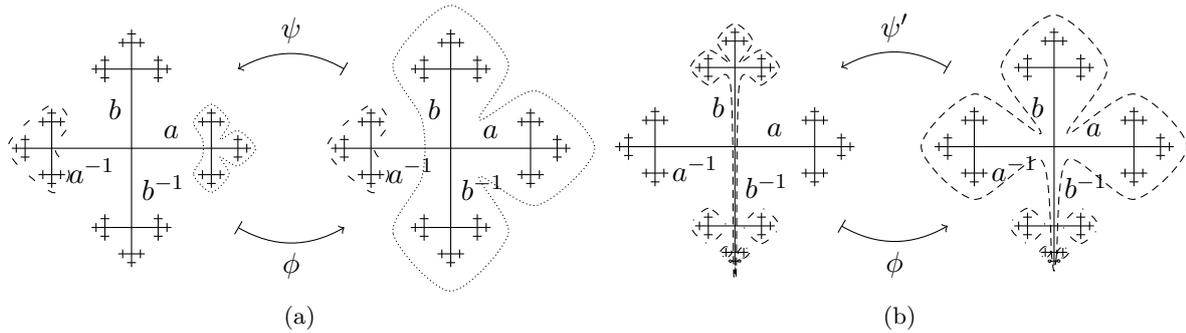

      \myfloatalign
      \begin{wide}
        \figureTreeTarskiFolnerTwoToOneMapPsiAndPsiPrime
        \caption{In each subfigure, the same part of the right $\setOf{a, b, a^{-1}, b^{-1}}$-Cayley graph of $F_2$ is depicted two times, where each time the vertex in the centre is the neutral element $e_{F_2}$; and the vertices of the graph on the right (or left) enclosed by a contour are mapped by $\psi$ or $\psi'$ (or $\phi$) to the vertices of the graph on the left (or right) enclosed by the contour with the same pattern and a similar shape. The map $\psi$ maps $\phi(A^{-})$ to $A^{-}$ and $\phi(A^{+})$ to $A^{+}$, the map $\psi'$ maps $\phi(B^{-})$ to $B^{-}$ and $\phi(B^{+})$ to $B^{+}$, and, as can be seen, the images of $\psi$ and $\psi'$ are disjoint and cover $M$.}
        \label{figure:tree:Tarski-Folner:2-to-1-map:psi-and-psi-prime}
      \end{wide}
    \end{figure}

    The set $E$ contains the neutral element $e_{F_2}$ and, for each finite subset $F$ of $M$, we have $\cardinalityOf{F \cdot E} \geq 2 \cardinalityOf{F}$, because $F = (F \cap A^{-}) \disjointUnionWith (F \cap A^{+}) \disjointUnionWith (F \cap B^{-}) \disjointUnionWith (F \cap B^{+})$ and hence $F \cdot E \supseteq ((F \cap A^{-}) \cdot e_{F_2}) \disjointUnionWith ((F \cap A^{-}) \cdot b) \disjointUnionWith ((F \cap A^{+}) \cdot a) \disjointUnionWith ((F \cap A^{+}) \cdot b) \disjointUnionWith ((F \cap B^{-}) \cdot e_{F_2}) \disjointUnionWith ((F \cap B^{-}) \cdot a) \disjointUnionWith ((F \cap B^{+}) \cdot a) \disjointUnionWith ((F \cap B^{+}) \cdot b)$ (see \cref{figure:tree:Tarski-Folner:finite-subset-E}). Note that $(F \cap A^{-}) \cdot e_{F_2} = \psi(F \cap A^{-})$, $(F \cap A^{-}) \cdot b = \psi'(F \cap A^{-})$, and so on. 
    \begin{figure}
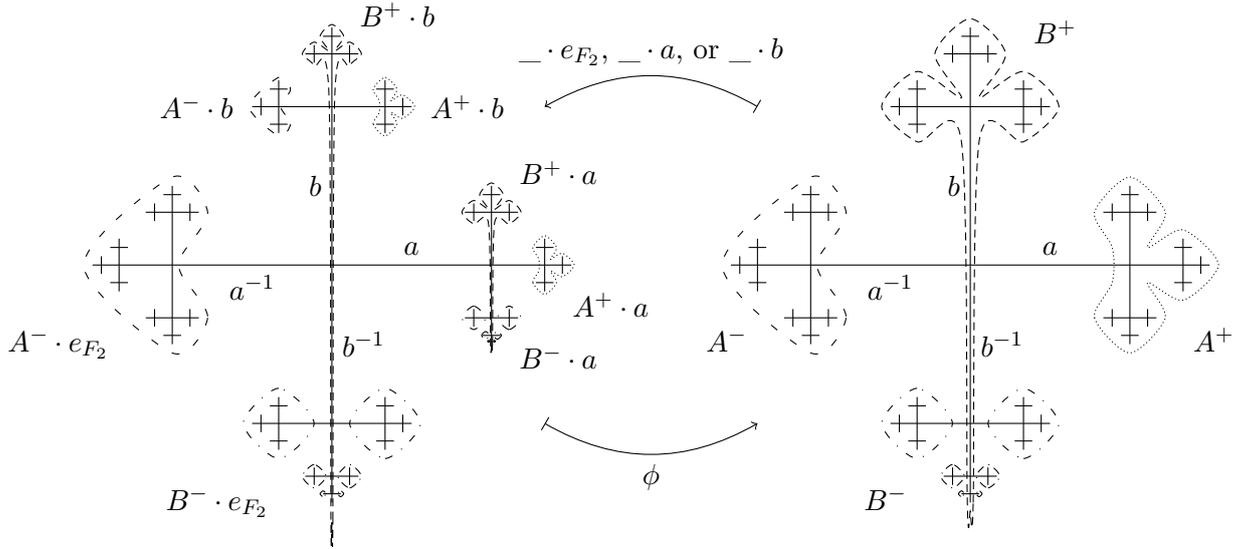

      \myfloatalign
      \begin{wide}
        \figureTreeTarskiFolnerFiniteSubsetE
        \caption{The same part of the right $\setOf{a, b, a^{-1}, b^{-1}}$-Cayley graph of $F_2$ is depicted two times, where each time the vertex in the centre is the neutral element $e_{F_2}$. The vertices of the graph on the right (or left) enclosed by a contour are mapped by $\blank \cdot e_{F_2}$, $\blank \cdot a$, or $\blank \cdot b$ (or $\phi$) to the vertices of the graph on the left (or right) enclosed by a contour with the same pattern and a similar shape. The images $A^{-} \cdot a$, $A^{+} \cdot e_{F_2}$, $B^{-} \cdot b$, and $B^{+} \cdot e_{F_2}$ are not depicted; all other images of $A^{-}$, $A^{+}$, $B^{-}$, and $A^{+}$ are depicted and, as can be seen, they are pairwise disjoint, cover $M$, and each kind of contour from the right graph occurs precisely twice in the left graph. It follows that $k$ vertices are mapped by $\blank \cdot E$ to at least $2k$ vertices, where $k \in \N_0$.} 
        \label{figure:tree:Tarski-Folner:finite-subset-E}
      \end{wide}
    \end{figure}
  \end{example}


  \begin{open-problem}
    According to section~4.5 in \cite{ceccherini-silberstein:coornaert:2010}, amenable groups are closed under taking subspaces, quotients, finite direct products, inductive limits, and arbitrary direct sums. Are right-a\-me\-na\-ble left-ho\-mo\-ge\-neous spaces closed under analogous notions?
  \end{open-problem}

  \section{From Left to Right Amenability}
  \label{section:left-versus-right}

  \paragraph{Introduction.} We begin by stating a general condition for when left amenability implies right amenability, then we take a look at ever more specific situations in which this condition is satisfied, and eventually we present a method to construct right-a\-me\-na\-ble cell spaces from left-a\-me\-na\-ble principal left-ho\-mo\-ge\-neous spaces using outer semi-direct products. The cell spaces yielded by this construction though are in the sense degenerated that a subgroup acts freely and transitively. A method to construct non-de\-gen\-er\-at\-ed right-a\-me\-na\-ble cell spaces is given in \cref{section:construction-of-example}.

  \begin{lemma}
  \label{lemma:general-sufficient-condition-for-left-implies-right-amenability}
    Let $\mathcal{R} = \ntuple{\ntuple{M, G, \actsOnPoint}, \ntuple{m_0, \family{g_{m_0, m}}_{m \in M}}}$ be a cell space and let $H$ be a subgroup of $G$ such that, for each element $\mathfrak{g} \in G \modulo G_0$, there is an element $h \in H$ such that the maps $\blank \isSemiActedUponBy \mathfrak{g}$ and $h \actsOnPoint \blank$ are inverse to each other.
    If $\ntuple{M, H, \actsOnPoint\restrictedTo_{H \times M}}$ is left amenable, then $\mathcal{R}$ is right amenable.
  \end{lemma}

  \begin{proof}
    Let $\mu \in \probabilityMeasuresOn(M)$. Furthermore, let $\mathfrak{g} \in G \modulo G_0$. There is an $h \in H$ such that $\blank \isSemiActedUponBy \mathfrak{g}$ and $h \actsOnPoint \blank$ are inverse to each other. 
    Moreover, let $A \subseteq M$. Because $\blank \isSemiActedUponBy \mathfrak{g} = (h \actsOnPoint \blank)^{-1} = h^{-1} \actsOnPoint \blank$, we have $A \isSemiActedUponBy \mathfrak{g} = h^{-1} \actsOnPoint A$. Therefore,
    \begin{align*}
      (\mu \measureIsKindOfSemiActedUponBy \mathfrak{g})(A)
      &= \mu(A \isSemiActedUponBy \mathfrak{g})\\
      &= \mu(h^{-1} \actsOnPoint A)\\
      &= (h \actsOnMeasure \mu)(A).
    \end{align*}
    Thus, $\mu \measureIsKindOfSemiActedUponBy \mathfrak{g} = h \actsOnMeasure \mu$. Hence, if $\mu$ is $\actsOnMeasure\restrictedTo_{H \times \closedInterval{0, 1}^{\powerSetOf(M)}}$-invariant, then $\mu$ is $\measureIsKindOfSemiActedUponBy$-semi-invariant. In conclusion, if $\ntuple{M, H, \actsOnPoint\restrictedTo_{H \times M}}$ is left amenable, then $\mathcal{R}$ is right amenable. 
  \end{proof}

  \begin{lemma}
  \label{lemma:less-general-sufficient-condition-for-left-implies-right-amenability}
    Let $\mathcal{R} = \ntuple{\ntuple{M, G, \actsOnPoint}, \ntuple{m_0, \family{g_{m_0, m}}_{m \in M}}}$ be a cell space and let $H$ be a subgroup of $G$ such that $G = G_0 H$, for each element $\mathfrak{g} \in G \modulo G_0$, the map $\blank \isSemiActedUponBy \mathfrak{g}$ is injective, 
    \begin{equation*}
      \ForEach h \in H \Holds \blank \isSemiActedUponBy h G_0 = h \actsOnPoint \blank,
    \end{equation*}
    and
    \begin{equation*}
      \ForEach h \in H \ForEach \mathfrak{g} \in G \modulo G_0 \Holds (\blank \isSemiActedUponBy h G_0) \isSemiActedUponBy \mathfrak{g} = \blank \isSemiActedUponBy h \cdot \mathfrak{g}.
    \end{equation*}
    If $\ntuple{M, H, \actsOnPoint\restrictedTo_{H \times M}}$ is left amenable, then $\mathcal{R}$ is right amenable.
  \end{lemma}

  \begin{proof}
    Let $g G_0 \in G \modulo G_0$. Because $g^{-1} \in G = G_0 H$, there is a $g_0 \in G_0$ and there is an $h \in H$ such that $g^{-1} = g_0 h$. Thus, $h = g_0^{-1} g^{-1} \in H$. Hence, for each $m \in M$,
    \begin{align*}
      ((\blank \isSemiActedUponBy g G_0) \after (h \actsOnPoint \blank))(m)
      &= (h \actsOnPoint m) \isSemiActedUponBy g G_0\\
      &= (m \isSemiActedUponBy h G_0) \isSemiActedUponBy g G_0\\
      &= m \isSemiActedUponBy h g G_0\\
      &= m \isSemiActedUponBy g_0^{-1} g^{-1} g G_0\\
      &= m \isSemiActedUponBy G_0\\
      &= m.
    \end{align*}
    Therefore, $h \actsOnPoint \blank$ is right inverse to $\blank \isSemiActedUponBy g G_0$. Hence, $\blank \isSemiActedUponBy g G_0$ is surjective and thus, because it is injective by precondition, bijective. Therefore, $\blank \isSemiActedUponBy g G_0$ and $h \actsOnPoint \blank$ are inverse to each other. In conclusion, according to \cref{lemma:general-sufficient-condition-for-left-implies-right-amenability}, if $\ntuple{M, H, \actsOnPoint\restrictedTo_{H \times M}}$ is left amenable, then $\mathcal{R}$ is right amenable.
  \end{proof}

  \begin{definition}
    Let $G$ be a group. The set
    \begin{equation*}
      \centreOf(G) = \setOf{z \in G \suchThat \ForEach g \in G \Holds z g = g z} \mathnote{center $\centreOf(G)$ of $G$}\index[symbols]{ZGcentre@$\centreOf(G)$}
    \end{equation*}
    is called \define{centre of $G$}.
  \end{definition}

  \begin{remark}
    The centre of $G$ is a subgroup of $G$.
  \end{remark}

  \begin{lemma}
  \label{lemma:sufficient-conditions-such-that-left-amenable-implies-right-amenable}
    Let $\mathcal{R} = \ntuple{\ntuple{M, G, \actsOnPoint}, \ntuple{m_0, \family{g_{m_0, m}}_{m \in M}}}$ be a cell space and let $H$ be a subgroup of $G$ such that $G$ is equal to $G_0 H$, $\actsOnPoint\restrictedTo_{H \times M}$ is free, and $\setOf{g_{m_0, m} \suchThat m \in M}$ is included in $\centreOf(H)$ (in particular, $\actsOnPoint\restrictedTo_{\centreOf(H) \times M}$ is transitive). If $\ntuple{M, H, \actsOnPoint\restrictedTo_{H \times M}}$ is left amenable, then $\mathcal{R}$ is right amenable.
  \end{lemma}

  \begin{proof}
    Let $g \in G$. For each $m \in M$,
    \begin{equation*}
      m \isSemiActedUponBy g G_0
      = g_{m_0, m} g \actsOnPoint m_0
      = g_{m_0, m} \actsOnPoint (g \actsOnPoint m_0).
    \end{equation*}
    Let $m \in M$. For each $m' \in M$, because $\actsOnPoint\restrictedTo_{\centreOf(H) \times M}$ is free and $g_{m_0, m}$, $g_{m_0, m'} \in \centreOf(H)$,
    \begin{align*}
      m' \isSemiActedUponBy g G_0 = m \isSemiActedUponBy g G_0
      &\ifAndOnlyIf g_{m_0, m'} = g_{m_0, m}\\
      &\ifAndOnlyIf m' = m.
    \end{align*}
    Therefore, $\blank \isSemiActedUponBy g G_0$ is injective. 

    Let $m \in M$ and let $h \in H$. Because $g_{m_0, m} \in \centreOf(G)$,
    \begin{align*}
      m \isSemiActedUponBy h G_0
      &= g_{m_0, m} h \actsOnPoint m_0\\
      &= h g_{m_0, m} \actsOnPoint m_0\\
      &= h \actsOnPoint m.
    \end{align*}
    Put $m' = m \isSemiActedUponBy h G_0$. Then,
    \begin{align*}
      g_{m_0, m} h \actsOnPoint m_0
      &= h g_{m_0, m} \actsOnPoint m_0\\
      &= h \actsOnPoint m\\
      &= m'.
    \end{align*}
    Hence, because $g_{m_0, m'} \actsOnPoint m_0 = m'$ also and $\actsOnPoint\restrictedTo_{H \times M}$ is free, $g_{m_0, m'} = g_{m_0, m} h$. Therefore,
    \begin{align*}
      (m \isSemiActedUponBy h G_0) \isSemiActedUponBy g G_0
      &= m' \isSemiActedUponBy g G_0\\
      &= g_{m_0, m'} g \actsOnPoint m_0\\
      &= g_{m_0, m} h g \actsOnPoint m_0\\
      &= m \isSemiActedUponBy h g G_0.
    \end{align*}
    In conclusion, according to \cref{lemma:less-general-sufficient-condition-for-left-implies-right-amenability}, if $\ntuple{M, H, \actsOnPoint\restrictedTo_{H \times M}}$ is left amenable, then $\mathcal{R}$ is right amenable.
  \end{proof}

  \begin{example} 
  \label{example:one-dimensional-affines-right-amenable}
    Let $M = \K$ be a field, let
    \begin{equation*}
      G = \setOf{f \from M \to M, x \mapsto a x + b \suchThat a, b \in M, a \neq 0} 
    \end{equation*}
    be the group of affine functions with composition as group multiplication, and let
    \begin{equation*}
      H = \setOf{f \from M \to M, x \mapsto x + b \suchThat b \in M}
    \end{equation*}
    be the group of translations also with composition as group multiplication. The group $H$ is an abelian subgroup of $G$, which in turn is a non-abelian subgroup of the symmetry group of $M$. Moreover, according to example~4.6.2 and theorem~4.6.3 in \cite{ceccherini-silberstein:coornaert:2010}, the group $G$ is left amenable and hence, according to proposition~4.5.1 in \cite{ceccherini-silberstein:coornaert:2010}, so is its subgroup $H$. Furthermore, the group $G$ acts transitively on $M$ by function application by $\actsOnPoint$ and so does $H$ by $\actsOnPoint\restrictedTo_{H \times M}$, even freely so. Because the groups $G$ and $H$ are left amenable, so are the left group sets $\ntuple{M, G, \actsOnPoint}$ and $\ntuple{M, H, \actsOnPoint\restrictedTo_{H \times M}}$. The stabiliser of $m_0 = 0$ is the group of dilations 
    \begin{equation*}
      G_0 = \setOf{f \from M \to M, x \mapsto a x \suchThat a \in M \smallsetminus \setOf{0}}.
    \end{equation*}
    We have $G = G_0 H$. For each $m \in M$, let
    \begin{align*}
      g_{m_0, m} \from M &\to M,\\
      x &\mapsto x + m,
    \end{align*} 
    be the translation by $m$. Then, $\setOf{g_{m_0, m} \suchThat m \in M}$ is included in $\centreOf(H) = H$. Hence, according to \cref{lemma:sufficient-conditions-such-that-left-amenable-implies-right-amenable}, the cell space $\mathcal{R} = \ntuple{\ntuple{M, G, \actsOnPoint}, \ntuple{m_0, \family{g_{m_0, m}}_{m \in M}}}$ is right amenable.
  \end{example}

  \begin{definition} 
    Let $H$ and $N$ be two groups, let $\phi$ be a group homomorphism from $H$ to $\automorphismsOf(N)$, let $G$ be the Cartesian product $N \times H$, and let
    \begin{align*}
      \cdot \from G \times G &\to G,\\
      ((n, h), (n', h')) &\mapsto (n \phi(h)(n'), h h').
    \end{align*}
    The tuple $(G, \cdot)$ is a group, called \graffito{outer semi-direct product $N \rtimes_\phi H$ of $H$ acting on $N$ by $\phi$}\define{outer semi-direct product of $H$ acting on $N$ by $\phi$}\index{semi-direct product!outer}\index[symbols]{NrtimesphiH@$N \rtimes_\phi H$}, and denoted by $N \rtimes_\phi H$.
  \end{definition}

  \begin{remark}
    Let $G$ be a semi-direct product of $H$ acting on $N$ by $\phi$. The neutral element of $G$ is $(e_N, e_H)$ and, for each element $(n, h) \in G$, the inverse of $(n, h)$ is $(\phi(h^{-1})(n^{-1}), h^{-1})$.
  \end{remark}

  \begin{lemma}
  \label{lemma:semi-direct-product-of-group-set-with-stabiliser}
    Let $\ntuple{M, H, \actsOnPoint_H}$ be a principal left-ho\-mo\-ge\-neous space. Furthermore, let $G_0$ be a group, let $\phi$ be a group homomorphism from $G_0$ to $\automorphismsOf(H)$, let $m_0$ be an element of $M$, for each element $m \in M$, let $h_{m_0,m}$ be the unique element of $H$ such that $h_{m_0,m} \actsOnPoint m_0 = m$, and let
    \begin{align*}
      \actsOnPoint_{G_0} \from G_0 \times M &\to M,\\
      (g_0, m) &\mapsto \phi(g_0)(h_{m_0,m}) \actsOnPoint_H m_0.
    \end{align*}
    Moreover, let $G$ be the outer semi-direct product of $G_0$ acting on $H$ by $\phi$, and let
    \begin{align*}
      \actsOnPoint \from G \times M &\to M,\\
      ((h, g_0), m) &\mapsto h \actsOnPoint_H (g_0 \actsOnPoint_{G_0} m).
    \end{align*}
    The triple $\ntuple{M, G_0, \actsOnPoint_{G_0}}$ is a left group set and the group $G_0$ is the stabiliser of $m_0$ under $\actsOnPoint_{G_0}$. Furthermore, the tuple $\mathcal{R} = \ntuple{\ntuple{M, G, \actsOnPoint}, \ntuple{m_0, \family{(h_{m_0, m}, e_{G_0})}_{m \in M}}}$ is a cell space and the group $\setOf{e_H} \times G_0$ is the stabiliser of $m_0$ under $\actsOnPoint$. Moreover, under the identification of $G_0$ with $\setOf{e_H} \times G_0$ and of $H$ with $H \times \setOf{e_{G_0}}$, the left group sets $\ntuple{M, G_0, \actsOnPoint_{G_0}}$ and $\ntuple{M, H, \actsOnPoint_H}$ are left group subsets of $\ntuple{M, G, \actsOnPoint}$. 
  \end{lemma}

  \begin{proof}
    Because $\phi(e_{G_0}) = \identityMap_{\automorphismsOf(H)}$, for each $m \in M$,
    \begin{align*}
      e_{G_0} \actsOnPoint_{G_0} m
      &= \phi(e_{G_0})(h_{m_0,m}) \actsOnPoint_H m_0\\
      &= h_{m_0,m} \actsOnPoint_H m_0\\
      &= m.
    \end{align*}
    Let $g_0$ and $g_0' \in G_0$, and let $m \in M$. Because $\actsOnPoint_H$ is free and $h_{m_0, \phi(g_0')(h_{m_0,m}) \actsOnPoint_H m_0} \actsOnPoint_H m_0 = \phi(g_0')(h_{m_0,m}) \actsOnPoint_H m_0$, we have $h_{m_0, \phi(g_0')(h_{m_0,m}) \actsOnPoint_H m_0} = \phi(g_0')(h_{m_0,m})$. Therefore,
    \begin{align*}
      g_0 g_0' \actsOnPoint_{G_0} m
      &= \phi(g_0 g_0')(h_{m_0,m}) \actsOnPoint_H m_0\\
      &= (\phi(g_0) \after \phi(g_0'))(h_{m_0,m}) \actsOnPoint_H m_0\\
      &= \phi(g_0)(\phi(g_0')(h_{m_0,m})) \actsOnPoint_H m_0\\
      &= \phi(g_0)(h_{m_0, \phi(g_0')(h_{m_0,m}) \actsOnPoint_H m_0}) \actsOnPoint_H m_0\\
      &= g_0 \actsOnPoint_{G_0} (\phi(g_0')(h_{m_0,m}) \actsOnPoint_H m_0)\\
      &= g_0 \actsOnPoint_{G_0} (g_0' \actsOnPoint_{G_0} m).
    \end{align*}
    In conclusion, $\ntuple{M, G_0, \actsOnPoint_{G_0}}$ is a left group set.

    Because $h_{m_0, m_0} = e_H$, for each $g_0 \in G_0$,
    \begin{align*}
      g_0 \actsOnPoint_{G_0} m_0
      &= \phi(g_0)(e_H) \actsOnPoint_H m_0\\
      &= e_H \actsOnPoint_H m_0\\
      &= m_0.
    \end{align*}
    In conclusion, $G_0$ is the stabiliser of $m_0$ under $\actsOnPoint_{G_0}$.

    For each $m \in M$,
    \begin{equation*}
      (e_H, e_{G_0}) \actsOnPoint m
      = e_H \actsOnPoint_H (e_{G_0} \actsOnPoint_{G_0} m)
      = m.
    \end{equation*}
    Let $g_0 \in G_0$, let $h \in H$, and let $m \in M$. Because $h h_{m_0,m} \actsOnPoint_H m_0 = h \actsOnPoint_H m$, we have $h h_{m_0,m} = h_{m_0, h \actsOnPoint_H m}$. Hence,
    \begin{align*}
      \phi(g_0)(h) \actsOnPoint_H (g_0 \actsOnPoint_{G_0} m)
      &= \phi(g_0)(h) \actsOnPoint_H (\phi(g_0)(h_{m_0,m}) \actsOnPoint_H m_0)\\
      &= \phi(g_0)(h) \phi(g_0)(h_{m_0,m}) \actsOnPoint_H m_0\\
      &= \phi(g_0)(h h_{m_0,m}) \actsOnPoint_H m_0\\
      &= \phi(g_0)(h_{m_0, h \actsOnPoint_H m}) \actsOnPoint_H m_0\\
      &= g_0 \actsOnPoint_{G_0} (h \actsOnPoint_H m).
    \end{align*}
    Therefore, for each $g_0 \in G_0$, each $g_0' \in G_0$, each $h \in H$, each $h' \in H$, and each $m \in M$,
    \begin{align*}
      (h, g_0) (h', g_0') \actsOnPoint m
      &= (h \phi(g_0)(h'), g_0 g_0') \actsOnPoint m\\
      &= h \phi(g_0)(h') \actsOnPoint_H (g_0 g_0' \actsOnPoint_{G_0} m)\\
      &= h \actsOnPoint_H \parens[\Big]{\phi(g_0)(h') \actsOnPoint_H \parens[\big]{g_0 \actsOnPoint_{G_0} (g_0' \actsOnPoint_{G_0} m)}}\\
      &= h \actsOnPoint_H \parens[\Big]{g_0 \actsOnPoint_{G_0} \parens[\big]{h' \actsOnPoint_H (g_0' \actsOnPoint_{G_0} m)}}\\
      &= (h, g_0) \actsOnPoint \parens[\big]{h' \actsOnPoint_H (g_0' \actsOnPoint_{G_0} m)}\\
      &= (h, g_0) \actsOnPoint \parens[\big]{(h', g_0') \actsOnPoint m}.
    \end{align*}
    In conclusion, $\ntuple{M, G, \actsOnPoint}$ is a left group action.

    Because $\actsOnPoint_H$ is transitive and, for each $h \in H$ and each $m \in M$, we have $(h, e_{G_0}) \actsOnPoint m = h \actsOnPoint m$, the left group action $\actsOnPoint$ is transitive and hence $\mathcal{M} = \ntuple{M, G, \actsOnPoint}$ is a left-ho\-mo\-ge\-neous space. Moreover, because, for each $m \in M$,
    \begin{align*}
      (h_{m_0, m}, e_{G_0}) \actsOnPoint m_0
      &= h_{m_0, m} \actsOnPoint_H (e_{G_0} \actsOnPoint_{G_0} m_0)\\
      &= h_{m_0, m} \actsOnPoint_H m_0\\
      &= m,
    \end{align*}
    the tuple $\mathcal{K} = \ntuple{m_0, \family{(h_{m_0, m}, e_{G_0})}_{m \in M}}$ is a coordinate system for $\mathcal{M}$. Therefore, $\mathcal{R} = \ntuple{\mathcal{M}, \mathcal{K}}$ is a cell space.

    Because $G_0$ is the stabiliser of $m_0$ under $\actsOnPoint_{G_0}$, for each $(h, g_0) \in G$, we have $(h, g_0) \actsOnPoint m_0 = h \actsOnPoint_H (g_0 \actsOnPoint m_0) = h \actsOnPoint m_0$. Because $\actsOnPoint_H$ is free, $\setOf{e_H} \times G_0$ is the stabiliser of $m_0$ under $\actsOnPoint$.

    Under the identification of $G_0$ with $\setOf{e_H} \times G_0$ and of $H$ with $H \times \setOf{e_{G_0}}$, we have $\actsOnPoint\restrictedTo_{G_0 \times M} = \actsOnPoint_{G_0}$ and $\actsOnPoint\restrictedTo_{H \times M} = \actsOnPoint_H$.  
  \end{proof}

  \begin{corollary}
  \label{corollary:sufficient-conditions-such-that-left-amenable-implies-right-amenable}
    In the situation of \cref{lemma:semi-direct-product-of-group-set-with-stabiliser}, let $H$ be abelian. The cell space $\mathcal{R}$ is right amenable.
  \end{corollary}

  \begin{proof} 
    According to theorem~4.6.1 in \cite{ceccherini-silberstein:coornaert:2010}, because $H$ is abelian, it is left amenable. Therefore, $\ntuple{M, H, \actsOnPoint_H}$ is left amenable. Identify $G_0$ with $\setOf{e_H} \times G_0$ and identify $H$ with $H \times \setOf{e_{G_0}}$. Then, $H$ is a subgroup of $G$, and $G = H G_0$, and $\actsOnPoint\restrictedTo_{H \times M} = \actsOnPoint_H$ is free, and, for each $m \in M$, we have $(h_{m_0, m}, e_{G_0}) \in H = \centreOf(H)$. Hence, according to \cref{lemma:sufficient-conditions-such-that-left-amenable-implies-right-amenable}, the cell space $\mathcal{R}$ is right amenable.
  \end{proof}

  \begin{example}[Euclidean Space (compare \cref{example:one-dimensional-affines-right-amenable})] 
  \label{example:euclidean-space}
    Let $d$ be a positive integer; let $E$ be the $d$-dimensional Euclidean group, that is, the symmetry group of the $d$-dimensional Euclidean space, in other words, the isometries of $\R^d$ with respect to the Euclidean metric with function composition; let $T$ be the $d$-dimensional translation group; and let $O$ be the $d$-dimensional orthogonal group. The group $T$ is abelian, a normal subgroup of $E$, and isomorphic to $\R^d$ with addition; the group $O$ is isomorphic to the quotient $E \modulo T$ and to the $(d \times d)$-dimensional orthogonal matrices with matrix multiplication; the group $E$ is isomorphic to the outer semi-direct product $T \rtimes_\iota O$, where $\iota \from O \to \automorphismsOf(\R^d)$ is the inclusion map. The groups $T$, $O$, and $E$ act on $\R^d$ on the left by function application, denoted by $\actsOnPoint_T$, $\actsOnPoint_O$, and $\actsOnPoint$ respectively; under the identification of $T$ with $\R^d$ by $t \mapsto [v \mapsto v + t]$, of $O$ with the orthogonal matrices of $\R^{d \times d}$ by $A \mapsto [v \mapsto A v]$, and of $E$ with $T \rtimes_\iota O$ by $(t, A) \mapsto [v \mapsto A v + t]$, we have 
    \begin{align*}
      \actsOnPoint_T \from T \times \R^d &\to \R^d,\\
      (t, v) &\mapsto v + t,
    \end{align*}
    and
    \begin{align*}
      \actsOnPoint_O \from O \times \R^d &\to \R^d,\\
      (A, v) &\mapsto A v,
    \end{align*}
    and
    \begin{align*}
      \actsOnPoint \from E \times \R^d &\to \R^d,\\
      ((t, A), v) &\mapsto A v + t,
    \end{align*}
    and
    \begin{align*}
      \iota \from O &\to \automorphismsOf(\R^d),\\
      A &\mapsto [v \mapsto A v].
    \end{align*}
    Hence, for each vector $v \in \R^d$, we have $v \actsOnPoint_T 0 = v$, therefore, $\actsOnPoint_O = [(A, v) \mapsto \iota(A)(v) \actsOnPoint_T 0]$, and thus $\actsOnPoint = [((t, A), v) \mapsto t \actsOnPoint_T (A \actsOnPoint_O v)]$. Moreover, because the group $(T, \after) \isIsomorphicTo (\R^d, +)$ is abelian, according to theorem~4.6.1 in \cite{ceccherini-silberstein:coornaert:2010}, it is left amenable and so is $\ntuple{\R^d, \R^d, +} \isIsomorphicTo \ntuple{\R^d, T, \actsOnPoint}$. In conclusion, according to \cref{corollary:sufficient-conditions-such-that-left-amenable-implies-right-amenable}, the cell space $\ntuple{\ntuple{\R^d, E, \actsOnPoint}, \ntuple{0, \family{v}_{v \in \R^d}}}$ is right amenable.
  \end{example} 

  \begin{example}[$d$-dimensional Lattice]
    Let $d$ be a positive integer, let $G$ be the symmetry group of $\Z^d$, and let $\actsOnPoint$ be the left group action of $G$ on $\Z^d$ by function application. The group $G$ is isomorphic to the outer semi-direct product $\Z^d \rtimes_\psi O$, where $O$ is the stabiliser of $0$ under $\actsOnPoint$ and $\psi$ is the group homomorphism $O \times \Z^d \to \Z^d$, $(f, v) \mapsto f(v)$. As in \cref{example:euclidean-space}, the cell space $\ntuple{\ntuple{\Z^d, G, \actsOnPoint}, \ntuple{0, \family{(v, \identityMap_{\Z^d})}_{v \in \Z^d}}}$ is right amenable.
  \end{example}

  \clearToOddPage
  \chapter{The Garden of Eden Theorem}
  \label{chapter:garden}

  \paragraph{Abstract.} We prove the Garden of Eden theorem for big-cellular automata with finite set of states and finite neighbourhood over right-a\-me\-na\-ble left-ho\-mo\-ge\-neous spaces with finite stabilisers. It states that the global transition function of such an automaton is surjective if and only if it is pre-injective. Pre-Injectivity means that two global configurations that differ at most on a finite subset and have the same image under the global transition function must be identical. The theorem is proven by showing that the global transition function of an automaton as above is surjective if and only if its image has maximal entropy and that its image has maximal entropy if and only if it is pre-injective. Entropy of a subset of global configurations measures the asymptotic growth rate of the number of finite patterns with growing domains that occur in the subset.

  \paragraph{Remark.} Most parts of this chapter appeared in the paper \enquote{\citetitle*{wacker:garden:2016}}\cite{wacker:garden:2016} and they generalise parts of chapter~5 of the monograph \enquote{\citetitle*{ceccherini-silberstein:coornaert:2010}}\cite{ceccherini-silberstein:coornaert:2010}.

  \paragraph{Summary.} For a right-a\-me\-na\-ble cell space with finite stabilisers we may choose a right Følner net $\mathcal{F} = \net{F_i}_{i \in I}$. The entropy of a subset $X$ of $Q^M$ with respect to $\mathcal{F}$, where $Q$ is a finite set, is, broadly speaking, the asymptotic growth rate of the number of finite patterns with domain $F_i$ that occur in $X$. For subsets $E$ and $E'$ of $G \modulo G_0$, an $\ntuple{E, E'}$-tiling is a subset $T$ of $M$ such that $\family{t \isSemiActedUponBy E}_{t \in T}$ is pairwise disjoint and $\family{t \isSemiActedUponBy E'}_{t \in T}$ is a cover of $M$. If for each point $t \in T$ not all patterns with domain $t \isSemiActedUponBy E$ occur in a subset of $Q^M$, then that subset does not have maximal entropy. 

  The global transition function $\Delta$ of a big-cellular automaton with finite set of states and finite neighbourhood over a right-a\-me\-na\-ble cell space with finite stabilisers, as introduced below, is surjective if and only if its image has maximal entropy. Indeed, if $\Delta$ is surjective, then its image is equal to the set of all global configurations, which has maximal entropy. And, if $\Delta$ is not surjective, then there is a global configuration that is not in its image; thus, because $\Delta$ is continuous and $Q^M$ is compact, where $Q^M$ is equipped with the prodiscrete topology, there is a finite pattern that does not occur in the global configurations of $\Delta(Q^M)$; and hence, $\Delta(Q^M)$ does not have maximal entropy. 

  The image of the global transition function $\Delta$ has maximal entropy if and only if it is pre-injective. Indeed, if $\Delta(Q^M)$ does not have maximal entropy, that is, the asymptotic growth rate of finite patterns in $\Delta(Q^M)$ is less than the one of $Q^M$, then there are two distinct finite patterns with the same domain that can be identically extended to global configurations with the same image under $\Delta$ and thus $\Delta$ is not pre-injective. And, if $\Delta$ is not pre-injective, then there are two distinct finite patterns $p$ and $p'$ with the same domain that have the same image under a restriction of $\Delta$; thus, $\Delta(Q^M)$ is equal to the image of the set $Y$ of all global configurations in which the pattern $p$ does not occur at the cells of a tiling, which is chosen such that its cells are far apart with respect to the domain of $p$; and hence, because the entropy of $\Delta(Y)$ is less than or equal to the one of $Y$ and the entropy of $Y$ is not maximal, the entropy of $\Delta(Q^M)$ is not maximal.

  The previous two paragraphs establish the Garden of Eden theorem, which states that a global transition function as above is surjective if and only if it is pre-injective. This answers a question posed by Sébastien Moriceau at the end of his paper \enquote{\citetitle*{moriceau:2011}}\cite{moriceau:2011}. The Garden of Eden theorem for cellular automata over $\Z^2$ is a famous theorem by Edward Forrest Moore and John R. Myhill from 1962 and 1963, see the papers \enquote{\citetitle*{moore:1962}}\cite{moore:1962} and \enquote{\citetitle*{myhill:1963}}\cite{myhill:1963}.


  \paragraph{Contents.} In \cref{section:patterns} we introduce patterns and blocks, and actions on these. In \cref{section:interiors-closures-and-boundaries} we introduce $E$-interiors, $E$-closures, and $E$-boundaries of subsets of $M$, and characterise right Følner nets using boundaries, which motivates the definition of right Erling nets and tractability. In \cref{section:tilings} we introduce $\ntuple{E, E'}$-tilings of cell spaces, show their existence, and relate them, interiors, and right Erling nets combinatorially. In \cref{section:entropies} we introduce entropies of subsets of $Q^M$, show that applications of global transition functions to subsets of $Q^M$ do not increase entropy, and show that subsets of $Q^M$ that miss a pattern at each cell of a tiling do not have maximal entropy. In \cref{section:gardens-of-Eden} we prove the Garden of Eden theorem by characterising surjectivity and pre-injectivity by maximality of the entropy of the image. And in \cref{section:construction-of-example} we construct non-de\-gen\-er\-at\-ed right-a\-me\-na\-ble left-ho\-mo\-ge\-neous spaces. 

  \paragraph{Preliminary Notions.} A \emph{left group set} is a triple $\ntuple{M, G, \actsOnPoint}$, where $M$ is a set, $G$ is a group, and $\actsOnPoint$ is a map from $G \times M$ to $M$, called \emph{left group action of $G$ on $M$}, such that $G \to \symmetricGroupOf(M)$, $g \mapsto [g \actsOnPoint \blank]$, is a group homomorphism. The action $\actsOnPoint$ is \emph{transitive} if $M$ is non-empty and for each $m \in M$ the map $\blank \actsOnPoint m$ is surjective; and \emph{free} if for each $m \in M$ the map $\blank \actsOnPoint m$ is injective. For each $m \in M$, the set $G \actsOnPoint m$ is the \emph{orbit of $m$}, the set $G_m = (\blank \actsOnPoint m)^{-1}(m)$ is the \emph{stabiliser of $m$}, and, for each $m' \in M$, the set $G_{m, m'} = (\blank \actsOnPoint m)^{-1}(m')$ is the \emph{transporter of $m$ to $m'$}.

  A \emph{left-ho\-mo\-ge\-neous space} is a left group set $\mathcal{M} = \ntuple{M, G, \actsOnPoint}$ such that $\actsOnPoint$ is transitive. A \emph{coordinate system for $\mathcal{M}$} is a tuple $\mathcal{K} = \ntuple{m_0, \family{g_{m_0, m}}_{m \in M}}$, where $m_0 \in M$ and for each $m \in M$ we have $g_{m_0, m} \actsOnPoint m_0 = m$. The stabiliser $G_{m_0}$ is denoted by $G_0$. The tuple $\mathcal{R} = \ntuple{\mathcal{M}, \mathcal{K}}$ is a \emph{cell space}. The set $\setOf{g G_0 \suchThat g \in G}$ of left cosets of $G_0$ in $G$ is denoted by $G \modulo G_0$. The map $\isSemiActedUponBy \from M \times G \modulo G_0 \to M$, $(m, g G_0) \mapsto g_{m_0, m} g g_{m_0, m}^{-1} \actsOnPoint m\ (= g_{m_0, m} g \actsOnPoint m_0)$ is a \emph{right semi-action of $G \modulo G_0$ on $M$ with defect $G_0$}, which means that
  \begin{equation*}
    \ForEach m \in M \Holds m \isSemiActedUponBy G_0 = m,
  \end{equation*}
  and
  \begin{multline*}
    \ForEach m \in M \ForEach g \in G \Exists g_0 \in G_0 \SuchThat \ForEach \mathfrak{g}' \in G \modulo G_0 \Holds\\
          m \isSemiActedUponBy g \cdot \mathfrak{g}' = (m \isSemiActedUponBy g G_0) \isSemiActedUponBy g_0 \cdot \mathfrak{g}'.
  \end{multline*}
  It is \emph{transitive}, which means that the set $M$ is non-empty and for each $m \in M$ the map $m \isSemiActedUponBy \blank$ is surjective; and \emph{free}, which means that for each $m \in M$ the map $m \isSemiActedUponBy \blank$ is injective; and \emph{semi-commutes with $\actsOnPoint$}, which means that
  \begin{multline*}
    \ForEach m \in M \ForEach g \in G \Exists g_0 \in G_0 \SuchThat \ForEach \mathfrak{g}' \in G \modulo G_0 \Holds\\
          (g \actsOnPoint m) \isSemiActedUponBy \mathfrak{g}' = g \actsOnPoint (m \isSemiActedUponBy g_0 \cdot \mathfrak{g}').
  \end{multline*}
  The maps $\iota \from M \to G \modulo G_0$, $m \mapsto G_{m_0, m}$, and $m_0 \isSemiActedUponBy \blank$ are inverse to each other. Under the identification of $M$ with $G \modulo G_0$ by either of these maps, we have $\isSemiActedUponBy \from (m, \mathfrak{g}) \mapsto g_{m_0, m} \actsOnPoint \mathfrak{g}$. (See \cref{chapter:automata}.)

  A left-ho\-mo\-ge\-neous space $\mathcal{M}$ is \emph{right amenable} if there is a coordinate system $\mathcal{K}$ for $\mathcal{M}$ and there is a finitely additive probability measure $\mu$ on $M$ such that 
  \begin{equation*}
    \ForEach \mathfrak{g} \in G \modulo G_0 \ForEach A \subseteq M \Holds \parens[\big]{(\blank \isSemiActedUponBy \mathfrak{g})\restrictedTo_A \text{ injective} \implies \mu(A \isSemiActedUponBy \mathfrak{g}) = \mu(A)},
  \end{equation*}
  in which case the cell space $\mathcal{R} = \ntuple{\mathcal{M}, \mathcal{K}}$ is called \emph{right amenable}. When the stabiliser $G_0$ is finite, that is the case if and only if there is a \emph{right Følner net in $\mathcal{R}$ indexed by $(I, \leq)$}, which is a net $\net{F_i}_{i \in I}$ in $\setOf{F \subseteq M \suchThat F \neq \emptyset, F \text{ finite}}$ such that
  \begin{equation*}
    \ForEach \mathfrak{g} \in G \modulo G_0 \Holds \lim_{i \in I} \frac{\cardinalityOf{F_i \smallsetminus (\blank \isSemiActedUponBy \mathfrak{g})^{-1}(F_i)}}{\cardinalityOf{F_i}} = 0.
  \end{equation*}
  If a net is a right Følner net for one coordinate system, then it is a right Følner net for each coordinate system. In particular, a left-ho\-mo\-ge\-neous space $\mathcal{M}$ with finite stabilisers is right amenable if and only if, for each coordinate system $\mathcal{K}$, the cell space $\ntuple{\mathcal{M}, \mathcal{K}}$ is right amenable. (See \cref{chapter:amenability}.)

  A cell space $\mathcal{R}$ is \emph{finitely and symmetrically right generated} if there is a finite subset $S$ of $G \modulo G_0$ with $G_0 \cdot S \subseteq S$ and $S^{-1} \subseteq S$, where $S^{-1} = \setOf{g^{-1} G_0 \suchThat s \in S, g \in s}$, such that
  \begin{multline*} 
    \ForEach m \in M \Exists k \in \N_0 \Exists \sequence{s_i}_{i \in \setOf{1, 2, \dotsc, k}} \text{ in } S \cup S^{-1} \SuchThat\\
        \parens[\Big]{\parens[\big]{(m_0 \isSemiActedUponBy s_1) \isSemiActedUponBy s_2} \isSemiActedUponBy \dotsb} \isSemiActedUponBy s_k = m.
  \end{multline*}
  The \emph{uncoloured $S$-Cayley graph of $\mathcal{R}$} is the symmetric and $2 \cardinalityOf{S}$-regular directed graph $\mathcal{G} = \ntuple{M, \setOf{(m, m \isSemiActedUponBy s) \suchThat m \in M, s \in S}}$, the \emph{$S$-metric on $\mathcal{R}$} is the distance $\distanceOf$ on $\mathcal{G}$, and the \emph{$S$-length on $\mathcal{R}$} is the map $\lengthOf{\blank} = \distanceOf(m_0, \blank)$. For each $m \in M$ and each $\rho \in \Z$, the \emph{$S$-ball of radius $\rho$ centred at $m$} is the set $\ball(m, \rho) = \setOf{m' \in M \suchThat \distanceOf(m, m') \leq \rho}$, the \emph{$S$-sphere of radius $\rho$ centred at $m$} is the set $\sphere(m, \rho) = \setOf{m' \in M \suchThat \distanceOf(m, m') = \rho}$, the ball $\ball(m_0, \rho)$ is denoted by $\ball(\rho)$, and the sphere $\sphere(m_0, \rho)$ by $\sphere(\rho)$. (See \cref{chapter:growth}.)

  A \emph{semi-cellular automaton} is a quadruple $\mathcal{C} = \ntuple{\mathcal{R}, Q, N, \delta}$, where $\mathcal{R}$ is a cell space; $Q$, called \emph{set of states}, is a set; $N$, called \emph{neighbourhood}, is a subset of $G \modulo G_0$ such that $G_0 \cdot N \subseteq N$; and $\delta$, called \emph{local transition function}, is a map from $Q^N$ to $Q$. A \emph{local configuration} is a map $\ell \in Q^N$ and a \emph{global configuration} is a map $c \in Q^M$. The stabiliser $G_0$ acts on $Q^N$ on the left by $\bullet \from G_0 \times Q^N \to Q^N$, $(g_0, \ell) \mapsto [n \mapsto \ell(g_0^{-1} \cdot n)]$, and the group $G$ acts on $Q^M$ on the left by $\actsOnMap \from G \times Q^M \to Q^M$, $(g, c) \mapsto [m \mapsto c(g^{-1} \actsOnPoint m)]$. The \emph{global transition function of $\mathcal{C}$} is the map $\Delta \from Q^M \to Q^M$, $c \mapsto [m \mapsto \delta(n \mapsto c(m \isSemiActedUponBy n))]$.

  A subgroup $H$ of $G$ is \emph{$\mathcal{K}$-big} if the set $\setOf{g_{m_0, m} \suchThat m \in M}$ is included in $H$. A \emph{big-cellular automaton} is a semi-cellular automaton $\mathcal{C} = \ntuple{\mathcal{R}, Q, N, \delta}$ such that, for some $\mathcal{K}$-big subgroup $H$ of $G$, the local transition function $\delta$ is \emph{$\bullet_{G_0 \cap H}$-invariant}, which means that, for each $h_0 \in G_0 \cap H$, we have $\delta(h_0 \bullet \blank) = \delta(\blank)$. Its global transition function is $\actsOnMap_H$-e\-qui\-var\-i\-ant, which means that, for each $h \in H$, we have $\Delta(h \actsOnMap \blank) = h \actsOnMap \Delta(\blank)$. Note that each $\mathcal{K}$-big subgroup of $G$ includes the subgroup of $G$ generated by $\setOf{g_{m_0, m} \suchThat m \in M}$ and that hence a semi-cellular automaton is a big-cellular automaton if and only if its local transition function is $\bullet_{G_0 \cap \groupGeneratedBy{g_{m_0, m} \suchThat m \in M}}$-invariant. (See \cref{chapter:automata}.)

  \paragraph{Context.} In this chapter, for each subset $A$ of $M$, let $\pi_A$ be the restriction map $Q^M \to Q^A$, $c \mapsto c\restrictedTo_A$.

  \section{Patterns}
  \label{section:patterns}

  In this section, let $\mathcal{R} = \ntuple{\ntuple{M, G, \actsOnPoint}, \ntuple{m_0, \family{g_{m_0, m}}_{m \in M}}}$ be a cell space and let $Q$ be a set.

  \paragraph{Contents.} A pattern is a map from a subset of cells to the set of states (see \cref{definition:pattern}), its size is the cardinality of its domain (see \cref{definition:pattern:size}), it is empty if its domain is empty (see \cref{definition:pattern:empty}), it is finite and called \emph{block} if its domain is finite (see \cref{definition:pattern:finite,definition:block}), and restricting its domain yields a subpattern (see \cref{definition:pattern:subpattern}). The group of symmetries acts on the set of all patterns on the left (see \cref{definition:induced-left-action-on-patterns}) and a pattern centred at the origin can be shifted to a new centre by sort of acting on the new centre (see \cref{definition:induced-right-semi-action}). A pattern occurs in a global configuration if a translation of it coincides with a part of the configuration (see \cref{definition:occurs}).

  \begin{definition}
  \label{definition:pattern}
    Let $A$ be a subset of $M$ and let $p$ be a map from $A$ to $Q$. The map $p$ is called \defineX{$A$-pattern}{pattern@$A$-pattern}\graffito{$A$-pattern $p$}\index[symbols]{pattern@$A$-pattern} and the set $\domainOf(p) = A$ is called \define{domain of $p$}\graffito{domain $\domainOf(p)$ of $p$}\index[symbols]{domp@$\domainOf(p)$}. 
  \end{definition}

  \begin{definition}
  \label{definition:pattern:size}
    Let $p$ be an $A$-pattern. The cardinal number $\sizeOf{p} = \cardinalityOf{A}$ is called \define{size of $p$}\graffito{size $\sizeOf{p}$ of $p$}\index[symbols]{absp@$\cardinalityOf{p}$}.
  \end{definition}


  \begin{definition}
  \label{definition:pattern:empty}
    The $\emptyset$-pattern is called \define{empty}\graffito{empty pattern $\varepsilon$}\index[symbols]{epsilon@$\varepsilon$}. 
  \end{definition}

  \begin{remark}
    The empty pattern is the only one of size $0$.
  \end{remark}

  \begin{definition}
  \label{definition:pattern:finite}
    Let $u$ be an $F$-pattern. It is called \define{finite}\graffito{finite pattern $u$}\index{pattern!finite} and \defineX{$F$-block}{block@$F$-block}\graffito{$F$-block $u$} if and only if its domain $F$ is finite.
  \end{definition}

  \begin{definition}
  \label{definition:block}
    The set $\setOf{u \in Q^F \suchThat F \subseteq M \text{ finite}}$\graffito{set $Q^*$ of all blocks} of all blocks is denoted by $Q^*$\index[symbols]{Qstar@$Q^*$}.
  \end{definition}

  \begin{definition} 
  \label{definition:pattern:subpattern}
    Let $p$ be an $A$- and let $p'$ be an $A'$-pattern. The pattern $p$ is called \define{subpattern of $p'$}\graffito{subpattern $p$ of $p'$}\index{pattern!sub-} if and only if $A \subseteq A'$ and $p = p'\restrictedTo_A$.
  \end{definition}




  \begin{definition}
  \label{definition:induced-left-action-on-patterns}
    The group $G$ acts on the set of patterns on the left by 
    \begin{align*}
      \actsOnMap \from G \times \bigcup_{A \subseteq M} Q^A &\to \bigcup_{A \subseteq M} Q^A, \mathnote{left group action $\actsOnMap$ of $G$ on $\bigcup_{A \subseteq M} Q^A$}\index[symbols]{arrowrightblack@$\actsOnMap$}\\
      (g, p) &\mapsto \left[
                        \begin{aligned}
                          g \actsOnPoint \domainOf(p) &\to Q,\\
                          m &\mapsto p(g^{-1} \actsOnPoint m). \qedhere
                        \end{aligned}
                      \right]
    \end{align*}
  \end{definition}


  \begin{remark}
    The restriction $\actsOnMap_{G \times Q^M \to Q^M}$ is the left group action on the set of global configurations that was introduced in \cref{definition:global-configuration} and that we so far also denoted by the symbol $\actsOnMap$.
  \end{remark}

  \begin{definition}
  \label{definition:induced-right-semi-action}
    Identify $M$ with $G \modulo G_0$ by $\iota \givenBy m \mapsto G_{m_0, m}$. The map
    \begin{align*}
      \actsByItsCoordinateOn \from M \times \bigcup_{A \subseteq M} Q^A &\to \bigcup_{A \subseteq M} Q^A, \mathnote{kind of left semi-action $\protect\actsByItsCoordinateOn$ of $M$ by coordinates on $\bigcup_{A \subseteq M} Q^A$}\index[symbols]{arrowleftblackunderscore@$\protect\actsByItsCoordinateOn$}\\
      (m, p) &\mapsto \left[
                         \begin{aligned}
                           m \isSemiActedUponBy \domainOf(p) &\to Q,\\ 
                           m \isSemiActedUponBy a &\mapsto p(a),
                         \end{aligned}
                       \right]
    \end{align*}
    broadly speaking, maps a point $m$ and a pattern $p$ that is centred at $m_0$ to the corresponding pattern centred at $m$ and, as we see in \cref{remark:induced-right-semi-action}, it is a kind of left semi-action of the set of cells by coordinates on the set of patterns.
  \end{definition}

  \begin{remark} 
  \label{remark:induced-right-semi-action}
    Let $p$ be an $A$-pattern and let $m$ be an element of $M$. Then, $m \actsByItsCoordinateOn p = g_{m_0, m} \actsOnMap p$ and $\domainOf(m \actsByItsCoordinateOn p) = m \isSemiActedUponBy \domainOf(p) = g_{m_0, m} \actsOnPoint \domainOf(p)$.
  \end{remark}

  \begin{remark}
    Each global transition function $\Delta$ of a big-cellular automaton is \defineX{$\actsByItsCoordinateOn$-e\-qui\-var\-i\-ant}{equivariant inducedrightsemiaction@$\protect\actsByItsCoordinateOn$-e\-qui\-var\-i\-ant}\graffito{$\actsByItsCoordinateOn$-e\-qui\-var\-i\-ant}, which means that
    \begin{equation*}
      \ForEach m \in M \ForEach c \in Q^M \Holds \Delta(m \actsByItsCoordinateOn c) = m \actsByItsCoordinateOn \Delta(c).
    \end{equation*}
    Note however that $\actsByItsCoordinateOn$-equivariance is only equivalent to $\actsOnMap_H$-equivariance if the subgroup $H$ of $G$ is equal to the subset $\setOf{g_{m_0, m} \suchThat m \in M}$ of $G$, which is in general not a subgroup of $G$.
  \end{remark}

  \begin{definition} 
  \label{definition:occurs}
    Identify $M$ with $G \modulo G_0$ by $\iota \givenBy m \mapsto G_{m_0, m}$, let $p$ be an $A$-pattern, let $c$ be a $M$-pattern, and let $m$ be an element of $M$. The pattern $p$ is said to \define{occur at $m$ in $c$}\graffito{$p$ occurs at $m$ in $c$} and we write $p \occursIn_m c$\graffito{$p \occursIn_m c$}\index[symbols]{psquaresubseteqmpprime@$p \occursIn_m p'$} if and only if $m \actsByItsCoordinateOn p = c\restrictedTo_{m \isSemiActedUponBy A}$.
  \end{definition}

  \section{Interiors, Closures, and Boundaries; Right Følner Nets and Right Erling Nets} 
  \sectionmark{Interiors, Closures, and Boundaries}
  \label{section:interiors-closures-and-boundaries}

  In this section, let $\mathcal{R} = \ntuple{\mathcal{M}, \mathcal{K}} = \ntuple{\ntuple{M, G, \actsOnPoint}, \ntuple{m_0, \family{g_{m_0, m}}_{m \in M}}}$ be a cell space.

  \paragraph{Contents.} In \cref{definition:interior-closure-boundary} we introduce $E$-interiors, $E$-closures, and (internal\slash external) $E$-boundaries of subsets $A$ of $M$, where the subset $E$ of $G \modulo G_0$ determines thickness and shape of the subtraction from, addition to, or boundary of $A$. In \cref{lemma:properties-of-interior-closure-and-boundary} we show essential properties of interiors, closures, and boundaries. In \cref{lemma:Delta-X-A-minus-plus-are-surjective} we define surjective restrictions $\Delta_{X, A}^-$ and $\Delta_A^-$ of global transition functions to patterns. In \cref{theorem:boundary-characterisation-of-Folner-net} we show that right Følner nets are those nets whose components are asymptotically invariant under taking finite boundaries. In \cref{definition:Erling-net,definition:tractable} we introduce right Erling nets and right tractability, which are weak variants of right Følner nets and right amenability. And in \cref{lemma:finitely-right-generated-is-tractable} we show that each finitely and symmetrically right-gen\-er\-at\-ed cell space is right tractable.

  \begin{definition} 
  \label{definition:interior-closure-boundary}
    Let $A$ be a subset of $M$ and let $E$ be a subset of $G \modulo G_0$.
    \begin{aenumerate}
      \item The set 
            \begin{equation*}
              A^{-E} = \setOf{m \in M \suchThat m \isSemiActedUponBy E \subseteq A} \quad \parens[\big]{= \bigcap_{e \in E} \bigcup_{a \in A} (\blank \isSemiActedUponBy e)^{-1}(a)} \mathnote{$E$-interior $A^{-E}$ of $A$}
            \end{equation*}
            is called \defineX{$E$-interior of $A$}{interior of $A$@$E$-interior of $A$}.
      \item The set
            \begin{equation*} 
              A^{+E} = \setOf{m \in M \suchThat (m \isSemiActedUponBy E) \cap A \neq \emptyset} \quad \parens[\big]{= \bigcup_{e \in E} \bigcup_{a \in A} (\blank \isSemiActedUponBy e)^{-1}(a)} \mathnote{$E$-closure $A^{+E}$ of $A$}
            \end{equation*}
            is called \defineX{$E$-closure of $A$}{closure of $A$@$E$-closure of $A$}.
      \item The set
            \begin{equation*}
              \boundaryOf_E A = A^{+E} \smallsetminus A^{-E}
              \mathnote{$E$-boundary $\boundaryOf_E A$ of $A$}
            \end{equation*}
            is called \defineX{$E$-boundary of $A$}{boundary of $A$@$E$-boundary of $A$}.
      \item The set
            \begin{equation*}
              \boundaryOf_E^- A = A \smallsetminus A^{-E}
              \mathnote{internal $E$-boundary $\boundaryOf_E^- A$ of $A$}
            \end{equation*}
            is called \define{internal $E$-boundary of $A$}\index{boundary of $A$!internal}.
      \item The set
            \begin{equation*}
              \boundaryOf_E^+ A = A^{+E} \smallsetminus A
              \mathnote{external $E$-boundary $\boundaryOf_E^+ A$ of $A$}
            \end{equation*}
            is called \define{external $E$-boundary of $A$}\index{boundary of $A$!external}. \qedhere 
    \end{aenumerate}
  \end{definition}

  \begin{remark}
  \label{remark:group:interior-closure-boundary}
    Let $\mathcal{R}$ be the cell space $\ntuple{\ntuple{G, G, \cdot}, \ntuple{e_G, \family{g}_{g \in G}}}$, where $G$ is a group and $e_G$ is its neutral element. Then, $G_0 = \setOf{e_G}$ and $\isSemiActedUponBy = \cdot$. Hence, the notions of $E$-interior, $E$-closure, and $E$-boundary are the same as the ones defined in paragraph~2 of section~5.4 in \cite{ceccherini-silberstein:coornaert:2010}.
  \end{remark}

  \begin{example}[Lattice] 
  \label{example:lattice:interior-closure-and-boundary}
    In the situation of \cref{example:lattice},
    let $A$ be the ball $\ball(2)$. The $\setOf{(1, 0)}$-interior of $A$ is the ball $A - (1, 0) = \ball((-1, 0), 2)$, which is not included in $A$; the $\setOf{(0, 0), (1, 0)}$-interior of $A$ is the set $A \cap (A - (1, 0)) = \ball_{\R^2}((-1/2, 0), 3/2) \cap M$, which is included in $A$ on the left; the $\setOf{(-1, 0), (0, 0), (1, 0)}$-interior of $A$ is the ball $(A + (1, 0)) \cap A \cap (A - (1, 0)) = \ball(1)$, which is included in $A$ in the middle; and the $E$-interior of $A$, where $E = \setOf{(-1, 0),\allowbreak (0, -1),\allowbreak (0, 0),\allowbreak (0, 1),\allowbreak (1, 0)} = \ball(1)$, is also the ball $\ball(1)$ (see \cref{figure:lattice:interior-closure-boundary}).

    The $\setOf{(1, 0)}$-closure of $A$ is the ball $A - (1, 0) = \ball((-1, 0), 2)$, which does not include $A$; the $\setOf{(0, 0), (1, 0)}$-closure of $A$ is the set $A \cup (A - (1, 0)) = \ball(2) \cup \ball((-1, 0), 2)$, which includes $A$ on the right; the $\setOf{(-1, 0), (0, 0), (1, 0)}$-closure of $A$ is the set $(A + (1, 0)) \cup A \cup (A - (1, 0)) = \ball((1, 0), 2) \cup \ball(2) \cup \ball((-1, 0), 2)$, which includes $A$ in the middle; and the $E$-interior of $A$ is the ball $\ball(3)$ (see \cref{figure:lattice:interior-closure-boundary}).

    The $\setOf{(1, 0)}$-boundary of $A$ is the empty set $\emptyset$, although $\setOf{(1, 0)}$ and $A$ are non-empty; the $\setOf{(0, 0), (1, 0)}$-boundary of $A$ is the set $(\ball_{\R^2}((-1/2, 0), 5/2) \smallsetminus \ball_{\R^2}((-1/2, 0), 3/2)) \cap M$, which includes $A$ on the right; the $\setOf{(-1, 0), (0, 0), (1, 0)}$-boundary of $A$ is the set $(\ball((1,\allowbreak 0), 2) \cup \ball(2) \cup \ball((-1, 0), 2)) \smallsetminus \ball(1)$, which includes $A$ in the middle; and the $E$-boundary of $A$ is the thickened sphere $\ball(3) \smallsetminus \ball(1)$ (see \cref{figure:lattice:interior-closure-boundary}).

    The above calculations suggest that the notions of $E$-interior, -closure and -boundary behave best if the set $E$ contains $(0, 0)$ and is invariant under $G_0$, which we also see in the forthcoming \cref{lemma:properties-of-interior-closure-and-boundary}.

    In general, for each cell $m \in M$ and for each pair $(\rho, \varrho) \in \N_0 \times \N_0$, we have $\ball(m, \rho)^{-\ball(\varrho)} = \ball(m, \rho - \varrho)$, $\ball(m, \rho)^{+\ball(\varrho)} = \ball(m, \rho + \varrho)$, and $\boundaryOf_{\ball(\varrho)} \ball(m, \rho) = \ball(m, \rho + \varrho) \smallsetminus \ball(m, \rho - \varrho)$. 
    \begin{figure}
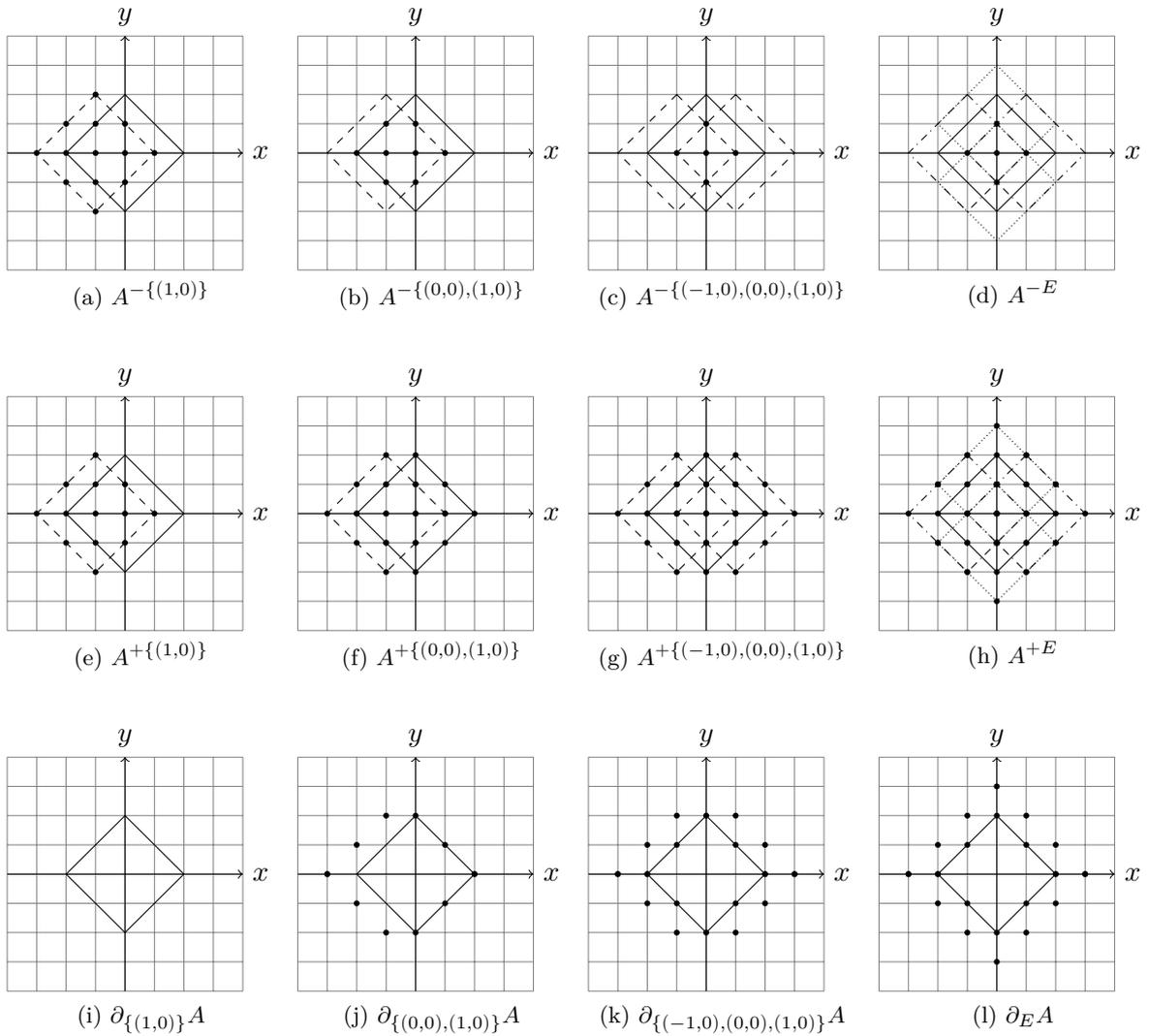

      \myfloatalign
      \figureLatticeInteriorClosureBoundary
      \caption{In each subfigure, the whole space is $\R^2$, the grid points are elements of $M = \Z^2$, the grid points in the region enclosed by the diamond with solid border are the elements of $A = \ball(2)$, the grid points in the regions enclosed by diamonds with dashed, dash-dotted, or dotted borders are the elements of the balls that occur in the calculation of the interior, closure, or boundary in \cref{example:lattice:interior-closure-and-boundary}, and the dots are the elements of the respective interior, closure, or boundary of $A$, where $E$ is the set $\setOf{(-1,0), (0,-1), (0,0), (0,1), (1,0)}$. The subfigures in the first row depict interiors, the ones in the second closures, and the ones in the third boundaries.}
      \label{figure:lattice:interior-closure-boundary}
    \end{figure}
  \end{example}

  \begin{example}[Tree] 
  \label{example:tree:interior-closure-and-boundary} 
    In the situation of \cref{example:tree}, let $A$ be the ball $\ball(2)$. The $\setOf{a}$-interior of $A$ is the set $A a^{-1}$, the $\setOf{a, b}$-interior of $A$ is the ball $A a^{-1} \cap A b^{-1} = \ball(1)$, and the $\setOf{a, a^{-1}}$- and $\setOf{a, b, a^{-1}, b^{-1}}$-interiors of $A$ are also the ball $\ball(1)$. The $\setOf{a}$-closure of $A$ is the set $A a^{-1}$, the $\setOf{a, b}$-closure of $A$ is the set $A a^{-1} \cup A b^{-1}$, the $\setOf{a, a^{-1}}$-closure of $A$ is the set $A a^{-1} \cup A a$, and the $\setOf{a, b, a^{-1}, b^{-1}}$-closure of $A$ is the ball $\ball(3)$. The $\setOf{a}$-boundary of $A$ is the empty set $\emptyset$, the $\setOf{a, b}$-boundary of $A$ is the symmetric difference $A a^{-1} \symmetricSetDifferenceOf A b^{-1}$, the $\setOf{a, a^{-1}}$-boundary of $A$ is the symmetric difference $A a^{-1} \symmetricSetDifferenceOf A a$, and the $\setOf{a, b, a^{-1}, b^{-1}}$-boundary of $A$ is the sphere $\ball(3) \smallsetminus \ball(1)$ (see \cref{figure:tree:interior-closure-boundary}).

    The above calculations suggest that the notions of $E$-interior, -closure and -boundary behave best if the set $E$ is invariant under $G_0$, which we also see in the forthcoming \cref{lemma:properties-of-interior-closure-and-boundary}.

    In general, for each cell $m \in M$ and for each pair $(\rho, \varrho) \in \N_0 \times \N_0$, we have $\ball(m, \rho)^{-\ball(\varrho)} = \ball(m, \rho - \varrho)$, $\ball(m, \rho)^{+\ball(\varrho)} = \ball(m, \rho + \varrho)$, and $\boundaryOf_{\ball(\varrho)} \ball(m, \rho) = \ball(m, \rho + \varrho) \smallsetminus \ball(m, \rho - \varrho)$. 
    \begin{figure}
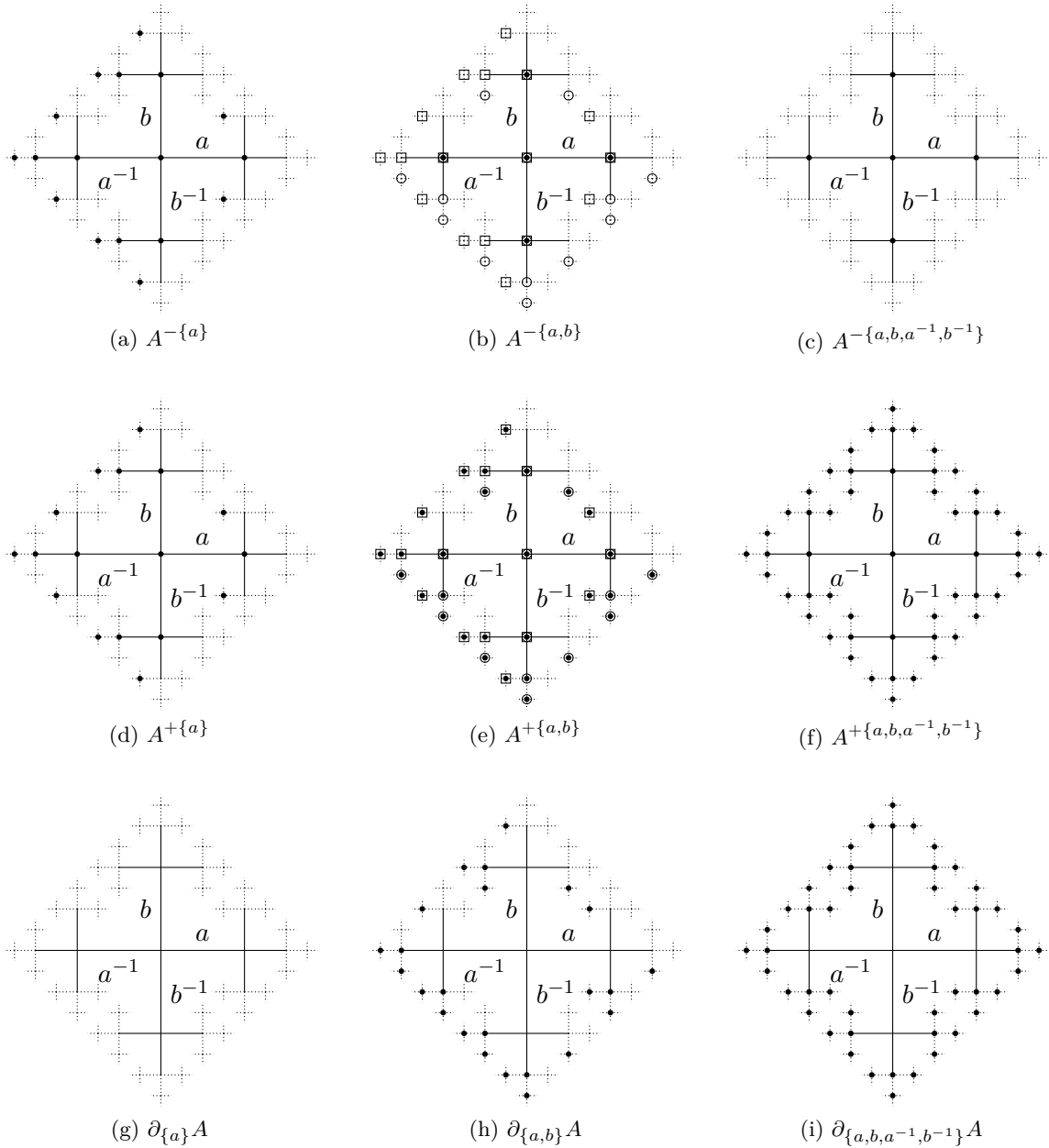
 
      \myfloatalign
      \figureTreeInteriorClosureBoundary
      \caption{In each subfigure, the same part of the $\setOf{a, b, a^{-1}, b^{-1}}$-Cayley graph of $F_2$ is depicted, the vertex in the centre is the neutral element $e_{F_2}$, the vertices that are adjacent to solid edges are the elements of $A = \ball(2)$, the dots are the elements of the respective
        interior, closure, or boundary of $A$, the squares are the elements of $A a^{-1}$, and the circles
        are the elements of $A b$.
        The subfigures in the first row depict interiors, the ones in the second closures, and the ones in the third boundaries.} 
      \label{figure:tree:interior-closure-boundary}
    \end{figure}
  \end{example}


  \begin{example}[Sphere] 
  \label{example:sphere:interior-closure-and-boundary}
    In the situation of \cref{example:sphere:liberation}, let $A$ be a curved circular disk of radius $3 \rho$ with the north pole $m_0$ at its centre, let $g$ be the rotation about an axis $a$ in the $(x,y)$-plane by $\rho$ radians, let $E$ be the set $\setOf{g_0 g G_0 \suchThat g_0 \in G_0}$, and, for each point $m \in M$, let $E_m$ be the set $m \isSemiActedUponBy E$. Because $G_0$ is the set of rotations about the $z$-axis and $m_0 \isSemiActedUponBy E = g_{m_0, m_0} G_0 g \actsOnPoint m_0 = G_0 \actsOnPoint (g \actsOnPoint m_0)$, the set $E_{m_0}$ is the boundary of a curved circular disk of radius $\rho$ with the north pole $m_0$ at its centre. And, for each point $m \in M$, because $m \isSemiActedUponBy E = g_{m_0, m} \actsOnPoint E_{m_0}$, the set $E_m$ is the boundary of a curved circular disk of radius $\rho$ with $m$ at its centre. 

    The $E$-interior of $A$ is the curved circular disk of radius $2 \rho$ with the north pole $m_0$ at its centre. The $E$-closure of $A$ is the curved circular disk of radius $4 \rho$ with the north pole $m_0$ at its centre. And the $E$-boundary of $A$ is the annulus bounded by the boundaries of the $E$-interior and the $E$-closure of $A$. 
  \end{example}

  Essential properties of and relations between interiors, closures, and boundaries are given in 

  \begin{lemma} 
  \label{lemma:properties-of-interior-closure-and-boundary}
    Let $A$ and $A'$ be two subsets of $M$, let $\family{A_i}_{i \in I}$ be a family of subsets of $M$, let $e$ be an element of $G \modulo G_0$, and let $E$ and $E'$ be two subsets of $G \modulo G_0$. Furthermore, for each element $e' \in E'$, let $e' \cdot E$ be the set $\setOf{g' \cdot e \suchThat e \in E, g' \in e'}$, let $E' \cdot E$ be the set $\bigcup_{e' \in E'} e' \cdot E\ (= \setOf{g' \cdot e \suchThat e \in E, e' \in E', g' \in e'})$, let $E^{-1}$ be the set $\setOf{g^{-1} G_0 \suchThat e \in E, g \in e}$, and let $(E')^{-1}$ be the set $\setOf{(g')^{-1} G_0 \suchThat e' \in E', g' \in e'}$. 
    \begin{aenumerate}
      \item \label{item:properties-of-interior-closure-and-boundary:only-neutral-element} 
            $A^{-\setOf{G_0}} = A$, $A^{+\setOf{G_0}} = A$, and $\boundaryOf_{\setOf{G_0}} A = \emptyset$.
      \item \label{item:properties-of-interior-closure-and-boundary:only-neutral-element-and-another} 
            $A^{-\setOf{G_0, e}} = A \cap (\blank \isSemiActedUponBy e)^{-1}(A)$, $A^{+\setOf{G_0, e}} = A \cup (\blank \isSemiActedUponBy e)^{-1}(A)$, and $\boundaryOf_{\setOf{G_0, e}} A = A \smallsetminus (\blank \isSemiActedUponBy e)^{-1}(A) \cup (\blank \isSemiActedUponBy e)^{-1}(A) \smallsetminus A$.
      \item \label{item:properties-of-interior-closure-and-boundary:complement} 
            $(M \smallsetminus A)^{-E} = M \smallsetminus A^{+E}$ and $(M \smallsetminus A)^{+E} = M \smallsetminus A^{-E}$.
      \item \label{item:properties-of-interior-closure-and-boundary:union}
            $\bigcup_{i \in I} A_i^{-E} \subseteq (\bigcup_{i \in I} A_i)^{-E}$ and $\bigcup_{i \in I} A_i^{+E} = (\bigcup_{i \in I} A_i)^{+E}$.
      \item \label{item:properties-of-interior-closure-and-boundary:intersection}
            $(\bigcap_{i \in I} A_i)^{-E} = \bigcap_{i \in I} A_i^{-E}$ and $(\bigcap_{i \in I} A_i)^{+E} \subseteq \bigcap_{i \in I} A_i^{+E}$.
      \item \label{item:properties-of-interior-closure-and-boundary:inclusions}
            Let $E \subseteq E'$. Then, $A^{-E} \supseteq A^{-E'}$, $A^{+E} \subseteq A^{+E'}$, and $\boundaryOf_E A \subseteq \boundaryOf_{E'} A$.
      \item \label{item:properties-of-interior-closure-and-boundary:subset}
            Let $A \subseteq A'$. Then, $A^{-E} \subseteq (A')^{-E}$ and $A^{+E} \subseteq (A')^{+E}$.
      \item \label{item:properties-of-interior-closure-and-boundary:neutral-element}
            Let $G_0 \in E$. Then, $A^{-E} \subseteq A \subseteq A^{+E}$.
      \item \label{item:properties-of-interior-closure-and-boundary:finite} 
            Let $G_0$, $A$, and $E$ be finite. Then, $A^{-E}$, $A^{+E}$, and $\boundaryOf_E A$ are finite. More precisely, $\cardinalityOf{A^{-E}} \leq \cardinalityOf{G_0} \cdot \cardinalityOf{A}$ and $\cardinalityOf{A^{+E}} \leq \cardinalityOf{G_0} \cdot \cardinalityOf{A} \cdot \cardinalityOf{E}$.
      \item \label{item:properties-of-interior-closure-and-boundary:commute}
            Let $g \in G$ and let $G_0 \cdot E \subseteq E$. Then, $g \actsOnPoint A^{-E} = (g \actsOnPoint A)^{-E}$, $g \actsOnPoint A^{+E} = (g \actsOnPoint A)^{+E}$, and $g \actsOnPoint \boundaryOf_E A = \boundaryOf_E (g \actsOnPoint A)$.
      \item \label{item:properties-of-interior-closure-and-boundary:commute-with-liberation}
            Let $m \in M$, let $G_0 \cdot E \subseteq E$, and let $\iota \from M \to G \modulo G_0$, $m \mapsto G_{m_0, m}$. Then, $m \isSemiActedUponBy \iota(A^{-E}) = (m \isSemiActedUponBy \iota(A))^{-E}$, $m \isSemiActedUponBy \iota(A^{+E}) = (m \isSemiActedUponBy \iota(A))^{+E}$, and $m \isSemiActedUponBy \iota(\boundaryOf_E A) = \boundaryOf_E (m \isSemiActedUponBy \iota(A))$.
      \item \label{item:properties-of-interior-closure-and-boundary:liberation}
            Let $G_0 \cdot E \subseteq E$ and let $E^{-1} \subseteq E$. Then, $A \isSemiActedUponBy E \subseteq A^{+E}$.
      \item \label{item:properties-of-interior-closure-and-boundary:same-repeated}
            Let $G_0 \cdot E \subseteq E$. Then, $(A^{-E})^{-E'} = A^{- E' \cdot E}$ and $(A^{+E})^{+E'} = A^{+ E' \cdot E}$.
      \item \label{item:properties-of-interior-closure-and-boundary:different-repeated}
            Let $G_0 \cdot E \subseteq E$. Then, $(A^{+E})^{-E'} = \bigcap_{e' \in E'} A^{+ e' \cdot E}$ and $(A^{-E})^{+E'} = \bigcup_{e' \in E'} A^{- e' \cdot E}$. And, if $(E')^{-1} \subseteq E$, then $A \subseteq (A^{+E})^{-E'}$ and $(A^{-E})^{+E'} \subseteq A$. \qedhere 
    \end{aenumerate}
  \end{lemma}

  \begin{proof}
    \begin{aenumerate}
      \item Because $\blank \isSemiActedUponBy G_0 = \identityMap_M$, this is a direct consequence of \cref{definition:interior-closure-boundary}.
      \item Because $(\blank \isSemiActedUponBy G_0)^{-1}(A) = A$, this is a direct consequence of \cref{definition:interior-closure-boundary}.
      \item For each $m \in M$,
            \begin{align*}
              m \in (M \smallsetminus A)^{-E}
              &\ifAndOnlyIf m \isSemiActedUponBy E \subseteq M \smallsetminus A\\
              &\ifAndOnlyIf (m \isSemiActedUponBy E) \cap A = \emptyset\\
              &\ifAndOnlyIf m \in M \smallsetminus A^{+E}.
            \end{align*}
            Hence, $(M \smallsetminus A)^{-E} = M \smallsetminus A^{+E}$. Therefore,
            \begin{align*}
              (M \smallsetminus A)^{+E}
              &= M \smallsetminus (M \smallsetminus (M \smallsetminus A)^{+E})\\
              &= M \smallsetminus (M \smallsetminus (M \smallsetminus A))^{-E}\\
              &= M \smallsetminus A^{-E}.
            \end{align*}
      \item For each $m \in M$,
            \begin{align*}
              m \in \bigcup_{i \in I} A_i^{-E} &\ifAndOnlyIf \Exists i \in I \SuchThat m \in A_i^{-E}\\
                                               &\ifAndOnlyIf \Exists i \in I \SuchThat m \isSemiActedUponBy E \subseteq A_i\\
                                               &\implies m \isSemiActedUponBy E \subseteq \bigcup_{i \in I} A_i\\ 
                                               &\ifAndOnlyIf m \in (\bigcup_{i \in I} A_i)^{-E}.
            \end{align*}
            Hence, $\bigcup_{i \in I} A_i^{-E} \subseteq (\bigcup_{i \in I} A_i)^{-E}$.

            Moreover, for each $m \in M$,
            \begin{align*}
              m \in \bigcup_{i \in I} A_i^{+E} &\ifAndOnlyIf \Exists i \in I \SuchThat m \in A_i^{+E}\\
                                               &\ifAndOnlyIf \Exists i \in I \SuchThat (m \isSemiActedUponBy E) \cap A_i \neq \emptyset\\
                                               &\ifAndOnlyIf (m \isSemiActedUponBy E) \cap (\bigcup_{i \in I} A_i) \neq \emptyset\\
                                               &\ifAndOnlyIf m \in (\bigcup_{i \in I} A_i)^{+E}.
            \end{align*}
            Therefore, $\bigcup_{i \in I} A_i^{+E} = (\bigcup_{i \in I} A_i)^{+E}$.
  %
      \item According to \cref{item:properties-of-interior-closure-and-boundary:complement} and \cref{item:properties-of-interior-closure-and-boundary:union},
            \begin{align*}
              (\bigcap_{i \in I} A_i)^{-E} &= M \smallsetminus (M \smallsetminus (\bigcap_{i \in I} A_i)^{-E})\\
                                           &= M \smallsetminus (M \smallsetminus (\bigcap_{i \in I} A_i))^{+E}\\
                                           &= M \smallsetminus (\bigcup_{i \in I} M \smallsetminus A_i)^{+E}\\
                                           &= M \smallsetminus (\bigcup_{i \in I} (M \smallsetminus A_i)^{+E})\\
                                           &= M \smallsetminus (\bigcup_{i \in I} M \smallsetminus A_i^{-E})\\
                                           &= M \smallsetminus (M \smallsetminus (\bigcap_{i \in I} A_i^{-E}))\\
                                           &= \bigcap_{i \in I} A_i^{-E}
            \end{align*}
            and
            \begin{align*}
              (\bigcap_{i \in I} A_i)^{+E} &=         M \smallsetminus (M \smallsetminus (\bigcap_{i \in I} A_i)^{+E})\\
                                           &=         M \smallsetminus (\bigcup_{i \in I} M \smallsetminus A_i)^{-E}\\
                                           &\subseteq M \smallsetminus (\bigcup_{i \in I} M \smallsetminus A_i^{+E})\\
                                           &=         M \smallsetminus (M \smallsetminus (\bigcap_{i \in I} A_i^{+E}))\\
                                           &=         \bigcap_{i \in I} A_i^{+E}.
            \end{align*}
      \item This is a direct consequence of \cref{definition:interior-closure-boundary}. 
      \item This is a direct consequence of \cref{definition:interior-closure-boundary}.
      \item This is a direct consequence of \cref{definition:interior-closure-boundary}.
      \item For each $e \in E$, according to \cref{lemma:liberation-preimage},
            \begin{equation*}
              \cardinalityOf{\bigcup_{a \in A} (\blank \isSemiActedUponBy e)^{-1}(a)}
              = \cardinalityOf{(\blank \isSemiActedUponBy e)^{-1}(A)}
              \leq \cardinalityOf{G_0} \cdot \cardinalityOf{A}.
            \end{equation*}
            Hence, because $A^{-E} = \bigcap_{e \in E} \bigcup_{a \in A} (\blank \isSemiActedUponBy e)^{-1}(a)$, we have $\cardinalityOf{A^{-E}} \leq \cardinalityOf{G_0} \cdot \cardinalityOf{A} < \infty$. And, because $A^{+E} = \bigcup_{e \in E} \bigcup_{a \in A} (\blank \isSemiActedUponBy e)^{-1}(a)$, we have $\cardinalityOf{A^{+E}} \leq \cardinalityOf{E} \cdot \cardinalityOf{G_0} \cdot \cardinalityOf{A} < \infty$.
            And, because $\boundaryOf_E A \subseteq A^{+E}$, we also have $\cardinalityOf{\boundaryOf_E A} < \infty$.
      \item Let $m \in M$. Because $\isSemiActedUponBy$ semi-commutes with $\actsOnPoint$, there is a $g_0 \in G_0$ such that $(g^{-1} \actsOnPoint m) \isSemiActedUponBy E = g^{-1} \actsOnPoint (m \isSemiActedUponBy g_0 \cdot E)$. And, because $G_0 \cdot E \subseteq E$, we have $g_0 \cdot E \subseteq E$ and $g_0^{-1} \cdot E \subseteq E$; hence $E = g_0 g_0^{-1} \cdot E = g_0 \cdot (g_0^{-1} \cdot E) \subseteq g_0 \cdot E$; thus $g_0 \cdot E = E$. Therefore, $(g^{-1} \actsOnPoint m) \isSemiActedUponBy E = g^{-1} \actsOnPoint (m \isSemiActedUponBy E)$. 

            Thus, for each $m \in M$,
            \begin{align*}
              m \in g \actsOnPoint A^{-E} &\ifAndOnlyIf \Exists m' \in A^{-E} \SuchThat g \actsOnPoint m' = m\\
                                         &\ifAndOnlyIf g^{-1} \actsOnPoint m \in A^{-E}\\
                                         &\ifAndOnlyIf (g^{-1} \actsOnPoint m) \isSemiActedUponBy E \subseteq A\\
                                         &\ifAndOnlyIf g^{-1} \actsOnPoint (m \isSemiActedUponBy E) \subseteq A\\
                                         &\ifAndOnlyIf m \isSemiActedUponBy E \subseteq g \actsOnPoint A\\
                                         &\ifAndOnlyIf m \in (g \actsOnPoint A)^{-E}.
            \end{align*}
            In conclusion, $g \actsOnPoint A^{-E} = (g \actsOnPoint A)^{-E}$. Moreover, for each $m \in M$,
            \begin{align*}
              m \in g \actsOnPoint A^{+E} &\ifAndOnlyIf g^{-1} \actsOnPoint m \in A^{+E}\\
                                         &\ifAndOnlyIf ((g^{-1} \actsOnPoint m) \isSemiActedUponBy E) \cap A \neq \emptyset\\
                                         &\ifAndOnlyIf (g^{-1} \actsOnPoint (m \isSemiActedUponBy E)) \cap A \neq \emptyset\\
                                         &\ifAndOnlyIf (m \isSemiActedUponBy E) \cap (g \actsOnPoint A) \neq \emptyset\\
                                         &\ifAndOnlyIf m \in (g \actsOnPoint A)^{+E}.
            \end{align*}
            In conclusion, $g \actsOnPoint A^{+E} = (g \actsOnPoint A)^{+E}$. Ultimately,
            \begin{align*}
              g \actsOnPoint \boundaryOf_E A &= g \actsOnPoint (A^{+E} \smallsetminus A^{-E})\\
                                          &= (g \actsOnPoint A^{+E}) \smallsetminus (g \actsOnPoint A^{-E})\\
                                          &= (g \actsOnPoint A)^{+E} \smallsetminus (g \actsOnPoint A)^{-E}\\
                                          &= \boundaryOf_E (g \actsOnPoint A).
            \end{align*}
      \item According to
            \cref{item:properties-of-interior-closure-and-boundary:commute},
            \begin{align*}
              m \isSemiActedUponBy \iota(A^{-E})
              &= g_{m_0, m} \actsOnPoint A^{-E}\\
              &= (g_{m_0, m} \actsOnPoint A)^{-E}\\
              &= (m \isSemiActedUponBy \iota(A))^{-E},
            \end{align*}
            and
            \begin{align*}
              m \isSemiActedUponBy \iota(A^{+E})
              &= g_{m_0, m} \actsOnPoint A^{+E}\\
              &= (g_{m_0, m} \actsOnPoint A)^{+E}\\
              &= (m \isSemiActedUponBy \iota(A))^{+E},
            \end{align*}
            and
            \begin{align*}
              m \isSemiActedUponBy \iota(\boundaryOf_E A)
              &= g_{m_0, m} \actsOnPoint \boundaryOf_E A\\
              &= \boundaryOf_E (g_{m_0, m} \actsOnPoint A)\\
              &= \boundaryOf_E (m \isSemiActedUponBy \iota(A)).
            \end{align*}
      \item Let $m \in A$ and let $e \in E$. Then, there is a $g \in G$ such that $g G_0 = e$. Hence, because $\isSemiActedUponBy$ is a semi-action with defect $G_0$, there is a $g_0 \in G_0$ such that, for each $\mathfrak{g} \in G \modulo G_0$, we have $(m \isSemiActedUponBy e) \isSemiActedUponBy g_0 \cdot \mathfrak{g} = m \isSemiActedUponBy g \cdot \mathfrak{g}$.

            Put $e' = g_0 \cdot g^{-1} G_0$. Then, because $G_0 \cdot E \subseteq E$ and $E^{-1} \subseteq E$, we have $e' \in E$. Moreover, $(m \isSemiActedUponBy e) \isSemiActedUponBy e' = m \isSemiActedUponBy g g^{-1} G_0 = m$. Therefore, $m \in (m \isSemiActedUponBy e) \isSemiActedUponBy E$ and hence $((m \isSemiActedUponBy e) \isSemiActedUponBy E) \cap A \neq \emptyset$. Thus, $m \isSemiActedUponBy e \in A^{+E}$. In conclusion, $A \isSemiActedUponBy E \subseteq A^{+E}$.
      \item For each $m \in M$, according to \cref{lemma:double-liberation-with-sets-to-one-liberation}, we have $(m \isSemiActedUponBy E') \isSemiActedUponBy E = m \isSemiActedUponBy E' \cdot E$. Therefore, for each $m \in M$,
            \begin{align*}
              m \in (A^{-E})^{-E'}
              &\ifAndOnlyIf m \isSemiActedUponBy E' \subseteq A^{-E}\\
              &\ifAndOnlyIf (m \isSemiActedUponBy E') \isSemiActedUponBy E \subseteq A\\
              &\ifAndOnlyIf m \isSemiActedUponBy E' \cdot E \subseteq A\\
              &\ifAndOnlyIf m \in A^{- E' \cdot E}.
            \end{align*}
            In conclusion, $(A^{-E})^{-E'} = A^{- E' \cdot E}$. Moreover, for each $m \in M$,
            \begin{align*}
              m \in (A^{+E})^{+E'}
              &\ifAndOnlyIf (m \isSemiActedUponBy E') \cap A^{+E} \neq \emptyset\\
              &\ifAndOnlyIf \Exists e' \in E' \SuchThat m \isSemiActedUponBy e' \in A^{+E}\\
              &\ifAndOnlyIf \Exists e' \in E' \SuchThat ((m \isSemiActedUponBy e') \isSemiActedUponBy E) \cap A \neq \emptyset\\
              &\ifAndOnlyIf ((m \isSemiActedUponBy E') \isSemiActedUponBy E) \cap A \neq \emptyset\\
              &\ifAndOnlyIf (m \isSemiActedUponBy E' \cdot E) \cap A \neq \emptyset\\
              &\ifAndOnlyIf m \in A^{+ E' \cdot E}.
            \end{align*}
            In conclusion, $(A^{+E})^{+E'} = A^{+ E' \cdot E}$.
      \item For each $m \in M$ and each $e' \in E'$, according to \cref{lemma:double-liberation-with-sets-to-one-liberation}, we have $(m \isSemiActedUponBy e') \isSemiActedUponBy E = m \isSemiActedUponBy e' \cdot E$. Therefore, for each $m \in M$,
            \begin{align*}
              m \in (A^{+E})^{-E'}
              &\ifAndOnlyIf m \isSemiActedUponBy E' \subseteq A^{+E}\\
              &\ifAndOnlyIf \ForEach e' \in E' \Holds m \isSemiActedUponBy e' \in A^{+E}\\
              &\ifAndOnlyIf \ForEach e' \in E' \Holds ((m \isSemiActedUponBy e') \isSemiActedUponBy E) \cap A \neq \emptyset\\
              &\ifAndOnlyIf \ForEach e' \in E' \Holds (m \isSemiActedUponBy e' \cdot E) \cap A \neq \emptyset\\
              &\ifAndOnlyIf \ForEach e' \in E' \Holds m \in A^{+ e' \cdot E}\\
              &\ifAndOnlyIf m \in \bigcap_{e' \in E'} A^{+ e' \cdot E}.
            \end{align*}
            In conclusion, $(A^{+E})^{-E'} = \bigcap_{e' \in E'} A^{+ e' \cdot E}$. Moreover, for each $m \in M$,
            \begin{align*}
              m \in (A^{-E})^{+E'}
              &\ifAndOnlyIf (m \isSemiActedUponBy E') \cap A^{-E} \neq \emptyset\\
              &\ifAndOnlyIf \Exists e' \in E' \SuchThat m \isSemiActedUponBy e' \in A^{-E}\\
              &\ifAndOnlyIf \Exists e' \in E' \SuchThat (m \isSemiActedUponBy e') \isSemiActedUponBy E \subseteq A\\
              &\ifAndOnlyIf \Exists e' \in E' \SuchThat m \isSemiActedUponBy e' \cdot E \subseteq A\\
              &\ifAndOnlyIf \Exists e' \in E' \SuchThat m \in A^{- e' \cdot E}\\
              &\ifAndOnlyIf m \in \bigcup_{e' \in E'} A^{- e' \cdot E}.
            \end{align*}
            In conclusion, $(A^{-E})^{+E'} = \bigcup_{e' \in E'} A^{- e' \cdot E}$.

            From now on, let $(E')^{-1} \subseteq E$. Then, for each $e' \in E'$, we have $G_0 \in e' \cdot E$ and hence, according to \cref{item:properties-of-interior-closure-and-boundary:neutral-element}, we have $A \subseteq A^{+ e' \cdot E}$ and $A^{- e' \cdot E} \subseteq A$. In conclusion, $A \subseteq (A^{+E})^{-E'}$ and $(A^{-E})^{+E'} \subseteq A$. \qedhere 
    \end{aenumerate}
  \end{proof} 

  The restriction $\Delta_{X, A}^-$ of $\Delta$ given in \cref{lemma:Delta-X-A-minus-plus-are-surjective} is well-defined according to the next lemma, which itself holds due to the locality of $\Delta$.

  \begin{lemma} 
  \label{lemma:global-transition-function-and-interior-closure}
    Let $\mathcal{C} = \ntuple{\mathcal{R}, Q, N, \delta}$ be a semi-cellular automaton, let $\Delta$ be the global transition function of $\mathcal{C}$, let $c$ and $c'$ be two global configurations of $\mathcal{C}$, and let $A$ be a subset of $M$.
    \begin{aenumerate}
      \item \label{item:global-transition-function-and-interior-closure:interior}
            If $c\restrictedTo_A = c'\restrictedTo_A$, then $\Delta(c)\restrictedTo_{A^{-N}} = \Delta(c')\restrictedTo_{A^{-N}}$.
      \item \label{item:global-transition-function-and-interior-closure:closure}
            If $c\restrictedTo_{M \smallsetminus A} = c'\restrictedTo_{M \smallsetminus A}$, then $\Delta(c)\restrictedTo_{M \smallsetminus A^{+N}} = \Delta(c')\restrictedTo_{M \smallsetminus A^{+N}}$.
      \item If $N^{-1} \subseteq N$ and $c\restrictedTo_{A^{+N}} = c'\restrictedTo_{A^{+N}}$, then $\Delta(c)\restrictedTo_A = \Delta(c')\restrictedTo_A$. \qedhere 
    \end{aenumerate}
  \end{lemma}

  \begin{proof}
    \begin{aenumerate}
      \item Let $c\restrictedTo_A = c'\restrictedTo_A$. Then, for each $m \in A^{-N}$, we have $m \isSemiActedUponBy N \subseteq A$ and hence $\Delta(c)(m) = \Delta(c')(m)$. In conclusion, $\Delta(c)\restrictedTo_{A^{-N}} = \Delta(c')\restrictedTo_{A^{-N}}$.
      \item This is a direct consequence of \cref{item:global-transition-function-and-interior-closure:interior} and \cref{item:properties-of-interior-closure-and-boundary:complement} of \cref{lemma:properties-of-interior-closure-and-boundary}.
      \item Let $N^{-1} \subseteq N$ and let $c\restrictedTo_{A^{+N}} = c'\restrictedTo_{A^{+N}}$. Then, for each $m \in A$, according to \cref{item:properties-of-interior-closure-and-boundary:liberation} of \cref{lemma:properties-of-interior-closure-and-boundary}, we have $m \isSemiActedUponBy N \subseteq A^{+N}$ and hence $\Delta(c)(m) = \Delta(c')(m)$. In conclusion, $\Delta(c)\restrictedTo_A = \Delta(c')\restrictedTo_A$. \qedhere 
    \end{aenumerate}
  \end{proof}

  \begin{lemma} 
  \label{lemma:Delta-X-A-minus-plus-are-surjective}
    Let $\mathcal{C} = \ntuple{\mathcal{R}, Q, N, \delta}$ be a semi-cellular automaton, let $\Delta$ be the global transition function of $\mathcal{C}$, let $X$ be a subset of $Q^M$, and let $A$ be a subset of $M$. The map
            \begin{align*}
              \Delta_{X, A}^- \from \pi_A(X) &\to \pi_{A^{-N}}(\Delta(X)), \mathnote{map $\Delta_{X, A}^-$}\index[symbols]{DeltaXAminus@$\Delta_{X, A}^-$}\\
              p &\mapsto \Delta(c)\restrictedTo_{A^{-N}}, \text{ where } c \in X \text{ such that } c\restrictedTo_A = p,
            \end{align*}
            is surjective. The map $\Delta_{Q^M, A}^-$ is denoted by $\Delta_A^-$\graffito{map $\Delta_A^-$}\index[symbols]{DeltaAminus@$\Delta_A^-$}.
  %
  \end{lemma}

  \begin{proof}
            Let $p' \in \pi_{A^{-N}}(\Delta(X))$. Then, there is a $c' \in \Delta(X)$ such that $c'\restrictedTo_{A^{-N}} = p'$. Moreover, there is a $c \in X$ such that $\Delta(c) = c'$. Put $p = c\restrictedTo_A \in \pi_A(X)$. Then, $\Delta_{X, A}^-(p) = \Delta(c)\restrictedTo_{A^{-N}} = c'\restrictedTo_{A^{-N}} = p'$. Hence, $\Delta_{X, A}^-$ is surjective.
  \end{proof}


  The restrictions $\Delta_{X, A}^-$ of global transition functions $\Delta$ of big-cellular automata are $\actsOnMap_H$-e\-qui\-var\-i\-ant in the sense given in

  \begin{lemma}
  \label{lemma:equivariance-of-restriction}
    Let $H$ be a $\mathcal{K}$-big subgroup of $G$, let $\mathcal{C} = \ntuple{\mathcal{R}, Q, N, \delta}$ be a semi-cellular automaton such that $\delta$ is $\bullet_{H_0}$-invariant, let $\Delta$ be the global transition function of $\mathcal{C}$, let $X$ be a subset of $Q^M$, and let $A$ be a subset of $M$. Then,
    \begin{equation*} 
      \ForEach h \in H \ForEach p \in \pi_A(X) \Holds \Delta_{h \actsOnMap X, h \actsOnPoint A}^-(h \actsOnMap p) = h \actsOnMap \Delta_{X, A}^-(p). \qedhere
    \end{equation*}
  \end{lemma}

  \begin{proof}
    Let $h \in H$ and let $p \in \pi_A(X)$. Then, the domain of $h \actsOnMap p$ is $h \actsOnPoint A$ and we have $h \actsOnMap p \in h \actsOnMap \pi_A(X) = \pi_{h \actsOnPoint A}(h \actsOnMap X)$. Hence, the term $\Delta_{h \actsOnMap X, h \actsOnPoint A}^-(h \actsOnMap p)$ is well-defined. Moreover, by definition of $p$, there is a $c \in X$ such that $c\restrictedTo_A = p$, and hence $(h \actsOnMap c)\restrictedTo_{h \actsOnPoint A} = h \actsOnMap p$. Therefore, by definition of $\Delta_{h \actsOnMap X, h \actsOnPoint A}^-$ and of $\Delta_{X, A}^-$, because $\Delta$ is $\actsOnMap_H$-e\-qui\-var\-i\-ant, we have $\Delta_{h \actsOnMap X, h \actsOnPoint A}^-(h \actsOnMap p) = \Delta(h \actsOnMap c)\restrictedTo_{h \actsOnPoint A^{-N}} = h \actsOnMap \Delta(c)\restrictedTo_{A^{-N}} = h \actsOnMap \Delta_{X, A}^-(p)$.
  \end{proof}

  \begin{corollary}
  \label{corollary:equivariance-of-restriction}
    In the situation of \cref{lemma:equivariance-of-restriction}, let $M$ be identified with $G \modulo G_0$ by $\iota \givenBy m \mapsto G_{m_0, m}$. Then,
    \begin{equation*} 
      \ForEach m \in M \ForEach p \in \pi_A(X) \Holds \Delta_{m \actsByItsCoordinateOn X, m \isSemiActedUponBy A}^-(m \actsByItsCoordinateOn p) = m \actsByItsCoordinateOn \Delta_{X, A}^-(p). \qedhere
    \end{equation*}
  \end{corollary}

  \begin{proof}
    This is a direct consequence of \cref{lemma:equivariance-of-restriction}, because $g_{m_0, m} \in H$, $m \actsByItsCoordinateOn \blank = g_{m_0, m} \actsOnMap \blank$, and $m \isSemiActedUponBy \blank = g_{m_0, m} \actsOnPoint \blank$.
  \end{proof}

  A net of non-empty and finite subsets of $M$ is a right Følner net if and only if these subsets are asymptotically invariant under the right semi-action induced by $\actsOnPoint$, which broadly speaking means that these subsets are asymptotically invariant under small perturbations. If the stabiliser of the origin under $\actsOnPoint$ is finite, that is the case if and only if those subsets are asymptotically invariant under taking finite boundaries, in other words, if and only if the finite boundaries of the subsets grow much slower than the subsets themselves. This is shown in

  \begin{theorem} 
  \label{theorem:boundary-characterisation-of-Folner-net} 
    Let $G_0$ be finite and let $\net{F_i}_{i \in I}$ be a net in $\setOf{F \subseteq M \suchThat F \neq \emptyset, F \text{ finite}}$ indexed by $(I, \leq)$. The net $\net{F_i}_{i \in I}$ is a right Følner net in $\mathcal{R}$ if and only if 
    \begin{equation*}
      \ForEach E \subseteq G \modulo G_0 \text{ finite} \Holds \lim_{i \in I} \frac{\cardinalityOf{\boundaryOf_E F_i}}{\cardinalityOf{F_i}} = 0. \qedhere
    \end{equation*}
  \end{theorem}

  \begin{proof}
    First, let $\net{F_i}_{i \in I}$ be a right Følner net in $\mathcal{R}$. Furthermore, let $E \subseteq G \modulo G_0$ be finite. Moreover, let $i \in I$. For each $e \in E$ and each $e' \in E$, put $A_{i,e,e'} = (\blank \isSemiActedUponBy e)^{-1}(F_i) \smallsetminus (\blank \isSemiActedUponBy e')^{-1}(F_i)$. For each $\mathfrak{g} \in G \modulo G_0$, put $B_{i,\mathfrak{g}} = F_i \smallsetminus (\blank \isSemiActedUponBy \mathfrak{g})^{-1}(F_i)$. 
    According to \cref{definition:interior-closure-boundary}, 
    \begin{align*}
      \boundaryOf_E F_i 
      &= \parens[\Big]{\bigcup_{e \in E} (\blank \isSemiActedUponBy e)^{-1}(F_i)} \smallsetminus \parens[\Big]{\bigcap_{e' \in E} (\blank \isSemiActedUponBy e')^{-1}(F_i)}\\
      &= \bigcup_{e, e' \in E} (\blank \isSemiActedUponBy e)^{-1}(F_i) \smallsetminus (\blank \isSemiActedUponBy e')^{-1}(F_i)\\
      &= \bigcup_{e, e' \in E} A_{i,e,e'}.
    \end{align*}
    Hence, $\cardinalityOf{\boundaryOf_E F_i} \leq \sum_{e, e' \in E} \cardinalityOf{A_{i,e,e'}}$.

    According to \cref{lemma:cardinality-of-inverse-image-of-liberation-minus-the-same-less-than-or-equal-to-whatever}, we have $\cardinalityOf{A_{i,e,e'}} \leq \cardinalityOf{G_0}^2 \cdot \max_{g \in e} \cardinalityOf{B_{i, g^{-1} \cdot e'}}$. Put $E' = \setOf{g^{-1} \cdot e' \suchThat e, e' \in E, g \in e}$. Because $E$ is finite, $G_0$ is finite, and, for each $e \in E$, we have $\cardinalityOf{e} = \cardinalityOf{G_0}$, the set $E'$ is finite. Therefore, 
    \begin{align*} 
      \frac{\cardinalityOf{\boundaryOf_E F_i}}{\cardinalityOf{F_i}}
      &\leq \frac{1}{\cardinalityOf{F_i}} \sum_{e, e' \in E} \cardinalityOf{A_{i,e,e'}}\\
      &\leq \frac{\cardinalityOf{G_0}^2}{\cardinalityOf{F_i}} \sum_{e, e' \in E} \max_{g \in e} \cardinalityOf{B_{i, g^{-1} \cdot e'}}\\
      &\leq \frac{\cardinalityOf{G_0}^2 \cdot \cardinalityOf{E}^2}{\cardinalityOf{F_i}} \max_{e' \in E'} \cardinalityOf{B_{i, e'}}\\
      &\leq \cardinalityOf{G_0}^2 \cdot \cardinalityOf{E}^2 \cdot \max_{e' \in E'} \frac{\cardinalityOf{F_i \smallsetminus (\blank \isSemiActedUponBy e')^{-1}(F_i)}}{\cardinalityOf{F_i}}\\ 
      &\underset{i \in I}{\to} 0.
    \end{align*}
    In conclusion, $\lim_{i \in I} \frac{\cardinalityOf{\boundaryOf_E F_i}}{\cardinalityOf{F_i}} = 0$.

    Secondly, for each finite $E \subseteq G \modulo G_0$, let $\lim_{i \in I} \frac{\cardinalityOf{\boundaryOf_E F_i}}{\cardinalityOf{F_i}} = 0$. Furthermore, let $i \in I$, let $e \in G \modulo G_0$, and put $E = \setOf{G_0, e}$. According to \cref{item:properties-of-interior-closure-and-boundary:only-neutral-element-and-another} of \cref{lemma:properties-of-interior-closure-and-boundary}, we have $F_i \smallsetminus (\blank \isSemiActedUponBy e)^{-1}(F_i) \subseteq \boundaryOf_E F_i$. Therefore,
    \begin{equation*}
      \frac{\cardinalityOf{F_i \smallsetminus (\blank \isSemiActedUponBy e)^{-1}(F_i)}}{\cardinalityOf{F_i}}
      \leq \frac{\cardinalityOf{\boundaryOf_E F_i}}{\cardinalityOf{F_i}}
      \underset{i \in I}{\to} 0.
    \end{equation*}
    In conclusion, $\net{F_i}_{i \in I}$ is a right Følner net in $\mathcal{R}$.
  \end{proof}

  \begin{example}[Lattice] 
  \label{example:lattice:Folner}
    In the situation of \cref{example:lattice:interior-closure-and-boundary}, the sequence $\sequence{\sphere(\rho)}_{\rho \in \N_+}$ grows linearly in size, more precisely, it grows like $\sequence{4 \rho}_{\rho \in \N_+}$; the sequence $\sequence{\ball(\rho)}_{\rho \in \N_0}$ grows polynomially in size, more precisely, it grows like $\sequence{\sum_{\varrho = 0}^\rho \cardinalityOf{\sphere(\varrho)}}_{\rho \in \N_0} = \sequence{2\rho(\rho + 1) + 1}_{\rho \in \N_0}$; the sequence $\sequence{\boundaryOf_{\ball(1)} \ball(\rho)}_{\rho \in \N_+}$ grows linearly in size, more precisely, it grows like $\sequence{\cardinalityOf{\ball(\rho + 1) \smallsetminus \ball(\rho - 1)}}_{\rho \in \N_+} = \sequence{4 (2\rho + 1)}_{\rho \in \N_+}$; in general, for each non-negative integer $\varrho$, we have $\sequence{\cardinalityOf{\boundaryOf_{\ball(\varrho)} \ball(\rho)}}_{\rho \in \N_+} = \sequence{4\varrho (2\rho + 1)}_{\rho \in \N_{\geq \varrho}}$. It follows that
    \begin{equation*}
      \ForEach \varrho \in \N_0 \Holds \lim_{\rho \to \infty} \frac{\cardinalityOf{\boundaryOf_{\ball(\varrho)} \ball(\rho)}}{\cardinalityOf{\ball(\rho)}} = 0.
    \end{equation*}
    Hence, according to \cref{theorem:boundary-characterisation-of-Folner-net}, the sequence $\sequence{\ball(\rho)}_{\rho \in \N_0}$ is a right Følner net in $\mathcal{R}$, which was also shown in \cref{example:lattice}. 
  \end{example}

  \begin{counterexample}[Tree] 
  \label{example:tree:Folner}
    In the situation of \cref{example:tree:interior-closure-and-boundary}, the sequence $\sequence{\sphere(\rho)}_{\rho \in \N_+}$ grows exponentially in size, more precisely, it grows like $\sequence{3^\rho + 3^{\rho - 1}}_{\rho \in \N_+}$; the sequence $\sequence{\ball(\rho)}_{\rho \in \N_0}$ also grows exponentially in size, more precisely, it grows like $\sequence{\sum_{\varrho = 0}^\rho \cardinalityOf{\sphere(\varrho)}}_{\rho \in \N_0} = \sequence{2 \cdot 3^\rho - 1}_{\rho \in \N_0}$; for each non-negative integer $\varrho$, the sequence $\sequence{\boundaryOf_{\ball(\varrho)} \ball(\rho)}_{\rho \in \N_{\geq \varrho}}$ also grows exponentially in size, more precisely, it grows like $\sequence{\cardinalityOf{\ball(\rho + \varrho) \smallsetminus \ball(\rho - \varrho)}}_{\rho \in \N_{\geq \varrho}} = \sequence{2(3^{\rho + \varrho} - 3^{\rho - \varrho})}_{\rho \in \N_{\geq \varrho}}$. It follows that
    \begin{equation*}
      \ForEach \varrho \in \N_0 \Holds \lim_{\rho \to \infty} \frac{\cardinalityOf{\boundaryOf_{\ball(\varrho)} \ball(\rho)}}{\cardinalityOf{\ball(\rho)}} = 3^\varrho - 3^{-\varrho} \neq 0.
    \end{equation*}
    Hence, according to \cref{theorem:boundary-characterisation-of-Folner-net}, each subsequence of $\sequence{\ball(\rho)}_{\rho \in \N_0}$ is not a right Følner net in $\mathcal{R}$, which was also shown in \cref{example:tree}. Actually, because we deduced in \cref{example:tree:Tarski-Folner} that the cell space $\mathcal{R}$ is not right amenable, there is no right Følner net in $\mathcal{R}$. 
  \end{counterexample}

  In \cref{definition:entropy} we use a net of non-empty and finite subsets of $M$ to define the entropy of a subset of global configurations. If that net is a right Følner net, then the global transition function of a big-cellular automaton is surjective if and only if the entropy of its image is maximal and the entropy of its image is maximal if and only if it is pre-injective. For the first equivalence it suffices that the net has a weaker property, which we define below; for the second equivalence though that weaker property does not suffice.

  \begin{definition} 
  \label{definition:Erling-net}
    Let $\net{F_i}_{i \in I}$ be a net in $\setOf{F \subseteq M \suchThat F \neq \emptyset, F \text{ finite}}$ indexed by $(I, \leq)$. It is called \define{right Erling net in $\mathcal{R}$ indexed by $(I, \leq)$}\graffito{right Erling net in $\mathcal{R}$ indexed by $(I, \leq)$}\index{Erling net@right Erling net}\index{net!Erling} if and only if
    \begin{equation*}
      \ForEach E \subseteq G \modulo G_0 \text{ finite} \Holds \limsup_{i \in I} \frac{\cardinalityOf{\boundaryOf_E^- F_i}}{\cardinalityOf{F_i}} < 1. \qedhere
    \end{equation*}
  \end{definition}

  \begin{remark}
    Regardless of whether $\net{F_i}_{i \in I}$ is an Erling net or not, the above limit superior is always less than or equal to $1$.
  \end{remark}

  \begin{definition}
  \label{definition:tractable} 
    The cell space $\mathcal{R}$ is called \define{right tractable}\graffito{right tractable}\index{tractable right@right tractable} if and only if there is a right Erling net in $\mathcal{R}$.
  \end{definition} 

  \begin{remark}
  \label{remark:amenable-implies-tractable}
    Because each right Følner net is a right Erling net, each right-a\-me\-na\-ble cell space with finite stabilisers is right tractable.
  \end{remark}

  The notions of being finitely and symmetrically right generated, of distance, and of balls will be introduced in \cref{chapter:growth}. A quick preview was given in the paragraph \textsc{preliminary notions} at the beginning of this chapter. A finitely and symmetrically right-gen\-er\-at\-ed cell space is right tractable, which is shown in

  \begin{lemma}
  \label{lemma:finitely-right-generated-is-tractable}
    Let $\mathcal{R}$ be finitely and symmetrically right generated. It is right tractable. And, for each symmetric and finite right-gen\-er\-at\-ing set $S$ of $\mathcal{R}$ such that $G_0 \in S$, the sequence $\sequence{\ball(\rho)}_{\rho \in \N_0}$ is a right Erling net in $\mathcal{R}$. 
  \end{lemma}

  \begin{proof}
    Let $S$ be a symmetric and finite right-gen\-er\-at\-ing set of $\mathcal{R}$ with $G_0 \in S$, let $\net{F_i}_{i \in I}$ be the sequence $\sequence{\ball(\rho)}_{\rho \in \N_0}$, and let $E$ be a finite subset of $G \modulo G_0$. There is a non-negative integer $\rho$ such that $m_0 \isSemiActedUponBy E \subseteq \ball(\rho)$ and there is a subset $E'$ of $G \modulo G_0$ such that $\ball(\rho) = m_0 \isSemiActedUponBy E'$. Note that, because $\isSemiActedUponBy$ is free, we have $E \subseteq E'$. 

    Let $i \in I$ such that $i \geq \rho + 1$. We have $\boundaryOf_E^- F_i \subseteq \boundaryOf_{E'}^- F_i$ and $\boundaryOf_{E'}^- F_i = F_i \smallsetminus F_{i - \rho}$. Moreover, because $G_0 \in S$, for each element $m \in F_i \smallsetminus F_{i - \rho}$, there is an element $m' \in F_{i - \rho} \smallsetminus F_{i - \rho - 1}$ and there is a family $\family{s_k}_{k \in \setOf{1, 2, \dotsc, \rho}}$ in $S$ such that $(((m' \isSemiActedUponBy s_1) \isSemiActedUponBy s_2) \isSemiActedUponBy \dotsb) \isSemiActedUponBy s_\rho = m$. Thus, 
    \begin{equation*} 
      \cardinalityOf{\boundaryOf_{E'}^- F_i}
      = \cardinalityOf{F_i \smallsetminus F_{i - \rho}}
      \leq \cardinalityOf{S}^\rho \cdot \cardinalityOf{F_{i - \rho} \smallsetminus F_{i - \rho - 1}}
      \leq \cardinalityOf{S}^\rho \cdot \cardinalityOf{F_{i - \rho}}. 
    \end{equation*}
    Furthermore, $\cardinalityOf{F_i} = \cardinalityOf{F_i \smallsetminus F_{i - \rho}} + \cardinalityOf{F_{i - \rho}} = \cardinalityOf{\boundaryOf_{E'}^- F_i} + \cardinalityOf{F_{i - \rho}}$. Hence,
    \begin{equation*}
      \frac{\cardinalityOf{F_i}}{\cardinalityOf{\boundaryOf_E^- F_i}} 
      \geq \frac{\cardinalityOf{F_i}}{\cardinalityOf{\boundaryOf_{E'}^- F_i}}
      =    1 + \frac{\cardinalityOf{F_{i - \rho}}}{\cardinalityOf{\boundaryOf_{E'}^- F_i}}
      \geq 1 + \frac{\cardinalityOf{F_{i - \rho}}}{\cardinalityOf{S}^\rho \cdot \cardinalityOf{F_{i - \rho}}}
      =    1 + \cardinalityOf{S}^{-\rho}.
    \end{equation*}
    Therefore,
    \begin{equation*}
      \limsup_{i \in I} \frac{\cardinalityOf{\boundaryOf_E^- F_i}}{\cardinalityOf{F_i}} \leq \frac{1}{1 + \cardinalityOf{S}^{-\rho}} < 1.
    \end{equation*} 
    In conclusion, $\net{F_i}_{i \in I}$ is a right Erling net and hence $\mathcal{R}$ is right tractable.
  \end{proof}

  \begin{remark}
    In the proof of \cref{lemma:finitely-right-generated-is-tractable}, we do not need the right-gen\-er\-at\-ing set to be symmetric. It is assumed merely for consistency, because we define balls solely for such generating sets. The reason is that only for those is the distance on Cayley graphs a metric and balls behave nicely.
  \end{remark}

  \begin{example}[Finitely] 
    Let $G$ be a finitely generated group, let $M$ be the vertices of a Cayley graph of $G$, and let $\actsOnPoint$ be the left group action of $G$ on $M$ by left multiplication. The cell space $\ntuple{\ntuple{M, G, \actsOnPoint}, \ntuple{e_G, \family{m}_{m \in M}}}$ is finitely and symmetrically right generated and hence right tractable.
  \end{example} 

  \begin{example}[Lattice]
  \label{example:lattice:erling}
    In the situation of \cref{example:lattice:Folner}, for each non-negative integer $\rho$ and each non-negative integer $\varrho$ such that $\rho \geq \varrho$, the interior $\ball(\varrho)$-boundary of $\ball(\rho)$ is equal to $\ball(\rho) \smallsetminus \ball(\rho - \varrho)$ and its cardinality is equal to $2\varrho (2\rho - \varrho + 1)$. Hence,
    \begin{equation*}
      \ForEach \varrho \in \N_0 \Holds \limsup_{\rho \to \infty} \frac{\cardinalityOf{\boundaryOf_{\ball(\varrho)}^- \ball(\rho)}}{\cardinalityOf{\ball(\rho)}} = 0.
    \end{equation*}
    Therefore, the sequence $\sequence{\ball(\rho)}_{\rho \in \N_0}$ is a right Erling net in $\mathcal{R}$ and hence the cell space $\mathcal{R}$ is right tractable.
  \end{example} 

  \begin{example}[Tree]
  \label{example:tree:erling}
    In the situation of \cref{example:tree:Folner}, for each non-negative integer $\rho$ and each non-negative integer $\varrho$ such that $\rho \geq \varrho$, the interior $\ball(\varrho)$-boundary of $\ball(\rho)$ is equal to $\ball(\rho) \smallsetminus \ball(\rho - \varrho)$ and its cardinality is equal to $2 \cdot 3^\rho (1 - 3^{-\varrho})$. Hence,
    \begin{equation*}
      \ForEach \varrho \in \N_0 \Holds \limsup_{\rho \to \infty} \frac{\cardinalityOf{\boundaryOf_{\ball(\varrho)}^- \ball(\rho)}}{\cardinalityOf{\ball(\rho)}} = 1 - 3^{-\varrho}.
    \end{equation*}
    Therefore, the sequence $\sequence{\ball(\rho)}_{\rho \in \N_0}$ is a right Erling net in $\mathcal{R}$ and hence the cell space $\mathcal{R}$ is right tractable. However, as we have seen in \cref{example:tree:Folner}, that sequence is not a right Følner net and that cell space is not right amenable.
  \end{example} 

  \section{Tilings} 
  \label{section:tilings}

  In this section, let $\mathcal{R} = \ntuple{\mathcal{M}, \mathcal{K}} = \ntuple{\ntuple{M, G, \actsOnPoint}, \ntuple{m_0, \family{g_{m_0, m}}_{m \in M}}}$ be a cell space.

  \paragraph{Contents.} In \cref{definition:tiling} we introduce the notion of $\ntuple{E,E'}$-tilings. In \cref{theorem:existence-of-tiling} we show using Zorn's lemma that, for each non-empty subset $E$ of $G \modulo G_0$, there is an $\ntuple{E, E'}$-tiling. And in lemma \ref{lemma:upper-bound-of-tiling-cap-Folner-net} we show that, for each right Erling net $\net{F_i}_{i \in I}$ and each $\ntuple{E, E'}$-tiling with finite sets $E$ and $E'$, the net $\net{\cardinalityOf{T \cap F_i^{-E}}}_{i \in I}$ is asymptotically not less than $\net{\cardinalityOf{F_i}}_{i \in I}$. 

  \begin{definition} 
    Let $\family{M_i}_{i \in I}$ be a family of subsets of $M$. The family $\family{M_i}_{i \in I}$ is called
    \begin{aenumerate}
      \item \defineX{pairwise disjoint}{pairwise disjoint!family}\graffito{pairwise disjoint} if and only if
            \begin{equation*}
              \ForEach i \in I \ForEach i' \in I \Holds (i \neq i' \implies M_i \cap M_{i'} = \emptyset);
            \end{equation*}
      \item \defineX{cover of $M$}{cover of $M$!family}\graffito{cover of $M$} if and only if $\bigcup_{i \in I} M_i = M$;
      \item \defineX{partition of $M$}{partition of $M$!family}\graffito{partition of $M$} if and only if it is pairwise disjoint and a cover of $M$. \qedhere
    \end{aenumerate}
  \end{definition}


  \begin{definition} 
  \label{definition:tiling}
    Let $T$ be a subset of $M$, and let $E$ and $E'$ be two subsets of $G \modulo G_0$. The set $T$ is called \defineX{$\ntuple{E, E'}$-tiling of $\mathcal{R}$}{tiling of $R$@$\ntuple{E, E'}$-tiling of $\mathcal{R}$}\graffito{$\ntuple{E, E'}$-tiling of $\mathcal{R}$} if and only if the family $\family{t \isSemiActedUponBy E}_{t \in T}$ is pairwise disjoint and the family $\family{t \isSemiActedUponBy E'}_{t \in T}$ is a cover of $M$. 
  \end{definition}

  \begin{remark}
    The set $T$ is an $\ntuple{E, E}$-tiling of $\mathcal{R}$ if and only if the family $\family{t \isSemiActedUponBy E}_{t \in T}$ is a partition of $M$.
  \end{remark}

  \begin{remark} 
  \label{remark:E-prime-of-tiling-non-empty}
    Let $T$ be an $\ntuple{E, E'}$-tiling of $\mathcal{R}$. Because $M$ is non-empty and $\family{t \isSemiActedUponBy E'}_{t \in T}$ is a cover of $M$, the set $E'$ is non-empty.
  \end{remark}

  \begin{remark} 
  \label{remark:E-and-E-prime-of-tiling-can-be-reduced-or-enlarged}
    Let $T$ be an $\ntuple{E, E'}$-tiling of $\mathcal{R}$. For each subset $F$ of $E$ and each superset $F'$ of $E'$ with $F' \subseteq G \modulo G_0$, the set $T$ is an $\ntuple{F, F'}$-tiling of $\mathcal{R}$. In particular, the set $T$ is an $\ntuple{E, E \cup E'}$-tiling of $\mathcal{R}$.
  \end{remark}


  \begin{remark}
  \label{remark:group:tilings}
    In the situation of \cref{remark:group:interior-closure-boundary}, the notion of $\ntuple{E, E'}$-tiling is the same as the one defined in paragraph~2 of section~5.6 in \cite{ceccherini-silberstein:coornaert:2010}.
  \end{remark}

  \begin{example}[Lattice]
  \label{example:lattice:tiling} 
    In the situation of \cref{example:lattice:erling}, let $E$ be the ball $\ball(1)$, let $E'$ be the set $\setOf{m - m' \suchThat m, m' \in E}$, which is the ball $\ball(2)$, and let $T$ be the set $\setOf{(z_1, z_2) \in M \suchThat (z_1, z_2) \in 2\Z^2, z_1 + z_2 \in 4\Z}$. The set $T$ is an $\ntuple{E, E'}$-tiling of $\mathcal{R}$ (see \cref{figure:lattice:tiling}). 
    \begin{figure}
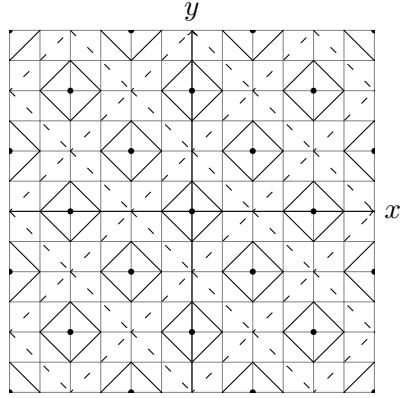

      \myfloatalign
      \figureLatticeTiling
      \caption{The grid points are elements of $M = \Z^2$; the dots are elements of the tiling $T$; for each element $t \in T$, the grid points in the region enclosed by the diamond with solid border about $t$ is the ball $t + E = \ball(t, 1)$ and the grid points in the region enclosed by the diamond with dashed border about $t$ is the ball $t + E' = \ball(t, 2)$.}
      \label{figure:lattice:tiling}
    \end{figure}
  \end{example}

  \begin{example}[Tree] 
  \label{example:tree:tiling} 
    In the situation of \cref{example:tree:erling}, recall that $\lengthOf{\blank}$ is the map $\distanceOf(m_0, \blank)$, let $E$ be the ball $\ball(1)$, let $E'$ be the set $\setOf{e (e')^{-1} \suchThat e, e' \in E}$, which is the ball $\ball(2)$, and let $T$ be the smallest subset of $M$ such that $m_0 \in T$ and, for each element $t \in T$, each element $x \in \setOf{a, b, a^{-1}, b^{-1}}$, and each element $y \in \setOf{a, b, a^{-1}, b^{-1}}$, we have $t x y^2 \in T$ if and only if $\lengthOf{t x y^2} = \lengthOf{t x} + 2$ and $\distanceOf(t x y^2, t) = 3$, in other words, if and only if $t = m_0$ and $y \neq x^{-1}$, or $t \neq m_0$, the last symbol of the reduced word that represents $t$ is not $x^{-1}$, and $y \neq x^{-1}$, or $t \neq m_0$, the last symbol is $x^{-1}$, and $y \notin \setOf{x, x^{-1}}$. The set $T$ is an $\ntuple{E, E}$-tiling of $\mathcal{R}$, in particular, because $E \subseteq E'$, it is an $\ntuple{E, E'}$-tiling of $\mathcal{R}$ (see \cref{figure:tree:tiling:balls}). 
    \begin{figure}
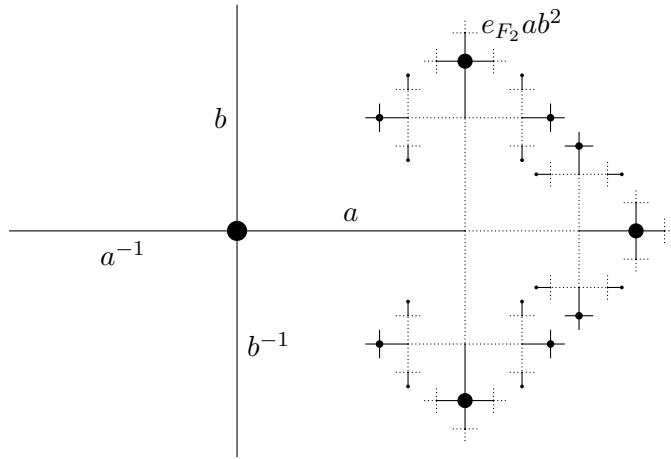

      \myfloatalign
      \figureTreeTilingBalls
      \caption{A part of the $\setOf{a, b, a^{-1}, b^{-1}}$-Cayley graph of $F_2$ is depicted, the leftmost dot is the neutral element $e_{F_2}$, the dots are elements of $T$, and the vertices that are adjacent to solid edges are elements of $\bigcup_{t \in T} t E$. The neutral element $e_{F_2}$ is contained in $T$ and, figuratively speaking, if you stand at an element $t \in T$, you take one step in a direction $x \in \setOf{a, b, a^{-1}, b^{-1}}$ to reach the leaf $t x$ of the cross $t E$ and you take another two steps in a direction $y \in \setOf{a, b, a^{-1}, b^{-1}}$ that leads away from $e_{F_2}$ and away from $t$, then you reach the element $t x y^2$, which is contained in $T$, and all elements of $T$ can be reached in that way.} 
      \label{figure:tree:tiling:balls}
    \end{figure}

    The set $\setOf{m a^{3z} \suchThat m \in M, z \in \Z, \lengthOf{m a^{\pm 1}} = \lengthOf{m} + 1}$ is a $\ntuple{\setOf{a^{-1}, e_{F_2}, a}, \setOf{a^{-1}, e_{F_2}, a}}$-tiling of $\mathcal{R}$ (see \cref{figure:tree:tiling:lines}). Note that, for each cell $m \in M$, we have $\lengthOf{m a^{\pm 1}} = \lengthOf{m} + 1$ if and only if the last symbol of the reduced word that represents $m$ is neither $a$ nor $a^{-1}$ but $b$ or $b^{-1}$, and that the set $\setOf{m a^{3z} \suchThat z \in \Z}$ contains every third element of the set $\setOf{m a^z \suchThat z \in \Z}$, which, in \cref{figure:tree:tiling:lines}, is a horizontal bi-infinite path in the $\setOf{a, b, a^{-1}, b^{-1}}$-Cayley graph of $F_2$. 
    \begin{figure}
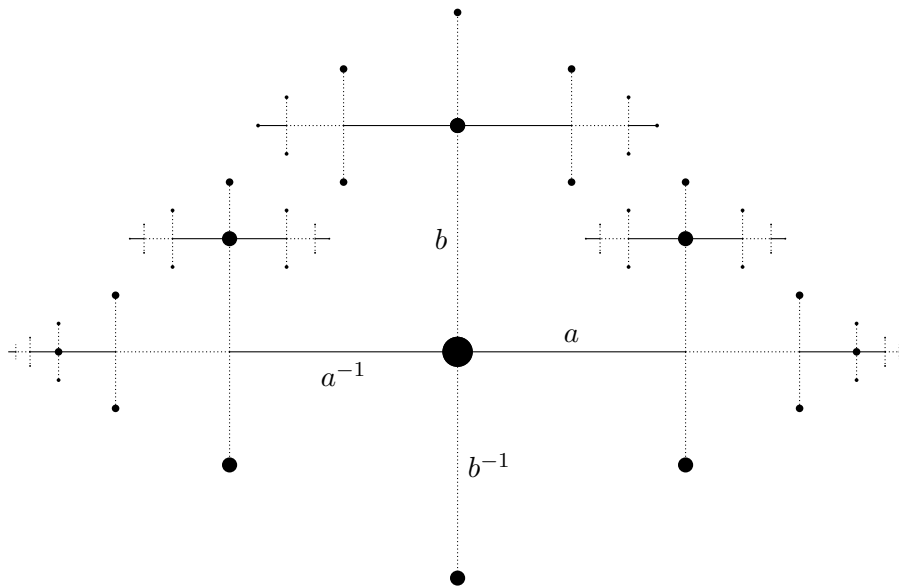

      \myfloatalign
      \figureTreeTilingLines
      \caption{A part of the $\setOf{a, b, a^{-1}, b^{-1}}$-Cayley graph of $F_2$ is depicted, the largest dot is the neutral element $e_{F_2}$, the dots are elements of $T$, and the vertices that are adjacent to solid edges are elements of $\bigcup_{t \in T} t E$.}
      \label{figure:tree:tiling:lines}
    \end{figure}
  \end{example}

  \begin{example}[Sphere]
  \label{example:sphere:tiling}
    In the situation of \cref{example:sphere:interior-closure-and-boundary}, let $E'$ be the set $\setOf{g (g')^{-1} G_0 \suchThat e, e' \in E, g \in e, g' \in e'}\ (= \setOf{g_0 g g_0' g^{-1} G_0 \suchThat g_0, g_0' \in G_0})$ and, for each point $m \in M$, let $E'_m = m \isSemiActedUponBy E'$. Because $g^{-1}$ is the rotation about the axis $a$ by $-\rho$ radians, the set $G_0 g^{-1} \actsOnPoint m_0$ is equal to $E_{m_0}$ and the set $g G_0 g^{-1} \actsOnPoint m_0$ is equal to $E_{g \actsOnPoint m_0}$. Because $m_0 \isSemiActedUponBy E' = g_{m_0, m_0} G_0 g G_0 g^{-1} \actsOnPoint m_0 = G_0 \actsOnPoint (g G_0 g^{-1} \actsOnPoint m_0) = G_0 \actsOnPoint E_{g \actsOnPoint m_0}$, the set $E'_{m_0}$ is the curved circular disk of radius $2 \rho$ with the north pole $m_0$ at its centre. And, for each point $m \in M$, because $m \isSemiActedUponBy E' = g_{m_0, m} \actsOnPoint E'_{m_0}$, the set $E'_m$ is the curved circular disk of radius $2 \rho$ with $m$ at its centre (see \cref{figure:tiling-of-sphere:general-alpha}). 
    \begin{figure} 
      \myfloatalign
      \includegraphics[trim = 75px 70px 50px 85px, clip]{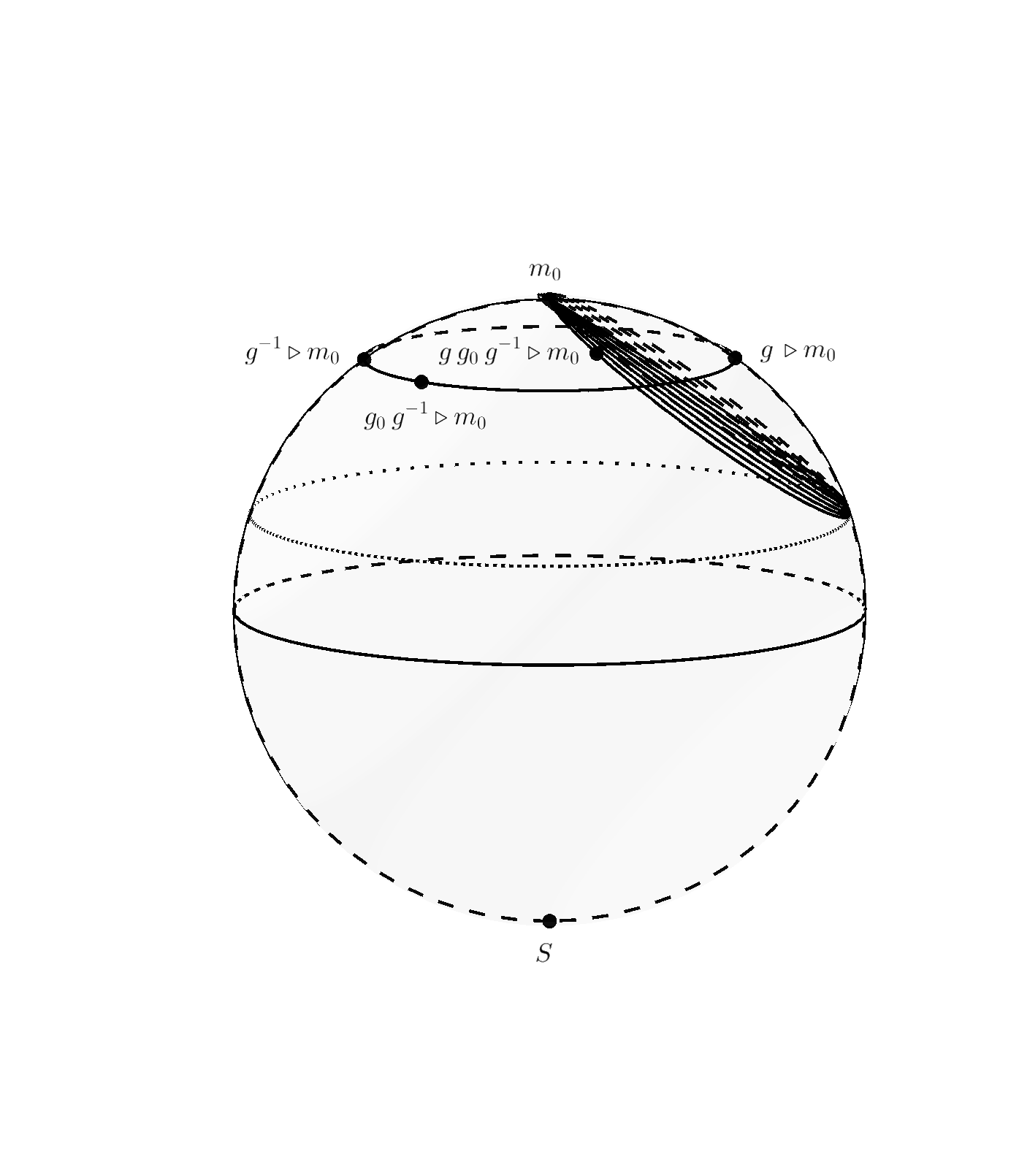}
      \caption{The circle $E_{m_0}$ is drawn solid; the boundary of the curved circular disk $E'_{m_0}$ is drawn dotted; the inclined circle on which $g g_0 g^{-1} \actsOnPoint m_0$ lies is the rotation $E_{g \actsOnPoint m_0}$ of $E_{m_0}$ by $\pi / 3$ about the axis $a$; and the other inclined circles are rotations $g_0 \actsOnPoint (E_{g \actsOnPoint m_0})$ of $E_{g \actsOnPoint m_0}$ about the $z$-axis, for a few $g_0 \in G_0$.}
      \label{figure:tiling-of-sphere:general-alpha}
    \end{figure}

    If the radius $\rho = \pi / 2$, then the circle $E_{m_0}$ is the equator and the curved circular disk $E'_{m_0}$ has radius $\pi$ and is thus the sphere $M$, and hence the set $T = \setOf{m_0}$ is an $\ntuple{E, E'}$-tiling of $\mathcal{R}$ (see \cref{figure:tiling-of-sphere:pi-half});
    \begin{figure}
      \myfloatalign
      \includegraphics[trim = 35px 70px 20px 85px, clip]{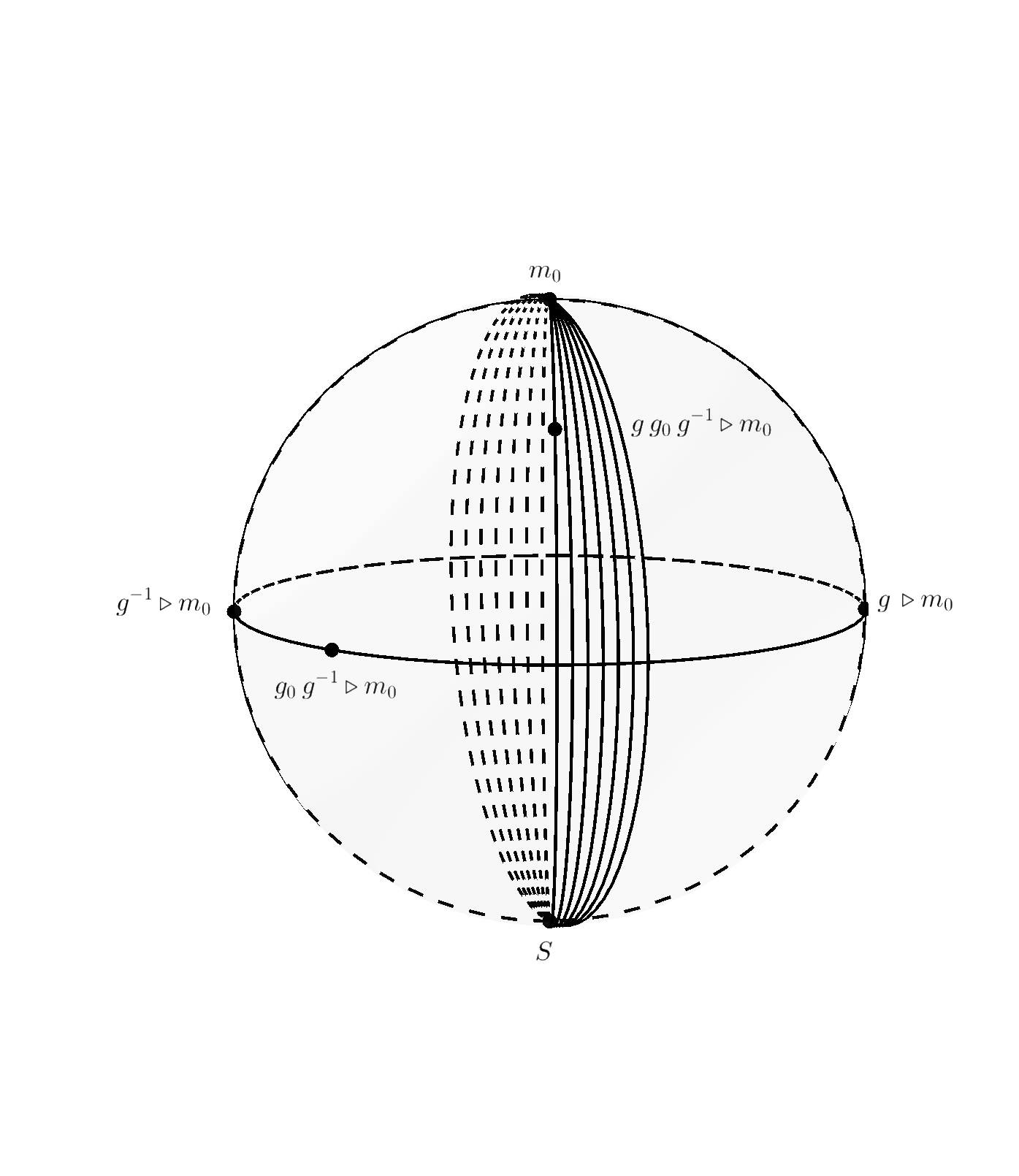}
      \caption{The set $\setOf{m_0}$ is an $\ntuple{E, E'}$-tiling of the sphere $\mathcal{R}$; the circle $E_{m_0}$ is the equator; the curved circular disk $E'_{m_0}$ is the sphere; the vertical circle on which $g g_0 g^{-1} \actsOnPoint m_0$ lies is the rotation $E_{g \actsOnPoint m_0}$ of $E_{m_0}$ by $\pi / 2$ about the axis $a$ and the other vertical circles are rotations $g_0 \actsOnPoint (E_{g \actsOnPoint m_0})$ of $E_{g \actsOnPoint m_0}$ about the $z$-axis, for a few $g_0 \in G_0$.}
      \label{figure:tiling-of-sphere:pi-half}
    \end{figure}
    if the radius $\rho = \pi / 4$, then the curved circular disks $E'_{m_0}$ and $E'_{S}$, where $S$ is the south pole, have radii $\pi / 2$, thus they are hemispheres, and hence the set $T = \setOf{m_0, S}$ is an $\ntuple{E, E'}$-tiling of $\mathcal{R}$ (see \cref{figure:tiling-of-sphere:pi-fourth});
    \begin{figure}
      \myfloatalign
      \includegraphics[trim = 65px 70px 45px 85px, clip]{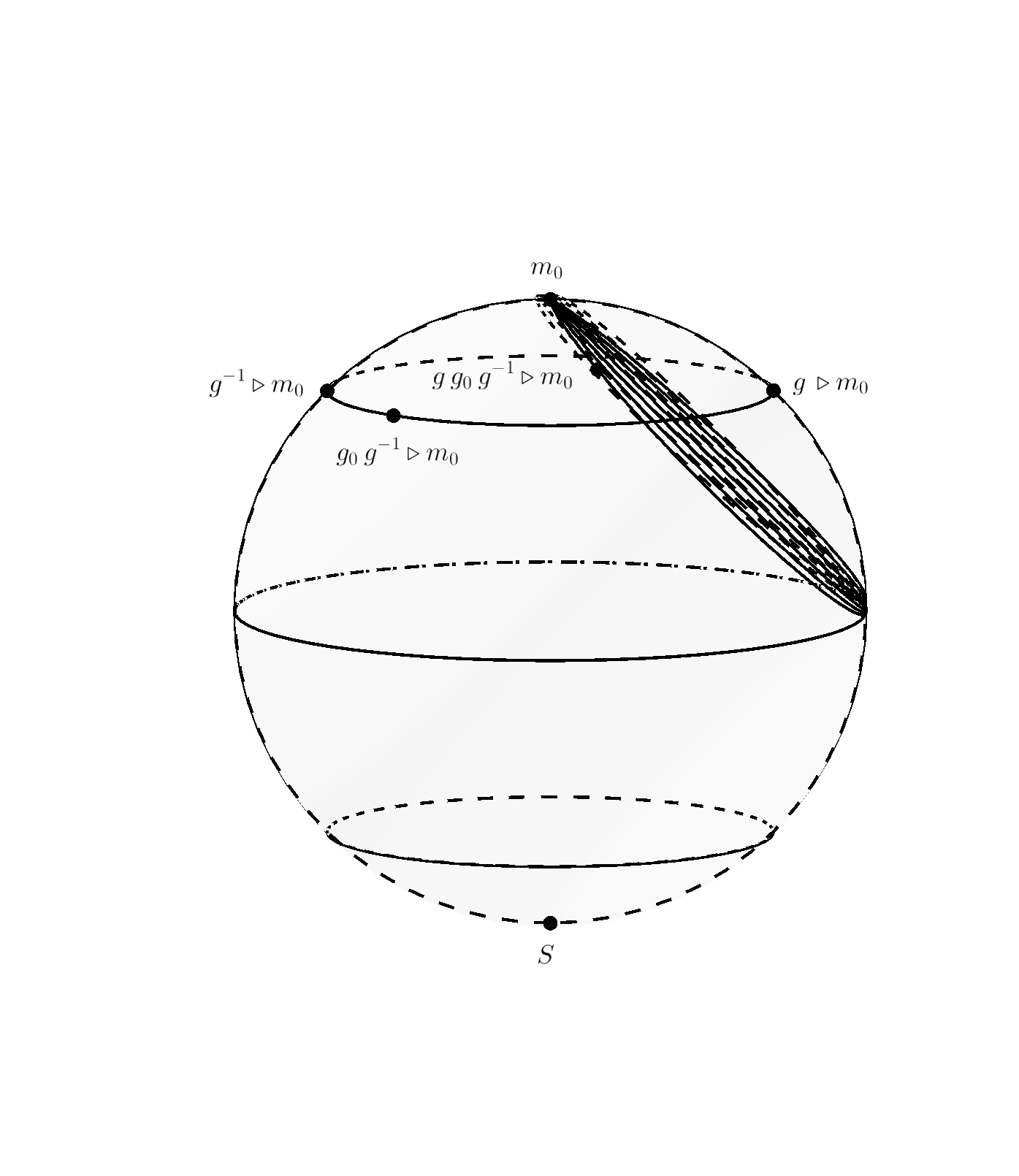}
      \caption{The set $\setOf{m_0, S}$ is an $\ntuple{E, E'}$-tiling of the sphere $\mathcal{R}$; the circles $E_{m_0}$ and $E_S$ are drawn solid; the boundaries of the curved circular disks $E'_{m_0}$ and $E'_S$ are the equator; the inclined circle on which $g g_0 g^{-1} \actsOnPoint m_0$ lies is the rotation $E_{g \actsOnPoint m_0}$ of $E_{m_0}$ by $\pi / 4$ about the axis $a$ and the other inclined circles are rotations $g_0 \actsOnPoint (E_{g \actsOnPoint m_0})$ of $E_{g \actsOnPoint m_0}$ about the $z$-axis, for a few $g_0 \in G_0$.}
      \label{figure:tiling-of-sphere:pi-fourth}
    \end{figure}
    if the radius $\rho = \pi / 8$, then the curved circular disks $E'_{m_0}$ and $E'_{S}$ have radii $\pi / 4$, and it can be shown with spherical geometry that the set $T$ --- consisting of the north pole $m_0$, the south pole $S$, four equidistant points $m_1$, $m_2$, $m_3$, and $m_4$ on the equator, and the circumcentres $c_1$, $c_2$, $\dotsc$, $c_8$ of the $8$ smallest spherical triangles with one vertex from $\setOf{m_0, S}$ and two vertices from $\setOf{m_1, m_2, m_3, m_4}$ --- is an $\ntuple{E, E'}$-tiling of $\mathcal{R}$ (see \cref{figure:tiling-of-sphere:pi-eighth}).
    \begin{figure} 
      \myfloatalign
      \figureTilingOfSpherePiEighth
      \caption{The points $m_0, S$, $m_1, m_2, m_3, m_4$, $c_1, c_2, \dotsc, c_8$ constitute an $\ntuple{E, E'}$-tiling of the sphere; the circles $E_m$ about these points are drawn solid; the boundaries of the curved circular disks $E'_m$ about these points are drawn dotted; the inclined circle about $g \actsOnPoint m_0$ is the rotation $E_{g \actsOnPoint m_0}$ of $E_{m_0}$ by $\pi / 8$ about the axis $a$; and the other inclined circles are rotations $g_0 \actsOnPoint (E_{g \actsOnPoint m_0})$ of $E_{g \actsOnPoint m_0}$ about the $z$-axis, for a few $g_0 \in G_0$.}
      \label{figure:tiling-of-sphere:pi-eighth}
    \end{figure}
  \end{example}

  Tilings exist as shown in

  \begin{theorem} 
  \label{theorem:existence-of-tiling}
    Let $E$ be a non-empty subset of $G \modulo G_0$. There is an $\ntuple{E, E'}$-tiling of $\mathcal{R}$, where $E' = \setOf{g (g')^{-1} G_0 \suchThat e, e' \in E, g \in e, g' \in e'}$.  
  \end{theorem}

  \begin{proof}
      Let
        $\mathcal{S} = \setOf{S \subseteq M \suchThat \family{s \isSemiActedUponBy E}_{s \in S} \text{ is pairwise disjoint}}$.
      Because $\setOf{m_0} \in \mathcal{S}$, the set $\mathcal{S}$ is non-empty. Moreover, it is preordered by inclusion.

      Let $\mathcal{C}$ be a chain in $\ntuple{\mathcal{S}, \subseteq}$. Then, $\bigcup_{S \in \mathcal{C}} S$ is an element of $\mathcal{S}$ and an upper bound of $\mathcal{C}$. According to Zorn's \cref{lemma:Zorns-lemma}, there is a maximal element $T$ in $\mathcal{S}$. By definition of $\mathcal{S}$, the family $\family{t \isSemiActedUponBy E}_{t \in T}$ is pairwise disjoint. 

      Let $m \in M$. Because $T$ is maximal and $m \isSemiActedUponBy E$ is non-empty, there is a $t \in T$ such that $(t \isSemiActedUponBy E) \cap (m \isSemiActedUponBy E) \neq \emptyset$. Hence, there are $e$, $e' \in E$ such that $t \isSemiActedUponBy e = m \isSemiActedUponBy e'$. According to \cref{lemma:rightsemiaction-can-be-undone}, there is a $g' \in e'$ such that $(m \isSemiActedUponBy e') \isSemiActedUponBy (g')^{-1} G_0 = m$, and there is a $g \in e$ such that $(t \isSemiActedUponBy e) \isSemiActedUponBy (g')^{-1} G_0 = t \isSemiActedUponBy g (g')^{-1} G_0$. Therefore, $m = t \isSemiActedUponBy g (g')^{-1} G_0$ (see \cref{figure:existence-of-tiling}). Because $g (g')^{-1} G_0 \in E'$, we have $m \in t \isSemiActedUponBy E'$. Thus, $\family{t \isSemiActedUponBy E'}_{t \in T}$ is a cover of $M$. 
      \begin{figure}
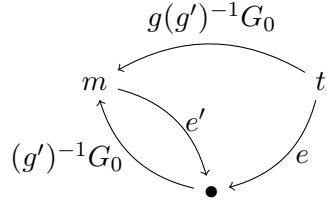

        \myfloatalign
        \figureExistenceOfTiling
        \caption{Schematic representation of a part of the proof of \cref{theorem:existence-of-tiling}.}
        \label{figure:existence-of-tiling}
      \end{figure}

      In conclusion, $T$ is an $\ntuple{E, E'}$-tiling of $\mathcal{R}$.
  \end{proof}


  Tilings have asymptotically many points in common with right Erling nets as shown in

  \begin{lemma} 
  \label{lemma:upper-bound-of-tiling-cap-Folner-net}
    Let $G_0$ be finite, let $\net{F_i}_{i \in I}$ be a right Erling net in $\mathcal{R}$ indexed by $(I, \leq)$, let $E$ and $E'$ be two finite subsets of $G \modulo G_0$, and let $T$ be an $\ntuple{E, E'}$-tiling of $\mathcal{R}$. There is a positive real number $\varepsilon \in \R_{> 0}$ and there is an index $i_0 \in I$ such that, for each index $i \in I$ with $i \geq i_0$, we have $\cardinalityOf{T \cap F_i^{-E}} \geq \varepsilon \cardinalityOf{F_i}$.
  \end{lemma}

  \begin{proof-sketch} 
    By the definition of internal boundaries, we have $\cardinalityOf{F_i} - \cardinalityOf{\boundaryOf_{E' \cdot E}^- F_i} \leq \cardinalityOf{F_i^{- E' \cdot E}}$. And, for nice $E$ and $E'$, we have $F_i^{- E' \cdot E} = (F_i^{-E})^{-E'}$. And, because $M = T \isSemiActedUponBy E'$, we have $\cardinalityOf{(F_i^{-E})^{-E'}} \leq \cardinalityOf{T \cap F_i^{-E}} \cdot \cardinalityOf{E'}$. Hence,
    \begin{equation*}
      \frac{\cardinalityOf{T \cap F_i^{-E}}}{\cardinalityOf{F_i}} \geq \frac{1}{\cardinalityOf{E'}} \cdot \parens*{1 - \frac{\cardinalityOf{\boundaryOf_{E' \cdot E}^- F_i}}{\cardinalityOf{F_i}}}.
    \end{equation*}
    For great enough indices $i \in I$, the right side is bounded below away from zero.
  \end{proof-sketch}

  \begin{proof}
    According to \cref{remark:E-and-E-prime-of-tiling-can-be-reduced-or-enlarged}, we may suppose, without loss of generality, that $E \subseteq E'$, $G_0 \cdot E' \subseteq E'$, and that $(E')^{-1} \subseteq E'$, where $(E')^{-1} = \setOf{(g')^{-1} G_0 \suchThat e' \in E', g' \in e'}$. 

    Let $i \in I$. Put
    \begin{equation*}
      T_i = T \cap F_i^{-E} = \setOf{t \in T \suchThat t \isSemiActedUponBy E \subseteq F_i}
    \end{equation*}
    and put
    \begin{equation*}
      T_i' = T \cap ((F_i^{-E})^{-E'})^{+E'} = \setOf{t \in T \suchThat (t \isSemiActedUponBy E') \cap (F_i^{-E})^{-E'} \neq \emptyset}
    \end{equation*}
    (see \cref{figure:upper-bound-of-tiling-cap-Folner-net}).
    \begin{figure}
      \myfloatalign
      \figureUpperBoundOfTilingCapFolnerNet
      \caption{Schematic representation of the set-up of the proof of \cref{lemma:upper-bound-of-tiling-cap-Folner-net}.}
      \label{figure:upper-bound-of-tiling-cap-Folner-net}
    \end{figure}
    Because the family $\family{t \isSemiActedUponBy E'}_{t \in T}$ is a cover of $M$,
    \begin{equation*}
      (F_i^{-E})^{-E'}
      = M \cap (F_i^{-E})^{-E'}
      = \bigcup_{t \in T} (t \isSemiActedUponBy E') \cap (F_i^{-E})^{-E'}.
    \end{equation*}
    And, for each element $t \in T \smallsetminus T_i'$, we have $(t \isSemiActedUponBy E') \cap (F_i^{-E})^{-E'} = \emptyset$. Hence,
    \begin{equation*}
      (F_i^{-E})^{-E'}
      = \bigcup_{t \in T_i'} (t \isSemiActedUponBy E') \cap (F_i^{-E})^{-E'}.
    \end{equation*}
    And, for each element $t \in T_i'$, we have $(t \isSemiActedUponBy E') \cap (F_i^{-E})^{-E'} \subseteq t \isSemiActedUponBy E'$. Thus,
    \begin{equation*}
      (F_i^{-E})^{-E'}
      \subseteq \bigcup_{t \in T_i'} t \isSemiActedUponBy E'.
    \end{equation*}
    And, according to \cref{item:properties-of-interior-closure-and-boundary:different-repeated} of \cref{lemma:properties-of-interior-closure-and-boundary}, we have $((F_i^{-E})^{-E'})^{+E'} \subseteq F_i^{-E}$ and hence $T_i' \subseteq T_i$. Therefore,
    \begin{equation*}
      (F_i^{-E})^{-E'}
      \subseteq \bigcup_{t \in T_i} t \isSemiActedUponBy E'.
    \end{equation*}
    And, because $\isSemiActedUponBy$ is free, for each $t \in T_i$, we have $\cardinalityOf{t \isSemiActedUponBy E'} = \cardinalityOf{E'}$. Hence,
    \begin{equation*}
      \cardinalityOf{(F_i^{-E})^{-E'}}
      \leq \cardinalityOf{T_i} \cdot \cardinalityOf{E'}.
    \end{equation*}
    And, according to Items~\cref{item:properties-of-interior-closure-and-boundary:inclusions} and~\ref{item:properties-of-interior-closure-and-boundary:same-repeated} of \cref{lemma:properties-of-interior-closure-and-boundary}, we have $(F_i^{-E})^{-E'} \supseteq (F_i^{- G_0 \cdot E})^{-E'} = F_i^{- E''}$, where $E'' = E' \cdot (G_0 \cdot E) = \setOf{g' \cdot (g_0 \cdot e) \suchThat e \in E, g_0 \in G_0, e' \in E', g' \in e'}$. And, because $\boundaryOf_{E''}^- F_i = F_i \smallsetminus (F_i^{-E''} \cap F_i)$,
    \begin{equation*}
      \cardinalityOf{F_i^{-E''}}
      \geq \cardinalityOf{F_i^{-E''} \cap F_i} 
      =    \cardinalityOf{F_i} - \cardinalityOf{\boundaryOf_{E''}^- F_i}.
    \end{equation*}
    Hence,
    \begin{equation*}
      \cardinalityOf{F_i} - \cardinalityOf{\boundaryOf_{E''}^- F_i} \leq \cardinalityOf{T_i} \cdot \cardinalityOf{E'}.
    \end{equation*}
    Therefore,
    \begin{equation*}
      \frac{\cardinalityOf{T_i}}{\cardinalityOf{F_i}}
      \geq \frac{1}{\cardinalityOf{E'}} \cdot \parens*{1 - \frac{\cardinalityOf{\boundaryOf_{E''}^- F_i}}{\cardinalityOf{F_i}}}.
    \end{equation*}
    Because $\net{F_i}_{i \in I}$ is a right Erling net, there is a real number $\xi \in \leftClosedAndRightOpenInterval{0, 1}$ and there is an index $i_0 \in I$ such that
    \begin{equation*} 
      \ForEach i \in I \Holds \parens*{i \geq i_0 \implies \frac{\cardinalityOf{\boundaryOf_{E''}^- F_i}}{\cardinalityOf{F_i}} \leq \xi}.
    \end{equation*}
    Put $\varepsilon = (1/\cardinalityOf{E''}) \cdot (1 - \xi)$. Then, for each $i \in I$ with $i \geq i_0$, we have $\cardinalityOf{T_i}/\cardinalityOf{F_i} \geq \varepsilon$.
  \end{proof}

  A tiling of a right-a\-me\-na\-ble and hence right-trac\-ta\-ble cell space is given in

  \begin{example}[Lattice]
  \label{example:lattice:upper-bound-of-tiling-cap-Folner-net} 
    In the situation of \cref{example:lattice:tiling}, for each non-negative integer $\rho$, one can see that $\cardinalityOf{T \cap \ball(\rho)} = \cardinalityOf{T \cap \ball(4 \floor{\rho/4})} = (2 \floor{\rho/4} + 1)^2$ and, for each positive integer $\rho$, recall that $\ball(\rho)^{-\ball(1)} = \ball(\rho - 1)$. 

    Let $\rho$ be a positive integer. The integer $\varrho = 4 \ceil{\rho/4}$ is the multiple of $4$ such that $\rho \leq \varrho < \rho + 4$. Thus, $\varrho - 1$ is the greatest integer such that $\floor{(\varrho - 1)/4} = \floor{(\rho - 1)/4}$, in particular, $\cardinalityOf{T \cap \ball(\rho - 1)} = \cardinalityOf{T \cap \ball(\varrho - 1)}$. Hence, because $\ball(\rho) \subseteq \ball(\varrho)$,
    \begin{equation*}
      \frac{\cardinalityOf{T \cap \ball(\rho - 1)}}{\cardinalityOf{\ball(\rho)}}
      \geq \frac{\cardinalityOf{T \cap \ball(\varrho - 1)}}{\cardinalityOf{\ball(\varrho)}}.
    \end{equation*} 
    Moreover, one can show that the sequence
    \begin{equation*}
      \sequence*{\frac{\cardinalityOf{T \cap \ball(\varsigma - 1)}}{\cardinalityOf{\ball(\varsigma)}}}_{\varsigma \in 4\N_+}
    \end{equation*}
    is increasing (and converges to $1/8$). Therefore, because $\varrho \in 4\N_+$,
    \begin{equation*}
      \frac{\cardinalityOf{T \cap \ball(\varrho - 1)}}{\cardinalityOf{\ball(\varrho)}}
      \geq \frac{\cardinalityOf{T \cap \ball(4 - 1)}}{\cardinalityOf{\ball(4)}}
      =    \frac{1}{41}.
    \end{equation*}
    In conclusion, $\cardinalityOf{T \cap \ball(\rho)^{-\ball(1)}} \geq (1/41) \cdot \cardinalityOf{\ball(\rho)}$ (actually, for each real number $\varepsilon \in \openInterval{0, 1/8}$, there is an index $\rho_0 \in \N_0$ such that, for each index $\rho \in \N_0$ with $\rho \geq \rho_0$, we have $\cardinalityOf{T \cap \ball(\rho)^{-\ball(1)}} \geq \varepsilon \cdot \cardinalityOf{\ball(\rho)}$).
  \end{example}

  A tiling of a right-trac\-ta\-ble but not right-a\-me\-na\-ble cell space is given in

  \begin{example}[Tree]
  \label{example:tree:upper-bound-of-tiling-cap-Folner-net}
    In the situation of \cref{example:tree:tiling}, by the construction of $T$, we have $m_0 \in T$; moreover, there are four pairwise distinct elements $x_{m_0, 1}$, $x_{m_0, 2}$, $x_{m_0, 3}$, and $x_{m_0, 4} \in \setOf{a, b, a^{-1}, b^{-1}}$ and, for each index $j \in \setOf{1, 2, 3, 4}$, there are three pairwise distinct elements $y_{m_0, 1}$, $y_{m_0, 2}$, and $y_{m_0, 3} \in \setOf{a, b, a^{-1}, b^{-1}}$ such that $x_{m_0, j} y_{m_0, k}^2 \in T$ and $\lengthOf{x_{m_0, j} y_{m_0, k}^2} = \lengthOf{x_{t, j}} + 2$ (note that $\lengthOf{x_{t, j}} + 2 = 3$); furthermore, for each element $t \in T \smallsetminus \setOf{m_0}$, there are three pairwise distinct elements $x_{t, 1}$, $x_{t, 2}$, and $x_{t, 3} \in \setOf{a, b, a^{-1}, b^{-1}}$ such that $\lengthOf{t x_{t, j}} = \lengthOf{t} + 1$, for $j \in \setOf{1, 2, 3}$, and, for each index $j \in \setOf{1, 2, 3}$, there are three pairwise distinct elements $y_{t, 1}$, $y_{t, 2}$, and $y_{t, 3} \in \setOf{a, b, a^{-1}, b^{-1}}$ such that $t x_{t, j} y_{t, k}^2 \in T$ and $\lengthOf{t x_{t, j} y_{t, k}^2} = \lengthOf{t x_{t, j}} + 2$ (note that $\lengthOf{t x_{t, j}} + 2 = \lengthOf{t} + 3$); and, there is an element $x_{t, 4} \in \setOf{a, b, a^{-1}, b^{-1}} \smallsetminus \setOf{x_{t, 1}, x_{t, 2}, x_{t, 3}}$ such that $\lengthOf{t x_{t, 4}} = \lengthOf{t} - 1$, and there are two distinct elements $y_{t, 4}$ and $y_{t, 4}' \in \setOf{a, b, a^{-1}, b^{-1}}$ such that $t x_{t, 4} y_{t, 4}^2$, $t x_{t, 4} (y_{t, 4}')^2 \in T$ and $\lengthOf{t x_{t, 4} y_{t, 4}^2}$ as well as $\lengthOf{t x_{t, 4} (y_{t, 4}')^2}$ is equal to $\lengthOf{t x_{t, 4}} + 2$ (note that $\lengthOf{t x_{t, 4}} + 2 = \lengthOf{t} + 1$); and, the elements $m_0$, $x_{m_0, j} y_{m_0, k}^2$, for $j \in \setOf{1, 2, 3, 4}$ and $k \in \setOf{1, 2, 3}$, $t x_{t, j} y_{t, k}^2$, $t x_{t, 4} y_{t, 4}^2$, $t x_{t, 4} (y_{t, 4}')^2$, for $t \in T \smallsetminus \setOf{m_0}$ and $j$, $k \in \setOf{1, 2, 3}$, are pairwise distinct and are the only elements of $T$.

    Therefore, broadly speaking, the origin $m_0$ contributes $4 \cdot 3 = 12$ elements to $T \cap \sphere(3)$; for each integer $\rho \in \Z_{\geq 4}$, each element $t \in T \cap \sphere(\rho - 3)$ contributes $3 \cdot 3 = 9$ elements to $T \cap \sphere(\rho)$; each element $t \in T \cap \sphere(\rho - 1)$ contributes $2$ elements to $T \cap \sphere(\rho)$; and, there are no other elements in $T \cap \sphere(\rho)$. And, because $\ball(\rho) = \ball(\rho - 1) \disjointUnionWith \sphere(\rho)$, the set $T \cap \ball(\rho)$ contains aside from the elements of $T \cap \ball(\rho - 1)$ only the elements of $T \cap \sphere(\rho)$. 

    For each non-negative integer $\rho$, let $\sigma_\rho = \cardinalityOf{T \cap \sphere(\rho)}$ and let $\beta_\rho = \cardinalityOf{T \cap \ball(\rho)}$. Then, $\sigma_0 = 1$, $\sigma_1 = 0$, $\sigma_2 = 0$, $\sigma_3 = 12$, and, for each integer $\rho \in \Z_{\geq 4}$, we have $\sigma_\rho = 9 \sigma_{\rho - 3} + 2 \sigma_{\rho - 1}$. Moreover, $\beta_0 = 1$ and, for each positive integer $\rho \in \N_+$, we have $\beta_\rho = \beta_{\rho - 1} + \sigma_\rho$. One can use mathematical software to determine a closed-form expression for $\beta_\rho$ and show that the sequence $\sequence{\beta_{\rho - 1}/\cardinalityOf{\ball(\rho)}}_{\rho \in \N_+}$
    converges to $1/15$ and hence is eventually bounded from below by, say, $1/15 - 1/30 = 1/30$. Hence, there is an index $\rho_0 \in \N_0$ such that, for each index $\rho \in \N_0$ with $\rho \geq \rho_0$, we have $\cardinalityOf{T \cap \ball(\rho)^{-\ball(1)}} \geq (1/30) \cdot \cardinalityOf{\ball(\rho)}$. 
  \end{example}

  \section{Entropies} 
  \label{section:entropies}

  In this section, let $\mathcal{R} = \ntuple{\mathcal{M}, \mathcal{K}} = \ntuple{\ntuple{M, G, \actsOnPoint}, \ntuple{m_0, \family{g_{m_0, m}}_{m \in M}}}$ be a cell space, let $\mathcal{C} = \ntuple{\mathcal{R}, Q, N, \delta}$ be a semi-cellular automaton, and let $\Delta$ be the global transition function of $\mathcal{C}$ such that the stabiliser $G_0$ of $m_0$ under $\actsOnPoint$, the set $Q$ of states, and the neighbourhood $N$ are finite, and the set $Q$ is non-empty. 

  \paragraph{Contents.} In \cref{definition:entropy} we introduce the entropy of a subset $X$ of $Q^M$ with respect to a net $\net{F_i}_{i \in I}$ of non-empty and finite subsets of $M$, which is the asymptotic growth rate of the number of finite patterns with domain $F_i$ that occur in $X$. In \cref{lemma:entropy-basic-facts} we show that $Q^M$ has entropy $\log\cardinalityOf{Q}$ and that entropy is non-decreasing. In \cref{theorem:entropy-does-non-increase} we show that applications of global transition functions of semi-cellular automata on subsets of $Q^M$ do not increase entropy. And in \cref{lemma:entropy-bounded-above-if-strange-tiling-exists} we show that if for each point $t$ of an $\ntuple{E,E'}$-tiling not all patterns with domain $t \isSemiActedUponBy E$ occur in a subset of $Q^M$, then that subset has less entropy than $Q^M$. 

  \begin{definition} 
  \label{definition:entropy}
    Let $X$ be a subset of $Q^M$ and let $\mathcal{F} = \net{F_i}_{i \in I}$ be a net in $\setOf{F \subseteq M \suchThat F \neq \emptyset, F \text{ finite}}$. The non-negative real number or negative infinity 
    \begin{equation*}
      \entropyOf_{\mathcal{F}}(X) = \limsup\limits_{i \in I} \frac{\log\cardinalityOf{\pi_{F_i}(X)}}{\cardinalityOf{F_i}}
      \mathnote{entropy $\entropyOf_{\mathcal{F}}(X)$ of $X$ with respect to $\mathcal{F}$}
    \end{equation*}
    is called \define{entropy of $X$ with respect to $\mathcal{F}$}.
  \end{definition}

  \begin{remark}
  \label{remark:group:entropy}
    In the situation of \cref{remark:group:interior-closure-boundary}, the notion of entropy is the same as the one defined in definition~5.7.1 in \cite{ceccherini-silberstein:coornaert:2010}.
  \end{remark}

  \begin{remark} 
    As already said, the entropy of $X$ with respect to the net $\mathcal{F}$ in $\mathcal{R}$ is the asymptotic growth rate of the number of finite patterns with domain $F_i$ that occur in $X$. In more precise terms but still hand-wavingly, this means that 
    \begin{equation*}
      \net{10^{\cardinalityOf{F_i} \cdot \entropyOf_{\mathcal{F}}(X)}}_{i \in I} \sim \net{\cardinalityOf{\pi_{F_i}(X)}}_{i \in I},
    \end{equation*}
    where $\sim$ is the binary relation, read \emph{asymptotic to}, given by
    \begin{equation*} 
      \ForEach \net{r_i}_{i \in I} \ForEach \net{r_i'}_{i \in I} \Holds \parens[\big]{\net{r_i}_{i \in I} \sim \net{r_i'}_{i \in I} \ifAndOnlyIf \lim_{i \in I} \frac{r_i}{r_i'} = 1}. \qedhere
    \end{equation*}
  \end{remark} 

  \begin{lemma}
  \label{lemma:entropy-basic-facts}
    Let $\mathcal{F} = \net{F_i}_{i \in I}$ be a net in $\setOf{F \subseteq M \suchThat F \neq \emptyset, F \text{ finite}}$. Then,
    \begin{aenumerate}
      \item \label{item:entropy-basic-facts:whole-space}
            $\entropyOf_{\mathcal{F}}(Q^M) = \log\cardinalityOf{Q}$;
      \item \label{item:entropy-basic-facts:monotonic}
            $\ForEach X \subseteq Q^M \ForEach X' \subseteq Q^M \Holds \parens[\big]{X \subseteq X' \implies \entropyOf_{\mathcal{F}}(X) \leq \entropyOf_{\mathcal{F}}(X')}$;
      \item \label{item:entropy-basic-facts:bound}
            $\ForEach X \subseteq Q^M \Holds \entropyOf_{\mathcal{F}}(X) \leq \log\cardinalityOf{Q}$. \qedhere
    \end{aenumerate}
  \end{lemma}

  \begin{proof}
    \begin{aenumerate}
      \item For each $i \in I$, we have $\pi_{F_i}(Q^M) = Q^{F_i}$ and hence
            \begin{equation*}
              \frac{\log\cardinalityOf{\pi_{F_i}(Q^M)}}{\cardinalityOf{F_i}}
              = \frac{\log\cardinalityOf{Q}^{\cardinalityOf{F_i}}}{\cardinalityOf{F_i}}
              = \frac{\cardinalityOf{F_i} \cdot \log\cardinalityOf{Q}}{\cardinalityOf{F_i}}
              = \log\cardinalityOf{Q}.
            \end{equation*}
            In conclusion, $\entropyOf_{\mathcal{F}}(Q^M) = \log\cardinalityOf{Q}$.
      \item Let $X$, $X' \subseteq Q^M$ such that $X \subseteq X'$. For each $i \in I$, we have $\pi_{F_i}(X) \subseteq \pi_{F_i}(X')$ and hence, because $\log$ is non-decreasing, $\log\cardinalityOf{\pi_{F_i}(X)} \leq \log\cardinalityOf{\pi_{F_i}(X')}$. In conclusion, $\entropyOf_{\mathcal{F}}(X) \leq \entropyOf_{\mathcal{F}}(X')$. 
      \item This is a direct consequence of Items~\ref{item:entropy-basic-facts:monotonic} and~\ref{item:entropy-basic-facts:whole-space}. \qedhere
    \end{aenumerate}
  \end{proof}

  Due to the locality of the global transition function $\Delta$ and the asymptotic invariance of a right Følner net $\net{F_i}_{i \in I}$ under taking finite boundaries, the information that flows into $F_i$ from the boundary under an application of $\Delta$ is asymptotically negligible and hence applications of $\Delta$ to subsets of global configurations do not increase entropy, which is shown in

  \begin{theorem} 
  \label{theorem:entropy-does-non-increase}
    Let $\mathcal{R}$ be right amenable, let $\mathcal{F} = \net{F_i}_{i \in I}$ be a right Følner net in $\mathcal{R}$ indexed by $(I, \leq)$, and let $X$ be a subset of $Q^M$. Then, $\entropyOf_{\mathcal{F}}(\Delta(X)) \leq \entropyOf_{\mathcal{F}}(X)$.
  \end{theorem} 

  \begin{proof}
    Suppose, without loss of generality, that $G_0 \in N$.
    Let $i \in I$. According to \cref{lemma:Delta-X-A-minus-plus-are-surjective}, the map
    $\Delta_{X, F_i}^- \from \pi_{F_i}(X) \to \pi_{F_i^{-N}}(\Delta(X))$
    is surjective. Therefore, $\cardinalityOf{\pi_{F_i^{-N}}(\Delta(X))} \leq \cardinalityOf{\pi_{F_i}(X)}$.
    Because $G_0 \in N$, according to \cref{item:properties-of-interior-closure-and-boundary:neutral-element} of \cref{lemma:properties-of-interior-closure-and-boundary}, we have $F_i^{-N} \subseteq F_i$. Thus, $\pi_{F_i}(\Delta(X)) \subseteq \pi_{F_i^{-N}}(\Delta(X)) \times Q^{F_i \smallsetminus F_i^{-N}}$. Hence, 
    \begin{align*} 
      \log\cardinalityOf{\pi_{F_i}(\Delta(X))}
      &\leq \log\cardinalityOf{\pi_{F_i^{-N}}(\Delta(X))} + \log\cardinalityOf{Q^{F_i \smallsetminus F_i^{-N}}}\\ 
      &\leq 
             \log\cardinalityOf{\pi_{F_i}(X)} + \cardinalityOf{F_i \smallsetminus F_i^{-N}} \cdot \log\cardinalityOf{Q}.
    \end{align*}
    Because $G_0 \in N$, according to \cref{item:properties-of-interior-closure-and-boundary:neutral-element} of \cref{lemma:properties-of-interior-closure-and-boundary}, we have $F_i \subseteq F_i^{+N}$. Therefore, $F_i \smallsetminus F_i^{-N} \subseteq F_i^{+N} \smallsetminus F_i^{-N} = \boundaryOf_N F_i$. Because $G_0$, $F_i$, and $N$ are finite, according to \cref{item:properties-of-interior-closure-and-boundary:finite} of \cref{lemma:properties-of-interior-closure-and-boundary}, the boundary $\boundaryOf_N F_i$ is finite. Hence,
    \begin{equation*}
      \frac{\log\cardinalityOf{\pi_{F_i}(\Delta(X))}}{\cardinalityOf{F_i}} \leq \frac{\log\cardinalityOf{\pi_{F_i}(X)}}{\cardinalityOf{F_i}} + \frac{\cardinalityOf{\boundaryOf_N F_i}}{\cardinalityOf{F_i}} \log\cardinalityOf{Q}.
    \end{equation*}
    Therefore, because $N$ is finite, according to \cref{theorem:boundary-characterisation-of-Folner-net},
    \begin{align*}
      \entropyOf_{\mathcal{F}}(\Delta(X))
      &\leq \limsup\limits_{i \in I} \frac{\log\cardinalityOf{\pi_{F_i}(X)}}{\cardinalityOf{F_i}} + \parens*{\lim_{i \in I} \frac{\cardinalityOf{\boundaryOf_N F_i}}{\cardinalityOf{F_i}}} \cdot \log\cardinalityOf{Q}\\
      &=    \entropyOf_{\mathcal{F}}(X). \qedhere
    \end{align*}
  \end{proof}

  \begin{counterexample}[Tree]
  \label{example:tree:entropy-does-increase}
    In the situation of \cref{example:tree:erling}, the cell space $\mathcal{R}$ is not right amenable and the sequence $\mathcal{F} = \sequence{\ball(\rho)}_{\rho \in \N_0}$ is not a right Følner net in $\mathcal{R}$. However, the cell space $\mathcal{R}$ is right tractable and the sequence $\mathcal{F}$ is a right Erling net in $\mathcal{R}$. Nevertheless, for the majority rule over $\mathcal{R}$, there is a subset $X$ of global configurations whose image has greater entropy than $X$ itself, as we show below. Broadly speaking, such a subset exists because the boundaries of components of $\mathcal{F}$ are so big that the information that flows in from them under applications of $\Delta$ is asymptotically significant. 

    Let $Q$ be the set $\setOf{0, 1}$, let $N$ be the ball $\ball(1)$, let $\delta$ be the $\bullet$-invariant map $Q^N \to Q$, $\ell \mapsto 0$, if $\sum_{n \in N} \ell(n) \leq \cardinalityOf{N}/2$, and $\ell \mapsto 1$, otherwise, which is known as \define{majority rule}\graffito{majority rule}\index{rule!majority}, and let $\mathcal{C}$ be the cellular automaton $\ntuple{\mathcal{R}, Q, N, \delta}$. Furthermore, for each non-negative integer $\rho$, let $Y_\rho$ be the set $\setOf{c \in Q^M \smallsetminus \setOf{\functionThatIsIdenticalToZero} \suchThat c\restrictedTo_{M \smallsetminus \sphere(\rho)} \equiv 0}$, let $X_{\rho + 1}$ be the set
    \begin{equation*} 
      \setOf{y \in Y_{\rho + 1} \suchThat \ForEach m' \in \sphere(\rho) \Exists q \in Q \SuchThat y\restrictedTo_{(m' \cdot N) \cap \sphere(\rho + 1)} \equiv q},
    \end{equation*} 
    let $Y$ be the set $\bigDisjointUnionOf_{\rho \in \N_0} Y_\rho$, and let $X$ be the set $\bigDisjointUnionOf_{\rho \in \N_0} X_{\rho + 1}$.

    The global transition function $\Delta$ of $\mathcal{C}$ maps $X$ bijectively onto $Y$; even more, for each non-negative integer $\rho$, it maps $X_{\rho + 1}$ bijectively onto $Y_\rho$; more precisely, for each global configuration $x \in X_{\rho + 1}$, the global configuration
    \begin{align*}
      y_x \from M &\to Q,\\
      m' &\mapsto \begin{dcases*}
                   0, &if $m' \in M \smallsetminus \sphere(\rho)$,\\
                   x(m), &if $m' \in \sphere(\rho)$, where $m \in (m' \cdot N) \cap \sphere(\rho + 1)$, 
                 \end{dcases*}
    \end{align*}
    is the unique element of $Y_\rho$ that satisfies $\Delta(x) = y_x$, and, for each global configuration $y \in Y_\rho$, the global configuration
    \begin{align*}
      x_y \from M &\to Q,\\
      m &\mapsto \begin{dcases*}
                   0, &if $m \in M \smallsetminus \sphere(\rho + 1)$,\\
                   y(m'), &if $m \in \sphere(\rho + 1)$,\\
                          &\phantom{if }where $m' \in \sphere(\rho)$ such that $m \in m' \cdot N$, 
                 \end{dcases*}
    \end{align*}
    is the unique element of $X_{\rho + 1}$ that satisfies $\Delta(x_y) = y$ (see \cref{figure:tree:entropy-does-increase:global-transition}).

    The entropy of $\Delta(X)$ with respect to $\mathcal{F}$ is greater than the entropy of $X$ with respect to $\mathcal{F}$. The reason is that, broadly speaking, the cardinality of $\pi_{\ball(i)}(\Delta(X))$ is approximately $2^{\cardinalityOf{\sphere(i)}}$, whereas the cardinality of $\pi_{\ball(i)}(X)$ is approximately $2^{\cardinalityOf{\sphere(i - 1)}}$, and the cardinality of $\ball(i)$ is approximately $\cardinalityOf{\sphere(i)}$. 
    \begin{figure}
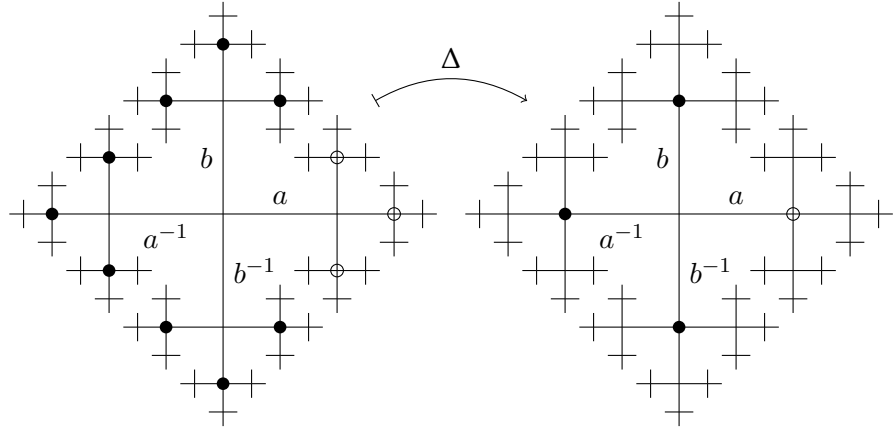

      \myfloatalign
      \figureTreeEntropyDoesIncreaseGlobalTransition
      \caption{The same part of two global configurations $x \in X_2$ and $\Delta(x) \in Y_1$ is depicted on the left and right in form of a $Q$-vertex-labelled $\setOf{a, b, a^{-1}, b^{-1}}$-Cayley graph of $F_2$, where the vertex in the centre is the neutral element $e_{F_2}$, the dotted vertices are in state $1$, the other vertices are in state $0$, the dotted and circled vertices of the left graph are the ones of $\sphere(2)$, and the ones of the right graph are the ones of $\sphere(1)$.}
      \label{figure:tree:entropy-does-increase:global-transition}
    \end{figure}
  \end{counterexample} 

  \begin{proof} 
    First, we prove that $\Delta(X) = Y$. Let $y$ be a global configuration of $Y$ and let $m'$ be a cell of $M$. Then,
    \begin{equation*}
      \Delta(x_y)(m')
      = \delta\parens[\big]{n \mapsto x_y(m' \cdot n)}
      = \begin{dcases*}
          0, &if $\sum_{n \in N} x_y(m' \cdot n) \leq \frac{5}{2}$,\\
          1, &otherwise.
        \end{dcases*}
    \end{equation*}
    Of the $5$ elements of $m' \cdot N$, the $4$ or $3$ elements of $(m' \cdot N) \cap \sphere(\lengthOf{m'} + 1)$ have the same state in $x_y$, namely $y(m')$. Thus, 
    \begin{equation*}
      \sum_{m \in m' \cdot N} x_y(m) \leq \frac{5}{2}
      \ifAndOnlyIf y(m') = 0.
    \end{equation*}
    Hence,
    \begin{equation*}
      \Delta(x_y)(m')
      = \left\{
          \begin{aligned}
            &0, &&\text{if $y(m') = 0$},\\
            &1, &&\text{otherwise},
          \end{aligned}
        \right\}
      = y(m').
    \end{equation*}
    Therefore, $\Delta(x_y) = y$. Moreover, for each $x \in X$, if $\Delta(x) = y$, then $y_x = y$, hence $x = x_{y_x} = x_y$, and therefore $x_y$ is unique. 

    Secondly, we prove that $\entropyOf_{\mathcal{F}}(\Delta(X)) > \entropyOf_{\mathcal{F}}(X)$. Recall that, for each positive integer $i$, we have $\cardinalityOf{\sphere(i)} = 3^i + 3^{i - 1}$ and $\cardinalityOf{\ball(i)} = 2 \cdot 3^i - 1$.

    Let $i$ be an integer such that $i \geq 2$. Then,
    \begin{equation*}
      \cardinalityOf{\pi_{\ball(i)}(X)}
      = \cardinalityOf{(\bigDisjointUnionOf_{\rho = 0}^{i - 1} \pi_{\ball(i)}(X_{\rho + 1})) \disjointUnionWith \setOf{\functionThatIsIdenticalToZero}} 
      = \sum_{\rho = 0}^{i - 1} (2^{\cardinalityOf{\sphere(\rho)}} - 1) + 1.
    \end{equation*}
    Thus, because $i \geq 1$,
    \begin{equation*}
      \cardinalityOf{\pi_{\ball(i)}(X)}
      \leq \sum_{\rho = 0}^{i - 1} (2^{\cardinalityOf{\sphere(i - 1)}} - 1) + 1
      =    i \cdot 2^{\cardinalityOf{\sphere(i - 1)}} - i + 1
      \leq i \cdot 2^{\cardinalityOf{\sphere(i - 1)}}.
    \end{equation*}
    Hence, because $i \geq 2$,
    \begin{align*}
      \log\cardinalityOf{\pi_{\ball(i)}(X)}
      &\leq \log(i) + \cardinalityOf{\sphere(i - 1)} \cdot \log(2)\\
      &=    \log(i) + (3^{i - 1} + 3^{i - 2}) \cdot \log(2).
    \end{align*}
    Therefore, because $\cardinalityOf{\ball(i)} \geq 2 \cdot 3^i - 3^i = 3^i$,
    \begin{align*}
      \frac{\log\cardinalityOf{\pi_{\ball(i)}(X)}}{\cardinalityOf{\ball(i)}}
      &\leq \frac{\log(i) + (3^{i - 1} + 3^{i - 2}) \cdot \log(2)}{3^i}\\
      &=    \frac{\log(i)}{3^i} + \frac{4}{9} \cdot \log(2).
    \end{align*}
    Thus, $\entropyOf_{\mathcal{F}}(X) \leq (4/9) \cdot \log(2)$.

    Let $i$ be an integer such that $i \geq 1$. Then,
    \begin{align*}
      \cardinalityOf{\pi_{\ball(i)}(Y)}
      &= \cardinalityOf{\bigDisjointUnionOf_{\rho \in \N_0} \pi_{\ball(i)}(Y_\rho) \disjointUnionWith \setOf{\functionThatIsIdenticalToZero}}\\
      &= \sum_{\rho = 0}^i (2^{\cardinalityOf{\sphere(\rho)}} - 1) + 1\\
      &= \sum_{\rho = 0}^i 2^{\cardinalityOf{\sphere(\rho)}} - i.
    \end{align*}
    Thus, because $\sum_{\rho = 0}^{i - 1} 2^{\cardinalityOf{\sphere(\rho)}} - i \geq 0$, 
    \begin{equation*}
      \log\cardinalityOf{\pi_{\ball(i)}(Y)}
      \geq \log 2^{\cardinalityOf{\sphere(i)}}
      =    \cardinalityOf{\sphere(i)} \cdot \log(2)
      =    (3^i + 3^{i - 1}) \cdot \log(2).
    \end{equation*}
    Hence, because $\cardinalityOf{\ball(i)} \leq 2 \cdot 3^i$,
    \begin{equation*}
      \frac{\log\cardinalityOf{\pi_{\ball(i)}(\Delta(X))}}{\cardinalityOf{\ball(i)}}
      \geq \frac{3^i + 3^{i - 1}}{2 \cdot 3^i} \cdot \log(2)
      =    \frac{6}{9} \cdot \log(2).
    \end{equation*}
    Therefore, $\entropyOf_{\mathcal{F}}(Y) \geq (6/9) \cdot \log(2)$.

    In conclusion, 
    \begin{equation*}
      \entropyOf_{\mathcal{F}}(\Delta(X))
      =    \entropyOf_{\mathcal{F}}(Y)
      \geq \frac{6}{9} \cdot \log(2)
      >    \frac{4}{9} \cdot \log(2)
      \geq \entropyOf_{\mathcal{F}}(X). \qedhere
    \end{equation*}
  \end{proof}

  For a right-trac\-ta\-ble cell space, if for each point $t$ of an $\ntuple{E,E'}$-tiling not all patterns with domain $t \isSemiActedUponBy E$ occur in a subset of $Q^M$, then that subset does not have maximal entropy, which is shown in

  \begin{lemma} 
  \label{lemma:entropy-bounded-above-if-strange-tiling-exists}
    Let $\mathcal{R}$ be right tractable, let $\mathcal{F} = \net{F_i}_{i \in I}$ be a right Erling net in $\mathcal{R}$ indexed by $(I, \leq)$, let $Q$ contain at least two elements, let $X$ be a subset of $Q^M$, let $E$ and $E'$ be two non-empty and finite subsets of $G \modulo G_0$, and let $T$ be an $\ntuple{E, E'}$-tiling of $\mathcal{R}$, such that, for each cell $t \in T$, we have $\pi_{t \isSemiActedUponBy E}(X) \subsetneqq Q^{t \isSemiActedUponBy E}$. Then, $\entropyOf_{\mathcal{F}}(X) < \log\cardinalityOf{Q}$.
  \end{lemma} 

  \begin{proof}
    For each $t \in T$, because $\pi_{t \isSemiActedUponBy E}(X) \subsetneqq Q^{t \isSemiActedUponBy E}$, $\cardinalityOf{Q} \geq 2$, and $\cardinalityOf{t \isSemiActedUponBy E} \geq 1$,
    \begin{equation*}
      \cardinalityOf{\pi_{t \isSemiActedUponBy E}(X)}
      \leq \cardinalityOf{Q^{t \isSemiActedUponBy E}} - 1
      =    \cardinalityOf{Q}^{\cardinalityOf{t \isSemiActedUponBy E}} - 1
      \geq 1.
    \end{equation*}
    Let $i \in I$. Put $T_i = T \cap F_i^{-E}$ and put $F_i^* = F_i \smallsetminus (\bigcup_{t \in T_i} t \isSemiActedUponBy E)$ (see \cref{figure:entropy-bounded-above-if-special-tiling-exists}).
    \begin{figure}
      \myfloatalign
      \figureEntropyBoundedAboveIfSpecialTilingExists
      \caption{Schematic representation of the set-up of the proof of \cref{lemma:entropy-bounded-above-if-strange-tiling-exists}.}
      \label{figure:entropy-bounded-above-if-special-tiling-exists}
    \end{figure}
    Because $\bigcup_{t \in T_i} t \isSemiActedUponBy E \subseteq F_i$ and $\family{t \isSemiActedUponBy E}_{t \in T}$ is pairwise disjoint,
    \begin{equation*}
      \pi_{F_i}(X) \subseteq \pi_{F_i^*}(X) \times \prod_{t \in T_i} \pi_{t \isSemiActedUponBy E}(X)
                   \subseteq Q^{F_i^*} \times \prod_{t \in T_i} \pi_{t \isSemiActedUponBy E}(X).
    \end{equation*}
    Therefore,
    \begin{align*}
      \log\cardinalityOf{\pi_{F_i}(X)}
      &\leq \log\cardinalityOf{Q}^{\cardinalityOf{F_i^*}} + \sum_{t \in T_i} \log\cardinalityOf{\pi_{t \isSemiActedUponBy E}(X)}\\ 
      &\leq \log\cardinalityOf{Q}^{\cardinalityOf{F_i^*}} + \sum_{t \in T_i} \log\parens[\big]{\cardinalityOf{Q}^{\cardinalityOf{t \isSemiActedUponBy E}} - 1}\\ 
      &=    \cardinalityOf{F_i^*} \cdot \log\cardinalityOf{Q} + \sum_{t \in T_i} \log\parens[\big]{\cardinalityOf{Q}^{\cardinalityOf{t \isSemiActedUponBy E}} (1 - \cardinalityOf{Q}^{- \cardinalityOf{t \isSemiActedUponBy E}})}\\
      &=    \begin{aligned}[t]
              \cardinalityOf{F_i^*} \cdot \log\cardinalityOf{Q} &+ \sum_{t \in T_i} \cardinalityOf{t \isSemiActedUponBy E} \cdot \log\cardinalityOf{Q}\\
                                            &+ \sum_{t \in T_i} \log\parens[\big]{1 - \cardinalityOf{Q}^{- \cardinalityOf{t \isSemiActedUponBy E}}}.
            \end{aligned}
    \end{align*}
    Moreover, for each $t \in T_i$, we have $t \isSemiActedUponBy E \subseteq F_i$. Thus,
    \begin{equation*}
      \cardinalityOf{F_i^*} = \cardinalityOf{F_i} - \sum_{t \in T_i} \cardinalityOf{t \isSemiActedUponBy E}.
    \end{equation*}
    And, because $\isSemiActedUponBy$ is free, we have $\cardinalityOf{t \isSemiActedUponBy E} = \cardinalityOf{E}$. Hence,
    \begin{equation*}
      \log\cardinalityOf{\pi_{F_i}(X)}
      \leq    \cardinalityOf{F_i} \cdot \log\cardinalityOf{Q} + \cardinalityOf{T_i} \cdot \log\parens[\big]{1 - \cardinalityOf{Q}^{- \cardinalityOf{E}}}. 
    \end{equation*}
    Put $c = - \log\parens[\big]{1 - \cardinalityOf{Q}^{- \cardinalityOf{E}}}$. Because $\cardinalityOf{Q} \geq 2$ and $\cardinalityOf{E} \geq 1$, we have $\cardinalityOf{Q}^{- \cardinalityOf{E}} \in \openInterval{0, 1}$ and hence $c > 0$. 
    According to \cref{lemma:upper-bound-of-tiling-cap-Folner-net}, there are $\varepsilon \in \R_{> 0}$ and $i_0 \in I$ such that, for each $i \in I$ with $i \geq i_0$, we have $\cardinalityOf{T_i} \geq \varepsilon \cardinalityOf{F_i}$. Therefore, for each such $i$, 
    \begin{equation*}
      \frac{\log\cardinalityOf{\pi_{F_i}(X)}}{\cardinalityOf{F_i}} \leq \log\cardinalityOf{Q} - c \varepsilon.
    \end{equation*}
    In conclusion,
    \begin{equation*}
      \entropyOf_{\mathcal{F}}(X)
      =    \limsup_{i \in I} \frac{\log\cardinalityOf{\pi_{F_i}(X)}}{\cardinalityOf{F_i}}
      \leq \log\cardinalityOf{Q} - c \varepsilon
      <    \log\cardinalityOf{Q}. \qedhere
    \end{equation*}
  \end{proof}

  For a right-trac\-ta\-ble cell space, if for a non-empty and finite subset $E$ of $G \modulo G_0$ not all patterns with domain $m_0 \isSemiActedUponBy E$ occur in a \emph{shift-invariant} subset of $Q^M$, then that subset does not have maximal entropy, which is shown in

  \begin{corollary} 
  \label{corollary:nice-properties-yield-less-entropy} 
    Let $\mathcal{R}$ be right tractable, let $\mathcal{F} = \net{F_i}_{i \in I}$ be a right Erling net in $\mathcal{R}$ indexed by $(I, \leq)$, let $Q$ contain at least two elements, let $H$ be a $\mathcal{K}$-big subgroup of $G$, let $X$ be a $\actsOnMap_H$-invariant subset of $Q^M$, and let $E$ be a non-empty and finite subset of $G \modulo G_0$, such that $\pi_{m_0 \isSemiActedUponBy E}(X) \subsetneqq Q^{m_0 \isSemiActedUponBy E}$. Then, $\entropyOf_{\mathcal{F}}(X) < \log\cardinalityOf{Q}$.
  \end{corollary}

  \begin{proof} 
    According to \cref{theorem:existence-of-tiling}, there is a subset $E'$ of $G \modulo G_0$ and an $\ntuple{E, E'}$-tiling $T$ of $\mathcal{R}$. Because $G_0$ and $E$ are finite, so is $E'$. Let $m \in M$.
    Put $h = g_{m_0, m_0} g_{m_0, m}^{-1}$. Then, because $H$ is $\mathcal{K}$-big, we have $h \in H$. And, $h \actsOnPoint (m \isSemiActedUponBy E) = m_0 \isSemiActedUponBy E$.
    Hence, because $X$ is $\actsOnMap_H$-invariant, 
    \begin{align*} 
      \pi_{m \isSemiActedUponBy E}(X)
      &= \pi_{m \isSemiActedUponBy E}(h^{-1} \actsOnMap X)\\
      &= h^{-1} \actsOnMap \pi_{h \actsOnPoint (m \isSemiActedUponBy E)}(X)\\
      &= h^{-1} \actsOnMap \pi_{m_0 \isSemiActedUponBy E}(X).
    \end{align*}
    And, because $\pi_{m_0 \isSemiActedUponBy E}(X) \subsetneqq Q^{m_0 \isSemiActedUponBy E}$, 
    \begin{equation*}
      h^{-1} \actsOnMap \pi_{m_0 \isSemiActedUponBy E}(X)
      \subsetneqq h^{-1} \actsOnMap Q^{m_0 \isSemiActedUponBy E}
      =           Q^{h^{-1} \actsOnPoint (m_0 \isSemiActedUponBy E)}
      =           Q^{m \isSemiActedUponBy E}.
    \end{equation*}
    Therefore, $\pi_{m \isSemiActedUponBy E}(X) \subsetneqq Q^{m \isSemiActedUponBy E}$. In conclusion, according to \cref{lemma:entropy-bounded-above-if-strange-tiling-exists}, we have $\entropyOf_{\mathcal{F}}(X) < \log\cardinalityOf{Q}$.
  \end{proof}

  \begin{example}
  \label{example:entropy-bounded-above-if-strange-tiling-exists}
    This example demonstrates that in \cref{lemma:entropy-bounded-above-if-strange-tiling-exists} it is necessary that $\mathcal{F}$, $E$, and $T$ are such that the limit superior of the net $\net{\cardinalityOf{T \cap F_i^{-E}} / \cardinalityOf{F_i}}_{i \in I}$ is greater than $0$, which, according to \cref{lemma:upper-bound-of-tiling-cap-Folner-net}, is the case if $\mathcal{F}$ is a right Erling net in $\mathcal{R}$. 

    Let $\mathcal{F} = \net{F_i}_{i \in I}$ be a net in $\setOf{F \subseteq M \suchThat F \neq \emptyset, F \text{ finite}}$ indexed by $(I, \leq)$, let $Q$ contain at least two elements, let $E$ and $E'$ be two non-empty and finite subsets of $G \modulo G_0$, and let $T$ be an $\ntuple{E, E'}$-tiling of $\mathcal{R}$, Moreover, let $q$ be an element of $Q$, let $p$ be the pattern of $Q^{m_0 \isSemiActedUponBy E}$ such that $p \equiv q$, and let $X$ be the subset $Q^M \smallsetminus \bigcup_{t \in T} \setOf{c \in Q^M \suchThat c\restrictedTo_{t \isSemiActedUponBy E} = t \actsByItsCoordinateOn p}$ of $Q^M$. Then, for each cell $t \in T$, we have $\pi_{t \isSemiActedUponBy E}(X) = Q^{t \isSemiActedUponBy E} \smallsetminus \setOf{t \actsByItsCoordinateOn p} \subsetneqq Q^{t \isSemiActedUponBy E}$. 

    Let $i$ be an index of $I$. Then, $\pi_{F_i}(X) = Q^{F_i} \smallsetminus \bigcup_{t \in T \cap F_i^{-E}} Y_t$, where $Y_t = \setOf{p' \in Q^{F_i} \suchThat p'\restrictedTo_{t \isSemiActedUponBy E} = t \actsByItsCoordinateOn p}$, for $t \in T \cap F_i^{-E}$. Indeed, we have $\pi_{F_i}(X) \subseteq Q^{F_i} \smallsetminus \bigcup_{t \in T \cap F_i^{-E}} Y_t$. To show the other inclusion, let $p' \in Q^{F_i} \smallsetminus \bigcup_{t \in T \cap F_i^{-E}} Y_t$. Then, there is a state $q' \in Q \smallsetminus \setOf{q}$ and there is a global configuration $c \in Q^M$ such that $c\restrictedTo_{F_i} = p'$ and $c\restrictedTo_{M \smallsetminus F_i} \equiv q'$. Hence, for each $t \in T \cap F_i^{-E}$, we have $c\restrictedTo_{t \isSemiActedUponBy E} = p'\restrictedTo_{t \isSemiActedUponBy E} \neq t \actsByItsCoordinateOn p$. And, for each $t \in T \smallsetminus F_i^{-E}$, there is an $e \in E$ such that $t \isSemiActedUponBy e \notin F_i$ and thus, by definition of $c$, we have $c(t \isSemiActedUponBy e) = q' \neq q$, and hence $c\restrictedTo_{t \isSemiActedUponBy E} \neq t \actsByItsCoordinateOn p$. Therefore, $c \in X$ and hence $p' = c\restrictedTo_{F_i} \in \pi_{F_i}(X)$. In conclusion, $Q^{F_i} \smallsetminus \bigcup_{t \in T \cap F_i^{-E}} Y_t \subseteq \pi_{F_i}(X)$. 

    For each subset $S$ of $T \cap F_i^{-E}$, we have $\cardinalityOf{\bigcap_{s \in S} Y_s} = \cardinalityOf{Q^{F_i \smallsetminus \bigDisjointUnionOf_{s \in S} s \isSemiActedUponBy E}} = \cardinalityOf{Q}^{\cardinalityOf{F_i} - \cardinalityOf{S} \cdot \cardinalityOf{E}}$, which only depends on the cardinality of $S$. Hence, according to a special case of the inclusion-exclusion principle and the binomial formula, 
    \begin{align*}
      \cardinalityOf{\pi_{F_i}(X)}
      &= \sum_{k = 0}^{\cardinalityOf{T \cap F_i^{-E}}} (-1)^k \binom{\cardinalityOf{T \cap F_i^{-E}}}{k} \cardinalityOf{Q}^{\cardinalityOf{F_i} - k \cardinalityOf{E}}\\
      &= \cardinalityOf{Q}^{\cardinalityOf{F_i}} \cdot \sum_{k = 0}^{\cardinalityOf{T \cap F_i^{-E}}} \binom{\cardinalityOf{T \cap F_i^{-E}}}{k} 1^{\cardinalityOf{T \cap F_i^{-E}} - k} (-\cardinalityOf{Q}^{-\cardinalityOf{E}})^k\\
      &= \cardinalityOf{Q}^{\cardinalityOf{F_i}} \cdot (1 - \cardinalityOf{Q}^{-\cardinalityOf{E}})^{\cardinalityOf{T \cap F_i^{-E}}}.
    \end{align*}
    Thus, 
    \begin{equation*}
      \frac{\log\cardinalityOf{\pi_{F_i}(X)}}{\cardinalityOf{F_i}}
      = \log\cardinalityOf{Q} + \frac{\cardinalityOf{T \cap F_i^{-E}}}{\cardinalityOf{F_i}} \cdot \log(1 - \cardinalityOf{Q}^{-\cardinalityOf{E}}).
    \end{equation*}
    Therefore,
    \begin{equation*}
      \entropyOf_{\mathcal{F}}(X) = \log\cardinalityOf{Q} + \limsup_{i \in I} \frac{\cardinalityOf{T \cap F_i^{-E}}}{\cardinalityOf{F_i}} \cdot \log(1 - \cardinalityOf{Q}^{-\cardinalityOf{E}}).
    \end{equation*}
    Hence, because $\log(1 - \cardinalityOf{Q}^{-\cardinalityOf{E}}) < 0$, we have $\entropyOf_{\mathcal{F}}(X) < \log\cardinalityOf{Q}$ if and only if $\limsup_{i \in I} \cardinalityOf{T \cap F_i^{-E}} / \cardinalityOf{F_i} > 0$.
  \end{example}

  \section{Gardens of Eden} 
  \label{section:gardens-of-Eden} 

  In this section, let $\mathcal{R} = \ntuple{\mathcal{M}, \mathcal{K}} = \ntuple{\ntuple{M, G, \actsOnPoint}, \ntuple{m_0, \family{g_{m_0, m}}_{m \in M}}}$ be a cell space and let $\mathcal{C} = \ntuple{\mathcal{R}, Q, N, \delta}$ be a semi-cellular automaton such that the stabiliser $G_0$ of $m_0$ under $\actsOnPoint$, the set $Q$ of states, and the neighbourhood $N$ are finite, and the set $Q$ is non-empty. Furthermore, let $\Delta$ be the global transition function of $\mathcal{C}$.

  \paragraph{Contents.} In \cref{theorem:not-surjective-implies-less-entropy} we show that if $\Delta$ is not surjective, then the entropy of its image is less than the entropy of $Q^M$. And the converse of that statement obviously holds. In \cref{theorem:less-entropy-implies-not-pre-injective} we show that if the entropy of the image of $\Delta$ is less than the entropy of $Q^M$, then $\Delta$ is not pre-injective. And in \cref{theorem:not-pre-injective-implies-less-entropy} we show the converse of that statement. These four statements establish the Garden of Eden theorem, which is \cref{theorem:garden-of-Eden}. 

  \paragraph{Body.} We first introduce the difference set of two global configurations and then use it to define pre-injectivity.

  \begin{definition} 
  \label{definition:difference}
    Let $c$ and $c'$ be two global configurations of $Q^M$. The set
    \begin{equation*}
      \differenceOf(c, c') = \setOf{m \in M \suchThat c(m) \neq c'(m)}
      \mathnote{difference $\differenceOf(c, c')$ of $c$ and $c'$}
    \end{equation*}
    is called \define{difference of $c$ and $c'$}.
  \end{definition}

  \begin{lemma}
  \label{lemma:difference-of-configurations-and-action}
    Let $c$ and $c'$ be two global configurations of $Q^M$. Then,
    \begin{equation*}
      \ForEach g \in G \Holds \differenceOf(g \actsOnMap c, g \actsOnMap c') = g \actsOnPoint \differenceOf(c, c'). \qedhere
    \end{equation*}
  \end{lemma}

  \begin{proof}
    Let $g \in G$. Then, for each $m \in M$,
    \begin{align*}
      m \in \differenceOf(g \actsOnMap c, \actsOnMap c')
      &\ifAndOnlyIf (g \actsOnMap c)(m) \neq (g \actsOnMap c')(m)\\
      &\ifAndOnlyIf c(g^{-1} \actsOnPoint m) \neq c'(g^{-1} \actsOnPoint m)\\
      &\ifAndOnlyIf g^{-1} \actsOnPoint m \in \differenceOf(c, c')\\
      &\ifAndOnlyIf m \in g \actsOnPoint \differenceOf(c, c').
    \end{align*}
    In conclusion, $\differenceOf(g \actsOnMap c, \actsOnMap c') = g \actsOnPoint \differenceOf(c, c')$.
  \end{proof}

  \begin{definition} 
  \label{definition:pre-injective}
    The map $\Delta$ is called \define{pre-injective}\graffito{pre-injective} if and only if, for each tuple $(c, c') \in Q^M \times Q^M$ such that $\differenceOf(c, c')$ is finite and $\Delta(c) = \Delta(c')$, we have $c = c'$. 
  \end{definition}

  \begin{remark}
    Each injective map is pre-injective. And, if $M$ is finite, then each pre-injective map is injective.
  \end{remark}

  In the proof of \cref{theorem:not-surjective-implies-less-entropy}, the existence of a Garden of Eden pattern, as stated in \cref{lemma:not-surjective-yields-garden-of-Eden-pattern}, is essential, which itself follows from the existence of a Garden of Eden configuration, the compactness of $Q^M$, and the continuity of $\Delta$. Garden of Eden configurations and patterns are introduced in

  \begin{definition}
    \begin{aenumerate}
      \item Let $c \from M \to Q$ be a global configuration. It is called \define{Garden of Eden configuration}\graffito{Garden of Eden configuration $c$ of $\mathcal{C}$} if and only if it is not contained in $\Delta(Q^M)$. 
      \item Let $p \from A \to Q$ be a pattern. It is called \define{Garden of Eden pattern}\graffito{Garden of Eden pattern $p$ of $\mathcal{C}$} if and only if, for each global configuration $c \in Q^M$, we have $\Delta(c)\restrictedTo_A \neq p$. \qedhere 
    \end{aenumerate}
  \end{definition}

  \begin{remark}
    \begin{aenumerate}
      \item The global transition function $\Delta$ is surjective if and only if there is no Garden of Eden configuration. 
      \item If $p \from A \to Q$ is a Garden of Eden pattern, then each global configuration $c \in Q^M$ with $c\restrictedTo_A = p$ is a Garden of Eden configuration. 
      \item If there is a Garden of Eden pattern, then $\Delta$ is not surjective. \qedhere 
    \end{aenumerate}
  \end{remark}

  \begin{lemma} 
  \label{lemma:not-surjective-yields-garden-of-Eden-pattern}
    Let $\Delta$ not be surjective. There is a Garden of Eden pattern with non-empty and finite domain. 
  \end{lemma} 

  \begin{proof}
    Because $\Delta$ is not surjective, there is a Garden of Eden configuration $c \in Q^M$. Equip $Q^M$ with the prodiscrete topology. According to \cref{lemma:image-of-phase-space-is-closed}, the image $\Delta(Q^M)$ is closed in $Q^M$. Hence, its complement $Q^M \smallsetminus \Delta(Q^M)$ is open. Therefore, because $c \in Q^M \smallsetminus \Delta(Q^M)$, according to \cref{remark:prodiscrete-topology}, there is a non-empty and finite subset $F$ of $M$ such that 
    \begin{equation*}
      \cylinder(c, F) = \setOf{c' \in Q^M \suchThat c'\restrictedTo_F = c\restrictedTo_F} \subseteq Q^M \smallsetminus \Delta(Q^M).
    \end{equation*}
    Hence, $c\restrictedTo_F$ is a Garden of Eden pattern with non-empty and finite domain.
  \end{proof}

  Under which assumptions the entropy of the image of a non-surjective global transition function is not maximal is shown in

  \begin{theorem} 
  \label{theorem:not-surjective-implies-less-entropy} 
    Let $\mathcal{R}$ be right tractable, let $\mathcal{F}$ be a right Erling net in $\mathcal{R}$, let $H$ be a $\mathcal{K}$-big subgroup of $G$, let $\delta$ be $\bullet_{H_0}$-invariant, let $Q$ contain at least two elements, and let $\Delta$ not be surjective. Then, $\entropyOf_{\mathcal{F}}(\Delta(Q^M)) < \log\cardinalityOf{Q}$.
  \end{theorem}

  \begin{proof} 
    According to \cref{lemma:not-surjective-yields-garden-of-Eden-pattern}, there is a Garden of Eden pattern $p \from F \to Q$ with non-empty and finite domain. Let $E = (m_0 \isSemiActedUponBy \blank)^{-1}(F)$. Then, $m_0 \isSemiActedUponBy E = F$ and, because $\isSemiActedUponBy$ is free, $\cardinalityOf{E} = \cardinalityOf{F} < \infty$. Because $p$ is a Garden of Eden pattern, $p \notin \pi_{m_0 \isSemiActedUponBy E}(\Delta(Q^M))$. Hence, $\pi_{m_0 \isSemiActedUponBy E}(\Delta(Q^M)) \subsetneqq Q^{m_0 \isSemiActedUponBy E}$. Moreover, according to \cref{theorem:local-invariance-versus-global-equivariance}, the map $\Delta$ is $\actsOnMap_H$-e\-qui\-var\-i\-ant. Hence, for each $h \in H$, we have $h \actsOnMap \Delta(Q^M) = \Delta(h \actsOnMap Q^M) = \Delta(Q^M)$. In other words, $\Delta(Q^M)$ is $\actsOnMap_H$-invariant. Thus, according to \cref{corollary:nice-properties-yield-less-entropy}, we have $\entropyOf_{\mathcal{F}}(\Delta(Q^M)) < \log\cardinalityOf{Q}$.
  \end{proof}

  This yields the characterisation of surjectivity by entropy that is given in

  \begin{corollary}
    Let $\mathcal{R}$ be right tractable, let $\mathcal{F}$ be a right Erling net in $\mathcal{R}$, let $H$ be a $\mathcal{K}$-big subgroup of $G$, let $\delta$ be $\bullet_{H_0}$-invariant, and let $Q$ contain at least two elements. Then, $\Delta$ is surjective if and only if $\entropyOf_{\mathcal{F}}(\Delta(Q^M)) = \log\cardinalityOf{Q}$.
  \end{corollary}

  \begin{proof}
    This is a direct consequence of \cref{theorem:not-surjective-implies-less-entropy}.
  \end{proof}

  In the remainder of this section, let $\mathcal{R}$ be right amenable and let $\mathcal{F} = \net{F_i}_{i \in I}$ be a right Følner net in $\mathcal{R}$ indexed by $(I, \leq)$. 

  In the proof of \cref{theorem:less-entropy-implies-not-pre-injective}, the fact that enlarging each element of $\mathcal{F}$ does not increase entropy, as stated in the next lemma, is essential.

  \begin{lemma}
  \label{lemma:entropy-invariant-under-closure-net-change}
    Let $X$ be a subset of $Q^M$ and let $E$ be a finite subset of $G \modulo G_0$ such that $G_0 \in E$. Then, $\entropyOf_{\net{F_i^{+E}}_{i \in I}}(X) \leq \entropyOf_{\mathcal{F}}(X)$. 
  \end{lemma}


  \begin{proof}
    Let $i \in I$. According to \cref{item:properties-of-interior-closure-and-boundary:neutral-element} of \cref{lemma:properties-of-interior-closure-and-boundary}, we have $F_i^{-E} \subseteq F_i \subseteq F_i^{+E}$. Hence, $\pi_{F_i^{+E}}(X) \subseteq \pi_{F_i}(X) \times Q^{F_i^{+E} \smallsetminus F_i}$ and $F_i^{+E} \smallsetminus F_i \subseteq \boundaryOf_E F_i$. Thus,
    \begin{align*}
      \log\cardinalityOf{\pi_{F_i^{+E}}(X)}
      &\leq \log\cardinalityOf{\pi_{F_i}(X)} + \cardinalityOf{F_i^{+E} \smallsetminus F_i} \cdot \log\cardinalityOf{Q}\\ 
      &\leq \log\cardinalityOf{\pi_{F_i}(X)} + \cardinalityOf{\boundaryOf_E F_i} \cdot \log\cardinalityOf{Q}.
    \end{align*}
    Therefore, according to \cref{theorem:boundary-characterisation-of-Folner-net},
    \begin{align*}
      \entropyOf_{\net{F_i^{+E}}_{i \in I}}(X)
      &\leq \limsup_{i \in I} \frac{\log\cardinalityOf{\pi_{F_i}(X)}}{\cardinalityOf{F_i}} + \parens*{\lim_{i \in I} \frac{\cardinalityOf{\boundaryOf_E F_i}}{\cardinalityOf{F_i}}} \cdot \log\cardinalityOf{Q}\\
      &=    \entropyOf_{\mathcal{F}}(X). \qedhere
    \end{align*}
  \end{proof}

  A global transition function whose image does not have maximal entropy is not pre-injective, which is shown in

  \begin{theorem} 
  \label{theorem:less-entropy-implies-not-pre-injective}
    Let $\entropyOf_{\mathcal{F}}(\Delta(Q^M)) < \log\cardinalityOf{Q}$. Then, $\Delta$ is not pre-injective.
  \end{theorem}

  \begin{proof-sketch}
    The asymptotic growth rate of finite patterns in $\Delta(Q^M)$ is less than the one of $Q^M$. Hence, there are at least two finite patterns that can be identically extended to global configurations that have the same image under $\Delta$. Therefore, $\Delta$ is not pre-injective. 
  \end{proof-sketch}

  \begin{proof}
    Suppose, without loss of generality, that $G_0 \in N$. Let $X = \Delta(Q^M)$. According to \cref{lemma:entropy-invariant-under-closure-net-change}, we have $\entropyOf_{\net{F_i^{+N}}_{i \in I}}(X) \leq \entropyOf_{\mathcal{F}}(X) < \log\cardinalityOf{Q}$. Hence, there is an $i \in I$ such that
    \begin{equation*}
      \frac{\log\cardinalityOf{\pi_{F_i^{+N}}(X)}}{\cardinalityOf{F_i}} < \log\cardinalityOf{Q}.
    \end{equation*}
    Thus, $\cardinalityOf{\pi_{F_i^{+N}}(X)} < \cardinalityOf{Q}^{\cardinalityOf{F_i}}$. Furthermore, let $q \in Q$ and let $X' = \setOf{c \in Q^M \suchThat c\restrictedTo_{M \smallsetminus F_i} \equiv q}$. Then, $\cardinalityOf{Q}^{\cardinalityOf{F_i}} = \cardinalityOf{X'}$. Hence,
      $\cardinalityOf{\pi_{F_i^{+N}}(X)} < \cardinalityOf{X'}$.
    Moreover, for each $(c, c') \in X' \times X'$, according to \cref{item:global-transition-function-and-interior-closure:closure} of \cref{lemma:global-transition-function-and-interior-closure}, we have $\Delta(c)\restrictedTo_{M \smallsetminus F_i^{+N}} = \Delta(c')\restrictedTo_{M \smallsetminus F_i^{+N}}$. Therefore, 
    \begin{align*}
      \cardinalityOf{\Delta(X')}
      &=    \cardinalityOf{\pi_{F_i^{+N}}(\Delta(X'))}\\
      &\leq \cardinalityOf{\pi_{F_i^{+N}}(\Delta(Q^M))}\\
      &=    \cardinalityOf{\pi_{F_i^{+N}}(X)}\\
      &<    \cardinalityOf{X'}.
    \end{align*}
    Hence, there are $c$, $c' \in X'$ such that $c \neq c'$ and $\Delta(c) = \Delta(c')$. Thus, because $\differenceOf(c, c') \subseteq F_i$ is finite, the map $\Delta$ is not pre-injective.
  \end{proof}

  \begin{counterexample}[Muller] 
  \label{example:non-maximal-entropy-but-pre-injective} 
    In this example we present a cellular automaton on a non-right-a\-me\-na\-ble but right-trac\-ta\-ble cell space that, although the image of its global transition function does not have maximal entropy with respect to a right Erling net, is pre-injective. It is Muller's counterexample to Myhill's theorem, see section~6, page~55, in \cite{machi:mignosi:1993}.

    Let $G$ be the group with presentation $\groupGeneratedBy{x, y, z \suchThat x^2, y^2, z^2}$ or, equivalently, the $3$-fold free product of the cyclic group of order $2$ with itself, let $Q$ be the $\F_2$-vector space $(\F_2)^2$, where $\F_2$ is the finite field of order $2$, let $N$ be the set $\setOf{x, y, z}$, let $\delta$ be the map $Q^N \to Q$, $\ell \mapsto (\ell(x)_1 + \ell(x)_2 + \ell(y)_1 + \ell(z)_2, 0)$, where, for each vector $v \in Q$, the first component of $v$ is denoted by $v_1$ and the second by $v_2$. The tuple $\mathcal{R} = \ntuple{\ntuple{G, G, \cdot}, \ntuple{e_G, \family{g}_{g \in G}}}$ is a cell space and the quadruple $\mathcal{C} = \ntuple{\mathcal{R}, Q, N, \delta}$ is a cellular automaton. 

    According to \cref{lemma:finitely-right-generated-is-tractable}, the cell space $\mathcal{R}$ is right tractable and the sequence $\mathcal{F} = \sequence{\ball(\rho)}_{\rho \in \N_0}$ is a right Erling net. Moreover, because point evaluation, projection, and addition are linear, the local transition function $\delta$ is linear and hence the global transition function $\Delta$ of $\mathcal{C}$ is linear. Furthermore, because the image of $\Delta$ is included in $(\F_2 \times \setOf{0})^G$, the global transition function $\Delta$ is not surjective. Hence, according to \cref{theorem:not-surjective-implies-less-entropy}, we have $\entropyOf_{\mathcal{F}}(\Delta(Q^G)) < \log\cardinalityOf{Q}$. However, as we show now, the global transition function $\Delta$ is pre-injective. 

    Let $c$ and $c'$ be two global configurations of $Q^G$ such that $\differenceOf(c, c')$ is finite and $\Delta(c) = \Delta(c')$. Then, because $\Delta$ is linear, we have $\Delta(c - c') \equiv 0$. Let $c''$ be the global configuration $c - c'$. Suppose that $c'' \not\equiv 0$. Then, because $c''$ has finite support, there is an element $g \in G$ such that $c''(g) \neq 0$ and $c''\restrictedTo_{M \smallsetminus \ball(\lengthOf{g})} \equiv 0$. And, there are two distinct elements $s$ and $s' \in \setOf{x, y, z}$ such that $g s$, $g s' \in \sphere(\lengthOf{g} + 1)$ (see \cref{figure:non-maximal-entropy-but-pre-injective}). Hence, because $c''\restrictedTo_{M \smallsetminus \ball(\lengthOf{g})} \equiv 0$, we have $c''(g s) = 0$ and $c''(g s') = 0$. And, for each element $s'' \in \setOf{s, s'}$, because $g s'' s'' = g$ and $g s'' s''' \in \sphere(\lengthOf{g} + 2)$, for $s''' \in \setOf{x, y, z} \smallsetminus \setOf{s''}$, 
    \begin{align*}
      \Delta(c'')(g s'')_1
      &= c''(g s'' x)_1 + c''(g s'' x)_2 + c''(g s'' y)_1 + c''(g s'' z)_2\\
      &= \begin{dcases*}
           c''(g)_1 + c''(g)_2, &if $s'' = x$,\\
           c''(g)_1, &if $s'' = y$,\\
           c''(g)_2, &if $s'' = z$.
         \end{dcases*}
    \end{align*}
    Because $c''(g) \neq 0$, for each element $s'' \in \setOf{s, s'}$, in two of the three cases we have $\Delta(c'')(g s'')_1 = 1$ and hence $\Delta(c'')(g s)_1 = 1$ or $\Delta(c'')(g s')_1 = 1$. Therefore, $\Delta(c'')(g s) \neq 0$ or $\Delta(c'')(g s') \neq 0$, which contradicts that $\Delta(c'') \equiv 0$. Thus, contrary to our supposition, we have $c'' \equiv 0$ and hence $c = c'$. In conclusion, the global transition function $\Delta$ is pre-injective.
    \begin{figure}
      \myfloatalign
      \figureNonMaximalEntropyButPreInjective
      \caption{A part of the $\setOf{x, y, z}$-Cayley graph of $G$.}
      \label{figure:non-maximal-entropy-but-pre-injective}
    \end{figure}
  \end{counterexample}

  %

  In the proof of \cref{theorem:not-pre-injective-implies-less-entropy}, the statement of \cref{lemma:exchanging-pattern-by-other-pattern-with-same-image-yields-configuration-with-same-image} is essential, which says that if two distinct patterns have the same image and we replace each occurrence of the first by the second in a global configuration, we get a new configuration in which the first pattern does not occur and that has the same image as the original one.

  How maps with disjoint domains are glued together, which is used to replace occurrences of patterns in global configurations by other patterns, is introduced in

  \begin{definition} 
  \label{definition:coproduct-of-maps}
            Let $I$ be a set and, for each index $i \in I$, let $A_i$ and $B_i$ be two sets and let $f_i$ be a map from $A_i$ to $B_i$, such that the sets $A_i$, for $i \in I$, are pairwise disjoint. The map\index[symbols]{coproductfiiinI@$\coprod_{i \in I} f_i$}
            \begin{align*}
              \coprod_{i \in I} f_i \from \bigDisjointUnionOf_{i \in I} A_i &\to \bigDisjointUnionOf_{i \in I} B_i, \mathnote{coproduct $\coprod_{i \in I} f_i$ of $\family{f_i}_{i \in I}$}\\
              x &\mapsto f_i(x), \text{ where } i \in I \text{ such that } x \in A_i,
            \end{align*}
            is called \define{coproduct of $\family{f_i}_{i \in I}$} and, if $I$ is the set $\setOf{1, 2, \dotsc, \cardinalityOf{I}}$, then it is also denoted by $f_1 \times f_2 \times \dotsb \times f_{\cardinalityOf{I}}$\graffito{$f_1 \times f_2 \times \dotsb \times f_{\cardinalityOf{I}}$}\index[symbols]{crossf1f2@$f_1 \times f_2 \times \dotsb \times f_{\cardinalityOf{I}}$}. 
  \end{definition}

  In the proof of \cref{lemma:exchanging-pattern-by-other-pattern-with-same-image-yields-configuration-with-same-image} we use the technical

  \begin{lemma}
  \label{lemma:m-liberation-A-is-either-outside-of-A-or-inside-of-closure-of-A}
    Let $A$ be a subset of $M$, and let $E$ and $E'$ be two subsets of $G \modulo G_0$ such that $\setOf{g^{-1} \cdot e' \suchThat e, e' \in E, g \in e} \subseteq E'$. Then, 
    \begin{equation*}
      \ForEach m \in M \Holds m \isSemiActedUponBy E \subseteq M \smallsetminus A \text{ or } m \isSemiActedUponBy E \subseteq A^{+E'}. \qedhere 
    \end{equation*}
  \end{lemma}

  \begin{proof}
    \begin{figure}
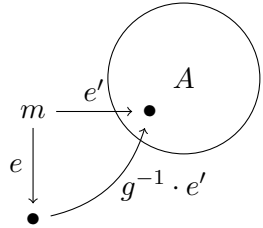

      \myfloatalign
      \figureMLiberationAIsEitherOutsideOfAOrInsideOfClosureOfA
      \caption{Schematic representation of the proof of \cref{lemma:m-liberation-A-is-either-outside-of-A-or-inside-of-closure-of-A}.}
      \label{figure:m-liberation-A-is-either-outside-of-A-or-inside-of-closure-of-A}
    \end{figure}
    Let $m \in M$ such that $m \isSemiActedUponBy E \nsubseteq M \smallsetminus A$. Then, $(m \isSemiActedUponBy E) \cap A \neq \emptyset$. Hence, there is an $e' \in E$ such that $m \isSemiActedUponBy e' \in A$. Let $e \in E$. According to \cref{lemma:rightsemiaction-can-be-undone}, there is a $g \in e$ such that $(m \isSemiActedUponBy e) \isSemiActedUponBy g^{-1} \cdot e' = m \isSemiActedUponBy e'$. Because $g^{-1} \cdot e' \in E'$ and $m \isSemiActedUponBy e' \in A$, we have $(m \isSemiActedUponBy e) \isSemiActedUponBy E' \cap A \neq \emptyset$ (see \cref{figure:m-liberation-A-is-either-outside-of-A-or-inside-of-closure-of-A}). Thus, $m \isSemiActedUponBy e \in A^{+E'}$. Therefore, $m \isSemiActedUponBy E \subseteq A^{+E'}$.  
  \end{proof}

  \begin{lemma} 
  \label{lemma:exchanging-pattern-by-other-pattern-with-same-image-yields-configuration-with-same-image} 
    Identify $M$ with $G \modulo G_0$ by $\iota \givenBy m \mapsto G_{m_0, m}$, let $H$ be a $\mathcal{K}$-big subgroup of $G$, let $\delta$ be $\bullet_{H_0}$-invariant, let $A$ be a subset of $M$, let $N'$ be the subset $\setOf{g^{-1} \cdot n' \suchThat n, n' \in N, g \in n}$ of $G \modulo G_0$, and let $p$ and $p'$ be two maps from $A^{+N'}$ to $Q$ such that $p\restrictedTo_{A^{+N'} \smallsetminus A} = p'\restrictedTo_{A^{+N'} \smallsetminus A}$ and $\Delta_{A^{+N'}}^-(p) = \Delta_{A^{+N'}}^-(p')$. Furthermore, let $c$ be a map from $M$ to $Q$ and let $S$ be a subset of $M$, such that the family $\family{s \isSemiActedUponBy A^{+N'}}_{s \in S}$ is pairwise disjoint and, for each cell $s \in S$, we have $p \occursIn_s c$. Put
    \begin{equation*} 
      c' = c\restrictedTo_{M \smallsetminus (\bigcup_{s \in S} s \isSemiActedUponBy A^{+N'})} \times \coprod_{s \in S} s \actsByItsCoordinateOn p'.
    \end{equation*}
    Then, for each cell $s \in S$, we have $p' \occursIn_s c'$, and $\Delta(c) = \Delta(c')$. In particular, if $p \neq p'$, then, for each cell $s \in S$, we have $p \not\occursIn_s c'$.
  \end{lemma} 

  \begin{proof} 
    For each $s \in S$, we have $\domainOf(s \actsByItsCoordinateOn p) = \domainOf(s \actsByItsCoordinateOn p') = s \isSemiActedUponBy A^{+N'}$. Hence, $c'$ is well-defined. Moreover, for each $s \in S$, we have $(s \actsByItsCoordinateOn p)\restrictedTo_{(s \isSemiActedUponBy A^{+N'}) \smallsetminus (s \isSemiActedUponBy A)} = (s \actsByItsCoordinateOn p')\restrictedTo_{(s \isSemiActedUponBy A^{+N'}) \smallsetminus (s \isSemiActedUponBy A)}$.

    Let $m \in M \smallsetminus (\bigcup_{s \in S} s \isSemiActedUponBy A)$. If $m \in M \smallsetminus (\bigcup_{s \in S} s \isSemiActedUponBy A^{+N'})$, then $c'(m) = c(m)$. And, if there is an $s \in S$ such that $m \in s \isSemiActedUponBy A^{+N'}$, then, because $m \notin s \isSemiActedUponBy A$, we have $c'(m) = (s \actsByItsCoordinateOn p')(m) = (s \actsByItsCoordinateOn p)(m) = c(m)$. Therefore,
    \begin{equation*}
      c' = c\restrictedTo_{M \smallsetminus (\bigcup_{s \in S} s \isSemiActedUponBy A)} \times \coprod_{s \in S} s \actsByItsCoordinateOn (p'\restrictedTo_A).
    \end{equation*}

    Let $m \in M$.
    \begin{description}
      \item[Case 1:] $m \isSemiActedUponBy N \subseteq M \smallsetminus (\bigcup_{s \in S} s \isSemiActedUponBy A)$. Then, $c'\restrictedTo_{m \isSemiActedUponBy N} = c\restrictedTo_{m \isSemiActedUponBy N}$. Hence, $\Delta(c')(m) = \Delta(c)(m)$.
      \item[Case 2:] $m \isSemiActedUponBy N \nsubseteq M \smallsetminus (\bigcup_{s \in S} s \isSemiActedUponBy A)$. Then, there is an $s \in S$ such that $m \isSemiActedUponBy N \nsubseteq M \smallsetminus (s \isSemiActedUponBy A)$. Thus, according to \cref{lemma:m-liberation-A-is-either-outside-of-A-or-inside-of-closure-of-A}, we have $m \isSemiActedUponBy N \subseteq (s \isSemiActedUponBy A)^{+N'}$. Hence, because $G_0 \cdot N' \subseteq N'$, according to \cref{item:properties-of-interior-closure-and-boundary:commute-with-liberation} of \cref{lemma:properties-of-interior-closure-and-boundary}, we have $m \isSemiActedUponBy N \subseteq s \isSemiActedUponBy A^{+N'}$ and hence $m \in (s \isSemiActedUponBy A^{+N'})^{-N}$. Therefore, because $c\restrictedTo_{s \isSemiActedUponBy A^{+N'}} = s \actsByItsCoordinateOn p$, $\Delta_{A^{+N'}}^-(p) = \Delta_{A^{+N'}}^-(p')$, and $c'\restrictedTo_{s \isSemiActedUponBy A^{+N'}} = s \actsByItsCoordinateOn p'$, according to \cref{corollary:equivariance-of-restriction}, 
            \begin{align*}
              \Delta(c)(m)
              &= \Delta_{s \isSemiActedUponBy A^{+N'}}^-(s \actsByItsCoordinateOn p)(m)\\
              &= (s \actsByItsCoordinateOn \Delta_{A^{+N'}}^-(p))(m)\\
              &= (s \actsByItsCoordinateOn \Delta_{A^{+N'}}^-(p'))(m)\\
              &= \Delta_{s \isSemiActedUponBy A^{+N'}}^-(s \actsByItsCoordinateOn p')(m)\\
              &= \Delta(c')(m).
            \end{align*}
    \end{description}
    In either case, $\Delta(c)(m) = \Delta(c')(m)$. Therefore, $\Delta(c) = \Delta(c')$.
  \end{proof}

  Under which assumptions the entropy of the image of a non-pre-injective global transition function is not maximal is shown in

  \begin{theorem} 
  \label{theorem:not-pre-injective-implies-less-entropy}
    Let $H$ be a $\mathcal{K}$-big subgroup of $G$, let $\delta$ be $\bullet_{H_0}$-invariant, let $Q$ contain at least two elements, and let $\Delta$ not be pre-injective. Then, $\entropyOf_{\mathcal{F}}(\Delta(Q^M)) < \log\cardinalityOf{Q}$.
  \end{theorem}

  \begin{proof-sketch} 
    For $N' = N^{-1} \cdot N$, there is a subset $A$ of $M$ and there are two distinct finite patterns $p$ and $p'$ with domain $A^{+N'}$ that have the same image under $\Delta_{A^{+N'}}^-$. The set $Y$ of all global configurations in which $p$ does not occur at the cells of a tiling has the same image under $\Delta$ as $Q^M$, because in a global configuration we may replace occurrences of $p$ by $p'$ without changing the image. It follows that $\entropyOf(\Delta(Q^M)) = \entropyOf(\Delta(Y)) \leq \entropyOf(Y)$; and, because $Y$ is missing the pattern $p$ at each cell of a tiling, we also have $\entropyOf(Y) < \entropyOf(Q^M) = \log\cardinalityOf{Q}$. In conclusion, $\entropyOf(\Delta(Q^M)) < \log\cardinalityOf{Q}$.
  \end{proof-sketch}

  \begin{proof} 
    Suppose, without loss of generality, that $G_0 \in N$. Identify $M$ with $G \modulo G_0$ by $\iota \givenBy m \mapsto G_{m_0, m}$.

    Because $\Delta$ is not pre-injective, there are $c$, $c' \in Q^M$ such that $\differenceOf(c, c')$ is finite, $\Delta(c) = \Delta(c')$, and $c \neq c'$. Put $A = \differenceOf(c, c')$, put $N' = \setOf{g^{-1} \cdot n' \suchThat n, n' \in N, g \in n}$, put $E = A^{+N'}$, and put $p = c\restrictedTo_E$ and $p' = c'\restrictedTo_E$. Because $\Delta(c) = \Delta(c')$, we have $\Delta_{A^{+N'}}^-(p) = \Delta_{A^{+N'}}^-(p')$. 

    Because $N$ is finite and, for each $n \in N$, we have $\cardinalityOf{n} = \cardinalityOf{G_0} < \infty$, the set $N'$ is finite. Moreover, $G_0 \cdot N' \subseteq N'$. According to \cref{item:properties-of-interior-closure-and-boundary:neutral-element} of \cref{lemma:properties-of-interior-closure-and-boundary}, because $G_0 \in N'$ and $A \neq \emptyset$, we have $E \supseteq A$ and hence $E$ is non-empty. According to \cref{item:properties-of-interior-closure-and-boundary:finite} of \cref{lemma:properties-of-interior-closure-and-boundary}, because $G_0$, $A$, and $N'$ are finite, so is $E$. Because $E$ is non-empty, according to \cref{theorem:existence-of-tiling}, there is a subset $E'$ of $G \modulo G_0$ and an $\ntuple{E, E'}$-tiling $T$ of $\mathcal{R}$. Because $G_0$ and $E$ are non-empty and finite, so is $E'$.

    Let $Y = \setOf{y \in Q^M \suchThat \ForEach t \in T \Holds p \not\occursIn_t y}$. For each $t \in T$, we have $t \actsByItsCoordinateOn p \notin \pi_{t \isSemiActedUponBy E}(Y)$ and therefore $\pi_{t \isSemiActedUponBy E}(Y) \subsetneqq Q^{t \isSemiActedUponBy E}$. According to \cref{lemma:entropy-bounded-above-if-strange-tiling-exists}, we have $\entropyOf_{\mathcal{F}}(Y) < \log\cardinalityOf{Q}$. Hence, according to \cref{theorem:entropy-does-non-increase}, we have $\entropyOf_{\mathcal{F}}(\Delta(Y)) < \log\cardinalityOf{Q}$.

    Let $x \in Q^M$. Put $S = \setOf{t \in T \suchThat p \occursIn_t x}$. According to \cref{lemma:exchanging-pattern-by-other-pattern-with-same-image-yields-configuration-with-same-image}, there is an $x' \in Q^M$ such that $x' \in Y$ and $\Delta(x) = \Delta(x')$. Therefore, $\Delta(Q^M) = \Delta(Y)$. In conclusion, $\entropyOf_{\mathcal{F}}(\Delta(Q^M)) < \log\cardinalityOf{Q}$.
  \end{proof}

  \begin{counterexample}[Tree]
  \label{example:tree:not-pre-injective-but-maximal-entropy}
    In the situation of counterexample \ref{example:tree:entropy-does-increase}, the global transition function $\Delta$ is not pre-injective but the entropy $\entropyOf_{\mathcal{F}}(\Delta(Q^M))$ is maximal as we show now. 

    First, let $c$ and $c'$ be the two global configurations of $Q^M$ such that $c \equiv 0$, $c'(m_0) = 1$, and $c'\restrictedTo_{M \smallsetminus \setOf{m_0}} \equiv 0$. Then, $\differenceOf(c, c')$ is finite and $\Delta(c) = \Delta(c')$ but $c \neq c'$. Hence, $\Delta$ is not pre-injective.

    Secondly, let $c$ be a global configuration of $Q^M$. And, let $c'$ be the global configuration of $Q^M$ such that $c'(m_0) = 0$ and, for each cell $m \in M$ and each generator $s \in \setOf{a, b, a^{-1}, b^{-1}}$ with $\lengthOf{m s} = \lengthOf{m} + 1$, we have $c'(m s) = c(m)$. For each generator $s \in \setOf{a, b, a^{-1}, b^{-1}}$, we have $\lengthOf{m_0 s} = \lengthOf{m_0} + 1$; and, for each cell $m \in M \smallsetminus \setOf{m_0}$, there are precisely three distinct generators $s_1$, $s_2$, and $s_3 \in \setOf{a, b, a^{-1}, b^{-1}}$ such that $\lengthOf{m s_k} = \lengthOf{m} + 1$, for $k \in \setOf{1, 2, 3}$. Hence, $\Delta(c') = c$. Therefore, $\Delta$ is surjective. In particular, $\entropyOf_{\mathcal{F}}(\Delta(Q^M)) = \log\cardinalityOf{Q}$. 
  \end{counterexample}

  The three preceding theorems yield

  \begin{main-theorem}(Garden of Eden Theorem; Edward Forrest Moore, 1962; John R. Myhill, 1963) 
  \label{theorem:garden-of-Eden} 
    Let $\mathcal{M} = \ntuple{M, G, \actsOnPoint}$ be a right-a\-me\-na\-ble left-ho\-mo\-ge\-neous space with finite stabilisers and let $\Delta$ be the global transition function of a big-cellular automaton over $\mathcal{M}$ with finite set of states and finite neighbourhood. The map $\Delta$ is surjective if and only if it is pre-injective. 
  \end{main-theorem} 

  \begin{proof}
    There is a coordinate system $\mathcal{K} = \ntuple{m_0, \family{g_{m_0, m}}_{m \in M}}$ such that there is a big-cellular automaton $\mathcal{C} = \ntuple{\mathcal{R}, Q, N, \delta}$ such that $Q$ and $N$ are finite and $\Delta$ is its global transition function. Moreover, because $\mathcal{C}$ is big, there is a $\mathcal{K}$-big subgroup $H$ of $G$ such that the local transition function $\delta$ is $\bullet_{H_0}$-invariant. And, because $G_0$ is finite, the cell space $\mathcal{R} = \ntuple{\mathcal{M}, \mathcal{K}}$ is right amenable.
    \begin{description}
      \item[Case 1:] $\cardinalityOf{Q} \leq 1$. If $\cardinalityOf{Q} = 0$, then, because $\cardinalityOf{M} \neq 0$, we have $\cardinalityOf{Q^M} = 0$. And, if $\cardinalityOf{Q} = 1$, then $\cardinalityOf{Q^M} = 1$. In either case, $\Delta$ is bijective, in particular, surjective and pre-injective.
      \item[Case 2:] $\cardinalityOf{Q} \geq 2$. According to \cref{theorem:not-surjective-implies-less-entropy} and \cref{item:entropy-basic-facts:whole-space} of \cref{lemma:entropy-basic-facts}, the map $\Delta$ is not surjective if and only if $\entropyOf_{\mathcal{F}}(\Delta(Q^M)) < \log\cardinalityOf{Q}$. And, according to \cref{theorem:less-entropy-implies-not-pre-injective} and \cref{theorem:not-pre-injective-implies-less-entropy}, we have $\entropyOf_{\mathcal{F}}(\Delta(Q^M)) < \log\cardinalityOf{Q}$ if and only if $\Delta$ is not pre-injective. Hence, $\Delta$ is not surjective if and only if it is not pre-injective. In conclusion, $\Delta$ is surjective if and only if it is pre-injective. \qedhere
    \end{description}
  \end{proof}

  \begin{remark}
    In the situation of \cref{remark:group:interior-closure-boundary}, \cref{theorem:garden-of-Eden} is theorem~5.3.1 in \cite{ceccherini-silberstein:coornaert:2010}.
  \end{remark}

  \begin{counterexample}[Tree] 
    The global transition function of the cellular automaton of \cref{example:tree:entropy-does-increase} is surjective but not pre-injective, which was shown in \cref{example:tree:not-pre-injective-but-maximal-entropy}.
  \end{counterexample}

  \begin{counterexample}[Muller] 
    The global transition function of the cellular automaton of \cref{example:non-maximal-entropy-but-pre-injective} is not surjective but pre-injective, which was shown there.
  \end{counterexample}

  \begin{example}[Exclusive Or]
    The global transition function of the elementary cellular automaton with Wolfram code $90$, whose local transition function combines the states of the left and right neighbours by exclusive or, is $4$-to-$1$ surjective and pre-injective but not injective. 
  \end{example}

  \begin{open-problem} 
    Laurent Bartholdi constructs in his papers \enquote{\citetitle*{bartholdi:2010}}\cite{bartholdi:2010} and~\enquote{\citetitle*{bartholdi:2016}}\cite{bartholdi:2016}, for each non-amenable group $G$, two finite sets $Q$ and $Q'$, and two cellular automata $\mathcal{C}$ and $\mathcal{C}'$ over $G$ such that the global transition function $\Delta \from Q^G \to Q^G$ of $\mathcal{C}$ is surjective but not pre-injective and the global transition function $\Delta' \from (Q')^G \to (Q')^G$ of $\mathcal{C}'$ is not surjective but pre-injective. Are similar constructions possible for non-right-a\-me\-na\-ble left-ho\-mo\-ge\-neous spaces with finite stabilisers?
  \end{open-problem}

  \section{Construction of Non-Degenerated Left Homogeneous Spaces}
  \label{section:construction-of-example}

  \paragraph{Introduction.} So far we have only seen examples of right-a\-me\-na\-ble left-ho\-mo\-ge\-neous spaces with finite stabilisers for which there is a subgroup that acts freely and transitively. The global transition function of a semi-cellular automaton on such a space is essentially the global transition function of a cellular automaton over a group. For those spaces, \cref{theorem:garden-of-Eden} states nothing new. A simple construction of right-a\-me\-na\-ble left-ho\-mo\-ge\-neous spaces with finite stabilisers for which there is no subgroup that acts freely and transitively goes like this: Act with the direct product of the automorphism groups of the coloured Cayley graph of a group and a vertex-transitive, finite, and directed non-Cayley graph on the direct product of the vertices of these graphs. Full details and concrete examples of this construction are given below.

  \paragraph{Body.} The Cartesian product of two cell spaces, and how its right quotient set semi-action and its notions of interior, closure, and boundary relate to the respective constructs and notions of its components is given and shown in

  \begin{lemma} 
  \label{lemma:interior-closure-and-boundary-of-Cartesian-product-of-cell-spaces}
    Let $\mathcal{R} = \ntuple{\ntuple{M, G, \actsOnPoint}, \ntuple{m_0, \family{g_{m_0, m}}_{m \in M}}}$ and $\mathcal{R}' = \ntuple{\ntuple{M', G', \actsOnPoint'}, \ntuple{m_0', \family{g_{m_0', m'}'}_{m' \in M'}}}$ be two cell spaces, and let $\mathcal{R}''$ be the cell space $\ntuple{\ntuple{M \times M', G \times G', \actsOnPoint \times \actsOnPoint'}, \ntuple{(m_0, m_0'), \family{(g_{m_0, m}, g_{m_0', m'}')}_{(m, m') \in M \times M'}}}$, where
    \begin{align*} 
      \actsOnPoint \times \actsOnPoint' \from (G \times G') \times (M \times M') &\to M \times M', \mathnote{left group action $\actsOnPoint \times \actsOnPoint'$ of $G \times G'$ on $M \times M'$}\\
      ((g, g'), (m, m')) &\mapsto (g \actsOnPoint m, g' \actsOnPoint' m').
    \end{align*}
    Furthermore, let $A$ be a subset of $M$, let $A'$ be a subset of $M'$, let $A''$ be a subset of $M \times M'$ such that $\setOf{m \in M \suchThat \Exists m' \in M' \SuchThat (m, m') \in A''} = A$, let $E$ be a subset of $G \modulo G_0$, and let $E'$ be a subset of $G' \modulo G_0'$. Then,
    \begin{aenumerate}
      \item\label{item:interior-closure-and-boundary-of-Cartesian-product-of-cell-spaces:rightsemiaction}
            $\isSemiActedUponBy'' = \isSemiActedUponBy \times \isSemiActedUponBy'$, where
            \begin{align*}
              \isSemiActedUponBy \times \isSemiActedUponBy' \from{} &(M \times M') \times ((G \times G') \modulo (G_0 \times G_0')) \to M \times M',\\
              &((m, m'), (g, g') (G_0 \times G_0')) \mapsto (m \isSemiActedUponBy g G_0, m' \isSemiActedUponBy' g' G_0').
            \end{align*}
      \item\label{item:interior-closure-and-boundary-of-Cartesian-product-of-cell-spaces:interior-closure-and-boundary}
            $(A \times A')^{-E \times E'} = A^{-E} \times (A')^{-E'}$, $(A \times A')^{+E \times E'} = A^{+E} \times (A')^{+E'}$, and $\boundaryOf_{E \times E'} (A \times A') = (\boundaryOf_E A \times (A')^{+E'}) \cup (A^{+E} \times \boundaryOf_{E'} A')$, where
            \begin{equation*}
              E \times E' = \setOf{(g, g') (G_0 \times G_0') \suchThat g G_0 \in E, g' G_0' \in E'}.
            \end{equation*}
      \item\label{item:interior-closure-and-boundary-of-Cartesian-product-of-cell-spaces:general-set-of-cells-special-thickness}
            $(A'')^{-E \times (G' \modulo G_0')} \subseteq A^{-E} \times M'$, $(A'')^{+E \times (G' \modulo G_0')} = A^{+E} \times M'$, and $\boundaryOf_{E \times (G' \modulo G_0')} A'' \supseteq \boundaryOf_E A \times M'$, where
            \begin{equation*} 
              E \times (G' \modulo G_0') = \setOf{(g, g') (G_0 \times G_0') \suchThat g G_0 \in E, g' G_0' \in G' \modulo G_0'}. \qedhere
            \end{equation*}
    \end{aenumerate}
  \end{lemma} 

  \begin{proof}
    Note that the stabiliser $G_0''$ of $(m_0, m_0')$ under $\actsOnPoint \times \actsOnPoint'$ is $G_0 \times G_0'$.
    \begin{aenumerate}
      \item For each $(m, m') \in M \times M'$ and each $(g, g') (G_0 \times G_0') \in (G \times G') \modulo (G_0 \times G_0')$,
            \begin{align*}
              &(m, m') \isSemiActedUponBy'' (g, g') (G_0 \times G_0')\\
              &{}= (g_{m_0, m}, g_{m_0', m'}') \cdot (g, g') \mathbin{\actsOnPoint \times \actsOnPoint'} (m_0, m_0')\\ 
              &{}= (g_{m_0, m} g \actsOnPoint m_0, g_{m_0', m'}' g' \actsOnPoint' m_0')\\
              &{}= (m \isSemiActedUponBy g G_0, m' \isSemiActedUponBy' g' G_0').
            \end{align*}
            Therefore, $\isSemiActedUponBy'' = \isSemiActedUponBy \times \isSemiActedUponBy'$.
      \item Let $(m, m') \in M \times M'$. Then, $(m, m') \isSemiActedUponBy'' E \times E' = (m \isSemiActedUponBy E) \times (m' \isSemiActedUponBy' E')$. Hence, if $E \neq \emptyset$ and $E' \neq \emptyset$, then $(m, m') \isSemiActedUponBy'' E \times E' \subseteq A \times A'$ if and only if $m \isSemiActedUponBy E \subseteq A$ and $m' \isSemiActedUponBy' E' \subseteq A'$. And, $((m, m') \isSemiActedUponBy'' E \times E') \cap (A \times A') \neq \emptyset$ if and only if $(m \isSemiActedUponBy E) \cap A \neq \emptyset$ and $(m' \isSemiActedUponBy' E') \cap A' \neq \emptyset$. Therefore, $(A \times A')^{-E \times E'} = A^{-E} \times (A')^{-E'}$ and $(A \times A')^{+E \times E'} = A^{+E} \times (A')^{+E'}$. Note that the first equality holds in the case that $E = \emptyset$ or $E' = \emptyset$, because in that case $(A \times A')^{-E \times E'} = \emptyset$, and $A^{-E} = \emptyset$ or $(A')^{-E'} = \emptyset$. Moreover, 
            \begin{align*}
              \boundaryOf_{E \times E'} (A \times A')
              &= \parens[\big]{A^{+E} \times (A')^{+E'}} \smallsetminus \parens[\big]{A^{-E} \times (A')^{-E'}}\\ 
              &= \begin{aligned}[t]
                   &\parens[\big]{(A^{+E} \smallsetminus A^{-E}) \times (A')^{+E'}}\\
                   &{}\cup \parens[\big]{A^{+E} \times ((A')^{+E'} \smallsetminus (A')^{-E'})}
                 \end{aligned}\\
              &= \parens[\big]{\boundaryOf_E A \times (A')^{+E'}} \cup \parens[\big]{A^{+E} \times \boundaryOf_{E'} A'}.
            \end{align*}
      \item Let $(m, m') \in M \times M'$. Then, $(m, m') \isSemiActedUponBy'' E \times (G' \modulo G_0') = (m \isSemiActedUponBy E) \times M'$. Hence, if $(m, m') \isSemiActedUponBy'' E \times (G' \modulo G_0') \subseteq A''$, then $m \isSemiActedUponBy E \subseteq A$. And, $((m, m') \isSemiActedUponBy'' E \times (G' \modulo G_0')) \cap A'' \neq \emptyset$ if and only if $(m \isSemiActedUponBy E) \cap A \neq \emptyset$. Therefore, $(A'')^{-E \times (G' \modulo G_0')} \subseteq A^{-E} \times M'$ and $(A'')^{+E \times (G' \modulo G_0')} = A^{+E} \times M'$. Moreover,
            \begin{align*}
              \boundaryOf_{E \times (G' \modulo G_0')} A''
              &\supseteq (A^{+E} \times M') \smallsetminus (A^{-E} \times M')\\
              &=         \begin{aligned}[t]
                           &\parens[\big]{(A^{+E} \smallsetminus A^{-E}) \times M'}\\
                           &{}\cup \parens[\big]{A^{+E} \times (M' \smallsetminus M')}
                         \end{aligned}\\
              &=         \boundaryOf_E A \times M'. \qedhere
            \end{align*}
    \end{aenumerate}
  \end{proof}

  The Cartesian product of a left-ho\-mo\-ge\-neous space with a finite one is right amenable if and only if the former one is right amenable, which is shown in

  \begin{lemma}
  \label{lemma:product-with-finite-right-amenable-if-and-only-if-first-component-right-amenable}
    Let $\mathcal{M} = \ntuple{M, G, \actsOnPoint}$ and $\mathcal{M}' = \ntuple{M', G', \actsOnPoint'}$ be two left-ho\-mo\-ge\-neous spaces with finite stabilisers such that $M'$ is finite, and let $\mathcal{M}''$ be the left-ho\-mo\-ge\-neous space $\ntuple{M \times M', G \times G', \actsOnPoint \times \actsOnPoint'}$. The space $\mathcal{M}''$ is right amenable if and only if the space $\mathcal{M}$ is right amenable.
  \end{lemma}

  \begin{proof} 
    Let $\mathcal{K} = \ntuple{m_0, \family{g_{m_0, m}}_{m \in M}}$ and $\mathcal{K}' = \ntuple{m_0', \family{g_{m_0', m'}'}_{m' \in M'}}$ be two coordinate systems for $\mathcal{M}$ and $\mathcal{M}'$ respectively, and let $\mathcal{K}''$ be the coordinate system $\ntuple{(m_0, m_0'), \family{(g_{m_0, m}, g_{m_0', m'}')}_{(m, m') \in M \times M'}}$ for $\mathcal{M}''$. Note that, because $\mathcal{M}$ and $\mathcal{M}'$ have finite stabilisers, so has $\mathcal{M}''$.

    First, let $\mathcal{M}''$ be right amenable. Then, because $\mathcal{M}''$ has finite stabilisers, there is a right Følner net $\net{F_i''}_{i \in I}$ in $\mathcal{R}'' = \ntuple{\mathcal{M}'', \mathcal{K}''}$. Put $F_i = \setOf{m \in M \suchThat \Exists m' \in M' \SuchThat (m, m') \in F_i''}$. Let $E$ be a finite subset of $G \modulo G_0$. Then, according to \cref{item:interior-closure-and-boundary-of-Cartesian-product-of-cell-spaces:general-set-of-cells-special-thickness} of \cref{lemma:interior-closure-and-boundary-of-Cartesian-product-of-cell-spaces} and because $F_i'' \subseteq F_i \times M'$,
    \begin{equation*} 
      \frac{\cardinalityOf{\boundaryOf_E F_i}}{\cardinalityOf{F_i}}
      \leq \frac{\cardinalityOf{M'}}{\cardinalityOf{M'}} \cdot \frac{\cardinalityOf{\boundaryOf_{E \times (G' \modulo G_0')} F_i''}}{\cardinalityOf{F_i''}}
      = \frac{\cardinalityOf{\boundaryOf_{E \times (G' \modulo G_0')} F_i''}}{\cardinalityOf{F_i''}}.
    \end{equation*}
    Hence, because $\net{F_i''}_{i \in I}$ is a right Følner net in $\mathcal{R''}$, the net $\net{F_i}_{i \in I}$ is a right Følner net in $\ntuple{\mathcal{M}, \mathcal{K}}$. In conclusion, $\mathcal{M}$ is amenable.

    Secondly, let $\mathcal{M}$ be amenable. Then, because $\mathcal{M}$ has finite stabilisers, there is a right Følner net $\net{F_i}_{i \in I}$ in $\mathcal{R} = \ntuple{\mathcal{M}, \mathcal{K}}$. Let $E''$ be a finite subset of $(G \times G') \modulo (G \times G')_0$. Put $E = \setOf{g G_0 \suchThat \Exists g' \in G' \SuchThat (g, g') (G \times G')_0 \in E''}$ and put $E' = G' \modulo G'_0$. Then, according to \cref{item:interior-closure-and-boundary-of-Cartesian-product-of-cell-spaces:interior-closure-and-boundary} of \cref{lemma:interior-closure-and-boundary-of-Cartesian-product-of-cell-spaces}, we have $\boundaryOf_{E''} (F_i \times M') \subseteq \boundaryOf_{E \times E'} (F_i \times M') \subseteq ((\boundaryOf_E F_i) \times M') \cup (F_i^{+E} \times (\boundaryOf_{E'} M')) = (\boundaryOf_E F_i) \times M'$. Hence, 
    \begin{equation*}
      \frac{\cardinalityOf{\boundaryOf_{E''} (F_i \times M')}}{\cardinalityOf{F_i \times M'}}
      \leq \frac{\cardinalityOf{\boundaryOf_E F_i} \cdot \cardinalityOf{M'}}{\cardinalityOf{F_i} \cdot \cardinalityOf{M'}}
      =    \frac{\cardinalityOf{\boundaryOf_E F_i}}{\cardinalityOf{F_i}}.
    \end{equation*}
    Therefore, because $\net{F_i}_{i \in I}$ is a right Følner net in $\mathcal{R}$, the net $\net{F_i \times M'}_{i \in I}$ is a right Følner net in $\mathcal{M}''$. In conclusion, $\mathcal{M}''$ is right amenable.
  \end{proof}

  Vertex-transitivity, the Cartesian product of two graphs, and an action on such a product are introduced in the forthcoming definitions.

  \begin{definition} 
    Let $\mathcal{G}$ be a directed graph. It is called \graffito{vertex-transitive}\define{vertex-transitive} if and only if its automorphism group acts transitively on its vertices by function application.
  \end{definition}

  \begin{remark} 
    Cayley graphs of groups are vertex-transitive.
  \end{remark}

  \begin{theorem}[Sabidussi's Theorem]
  \label{theorem:Sabidussi-theorem}
    Let $\mathcal{G}$ be a directed graph. It is a Cayley graph if and only if a subgroup of its automorphism group acts freely and transitively on its vertices by function application.
  \end{theorem}

  \begin{proof}
    See proposition~3.1(b) in \cite{babai:1996}. 
  \end{proof}

  \begin{definition}
    Let $\mathcal{G} = \ntuple{V, E}$ and $\mathcal{G}' = \ntuple{V', E'}$ be two directed graphs. The graph $\mathcal{G} \graphProductOf \mathcal{G}' = \ntuple{V \times V', \setOf{((v_1, v_1'), (v_2, v_2')) \in (V \times V') \times (V \times V') \suchThat v_1 = v_2 \land (v_1', v_2') \in E' \text{ or } v_1' = v_2' \land (v_1, v_2) \in E}}$ is called \define{Cartesian product of $\mathcal{G}$ and $\mathcal{G}'$}\graffito{Cartesian product $\mathcal{G} \graphProductOf \mathcal{G}'$ of $\mathcal{G}$ and $\mathcal{G}'$}.
  \end{definition}

  \begin{definition}
  \label{definition:group-action-of-automorphisms-on-cartesian-product-of-graphs}
    Let $\mathcal{G} \graphProductOf \mathcal{G}'$ be the Cartesian product of $\mathcal{G} = \ntuple{V, E}$ and $\mathcal{G}' = \ntuple{V', E'}$, and let $\actsOnPoint$ and $\actsOnPoint'$ be the left group actions of $\automorphismsOf(\mathcal{G})$ on $V$ and of $\automorphismsOf(\mathcal{G}')$ on $V'$ by function application. The direct product $\automorphismsOf(\mathcal{G}) \times \automorphismsOf(\mathcal{G}')$ acts on $V \times V'$ on the left by
    \begin{align*}
      \actsOnPoint \graphProductOf \actsOnPoint' \from (\automorphismsOf(\mathcal{G}) \times \automorphismsOf(\mathcal{G}')) \times (V \times V') &\to V \times V', \mathnote{left group action $\actsOnPoint \graphProductOf \actsOnPoint'$ of $\automorphismsOf(\mathcal{G}) \times \automorphismsOf(\mathcal{G}')$ on $V \times V'$}\\
      ((\varphi, \varphi'), (v, v')) &\mapsto (\varphi \actsOnPoint v, \varphi' \actsOnPoint' v'). \qedhere
    \end{align*}
  \end{definition}

  \begin{remark}
    The map $\blank \mathbin{(\actsOnPoint \graphProductOf \actsOnPoint')} \from \automorphismsOf(\mathcal{G}) \times \automorphismsOf(\mathcal{G}') \to \automorphismsOf(\mathcal{G} \graphProductOf \mathcal{G}')$
    is an injective group homomorphism.
  \end{remark} 

  \begin{remark} 
    If $\mathcal{G}$ and $\mathcal{G}'$ are vertex-transitive, then the left group action $\actsOnPoint \graphProductOf \actsOnPoint'$ is transitive.
  \end{remark}

  When a subgroup acts freely and transitively on the Cartesian product of two graphs is characterised in

  \begin{lemma}
  \label{lemma:free-and-transitive-action-on-cartesian-product-of-graphs}
    Let $\mathcal{G} \graphProductOf \mathcal{G}'$ be the Cartesian product of $\mathcal{G} = \ntuple{V, E}$ and $\mathcal{G}' = \ntuple{V', E'}$. There is a subgroup $H''$ of $\automorphismsOf(\mathcal{G}) \times \automorphismsOf(\mathcal{G}')$ that acts freely and transitively on $V \times V'$ by $\actsOnPoint \graphProductOf \actsOnPoint'$ if and only if there is a subgroup $H$ of $\automorphismsOf(\mathcal{G})$ and there is a subgroup $H'$ of $\automorphismsOf(\mathcal{G}')$ such that $H$ acts freely and transitively on $V$ by $\actsOnPoint$ and $H'$ acts freely and transitively on $V'$ by $\actsOnPoint'$.
  \end{lemma} 

  \begin{proof} 
    First, let $H''$ be a subgroup of $\automorphismsOf(\mathcal{G}) \times \automorphismsOf(\mathcal{G}')$ that acts freely and transitively on $V \times V'$ by $\actsOnPoint \graphProductOf \actsOnPoint'$. Furthermore, let $v'$ be a vertex of $V'$. The set $F = \setOf{h'' \in H'' \suchThat h'' \mathbin{\actsOnPoint \graphProductOf \actsOnPoint'} V \times \setOf{v'} \subseteq V \times \setOf{v'}}$ is a subgroup of $H''$ and the left group action $(\actsOnPoint \graphProductOf \actsOnPoint')\restrictedTo_{F \times (V \times \setOf{v'}) \to V \times \setOf{v'}}$ is free and transitive. The set $H = \pi_1(F)$ is a subgroup of $\automorphismsOf(\mathcal{G})$, which acts freely and transitively on $V$ by $\actsOnPoint$, where $\pi_1$ is the projection homomorphism from $\automorphismsOf(\mathcal{G}) \times \automorphismsOf(\mathcal{G}')$ onto $\automorphismsOf(\mathcal{G})$. Analogously, a subgroup $H'$ of $\automorphismsOf(\mathcal{G}')$ can be constructed that acts freely and transitively on $V'$ by $\actsOnPoint'$.

    Secondly, let $H$ be a subgroup of $\automorphismsOf(\mathcal{G})$ and let $H'$ be a subgroup of $\automorphismsOf(\mathcal{G}')$ such that $H$ acts freely and transitively on $V$ by $\actsOnPoint$ and $H'$ acts freely and transitively on $V'$ by $\actsOnPoint'$. The direct product $H'' = H \times H'$ acts freely and transitively on $V \times V'$ by $\actsOnPoint \graphProductOf \actsOnPoint'$.
  \end{proof}

  It follows that there is no subgroup that acts freely and transitively on the Cartesian product of a Cayley and a non-Cayley graph, which is stated in

  \begin{corollary}
  \label{corollary:free-and-transitive-action-on-cartesian-product-of-graphs}
    Let $\mathcal{G} \graphProductOf \mathcal{G}'$ be the Cartesian product of $\mathcal{G} = \ntuple{V, E}$ and $\mathcal{G}' = \ntuple{V', E'}$, where $\mathcal{G}$ or $\mathcal{G}'$ is not a Cayley graph. There is no subgroup of $\automorphismsOf(\mathcal{G}) \times \automorphismsOf(\mathcal{G}')$ that acts freely and transitively on $V \times V'$ by $\actsOnPoint \graphProductOf \actsOnPoint'$.
  \end{corollary}

  \begin{proof}
    This is a direct consequence of \cref{theorem:Sabidussi-theorem,lemma:free-and-transitive-action-on-cartesian-product-of-graphs}.
  \end{proof}

  \begin{example}
    Let $G$ be a group, let $S$ be a generating set of $G$, let $\mathcal{G} = \ntuple{V, E}$ be the coloured $S$-Cayley graph of $G$, let $\actsOnPoint$ be the left group action of $\automorphismsOf(\mathcal{G})$ on $V$ by function application, let $\mathcal{G}' = \ntuple{V', E'}$ be a vertex-transitive, finite, and directed non-Cayley graph, and let $\actsOnPoint'$ be the left group action of $\automorphismsOf(\mathcal{G}')$ on $V'$ by function application, and let $\mathcal{M}''$ be the left-ho\-mo\-ge\-neous space $\ntuple{V \times V', \automorphismsOf(\mathcal{G}) \times \automorphismsOf(\mathcal{G}'), \actsOnPoint \graphProductOf \actsOnPoint'}$.

    Note that $V = G$, that $G \simeq \automorphismsOf(\mathcal{G})$ by $g \mapsto g \cdot \blank$, that under this identification the map $\actsOnPoint$ is the group multiplication $\cdot$ and the map $\actsOnPoint \graphProductOf \actsOnPoint'$ is $(G \times \automorphismsOf(\mathcal{G}')) \times (G \times V') \to G \times V'$, $((g, \varphi'), (v, v')) \mapsto (g \cdot v, \varphi'(v'))$, and that $\ntuple{V, \automorphismsOf(\mathcal{G}), \actsOnPoint}$ is right amenable if and only if $G$ is amenable.

    According to \cref{corollary:free-and-transitive-action-on-cartesian-product-of-graphs}, there is no subgroup of $\automorphismsOf(\mathcal{G}) \times \automorphismsOf(\mathcal{G}')$ that acts freely and transitively on $V \times V'$ by $\actsOnPoint \graphProductOf \actsOnPoint'$. And, according to \cref{lemma:product-with-finite-right-amenable-if-and-only-if-first-component-right-amenable} and the note above, the space $\mathcal{M}''$ is right amenable if and only if the group $G$ is amenable. For example:
    \begin{aenumerate} 
      \item If $G$ is the integer lattice $\Z$ and $\mathcal{G}'$ is the Petersen graph, then $\mathcal{M}''$ is right amenable. 
      \item If $G$ is the free group $F_2$ over $\setOf{a, b}$, where $a \neq b$ and $\mathcal{G}'$ is the Coxeter graph, then $\mathcal{M}''$ is not right amenable. \qedhere 
    \end{aenumerate}
  \end{example}

  \clearToOddPage
  \chapter{Finitely Right Generated, Growth, and Right Amenability}
  \label{chapter:growth}

  \paragraph{Abstract.} We introduce right-gen\-er\-at\-ing sets, left Cayley graphs, growth functions, types and rates, and isoperimetric constants for cell spaces and some of these notions also for left-ho\-mo\-ge\-neous spaces; characterise right-a\-me\-na\-ble finitely right-gen\-er\-at\-ed cell spaces with finite stabilisers as those whose isoperimetric constant is $0$; and prove that finitely right-gen\-er\-at\-ed left-ho\-mo\-ge\-neous spaces with finite stabilisers of sub-exponential growth are right amenable, in particular, quotient sets of groups of sub-exponential growth by finite subgroups acted upon by left multiplication are right amenable.

  \paragraph{Remark.} Most parts of this chapter appeared in the preprint \enquote{\citetitle*{wacker:growth:2017}}\cite{wacker:growth:2017} and they generalise parts of chapter~6 of the monograph \enquote{\citetitle*{ceccherini-silberstein:coornaert:2010}}\cite{ceccherini-silberstein:coornaert:2010}.

  \paragraph{Summary.} A cell space $\mathcal{R}$ is \emph{finitely right generated} if there is a finite subset $S$ of $G \modulo G_0$ with $G_0 \cdot S \subseteq S$ such that, for each point $m \in M$, there is a sequence $\sequence{s_i}_{i \in \setOf{1, 2, \dotsc, k}}$ of elements in $S \cup S^{-1}$ such that $m = (((m_0 \isSemiActedUponBy s_1) \isSemiActedUponBy s_2) \isSemiActedUponBy \dotsb) \isSemiActedUponBy s_k$. The finite right-gen\-er\-at\-ing set $S$ induces the \emph{left $S$-Cayley graph} structure on $M$ given by: For each point $m \in M$ and each generator $s \in S$, there is an edge from $m$ to $m \isSemiActedUponBy s$. The length of the shortest path between two points of $M$ yields the \emph{$S$-metric}. The ball of radius $\rho \in \N_0$ centred at $m \in M$, denoted by $\ball_S(m, \rho)$, is the set of all points whose distance to $m$ is not greater than $\rho$. The \emph{$S$-growth function} is the map $\gamma_S \from \N_0 \to \N_0$, $k \mapsto \cardinalityOf{\ball_S(m_0, k)}$; the \emph{growth type of $\mathcal{R}$}, which does not depend on $S$, is the equivalence class $\equivalenceClassOf{\gamma_S}_\isEquivalentTo$, where two growth functions are equivalent if they dominate each other; and the \emph{$S$-growth rate} is the limit point of the sequence $\sequence{\sqrt[k]{\gamma_S(k)}}_{k \in \N_0}$.

  A finitely right-gen\-er\-at\-ed cell space $\mathcal{R}$ is said to have \emph{sub-exponential growth} if its growth type is not $\equivalenceClassOf{\exp}_\isEquivalentTo$, which is the case if and only if its growth rates are $1$. The \emph{$S$-isoperimetric constant} is a real number between $0$ and $1$ that measures, broadly speaking, the invariance under $\isSemiActedUponBy\restrictedTo_{M \times S}$ that a finite subset of $M$ can have, where $0$ means maximally and $1$ minimally invariant. In the case that $G_0$ is finite, this constant is $0$ if and only if $\mathcal{R}$ is right amenable; and, if $\mathcal{R}$ has sub-exponential growth, then it is right amenable; and, if $G$ has sub-exponential growth, then so has $\mathcal{R}$.

  Cayley graphs were introduced by Arthur Cayley in his paper \enquote{\citetitle*{cayley:1878}}\cite{cayley:1878}. The notion of growth was introduced by Vadim Arsenyevich Efremovich and Albert S. Švarc in their papers \enquote{\citetitle*{efremovich:1952}}\cite{efremovich:1952} and \enquote{\citetitle*{svarc:1955}}\cite{svarc:1955}. Mikhail Leonidovich Gromov was the first to study groups through their word metrics, see for example his paper \enquote{\citetitle*{gromov:1984}}\cite{gromov:1984}.

  \paragraph{Contents.} In \cref{section:gen-set} we introduce right-gen\-er\-at\-ing sets. In \cref{section:multigraphs} we recapitulate directed multigraphs. In \cref{section:Cayley} we introduce left Cayley graphs induced by right-gen\-er\-at\-ing sets. In \cref{section:metrics} we introduce metrics and lengths induced by left Cayley graphs. In \cref{section:balls-and-spheres} we consider balls and spheres induced by metrics. In \cref{section:interiors-closures-boundaries} we consider interiors, closures, and boundaries of any thickness of sets. In \cref{section:growth-functions-and-types} we recapitulate growth functions and types. In \cref{section:cell-spaces-growth-functions-and-types} we introduce growth functions and types of cell spaces. In \cref{section:growth-rates} we introduce growth rates of cell spaces. In \cref{section:Folner-conditions} we prove that right amenability and having isoperimetric constant $0$ are equivalent, and we characterise right Følner nets. And in \cref{section:sub-exponential-growth-and-amenability} we prove that having sub-exponential growth implies right amenability.

  \paragraph{Preliminary Notions.} A \emph{left group set} is a triple $\ntuple{M, G, \actsOnPoint}$, where $M$ is a set, $G$ is a group, and $\actsOnPoint$ is a map from $G \times M$ to $M$, called \emph{left group action of $G$ on $M$}, such that $G \to \symmetricGroupOf(M)$, $g \mapsto [g \actsOnPoint \blank]$, is a group homomorphism. The action $\actsOnPoint$ is \emph{transitive} if $M$ is non-empty and for each $m \in M$ the map $\blank \actsOnPoint m$ is surjective; and \emph{free} if for each $m \in M$ the map $\blank \actsOnPoint m$ is injective. For each $m \in M$, the set $G \actsOnPoint m$ is the \emph{orbit of $m$}, the set $G_m = (\blank \actsOnPoint m)^{-1}(m)$ is the \emph{stabiliser of $m$}, and, for each $m' \in M$, the set $G_{m, m'} = (\blank \actsOnPoint m)^{-1}(m')$ is the \emph{transporter of $m$ to $m'$}.

  A \emph{left-ho\-mo\-ge\-neous space} is a left group set $\mathcal{M} = \ntuple{M, G, \actsOnPoint}$ such that $\actsOnPoint$ is transitive. A \emph{coordinate system for $\mathcal{M}$} is a tuple $\mathcal{K} = \ntuple{m_0, \family{g_{m_0, m}}_{m \in M}}$, where $m_0 \in M$ and for each $m \in M$ we have $g_{m_0, m} \actsOnPoint m_0 = m$. The stabiliser $G_{m_0}$ is denoted by $G_0$. The tuple $\mathcal{R} = \ntuple{\mathcal{M}, \mathcal{K}}$ is a \emph{cell space}. The map $\isSemiActedUponBy \from M \times G \modulo G_0 \to M$, $(m, g G_0) \mapsto g_{m_0, m} g g_{m_0, m}^{-1} \actsOnPoint m\ (= g_{m_0, m} g \actsOnPoint m_0)$ is a \emph{right semi-action of $G \modulo G_0$ on $M$ with defect $G_0$}, which means that
  \begin{equation*}
    \ForEach m \in M \Holds m \isSemiActedUponBy G_0 = m,
  \end{equation*}
  and
  \begin{multline*}
    \ForEach m \in M \ForEach g \in G \Exists g_0 \in G_0 \SuchThat \ForEach \mathfrak{g}' \in G \modulo G_0 \Holds\\
          m \isSemiActedUponBy g \cdot \mathfrak{g}' = (m \isSemiActedUponBy g G_0) \isSemiActedUponBy g_0 \cdot \mathfrak{g}'.
  \end{multline*}
  It is \emph{transitive}, which means that the set $M$ is non-empty and for each $m \in M$ the map $m \isSemiActedUponBy \blank$ is surjective; and \emph{free}, which means that for each $m \in M$ the map $m \isSemiActedUponBy \blank$ is injective; and \emph{semi-commutes with $\actsOnPoint$}, which means that
  \begin{multline*}
    \ForEach m \in M \ForEach g \in G \Exists g_0 \in G_0 \SuchThat \ForEach \mathfrak{g}' \in G \modulo G_0 \Holds\\
          (g \actsOnPoint m) \isSemiActedUponBy \mathfrak{g}' = g \actsOnPoint (m \isSemiActedUponBy g_0 \cdot \mathfrak{g}').
  \end{multline*}
  The maps $\iota \from M \to G \modulo G_0$, $m \mapsto G_{m_0, m}$, and $m_0 \isSemiActedUponBy \blank$ are inverse to each other. Under the identification of $M$ with $G \modulo G_0$ by either of these maps, we have $\isSemiActedUponBy \from (m, \mathfrak{g}) \mapsto g_{m_0, m} \actsOnPoint \mathfrak{g}$.

  A left-ho\-mo\-ge\-neous space $\mathcal{M}$ is \emph{right amenable} if there is a coordinate system $\mathcal{K}$ for $\mathcal{M}$ and there is a finitely additive probability measure $\mu$ on $M$ such that 
  \begin{equation*}
    \ForEach \mathfrak{g} \in G \modulo G_0 \ForEach A \subseteq M \Holds \parens[\big]{(\blank \isSemiActedUponBy \mathfrak{g})\restrictedTo_A \text{ injective} \implies \mu(A \isSemiActedUponBy \mathfrak{g}) = \mu(A)},
  \end{equation*}
  in which case the cell space $\mathcal{R} = \ntuple{\mathcal{M}, \mathcal{K}}$ is called \emph{right amenable}. When the stabiliser $G_0$ is finite, that is the case if and only if there is a \emph{right Følner net in $\mathcal{R}$ indexed by $(I, \leq)$}, which is a net $\net{F_i}_{i \in I}$ in $\setOf{F \subseteq M \suchThat F \neq \emptyset, F \text{ finite}}$ such that
  \begin{equation*}
    \ForEach \mathfrak{g} \in G \modulo G_0 \Holds \lim_{i \in I} \frac{\cardinalityOf{F_i \smallsetminus (\blank \isSemiActedUponBy \mathfrak{g})^{-1}(F_i)}}{\cardinalityOf{F_i}} = 0.
  \end{equation*}

  \section{Right Generating Sets}
  \label{section:gen-set}

  In this section, let $\mathcal{R} = \ntuple{\mathcal{M}, \mathcal{K}} = \ntuple{\ntuple{M, G, \actsOnPoint}, \ntuple{m_0, \family{g_{m_0, m}}_{m \in M}}}$ be a cell space.

  \paragraph{Contents.} In \cref{definition:right-generating-set} we introduce right-gen\-er\-at\-ing sets of $\mathcal{R}$. And in \cref{lemma:gen-set-of-group-induces-right-gen-set-of-action} we show how generating sets of $G$ induce right ones of $\mathcal{R}$.

  \begin{definition}
  \label{definition:right-generating-set}
    Let $S$ be a subset of $G \modulo G_0$ such that $G_0 \cdot S \subseteq S$.
    \begin{aenumerate}
      \item The set $\setOf{g^{-1} G_0 \suchThat s \in S, g \in s}$ is denoted by $S^{-1}$\graffito{$S^{-1}$}\index[symbols]{Sminus1@$S^{-1}$}.
      \item For each element $m \in M$, each non-negative integer $k$, and each finite sequence $\sequence{s_i}_{i \in \setOf{1, 2, \dotsc, k}}$ of elements in $S \cup S^{-1}$, the element
            \begin{equation*}
              \Big(\big((m \isSemiActedUponBy s_1) \isSemiActedUponBy s_2\big) \isSemiActedUponBy \dotsb\Big) \isSemiActedUponBy s_k
            \end{equation*}
            is denoted by $m \isSemiActedUponBy \sequence{s_i}_{i \in \setOf{1, 2, \dotsc, k}}$\graffito{$m \isSemiActedUponBy \sequence{s_i}_{i \in \setOf{1, 2, \dotsc, k}}$}\index[symbols]{arrowleftunderscorem@$m \isSemiActedUponBy \sequence{s_i}_{i \in \setOf{1, 2, \dotsc, k}}$}, where, in the case that $k = 0$, the sequence $\sequence{s_i}_{i \in \setOf{1, 2, \dotsc, k}}$ is the empty sequence $()$ and $m \isSemiActedUponBy \sequence{s_i}_{i \in \setOf{1, 2, \dotsc, k}}$ is equal to $m$.
      \item The set $S$ is said to \define{right generate $\mathcal{R}$}\graffito{$S$ right generates $\mathcal{R}$}\index{generate right R@right generate $\mathcal{R}$}, called \define{right-gen\-er\-at\-ing set of $\mathcal{R}$}\graffito{right-gen\-er\-at\-ing set $S$ of $\mathcal{R}$}\index{generating set right@right-generating set of $\mathcal{R}$}, and each element $s \in S$ is called \define{right generator}\graffito{right generator $s$}\index{generator right@right generator} if and only if
            \begin{equation*}
              \ForEach m \in M \Exists \sequence{s_i}_{i \in \setOf{1, 2, \dotsc, k}} \text{ in } S \cup S^{-1} \SuchThat m_0 \isSemiActedUponBy \sequence{s_i}_{i \in \setOf{1, 2, \dotsc, k}} = m.
            \end{equation*}
      \item The set $S$ is called \defineX{symmetric}{symmetric!right-gen\-er\-at\-ing set}\graffito{symmetric} if and only if $S^{-1} \subseteq S$. \qedhere
    \end{aenumerate}
  \end{definition}

  \begin{remark}
    If $S$ is a right-gen\-er\-at\-ing set of $\mathcal{R}$, then $S \cup S^{-1}$ is a symmetric one; and, if $S$ is also finite and $G_0$ is finite, then $S \cup S^{-1}$ is finite. 
  \end{remark}

  \begin{definition}
    Let $P$ be an adjective. The cell space $\mathcal{R}$ is called \define{$P$ly right generated}\graffito{$P$ly right-gen\-er\-at\-ed cell space $\mathcal{R}$} if and only if there is a right-gen\-er\-at\-ing set of $\mathcal{R}$ that is $P$.
  \end{definition}

  \begin{remark}
    If $\mathcal{R}$ is finitely right generated, then $M$ is countably infinite.
  \end{remark}

  \begin{remark}
  \label{remark:independence-of-being-right-generated-of-coordinate-system}
    For each right-gen\-er\-at\-ing set $S$ of $\mathcal{R}$, each coordinate system $\mathcal{K}' = \ntuple{m_0', \family{g_{m_0', m}'}_{m \in M}}$ for $\mathcal{M}$, and each element $g \in G$ such that $g \actsOnPoint m_0 = m_0'$, according to \cref{lemma:liberation-and-coordinate-systems}, because $G_0 \cdot S \subseteq S$, the subset $g \conjugates S$ of $G \modulo G_{m_0'}$ is a right-gen\-er\-at\-ing set of $\ntuple{\mathcal{M}, \mathcal{K}'}$, which is finite or symmetric if $S$ has the respective property. In particular, being finitely and/or symmetrically right generated does not depend on the coordinate system.
  \end{remark}

  \begin{definition}
    The left-ho\-mo\-ge\-neous space $\mathcal{M}$ is called \graffito{finitely and/or symmetrically right-gen\-er\-at\-ed left-ho\-mo\-ge\-neous space}\define{finitely and/or symmetrically right generated} if and only if, there is a (for each) coordinate system $\mathcal{K}$ for $\mathcal{M}$, the cell space $\ntuple{\mathcal{M}, \mathcal{K}}$ is finitely and/or symmetrically right generated.
  \end{definition}

  \begin{lemma}
  \label{lemma:gen-set-of-group-induces-right-gen-set-of-action}
    Let $T$ be a generating set of $G$. The set $S = \setOf{g_0 \cdot t G_0 \suchThat g_0 \in G_0, t \in T}$ is a right-gen\-er\-at\-ing set of $\mathcal{R}$. And, if $T$ is symmetric, then so is $S$. And, if $T$ and $G_0$ are finite, then so is $S$.
  \end{lemma}

  \begin{proof}
    First, let $m \in M$. Then, because $\isSemiActedUponBy$ is transitive, there is a $g \in G$ such that $m_0 \isSemiActedUponBy g G_0 = m$. And, because $T$ generates $G$, there is a finite sequence $\sequence{t_i}_{i \in \setOf{1, 2, \dotsc, k}}$ in $T \cup T^{-1}$ such that $t_1 t_2 \dotsb t_k = g$. And, because $\isSemiActedUponBy$ is a semi-action, there is a finite sequence $\sequence{g_{i,0}}_{i \in \setOf{2,3,\dotsc,k}}$ in $G_0$ such that
    \begin{equation*}
      m_0 \isSemiActedUponBy \sequence{t_1 G_0, g_{2,0} t_2 G_0, \dotsc, g_{k,0} t_k G_0}
      = m_0 \isSemiActedUponBy t_1 t_2 \dotsb t_k G_0
      = m.
    \end{equation*}
    In conclusion, because $t_1 G_0 \in S \cup S^{-1}$ and $g_{i,0} t_i G_0 \in S \cup S^{-1}$, for $i \in \setOf{2,3,\dotsc,k}$, the set $S$ is a right-gen\-er\-at\-ing set of $\mathcal{R}$.

    Secondly, let $T$ be symmetric. Furthermore, let $s \in S$ and let $g \in s$. Then, there is a $g_0 \in G_0$, there is a $t \in T$, and there is a $g_0' \in G_0$ such that $g_0 \cdot t G_0 = s$ and $g_0 t g_0' = g$. Hence, because $(g_0')^{-1} \in G_0$ and $t^{-1} \in T$,
    \begin{equation*}
      g^{-1} G_0 = (g_0')^{-1} t^{-1} g_0^{-1} G_0 = (g_0')^{-1} \cdot t^{-1} G_0 \in S.
    \end{equation*}
    Therefore, $S^{-1} \subseteq S$. In conclusion, $S$ is symmetric.

    Lastly, let $T$ and $G_0$ be finite. Then, because $\cardinalityOf{S} \leq \cardinalityOf{G_0} \cdot \cardinalityOf{T}$, the set $S$ is finite.
  \end{proof}

  \section{Directed Multigraphs}
  \label{section:multigraphs}

  \begin{definition}
    Let $V$ and $E$ be two sets, and let $\headOf$ and $\tailOf$ be two maps from $E$ to $V$. The quadruple $\mathcal{G} = \ntuple{V, E, \headOf, \tailOf}$ is called \graffito{directed multigraph $\mathcal{G}$}\define{directed multigraph}\index{multigraph directed@directed multigraph}\index[symbols]{Gcalligraphic@$\mathcal{G}$}; each element $v \in V$ is called \define{vertex}\graffito{vertex $v$}\index[symbols]{v@$v$}; each element $e \in E$ is called \define{edge from $\headOf(e)$ to $\tailOf(e)$}\graffito{edge $e$ from $\headOf(e)$ to $\tailOf(e)$}\index[symbols]{e@$e$}; for each element $e \in E$, the vertex $\headOf(e)$ is called \define{head of $e$}\graffito{head $\headOf(e)$ of $e$}\index[symbols]{sigmae@$\headOf(e)$} and the vertex $\tailOf(e)$ is called \define{tail of $e$}\graffito{tail $\tailOf(e)$ of $e$}\index[symbols]{taue@$\tailOf(e)$}. 
  \end{definition}

  \begin{remark}
    In the case that the map $\eta \from E \to V \times V$, $e \mapsto (\headOf(e), \tailOf(e))$, is injective, the directed multigraph $\mathcal{G}$ is a directed graph and is identified with $\ntuple{V, \eta(E)}$.
  \end{remark}

  In the remainder of this section, let $\mathcal{G} = \ntuple{V, E, \headOf, \tailOf}$ be a directed multigraph.

  \begin{definition}
    Let $e$ be an edge of $\mathcal{G}$. The edge $e$ is called \define{loop}\graffito{loop $e$} if and only if $\tailOf(e) = \headOf(e)$.
  \end{definition}

  \begin{definition}
    Let $v$ be a vertex of $\mathcal{G}$.
    \begin{aenumerate}
      \item The cardinal number
            \begin{equation*}
              \degreeOf^+(v) = \cardinalityOf{\setOf{e \in E \suchThat \headOf(e) = v}}
              \mathnote{out-degree $\degreeOf^+(v)$ of $v$}
              \index[symbols]{degplusv@$\degreeOf^+(v)$}
            \end{equation*}
            is called \define{out-degree of $v$}\index{degree!out-}.
      \item The cardinal number
            \begin{equation*}
              \degreeOf^-(v) = \cardinalityOf{\setOf{e \in E \suchThat \tailOf(e) = v}}
              \mathnote{in-degree $\degreeOf^-(v)$ of $v$}
              \index[symbols]{degminusv@$\degreeOf^-(v)$}
            \end{equation*}
            is called \define{in-degree of $v$}\index{degree!in-}.
      \item The cardinal number
            \begin{equation*}
              \degreeOf(v) = \degreeOf^+(v) + \degreeOf^-(v)
              \mathnote{degree $\degreeOf(v)$ of $v$}
              \index[symbols]{degreev@$\degreeOf(v)$}
            \end{equation*}
            is called \defineX{degree of $v$}{degree@degree of $v$}. \qedhere
    \end{aenumerate}
  \end{definition}

  \begin{remark}
    In the degree of $v$ loops are counted twice.
  \end{remark}

  \begin{definition}
    Let $v$ and $v'$ be two vertices of $\mathcal{G}$. They are called \define{adjacent}\graffito{adjacent vertices $v$ and $v'$} if and only if there is an edge from $v$ to $v'$ or one from $v'$ to $v$.
  \end{definition}

  \begin{definition} 
    Let $p = \sequence{e_i}_{i \in \setOf{1, 2, \dotsc, k}}$ be a finite sequence of edges of $\mathcal{G}$. It is called \define{path from $\headOf(e_1)$ to $\tailOf(e_k)$}\graffito{path $p$ from $\headOf(e_1)$ to $\tailOf(e_k)$}\index[symbols]{p@$p$} if and only if, for each index $i \in \setOf{1, 2, \dotsc, k-1}$, we have $\tailOf(e_i) = \headOf(e_{i + 1})$.
  \end{definition}

  \begin{definition} 
    Let $p = \sequence{e_i}_{i \in \setOf{1, 2, \dotsc, k}}$ be a path in $\mathcal{G}$. The non-negative integer $\lengthOfPath{p} = k$ is called \define{length of $p$}\graffito{length $\lengthOfPath{p}$ of $p$}\index[symbols]{absp@$\lengthOfPath{p}$}.
  \end{definition}

  \begin{definition}
    The directed multigraph $\mathcal{G}$ is called
    \begin{aenumerate}
      \item \define{symmetric}\graffito{symmetric directed multigraph} if and only if, for each edge $e \in E$, there is an edge $e' \in E$ such that $\headOf(e') = \tailOf(e)$ and $\tailOf(e') = \headOf(e)$; 
      \item \define{strongly connected}\graffito{strongly connected directed multigraph} if and only if, for each vertex $v \in V$ and each vertex $v' \in V$, there is a path $p$ from $v$ to $v'$; 
      \item \define{regular}\graffito{regular directed multigraph} if and only if all vertices of $\mathcal{G}$ have the same degree and, for each vertex $v \in V$, we have $\degreeOf^-(v) = \degreeOf^+(v)$. \qedhere 
    \end{aenumerate}
  \end{definition}

  \begin{definition}
    Let $W$ be a subset of $V$ and let $F$ be a subset of $E$ such that $\headOf(F) \subseteq W$ and $\tailOf(F) \subseteq W$. The directed multigraph $\ntuple{W, F, \headOf\restrictedTo_{F \to W}, \tailOf\restrictedTo_{F \to W}}$ is called \define{submultigraph of $\mathcal{G}$}\graffito{submultigraph of $\mathcal{G}$}.
  \end{definition}

  \begin{definition}
    Let $W$ be a subset of $V$, let $F$ be the set $\setOf{e \in E \suchThat \headOf(e), \tailOf(e) \in W}$, let $\varsigma$ be the map $\headOf\restrictedTo_{F \to W}$, and let $\upsilon$ be the map $\tailOf\restrictedTo_{F \to W}$. The submultigraph $\mathcal{G}[W] = \ntuple{W, F, \varsigma, \upsilon}$ of $\mathcal{G}$ is called \graffito{submultigraph $\mathcal{G}[W]$ of $\mathcal{G}$ induced by $W$}\define{induced by $W$}\index[symbols]{GWcalligraphicbrackets@$\mathcal{G}[W]$}.
  \end{definition}

  \begin{definition} 
    Let $\mathcal{G} = \ntuple{V, E, \headOf, \tailOf}$ be symmetric and strongly connected. The map
    \begin{align*}
      \distanceOf \from V \times V &\to \N_0, \mathnote{distance $\distanceOf$ on $\mathcal{G}$}\index[symbols]{d@$\distanceOf$}\\
      (v, v') &\mapsto \min\setOf{\lengthOfPath{p} \suchThat \text{$p$ path from $v$ to $v'$}},
    \end{align*}
    is a metric on $V$ and called \define{distance on $\mathcal{G}$}.
  \end{definition}

  \begin{definition}
    Let $\Lambda$ be a set and let $\lambda$ be a map from $E$ to $\Lambda$. The map $\lambda$ is called \defineX{$\Lambda$-edge-labelling of $\mathcal{G}$}{edge-labelling of $\mathcal{G}$@$\Lambda$-edge-labelling of $\mathcal{G}$}\graffito{$\Lambda$-edge-labelling $\lambda$ of $\mathcal{G}$}\index[symbols]{lambda@$\lambda$} and the multigraph $\mathcal{G}$ equipped with $\lambda$ is called \defineX{$\Lambda$-edge-labelled}{edge-labelled@$\Lambda$-edge-labelled}\graffito{$\Lambda$-edge-labelled directed multigraph}. 
  \end{definition}

  \begin{remark}
    In the case that $\mathcal{G}$ is a directed graph, the $\Lambda$-edge-labelled graph $\mathcal{G}$ is identified with $\ntuple{V, \setOf{(\sigma(e), \lambda(e), \tau(e)) \suchThat e \in E}}$.
  \end{remark}

  \section{Coloured and Uncoloured Left Cayley Graphs}
  \label{section:Cayley}

  In this section, let $\mathcal{R} = \ntuple{\mathcal{M}, \mathcal{K}} = \ntuple{\ntuple{M, G, \actsOnPoint}, \ntuple{m_0, \family{g_{m_0, m}}_{m \in M}}}$ be a cell space and let $S$ be a right-gen\-er\-at\-ing set of $\mathcal{R}$. 

  \begin{definition}
    \begin{aenumerate}
      \item The directed graph
            \begin{equation*}
              \ntuple{M, \setOf{(m, m \isSemiActedUponBy s) \suchThat m \in M, s \in S}}
            \end{equation*}
            is called \define{uncoloured (left) $S$-Cayley graph of $\mathcal{R}$}\graffito{uncoloured left $S$-Cayley graph of $\mathcal{R}$}\index{left $S$-Cayley graph of $\mathcal{R}$!uncoloured}\index{Cayley graph of $\mathcal{R}$@left $S$-Cayley graph of $\mathcal{R}$}\index{Cayley graph of $\mathcal{R}$!uncoloured}.
      \item The $S$-edge-labelled directed graph
            \begin{equation*}
              \ntuple{M, \setOf{(m, s, m \isSemiActedUponBy s) \suchThat m \in M, s \in S}}
            \end{equation*}
            is called \define{coloured (left) $S$-Cayley graph of $\mathcal{R}$}\graffito{coloured left $S$-Cayley graph of $\mathcal{R}$}\index{left $S$-Cayley graph of $\mathcal{R}$!coloured}\index{Cayley graph of $\mathcal{R}$!coloured}. \qedhere
    \end{aenumerate}
  \end{definition}

  \begin{example}
    Examples of coloured Cayley graphs of groups that illustrate their dependence on the generating set are given in examples 6.3.2 in \cite{ceccherini-silberstein:coornaert:2010}. 
  \end{example}

  \begin{remark}
    Because the semi-action $\isSemiActedUponBy$ is free, the uncoloured and coloured $S$-Cayley graphs are isomorphic (if we ignore labels). Moreover, each automorphism of the coloured $S$-Cayley graph is an automorphism of the uncoloured one. However, because automorphisms of labelled graphs are label-preserving, the converse is not true in general.
  \end{remark}

  \begin{remark}
  \label{remark:independence-of-uncoloured-Cayley-of-coordinate-system}
    According to \cref{lemma:liberation-and-coordinate-systems}, because $G_0 \cdot S \subseteq S$, the uncoloured $S$-Cayley graph of $\mathcal{R}$ does not depend on the coordinates $\family{g_{m_0, m}}_{m \in M}$; and, for each element $g \in G$, it is equal to the uncoloured $(g \conjugates S)$-Cayley graph of $\ntuple{\mathcal{M}, \ntuple{g \actsOnPoint m_0, \family{g_{m_0, m} g^{-1}}_{m \in M}}}$.
  \end{remark}

  In the remainder of this section, let $\mathcal{G}$ be the uncoloured or coloured $S$-Cayley graph of $\mathcal{R}$.

  \begin{remark}
    Because the element $G_0 \in G \modulo G_0$ is the only one that acts trivially by $\isSemiActedUponBy$, the following three statements are equivalent:
    \begin{aenumerate}
      \item $G_0 \in S$;
      \item At least one vertex of $\mathcal{G}$ has a loop;
      \item All vertices of $\mathcal{G}$ have a loop. \qedhere
    \end{aenumerate}
  \end{remark}

  \begin{remark} 
    Because $\isSemiActedUponBy$ is a semi-action and $G_0 \cdot S \subseteq S$, if $S$ is symmetric, then $\mathcal{G}$ is symmetric and strongly connected.
  \end{remark} 

  \begin{remark}
    Let $m$ be a vertex of $\mathcal{G}$. The map $S \to m \isSemiActedUponBy S$, $s \mapsto m \isSemiActedUponBy s$, is a bijection onto the out-neighbourhood of $m$. It is injective, because $\isSemiActedUponBy$ is free, and it is surjective, by definition. Therefore, if $S$ is symmetric, then the degree of $m$ is $2 \cardinalityOf{S}$ in cardinal arithmetic and the graph $\mathcal{G}$ is regular. 
  \end{remark}

  \section{Metrics and Lengths}
  \label{section:metrics}

  In this section, let $\mathcal{R} = \ntuple{\mathcal{M}, \mathcal{K}} = \ntuple{\ntuple{M, G, \actsOnPoint}, \ntuple{m_0, \family{g_{m_0, m}}_{m \in M}}}$ be a cell space and let $S$ be a symmetric right-gen\-er\-at\-ing set of $\mathcal{R}$.

  \paragraph{Contents.} In \cref{definition:metric,definition:length} we introduce the $S$-metric $\distanceOf_S$ and the $S$-length $\lengthOf{\blank}_S$ on $\mathcal{R}$ induced by the uncoloured $S$-Cayley graph. And in \cref{lemma:metric-and-liberation,lemma:metric-invariant-under-left-action} we show how the $S$-metric relates to the left group action $\actsOnPoint$ and the right quotient set semi-action $\isSemiActedUponBy$.


  \begin{definition} 
  \label{definition:metric}
    The distance on the uncoloured $S$-Cayley graph of $\mathcal{R}$ is called \defineX{$S$-metric on $\mathcal{R}$}{metric on $\mathcal{R}$@$S$-metric on $\mathcal{R}$}\graffito{$S$-metric $\distanceOf_S$ on $\mathcal{R}$} and denoted by $\distanceOf_S$\index[symbols]{dS@$\distanceOf_S$}.
  \end{definition}

  \begin{remark}
    The $S$-metric on $\mathcal{R}$ is the map
    \begin{align*}
      \distanceOf_S \from M \times M &\to \N_0,\\ 
      (m, m') &\mapsto \min\{\begin{aligned}[t]
        k \in \N_0 \suchThat{} &\Exists \sequence{s_i}_{i \in \setOf{1, 2, \dotsc, k}} \text{ in } S \SuchThat\\
        &\;m \isSemiActedUponBy \sequence{s_i}_{i \in \setOf{1, 2, \dotsc, k}} = m'\}. \qedhere
      \end{aligned}
    \end{align*}
  \end{remark}

  \begin{remark}
  \label{remark:independence-of-uncoloured-cayley-of-distance}
    According to \cref{remark:independence-of-uncoloured-Cayley-of-coordinate-system}, the $S$-metric on $\mathcal{R}$ does not depend on the coordinates $\family{g_{m_0, m}}_{m \in M}$ and is identical to the $(g \conjugates S)$-metric on $\ntuple{\mathcal{M}, \ntuple{g \actsOnPoint m_0, \family{g_{m_0, m} g^{-1}}_{m \in M}}}$.
  \end{remark}

  \begin{lemma}
  \label{lemma:metric-and-liberation}
    Let $m$ and $m'$ be two elements of $M$, and let $s$ be an element of $S$. Then, $\distanceOf_S(m, m' \isSemiActedUponBy s) \leq \distanceOf_S(m, m') + 1$.
  \end{lemma}

  \begin{proof}
    Let $k = \distanceOf_S(m, m')$. Then, there is a finite sequence $\sequence{s_i}_{i \in \setOf{1, 2, \dotsc, k}}$ in $S$ such that $m \isSemiActedUponBy \sequence{s_i}_{i \in \setOf{1, 2, \dotsc, k}} = m'$. Hence, $m \isSemiActedUponBy \sequence{s_1, s_2, \dotsc, s_k, s} = m' \isSemiActedUponBy s$. Therefore, $\distanceOf_S(m, m' \isSemiActedUponBy s) \leq \distanceOf_S(m, m') + 1$.
  \end{proof}

  \begin{lemma} 
  \label{lemma:metric-invariant-under-left-action}
    Let $m$ and $m'$ be two elements of $M$, and let $g$ be an element of $G$. Then, $\distanceOf_S(g \actsOnPoint m, g \actsOnPoint m') = \distanceOf_S(m, m')$.
  \end{lemma}

  \begin{proof}
    Let $k = \distanceOf_S(g \actsOnPoint m, g \actsOnPoint m')$. Then, there is a finite sequence $\sequence{s_i}_{i \in \setOf{1, 2, \dotsc, k}}$ in $S$ such that $(g \actsOnPoint m) \isSemiActedUponBy \sequence{s_i}_{i \in \setOf{1, 2, \dotsc, k}} = g \actsOnPoint m'$. And, because $\isSemiActedUponBy$ semi-commutes with $\actsOnPoint$, there is a finite sequence $\sequence{g_{i,0}}_{i \in \setOf{1, 2, \dotsc, k}}$ in $G_0$ such that $(g \actsOnPoint m) \isSemiActedUponBy \sequence{s_i}_{i \in \setOf{1, 2, \dotsc, k}} = g \actsOnPoint (m \isSemiActedUponBy \sequence{g_{i,0} \cdot s_i}_{i \in \setOf{1, 2, \dotsc, k}})$.
    Hence, $g \actsOnPoint (m \isSemiActedUponBy \sequence{g_{i,0} \cdot s_i}_{i \in \setOf{1, 2, \dotsc, k}}) = m'$. Therefore, $\distanceOf_S(m, m') \leq k = \distanceOf_S(g \actsOnPoint m, g \actsOnPoint m')$. Analogously, $\distanceOf_S(g \actsOnPoint m, g \actsOnPoint m') \leq \distanceOf_S(g^{-1} \actsOnPoint (g \actsOnPoint m), g^{-1} \actsOnPoint (g \actsOnPoint m')) = \distanceOf_S(m, m')$. In conclusion, $\distanceOf_S(g \actsOnPoint m, g \actsOnPoint m') = \distanceOf_S(m, m')$. 
  \end{proof}

  \begin{lemma}
  \label{lemma:truncated-minimal-path-yields-minimal-path}
    Let $m$ and $m'$ be two elements of $M$, let $\sequence{s_i}_{i \in \setOf{1, 2, \dotsc, \distanceOf_S(m, m')}}$ be a finite sequence in $S$ such that $m' = m \isSemiActedUponBy \sequence{s_i}_{i \in \setOf{1, 2, \dotsc, \distanceOf_S(m, m')}}$, let $i$ be an element of $\setOf{0, 1, 2, \dotsc, \distanceOf_S(m, m')}$, and let $m_i = m \isSemiActedUponBy \sequence{s_i}_{i \in \setOf{1, 2, \dotsc, i}}$. Then, $\distanceOf_S(m, m_i) = i$.
  \end{lemma}

  \begin{proof}
    By definition of $m_i$, we have $\distanceOf_S(m, m_i) \leq i$ and $\distanceOf_S(m_i, m') \leq \distanceOf_S(m, m') - i$. Therefore, because $\distanceOf_S(m, m') \leq \distanceOf_S(m, m_i) + \distanceOf_S(m_i, m')$,
    \begin{align*}
      \distanceOf_S(m, m_i)
      &\geq \distanceOf_S(m, m') - \distanceOf_S(m_i, m')\\
      &\geq \distanceOf_S(m, m') - \parens[\big]{\distanceOf_S(m, m') - i}\\
      &= i
    \end{align*}
    In conclusion, $\distanceOf_S(m, m_i) = i$.
  \end{proof}

  \begin{definition}
  \label{definition:length}
    The map
    \begin{align*}
      \lengthOf{\blank}_S \from M &\to \N_0, \mathnote{$S$-length $\lengthOf{\blank}_S$ on $\mathcal{R}$}\\ 
      m &\mapsto \distanceOf_S(m_0, m),
    \end{align*}
    is called \defineX{$S$-length on $\mathcal{R}$}{length on $\mathcal{R}$@$S$-length on $\mathcal{R}$}. 
  \end{definition}

  \begin{remark}
    For each element $m \in M$, we have $\lengthOf{m}_S = 0$ if and only if $m = m_0$.
  \end{remark}


  \section{Balls and Spheres}
  \label{section:balls-and-spheres}

  In this section, let $\mathcal{R} = \ntuple{\mathcal{M}, \mathcal{K}} = \ntuple{\ntuple{M, G, \actsOnPoint}, \ntuple{m_0, \family{g_{m_0, m}}_{m \in M}}}$ be a cell space and let $S$ be a symmetric right-gen\-er\-at\-ing set of $\mathcal{R}$.

  \paragraph{Contents.} In \cref{definition:balls-and-spheres} we introduce balls $\ball_S$ and spheres $\sphere_S$ in the $S$-metric on $\mathcal{R}$. In \cref{subsection:acting-on-spheres-and-balls} we show how the left group action $\actsOnPoint$ and the right quotient set semi-action $\isSemiActedUponBy$ act on balls and spheres. And in \cref{subsection:distances-of-spheres-and-balls} we calculate and approximate the distances of balls and spheres in various circumstances.

  \begin{definition}
  \label{definition:balls-and-spheres}
    Let $m$ be an element of $M$ and let $\rho$ be an integer.
    \begin{aenumerate}
      \item The set
            \begin{equation*}
              \ball_S(m, \rho) = \setOf{m' \in M \suchThat \distanceOf_S(m, m') \leq \rho}
              \mathnote{$S$-ball $\ball_S(m, \rho)$ of radius $\rho$ centred at $m$}
              \index[symbols]{BSmrho@$\ball_S(m, \rho)$}
            \end{equation*}
            is called \defineX{$S$-ball of radius $\rho$ centred at $m$}{ball of radius $\rho$ centred at $m$@$S$-ball of radius $\rho$ centred at $m$}. The ball of radius $\rho$ centred at $m_0$ is denoted by $\ball_S(\rho)$\index[symbols]{BSrho@$\ball_S(\rho)$}.
      \item The set
            \begin{equation*}
              \sphere_S(m, \rho) = \setOf{m' \in M \suchThat \distanceOf_S(m, m') = \rho}
              \mathnote{$S$-sphere $\sphere_S(m, \rho)$ of radius $\rho$ centred at $m$}
              \index[symbols]{SSmrho@$\sphere_S(m, \rho)$}
            \end{equation*}
            is called \defineX{$S$-sphere of radius $\rho$ centred at $m$}{sphere of radius $\rho$ centred at $m$@$S$-sphere of radius $\rho$ centred at $m$}. The sphere of radius $\rho$ centred at $m_0$ is denoted by $\sphere_S(\rho)$\index[symbols]{SSrho@$\sphere_S(\rho)$}. \qedhere
    \end{aenumerate}
  \end{definition}

  \begin{remark}
  \label{remark:independence-of-uncoloured-cayley-of-balls-and-spheres}
    According to \cref{remark:independence-of-uncoloured-cayley-of-distance}, the $S$-balls and~spheres of $\mathcal{R}$ do not depend on the coordinates $\family{g_{m_0, m}}_{m \in M}$ and are identical to the $(g \conjugates S)$-balls and~spheres of $\ntuple{\mathcal{M}, \ntuple{g \actsOnPoint m_0, \family{g_{m_0, m} g^{-1}}_{m \in M}}}$.
  \end{remark}

  \begin{remark}
    For each negative integer $\rho$,
    \begin{equation*}
      \sphere_S(m, \rho) = \ball_S(m, \rho) = \emptyset.
    \end{equation*}
    And, $\sphere_S(m, 0) = \ball_S(m, 0) = \setOf{m}$.
  \end{remark}

  \begin{remark}
  \label{remark:spheres-expressed-in-terms-of-balls}
    For each integer $\rho$,
    \begin{equation*}
      \sphere_S(m, \rho) = \ball_S(m, \rho) \smallsetminus \ball_S(m, \rho - 1). \qedhere
    \end{equation*}
  \end{remark}

  \begin{remark}
  \label{remark:m-in-ball-about-mprime-if-and-only-if-converse-holds}
    Because the metric $\distanceOf_S$ is symmetric, for each integer $\rho$, each element $m \in M$, and each element $m' \in M$,
    \begin{equation*}
      m' \in \ball_S(m, \rho) \ifAndOnlyIf m \in \ball_S(m', \rho)
    \end{equation*}
    and
    \begin{equation*}
      m' \in \sphere_S(m, \rho) \ifAndOnlyIf m \in \sphere_S(m', \rho). \qedhere
    \end{equation*}
  \end{remark}

  \begin{remark}
    For each integer $\rho$,
    \begin{equation*}
      \ball_S(\rho) = \setOf{m \in M \suchThat \lengthOf{m}_S \leq \rho}
    \end{equation*}
    and
    \begin{equation*}
      \sphere_S(\rho) = \setOf{m \in M \suchThat \lengthOf{m}_S = \rho}. \qedhere
    \end{equation*}
  \end{remark}

  \begin{definition}
    Let $\sequence{A_k}_{k \in \N_0}$ be a sequence of sets. It is called
    \begin{aenumerate}
      \item \define{non-decreasing}\graffito{non-decreasing sequence}\index{decreasingnon@non-decreasing} if and only if
            \begin{equation*}
              \ForEach k \in \N_0 \Holds A_k \subseteq A_{k + 1}.
            \end{equation*}
      \item \define{non-increasing}\graffito{non-increasing sequence}\index{increasingnon@non-increasing} if and only if
            \begin{equation*}
              \ForEach k \in \N_0 \Holds A_{k + 1} \subseteq A_k. \qedhere
            \end{equation*}
    \end{aenumerate}
  \end{definition}

  \begin{definition} 
    Let $\sequence{A_k}_{k \in \N_0}$ be a sequence of sets.
    \begin{aenumerate}
      \item The set
            \begin{equation*}
              \liminf_{k \to \infty} A_k = \bigcup_{k \in \N_0} \bigcap_{\substack{j \in \N_0 \\ j \geq k}} A_j
              \mathnote{limit inferior $\liminf_{k \to \infty} A_k$ of $\sequence{A_k}_{k \in \N_0}$}
              \index[symbols]{liminfktoinfinityAk@$\liminf_{k \to \infty} A_k$}
            \end{equation*}
            is called \define{limit inferior of $\sequence{A_k}_{k \in \N_0}$}.
      \item The set
            \begin{equation*}
              \limsup_{k \to \infty} A_k = \bigcap_{k \in \N_0} \bigcup_{\substack{j \in \N_0 \\ j \geq k}} A_j
              \mathnote{limit superior $\limsup_{k \to \infty} A_k$ of $\sequence{A_k}_{k \in \N_0}$}
              \index[symbols]{limsupktoinfinityAk@$\limsup_{k \to \infty} A_k$}
            \end{equation*}
            is called \define{limit superior of $\sequence{A_k}_{k \in \N_0}$}.
      \item Let $A$ be a set. The sequence $\sequence{A_k}_{k \in \N_0}$ is said to \graffito{$\sequence{A_k}_{k \in \N_0}$ converges to $A$}\define{converge to $A$}, the set $A$ is called \define{limit set of $\sequence{A_k}_{k \in \N_0}$}\graffito{limit set $\lim_{k \to \infty} A_k$ of $\sequence{A_k}_{k \in \N_0}$}, and $A$ is denoted by $\lim_{k \to \infty} A_k$\index[symbols]{limktoinfinityAk@$\lim_{k \to \infty} A_k$} if and only if
            \begin{equation*}
              \liminf_{k \to \infty} A_k = \limsup_{k \to \infty} A_k = A.
            \end{equation*}
      \item The sequence $\sequence{A_k}_{k \in \N_0}$ is called \define{convergent}\graffito{convergent sequence $\sequence{A_k}_{k \in \N_0}$} if and only if
            \begin{equation*}
              \liminf_{k \to \infty} A_k = \limsup_{k \to \infty} A_k. \qedhere
            \end{equation*}
    \end{aenumerate}
  \end{definition}

  \begin{remark} 
    Let $\sequence{A_k}_{k \in \N_0}$ be a non-decreasing or non-increasing sequence of sets with respect to inclusion. It converges to $\bigcup_{k \in \N_0} A_k$ or $\bigcap_{k \in \N_0} A_k$ respectively.
  \end{remark}

  \begin{remark}
  \label{remark:ball-of-radius-0-contains-one-element-and-sequence-of-balls-is-monotonic}
    Let $m$ be an element of $M$. Then, $\ball_S(m, 0) = \setOf{m}$ and the sequence $\sequence{\ball_S(m, \rho)}_{\rho \in \N_0}$ is non-decreasing with respect to inclusion and converges to $M$. In particular, for each non-negative integer $\rho$, 
    \begin{equation*}
      \bigcup_{\substack{\rho' \in \N_0 \\ \rho' \geq \rho}} \ball_S(m, \rho') = M.
    \end{equation*}
    Indeed, for each $m \in M$, each $\rho \in \N_0$, and each $m' \in M$, we have $m' \in \ball_S(m, \rho')$, where $\rho' = \max\setOf{\rho, \distanceOf_S(m, m')}$.
  \end{remark}

  \begin{remark}
  \label{remark:upper-bound-for-cardinality-of-balls}
    For each element $m \in M$ and each integer $\rho$, in cardinal arithmetic,
    \begin{equation*}
      \cardinalityOf{\ball_S(m, \rho)} \leq (1 + \cardinalityOf{S})^\rho,
    \end{equation*}
    because the map
    \begin{align*}
      (\setOf{G_0} \cup S)^\rho &\to \ball_S(m, \rho),\\
      \sequence{s_i}_{i \in \setOf{1, 2, \dotsc, \rho}} &\mapsto m \isSemiActedUponBy \sequence{s_i}_{i \in \setOf{1, 2, \dotsc, \rho}},
    \end{align*}
    is surjective and $\cardinalityOf{\setOf{G_0} \cup S}^\rho \leq (1 + \cardinalityOf{S})^\rho$.
  \end{remark}

  \begin{lemma} 
  \label{lemma:finite-set-contained-in-ball-if-gen-set-contains-neutral-element}
    Let $A$ be a finite subset of $M$ and let $S'$ be the set $\setOf{G_0} \cup S$. Then, there is a non-negative integer $k$ such that
    \begin{equation*}
      A \subseteq \setOf{m \in M \suchThat \Exists \sequence{s_i'}_{i \in \setOf{1, 2, \dotsc, k}} \text{ in } S' \SuchThat m_0 \isSemiActedUponBy \sequence{s_i'}_{i \in \setOf{1, 2, \dotsc, k}} = m}. \qedhere
    \end{equation*}
  \end{lemma}

  \begin{proof}
    If $A$ is empty, then any $k \in \N_0$ works. Otherwise, let $k = \sup_{a \in A} \distanceOf_S(m_0, a)$. Then, because $A$ is non-empty and finite, we have $k \in \N_0$. And, by the choice of $k$, we have $A \subseteq \ball(m_0, k)$. And, because $G_0 \in S'$ and $\blank \isSemiActedUponBy G_0 = \identityMap_M$, we have $\ball(m_0, k) = \setOf{m \in M \suchThat \Exists \sequence{s_i'}_{i \in \setOf{1, 2, \dotsc, k}} \text{ in } S' \SuchThat m_0 \isSemiActedUponBy \sequence{s_i'}_{i \in \setOf{1, 2, \dotsc, k}} = m}$. In conclusion, the stated inclusion holds.
  \end{proof}

  \subsection{Acting on Balls and Spheres}
  \label{subsection:acting-on-spheres-and-balls}

  \begin{lemma}
  \label{lemma:ball-liberation-included-in-ball-one-larger}
    Let $m$ be an element of $M$, let $\rho$ be a non-negative integer, and let $s$ be an element of $S$. Then, $\ball_S(m, \rho) \isSemiActedUponBy s \subseteq \ball_S(m, \rho + 1)$.
  \end{lemma} 

  \begin{proof}
    Let $m' \in \ball_S(m, \rho) \isSemiActedUponBy s$. Then, there is an $m'' \in \ball_S(m, \rho)$ such that $m'' \isSemiActedUponBy s = m'$. Hence, according to \cref{lemma:metric-and-liberation}, we have $\distanceOf_S(m, m') = \distanceOf_S(m, m'' \isSemiActedUponBy s) \leq \distanceOf_S(m, m'') + 1 \leq \rho + 1$. Therefore, $m' \in \ball_S(m, \rho + 1)$. In conclusion, $\ball_S(m, \rho) \isSemiActedUponBy s \subseteq \ball_S(m, \rho + 1)$.
  \end{proof}

  %

  \begin{lemma}
  \label{lemma:left-action-and-balls}
    Let $m$ be an element of $M$, let $\rho$ be a non-negative integer, and let $g$ be an element of $G$. Then, $g \actsOnPoint \ball_S(m, \rho) = \ball_S(g \actsOnPoint m, \rho)$.
  \end{lemma}

  \begin{proof} 
      First, let $m' \in g \actsOnPoint \ball_S(m, \rho)$. Then, $g^{-1} \actsOnPoint m' \in \ball_S(m, \rho)$ and thus $\distanceOf_S(m, g^{-1} \actsOnPoint m') \leq \rho$. Hence, according to \cref{lemma:metric-invariant-under-left-action},
      \begin{align*}
        \distanceOf_S(g \actsOnPoint m, m') &= \distanceOf_S(g^{-1} \actsOnPoint (g \actsOnPoint m), g^{-1} \actsOnPoint m')\\
                                 &= \distanceOf_S(m, g^{-1} \actsOnPoint m')\\
                                 &\leq \rho.
      \end{align*}
      Therefore, $m' \in \ball_S(g \actsOnPoint m, \rho)$. In conclusion, $g \actsOnPoint \ball_S(m, \rho) \subseteq \ball_S(g \actsOnPoint m, \rho)$.

      Secondly, let $m' \in \ball_S(g \actsOnPoint m, \rho)$. Then, $\distanceOf_S(g \actsOnPoint m, m') \leq \rho$. Thus, according to \cref{lemma:metric-invariant-under-left-action},
      \begin{align*}
        \distanceOf_S(m, g^{-1} \actsOnPoint m') &= \distanceOf_S(g \actsOnPoint m, g \actsOnPoint (g^{-1} \actsOnPoint m'))\\
                                      &= \distanceOf_S(g \actsOnPoint m, m')\\
                                      &\leq \rho.
      \end{align*}
      Hence, $g^{-1} \actsOnPoint m' \in \ball_S(m, \rho)$. Therefore, $m' \in g \actsOnPoint \ball_S(m, \rho)$. In conclusion, $\ball_S(g \actsOnPoint m, \rho) \subseteq g \actsOnPoint \ball_S(m, \rho)$.
  \end{proof}

  \begin{corollary}
    Let $m$ be an element of $M$, let $\rho$ be a non-negative integer, and let $g_m$ be an element of $G_m$. Then, $g_m \actsOnPoint \ball_S(m, \rho) = \ball_S(m, \rho)$. In particular, $G_m \actsOnPoint \ball_S(m, \rho) = \ball_S(m, \rho)$.
  \end{corollary}

  \begin{proof}
    This is a direct consequence of \cref{lemma:left-action-and-balls}, because $g_m \actsOnPoint m = m$.
  %
  \end{proof}

  \begin{corollary}
  \label{corollary:balls-of-equal-radius-have-same-number-of-elements}
    Let $m$ and $m'$ be two elements of $M$, and let $\rho$ be a non-negative integer. Then, $\cardinalityOf{\ball_S(m, \rho)} = \cardinalityOf{\ball_S(m', \rho)}$.
  \end{corollary}

  \begin{proof}
    This is a direct consequence of \cref{lemma:left-action-and-balls}, because there is a $g \in G$ such that $g \actsOnPoint m = m'$, and $g \actsOnPoint \blank$ is injective
  \end{proof}

  \begin{remark}
  \label{remark:liberation-under-identification-of-quotient-set-with-M}
    Recall from \cref{lemma:identification-of-G-quotient-Gzero-with-M-by-right-semiaction} that, under the identification of $M$ with $G \modulo G_0$ by $\iota \givenBy m \mapsto G_{m_0, m}$,
    \begin{equation*}
      \ForEach g \in G \ForEach m \in M \Holds g \cdot m = g \actsOnPoint m,
    \end{equation*}
    and
    \begin{equation*}
      \ForEach m \in M \ForEach m' \in M \Holds m \isSemiActedUponBy m' = g_{m_0, m} \actsOnPoint m'. \qedhere
    \end{equation*}
  \end{remark}

  \begin{lemma}
  \label{lemma:almost-associativity-of-liberation-under-identification}
    Let $m$, $m'$, and $m''$ be three elements of $M$ and identify $M$ with $G \modulo G_0$ by $\iota \givenBy m \mapsto G_{m_0, m}$. Then, there is an element $g_0 \in G_0$ such that $(m \isSemiActedUponBy m') \isSemiActedUponBy m'' = m \isSemiActedUponBy (m' \isSemiActedUponBy (g_0 \actsOnPoint m''))$.
  \end{lemma}

  \begin{proof}
    Because $\isSemiActedUponBy$ is a right semi-action, there is an element $g_0 \in G_0$ such that $m \isSemiActedUponBy g_{m_0, m'} \cdot g_0 \cdot G_{m_0, m''} = (m \isSemiActedUponBy g_{m_0, m'} G_0) \isSemiActedUponBy G_{m_0, m''}$. And, according to \cref{remark:liberation-under-identification-of-quotient-set-with-M}, we have $G_{m_0, m''} = m''$, $g_{m_0, m'} G_0 = G_{m_0, m'} = m'$, and $g_{m_0, m'} \cdot g_0 \cdot G_{m_0, m''} = m' \isSemiActedUponBy (g_0 \actsOnPoint m'')$. Therefore, $m \isSemiActedUponBy (m' \isSemiActedUponBy (g_0 \actsOnPoint m'')) = (m \isSemiActedUponBy m') \isSemiActedUponBy m''$.
  \end{proof}

  \begin{corollary}
  \label{corollary:ball-centred-at-mzero-to-m}
    Let $m$ be an element of $M$, let $\rho$ be a non-negative integer, and identify $M$ with $G \modulo G_0$ by $\iota \givenBy m \mapsto G_{m_0, m}$. Then, $m \isSemiActedUponBy \ball_S(\rho) = \ball_S(m, \rho)$.
  \end{corollary}

  \begin{proof}
    According to \cref{remark:liberation-under-identification-of-quotient-set-with-M} and \cref{lemma:left-action-and-balls},
    \begin{align*}
      m \isSemiActedUponBy \ball_S(\rho)
      &= g_{m_0, m} \actsOnPoint \ball_S(\rho)\\
      &= \ball_S(g_{m_0, m} \actsOnPoint m_0, \rho)\\
      &= \ball_S(m, \rho). \qedhere
    \end{align*}
  \end{proof}

  \begin{corollary}
  \label{corollary:liberation-of-balls-yields-bigger-one}
    Let $m$ be an element of $M$, let $\rho$ and $\rho'$ be two non-negative integers, and identify $M$ with $G \modulo G_0$ by $\iota \givenBy m \mapsto G_{m_0, m}$. Then, $\ball_S(m, \rho) \isSemiActedUponBy \ball_S(\rho') = \ball_S(m, \rho + \rho')$.
  \end{corollary}

  \begin{proof}
    First, let $m' \in \ball_S(m, \rho) \isSemiActedUponBy \ball_S(\rho')$. Then, there is an $m'' \in \ball_S(m, \rho)$ such that $m' \in m'' \isSemiActedUponBy \ball_S(\rho')$. And, according to \cref{corollary:ball-centred-at-mzero-to-m}, we have $m'' \isSemiActedUponBy \ball_S(\rho') = \ball_S(m'', \rho')$. Hence, because $\distanceOf_S$ is subadditive, we have $\distanceOf_S(m, m') \leq \distanceOf_S(m, m'') + \distanceOf_S(m'', m') \leq \rho + \rho'$. Therefore, $m' \in \ball_S(m, \rho + \rho')$. In conclusion, $\ball_S(m, \rho) \isSemiActedUponBy \ball_S(\rho') \subseteq \ball_S(m, \rho + \rho')$.

    Secondly, let $m' \in \ball_S(m, \rho + \rho')$.
    \begin{description}
      \item[Case 1:] $m' \in \ball_S(m, \rho)$. Then, because $m_0 \in \ball_S(\rho')$, we have $m' = m' \isSemiActedUponBy m_0 \in \ball_S(m, \rho) \isSemiActedUponBy \ball_S(\rho')$.
      \item[Case 2:] $m' \notin \ball_S(m, \rho)$. Then, there is a $j \in \setOf{1, 2, \dotsc, \rho'}$ and there is a finite sequence $\sequence{s_i}_{i \in \setOf{1, 2, \dotsc, \rho + j}}$ in $S$ such that $m'' \isSemiActedUponBy \sequence{s_i}_{i \in \setOf{\rho + 1, \rho + 2, \dotsc, \rho + j}} = m'$, where $m'' = m \isSemiActedUponBy \sequence{s_i}_{i \in \setOf{1, 2, \dotsc, \rho}} \in \ball_S(m, \rho)$. Hence, $m' \in \ball_S(m'', j) \subseteq \ball_S(m'', \rho') = m'' \isSemiActedUponBy \ball_S(\rho') \subseteq \ball_S(m, \rho) \isSemiActedUponBy \ball_S(\rho')$.
    \end{description}
    In either case, $m' \in \ball_S(m, \rho) \isSemiActedUponBy \ball_S(\rho')$. In conclusion, $\ball_S(m, \rho + \rho') \subseteq \ball_S(m, \rho) \isSemiActedUponBy \ball_S(\rho')$.
  \end{proof}

  \subsection{Distances of Balls and Spheres}
  \label{subsection:distances-of-spheres-and-balls}

  \begin{definition}
    Let $A$ and $A'$ be two subsets of $M$. The non-negative number or infinity
    \begin{equation*}
      \distanceOf_S(A, A') = \min\setOf{\distanceOf_S(a, a') \suchThat a \in A, a' \in A'}
      \mathnote{distance $\distanceOf_S(A, A')$ of $A$ and $A'$}
      \index[symbols]{dSAAprime@$\distanceOf_S(A, A')$}
    \end{equation*}
    is called \define{distance of $A$ and $A'$}, where we put $\min\emptyset = \infty$. In the case that $A = \setOf{a}$, we write $\distanceOf_S(a, A')$ in place of $\distanceOf_S(\setOf{a}, A')$; and in the case that $A' = \setOf{a'}$, we write $\distanceOf_S(A, a')$ in place of $\distanceOf_S(A, \setOf{a'})$.
  \end{definition}

  \begin{lemma}
  \label{lemma:distance-of-sphere-and-point}
    Let $m$ and $m'$ be two elements of $M$, and let $\rho$ be a non-negative integer such that $\rho \leq \distanceOf_S(m, m')$. Then, $\distanceOf_S(\sphere_S(m, \rho), m') = \distanceOf_S(m, m') - \rho$.
  \end{lemma}

  \begin{proof}
    Let $\rho' = \distanceOf_S(m, m')$. Then, there is a finite sequence $\sequence{s_i}_{i \in \setOf{1,2,\dotsc,\rho'}}$ in $S$ such that $m \isSemiActedUponBy \sequence{s_i}_{i \in \setOf{1,2,\dotsc,\rho'}} = m'$. Let $m'' = m \isSemiActedUponBy \sequence{s_i}_{i \in \setOf{1,2,\dotsc,\rho}}$. Then, $m'' \isSemiActedUponBy \sequence{s_i}_{i \in \setOf{\rho + 1, \rho + 2, \dotsc, \rho'}} = m'$. And, according to \cref{lemma:truncated-minimal-path-yields-minimal-path}, we have $m'' \in \sphere_S(m, \rho)$. Thus, $\distanceOf_S(\sphere_S(m, \rho), m') \leq \distanceOf_S(m'', m') \leq \rho' - \rho$.

    Suppose that $\distanceOf_S(\sphere_S(m, \rho), m') < \rho' - \rho$. Then, there is an $m'' \in \sphere_S(m, \rho)$ such that $\distanceOf_S(m'', m') < \rho' - \rho$. Hence, $\distanceOf_S(m, m') \leq \distanceOf_S(m, m'') + \distanceOf_S(m'', m') < \rho + (\rho' - \rho) = \rho'$, which contradicts $\distanceOf_S(m, m') = \rho'$. Therefore, $\distanceOf_S(\sphere_S(m, \rho), m') \geq \rho' - \rho$. In conclusion, $\distanceOf_S(\sphere_S(m, \rho), m') = \rho' - \rho = \distanceOf_S(m, m') - \rho$.
  \end{proof}

  \begin{corollary}
  \label{corollary:distance-of-spheres}
    Let $m$ be an element of $M$, and let $\rho$ and $\rho'$ be two non-negative integers such that the spheres $\sphere_S(m, \rho)$ and $\sphere_S(m, \rho')$ are non-empty. Then, $\distanceOf_S(\sphere_S(m, \rho), \sphere_S(m, \rho')) = \absoluteValueOf{\rho - \rho'}$.
  \end{corollary}

  \begin{proof}
    Without loss of generality, let $\rho \leq \rho'$. Then, for each $m' \in \sphere_S(m, \rho')$, according to \cref{lemma:distance-of-sphere-and-point}, we have $\distanceOf_S(\sphere_S(m, \rho), m') = \rho' - \rho$. In conclusion, $\distanceOf_S(\sphere_S(m, \rho), \sphere_S(m, \rho')) = \rho' - \rho = \absoluteValueOf{\rho - \rho'}$.
  \end{proof}

  \begin{corollary} 
    Let $m$ and $m'$ be two elements of $M$, and let $\rho$ be a non-negative integer. Then, $\distanceOf_S(\sphere_S(m, \rho), m') \geq \absoluteValueOf{\distanceOf_S(m, m') - \rho}$.
  \end{corollary}

  \begin{proof}
    If $\sphere_S(m, \rho) = \emptyset$, then $\distanceOf_S(\sphere_S(m, \rho), m') = \infty \geq \absoluteValueOf{\distanceOf_S(m, m') - \rho}$. Otherwise, let $\rho' = \distanceOf_S(m, m')$. Then, $m' \in \sphere_S(m, \rho') \neq \emptyset$. Hence, according to \cref{corollary:distance-of-spheres},
    \begin{align*}
      \distanceOf_S(\sphere_S(m, \rho), m')
      &\geq \distanceOf_S(\sphere_S(m, \rho), \sphere_S(m, \rho'))\\
      &= \absoluteValueOf{\rho - \rho'}\\
      &= \absoluteValueOf{\distanceOf_S(m, m') - \rho}. \qedhere
    \end{align*}
  \end{proof}

  \begin{lemma}
  \label{lemma:distance-of-balls}
    Let $m$ and $m'$ be two elements of $M$, and let $\rho$ and $\rho'$ be two non-negative integers such that $\rho + \rho' \leq \distanceOf_S(m, m')$. Then, $\distanceOf_S(\ball_S(m, \rho), \ball_S(m', \rho')) = \distanceOf_S(m, m') - (\rho + \rho')$.
  \end{lemma}

  \begin{proof}
    First, for each $m_\rho \in \ball_S(m, \rho)$ and each $m'_{\rho'} \in \ball_S(m', \rho')$, because $\distanceOf_S$ is subadditive,
    \begin{align*}
      \distanceOf_S(m, m') &\leq \distanceOf_S(m, m_\rho) + \distanceOf_S(m_\rho, m'_{\rho'}) + \distanceOf_S(m'_{\rho'}, m')\\
                 &\leq \rho + \distanceOf_S(m_\rho, m'_{\rho'}) + \rho',
    \end{align*}
    and hence $\distanceOf_S(m_\rho, m'_{\rho'}) \geq \distanceOf_S(m, m') - (\rho + \rho')$. In conclusion,
    \begin{equation*}
      \distanceOf_S(\ball_S(m, \rho), \ball_S(m', \rho')) \geq \distanceOf_S(m, m') - (\rho + \rho').
    \end{equation*}

    Secondly, there is a finite sequence $\sequence{s_i}_{i \in \setOf{1,2,\dotsc,\distanceOf_S(m, m')}}$ in $S$ such that $m \isSemiActedUponBy \sequence{s_i}_{i \in \setOf{1,2,\dotsc,\distanceOf_S(m, m')}} = m'$. Let $m_\rho = m \isSemiActedUponBy \sequence{s_i}_{i \in \setOf{1,2,\dotsc,\rho}}$ and let $m'_{\rho'} = m_\rho \isSemiActedUponBy \sequence{s_i}_{i \in \setOf{\rho + 1, \rho + 2, \dotsc, \distanceOf_S(m, m') - \rho'}}$. Then, $m_\rho \in \ball_S(m, \rho)$. And, because
    \begin{equation*}
      m'_{\rho'} \isSemiActedUponBy \sequence{s_i}_{i \in \setOf{\distanceOf_S(m, m') - \rho' + 1, \distanceOf_S(m, m') - \rho' + 2, \dotsc, \distanceOf_S(m, m')}} = m',
    \end{equation*}
    we have $m' \in \ball_S(m'_{\rho'}, \rho')$ and hence $m'_{\rho'} \in \ball_S(m', \rho')$. And, $\distanceOf_S(m_\rho, m'_{\rho'}) \leq \distanceOf_S(m, m') - \rho' - \rho$. In conclusion,
    \begin{equation*}
      \distanceOf_S(\ball_S(m, \rho), \ball_S(m', \rho')) \leq \distanceOf_S(m, m') - (\rho + \rho'). \qedhere
    \end{equation*}
  \end{proof}

  \section{Interiors, Closures, and Boundaries}
  \label{section:interiors-closures-boundaries}

  In this section, let $\mathcal{R} = \ntuple{\ntuple{M, G, \actsOnPoint}, \ntuple{m_0, \family{g_{m_0, m}}_{m \in M}}}$ be a cell space, let $S$ be a symmetric right-gen\-er\-at\-ing set of $\mathcal{R}$, let us omit the subscript $S$ of $\distanceOf_S$, $\lengthOf{\blank}_S$, $\ball_S$, and $\sphere_S$, and let $M$ be identified with $G \modulo G_0$ by $\iota \givenBy m \mapsto G_{m_0, m}$.

  \paragraph{Contents.} In \cref{definition:k-interior-closure-and-boundary} we introduce $\theta$-interiors $A^{-\theta}$, $\theta$-clo\-sures $A^{+\theta}$, and (internal/external) $\theta$-boundaries $\boundaryOf_\theta A$, $\boundaryOf_\theta^- A$, or $\boundaryOf_\theta^+ A$. And in the lemmata and corollaries of this section we characterise them and show how they and the $S$-metric relate to each other.

  \begin{definition}
  \label{definition:k-interior-closure-and-boundary}
    Let $A$ be a subset of $M$, let $\theta$ be an integer.
    \begin{aenumerate}
      \item The set
            \begin{equation*}
              A^{-\theta} = A^{-\ball(\theta)} \quad \big(= \setOf{m \in M \suchThat m \isSemiActedUponBy \ball(\theta) \subseteq A}\big)
              \mathnote{$\theta$-interior $A^{-\theta}$ of $A$}
              \index[symbols]{A-theta@$A^{-\theta}$}
            \end{equation*}
            is called \defineX{$\theta$-interior of $A$}{interior of $A$ theta@$\theta$-interior of $A$}.
      \item The set
            \begin{equation*}
              A^{+\theta} = A^{+\ball(\theta)} \quad \big(= \setOf{m \in M \suchThat (m \isSemiActedUponBy \ball(\theta)) \cap A \neq \emptyset}\big)
              \mathnote{$\theta$-closure $A^{+\theta}$ of $A$}
              \index[symbols]{A+theta@$A^{+\theta}$}
            \end{equation*}
            is called \defineX{$\theta$-closure of $A$}{closure of $A$ theta@$\theta$-closure of $A$}.
      \item The set
            \begin{equation*}
              \boundaryOf_\theta A = A^{+\theta} \smallsetminus A^{-\theta} \quad \big(= A^{+\ball(\theta)} \smallsetminus A^{-\ball(\theta)} = \boundaryOf_{\ball(\theta)} A\big)
              \mathnote{$\theta$-boundary $\boundaryOf_\theta A$ of $A$}
              \index[symbols]{partialthetaA@$\boundaryOf_\theta A$}
            \end{equation*}
            is called \defineX{$\theta$-boundary of $A$}{boundary of $A$ theta@$\theta$-boundary of $A$}.
      \item The set
            \begin{equation*}
              \boundaryOf_\theta^- A = A \smallsetminus A^{-\theta} \quad \big(= A \smallsetminus A^{-\ball(\theta)} = \boundaryOf_{\ball(\theta)}^- A\big)
              \mathnote{internal $\theta$-boundary $\boundaryOf_\theta^- A$ of $A$}
              \index[symbols]{partialtheta-A@$\boundaryOf_\theta^- A$}
            \end{equation*}
            is called \define{internal $\theta$-boundary of $A$}\index{boundary of $A$ theta!internal}.
      \item The set
            \begin{equation*}
              \boundaryOf_\theta^+ A = A^{+\theta} \smallsetminus A \quad \big(= A^{+\ball(\theta)} \smallsetminus A = \boundaryOf_{\ball(\theta)}^+ A\big)
              \mathnote{external $\theta$-boundary $\boundaryOf_\theta^+ A$ of $A$}
              \index[symbols]{partialtheta+A@$\boundaryOf_\theta^+ A$}
            \end{equation*}
            is called \define{external $\theta$-boundary of $A$}\index{boundary of $A$ theta!external}. \qedhere
    \end{aenumerate}
  \end{definition}

  \begin{remark}
    According to \cref{remark:independence-of-uncoloured-cayley-of-balls-and-spheres}, the above notions in $\mathcal{R}$ do not depend on the coordinates $\family{g_{m_0, m}}_{m \in M}$ and are identical to the respective notions with respect to the symmetric right-gen\-er\-at\-ing set $g \conjugates S$ in $\ntuple{\mathcal{M}, \ntuple{g \actsOnPoint m_0, \family{g_{m_0, m} g^{-1}}_{m \in M}}}$.
  \end{remark}

  \begin{remark}
    For each negative integer $\theta$, we have $A^{-\theta} = A$, and $A^{+\theta} = \emptyset$, and $\boundaryOf_\theta A = \boundaryOf_\theta^- A = \boundaryOf_\theta^+ A = \emptyset$. Moreover, $A^{-0} = A^{+0} = A$ and $\boundaryOf_0 A = \boundaryOf_0^- A = \boundaryOf_0^+ A = \emptyset$. Furthermore, for each non-negative integer $\theta$, we have $M^{-\theta} = M^{+\theta} = M$ and $\boundaryOf_\theta M = \boundaryOf_\theta^- M = \boundaryOf_\theta^+ M = \emptyset$.
  \end{remark}

  \begin{lemma}
  \label{lemma:characterisation-of-k-closure-and-interior}
    Let $A$ be a subset of $M$ and let $\theta$ be a non-negative integer. Then,
    \begin{aenumerate}
      \item\label{item:characterisation-of-k-closure-and-interior:interior}
            $\displaystyle A^{-\theta} = \setOf{m \in A \suchThat \ball(m, \theta) \subseteq A}$;
      \item\label{item:characterisation-of-k-closure-and-interior:closure}
            $\displaystyle A^{+\theta} = \bigcup_{m \in A} \ball(m, \theta) = A \isSemiActedUponBy \ball(\theta)$. \qedhere
    \end{aenumerate}
  \end{lemma}

  \begin{proof}
    \begin{aenumerate}
      \item According to \cref{corollary:ball-centred-at-mzero-to-m} and because $m \in \ball(m, \theta)$, for $m \in M$,
            \begin{align*}
              A^{-\theta}
              &= \setOf{m \in M \suchThat \ball(m, \theta) \subseteq A}\\
              &= \setOf{m \in A \suchThat \ball(m, \theta) \subseteq A}.
            \end{align*}
      \item According to \cref{corollary:ball-centred-at-mzero-to-m} and \cref{remark:m-in-ball-about-mprime-if-and-only-if-converse-holds},
            \begin{align*}
              A^{+\theta}
              &= \setOf{m \in M \suchThat \ball(m, \theta) \cap A \neq \emptyset}\\
              &= \setOf{m \in M \suchThat \Exists m' \in A \SuchThat m' \in \ball(m, \theta)}\\
              &= \setOf{m \in M \suchThat \Exists m' \in A \SuchThat m \in \ball(m', \theta)}\\
              &= \bigcup_{m' \in A} \ball(m', \theta)\\
              &= \bigcup_{m' \in A} m' \isSemiActedUponBy \ball(\theta)\\
              &= A \isSemiActedUponBy \ball(\theta). \qedhere
            \end{align*}
    \end{aenumerate}
  \end{proof}

  \begin{corollary}
  \label{corollary:characterisation-of-k-closure-and-interior-of-balls}
    Let $m$ be an element of $M$, let $\rho$ be a non-negative integer, and let $\theta$ be a non-negative integer. Then,
    \begin{aenumerate}
      \item\label{item:characterisation-of-k-closure-and-interior-of-balls:interior}
            $\displaystyle \ball(m, \rho)^{-\theta} \supseteq \ball(m, \rho - \theta)$;
      \item\label{item:characterisation-of-k-closure-and-interior-of-balls:closure}
            $\displaystyle \ball(m, \rho)^{+\theta} = \ball(m, \rho + \theta)$;
      \item\label{item:characterisation-of-k-closure-and-interior-of-balls:boundary}
            $\displaystyle \boundaryOf_\theta \ball(m, \rho) \subseteq \ball(m, \rho + \theta) \smallsetminus \ball(m, \rho - \theta)$. \qedhere
    \end{aenumerate}
  \end{corollary}

  \begin{proof}
    \begin{aenumerate}
      \item If $\rho < \theta$, then $\ball(m, \rho - \theta) = \emptyset$ and hence $\ball(m, \rho - \theta) \subseteq \ball(m, \rho)^{-\theta}$. Otherwise, according to \cref{corollary:liberation-of-balls-yields-bigger-one}, we have $\ball(m, \rho - \theta) \isSemiActedUponBy \ball(\theta) \subseteq \ball(m, \rho)$ and hence, according to \cref{definition:k-interior-closure-and-boundary}, we have $\ball(m, \rho - \theta) \subseteq \ball(m, \rho)^{-\theta}$. 
      \item According to \cref{item:characterisation-of-k-closure-and-interior:closure} of \cref{lemma:characterisation-of-k-closure-and-interior} and \cref{corollary:liberation-of-balls-yields-bigger-one}, we have $\ball(m, \rho)^{+\theta} = \ball(m, \rho) \isSemiActedUponBy \ball(\theta) = \ball(m, \rho + \theta)$.
      \item This is a direct consequence of \cref{item:characterisation-of-k-closure-and-interior-of-balls:interior,item:characterisation-of-k-closure-and-interior-of-balls:closure}. \qedhere
    \end{aenumerate}
  \end{proof}

  \begin{remark} 
    Let $M$ be finite, let $\rho$ be the least integer such that $\ball(m, \rho) = M$, which is non-negative because $M \neq \emptyset$, and let $\theta$ be a positive integer. Then, $\ball(m, \rho)^{-\theta} = M \neq \ball(m, \rho - \theta)$.
  \end{remark}

  \begin{lemma} 
  \label{lemma:repeated-k-boundaries-etc}
    Let $A$ be a subset of $M$, and let $\theta$ and $\theta'$ be two non-negative integers. The following five statements hold:
    \begin{aenumerate}
      \item \label{item:repeated-k-boundaries-etc:interior}
            $\displaystyle (A^{-\theta})^{-\theta'} = A^{-(\theta + \theta')}$;
      \item $\displaystyle \boundaryOf_{\theta'}^- A^{-\theta} = A^{-\theta} \smallsetminus A^{-(\theta + \theta')}$;
      \item \label{item:repeated-k-boundaries-etc:closure}
            $\displaystyle (A^{+\theta})^{+\theta'} = A^{+(\theta + \theta')}$;
      \item \label{item:repeated-k-boundaries-etc:external-boundary}
            $\displaystyle \boundaryOf_{\theta'}^+ A^{+\theta} = A^{+(\theta + \theta')} \smallsetminus A^{+\theta}$;
      \item \label{item:repeated-k-boundaries-etc:closure-interior}
            If $\theta' \leq \theta$, then $A^{+(\theta - \theta')} \subseteq (A^{+\theta})^{-\theta'}$ and $(A^{-\theta})^{+\theta'} \subseteq A^{-(\theta - \theta')}$. \qedhere
    \end{aenumerate}
  \end{lemma}

  \begin{proof}
    \begin{aenumerate}
      \item For each $m' \in A$, according to \cref{corollary:ball-centred-at-mzero-to-m} and \cref{item:characterisation-of-k-closure-and-interior:interior} of \cref{lemma:characterisation-of-k-closure-and-interior}, we have $m' \in A^{-\theta}$ if and only if $m' \isSemiActedUponBy \ball(\theta) = \ball(m', \theta) \subseteq A$. Hence, according to \cref{corollary:liberation-of-balls-yields-bigger-one},
            \begin{align*}
              (A^{-\theta})^{-\theta'}
              &= \setOf{m' \in A \suchThat \ball(m', \theta') \subseteq A^{-\theta}}\\
              &= \setOf{m' \in A \suchThat \ball(m', \theta') \isSemiActedUponBy \ball(\theta) \subseteq A}\\
              &= \setOf{m' \in A \suchThat \ball(m', \theta + \theta') \subseteq A}\\
              &= A^{-(\theta + \theta')}.
            \end{align*}
      \item This is a direct consequence of \cref{item:repeated-k-boundaries-etc:interior}.
      \item According to \cref{item:characterisation-of-k-closure-and-interior:closure} of \cref{lemma:characterisation-of-k-closure-and-interior} and \cref{corollary:liberation-of-balls-yields-bigger-one},
            \begin{align*}
              (A^{+\theta})^{+\theta'}
              &= A^{+\theta} \isSemiActedUponBy \ball(\theta')\\
              &= \parens*{\bigcup_{m \in A} \ball(m, \theta)} \isSemiActedUponBy \ball(\theta')\\
              &= \bigcup_{m \in A} \ball(m, \theta) \isSemiActedUponBy \ball(\theta')\\
              &= \bigcup_{m \in A} \ball(m, \theta + \theta')\\
              &= A^{+(\theta + \theta')}. 
            \end{align*}
      \item This is a direct consequence of \cref{item:repeated-k-boundaries-etc:closure}.
      \item Let $\theta' \leq \theta$. Then, according to \cref{item:characterisation-of-k-closure-and-interior:closure} of \cref{lemma:characterisation-of-k-closure-and-interior} and \cref{item:repeated-k-boundaries-etc:closure}, 
            \begin{align*}
              A^{+(\theta - \theta')} \isSemiActedUponBy \ball(\theta')
              &= (A^{+(\theta - \theta')})^{+\theta'}\\
              &= A^{+((\theta - \theta') + \theta')}\\
              &= A^{+\theta}.
            \end{align*}
            Therefore, according to \cref{definition:k-interior-closure-and-boundary}, we have $A^{+(\theta - \theta')} \subseteq (A^{+\theta})^{-\theta'}$.

            Moreover, according to \cref{item:repeated-k-boundaries-etc:closure}, \cref{item:characterisation-of-k-closure-and-interior:closure} of \cref{lemma:characterisation-of-k-closure-and-interior}, and \cref{definition:k-interior-closure-and-boundary}, 
            \begin{align*}
              (A^{-\theta})^{+\theta'} \isSemiActedUponBy \ball(\theta - \theta')
              &= ((A^{-\theta})^{+\theta'})^{+(\theta - \theta')}\\
              &= (A^{-\theta})^{+\theta' + (\theta - \theta')}\\
              &= (A^{-\theta})^{+\theta}\\
              &= A^{-\theta} \isSemiActedUponBy \ball(\theta)\\
              &\subseteq A.
            \end{align*}
            Therefore, according to \cref{definition:k-interior-closure-and-boundary}, we have $(A^{-\theta})^{+\theta'} \subseteq A^{-(\theta - \theta')}$. \qedhere
    \end{aenumerate}
  \end{proof}

  \begin{remark}
    Let $M$ be finite, let $\cardinalityOf{M} \geq 2$, let $A = \setOf{m_0}$, let $\theta$ be the least integer such that $A^{+\theta} = M$, which is positive because $\cardinalityOf{M} \geq 2$, and let $\theta'$ be a positive integer such that $\theta' \leq \theta$. Then, $A^{+(\theta - \theta')} \neq M = (A^{+\theta})^{-\theta'}$.
  \end{remark}

  \begin{lemma}
  \label{lemma:distance-of-set-minus-closure-to-set}
    Let $\theta$ be a non-negative integer, and let $A$ and $A'$ be two subsets of $M$. Then, $\distanceOf(A, A' \smallsetminus A^{+\theta}) \geq \theta + 1$.
  \end{lemma}

  \begin{proof}
    If $A$ or $A' \smallsetminus A^{+\theta}$ is empty, then $\distanceOf(A, A' \smallsetminus A^{+\theta}) = \infty \geq \theta + 1$. Otherwise, let $m' \in A' \smallsetminus A^{+\theta}$. Then, according to \cref{item:properties-of-interior-closure-and-boundary:complement} of \cref{lemma:properties-of-interior-closure-and-boundary}, we have $A' \smallsetminus A^{+\theta} = (A' \smallsetminus A)^{-\theta}$. Hence, according to \cref{item:characterisation-of-k-closure-and-interior:interior} of \cref{lemma:characterisation-of-k-closure-and-interior}, we have $\ball(m', \theta) \subseteq A' \smallsetminus A$. Therefore, for each $m \in A$, we have $m \notin \ball(m', \theta)$ and hence $\distanceOf(m, m') \geq \theta + 1$. Thus, $\distanceOf(A, m') \geq \theta + 1$. In conclusion, $\distanceOf(A, A' \smallsetminus A^{+\theta}) \geq \theta + 1$.
  \end{proof}

  \begin{corollary}
  \label{corollary:distance-of-closure-boundary-of-closure-to-set-greater-than-closure-plus-one}
    Let $\theta$ and $\theta'$ be two non-negative integers, and let $A$ be a subset of $M$. Then, $\distanceOf(A, \boundaryOf_{\theta'}^+ A^{+\theta}) \geq \theta + 1$.
  \end{corollary}

  \begin{proof}
    This is a direct consequence of \cref{lemma:distance-of-set-minus-closure-to-set}, because $\boundaryOf_{\theta'}^+ A^{+\theta} = (A^{+\theta})^{+\theta'} \smallsetminus A^{+\theta}$.
  \end{proof}

  \section{Growth Functions and Types}
  \label{section:growth-functions-and-types}

  In this section we recapitulate growth functions and types, more or less as presented in section~6.4 in \cite{ceccherini-silberstein:coornaert:2010}.

  \begin{definition}
    Let $\gamma$ be a map from $\N_0$ to $\R_{\geq 0}$. It is called \graffito{growth function $\gamma$}\define{growth function}\index[symbols]{gamma@$\gamma$} if and only if it is \define{non-decreasing}\graffito{non-decreasing map}\index{decreasing non@non-decreasing}, that is to say, that
    \begin{equation*}
      \ForEach k \in \N_0 \ForEach k' \in \N_0 \Holds \big(k \leq k' \implies \gamma(k) \leq \gamma(k')\big). \qedhere
    \end{equation*}
  \end{definition}

  \begin{definition} 
    Let $\gamma$ and $\gamma'$ be two growth functions. The growth function $\gamma$ is said to \define{dominate}\graffito{$\gamma$ dominates $\gamma'$} $\gamma'$ and we write $\gamma \dominates \gamma'$\graffito{$\gamma \dominates \gamma'$}\index[symbols]{greaterthanorequaltocurly@$\gamma \dominates \gamma'$} if and only if
    \begin{equation*}
      \Exists \alpha \in \N_+ \SuchThat \ForEach k \in \N_+ \Holds \alpha \cdot \gamma(\alpha \cdot k) \geq \gamma'(k). \qedhere
    \end{equation*}
  \end{definition}

  \begin{definition}
    Let $\gamma$ and $\gamma'$ be two growth functions. They are called \define{equivalent}\graffito{$\gamma$ and $\gamma'$ are equivalent} and we write $\gamma \isEquivalentTo \gamma'$\graffito{$\gamma \isEquivalentTo \gamma'$}\index[symbols]{tilde@$\isEquivalentTo$} if and only if $\gamma \dominates \gamma'$ and $\gamma' \dominates \gamma$.
  \end{definition}

  \begin{lemma}
  \label{lemma:dominate-equivalence-relations}
    \begin{aenumerate}
      \item The binary relation $\dominates$ is reflexive and transitive.
      \item The binary relation $\isEquivalentTo$ is an equivalence relation.
      \item\label{item:dominate-equivalence-relations:interplay}
            If $\gamma_1 \isEquivalentTo \gamma_2$ and $\gamma_1' \isEquivalentTo \gamma_2'$, then $\gamma_1 \dominates \gamma_1'$ implies $\gamma_2 \dominates \gamma_2'$. \qedhere
    \end{aenumerate}
  \end{lemma}

  \begin{proof}
    See proposition~6.4.3 in \cite{ceccherini-silberstein:coornaert:2010}.
  \end{proof}

  \begin{definition}
    For each growth function $\gamma$, the \graffito{equivalence class $\equivalenceClassOf{\gamma}_\isEquivalentTo$ of $\gamma$ with respect to $\isEquivalentTo$}equivalence class of $\gamma$ with respect to $\isEquivalentTo$ is denoted by $\equivalenceClassOf{\gamma}_\isEquivalentTo$\index[symbols]{gammabracketstilde@$\equivalenceClassOf{\gamma}_\isEquivalentTo$}. And, each equivalence class $\Gamma$ with respect to $\isEquivalentTo$ is called \define{growth type}\graffito{growth type $\Gamma$}\index[symbols]{Gamma@$\Gamma$}.
  \end{definition}

  \begin{definition}
    Let $\Gamma$ and $\Gamma'$ be two growth types. The growth type $\Gamma$ is said to \define{dominate}\graffito{$\Gamma$ dominates $\Gamma'$} $\Gamma'$ and we write $\Gamma \dominates \Gamma'$\graffito{$\Gamma \dominates \Gamma'$}\index[symbols]{greaterthanorequaltocurly@$\Gamma \dominates \Gamma'$} if and only if 
    \begin{equation*}
      \Exists \gamma \in \Gamma \Exists \gamma' \in \Gamma' \SuchThat \gamma \dominates \gamma'. \qedhere
    \end{equation*} 
  \end{definition}

  \begin{remark}
    According to \cref{item:dominate-equivalence-relations:interplay} of \cref{lemma:dominate-equivalence-relations}, the growth type $\Gamma$ dominates $\Gamma'$ if and only if
    \begin{equation*}
      \ForEach \gamma \in \Gamma \ForEach \gamma' \in \Gamma' \Holds \gamma \dominates \gamma'. \qedhere
    \end{equation*}
  \end{remark}

  \begin{example}[{\cite[Examples~6.4.4, except for the first item]{ceccherini-silberstein:coornaert:2010}}]\leavevmode
  \label{example:growth-functions}
    \begin{aenumerate}
      \item \label{item:growth-functions:id-versus-unity}
            The growth function $[k \mapsto k]$ dominates $\functionThatIsIdenticalToOne$ but they are not equivalent.

            \begin{proof}
              For each $k \in \N_+$, we have $k \geq 1 = \functionThatIsIdenticalToOne(k)$. However, for each $\alpha \in \N_+$, there is a $k \in \N_+$, for example $k = \alpha + 1$, such that $\alpha \functionThatIsIdenticalToOne(\alpha k) = \alpha < k$.
            \end{proof}
      \item Let $r$ and $s$ be two non-negative real numbers. Then, $[k \mapsto k^r] \dominates [k \mapsto k^s]$ if and only if $r \geq s$. And, $[k \mapsto k^r] \isEquivalentTo [k \mapsto k^s]$ if and only if $r = s$. 
      \item Let $\gamma$ be a growth function such that it is a polynomial function of degree $d \in \N_0$. Then, $\gamma \isEquivalentTo [k \mapsto k^d]$. 
      \item \label{item:growth-functions:exponential-growth}
            Let $r$ and $s$ be two elements of $\R_{> 1}$. Then, $[k \mapsto r^k] \isEquivalentTo [k \mapsto s^k]$. In particular, $[k \mapsto r^k] \isEquivalentTo \exp$.

            \begin{proof}
              Without loss of generality, suppose that $r \leq s$. Then, for each $k \in \N_+$, we have $r^k \leq s^k$. Hence, $[k \mapsto r^k] \isDominatedBy [k \mapsto s^k]$. Moreover, let $\alpha = \ceil{\log_r s} \in \N_+$. Then, for each $k \in \N_+$,
              \begin{equation*}
                s^k = (r^{\log_r s})^k = r^{(\log_r s) k} \leq r^{\alpha k} \leq \alpha r^{\alpha k}.
              \end{equation*}
              Hence, $[k \mapsto r^k] \dominates [k \mapsto s^k]$. In conclusion, $[k \mapsto r^k] \isEquivalentTo [k \mapsto s^k]$.
            \end{proof}
      \item\label{item:growth-functions:exp-dominates-polynomials}
            Let $d$ be a non-negative integer. Then, $\exp \dominates [k \mapsto k^d]$ and $\exp \isNotEquivalentTo [k \mapsto k^d]$.

            \begin{proof}
              See examples~6.4.4 (d) in \cite{ceccherini-silberstein:coornaert:2010}. \qedhere
            \end{proof}
    \end{aenumerate}\let\qed\relax
  \end{example}

  \begin{lemma}
  \label{lemma:growth-function-dominated-by-poly-is-dominated-by-exp-and-inequivalent}
    Let $\gamma$ be a growth function and let $d$ be a non-negative integer such that $[k \mapsto k^d] \dominates \gamma$. Then, $\exp \dominates \gamma$ and $\exp \isNotEquivalentTo \gamma$.
  \end{lemma}

  \begin{proof}
    According to \cref{item:growth-functions:exp-dominates-polynomials} of \cref{example:growth-functions}, we have $\exp \dominates [k \mapsto k^d]$ and $\exp \isNotEquivalentTo [k \mapsto k^d]$. Hence, because $\dominates$ is transitive and $[k \mapsto k^d] \dominates \gamma$, we have $\exp \dominates \gamma$ and $\exp \isNotEquivalentTo \gamma$.
  \end{proof}

  \section{Spaces' Growth Functions and Types}
  \label{section:cell-spaces-growth-functions-and-types}

  In this section, let $\mathcal{R} = \ntuple{\mathcal{M}, \mathcal{K}} = \ntuple{\ntuple{M, G, \actsOnPoint}, \ntuple{m_0, \family{g_{m_0, m}}_{m \in M}}}$ be a cell space such that there is a finite and symmetric right-gen\-er\-at\-ing set $S$ of $\mathcal{R}$.

  \paragraph{Contents.} In \cref{definition:growth-function-of-cell-space} we introduce the $S$-growth function $\gamma_S$ of $\mathcal{R}$. In \cref{lemma:metric-and-generating-set} and its corollaries we show that $\gamma_S$ is dominated by $\exp$ and that the $\isEquivalentTo$-equivalence class $\equivalenceClassOf{\gamma_S}_\isEquivalentTo$ does not depend on $S$. In \cref{definition:growth-type-of-cell-space} we introduce the growth type $\gamma(\mathcal{M})$ of $\mathcal{M}$ as that equivalence class. In \cref{lemma:ball-sequence-strictly-increasing-or-eventually-constant} and its corollary we relate the inclusion-behaviour of the sequence of balls to the cardinality of $M$. And in \cref{definition:growth} we say what exponential, sub-exponential, polynomial, and intermediate growth of $\mathcal{R}$ mean.

  \begin{definition}
  \label{definition:growth-function-of-cell-space}
    The map
    \begin{align*}
      \gamma_S \from \N_0 &\to \N_0, \mathnote{$S$-growth function $\gamma_S$ of $\mathcal{R}$}\index[symbols]{gammaS@$\gamma_S$}\\
      k &\mapsto \cardinalityOf{\ball_S(k)},
    \end{align*}
    is called \defineX{$S$-growth function of $\mathcal{R}$}{growth function of $\mathcal{R}$@$S$-growth function of $\mathcal{R}$}.
  \end{definition}

  \begin{remark}
  \label{remark:independence-of-uncoloured-cayley-of-growth-function}
    According to \cref{remark:independence-of-uncoloured-cayley-of-balls-and-spheres}, the $S$-growth function of $\mathcal{R}$ does not depend on the coordinates $\family{g_{m_0, m}}_{m \in M}$ and is identical to the $(g \conjugates S)$-growth function of $\ntuple{\mathcal{M}, \ntuple{g \actsOnPoint m_0, \family{g_{m_0, m} g^{-1}}_{m \in M}}}$.
  \end{remark}

  \begin{remark}
  \label{remark:growth-function-of-cell-space}
    According to \cref{remark:ball-of-radius-0-contains-one-element-and-sequence-of-balls-is-monotonic}, we have $\gamma_S(0) = 1$ and the sequence $\sequence{\gamma_S(k)}_{k \in \N_0}$ is non-decreasing with respect to the partial order $\leq$. Moreover, according to \cref{remark:upper-bound-for-cardinality-of-balls}, for each non-negative integer $k$, we have $\gamma_S(k) \leq (1 + \cardinalityOf{S})^k$.
  \end{remark}

  \begin{lemma} 
  \label{lemma:metric-and-generating-set}
    Let $S'$ be a finite and symmetric right-gen\-er\-at\-ing set of $\mathcal{R}$ and let $\alpha$ be the non-negative integer $\min\setOf{k \in \N_0 \suchThat \ball_S(1) \subseteq \ball_{S'}(k)}$. Then, 
    \begin{equation*}
      \ForEach m \in M \ForEach m' \in M \Holds \distanceOf_{S'}(m, m') \leq \alpha \cdot \distanceOf_S(m, m'),
    \end{equation*}
    in particular, 
    \begin{equation*}
      \ForEach m \in M \Holds \lengthOf{m}_{S'} \leq \alpha \cdot \lengthOf{m}_S. \qedhere
    \end{equation*}
  \end{lemma}

  \begin{proof}
    For each $m \in M$, let $\alpha_m = \min\setOf{k \in \N_0 \suchThat \ball_S(m, 1) \subseteq \ball_{S'}(m, k)}$, in particular, $\alpha_{m_0} = \alpha$.

    First, let $m \in M$. Furthermore, let $k \in \N_0$ and let $g \in G_{m_0, m}$. Then, because $g \actsOnPoint \blank$ is bijective, we have $\ball_S(1) \subseteq \ball_{S'}(k)$ if and only if $g \actsOnPoint \ball_S(1) \subseteq g \actsOnPoint \ball_{S'}(k)$. And, according to \cref{lemma:left-action-and-balls}, we have $g \actsOnPoint \ball_S(1) = \ball_S(m, 1)$ and $g \actsOnPoint \ball_{S'}(k) = \ball_{S'}(m, k)$. Hence, $\ball_S(1) \subseteq \ball_{S'}(k)$ if and only if $\ball_S(m, 1) \subseteq \ball_{S'}(m, k)$. Therefore, $\alpha_m = \alpha$. In conclusion, for each $m \in M$, we have $\alpha_m = \alpha$.

    Secondly, we prove by induction on the distance $k$ that
    \begin{multline*}
      \ForEach k \in \N_0 \ForEach m \in M \ForEach m' \in M \Holds\\
          \parens[\big]{\distanceOf_S(m, m') = k \implies \distanceOf_{S'}(m, m') \leq \alpha \cdot k}.
    \end{multline*}

    \proofPart{Base Case}
      Let $k = 0$. Furthermore, let $m$ and $m' \in M$ such that $\distanceOf_S(m, m') = k$. Then, $m = m'$. Hence, $\distanceOf_{S'}(m, m') = 0$. Therefore, $\distanceOf_{S'}(m, m') \leq \alpha \cdot k$.

    \proofPart{Inductive Step}
      Let $k \in \N_0$ such that
      \begin{equation*}
        \ForEach m \in M \ForEach m' \in M \Holds \parens[\big]{\distanceOf_S(m, m') = k \implies \distanceOf_{S'}(m, m') \leq \alpha \cdot k}.
      \end{equation*}
      Furthermore, let $m$ and $m'' \in M$ such that $\distanceOf_S(m, m'') = k + 1$. Then, there is a finite sequence $\sequence{s_i}_{i \in \setOf{1, 2, \dotsc, k + 1}}$ in $S$ such that $m' \isSemiActedUponBy s_{k + 1} = m''$, where $m' = m \isSemiActedUponBy \sequence{s_i}_{i \in \setOf{1, 2, \dotsc, k}}$. And, according to \cref{lemma:truncated-minimal-path-yields-minimal-path}, we have $\distanceOf_S(m, m') = k$. Therefore, according to the inductive hypothesis, $\distanceOf_{S'}(m, m') \leq \alpha \cdot k$. Moreover, by definition of $\alpha_{m'}$, we have $m'' = m' \isSemiActedUponBy s_{k + 1} \in \ball_S(m', 1) \subseteq \ball_{S'}(m', \alpha_{m'})$. Hence, because $\alpha_{m'} = \alpha$, we have $\distanceOf_{S'}(m', m'') \leq \alpha_{m'} = \alpha$. In conclusion, because $\distanceOf_{S'}$ is subadditive, we have $\distanceOf_{S'}(m, m'') \leq \distanceOf_{S'}(m, m') + \distanceOf_{S'}(m', m'') \leq \alpha \cdot k + \alpha = \alpha \cdot (k + 1)$.
  \end{proof}

  \begin{corollary}
  \label{corollary:balls-and-generating-set}
    In the situation of \cref{lemma:metric-and-generating-set}, for each element $m \in M$ and each non-negative integer $k$, we have $\ball_S(m, k) \subseteq \ball_{S'}(m, \alpha \cdot k)$.
  \end{corollary}

  \begin{proof}
    This is a direct consequence of \cref{lemma:metric-and-generating-set}, because for each element $m \in M$, each non-negative integer $k$, and each element $m' \in M$, if $\distanceOf_S(m, m') \leq k$, then $\distanceOf_{S'}(m, m') \leq \alpha \cdot k$.
  \end{proof}

  \begin{corollary}
  \label{corollary:growth-functions-and-generating-set}
    In the situation of \cref{lemma:metric-and-generating-set}, for each non-negative integer $k$, we have $\gamma_S(k) \leq \gamma_{S'}(\alpha \cdot k)$.
  \end{corollary}

  \begin{proof}
    This is a direct consequence of \cref{corollary:balls-and-generating-set}.
  \end{proof}


  \begin{corollary} 
  \label{corollary:Lipschitz-equivalent-metrics}
    Let $S'$ be a finite and symmetric right-gen\-er\-at\-ing set of $\mathcal{R}$. The metrics $\distanceOf_S$ and $\distanceOf_{S'}$ are \define{Lipschitz equivalent}\graffito{$\distanceOf_S$ and $\distanceOf_{S'}$ are Lipschitz equivalent}\index{equivalent!Lipschitz}, that is to say, that there are two positive real numbers $\kappa$ and $\varkappa$ such that $\kappa \cdot \distanceOf_{S'} \leq \distanceOf_S \leq \varkappa \cdot \distanceOf_{S'}$, where the scalar multiplication $\cdot$ and the partial order $\leq$ are pointwise.
  \end{corollary}

  \begin{proof}
    Let $\alpha = \min\setOf{k \in \N_0 \suchThat \ball_S(1) \subseteq \ball_{S'}(k)}$ and let $\alpha' = \min\setOf{k \in \N_0 \suchThat \ball_{S'}(1) \subseteq \ball_S(k)}$. If $\alpha = 0$ or $\alpha' = 0$, then $M = \setOf{m_0}$, hence $\distanceOf_S = \functionThatIsIdenticalToZero = \distanceOf_{S'}$, and therefore $\distanceOf_S \leq \distanceOf_{S'} \leq \distanceOf_S$. Otherwise, according to \cref{lemma:metric-and-generating-set}, we have $\frac{1}{\alpha} \cdot \distanceOf_{S'} \leq \distanceOf_S \leq \alpha' \cdot \distanceOf_{S'}$.
  \end{proof}

  \begin{corollary} 
  \label{corollary:growth-functions-independent-of-generating-set-and-gamma-dominated-by-exp}
    Let $S'$ be a finite and symmetric right-gen\-er\-at\-ing set of $\mathcal{R}$. The $S$-growth function $\gamma_S$ of $\mathcal{R}$ and the $S'$-growth function $\gamma_{S'}$ of $\mathcal{R}$ are equivalent.
  \end{corollary}

  \begin{proof}
    According to \cref{corollary:growth-functions-and-generating-set}, there is a $\alpha \in \N_0$ such that, for each $k \in \N_0$, we have $\gamma_S(k) \leq \gamma_{S'}(\alpha \cdot k)$. Hence, according to \cref{remark:growth-function-of-cell-space}, for each $k \in \N_0$, we have $\gamma_S(k) \leq (\alpha + 1) \gamma_{S'}((\alpha + 1) \cdot k)$. Therefore, $\gamma_S$ is dominated by $\gamma_{S'}$. Switching roles of $S$ and $S'$ yields that $\gamma_{S'}$ is dominated by $\gamma_S$. In conclusion, $\gamma_S$ and $\gamma_{S'}$ are equivalent.
  \end{proof}

  \begin{corollary}
  \label{corollary:growth-function-dominated-by-exp}
    The $S$-growth function $\gamma_S$ of $\mathcal{R}$ is dominated by $\exp$.
  \end{corollary}

  \begin{proof}
    According to \cref{remark:growth-function-of-cell-space}, for each $k \in \N_0$, we have $\gamma_S(k) \leq r^k$, where $r = 1 + \cardinalityOf{S}$. Hence, $\gamma_S \isDominatedBy [k \mapsto r^k]$. Moreover, according to \cref{item:growth-functions:exponential-growth} of \cref{example:growth-functions}, we have $[k \mapsto r^k] \isEquivalentTo \exp$. In conclusion, $\gamma_S \isDominatedBy \exp$.
  \end{proof}

  \begin{definition}
  \label{definition:growth-type-of-cell-space}
    The equivalence class $\gamma(\mathcal{M}) = \equivalenceClassOf{\gamma_S}_\isEquivalentTo$ is called \define{growth type of $\mathcal{M}$}\graffito{growth type $\gamma(\mathcal{M})$ of $\mathcal{M}$}\index[symbols]{gammaMcalligraphic@$\gamma(\mathcal{M})$}.
  \end{definition}

  \begin{remark}
    Note that, according to \cref{remark:independence-of-uncoloured-cayley-of-growth-function} and corollary \ref{corollary:growth-functions-independent-of-generating-set-and-gamma-dominated-by-exp}, the equivalence class $\equivalenceClassOf{\gamma_S}_\isEquivalentTo$ does neither depend on the right-gen\-er\-at\-ing set $S$ nor on the coordinate system $\mathcal{K}$. 
  \end{remark}

  \begin{example}[Groups {\cite[Examples~6.4.11]{ceccherini-silberstein:coornaert:2010}}]
    The growth type of the group $\Z$ of integers, of the direct product of the groups $\Z$ and $\Z \modulo 2\Z$, and of the infinite dihedral group $\dihedralGroup_\infty$ is $[k \mapsto k]_\isEquivalentTo$. The growth type of the group $\Z^2$ is $[k \mapsto k^2]_\isEquivalentTo$. And, the growth type of each free group of finite rank is $[\exp]_\isEquivalentTo$.
  \end{example}

  \begin{lemma} 
  \label{lemma:growth-function-equivalent-to-unity-if-and-only-if-bounded}
    Let $\gamma$ be a growth function such that $\gamma(0) > 0$. Then, $\gamma$ is equivalent to $\functionThatIsIdenticalToOne$ if and only if $\gamma$ is bounded.
  \end{lemma}

  \begin{proof}
    See proposition~6.4.6 in \cite{ceccherini-silberstein:coornaert:2010}.
  \end{proof}

  \begin{corollary} 
  \label{corollary:m-finite-if-and-only-if-growth-type-of-m-constant}
    The set $M$ is finite if and only if the growth types $\gamma(\mathcal{M})$ and $\equivalenceClassOf{\functionThatIsIdenticalToOne}_\isEquivalentTo$ are equal.
  \end{corollary}

  \begin{proof}
    First, let $M$ be finite. Then, for each $k \in \N_0$, we have $\gamma_S(k) \leq \cardinalityOf{M}$. Hence, according to \cref{lemma:growth-function-equivalent-to-unity-if-and-only-if-bounded}, we have $\gamma_S \isEquivalentTo \functionThatIsIdenticalToOne$. In conclusion, $\gamma(\mathcal{M}) = \equivalenceClassOf{\functionThatIsIdenticalToOne}_\isEquivalentTo$.

    Secondly, let $\gamma(\mathcal{M}) = \equivalenceClassOf{\functionThatIsIdenticalToOne}_\isEquivalentTo$. Then, according to \cref{lemma:growth-function-equivalent-to-unity-if-and-only-if-bounded}, the growth function $\gamma_S$ is bounded by some $\xi \in \R_{> 0}$. And, according to \cref{remark:ball-of-radius-0-contains-one-element-and-sequence-of-balls-is-monotonic}, we have $M = \bigcup_{k \in \N_0} \ball_S(k)$ and $\sequence{\ball_S(k)}_{k \in \N_0}$ is non-decreasing with respect to $\subseteq$. Therefore, because $\sequence{\gamma_S(k)}_{k \in \N_0} = \sequence{\cardinalityOf{\ball_S(k)}}_{k \in \N_0}$, we have $\cardinalityOf{M} \leq \sup_{k \in \N_0} \gamma_S(k) \leq \xi$. In conclusion, $M$ is finite.
  \end{proof}

  \begin{lemma}
  \label{lemma:ball-sequence-strictly-increasing-or-eventually-constant}
    Either the sequence $\sequence{\ball_S(k)}_{k \in \N_0}$ is strictly increasing with respect to $\subseteq$ or \define{eventually constant}\graffito{eventually constant sequence}, that is to say, that there is a non-negative integer $k$ such that, for each non-negative integer $k'$ with $k' \geq k$, we have $\ball_S(k') = \ball_S(k)$.
  \end{lemma}

  \begin{proof}
    According to \cref{remark:ball-of-radius-0-contains-one-element-and-sequence-of-balls-is-monotonic}, the sequence $\sequence{\ball_S(k)}_{k \in \N_0}$ is non-decreasing with respect to $\subseteq$. If it is strictly increasing with respect to $\subseteq$, it is not eventually constant. Otherwise, there is a $k \in \N_0$ such that $\ball_S(k) = \ball_S(k + 1)$. We prove by induction on $k'$ that, for each $k' \in \N_0$ with $k' \geq k$, we have $\ball_S(k') = \ball_S(k)$.

    \proofPart{Base Case}
      Let $k' = k$. Then, $\ball_S(k') = \ball_S(k)$.

    \proofPart{Inductive Step}
      Let $k' \in \N_0$ with $k' \geq k$ such that $\ball_S(k') = \ball_S(k)$. Furthermore, let $m \in \ball_S(k' + 1)$.
      \begin{description}
        \item[Case 1:] $m \in \ball_S(k')$. Then, according to the inductive hypothesis, $m \in \ball_S(k)$.
        \item[Case 2:] $m \notin \ball_S(k')$. Then, there is a finite sequence $\sequence{s_i}_{i \in \setOf{1, 2, \dotsc, k' + 1}}$ in $S$ such that $m' \isSemiActedUponBy s_{k' + 1} = m$, where $m' = m_0 \isSemiActedUponBy \sequence{s_i}_{i \in \setOf{1, 2, \dotsc, k'}}$. Hence, $m' \in \ball_S(k')$ and thus, according to the inductive hypothesis, $m' \in \ball_S(k)$. Therefore, according to \cref{lemma:ball-liberation-included-in-ball-one-larger}, we have $m \in \ball_S(k + 1)$. Thus, because $\ball_S(k + 1) = \ball_S(k)$, we have $m \in \ball_S(k)$.
      \end{description}
      In either case, $m \in \ball_S(k)$. Therefore, $\ball_S(k' + 1) \subseteq \ball_S(k) \subseteq \ball_S(k') \subseteq \ball_S(k' + 1)$. In conclusion, $\ball_S(k' + 1) = \ball_S(k)$.
  %
  \end{proof}

  \begin{corollary}
  \label{corollary:finite-if-and-only-if-eventually-equal-to-m}
    The set $M$ is finite if and only if the sequence $\sequence{\ball_S(k)}_{k \in \N_0}$ is \define{eventually equal to $M$}\graffito{eventually equal to $M$ sequence}, that is to say, that there is a non-negative integer $k$ such that, for each non-negative integer $k'$ with $k' \geq k$, we have $\ball_S(k') = M$.
  \end{corollary}

  \begin{proof}
    First, let $M$ be finite. Then, $\sequence{\ball_S(k)}_{k \in \N_0}$ is not strictly increasing. Hence, according to \cref{lemma:ball-sequence-strictly-increasing-or-eventually-constant}, it is eventually constant. And, according to \cref{remark:ball-of-radius-0-contains-one-element-and-sequence-of-balls-is-monotonic}, it converges to $M$. In conclusion, it is eventually equal to $M$.

    Secondly, let $\sequence{\ball_S(k)}_{k \in \N_0}$ be eventually equal to $M$. Then, there is a $k \in \N_0$ such that $\ball_S(k) = M$. In conclusion, according to \cref{remark:upper-bound-for-cardinality-of-balls}, the set $M$ is finite.
  \end{proof}

  \begin{corollary}
  \label{corollary:m-infinite-if-and-only-if-balls-strictly-increasing}
    The set $M$ is infinite if and only if the sequence $\sequence{\ball_S(k)}_{k \in \N_0}$ is strictly increasing with respect to $\subseteq$.
  \end{corollary}

  \begin{proof}
    First, let $M$ be infinite. Then, because $\ball_S(k)$ is finite, for $k \in \N_0$, and $\sequence{\ball_S(k)}_{k \in \N_0}$ converges to $M$, the sequence $\sequence{\ball_S(k)}_{k \in \N_0}$ is not eventually constant. In conclusion, according to \cref{lemma:ball-sequence-strictly-increasing-or-eventually-constant}, it is strictly increasing with respect to $\subseteq$.

    Secondly, let $\sequence{\ball_S(k)}_{k \in \N_0}$ be strictly increasing with respect to $\subseteq$. Then, because it converges to $M$, the set $M$ is infinite.
  \end{proof}

  \begin{corollary}
  \label{corollary:m-infinite-if-and-only-if-spheres-are-non-empty}
    The set $M$ is infinite if and only if
    \begin{equation}
    \label{equation:m-infinite-if-and-only-if-spheres-are-non-empty}
      \ForEach \rho \in \N_0 \Holds \sphere_S(\rho) \neq \emptyset. \qedhere
    \end{equation}
  \end{corollary}

  \begin{proof}
    We have $\sphere_S(0) = \setOf{m_0} \neq \emptyset$. And, according to \cref{remark:spheres-expressed-in-terms-of-balls}, for each $\rho \in \N_+$, we have $\sphere_S(\rho) = \ball_S(\rho) \smallsetminus \ball_S(\rho - 1)$. Hence, the sequence $\sequence{\ball_S(k)}_{k \in \N_0}$ is strictly increasing with respect to $\subseteq$ if and only if \cref{equation:m-infinite-if-and-only-if-spheres-are-non-empty} holds. Therefore, according to \cref{corollary:m-infinite-if-and-only-if-balls-strictly-increasing}, the set $M$ is infinite if and only if \cref{equation:m-infinite-if-and-only-if-spheres-are-non-empty} holds.
  \end{proof}

  \begin{lemma} 
    The set $M$ is infinite if and only if the growth type of $\mathcal{M}$ dominates $\equivalenceClassOf{k \mapsto k}_\isEquivalentTo$.
  \end{lemma}

  \begin{proof}
    First, let $M$ be infinite. Then, according to \cref{corollary:m-infinite-if-and-only-if-balls-strictly-increasing}, the sequence $\sequence{\ball_S(k)}_{k \in \N_0}$ is strictly increasing with respect to $\subseteq$. Hence, because $\ball_S(0) = \setOf{m_0}$, for each $k \in \N_0$, we have $\gamma_S(k) = \cardinalityOf{\ball_S(k)} \geq k + 1$. In conclusion, $\gamma_S$ dominates $[k \mapsto k]$ and hence $\gamma(\mathcal{M})$ dominates $\equivalenceClassOf{k \mapsto k}_\isEquivalentTo$.

    Secondly, let $M$ be finite. Then, according to \cref{corollary:m-finite-if-and-only-if-growth-type-of-m-constant}, we have $\gamma(\mathcal{M}) = \equivalenceClassOf{\functionThatIsIdenticalToOne}_\isEquivalentTo$. Hence, according to \cref{item:growth-functions:id-versus-unity} of \cref{example:growth-functions}, the growth type $\gamma(\mathcal{M})$ does not dominate $\equivalenceClassOf{k \mapsto k}_\isEquivalentTo$.
  \end{proof}

  Because Cayley graphs of $\mathcal{R}$ are in a sense quotients of Cayley graphs of $G$ by the finite subgroup $G_0$, the distances on these graphs are related by the multiplicative constant $\cardinalityOf{G_0}$, which is implicitly used in

  \begin{lemma} 
  \label{lemma:group-dominates-cell-space-growth-rate}
    Let the group $G$ be finitely generated and let the stabiliser $G_0$ be finite. The growth types of $G$ and $\mathcal{M}$ are equal.
  \end{lemma}

  \begin{proof} 
    Because the group $G_0$ is finite, there is a finite and symmetric generating set $T$ of $G$ such that $G_0 T \subseteq T$. And, according to \cref{lemma:gen-set-of-group-induces-right-gen-set-of-action}, the set $S = \setOf{t G_0 \suchThat t \in T} = \setOf{g_0 \cdot t G_0 \suchThat g_0 \in G_0, t \in T}$ is a finite and symmetric right-gen\-er\-at\-ing set of $\mathcal{R}$.

    Let $k$ be a non-negative integer. Furthermore, let $m$ be an element of $\ball_S^{\mathcal{R}}(k)$. Then, there is a non-negative integer $j \in \setOf{0,1, 2, \dotsc, k}$ and there is a finite sequence $\sequence{s_i}_{i \in \setOf{1, 2, \dotsc, j}}$ of elements in $S$ such that
    \begin{equation*}
      m = m_0 \isSemiActedUponBy \sequence{s_i}_{i \in \setOf{1, 2, \dotsc, j}}.
    \end{equation*}
    And, by the definition of $S$, there is a finite sequence $\sequence{t_i}_{i \in \setOf{1, 2, \dotsc, j}}$ of elements in $T$ such that $\sequence{t_i G_0}_{i \in \setOf{1, 2, \dotsc, j}} = \sequence{s_i}_{i \in \setOf{1, 2, \dotsc, j}}$. And, because $\isSemiActedUponBy$ is a semi-action, there is a finite sequence $\sequence{g_{i,0}}_{i \in \setOf{1, 2, \dotsc, j}}$ of elements in $G_0$ such that
    \begin{align*}
      m
      &= m_0 \isSemiActedUponBy \sequence{s_i}_{i \in \setOf{1, 2, \dotsc, j}}\\
      &= m_0 \isSemiActedUponBy t_1 g_{2,0} t_2 g_{3,0} t_3 \dotsb g_{j,0} t_j G_0\\
      &= g_{1,0} t_1 g_{2,0} t_2 g_{3,0} t_3 \dotsb g_{j,0} t_j \actsOnPoint m_0,
    \end{align*}
    where $g_{1,0} = e_G$. And, because $G_0 T \subseteq T$, the sequence $\sequence{g_{i,0} t_i}_{i \in \setOf{1, 2, \dotsc, j}}$ is one of elements in $T$. Hence, the element $m$ is contained in $\ball_T^G(j) \actsOnPoint m_0$, which is included in $\ball_T^G(k) \actsOnPoint m_0$. Therefore, the set $\ball_S^{\mathcal{R}}(k)$ is included in $\ball_T^G(k) \actsOnPoint m_0$. Analogously, one can show that the set $\ball_T^G(k) \actsOnPoint m_0$ is included in $\ball_S^{\mathcal{R}}(k)$. It follows that $\ball_S^{\mathcal{R}}(k)$ is equal to $\ball_T^G(k) \actsOnPoint m_0$.

    Hence, $\cardinalityOf{\ball_S^{\mathcal{R}}(k)} \leq \cardinalityOf{\ball_T^G(k)}$. And, because the map $\blank \actsOnPoint m_0$ from $G$ to $M$ is $\cardinalityOf{G_0}$-to-$1$ surjective, we have $\cardinalityOf{\ball_T^G(k)} \leq \cardinalityOf{G_0} \cdot \cardinalityOf{\ball_S^{\mathcal{R}}(k)} \leq \cardinalityOf{G_0} \cdot \cardinalityOf{\ball_S^{\mathcal{R}}(\cardinalityOf{G_0} \cdot k)}$. Therefore, the growth function $\gamma_T^G$ dominates $\gamma_S^{\mathcal{R}}$ and vice versa. In conclusion, the growth types $\gamma(G)$ and $\gamma(\mathcal{M})$ are equal.
  \end{proof}


  \begin{definition} 
  \label{definition:growth}
    The left-ho\-mo\-ge\-neous space $\mathcal{M}$ is said to have
    \begin{aenumerate} 
      \item \define{exponential growth}\graffito{exponential growth}\index{growth!exponential} if and only if its growth type is $\equivalenceClassOf{\exp}_\isEquivalentTo$;
      \item \define{sub-exponential growth}\graffito{sub-exponential growth}\index{growth!sub-exponential} if and only if it does not have exponential growth;
      \item \define{polynomial growth}\graffito{polynomial growth}\index{growth!polynomial} if and only if there is a non-negative integer $d \in \N_0$ such that $\equivalenceClassOf{k \mapsto k^d}_\isEquivalentTo$ dominates $\gamma(\mathcal{M})$.
      \item \define{intermediate growth}\graffito{intermediate growth}\index{growth!intermediate} if and only if it has sub-exponential growth but not polynomial growth. \qedhere
    \end{aenumerate}
  \end{definition}

  \begin{example}[Free Groups of Finite Rank]
    Each free group of finite rank has exponential growth (see examples~6.4.11 (g) in \cite{ceccherini-silberstein:coornaert:2010}).
  \end{example}

  \begin{example}[Virtually Nilpotent Groups]
  \label{example:virtually-nilpotent-group}
    Each virtually nilpotent finitely generated group has polynomial growth (see corollary~6.8.5 in \cite{ceccherini-silberstein:coornaert:2010}. Examples of such groups are abelian groups, like the group of integers under addition, nilpotent but non-abelian groups, like the discrete Heisenberg group, and virtually nilpotent but non-nilpotent groups, like the infinite dihedral group. 
  \end{example}

  \begin{example}[Grigorchuk Group]
  \label{example:Grigorchuck-group}
    Rostislav Ivanovich Grigorchuk was the first to construct a finitely generated group of intermediate growth in 1984. This group is known as \define{Grigorchuk group}\graffito{Grigorchuk group}. The original construction and proofs were published in the paper \enquote{\citetitle*{grigorchuk:1985}}\cite{grigorchuk:1985}; a more accessible exposition can be found in section~6.9 in \cite{ceccherini-silberstein:coornaert:2010}.
  \end{example}

  \begin{example}[Left Homogeneous Spaces]
    Each quotient set of any group of the previous examples by a finite subgroup acted upon by left multiplication is, according to \cref{lemma:group-dominates-cell-space-growth-rate}, a left-ho\-mo\-ge\-neous space of the respective growth.
  \end{example}

  \begin{lemma} 
    Let $\mathcal{M}$ have polynomial growth. It has sub-exponential growth.
  \end{lemma}

  \begin{proof}
    There is a $d \in \N_0$ such that $\equivalenceClassOf{k \mapsto k^d}_\isEquivalentTo \dominates \gamma(\mathcal{M})$. Hence, $[k \mapsto k^d] \dominates \gamma_S$. Therefore, according to \cref{lemma:growth-function-dominated-by-poly-is-dominated-by-exp-and-inequivalent}, we have $\exp \isNotEquivalentTo \gamma_S$. In conclusion, $\gamma(\mathcal{M}) \neq \equivalenceClassOf{\exp}_\isEquivalentTo$.
  \end{proof}

  \section{Growth Rates}
  \label{section:growth-rates}

  In this section, let $\mathcal{R} = \ntuple{\mathcal{M}, \mathcal{K}} = \ntuple{\ntuple{M, G, \actsOnPoint}, \ntuple{m_0, \family{g_{m_0, m}}_{m \in M}}}$ be a cell space such that there is a finite and symmetric right-gen\-er\-at\-ing set $S$ of $\mathcal{R}$.

  \paragraph{Contents.} In \cref{definition:growth-rate} we introduce the $S$-growth rate of $\mathcal{R}$. And in \cref{theorem:growth-rate-greater-than-one-if-and-only-if-exp-growth} show how that growth rate and exponential growth relate to each other.

  \begin{lemma} 
  \label{lemma:growth-rate-exists-and-is-greater-than-or-equal-to-one}
    The sequence $\sequence{\sqrt[k]{\gamma_S(k)}}_{k \in \N_0}$ converges to $\inf_{k \in \N_0} \sqrt[k]{\gamma_S(k)} \in \R_{\geq 1}$.
  \end{lemma}

  \begin{proof}
    According to \cref{corollary:liberation-of-balls-yields-bigger-one},
    \begin{align*}
      \gamma_S(k + k') &=    \cardinalityOf{\ball_S(k + k')}\\
                       &=    \cardinalityOf{\ball_S(k) \isSemiActedUponBy \ball_S(k')}\\
                       &\leq \cardinalityOf{\ball_S(k)} \cdot \cardinalityOf{\ball_S(k')}\\
                       &=    \gamma_S(k) \cdot \gamma_S(k').
    \end{align*}
    Hence, according to lemma~6.5.1 in \cite{ceccherini-silberstein:coornaert:2010}, the sequence $\sequence{\sqrt[k]{\gamma_S(k)}}_{k \in \N_0}$ converges to $\inf_{k \in \N_0} \sqrt[k]{\gamma_S(k)}$. Moreover, because, for each $k \in \N_0$, we have $\gamma_S(k) \geq 1$, that limit point is in $\R_{\geq 1}$.
  \end{proof}

  \begin{definition} 
  \label{definition:growth-rate}
    The limit point $\lambda_S = \lim_{k \to \infty} \sqrt[k]{\gamma_S(k)}$ is called \defineX{$S$-growth rate of $\mathcal{R}$}{growth rate of $\mathcal{R}$@$S$-growth rate of $\mathcal{R}$}\graffito{$S$-growth rate $\lambda_S$ of $\mathcal{R}$}\index[symbols]{lambdaS@$\lambda_S$}.
  \end{definition}

  \begin{example}[Lattice]
    In the situation of \cref{example:lattice:Folner}, the $\setOf{(-1, 0),\allowbreak (0, -1),\allowbreak (0, 1),\allowbreak (1, 0)}$-growth rate $\lim_{k \to \infty} \sqrt[k]{\cardinalityOf{\ball(k)}}$ of $\mathcal{R}$ is equal to $1$.
  \end{example}

  \begin{example}[Tree]
    In the situation of \cref{example:tree:Folner}, the $\setOf{a, b, a^{-1}, b^{-1}}$-growth rate $\lim_{k \to \infty} \sqrt[k]{\cardinalityOf{\ball(k)}}$ of $\mathcal{R}$ is equal to $3$.
  \end{example}

  \begin{theorem} 
  \label{theorem:growth-rate-greater-than-one-if-and-only-if-exp-growth}
    The $S$-growth rate of $\mathcal{R}$ is greater than $1$ if and only if the left-ho\-mo\-ge\-neous space $\mathcal{M}$ has exponential growth. 
  \end{theorem}

  \begin{proof} 
    First, let $\lambda_S > 1$. According to \cref{lemma:growth-rate-exists-and-is-greater-than-or-equal-to-one}, for each $k \in \N_0$, we have $\sqrt[k]{\gamma_S(k)} \geq \lambda_S$ and hence $\gamma_S(k) \geq \lambda_S^k$. Hence, $\gamma_S$ dominates $\lambda_S^{(\blank)}$. And, because $\lambda_S > 1$, according to \cref{item:growth-functions:exponential-growth} of \cref{example:growth-functions}, the growth function $\lambda_S^{(\blank)}$ is equivalent to $\exp$. Therefore, $\gamma_S$ dominates $\exp$. Moreover, according to \cref{corollary:growth-function-dominated-by-exp}, the growth function $\gamma_S$ is dominated by $\exp$. Altogether, $\gamma_S$ and $\exp$ are equivalent. In conclusion, $\gamma(\mathcal{M}) = \equivalenceClassOf{\gamma_S}_\isEquivalentTo = \equivalenceClassOf{\exp}_\isEquivalentTo$.

    Secondly, let $\gamma(\mathcal{M}) = \equivalenceClassOf{\exp}_\isEquivalentTo$. Then, $\gamma_S$ and $\exp$ are equivalent. In particular, $\gamma_S$ dominates $\exp$. Hence, there is a $\alpha \in \N_+$ such that, for each $k \in \N_+$, we have $\alpha \cdot \gamma_S(\alpha \cdot k) \geq \exp(k)$. Therefore, for each $k \in \N_+$, 
    \begin{equation*}
      \sqrt[\alpha k]{\alpha} \cdot \sqrt[\alpha k]{\gamma_S(\alpha \cdot k)}
      =    \sqrt[\alpha k]{\alpha \cdot \gamma_S(\alpha \cdot k)}
      \geq \sqrt[\alpha k]{\exp(k)}
      =    \sqrt[\alpha]{\EulersNumber}.
    \end{equation*}
    Thus, because $\sequence{\sqrt[\alpha k]{\alpha}}_{k \in \N_+}$ converges to $1$ and $\sequence{\sqrt[\alpha k]{\gamma_S(\alpha \cdot k)}}_{k \in \N_+}$, as subsequence of $\sequence{\sqrt[k]{\gamma_S(k)}}_{k \in \N_0}$, converges to $\lambda_S$, we conclude that $\lambda_S \geq \sqrt[\alpha]{\EulersNumber} > 1$.
  \end{proof}

  \begin{corollary} 
  \label{corollary:growth-rate-equal-to-one-if-and-only-if-subexp-growth}
    The $S$-growth rate of $\mathcal{R}$ is equal to $1$ if and only if the left-ho\-mo\-ge\-neous space $\mathcal{M}$ has sub-exponential growth.
  \end{corollary}

  \begin{proof}
    This is a direct consequence of \cref{theorem:growth-rate-greater-than-one-if-and-only-if-exp-growth}.
  \end{proof}

  \begin{corollary} 
    Let $S'$ be a finite and symmetric right-gen\-er\-at\-ing set of $\mathcal{R}$. The $S$-growth rate of $\mathcal{R}$ is equal to $1$ or greater than $1$ if and only if the $S'$-growth rate of $\mathcal{R}$ is equal to $1$ or greater than $1$ respectively.
  \end{corollary}

  \begin{proof} 
    This is a direct consequence of \cref{corollary:growth-rate-equal-to-one-if-and-only-if-subexp-growth} and theorem \ref{theorem:growth-rate-greater-than-one-if-and-only-if-exp-growth}. 
  \end{proof}

  \section{Amenability, Følner Conditions/Nets, and Isoperimetric Constants} 
  \label{section:Folner-conditions}

  In this section, let $\mathcal{R} = \ntuple{\ntuple{M, G, \actsOnPoint}, \ntuple{m_0, \family{g_{m_0, m}}_{m \in M}}}$ be a finitely right-gen\-er\-at\-ed cell space such that the stabiliser $G_0$ is finite, and let $S$ be a finite and symmetric right-gen\-er\-at\-ing set of $\mathcal{R}$ (note that such a set $S$ exists due to the assumptions on $\mathcal{R}$).

  \paragraph{Contents.} In \cref{definition:isoperimetric-constant} we introduce the $S$-isoperimetric constant of $\mathcal{R}$, which measures, broadly speaking, the invariance under $\isSemiActedUponBy\restrictedTo_{M \times S}$ that a finite subset of $M$ can have, where $0$ means maximally and $1$ minimally invariant. In \cref{theorem:Folner-condition-for-finitely-right-generated-group-sets} we show that $\mathcal{R}$ is right amenable if and only if a kind of Følner condition holds, which in turn holds if and only if the $S$-isoperimetric constant is $0$. And in \cref{theorem:k-boundary-characterisation-of-Folner-net} we characterise right Følner nets using $\rho$-boundaries.

  \begin{definition}
  \label{definition:isoperimetric-constant}
    Let $E$ be a subset of $G \modulo G_0$ and let $\mathcal{F}$ be the set $\setOf{F \subseteq M \suchThat F \neq \emptyset, F \text{ finite}}$. The non-negative real number
    \begin{equation*} 
      \iota_E(\mathcal{R}) = \inf_{F \in \mathcal{F}} \frac{\cardinalityOf{\bigcup_{e \in E} F \smallsetminus (\blank \isSemiActedUponBy e)^{-1}(F)}}{\cardinalityOf{F}} \in \closedInterval{0, 1}
      \mathnote{$E$-isoperimetric constant $\iota_E(\mathcal{R})$ of $\mathcal{R}$}
      \index[symbols]{iotaERcalligraphic@$\iota_E(\mathcal{R})$}
    \end{equation*}
    is called \defineX{$E$-isoperimetric constant of $\mathcal{R}$}{isoperimetric constant of $\mathcal{R}$@$E$-isoperimetric constant of $\mathcal{R}$}.
  \end{definition}

  \begin{example}[Lattice] 
    In the situation of \cref{example:lattice:Folner}, because the balls $\ball(\rho)$, for $\rho \in \N_0$, are non-empty and finite, and the sequence $\sequence{\ball(\rho)}_{\rho \in \N_0}$ is a right Følner net in $\mathcal{R}$, for each finite subset $E$ of $\Z^2$, the $E$-isoperimetric constant of $\mathcal{R}$ is equal to $0$.
  \end{example}

  \begin{remark}
  \label{remark:liberation-of-setminus-is-liberation-of-setminus-of-liberation-and-}
    Let $\mathfrak{g}$ and $\mathfrak{g}'$ be two elements of $G \modulo G_0$, and let $A$, $B$, and $C$ be three sets. Then,
    \begin{aenumerate}
      \item\label{item:liberation-of-setminus-is-liberation-of-setminus-of-liberation-and-:one}
            $\displaystyle \big((\blank \isSemiActedUponBy \mathfrak{g}) \isSemiActedUponBy \mathfrak{g}'\big)^{-1}(A) = (\blank \isSemiActedUponBy \mathfrak{g})^{-1}\big((\blank \isSemiActedUponBy \mathfrak{g}')^{-1}(A)\big)$;
      \item\label{item:liberation-of-setminus-is-liberation-of-setminus-of-liberation-and-:two} 
            $\displaystyle (\blank \isSemiActedUponBy \mathfrak{g})^{-1}(A \smallsetminus B) = (\blank \isSemiActedUponBy \mathfrak{g})^{-1}(A) \smallsetminus (\blank \isSemiActedUponBy \mathfrak{g})^{-1}(B)$. \qedhere
    \end{aenumerate}
  \end{remark}

  \begin{remark} 
  \label{remark:a-setminus-b-leq-a-setminus-c-plus-c-setminus-b}
    Let $A$, $B$, and $C$ be three finite sets. Then,
    \begin{equation*}
      \cardinalityOf{A \smallsetminus B} \leq \cardinalityOf{A \smallsetminus C} + \cardinalityOf{C \smallsetminus B}. \qedhere 
    \end{equation*}
  \end{remark}

  \begin{lemma}
  \label{lemma:liberation-minus-liberation-lib-s-subseteq-bigcup-liberation-of-setminus-liberation-gzero-s}
    Let $A$ be a subset of $M$, let $\mathfrak{g}$ and $\mathfrak{g}'$ be two elements of $G \modulo G_0$, and identify $M$ with $G \modulo G_0$ by $\iota \givenBy m \mapsto G_{m_0, m}$. Then,
    \begin{multline*}
      (\blank \isSemiActedUponBy \mathfrak{g})^{-1}(A) \smallsetminus \big(\blank \isSemiActedUponBy (\mathfrak{g} \isSemiActedUponBy \mathfrak{g}')\big)^{-1}(A)\\
      \subseteq \bigcup_{g_0 \in G_0} (\blank \isSemiActedUponBy \mathfrak{g})^{-1}\big(A \smallsetminus (\blank \isSemiActedUponBy g_0 \cdot \mathfrak{g}')^{-1}(A)\big). \qedhere
    \end{multline*}
  \end{lemma}

  \begin{proof}
    If $(\blank \isSemiActedUponBy \mathfrak{g})^{-1}(A) \smallsetminus (\blank \isSemiActedUponBy (\mathfrak{g} \isSemiActedUponBy \mathfrak{g}'))^{-1}(A)$ is empty, then there is nothing to show. Otherwise, let $m \in (\blank \isSemiActedUponBy \mathfrak{g})^{-1}(A) \smallsetminus (\blank \isSemiActedUponBy (\mathfrak{g} \isSemiActedUponBy \mathfrak{g}'))^{-1}(A)$. Then, according to \cref{lemma:almost-associativity-of-liberation-under-identification}, there is a $g_0 \in G_0$ such that $(m \isSemiActedUponBy \mathfrak{g}) \isSemiActedUponBy g_0 \cdot \mathfrak{g}' = m \isSemiActedUponBy (\mathfrak{g} \isSemiActedUponBy \mathfrak{g}') \notin A$. Therefore, $m \notin ((\blank \isSemiActedUponBy \mathfrak{g}) \isSemiActedUponBy g_0 \cdot \mathfrak{g}')^{-1}(A)$ and hence $m \in (\blank \isSemiActedUponBy \mathfrak{g})^{-1}(A) \smallsetminus ((\blank \isSemiActedUponBy \mathfrak{g}) \isSemiActedUponBy g_0 \cdot \mathfrak{g}')^{-1}(A)$.

    Moreover, according to \cref{item:liberation-of-setminus-is-liberation-of-setminus-of-liberation-and-:one} of \cref{remark:liberation-of-setminus-is-liberation-of-setminus-of-liberation-and-}, we have $((\blank \isSemiActedUponBy \mathfrak{g}) \isSemiActedUponBy g_0 \cdot \mathfrak{g}')^{-1}(A) = (\blank \isSemiActedUponBy \mathfrak{g})^{-1}((\blank \isSemiActedUponBy g_0 \cdot \mathfrak{g}')^{-1}(A))$. Hence, according to \cref{item:liberation-of-setminus-is-liberation-of-setminus-of-liberation-and-:two} of \cref{remark:liberation-of-setminus-is-liberation-of-setminus-of-liberation-and-}, we have $(\blank \isSemiActedUponBy \mathfrak{g})^{-1}(A) \smallsetminus ((\blank \isSemiActedUponBy \mathfrak{g}) \isSemiActedUponBy g_0 \cdot \mathfrak{g}')^{-1}(A) = (\blank \isSemiActedUponBy \mathfrak{g})^{-1}(A \smallsetminus (\blank \isSemiActedUponBy g_0 \cdot \mathfrak{g}')^{-1}(A))$. Therefore, $m \in \bigcup_{g_0 \in G_0} (\blank \isSemiActedUponBy \mathfrak{g})^{-1}(A \smallsetminus (\blank \isSemiActedUponBy g_0 \cdot \mathfrak{g}')^{-1}(A))$. In conclusion, the stated inclusion holds. 
  \end{proof}

  \begin{lemma}
  \label{lemma:one-Folner-condition-for-finitely-right-generated-group-sets:Folner-condition}
    Let $\mathcal{F}$ be a subset of $\setOf{F \subseteq M \suchThat F \neq \emptyset, F \text{ finite}}$. The following two statements are equivalent:
    \begin{aenumerate}
      \item\label{item:Folner-condition-for-finitely-right-generated-group-sets:Folner-condition-new}
            For each positive real number $\varepsilon \in \R_{> 0}$, there is an element $F \in \mathcal{F}$ such that
            \begin{equation}
            \label{equation:one-Folner-condition-for-finitely-right-generated-group-sets:Folner-condition-new}
              \ForEach s \in S \Holds \frac{\cardinalityOf{F \smallsetminus (\blank \isSemiActedUponBy s)^{-1}(F)}}{\cardinalityOf{F}} < \varepsilon;
            \end{equation}
      \item\label{item:Folner-condition-for-finitely-right-generated-group-sets:Folner-condition-old}
            For each finite subset $E$ of $G \modulo G_0$ and each positive real number $\varepsilon \in \R_{> 0}$, there is an element $F \in \mathcal{F}$ such that
            \begin{equation}
            \label{equation:one-Folner-condition-for-finitely-right-generated-group-sets:Folner-condition-old}
              \ForEach e \in E \Holds \frac{\cardinalityOf{F \smallsetminus (\blank \isSemiActedUponBy e)^{-1}(F)}}{\cardinalityOf{F}} < \varepsilon. \qedhere
            \end{equation}
    \end{aenumerate}
  \end{lemma}

  \begin{proof}
    First, because $S$ is a finite subset of $G \modulo G_0$, \cref{item:Folner-condition-for-finitely-right-generated-group-sets:Folner-condition-old} implies \cref{item:Folner-condition-for-finitely-right-generated-group-sets:Folner-condition-new}.

    Secondly, let \cref{item:Folner-condition-for-finitely-right-generated-group-sets:Folner-condition-new} hold. Furthermore, let $\varepsilon' \in \R_{> 0}$, let $E$ be a finite subset of $G \modulo G_0$, and identify $M$ with $G \modulo G_0$ by $\iota \givenBy m \mapsto G_{m_0, m}$. Then, according to \cref{lemma:finite-set-contained-in-ball-if-gen-set-contains-neutral-element}, there is a $k \in \N_0$ such that
    \begin{equation*}
      E \subseteq \setOf{m \in M \suchThat \Exists \sequence{s_i'}_{i \in \setOf{1, 2, \dotsc, k}} \text{ in } S' \SuchThat m_0 \isSemiActedUponBy \sequence{s_i'}_{i \in \setOf{1, 2, \dotsc, k}}},
    \end{equation*}
    where $S' = \setOf{G_0} \cup S$. Let $\varepsilon = \varepsilon' / (\cardinalityOf{G_0}^2 \cdot k)$ and let $F \in \mathcal{F}$ such that \cref{equation:Folner-condition-for-finitely-right-generated-group-sets:Folner-condition} holds. Furthermore, let $e \in E$. Then, there is a finite sequence $\sequence{s_i'}_{i \in \setOf{1, 2, \dotsc, k}}$ in $S'$ such that $m_0 \isSemiActedUponBy \sequence{s_i'}_{i \in \setOf{1, 2, \dotsc, k}} = e$. For each $i \in \setOf{0, 1, \dotsc, k}$, let $m_i = m_0 \isSemiActedUponBy \sequence{s_j'}_{j \in \setOf{1, 2, \dotsc, i}}$ and let $F_i = (\blank \isSemiActedUponBy m_i)^{-1}(F)$. Note that $m_k = e$ and that $F_0 = F$. Then, according to \cref{remark:a-setminus-b-leq-a-setminus-c-plus-c-setminus-b},
    \begin{align*}
      \cardinalityOf{F \smallsetminus F_k}
      &=    \cardinalityOf{F_0 \smallsetminus F_k}\\
      &\leq \cardinalityOf{F_0 \smallsetminus F_1} + \cardinalityOf{F_1 \smallsetminus F_k}\\
      &\leq \cardinalityOf{F_0 \smallsetminus F_1} + \cardinalityOf{F_1 \smallsetminus F_2} + \cardinalityOf{F_2 \smallsetminus F_k}\\
      &\leq \dotso\\
      &\leq \sum_{i = 1}^k \cardinalityOf{F_{i - 1} \smallsetminus F_i}.
    \end{align*}
    Let $i \in \setOf{1, 2, \dotsc, k}$. Then, because $m_i = m_{i - 1} \isSemiActedUponBy s_i'$, we have $F_{i - 1} \smallsetminus F_i = (\blank \isSemiActedUponBy m_{i - 1})^{-1}(F) \smallsetminus (\blank \isSemiActedUponBy (m_{i - 1} \isSemiActedUponBy s_i'))^{-1}(F)$. Hence, according to \cref{lemma:liberation-minus-liberation-lib-s-subseteq-bigcup-liberation-of-setminus-liberation-gzero-s}, we have $F_{i - 1} \smallsetminus F_i \subseteq \bigcup_{g_0 \in G_0} (\blank \isSemiActedUponBy m_{i - 1})^{-1}(F \smallsetminus (\blank \isSemiActedUponBy g_0 \cdot s_i')^{-1}(F))$. Therefore,
    \begin{equation*}
      \cardinalityOf{F_{i - 1} \smallsetminus F_i}
      \leq \sum_{g_0 \in G_0} \cardinalityOf{(\blank \isSemiActedUponBy m_{i - 1})^{-1}(F \smallsetminus (\blank \isSemiActedUponBy g_0 \cdot s_i')^{-1}(F))}.
    \end{equation*}
    Thus, according to \cref{lemma:liberation-preimage},
    \begin{equation*}
      \cardinalityOf{F_{i - 1} \smallsetminus F_i} \leq \sum_{g_0 \in G_0} \cardinalityOf{G_0} \cdot \cardinalityOf{F \smallsetminus (\blank \isSemiActedUponBy g_0 \cdot s_i')^{-1}(F)}.
    \end{equation*}
    Hence, because $G_0 \cdot S' \subseteq S'$, $F \smallsetminus (\blank \isSemiActedUponBy G_0)^{-1}(F) = F \smallsetminus F = \emptyset$, and \cref{equation:Folner-condition-for-finitely-right-generated-group-sets:Folner-condition} holds,
    \begin{align*}
      \cardinalityOf{F_{i - 1} \smallsetminus F_i}
      &< \sum_{g_0 \in G_0} \cardinalityOf{G_0} \cdot \varepsilon \cdot \cardinalityOf{F}\\
      &= \cardinalityOf{G_0}^2 \cdot \varepsilon \cdot \cardinalityOf{F}\\
      &= \cardinalityOf{G_0}^2 \cdot \frac{\varepsilon'}{\cardinalityOf{G_0}^2 \cdot k} \cdot \cardinalityOf{F}\\
      &= \frac{\varepsilon'}{k} \cdot \cardinalityOf{F}.
    \end{align*}
    Therefore,
    \begin{equation*}
      \cardinalityOf{F \smallsetminus F_k}
      \leq \sum_{i = 1}^k \cardinalityOf{F_{i - 1} \smallsetminus F_i}
      <    k \cdot \frac{\varepsilon'}{k} \cdot \cardinalityOf{F}
      =    \varepsilon' \cdot \cardinalityOf{F}.
    \end{equation*}
    Thus, because $F_k = (\blank \isSemiActedUponBy e)^{-1}(F)$,
    \begin{equation*}
      \frac{\cardinalityOf{F \smallsetminus (\blank \isSemiActedUponBy e)^{-1}(F)}}{\cardinalityOf{F}} < \varepsilon'. 
    \end{equation*}
    Hence, \cref{equation:one-Folner-condition-for-finitely-right-generated-group-sets:Folner-condition-old} holds. In conclusion, \cref{item:Folner-condition-for-finitely-right-generated-group-sets:Folner-condition-old} holds.
  \end{proof}

  \begin{theorem} 
  \label{theorem:Folner-condition-for-finitely-right-generated-group-sets}
    The following three statements are equivalent:
    \begin{aenumerate}
      \item \label{item:Folner-condition-for-finitely-right-generated-group-sets:right-amenable}
            The cell space $\mathcal{R}$ is right amenable;
      \item \label{item:Folner-condition-for-finitely-right-generated-group-sets:Folner-condition}
            For each positive real number $\varepsilon$, there is a non-empty and finite subset $F$ of $M$ such that
            \begin{equation}
            \label{equation:Folner-condition-for-finitely-right-generated-group-sets:Folner-condition}
              \ForEach s \in S \Holds \frac{\cardinalityOf{F \smallsetminus (\blank \isSemiActedUponBy s)^{-1}(F)}}{\cardinalityOf{F}} < \varepsilon;
            \end{equation}
      \item \label{item:Folner-condition-for-finitely-right-generated-group-sets:isoperimetric-constant}
            The isoperimetric constant $\iota_S(\mathcal{R})$ is $0$. \qedhere
    \end{aenumerate}
  \end{theorem}

  \begin{proof}
    \proofPart{\ref{item:Folner-condition-for-finitely-right-generated-group-sets:right-amenable} $\implies$ \ref{item:Folner-condition-for-finitely-right-generated-group-sets:Folner-condition}}
      Let $\mathcal{R}$ be right amenable. Then, according to \cref{theorem:Tarski-Folner}, there is a right Følner net in $\mathcal{R}$. Hence, according to \cref{theorem:epsilon-characterisation-of-Folner-net}, for each $\varepsilon \in \R_{> 0}$, there is a non-empty and finite subset $F$ of $M$ such that \cref{equation:Folner-condition-for-finitely-right-generated-group-sets:Folner-condition} holds. 

    \proofPart{\ref{item:Folner-condition-for-finitely-right-generated-group-sets:Folner-condition} $\implies$ \ref{item:Folner-condition-for-finitely-right-generated-group-sets:right-amenable}}
      Let \cref{item:Folner-condition-for-finitely-right-generated-group-sets:Folner-condition} hold. Then, according to \cref{lemma:one-Folner-condition-for-finitely-right-generated-group-sets:Folner-condition} and \cref{theorem:epsilon-characterisation-of-Folner-net}, there is a right Følner net in $\mathcal{R}$. In conclusion, according to \cref{theorem:Tarski-Folner}, the cell space $\mathcal{R}$ is right amenable.

    \proofPart{\ref{item:Folner-condition-for-finitely-right-generated-group-sets:Folner-condition} $\implies$ \ref{item:Folner-condition-for-finitely-right-generated-group-sets:isoperimetric-constant}}
      Let \cref{item:Folner-condition-for-finitely-right-generated-group-sets:Folner-condition} hold. Furthermore, let $\varepsilon' \in \R_{> 0}$ and let $\varepsilon = \varepsilon' / \cardinalityOf{S}$. Then, there is a non-empty and finite subset $F$ of $M$ such that \cref{equation:Folner-condition-for-finitely-right-generated-group-sets:Folner-condition} holds. Therefore,
      \begin{equation*}
        \frac{\cardinalityOf{\bigcup_{s \in S} F \smallsetminus (\blank \isSemiActedUponBy s)^{-1}(F)}}{\cardinalityOf{F}}
        \leq \sum_{s \in S} \frac{\cardinalityOf{F \smallsetminus (\blank \isSemiActedUponBy s)^{-1}(F)}}{\cardinalityOf{F}}
        <    \cardinalityOf{S} \cdot \varepsilon
        =    \varepsilon'.
      \end{equation*}
      In conclusion, $\iota_S(\mathcal{R}) = 0$.

    \proofPart{\ref{item:Folner-condition-for-finitely-right-generated-group-sets:isoperimetric-constant} $\implies$ \ref{item:Folner-condition-for-finitely-right-generated-group-sets:Folner-condition}}
      Let $\iota_S(\mathcal{R}) = 0$. Furthermore, let $\varepsilon \in \R_{> 0}$. Then, because $\iota_S(\mathcal{R}) = 0$, there is a non-empty and finite subset $F$ of $M$ such that
      \begin{equation*}
        \frac{\cardinalityOf{\bigcup_{s \in S} F \smallsetminus (\blank \isSemiActedUponBy s)^{-1}(F)}}{\cardinalityOf{F}} < \varepsilon.
      \end{equation*}
      Hence, for each $s \in S$, because $F \smallsetminus (\blank \isSemiActedUponBy s)^{-1}(F) \subseteq \bigcup_{s' \in S} F \smallsetminus (\blank \isSemiActedUponBy s')^{-1}(F)$,
      \begin{equation*}
        \frac{\cardinalityOf{F \smallsetminus (\blank \isSemiActedUponBy s)^{-1}(F)}}{\cardinalityOf{F}} < \varepsilon.
      \end{equation*}
      In conclusion, \cref{equation:Folner-condition-for-finitely-right-generated-group-sets:Folner-condition} holds.
  \end{proof}

  \begin{remark}
    In the case that $M = G$ and $\actsOnPoint$ is the group multiplication of $G$, \cref{theorem:Folner-condition-for-finitely-right-generated-group-sets} is proposition~6.10.2 in \cite{ceccherini-silberstein:coornaert:2010}.
  \end{remark}

  \begin{theorem} 
  \label{theorem:k-boundary-characterisation-of-Folner-net}
    Let $\net{F_i}_{i \in I}$ be a net in $\setOf{F \subseteq M \suchThat F \neq \emptyset, F \text{ finite}}$ indexed by $(I, \leq)$. It is a right Følner net in $\mathcal{R}$ if and only if
    \begin{equation}
    \label{equation:k-boundary-characterisation-of-Folner-net}
      \ForEach \rho \in \N_0 \Holds \lim_{i \in I} \frac{\cardinalityOf{\boundaryOf_\rho F_i}}{\cardinalityOf{F_i}} = 0. \qedhere
    \end{equation}
  \end{theorem}

  \begin{proof}
    First, let $\net{F_i}_{i \in I}$ be a right Følner net in $\mathcal{R}$. Furthermore, let $\rho$ be a non-negative integer. Then, for each index $i \in I$, we have $\boundaryOf_\rho F_i = \boundaryOf_{\ball_S(\rho)} F_i$. And, according to \cref{remark:upper-bound-for-cardinality-of-balls}, the ball $\ball_S(\rho)$ is finite. Hence, according to \cref{theorem:boundary-characterisation-of-Folner-net},
    \begin{equation*}
      \lim_{i \in I} \frac{\cardinalityOf{\boundaryOf_{\ball_S(\rho)} F_i}}{\cardinalityOf{F_i}} = 0.
    \end{equation*}
    In conclusion, \cref{equation:k-boundary-characterisation-of-Folner-net} holds.

    Secondly, let \cref{equation:k-boundary-characterisation-of-Folner-net} hold. Furthermore, let $N$ be a finite subset of $G \modulo G_0$. Then, according to \cref{remark:ball-of-radius-0-contains-one-element-and-sequence-of-balls-is-monotonic}, there is a non-negative integer $\rho$ such that $N \subseteq \ball_S(\rho)$. Hence, for each index $i \in I$, according to \cref{item:properties-of-interior-closure-and-boundary:inclusions} of \cref{lemma:properties-of-interior-closure-and-boundary},
    we have $\boundaryOf_N F_i \subseteq \boundaryOf_{\ball_S(\rho)} F_i = \boundaryOf_\rho F_i$. Therefore,
    \begin{equation*}
      \lim_{i \in I} \frac{\cardinalityOf{\boundaryOf_N F_i}}{\cardinalityOf{F_i}} = 0.
    \end{equation*}
    In conclusion, according to \cref{theorem:boundary-characterisation-of-Folner-net}, the net $\net{F_i}_{i \in I}$ is a right Følner net in $\mathcal{R}$.
  \end{proof}

  \section{Sub-Exponential Growth and Amenability}
  \label{section:sub-exponential-growth-and-amenability}

  In this section, let $\mathcal{M} = \ntuple{M, G, \actsOnPoint}$ be a finitely right-gen\-er\-at\-ed left-ho\-mo\-ge\-neous space with finite stabilisers. Note that it is even finitely and symmetrically right generated. 

  \paragraph{Contents.} In \cref{theorem:sub-exp-implies-amenable} we show that if $\mathcal{M}$ has sub-exponential growth, then it is right amenable and a sequence of balls is a right Følner net. And in \cref{corollary:group-sub-exp-implies-cell-space-sub-exp} we show that if $G$ has sub-exponential growth, then so has $\mathcal{M}$.

  \begin{lemma}
  \label{lemma:liminf-of-frac-is-liminf-of-sqrt}
    Let $\sequence{r_k}_{k \in \N_0}$ be a sequence of positive real numbers. Then,
    \begin{equation*}
      \liminf_{k \to \infty} \frac{r_{k + 1}}{r_k} \leq \liminf_{k \to \infty} \sqrt[k]{r_k}. \qedhere 
    \end{equation*}
  \end{lemma}

  \begin{proof}
    See lemma~6.11.1 in \cite{ceccherini-silberstein:coornaert:2010}.
  \end{proof}

  \begin{main-theorem} 
  \label{theorem:sub-exp-implies-amenable}
    Let the space $\mathcal{M}$ have sub-exponential growth. Then, it is right amenable. And, for each coordinate system $\mathcal{K}$ for $\mathcal{M}$ and each finite and symmetric right-gen\-er\-at\-ing set $S$ of $\mathcal{R} = \ntuple{\mathcal{M}, \mathcal{K}}$, there is a subsequence of $\sequence{\ball_S(\rho)}_{\rho \in \N_0}$ that is a right Følner net in $\mathcal{R}$.
  \end{main-theorem}

  \begin{proof}
    Let $\mathcal{K}$ be a coordinate system for $\mathcal{M}$ and let $S$ be a finite and symmetric right-gen\-er\-at\-ing set of $\mathcal{R} = \ntuple{\mathcal{M}, \mathcal{K}}$. Then, according to \cref{lemma:liminf-of-frac-is-liminf-of-sqrt} and \cref{corollary:growth-rate-equal-to-one-if-and-only-if-subexp-growth},
    \begin{equation*}
      1 \leq \liminf_{k \to \infty} \frac{\gamma_S(k + 1)}{\gamma_S(k)}
        \leq \lim_{k \to \infty} \sqrt[k]{\gamma_S(k)}
        =    \lambda_S
        =    1.
    \end{equation*}
    Hence, $\liminf_{k \to \infty} \frac{\gamma_S(k + 1)}{\gamma_S(k)} = 1$. 
    Let $\varepsilon \in \R_{> 0}$. Then, there is a $k \in \N_+$ such that $\frac{\gamma_S(k)}{\gamma_S(k - 1)} < 1 + \varepsilon$. Thus, $\gamma_S(k) - \gamma_S(k - 1) < \varepsilon \cdot \gamma_S(k - 1)$. 
    Let $s \in S$. Then, according to \cref{lemma:ball-liberation-included-in-ball-one-larger}, we have $\ball_S(k - 1) \subseteq (\blank \isSemiActedUponBy s)^{-1}(\ball_S(k))$. Hence, because $\ball_S(k - 1) \subseteq \ball_S(k)$ and $\gamma_S(k - 1) \leq \gamma_S(k)$,
    \begin{align*}
      \cardinalityOf{\ball_S(k) \smallsetminus (\blank \isSemiActedUponBy s)^{-1}(\ball_S(k))}
      &\leq \cardinalityOf{\ball_S(k) \smallsetminus \ball_S(k - 1)}\\
      &=    \cardinalityOf{\ball_S(k)} - \cardinalityOf{\ball_S(k - 1)}\\
      &=    \gamma_S(k) - \gamma_S(k - 1)\\
      &<    \varepsilon \cdot \gamma_S(k - 1)\\
      &\leq \varepsilon \cdot \gamma_S(k)\\
      &=    \varepsilon \cdot \cardinalityOf{\ball_S(k)}.
    \end{align*}
    Therefore, for each $\varepsilon \in \R_{> 0}$, there is a $k \in \N_+$ such that
    \begin{equation}
    \label{equation:liminf-of-frac-is-liminf-of-sqrt:varepsilon}
      \ForEach s \in S \Holds \frac{\cardinalityOf{\ball_S(k) \smallsetminus (\blank \isSemiActedUponBy s)^{-1}(\ball_S(k))}}{\cardinalityOf{\ball_S(k)}} < \varepsilon.
    \end{equation}
    In conclusion, according to \cref{theorem:Folner-condition-for-finitely-right-generated-group-sets}, the cell space $\mathcal{R}$ is right amenable and hence so is $\mathcal{M}$. And, by going through the proof of \cref{theorem:Folner-condition-for-finitely-right-generated-group-sets} and the proofs of the lemmata, corollaries, and theorems it uses, one can see that a subsequence of $\sequence{\ball_S(\rho)}_{\rho \in \N_0}$ is a right Følner net in $\mathcal{R}$. As this is rather tedious, we construct such a subsequence in the following.

    If $M$ is finite, then, according to \cref{corollary:finite-if-and-only-if-eventually-equal-to-m}, the sequence $\sequence{\ball(\rho)}_{\rho \in \N_0}$ is eventually equal to $M$. Hence, because
    \begin{equation*}
      \ForEach \mathfrak{g} \in G \modulo G_0 \Holds \frac{\cardinalityOf{M \smallsetminus (\blank \isSemiActedUponBy \mathfrak{g})^{-1}(M)}}{\cardinalityOf{M}} = 0,
    \end{equation*}
    the sequence $\sequence{\ball(\rho)}_{\rho \in \N_0}$ is a right Følner net in $\mathcal{R}$.

    From now on, let $M$ be infinite. We showed above that, for each $\varepsilon \in \R_{> 0}$, there is a $k \in \N_+$ such that \cref{equation:liminf-of-frac-is-liminf-of-sqrt:varepsilon} holds. Hence, according to \cref{lemma:one-Folner-condition-for-finitely-right-generated-group-sets:Folner-condition}, under the identification of $M$ with $G \modulo G_0$, for each $j \in \N_+$ and each $\varepsilon \in \R_{> 0}$, there is a $k \in \N_+$ such that
    \begin{equation*}
      \ForEach e \in \ball_S(j) \Holds \frac{\cardinalityOf{\ball_S(k) \smallsetminus (\blank \isSemiActedUponBy e)^{-1}(\ball_S(k))}}{\cardinalityOf{\ball_S(k)}} < \varepsilon.
    \end{equation*}
    Therefore, for each $n \in \N_+$, there is a least $k_n \in \N_+$ such that
    \begin{equation*}
      \ForEach e \in \ball_S(n) \Holds \frac{\cardinalityOf{\ball_S(k_n) \smallsetminus (\blank \isSemiActedUponBy e)^{-1}(\ball_S(k_n))}}{\cardinalityOf{\ball_S(k_n)}} < \frac{1}{n}.
    \end{equation*}
    By the choices of $k_n$, for $n \in \N_+$, the sequence $\sequence{k_n}_{n \in \N_+}$ is non-decreasing. If it is not eventually constant, then, by skipping duplicate entries, we get an increasing subsequence $\sequence{k_{n_i}}_{i \in \N_+}$.

    Otherwise, there is an $n \in \N_+$ such that, for each $e \in \ball_S(1)$, we have $\ball_S(k_n) \smallsetminus (\blank \isSemiActedUponBy e)^{-1}(\ball_S(k_n)) = \emptyset$ and thus $\ball_S(k_n) \isSemiActedUponBy e \subseteq \ball_S(k_n)$. Hence, $\ball_S(k_n + 1) = \ball_S(k_n) \isSemiActedUponBy \ball_S(1) \subseteq \ball_S(k_n)$. Therefore, according to \cref{corollary:m-infinite-if-and-only-if-balls-strictly-increasing}, the set $M$ is finite, which contradicts our assumption that $M$ is infinite. Hence, the case that $\sequence{k_n}_{n \in \N_+}$ is eventually constant does not occur.

    In the case that does occur, for each $\mathfrak{g} \in G \modulo G_0$ and each $\varepsilon \in \R_{> 0}$, there is an $i_0 \in \N_+$ such that $\mathfrak{g} \in \ball_S(n_{i_0})$ and $1/n_{i_0} \leq \varepsilon$, and hence, for each $i \in \N_+$ with $i \geq i_0$,
    \begin{equation*}
      \frac{\cardinalityOf{\ball_S(k_{n_i}) \smallsetminus (\blank \isSemiActedUponBy \mathfrak{g})^{-1}(\ball_S(k_{n_i}))}}{\cardinalityOf{\ball_S(k_{n_i})}} < \frac{1}{n_i} \leq \frac{1}{n_{i_0}} \leq \varepsilon.
    \end{equation*}
    Therefore,
    \begin{equation*}
      \ForEach \mathfrak{g} \in G \modulo G_0 \Holds
          \lim_{i \to \infty} \frac{\cardinalityOf{\ball_S(k_{n_i}) \smallsetminus (\blank \isSemiActedUponBy \mathfrak{g})^{-1}(\ball_S(k_{n_i}))}}{\cardinalityOf{\ball_S(k_{n_i})}} = 0.
    \end{equation*}
    In conclusion, the sequence $\sequence{\ball_S(k_{n_i})}_{i \in \N_+}$ is a right Følner net in $\mathcal{R}$.
  \end{proof} 

  \begin{remark}
    In the case that $M = G$ and $\actsOnPoint$ is the group multiplication of $G$, the first part of \cref{theorem:sub-exp-implies-amenable} is theorem~6.11.2 in \cite{ceccherini-silberstein:coornaert:2010}.
  \end{remark}

  \begin{corollary}
  \label{corollary:group-sub-exp-implies-cell-space-sub-exp}
    Let the group $G$ be finitely generated and let it have sub-exponential growth. The space $\mathcal{M}$ has sub-exponential growth and is right amenable.
  \end{corollary}

  \begin{proof}
    This is a direct consequence of \cref{lemma:group-dominates-cell-space-growth-rate} and \cref{theorem:sub-exp-implies-amenable}.
  \end{proof}

  \begin{example} 
    Each quotient set of any virtually nilpotent group (see \cref{example:virtually-nilpotent-group}) or the Grigorchuk group (see \cref{example:Grigorchuck-group}) by a finite subgroup acted upon by left multiplication is a right-a\-me\-na\-ble left-ho\-mo\-ge\-neous space.
  \end{example}


  \clearToOddPage
  \chapter{Shift Spaces and the Moore and the Myhill Properties}
  \label{chapter:Moore}

  \paragraph{Abstract.} We prove the Moore and the Myhill property for strongly irreducible subshifts over right-a\-me\-na\-ble and finitely right-gen\-er\-at\-ed left-ho\-mo\-ge\-neous spaces with finite stabilisers. Both properties together mean that the global transition function of each big-cellular automaton with finite set of states and finite neighbourhood over such a subshift is surjective if and only if it is pre-injective. This statement is known as Garden of Eden theorem. Pre-Injectivity means that two global configurations that differ at most on a finite subset and have the same image under the global transition function must be identical.

  \paragraph{Remark.} This chapter generalises parts of the paper \enquote{\citetitle*{fiorenzi:2003}}\cite{fiorenzi:2003}. 

  \paragraph{Summary.} A subset $X$ of the phase space $Q^M$, where $Q$ is a finite set of states, is a \emph{shift space of finite type} if it is generated by a finite set of forbidden blocks. Such a space $X$ shift-invariant and compact. And it is \emph{strongly irreducible} if each pair of finite patterns that are allowed in $X$ and at least some fixed positive integer apart, are embedded in a point of $X$. A map $\Delta$ from a shift space $X$ to a shift space $Y$ is \emph{local} if the state $\Delta(x)(m)$ is uniformly and locally determined in $m$, in other words, if the map $\Delta$ is the restriction of the global transition function of a big-cellular automaton with finite neighbourhood to the domain $X$ and the codomain $Y$. 

  For a right-a\-me\-na\-ble and finitely right-gen\-er\-at\-ed cell space with finite stabilisers we may choose a right Følner net $\mathcal{F} = \family{F_i}_{i \in I}$. The \emph{entropy} of a subset $X$ of $Q^M$ with respect to $\mathcal{F}$, where $Q$ is a finite set, is, broadly speaking, the asymptotic growth rate of the number of finite patterns with domain $F_i$ that occur in $X$. For non-negative integers $\theta$, $\kappa$, and $\theta'$, a \emph{$\ntuple{\theta, \kappa, \theta'}$-tiling} is a subset $T$ of $M$ such that $\family{\ball(t, \theta)}_{t \in T}$ is pairwise at least $\kappa + 1$ apart and $\family{\ball(t, \theta')}_{t \in T}$ is a cover of $M$. If for each point $t \in T$ not all patterns with domain $\ball(t, \theta)$ occur in a subset of $Q^M$, then that subset does not have maximal entropy. 

  A local map from a non-empty strongly irreducible shift space of finite type to a strongly irreducible shift space with the same entropy over a right-a\-me\-na\-ble and finitely right-gen\-er\-at\-ed cell space with finite stabilisers is surjective if and only if its image has maximal entropy and its image has maximal entropy if and only if it is pre-injective. This establishes the Garden of Eden theorem, which states that a local map as above is surjective if and only if it is pre-injective. This answers a question posed by Sébastien Moriceau at the end of his paper \enquote{\citetitle*{moriceau:2011}}\cite{moriceau:2011}. And it follows that strongly irreducible shift spaces of finite type over right-a\-me\-na\-ble and finitely right-gen\-er\-at\-ed cell spaces have the Moore and the Myhill property.

  The Garden of Eden theorem for cellular automata over $\Z^2$ is a famous theorem by Edward Forrest Moore and John R. Myhill from 1962 and 1963, which was proved in their papers \enquote{\citetitle*{moore:1962}}\cite{moore:1962} and \enquote{\citetitle*{myhill:1963}}\cite{myhill:1963}. That theorem also holds for cellular automata over amenable finitely generated groups, which was proved by Tullio Ceccherini-Silberstein, Antonio Machi, and Fabio Scarabotti in their paper \enquote{\citetitle*{ceccherini-silberstein:machi:scarabotti:1999}}\cite{ceccherini-silberstein:machi:scarabotti:1999}. It even holds for such automata on strongly irreducible shifts of finite type, which was proved by Francesca Fiorenzi in her paper \enquote{\citetitle*{fiorenzi:2003}}\cite{fiorenzi:2003}.

  \paragraph{Contents.} In \cref{section:shift-spaces} we introduce full shifts, shift-invariance, shift spaces or subshifts (of finite type), strong irreducibility, bounded propagation, local maps, conjugacies, and the Moore and the Myhill property. In \cref{section:apart-tilings} we introduce tilings, prove their existence, and relate them to entropies. And in \cref{section:gardens-of-Eden-on-shifts} we prove the Garden of Eden theorem, from which we deduce that both the Moore and the Myhill property hold.

  \paragraph{Preliminary Notions.} A \emph{left group set} is a triple $\ntuple{M, G, \actsOnPoint}$, where $M$ is a set, $G$ is a group, and $\actsOnPoint$ is a map from $G \times M$ to $M$, called \emph{left group action of $G$ on $M$}, such that $G \to \symmetricGroupOf(M)$, $g \mapsto [g \actsOnPoint \blank]$, is a group homomorphism. The action $\actsOnPoint$ is \emph{transitive} if $M$ is non-empty and for each $m \in M$ the map $\blank \actsOnPoint m$ is surjective; and \emph{free} if for each $m \in M$ the map $\blank \actsOnPoint m$ is injective. For each $m \in M$, the set $G \actsOnPoint m$ is the \emph{orbit of $m$}, the set $G_m = (\blank \actsOnPoint m)^{-1}(m)$ is the \emph{stabiliser of $m$}, and, for each $m' \in M$, the set $G_{m, m'} = (\blank \actsOnPoint m)^{-1}(m')$ is the \emph{transporter of $m$ to $m'$}.

  A \emph{left-ho\-mo\-ge\-neous space} is a left group set $\mathcal{M} = \ntuple{M, G, \actsOnPoint}$ such that $\actsOnPoint$ is transitive. A \emph{coordinate system for $\mathcal{M}$} is a tuple $\mathcal{K} = \ntuple{m_0, \family{g_{m_0, m}}_{m \in M}}$, where $m_0 \in M$ and for each $m \in M$ we have $g_{m_0, m} \actsOnPoint m_0 = m$. The stabiliser $G_{m_0}$ is denoted by $G_0$. The tuple $\mathcal{R} = \ntuple{\mathcal{M}, \mathcal{K}}$ is a \emph{cell space}. The map $\isSemiActedUponBy \from M \times G \modulo G_0 \to M$, $(m, g G_0) \mapsto g_{m_0, m} g g_{m_0, m}^{-1} \actsOnPoint m\ (= g_{m_0, m} g \actsOnPoint m_0)$ is a \emph{right semi-action of $G \modulo G_0$ on $M$ with defect $G_0$}, which means that
  \begin{equation*}
    \ForEach m \in M \Holds m \isSemiActedUponBy G_0 = m,
  \end{equation*}
  and
  \begin{multline*}
    \ForEach m \in M \ForEach g \in G \Exists g_0 \in G_0 \SuchThat \ForEach \mathfrak{g}' \in G \modulo G_0 \Holds\\
          m \isSemiActedUponBy g \cdot \mathfrak{g}' = (m \isSemiActedUponBy g G_0) \isSemiActedUponBy g_0 \cdot \mathfrak{g}'.
  \end{multline*}
  It is \emph{transitive}, which means that the set $M$ is non-empty and for each $m \in M$ the map $m \isSemiActedUponBy \blank$ is surjective; and \emph{free}, which means that for each $m \in M$ the map $m \isSemiActedUponBy \blank$ is injective; and \emph{semi-commutes with $\actsOnPoint$}, which means that
  \begin{multline*}
    \ForEach m \in M \ForEach g \in G \Exists g_0 \in G_0 \SuchThat \ForEach \mathfrak{g}' \in G \modulo G_0 \Holds\\
          (g \actsOnPoint m) \isSemiActedUponBy \mathfrak{g}' = g \actsOnPoint (m \isSemiActedUponBy g_0 \cdot \mathfrak{g}').
  \end{multline*}
  The maps $\iota \from M \to G \modulo G_0$, $m \mapsto G_{m_0, m}$, and $m_0 \isSemiActedUponBy \blank$ are inverse to each other. Under the identification of $M$ with $G \modulo G_0$ by either of these maps, we have $\isSemiActedUponBy \from (m, \mathfrak{g}) \mapsto g_{m_0, m} \actsOnPoint \mathfrak{g}$. (See \cref{chapter:automata}.)

  A left-ho\-mo\-ge\-neous space $\mathcal{M}$ is \emph{right amenable} if there is a coordinate system $\mathcal{K}$ for $\mathcal{M}$ and there is a finitely additive probability measure $\mu$ on $M$ such that 
  \begin{equation*}
    \ForEach \mathfrak{g} \in G \modulo G_0 \ForEach A \subseteq M \Holds \parens[\big]{(\blank \isSemiActedUponBy \mathfrak{g})\restrictedTo_A \text{ injective} \implies \mu(A \isSemiActedUponBy \mathfrak{g}) = \mu(A)},
  \end{equation*}
  in which case the cell space $\mathcal{R} = \ntuple{\mathcal{M}, \mathcal{K}}$ is called \emph{right amenable}. When the stabiliser $G_0$ is finite, that is the case if and only if there is a \emph{right Følner net in $\mathcal{R}$ indexed by $(I, \leq)$}, which is a net $\net{F_i}_{i \in I}$ in $\setOf{F \subseteq M \suchThat F \neq \emptyset, F \text{ finite}}$ such that
  \begin{equation*} 
    \ForEach \rho \in \N_0 \Holds \lim_{i \in I} \frac{\cardinalityOf{\boundaryOf_\rho F_i}}{\cardinalityOf{F_i}} = 0.
  \end{equation*}
  If a net is a right Følner net for one coordinate system, then it is a right Følner net for each coordinate system. In particular, a left-ho\-mo\-ge\-neous space $\mathcal{M}$ with finite stabilisers is right amenable if and only if, for each coordinate system $\mathcal{K}$ for $\mathcal{M}$, the cell space $\ntuple{\mathcal{M}, \mathcal{K}}$ is right amenable. (See \cref{chapter:amenability,chapter:growth}.)

  A left-ho\-mo\-ge\-neous space $\mathcal{M}$ is \emph{finitely right generated} if there is a coordinate system $\mathcal{K}$ for $\mathcal{M}$ and there is a finite subset $S$ of $G \modulo G_0$ such that $G_0 \cdot S \subseteq S$ and, for each $m \in M$, there is a $k \in \N_0$ and there is a $\sequence{s_i}_{i \in \setOf{1, 2, \dotsc, k}}$ in $S \cup S^{-1}$, where $S^{-1} = \setOf{g^{-1} G_0 \suchThat s \in S, g \in s}$, such that
  \begin{equation*}
    m = \parens[\Big]{\parens[\big]{(m_0 \isSemiActedUponBy s_1) \isSemiActedUponBy s_2} \dotsb} \isSemiActedUponBy s_k.
  \end{equation*}
  in which case the cell space $\mathcal{R} = \ntuple{\mathcal{M}, \mathcal{K}}$ is called \emph{finitely right generated}. The left-ho\-mo\-ge\-neous space $\mathcal{M}$ is finitely right generated if and only if, for each coordinate system $\mathcal{K}$ for $\mathcal{M}$, the cell space $\ntuple{\mathcal{M}, \mathcal{K}}$ is finitely right generated. The right-gen\-er\-at\-ing set $S$ is \emph{symmetric} if $S^{-1} \subseteq S$. The $S$-edge-labelled directed graph $\ntuple{M, E}$, where $E = \setOf{(m, s, m \isSemiActedUponBy s) \suchThat m \in M, s \in S}$, is the \emph{coloured $S$-Cayley graph}. The distance $\distanceOf_S$ on that graph is the \emph{$S$-metric} and the map $\lengthOf{\blank}_S = \distanceOf_S(m_0, \blank)$ is the \emph{$S$-length}. For each $m \in M$ and each $\rho \in \Z$, the sets
  \begin{align*}
    \ball_S(m, \rho) &= \setOf{m' \in M \suchThat \distanceOf_S(m, m') \leq \rho},\\
    \sphere_S(m, \rho) &= \setOf{m' \in M \suchThat \distanceOf_S(m, m') = \rho}
  \end{align*}
  are the \emph{ball/sphere of radius $\rho$ centred at $m$}, the ball $\ball_S(m_0, \rho)$ is denoted by $\ball_S(\rho)$, and the sphere $\sphere_S(m_0, \rho)$ by $\sphere_S(\rho)$. For each $A \subseteq M$, each $\theta \in \N_0$, the set $A^{-\theta} = \setOf{m \in A \suchThat \ball_S(m, \theta) \subseteq A}$ is the \emph{$\theta$-interior of $A$}, the set $A^{+\theta} = \setOf{m \in M \suchThat \ball_S(m, \theta) \cap A \neq \emptyset}$ is the \emph{$\theta$-closure of $A$}, the set $\boundaryOf_\theta A = A^{+\theta} \smallsetminus A^{-\theta}$ is the \emph{$\theta$-boundary of $A$}, the set $\boundaryOf_\theta^- A = A \smallsetminus A^{-\theta}$ is the \emph{internal $\theta$-boundary of $A$}, and the set $\boundaryOf_\theta^+ A = A^{+\theta} \smallsetminus A$ is the \emph{external $\theta$-boundary of $A$}. (See \cref{chapter:growth}.)

  A \emph{semi-cellular automaton} is a quadruple $\mathcal{C} = \ntuple{\mathcal{R}, Q, N, \delta}$, where $\mathcal{R}$ is a cell space; $Q$, called \emph{set of states}, is a set; $N$, called \emph{neighbourhood}, is a subset of $G \modulo G_0$ such that $G_0 \cdot N \subseteq N$; and $\delta$, called \emph{local transition function}, is a map from $Q^N$ to $Q$. A \emph{local configuration} is a map $\ell \in Q^N$, a \emph{global configuration} is a map $c \in Q^M$, an \emph{$A$-pattern} is a map $p \in Q^A$, where $A$ is a subset of $M$, the number $\sizeOf{p} = \cardinalityOf{A}$ is the \emph{size of $p$}, a finite pattern is a \emph{block}, and the set of blocks is denoted by $Q^*$. The stabiliser $G_0$ acts on $Q^N$ on the left by $\bullet \from G_0 \times Q^N \to Q^N$, $(g_0, \ell) \mapsto [n \mapsto \ell(g_0^{-1} \cdot n)]$, and the group $G$ acts on the set of patterns on the left by
  \begin{align*}
    \actsOnMap \from G \times \bigcup_{A \subseteq M} Q^A &\to \bigcup_{A \subseteq M} Q^A,\\
    (g, p) &\mapsto \left[
                      \begin{aligned}
                        g \actsOnPoint \domainOf(p) &\to Q,\\
                        m &\mapsto p(g^{-1} \actsOnPoint m).
                      \end{aligned}
                    \right]
  \end{align*}
  The \emph{global transition function of $\mathcal{C}$} is the map $\Delta \from Q^M \to Q^M$, $c \mapsto [m \mapsto \delta(n \mapsto c(m \isSemiActedUponBy n))]$.

  A subgroup $H$ of $G$ is \emph{$\mathcal{K}$-big} if the set $\setOf{g_{m_0, m} \suchThat m \in M}$ is included in $H$. A \emph{big-cellular automaton} is a semi-cellular automaton $\mathcal{C} = \ntuple{\mathcal{R}, Q, N, \delta}$ such that, for some $\mathcal{K}$-big subgroup $H$ of $G$, the local transition function $\delta$ is \emph{$\bullet_{G_0 \cap H}$-invariant}, which means that, for each $h_0 \in G_0 \cap H$, we have $\delta(h_0 \bullet \blank) = \delta(\blank)$. Its global transition function is $\actsOnMap_H$-e\-qui\-var\-i\-ant, which means that, for each $h \in H$, we have $\Delta(h \actsOnMap \blank) = h \actsOnMap \Delta(\blank)$.

  Under the identification of $M$ with $G \modulo G_0$ by $\iota \givenBy m \mapsto G_{m_0, m}$, the map
  \begin{align*}
    \actsByItsCoordinateOn \from M \times \bigcup_{A \subseteq M} Q^A &\to \bigcup_{A \subseteq M} Q^A,\\
    (m, p) &\mapsto \left[
                      \begin{aligned}
                        m \isSemiActedUponBy \domainOf(p) &\to Q,\\
                        m \isSemiActedUponBy a &\mapsto p(a),
                      \end{aligned}
                    \right]
  \end{align*}
  broadly speaking, maps a point $m$ and a pattern $p$ that is centred at $m_0$ to the corresponding pattern centred at $m$. For each cell $m \in M$, each subset $A$ of $M$, and each pattern $p \in Q^A$, we have $m \actsByItsCoordinateOn p = g_{m_0, m} \actsOnMap p$. It follows that the global transition function $\Delta$ of a big-cellular automaton is $\actsByItsCoordinateOn$-e\-qui\-var\-i\-ant, which means that, for each $m \in M$, we have $\Delta(m \actsByItsCoordinateOn \blank) = m \actsByItsCoordinateOn \Delta(\blank)$. (See \cref{chapter:automata,chapter:garden}.)

  \paragraph{Context.} In this chapter, let $\mathcal{R} = \ntuple{\mathcal{M}, \mathcal{K}} = \ntuple{\ntuple{M, G, \actsOnPoint}, \ntuple{m_0, \family{g_{m_0, m}}_{m \in M}}}$ be a finitely right-gen\-er\-at\-ed cell space such that the stabiliser $G_0$ of $m_0$ under $\actsOnPoint$ is finite; let $S$ be a finite and symmetric right-gen\-er\-at\-ing set of $\mathcal{R}$; let $H$ be a $\mathcal{K}$-big subgroup of $G$; let $H_0$ be the stabiliser of $m_0$ under $\actsOnPoint\restrictedTo_{H \times M}$, which is $H \cap G_0$; for each cell $m \in M$, let $H_{m_0, m}$ be the transporter of $m_0$ to $m$ under $\actsOnPoint\restrictedTo_{H \times M}$; let $Q$ be a finite set; let $Q^M$ be equipped with the prodiscrete topology; for each subset $A$ of $M$, let $\pi_A$ be the restriction map $Q^M \to Q^A$, $c \mapsto c\restrictedTo_A$; and identify $M$ with $G \modulo G_0$ by $\iota \givenBy m \mapsto G_{m_0, m}$. Moreover, we omit the subscript $S$, in particular, instead of $\distanceOf_S$ we write $\distanceOf$, instead of $\lengthOf{\blank}_S$ we write $\lengthOf{\blank}$, instead of $\ball_S$ we write $\ball$, and instead of $\sphere_S$ we write $\sphere$. 

  \section{Shift Spaces} 
  \label{section:shift-spaces}

  \paragraph{Contents.} The full shift is the set of global configurations, the points of the full shift (see \cref{definition:full-shift}). A subset of the full shift is shift-invariant if it is invariant under a group that contains the coordinates (see \cref{definition:shift-invariant}). For example the full shift is shift-invariant. A pattern semi-occurs in another pattern if a rotation of it occurs in the other pattern (see \cref{definition:occurs-semi-occurs,remark:semi-occurs-at,remark:semi-occurs}). It is allowed in a subset of the full shift if it semi-occurs in one of its points and it is forbidden otherwise (see \cref{definition:allowed-forbidden}). The set of points of the full shift in which each block of a given set is forbidden is generated by that set (see \cref{definition:generated-by-forbidden-blocks}) and it is a shift space (see \cref{definition:shift-space}), which is shift-invariant (see \cref{lemma:shift-space-is-shift-invariant}), closed (see \cref{lemma:if-restrictions-to-balls-are-in-shift-then-so-is-pattern}), and compact (see \cref{lemma:shift-space-if-and-only-if-invariant-and-compact}).

  If there is a finite generating set, then the shift space is of finite type (see \cref{definition:of-finite-type}), the radius $\kappa$ of a ball that includes the domains of the patterns of the generating set is its memory and it itself is called $\kappa$-step (see \cref{definition:step-and-memory,lemma:of-finite-type-implies-step}), and its points are characterised by restrictions to balls with its memory as radius (see \cref{lemma:characterisation-of-kappa-step-subshifts}). Finitely many points of shift spaces of finite type can be, in various ways, cut into pieces and glued together to construct new points, as long as the pieces agree on a big enough boundary of the cuts (see \cref{lemma:overlapping-global-configurations-can-be-glued,corollary:overlapping-global-configurations-can-be-glued,corollary:overlapping-patterns-can-be-glued}). A shift space is strongly irreducible if allowed finite patterns that are at least a certain distance apart can be embedded in the same point (see \cref{definition:kappa-strongly-irreducible,definition:strongly-irreducible}). It has bounded propagation if finite patterns are allowed whenever their restrictions to balls of a certain radius are allowed (see \cref{definition:rho-bounded-propagation,definition:bounded-propagation}). Such spaces are strongly irreducible and of finite type (see \cref{lemma:bounded-propagation-implies-strong-irreducibility-and-finite-type}).

  A map from a shift space to another one is local if it is uniformly and locally determined in each cell (see \cref{definition:kappa-local-map,definition:local-map}). Such maps are global transition functions of big-cellular automata with finite neighbourhoods (see \cref{remark:local-map-if-and-only-if-cellular-automaton}), their domains and codomains can be simultaneously restricted to a subset of cells and its interior (see \cref{definition:restriction-of-local-map}), and their images are shift spaces (see \cref{lemma:image-of-local-map-is-subshift}). The difference of two points of the full shift is the set of cells in which they differ (see \cref{definition:difference-of-two-points-of-the-full-shift}) and a local map is pre-injective if it is injective on points with finite support (see \cref{definition:pre-injective-for-local-maps}). A local map that has a local inverse is a conjugacy (see \cref{definition:conjugacy}), and its domain and codomain are conjugate (see \cref{definition:conjugate}). Entropy is invariant under conjugacy (see \cref{lemma:entropy-invariant-under-conjugacy}). A subshift has the Moore property if each surjective local map is pre-injective, and the Myhill property if the converse holds (see \cref{definition:Moore-and-Myhill-properties}). Both these properties are invariant under conjugacy (see \cref{remark:Moore-and-Myhill-are-invariant-under-conjugacy}).


  \begin{definition}
  \label{definition:full-shift}
    The set $Q^M$ is called \define{full shift}\graffito{full shift $Q^M$}\index{shift!full} and each element $c \in Q^M$ is called \define{point}\graffito{point $c$}.
  \end{definition}

  \begin{example}[{\cite[Definition~1.1.1]{lind:marcus:1995}}]
  \label{example:full-shift-over-integers}
    Let $\mathcal{M}$ be the left-ho\-mo\-ge\-neous space $\ntuple{\Z, \Z, +}$, let $\mathcal{K}$ be the coordinate system $\ntuple{0, \family{z}_{z \in \Z}}$, let $\mathcal{R}$ be the cell space $\ntuple{\mathcal{M}, \mathcal{K}}$, let $S$ be the set $\setOf{-1, 1}$, let $H$ be the only $\mathcal{K}$-big subgroup $\Z$ of $\Z$, and let $Q$ be the binary set $\setOf{0, 1}$. The stabiliser $\Z_0$ of $0$ under $+$ is the singleton set $\setOf{0}$; under the identification of $\Z$ with $\Z \modulo \Z_0$ by $z \mapsto z + \Z_0$, the right semi-action of $\Z \modulo \Z_0$ on $\Z$ is but $+$; and the set $S$ is a finite and symmetric right-gen\-er\-at\-ing set of $\mathcal{R}$. The full shift $Q^\Z$ is the usual full $2$-shift considered in symbolic dynamics and its points are called \define{bi-infinite binary sequences}\graffito{bi-infinite binary sequences $Q^\Z$}.
  \end{example}

  \begin{remark}
    There is a bijective map $\varphi$ from $Q$ to $\Z_{\cardinalityOf{Q}}\ (= \setOf{0, 1, \dotsc,\allowbreak \cardinalityOf{Q} - 1})$. It induces the bijective map
    \begin{align*}
      \Phi \from Q^M &\to \Z_{\cardinalityOf{Q}}^M,\\
      c &\mapsto \big[m \mapsto \varphi\parens[\big]{c(m)}\big]. \qedhere
    \end{align*}
  \end{remark}

  \begin{definition}
  \label{definition:shift-invariant}
    Let $X$ be a subset of $Q^M$. It is called \define{shift-invariant}\graffito{shift-invariant} if and only if
    \begin{equation*}
      \ForEach h \in H \Holds h \actsOnMap X = X. \qedhere
    \end{equation*}
  \end{definition}

  \begin{remark}
    Shift-invariance is the same as $\actsOnMap_H$-invariance.
  \end{remark}

  \begin{remark}
  \label{remark:subseteq-sufficient-for-shift-invariance}
    The set $X$ is shift-invariant if and only if
    \begin{equation*}
      \ForEach h \in H \Holds h \actsOnMap X \subseteq X. \qedhere
    \end{equation*}
  \end{remark}

  \begin{remark}
  \label{remark:shift-invariance-induces-invariance-under-right-semi-action}
    Let $X$ be a shift-invariant subset of $Q^M$. Then,
    \begin{equation*}
      \ForEach m \in M \Holds m \actsByItsCoordinateOn X = X. \qedhere
    \end{equation*}
  \end{remark}

  A pattern occurs in another pattern if a translation of it coincides with a subpattern the other pattern and it semi-occurs in another pattern if a rotation and translation of it coincides with a subpattern of the other pattern, as defined in

  \begin{definition}
  \label{definition:occurs-semi-occurs}
    Let $p$ be an $A$-pattern and let $p'$ be an $A'$-pattern.
    \begin{enumerate} 
      \item Let $m$ be an element of $M$. The pattern $p$ is said to \define{occur at $m$ in $p'$}\graffito{$p$ occurs at $m$ in $p'$} and we write $p \occursIn_m p'$\graffito{$p \occursIn_m p'$}\index[symbols]{psquaresubseteqmpprime@$p \occursIn_m p'$} if and only if
            \begin{equation*}
              m \isSemiActedUponBy A \subseteq A' \text{ and } m \actsByItsCoordinateOn p = p'\restrictedTo_{m \isSemiActedUponBy A}.
            \end{equation*}
            And it is said to \define{semi-occur at $m$ in $p'$}\graffito{$p$ semi-occurs at $m$ in $p'$}\index{occur at $m$ in $p'$!semi-} and we write $p \semiOccursIn_m p'$\graffito{$p \semiOccursIn_m p'$}\index[symbols]{psquaresubseteqmzpprime@$p \semiOccursIn_m p'$} if and only if
            \begin{equation*}
              \Exists h \in H_{m_0, m} \Holds h \actsOnPoint A \subseteq A' \text{ and } h \actsOnMap p = p'\restrictedTo_{h \actsOnPoint A}.
            \end{equation*}
      \item The pattern $p$ is said to \define{occur in $p'$}\graffito{$p$ occurs in $p'$} and we write $p \occursIn p'$\graffito{$p \occursIn p'$}\index[symbols]{psquaresubseteqpprime@$p \occursIn p'$} if and only if 
            \begin{equation*}
              \Exists m \in M \SuchThat p \occursIn_m p'.
            \end{equation*}
            And it is said to \define{semi-occur in $p'$}\graffito{$p$ semi-occurs in $p'$}\index{occur in $p'$!semi-} and we write $p \semiOccursIn p'$\graffito{$p \semiOccursIn p'$}\index[symbols]{psquaresubseteqzpprime@$p \semiOccursIn p'$} if and only if 
            \begin{equation*}
              \Exists m \in M \SuchThat p \semiOccursIn_m p'. \qedhere
            \end{equation*}
    \end{enumerate}
  \end{definition}

  \begin{remark} 
  \label{remark:groups:occurs}
    Let $\mathcal{R}$ be the cell space $\ntuple{\ntuple{G, G, \cdot}, \ntuple{e_G, \family{g}_{g \in G}}}$, where $G$ is a group, $\cdot$ is its multiplication, and $e_G$ is its neutral element. Then, $G_0 = \setOf{e_G}$, $\isSemiActedUponBy = \cdot$, $\actsByItsCoordinateOn = \actsOnMap$, and, for each element $g \in G$, we have $H_{e_G, g} = \setOf{g}$. Hence, the notions \emph{occurs} and \emph{semi-occurs} are identical, and they are the common notion of \emph{occurrence} as used in \cite{fiorenzi:2003}.
  \end{remark}

  Semi-occurrence can be characterised in many ways, each illuminating a different aspect, some of which are given in

  \begin{remark}
  \label{remark:semi-occurs-at}
    Let $p$ be an $A$-pattern, let $p'$ be an $A'$-pattern, and let $m$ be an element of $M$. The following statements are equivalent:
    \begin{enumerate}
      \item $\displaystyle p \semiOccursIn_m p'$; 
      \item $\displaystyle \Exists h_0 \in H_0 \SuchThat h_0 \actsOnMap p \occursIn_m p'$;
      \item $\displaystyle \Exists h_0 \in H_0 \SuchThat g_{m_0, m} h_0 \actsOnPoint A \subseteq A' \text{ and } g_{m_0, m} h_0 \actsOnMap p = p'\restrictedTo_{g_{m_0, m} h_0 \actsOnPoint A}$;
      \item $\displaystyle \Exists h \in H_{m, m_0} \SuchThat A \subseteq h \actsOnPoint A' \text{ and } (h \actsOnMap p')\restrictedTo_A = p$;
      \item $\displaystyle \Exists h_0 \in H_0 \SuchThat A \subseteq h_0 g_{m_0, m}^{-1} \actsOnPoint A' \text{ and } (h_0 g_{m_0, m}^{-1} \actsOnMap p')\restrictedTo_A = p$. \qedhere
    \end{enumerate}
  \end{remark}

  \begin{remark}
  \label{remark:semi-occurs}
    Let $p$ be an $A$-pattern and let $p'$ be an $A'$-pattern. The following statements are equivalent:
    \begin{enumerate}
      \item $\displaystyle p \semiOccursIn p'$; 
      \item $\displaystyle \Exists h_0 \in H_0 \SuchThat h_0 \actsOnMap p \occursIn p'$;
      \item $\displaystyle \Exists h \in H \SuchThat h \actsOnPoint A \subseteq A' \text{ and } h \actsOnMap p = p'\restrictedTo_{h \actsOnPoint A}$;
      \item $\displaystyle \Exists h \in H \SuchThat A \subseteq h \actsOnPoint A' \text{ and } (h \actsOnMap p')\restrictedTo_A = p$. \qedhere
    \end{enumerate}
  \end{remark}

  %


  A pattern is allowed in a subset of the full shift if it semi-occurs in one of its points and forbidden otherwise, as defined in

  \begin{definition}
  \label{definition:allowed-forbidden}
    Let $X$ be a subset of $Q^M$ and let $p$ be an $A$-pattern. The pattern $p$ is called
    \begin{enumerate}
      \item \define{allowed in $X$}\graffito{pattern $p$ allowed in $X$}\index{pattern!allowed in} if and only if 
            \begin{equation*}
              \Exists x \in X \SuchThat p \semiOccursIn x;
            \end{equation*}
      \item \define{forbidden in $X$}\graffito{pattern $p$ forbidden in $X$}\index{pattern!forbidden in} if and only if 
            \begin{equation*}
              \ForEach x \in X \Holds p \not\semiOccursIn x. \qedhere
            \end{equation*}
    \end{enumerate}
  \end{definition}

  The greatest subset of the full shift with respect to inclusion in which each block of a given set of blocks is forbidden is said to be generated by the set of blocks, as defined in

  \begin{definition}
  \label{definition:generated-by-forbidden-blocks}
    Let $\mathfrak{F}$ be a subset of $Q^*$. The set
    \begin{equation*}
      \subshiftGeneratedBy{\mathfrak{F}} = \setOf{c \in Q^M \suchThat \ForEach \mathfrak{f} \in \mathfrak{F} \Holds \mathfrak{f} \not\semiOccursIn c} \mathnote{set $\subshiftGeneratedBy{\mathfrak{F}}$ generated by $\mathfrak{F}$}
      \index[symbols]{angleFanglefraktur@$\subshiftGeneratedBy{\mathfrak{F}}$}
    \end{equation*}
    is said to be \define{generated by $\mathfrak{F}$}\graffito{$\subshiftGeneratedBy{\mathfrak{F}}$ generated by $\mathfrak{F}$}. 
  \end{definition}


  A shift space is a subset of the full shift that is generated by a set of blocks, as defined in

  \begin{definition} 
  \label{definition:shift-space}
    Let $X$ be a subset of $Q^M$. It is called \define{shift space}\graffito{shift space $X$}\index{shift!space} and \define{subshift of $Q^M$}\graffito{subshift $X$ of $Q^M$}\index{shift!sub-} if and only if there is a subset $\mathfrak{F}$ of $Q^*$ such that $\subshiftGeneratedBy{\mathfrak{F}} = X$. 
  \end{definition}

  \begin{example}[Full Shift] 
  \label{example:shift:full}
    Because $\subshiftGeneratedBy{\emptyset} = Q^M$, the set $Q^M$ is a shift space.
  \end{example}

  \begin{example}[Empty Shift] 
  \label{example:shift:empty}
    If we identify each $q \in Q$ with the $\ball(m_0)$-block $[m_0 \mapsto q]$, then $\subshiftGeneratedBy{Q} = \emptyset$ and hence the set $\emptyset$ is a shift space.
  \end{example}

  \begin{example}[{Golden Mean Shift \cite[Example~1.2.3]{lind:marcus:1995}}]
  \label{example:shift:golden-mean}
    In the situation of \cref{example:full-shift-over-integers}, the set $X$ of all bi-infinite binary sequences with no two $1$'s next to each other is the shift space known as \graffito{golden mean shift}\define{golden mean shift}\index{shift!golden mean}. It is for example generated by the forbidden block $11$. 
  \end{example}

  \begin{example}[{Even Shift \cite[Example~1.2.4]{lind:marcus:1995}}]
  \label{example:shift:even}
    In the situation of \cref{example:full-shift-over-integers}, the set $X$ of all bi-infinite binary sequences such that, between any two occurrences of $1$'s, there are an even number of $0$'s, is the shift space known as \define{even shift}\graffito{even shift}\index{shift!even}. It is for example generated by the forbidden blocks $1 0^{2k + 1} 1$, for $k \in \N_0$.
  \end{example}


  Each shift space is shift-invariant, which is shown in

  \begin{lemma}
  \label{lemma:shift-space-is-shift-invariant}
    Let $X$ be a subshift of $Q^M$. It is shift-invariant.
  \end{lemma}

  \begin{proof} 
    There is a subset $\mathfrak{F}$ of $Q^*$ such that $\subshiftGeneratedBy{\mathfrak{F}} = X$. Let $x \in X$ and let $h \in H$. Suppose that there is an $\mathfrak{f} \in \mathfrak{F}$ such that $\mathfrak{f} \semiOccursIn h \actsOnMap x$. Then, there is an $h' \in H$ such that $h' \actsOnMap \mathfrak{f} = (h \actsOnMap x)\restrictedTo_{h' \actsOnPoint \domainOf(\mathfrak{f})}$. Thus, because $(h \actsOnMap x)\restrictedTo_{h' \actsOnPoint \domainOf(\mathfrak{f})} = h \actsOnMap (x\restrictedTo_{h^{-1} h' \actsOnPoint \domainOf(\mathfrak{f})})$, we have $h^{-1} h' \actsOnMap \mathfrak{f} = x\restrictedTo_{h^{-1} h' \actsOnPoint \domainOf(\mathfrak{f})}$. Hence, because $h^{-1} h' \in H$, we have $\mathfrak{f} \semiOccursIn x$, which contradicts that $x \in \subshiftGeneratedBy{\mathfrak{F}}$. Therefore, contrary to the supposition, for each $\mathfrak{f} \in \mathfrak{F}$, we have $\mathfrak{f} \not\semiOccursIn h \actsOnMap x$. Hence, $h \actsOnMap x \in X$. Therefore, $h \actsOnMap X \subseteq X$. In conclusion, according to \cref{remark:subseteq-sufficient-for-shift-invariance}, the subshift $X$ is shift-invariant.
  \end{proof}

  \begin{example}[{Shift-Invariant Non-Shift \cite[Example~1.2.10]{lind:marcus:1995}}]
  \label{example:shift-invariant-non-shift}
    In the situation of \cref{example:full-shift-over-integers}, the set $X$ of all bi-infinite binary sequences in which the symbol $1$ occurs exactly once is shift-invariant, but it is not a shift space.
  \end{example}

  Patterns with the same domain can all be restricted to some subdomain, as done in

  \begin{definition}
    Let $A$ be a subset of $M$, let $B$ be a subset of $A$, and let $P$ be a subset of $Q^A$. The set
    \begin{equation*}
      P_B = \setOf{p\restrictedTo_B \suchThat p \in P} \mathnote{$P_B$} 
    \end{equation*}
    is the set of all $B$-subpatterns of patterns of $P$.
  \end{definition}

  Because subshifts are shift-invariant, restrictions to patterns behave nicely with translations and rotations, as remarked in

  \begin{remark}
  \label{remark:pattern-belongs-to-shift-if-and-only-if-translated-pattern-belongs}
    Let $X$ be a subshift of $Q^M$, let $A$ be a subset of $M$, let $h$ be an element of $H$, and let $m$ be an element of $M$. For each $A$-pattern $p$, we have $p \in X_A$ if and only if $h \actsOnMap p \in X_{h \actsOnPoint A}$, and $p \in X_A$ if and only if $m \actsByItsCoordinateOn p \in X_{m \isSemiActedUponBy A}$. In particular, $h \actsOnMap X_A = X_{h \actsOnPoint A}$ and $m \actsByItsCoordinateOn X_A = X_{m \isSemiActedUponBy A}$. And, if $A$ is finite, then $\cardinalityOf{X_A} = \cardinalityOf{X_{h \actsOnPoint A}}$ and $\cardinalityOf{X_A} = \cardinalityOf{X_{m \isSemiActedUponBy A}}$. 
  \end{remark}

  Because shift spaces are generated by forbidden blocks, which have finite domains, a point of the full shift belongs to a shift space if and only if its subpatterns do, which is shown in 

  \begin{lemma}
  \label{lemma:if-restrictions-to-balls-are-in-shift-then-so-is-pattern}
    Let $X$ be a subshift of $Q^M$ and let $c$ be a point of $Q^M$. Then, $c \in X$ if and only if 
    \begin{equation}
    \label{equation:if-restrictions-to-balls-are-in-shift-then-so-is-pattern}
      \ForEach \rho \in \N_0 \Holds c\restrictedTo_{\ball(\rho)} \in X_{\ball(\rho)}. \qedhere
    \end{equation}
  \end{lemma}

  \begin{proof}
    If $c \in X$, then \cref{equation:if-restrictions-to-balls-are-in-shift-then-so-is-pattern} holds. From now on, let \cref{equation:if-restrictions-to-balls-are-in-shift-then-so-is-pattern} hold. Because $X$ is a subshift, there is a subset $\mathfrak{F}$ of $Q^*$ such that $\subshiftGeneratedBy{\mathfrak{F}} = X$. Let $\mathfrak{f} \in \mathfrak{F}$. Suppose that $\mathfrak{f}$ semi-occurs in $c$. Then, because $\sizeOf{\mathfrak{f}} < \infty$, according to \cref{remark:ball-of-radius-0-contains-one-element-and-sequence-of-balls-is-monotonic}, there is a $\rho \in \N_0$ such that $\mathfrak{f}$ semi-occurs in $c\restrictedTo_{\ball(\rho)}$. Hence, $c\restrictedTo_{\ball(\rho)} \notin X_{\ball(\rho)}$, which contradicts \cref{equation:if-restrictions-to-balls-are-in-shift-then-so-is-pattern}. Therefore, $\mathfrak{f}$ does not semi-occur in $c$. In conclusion, $c \in X$.
  \end{proof}

  Shift spaces are characterised by shift-invariance and compactness, which is shown in

  \begin{lemma} 
  \label{lemma:shift-space-if-and-only-if-invariant-and-compact}
    Let $X$ be a subset of $Q^M$. It is a shift space if and only if it is shift-invariant and compact.
  \end{lemma}

  \begin{proof}
    First, let $X$ be a shift space. Then, according to \cref{lemma:shift-space-is-shift-invariant}, it is shift-invariant. Moreover, let $\sequence{x_k}_{k \in \N_+}$ be a sequence in $X$ that converges to a point $c \in Q^M$. Then, for each $\rho \in \N_0$, there is a $k \in \N_+$ such that $c\restrictedTo_{\ball(\rho)} = x_k\restrictedTo_{\ball(\rho)} \in X_{\ball(\rho)}$. Thus, according to \cref{lemma:if-restrictions-to-balls-are-in-shift-then-so-is-pattern}, we have $c \in X$. Hence, $X$ is closed. And, according to the first paragraph of section~1.8 in \cite{ceccherini-silberstein:coornaert:2010}, the set $Q^M$ is compact. Therefore, $X$ is compact. 

    Secondly, let $X$ be shift-invariant and compact. Then, $X$ is closed and $Q^M \smallsetminus X$ is open. Hence, for each $c \in Q^M \smallsetminus X$, there is a $\rho_c \in \N_0$ such that $\cylinder(c\restrictedTo_{\ball(\rho_c)}) \subseteq Q^M \smallsetminus X$. Put $\mathfrak{F} = \setOf{c\restrictedTo_{\ball(\rho_c)} \suchThat c \in Q^M \smallsetminus X}$.

    Let $x \in X$. Suppose that there is an $\mathfrak{f} \in \mathfrak{F}$ such that $\mathfrak{f} \semiOccursIn x$. Then, according to \cref{remark:semi-occurs}, there is an $h \in H$ such that $(h \actsOnMap x)\restrictedTo_{\domainOf(\mathfrak{f})} = \mathfrak{f}$. Thus, $h \actsOnMap x \in \cylinder(\mathfrak{f}) \subseteq Q^M \smallsetminus X$. And, because $X$ is shift-invariant, we also have $h \actsOnMap x \in X$, which contradicts that $h \actsOnMap x \in Q^M \smallsetminus X$. Hence, contrary to the supposition, for each $\mathfrak{f} \in \mathfrak{F}$ we have $\mathfrak{f} \not\semiOccursIn x$. Therefore, $X \subseteq \subshiftGeneratedBy{\mathfrak{F}}$.

    Let $c \in \subshiftGeneratedBy{\mathfrak{F}}$. Suppose that $c \notin X$. Then, $c\restrictedTo_{\ball(\rho_c)} \in \mathfrak{F}$. Thus $c \notin \subshiftGeneratedBy{\mathfrak{F}}$, which contradicts that $c \in \subshiftGeneratedBy{\mathfrak{F}}$. Hence, $c \in X$. Therefore, $\subshiftGeneratedBy{\mathfrak{F}} \subseteq X$.

    Altogether, $\subshiftGeneratedBy{\mathfrak{F}} = X$. In conclusion, $X$ is a shift space.
  \end{proof}

  \begin{remark}
    Compactness cannot be omitted in the equivalence. For example, the subset of $\setOf{0, 1}^\Z$ from \cref{example:shift-invariant-non-shift} is shift-invariant but not a shift space.
  \end{remark}

  A shift space of finite type is one that is generated by finitely many forbidden blocks, as defined in

  \begin{definition} 
  \label{definition:of-finite-type}
    Let $X$ be a subshift of $Q^M$. It is said to be \graffito{of finite type}\define{of finite type}\index{shift!of finite type}\index{finite type@of finite type} if and only if there is a finite subset $\mathfrak{F}$ of $Q^*$ such that $\subshiftGeneratedBy{\mathfrak{F}} = X$.
  \end{definition}

  \begin{example}[Of Finite Type or Not]
  \label{example:of-finite-type-or-not}
    As is apparent from their definition, the full shift (\cref{example:shift:full}), the empty shift (\cref{example:shift:empty}), and the golden mean shift (\cref{example:shift:golden-mean}) are of finite type. According to example~2.1.5 in \cite{lind:marcus:1995}, the even shift (\cref{example:shift:even}) is \emph{not} of finite type.
  \end{example}

  A $\kappa$-step shift space is one that is generated by forbidden $\ball(\kappa)$-blocks, as defined in

  \begin{definition} 
  \label{definition:step-and-memory}
    Let $X$ be a subshift of $Q^M$ and let $\kappa$ be a non-negative integer. The subshift $X$ is called \defineX{$\kappa$-step}{step@$\kappa$-step}\graffito{$\kappa$-step} and the integer $\kappa$ is called \define{memory of $X$}\graffito{memory $\kappa$ of $X$} if and only if there is a subset $\mathfrak{F}$ of $Q^{\ball(\kappa)}$ such that $\subshiftGeneratedBy{\mathfrak{F}} = X$.
  \end{definition}

  \begin{remark}
  \label{remark:k-step-for-greater-k}
    Let $X$ be a $\kappa$-step subshift of $Q^M$. Because the set $Q^{\ball(\kappa)}$ is finite, the subshift $X$ is of finite type. And, for each non-negative integer $\kappa'$ such that $\kappa' \geq \kappa$, the subshift $X$ is $\kappa'$-step.
  \end{remark}

  A shift space of finite type is $\kappa$-step, where $\kappa$ is the radius of a ball that includes all domains of a finite generating set of the shift space, which is shown in

  \begin{lemma}
  \label{lemma:of-finite-type-implies-step}
    Let $X$ be a subshift of $Q^M$ of finite type. There is a non-negative integer $\kappa$ such that $X$ is $\kappa$-step.
  \end{lemma}

  \begin{proof}
    Because $X$ is of finite type, there is a finite subset $\mathfrak{F}$ of $Q^*$ such that $\subshiftGeneratedBy{\mathfrak{F}} = X$. And, because the set $\mathfrak{F}$ is finite, according to \cref{remark:ball-of-radius-0-contains-one-element-and-sequence-of-balls-is-monotonic}, there is a non-negative integer $\kappa$ such that, for each $\mathfrak{f} \in \mathfrak{F}$, we have $\domainOf(\mathfrak{f}) \subseteq \ball(\kappa)$. Let $\mathfrak{F}'$ be the set $\setOf{p \in Q^{\ball(\kappa)} \suchThat \Exists \mathfrak{f} \in \mathfrak{F} \SuchThat p\restrictedTo_{\domainOf(\mathfrak{f})} = \mathfrak{f}}$. Then, $\subshiftGeneratedBy{\mathfrak{F}'} = \subshiftGeneratedBy{\mathfrak{F}} = X$. In conclusion, $X$ is $\kappa$-step. 
  \end{proof}

  A shift space is $\kappa$-step if and only if it contains each point of the full shift whose restrictions to the balls of radius $\kappa$ are allowed patterns, which is shown in

  \begin{lemma} 
  \label{lemma:characterisation-of-kappa-step-subshifts}
    Let $X$ be a subshift of $Q^M$ and let $\kappa$ be a non-negative integer. The subshift $X$ is $\kappa$-step if and only if 
    \begin{equation}
    \label{equation:characterisation-of-kappa-step-subshifts:char}
      \ForEach c \in Q^M \Holds \parens[\Big]{\parens[\big]{\ForEach m \in M \Holds c\restrictedTo_{\ball(m, \kappa)} \in X_{\ball(m, \kappa)}} \implies c \in X}. \qedhere
    \end{equation}
  \end{lemma}

  \begin{proof}
    First, let $X$ be $\kappa$-step. Then, there is an $\mathfrak{F} \subseteq Q^{\ball(\kappa)}$ such that $\subshiftGeneratedBy{\mathfrak{F}} = X$. Furthermore, let $c \in Q^M$ such that
    \begin{equation}
    \label{equation:characterisation-of-kappa-step-subshifts:antecedent}
      \ForEach m \in M \Holds c\restrictedTo_{\ball(m, \kappa)} \in X_{\ball(m, \kappa)}.
    \end{equation}
    Suppose that there is an $\mathfrak{f} \in \mathfrak{F}$ such that $\mathfrak{f} \semiOccursIn c$. Then, there is an $m \in M$ such that $\mathfrak{f} \semiOccursIn_m c$. Thus, there is an $h \in H_{m_0, m}$ such that $h \actsOnMap \mathfrak{f} = c\restrictedTo_{h \actsOnPoint \ball(\kappa)}$. Moreover, because $h \actsOnPoint m_0 = m$, according to \cref{lemma:left-action-and-balls}, we have $h \actsOnPoint \ball(\kappa) = \ball(m, \kappa)$. Hence, according to \cref{equation:characterisation-of-kappa-step-subshifts:antecedent}, we have $h \actsOnMap \mathfrak{f} = c\restrictedTo_{\ball(m, \kappa)} \in X_{\ball(m, \kappa)} = X_{h \actsOnPoint \ball(\kappa)}$. Therefore, there is an $x \in X$ such that $\mathfrak{f} \semiOccursIn x$, which contradicts that $\subshiftGeneratedBy{\mathfrak{F}} = X$. In conclusion, for each $\mathfrak{f} \in \mathfrak{F}$, we have $\mathfrak{f} \not\semiOccursIn c$, and hence $c \in X$.

    Secondly, let \cref{equation:characterisation-of-kappa-step-subshifts:char} hold. Furthermore, let $\mathfrak{F} = Q^{\ball(\kappa)} \smallsetminus X_{\ball(\kappa)}$. We show below that $X \subseteq \subshiftGeneratedBy{\mathfrak{F}}$ and $\subshiftGeneratedBy{\mathfrak{F}} \subseteq X$. Hence, $\subshiftGeneratedBy{\mathfrak{F}} = X$. In conclusion, $X$ is $\kappa$-step.

    \proofPart{Subproof of: $X \subseteq \subshiftGeneratedBy{\mathfrak{F}}$}
    Let $c \in Q^M \smallsetminus \subshiftGeneratedBy{\mathfrak{F}}$. Then, there is an $\mathfrak{f} \in \mathfrak{F}$ such that $\mathfrak{f} \semiOccursIn c$. Thus, there is an $m \in M$ and there is an $h \in H_{m_0, m}$ such that $h \actsOnMap \mathfrak{f} = c\restrictedTo_{\ball(m, \kappa)}$. Therefore, because $h \actsOnMap \blank$ is bijective, according to \cref{remark:pattern-belongs-to-shift-if-and-only-if-translated-pattern-belongs} and \cref{lemma:left-action-and-balls}, we have $c\restrictedTo_{\ball(m, \kappa)} \in h \actsOnMap \mathfrak{F} = h \actsOnMap (Q^{\ball(\kappa)} \smallsetminus X_{\ball(\kappa)}) = (h \actsOnMap Q^{\ball(\kappa)}) \smallsetminus (h \actsOnMap X_{\ball(\kappa)}) = Q^{\ball(m, \kappa)} \smallsetminus X_{\ball(m, \kappa)}$. Hence, $c \in Q^M \smallsetminus X$. In conclusion, $X \subseteq \subshiftGeneratedBy{\mathfrak{F}}$.

    \proofPart{Subproof of: $\subshiftGeneratedBy{\mathfrak{F}} \subseteq X$}
    Let $c \in Q^M \smallsetminus X$. Then, according to \cref{equation:characterisation-of-kappa-step-subshifts:char}, there is an $m \in M$ such that $c\restrictedTo_{\ball(m, \kappa)} \in Q^{\ball(m, \kappa)} \smallsetminus X_{\ball(m, \kappa)}$. Furthermore, let $h = g_{m_0, m}$. Then, similar as above, $h^{-1} \actsOnMap c\restrictedTo_{\ball(m, \kappa)} \in Q^{\ball(\kappa)} \smallsetminus X_{\ball(\kappa)} = \mathfrak{F}$. Thus, there is an $\mathfrak{f} \in \mathfrak{F}$ such that $h \actsOnMap \mathfrak{f} = c\restrictedTo_{\ball(m, \kappa)}$. Hence, because $h \in H$, we have $\mathfrak{f} \semiOccursIn c$. Therefore, $c \in Q^M \smallsetminus \subshiftGeneratedBy{\mathfrak{F}}$. In conclusion, $\subshiftGeneratedBy{\mathfrak{F}} \subseteq X$.
  \end{proof} 

  If we cut holes in a point of a shift space of finite type that are far enough apart and fill these holes with pieces from other points of the shift space that agree on big enough boundaries of the holes with the holey point, then we still have a point of the shift space, which is shown in

  \begin{lemma} 
  \label{lemma:overlapping-global-configurations-can-be-glued}
    Let $X$ be a $\kappa$-step subshift of $Q^M$, let $x$ be a point of $X$, let $\family{A_i}_{i \in I}$ be a family of subsets of $M$ such that the family $\family{A_i^{+2\kappa}}_{i \in I}$ is pairwise disjoint, and let $\family{x_i}_{i \in I}$ be a family of points of $X$ such that, for each index $i \in I$, we have $x_i\restrictedTo_{\boundaryOf_{2\kappa}^+ A_i} = x\restrictedTo_{\boundaryOf_{2\kappa}^+ A_i}$. The map $x\restrictedTo_{M \smallsetminus (\bigcup_{i \in I} A_i)} \times \coprod_{i \in I} x_i\restrictedTo_{A_i}$ is identical to $x\restrictedTo_{M \smallsetminus (\bigcup_{i \in I} A_i^{+2\kappa})} \times \coprod_{i \in I} x_i\restrictedTo_{A_i^{+2\kappa}}$ and a point of $X$. Recall that such coproducts were introduced in \cref{definition:coproduct-of-maps}.
  \end{lemma}

  \begin{proof}
    Because $\family{A_i^{+2\kappa}}_{i \in I}$ is pairwise disjoint, so is $\family{A_i}_{i \in I}$. Hence, $x' = x\restrictedTo_{M \smallsetminus (\bigcup_{i \in I} A_i)} \times \coprod_{i \in I} x_i\restrictedTo_{A_i}$ is well-defined. Furthermore, for each $i \in I$, we have $x_i\restrictedTo_{A_i^{+2\kappa} \smallsetminus A_i} = x\restrictedTo_{A_i^{+2\kappa} \smallsetminus A_i}$. Hence, $x' = x\restrictedTo_{M \smallsetminus (\bigcup_{i \in I} A_i^{+2\kappa})} \times \coprod_{i \in I} x_i\restrictedTo_{A_i^{+2\kappa}}$. Moreover, because $X$ is $\kappa$-step, there is an $\mathfrak{F} \subseteq Q^{\ball(\kappa)}$ such that $\subshiftGeneratedBy{\mathfrak{F}} = X$.

    Let $u \in Q^{\ball(\kappa)}$ semi-occur in $x'$. Then, there is an $m \in M$ and there is an $h \in H_{m_0, m}$ such that $h \actsOnMap u = x'\restrictedTo_{h \actsOnPoint \ball(\kappa)}$. Moreover, according to \cref{lemma:left-action-and-balls} and \cref{item:characterisation-of-k-closure-and-interior:closure} of \cref{lemma:characterisation-of-k-closure-and-interior}, we have $h \actsOnPoint \ball(\kappa) = \ball(m, \kappa) = \setOf{m}^{+\kappa}$.
    \begin{description}
      \item[Case 1:] $\Exists i \in I \SuchThat m \in A_i^{+\kappa}$. Then, according to \cref{item:properties-of-interior-closure-and-boundary:subset} of \cref{lemma:properties-of-interior-closure-and-boundary} and \cref{item:repeated-k-boundaries-etc:closure} of \cref{lemma:repeated-k-boundaries-etc}, we have $\setOf{m}^{+\kappa} \subseteq (A_i^{+\kappa})^{+\kappa} \subseteq A_i^{+2\kappa}$. Hence, $h \actsOnPoint \ball(\kappa) \subseteq A_i^{+2\kappa}$. Thus, $x'\restrictedTo_{h \actsOnPoint \ball(\kappa)} = x_i\restrictedTo_{h \actsOnPoint \ball(\kappa)}$. Therefore, $u$ semi-occurs in $x_i$. Hence, $u \notin \mathfrak{F}$.
      \item[Case 2:] $m \in M \smallsetminus (\bigcup_{i \in I} A_i^{+\kappa})$. Then, according to \cref{item:properties-of-interior-closure-and-boundary:subset}, \cref{item:properties-of-interior-closure-and-boundary:complement}, and \cref{item:properties-of-interior-closure-and-boundary:union} of \cref{lemma:properties-of-interior-closure-and-boundary} and \cref{item:repeated-k-boundaries-etc:closure-interior} of \cref{lemma:repeated-k-boundaries-etc}, we have
            \begin{align*}
              \setOf{m}^{+\kappa}
              &\subseteq (M \smallsetminus (\bigcup_{i \in I} A_i^{+\kappa}))^{+\kappa}\\
              &\subseteq M \smallsetminus (\bigcup_{i \in I} (A_i^{+\kappa})^{-\kappa})\\
              &\subseteq M \smallsetminus (\bigcup_{i \in I} A_i^{+(\kappa - \kappa)})\\
              &= M \smallsetminus (\bigcup_{i \in I} A_i).
            \end{align*}
            Hence, $h \actsOnPoint \ball(\kappa) \subseteq M \smallsetminus (\bigcup_{i \in I} A_i)$. Thus, $x'\restrictedTo_{h \actsOnPoint \ball(\kappa)} = x\restrictedTo_{h \actsOnPoint \ball(\kappa)}$. Therefore, $u$ semi-occurs in $x$. Hence, $u \notin \mathfrak{F}$.
    \end{description}
    In either case, $u \notin \mathfrak{F}$. In conclusion, $x' \in X$.
  \end{proof}

  If we sew one part and respectively the other part of two sufficiently overlapping points of a shift space of finite type together, then we get a point of the shift space, which is shown in

  \begin{corollary} 
  \label{corollary:overlapping-global-configurations-can-be-glued}
    Let $X$ be a $\kappa$-step subshift of $Q^M$, let $A$ be a subset of $M$, and let $x$ and $x'$ be two points of $X$ such that $x\restrictedTo_{\boundaryOf_{2\kappa}^+ A} = x'\restrictedTo_{\boundaryOf_{2\kappa}^+ A}$. The map $x\restrictedTo_A \times x'\restrictedTo_{M \smallsetminus A}$ is identical to $x\restrictedTo_{A^{+2\kappa}} \times x'\restrictedTo_{M \smallsetminus A^{+2\kappa}}$ and a point of $X$.
  \end{corollary}

  \begin{proof}
    This is a direct consequence of \cref{lemma:overlapping-global-configurations-can-be-glued} with $\family{x_i}_{i \in I} = \family{x'}$.
  \end{proof}

  If two patterns with the same domain, that are allowed in a shift space of finite type, agree on a big enough boundary, then they can be identically extended to points of the shift space, which is shown in

  \begin{corollary} 
  \label{corollary:overlapping-patterns-can-be-glued}
    Let $X$ be a $\kappa$-step subshift of $Q^M$, let $A$ be a subset of $M$, let $p$ and $p'$ be two patterns of $X_{A^{+2\kappa}}$ such that $p\restrictedTo_{\boundaryOf_{2\kappa}^+ A} = p'\restrictedTo_{\boundaryOf_{2\kappa}^+ A}$. There are two points $x$ and $x'$ of $X$ such that $x\restrictedTo_{A^{+2\kappa}} = p$, $x'\restrictedTo_{A^{+2\kappa}} = p'$, and $x\restrictedTo_{M \smallsetminus A} = x'\restrictedTo_{M \smallsetminus A}$.
  \end{corollary}

  \begin{proof}
    Because $p$, $p' \in X_{A^{+2\kappa}}$, there are $x''$, $x' \in X$ such that $x''\restrictedTo_{A^{+2\kappa}} = p$ and $x'\restrictedTo_{A^{+2\kappa}} = p'$. Because $p\restrictedTo_{\boundaryOf_{2\kappa}^+ A} = p'\restrictedTo_{\boundaryOf_{2\kappa}^+ A}$, we have $x''\restrictedTo_{\boundaryOf_{2\kappa}^+ A} = x'\restrictedTo_{\boundaryOf_{2\kappa}^+ A}$. Hence, according to \cref{lemma:overlapping-global-configurations-can-be-glued}, we have $x = x''\restrictedTo_A \times x'\restrictedTo_{M \smallsetminus A} \in X$. Moreover, $x\restrictedTo_{M \smallsetminus A} = x'\restrictedTo_{M \smallsetminus A}$. And, because $x = x''\restrictedTo_{A^{+2\kappa}} \times x'\restrictedTo_{M \smallsetminus A^{+2\kappa}}$, we have $x\restrictedTo_{A^{+2\kappa}} = p$.
  \end{proof}



  A shift space is strongly irreducible if two allowed finite patterns that are far enough apart are embedded in a point of the shift space, which is defined in

  \begin{definition} 
  \label{definition:kappa-strongly-irreducible}
    Let $X$ be a subshift of $Q^M$ and let $\kappa$ be a non-negative integer. The subshift $X$ is called \defineX{$\kappa$-strongly irreducible}{strongly irreducible kappa@$\kappa$-strongly irreducible}\graffito{$\kappa$-strongly irreducible} if and only if, for each tuple $(p, p')$ of finite patterns allowed in $X$ such that $\distanceOf(\domainOf(p), \domainOf(p')) \geq \kappa + 1$, there is a point $x \in X$ such that $x\restrictedTo_{\domainOf(p)} = p$ and $x\restrictedTo_{\domainOf(p')} = p'$.
  \end{definition}

  \begin{remark}
  \label{remark:k-strongly-irreducible-for-greater-k}
    Let $X$ be a $\kappa$-strongly irreducible subshift of $Q^M$. For each non-negative integer $\kappa'$ such that $\kappa' \geq \kappa$, the subshift $X$ is $\kappa'$-strongly irreducible.
  \end{remark}

  \begin{definition}
  \label{definition:strongly-irreducible}
    Let $X$ be a subshift of $Q^M$. It is called \graffito{strongly irreducible}\define{strongly irreducible} if and only if there is a non-negative integer $\kappa$ such that it is $\kappa$-strongly irreducible.
  \end{definition} 

  \begin{example}[Strongly Irreducible]
  \label{example:strongly-irreducible} 
    The full shift (\cref{example:shift:full}) and the empty shift (\cref{example:shift:empty}) are $0$-strongly irreducible. The golden mean shift is $1$-strongly irreducible. The even shift is $2$-strongly irreducible. The generalised golden mean shifts (\cref{example:generalised-golden-mean-shift}) and the shift space of \cref{example:strongly-irreducible-of-finite-type-without-bounded-propagation} are strongly irreducible. Of these examples, according to \cref{example:of-finite-type-or-not}, all but the even shift are of finite type.
  \end{example}

  \begin{example}[{Not Strongly Irreducible \cite[Example~4.6]{ceccherini-silberstein:coornaert:2012}}] 
  \label{example:not-strongly-irreducible}
    In the situation of \cref{example:full-shift-over-integers}, the set $X$ of all bi-infinite binary sequences with no two $0$'s and no two $1$'s next to each other is a shift space. It is for example generated by the forbidden blocks $00$ and $11$, in particular, it is of finite type. It consists of the two bi-infinite binary sequences with alternating $0$'s and $1$'s. And, it is not strongly irreducible; indeed, for each even and non-negative integer $\kappa$, the finite patterns $01$ and $10$ are allowed in $X$, the allowed patterns of size $\kappa$ are of the form $(01)^{\kappa/2}$ or $(10)^{\kappa/2}$, but the patterns $01(01)^{\kappa/2}10$ and $01(10)^{\kappa/2}10$ are not allowed and hence not embedded in a point of $X$.
  \end{example}


  A shift space has bounded propagation if a finite pattern is allowed whenever all restrictions of it to balls of a fixed radius are allowed, as defined in

  \begin{definition} 
  \label{definition:rho-bounded-propagation}
    Let $X$ be a subshift of $Q^M$ and let $\rho$ be a non-negative integer. The subshift $X$ is said to have \defineX{$\rho$-bounded propagation}{bounded propagation@$\rho$-bounded propagation}\graffito{$\rho$-bounded propagation} if and only if
    \begin{multline*}
      \ForEach F \subseteq M \text{ finite} \ForEach p \in Q^F \Holds\\
          \parens[\Big]{ \parens[\big]{\ForEach f \in F \Holds p\restrictedTo_{\ball(f, \rho) \cap F} \in X_{\ball(f, \rho) \cap F}} \implies p \in X_F }. \qedhere
    \end{multline*}
  \end{definition}

  \begin{remark}
    Let $X$ be a subshift of $Q^M$ with $\rho$-bounded propagation. For each non-negative integer $\rho'$ such that $\rho' \geq \rho$, the subshift $X$ has $\rho'$-bounded propagation.
  \end{remark}

  \begin{definition}
  \label{definition:bounded-propagation}
    Let $X$ be a subshift of $Q^M$. It is said to have \define{bounded propagation}\graffito{bounded propagation} if and only if there is a non-negative integer $\rho$ such that it has $\rho$-bounded propagation.
  \end{definition}

  \begin{remark}
    The notion of bounded propagation was first introduced by Mikhail Leonidovich Gromov in paragraph~7.A in \cite{gromov:1999}.
  \end{remark}

  A subshift with bounded propagation is strongly irreducible and of finite type, which is shown in

  \begin{lemma} 
  \label{lemma:bounded-propagation-implies-strong-irreducibility-and-finite-type}
    Let $X$ be a subshift of $Q^M$ with $\rho$-bounded propagation. It is $\rho$-strongly irreducible and $\rho$-step.
  \end{lemma}

  \begin{proof}
    Let $(p, p')$ be a tuple of finite patterns allowed in $X$ such that $\distanceOf(\domainOf(p), \domainOf(p')) \geq \rho + 1$. Moreover, let $F = \domainOf(p) \cup \domainOf(p')$ and let $p'' = p \times p' \in Q^F$. Furthermore, let $f \in F$. Then, if $f \in \domainOf(p)$, then $\ball(f, \rho) \cap F \subseteq \domainOf(p)$; and, if $f \in \domainOf(p')$, then $\ball(f, \rho) \cap F \subseteq \domainOf(p')$. Thus, in either case, $p''\restrictedTo_{\ball(f, \rho) \cap F} \in X_{\ball(f, \rho) \cap F}$. Hence, because $X$ has $\rho$-bounded propagation, we have $p'' \in X_F$. Therefore, there is an $x \in X$ such that $x\restrictedTo_F = p''$, in particular, $x\restrictedTo_{\domainOf(p)} = p$ and $x\restrictedTo_{\domainOf(p')} = p'$. In conclusion, $X$ is $\rho$-strongly irreducible.

    Let $c \in Q^M$ such that, for each $m \in M$, we have $c\restrictedTo_{\ball(m, \rho)} \in X_{\ball(m, \rho)}$. Moreover, let $\rho' \in \N_0$. Then, for each $f \in \ball(\rho')$, we have
    \begin{multline*}
      (c\restrictedTo_{\ball(\rho')})\restrictedTo_{\ball(f, \rho) \cap \ball(\rho')}
      = (c\restrictedTo_{\ball(f, \rho)})\restrictedTo_{\ball(\rho') \cap \ball(f, \rho)}\\
      \in (X_{\ball(f, \rho)})_{\ball(\rho') \cap \ball(f, \rho)}
      = X_{\ball(f, \rho) \cap \ball(\rho')}.
    \end{multline*}
    Thus, because $X$ has $\rho$-bounded propagation, we have $c\restrictedTo_{\ball(\rho')} \in X_{\ball(\rho')}$. Hence, according to \cref{lemma:if-restrictions-to-balls-are-in-shift-then-so-is-pattern}, we have $c \in X$. In conclusion, according to \cref{lemma:characterisation-of-kappa-step-subshifts}, the subshift $X$ is $\rho$-step.
  \end{proof}

  \begin{example}[Has Bounded Propagation or Not]
  \label{example:bounded-propagation}
    The full shift (\cref{example:shift:full}) and the empty shift (\cref{example:shift:empty}) have $0$-bounded propagation. The golden mean shift (\cref{example:shift:golden-mean}) has $1$-bounded propagation. The even shift (\cref{example:shift:even}) is, according to example \ref{example:of-finite-type-or-not}, not of finite type and hence it is, according to \cref{lemma:bounded-propagation-implies-strong-irreducibility-and-finite-type}, does not have bounded propagation. The shift of \cref{example:not-strongly-irreducible} is not strongly irreducible and hence it is, according to \cref{lemma:bounded-propagation-implies-strong-irreducibility-and-finite-type}, does not have bounded propagation. 
  \end{example}

  \begin{example}[{Generalised Golden Mean Shifts \cite[Example~2.8]{ceccherini-silberstein:coornaert:2012}}] 
  \label{example:generalised-golden-mean-shift}
    Let $q$ be a positive integer, let $Q$ be the set $\setOf{0, 1, \dotsc, q}$, let $k$ be a positive integer, let $\family{F_i}_{i \in \setOf{1, 2, \dotsc, k}}$ be a family of finite subsets of $M$ that contain $m_0$, let $\rho$ be the least non-negative integer such that $\bigcup_{i \in \setOf{1, 2, \dotsc, k}} F_i \subseteq \ball(\rho)$, and let $X$ be the $\rho$-step subshift $\subshiftGeneratedBy{p \in Q^{\ball(\rho)} \suchThat \ForEach i \in \setOf{1, 2, \dotsc, k} \ForEach m \in F_i \Holds p(m) \neq 0}$ of $Q^M$. The subshift $X$ is called \define{generalised golden mean shift}\graffito{generalised golden mean shift}\index{golden mean shift!generalised}\index{shift!generalised golden mean}, it is equal to $\setOf{c \in Q^M \suchThat \ForEach i \in \setOf{1, 2, \dotsc, k} \ForEach h \in H \Exists m \in h \actsOnPoint F_i \SuchThat c(m) = 0}$, and it has $\rho$-bounded propagation, in particular, according to \cref{lemma:bounded-propagation-implies-strong-irreducibility-and-finite-type}, it is $\rho$-strongly irreducible.

    In the case that $\mathcal{R}$ is the cell space from \cref{example:full-shift-over-integers}, $q = 1$, $k = 1$, and $F_1 = \setOf{0, 1}$, the non-negative integer $\rho$ is equal to $1$ and the generalised golden mean shift $X$ is equal to the golden mean shift from \cref{example:shift:golden-mean}.
  \end{example}

  \begin{proof}[Proof of Bounded Propagation]
    Let $F$ be a finite subset of $M$, let $p$ be a pattern of $Q^F$ such that
    \begin{equation*}
      \ForEach f \in F \Holds p\restrictedTo_{\ball(f, \rho) \cap F} \in X_{\ball(f, \rho) \cap F},
    \end{equation*}
    and let $c$ be the point of $Q^M$ that is equal to $p$ on $F$ and identically $0$ on $M \smallsetminus F$. Moreover, let $i \in I$ and let $h \in H$. If $h \actsOnPoint F_i \nsubseteq F$, then, there is an $m \in (h \actsOnPoint F_i) \smallsetminus F$, and, by definition of $c$, we have $c(m) = 0$. Otherwise, if $h \actsOnPoint F_i \subseteq F$, then, because $m_0 \in F_i$, we have $h \actsOnPoint m_0 \in F$; thus, because $p\restrictedTo_{\ball(h \actsOnPoint m_0, \rho) \cap F} \in X_{\ball(h \actsOnPoint m_0, \rho) \cap F}$, there is a point $x \in X$ such that $p\restrictedTo_{\ball(h \actsOnPoint m_0, \rho) \cap F} = x\restrictedTo_{\ball(h \actsOnPoint m_0, \rho) \cap F}$; hence, by the characterisation of $X$, there is a cell $m \in h \actsOnPoint F_i$ such that $x(m) = 0$; and therefore, because $h \actsOnPoint F_i \subseteq (h \actsOnPoint \ball(\rho)) \cap F = \ball(h \actsOnPoint m_0, \rho) \cap F$, we have $c(m) = p(m) = x(m) = 0$. Hence, in either case, there is an $m \in h \actsOnPoint F_i$ such that $c(m) = 0$. Therefore, by the characterisation of $X$, we have $c \in X$ and hence $p = c\restrictedTo_F \in X_F$. In conclusion, $X$ has $\rho$-bounded propagation.
  \end{proof}

  \begin{example}[{\cite[Section 4, at the very end]{fiorenzi:2003}}]
  \label{example:strongly-irreducible-of-finite-type-without-bounded-propagation}
    In the situation of \cref{example:full-shift-over-integers}, the subshift $\subshiftGeneratedBy{010, 111}$ of $\setOf{0, 1}^\Z$ is strongly irreducible and of finite type but does not have bounded propagation. 
  \end{example}

  A map from a shift space to another shift space is local if the image of a point is uniformly and locally determined in each cell, as defined in

  \begin{definition} 
  \label{definition:kappa-local-map}
    Let $X$ and $Y$ be two subshifts of $Q^M$, let $\Delta$ be a map from $X$ to $Y$, let $\kappa$ be a non-negative integer, let $N$ be a subset of $\ball(\kappa)$ such that $G_0 \actsOnPoint N \subseteq N$, and let $\bullet_{H_0}$ be the left group action $\actsOnMap\restrictedTo_{H_0 \times X_N \to X_N}$ of $H_0$ on $X_N$. The map $\Delta$ is called \defineX{$\kappa$-local}{local@$\kappa$-local}\graffito{$\kappa$-local} if and only if there is a $\bullet_{H_0}$-invariant map $\delta \from X_N \to Q$ such that 
    \begin{equation*}
      \ForEach x \in X \ForEach m \in M \Holds \Delta(x)(m) = \delta(n \mapsto x(m \isSemiActedUponBy n)). \qedhere %
    \end{equation*}
  \end{definition}

  \begin{remark} 
    For each point $x \in X$ and each cell $m \in M$, we have $\Delta(x)(m) = \delta((g_{m_0, m}^{-1} \actsOnMap x)\restrictedTo_N)$.
  \end{remark}

  \begin{remark}
  \label{remark:k-local-for-greater-k}
    Let $\Delta$ be a $\kappa$-local map from $X$ to $Y$. For each non-negative integer $\kappa'$ such that $\kappa' \geq \kappa$, the map $\Delta$ is $\kappa'$-local.
  \end{remark}

  \begin{definition}
  \label{definition:local-map}
    Let $X$ and $Y$ be two subshifts of $Q^M$ and let $\Delta$ be a map from $X$ to $Y$. The map $\Delta$ is called \define{local}\graffito{local} if and only if there is a non-negative integer $\kappa$ such that it is $\kappa$-local.
  \end{definition}

  \begin{remark}
  \label{remark:local-map-if-and-only-if-cellular-automaton}
    Let $X$ and $Y$ be two subshifts of $Q^M$ and let $\Delta$ be a map from $X$ to $Y$. The map $\Delta$ is local if and only if it is the restriction to $X \to Y$ of the global transition function of a big-cellular automaton over $\mathcal{R}$ with set of states $Q$ and finite neighbourhood. 
  \end{remark}

  \begin{example}[{Sliding Block Codes \cite[Definition~1.5.1]{lind:marcus:1995}}]
    In the situation of \cref{example:full-shift-over-integers}, local maps are but sliding block codes.
  \end{example}

  \begin{example}[{From Golden Mean to Even Shift \cite[Example~1.5.6]{lind:marcus:1995}}]
  \label{example:from-golden-mean-to-even-shift-local-map}
    Let $X$ be the golden mean shift (\cref{example:shift:golden-mean}), let $Y$ be the even shift (\cref{example:shift:even}), let $\delta$ be the map from $X_{\setOf{0, 1}}$ to $\setOf{0, 1}$ given by $00 \mapsto 1$, $01 \mapsto 0$, and $10 \mapsto 0$, and let $\Delta$ be the map from $X$ to $Y$ given by $x \mapsto [z \mapsto \delta(n \mapsto x(z + n))]$. The map $\Delta$ is, by definition, local and it is, according to example~1.5.6 in \cite{lind:marcus:1995}, surjective.
  \end{example}

  Domain and codomain of a local map can be restricted simultaneously to a subset of cells and its interior, as is done in

  \begin{definition} 
  \label{definition:restriction-of-local-map}
    Let $\Delta$ be a $\kappa$-local map from $X$ to $Y$ and let $A$ be a subset of $M$. The map
    \begin{align*}
      \Delta_A^- \from X_A &\to Y_{A^{-\kappa}}, \mathnote{restriction $\Delta_A^-$ of $\Delta$ to $A$}\index[symbols]{DeltaAminus@$\Delta_A^-$}\\
      p &\mapsto \Delta(c)\restrictedTo_{A^{-\kappa}}, \text{ where } c \in X \text{ such that } c\restrictedTo_A = p,
    \end{align*}
    is called \define{restriction of $\Delta$ to $A$}.
  \end{definition}

  The image of a local map is a shift space, which is shown in

  \begin{lemma} 
  \label{lemma:image-of-local-map-is-subshift}
    Let $\Delta$ be a local map from $X$ to $Y$. Its image $\Delta(X)$ is a subshift of $Q^M$.
  \end{lemma}

  \begin{proof} 
    According to \cref{lemma:shift-space-if-and-only-if-invariant-and-compact}, the shift space $X$ is compact. And, according to \cref{remark:local-map-if-and-only-if-cellular-automaton} and \cref{corollary:Curtis-Hedlund-Lyndon}, the map $\Delta$ is continuous. Therefore, the topological space $\Delta(X)$ is compact and hence closed. Moreover, because $h \actsOnMap \Delta(X) = \Delta(h \actsOnMap X) = \Delta(X)$, the topological space $\Delta(X)$ is shift-invariant. Therefore, according to \cref{lemma:shift-space-if-and-only-if-invariant-and-compact}, the topological space $\Delta(X)$ is a subshift of $Q^M$.
  \end{proof}

  The difference of two points of the full shift is the set of cells in which they differ, as defined in

  \begin{definition} 
  \label{definition:difference-of-two-points-of-the-full-shift}
    Let $c$ and $c'$ be two points of $Q^M$. The set
    \begin{equation*}
      \differenceOf(c, c') = \setOf{m \in M \suchThat c(m) \neq c'(m)}
      \mathnote{difference $\differenceOf(c, c')$ of $c$ and $c'$}
    \end{equation*}
    is called \define{difference of $c$ and $c'$}.
  \end{definition}

  A local map is pre-injective if it is injective on points with finite support, as defined in

  \begin{definition} 
  \label{definition:pre-injective-for-local-maps}
    Let $\Delta$ be a local map from $X$ to $Y$. It is called \define{pre-injective}\graffito{pre-injective} if and only if, for each tuple $(x, x') \in X \times X$ such that $\differenceOf(x, x')$ is finite and $\Delta(x) = \Delta(x')$, we have $x = x'$.
  \end{definition}

  A bijective local map with local inverse is a conjugacy, and its domain and codomain are conjugate, which is defined in

  \begin{definition} 
  \label{definition:conjugacy}
    Let $\Delta$ be a local map from $X$ to $Y$. It is called \define{conjugacy}\graffito{conjugacy} if and only if it is bijective and its inverse is local.
  \end{definition}

  \begin{definition} 
  \label{definition:conjugate}
    Let $X$ and $Y$ be two shift spaces. They are called \define{conjugate}\graffito{conjugate} if and only if there is a conjugacy from $X$ to $Y$.
  \end{definition}

  Entropy of shift spaces is invariant under conjugacy, which is shown in

  \begin{lemma}
  \label{lemma:entropy-invariant-under-conjugacy}
    Let $\mathcal{R}$ be right amenable, let $\mathcal{F}$ be a right Følner net in $\mathcal{R}$, and let $X$ and $Y$ be two conjugate subshifts of $Q^M$. Then, $\entropyOf_{\mathcal{F}}(X) = \entropyOf_{\mathcal{F}}(Y)$.
  \end{lemma}

  \begin{proof}
    This is a direct consequence of \cref{remark:local-map-if-and-only-if-cellular-automaton} and \cref{theorem:entropy-does-non-increase}.
  \end{proof}

  Both the Moore and the Myhill property are introduced in

  \begin{definition} 
  \label{definition:Moore-and-Myhill-properties}
    Let $X$ be a subshift of $Q^M$. It is said to have the
    \begin{enumerate}
      \item \define{Moore property}\graffito{Moore property}\index{property!Moore} if and only if each surjective local map from $X$ to $X$ is pre-injective;
      \item \define{Myhill property}\graffito{Myhill property}\index{property!Myhill} if and only if each pre-injective local map from $X$ to $X$ is surjective. \qedhere
    \end{enumerate}
  \end{definition}

  \begin{remark}
  \label{remark:Moore-and-Myhill-are-invariant-under-conjugacy}
    Both the Moore and the Myhill property is invariant under conjugacy.
  \end{remark}

  \section{Tilings}
  \label{section:apart-tilings} 

  \paragraph{Contents.} A $\ntuple{\theta, \kappa, \theta'}$-tiling is a subset of cells such that the balls of radius $\theta$ about those cells are pairwise at least $\kappa + 1$ apart and the balls of radius $\theta'$ about those cells cover all cells (see \cref{definition:pairwise-at-least-apart,definition:theta-kappa-theta-prime-tiling}). If there are infinitely many cells, then, for each $\theta$ and $\kappa$, there is a $\ntuple{\theta, \kappa, 4 \theta + 2 \kappa}$-tiling (see \cref{theorem:there-is-a-rho-kappa-theta-tiling}). And, a subset of a strongly irreducible shift space has less entropy than that space if about each point of a $\ntuple{\theta, \kappa, \theta'}$-tiling the subset has fewer patterns with ball-shaped domains of radius $\theta$ than the space (see \cref{theorem:entropy-less-if-real-subsets-for-rho-kappa-theta-tiling}); in the proof of that statement we use \cref{lemma:technical-inequality-for-theorem-entropy-less-if-real-subsets-for-rho-kappa-theta-tiling,lemma:set-contained-in-union-of-balls-and-of-interior-of-set,corollary:subshifts-upper-bound-of-tiling-cap-Folner-net,lemma:subshift-with-at-least-two-points-has-at-least-two-patterns-for-each-domain}.

  \begin{definition} 
  \label{definition:pairwise-at-least-apart}
    Let $\family{A_j}_{j \in J}$ be a family of subsets of $M$ and let $\kappa$ be a non-negative integer. The family $\family{A_j}_{j \in J}$ is called \graffito{pairwise at least $\kappa + 1$ apart}\define{pairwise at least $\kappa + 1$ apart}\index{apart!pairwise at least $\kappa + 1$} if and only if 
    \begin{equation*}
      \ForEach j \in J \ForEach j' \in J \Holds \parens[\big]{j \neq j' \implies \distanceOf(A_j, A_{j'}) \geq \kappa + 1}. \qedhere
    \end{equation*}
  \end{definition}

  \begin{remark}
  \label{remark:pairwise-kappa-apart-versus-pairwise-disjoint}
    Each pairwise at least $\kappa + 1$ apart family is pairwise disjoint. And each pairwise disjoint family is pairwise at least $0 + 1$ apart.
  \end{remark}

  \begin{definition} 
  \label{definition:theta-kappa-theta-prime-tiling}
    Let $T$ be a subset of $M$, and let $\theta$, $\kappa$, and $\theta'$ be three non-negative integers. The set $T$ is called \defineX{$\ntuple{\theta, \kappa, \theta'}$-tiling of $\mathcal{R}$}{tiling of R@$\ntuple{\theta, \kappa, \theta'}$-tiling of $\mathcal{R}$}\graffito{$\ntuple{\theta, \kappa, \theta'}$-tiling $T$ of $\mathcal{R}$}\index{tilingthetakappathetaprime@$\ntuple{\theta, \kappa, \theta'}$-tiling of $\mathcal{R}$} if and only if the family $\family{\ball(t, \theta)}_{t \in T}$ is pairwise at least $\kappa + 1$ apart and the family $\family{\ball(t, \theta')}_{t \in T}$ is a cover of $M$.
  \end{definition}

  \begin{remark}
    According to \cref{remark:pairwise-kappa-apart-versus-pairwise-disjoint}, each $\ntuple{\theta, \kappa, \theta'}$-tiling of $\mathcal{R}$ is a $\ntuple{\ball(\theta), \ball(\theta')}$-tiling of $\mathcal{R}$, see \cref{definition:tiling}; and each $\ntuple{\ball(\theta), \ball(\theta')}$-tiling of $\mathcal{R}$ is a $\ntuple{\theta, 0, \theta'}$-tiling of $\mathcal{R}$.
  \end{remark}

  \begin{example}[Lattice]
    The $\ntuple{\ball(1), \ball(2)}$-tiling from \cref{example:lattice:tiling} is a $\ntuple{1, 1, 2}$-tiling of $\mathcal{R}$.
  \end{example}

  \begin{example}[Tree]
    The $\ntuple{\ball(1), \ball(2)}$-tiling from \cref{example:tree:tiling} is a $\ntuple{1, 0, 2}$-tiling of $\mathcal{R}$.
  \end{example}

  Greedily picking elements that are pairwise far enough apart yields a tiling, which we show in

  \begin{theorem}[Analogue of \cref{theorem:existence-of-tiling}] 
  \label{theorem:there-is-a-rho-kappa-theta-tiling}
    Let $M$ be infinite, and let $\theta$ and $\kappa$ be two non-negative integers. There is a countably infinite $\ntuple{\theta, \kappa, \theta'}$-tiling $T$ of $\mathcal{R}$, where $\theta' = 4 \theta + 2 \kappa$.
  \end{theorem}

  \begin{usage-note}
    In the proof, infiniteness of $M$ is used to deduce that spheres of arbitrarily big radii are non-empty.
  \end{usage-note}

  \begin{proof-sketch}
    From each of the spheres $\sphere(i (2 \theta + \kappa + 1))$, for $i \in \N_0$, pick elements that are pairwise at least $2 \theta + \kappa + 1$ apart and whose $(2 \theta + \kappa)$-closure covers the sphere --- they constitute a set $T$. The family $\family{\ball(t, \theta)}_{t \in T}$ is pairwise at least $\kappa + 1$ apart and the family $\family{\ball(t, 4 \theta + 2 \kappa)}_{t \in T}$ is a cover of $M$. See \cref{figure:there-is-a-rho-kappa-theta-tiling} for a schematic representation.
    \begin{figure}
      \myfloatalign
      \figureThereIsARhoKappaThetaTiling
      \caption{%
        Schematic representation of the set-up of the proof of \cref{theorem:there-is-a-rho-kappa-theta-tiling}. 
        The whole space is $M$; the dot in the centre is $m_0$; the smaller solid circle is $M_{i, 1} = \sphere(i (2 \theta + \kappa + 1))$ and the larger solid circle is $M_{i + 1, 1} = \sphere((i + 1) (2 \theta + \kappa + 1))$; the dots on the smaller solid circle are the elements of $M_i$; the region enclosed by the dotted circle about $m_{i,j}$ is $\ball(m_{i,j}, \theta)$; the region enclosed by the dash-dotted circle about $m_{i,j}$ is $\ball(m_{i,j}, \theta + \kappa)$; the region enclosed by the dashed circle about $m_{i,j}$ is $\ball(m_{i,j}, 2 \theta + \kappa)$; the dots labelled $m$, $m'$, and $m''$ are the respective elements from the last part of the proof.
      }
      \label{figure:there-is-a-rho-kappa-theta-tiling}
    \end{figure}
  \end{proof-sketch}

  \begin{proof} 
    Let $i \in \N_0$. Furthermore, let $M_{i,1} = \sphere(i (2 \theta + \kappa + 1))$. Then, because $M$ is infinite, according to \cref{corollary:m-infinite-if-and-only-if-spheres-are-non-empty}, the set $M_{i,1}$ is non-empty and finite. For $j \in \N_+$ in increasing order, if $M_{i,j} \smallsetminus \setOf{m_{i,1}, m_{i,2}, \dotsc, m_{i,j-1}} \neq \emptyset$, then choose $m_{i,j} \in M_{i,j} \smallsetminus \setOf{m_{i,1},\allowbreak m_{i,2},\allowbreak \dotsc,\allowbreak m_{i,j-1}}$ and put
    \begin{align*}
      M_{i,j+1} &= \setOf{m_{i,j}} \cup M_{i,j} \smallsetminus \ball(m_{i,j}, 2 \theta + \kappa)\\
                &= \setOf{m \in M_{i,j} \suchThat \distanceOf(m, m_{i,j}) = 0 \text{ or } \distanceOf(m, m_{i,j}) \geq 2 \theta + \kappa + 1};
    \end{align*}
    otherwise stop, put $j_i = j$ and put $M_i = M_{i, j_i}$ (see \cref{figure:there-is-a-rho-kappa-theta-tiling}).

    By construction,
    \begin{equation}
    \label{equation:there-is-a-rho-kappa-theta-tiling:apart}
      \ForEach m \in M_i \ForEach m' \in M_i \Holds \parens[\big]{m \neq m' \implies \distanceOf(m, m') \geq 2 \theta + \kappa + 1},
    \end{equation}
    and
    \begin{equation}
    \label{equation:there-is-a-rho-kappa-theta-tiling:base-point}
      \ForEach m \in M_{i,1} \Exists m' \in M_i \SuchThat \distanceOf(m, m') \leq 2 \theta + \kappa,
    \end{equation}
    and, for each $i' \in \N_0$ with $i' \neq i$, because $M_i \subseteq M_{i,1}$ and $M_{i'} \subseteq M_{i',1}$, and $M$ is infinite, according to \cref{corollary:distance-of-spheres},
    \begin{equation}
    \label{equation:there-is-a-rho-kappa-theta-tiling:sets-apart}
      \distanceOf(M_i, M_{i'}) \geq \distanceOf(M_{i,1}, M_{i',1}) \geq 2 \theta + \kappa + 1.
    \end{equation}

    Let $T = \bigcup_{i \in \N_0} M_i$. Because, for each $i \in \N_0$, the set $M_i$ is finite, the set $T$ is countable. And, because $\sequence{M_{i,1}}_{i \in \N_0}$ is pairwise disjoint, so is $\sequence{M_i}_{i \in \N_0}$ and hence $T$ is infinite. 
    Thus, $T$ is countably infinite.

    \proofPart{Subproof of: $\family{\ball(t, \theta)}_{t \in T}$ is pairwise at least $\kappa + 1$ apart}
    Let $t$, $t' \in T$ such that $t \neq t'$. If there is an $i \in \N_0$ such that $t$, $t' \in M_i$, then, according to \cref{equation:there-is-a-rho-kappa-theta-tiling:apart}, we have $\distanceOf(t, t') \geq 2 \theta + \kappa + 1$. Otherwise, there are $i$, $i' \in \N_0$ with $i \neq i'$ such that $t \in M_i$ and $t' \in M_{i'}$, and then, according to \cref{equation:there-is-a-rho-kappa-theta-tiling:sets-apart}, we have $\distanceOf(t, t') \geq \distanceOf(M_i, M_{i'}) \geq 2 \theta + \kappa + 1$. In conclusion, in both cases, according to \cref{lemma:distance-of-balls}, we have $\distanceOf(\ball(t, \theta), \ball(t', \theta)) \geq \kappa + 1$.

    \proofPart{Subproof of: $\family{\ball(t, 4 \theta + 2 \kappa)}_{t \in T}$ is a cover of $M$ (see \cref{figure:there-is-a-rho-kappa-theta-tiling})}
    Let $m \in M$. Then, there is an $i \in \N_0$ such that $i (2 \theta + \kappa + 1) \leq \lengthOf{m} < (i + 1) (2 \theta + \kappa + 1)$. Hence, according to \cref{lemma:distance-of-sphere-and-point},
    \begin{align*}
      \distanceOf(m, M_{i,1})
      &= \distanceOf(M_{i,1}, m)\\
      &\leq \distanceOf(m_0, m) - i (2 \theta + \kappa + 1)\\
      &< (i + 1) (2 \theta + \kappa + 1) - i (2 \theta + \kappa + 1)\\
      &= 2 \theta + \kappa + 1.
    \end{align*}
    Thus, $\distanceOf(m, M_{i,1}) \leq 2 \theta + \kappa$. Therefore, there is an $m' \in M_{i,1}$ such that $\distanceOf(m, m') \leq 2 \theta + \kappa$. Moreover, according to \cref{equation:there-is-a-rho-kappa-theta-tiling:base-point}, there is an $m'' \in M_i$ such that $\distanceOf(m', m'') \leq 2 \theta + \kappa$. Thus,
    \begin{equation*}
      \distanceOf(m, m'') \leq \distanceOf(m, m') + \distanceOf(m', m'') \leq 4 \theta + 2 \kappa.
    \end{equation*}
    Hence, $m \in \ball(m'', 4 \theta + 2 \kappa)$. Therefore, because $m'' \in T$, we have $m \in \bigcup_{t \in T} \ball(t, 4 \theta + 2 \kappa)$. In conclusion, $\bigcup_{t \in T} \ball(t, 4 \theta + 2 \kappa) = M$. 
  \end{proof}

  One by one banning subpatterns with similar domains that are far enough apart in a set of finite patterns, decreases its size by at least a multiplicative constant between $0$ and $1$ in each step. In other words, the number of finite patterns with a fixed domain, excluding those in which some subpatterns with similar domains that are far enough apart are embedded, is bounded above by some constant between $0$ and $1$ raised to the power of the number of forbidden subpatterns times the number of all finite patterns with the fixed domain, which is shown in

  \begin{lemma} 
  \label{lemma:technical-inequality-for-theorem-entropy-less-if-real-subsets-for-rho-kappa-theta-tiling}
    Let $X$ be a non-empty and $\kappa$-strongly irreducible subshift of $Q^M$, let $F$ be a finite subset of $M$, let $\theta$ be a non-negative integer, let $T$ be a subset of $M$ such that the family $\family{\ball(t, \theta)}_{t \in T}$ is pairwise at least $\kappa + 1$ apart, and, for each element $t \in T$, let $p_t$ be a pattern of $X_{\ball(t, \theta)}$. Furthermore, let $\xi$ be the positive integer $\cardinalityOf{X_{\ball(\theta)^{+\kappa}}}$, let $S$ be the finite set $T \cap F^{-(\theta + \kappa)}\ (= \setOf{t \in T \suchThat \ball(t, \theta)^{+\kappa} \subseteq F})$, and, for each element $s \in S$, let $\pi_s$ be the map $X_F \to X_{\ball(s, \theta)}$, $p \mapsto p\restrictedTo_{\ball(s, \theta)}$ (see \cref{figure:technical-inequality-for-theorem-entropy-less-if-real-subsets-for-rho-kappa-theta-tiling}).
    \begin{figure}
      \myfloatalign
      \figureTechnicalInequalityForTheoremEntropyLessIfRealSubsetsForRhoKappaThetaTiling
      \caption{Schematic representation of the set-up of \cref{lemma:technical-inequality-for-theorem-entropy-less-if-real-subsets-for-rho-kappa-theta-tiling}.}
      \label{figure:technical-inequality-for-theorem-entropy-less-if-real-subsets-for-rho-kappa-theta-tiling}
    \end{figure}
    Then, 
    \begin{equation*}
      \cardinalityOf{X_F \smallsetminus \bigcup_{s \in S} \pi_s^{-1}(p_s)}
      \leq
      (1 - \xi^{-1})^{\cardinalityOf{S}} \cdot \cardinalityOf{X_F}. \qedhere
    \end{equation*}
  \end{lemma}

  \begin{usage-note}
    In the proof, $\kappa$-strong irreducibility is used to extend an in $X$ allowed $F \smallsetminus \ball(s, \theta)^{+\kappa}$-pattern by the $\ball(s, \theta)$-pattern $p_s$ and a $\boundaryOf_\kappa^+ \ball(s, \theta)$-pattern to an in $X$ allowed $F$-pattern. 
  \end{usage-note}

  \begin{proof-sketch}
    Let $\family{s_j}_{j \in \setOf{1, 2, \dotsc, \cardinalityOf{S}}}$ be an enumeration of $S$, let $Z_0 = X_F$, and, for each $\vartheta \in \setOf{0, 1, \dotsc, \cardinalityOf{S} - 1}$, let $Z_{\vartheta + 1} = Z_\vartheta \smallsetminus \parens{\pi_{s_{\vartheta + 1}}^{-1}(p_{s_{\vartheta + 1}}) \cap Z_\vartheta}$. Furthermore, let $\vartheta \in \setOf{0, 1, \dotsc, \cardinalityOf{S} - 1}$. Then, $\cardinalityOf{Z_\vartheta} \leq \cardinalityOf{X_{\ball(s_{\vartheta + 1}, \theta)^{+ \kappa}}} \cdot \cardinalityOf{(Z_\vartheta)_{F \smallsetminus \ball(s_{\vartheta + 1}, \theta)^{+ \kappa}}} = \xi \cdot \cardinalityOf{(Z_\vartheta)_{F \smallsetminus \ball(s_{\vartheta + 1}, \theta)^{+ \kappa}}}$. And, because $X$ is $\kappa$-strongly irreducible, each pattern of $(Z_\vartheta)_{F \smallsetminus \ball(s_{\vartheta + 1}, \theta)^{+ \kappa}}$ can be extended by $p_{s_{\vartheta + 1}}$ and a pattern with domain $\ball(s_{\vartheta + 1}, \theta)^{+ \kappa} \smallsetminus \ball(s_{\vartheta + 1}, \theta)$ to a pattern of $\pi_{s_{\vartheta + 1}}^{-1}(p_{s_{\vartheta + 1}}) \cap Z_\vartheta$ and thus $\cardinalityOf{(Z_\vartheta)_{F \smallsetminus \ball(s_{\vartheta + 1}, \theta)^{+ \kappa}}} \leq \cardinalityOf{\pi_{s_{\vartheta + 1}}^{-1}(p_{s_{\vartheta + 1}}) \cap Z_\vartheta}$. Hence, $\cardinalityOf{\pi_{s_{\vartheta + 1}}^{-1}(p_{s_{\vartheta + 1}}) \cap Z_\vartheta} \geq \xi^{-1} \cdot \cardinalityOf{Z_\vartheta}$. Therefore,
      $\cardinalityOf{Z_{\vartheta + 1}}
      =    \cardinalityOf{Z_\vartheta} - \cardinalityOf{\pi_{s_{\vartheta + 1}}^{-1}(p_{s_{\vartheta + 1}}) \cap Z_\vartheta}
      \leq (1 - \xi^{-1}) \cdot \cardinalityOf{Z_\vartheta}$.
    The statement follows by induction.
  \end{proof-sketch}

  \begin{proof}
    As claimed, because $X \neq \emptyset$, the integer $\xi$ is positive; because, according to \cref{corollary:ball-centred-at-mzero-to-m} and \cref{item:characterisation-of-k-closure-and-interior-of-balls:closure} of \cref{corollary:characterisation-of-k-closure-and-interior-of-balls}, $s \isSemiActedUponBy \ball(\theta + \kappa) = \ball(s, \theta)^{+\kappa}$, we have $S = T \cap F^{-(\theta + \kappa)} = \setOf{t \in T \suchThat \ball(t, \theta)^{+\kappa} \subseteq F}$; and, because $S \subseteq \bigcup_{s \in S} \ball(s, \theta)^{+\kappa} \subseteq F$ and $F$ is finite, the set $S$ is finite.

    Let $\family{B_s}_{s \in S} = \family{\ball(s, \theta)}_{s \in S}$, let $\family{s_j}_{j \in \setOf{1, 2, \dotsc, \cardinalityOf{S}}}$ be an enumeration of $S$, and, for each $\vartheta \in \setOf{0, 1, \dotsc, \cardinalityOf{S}}$, let 
    \begin{align*}
      Z_\vartheta
      &= X_F \smallsetminus \bigcup_{j = 1}^\vartheta \pi_{s_j}^{-1}(p_{s_j})\\
      \big(&= \setOf{p \in X_F \suchThat \ForEach j \in \setOf{1, 2, \dotsc, \vartheta} \Holds p\restrictedTo_{B_{s_j}} \neq p_{s_j}}\big).
    \end{align*}
    To establish the claim, we prove by induction on $\vartheta$, that, for each $\vartheta \in \setOf{0, 1, \dotsc, \cardinalityOf{S}}$,
    \begin{equation}
    \label{equation:entropy-less-if-real-subsets-for-rho-kappa-theta-tiling:inductive-hypothesis}
      \cardinalityOf{Z_\vartheta}
      \leq
      (1 - \xi^{-1})^\vartheta \cdot \cardinalityOf{X_F}.
    \end{equation}

    \proofPart{Base Case}
      Let $\vartheta = 0$. Then, because $\bigcup_{j = 1}^0 \pi_{s_j}^{-1}(p_{s_j}) = \emptyset$ and $(1 - \xi^{-1})^0 = 1$, \cref{equation:entropy-less-if-real-subsets-for-rho-kappa-theta-tiling:inductive-hypothesis} holds. Note that $0^0 = 1$.

    \proofPart{Inductive Step}
      Let $\vartheta \in \setOf{0, 1, \dotsc, \cardinalityOf{S} - 1}$ such that \cref{equation:entropy-less-if-real-subsets-for-rho-kappa-theta-tiling:inductive-hypothesis}, called \emph{inductive hypothesis}, holds. Furthermore, let $Z = Z_\vartheta$. Because $B_{s_{\vartheta + 1}}^{+\kappa} \subseteq F$, 
      \begin{equation*}
        Z \subseteq Z_{B_{s_{\vartheta + 1}}^{+\kappa}} \times Z_{F \smallsetminus B_{s_{\vartheta + 1}}^{+\kappa}}
          \subseteq X_{B_{s_{\vartheta + 1}}^{+\kappa}} \times Z_{F \smallsetminus B_{s_{\vartheta + 1}}^{+\kappa}}.
      \end{equation*}
      Hence, $\cardinalityOf{Z} \leq \cardinalityOf{X_{B_{s_{\vartheta + 1}}^{+\kappa}}} \cdot \cardinalityOf{Z_{F \smallsetminus B_{s_{\vartheta + 1}}^{+\kappa}}}$. Moreover, according to \cref{remark:pattern-belongs-to-shift-if-and-only-if-translated-pattern-belongs}, we have $\cardinalityOf{X_{B_{s_{\vartheta + 1}}^{+\kappa}}} = \cardinalityOf{X_{\ball(\theta)^{+\kappa}}} = \xi$ (where we used that $B_{s_{\vartheta + 1}}^{+\kappa} = s_{\vartheta + 1} \isSemiActedUponBy \ball(\theta)^{+\kappa}$, which holds according to \cref{corollary:ball-centred-at-mzero-to-m} and \cref{item:characterisation-of-k-closure-and-interior-of-balls:closure} of \cref{corollary:characterisation-of-k-closure-and-interior-of-balls}). Therefore, $\cardinalityOf{Z} \leq \xi \cdot \cardinalityOf{Z_{F \smallsetminus B_{s_{\vartheta + 1}}^{+\kappa}}}$. 

      Let $p \in Z_{F \smallsetminus B_{s_{\vartheta + 1}}^{+\kappa}}$. Then, $p \in X_{F \smallsetminus B_{s_{\vartheta + 1}}^{+\kappa}}$. Moreover, according to \cref{lemma:distance-of-set-minus-closure-to-set}, we have $\distanceOf(B_{s_{\vartheta + 1}}, F \smallsetminus B_{s_{\vartheta + 1}}^{+\kappa}) \geq \kappa + 1$. Hence, because $X$ is $\kappa$-strongly irreducible, there is a $p'' \in X_F$ such that $p''\restrictedTo_{\domainOf(p)} = p$ and $p''\restrictedTo_{\domainOf(p_{s_{\vartheta + 1}})} = p_{s_{\vartheta + 1}}$.
      Furthermore, because $\family{\ball(t, \theta)}_{t \in T}$ is pairwise at least $\kappa + 1$ apart, for each $j \in \setOf{1, 2, \dotsc, \vartheta}$, we have $B_{s_j} \subseteq F \smallsetminus B_{s_{\vartheta + 1}}^{+\kappa}$. Therefore, for each $j \in \setOf{1, 2, \dotsc, \vartheta}$, we have $p''\restrictedTo_{B_{s_j}} = p\restrictedTo_{B_{s_j}} \neq p_{s_j}$ and hence $p'' \notin \pi_{s_j}^{-1}(p_{s_j})$. Thus, $p'' \in Z$. Moreover, $p'' \in \pi_{s_{\vartheta + 1}}^{-1}(p_{s_{\vartheta + 1}})$. Therefore, $\cardinalityOf{Z_{F \smallsetminus B_{s_{\vartheta + 1}}^{+\kappa}}} \leq \cardinalityOf{\pi_{s_{\vartheta + 1}}^{-1}(p_{s_{\vartheta + 1}}) \cap Z}$. 

      Together,
      \begin{equation*}
        \cardinalityOf{Z} \leq \xi \cdot \cardinalityOf{\pi_{s_{\vartheta + 1}}^{-1}(p_{s_{\vartheta + 1}}) \cap Z}.
      \end{equation*}
      Because $Z_{\vartheta + 1} = Z \smallsetminus \pi_{s_{\vartheta + 1}}^{-1}(p_{s_{\vartheta + 1}})$,
      \begin{align*}
        \cardinalityOf{Z_{\vartheta + 1}}
        &=    \cardinalityOf{Z \smallsetminus \pi_{s_{\vartheta + 1}}^{-1}(p_{s_{\vartheta + 1}})}\\
        &=    \cardinalityOf{Z \smallsetminus (\pi_{s_{\vartheta + 1}}^{-1}(p_{s_{\vartheta + 1}}) \cap Z)}\\ 
        &=    \cardinalityOf{Z} - \cardinalityOf{\pi_{s_{\vartheta + 1}}^{-1}(p_{s_{\vartheta + 1}}) \cap Z}\\
        &\leq \cardinalityOf{Z} - \xi^{-1} \cdot \cardinalityOf{Z}\\
        &=    (1 - \xi^{-1}) \cdot \cardinalityOf{Z}.
      \end{align*}
      Hence, according to the inductive hypothesis,
      \begin{align*}
        \cardinalityOf{Z_{\vartheta + 1}}
        &\leq (1 - \xi^{-1}) \cdot (1 - \xi^{-1})^\vartheta \cdot \cardinalityOf{X_F}\\
        &=    (1 - \xi^{-1})^{\vartheta + 1} \cdot \cardinalityOf{X_F}.
      \end{align*}

    In conclusion, according to the principle of mathematical induction, for each $\vartheta \in \setOf{0, 1, \dotsc, \cardinalityOf{S}}$, \cref{equation:entropy-less-if-real-subsets-for-rho-kappa-theta-tiling:inductive-hypothesis} holds.
  \end{proof}

  The number of elements in a finite set is bounded above by the number of elements its interior shares with a tiling times the number of elements of a big enough ball plus the number of elements of a big enough boundary of the finite set, which is shown in

  \begin{lemma} 
  \label{lemma:set-contained-in-union-of-balls-and-of-interior-of-set}
    Let $F$ be a finite subset of $M$, let $\theta$, $\kappa$, and $\theta'$ be three non-negative integers, let $T$ be a subset of $M$ such that $\family{\ball(t, \theta')}_{t \in T}$ is a cover of $M$, and let $S$ be the finite set $T \cap F^{-(\theta + \kappa)}\ (= \setOf{t \in T \suchThat \ball(t, \theta)^{+\kappa} \subseteq F})$ (see \cref{figure:set-contained-in-union-of-balls-and-of-interior-of-set}).
    \begin{figure}
      \myfloatalign
      \figureSetContainedInUnionOfBallsAndOfInteriorOfSet
      \caption{Schematic representation of the set-up of \cref{lemma:set-contained-in-union-of-balls-and-of-interior-of-set}.}
      \label{figure:set-contained-in-union-of-balls-and-of-interior-of-set}
    \end{figure}
    Then,
    \begin{equation*}
      \cardinalityOf{F} \leq \cardinalityOf{S} \cdot \cardinalityOf{\ball(\theta')} + \cardinalityOf{\boundaryOf_{\theta + \kappa + \theta'}^- F}. \qedhere
    \end{equation*}%
  \end{lemma}

  \begin{proof-sketch}
    For each $m \in F$, if $m \notin S^{+ \theta'}$, then $m \notin F^{-(\theta + \kappa + \theta')}$. Thus, $F \subseteq S^{+ \theta'} \cup F \smallsetminus F^{-(\theta + \kappa + \theta')} = \parens{\bigcup_{s \in S} \ball(s, \theta')} \cup \boundaryOf_{\theta + \kappa + \theta'}^- F$. Hence, $\cardinalityOf{F} \leq \cardinalityOf{S} \cdot \cardinalityOf{\ball(\theta')} + \cardinalityOf{\boundaryOf_{\theta + \kappa + \theta'}^- F}$.
  \end{proof-sketch}

  \begin{proof}
    Let $m \in F \smallsetminus \bigcup_{s \in S} \ball(s, \theta')$. Because $\family{\ball(t, \theta')}_{t \in T}$ is a cover of $M$, there is a $t \in T$ such that $m \in \ball(t, \theta')$. Because $m \notin \bigcup_{s \in S} \ball(s, \theta')$, we have $t \notin S$ and hence $\ball(t, \theta)^{+\kappa} \nsubseteq F$. Because $m \in \ball(t, \theta')$, we have $\distanceOf(m, t) \leq \theta'$ and hence $t \in \ball(m, \theta') = \setOf{m}^{+\theta'}$. Suppose that $m \in F^{-(\theta + \kappa + \theta')}$. Then, according to \cref{item:characterisation-of-k-closure-and-interior:closure} of \cref{lemma:characterisation-of-k-closure-and-interior}, \cref{item:properties-of-interior-closure-and-boundary:subset} of \cref{lemma:properties-of-interior-closure-and-boundary}, and \cref{item:repeated-k-boundaries-etc:closure-interior} of \cref{lemma:repeated-k-boundaries-etc},
    \begin{align*} 
      \ball(t, \theta)^{+\kappa}
      &= (\setOf{t}^{+\theta})^{+\kappa}\\
      &\subseteq \parens[\Big]{\parens[\big]{\setOf{m}^{+\theta'}}^{+\theta}}^{+\kappa}\\
      &= \setOf{m}^{+(\theta + \kappa + \theta')}\\
      &\subseteq (F^{-(\theta + \kappa + \theta')})^{+(\theta + \kappa + \theta')}\\
      &= F,
    \end{align*}
    which contradicts that $\ball(t, \theta)^{+\kappa} \nsubseteq F$. Hence, $m \in F \smallsetminus F^{-(\theta + \kappa + \theta')} = \boundaryOf_{\theta + \kappa + \theta'}^- F$. Therefore,
    \begin{equation*}
      F \subseteq \parens[\big]{\bigcup_{s \in S} \ball(s, \theta')} \cup \boundaryOf_{\theta + \kappa + \theta'}^- F.
    \end{equation*}
    Moreover, because $S \subseteq \bigcup_{s \in S} \ball(s, \theta)^{+\kappa} \subseteq F$ and $F$ is finite, the set $S$ is finite. And, for each $s \in S$, according to \cref{corollary:balls-of-equal-radius-have-same-number-of-elements}, we have $\cardinalityOf{\ball(s, \theta')} = \cardinalityOf{\ball(\theta')}$. In conclusion,
    \begin{align*}
      \cardinalityOf{F} &\leq \sum_{s \in S} \cardinalityOf{\ball(\theta')} + \cardinalityOf{\boundaryOf_{\theta + \kappa + \theta'}^- F}\\
              &=    \cardinalityOf{S} \cdot \cardinalityOf{\ball(\theta')} + \cardinalityOf{\boundaryOf_{\theta + \kappa + \theta'}^- F}. \qedhere
    \end{align*}
  \end{proof}

  The number of elements that the components of a right Følner net share with a tiling is asymptotically bounded below away from zero, which is shown in

  \begin{corollary}[Analogue of \cref{lemma:upper-bound-of-tiling-cap-Folner-net}] 
  \label{corollary:subshifts-upper-bound-of-tiling-cap-Folner-net} 
    Let $\net{F_i}_{i \in I}$ be a right Erling net in $\mathcal{R}$ (which exists according to \cref{lemma:finitely-right-generated-is-tractable}), let $\theta$, $\kappa$, and $\theta'$ be three non-negative integers, let $T$ be a subset of $M$ such that $\family{\ball(t, \theta')}_{t \in T}$ is a cover of $M$. There is a positive real number $\varepsilon \in \R_{> 0}$ and there is an index $i_0 \in I$ such that, for each index $i \in I$ with $i \geq i_0$, we have $\cardinalityOf{T \cap F_i^{-(\theta + \kappa)}} \geq \varepsilon \cardinalityOf{F_i}$.
  \end{corollary}

  \begin{proof}
    Let $i \in I$ and let $T_i = T \cap F_i^{-(\theta + \kappa)}$. Then, according to \cref{lemma:set-contained-in-union-of-balls-and-of-interior-of-set}, we have $\cardinalityOf{F_i} \leq \cardinalityOf{T_i} \cdot \cardinalityOf{\ball(\theta')} + \cardinalityOf{\boundaryOf_{\theta + \kappa + \theta'}^- F_i}$. Hence,
    \begin{equation*}
      \frac{\cardinalityOf{T_i}}{\cardinalityOf{F_i}} \geq \frac{1}{\cardinalityOf{\ball(\theta')}} \cdot \parens*{1 - \frac{\cardinalityOf{\boundaryOf_{\theta + \kappa + \theta'}^- F_i}}{\cardinalityOf{F_i}}}.
    \end{equation*}
    Moreover, because $\net{F_i}_{i \in I}$ is a right Erling net, there is a $\xi \in \leftClosedAndRightOpenInterval{0, 1}$ and there is an $i_0 \in I$ such that
    \begin{equation*}
      \ForEach i \in I \Holds \parens*{i \geq i_0 \implies \frac{\cardinalityOf{\boundaryOf_{\theta + \kappa + \theta'}^- F_i}}{\cardinalityOf{F_i}} \leq \xi}.
    \end{equation*}
    Let $\varepsilon = (1/\ball(\theta')) \cdot (1 - \xi)$. Then, for each $i \in I$ with $i \geq i_0$, we have $\cardinalityOf{T_i}/\cardinalityOf{F_i} \geq \varepsilon$.
  \end{proof}

  If a shift space has at least two points, then, for each non-empty domain, it has at least two patterns.

  \begin{lemma}
  \label{lemma:subshift-with-at-least-two-points-has-at-least-two-patterns-for-each-domain}
    Let $X$ be a subshift of $Q^M$ such that $\cardinalityOf{X} \geq 2$ and let $A$ be a non-empty subset of $M$. Then, $\cardinalityOf{X_A} \geq 2$.
  \end{lemma}

  \begin{proof}
    Because $\cardinalityOf{X} \geq 2$, there are $x$ and $x'$ in $X$ such that $x \neq x'$. Thus, there is an $m \in M$ such that $x(m) \neq x'(m)$. And, because $A$ is non-empty, there is an $a \in A$. The element $h = g_{m_0, a} g_{m_0, m}^{-1}$ is contained in $H$ and satisfies $h^{-1} \actsOnPoint a = m$. Hence, $(h \actsOnMap x)(a) \neq (h \actsOnMap x')(a)$. Therefore, $(h \actsOnMap x)\restrictedTo_A$ and $(h \actsOnMap x')\restrictedTo_A$ are distinct and are contained in $X_A$. In conclusion, $\cardinalityOf{X_A} \geq 2$.
  \end{proof}

  A subset of a strongly irreducible shift space has less entropy than that space if about each point of a tiling the subset has fewer patterns of a certain radius than the space, which is shown in

  \begin{theorem}[Analogue of \cref{lemma:entropy-bounded-above-if-strange-tiling-exists}] 
  \label{theorem:entropy-less-if-real-subsets-for-rho-kappa-theta-tiling}
    Let $\mathcal{F} = \net{F_i}_{i \in I}$ be a right Erling net in $\mathcal{R}$ (which exists according to \cref{lemma:finitely-right-generated-is-tractable}), let $X$ be a $\kappa$-strongly irreducible subshift of $Q^M$ such that $\cardinalityOf{X} \geq 2$, let $Y$ be a subset of $X$, and let $T$ be a $\ntuple{\theta, \kappa, \theta'}$-tiling of $\mathcal{R}$ such that, for each element $t \in T$, we have $Y_{\ball(t, \theta)} \subsetneqq X_{\ball(t, \theta)}$. Then, $\entropyOf_{\mathcal{F}}(Y) < \entropyOf_{\mathcal{F}}(X)$.
  \end{theorem} 

  \begin{usage-note}
    In the proof, $\kappa$-strong irreducibility is used to apply \cref{lemma:technical-inequality-for-theorem-entropy-less-if-real-subsets-for-rho-kappa-theta-tiling} yielding the inequality $\cardinalityOf{X_{F_i} \smallsetminus \bigcup_{t \in T_i} \pi_{i, t}^{-1}(p_t)} \leq (1 - \xi^{-1})^{\cardinalityOf{T_i}} \cdot \cardinalityOf{X_{F_i}}$.
  \end{usage-note}

  \begin{proof-sketch} 
    Let $p_t \in X_{\ball(t, \theta)} \smallsetminus Y_{\ball(t, \theta)}$ and let $T_i = T \cap F_i^{-(\theta + \kappa)}$. Then, $Y_{F_i} \subseteq X_{F_i} \smallsetminus \bigcup_{t \in T_i} \pi_{i, t}^{-1}(p_t)$. Hence, $\cardinalityOf{Y_{F_i}} \leq (1 - \xi^{-1})^{\cardinalityOf{T_i}} \cdot \cardinalityOf{X_{F_i}}$. Therefore, $\log\cardinalityOf{Y_{F_i}}/\cardinalityOf{F_i} \leq \log(1 - \xi^{-1}) \cdot \cardinalityOf{T_i}/\cardinalityOf{F_i} + \log\cardinalityOf{X_{F_i}}/\cardinalityOf{F_i}$. In conclusion, because $\log(1 - \xi^{-1}) < 0$ and $\net{\cardinalityOf{T_i}/\cardinalityOf{F_i}}_{i \in I}$ is eventually bounded below away from zero, we have $\entropyOf_{\mathcal{F}}(Y) < \entropyOf_{\mathcal{F}}(X)$.
  \end{proof-sketch}

  \begin{proof} 
    For each $t \in T$, because $Y_{\ball(t, \theta)} \subsetneqq X_{\ball(t, \theta)}$, we have $X_{\ball(t, \theta)} \smallsetminus Y_{\ball(t, \theta)} \neq \emptyset$. Let $\family{p_t}_{t \in T}$ be a transversal of $\family{X_{\ball(t, \theta)} \smallsetminus Y_{\ball(t, \theta)}}_{t \in T}$ and let $\xi = \cardinalityOf{X_{\ball(\theta)^{+\kappa}}}$. Furthermore, let $i \in I$, let $T_i = T \cap F_i^{-(\theta + \kappa)}\ (= \setOf{t \in T \suchThat \ball(t, \theta)^{+\kappa} \subseteq F_i})$ and, for each $t \in T_i$, let $\pi_{i, t} \from X_{F_i} \to X_{\ball(t, \theta)}$, $p \mapsto p\restrictedTo_{\ball(t, \theta)}$. Note that, because $\cardinalityOf{X} \geq 2$, according to \cref{lemma:subshift-with-at-least-two-points-has-at-least-two-patterns-for-each-domain}, we have $\xi \geq 2$ and hence $1 - \xi^{-1} > 0$. 

    Because $\family{\ball(t, \theta)}_{t \in T}$ is pairwise at least $\kappa + 1$ apart, according to \cref{lemma:technical-inequality-for-theorem-entropy-less-if-real-subsets-for-rho-kappa-theta-tiling},
    \begin{equation*}
      \cardinalityOf{X_{F_i} \smallsetminus \bigcup_{t \in T_i} \pi_{t}^{-1}(p_{t})}
      \leq
      (1 - \xi^{-1})^{\cardinalityOf{T_i}} \cdot \cardinalityOf{X_{F_i}}.
    \end{equation*}
    For each $t \in T_i$, because $p_t \notin Y_{\ball(t, \theta)}$, we have $\pi_{i, t}^{-1}(p_t) \cap Y_{F_i} = \emptyset$. Hence, $\parens*{\bigcup_{t \in T_i} \pi_{i, t}^{-1}(p_t)} \cap Y_{F_i} = \emptyset$. Therefore,
    \begin{align*}
      \cardinalityOf{Y_{F_i}}
      &=    \cardinalityOf{Y_{F_i} \smallsetminus \bigcup_{t \in T_i} \pi_{i, t}^{-1}(p_t)}\\
      &\leq \cardinalityOf{X_{F_i} \smallsetminus \bigcup_{t \in T_i} \pi_{i, t}^{-1}(p_t)}\\
      &\leq (1 - \xi^{-1})^{\cardinalityOf{T_i}} \cdot \cardinalityOf{X_{F_i}}.
    \end{align*}
    Thus,
    \begin{equation*}
      \frac{\log\cardinalityOf{Y_{F_i}}}{\cardinalityOf{F_i}}
      \leq \frac{\cardinalityOf{T_i}}{\cardinalityOf{F_i}} \cdot \log(1 - \xi^{-1})
           + \frac{\log\cardinalityOf{X_{F_i}}}{\cardinalityOf{F_i}}.
    \end{equation*}
    Because $\family{\ball(t, \theta')}_{t \in T}$ is a cover of $M$, according to \cref{corollary:subshifts-upper-bound-of-tiling-cap-Folner-net}, there is an $\varepsilon \in \R_{> 0}$ and there is an $i_0 \in I$ such that
    \begin{equation*}
      \ForEach i \in I \Holds \parens*{i \geq i_0 \implies \frac{\cardinalityOf{T_i}}{\cardinalityOf{F_i}} \geq \varepsilon}.
    \end{equation*}
    Hence, because $\log(1 - \xi^{-1}) < 0$,
    \begin{equation*}
      \entropyOf_{\mathcal{F}}(Y) \leq \varepsilon \cdot \log(1 - \xi^{-1}) + \entropyOf_{\mathcal{F}}(X)
                                <    \entropyOf_{\mathcal{F}}(X). \qedhere
    \end{equation*}
  \end{proof}

  \section{The Moore and the Myhill Properties} 
  \label{section:gardens-of-Eden-on-shifts}

  \paragraph{Contents.} The image of a local map to a strongly irreducible shift space that is not surjective does not have maximal entropy (see \cref{theorem:subshift-not-surjective-implies-less-entropy}). And the converse of that statement obviously holds. Moreover, a local map from a strongly irreducible shift space of finite type whose image has less entropy than its domain is not pre-injective (see \cref{theorem:subshift-less-entropy-implies-not-pre-injective}). And the converse of that statement also holds (see \cref{theorem:subshift-not-pre-injective-implies-less-entropy}). These four statements establish the Garden of Eden theorem (see Main \cref{theorem:subshift-garden-of-Eden}). It follows that strongly irreducible shift spaces of finite type have the Moore and the Myhill property (see \cref{corollary:Moore-and-Myhill}). 


  \paragraph{Body.} Because a local map that is not surjective has a Garden of Eden pattern, the entropy of its image is not maximal, which is shown in

  \begin{theorem}[Analogue of \cref{theorem:not-surjective-implies-less-entropy}] 
  \label{theorem:subshift-not-surjective-implies-less-entropy}
    Let $M$ be infinite, let $X$ be a non-empty subshift of $Q^M$, let $Y$ be a strongly irreducible subshift of $Q^M$, let $\Delta$ be a local map from $X$ to $Y$ that is not surjective, and let $\mathcal{F}$ be a right Erling net in $\mathcal{R}$ (which exists according to \cref{lemma:finitely-right-generated-is-tractable}). Then, $\entropyOf_{\mathcal{F}}(\Delta(X)) < \entropyOf_{\mathcal{F}}(Y)$.
  \end{theorem}

  \begin{usage-note}
    In the proof, infiniteness of $M$ is used to apply \cref{theorem:there-is-a-rho-kappa-theta-tiling} yielding a tiling, locality of $\Delta$ is used to apply \cref{lemma:image-of-local-map-is-subshift} yielding that $\Delta(X)$ is a subshift of $Q^M$, and strong irreducibility of $Y$ is used to apply \cref{theorem:entropy-less-if-real-subsets-for-rho-kappa-theta-tiling} yielding a strict inequality for entropies.
  \end{usage-note}

  \begin{proof-sketch} 
    Because $\Delta$ is not surjective, there is a Garden of Eden configuration. Thus, because $\Delta$ is local, there is a Garden of Eden pattern. Hence, there are too many Garden of Eden configurations for the entropy to be maximal. 
  \end{proof-sketch}

  \begin{proof}
    Because $Y$ is strongly irreducible, there is a $\kappa \in \N_0$ such that $Y$ is $\kappa$-strongly irreducible. And, because $\Delta$ is not surjective, there is a $y \in Y \smallsetminus \Delta(X)$. Hence, according to \cref{lemma:if-restrictions-to-balls-are-in-shift-then-so-is-pattern}, there is a $\rho \in \N_0$ such that $y\restrictedTo_{\ball(\rho)} \notin (\Delta(X))_{\ball(\rho)}$ and thus $y\restrictedTo_{\ball(\rho)} \in Y_{\ball(\rho)} \smallsetminus (\Delta(X))_{\ball(\rho)}$. And, because $M$ is infinite, according to \cref{theorem:there-is-a-rho-kappa-theta-tiling}, there is a $\ntuple{\rho, \kappa, \theta'}$-tiling $T$ of $\mathcal{R}$.

    According to \cref{lemma:image-of-local-map-is-subshift}, the set $\Delta(X)$ is a subshift of $Q^M$. And, for each $t \in T$, according to \cref{corollary:ball-centred-at-mzero-to-m}, we have $t \isSemiActedUponBy \ball(\rho) = \ball(t, \rho)$. Therefore, for each $t \in T$, because $t \actsByItsCoordinateOn \blank$ is bijective and according to \cref{remark:pattern-belongs-to-shift-if-and-only-if-translated-pattern-belongs}, we have $t \actsByItsCoordinateOn (y\restrictedTo_{\ball(\rho)}) \in Y_{\ball(t, \rho)} \smallsetminus (\Delta(X))_{\ball(t, \rho)}$ and thus $(\Delta(X))_{\ball(t, \rho)} \subsetneqq Y_{\ball(t, \rho)}$. 

    Because $X$ is non-empty and $\Delta$ is not surjective, we have $\cardinalityOf{Y} \geq 2$. In conclusion, because $Y$ is $\kappa$-strongly irreducible, according to \cref{theorem:entropy-less-if-real-subsets-for-rho-kappa-theta-tiling}, we have $\entropyOf_{\mathcal{F}}(\Delta(X)) < \entropyOf_{\mathcal{F}}(Y)$.
  \end{proof}

  If there are less patterns in the codomain of a local map than in its domain, at least two patterns have the same image, which is shown in

  \begin{lemma} 
  \label{lemma:there-are-patterns-that-agree-on-boundary-and-have-the-same-image-under-Delta-restriction-k}
    Let $X$ be a $\kappa$-strongly irreducible subshift of $Q^M$, let $Y$ be a subshift of $Q^M$, let $\Delta$ be a $\kappa$-local map from $X$ to $Y$, and let $F$ be a finite subset of $M$ such that $\cardinalityOf{Y_{F^{+2\kappa}}} < \cardinalityOf{X_F}$. There are two patterns $p$ and $p'$ in $X_{F^{+3\kappa}}$ such that $p \neq p'$, $p\restrictedTo_{\boundaryOf_{2\kappa}^+ F^{+\kappa}} = p'\restrictedTo_{\boundaryOf_{2\kappa}^+ F^{+\kappa}}$, and $\Delta_{F^{+3\kappa}}^-(p) = \Delta_{F^{+3\kappa}}^-(p')$. 
  \end{lemma}

  \begin{usage-note}
    In the proof, strong irreducibility of $X$ is used to extend an in $X$ allowed $F$-pattern by an in $X$ allowed $\boundaryOf_{2\kappa}^+ F^{+\kappa}$-pattern and an $\boundaryOf_\kappa^+ F$-pattern to an in $X$ allowed $F^{+3\kappa}$-pattern; and $\kappa$-locality of $\Delta$ is used to restrict it to a map from $X_{F^{+3\kappa}}$ to $Y_{F^{+2\kappa}}$.
  \end{usage-note}

  \begin{proof}
    Because $\cardinalityOf{Y_{F^{+2\kappa}}} < \cardinalityOf{X_F}$, we have $\cardinalityOf{X_F} > 0$, thus $X_F \neq \emptyset$, and hence $X \neq \emptyset$. 
    Therefore, there is a $v \in X_{\boundaryOf_{2\kappa}^+ F^{+\kappa}}$. Let $P_v = \setOf{p \in X_{F^{+3\kappa}} \suchThat p\restrictedTo_{\domainOf(v)} = v}$. Note that, according to \cref{item:repeated-k-boundaries-etc:external-boundary} of \cref{lemma:repeated-k-boundaries-etc}, we have $\boundaryOf_{2\kappa}^+ F^{+\kappa} = F^{+3\kappa} \smallsetminus F^{+\kappa}$.

      Let $u \in X_F$. According to \cref{corollary:distance-of-closure-boundary-of-closure-to-set-greater-than-closure-plus-one}, we have $\distanceOf(F, \domainOf(v)) \geq \kappa + 1$. Hence, because $X$ is $\kappa$-strongly irreducible, there is an $x \in X$ such that $x\restrictedTo_F = u$ and $x\restrictedTo_{\domainOf(v)} = v$. Let $p = x\restrictedTo_{F^{+3\kappa}}$. Then, $p\restrictedTo_F = u$ and $p \in P_v$.

    Therefore, $\cardinalityOf{P_v} \geq \cardinalityOf{X_F} > \cardinalityOf{Y_{F^{+2\kappa}}}$. The restriction $\phi$ of $\Delta_{F^{+3\kappa}}^- \from X_{F^{+3\kappa}} \to Y_{F^{+2\kappa}}$ to $P_v \to \Delta_{F^{+3\kappa}}^-(P_v)$ is surjective. Note that, because $\Delta$ is $\kappa$-local and, according to \cref{item:repeated-k-boundaries-etc:closure-interior} of \cref{lemma:repeated-k-boundaries-etc}, we have $(F^{+3\kappa})^{-\kappa} \supseteq F^{+2\kappa}$, we can choose $Y_{F^{+2\kappa}}$ as the codomain of $\Delta_{F^{+3\kappa}}^-$. If $\phi$ were injective, then $\cardinalityOf{P_v} = \cardinalityOf{\phi(P_v)} \leq \cardinalityOf{Y_{F^{+2\kappa}}}$, which contradicts that $\cardinalityOf{P_v} > \cardinalityOf{Y_{F^{+2\kappa}}}$. Hence, $\phi$ is not injective. In conclusion, there are $p$, $p' \in P_v$ such that $p \neq p'$ and $\phi(p) = \phi(p')$.
  \end{proof}

  Because a local map, that has an image whose entropy is less than the entropy of its domain, maps at least two finite patterns to the same pattern, it is not pre-injective, which is shown in

  \begin{theorem}[Analogue of \cref{theorem:less-entropy-implies-not-pre-injective}] 
  \label{theorem:subshift-less-entropy-implies-not-pre-injective}
    Let $\mathcal{R}$ be right amenable, let $\mathcal{F} = \net{F_i}_{i \in I}$ be a right Følner net in $\mathcal{R}$, let $X$ be a strongly irreducible subshift of $Q^M$ of finite type, let $Y$ be a subshift of $Q^M$, and let $\Delta$ be a local map from $X$ to $Y$ such that $\entropyOf_{\mathcal{F}}(\Delta(X)) < \entropyOf_{\mathcal{F}}(X)$. The map $\Delta$ is not pre-injective.
  \end{theorem}

  \begin{usage-note}
    In the proof, strong irreducibility of $X$ and locality of $\Delta$ are used to apply \cref{lemma:there-are-patterns-that-agree-on-boundary-and-have-the-same-image-under-Delta-restriction-k} yielding two distinct finite patterns with the same domain, identical boundaries, and identical images; and of finite typeness of $X$ is used to apply \cref{corollary:overlapping-patterns-can-be-glued} to identically extend these patterns to points of $X$.
  \end{usage-note}

  \begin{proof-sketch}
    Because the entropy of $\Delta(X)$ is less than the one of $X$, the number of finite patterns in $\Delta(X)$ grows slower than in $X$. Hence, there are two distinct finite patterns in $X$ that have the same image and these can be identically extended to two distinct points of $X$ that have the same image. Therefore, the map $\Delta$ is not pre-injective.
  \end{proof-sketch}

  \begin{proof}
    According to \cref{remark:k-strongly-irreducible-for-greater-k}, \cref{lemma:of-finite-type-implies-step} and \cref{remark:k-step-for-greater-k}, and \cref{remark:k-local-for-greater-k}, there is a $\kappa \in \N_0$ such that $X$ is $\kappa$-strongly irreducible, $X$ is $\kappa$-step, and $\Delta$ is $\kappa$-local.

    Let $Y = \Delta(X)$. According to \cref{remark:local-map-if-and-only-if-cellular-automaton} and \cref{lemma:entropy-invariant-under-closure-net-change} and the precondition $\entropyOf_{\mathcal{F}}(Y) < \entropyOf_{\mathcal{F}}(X)$, we have $\entropyOf_{\net{F_i^{+2\kappa}}_{i \in I}}(Y) \leq \entropyOf_{\mathcal{F}}(Y) < \entropyOf_{\mathcal{F}}(X)$. Hence, there is an $i \in I$ such that
    \begin{equation*}
      \frac{\log\cardinalityOf{Y_{F_i^{+2\kappa}}}}{\cardinalityOf{F_i}} < \frac{\log\cardinalityOf{X_{F_i}}}{\cardinalityOf{F_i}}.
    \end{equation*}
    Thus, $\log\cardinalityOf{Y_{F_i^{+2\kappa}}} < \log\cardinalityOf{X_{F_i}}$ and thus $\cardinalityOf{Y_{F_i^{+2\kappa}}} < \cardinalityOf{X_{F_i}}$.
    Therefore, because $X$ is $\kappa$-strongly irreducible and $\Delta$ is $\kappa$-local, according to \cref{lemma:there-are-patterns-that-agree-on-boundary-and-have-the-same-image-under-Delta-restriction-k}, there are $p$ and $p'$ in $X_{F_i^{+3\kappa}}$ such that $p \neq p'$, $p\restrictedTo_{\boundaryOf_{2\kappa}^+ F_i^{+\kappa}} = p'\restrictedTo_{\boundaryOf_{2\kappa}^+ F_i^{+\kappa}}$, and $\Delta_{F_i^{+3\kappa}}^-(p) = \Delta_{F_i^{+3\kappa}}^-(p')$.

    Hence, because $X$ is $\kappa$-step, according to \cref{corollary:overlapping-patterns-can-be-glued}, there are $x$ and $x'$ in $X$ such that $x\restrictedTo_{\domainOf(p)} = p$, $x'\restrictedTo_{\domainOf(p')} = p'$, and $x\restrictedTo_{M \smallsetminus F_i^{+\kappa}} = x'\restrictedTo_{M \smallsetminus F_i^{+\kappa}}$. In particular, because $p \neq p'$, we have $x \neq x'$ and, because $F_i^{+\kappa}$ is finite, the set $\differenceOf(x, x')$ is finite.

    Moreover, $\Delta(x)\restrictedTo_{F_i^{+2\kappa}} = \Delta_{F_i^{+3\kappa}}^-(p) = \Delta_{F_i^{+3\kappa}}^-(p') = \Delta(x')\restrictedTo_{F_i^{+2\kappa}}$. And, according to \cref{remark:local-map-if-and-only-if-cellular-automaton} and \cref{lemma:global-transition-function-and-interior-closure}, we have $\Delta(x)\restrictedTo_{M \smallsetminus F_i^{+2\kappa}} = \Delta(x')\restrictedTo_{M \smallsetminus F_i^{+2\kappa}}$. Therefore, $\Delta(x) = \Delta(x')$. In conclusion, $\Delta$ is not pre-injective.
  \end{proof}

  If in a point of a shift space we replace all occurrences of a pattern by another pattern with the same image that agree on a big enough boundary, we get a new point of the shift space in which the first pattern does not occur that has the same image as the original point, which is shown in

  \begin{lemma}[Analogue of \cref{lemma:exchanging-pattern-by-other-pattern-with-same-image-yields-configuration-with-same-image}] 
  \label{lemma:exchanging-pattern-by-other-pattern-with-same-image-yields-configuration-with-same-image-and-we-stay-in-subshift}
    Let $X$ be a $\kappa$-step subshift of $Q^M$, let $Y$ be a subshift of $Q^M$, let $\Delta$ be a $\kappa$-local map from $X$ to $Y$, let $A$ be a subset of $M$, let $p$ and $p'$ be two patterns in $X_{A^{+2\kappa}}$ such that $p\restrictedTo_{\boundaryOf_{2\kappa}^+ A} = p'\restrictedTo_{\boundaryOf_{2\kappa}^+ A}$ and $\Delta_{A^{+2\kappa}}^-(p) = \Delta_{A^{+2\kappa}}^-(p')$. Furthermore, let $c$ be a point of $X$ and let $T$ be a subset of $M$ such that the family $\family{t \isSemiActedUponBy A^{+2\kappa}}_{t \in T}$ is pairwise disjoint and that, for each element $t \in T$, we have $p \occursIn_t c$. Put 
    \begin{equation*} 
      c' = c\restrictedTo_{M \smallsetminus (\bigcup_{t \in T} t \isSemiActedUponBy A^{+2\kappa})} \times \coprod_{t \in T} t \actsByItsCoordinateOn p'.
    \end{equation*}
    Then, for each element $t \in T$, we have $p' \occursIn_t c'$, and $c' \in X$, and $\Delta(c) = \Delta(c')$. In particular, if $p \neq p'$, then, for each element $t \in T$, we have $p \not\occursIn_t c'$.
  \end{lemma}

  \begin{usage-note}
    In the proof, $\kappa$-stepness of $X$ is used to apply lemma \ref{lemma:overlapping-global-configurations-can-be-glued} to deduce that $c'$ is a point of $X$; and locality of $\Delta$ is used to deduce that $\Delta(c) = \Delta(c')$. 
  \end{usage-note}

  \begin{proof}
    There are $x$ and $x'$ in $X$ such that $x\restrictedTo_{A^{+2\kappa}} = p$ and $x'\restrictedTo_{A^{+2\kappa}} = p'$. Thus, for each $t \in T$, we have $(t \actsByItsCoordinateOn x')\restrictedTo_{t \isSemiActedUponBy A^{+2\kappa}} = t \actsByItsCoordinateOn p'$. Hence,
    \begin{align*}
      c' = c\restrictedTo_{M \smallsetminus (\bigcup_{t \in T} t \isSemiActedUponBy A^{+2\kappa})} \times \coprod_{t \in T} (t \actsByItsCoordinateOn x')\restrictedTo_{t \isSemiActedUponBy A^{+2\kappa}}
    \end{align*}
    Moreover, for each $t \in T$, according to \cref{remark:shift-invariance-induces-invariance-under-right-semi-action}, we have $t \actsByItsCoordinateOn x' \in X$. And, by precondition, $\family{(t \isSemiActedUponBy A)^{+2\kappa}}_{t \in T}$ is pairwise disjoint (where we used that $t \isSemiActedUponBy A^{+2\kappa} = (t \isSemiActedUponBy A)^{+2\kappa}$, which holds according to \cref{item:properties-of-interior-closure-and-boundary:commute-with-liberation} of \cref{lemma:properties-of-interior-closure-and-boundary}). And, for each $t \in T$, we have $(t \actsByItsCoordinateOn x')\restrictedTo_{\boundaryOf_{2\kappa}^+(t \isSemiActedUponBy A)} = (t \actsByItsCoordinateOn p')\restrictedTo_{\boundaryOf_{2\kappa}^+(t \isSemiActedUponBy A)} = (t \actsByItsCoordinateOn p)\restrictedTo_{\boundaryOf_{2\kappa}^+(t \isSemiActedUponBy A)} = c\restrictedTo_{\boundaryOf_{2\kappa}^+(t \isSemiActedUponBy A)}$ (where we used that $\boundaryOf_{2\kappa}^+(t \isSemiActedUponBy A) = t \isSemiActedUponBy \boundaryOf_{2\kappa}^+ A$, which holds according to \cref{item:properties-of-interior-closure-and-boundary:commute-with-liberation} of \cref{lemma:properties-of-interior-closure-and-boundary}). Therefore, because $X$ is $\kappa$-step, according to \cref{lemma:overlapping-global-configurations-can-be-glued}, we have $c' \in X$.

    Let $m \in M$.
    \begin{description}
      \item[Case 1:] $\Exists t \in T \SuchThat m \in t \isSemiActedUponBy A^{+\kappa}$. Then, according to \cref{item:characterisation-of-k-closure-and-interior:closure} of lemma \ref{lemma:characterisation-of-k-closure-and-interior} and \cref{item:repeated-k-boundaries-etc:closure} of \cref{lemma:repeated-k-boundaries-etc} 
                     \begin{align*}
                       m \isSemiActedUponBy \ball(\kappa)
                       &\subseteq (t \isSemiActedUponBy A^{+\kappa}) \isSemiActedUponBy \ball(\kappa)\\
                       &=         (t \isSemiActedUponBy A^{+\kappa})^{+\kappa}\\
                       &=          t \isSemiActedUponBy A^{+2\kappa}.
                     \end{align*}
                     Hence, because $\Delta$ is $\kappa$-local, according to \cref{remark:local-map-if-and-only-if-cellular-automaton} and \cref{corollary:equivariance-of-restriction},
                     \begin{align*}
                       \Delta(c')(m)
                       &= \Delta_{t \isSemiActedUponBy A^{+2\kappa}}^-(t \actsByItsCoordinateOn p')\\
                       &= t \actsByItsCoordinateOn \Delta_{A^{+2\kappa}}^-(p')\\
                       &= t \actsByItsCoordinateOn \Delta_{A^{+2\kappa}}^-(p)\\
                       &= \Delta_{t \isSemiActedUponBy A^{+2\kappa}}^-(t \actsByItsCoordinateOn p)\\
                       &= \Delta(c)(m).
                     \end{align*}
      \item[Case 2:] $\ForEach t \in T \Holds m \notin t \isSemiActedUponBy A^{+\kappa}$. Then, $m \in M \smallsetminus \bigcup_{t \in T} t \isSemiActedUponBy A^{+\kappa}$. Hence, according to \cref{item:characterisation-of-k-closure-and-interior:closure} of \cref{lemma:characterisation-of-k-closure-and-interior}, \cref{item:properties-of-interior-closure-and-boundary:complement} of \cref{lemma:properties-of-interior-closure-and-boundary}, \cref{item:properties-of-interior-closure-and-boundary:union} of \cref{lemma:properties-of-interior-closure-and-boundary}, and \cref{item:repeated-k-boundaries-etc:closure-interior} of \cref{lemma:repeated-k-boundaries-etc},
                     \begin{align*}
                       m \isSemiActedUponBy \ball(\kappa)
                       &\subseteq \parens[\big]{M \smallsetminus \bigcup_{t \in T} t \isSemiActedUponBy A^{+\kappa}} \isSemiActedUponBy \ball(\kappa)\\
                       &=         \parens[\big]{M \smallsetminus \bigcup_{t \in T} t \isSemiActedUponBy A^{+\kappa}}^{+\kappa}\\
                       &=          M \smallsetminus \parens[\big]{\bigcup_{t \in T} t \isSemiActedUponBy A^{+\kappa}}^{-\kappa}\\
                       &\subseteq  M \smallsetminus \bigcup_{t \in T} (t \isSemiActedUponBy A^{+\kappa})^{-\kappa}\\
                       &\subseteq  M \smallsetminus \bigcup_{t \in T} t \isSemiActedUponBy A^{+ 0}\\
                       &=          M \smallsetminus \bigcup_{t \in T} t \isSemiActedUponBy A.
                     \end{align*}
                     Therefore, because $\Delta$ is $\kappa$-local and $c'\restrictedTo_{M \smallsetminus \bigcup_{t \in T} t \isSemiActedUponBy A} = c\restrictedTo_{M \smallsetminus \bigcup_{t \in T} t \isSemiActedUponBy A}$, we have $\Delta(c')(m) = \Delta(c)(m)$.
    \end{description}
    In either case, $\Delta(c')(m) = \Delta(c)(m)$. Therefore, $\Delta(c') = \Delta(c)$.
  \end{proof}

  Because a local map that is not pre-injective maps at least two finite patterns to the same pattern, the entropy of its image is less than the entropy of its domain, which is shown in

  \begin{theorem}[Analogue of \cref{theorem:not-pre-injective-implies-less-entropy}] 
  \label{theorem:subshift-not-pre-injective-implies-less-entropy}
    Let $M$ be infinite, let $X$ be a strongly irreducible subshift of $Q^M$ of finite type, let $Y$ be a subshift of $Q^M$, let $\Delta$ be a local map from $X$ to $Y$ that is not pre-injective, and let $\mathcal{F}$ be a right Erling net in $\mathcal{R}$ (which exists according to \cref{lemma:finitely-right-generated-is-tractable}). Then, $\entropyOf_{\mathcal{F}}(\Delta(X)) < \entropyOf_{\mathcal{F}}(X)$.
  \end{theorem}

  \begin{usage-note}
    In the proof, infiniteness of $M$ is used to apply \cref{theorem:there-is-a-rho-kappa-theta-tiling} yielding a tiling, strong irreducibility of $X$ is used to apply \cref{theorem:entropy-less-if-real-subsets-for-rho-kappa-theta-tiling} yielding a strict inequality for entropies, finite typeness of $X$ and locality of $\Delta$ is used to apply \cref{lemma:exchanging-pattern-by-other-pattern-with-same-image-yields-configuration-with-same-image-and-we-stay-in-subshift} yielding that the image of all points of $X$ in which a certain pattern does not occur at points of a tiling is the same as the image of $X$.
  \end{usage-note}

  \begin{proof-sketch}
    Because $\Delta$ is not pre-injective, there are two distinct points of $X$ with the same image that differ only in finitely many cells. Thus, there are two distinct finite patterns, say $p$ and $p'$, with the same image. Hence, the image of $X$ is equal to the image of the set $Z$ of all points of $X$ in which the pattern $p$ does not occur. Because there are too many points not in $Z$, this set does have less entropy than $X$.
  \end{proof-sketch}

  \begin{proof}
    According to \cref{remark:k-strongly-irreducible-for-greater-k}, \cref{lemma:of-finite-type-implies-step} and \cref{remark:k-step-for-greater-k}, and \cref{remark:k-local-for-greater-k}, there is a $\kappa \in \N_0$ such that $X$ is $\kappa$-strongly irreducible, $X$ is $\kappa$-step, and $\Delta$ is $\kappa$-local.

    Because $\Delta$ is not pre-injective, there are $c$ and $c'$ in $X$ such that $\differenceOf(c, c')$ is finite, $\Delta(c) = \Delta(c')$, and $c \neq c'$; in particular, $\cardinalityOf{X} \geq 2$. Hence, there is a $\rho \in \N_0$ such that $\differenceOf(c, c') \subseteq \ball(\rho)$. Let $p = c\restrictedTo_{\ball(\rho)^{+2\kappa}}$ and let $p' = c'\restrictedTo_{\ball(\rho)^{+2\kappa}}$. Then, $p \neq p'$; $m_0 \in \ball(\rho)$; $p$, $p' \in X_{\ball(\rho)^{+2\kappa}}$; $p\restrictedTo_{\boundaryOf_{2\kappa}^+ \ball(\rho)} = p'\restrictedTo_{\boundaryOf_{2\kappa}^+ \ball(\rho)}$; and, because $\Delta(c) = \Delta(c')$, we have $\Delta_{\ball(\rho)^{+2\kappa}}^-(p) = \Delta_{\ball(\rho)^{+2\kappa}}^-(p')$.

    Because $M$ is infinite, according to \cref{theorem:there-is-a-rho-kappa-theta-tiling}, there is a $\ntuple{\rho + 2 \kappa, \kappa, \theta'}$-tiling $T$ of $\mathcal{R}$. Let
    \begin{equation*}
      Z = \setOf{x \in X \suchThat \ForEach t \in T \Holds p \not\occursIn_t x}.
    \end{equation*}
    For each $t \in T$, according to \cref{remark:pattern-belongs-to-shift-if-and-only-if-translated-pattern-belongs}, we have $t \actsByItsCoordinateOn p \in X_{t \isSemiActedUponBy \ball(\rho)^{+2\kappa}} \smallsetminus Z_{t \isSemiActedUponBy \ball(\rho)^{+2\kappa}}$ and hence $Z_{t \isSemiActedUponBy \ball(\rho)^{+2\kappa}} \subsetneqq X_{t \isSemiActedUponBy \ball(\rho)^{+2\kappa}}$. Moreover, for each $t \in T$, according to \cref{item:characterisation-of-k-closure-and-interior-of-balls:closure} of \cref{corollary:characterisation-of-k-closure-and-interior-of-balls} and \cref{corollary:ball-centred-at-mzero-to-m}, we have $t \isSemiActedUponBy \ball(\rho)^{+2\kappa} = \ball(t, \rho + 2 \kappa)$.
    Therefore, because $X$ is $\kappa$-strongly irreducible and $\cardinalityOf{X} \geq 2$, according to \cref{theorem:entropy-less-if-real-subsets-for-rho-kappa-theta-tiling}, we have $\entropyOf_{\mathcal{F}}(Z) < \entropyOf_{\mathcal{F}}(X)$. Hence, according to \cref{theorem:entropy-does-non-increase}, we have $\entropyOf_{\mathcal{F}}(\Delta(Z)) < \entropyOf_{\mathcal{F}}(X)$.

    Let $x \in X$. Put $U = \setOf{t \in T \suchThat p \occursIn_t x}$. Because $X$ is $\kappa$-step and $\Delta$ is $\kappa$-local, according to \cref{lemma:exchanging-pattern-by-other-pattern-with-same-image-yields-configuration-with-same-image-and-we-stay-in-subshift}, there is an $x' \in X$ such that $x' \in Z$ and $\Delta(x) = \Delta(x')$. Therefore, $\Delta(X) = \Delta(Z)$. In conclusion, $\entropyOf_{\mathcal{F}}(\Delta(X)) < \entropyOf_{\mathcal{F}}(X)$. 
  \end{proof}

  Because a right Følner net in a finite cell space is eventually equal to the set of cells, the entropy of a subset of the full shift is a function of the cardinality of that set, which is shown in

  \begin{lemma} 
  \label{lemma:if-M-is-finite-abs-X-determined-by-entropy}
    Let $\mathcal{R}$ be right amenable, let $\mathcal{F} = \net{F_i}_{i \in I}$ be a right Følner net in $\mathcal{R}$, let $M$ be finite, and let $X$ be a subset of $Q^M$. Then,
    \begin{equation*}
      \cardinalityOf{X} = \exp\parens[\big]{\cardinalityOf{M} \cdot \entropyOf_{\mathcal{F}}(X)}. \qedhere
    \end{equation*}
  \end{lemma}

  \begin{proof} 
    Let $F$ be a non-empty and finite subset of $M$ such that $F \neq M$. Then, because $\isSemiActedUponBy$ is transitive, there is a $\mathfrak{g} \in G \modulo G_0$ such that $(F \isSemiActedUponBy \mathfrak{g}) \cap (M \smallsetminus F) \neq \emptyset$. Hence, $F \isSemiActedUponBy \mathfrak{g} \nsubseteq F$, thus $F \nsubseteq (\blank \isSemiActedUponBy \mathfrak{g})^{-1}(F)$, thus $F \smallsetminus (\blank \isSemiActedUponBy \mathfrak{g})^{-1}(F) \neq \emptyset$, and therefore $\cardinalityOf{F \smallsetminus (\blank \isSemiActedUponBy \mathfrak{g})^{-1}(F)} \neq 0$. On the other hand, $\cardinalityOf{M \smallsetminus (\blank \isSemiActedUponBy \mathfrak{g})^{-1}(M)} = 0$. Moreover, because $M$ is finite, the set $\setOf{F \subseteq M \suchThat F \neq \emptyset, F \text{ finite}}$ is finite and hence its subset $\setOf{F_i \suchThat i \in I}$ is finite too. Furthermore, because $\mathcal{F}$ is a right Følner net,
    \begin{equation*}
      \ForEach \mathfrak{g} \in G \modulo G_0 \Holds \lim_{i \in I} \frac{\cardinalityOf{F_i \smallsetminus (\blank \isSemiActedUponBy \mathfrak{g})^{-1}(F_i)}}{\cardinalityOf{F_i}} = 0.
    \end{equation*}
    Altogether, $\mathcal{F}$ is eventually equal to $M$. Therefore,
    \begin{equation*}
      \entropyOf_{\mathcal{F}}(X)
      = \frac{\log\cardinalityOf{\pi_M(X)}}{\cardinalityOf{M}}
      = \frac{\log\cardinalityOf{X}}{\cardinalityOf{M}}.
    \end{equation*}
    In conclusion, $\cardinalityOf{X} = \exp(\cardinalityOf{M} \cdot \entropyOf_{\mathcal{F}}(X))$.
  \end{proof}

  Because surjectivity as well as pre-injectivity of a local map is characterised by maximal entropy of its image with respect to its codomain or domain, if both domains have the same entropy, then a local map is surjective if and only if it is pre-injective, which is shown in

  \begin{main-theorem}[Garden of Eden Theorem; Edward Forrest Moore, 1962; John R. Myhill, 1963; Analogue of \cref{theorem:garden-of-Eden}] 
  \label{theorem:subshift-garden-of-Eden}
    Let $\mathcal{R}$ be right amenable, let $\mathcal{F}$ be a right Følner net in $\mathcal{R}$, let $X$ be a non-empty strongly irreducible subshift of $Q^M$ of finite type, let $Y$ be a strongly irreducible subshift of $Q^M$ such that $\entropyOf_{\mathcal{F}}(X) = \entropyOf_{\mathcal{F}}(Y)$, and let $\Delta$ be a local map from $X$ to $Y$. The map $\Delta$ is surjective if and only if it is pre-injective.
  \end{main-theorem}

  \begin{usage-note}
    In the proof, non-emptiness of $X$, strong irreducibility of $Y$, and locality of $\Delta$ are used to apply \cref{theorem:subshift-not-surjective-implies-less-entropy} yielding a characterisation of surjectivity; and strong irreducibility and finite typeness of $X$, and locality of $\Delta$ are used to apply \cref{theorem:subshift-less-entropy-implies-not-pre-injective,theorem:subshift-not-pre-injective-implies-less-entropy} yielding a characterisation of pre-injectivity.
  \end{usage-note}

  \begin{proof}
    First, let $M$ be finite. Then, $Q^M$ is finite, and thus $X$ and $Y$ are finite. Hence, because $\entropyOf_{\mathcal{F}}(X) = \entropyOf_{\mathcal{F}}(Y)$, according to \cref{lemma:if-M-is-finite-abs-X-determined-by-entropy}, we have $\cardinalityOf{X} = \cardinalityOf{Y}$. Therefore, $\Delta$ is surjective if and only if it is injective. Moreover, because $M$ is finite, the map $\Delta$ is pre-injective if and only if it is injective. In conclusion, $\Delta$ is surjective if and only if it is pre-injective.

    Secondly, let $M$ be infinite. According to \cref{theorem:subshift-not-surjective-implies-less-entropy}, the map $\Delta$ is not surjective if and only if $\entropyOf_{\mathcal{F}}(\Delta(X)) < \entropyOf_{\mathcal{F}}(Y)$. And, according to \cref{theorem:subshift-less-entropy-implies-not-pre-injective} and \cref{theorem:subshift-not-pre-injective-implies-less-entropy}, because $\entropyOf_{\mathcal{F}}(X) = \entropyOf_{\mathcal{F}}(Y)$, we have $\entropyOf_{\mathcal{F}}(\Delta(X)) < \entropyOf_{\mathcal{F}}(Y)$ if and only if $\Delta$ is not pre-injective. Hence, $\Delta$ is not surjective if and only if it is not pre-injective. In conclusion, $\Delta$ is surjective if and only if it is pre-injective.
  \end{proof}

  \begin{corollary} 
  \label{corollary:Moore-and-Myhill} 
    Let $\mathcal{R}$ be right amenable. Each strongly irreducible subshift of $Q^M$ of finite type has the Moore and the Myhill property.
  \end{corollary}

  \begin{proof}
    This is a direct consequence of \cref{theorem:subshift-garden-of-Eden} and the fact that the empty strongly irreducible subshift of $Q^M$ has the Moore and the Myhill property.
  \end{proof}

  \begin{corollary}
  \label{corollary:Moore-and-Myhill-for-left-homogeneous-spaces}
    Let $\mathcal{M} = \ntuple{M, G, \actsOnPoint}$ be a right-a\-me\-na\-ble and finitely right-gen\-er\-at\-ed left-ho\-mo\-ge\-neous space with finite stabilisers and let $Q$ be a finite set. For each coordinate system $\mathcal{K}$ for $\mathcal{M}$ and each $\mathcal{K}$-big subgroup $H$ of $G$, each strongly irreducible subshift of $Q^M$ of finite type with respect to $\mathcal{R} = \ntuple{\mathcal{M}, \mathcal{K}}$ and $H$ has the Moore and the Myhill property with respect to $\mathcal{R}$ and $H$.
  \end{corollary}

  \begin{proof}
    This is a direct consequence of \cref{theorem:subshift-garden-of-Eden}.
  \end{proof}

  \begin{remark}
    Note that in \cref{corollary:Moore-and-Myhill-for-left-homogeneous-spaces} we do not have to choose a finite and symmetric right-gen\-er\-at\-ing set $S$, because being a subshift, being strongly irreducible, being of finite type, being local, being surjective, being pre-injective, having the Moore property, and having the Myhill property, do not depend on the choice of a finite and symmetric right-gen\-er\-at\-ing set of $\mathcal{R}$; the reason for the properties that depend on the metric induced by such a right-gen\-er\-at\-ing set is that those metrics are, according to \cref{corollary:Lipschitz-equivalent-metrics}, pairwise Lipschitz equivalent.
  \end{remark}

  \begin{example}[From Golden Mean to Even Shift]
    In the situation of \cref{example:from-golden-mean-to-even-shift-local-map}, let $\mathcal{F}$ be the right Følner net $\sequence{\setOf{1, 2, \dotsc, n}}_{n \in \N_+}$. Recall that the golden mean shift $X$ is non-empty (\cref{example:shift:golden-mean}), strongly irreducible (\cref{example:strongly-irreducible}), and of finite type (\cref{example:of-finite-type-or-not}); and that the even shift $Y$ is strongly irreducible (\cref{example:strongly-irreducible}) but \emph{not} of finite type (\cref{example:of-finite-type-or-not}). And, according to example~4.1.4 in \cite{lind:marcus:1995} and Example~4.1.6 in \cite{lind:marcus:1995}, the entropy of $X$ with respect to $\mathcal{F}$ and the one of $Y$ are both the golden mean $(1 + \sqrt{5})/2$. Therefore, according to \cref{theorem:subshift-garden-of-Eden}, because the local map $\Delta$ from $X$ to $Y$ is surjective (\cref{example:from-golden-mean-to-even-shift-local-map}), it is also pre-injective. However, it is not injective, because the two points of $X$ with alternating $0$'s and $1$'s, that is, those of the form $\dotso 010101 \dotso$, are both mapped to the point of $Y$ with only $0$'s, that is, the one of the form $\dotso 000 \dotso$. 
  \end{example}

  \clearToOddPage
  \chapter{A Quasi-Solution of the Firing Mob Synchronisation Problem} 
  \label{chapter:fssp}

  \paragraph{Abstract.} We construct a time-optimal quasi-solution of the firing mob synchronisation problem over finite, connected, and undirected multigraphs whose maximum degrees are uniformly bounded by a constant. It is only a quasi-solution because its number of states depends on the graph or, from another perspective, does not depend on the graph but is countably infinite. To construct this quasi-solution we introduce signal machines over continuum representations of such multigraphs and construct a signal machine whose discretisation is a cellular automaton that quasi-solves the problem. This automaton uses a time-optimal solution of the firing squad synchronisation problem in dimension one with one general at one end to synchronise edges, and freezes and thaws the synchronisation of edges in such a way that all edges synchronise at the same time.

  \paragraph{Introduction.} The firing squad synchronisation problem in dimension one with one general at one end is to synchronise each finite one-dimensional array of cells starting from one end of the array and the cell at this end is called \emph{general}. It was proposed by John R. Myhill in 1957, solved by John McCarthy and Marvin Lee Minsky, and published by Edward Forrest Moore in 1962 (see \cite{moore:1964}). The first time-optimal several-thousand-states solution was found by Eiichi Goto in 1962 (see \cite{goto:1962}), reduced to $16$ states by Abraham Waksman in 1966 (see \cite{waksman:1966}), and reduced to $8$ states by Robert Balzer in 1967 (see \cite{balzer:1967}). Hein D. Gerken found another time-optimal $7$-states solution in 1987 (see \cite{gerken:1987}) and Jacques Mazoyer found a time-optimal $6$-states solution also in 1987 (see \cite{mazoyer:1987}). It is unknown whether there is a time-optimal $5$-states solution but it is known that there is no time-optimal $4$-states solution, a result due to Robert Balzer and Peter Sanders (see \cite{balzer:1967,sanders:1994}). 

  The firing mob synchronisation problem is to synchronise each finite, connected, and undirected graph whose maximum degree is bounded by a fixed constant starting from any vertex and this vertex is called \emph{general}. It was solved by P. Rosenstiehl, J.\;R. Fiksel, and A. Holliger in 1972 (see \cite{rosenstiehl:1972}) and also by Francesco Romani in 1976 (see \cite{romani:1976}), where the latter solution achieves better running times than the former. The problem for specific classes of graphs were for example studied by Kojiro Kobayashi in 1977 and 1978 (see \cite{kobayashi:1977,kobayashi2:1978,kobayashi:1978}) and by Zsuzsanna Róka in 2000 (see \cite{roka:2000}). Karel Culik II and Simant Dube presented a solution of the general case in 1991 (see \cite{culik-dube:1991}). It needs $3.5 r$-many steps, where $r$ is the maximal distance of the general to a vertex and is called \emph{radius of the graph with respect to the general}. By using more and more states, the solution can be adjusted such that the number of steps it needs approaches $3 r$.

  It was shown that $r + d$ is a lower bound for the number of steps that solutions of the firing mob synchronisation problem need by John J. Grefenstette in 1983 (see \cite{grefenstette:1983}), where $d$ is the maximal distance between two vertices of the graph and is called \emph{diameter of the graph}. Because there are graphs and choices of generals such that the diameter is $2 r$, the solutions by Karel Culik II and Simant Dube approach the optimal number of steps, namely $3 r$, if $r$ is taken as problem size. However, if $r + d$ is taken as problem size, then their solutions do not approach the optimal number of steps.

  In the present chapter we construct a time-optimal quasi-solution that needs exactly $r + d$ steps but whose number of states depends on the graph or, from another perspective, does not depend on the graph but is countably infinite (this is why we call it a quasi-solution). It can also be turned into a time-optimal quasi-solution of the firing squad synchronisation problem for any region in any dimension with one general at any position by regarding each region to be synchronised as a graph, where cells in the region are vertices and edges are neighbourhood relationships. 

  However, restricted to specific classes of problems, the quasi-solution may not be time-optimal. For example, restricted to rectangular regions with one general at one corner, the quasi-solution needs $2 (k + \ell - 2)$-many steps whereas $(k + \ell + \max\setOf{k, \ell} - 3)$-many steps is optimal (see for example \cite{umeo:2009}), where $k$ and $\ell$ are the side lengths of the rectangle. Nevertheless, because the quasi-solution is (trivially) embeddable in the sense of \cite{grefenstette:1983}, according to theorem~1 in \cite{grefenstette:1983}, it can be combined with finitely many embeddable time-optimal solutions for specific classes of problems to get one quasi-solution that is also time-optimal for those classes. Examples of solutions for specific classes, like rectangular regions with one general at the upper left corner, are given in sections~5 and~6 in \cite{grefenstette:1983}.

  To design, explain, and draw solutions of firing squad/mob synchronisation problems, it is convenient to think about, talk about, and draw continuous space-time diagrams of different kinds of signals that move across the cell space, vanish or give rise to new signals upon reaching boundaries or junctions of the space or upon colliding with each other. This is mostly done in an informal way, but the idea of signals has also been formalised for one-dimensional cellular automata by Jérôme Olivier Durand-Lose in 2005 (see \cite{durand-lose:2005}).

  This formalisation however does not handle accumulations of events like collisions and does not allow infinitely many signals of different speeds, which naturally occur and are necessary in descriptions of many solutions of the firing squad synchronisation problem by signals. For example, collisions accumulate at the time synchronisation finishes and infinitely many signals of different speeds may originate from the general. In the time evolutions of the actual cellular automata, the accumulations of collisions disappear due to the discreteness of space and time, and the infinitely many signals are cleverly produced by finitely many states (see for example \cite{mazoyer:1987}).

  Because we want to describe our quasi-solution in terms of signals in a formal way, we first introduce continuum representations of finite and connected multigraphs (without self-loops), we secondly introduce signal machines over such representations that allow infinitely many signals of different speeds and seamlessly handle accumulations of events and accumulations of accumulations of events and so forth, and we thirdly construct a signal machine for the continuous firing mob synchronisation problem over such representations and shortly note how to discretise it to get a cellular automaton quasi-solution of the firing mob synchronisation problem.

  \paragraph{Contents.} In \cref{section:firing-squad-problem} we state the firing squad and the firing mob synchronisation problems. In \cref{section:undirected-multigraphs} we introduce undirected multigraphs (without self-loops) and direction-preserving paths in such graphs, which are paths that do not make U-turns. In \cref{section:continuum-representation} we introduce continuum representations of undirected multigraphs, which are in a sense drawings of graphs in a high-dimensional Euclidean space. In \cref{section:signal-machines} we introduce signal machines, which can be studied in their own right, but which can also be thought of as high-level views of time evolutions of cellular automata over graphs, like cellular automata over finitely right-gen\-er\-at\-ed cell spaces, that are restricted to configurations with a fixed finite support. In \cref{section:firing-squad-solution} we construct a signal machine whose discretisation is a cellular automaton that quasi-solves the firing mob synchronisation problem in $(r + d)$-many steps. And in \cref{section:proof-sketch} we sketch a proof for that statement. The impatient may right now have a look at the continuous space-time diagrams of the synchronisations of small trees as performed by the quasi-solution: See \cref{figure:MazoyerWithOneAndTwoEdges,figure:fsspWithTwoEdgesAndGeneralInBetweenAndAtTheLeft,figure:fsspWithThreeEdgesInARowAndGeneralAtTheSecondVertexFromTheLeft,figure:fsspWithThreeEdgesIncidentToTheSameVertexAndGeneralAtTheVertexAndAtTheLeafOfTheShortestEdge} on \cpageref{figure:MazoyerWithOneAndTwoEdges,figure:fsspWithTwoEdgesAndGeneralInBetweenAndAtTheLeft,figure:fsspWithThreeEdgesInARowAndGeneralAtTheSecondVertexFromTheLeft,figure:fsspWithThreeEdgesIncidentToTheSameVertexAndGeneralAtTheVertexAndAtTheLeafOfTheShortestEdge}.

  \paragraph{Preliminary Notions.} The affinely extended real numbers\graffito{affinely extended real numbers $\overline{\R}$} $\R \cup \setOf{-\infty, +\infty}$ are denoted by $\overline{\R}$\index[symbols]{Roverline@$\overline{\R}$}. For each tuple $(r, r') \in \overline{\R} \times \overline{\R}$ such that $r \leq r'$, the closed, open, and the two half-open extended real intervals\graffito{extended real intervals $\closedInterval{r, r'}$, $\openInterval{r, r'}$, and $[r, r'[$ and $]r, r']$} with the endpoints $r$ and $r'$ are denoted by $\closedInterval{r, r'}$, $\openInterval{r, r'}$, and $[r, r'[$ and $]r, r']$ respectively. And, for each tuple $(z, z') \in \Z \times \Z$ such that $z \leq z'$, the closed integer-valued interval\graffito{closed integer-valued interval $\discreteInterval{z}{z'}$} with the endpoints $z$ and $z'$, namely $\setOf{z, z + 1, \dotsc, z'}$ or equivalently $\closedInterval{z, z'} \cap \Z$, is denoted by $\discreteInterval{z}{z'}$. 

  \section{The Firing Squad/Mob Synchronisation Problems}
  \label{section:firing-squad-problem}

  In this section, let $\mathcal{R} = \ntuple{\ntuple{M, G, \actsOnPoint}, \ntuple{m_0, \family{g_{m_0, m}}_{m \in M}}}$ be a finitely right-gen\-er\-at\-ed cell space, let $N$ be a finite right-gen\-er\-at\-ing set of $\mathcal{R}$ that contains $G_0$, where $G_0$ is the stabiliser of $m_0$ under $\actsOnPoint$, let $\mathcal{G} = \ntuple{M, E}$ be the coloured $N$-Cayley graph of $\mathcal{R}$, let $\mathcal{C}$ be a semi-cellular automaton over $\mathcal{R}$ with state set $Q$, neighbourhood $N$, and local transition function $\delta$, and let $\Delta$ be the global transition function of $\mathcal{C}$. 

  To state the problems succinctly we introduce the notions of passive subsets of states, dead states, supports of global configurations with respect to a distinguished dead state, and what it means for a global configuration to be of the form of a pattern in the following four definitions.

  \begin{definition}
    Let $P$ be a subset of $Q$. It is called \define{passive}\graffito{passive set of states} if and only if, for each local configuration $\ell \in Q^N$ with $\imageOf(\ell) \subseteq P$, we have $\delta(\ell) = \ell(G_0)$.
  \end{definition}

  \begin{definition}
    Let $q$ be a state of $Q$. It is called \define{dead}\graffito{dead state} if and only if, for each local configuration $\ell \in Q^N$ with $\ell(G_0) = q$, we have $\delta(\ell) = q$.
  \end{definition}

  In the remainder of this section, let $Q$ contain a distinguished dead state named $\deadState$.

  \begin{definition}
    Let $c$ be a global configuration of $Q^M$. The set $\supportOf(c) = M \smallsetminus c^{-1}(\deadState)$ is called \define{support of $c$}\graffito{support $\supportOf(c)$ of $c$}\index[symbols]{suppc@$\supportOf(c)$}.
  \end{definition}

  \begin{definition} 
    Let $A$ be a subset of $M$, let $p$ be a pattern of $Q^A$, and let $c$ be a global configuration of $Q^M$. The global configuration $c$ is said to \define{be of the form $p$}\graffito{global configuration is of the form $p$} if and only if there is an element $g \in G$ such that $c\restrictedTo_{g \actsOnPoint A} = g \actsOnMap p$ and $c\restrictedTo_{M \smallsetminus (g \actsOnPoint A)} \equiv \deadState$. 
  \end{definition}


  We state the firing squad synchronisation problem in

  \begin{definition} 
    Let $\deadState$, $\generalState$, $\soldierState$, and $\fireState$ be four distinct states, and let $Q'$ be the set that consists of those states. A solution of the \graffito{firing squad synchronisation problem in dimension one with one general at the left end}\define{firing squad synchronisation problem in dimension one with one general at the left end} is a cellular automaton $\mathcal{C}$ over $\ntuple{\ntuple{\Z, \Z, +}, \ntuple{0, \family{z}_{z \in \Z}}}$ with neighbourhood $\setOf{-1, 0, 1}$ and finite set of states that includes $Q'$ such that the state $\deadState$ is dead and the set $\setOf{\deadState, \soldierState}$ is passive, and whose global transition function $\Delta$ has the following property:

    For each global configuration $c$ with finite support of the form $\generalState \soldierState \soldierState \dotsb \soldierState$, there is a non-negative integer $k$ such that the global configuration $\Delta^k(c)$ is of the form $\fireState \fireState \dotsb \fireState$ and has the same support as $c$, and such that the state $\fireState$ does not occur in any of the global configurations $\Delta^j(c)$, for $j \in \N_0$ with $j < k$.
  \end{definition}

  \begin{remark}
    Let $\mathcal{C}$ be a solution of the above problem, let $c$ be a global configurations of the form $\generalState \soldierState \soldierState \dotsb \soldierState$, and let $k$ be the non-negative integer from the problem definition. Then, because the state $\deadState$ is dead and the support of $\Delta^k(c)$ is the same as the one of $c$, for each non-negative integer $j$ with $j \leq k$, the support of $\Delta^j(c)$ is the same as the one of $c$. Broadly speaking, in the time evolution of solutions, the support of initial configurations can neither shrink nor grow before synchronisation is finished. Moreover, because the set $\setOf{\deadState, \soldierState}$ is passive, if the support of $c$ consists of at least $3$ cells, then $\Delta(c)$ cannot be of the form $\fireState \fireState \dotsb \fireState$. Broadly speaking, the problem cannot be solved trivially.
  \end{remark}

  \begin{remark}
  \label{rem:one-dimensional-array}
    As mentioned above, for each global configuration $c$ of the form $\generalState \soldierState \soldierState \dotsb \soldierState$, the supports of the global configurations that are observable in the time evolutions that begin in the configuration $c$ of cellular automata that solve the above problem, are included in the support of $c$. Hence, we can regard such cellular automata as automata over one-dimensional arrays with one dummy neighbour in the state $\deadState$ at each end.
  \end{remark}

  \begin{remark}
    The above problem can be generalised in many ways. For example, by allowing the general to be placed anywhere or by allowing more than one general.
  \end{remark}

  We state the firing mob synchronisation problem in

  \begin{definition} 
    Let $\deadState$, $\generalState$, $\soldierState$, and $\fireState$ be four distinct states, and let $Q'$ be the set that consists of those states. A solution of the \graffito{firing mob synchronisation problem in $\mathcal{R}$ with respect to $S$}\define{firing mob synchronisation problem in $\mathcal{R}$ with respect to $S$} is a semi-cellular automaton over $\mathcal{R}$ with neighbourhood $S$ and finite set of states that includes $Q'$ such that the state $\deadState$ is dead and the set $\setOf{\deadState, \soldierState}$ is passive, and whose global transition function $\Delta$ has the following property:

    For each finite subset $A$ of $M$ such that the subgraph $\mathcal{G}[A]$ of $\mathcal{G}$ induced by $A$ is connected, each element $a \in A$, each pattern $p \in Q^A$ such that $p(a) = \generalState$ and $p\restrictedTo_{A \smallsetminus \setOf{a}} \equiv \soldierState$, and each global configuration $c$ of the form $p$, there is a non-negative integer $k$ such that the global configuration $\Delta^k(c)$ is of the form $A \to Q$, $a \mapsto \fireState$, and such that the state $\fireState$ does not occur in any of the global configurations $\Delta^j(c)$, for $j \in \N_0$ with $j < k$.
  \end{definition}

  \begin{remark}
    The firing squad synchronisation problem with one general at an arbitrary position is the firing mob synchronisation problem in $\ntuple{\ntuple{\Z, \Z, +}, \ntuple{0, \family{z}_{z \in \Z}}}$ with respect to $\setOf{-1, 0, 1}$. Note that the notions of semi-cellular and cellular automata are identical over $\ntuple{\Z, \Z, +}$.
  \end{remark}

  \begin{remark}
    Each semi-cellular automaton over $\mathcal{R}$ with neighbourhood $S$ is equivalent to a cellular automaton over the coloured $S$-Cayley graph of $\mathcal{R}$ acted upon by its automorphism group, in the sense that, for each of the former kind of automata, there is one of the latter kind with the same global transition function, and vice versa. Note that the stabilisers of coloured $S$-Cayley graphs of $\mathcal{R}$ are trivial, and hence the notions of semi-cellular and cellular automata are identical over such graphs.
  \end{remark}

  \begin{remark}
    We can regard semi-cellular automata that solve the above problem as semi-cellular automata over subgraphs of $\mathcal{G}$ that are induced by finite subsets of $M$ with one dummy neighbour in the state $\deadState$ at each edge that leads out of the subgraph. Note that, because the graph $\mathcal{G}$ is of bounded degree, the maximum degrees of the subgraphs it induces are uniformly bounded by a constant.
  \end{remark}

  \begin{remark}
    Ideally we would like an abstract description of a semi-cellular automaton that does not depend on any specifics of $\mathcal{R}$ and $S$ and that yields a solution for each choice of $\mathcal{R}$ and $S$ or at least for as huge a class of such choices as possible.
  \end{remark}

  \section{Undirected Multigraphs}
  \label{section:undirected-multigraphs}

  Undirected multigraphs without self-loops are introduced in

  \begin{definition} 
    Let $\Vertices$ and $\Edges$ be two disjoint sets, and let $\eendsOf$ be a map from $\Edges$ to $\setOf{\setOf{v, v'} \subseteq \Vertices \suchThat v \neq v'}$. The triple $\Graph = \ntuple{\Vertices, \Edges, \eendsOf}$ is called \define{undirected multigraph}\graffito{undirected multigraph $\Graph = \ntuple{\Vertices, \Edges, \eendsOf}$}\index[symbols]{Gcalligraphic@$\Graph$}\index[symbols]{V@$\Vertices$}\index[symbols]{E@$\Edges$}; each element $v \in \Vertices$ is called \define{vertex}\graffito{vertex $v$}\index[symbols]{v@$v$}; each element $e \in \Edges$ is called \define{edge}\graffito{edge $e$}\index[symbols]{e@$e$}; and, for each edge $e \in \Edges$, each vertex of $\eendsOf(e)$ is called \define{end of $e$}\graffito{ends $\eendsOf(e)$ of $e$}\index[symbols]{epsilonvare@$\eendsOf(e)$}.
  \end{definition}

  \begin{remark}
    Because each set in the codomain of $\eendsOf$ consists of exactly two distinct vertices, there are no self-loops in the undirected multigraph $\Graph$. With minor modifications the theory and the automata presented in this chapter also work if there are self-loops. They were merely excluded to make the presentation a little simpler. 
  \end{remark}

  In the remainder of this section, let $\Graph = \ntuple{\Vertices, \Edges, \eendsOf}$ be an undirected multigraph.

  What being finite means for multigraphs is said in

  \begin{definition}
    The multigraph $\Graph$ is called \define{finite}\graffito{finite multigraph} if and only if the sets $\Vertices$ and $\Edges$ are both finite.
  \end{definition}

  Isolated vertices are the ones without incident edges as introduced in

  \begin{definition}
    Let $v$ be a vertex of $\Graph$. It is called \define{isolated}\graffito{isolated vertex} if and only if, for each edge $e \in \Edges$, we have $v \notin \eendsOf(e)$.
  \end{definition}

  Directed edges are edges with distinguished source and target vertices as introduced in

  \begin{definition} 
    Let $e$ be an edge of $\Graph$, and let $v_1$ and $v_2$ be two vertices of $\Graph$ such that $\setOf{v_1, v_2} = \eendsOf(e)$. The triple $\direct{e} = (v_1, e, v_2)$ is called \define{directed edge from $v_1$ through $e$ to $v_2$}\graffito{directed edge $\direct{e}$ from $v_1$ through $e$ to $v_2$}\index{edge!directed}\index[symbols]{earrowontop@$\direct{e}$}; the vertex $\sourceOf(\direct{e}) = v_1$ is called \define{source of $\direct{e}$}\graffito{source $\sourceOf(\direct{e})$ of $\direct{e}$}\index[symbols]{sigmaearrowontop@$\sourceOf(\direct{e})$}; the edge $\bedOf(\direct{e}) = e$ is called \define{bed of $\direct{e}$}\graffito{bed $\bedOf(\direct{e})$ of $\direct{e}$}\index[symbols]{betaearrowontop@$\bedOf(\direct{e})$}; and the vertex $\targetOf(\direct{e}) = v_2$ is called \define{target of $\direct{e}$}\graffito{target $\targetOf(\direct{e})$ of $\direct{e}$}\index[symbols]{tauearrowontop@$\targetOf(\direct{e})$}.
  \end{definition}

  At each vertex there is an empty path that starts and ends at the vertex, and non-empty paths are concatenations of directed edges with matching source and target vertices as introduced in

  \begin{definition}
    \begin{aenumerate}
      \item Let $v$ be a vertex of $\Graph$. The singleton $p = (v)$ is called \define{empty path in $v$}\graffito{empty path $(v)$ in $v$}\index{path!empty}\index[symbols]{vleftparenrightparen@$(v)$}, the vertex $\sourceOf(p) = \targetOf(p) = v$ is called \defineX{source}{source of $p$}\index[symbols]{sigmap@$\sourceOf(p)$} and \define{target of $p$}\graffito{source $\sourceOf((v))$ and target $\targetOf((v))$ of $(v)$}\index[symbols]{taup@$\targetOf(p)$}, and the non-negative integer $\lengthOfPath{p} = 0$ is called \define{length of $p$}\graffito{length $\lengthOfPath{(v)}$ of $(v)$}\index[symbols]{absolutep@$\lengthOfPath{p}$}.
      \item Let $n$ be a positive integer and, for each index $i \in \discreteInterval{1}{n}$, let $\direct{e}_i$ be a directed edge of $\Graph$ such that, if $i \neq 1$, then $\sourceOf(\direct{e}_i) = \targetOf(\direct{e}_{i - 1})$. The $(2n + 1)$-tuple $p = (\sourceOf(\direct{e}_1), \bedOf(\direct{e}_1), \targetOf(\direct{e}_1), \dotsc, \bedOf(\direct{e}_n), \targetOf(\direct{e}_n))$ is called \define{path from $\sourceOf(\direct{e}_1)$ to $\targetOf(\direct{e}_n)$}\graffito{path $p$ from $\sourceOf(\direct{e}_1)$ to $\targetOf(\direct{e}_n)$}\index[symbols]{p@$p$}; the vertex $\sourceOf(p) = \sourceOf(\direct{e}_1)$ is called \define{source of $p$}\graffito{source $\sourceOf(p)$ of $p$}\index[symbols]{sigmap@$\sourceOf(p)$}; the vertex $\targetOf(p) = \targetOf(\direct{e}_n)$ is called \define{target of $p$}\graffito{target $\targetOf(p)$ of $p$}\index[symbols]{taup@$\targetOf(p)$}; and the positive integer $\lengthOfPath{p} = n$ is called \define{length of $p$}\graffito{length $\lengthOfPath{p}$ of $p$}\index[symbols]{absolutep@$\lengthOfPath{p}$}.
      \item The set of paths is denoted by $\Paths$\graffito{set $\Paths$ of paths}\index[symbols]{P aths@$\Paths$}. \qedhere
    \end{aenumerate}
  \end{definition}

  \begin{remark}
    Each directed edge is a path of length $1$.
  \end{remark}

  Subpaths are connected parts of paths as introduced in

  \begin{definition}
    Let $p = (v_0, e_1, v_1, \dotsc, e_n, v_n)$ be a path of $\Graph$, and let $k$ and $\ell$ be two indices of $\discreteInterval{0}{n}$ such that $k \leq \ell$. The path $(v_k)$, if $k = \ell$, or $(v_k, e_{k + 1}, v_{k + 1}, \dotsc, e_\ell, v_\ell)$, otherwise, is called \define{subpath of $p$}\graffito{subpath $p_{\discreteInterval{k}{\ell}}$ of $p$} and is denoted by $p_{\discreteInterval{k}{\ell}}$\index[symbols]{pklsubscript@$p_{\discreteInterval{k}{\ell}}$}.
  \end{definition}

  Direction-preserving paths are the ones without U-turns as introduced in

  \begin{definition}
    Let $p = (v_0, e_1, v_1, \dotsc, e_n, v_n)$ be a path of $\Graph$. It is called \define{direction-preserving}\graffito{direction-preserving} if and only if, for each index $i \in \discreteInterval{1}{n - 1}$, we have $e_i \neq e_{i + 1}$. The \graffito{set $\Paths_{\directionPreserving}$ of direction-preserving paths}set of direction-preserving paths is denoted by $\Paths_{\directionPreserving}$\index[symbols]{Parrowrightshortsubscript@$\Paths_{\directionPreserving}$}.
  \end{definition}


  Two paths with matching target and source vertices can be concatenated as introduced in

  \begin{definition}
    Let $p = (v_0, e_1, v_1, \dotsc, e_n, v_n)$ and $p' = (v_0', e_1', v_1',\allowbreak \dotsc, e_{n'}', v_{n'}')$ be two paths of $\Graph$ such that $v_n = v_0'$. The path $p \concat p' = (v_0, e_1, v_1, \dotsc, e_n, v_n, e_1', v_1', \dotsc, e_{n'}', v_{n'}')$ is called \graffito{concatenation of $p$ and $p'$}\define{concatenation $p \concat p'$ of $p$ and $p'$}\index[symbols]{bullet@$\bullet$}.
  \end{definition}

  \begin{remark}
    Each non-empty path is the concatenation of directed edges.
  \end{remark}

  What being connected means for multigraphs is said in

  \begin{definition}
    The multigraph $\Graph$ is called \define{connected}\graffito{connected multigraph} if and only if, for each tuple $(v, v') \in \Vertices \times \Vertices$, there is a path from $v$ to $v'$.
  \end{definition}

  \begin{remark}
    The multigraph $\Graph$ is connected if and only if, for each tuple $(v, v') \in \Vertices \times \Vertices$, there is a direction-preserving path from $v$ to $v'$.
  \end{remark}


  We can assign weights to edges as done in

  \begin{definition} 
    Let $\weightOf$ be a map from $\Edges$ to $\R_{> 0}$. It is called \graffito{edge weighting $\weightOf$ of $\Graph$}\define{edge weighting of $\Graph$}\index[symbols]{omega@$\weightOf$}, and, for each edge $e \in E$, the element $\weightOf(e)$ is called \define{edge weight of $e$}\graffito{edge weight $\weightOf(e)$ of $e$}\index[symbols]{omegae@$\weightOf(e)$}.
  \end{definition}


  Edge weights induce weights of paths as introduced in

  \begin{definition}
    Let $\weightOf$ be an edge weighting of $\Graph$ and let $p = (v_0, e_1, v_1, \dotsc, e_n, v_n)$ be a path of $\Graph$. The sum $\weightOf(p) = \sum_{i = 1}^n \weightOf(e_i)$ is called \define{weight of $p$}\graffito{weight $\weightOf(p)$ of $p$}\index[symbols]{omegap@$\weightOf(p)$}. 
  \end{definition}

  \begin{remark}
    Each empty path has weight $0$.
  \end{remark}

  \begin{remark}
    Each directed edge has the same weight as its bed.
  \end{remark}

  \section{Continuum Representation}
  \label{section:continuum-representation}

  In this section, let $\Graph = \ntuple{\Vertices, \Edges, \eendsOf}$ be an undirected multigraph without isolated vertices and let $\weightOf$ be an edge weighting of $\Graph$. 

  An orientation is a choice of source and target vertices for each edge as introduced in

  \begin{definition}
    Let $\sigma$ and $\tau$ be two maps from $\Edges$ to $\Vertices$ such that, for each edge $e \in \Edges$, we have $\setOf{\sigma(e), \tau(e)} = \eendsOf(e)$. The tuple $(\sigma, \tau)$ is called \define{orientation of $\Graph$}\graffito{orientation $(\sigma, \tau)$ of $\Graph$}\index[symbols]{sigmatautuple@$(\sigma, \tau)$}.
  \end{definition}

  Realising weighted edges as disjoint intervals and gluing these intervals together at shared ends yields a continuum representation of $\Graph$ and is done in

  \begin{definition}
  \label{definition:continuum-representation}
    Let $(\sigma, \tau)$ be an orientation of $\Graph$, let 
    \begin{align*}
      \zeta \from \R \smallsetminus \setOf{0} &\to \setOf{\sigma, \tau}, \mathnote{map $\zeta$ from $\R \smallsetminus \setOf{0}$ to $\setOf{\sigma, \tau}$}\index[symbols]{zeta@$\zeta$}\\
      r &\mapsto \begin{dcases*}
                   \sigma, &if $r < 0$,\\
                   \tau,   &if $r > 0$,
                 \end{dcases*}
    \end{align*}
    let
    \begin{equation*}
      \left\{
        \begin{aligned}
          \phi \from \Edges &\to \R_{< 0}, \index[symbols]{phi@$\phi$}\\
          e &\mapsto - \frac{\weightOf(e)}{2},
        \end{aligned}
      \right\}
      \text{ and }
      \left\{
        \begin{aligned}
          \psi \from \Edges &\to \R_{> 0}, \index[symbols]{psi@$\psi$}\\
          e &\mapsto \frac{\weightOf(e)}{2},
        \end{aligned}
      \right\}
      \mathnote{maps $\phi$ and $\psi$ from $\Edges$ to $\R_{< 0}$ and $\R_{> 0}$}
    \end{equation*}
    and let $\sim$\graffito{equivalence relation $\sim$ on $\R \times \Edges$}\index[symbols]{tilde@$\sim$} be the equivalence relation on $\R \times \Edges$ such that, for each tuple $(r, e) \in \R \times \Edges$ and each tuple $(r', e') \in \R \times \Edges$, 
    \begin{align*}
      (r, e) \sim (r', e') \ifAndOnlyIf
      &r \in \setOf{\phi(e), \psi(e)}\\
      &\land r' \in \setOf{\phi(e'), \psi(e')}\\
      &\land \zeta(r)(e) = \zeta(r')(e').
    \end{align*}
    The set $\continuumRepresentationOf{\Graph} = (\bigcup_{e \in \Edges} \closedInterval{\phi(e), \psi(e)} \times \setOf{e}) \modulo {\sim}$ is called \graffito{continuum representation $\continuumRepresentationOf{\Graph}$ of $\Graph$}\define{continuum representation of $\Graph$}\index[symbols]{Gcalligraphicbar@$\continuumRepresentationOf{\Graph}$}.
  \end{definition} 

  \begin{remark}
    Each weighted edge is realised as a closed interval whose length is the edge's weight. These intervals are made disjoint by taking the Cartesian product with the respective edge. And they are glued together at shared ends by taking the set of all these disjoint intervals modulo the equivalence relation $\sim$. The vertices are implicitly realised as end points or junctions of the glued disjoint intervals.
  \end{remark}

  \begin{remark}
    If the graph $\Graph$ contained isolated vertices, then they would not be represented in $\continuumRepresentationOf{\Graph}$. With minor modifications the theory presented in this chapter also works if there are isolated vertices. They were merely excluded to make the presentation a little simpler.
  \end{remark}

  In the remainder of this section, let $\continuumGraph$ be a continuum representation of $\Graph$ with respect to an orientation $(\sigma, \tau)$, and let $\zeta$, $\phi$, $\psi$, and $\sim$ be the maps and the equivalence relation from \cref{definition:continuum-representation}.

  Vertices are canonically embedded into $\continuumGraph$ as is done in

  \begin{definition}
    The map
    \begin{align*}
      \continuumRepresentationOf{\phantom{v}} \from V &\to \continuumGraph, \mathnote{vertex embedding $\continuumRepresentationOf{\phantom{v}}$}\index[symbols]{vbar@$\continuumRepresentationOf{v}$}\\
      \sourceOf(e) &\mapsto \equivalenceClassOf{(\phi(e), e)}_\sim,\\
      \targetOf(e) &\mapsto \equivalenceClassOf{(\psi(e), e)}_\sim,
    \end{align*}
    embeds vertices of $\Graph$ into $\continuumGraph$. Its image is denoted by $\continuumVertices$\graffito{vertices $\continuumVertices$}\index[symbols]{Vfraktur@$\continuumVertices$} and each element $\mathfrak{v} \in \continuumVertices$ is called \define{vertex}\graffito{vertex $\mathfrak{v}$}\index[symbols]{vfraktur@$\mathfrak{v}$}. 
  \end{definition}

  \begin{remark}
    The embedding is well-defined due to the definition of the equivalence relation $\sim$. 
  \end{remark}

  Edges are canonically embedded into the power set of $\continuumGraph$ as is done in

  \begin{definition}
    The map
    \begin{align*}
      \continuumRepresentationOf{\phantom{e}} \from E &\to \powerSetOf(\continuumGraph), \mathnote{edge embedding $\continuumRepresentationOf{\phantom{e}}$}\index[symbols]{ebar@$\continuumRepresentationOf{e}$}\\
      e &\mapsto ([\phi(e), \psi(e)] \times \setOf{e}) \modulo {\sim},
    \end{align*}
    embeds edges of $\Graph$ into $\continuumGraph$. Its image is denoted by $\continuumEdges$\graffito{edges $\continuumEdges$}\index[symbols]{Efraktur@$\continuumEdges$} and each element $\mathfrak{e} \in \continuumEdges$ is called \define{edge}\graffito{edge $\mathfrak{e}$}\index[symbols]{efraktur@$\mathfrak{e}$}. 
  \end{definition}

  At each point of $\continuumGraph$ there is at least one direction to move: In a vertex of degree $k$, there are $k$ directions; and on an edge but not in one of its endpoints, there are $2$ directions. An inefficient but immediate way to represent these directions is as in

  \begin{definition}
    The set $\setOf{-1, 1} \times \Edges$ is denoted by $\Directions$\graffito{set $\Directions$ of directions}\index[symbols]{Dir@$\Directions$}, each element $d = (o, e) \in \Directions$ is called \defineX{direction on $e$}{direction}\graffito{direction $d$ on $e$}\index[symbols]{direction@$d$}, the element $o$ is called \define{orientation of $d$}\graffito{orientation $o$ of $d$}\index[symbols]{orientation@$o$}, the involution
    \begin{align*}
      \reverse \from \Directions &\to \Directions, \mathnote{orientation reversing involution $\reverse$}\index[symbols]{minus@$\reverse$}\\
      (o, e) &\mapsto (-o, e),
    \end{align*}
    reverses the orientation of directions, and the map
    \begin{align*}
      \directionOf \from \continuumGraph &\to \powerSetOf(\Directions), \mathnote{map $\directionOf$ that assigns directions}\index[symbols]{direction@$\directionOf$}\\
      \equivalenceClassOf{(r, e)}_\sim &\mapsto
        \left\{\begin{aligned} 
          &\setOf{(-1, e), (1, e)}, \text{ if $r \in \openInterval{\phi(e), \psi(e)}$},\\ 
          &\setOf{(- \signOf(r'), e') \suchThat (r', e') \in \equivalenceClassOf{(r, e)}_\sim}, \text{ otherwise}, 
        \end{aligned}\right.
    \end{align*}
    assigns to each point in $\continuumGraph$ the set of possible directions in which someone standing on that point can move.
  \end{definition}

  \begin{remark}
  \label{remark:efficient-representation-of-directions}
    This representation of directions is inefficient in the following sense: If we stand on an edge but not on one of its endpoints, then the orientation is enough directional information; and if we stand on a vertex, then the edge is enough directional information, because the orientation is implicit in the fact that we can only move onto the edge but not off it since we are in one end of the edge. On a vertex we do not even need the edge itself but only an identifier for the edge that is locally unique; for example, we could colour the edges such that no two edges of the same colour are incident to the same vertex and use this colour instead.
  \end{remark} 

  Like vertices, paths of $\Graph$ are also canonically embedded into $\continuumGraph$ and each embedding can be unit-speed parametrised by the interval from $0$ to the path's weight as is inductively done in

  \begin{definition}
    The map 
    \begin{align*}
      \continuumRepresentationOf{\phantom{p}} \from \Paths &\to \continuumGraph^{\setOf{\closedInterval{0, r} \suchThat r \in \R_{\geq 0}}}, \mathnote{path embedding $\continuumRepresentationOf{\phantom{p}}$}\index[symbols]{pbar@$\continuumRepresentationOf{p}$}\\
      (v) &\mapsto \left[
                     \begin{aligned}
                       \closedInterval{0, 0} &\to \continuumGraph,\\
                       r &\mapsto \continuumRepresentationOf{v}
                     \end{aligned}
                   \right]\\
      (\sourceOf(e), e, \targetOf(e)) &\mapsto \left[
                                             \begin{aligned}
                                               \closedInterval{0, \weightOf(e)} &\to \continuumGraph,\\
                                               r &\mapsto \equivalenceClassOf{(\phi(e) + r, e)}_\sim,
                                             \end{aligned}
                                           \right]\\
      (\targetOf(e), e, \sourceOf(e)) &\mapsto \left[
                                             \begin{aligned}
                                               \closedInterval{0, \weightOf(e)} &\to \continuumGraph,\\
                                               r &\mapsto \equivalenceClassOf{(\psi(e) - r, e)}_\sim,
                                             \end{aligned}
                                           \right]\\
      (v_0, e_1, v_1) \concat p' &\mapsto \left[ 
                                            \begin{aligned}
                                              &\closedInterval{0, \weightOf((v_0, e_1, v_1) \concat p')} \to \continuumGraph,\\
                                              &r \mapsto \begin{dcases*}
                                                           \continuumRepresentationOf{(v_0, e_1, v_1)}(r), &if $r \leq \weightOf(e_1)$,\\
                                                           \continuumRepresentationOf{p'}(r - \weightOf(e_1)), &otherwise,
                                                         \end{dcases*} 
                                            \end{aligned}
                                          \right]
    \end{align*}
    maps paths of $\Graph$ to unit-speed parametrisations of them in $\continuumGraph$.
  \end{definition}

  \begin{remark}
    The base cases of the inductive definition do not overlap because there are no self-loops, and the inductive step is well-defined because $\weightOf((v_0, e_1, v_1) \concat p') = \weightOf(e_1) + \weightOf(p')$.
  \end{remark}

  The images $\continuumRepresentationOf{\Paths}$ and $\continuumRepresentationOf{\Paths_{\directionPreserving}}$ consist broadly speaking of paths and direction-preserving paths in $\continuumGraph$ from vertices to vertices that only change direction at vertices. Restricting the parametrisation intervals of paths in $\continuumRepresentationOf{\Paths_{\directionPreserving}}$ to subintervals and doing a reparametrisation such that the new parametrisation starts at $0$ yields all direction-preserving paths in $\continuumGraph$ and is done in

  \begin{definition}
    The set
    \begin{equation*}
      \setOf{\continuumRepresentationOf{p}\restrictedTo_{\closedInterval{r, s}}(\blank + r) \suchThat p \in \Paths_{\directionPreserving} \text{ and } r, s \in \closedInterval{0, \weightOf(p)} \text{ with } r \leq s}
    \end{equation*} 
    is denoted by $\continuumPaths_{\directionPreserving}$\graffito{set $\continuumPaths_{\directionPreserving}$ of direction-preserving paths $\mathfrak{p}$}\index[symbols]{Parrowrightsubscriptfraktur@$\continuumPaths_{\directionPreserving}$}; each element $\mathfrak{p} \in \continuumPaths_{\directionPreserving}$ is called \define{direction-preserving path}\index[symbols]{pfraktur@$\mathfrak{p}$}, the length of the interval $\domainOf(\mathfrak{p})$ is called \define{length of $\mathfrak{p}$}\graffito{length $\length(\mathfrak{p})$ of $\mathfrak{p}$} and is denoted by $\length(\mathfrak{p})$\index[symbols]{omegapfraktur@$\length(\mathfrak{p})$}, the point $\sourceOf(\mathfrak{p}) = \mathfrak{p}(0)$ is called \define{source of $\mathfrak{p}$}\graffito{source $\sourceOf(\mathfrak{p})$ of $\mathfrak{p}$}\index[symbols]{sigmapfraktur@$\sourceOf(\mathfrak{p})$}, the point $\targetOf(\mathfrak{p}) = \mathfrak{p}(\length(\mathfrak{p}))$ is called \define{target of $\mathfrak{p}$}\graffito{target $\targetOf(\mathfrak{p})$ of $\mathfrak{p}$}\index[symbols]{taupfraktur@$\targetOf(\mathfrak{p})$}, and the path $\mathfrak{p}$ is called \define{empty}\graffito{empty path} if and only if $\length(\mathfrak{p}) = 0$.
  \end{definition} 

  \begin{remark}
    Doing the same with the paths in $\continuumRepresentationOf{\Paths}$ does not yield all paths in $\continuumGraph$ but only those that change direction at vertices and not on edges. Because we only need direction-preserving paths in what is to come, we do not define what a general path on $\continuumGraph$ is.
  \end{remark}

  \begin{remark}
    Sources and targets of direction-preserving paths in $\continuumGraph$ are in general not vertices.
  \end{remark}

  The distance between two points is the length of the shortest path between the points as introduced in

  \begin{definition}
    The map
    \begin{align*}
      \distanceOf \from \continuumGraph \times \continuumGraph &\to \R_{\geq 0} \cup \setOf{\infty}, \mathnote{distance $\distanceOf$}\index[symbols]{d@$\distanceOf$}\\
      (\mathfrak{m}, \mathfrak{m}') &\mapsto \inf\setOf{\length(\mathfrak{p}) \suchThat \mathfrak{p} \in \continuumPaths_{\directionPreserving} \text{ with } \sourceOf(\mathfrak{p}) = \mathfrak{m} \text{ and } \targetOf(\mathfrak{p}) = \mathfrak{m}'} 
    \end{align*}
    is called \define{distance}, where the infimum of the empty set is infinity.
  \end{definition}

  \begin{remark}
    If the graph $\Graph$ is finite and connected, then the distance map $\distanceOf$ is a metric. Otherwise, it may not be a metric. For example, if there are two distinct vertices $v$, $v' \in \Vertices$ such that, for each $n \in \N_+$, there is an edge $e \in \Edges$ whose weight is $1/n$, then the distance of $\continuumRepresentationOf{v}$ and $\continuumRepresentationOf{v'}$ is $0$ although $\continuumRepresentationOf{v} \neq \continuumRepresentationOf{v'}$. Or, if the graph $\Graph$ is not connected, then there are two points $\mathfrak{m}$, $\mathfrak{m}' \in \continuumGraph$ whose distance is $\infty$.
  \end{remark}

  Each non-zero vector of a vector space is uniquely determined by its magnitude and its direction, and the zero vector is already uniquely determined by its magnitude, which is $0$, and can be thought of as pointing in every direction, which can be represented by the set of directions. A generalisation of vector spaces is given in

  \begin{definition} 
    Let $\every$ be the set $\Directions$. The set
    \begin{equation*}
      \Arrows = \setOf{(0, \every)} \cup (\R_{> 0} \times \Directions) \mathnote{arrow space $\Arrows$}\index[symbols]{Arr@$\Arrows$}
    \end{equation*}
    is called \define{arrow space}; each element $a \in \Arrows$ is called \define{arrow}\graffito{arrow $a$}\index[symbols]{a@$a$}; the set $\every$ is called \define{semi-direction}\graffito{semi-direction $\every$}\index{direction!semi-}\index[symbols]{vry@$\every$}; for each element $a = (r, d) \in \Arrows$, the real number $\normOf{a} = r$ is called \define{magnitude of $a$}\graffito{magnitude $\normOf{a}$ of $a$}\index[symbols]{norma@$\normOf{a}$}, and the (semi-)direction $\directionOf(a) = d$ is called \defineX{(semi-)direction of $a$}{direction of $a$ semi@(semi-)direction of $a$}\graffito{(semi-)direction $\directionOf(a)$ of $a$}\index{semi-direction of $a$@(semi-)direction of $a$}\index[symbols]{dira@$\directionOf(a)$}.
  \end{definition}

  Arrow spaces will be used to represent both velocities, which are directed speeds, and directed distances. Multiplying a velocity by a time yields a directed distance. This scalar multiplication is introduced in

  \begin{definition}
    The map
    \begin{align*}
      \multipli \from \Arrows \times \R_{\geq 0} &\to \Arrows, \mathnote{scalar multiplication $\multipli$}\index[symbols]{dotcentre@$\multipli$}\index[symbols]{centredot@$\multipli$}\\
      ((r, d), s) &\mapsto \begin{dcases*}
                             (0, \every), &if $s = 0$,\\
                             (r \multipli s, d), &if $s > 0$,
                           \end{dcases*}
    \end{align*}
    is called \define{scalar multiplication}.
  \end{definition}

  When we stand at the beginning of a non-empty direction-preserving path and walk along it until we reach its end, we start our walk on the first edge of the path in a certain direction and we end it on the last edge of the path in a certain direction. These directions are introduced in

  \begin{definition}
    Let $\mathfrak{p}$ be a direction-preserving path of $\continuumPaths_{\directionPreserving}$. If $\mathfrak{p}$ is empty, let $\directionOf_\sourceOf(\mathfrak{p}) = \every$ and let $\directionOf_\targetOf(\mathfrak{p}) = \every$.

    Otherwise, there are two edges $e_\sourceOf$, $e_\targetOf \in E$, which may be the same, and there are two positive real numbers $\xi_\sourceOf$, $\xi_\targetOf \in \leftOpenAndRightClosedInterval{0, \length(\mathfrak{p})}$ such that $\mathfrak{p}(\closedInterval{0, \xi_\sourceOf}) \subseteq \continuumRepresentationOf{e}_\sourceOf$ and $\mathfrak{p}(\closedInterval{\length(\mathfrak{p}) - \xi_\targetOf, \length(\mathfrak{p})}) \subseteq \continuumRepresentationOf{e}_\targetOf$. Moreover, there are four real numbers $r_\sourceOf$, $r_\sourceOf'$, $r_\targetOf$, and $r_\targetOf'$ such that $\equivalenceClassOf{(r_\sourceOf, e_\sourceOf)}_\sim = \mathfrak{p}(0)$, $\equivalenceClassOf{(r_\sourceOf', e_\sourceOf)}_\sim = \mathfrak{p}(\xi_\sourceOf)$, $\equivalenceClassOf{(r_\targetOf, e_\targetOf)}_\sim = \mathfrak{p}(\length(\mathfrak{p}) - \xi_\targetOf)$, and $\equivalenceClassOf{(r_\targetOf', e_\targetOf)}_\sim = \mathfrak{p}(\length(\mathfrak{p}))$. Let $\directionOf_\sourceOf(\mathfrak{p}) = (\signOf(r_\sourceOf' - r_\sourceOf), e_\sourceOf)$ and let $\directionOf_\targetOf(\mathfrak{p}) = (\signOf(r_\targetOf' - r_\targetOf), e_\targetOf)$.

    In both cases, the (semi-)direction $\directionOf_\sourceOf(\mathfrak{p})$ is called \graffito{source direction $\directionOf_\sourceOf(\mathfrak{p})$ of $\mathfrak{p}$}\define{source direction of $\mathfrak{p}$}\index[symbols]{dirsigmapfraktur@$\directionOf_\sourceOf(\mathfrak{p})$} and the (semi-)direction $\directionOf_\targetOf(\mathfrak{p})$ is called \define{target direction of $\mathfrak{p}$}\graffito{target direction $\directionOf_\targetOf(\mathfrak{p})$ of $\mathfrak{p}$}\index[symbols]{dirtaupfraktur@$\directionOf_\targetOf(\mathfrak{p})$}.
  \end{definition} 

  \begin{remark}
    The edge $e_\sourceOf$ is the first edge of the path $\mathfrak{p}$ and the edge $e_\targetOf$ is its last edge. The numbers $\xi_\sourceOf$ and $\xi_\targetOf$ are two positive real numbers such that the first $\xi_\sourceOf$ length units of the path run on its first edge and the last $\xi_\targetOf$ length units of the path run on its last edge. The numbers $r_\sourceOf$ and $r_\sourceOf'$ are the positions of the path on its first edge at its very beginning and after $\xi_\sourceOf$ length units, and the numbers $r_\targetOf$, and $r_\targetOf'$ are the positions of the path on its last edge $\xi_\targetOf$ length units before its end and at its very end. Therefore, the signum of $r_\sourceOf' - r_\sourceOf$ is the start direction on the first edge of the path and the signum of $r_\targetOf' - r_\targetOf$ is the end direction on the last edge of the path.
  \end{remark}

  \begin{remark}
    For each non-empty path $\mathfrak{p} \in \continuumPaths_{\directionPreserving}$, we have $\directionOf_\sourceOf(\mathfrak{p}) \in \directionOf(\sourceOf(\mathfrak{p}))$ and $\directionOf_\targetOf(\mathfrak{p}) \in \reverse \directionOf(\targetOf(\mathfrak{p}))$.
  \end{remark}

  \section{Signal Machines}
  \label{section:signal-machines}

  In this section, let $\Graph = \ntuple{\Vertices, \Edges, \eendsOf}$ be a non-trivial, finite, and connected undirected multigraph, let $\weightOf$ be an edge weighting of $\Graph$, and let $M$ be a continuum representation of $\Graph$. Recall that, according to our definition of undirected multigraphs, there are no self-loops in $\Graph$. To motivate the definitions in this section, we talk as if there were a signal machine in front of us whose time evolution we can observe, although this evolution is not completely defined until the end of this section. 

  If you observe the time evolution of a signal machine on the graph $M$, you see signals of different kinds and various speeds each carrying some data move along edges. When signals collide, they may be reflected, removed, new signals may be created, and so on. Similarly, when signals reach a vertex, they may be removed, copies of them may be sent onto all incident edges, new signals may be created, and so on. You may also see stationary signals and signals that travel side-by-side at the same speed. What happens when signals collide or reach a vertex is decided by two local transition functions, one that handles such \defineX{events}{event}\graffito{event} in vertices and one that handles them on edges.

  The only events on edges are collisions. In each collision on an edge, there are at least two signals involved, the involved signals are either stationary or they move in one of the two possible directions, and at least two of the signals collide head-on or rear-end. Such a collision results in a set of signals that are either stationary or move in one of the two possible directions.

  A vertex may be reached by just one signal or multiple signals may collide in it. In both cases, there is at least one signal involved, the involved signals are either stationary or they moved towards the vertex just before the event, and at least one signal is moving. Such an event results in a set of signals that are either stationary or move away from the vertex along incident edges.

  \begin{definition}
    Let $\Kinds$\graffito{set $\Kinds$}\index[symbols]{Knd@$\Kinds$} be a set, let $\speedOf$\graffito{map $\speedOf$} be a map from $\Kinds$ to $\R_{\geq 0}$, let $\family{\Data_k}_{k \in \Kinds}$\graffito{family $\family{\Data_k}_{k \in \Kinds}$} be a family of sets, let\graffito{set $\Signals$}\graffito{set $\Directions_e$} 
    \begin{equation*}
      \Signals = \setOf{(k, d, u) \suchThat k \in \Kinds \text{, } (\speedOf(k), d) \in \Arrows \text{, and } u \in \Data_k}, \index[symbols]{Sgnl@$\Signals$}
    \end{equation*}
    let $\Directions_e = \setOf{\setOf{d, \reverse d} \suchThat d \in \Directions}$\index[symbols]{Dire@$\Directions_e$}, let
    \begin{align*}
      \domainOf(\localTransitionFunction_e) = \{&S \in \powerSetOf(\Signals) \suchThat
          \Exists D \in \Directions_e \SuchThat
              \cardinalityOf{S} \geq 2 \text{ and }\mathnote{set $\domainOf(\localTransitionFunction_e)$}\index[symbols]{domdeltae@$\domainOf(\localTransitionFunction_e)$}\\ 
              &\ForEach (k, d, u) \in S \Holds d \in \setOf{\every} \cup D \text{ and }\\ 
              &\Exists (k, d, u) \in S \Exists (k', d', u') \in S \SuchThat
                  \begin{aligned}[t]
                    &d \neq d' \text{ or }\\
                    &\speedOf(k) \neq \speedOf(k')\},
                  \end{aligned} 
    \end{align*}
    let $\localTransitionFunction_e$\graffito{map $\localTransitionFunction_e$}\index[symbols]{deltae@$\localTransitionFunction_e$} be a map from $\domainOf(\localTransitionFunction_e)$ to $\powerSetOf(\Signals)$ such that
    \begin{equation*}
      \ForEach S \in \domainOf(\localTransitionFunction_e) \Exists D \in \Directions_e \SuchThat \ForEach (k, d, u) \in S \cup \localTransitionFunction_e(S) \Holds d \in \setOf{\every} \cup D,
    \end{equation*}
    let $\Directions_v = \directionOf(\continuumVertices)$\graffito{set $\Directions_v$}\index[symbols]{Dirv@$\Directions_v$}, let
    \begin{align*}
      \domainOf(\localTransitionFunction_v) = \setOf{&(D, S) \in \Directions_v \times \powerSetOf(\Signals) \suchThat
          \cardinalityOf{S} \geq 1 \text{ and }\mathnote{set $\domainOf(\localTransitionFunction_v)$}\index[symbols]{domdeltav@$\domainOf(\localTransitionFunction_v)$}\\
          &\ForEach (k, d, u) \in S \Holds d \in \setOf{\every} \cup \reverse D \text{ and }\\ 
          &\Exists (k, d, u) \in S \SuchThat \speedOf(k) > 0},
    \end{align*}
    and let $\localTransitionFunction_v$\graffito{map $\localTransitionFunction_v$}\index[symbols]{deltav@$\localTransitionFunction_v$} be a map from $\domainOf(\localTransitionFunction_v)$ to $\powerSetOf(\Signals)$ such that
    \begin{equation*}
      \ForEach (D, S) \in \domainOf(\localTransitionFunction_v) \ForEach (k, d, u) \in \localTransitionFunction_v(D, S) \Holds d \in \setOf{\every} \cup D.
    \end{equation*}

    The quadruple $\mathcal{S} = \ntuple{\Kinds, \speedOf, \family{\Data_k}_{k \in \Kinds}, (\localTransitionFunction_e, \localTransitionFunction_v)}$\graffito{signal machine $\mathcal{S}$}\index[symbols]{Scalligraphic@$\mathcal{S}$} is called \define{signal machine}; each element $k \in \Kinds$ is called \define{kind}\graffito{kind $k$}\index[symbols]{k@$k$}; for each kind $k \in \Kinds$, the non-negative real number $\speedOf(k)$ is called \define{speed of $k$}\graffito{speed $\speedOf(k)$ of $k$}\index[symbols]{spdk@$\speedOf(k)$} and the set $\Data_k$ is called \define{data set of $k$}\graffito{data set $\Data_k$ of $k$}\index[symbols]{Dtk@$\Data_k$}; each element $s \in \Signals$ is called \define{signal}\graffito{signal $s$}\index[symbols]{s@$s$}; and the maps $\localTransitionFunction_e$ and $\localTransitionFunction_v$ are called \define{local transition function on edges} and \define{in vertices}\graffito{local transition functions $\localTransitionFunction_e$ and $\localTransitionFunction_v$ on edges and in vertices} respectively.
  \end{definition} 

  \begin{remark}
    The local transition function $\localTransitionFunction_e$ is used to handle events on edges but not in their endpoints. It gets the signals that are involved in the event and returns the resulting signals. Because in each event at least one moving signal is involved, the direction of this signal and the map that reverses orientation can be used by $\localTransitionFunction_e$ to determine the two possible directions the resulting signals may have.

    The local transition function $\localTransitionFunction_v$ is used to handle events in vertices. It gets the directions signals may take at the respective vertex and the signals that are involved in the event and returns the resulting signals.

    The local transition functions $\localTransitionFunction_e$ and $\localTransitionFunction_v$ are supposed to regard directions as black boxes that can merely be distinguished and whose orientation can be reversed. They must not determine edges or vertices by deconstructing directions, which is possible with the chosen representation of directions. If they did something like that, the signal machine would not be uniform.
  \end{remark} 

  \begin{remark}
    At the beginning of this section we fixed a general multigraph. This multigraph should be regarded as the blueprint of a multigraph. A signal machine depends only on that blueprint and not on any specific properties that a concrete choice of a multigraph may have. So, one and the same signal machine can be instantiated for any multigraph, each instantiation results in a machine on a concrete multigraph, and these instantiations are uniform in the chosen multigraphs. In other words, a signal machine is a map from the set of all multigraphs to the set of quadruples that describe instantiations of the machine on concrete multigraphs and this map depends only in a trivial way on its argument. 

    The quadruple that describes signal machines could be made independent of multigraphs by choosing a different representation for directions. They could for example be represented by integers or vectors or colours equipped with an involution to switch the orientation of directions.

    Even the global transition function, which is introduced below and describes the time evolution of a signal machine, could be made independent of multigraphs by representing them as patterns in high-dimensional Euclidean spaces (think of the drawing of a graph on a piece of paper). The directions are then vectors that are tangential to edges.

    In classical solutions of the firing squad synchronisation problem the regions to be synchronised are actually represented as patterns in integer lattices: The cells outside the region are in the same state, say $0$, and all cells inside the region are not in state $0$, more precisely, one cell inside the region is in a state that distinguishes it as the general, say $1$, and the other cells inside the region are all in the same state, say $2$.
  \end{remark}

  In the remainder of this section, let $\mathcal{S} = \ntuple{\Kinds, \speedOf, \family{\Data_k}_{k \in \Kinds}, (\localTransitionFunction_e, \localTransitionFunction_v)}$ be a signal machine.

  To describe the time evolution of our signal machine, the following notions are convenient.

  \begin{definition}
    Let $\Times$\graffito{set $\Times$}\index[symbols]{T@$\Times$} be the set $\R_{\geq 0}$, let $\extendedTimes$\graffito{set $\extendedTimes$}\index[symbols]{Tbar@$\extendedTimes$} be the set $\Times \cup \setOf{\infty}$, let $\States$\graffito{set $\States$}\index[symbols]{Q@$\States$} be the set $\powerSetOf(\Signals)$, and let $\Configurations$\graffito{set $\Configurations$}\index[symbols]{Cnf@$\Configurations$} be the set $Q^M$. Each element $t \in \extendedTimes$ is called \define{time}\graffito{time $t$}\index[symbols]{t@$t$} and the time $\infty$ is called \define{improper}\graffito{improper time}\index[symbols]{infinity@$\infty$}, each element $q \in \States$ is called \define{state}\graffito{state $q$}\index[symbols]{q@$q$}, and each element $c \in \Configurations$ is called \define{configuration}\graffito{configuration $c$}\index[symbols]{c@$c$}.
  \end{definition}

  The components of signals and some compounds of them are named in

  \begin{definition}
    Let $s = (k, d, u)$ be a signal of $\Signals$. The kind $k$ of $s$ is denoted by $\kindOf(s)$\graffito{kind $\kindOf(s)$ of $s$}\index[symbols]{knds@$\kindOf(s)$}; the speed $\speedOf(k)$ of $s$ is denoted by $\speedOf(s)$\graffito{speed $\speedOf(s)$ of $s$}\index[symbols]{spds@$\speedOf(s)$}; the (semi-)direction $d$ of $s$ is denoted by $\directionOf(s)$\graffito{(semi-)direction $\directionOf(s)$ of $s$}\index[symbols]{dirs@$\directionOf(s)$}; the velocity $(\speedOf(k), d)$ of $s$ is denoted by $\velocityOf(s)$\graffito{velocity $\velocityOf(s)$ of $s$}\index[symbols]{vels@$\velocityOf(s)$}; and the datum $u$ of $s$ is denoted by $\datumOf(s)$\graffito{datum $\datumOf(s)$ of $s$}\index[symbols]{dts@$\datumOf(s)$}.
  \end{definition}

  \begin{remark}
    For each time $t \in \Times$, the arrow $\velocityOf(s) \multipli t$ is equal to the arrow $(\speedOf(k) \multipli t, d)$, which can be interpreted as a directed distance.
  \end{remark}

  Signals of speed $0$ are named in

  \begin{definition}
    Let $s$ be a signal of $\Signals$. It is called \define{stationary}\graffito{stationary signal} if and only if $\speedOf(s) = 0$.
  \end{definition}

  When we stand on a point facing in a direction and from there we walk a fixed distance without making U-turns but otherwise making arbitrary choices at each vertex, then there is a finite number of points that we may reach and we reach them walking in some direction. The set of these points with and without target-directions is introduced in

  \begin{definition}
    Let $m$ be a point of $M$ and let $(\ell, d)$ be an arrow. The set of points with target-directions that can be reached from $m$ by a direction-preserving path of length $\ell$ with source-direction $d$ is
    \begin{equation*}
      \ReachOf_{(\ell, d)}^{\directed}(m) =
          \begin{aligned}[t]
            \setOf{(\targetOf(\mathfrak{p}), \directionOf_{\xend}(\mathfrak{p})) \suchThat{}
                 &\mathfrak{p} \in \continuumPaths_{\directionPreserving} \text{, } \length(\mathfrak{p}) = \ell \text{, }\\
                 &\directionOf_{\start}(\mathfrak{p}) = d \text{, and } \sourceOf(\mathfrak{p}) = m}
          \end{aligned}
      \mathnote{$\ReachOf_{(\ell, d)}^{\directed}(m)$}
    \end{equation*}
    and without target-directions it is
    \begin{equation*}
      \ReachOf_{(\ell, d)}(m) = \setOf{m' \suchThat \Exists d \in \Directions \SuchThat (m', d) \in \ReachOf_{(\ell, d)}^{\directed}(m)}.\mathnote{$\ReachOf_{(\ell, d)}(m)$} \qedhere 
    \end{equation*}
  \end{definition}

  An event occurs when a signal reaches a vertex or two signals coming from different points collide. The time of the next event is given a name in

  \begin{definition}
    Let $c$ be a configuration of $\Configurations$. The minimum time until a signal in $c$ reaches a vertex is
    \begin{equation*}
      t' = \inf_{m \in M} \inf_{\substack{s \in c(m)\\ \speedOf(s) > 0}} \inf \setOf{t \in \R_{> 0} \suchThat \ReachOf_{\velocityOf(s) \multipli t}(m) \cap \continuumVertices \neq \emptyset}.
    \end{equation*}
    The minimum time until at least two signals in $c$ collide is
    \begin{equation*}
      t'' = \inf_{\substack{m, m' \in M\\ m \neq m'}} \inf_{\substack{s \in c(m)\\ s' \in c(m')}} \inf \setOf{t \in \R_{> 0} \suchThat \ReachOf_{\velocityOf(s) \multipli t}(m) \cap \ReachOf_{\velocityOf(s') \multipli t}(m') \neq \emptyset}.
    \end{equation*}
    The minimum time until the next event(s) in $c$ occurs is \graffito{next event(s) time $t_0(c)$}$t_0(c) = \min\setOf{t', t''}$\index[symbols]{t0c@$t_0(c)$}.
  \end{definition}

  \begin{remark}
    A stationary signal at a vertex does never reach the vertex (it is already there) and two signals that already are at the same vertex do never collide (they may have collided when they got there but now they just are at the same vertex; if they have a non-zero velocity, then they will leave the vertex without interfering each other, and, if they have the same positive velocity, then they will travel alongside each other).
  \end{remark}

  \begin{remark}
    The next event time may be $0$, which means that events accumulate at time $0$, or $\infty$, which means that there are no more events in the future; note that $\inf\emptyset = \infty$. It is for example $0$ if there is a sequence of signals moving at the same velocity towards the same vertex each one being already a little closer to the vertex than the previous one. And it is for example $\infty$ if there are no signals at all or there are only stationary signals.
  \end{remark}

  If each signal moves with its velocity and upon reaching a vertex propagates to all incident edges except the one it came from (making copies of itself if necessary) and if collisions of signals are ignored, then the set of points an event occurs in at a time in the future is given a name in 

  \begin{definition} 
    Let $c$ be a configuration of $\Configurations$ and let $t$ be a time of $\Times$. The set of vertices that signals in $c$ reach at time $t$ (under the assumptions given in the introduction to the present definition) is
    \begin{equation*}
      M' = \begin{dcases*}
             \emptyset, &if $t = 0$,\\ 
             \bigcup_{\substack{m \in M\\ s \in c(m)\\ \speedOf(s) > 0}} \ReachOf_{\velocityOf(s) \multipli t}(m) \cap \continuumVertices, &otherwise.
           \end{dcases*}
    \end{equation*}
    The set of points that signals in $c$ collide in at time $t$ is
    \begin{equation*}
      M'' = \bigcup_{\substack{m, m' \in M\\ m \neq m'}} \bigcup_{\substack{s \in c(m)\\ s' \in c(m')}} \ReachOf_{\velocityOf(s) \multipli t}(m) \cap \ReachOf_{\velocityOf(s') \multipli t}(m').
    \end{equation*}
    The set of points that signals in $c$ are involved in an event in at time $t$ is $M_t(c) = M' \cup M''$\graffito{set $M_t(c)$ of points that signals in $c$ are involved in an event in at time $t$}\index[symbols]{Mtc@$M_t(c)$}.
  \end{definition}

  \begin{remark}
    For times $t$ before and including the time $t_0(c)$ of the next event, the definition of $M_t(c)$ is natural. And for other times, it is plausible with the explanation given before the definition and it is used to handle accumulations of events and accumulations of accumulations of events and so on.
  \end{remark}

  As above, if each signal moves with its velocity and upon reaching a vertex propagates to all incident edges except the one it came from (making copies of itself if necessary) and if collisions of signals are ignored, then starting our signal machine in a configuration $c$ and letting it run for a time $t$ without handling propagation of signals in vertices at time $t$ yields a new configuration $\gtfIgnoreCollisions(t)(c)$ as defined in 

  \begin{definition}
    \begin{align*}
      \gtfIgnoreCollisions \from \Times &\to (\Configurations \to \Configurations), \mathnote{map $\gtfIgnoreCollisions$ from $\Times$ to $\Configurations \to \Configurations$}\index[symbols]{boxplus@$\gtfIgnoreCollisions$}\\
      0 &\mapsto \identityMap_{\Configurations},\\ 
      t &\mapsto [c \mapsto [m \mapsto \begin{aligned}[t]
                                         \setOf{s \in \Signals \suchThat{} &\Exists m' \in M \Exists s' \in c(m') \SuchThat\\
                                                                           &(m, \directionOf(s)) \in \ReachOf_{\velocityOf(s') \multipli t}^{\directed}(m') \text{, }\\
                                                                           &\kindOf(s) = \kindOf(s') \text{, and }\\
                                                                           &\datumOf(s) = \datumOf(s')}]]. \qedhere
                                       \end{aligned}
    \end{align*}
  \end{definition}

  \begin{remark}
    For each time $t \in \Times$ and each configuration $c \in \Configurations$, if there is a signal in $c$ that reaches a vertex in time $t$ (or one of its duplicates does), then in the configuration $\gtfIgnoreCollisions(t)(c)$ that signal is at the vertex and its velocity is the one it had just before reaching the vertex. The direction of that velocity points away from the edge the signal came from and it does not point to any edge that is incident to the vertex. It is then up to the local transition function to decide what to do with the signal and, if it is not removed, what its direction shall be.
  \end{remark}

  \begin{remark}
    The map $\gtfIgnoreCollisions$ is used to determine future configurations before or right until the next event occurs and needs to be handled, and also to make crude predictions of future configurations beyond the next event time by propagating at vertices and ignoring collisions as explained above. These predictions will be used to handle accumulations of events and accumulations of accumulations of events and so on.
  \end{remark}

  Until the first events occur, signals move along edges without colliding. At the time the first events occur, a signal reached a vertex or two signals collided (on a vertex or edge) or multiple such events happened. An event in a vertex is handled by the local transition function $\localTransitionFunction_v$ and on an edge by $\localTransitionFunction_e$. This global behaviour can be described by a map that maps a configuration to the configuration right after the first events occurred and have been handled. This map is given in 

  \begin{definition}
      \index[symbols]{Deltadot0@$\gtfJump_0$}\begin{align*}
        \gtfJump_0 \from \Configurations &\to \Configurations, \mathnote{map $\gtfJump_0$ from $\Configurations$ to $\Configurations$}\\
        c &\mapsto
          \left\{
            \begin{aligned}
              &c, \text{ if $t_0(c) \in \setOf{0, \infty}$},\\
              &[m \mapsto
                \left\{
                  \begin{aligned}
                    &c'(m), &&\text{if $m \notin M_{t_0(c)}(c)$,}\\
                    &\localTransitionFunction_v(\directionOf(m), c'(m)), &&\text{if $m \in M_{t_0(c)}(c) \cap \continuumVertices$,}\\
                    &\localTransitionFunction_e(c'(m)), &&\text{if $m \in M_{t_0(c)}(c) \smallsetminus \continuumVertices$,}
                  \end{aligned}
                \right\}
              ],\\
              &&\llap{otherwise,\quad}
            \end{aligned}
          \right.\\
          &\phantom{\mapsto{}}\text{ where } c' = \gtfIgnoreCollisions(t_0(c))(c). \qedhere
      \end{align*}
  \end{definition}

  \begin{remark}
    The map $\gtfJump_0$ maps a configuration to itself if the next event time is $0$, which means that event times accumulated at $0$, or $\infty$, which means that there is no next event. And it maps a configuration to the configuration that is reached after the first events have been handled, by first using $\gtfIgnoreCollisions$ to determine the configuration in which the events occur and then handling all occurring events with $\localTransitionFunction_v$ and $\localTransitionFunction_e$.
  \end{remark}

  If the next event time is $0$, which we call \index{singularity}\graffito{singularity of order $-1$}\defineX{singularity of order $-1$}{singularity!of order $-1$} (see \cref{figure:singularity-of-order-minus-one}), then the machine is sometimes stuck, in the sense that there is no natural way to define what configuration the machine is in at any time in the future, and sometimes the machine can go forward in time, in the sense that there is a natural way to define what configuration the machine is in at least until some time in the future; the latter case is handled later and largely ignored for now. If the next event time is $\infty$, then the machine does nothing for eternity.

  Otherwise, the machine can at least proceed until the next event time and handle the occurring events, which we call \graffito{singularity of order $0$}\defineX{singularities of order $0$}{singularity!of order $0$}, and then the next event time may again be $0$, $\infty$, or something in between. It may happen that the next event times are never $0$ or $\infty$ but accumulate at some time in the future, which we call \graffito{singularity of order $j$, for $j \in \N_+$}\defineX{singularity of order $1$}{singularity!of order $j$, for $j \in \N_+$} (see \cref{figure:singularity-of-order-one}). In that case repeated applications of $\gtfJump_0$ never reach a configuration at that future time or a time beyond. But we can in a sense take the limit of the sequence of configurations that repeated applications of $\gtfJump_0$ yield. Yet, it may even happen that singularities of order $1$ accumulate at some time in the future, which we call \defineX{singularity of order $2$}{singularity!of order $j$, for $j \in \N_+$}. Again, we can in the same sense as before take the limit of the sequence of configurations at these singularities. It may continue this way ad infinitum. In precise terms this is done in

  \begin{figure}
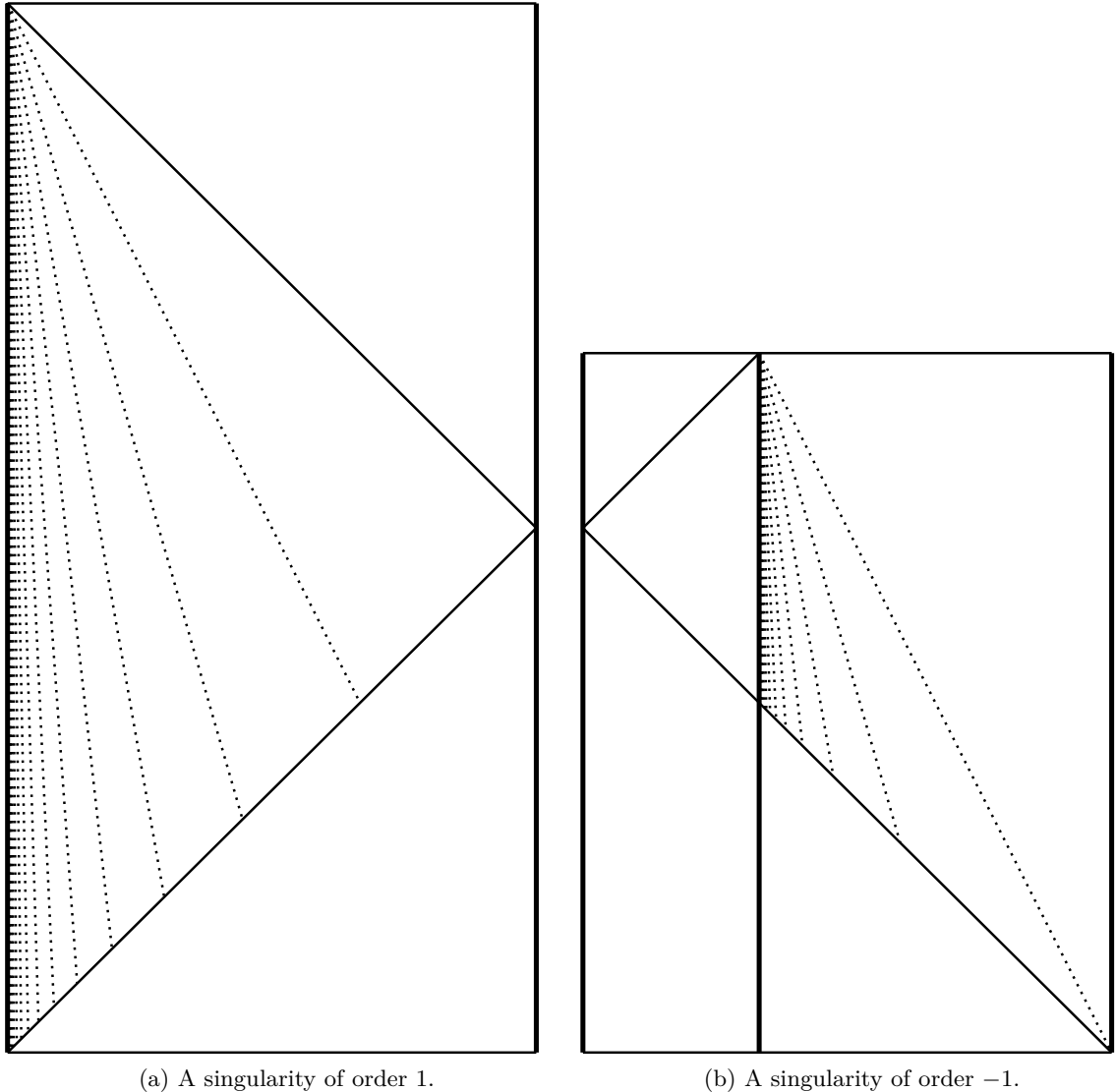

    \myfloatalign
    \begin{wide}
      \figureSingularities
      \caption{Both subfigures depict a space-time diagram of the time evolution of an unspecified signal machine, where space is drawn on the horizontal axis and time on the vertical axis evolving from top to bottom. In \cref{figure:singularity-of-order-one}, there are infinitely many signals of various speeds arbitrarily close to $0$ and slower than $1$ emanating from the left vertex, the fastest signal, which is the one of speed $1$, is reflected at the right vertex, and the reflected signal collides with all other signals at shorter and shorter time spans between collisions, resulting in a singularity of order $1$ at the last depicted time. In \cref{figure:singularity-of-order-minus-one}, there are infinitely many signals of various speeds arbitrarily close to $0$ and slower than $1$ emanating to the right from the middle vertex, and there is one signal of speed $1$ emanating to the left from the middle vertex, this signal is reflected at the left vertex, and at the time the reflected signal reaches the middle vertex is a singularity of order $-1$ because it collides with all signals that emanated to the right and are arbitrarily close to the vertex.}
      \label{figure:singularities}
    \end{wide}
  \end{figure}

  \begin{definition}
    The sequence $\sequence{t_{j - 1}^n}_{n \in \N_0}$, where the $n$ in $t_{j - 1}^n$ is an upper index and does not stand for exponentiation, the map $t_j$, and the map $\gtfJump_j$, for $j \in \N_+$, are defined by mutual induction as follows: The maps $t_0$ and $\gtfJump_0$ have already been defined and, for each positive integer $j$, let
    \begin{equation*}
      \sequence*{
        \begin{aligned}
          t_{j - 1}^n \from \Configurations &\to \extendedTimes, \index[symbols]{tj-1n@$t_{j - 1}^n$}\\ 
          c &\mapsto \sum_{i = 0}^{n - 1} t_{j - 1}(\gtfJump_{j - 1}^i(c)),
        \end{aligned}
      }_{n \in \N_0} \mathnote{map $t_{j - 1}^n$ from $\Configurations$ to $\extendedTimes$}
    \end{equation*} 
    (note that $t_{j - 1}^0 = 0$), let
    \begin{align*}
      t_j \from \Configurations &\to \extendedTimes, \mathnote{map $t_j$ from $\Configurations$ to $\extendedTimes$}\index[symbols]{tjsubscript@$t_j$}\\ 
      c &\mapsto \lim_{n \to \infty} t_{j - 1}^n(c),
    \end{align*}
    and let\index[symbols]{Deltadotj@$\gtfJump_j$}
    \begin{align*}
      \gtfJump_j \from \Configurations &\to \Configurations, \mathnote{map $\gtfJump_j$ from $\Configurations$ to $\Configurations$}\\
      c &\mapsto
        \left\{
          \begin{aligned}
            &c, \text{ if $t_j(c) \in \setOf{0, \infty}$},\\
            &[m \mapsto
              \left\{
                \begin{aligned}
                  &c'(m), &&\text{if $m \notin M_j(c)$,}\\
                  &\localTransitionFunction_v(\directionOf(m), c'(m)), &&\text{if $m \in M_j(c) \cap \continuumVertices$,}\\
                  &\localTransitionFunction_e(c'(m)), &&\text{if $m \in M_j(c) \smallsetminus \continuumVertices$,}
                \end{aligned}
              \right\}
            ],\\
            &&\llap{otherwise,\quad}
          \end{aligned}
        \right.\\
        &\phantom{\mapsto{}}\text{ where } c' = \liminf_{n \to \infty} \gtfIgnoreCollisions(t_j(c) - t_{j - 1}^n(c))(\gtfJump_{j - 1}^n(c)),\\
        &\phantom{\mapsto{}}\text{ and } M_j(c) = \liminf_{n \to \infty} M_{t_j(c) - t_{j - 1}^n(c)}(\gtfJump_{j - 1}^n(c)),
    \end{align*}
    where the first limit inferior is the pointwise limit inferior of sequences of set-valued maps and the second limit inferior is the limit inferior of sequences of sets. In greater detail, for each sequence $\sequence{c_n}_{n \in \N_0}$ of set-valued maps from $M$ to $Q$, the pointwise limit inferior of $\sequence{c_n}_{n \in \N_0}$ is the map $c \from M \to Q$, $m \mapsto \liminf_{n \to \infty} c_n(m)$ and is denoted by $\liminf_{n \to \infty} c_n$; and, for each sequence $\sequence{A_n}_{n \in \N_0}$ of subsets of $M$, the limit inferior of $\sequence{A_n}_{n \in \N_0}$ is the subset $\bigcup_{n \in \N_0} \bigcap_{k \geq n} A_k$ of $M$ and is denoted by $\liminf_{n \to \infty} A_n$.
  \end{definition} 

  \begin{remark}
    Let $c$ be a configuration of $\Configurations$. Then, for each positive integer $j$, we have $t_{j - 1}^0(c) = 0$ and $t_{j - 1}^1(c) = t_{j - 1}(c)$, and the sequence $\sequence{t_{j - 1}^n(c)}_{n \in \N_0}$ in $\extendedTimes$ is non-decreasing and hence converges in $\extendedTimes$. And, the sequence $\sequence{t_j(c)}_{j \in \N_0}$ in $\extendedTimes$ is non-decreasing and hence converges in $\extendedTimes$.
  \end{remark}

  \begin{remark} 
    Let the signal machine be in a configuration $c$ at time $0$ and let there be no future configuration whose next event time is $0$ or $\infty$. The latter is the case if and only if the sequences $\sequence{t_{j - 1}^n(c)}_{n \in \N_0}$, for $j \in \N_+$, and hence also the sequence $\sequence{t_j(c)}_{j \in \N_0}$, are strictly increasing sequences in $\Times$.

    Then, for each positive integer $j$ and each non-negative integer $n$, at time $t_{j - 1}^n(c)$ the machine is in configuration $\gtfJump_{j - 1}^n(c)$. And, the time $t_0^n(c)$ is the $n$-th time an event occurs (singularity of order $0$), the time $t_1^n(c)$ is the $n$-th time an accumulation of events occurs (singularity of order $1$), the time $t_2^n(c)$ is the $n$-th time an accumulation of accumulation of events occurs (singularity of order $2$), and so forth.

    Moreover, for each non-negative integer $j$, at time $t_j(c)$ the machine is in configuration $\gtfJump_j(c)$. And, the time $t_0(c)$ is the next time an event occurs, the time $t_1(c)$ is the next time an accumulation of events occurs, the time $t_2(c)$ is the next time an accumulation of accumulations of events occurs, and so forth. 

    Furthermore, for each positive integer $j$, the map $\gtfJump_j$ maps the configuration $c$ to the configuration that is reached after an accumulation of singularities of order $j - 1$, which is a singularity of order $j$. First, it calculates the configurations that accumulate, namely $\gtfJump_{j - 1}^n(c)$; secondly, for each of these configurations, it uses $\gtfIgnoreCollisions$ to determine the configuration that would be reached at the accumulation time if there were no further events, which is a crude prediction of the future that becomes better the greater $n$ is; thirdly, it calculates the pointwise limit inferior of these configurations, which is essentially the configuration that contains the signals that all but finitely many of the configurations have in common (in particular, if for a point $m$ the sequence of signals at $m$ become constant, then the limit at $m$ is that set of signals); lastly, it handles collisions.
  \end{remark}

  \begin{remark}
  \label{remark:limit-of-jumps-is-singularity-of-order-minus-1}
    Let the signal machine be in a configuration $c$ at time $0$ and let there be a future configuration whose next event time is $0$. Then, there is a least positive integer $j$ such that the sequence $\sequence{t_{j - 1}^n(c)}_{n \in \N_0}$ is eventually constant. And, there is a least non-negative integer $n$ such that $t_{j - 1}^n(c) = t_{j - 1}^{n + 1}(c)$. The time $t' = t_{j - 1}^n(c)$ is the first time at which the signal machine is in a configuration whose next event time is $0$ and this configuration is $c' = \gtfJump_{j - 1}^n(c)$.

    For each non-negative integer $n'$ such that $n' \geq n$, we have $t_{j - 1}^{n'}(c) = t'$ and $\gtfJump_{j - 1}^{n'}(c) = c'$. And, for each positive integer $j'$ such that $j' \geq j$, the time $t_{j'}(c)$ is equal to $t'$ and the configuration $\gtfJump_{j'}(c)$ is equal to $c'$, and the sequences $\sequence{t_{j'}^n(c)}_{n \in \N_0}$ and $\sequence{\gtfJump_{j'}^n(c)}_{n \in \N_0}$ are the constant sequences $\sequence{t'}_{n \in \N_0}$ and $\sequence{c'}_{n \in \N_0}$. In particular, the limit $\lim_{j \to \infty} t_j(c)$ is equal to $t'$ and the limit inferior $\liminf_{j \to \infty} \gtfJump_j(c)$ is equal to $c'$.
  \end{remark}

  \begin{remark}
    The limit of sequences of configurations and of sets of points does in general not exist. However, the limit inferior and the limit superior always exist. We decided not to use the limit, to avoid case distinctions that would have to be made. Instead, we decided to use the limit inferior; we could as well have decided to use the limit superior. Which of the two has the desired outcome depends on the specific use case.

    For the signal machine that solves the firing squad synchronisation problem that we construct in the next section and the configurations it is initialised with and the configurations it encounters during its time evolution, the encountered limit inferiors and superiors are actually always the same, which means that the limits exist, and hence the choice of limit inferior or superior is irrelevant in that use case.
  \end{remark}

  If the machine never assumes a configuration in which events accumulate at time $0$ and if non-negative singularities of ever higher orders ad infinitum do not accumulate, then the machine can be observed for eternity. Otherwise, it can for now only be observed for all times before $\lim_{j \to \infty} t_j(c)$, where $c$ is the initial configuration of the machine. This time is given a name in 

  \begin{definition}
    For each configuration $c \in \Configurations$, the non-negative real number or infinity $t_\infty(c) = \lim_{j \to \infty} t_j(c)$ is called \defineX{$\infty$-existence time of $c$}{existence time infinity@$\infty$-existence time of $c$}\graffito{$\infty$-existence time $t_\infty(c)$ of $c$}\index[symbols]{tinftyc@$t_\infty(c)$}, and the closed interval $\closedInterval{0, t_\infty(c)}$ is called \defineX{$\infty$-existence interval of $c$}{existence interval infinity@$\infty$-existence interval of $c$}\graffito{$\infty$-existence interval $\closedInterval{0, t_\infty(c)}$ of $c$}.
  \end{definition}

  Repeated applications of powers of the maps $\gtfJump_j$, for decreasing $j \in \N_0$, let us jump to and from configurations right after singularities of decreasing orders. The resulting map is given a name in

  \begin{definition}
    For each non-negative integer $j$, each non-negative integer $k$, and each finite sequence $\sequence{n_i}_{i = j}^k$ of non-negative integers, let
    \begin{equation*} 
      \gtfJump_{\sequence{n_i}_{i = j}^k}
      = \begin{dcases*}
          \identityMap_{\Configurations}, &if $j > k$,\\
          \gtfJump_j^{n_j} \after \gtfJump_{j + 1}^{n_{j + 1}} \after \dotsb \after \gtfJump_k^{n_k}, &otherwise, 
        \end{dcases*}
      \mathnote{map $\gtfJump_{\sequence{n_i}_{i = j}^k}$ from $\Configurations$ to $\Configurations$}
      \index[symbols]{Deltadotniijksequence@$\gtfJump_{\sequence{n_i}_{i = j}^k}$}
    \end{equation*}
    and let
    \begin{equation*}
      t_{\sequence{n_i}_{i = j}^k} = \sum_{i = j}^k t_i^{n_i} \after \gtfJump_{\sequence{n_\ell}_{\ell = i + 1}^k}.
      \mathnote{map $t_{\sequence{n_i}_{i = j}^k}$ from $\Configurations$ to $\extendedTimes$}
      \index[symbols]{tniijksequence@$t_{\sequence{n_i}_{i = j}^k}$} \qedhere
    \end{equation*}
  \end{definition}

  \begin{remark}
    Let the signal machine be in a configuration $c$ at time $0$ and let there be no future configuration whose next event time is $0$ or $\infty$. The map $\gtfJump_{\sequence{n_i}_{i = j}^k}$ applied to $c$, first applies $\gtfJump_k^{n_k}$ to jump from $c$ to the configuration right after the $n_k$-th time a singularity of order $k$ occurs, secondly it applies $\gtfJump_{k - 1}^{n_{k - 1}}$ to jump from that configuration, namely $\gtfJump_{\sequence{n_i}_{i = k}^k}(c)$, to the configuration right after the $n_{k - 1}$-th time a singularity of order $k - 1$ occurs (counting from the time at which the machine is in configuration $\gtfJump_{\sequence{n_i}_{i = k}^k}(c)$), and so forth until it finally applies $\gtfJump_j^{n_j}$ to jump from the configuration $\gtfJump_{\sequence{n_i}_{i = j + 1}^k}(c)$ to the configuration right after the $n_j$-th time a singularity of order $j$ occurs (counting from the time at which the machine is in configuration $\gtfJump_{\sequence{n_i}_{i = j + 1}^k}(c)$), where in the case that $j = 0$, a singularity of order $0$ is nothing but an event.

    The time it takes the machine to get from $c$ to the configuration $\gtfJump_k^{n_k}(c)$ is $t_k^{n_k}(c)$, the time it takes to get from $\gtfJump_k^{n_k}(c)$ to $\gtfJump_{k - 1}^{n_{k - 1}}(\gtfJump_k^{n_k}(c))$ is $t_{k - 1}^{n_{k - 1}}(\gtfJump_k^{n_k}(c))$, and so forth; in total, the time it takes to get from $c$ to $\gtfJump_{\sequence{n_i}_{i = j}^k}(c)$ is $t_{\sequence{n_i}_{i = j}^k}(c)$.
  \end{remark}

  To compute the configuration the machine is in at the $\infty$-existence time of the initial configuration, we can use the maps $\gtfJump_j$ for increasing $j$ to jump from singularities of non-negative lower orders to singularities of ever higher orders, which in the case there are any singularities of order $-1$ comes to a halt at the first such singularity at the $\infty$-existence time and in the other case yields in a sense the limit configuration. And, to compute the configuration at time $t$ before the $\infty$-existence time, we can use one of the maps $\gtfJump_{\sequence{n_i}_{i = 0}^k}$ to jump to the configuration right after the last event before time $t$ and then we can use the map $\gtfIgnoreCollisions$ to jump from there to time $t$. The resulting map describes the time evolution of the signal machine before or at $\infty$-existence times and it is given in

  \begin{definition}
    For each time $t \in \extendedTimes$, the set of configurations whose $\infty$-existence interval contains $t$ is
    \begin{equation*}
      \Configurations_t^\infty = \setOf{c \in \Configurations \suchThat t \leq t_\infty(c)}.
      \mathnote{set $\Configurations_t^\infty$} 
      \index[symbols]{Cnftinfinity@$\Configurations_t^\infty$}
    \end{equation*}
    Let
    \begin{align*} 
      \gtfNonNegativeSingularities \from \extendedTimes &\to \bigcup_{t \in \extendedTimes} (\Configurations_t^\infty \to \Configurations), \mathnote{map $\gtfNonNegativeSingularities$ from $\extendedTimes$ to $\bigcup_{t \in \extendedTimes} (\Configurations_t^\infty \to \Configurations)$}\index[symbols]{boxminus@$\gtfNonNegativeSingularities$}\\ 
      t &\mapsto
        \left[
          \begin{aligned}
            \Configurations_t^\infty &\to \Configurations,\\
            c &\mapsto
              \begin{dcases*}
                \liminf_{j \to \infty} \gtfJump_j(c), &if $t = t_\infty(c)$,\\ 
                \gtfIgnoreCollisions(t - t_{\sequence{n_i}_{i = 0}^k}(c))(\gtfJump_{\sequence{n_i}_{i = 0}^k}(c)), &otherwise,
              \end{dcases*}\\
              &\phantom{\mapsto{}}\text{ for the least } k \in \N_0 \text{ and the } \sequence{n_i}_{i = 0}^k \in \N_0^{k + 1}\\
              &\phantom{\mapsto{}}\text{ with } t \in \leftClosedAndRightOpenInterval{t_{\sequence{n_i}_{i = 0}^k}(c), t_{(n_0 + 1, n_1, n_2, \dotsc, n_k)}(c)}.
          \end{aligned}
        \right]
    \end{align*} 
    Note that the least index $k$ and the finite sequence $\sequence{n_i}_{i = 0}^k$ that occur above are uniquely determined by the time $t$ and the configuration $c$. 
  \end{definition}

  \begin{remark}
    In the second case in the definition of $\gtfNonNegativeSingularities$, we have $t \in \leftClosedAndRightOpenInterval{t_k^{n_k}(c), t_k^{n_k + 1}(c)}$, and $t - t_k^{n_k}(c) \in \leftClosedAndRightOpenInterval{t_{k - 1}^{n_{k - 1}}(\gtfJump_k^{n_k}(c)), t_{k - 1}^{n_{k - 1} + 1}(\gtfJump_k^{n_k}(c))}$, and so forth.
  \end{remark}

  \begin{remark}
    For each configuration $c \in \Configurations$, the map $\gtfNonNegativeSingularities(\blank)(c)$ is defined on the closed interval $\closedInterval{0, t_\infty(c)}$.
  \end{remark}

  \begin{remark}
    Let $t$ be a time of $\Times$ and let $c$ be a configuration of $\Configurations_t^\infty$ such that $t \neq t_\infty(c)$. If events do not accumulate for the initial configuration $c$, then $\gtfNonNegativeSingularities(t)(c) = \gtfIgnoreCollisions(t - t_0^{n_0}(c))(\gtfJump_0^{n_0}(c))$, for the $n_0 \in \N_0$ with $t \in \leftClosedAndRightOpenInterval{t_0^{n_0}(c), t_0^{n_0 + 1}(c)}$; in words, we apply $\gtfJump_0$ repeatedly, jumping from event to event, until we reach the configuration $\gtfJump_0^{n_0}(c)$ at time $t_0^{n_0}(c)$ with the property that the next event (if there even is one) occurs after $t$, at which point we use $\gtfIgnoreCollisions$ to move signals along edges for the remaining time $t - t_0^{n_0}(c)$.

    If events do accumulate for $c$ but singularities of order $1$ do not accumulate, then we apply $\gtfJump_1$ repeatedly, jumping from singularity to singularity, until we reach the configuration $\gtfJump_1^{n_1}(c)$ at time $t_1^{n_1}(c)$ with the property that the next singularity (if there even is one) occurs after $t$, at which point we apply $\gtfJump_0$ repeatedly, jumping from event to event, until we reach the configuration $\gtfJump_0^{n_0}(\gtfJump_1^{n_1}(c))$ at time $t_0^{n_0}(c) + t_1^{n_1}(c)$ with the property that the next event (if there even is one) occurs after $t$, at which point we use $\gtfIgnoreCollisions$ to move signals along edges for the remaining time $t - t_0^{n_0}(c) - t_1^{n_1}(c)$.

    And so forth.
  \end{remark}

  \begin{remark}
  \label{remark:gtf-at-infinity-existence-time-sometimes-yields-the-configuration-at-a-singulairty-of-order-minus-1}
    Let the signal machine be in a configuration $c$ at time $0$. If there is a singularity of order $-1$ in the future, then, according to \cref{remark:limit-of-jumps-is-singularity-of-order-minus-1}, the time $t_\infty(c)$ is the first time at which the signal machine is at such a singularity and the corresponding configuration is $\gtfNonNegativeSingularities(t_\infty(c))(c)$. Otherwise, the time $t_\infty(c)$, which may be the improper time $\infty$, is the time just after all singularities and the corresponding configuration is $\gtfNonNegativeSingularities(t_\infty(c))(c)$. In either case, in what is to come, for simplicity, if $t_\infty(c)$ is finite, then we talk as if there is a singularity of order $-1$ at time $t_\infty(c)$. 
  \end{remark}


  At the time of a singularity of order $1$, there are infinitely many events that occur just \emph{before} that time and arbitrarily close to it, and the problem is to define the configuration at the time of the singularity. At the time of a singularity of order $-1$, there are infinitely many events that occur just \emph{after} that time and arbitrarily close to it, and the problem is to define the configurations at all times after the singularity. Analogous problems exist for accumulations of singularities of order $1$ or $-1$, and accumulations of accumulations of singularities of order $1$ or $-1$, and so forth. For singularities of positive orders these problems have been solved above but not for singularities of order $-1$ and its accumulations.

  At a singularity of order $-1$ at time $0$, to compute the configuration at a small enough time $t$, first, we make crude predictions of the future with $\gtfIgnoreCollisions$ by jumping past the singularity to future times $\varepsilon$ ignoring events, secondly, we extrapolate these predictions to the time $t$ with $\gtfNonNegativeSingularities$ by letting the machine evolve them until the time $t$, and, lastly, we take the limit inferior of these predictions as $\varepsilon$ tends to $0$. This does not work for all singularities of order $-1$ regardless of how small we choose $t$ (for example if at each point in time of a time span after and including $0$ an event occurs), but it does work for the singularities of order $-1$ that occur in our quasi-solution of the firing squad synchronisation problem. 

  If the $\infty$-existence time of the current configuration is $\infty$, then the machine can be observed for eternity. Otherwise, it can be observed until the $\infty$-existence time using $\gtfNonNegativeSingularities$ and from there at least until the least time until which all the crude predictions mentioned above can be observed; this time span is named in

  \begin{definition}
    \begin{align*}
      t_{-1} \from \Configurations &\to \extendedTimes, \mathnote{map $t_{-1}$ from $\Configurations$ to $\extendedTimes$}\index[symbols]{t-1@$t_{-1}$}\\
      c &\mapsto
        \begin{dcases*}
          \infty, &if $t_\infty(c) = \infty$,\\
          \inf_{\varepsilon \in \R_{> 0}}((\varepsilon + t_\infty(c)) + t_\infty(c_\varepsilon)), &otherwise,
        \end{dcases*}\\
        &\phantom{\mapsto{}}\text{ where } c_\varepsilon = \gtfIgnoreCollisions(\varepsilon)(\gtfNonNegativeSingularities(t_\infty(c))(c)). \qedhere
    \end{align*}
  \end{definition}

  How the machine evolves until and beyond $t_\infty$ and at most until $t_{-1}$ is given in

  \begin{definition}
    For each time $t \in \extendedTimes$, let
    \begin{equation*}
      \Configurations_t^{-1} = \setOf{c \in \Configurations \suchThat t \leq t_{-1}(c)}
      \mathnote{set $\Configurations_t^{-1}$} 
      \index[symbols]{Cnft-1@$\Configurations_t^{-1}$}
    \end{equation*}
    and let
    \begin{align*}
      \gtfNegativeSingularityOfOrderMinusOne \from \extendedTimes &\to \bigcup_{t \in \extendedTimes} (\Configurations_t^{-1} \to \Configurations), \mathnote{map $\gtfNegativeSingularityOfOrderMinusOne$ from $\extendedTimes$ to $\bigcup_{t \in \extendedTimes} (\Configurations_t^{-1} \to \Configurations)$}\index[symbols]{boxast@$\gtfNegativeSingularityOfOrderMinusOne$}\\
      t &\mapsto
        \left[
          \begin{aligned}
            \Configurations_t^{-1} &\to \Configurations,\\
            c &\mapsto
              \begin{dcases*}
                \gtfNonNegativeSingularities(t)(c), &if $t \leq t_\infty(c)$,\\
                \liminf_{\varepsilon \downTo 0} \gtfNonNegativeSingularities(t - (\varepsilon + t_\infty(c)))(c_\varepsilon), &otherwise,
              \end{dcases*}\\
              &\phantom{\mapsto{}}\text{ where } c_\varepsilon = \gtfIgnoreCollisions(\varepsilon)(\gtfNonNegativeSingularities(t_\infty(c))(c)),
          \end{aligned}
        \right]
    \end{align*}
    where the limit inferior is the pointwise limit inferior of set-valued maps, in greater detail, for each family $\family{c_\varepsilon}_{\varepsilon \in \R_{> 0}}$ of set-valued maps from $M$ to $Q$, the pointwise limit inferior $\liminf_{\varepsilon \downTo 0} c_\varepsilon$ is the map $c \from M \to Q$, $m \mapsto \bigcup_{r \in \R_{> 0}} \bigcap_{s \geq r} c_{1/s}(m)$.
  \end{definition}

  \begin{remark}
    For each configuration $c \in \Configurations$, we have $t_\infty(c) \leq t_{-1}(c)$; for each time $t \in \extendedTimes$, we have $\Configurations_t^\infty \subseteq \Configurations_t^{-1}$; and for each configuration $c \in \Configurations$ and each positive real number $\varepsilon$, we have $t_{-1}(c) - (\varepsilon + t_\infty(c)) \leq t_\infty(c_\varepsilon)$.
  \end{remark}

  \begin{remark}
    Let the signal machine be in a configuration $c$ at time $0$. If $t_\infty(\gtfNonNegativeSingularities(t_\infty(c))(c)) = 0$, then there is a singularity of order $-1$ at time $t_\infty(c)$. Otherwise, there is not. In the latter case, as one would hope, $t_{-1}(c) = t_\infty(\gtfNonNegativeSingularities(t_\infty(c))(c))$, and, for each time $t \in \closedInterval{t_\infty(c), t_{-1}(c)}$, we have $\gtfNegativeSingularityOfOrderMinusOne(t)(c) = \gtfNonNegativeSingularities(t - t_\infty(c))(\gtfNonNegativeSingularities(t_\infty(c))(c))$. 
  \end{remark}

  Unlike for a singularity of a non-negative order, for a singularity of order $-1$ it is not possible to jump to the time right after the singularity as there is no such time. However, we may jump over the singularity at time $t_\infty$ to the time $t_{-1}$. This is done by the map 

  \begin{definition}
    \index[symbols]{Deltadot-1@$\gtfJump_{-1}$}\begin{align*}
      \gtfJump_{-1} \from \Configurations &\to \Configurations, \mathnote{map $\gtfJump_{-1}$ from $\Configurations$ to $\Configurations$}\\
      c &\mapsto
        \begin{dcases*}
          c, &if $t_\infty(c) = \infty$,\\
          \gtfNegativeSingularityOfOrderMinusOne(t_{-1}(c))(c), &otherwise. \qedhere
        \end{dcases*}
    \end{align*}
  \end{definition}

  Like for singularities of non-negative orders, such of order $-1$ may accumulate, which is a \index{singularity}\defineX{singularity of order $-2$}{singularity!of order $-j$, for $j \in \Z_{\geq 2}$}\graffito{singularity of order $-j$, for $j \in \Z_{\geq 2}$}, such of order $-2$ my accumulate, which is a \defineX{singularity of order $-3$}{singularity!of order $-j$, for $j \in \Z_{\geq 2}$}, and so forth. The map to jump over singularities of order $-1$ has already been given and the maps to jump over singularities of smaller orders are introduced in

  \begin{definition}
    The sequence $\sequence{t_{-(j - 1)}^n}_{n \in \N_0}$, where the $n$ in $t_{-(j - 1)}^n$ is an upper index and does not stand for exponentiation, the map $t_{-j}$, and the map $\gtfJump_{-j}$, for $j \in \Z_{\geq 2}$, are defined by mutual induction as follows: The maps $t_{-1}$ and $\gtfJump_{-1}$ have already been defined and, for each integer $j$ such that $j \geq 2$, let
    \begin{equation*}
      \sequence*{
        \begin{aligned}
          t_{-(j - 1)}^n \from \Configurations &\to \extendedTimes, \index[symbols]{tj-1nz@$t_{-(j - 1)}^n$}\\
          c &\mapsto \sum_{i = 0}^{n - 1} t_{-(j - 1)}(\gtfJump_{-(j - 1)}^i(c)),
        \end{aligned}
      }_{n \in \N_0} \mathnote{map $t_{-(j - 1)}^n$ from $\Configurations$ to $\extendedTimes$}
    \end{equation*}
    (note that $t_{-(j - 1)}^0 = 0$), let
    \begin{align*}
      t_{-j} \from \Configurations &\to \extendedTimes, \mathnote{map $t_{-j}$ from $\Configurations$ to $\extendedTimes$}\index[symbols]{tjsubscriptz@$t_{-j}$}\\
      c &\mapsto \lim_{n \to \infty} t_{-(j - 1)}^n(c),
    \end{align*}
    and let\index[symbols]{Deltadotj-@$\gtfJump_{-j}$}
    \begin{align*}
      \gtfJump_{-j} \from \Configurations &\to \Configurations, \mathnote{map $\gtfJump_{-j}$ from $\Configurations$ to $\Configurations$}\\
      c &\mapsto
        \left\{
          \begin{aligned}
            &c, \text{ if $t_{-j}(c) = \infty$},\\
            &\liminf_{n \to \infty} \gtfIgnoreCollisions(t_{-j}(c) - t_{-(j - 1)}^n(c))(\gtfJump_{-(j - 1)}^n(c)),\\
            &&\llap{otherwise.} \qedhere
          \end{aligned}
        \right.
    \end{align*}
  \end{definition}

  The machine can be observed for all times before $\lim_{j \to \infty} t_{-j}(c)$, where $c$ is the initial configuration of the machine. This time is given a name in

  \begin{definition}
    For each configuration $c \in \Configurations$, the non-negative real number or infinity $t_{-\infty}(c) = \lim_{j \to \infty} t_{-j}(c)$ is called \graffito{$(-\infty)$-existence time $t_{-\infty}(c)$ of $c$}\defineX{$(-\infty)$-existence time of $c$}{existence time infinity minus@$(-\infty)$-existence time of $c$}\index[symbols]{tinftycminus@$t_{-\infty}(c)$}, and the closed interval $\closedInterval{0, t_{-\infty}(c)}$ is called \graffito{$(-\infty)$-existence interval $\closedInterval{0, t_{-\infty}(c)}$ of $c$}\defineX{$(-\infty)$-existence interval of $c$}{existence interval infinity minus@$(-\infty)$-existence interval of $c$}.
  \end{definition}

  Like for singularities of non-negative orders, repeated applications of powers of the maps $\gtfJump_{-j}$, for decreasing $j \in \N_+$, let us jump over singularities of decreasing negative orders down to order $-1$. The resulting map is given a name in

  \begin{definition}
    For each positive integer $j$, each non-negative integer $k$, and each finite sequence $\sequence{n_i}_{i = -j}^{-k}$ of non-negative integers that is indexed from $-j$ down to $-k$, let
    \begin{equation*}
      \gtfJump_{\sequence{n_i}_{i = -j}^{-k}}
      = \begin{dcases*}
          \identityMap_{\Configurations}, &if $-j < -k$,\\
          \gtfJump_{-j}^{n_{-j}} \after \gtfJump_{-j - 1}^{n_{-j - 1}} \after \dotsb \after \gtfJump_{-k}^{n_{-k}}, &otherwise, 
        \end{dcases*}
      \mathnote{map $\gtfJump_{\sequence{n_i}_{i = -j}^{-k}}$ from $\Configurations$ to $\Configurations$}
      \index[symbols]{Deltadotniijksequenceminus@$\gtfJump_{\sequence{n_i}_{i = -j}^{-k}}$}
    \end{equation*}
    and let
    \begin{equation*}
      t_{\sequence{n_i}_{i = -j}^{-k}} = \sum_{i = -j}^{-k} t_i^{n_i} \after \gtfJump_{\sequence{n_\ell}_{\ell = i - 1}^{-k}}. \qedhere
      \mathnote{map $t_{\sequence{n_i}_{i = -j}^{-k}}$ from $\Configurations$ to $\extendedTimes$}
      \index[symbols]{tniijksequenceminus@$t_{\sequence{n_i}_{i = -j}^{-k}}$}
    \end{equation*}
  \end{definition}

  Like for singularities of non-negative orders, to compute the configuration the machine is in at the $(-\infty)$-existence time of the initial configuration, we can use the maps $\gtfJump_{-j}$ for increasing $j$ to jump from singularities of negative lower orders to singularities of ever higher orders ad infinitum. And, to compute the configuration at time $t$ before the $(-\infty)$-existence time, we can use one of the maps $\gtfJump_{\sequence{n_i}_{i = -1}^{-k}}$ to jump over all the singularities before time $t$ such that the next jump over a singularity of order $-1$ would be beyond time $t$ and then we can use the map $\gtfNegativeSingularityOfOrderMinusOne$ to jump from there to time $t$ whereby we may cross a singularity of order $-1$. The resulting map describes the time evolution of the signal machine before or at $(-\infty)$-existence times, which for the purposes of this treatise is the complete time evolution, and this map is given in

  \begin{definition}
    For each time $t \in \extendedTimes$, the set of configurations whose $(-\infty)$-existence interval contains $t$ is
    \begin{equation*}
      \Configurations_t^{-\infty} = \setOf{c \in \Configurations \suchThat t \leq t_{-\infty}(c)}.
      \mathnote{set $\Configurations_t^{-\infty}$} 
      \index[symbols]{Cnftinfinityminus@$\Configurations_t^{-\infty}$}
    \end{equation*}
    The map
    \begin{align*}
      \globalTransitionFunction \from \extendedTimes &\to \bigcup_{t \in \extendedTimes} (\Configurations_t^{-\infty} \to \Configurations), \mathnote{global transition function $\globalTransitionFunction$ from $\extendedTimes$ to $\bigcup_{t \in \extendedTimes} (\Configurations_t^{-\infty} \to \Configurations)$}\index[symbols]{boxdot@$\globalTransitionFunction$}\\ 
      t &\mapsto
        \left[
          \begin{aligned}
            \Configurations_t^{-\infty} &\to \Configurations,\\
            c &\mapsto
              \left\{
                \begin{aligned}
                  &\liminf_{j \to \infty} \gtfJump_{-j}(c), \text{ if $t = t_{-\infty}(c)$},\\
                  &\gtfNegativeSingularityOfOrderMinusOne(t - t_{\sequence{n_i}_{i = -1}^{-k}}(c))(\gtfJump_{\sequence{n_i}_{i = -1}^{-k}}(c)), \text{ otherwise},
                \end{aligned}
              \right.\\
              &\phantom{\mapsto{}}\text{ for the least } k \in \N_+ \text{ and the } \sequence{n_i}_{i = -1}^{-k} \in \N_0^k\\
              &\phantom{\mapsto{}}\text{ with } t \in \leftClosedAndRightOpenInterval{t_{\sequence{n_i}_{i = -1}^{-k}}(c), t_{(n_{-1} + 1, n_{-2}, n_{-3}, \dotsc, n_{-k})}(c)},
          \end{aligned}
        \right]
    \end{align*}
    is called \define{global transition function}.
  \end{definition} 

  \begin{remark}
    Jérôme Olivier Durand-Lose introduces and studies another notion of signal machines in his paper \enquote{\citetitle*{durand-lose:2008}}\cite{durand-lose:2008}. The notable differences are the following: While our machines are defined over continuum graphs, his machines are defined over the real number line; while our machines may have infinitely many different kinds of signals, his machines may only have signals of a finite number of kinds (which he calls meta-signals); while in a configuration of our machines there may exist infinitely many signals (even at the same point), in configurations of his machines there may only exist finitely many signals; while the time evolution of our machines can be observed beyond singularities of any order, the time evolution of his machines already stops before singularities of order $1$, that is, before accumulations of collisions (as there are no vertices, other events do not exist for his machines).
  \end{remark}


  \section{Firing Squad Signal Machines}
  \label{section:firing-squad-solution}

  In this section, let $\Graph = \ntuple{\Vertices, \Edges, \eendsOf}$ be a non-trivial, finite, and connected undirected multigraph, let $\weightOf$ be an edge weighting of $\Graph$, let $M$ be a continuum representation of $\Graph$, and let $\general$ be a vertex of $M$, which we call \define{general}\graffito{general vertex $\general$}. When we say \emph{path}, we either mean \emph{path in $\Graph$} or \emph{path in $M$}; when we say \emph{longest path}, we either mean \emph{maximum-weight path in $\Graph$} or \emph{longest path in $M$}; and, when we say \emph{path}, we often mean \emph{non-empty direction-preserving path from vertex to vertex}; in all cases it should be clear from the context what is meant.

  From a broad perspective, the signal machine we construct in this section performs the following tasks: It cuts the graph such that the graph turns into a virtual tree; it starts synchronisation of edges as soon as possible and freezes it as late as possible; it determines the midpoints of all non-empty direction-preserving paths from vertices to vertices; it determines which midpoints are the ones of the longest paths; starting from the midpoints of the longest paths, it traverses midpoints of shorter and shorter paths and upon reaching midpoints of edges, it thaws synchronisation of the respective edges; all edges finish synchronisation at the same time with the creation of fire signals that lie dense in the graph (see \cref{figure:determinedMidpointOfOneEdgeOrOfTwoEdges,figure:MazoyerWithOneAndTwoEdges,figure:fsspWithTwoEdgesAndGeneralInBetweenAndAtTheLeft,figure:fsspWithThreeEdgesInARowAndGeneralAtTheSecondVertexFromTheLeft,figure:fsspWithThreeEdgesInARowAndGeneralAtTheSecondVertexFromTheLeft,figure:fsspWithThreeEdgesIncidentToTheSameVertexAndGeneralAtTheVertexAndAtTheLeafOfTheShortestEdge}). A more detailed account is given in

  \begin{figure}
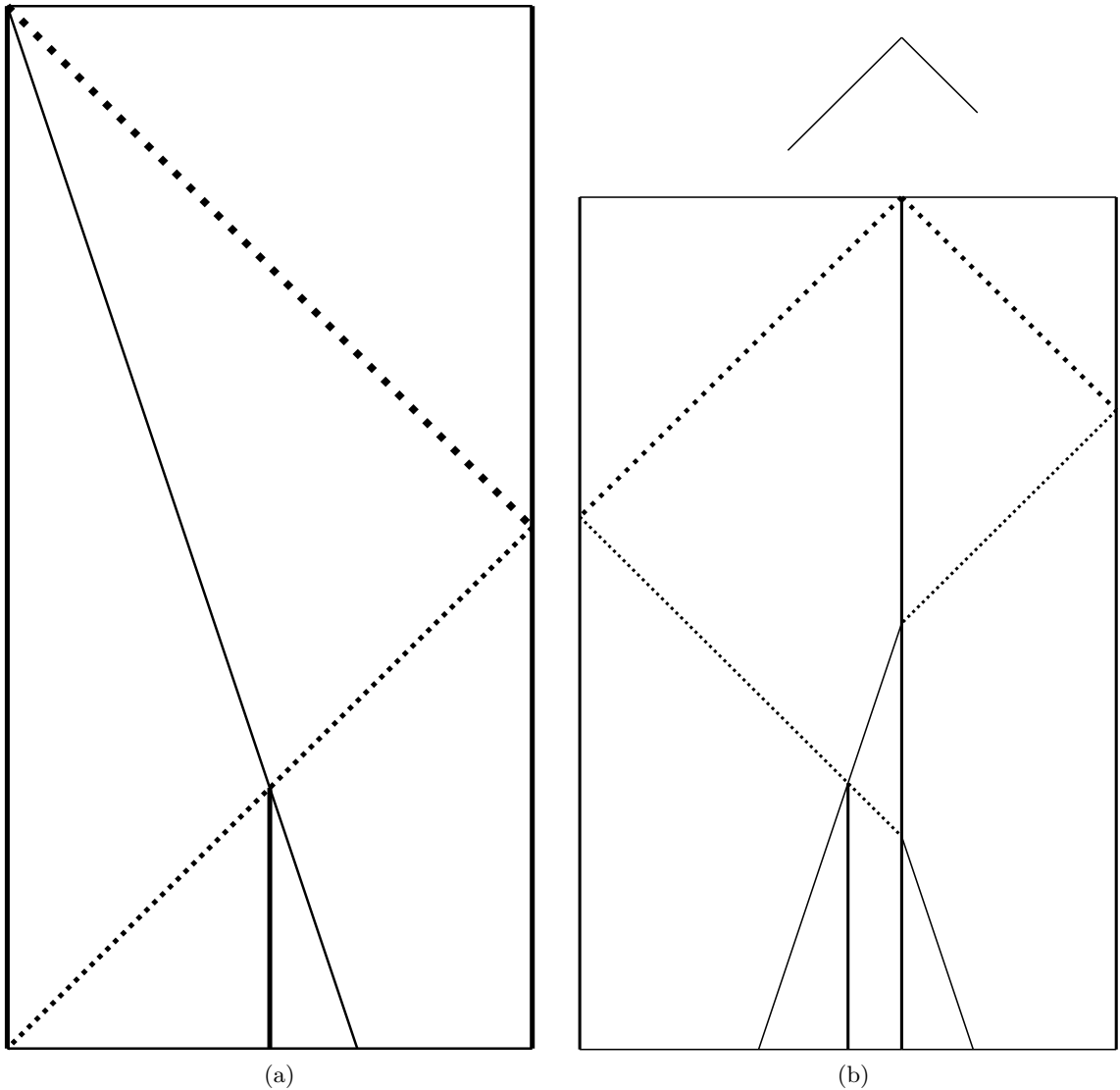

    \begin{wide}
      \mbox{
        \subfloat[]{
          \resizebox{(\linewidth-1em)/2}{!}{\figureMidpointOneEdge}
          \label{figure:determineMidpointOfOneEdge}
        }
        \subfloat[]{
          \resizebox{(\linewidth-1em)/2}{!}{\figureMidpointTwoEdges}
          \label{figure:determineMidpointOfTwoEdges}
        }
      }
      \caption{Both subfigures depict a space-time diagram of the time evolution of the signal machine that we construct in \cref{section:firing-squad-solution}. The diagram in \cref{figure:determineMidpointOfOneEdge} illustrates how, beginning from one end of an edge, the midpoint of the edge is found, which works analogously for paths instead of edges, and the diagram in \cref{figure:determineMidpointOfTwoEdges} illustrates how, beginning from the inner vertex of a path consisting of two edges, the midpoint of the path is found, which works analogously for longer paths; above the right space-time diagram is a scaled-down depiction of the two-edged path. Vertices are thick and solid lines; find-midpoint signals of speed $1$ are thick and dotted lines; reflected find-midpoint signals of speed $1$ are densely dotted lines; slowed-down find-midpoint signals of speed $1/3$ are solid lines; stationary midpoint signals are thick and solid lines.} 
      \label{figure:determinedMidpointOfOneEdgeOrOfTwoEdges}
    \end{wide}
  \end{figure}

  \begin{figure}
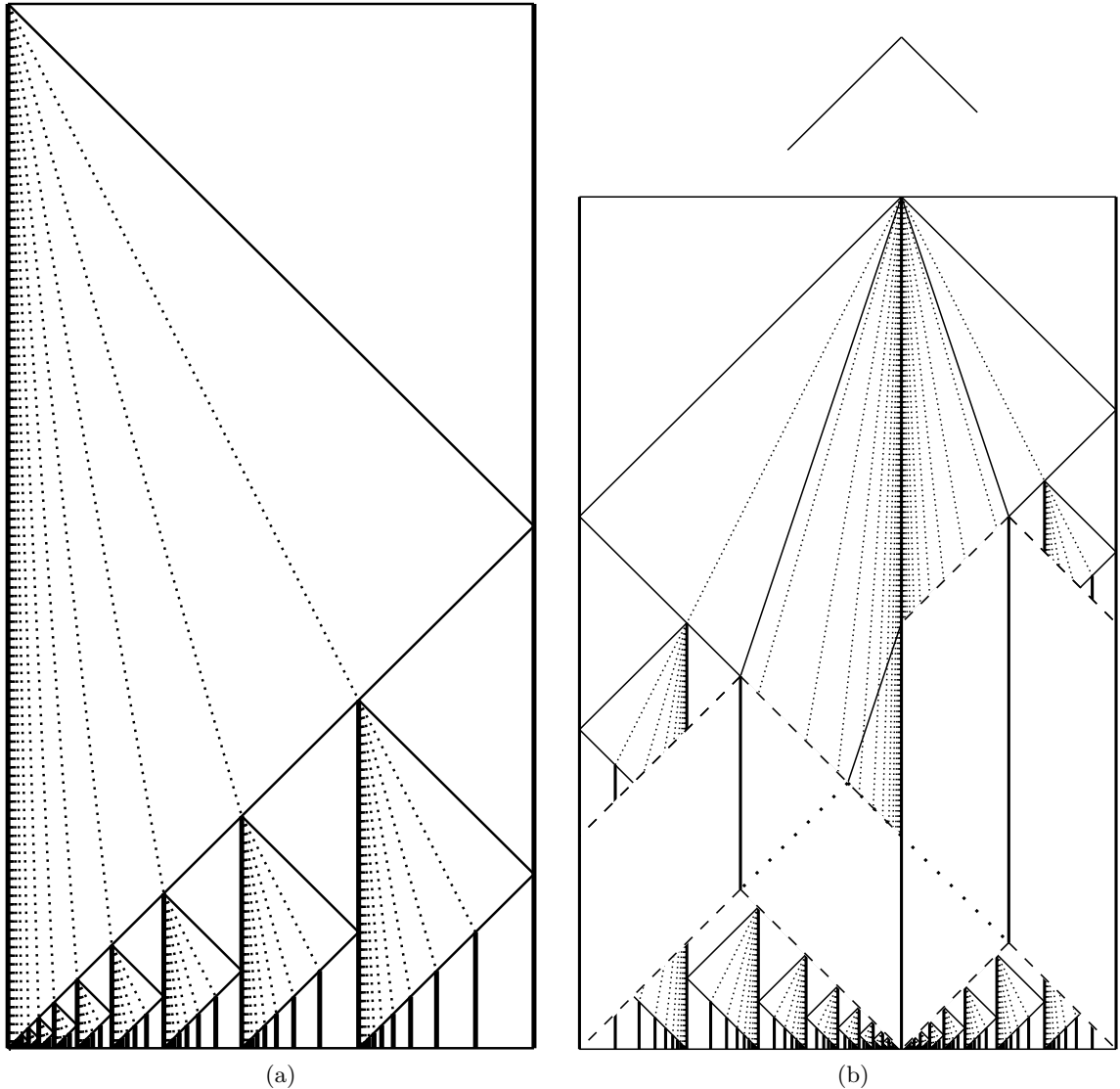

    \begin{wide}
      \mbox{
        \subfloat[]{
          \resizebox{(\linewidth-1em)/2}{!}{\figureMazoyer}
          \label{figure:Mazoyer}
        }
        \subfloat[]{
          \resizebox{(\linewidth-1em)/2}{!}{\figureFSSPMazoyerWithTwoEdges}
          \label{figure:MazoyerWithTwoEdges}
        }
      }
      \caption{Both subfigures depict a space-time diagram of the time evolution of the signal machine that we construct in \cref{section:firing-squad-solution}. The diagram in \cref{figure:Mazoyer} illustrates how, beginning from one end of an edge, the edge is synchronised by recursively dividing it into two parts, where one part is two-thirds the superpart's length and the other part is one-third the superpart's length, and creating fire signals when reflected divide signals reach the boundaries they originated at, which happens for all at the same time and these fire signals lie dense in the edge. And, the diagram in \cref{figure:MazoyerWithTwoEdges} illustrates how, beginning from the inner vertex of a path consisting of two edges, the path is synchronised by synchronising both edges, freezing the synchronisation of individual edges as late as possible and thawing it as early as possible such that both edges finish at the same time. Note that in both subfigures, at each boundary, there is a singularity of order $1$ at the last depicted time, and in \cref{figure:MazoyerWithTwoEdges}, there are additional singularities of order $1$ at time twice the longer edge's length and at time twice the shorter edge's length, and of order $-1$ at various points of the freeze and thaw signals. Only the most relevant signals are depicted and these only for the most relevant time spans. Initiate signals, divide signals of type $0$, and find-midpoint signals, all of speed $1$, are solid lines; reflected divide signals and reflected find-midpoint signals, both of speed $1$, are solid lines; divide signals of type $n \in \N_+$, which have speed $(2/3)^n / (2 - (2/3)^n)$, are densely dotted lines; slowed-down find-midpoint signals of speed $1/3$ are solid lines; stationary boundary and midpoint signals are thick and solid lines; freeze signals of speed $1$ are dashed lines; thaw signals of speed $1$ that do not thaw synchronisation of the edge they are on are thick and loosely dotted lines, and the other thaw signals of speed $1$ are dashed lines.} 
      \label{figure:MazoyerWithOneAndTwoEdges}
    \end{wide}
  \end{figure}

  \begin{figure}
    \begin{wide}
      \mbox{
        \subfloat[]{
          \resizebox{(\linewidth-1em)/2}{!}{\figureFSSPWithTwoEdgesAndGeneralInBetween}
          \label{figure:fsspWithTwoEdgesAndGeneralInBetween}
        }
        \subfloat[]{
          \resizebox{(\linewidth-1em)/2}{!}{\figureFSSPWithTwoEdgesAndGeneralAtTheLeft}
          \label{figure:fsspWithTwoEdgesAndGeneralAtTheLeft}
        }
      }
      \caption{Both subfigures depict the space-time diagram of the time evolution of the signal machine that we construct in \cref{section:firing-squad-solution} for the scaled-down depicted tree with two edges, beginning with the initial configuration for the firing squad synchronisation problem --- which is the configuration in which initiate signals, find-midpoint signals, and slowed-down find-midpoint signals emanate from the general onto all incident edges --- and ending with the final configuration of the problem --- which is the configuration in which the fire signals lie dense in the multigraph ---, where in \cref{figure:fsspWithTwoEdgesAndGeneralInBetween} the general is the vertex that is incident to both edges and in \cref{figure:fsspWithTwoEdgesAndGeneralAtTheLeft} it is the leaf of the longer edge. Only the most relevant signals are depicted and these only for the most relevant time spans. Initiate signals and find-midpoint signals, which in the depicted cases always travel alongside each other, are thick and dotted lines; reflected find-midpoint signals are densely dotted lines; slowed-down find-midpoint signals are solid lines; midpoint signals are thick and solid lines; freeze signals are dashed lines; thaw signals that do not thaw synchronisation of the edge they are on are thick and loosely dotted lines, and the other thaw signals are dashed lines. The synchronisation of individual edges, before it is frozen and after it is thawed, is schematically represented by hatch patterns.} 
      \label{figure:fsspWithTwoEdgesAndGeneralInBetweenAndAtTheLeft}
    \end{wide}
  \end{figure}

  \begin{figure}
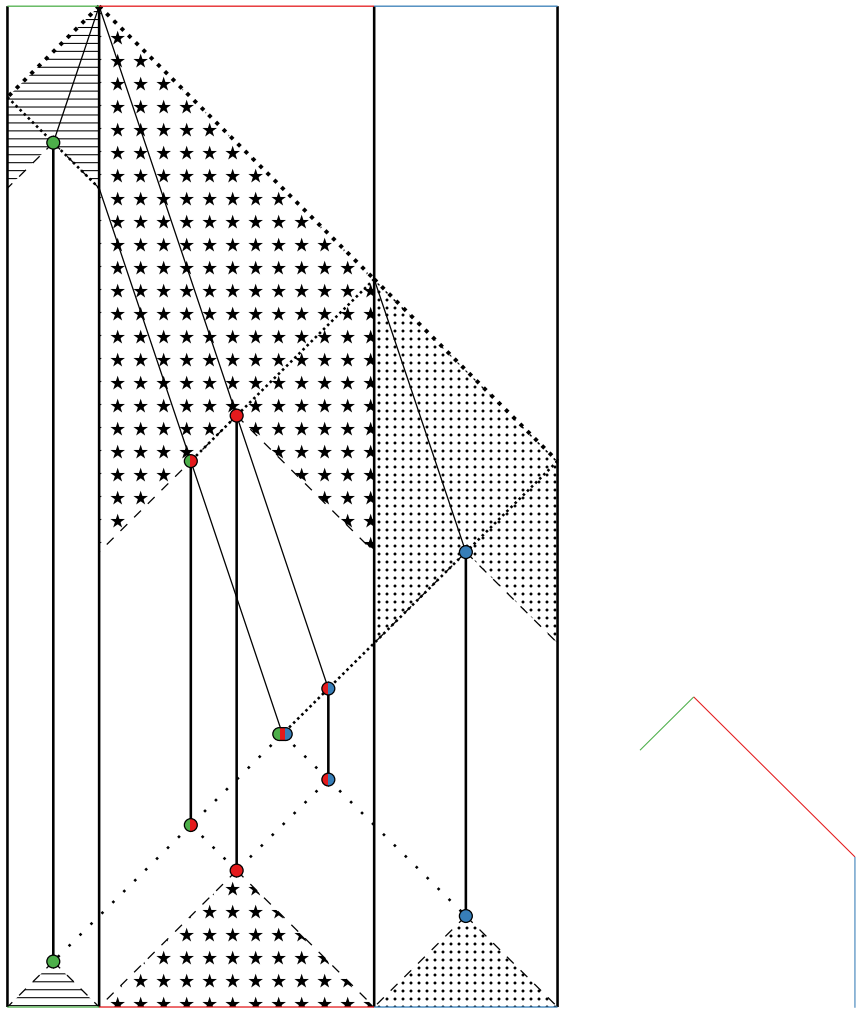

    \begin{wide}
      \resizebox{(\linewidth-1em)/2}{!}{\figureFSSPWithThreeEdgesInARowAndGeneralAtTheSecondVertexFromTheLeft}
      \qquad
      \figureTheTreeThatIsSynchronised
      \caption{The figure depicts the space-time diagram of the time evolution of the signal machine that we construct in \cref{section:firing-squad-solution} for the tree with three edges depicted on the right, beginning with the initial configuration for the firing squad synchronisation problem --- which is the configuration in which initiate signals, find-midpoint signals, and slowed-down find-midpoint signals emanate from the general onto all incident edges --- and ending with the final configuration of the problem --- which is the configuration in which the fire signals lie dense in the multigraph ---, where the general is the second vertex from the left. Only the most relevant signals are depicted and these only for the most relevant time spans. Initiate signals and find-midpoint signals are thick and dotted lines; reflected find-midpoint signals are densely dotted lines; slowed-down find-midpoint signals are solid lines; midpoint signals are thick and solid lines; freeze signals are dashed lines; thaw signals that do not thaw synchronisation of the edge they are on are thick and loosely dotted lines, and the other thaw signals are dashed lines. When a midpoint is found and when the corresponding midpoint signal collides with a matching thaw signal, the path it is the midpoint of is represented by the colours of the path's edges in a little disk or rounded rectangle. The synchronisation of individual edges, before it is frozen and after it is thawed, is schematically represented by hatch patterns.}
      \label{figure:fsspWithThreeEdgesInARowAndGeneralAtTheSecondVertexFromTheLeft}
    \end{wide}
  \end{figure}

  \begin{figure}
    \begin{wide}
      \mbox{
        \subfloat[]{
          \resizebox{(\linewidth-1em)/2}{!}{\figureFSSPWithThreeEdgesIncidentToTheGeneralVertex} 
          \label{figure:fsspWithThreeEdgesIncidentToTheSameVertexAndGeneralAtTheVertex}
        }
        \subfloat[]{
          \resizebox{(\linewidth-1em)/2}{!}{\figureFSSPWithThreeEdgesIncidentToTheSameVertexAndTheGeneralIsNotOnTheLongestPath}
          \label{figure:fsspWithThreeEdgesIncidentToTheSameVertexAndGeneralAtTheShortestEdge}
        }
      }
      \caption{Both subfigures depict the space-time diagram of the time evolution of the signal machine that we construct in \cref{section:firing-squad-solution} for the scaled-down depicted tree with three edges, beginning with the initial configuration for the firing squad synchronisation problem --- which is the configuration in which initiate signals, find-midpoint signals, and slowed-down find-midpoint signals emanate from the general onto all incident edges --- and ending with the final configuration of the problem --- which is the configuration in which the fire signals lie dense in the multigraph. In \cref{figure:fsspWithThreeEdgesIncidentToTheSameVertexAndGeneralAtTheVertex} the general is the one that is incident to all edges; and in \cref{figure:fsspWithThreeEdgesIncidentToTheSameVertexAndGeneralAtTheShortestEdge} it is the leaf of the shortest edge. Only the most relevant signals are depicted and these only for the most relevant time spans. Each edge has a colour, leaves are coloured according to the edge they are incident to, and signals and hatch patterns are coloured according to the edge they are on. The shortest edge is depicted twice, once coloured red and overlaying the green edge, and once coloured violet and overlaying the blue edge; this makes it easier to follow signals that travel from or onto the shortest edge. Initiate signals and find-midpoint signals, which in the depicted cases always travel alongside each other, are thick and dotted lines; reflected find-midpoint signals are densely dotted lines; slowed-down find-midpoint signals are solid lines; midpoint signals are thick and solid lines; freeze signals are dashed lines; thaw signals that do not thaw synchronisation of the edge they are on are thick and loosely dotted lines, and the other thaw signals are dashed lines. When a midpoint is found and when the corresponding midpoint signal collides with a matching thaw signal, the path it is the midpoint of is represented by the colours of the path's edges in a little disk. The synchronisation of individual edges, before it is frozen and after it is thawed, is schematically represented by hatch patterns.}
      \label{figure:fsspWithThreeEdgesIncidentToTheSameVertexAndGeneralAtTheVertexAndAtTheLeafOfTheShortestEdge}
    \end{wide}
  \end{figure}

  \begin{remark} 
    \begin{aenumerate}
      \item Turn the graph into a tree: Initiate signals of speed $1$ spread from the general throughout the graph without making U-turns and vanish at leaves. When they collide on an edge but not in one of its ends, the edge is cut in two at the point of collision by two stationary leaf signals (also called \emph{virtual leaves}), one for each of the two new ends. When they collide in a vertex coming from all incident edges, all edges are cut off by leaf signals; and when the do not come from all incident edges, all but one of the edges the signals come from are cut off. In the latter case, the initiate signal that comes from the edge that is not cut off spreads onto the edges from which no initiate signals came.

            In this way all cycles are eventually broken up and the graph is turned into a virtual tree. Because the virtual leaves are created as soon as possible and because they are treated just like normal leaves, we may and will assume in the description of the other tasks that the graph is a tree.
      \item Start and freeze synchronisation of edges (see \cref{figure:MazoyerWithOneAndTwoEdges}): When initiate signals reach a vertex, for each incident edge, synchronisation of the edge is started from the vertex immediately and freezing of this synchronisation is started from the midpoint of the edge after $3/2$ times the edge's length time units by sending freeze signals of speed $1$ to both ends of the edge and finishes after twice times the edge's length time units.

            Each edge is synchronised by recursively dividing it into two parts, one having two-third its length and the other having one-third its length. This division procedure becomes finer and finer, in other words, divisions accumulate, the closer the time evolution gets to the time $2$ times the edge's length after of the edge's synchronisation started.

            The division of one part is performed by sending divide signals of speed $(2/3)^n / (2 - (2/3)^n)$, for $n \in \N_0$, from one boundary onto the part, reflecting the divide signal of speed $1$ at the other boundary, and, upon collision of the reflected divide signal with a divide signal, letting the reflected divide signal move on, removing the involved divide signal, creating a stationary boundary signal and sending divide signals of the above speeds onto the newly created part the reflected divide signal comes from. 
      \item Determine the midpoints of all direction-preserving paths from vertices to vertices (see \cref{figure:determinedMidpointOfOneEdgeOrOfTwoEdges}): The general and initiate signals that reach vertices send find-midpoint signals of speed $1$ and slowed-down find-midpoint signals of speed $1/3$ in all directions to determine the midpoints of all paths that contain the general and the reached vertex respectively. These signals spread throughout the graph without making U-turns \emph{memorising the paths they take} and (fast) find-midpoint signals are additionally reflected. 

            Reflected find-midpoint signals of speed $1$ take the paths back they took before being reflected \emph{memorising the path to the vertex they were reflected at and the remaining path they have to take}. Upon finishing their path at the vertex they originated at, slowed-down find-midpoint signals of speed $1/3$ are sent onto all edges except the one the reflected signal comes from, \emph{they memorise the path from the origin vertex to the reflection vertex of the reflected find-midpoint signal}, and they spread throughout the graph without making U-turns \emph{memorizing the paths they take} and vanish at leaves.

            When a reflected find-midpoint signal collides with a slowed-down find-midpoint signal that originated at the same vertex, the point of collision is the midpoint of the concatenation of the paths the two signals took after being reflected, that is, the path from the vertex the reflected find-midpoint signal was reflected at over the point of collision over the vertex both signals originated at to the vertex the reflected find-midpoint signal that spawned the slowed-down find-midpoint signal was reflected at. Each midpoint is designated by a stationary midpoint signal that \emph{memorises the path it is the midpoint of along with its position on the path}.

            When a reflected find-midpoint signal collides with a reflected find-midpoint signal that originated at the same vertex and the point of collision is the origin vertex itself, both reflection vertices have the same distance to the origin vertex and this vertex is the midpoint of the concatenation of the paths the two signals took after being reflected. 

            To get a clearer picture of how midpoints are determined, let us focus on only a few signals and let us ignore boundary cases: When an initiate signal reaches a vertex, one find-midpoint signal is sent along one incident edge, and another find-midpoint signal is sent along another incident edge. Both signals travel along edges and upon reaching a vertex they are either reflected and travel back or they take one of the incident edges that leads them further away. At the latest, they are reflected upon reaching a leaf. One of the reflected find-midpoint signals returns first to the vertex it originated at, is slowed down there, and the slowed down signal travels towards the other (reflected) find-midpoint signal. At some time in the future, the slowed-down find-midpoint signal collides with the other, now reflected, find-midpoint signal and the point of collision is the midpoint of the concatenation of the paths the two signals took after being reflected. 
      \item Determine the midpoints of the longest direction-preserving paths (see \cref{figure:fsspWithTwoEdgesAndGeneralInBetweenAndAtTheLeft,figure:fsspWithThreeEdgesInARowAndGeneralAtTheSecondVertexFromTheLeft}): The midpoints of the longest paths are eventually found. But this is not sufficient, they also need to be recognised as such. To that end, each reflected find-midpoint signal and each slowed-down find-midpoint signal carries a boolean that indicates whether the path it took from the vertex it was reflected at may be the subpath of a longest path that has the same source (or target), and whether the signal would be the first one to find the longest path's midpoint (recall that a slowed-down find-midpoint signal was either spawned by the general vertex or an initiate signal that reached a vertex, in which case we can think of this vertex as the one the slowed-down find-midpoint signal was reflected at; or it was originally a reflected find-midpoint signal and knows where that signal was reflected at). Let us call signals that carry the boolean $\yes$ \emph{marked} and the others not.

            At each vertex that is not a leaf, a stationary count signal memorises the directions from which marked reflected find-midpoint signals that originated at the vertex or from which marked slowed-down find-midpoint signals with any origin have already returned. Because longest paths always end at leaves, for each leaf, find-midpoint signals are marked when they are the first ones to be reflected at the leaf and not otherwise. When a marked signal reaches a vertex, it stays marked, if, including itself, from all but one direction have marked signals already returned (which is the case if and only if it is the last signal to return from its direction and the penultimate signal to return among such signals from either direction), and it is unmarked, otherwise.

            This means that a marked reflected find-midpoint signal is unmarked, if there are at least two other signals that started out at the same time at the same vertex but take longer to return (because they have a longer way; the combined paths that two of the other signals take may be a longest path), and it stays marked, if all other signals except for one that all started out at the same time at the same vertex returned earlier (because they had a shorter way).

            And it means that a marked slowed-down find-midpoint signal is unmarked, if it is unclear which of the incident edges belong to longest paths, and it stays marked, otherwise, which is the case if from precisely one direction no slowest reflected find-midpoint signal that originated at the vertex has returned yet. This is not only pessimistic, meaning that we do not falsely consider midpoints of paths that are not among the longest as such, but also correct, meaning that we still find the midpoints of all longest paths in time (see \cref{subsection:midpoints-of-maximum-weight-paths-are-recognised}).

            When a marked reflected find-midpoint signal reaches its origin and is the penultimate such signal to do so, it turns into a marked slowed-down find-midpoint signal that travels in the one direction from which no marked signal has returned yet and the point at which this signal collides with the one signal that has not yet returned will be the midpoint of a longest path, if both signals are still marked at the time of collision.
      \item Traverse midpoints, thaw synchronisation of edges, and fire (see \cref{figure:fsspWithThreeEdgesInARowAndGeneralAtTheSecondVertexFromTheLeft,figure:fsspWithThreeEdgesIncidentToTheSameVertexAndGeneralAtTheVertexAndAtTheLeafOfTheShortestEdge}): The midpoints of the longest paths are found at the same time, at which, from each such midpoint, two thaw signals of speed $1$ are sent that travel along the midpoint's path towards both of its ends. When a thaw signal collides with the midpoint of a path such that one of the two subpaths from the midpoint to either end of the path coincides with the remaining path the thaw signal travels along, an additional thaw signal is created that travels along the other subpath. On its way from the midpoint of the last edge of the path a thaw signal travels on to the end of the path, the thaw signal thaws all frozen signals it collides with. All thaw signals reach the ends of their paths at the same time, which is also the time all edges finish synchronisation with the creation of stationary fire signals that lie dense in the graph. \qedhere
    \end{aenumerate} 
  \end{remark}

  We introduce a typographic convention in

  \begin{definition}
    Each word of letters of the Latin alphabet that is written in typewriter font shall denote the word itself and shall for example not be the name of a variable.
  \end{definition}

  We introduce a boolean algebra in

  \begin{definition} 
    Let $\booleans = \setOf{\no, \yes}$\graffito{booleans $\booleans = \setOf{\no, \yes}$}\index[symbols]{Bblackboard@$\booleans$}. Each element $b \in \booleans$ is called \define{boolean}\index[symbols]{b@$b$}. The map
    \begin{align*}
      \neg \from \booleans &\to \booleans, \mathnote{negation $\neg$}\index[symbols]{negate@$\neg$}\\
      \no &\mapsto \yes,\\
      \yes &\mapsto \no,
    \end{align*}
    is called \define{negation}. The map
    \begin{align*}
      \land \from \booleans \times \booleans &\to \booleans, \mathnote{conjunction $\land$}\index[symbols]{wedge@$\land$}\\
      (b, b') &\mapsto \begin{dcases*}
        \no, &if $\no \in \setOf{b, b'}$,\\
        \yes, &otherwise,
      \end{dcases*}
    \end{align*}
    is called \define{conjunction}.
  \end{definition}

  We introduce finite lists of directions in

  \begin{definition}
    Let $\Directions^*$ be the set $\setOf{w \from \discreteInterval{1}{n} \to \Directions \suchThat n \in \N_0}$\graffito{set $\Directions^*$}\index[symbols]{Dirstarsuperscript@$\Directions^*$}. Each element $w \in \Directions^*$ is called \define{word over $\Directions$}\graffito{word $w$ over $\Directions$}\index[symbols]{w@$w$}; for each word $w \in \Directions^*$, the non-negative integer $\lengthOf{w} = \cardinalityOf{\domainOf(w)}$ is called \define{length of $w$}\graffito{length $\lengthOf{w}$ of $w$}\index[symbols]{absolutew@$\lengthOf{w}$}; the word $\emptyWord \from \emptyset \to \Directions$ is called \define{empty}\graffito{empty word $\emptyWord$}\index[symbols]{lambda@$\emptyWord$}; the map
    \begin{align*}
      \concat \from \Directions^* \times \Directions^* &\to \Directions^*, \mathnote{concatenation $\bullet$}\index[symbols]{bullet@$\concat$}\\
      (w, w') &\mapsto \left[
                         \begin{aligned}
                           \discreteInterval{1}{\lengthOf{w} + \lengthOf{w'}} &\to \Directions,\\
                           i &\mapsto \begin{dcases*}
                                        w(i), &if $i \leq \lengthOf{w}$,\\
                                        w'(i - \lengthOf{w}), &otherwise,
                                      \end{dcases*}
                         \end{aligned}
                       \right]
    \end{align*}
    is called \define{concatenation}.
  \end{definition}

  \begin{remark}
    The empty word is the only word of length $0$ and it is the neutral element of $\concat$.
  \end{remark}

  The signal machine we construct in this section has infinitely many kinds of signals, which are explained in \cref{remark:explanation-of-signals}, and given names, speeds, and data sets in

  \begin{definition}
    Let
    \begin{equation*}
      \Kinds = \setOf{\initiateKind, \leafKind, \countKind, \midpointKind, \findMidpointKind, \reflectedFindMidpointKind, \slowedDownFindMidpointKind, \freezeKind, \thawKind} \cup \parens[\big]{\bigcup_{n \in \N_0} \setOf{\divideKind_n, \frozenDivideKind_n}} \cup \setOf{\reflectedDivideKind, \boundaryKind, \fireKind, \frozenFireKind}, \mathnote{set $\Kinds$} 
    \end{equation*}
    let
    \begin{align*}
      \speedOf \from \Kinds &\to \R_{\geq 0}, \mathnote{map $\speedOf$}\\
            \initiateKind &\mapsto 1,\\
                \leafKind &\mapsto 0,\\
               \countKind &\mapsto 0,\\
            \midpointKind &\mapsto 0,\\
            \findMidpointKind &\mapsto 1,\\
        \reflectedFindMidpointKind &\mapsto 1,\\
            \slowedDownFindMidpointKind &\mapsto \frac{1}{3},\\
              \freezeKind &\mapsto 1,\\
                \thawKind &\mapsto 1,\\
            \divideKind_n &\mapsto \frac{\parens*{\frac{2}{3}}^n}{2 - \parens*{\frac{2}{3}}^n}, \text{ for } n \in \N_0,\\
      \frozenDivideKind_n &\mapsto 0, \text{ for } n \in \N_+,\\
          \reflectedDivideKind &\mapsto 1,\\
            \boundaryKind &\mapsto 0,\\
                \fireKind &\mapsto 0,\\
          \frozenFireKind &\mapsto 0,
    \end{align*}
    and let
    \begin{gather*}
      \NoData = \setOf{0}, \mathnote{set $\NoData$}\\ 
      \Data_{\initiateKind} = \NoData, \mathnote{family $\family{\Data_k}_{k \in \Kinds}$}\\
      \Data_{\leafKind} = \Directions,\\
      \Data_{\countKind} = \powerSetOf(\Directions),\\
      \Data_{\midpointKind} = \setOf{\setOf{w, w'} \subseteq \Directions^* \suchThat w \neq w'},\\
      \Data_{\findMidpointKind} = \Directions^*,\\
      \Data_{\reflectedFindMidpointKind} = \Directions^* \times \Directions^* \times \booleans,\\
      \Data_{\slowedDownFindMidpointKind} = \Directions^* \times \Directions^* \times \booleans,\\
      \Data_{\freezeKind} = \NoData,\\
      \Data_{\thawKind} = \Directions^* \times \booleans,\\
      \Data_{\divideKind_n} = \NoData, \text{ for } n \in \N_0,\\
      \Data_{\frozenDivideKind_n} = \Directions, \text{ for } n \in \N_+,\\
      \Data_{\reflectedDivideKind} = \NoData,\\
      \Data_{\boundaryKind} = \NoData,\\
      \Data_{\fireKind} = \NoData,\\
      \Data_{\frozenFireKind} = \NoData. \qedhere
    \end{gather*}
  \end{definition}

  The kinds together with their speeds and data sets determine the possible signals, which are recalled and given abbreviations in

  \begin{definition}
    The set of signals is
    \begin{equation*}
      \Signals = \setOf{(k, d, u) \suchThat k \in \Kinds \text{, } (\speedOf(k), d) \in \Arrows \text{, and } u \in \Data_k}. \mathnote{set $\Signals$}
    \end{equation*}
    Let
    \begin{gather*}
      \initiateSignal{d} = (\initiateKind, d, 0), \text{ for } d \in \Directions, \mathnote{abbreviations of signals like $\initiateSignal{d}$, $\reflectedFindMidpointSignal{w_o}{d}{w_r}{b}$ and $\frozenFireSignal$}\\
      \leafSignal{d} = (\leafKind, \every, d), \text{ for } d \in \Directions,\\
      \countSignal{D} = (\countKind, \every, D), \text{ for } D \in \Data_{\countKind},\\
      \midpointSignal{w}{w'} = (\midpointKind, \every, \setOf{w, w'}), \text{ for } \setOf{w, w'} \in \Data_{\midpointKind},\\
      \findMidpointSignal{w_o}{d} = (\findMidpointKind, d, w_o), \text{ for } w_o \in \Data_{\findMidpointKind},\\
      \reflectedFindMidpointSignal{w_o}{d}{w_r}{b} = (\reflectedFindMidpointKind, d, (w_o, w_r, b)), \text{ for } d \in \Directions \text{ and } (w_o, w_r, b) \in \Data_{\reflectedFindMidpointKind},\\
      \slowedDownFindMidpointSignal{d}{w_o}{w_r}{b} = (\slowedDownFindMidpointKind, d, (w_o, w_r, b)), \text{ for } d \in \Directions \text{ and } (w_o, w_r, b) \in \Data_{\slowedDownFindMidpointKind},\\
      \freezeSignal{d} = (\freezeKind, d, 0), \text{ for } d \in \Directions,\\
      \thawSignal{d}{w}{b} = (\thawKind, d, (w, b)), \text{ for } d \in \Directions \text{ and } (w, b) \in \Data_{\thawKind},\\
      \divideSignal{n}{d} = (\divideKind_n, d, 0), \text{ for } n \in \N_0 \text{ and } d \in \Directions,\\
      \frozenDivideSignal{n}{d} = (\frozenDivideKind_n, \every, d), \text{ for } n \in \N_+ \text{ and } d \in \Directions,\\
      \reflectedDivideSignal{d} = (\reflectedDivideKind, d, 0), \text{ for } d \in \Directions,\\
      \boundarySignal = (\boundaryKind, \every, 0),\\
      \fireSignal = (\fireKind, \every, 0),\\
      \frozenFireSignal = (\frozenFireKind, \every, 0). \qedhere
    \end{gather*}
  \end{definition}

  What signals of various kinds do when they reach a vertex or collide with one another, what the data they carry means, and what we call them is given a glimpse at in

  \begin{remark} 
  \label{remark:explanation-of-signals}
    \begin{aenumerate}
      \item Each signal of kind $\initiateKind$ has speed $1$; at each vertex it reaches it spreads in all directions that lead away from where it comes from, it initiates synchronisation of all incident edges except the one it comes from by sending divide signals onto them, it initiates the search for midpoints of paths that contain the vertex by sending (slowed-down) find-midpoint signals in all directions, and it initiates one component of the search for the longest paths of the graph by marking slowed-down find-midpoint signals if the vertex is a leaf; it carries no data; and it is called \graffito{initiate signal}\define{initiate signal}\index{signal!initiate}. The very first initiate signals spread from the general in all directions.
      \item Each signal of kind $\leafKind$ is stationary, designates a virtual leaf, carries the direction that leads onto the edge that is incident to the virtual leaf, and is called \define{leaf signal}\index{signal!leaf}\graffito{leaf signal}. Such signals are created when initiate signals collide in a vertex (or on an edge), which means that there is a cycle in the graph, and this cycle is broken up by virtually terminating the involved edge(s) with leaf signals. Each leaf signal is treated like a leaf in the following way: When signals collide with each other and with leaf signals, for each involved leaf signal, the collision of the signals that move in the opposite direction than the one the leaf signal carries is handled as if those signals collided in a leaf. Because leaf signals are created at points at the same time or before any other signal reaches them, the graph looks like a tree for all other signals.
      \item Each signal of kind $\countKind$ is stationary; is positioned at a vertex that is not a leaf; memorises the directions from which find-midpoint signals that originated at the vertex, were reflected, and may be on longest paths and would be the first ones to find their midpoints have already returned, in other words, it memorises the directions from which the slowest find-midpoint signals that originated at the vertex and travelled alongside initiate signals before they were reflected at a leaf have already returned; and is called \define{count signal}\index{signal!count}\graffito{count signal}. When an initiate signal reaches a vertex that is not a leaf, a count signal is created. 

            Note that, for the data set of count signals, instead of the infinite set $\powerSetOf(\Directions)$, we could have used the finite set $\powerSetOf(\setOf{1, 2, \dotsc, k})$, where $k$ is the maximum degree of the graph or the upper bound of the maximum degrees of the graphs to be considered and the numbers represent directions that lead away from vertices. The first choice of $k$ would make the signal machine depend on the graph, which is unconventional, whereas the latter would not and would also fit to the fact that solutions of firing mob synchronisation problems are usually considered for graphs whose maximum degrees are uniformly bounded by a constant.
      \item Each signal of kind $\midpointKind$ is stationary, designates the midpoint of a path, carries the directions that lead from its position to both ends of the path, and is called \define{midpoint signal}\index{signal!midpoint}\graffito{midpoint signal}. Such a signal is created when a reflected find-midpoint signal collides with a slowed-down find-midpoint signal that originated at the same vertex, or when two reflected find-midpoint signals that originated at the same vertex collide with each other, which only happens at the origin vertex itself. See \cref{figure:determinedMidpointOfOneEdgeOrOfTwoEdges}.
      \item Each signal of kind $\findMidpointKind$ has speed $1$, at each vertex it reaches it spreads in all directions that lead away from where it comes from (in the sense that, in each such direction, a signal of its kind is sent) and it is also reflected (in the sense that a reflected find-midpoint signal is sent in the direction from where it comes from), it carries the directions that lead from its position to the vertex the signal originated at, and it is called \define{find-midpoint signal}\index{signal!find-midpoint}\graffito{find-midpoint signal}. When an initiate signal reaches a vertex, for each incident edge, a find-midpoint signal whose origin is the vertex is created that travels onto the edge.
      \item Each signal of kind $\reflectedFindMidpointKind$ has speed $1$; is the reflection of a find-midpoint signal at a vertex, travels back along the path this signal took before it was reflected and slows down when it reaches the vertex the find-midpoint signal originated at; carries the directions that lead from its position to the vertex the find-midpoint signal originated at, the directions that lead from its position to the vertex the find-midpoint signal was reflected at, and a boolean that indicates whether the path described by its position and both directions, which leads from the reflection vertex to the origin vertex, may be the subpath of a longest path that has the same source (or target), and the boolean also indicates whether the signal would be the first one to find the longest path's midpoint; and is called \define{reflected find-midpoint signal}\index{signal!reflected find-midpoint}\index{signal!reflected find-midpoint!marked}\graffito{(marked) reflected find-midpoint signal} and, if the boolean it carries is $\yes$, it is called \define{marked}.

            As has already been pointed at, when a find-midpoint signal reaches a vertex, a reflected find-midpoint signal is created that travels onto the edge the find-midpoint signal comes from. If the vertex is a leaf and the find-midpoint signal is one of the first signals to reach it, which is precisely the case if the signal reaches the leaf together with an initiate signal, then its reflection is marked, and otherwise, not. The reasons are that both ends of a longest path are leaves and that a find-midpoint signal that is not among the first signals to reach one end of a longest path would not find its midpoint after another one has already found it.

            When a marked reflected find-midpoint signal reaches a vertex that is not a leaf, the count signal at the vertex memorises the direction the marked signal comes from, and the signal stays marked, if the memory of the count signal contains each but one direction that leads away from the vertex, and it is unmarked, otherwise. Why is that? Each vertex that is not the general is reached precisely once by an initiate signal, at which point find-midpoint signals are sent in all directions; for each direction, the marked reflected find-midpoint signal to return from that direction is memorised, which is the slowest one or, in other words, the last one or the one that had the longest way (note that although only one find-midpoint signal is sent in a direction, multiple reflected find-midpoint signals may return from that direction); the penultimate marked reflected find-midpoint signal to return may come from one edge of a longest path that runs through the vertex and hence it stays marked (note that the other edge of the longest path that is incident to the vertex would be the one from which the marked signal has not yet returned); the signals that return before the penultimate one are too fast to be on a longest path and the last signal to return has already collided with the slowed-down penultimate signal that returned before it (if they do not return at the same time) and hence they are unmarked.
      \item Each signal of kind $\slowedDownFindMidpointKind$ has speed $1/3$; is the slow-down of a reflected find-midpoint signal at the vertex the find-midpoint signal originated at; at each vertex it reaches it spreads in all directions that lead away from where it comes from; it carries the directions that lead from its position to the vertex the find-midpoint signal originated at, the directions that lead from the origin vertex to the vertex the find-midpoint signal was reflected at, and a boolean that indicates whether the path described by its position and both directions, which leads from the reflection vertex over the origin vertex to its position, may be the subpath of a longest path that has the same source (or target), and the boolean also indicates whether the signal would be the first one to find the longest path's midpoint; and it is called \define{slowed-down find-midpoint signal}\index{signal!slowed-down find-midpoint}\index{signal!slowed-down find-midpoint!marked}\graffito{(marked) slowed-down find-midpoint signal} and, if the boolean it carries is $\yes$, it is called \define{marked}.

            As has already been pointed at, when a reflected find-midpoint signal reaches the vertex the find-midpoint signal originated at, for each incident edge except the one the reflected find-midpoint signal comes from, a slowed-down find-midpoint signal is created that travels onto the edge. Additionally, when an initiate signal reaches a vertex, for each incident edge, a slowed-down find-midpoint signal is created that travels onto the edge.

            The latter case is in the following senses the boundary or limiting case of the former: Imagine that the vertex the initiate signal reaches is in fact two vertices that are infinitesimally close; then a find-midpoint signal is created at one vertex, this signal immediately reaches the infinitesimally close other vertex, there it is reflected, the reflected find-midpoint signal immediately reaches the infinitesimally close other vertex, and there it is slowed down. Or, analogously, imagine the limit of the cascade of the creation of a find-midpoint signal, its reflection, and slow-down for shorter and shorter distances between the vertex the find-midpoint signal originates at and the one it is reflected at; then in the limit the find-midpoint signal and its reflection vanish and only the slowed-down find-midpoint signal remains.

            When a marked slowed-down find-midpoint signal reaches a vertex that is not a leaf, the count signal at the vertex memorises the direction the marked signal comes from, and the signal stays marked, if the memory of the count signal contains each but one direction that leads away from the vertex, and it is unmarked, otherwise. Why is that? The slowed-down signal reaches the vertex at the same time and from the same direction as the slowest reflected find-midpoint signal from that direction that originated at the vertex. The latter signal is however not marked, because it did not travel alongside initiate signals before it was reflected at a leaf (the reason is that if it had travelled alongside initiate signals, then it would have been reflected at the same time as the find-midpoint signal whose reflection turned into the marked slowed-down find-midpoint signal and hence, because the paths from the reflection leaves to the vertex they reach together have the same lengths, the find-midpoint signal that reaches it slowed down would have taken longer, and therefore the signals would not reach the vertex at the same time).

            Therefore, the memory of the count signal contains each but one direction if and only if from each but one direction the slowest reflected find-midpoint signals that originated at the vertex have already returned. If this is the case, then the incident edge belonging to that direction may be the edge of a longest path that runs through the vertex and the slowed-down find-midpoint signal may be the one to collide with the not yet returned signal somewhere on or beyond the edge precisely at the midpoint of the longest path. If from more than one direction signals are overdue, then the paths running through each pair of these directions are longer than the paths running through any of these directions and the direction the slowed-down find-midpoint signal comes from. And, if all signals have already returned, then they have already collided with the slowed-down signal and found the midpoints of the longest paths whose determination involves the slowed-down signal if there are any. 
      \item Each signal of kind $\freezeKind$ has speed $1$, is created at the midpoint of an edge, moves towards one end of the edge, and freezes synchronisation of the edge, carries no data, and is called \define{freeze signal}\index{signal!freeze}\graffito{freeze signal}. When an initiate signal reaches a vertex, for each incident edge, a find-midpoint signal and a slowed-down find-midpoint signal are created that travel onto the edge, the former is reflected at the other end of the edge and collides with the latter at the midpoint of the edge, at which point two freeze signals are created that travel to both ends of the edge. See \cref{figure:fsspWithTwoEdgesAndGeneralInBetweenAndAtTheLeft}
      \item Each signal of kind $\thawKind$ has speed $1$, is created at the midpoint of a path, travels along the path towards one of end of the path, creates a new signal of its kind when it collides with the midpoint signal that designates the midpoint of a path such that one of the two subpaths from the midpoint to either end of the path (a \graffito{half-path}\define{half-path}\index{path!half-}) coincides with the path it takes itself and the new signal travels along the other half-path, and thaws synchronisation of an edge if it collided with or was created at the midpoint signal that designates the midpoint of the last edge of the path it takes, carries the directions of the path it takes and a boolean that indicates whether it thaws synchronisation of the edge it is on or not, and is called \define{thaw signal}\index{signal!thaw}\graffito{thaw signal}.

            The first thaw signals are created simultaneously at the midpoints of longest paths. For each such midpoint, two thaw signals are created, one that travels along the path to one end of the path and the other that travels to the other end of the path. When a thaw signal collides with the midpoint of a path whose one half-path coincides with the remaining path the thaw signal travels along, an additional thaw signal is created that travels along the other half-path. On its way from the midpoint of the last edge of the path a thaw signal travels on to the end of the path, the thaw signal thaws all frozen signals it collides with.

            In this way, starting at the midpoints of longest paths, thaw signals traverse the midpoints of shorter and shorter paths, reach the ends of their paths at the same time, and thaw synchronisation of edges, which finishes at the same time. See \cref{figure:fsspWithTwoEdgesAndGeneralInBetweenAndAtTheLeft,figure:fsspWithThreeEdgesInARowAndGeneralAtTheSecondVertexFromTheLeft,figure:fsspWithThreeEdgesIncidentToTheSameVertexAndGeneralAtTheVertexAndAtTheLeafOfTheShortestEdge}
      \item Each signal of kind $\divideKind_0$ has speed $1$, moves from one boundary (which may be one end of an edge or a boundary signal) to the next boundary (which may be the other end of the edge or a boundary signal) and is reflected there, carries no data, and is called \define{divide signal of type $0$}\graffito{divide signal of type $0$}. When an initiate signal reaches a vertex, for each incident edge, a divide signal of type $0$ is created that travels onto the edge. And, when a divide signal of any type collides with a reflected divide signal, a divide signal of type $0$ is created that travels in the same direction as the (non-reflected) divide signal.
      \item Each signal of kind $\divideKind_n$, for $n \in \N_+$, has speed $(2/3)^n / (2 - (2/3)^n)$, moves from one boundary (which may be one end of an edge or a boundary signal) towards the next boundary (which may be the other end of the edge or a boundary signal) but never reaches it and can be frozen, carries no data, and is called \graffito{divide signal of type $n$}\define{divide signal of type $n$}. When an initiate signal reaches a vertex, for each incident edge, and for each $n \in \N_+$, a divide signal of type $n$ is created that travels onto the edge. And, when a divide signal of any type collides with a reflected divide signal, for each $n \in \N_+$, a divide signal of type $n$ is created that travels in the same direction as the (non-reflected) divide signal.

            Note that although $\divideKind_n$, for $n \in \N_+$, are different kinds, events that involve signals of these kinds are handled the same way, in other words, signals of these kinds are not differentiated by the two local transition functions of the signal machine. The only reason they are different kinds is because we need them to have different speeds and by definition all signals of the same kind have the same speed.
      \item Each signal of kind $\frozenDivideKind_n$, for $n \in \N_+$, has speed $0$, is a frozen divide signal of type $n$, carries the direction the non-frozen divide signal had, and is called \define{frozen divide signal of type $n$}\graffito{frozen divide signal of type $n$}. When a freeze signal collides with or is created at the same time as a divide signal of type $n \in \N_+$, the divide signal is frozen.
      \item Each signal of kind $\reflectedDivideKind$ has speed $1$, is the reflection of a divide signal of type $0$, creates a boundary signal when it collides with a divide signal of type $n$, creates a fire signal when it reaches the end of the edge it traverses, carries no data, and is called \graffito{reflected divide signal}\define{reflected divide signal}\index{signal!reflected divide}. When a divide signal of type $0$ reaches a boundary (which may be a vertex or a boundary signal), a reflected divide signal is created that travels in the opposite direction.
      \item Each signal of kind $\boundaryKind$ is stationary, designates a boundary for the synchronisation of an edge, carries no data, and is called \graffito{boundary signal}\define{boundary signal}\index{signal!boundary}. Such signals are created when divide signals collide with reflected divide signals.

            On each edge, the interplay of divide signals, reflected divide signals, and boundary signals has the following effect: At first a divide signal of type $1$ collides with a reflected divide signal that originated at the same end of the edge. This collision results in the creation of a boundary signal that divides the edge into two parts. The length of the part from the origin vertex to the boundary signal is $2/3$ times the length of the edge and the length of the part from the boundary signal to the other end of the edge is $1/3$ times the length of the edge.

            In the same manner as the edge itself, the $(1/3)$-part is recursively divided further and further. In the $(2/3)$-part, a signal of type $2$ collides with the reflected divide signal from before. This collision results in the creation of a boundary signal that divides the $(2/3)$-part into two subparts. One has $(2/3) \cdot (2/3)$ the length of the edge and the other has $(1/3) \cdot (2/3)$ the length of the edge. The $((1/3) \cdot (2/3))$-part is recursively divided further and further. The $((2/3) \cdot (2/3))$-part is divided into a $((2/3) \cdot (2/3) \cdot (2/3))$-part and a $((1/3) \cdot (2/3) \cdot (2/3))$-part and so forth.

            If there were no freeze and thaw signals, after twice the time the edge is long --- which is precisely the time it took the divide signal of type $0$ to reach the other end of the edge, to be reflected there, and to return to the end it originated at --- the boundary signals together are dense on the edge, which means that each point on the edge is arbitrarily close to a boundary signal, and, at this point in time, each boundary signal collides with a reflected divide signal, which results in the creation of fire signals that designate that synchronisation has finished. 

            However, because the synchronisation of each edge is started at different times and takes different times depending on how far away the edge is from the general and how long the edge is, synchronisation of each edge is frozen at the last possible moment --- the freezing starts from the midpoint of the edge $3/2$ times the edge's length many time units after synchronisation of the edge was initiated --- and it is thawed such that all edges finish synchronisation at the same time --- the thawing starts from the midpoint of the edge $1/3$ times the edge's length many time units before the total synchronisation finishes, which is the sum of the radius of the graph and its diameter. Recall that the radius is the longest distance from the general to another vertex and that the diameter is the longest distance between two vertices.

            Note that, for each edge, collisions of divide signals with reflected divide signals and with boundary signals accumulate at the times the two freeze signals and the two thaw signals reach the ends of the edge. These accumulations are singularities of order $1$. See \cref{figure:MazoyerWithOneAndTwoEdges}.
      \item Each signal of kind $\fireKind$ is stationary, designates that synchronisation has finished and can be frozen, carries no data, and is called \define{fire signal}\index{signal!fire}\graffito{fire signal}. Such signals are created when reflected divide signals collide with boundary signals or reach vertices.
      \item Each signal of kind $\frozenFireKind$ is stationary, is a frozen fire signal, carries no data, and is called \define{frozen fire signal}\index{signal!frozen fire}\graffito{frozen fire signal}. On each edge, one of the two freeze signals reaches an end of the edge at the same time as the reflected divide signal, at which point a frozen fire signal is created; it is thawed at the same time at which all other fire signals are created, which happens on all edges at the same time and the fire signals lie dense in the multigraph. \qedhere 
    \end{aenumerate}
  \end{remark}

  In the forthcoming definitions of maps we make extensive use of pattern matching. To make the exposition concise and readable we introduce some pattern matching conventions in

  \begin{definition}
    \begin{aenumerate}
      \item In\graffito{order matters} the case that patterns of multiple rules overlap, the rule that occurs first is the one to use. For example, the map $f \from \Z \to \Z$, $0 \mapsto 1$, $1 \mapsto 0$, $z \mapsto z$, maps $0$ to $1$ and $1$ to $0$ and each other integer to itself.
      \item In\graffito{wildcard $\blank$}\index[symbols]{underscore@$\blank$} the case that we do not care to name some part of a matched structure, we write $\blank$ instead of a name for the part. For example, the pattern $(E, \blank, \initiateSignal{d})$ matches each triple whose last component is an initiate signal, gives the first component the name $E$, does not give a name to the second component, and gives the direction of the initiate signal the name $d$. And, the pattern $\findMidpointSignal{w_o}{\blank}$ matches each find-midpoint signal, gives the directions to the vertex the signal originated at the name $w_o$, but gives no name to the direction of the signal.
      \item To\graffito{same names express equality} express equality of different parts of a matched structure, we give those parts the same name. For example, the pattern $\setOf{\reflectedDivideSignal{\reverse d}, \divideSignal{n}{d}}$ matches each set that consists of a reflected divide signal and a divide signal of any type such that the direction of the reflected divide signal is the opposite of the direction of the divide signal, gives the direction of the divide signal the name $d$, and gives its type the name $n$.
      \item To\graffito{$@$-notation} name both a structure and its parts, we use a Haskell-like $@$-notation. For example, the pattern $s@\reflectedFindMidpointSignal{\emptyWord}{d}{\emptyWord}{b}$ matches each reflected find-midpoint signal whose directions to the vertex the signal originated at and to the vertex the signal was reflected at are empty, gives the direction of the signal the name $d$, gives the boolean that indicates whether the signal may be on a longest path and would be the first to find its midpoint the name $b$, and gives the signal itself the name $s$. \qedhere
    \end{aenumerate}
  \end{definition}

  Some of the maps we define below are actually partial maps. We represent them by (total) maps as specified in

  \begin{definition}
    Let\graffito{representations of partial maps using the bottom symbol $\bot$} $X$ and $X'$ be two sets, let $Y$ be a subset of $X$, let $\bot$ be an element that is not in $X'$, which we call \define{bottom}, and let $f$ be a map from $X$ to $X' \cup \setOf{\bot}$ such that, for each element $x \in X$, we have $f(x) = \bot$ if and only if $x \notin Y$. The map $f$ represents a partial map whose domain of definition is $Y$, whose domain is $X$, and whose codomain is $X'$.

    In the following, for maps like $f$, we do not explicitly specify the domain $Y$ of definition (it is the set $X \smallsetminus f^{-1}(\bot)$) and we implicitly assume that $\bot$ does not occur in the codomain $X'$ of the represented partial map. 
  \end{definition}

  To define the local transition functions, we begin with definitions for special cases and use those to gradually arrive at definitions for the general case. For trees and without freezing and thawing, the map $\localTransitionFunction_{v, 1}^{\tree}$ handles the event that precisely one signal reaches a vertex, and the maps $\localTransitionFunction_{v, 2}^{\tree}$ and $\localTransitionFunction_{e, 2}^{\tree}$ handle the event that precisely two signals collide in a vertex and an edge respectively (see \cref{definition:for-trees:single-signals-reach-a-vertex-and-pairs-collide}). For trees, the maps $\localTransitionFunction_v^{\tree}$ and $\localTransitionFunction_e^{\tree}$ handle events involving arbitrarily many signals by considering unordered pairs of signals and applying $\localTransitionFunction_{v, 1}^{\tree}$, $\localTransitionFunction_{v, 2}^{\tree}$, and $\localTransitionFunction_{e, 2}^{\tree}$, and by also freezing and thawing signals if needed (see \cref{definition:for-trees:local-transition-functions}). For virtual trees, which means that edges of the graph have been virtually cut by leaf signals to remove circles, the maps $\localTransitionFunction_v^{\virtualTree}$ and $\localTransitionFunction_e^{\virtualTree}$ handle events by partitioning signals at virtual cuts into those belonging to one or the other leaf and applying $\localTransitionFunction_v^{\tree}$ and $\localTransitionFunction_e^{\tree}$ (see \cref{definition:for-virtual-trees:local-transition-functions}). For general graphs, the maps $\localTransitionFunction_v$ and $\localTransitionFunction_e$ handle events by virtually cutting the graph, which eventually creates a virtual tree, and applying $\localTransitionFunction_v^{\virtualTree}$ and $\localTransitionFunction_e^{\virtualTree}$ (see \cref{definition:for-graphs:local-transition-functions}).

  Most of the forthcoming definitions and parts of them are annotated with intuitive explanations of what they mean. For example, after each rule that handles a specific kind of event, it is explained what kind of event in the time evolution of the signal machine is handled, how it is handled, and sometimes why.

  How events for trees and without freezing and thawing, with one or two signals involved are handled is given in

  \begin{definition}
  \label{definition:for-trees:single-signals-reach-a-vertex-and-pairs-collide}
    The following map tells whether two words of directions are both empty:
    \begin{align*}
      \areEmpty \from \Directions^* \times \Directions^* &\to \booleans, \mathnote{map $\areEmpty$}\\
      (w, w') &\mapsto \begin{dcases*}
                         \yes, &if $\lengthOf{w} = 0$ and $\lengthOf{w'} = 0$,\\
                         \no, &otherwise.
                       \end{dcases*}
    \end{align*}

    For trees and without freezing and thawing, the case that precisely two signals collide on an edge but not in one of its ends is handled by the following map, which maps colliding unordered pairs of signals for which a collision rule is specified to the resulting signals and all other sets of signals to $\bot$:
    \begin{align*} 
      \localTransitionFunction_{e, 2}^{\tree} \from \powerSetOf(\Signals) &\to \powerSetOf(\Signals) \cup \setOf{\bot}, \mathnote{map $\localTransitionFunction_{e, 2}^{\tree}$}\\
      \setOf{\divideSignal{0}{d}, \boundaryKind}
          &\mapsto \setOf{\reflectedDivideSignal{\reverse d}, \boundarySignal}, 
          \intertext{(If a divide signal of type $0$ collides with a boundary signal, then reflect the divide signal.)}
      \setOf{\reflectedDivideSignal{\reverse d}, \divideSignal{n}{d}}
          &\mapsto \setOf{\boundarySignal} \cup \setOf{\divideSignal{n'}{d} \suchThat n' \in \N_0},
          \intertext{(If a reflected divide signal collides with a divide signal of any type, then create a boundary signal and send divide signals of all types in the direction of the original divide signal.)}
      \setOf{s@\reflectedFindMidpointSignal{\emptyWord}{d}{\emptyWord}{b}, s'@\slowedDownFindMidpointSignal{\reverse d}{\emptyWord}{\emptyWord}{b'}}
          &\mapsto\\
                   &&\llap{$\begin{dcases*}
                     \setOf{s, s', \midpointSignal{\reverse d}{d}, \freezeSignal{\reverse d}, \freezeSignal{d}}, &if $b = \no$ or $b' = \no$,\\
                     \emptyset, &otherwise,
                   \end{dcases*}$}
          \intertext{(If a reflected find-midpoint signal collides with a slowed-down find-midpoint signal that originated at the same end of an edge and only travelled on this edge, then the point of collision is the midpoint of the edge and, if the graph has at least two edges, then let the signals move on, designate the point by a midpoint signal, and send freeze signals to both ends of the edge to freeze synchronisation of the edge, and otherwise, do not create and send any signals.)}
      \setOf{s@\reflectedFindMidpointSignal{w_o}{d}{w_r}{b}, s'@\slowedDownFindMidpointSignal{\reverse d}{w_o}{w_r'}{b'}}
          &\mapsto\\
                   &&\llap{$\begin{dcases*}
                     \setOf{s, s', \midpointSignal{(\reverse d) \concat w_r}{d \concat w_o \concat w_r'}}, &if $b = \no$ or $b' = \no$,\\
                     \setOf{\thawSignal{\reverse d}{w_r}{\no}, \thawSignal{d}{w_o \concat w_r'}{\no}}, &otherwise,
                   \end{dcases*}$}
          \intertext{(If a reflected find-midpoint signal collides with a slowed-down find-midpoint signal that originated at the same vertex, then the point of collision is the midpoint of the shortest path in the virtual tree from the vertex the reflected find-midpoint signal was reflected at to the vertex the signal that spawned the slowed-down find-midpoint signal was reflected at and, if this midpoint is not the midpoint of a longest path, then let the signals move on and designate the point by a midpoint signal, and otherwise, send thaw signals along that longest path to both its ends.)}
      \setOf{\thawSignal{d}{w}{\no}, \midpointSignal{d \concat w}{d' \concat w'}}
          &\mapsto \setOf{\thawSignal{d}{w}{\areEmpty(w, w')}, \thawSignal{d'}{w'}{\areEmpty(w, w')}},
          \intertext{(If a thaw signal collides with a midpoint signal that designates the midpoint of a path whose one half-path coincides with the path the thaw signal is going to take, then send an additional thaw signal along the other half-path and, if this path is just a directed edge, then make the thaw signals thaw synchronisation of the edge.)}
      \blank &\mapsto \bot.
          \shortintertext{(If none of the above happened, then indicate that by returning bottom.)}
    \end{align*}

    A signal that reaches a vertex is in a leaf if and only if there is precisely one direction that leads away from the vertex. And a signal that reaches a vertex is the penultimate one to do so if and only if the number of signals that have already returned including the signal itself is one less than the number of directions that lead away from the vertex. The two maps that express this in an abstract way using booleans are
    \begin{align*} 
      \isInLeaf \from \powerSetOf(\Directions) &\to \booleans, \mathnote{map $\isInLeaf$}\\
      E &\mapsto \begin{dcases*}
                   \no,  &if $\cardinalityOf{E} \neq 1$,\\
                   \yes, &otherwise,
                 \end{dcases*}
    \end{align*}
    and
    \begin{align*}
      \isPenultimate \from \powerSetOf(\Directions) \times \N_0 &\to \booleans, \mathnote{map $\isPenultimate$}\\
      (E, n) &\mapsto \begin{dcases*}
                        \no,  &if $\cardinalityOf{E} - 1 \neq n$,\\
                        \yes, &otherwise.
                      \end{dcases*}
    \end{align*}

    For trees and without freezing and thawing, the case that precisely one signal reaches a vertex is handled by the following map, which maps quadruples --- consisting of first, the set of directions that lead away from the vertex; secondly, the number of the directions from which the slowest reflected find-midpoint signals that originated at the vertex have already returned or have just arrived; thirdly, a boolean that indicates whether the signal is among the first ones to reach the vertex; and lastly, the signal that reaches the vertex --- to the resulting signals:
    \begin{align*}
      \localTransitionFunction_{v, 1}^{\tree} \from \powerSetOf(\Directions) \times \N_0 \times \booleans \times \Signals &\to \powerSetOf(\Signals), \mathnote{map $\localTransitionFunction_{v, 1}^{\tree}$}\\
      (\blank, \blank, \blank, \divideSignal{0}{d})
          &\mapsto \setOf{\reflectedDivideSignal{\reverse d}},
          \intertext{(If a divide signal of type $0$ reaches a vertex, then reflect it.)}
      (\blank, \blank, \blank, \reflectedDivideSignal{d})
          &\mapsto \setOf{\fireSignal},
          \intertext{(If a reflected divide signal reaches a vertex, then create a fire signal.)}
      (E, \blank, \blank, \initiateSignal{d})
          &\mapsto \begin{aligned}[t]
                     &\parens[\big]{\bigcup_{e \in E \smallsetminus \setOf{\reverse d}} \setOf{\initiateSignal{e}}}\\
                     &\cup \parens[\big]{\bigcup_{e \in E \smallsetminus \setOf{\reverse d}} \setOf{\divideSignal{n}{e} \suchThat n \in \N_0}}\\
                     &\cup \parens[\big]{\bigcup_{e \in E} \setOf{\findMidpointSignal{\emptyWord}{e}, \slowedDownFindMidpointSignal{e}{\emptyWord}{\emptyWord}{\isInLeaf(E)}}},
                   \end{aligned}
          \intertext{(If an initiate signal reaches a vertex, then send initiate signals onto all incident edges except the one the original initiate signal comes from, send divide signals of all types onto all incident edges except the one the original initiate signal comes from, and send find-midpoint and slowed-down find-midpoint signals onto all incident edges to find all midpoints of paths that contain the vertex, where, in the case that the vertex is a leaf, the slowed-down find-midpoint signal is marked, where the mark means that it may be on a longest path and would be the first to find its midpoint.)}
      (E, \blank, b, \findMidpointSignal{w_o}{d})
          &\mapsto \begin{aligned}[t]
                     &\setOf{\reflectedFindMidpointSignal{w_o}{\reverse d}{\emptyWord}{b \land \isInLeaf(E)}}\\
                     &\cup \parens[\big]{\bigcup_{e \in E \smallsetminus \setOf{\reverse d}} \setOf{\findMidpointSignal{(\reverse d) \concat w_o}{e}}},
                   \end{aligned}
          \intertext{(If a find-midpoint signal reaches a leaf, then reflect it and, if it is one of the first signals to reach the leaf, then also mark it as a signal that may be on a longest path and would be the first to find its midpoint. And, if a find-midpoint signal reaches a vertex that is not a leaf, then reflect it and send find-midpoint signals onto all incident edges except the one the original signal comes from.)}
      (E, n, \blank, \reflectedFindMidpointSignal{\emptyWord}{d}{w_r}{b})
          &\mapsto \bigcup_{e \in E \smallsetminus \setOf{\reverse d}} \setOf{\slowedDownFindMidpointSignal{e}{\emptyWord}{(\reverse d) \concat w_r}{b \land \isPenultimate(E, n)}}, 
          \intertext{(If a reflected find-midpoint signal reaches the vertex it originated at, then send slowed-down find-midpoint signals onto all incident edges except the one the original signal comes from and, if the original signal is marked and is the penultimate marked signal that originated at and has returned to the vertex, then also mark the slowed-down signals as signals that may be on a longest path and would be the first to find their midpoints.)}
      (E, n, \blank, \reflectedFindMidpointSignal{e \concat w_o}{d}{w_r}{b})
          &\mapsto \setOf{\reflectedFindMidpointSignal{w_o}{e}{(\reverse d) \concat w_r}{b \land \isPenultimate(E, n)}},
          \intertext{(If a reflected find-midpoint signal reaches a vertex, then it takes the way back it took before it was reflected and, if it is not one of the penultimate marked signals that reaches the vertex, then it is unmarked.)}
      (E, n, \blank, \slowedDownFindMidpointSignal{d}{w_o}{w_r}{b})
          &\mapsto \bigcup_{e \in E \smallsetminus \setOf{\reverse d}} \setOf{\slowedDownFindMidpointSignal{e}{(\reverse d) \concat w_o}{w_r}{b \land \isPenultimate(E, n)}},
          \intertext{(If a slowed-down find-midpoint signal reaches a vertex, then send slowed-down find-midpoint signals onto all incident edges except the one the original signal comes from, and mark these signals, if the original signal is marked --- in which case it arrives at the same time and from the same direction as the slowest reflected find-midpoint signal that originated at the vertex, which however is not marked because it arrived too late at the leaf it was reflected at --- and from exactly one direction the slowest find-midpoint signal that originated at the vertex has not returned yet.)}
      (\blank, \blank, \blank, \freezeSignal{d})
          &\mapsto \emptyset,
          \intertext{(If a freeze signal reaches the end of the edge it freezes, then it vanishes.)}
      (\blank, \blank, \blank, \thawSignal{d}{\emptyWord}{\yes})
          &\mapsto \emptyset,
          \intertext{(If a thaw signal reaches the end of its path, then it vanishes.)}
      (\blank, \blank, \blank, \thawSignal{d}{e \concat w}{\no})
          &\mapsto \setOf{\thawSignal{e}{w}{\no}},
          \intertext{(If a thaw signal reaches a vertex, then it takes the direction that makes it stay on its path.)}
      (E, \blank, \blank, (s, d, y)) 
          &\mapsto \left\{
                     \begin{aligned}
                       &\setOf{(s, d, y)}, \text{ if $\speedOf(s) = 0$},\\
                       &\bigcup_{e \in E \smallsetminus \setOf{\reverse d}} \setOf{(s, e, y)}, \text{ otherwise}. 
                     \end{aligned}
                   \right.
          \shortintertext{(A stationary signal in a vertex stays there, and if a non-stationary signal reaches a vertex, then copies of it are sent onto each edge except the one the signal comes from.)}
    \end{align*}

    For trees and without freezing and thawing, the case that precisely two signals collide in a vertex is handled by the following map, which maps triples --- consisting of first, the set of directions that lead away from the vertex; secondly, the number of the directions from which the slowest reflected find-midpoint signals that originated at the vertex have already returned or have just arrived; and lastly, the set of colliding signals that is supposed to consist of precisely two signals --- to the resulting signals, if a collision rule is specified, and to $\bot$, otherwise:
      \begin{align*}
        \localTransitionFunction_{v, 2}^{\tree} \from \powerSetOf(\Directions) \times \N_0 \times \powerSetOf(\Signals) &\to \powerSetOf(\Signals) \cup \setOf{\bot}, \mathnote{map $\localTransitionFunction_{v, 2}^{\tree}$}\\
        (\blank, \blank, \setOf{\reflectedDivideSignal{d}, \boundarySignal})
            &\mapsto \setOf{\fireSignal},
            \intertext{(If a reflected divide signal collides with a boundary signal, then create a fire signal.)}
        (E, n, \setOf{s@\reflectedFindMidpointSignal{(\reverse d') \concat w_o'}{d}{w_r}{b}, s'@\slowedDownFindMidpointSignal{d'}{w_o'}{w_r'}{b'}})
            &\mapsto\\
                    &&\llap{$\begin{dcases*}
                       \left.\begin{aligned}
                         &\localTransitionFunction_{v, 1}^{\tree}(E, n, \no, s)\\
                         &\quad\cup \localTransitionFunction_{v, 1}^{\tree}(E, n, \no, s')\\
                         &\quad\cup \setOf{\midpointSignal{(\reverse d) \concat w_r}{(\reverse d') \concat w_o' \concat w_r'}},
                       \end{aligned}\right\}
                       &if $b = \no$ or $b' = \no$,\\ 
                       \setOf{\thawSignal{\reverse d}{w_r}{\no}, \thawSignal{\reverse d'}{w_o' \concat w_r'}{\no}}, &otherwise,
                     \end{dcases*}$}
            \intertext{(If a reflected find-midpoint signal collides with a slowed-down find-midpoint signal that originated at the same vertex, then the vertex of collision is the midpoint of the shortest path in the virtual tree from the vertex the reflected find-midpoint signal was reflected at to the vertex the signal that spawned the slowed-down find-midpoint signal was reflected at and, if this midpoint is not the midpoint of a longest path, then treat the original signals as if they reached the vertex alone and designate the point by a midpoint signal, and otherwise, send thaw signals along that longest path to both its ends.)}
        (E, n, \setOf{s@\reflectedFindMidpointSignal{\emptyWord}{d}{w_r}{b}, s'@\reflectedFindMidpointSignal{\emptyWord}{d'}{w_r'}{b'}})
            &\mapsto\\
                     &&\llap{$\begin{dcases*}
                       \left.\begin{aligned}
                         &\localTransitionFunction_{v, 1}^{\tree}(E, n, \no, s)\\
                         &\quad\cup \localTransitionFunction_{v, 1}^{\tree}(E, n, \no, s')\\
                         &\quad\cup \setOf{\midpointSignal{(\reverse d) \concat w_r}{(\reverse d') \concat w_r'}},
                       \end{aligned}\right\}
                       &if $b = \no$ or $b' = \no$ or $n \neq \cardinalityOf{E}$,\\ 
                       \setOf{\thawSignal{\reverse d}{w_r}{\no}, \thawSignal{\reverse d'}{w_r'}{\no}}, &otherwise,
                     \end{dcases*}$}
            \intertext{(If two reflected find-midpoint signals that originated at the same vertex collide with each other, then the vertex of collision is the vertex the signals originated at and it is the midpoint of the shortest path in the virtual tree between the vertices the signals were reflected at and, if this midpoint is not the midpoint of a longest path, then treat the original signals as if they reached the vertex alone and designate the point by a midpoint signal, and otherwise, send thaw signals along that longest path to both its ends.)}
        (\blank, \blank, \setOf{\thawSignal{d}{w}{\no}, \midpointSignal{w@(d' \concat w')}{d'' \concat w''}})
            &\mapsto \setOf{\thawSignal{d'}{w'}{\no}, \thawSignal{d''}{w''}{\no}}, 
            \intertext{(If a thaw signal collides with a midpoint signal that designates the midpoint of a path whose one half-path coincides with the path the thaw signal is going to take, then send an additional thaw signal along the other half-path.)}
        \blank &\mapsto \bot.
            \shortintertext{(If none of the above happened, then indicate that by returning bottom.) \qedhere}
      \end{align*}
  \end{definition}

  How events for trees are handled is given in

  \begin{definition}
  \label{definition:for-trees:local-transition-functions}
    The maps
    \begin{equation*}
      \left\{
        \begin{aligned}
          \xi \from \Signals &\to \Signals,\\
          \divideSignal{n}{d} &\mapsto \frozenDivideSignal{n}{d},\\
          \fireSignal &\mapsto \frozenFireSignal,\\
          s &\mapsto s,
        \end{aligned}
      \right\}
      \text{ and }
      \left\{
        \begin{aligned} 
          \chi \from \Signals &\to \Signals,\\
          \frozenDivideSignal{n}{d} &\mapsto \divideSignal{n}{d},\\
          \frozenFireSignal &\mapsto \fireSignal,\\
          s &\mapsto s,
        \end{aligned}
      \right\}
      \mathnote{maps $\xi$ and $\chi$}
    \end{equation*}
    freeze and thaw signals that can be frozen and thawed respectively.

    The map
    \begin{align*}
      \nu \from \powerSetOf(\Signals) \times \powerSetOf(\Signals) &\to \powerSetOf(\Signals), \mathnote{map $\nu$}\\
      (S, S') &\mapsto \begin{dcases*}
                         \xi(S'), &if $\Exists \freezeSignal{d} \in S \cup S'$\\
                                  &and $\notExists \thawSignal{d}{\emptyWord}{\yes} \in S \cup S'$,\\
                         \chi(S'), &if $\Exists \thawSignal{d}{\emptyWord}{\yes} \in S \cup S'$\\
                                   &and $\notExists \freezeSignal{d} \in S \cup S'$,\\
                         S', &otherwise,
                       \end{dcases*}
    \end{align*}
    takes a set of old signals and a set of new signals and freezes the new signals, if the old or new signals contain a freeze signal but not a thaw signal that thaws the synchronisation of an edge; thaws the new signals, if the old or new signals contain a thaw signal that thaws the synchronisation of an edge but not a freeze signal; and does nothing, otherwise. 

    The maps
    \begin{align*}
      \zeta_2 \from \powerSetOf(\Signals) &\to \powerSetOf(\powerSetOf(\Signals)), \mathnote{map $\zeta_2$}\\
      S &\mapsto \setOf{\setOf{s, s'} \subseteq S \suchThat s \neq s' \text{ and } \localTransitionFunction_{e, 2}^{\tree}(\setOf{s, s'}) \neq \bot},
    \end{align*}
    and
    \begin{align*}
      \eta_2 \from \powerSetOf(\Directions) \times \powerSetOf(\Signals) &\to \powerSetOf(\powerSetOf(\Signals)), \mathnote{map $\eta_2$}\\
      (D, S) &\mapsto \begin{aligned}[t]
                        \{
                          \setOf{s, s'} \subseteq S \suchThat{} &s \neq s' \text{ and }\\
                                                                &\localTransitionFunction_{v, 2}^{\tree}(D, x, \setOf{s, s'}) \neq \bot
                        \},
                      \end{aligned}
    \end{align*}
    both take a set of signals and return the set of unordered pairs of distinct signals from the given set for which a collision rule is specified in $\localTransitionFunction_{e, 2}^{\tree}$ and $\localTransitionFunction_{v, 2}^{\tree}$ respectively.

    %

    The map
    \begin{align*} 
      \localTransitionFunction_e^{\tree} \from \powerSetOf(\Signals) &\to \powerSetOf(\Signals), \mathnote{map $\localTransitionFunction_e^{\tree}$}\\
      S &\mapsto \nu(S,
                   \parens[\big]{ S \smallsetminus \bigcup_{P \in \zeta_2(S)} P } \cup
                   \parens[\big]{ \bigcup_{P \in \zeta_2(S)} \localTransitionFunction_{e, 2}^{\tree}(P) }
                 ),
    \end{align*}
    handles collisions of signals on edges by leaving signals for which no pairwise collision rule with any other signal is specified in $\localTransitionFunction_{e, 2}^{\tree}$ as is, by applying $\localTransitionFunction_{e, 2}^{\tree}$ to each unordered pair of distinct signals for which a collision rule is specified, and by applying $\nu$ to freeze or thaw signals if there are or were any freeze or thaw signals.

    The map
    \begin{align*} 
      \kappa \from \powerSetOf(\Signals) &\to \powerSetOf(\Directions), \mathnote{map $\kappa$}\\
      S &\mapsto \begin{aligned}[t]
                   &\parens{\bigcup_{\countSignal{D} \in S} D}\\
                   &\cup \setOf{d \in \Directions \suchThat \Exists w_r \in \Directions^* \SuchThat \reflectedFindMidpointSignal{\emptyWord}{d}{w_r}{\yes} \in S}\\
                   &\cup \setOf{d \in \Directions \suchThat \Exists w_o \in \Directions^* \Exists w_r \in \Directions^* \SuchThat \slowedDownFindMidpointSignal{d}{w_o}{w_r}{\yes} \in S}\\
                 \end{aligned}
    \end{align*} 
    takes a set of signals that are at a vertex and returns the directions from which the slowest reflected find-midpoint signals that originated at the vertex have already returned or have just arrived. The directions from which signals have already returned is memorised by a count signal, and the other directions are the ones from which marked reflected find-midpoint signals that originated at the vertex or marked slowed-down find-midpoint signals with any origin have just arrived (the slowest ones of the latter kind always arrive at the same time and from the same direction as the slowest but unmarked reflected find-midpoint signal that originated at the vertex arrive from that direction). Note that although we take the union of the memories of all count signals, there is actually gonna be no such signal in leaves (in which case the union is $\emptyset$) and precisely one such signal in each vertex that is not a leaf (in which case the union is the memory stored by this signal). 

    The map
    \begin{align*} 
      \varkappa \from \powerSetOf(\Signals) &\to \booleans, \mathnote{map $\varkappa$}\\
      S &\mapsto \begin{dcases*}
                   \yes, &if $\Exists d \in \Directions \SuchThat \initiateSignal{d} \in S$,\\
                   \no, &otherwise,
                 \end{dcases*}
    \end{align*}
    tells whether a set of signals contains any initiate signals or not. It is used by $\localTransitionFunction_v^{\tree}$ to tell whether a find-midpoint signal that reaches a leaf is among the first ones to do so, which is the case if and only if the signal travels alongside an initiate signal.

    The map
    \begin{align*}
      \localTransitionFunction_v^{\tree} \from{} &\powerSetOf(\Directions) \times \powerSetOf(\Signals) \to \powerSetOf(\Signals), \mathnote{map $\localTransitionFunction_v^{\tree}$}\\
      &(D, S) \begin{aligned}[t]
                {}\mapsto \nu(S,
                  &\left\{\begin{aligned} 
                    &\emptyset, &&\text{if $\isInLeaf(D) = \yes$},\\ 
                    &\setOf{\countSignal{\kappa(S)}}, &&\text{otherwise},
                  \end{aligned}\right\}\\
                  &\cup \localTransitionFunction_{v, 1}^{\tree}\parens[\big]{D, \cardinalityOf{\kappa(S)}, \varkappa(S), S' \smallsetminus \bigcup_{P \in \eta_2(D, S')} P}\\
                  &\cup \bigcup_{P \in \eta_2(D, S')} \localTransitionFunction_{v, 2}^{\tree}(D, \cardinalityOf{\kappa(S)}, P)
                ),\\
                &\text{ where } S' = S \smallsetminus \setOf{\countSignal{D} \in S \suchThat D \subseteq \Directions}.
              \end{aligned}
    \end{align*}
    handles events in vertices by updating the memory of the directions from which the slowest reflected find-midpoint signals that originated at the vertex have already returned or have just arrived, by applying $\localTransitionFunction_{v, 1}^{\tree}$ to non-count signals for which no pairwise collision rule with any other non-count signal is specified in $\localTransitionFunction_{v, 2}^{\tree}$, by applying $\localTransitionFunction_{v, 2}^{\tree}$ to each unordered pair of distinct non-count signals for which a collision rule is specified, and by applying $\nu$ to freeze or thaw signals if there are or were any freeze or thaw signals.
  \end{definition}

  How events for virtual trees, where virtual leaves already exist, are handled is given in

  \begin{definition}
  \label{definition:for-virtual-trees:local-transition-functions}
    When colliding signals in the graph and the involved directions are partitioned with respect to the virtual tree, some of the components may be degenerated and applying $\localTransitionFunction_e^{\tree}$ and $\localTransitionFunction_v^{\tree}$ to them may have unwanted effects. Such boundary cases are properly handled by the maps 
    \begin{align*}
      \localTransitionFunction_e^{\boundaryCases} \from \powerSetOf(\Signals) &\to \powerSetOf(\Signals), \mathnote{map $\localTransitionFunction_e^{\boundaryCases}$}\\
      S &\mapsto \begin{dcases*}
                   \emptyset, &if $\cardinalityOf{S} \leq 1$,\\
                   \localTransitionFunction_e^{\tree}(S), &otherwise,
                 \end{dcases*}
    \end{align*}
    and
    \begin{align*}
      \localTransitionFunction_v^{\boundaryCases} \from \powerSetOf(\Directions) \times \powerSetOf(\Signals) &\to \powerSetOf(\Signals), \mathnote{map $\localTransitionFunction_v^{\boundaryCases}$}\\
      (D, S) &\mapsto \begin{dcases*}
                        \emptyset, &if $D = \emptyset$ or $S = \emptyset$,\\
                        \localTransitionFunction_v^{\tree}(D, S), &otherwise.
                      \end{dcases*}
    \end{align*}

    To handle an event involving the signals $S$ in a virtual tree we do the following: We partition the signals into leaf signals, namely $S_{\leafKind}$, for each leaf signal $\leafSignal{d} \in S_{\leafKind}$, the signals coming from the direction $\reverse d$, namely $S_d$, and all other signals, namely $S_o$. And we denote the set of directions that no leaf signal has by $D'$. Intuitively, $S_{\leafKind}$ is the set of virtual leaves, $S_d$ is the set of signals that reach the virtual leaf $\leafSignal{d}$, and $D'$ is the set of directions that do not lead away from any virtual leaf. In the configurations we will encounter, in the case of a collision on an edge, there are either no leaf signals, in which case $S_o = S$, or there are two leaf signals, in which case $S_o = \emptyset$ (because there are no stationary signals in virtual leaves besides the leaf signals themselves) and hence either the signals simply collide, or one subset of signals reaches one virtual leaf and the other subset reaches the other virtual leaf. And, in the case of a collision in a vertex, signals coming from some edges may reach a virtual leaf and signals coming from other edges may reach the virtual vertex (we call it virtual because some of its original edges have been cut off; the directions onto the ones that have not been cut off are those in the set $D'$). Note that some collisions are not collisions in the virtual tree, because the signals came from different directions of virtual cuts. The maps that do what we just explained are 
    \begin{align*}
      \localTransitionFunction_e^{\virtualTree} \from \powerSetOf(\Signals) &\to \powerSetOf(\Signals), \mathnote{map $\localTransitionFunction_e^{\virtualTree}$}\\
      S &\mapsto \begin{aligned}[t]
                   S_{\leafKind}
                   &\cup \localTransitionFunction_e^{\boundaryCases}(S_o)\\
                   &\cup \parens[\big]{\bigcup_{d \in X} \localTransitionFunction_v^{\boundaryCases}(\setOf{d}, S_d)},
                 \end{aligned}
    \end{align*}
    and
    \begin{align*}
      \localTransitionFunction_v^{\virtualTree} \from \powerSetOf(\Directions) \times \powerSetOf(\Signals) &\to \powerSetOf(\Signals), \mathnote{map $\localTransitionFunction_v^{\virtualTree}$}\\
      (D, S) &\mapsto \begin{aligned}[t]
                        S_{\leafKind}
                        &\cup \localTransitionFunction_v^{\boundaryCases}(D', S_o)\\
                        &\cup \parens[\big]{\bigcup_{d \in X} \localTransitionFunction_v^{\boundaryCases}(\setOf{d}, S_d)},
                      \end{aligned}
    \end{align*}
    where $X = \setOf{d \in \Directions \suchThat \leafSignal{d} \in S}$, the set of directions that lead away from virtual leaves, $S_{\leafKind} = \setOf{\leafSignal{d} \suchThat d \in X}$, the set of virtual leaves, $\family{S_d}_{d \in X} = \family{\setOf{s \in S \suchThat \directionOf(s) = \reverse d}}_{d \in X}$, for each direction of a virtual leaf, the set of signals that reach the virtual leaf corresponding to the direction moving towards it, $S_o = S \smallsetminus (S_{\leafKind} \cup (\bigcup_{d \in X} S_d))$, the signals that are not virtual leaves and that do not reach a virtual leaf, and $D' = D \smallsetminus X$, the set of directions that do not lead away from virtual leaves. 
  \end{definition}

  How events for graphs are handled is given in

  \begin{definition}
  \label{definition:for-graphs:local-transition-functions}
    The map
    \begin{align*}
      \mu \from \powerSetOf(\Signals) &\to \powerSetOf(\Signals), \mathnote{map $\mu$}\\
      S &\mapsto \reverse \setOf{d \in \Directions \suchThat \initiateSignal{d} \in S},
    \end{align*}
    takes a set of signals and returns the set of the reverses of the directions that initiate signals have.

    The map
    \begin{align*}
      \varphi_{\leafKind}^e \from \powerSetOf(\Signals) &\to \powerSetOf(\Signals), \mathnote{map $\varphi_{\leafKind}^e$}\\
      S &\mapsto \begin{dcases*}
                   \emptyset, &if $\cardinalityOf{\mu(S)} \leq 1$,\\
                   \setOf{\leafSignal{d} \suchThat d \in \mu(S)}, &otherwise,
                 \end{dcases*}
    \end{align*}
    takes a set of signals and returns the empty set, if there is at most one initiate signal, and the set that consists of a virtual leaf for each initiate signal, otherwise.

    The map
    \begin{align*}
      &\varphi_{\leafKind}^v \from \powerSetOf(\Directions) \times \powerSetOf(\Signals) \to \powerSetOf(\Signals), \mathnote{map $\varphi_{\leafKind}^v$}\\
      &(D, S) \mapsto \left\{
                        \begin{aligned} 
                          &\emptyset, \text{ if $\cardinalityOf{\mu(S)} \leq 1$},\\
                          &\setOf{\leafSignal{d} \suchThat d \in D}, \text{ if $\cardinalityOf{\mu(S)} \geq 2$ and $\mu(S) = D$,}\\
                              &&\llap{and $\isInLeaf(D) = \no$},\\ 
                          &\setOf{\leafSignal{d} \suchThat d \in D \smallsetminus \setOf{d_x}}, \text{ if $\cardinalityOf{\mu(S)} \geq 2$ and $\mu(S) \neq D$},\\
                              &&\llap{for some $d_x \in D \smallsetminus \mu(S)$}, 
                        \end{aligned}
                      \right.
    \end{align*}
    takes a set of directions and a set of signals and returns the empty set, if there is at most one initiate signal, or the set that consists of a virtual leaf for each initiate signal, if there are at least two initiate signals and initiate signals reached the vertex from all incident edges, or the set that consists of a virtual leaf for each initiate signal but one, otherwise. Note that right now the choice of $d_x$ is non-deterministic; however, if the finite set of directions carried a total order, then we could deterministically choose for example the smallest direction; or if continuum representations of graphs were embedded in high-dimensional Euclidean spaces and a Cartesian coordinate system was chosen such that the occurring directions are unit vectors, then the lexicographic order is a total order on the set of directions.

    The map
    \begin{align*}
      \localTransitionFunction_e \from \domainOf(\localTransitionFunction_e) &\to \powerSetOf(\Signals), \mathnote{map $\localTransitionFunction_e$}\\
      S &\mapsto \localTransitionFunction_e^{\virtualTree}(S \cup \varphi_{\leafKind}^e(S))
    \end{align*}
    handles collisions on edges by creating two virtual leaves, if two initiate signals collide, which cuts the edge virtually, and then applying $\localTransitionFunction_e^{\virtualTree}$ to the maybe new set of signals.

    The map
    \begin{align*}
      \localTransitionFunction_v \from \domainOf(\localTransitionFunction_v) &\to \powerSetOf(\Signals), \mathnote{map $\localTransitionFunction_v$}\\
      (D, S) &\mapsto \localTransitionFunction_v^{\virtualTree}(D, S \cup \varphi_{\leafKind}^v(D, S))
    \end{align*}
    handles events in vertices by creating a virtual leaf for each incident edge, if the vertex is not a leaf and initiate signals reached the vertex from all incident edges, or by creating a virtual leaf for each incident edge from which an initiate signal arrived except for one such edge, otherwise, and then applying $\localTransitionFunction_v^{\virtualTree}$ to the maybe new set of signals.

    Intuitively, initiate signals are used to turn the graph into a virtual tree by cutting edges at points where such signals collide. These cuts create virtual leaves which are represented by leaf signals. More precisely: When two initiate signals collide on an edge, it is cut by two leaf signals, one for each of the two directions. And when at least two initiate signals collide in a vertex, there are two cases: If initiate signals arrive from all directions, then each incident edge is cut by a leaf signal; otherwise, the incident edges from which initiate signals arrive are cut except for one such edge --- the initiate signal from this excluded edge will spread to all edges that have not been cut. 
  \end{definition}

  \begin{main-theorem} 
  \label{theorem:signal-machine-is-quasi-solution-of-the-firing-mob-synchronisation-problem}
    The signal machine \graffito{signal machine $\mathcal{S}$}$\mathcal{S} = \ntuple{\Kinds, \speedOf, \family{\Data_k}_{k \in \Kinds}, (\localTransitionFunction_e, \localTransitionFunction_v)}$ is a time-optimal quasi-solution of the firing mob synchronisation problem over continuum representations of weighted, non-trivial, finite, and connected undirected multigraphs in the following sense: For each representation $M$ of such a graph, each vertex $\general \in M$, for the time $t = r + d$, where $r = \sup_{m \in M} \distanceOf(\general, m)$ is the radius of $M$ with respect to $\general$ and $d = \sup_{m, m' \in M} \distanceOf(m, m')$ is the diameter of $M$, for the instantiation of $\mathcal{S}$ for $M$, for the configuration $c \in \Configurations$ such that $c(\general) = \bigcup_{d \in \directionOf(\general)} \setOf{\initiateSignal{d}} \cup \setOf{\divideSignal{n}{d} \suchThat n \in \N_0} \cup \setOf{\findMidpointSignal{\emptyWord}{d}, \slowedDownFindMidpointSignal{d}{\emptyWord}{\emptyWord}{\isInLeaf(\directionOf(\general))}}$ and $c\restrictedTo_{M \smallsetminus \setOf{\general}} \equiv \emptyset$, the points in the configuration $\globalTransitionFunction(t)(c)$ at which a fire signal occurs lies dense in $M$ with respect to the metric $\distanceOf$, and no fire signals occur in any of the configurations $\globalTransitionFunction(s)(c)$, for $s \in \R_{\geq 0}$ with $s < t$.
  \end{main-theorem}

  \begin{proof-sketch}
    A proof is sketched in \cref{section:proof-sketch}.
  \end{proof-sketch}

  \begin{remark}
  \label{remark:data-sets-of-kinds-can-be-chosen-to-be-finite}
    For each positive integer $k$, under the restriction to multigraphs whose maximum degree is bounded by $k$, for each such multigraph, because directions only need to be locally unique (compare \cref{remark:efficient-representation-of-directions}), the set $\setOf{1, 2, \dotsc, 2 k}$ can be chosen as the set of directions, which makes the data sets of the kinds $\leafKind$ and $\frozenDivideKind$ finite and independent of the multigraph, the finite set $\powerSetOf(\setOf{1, 2, \dotsc, k})$ can be chosen as the data set of the kind $\countKind$, and, depending on the diameter $d$ of the multigraph, the finite set of words over $\Directions$ with maximum length $d$ can be chosen as the sets of words over $\Directions$ that occur in the data sets of the kinds $\midpointKind$, $\findMidpointKind$, $\reflectedFindMidpointKind$, $\slowedDownFindMidpointKind$, and $\thawKind$ --- altogether, the data sets of all kinds can be chosen to be finite but some depend on the multigraph.
  \end{remark}

  \begin{corollary}
    A discretisation of the signal machine $\mathcal{S}$ is a time-optimal cellular automaton quasi-solution of the firing mob synchronisation problem over non-trivial, finite, and connected undirected multigraphs.
  \end{corollary}

  \begin{proof-sketch}
    Let $\Graph$ be a non-trivial, finite, and connected undirected multigraph. It is made up of paths whose source and target vertices are not of degree $2$ and whose other vertices are of degree $2$. Each such path together with its inverse can be regarded as an undirected \emph{uber-edge} whose weight is the length of the path and whose ends are the source and target vertices of the path or its inverse. We call vertices that are not of degree $2$ \emph{uber-vertices} and vertices that are of degree $2$ \emph{under-vertices}.

    Signals jump from vertices to vertices along edges. They \emph{collide} when they jump simultaneously from different vertices onto the same vertex or along the same edge but in different directions, or when signals jump onto vertices on which stationary signals reside. Collisions in uber-vertices are handled as collisions in vertices, collisions in under-vertices are handled as collisions on edges (namely on the uber-edges that contain the under-vertices), and collisions in the midst of edges are handled as collisions on edges (namely on the uber-edges that contain the edges). Signals \emph{reach a vertex} when they jump onto an uber-vertex, but not when they jump onto an under-vertex (because the latter just means that they travel along an uber-edge).

    When signals collide in a vertex, be it a uber- or under-vertex, the resulting signals are on the vertex. But when signals collide in the midst of an edge, the resulting signals must be distributed onto both its ends depending on their directions. This last case makes it rather cumbersome to write down the local transition functions explicitly. Vertices must be virtually divided into multiple parts: One part that plays the role of the vertex itself and, for each incident edge, an additional part that plays the role of the midpoint of the edge together with the corresponding part of the other end of the edge. And signals must be cleverly distributed onto these parts depending on their direction and how they came into being, and collisions of signals must also be cleverly handled taking the parts the involved signals came from and are on into account.

    This discretisation of the signal machine is actually a cellular automaton over the multigraph with appropriate dummy neighbours that are in a dead state (think for example of the multigraph as being embedded in a coloured $S$-Cayley graph with sufficient maximum degree and of the vertices that do not belong the multigraph as being in a dead state).
  \end{proof-sketch}

  \begin{remark}
    Jacques Mazoyer showed in 1987 that all infinitely many divide signals of type $n$, for $n \in \N_0$, that emanate from the same point can be generated by a cellular automaton with only finitely many states (see \cite{mazoyer:1987}). And, as illustrated in \cref{remark:data-sets-of-kinds-can-be-chosen-to-be-finite}, under the restriction to multigraphs whose maximum degrees are uniformly bounded by a constant, the data sets of all kinds can be chosen to be finite but some depend on the multigraph. Therefore, depending on the multigraph, the discretisation of $\mathcal{S}$ is a cellular automaton with a finite number of states.
  \end{remark}

  \begin{open-problem}
    Are there time-optimal signal machine and cellular automaton solutions of the firing mob synchronisation problem over non-trivial, finite, and connected undirected multigraphs whose maximum degrees are uniformly bounded by a constant? Or, more specifically, is it possible to adapt the signal machine $\mathcal{S}$ (and thereby its discretisation) such that the data sets of all kinds can be chosen to be finite and independent of the multigraph (and thereby making its discretisation have a finite set of states), for example by reducing the number of midpoints that are and need to be determined?
  \end{open-problem}

  \section{Proof Sketch of the Main Theorem}
  \label{section:proof-sketch}

  In this section, we sketch a proof of \cref{theorem:signal-machine-is-quasi-solution-of-the-firing-mob-synchronisation-problem}. To that end, let $\Graph = \ntuple{\Vertices, \Edges, \eendsOf}$ be a non-trivial, finite, and connected undirected multigraph, let $\weightOf$ be an edge weighting of $\Graph$, let $M$ be a continuum representation of $\Graph$, identify vertices of $\Graph$ and $M$ and direction-preserving paths from vertices to vertices of $\Graph$ and $M$, and let $\general$ be a vertex of $M$, which we call \define{general}\graffito{general vertex $\general$}. Furthermore, let $\mathcal{S}$ be the signal machine and let $c$ be the initial configuration of the firing mob synchronisation problem from \cref{theorem:signal-machine-is-quasi-solution-of-the-firing-mob-synchronisation-problem}, and, whenever we talk about time evolution, we mean the one of $\mathcal{S}$ that is in the configuration $c$ at time $0$, for example, \emph{at time $t$} either means \emph{in configuration $\globalTransitionFunction(t)(c)$} or \emph{essentially in configuration $\globalTransitionFunction(t)(c)$ but before events have been handled}.

  To proof \cref{theorem:signal-machine-is-quasi-solution-of-the-firing-mob-synchronisation-problem}, we need to ascertain that the signal machine performs the following tasks: First, it cuts the multigraph such that the multigraph turns into a virtual tree and looks like a tree to all other tasks; secondly, it starts synchronisation of edges and freezes it in time; thirdly, it determines the midpoints of all non-empty direction-preserving paths from vertices to vertices in time; fourthly, it determines which midpoints are the ones of the longest paths; fifthly, starting from the midpoints of the longest paths, it traverses midpoints of shorter and shorter paths and upon reaching midpoints of edges, it thaws synchronisation of the respective edges; sixthly, all edges finish synchronisation at time $r + d$ with the creation of fire signals that lie dense in the graph, where $r$ is the radius of the graph with respect to the general and $d$ is the diameter of the graph.

  That the first task is performed is evident from the definitions of $\localTransitionFunction_e$, $\localTransitionFunction_v$, $\localTransitionFunction_e^{\virtualTree}$, and $\localTransitionFunction_v^{\virtualTree}$. The only subtlety here is that besides leaf signals there cannot be stationary signals in virtual leaves or, more precisely, at points that are virtually cut by leaf signals, because stationary signals carry the semi-direction $\every$, which is insufficient to associate them with one or the other leaf signal as is done for non-stationary signals. This is no problem because the other tasks do not place stationary signals in (virtual) leaves. Therefore, we assume from now on, without loss of generality, that the multigraph $\Graph$ is a \emph{tree}.

  That the second task is performed can be seen from a careful examination of the definitions of $\localTransitionFunction_e^{\tree}$, $\localTransitionFunction_v^{\tree}$, $\localTransitionFunction_{e, 2}^{\tree}$, $\localTransitionFunction_{v, 1}^{\tree}$, and $\localTransitionFunction_{v, 2}^{\tree}$, where from the third task it is used that midpoints of edges are found in time to start the freezing process. 

  That the third, fourth, and fifth and sixth parts are performed is proven in \cref{subsection:midpoints-are-determined}, \cref{subsection:midpoints-of-maximum-weight-paths-are-recognised}, and \cref{subsection:thaw-signals-traverse-midpoints-and-thaw-synchronisation-of-edges-just-in-time}. 

  \subsection{Midpoints are Determined}
  \label{subsection:midpoints-are-determined}

  The midpoint of a path in a multigraph is the midpoint of its embedding in the continuum representation of the multigraph as introduced in

  \begin{definition} 
    Let $p$ be a path in $\Graph$. The point $\midpoint_p = \continuumRepresentationOf{p}(\weightOf(p) / 2)$ is called \define{midpoint of $p$}\graffito{midpoint $\midpoint_p$ of $p$}\index[symbols]{mpfraktur@$\midpoint_p$}.
  \end{definition}

  \begin{remark}
    The midpoint of the empty path in $v$ is the vertex $v$ itself.
  \end{remark}

  \begin{remark}
    Let $p$, $q$, and $q'$ be three paths such that $\targetOf(p) = \sourceOf(q) = \sourceOf(q')$, $\midpoint_{p \concat q} \in \imageOf p$, and $\weightOf(q) \geq \weightOf(q')$. Then, $\midpoint_{p \concat q'} \in \imageOf p$ and
    \begin{equation*}
      \distanceOf(\midpoint_{p \concat q}, \midpoint_{p \concat q'})
      = \weightOf(p \concat q) / 2 - \weightOf(p \concat q') / 2
      = \weightOf(q) / 2 - \weightOf(q') / 2.
    \end{equation*}

    Analogously, let $q$, $q'$, and $p$ be three paths such that $\targetOf(q) = \targetOf(q') = \sourceOf(p)$, $\midpoint_{q \concat p} \in \imageOf p$, and $\weightOf(q) \geq \weightOf(q')$. Then, $\midpoint_{q' \concat p} \in \imageOf p$ and
    \begin{equation*}
      \distanceOf(\midpoint_{q \concat p}, \midpoint_{q' \concat p})
      = \weightOf(q \concat p) / 2 - \weightOf(q' \concat p) / 2
      = \weightOf(q) / 2 - \weightOf(q') / 2. \qedhere
    \end{equation*} 
  \end{remark}

  When the midpoints of non-empty direction-preserving paths are found is stated in

  \begin{lemma}
  \label{lemma:the-time-at-which-the-midpoint-of-a-path-is-found}
    Let $p$ be a non-empty direction-preserving path in $\Graph$. The midpoint signal that designates the midpoint of $p$ is created at $\midpoint_p$ at time $t_p = \max\setOf{\distanceOf(\general, \sourceOf(p)),\allowbreak \distanceOf(\general, \targetOf(p))} + \weightOf(p)/2$.
  \end{lemma}

  \begin{proof-sketch}
    First, let the general $\general$ be the source or target of $p$. Then, a find-midpoint signal with origin $\general$ of speed $1$ and a slowed-down find-midpoint signal with origin $\general$ of speed $1/3$ travel from $\general$ towards the other end of $p$. The find-midpoint signal is reflected at the other end at time $(2/2) \cdot \weightOf(p)$ and this reflection collides with the slowed-down find-midpoint signal at the midpoint of $p$ at time $(3/2) \cdot \weightOf(p)$ creating a midpoint signal for $p$ (because both signals have the same origin). Note that the time of collision is equal to $t_p$ (see \cref{figure:determineMidpointOfOneEdge}).

    Secondly, let the general $\general$ lie on $p$ without being its source or target. Then, two find-midpoint signals with origin $\general$ travel from $\general$ to the ends of $p$; the source of $p$ is reached at time $\distanceOf(\general, \sourceOf(p))$ and the target at time $\distanceOf(\general, \targetOf(p))$. When such a signal reaches its end, it is reflected and travels back. If $\distanceOf(\general, \sourceOf(p)) = \distanceOf(\general, \targetOf(p))$, then this distance is equal to $\weightOf(p)/2$, the reflected signals collide at time $\weightOf(p)$ at the midpoint of $p$ creating a midpoint signal for $p$; note that the time of collision is equal to $t_p$. Otherwise, the reflected signal that is nearer to $\general$ reaches this vertex first, where the signal is slowed down and travels towards the other reflected signal with which it collides at the midpoint of $p$ at time $t_p$ (see \cref{figure:determineMidpointOfTwoEdges}).

    Lastly, let the general $\general$ not lie on $p$. Then, an initiate signal travels from $\general$ to the nearest vertex $v$ on $p$, where it creates two find-midpoint signals with origin $v$ that travel to the ends of $p$, are reflected at these ends and travel back, one is slowed-down upon reaching $v$, and the slowed-down signal collides with the reflected signal in the midpoint of $p$ at time $\distanceOf(\general, v) + \max\setOf{\distanceOf(v, \sourceOf(p)),\allowbreak \distanceOf(v, \targetOf(p))} + \weightOf(p)/2$ creating a midpoint signal for $p$. Note that the time of collision is equal to $t_p$.
  \end{proof-sketch}

  \begin{remark} 
    When reflected/slowed-down find-midpoint signals with different origins collide, nothing happens, the signals just move on. Hence, no points are falsely found to be midpoints.
  \end{remark}

  That the midpoints of maximum-weight direction-preserving paths are identical and found at time $r + d$ is shown in

  \begin{lemma}
  \label{lemma:the-time-at-which-the-midpoints-of-longest-paths-are-found}
    All maximum-weight direction-preserving paths in $\Graph$ have the same midpoint $\hat{\midpoint}$ and the midpoint signals that designate the midpoints of such paths are created at $\hat{\midpoint}$ at time $r + d/2$, where $r = \max_{v \in \Vertices} \distanceOf(\general, v)$ is the radius of $\Graph$ with respect to $\general$ and $d = \max_{v, v' \in \Vertices} \distanceOf(v, v')$ is the diameter of $\Graph$. Note that $r$ is equal to the radius $\sup_{m \in M} \distanceOf(\general, m)$ of $M$ with respect to $\general$ and $d$ is equal to the diameter $\sup_{m, m' \in M} \distanceOf(m, m')$ of $M$.
  \end{lemma}

  \begin{proof} 
    We prove both statements by contradiction.

    First, suppose that there are two maximum-weight direction-pre\-serv\-ing paths $\hat{p}$ and $\hat{p}'$ in $\Graph$ that do not have the same midpoint. Then, there is a non-empty direction-preserving path $\mathfrak{p}_{\midpoint}$ in $M$ from the midpoint of $\hat{p}$ to the one of $\hat{p}'$. And, there is a direction-preserving subpath $\mathfrak{p}$ of $\continuumRepresentationOf{\hat{p}}$ in $M$ from one end of $\hat{p}$ to its midpoint whose target-direction is not the reverse of the source-direction of $p_{\midpoint}$. And, there is a direction-preserving subpath $\mathfrak{p'}$ of $\continuumRepresentationOf{\hat{p}'}$ in $M$ from the midpoint of $\hat{p}'$ to one of its ends whose source-direction is not the reverse of the target-direction of $p_{\midpoint}$. The concatenation of $\mathfrak{p}$, $\mathfrak{p}_{\midpoint}$, and $\mathfrak{p}'$ is a direction-preserving path from vertex to vertex in $M$ whose length is equal to $d/2 + \length(\mathfrak{p}_{\midpoint}) + d/2 > d$. It corresponds to a direction-preserving path $p$ in $\Graph$ whose weight is greater than $d$, which contradicts that $d$ is the diameter of $\Graph$. Therefore, all maximum-weight direction-preserving paths in $\Graph$ have the same midpoint, which we denote by $\hat{\midpoint}$.

    Secondly, suppose that there is a maximum-weight direction-pre\-serv\-ing path $\hat{p}$ in $\Graph$ such that $\max\setOf{\distanceOf(\general, \sourceOf(\hat{p})),\allowbreak \distanceOf(\general, \targetOf(\hat{p}))} < r$. Let $v$ be a vertex of $\Graph$ such that $\distanceOf(\general, v) = r$, let $p_v$, $p_{\sourceOf}$, and $p_{\targetOf}$ be the direction-preserving paths in $\Graph$ from $v$, $\sourceOf(\hat{p})$, and $\targetOf(\hat{p})$ to $\general$, let $v'$ be the vertex on $p_v$ and $p_{\sourceOf}$ or on $p_v$ and $p_{\targetOf}$ that is the furthest from $\general$, and let $v''$ be the vertex on $\hat{p}$ that is the nearest to $\general$. Then, the weight of $p_v$ is $r$, the one of $p_{\sourceOf}$ is $\distanceOf(\general, \sourceOf(\hat{p})) < r$, and the one of $p_{\targetOf}$ is $\distanceOf(\general, \targetOf(\hat{p})) < r$. If $v'$ lies on the subpath of $p_{\sourceOf}$ from $\sourceOf(\hat{p})$ to $v''$ (which is equal to the subpath of $\hat{p}$ with the same ends), then let $p$ be the direction-preserving path from $v$ over $v'$ over $v''$ to $\targetOf(\hat{p})$; if $v'$ lies on the subpath of $p_{\targetOf}$ from $\targetOf(\hat{p})$ to $v''$ (which is equal to the subpath of the inverse of $\hat{p}$ with the same ends), then let $p$ be the direction-preserving path from $v$ over $v'$ over $v''$ to $\sourceOf(\hat{p})$; and otherwise, let $p$ be the direction-preserving path from $v$ over $v'$ over $v''$ to $\targetOf(\hat{p})$ (we could have chosen $\sourceOf(\hat{p})$ as well). See \cref{figure:the-time-at-which-the-midpoints-of-longest-paths-are-found} for a schematic representation of the three cases.
    \begin{figure}
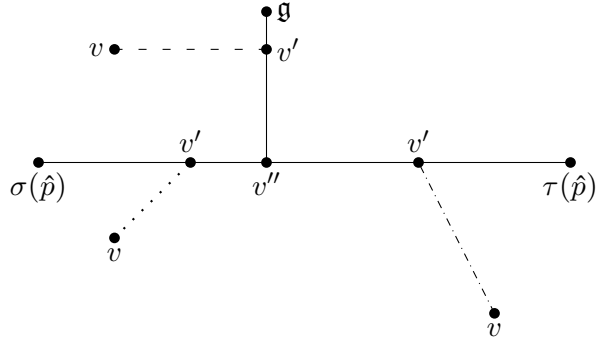

      \myfloatalign
      \figureTheTimeAtWhichTheMidpointsOfLongestPathsAreFound
      \caption{Schematic representation of the set-up of the proof of \cref{lemma:the-time-at-which-the-midpoints-of-longest-paths-are-found} with the three cases of where $v'$ may be located. The first case corresponds to the vertices $v$ and $v'$ that are incident to the dotted path, the second case to the dash-dotted path, and the third case to the dashed path.}
      \label{figure:the-time-at-which-the-midpoints-of-longest-paths-are-found}
    \end{figure}

    In the first case, because the subpaths of $p_v$ and $p_{\sourceOf}$ from $v''$ to $\general$ coincide and the weight of $p_v$ is greater than the weight of $p_{\sourceOf}$, the weight of the subpath of $p_v$ from $v$ over $v'$ to $v''$ (which is equal to the subpath of $p$ from $v$ over $v'$ to $v''$) is greater than the weight of the subpath of $p_{\sourceOf}$ from $\sourceOf(\hat{p})$ over $v'$ to $v''$ (which is equal to the subpath of $\hat{p}$ from $\sourceOf(\hat{p})$ over $v'$ to $v''$) and hence, because the subpaths of $p$ and $\hat{p}$ from $v''$ to $\targetOf(\hat{p})$ coincide, the weight of $p$ is greater than the weight of $\hat{p}$. In the second case, it follows analogously that the weight of $p$ is greater than the weight of $\hat{p}$. And in the third case, because the subpaths of $p_v$ and $p_{\sourceOf}$ from $v'$ to $\general$ coincide and the weight of $p_v$ is greater than the weight of $p_{\sourceOf}$, the weight of the subpath of $p_v$ from $v$ to $v'$ is greater than the weight of the subpath of $p_{\sourceOf}$ from $\sourceOf(\hat{p})$ over $v''$ to $v'$, hence the weight of the subpath of $p$ from $v$ over $v'$ to $v''$ is greater than the weight of the subpath of $\hat{p}$ from $\sourceOf(\hat{p})$ to $v''$, and therefore, the weight of $p$ is greater than the weight of $\hat{p}$.

    In either case, the inequality $\weightOf(p) > \weightOf(\hat{p})$ contradicts that $\hat{p}$ is a maximum-weight path. Therefore, for each maximum-weight direction-preserving path in $\Graph$, we have $\max\setOf{\distanceOf(\general, \sourceOf(\hat{p})),\allowbreak \distanceOf(\general, \targetOf(\hat{p}))} = r$. It follows from \cref{lemma:the-time-at-which-the-midpoint-of-a-path-is-found} that the midpoint signals that designate the midpoints of maximum-weight direction-preserving paths in $\Graph$ are created at $\hat{\midpoint}$ at time $r + d/2$.
  \end{proof}

  It follows that the midpoints of non-empty non-maximum-weight direction-preserving paths are found before the ones of maximum-weight direction-preserving paths as shown in

  \begin{corollary}
    Let $p$ be a non-empty direction-preserving path in $\Graph$. The midpoint signal that designates the midpoint of $p$ is created at $\midpoint_p$ before time $r + d/2$.
  \end{corollary}

  \begin{proof}
    This is a direct consequence of \cref{lemma:the-time-at-which-the-midpoint-of-a-path-is-found,lemma:the-time-at-which-the-midpoints-of-longest-paths-are-found}, because $\distanceOf(\general, \sourceOf(p)) \leq r$, $\distanceOf(\general, \targetOf(p)) \leq r$, and $\weightOf(p) < d$.
  \end{proof}

  \subsection{Midpoints of Maximum-Weight Paths are Recognised as Such}
  \label{subsection:midpoints-of-maximum-weight-paths-are-recognised}

  That the midpoints of maximum-weight direction-preserving paths are recognised as such is sketched in

  \begin{remark}
  \label{remark:midpoints-of-maximum-weight-paths-are-recognised}
    The first two reflected find-midpoint signals, or the first two reflected and slowed-down find-midpoint signals to collide that originated at the same vertex and were reflected at the ends of a maximum-weight direction-preserving path are marked at the time of collision and hence recognise that the midpoint of the path they collide at is the one of a maximum-weight direction-preserving path.
  \end{remark}

  \begin{proof-sketch} 
    Let $\hat{p}$ be a maximum-weight direction-preserving path in $\Graph$ (see \cref{figure:remark:midpoints-of-maximum-weight-paths-are-recognised}). Then, the ends $\hat{v}_1$ and $\hat{v}_2$ of $\hat{p}$ are leaves. And, among the first find-midpoint signals to reach the ends of $\hat{p}$ are the two that originated at the vertex $\hat{v}$ on $\hat{p}$ that is nearest to $\general$, and, because they travel alongside initiate signals, their reflections $s$ and $s'$ at the ends of $\hat{p}$ are marked. When one of them reaches a vertex on its way back it stays marked, because from all directions excluding from the direction it is headed but including the direction it is coming from have the marked reflected find-midpoint signals that originated at the vertex just or already returned, which is memorised by a count signal that is located at the vertex; the reason that such marked signals have already returned is that otherwise there would be a direction-preserving path with more weight than $\hat{p}$ that would be the concatenation the maximal subpath of $\hat{p}$ that lies in the direction the signal is headed and a path that begins with one of the edges from which no marked signal has returned yet. Note that the signals $s$ and $s'$ travel back alongside the marked reflected find-midpoint signals that originated at the vertices the signals $s$ and $s'$ passed by before they were reflected and that therefore, when they reach vertices on their way back, the count signals are just updated and hence up-to-date.
    \begin{figure}
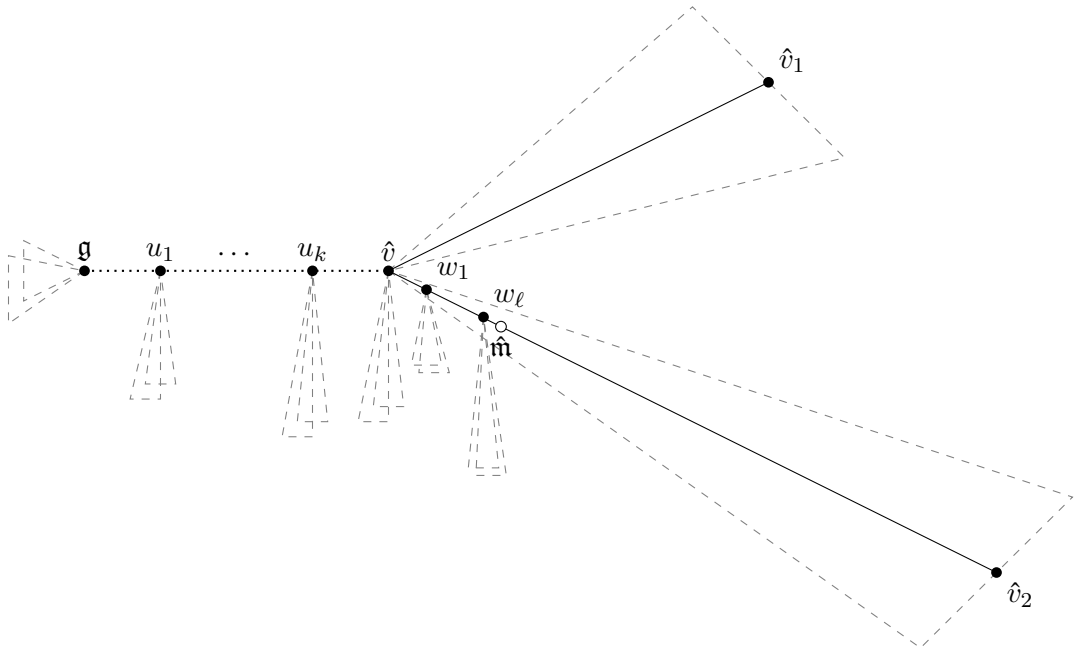

      \begin{wide}
        \figureRemarkMidpointsOfMaximumWeightPathsAreRecognised
        \caption{The path $\hat{p}$ is drawn solid, the direction-preserving path from $\mathfrak{g}$ to $\hat{v}$ is drawn dotted, the two subtrees in $\subtreesOf(\hat{v})$ that contain $\hat{v}_1$ and $\hat{v}_2$ are represented by dashed triangles, and, for each depicted vertex $v$, the possible existence of subtrees in $\subtreesOf(v)$ that correspond to non-depicted edges that are incident to $v$ is hinted at by dashed triangles.}
        \label{figure:remark:midpoints-of-maximum-weight-paths-are-recognised}
      \end{wide}
    \end{figure}

    When $s$ and $s'$ reach the vertex they originated at at the same time, then this vertex is found to be the midpoint of the path $\hat{p}$ and, because $s$ and $s'$ are marked, it is recognised as the midpoint of a maximum-weight path. When one of the signals, say $s$, reaches the vertex it originated at, namely $\hat{v}$, first, then it is slowed down and spreads throughout the graph away from the edge it comes from, in particular, one slowed-down find-midpoint signal, let us denote it by $s''$, travels towards $s'$. When $s''$ reaches a vertex on its way towards $s'$ it stays marked, because from all directions excluding from the direction it is headed have the \emph{slowest} reflected find-midpoint signals that originated at the vertex just or already returned, more precisely, from all directions excluding from the direction it is \emph{coming from and the one it is headed} have the marked reflected find-midpoint signals that originated at the vertex just or already returned and from the direction it is coming from has the slowest but unmarked reflected find-midpoint signal that originated at the vertex and the first marked slowed-down find-midpoint signal just returned. When $s''$ and $s'$ collide, the midpoint of $\hat{p}$ is found and, because $s''$ and $s'$ are marked, it is recognised as the midpoint of a maximum-weight path. 

    A detailed proof is given in the remainder of the present subsection.
  \end{proof-sketch}

  Symbolic notations for the predicates \emph{is a vertex of} and \emph{is an edge of} are introduced in

  \begin{definition}
    Let $\Tree = \ntuple{W, F}$ be a subtree of $\Graph$, let $v$ be a vertex of $\Graph$, and let $e$ be an edge of $\Graph$. We write $v \in \Tree$\graffito{$v \in \Tree$}\index[symbols]{vinTcalligraphic@$v \in \Tree$} instead of $v \in W$ and $e \in \Tree$\graffito{$e \in \Tree$}\index[symbols]{einTcalligraphic@$e \in \Tree$} instead of $e \in F$.
  \end{definition}

  The set of greatest subtrees of a vertex $v$ that correspond to its incident edges is named in

  \begin{definition}
    Let $v$ be a vertex of $\Graph$.
    \begin{aenumerate}
      \item Let $\Edges_v$ be the set of edges that are incident to $v$, and, for each edge $e \in \Edges_v$, let $\mathcal{T}_{v, e}$ be the greatest subtree of $\Graph$ that is rooted at $v$, contains the edge $e$, and does not contain any other edge of $\Edges_v$ (note that by \emph{greatest} we mean greatest with respect to the number of vertices). The set $\setOf{\mathcal{T}_{v, e} \suchThat e \in \Edges_v}$ is denoted by $\subtreesOf(v)$\graffito{set $\subtreesOf(v)$}\index[symbols]{sbtrsv@$\subtreesOf(v)$}.
      \item The set
            \begin{equation*}
              \maxSubtreesOf(v) = \argmax_{\Tree \in \subtreesOf(v)}\radiusOf{\Tree}_v
              \mathnote{set $\maxSubtreesOf(v)$ of trees}
              \index[symbols]{mxsbtrsv@$\maxSubtreesOf(v)$}
            \end{equation*}
            is the set of maximum-radius trees of $\subtreesOf(v)$. In the case that it is a singleton set, we denote its one and only element by $\maxTree_v$\graffito{tree $\maxTree_v$}\index[symbols]{Tvcalligraphichat@$\maxTree_v$}. 
      \item The set
            \begin{equation*}
              \secondMaxSubtreesOf(v) = \argmax_{\Tree \in \subtreesOf(v) \smallsetminus \maxSubtreesOf(v)}\radiusOf{\Tree}_v
              \mathnote{set $\secondMaxSubtreesOf(v)$ of trees}
              \index[symbols]{sbtrsvclosure@$\secondMaxSubtreesOf(v)$}
            \end{equation*}
            is the set of second-maximum-radius trees of $\subtreesOf(v)$. \qedhere
    \end{aenumerate}
  \end{definition}

  \begin{remark}
    For each leaf $v$ of $\Graph$, the set $\maxSubtreesOf(v)$ is a singleton and it is equal to $\subtreesOf(v)$, and the set $\secondMaxSubtreesOf(v)$ is empty.
  \end{remark}

  Things associated with a greatest subtree of a vertex are introduced in

  \begin{definition}
    Let $v$ be a vertex of $\Graph$ and let $\Tree$ be a tree of $\subtreesOf(v)$.
    \begin{aenumerate}
      \item Let $e$ be the edge of $\Tree$ that is incident to $v$. The direction that leads from $v$ onto $e$ is uniquely determined by $\Tree$ and is denoted by $\directionOf_v(\Tree)$\graffito{direction $\directionOf_v(\Tree)$}\index[symbols]{dirvTcalligraphic@$\directionOf_v(\Tree)$}. 
      \item The non-negative integer $\radiusOf{\Tree}_v = \max_{v' \in \Tree} \distanceOf(v, v')$ is called \graffito{radius $\radiusOf{\Tree}_v$ of $\Tree$ with respect to $v$}\define{radius of $\Tree$ with respect to $v$}\index[symbols]{normvsubscript@$\radiusOf{\blank}_v$}.
      \item Let $\Paths_{\Tree}$ be the set of direction-preserving paths in $\Tree$. The set
            \begin{equation*}
              \maxPaths_v(\Tree) = \setOf{p \in \Paths_{\Tree} \suchThat \sourceOf(p) = v \text{ and } \weightOf(p) = \radiusOf{\Tree}_v}
              \mathnote{set $\maxPaths_v(\Tree)$ of paths}
              \index[symbols]{mxpthsvTcalligraphic@$\maxPaths_v(\Tree)$}
            \end{equation*}
            is the set of maximum-weight direction-preserving paths from $v$ in $\Tree$. 
      \item The set
            \begin{equation*}
              \maxVertices_v(\Tree) = \setOf{\targetOf(p) \suchThat p \in \maxPaths_v(\Tree)}
              \mathnote{set $\maxVertices_v(\Tree)$ of vertices}
              \index[symbols]{mxvrtcsvTcalligraphic@$\maxVertices_v(\Tree)$}
            \end{equation*}
            is the set of vertices in $\Tree$ that are furthest away from $v$. \qedhere 
    \end{aenumerate}
  \end{definition}

  \begin{remark}
    We have $\cardinalityOf{\subtreesOf(v)} = \degreeOf(v)$ and $\directionOf_v(\subtreesOf(v)) = \directionOf(v)$.
  \end{remark}

  \begin{remark}
    Each vertex of $\maxVertices_v(\Tree)$ is a leaf.
  \end{remark}

  The unique direction-preserving path from one vertex to another is named in

  \begin{definition}
    Let $v$ and $v'$ be two vertices of $\Graph$. The direction-preserving path in $\Graph$ from $v$ to $v'$ is denoted by $p_{v, v'}$\graffito{path $p_{v, v'}$}\index[symbols]{pvvprime@$p_{v, v'}$} and the vertices on this path are denoted by $\Vertices_{v, v'}$\graffito{set $\Vertices_{v, v'}$ of vertices}\index[symbols]{Vvvprime@$\Vertices_{v, v'}$}.
  \end{definition}

  When and why count signals, initiate signals, and (maybe-marked slowed-down/reflected) find-midpoint signals are created and how they spread throughout the tree is said in

  \begin{remark} 
  \label{remark:how-initiate-and-find-midpoint-signals-spread}
    Let the signal machine $\mathcal{S}$ be in the configuration $c_{\initiateKind}$ at time $0$.
    \begin{aenumerate}
      \item At time $0$, a count signal with empty memory is created at $\general$, and initiate signals, find-midpoint signals with origin $\general$, and maybe-marked slowed-down find-midpoint signals with origin $\general$ and reflection vertex $\general$ are sent from $\general$ in all directions, where the slowed-down signal is marked if and only if $\general$ is a leaf. For each vertex $v \in \Vertices \smallsetminus \setOf{\general}$, at time $\distanceOf(\general, v)$, an initiate signal reaches $v$ from the direction towards $\general$, whereupon
            \begin{aenumerate}
              \item a count signal with empty memory is created at $v$,
              \item initiate signals are sent from $v$ in all directions away from $\general$, and
              \item find-midpoint signals with origin $v$, and maybe-marked slowed-down find-midpoint signals with origin $v$ and reflection vertex $v$ are sent from $v$ in all directions,
            \end{aenumerate}
            where the slowed-down signal is marked if and only if $v$ is a leaf.

            In short, initiate signals spread from $\general$ to all leaves, where they vanish, and they initiate the search for midpoints at all vertices. Note that to simplify the exposition of the forthcoming proofs, we also create count signals at leaves. 
      \item Let $v$ be a vertex of $\Graph$. As said above, at time $\distanceOf(\general, v)$, find-midpoint signals with origin $v$, and maybe-marked slowed-down find-midpoint signals with origin $v$ and reflection vertex $v$ are sent from $v$ in all directions, where the slowed-down signal is marked if and only if $v$ is a leaf. For each vertex $v' \in \Vertices \smallsetminus \setOf{v}$,
            \begin{aenumerate}
              \item at time $\distanceOf(\general, v) + \distanceOf(v, v')$, a find-midpoint signal with origin $v$ reaches $v'$ from the direction towards $v$, whereupon a maybe-marked reflected find-midpoint signal with origin $v$ and reflection vertex $v'$ is sent from $v'$ towards $v$ and find-midpoint signals with origin $v$ are sent from $v'$ in all directions away from $v$, where the reflected signal is marked if and only if $v$ is a leaf and an initiate signal just reached $v'$ (the latter is the case if and only if $\distanceOf(\general, v) + \distanceOf(v, v') = \distanceOf(\general, v')$), and,
              \item at time $\distanceOf(\general, v) + 3 \cdot \distanceOf(v, v')$, a maybe-marked slowed-down find-midpoint signal with origin $v$ and reflection vertex $v$ reaches $v'$ from the direction towards $v$, whereupon maybe-marked slowed-down find-midpoint signals with origin $v$ and reflection vertex $v$ are sent from $v'$ in all directions away from $v$. Note that even if the slowed-down signal that is sent from $v$ at time $\distanceOf(\general, v)$ is marked, the slowed-down signal that reaches $v'$ may be unmarked, and even if the latter signal is marked, the slowed-down signals that are sent from $v'$ may be unmarked.
            \end{aenumerate}
            In short, find-midpoint signals with origin $v$ and maybe-marked slowed-down find-midpoint signals with origin $v$ and reflection vertex $v$ spread from $v$ to all leaves, and whenever one of the former signals reaches a vertex, it is \emph{also} reflected. Note that to simplify the exposition of the forthcoming proofs, we talk as if maybe-marked slowed-down find-midpoint signals only vanish at leaves, although when a marked slowed-down find-midpoint signal collides with a marked reflected find-midpoint signal with the same origin, both signals vanish.
      \item Let $v$ and $v'$ be two vertices of $\Graph$ such that $v \neq v'$. As said above, at time $\distanceOf(\general, v) + \distanceOf(v, v')$ a maybe-marked reflected find-midpoint signal with origin $v$ and reflection vertex $v'$ is sent from $v'$ towards $v$. For each vertex $w$ on the direction-preserving path from $v'$ to $v$ except for $v'$, at time $\distanceOf(\general, v) + \distanceOf(v, v') + \distanceOf(v', w)$, the signal reaches $w$ from the direction towards $v'$, whereupon, if $w \neq v$, the signal is sent from $w$ towards $v$, and otherwise, maybe-marked slowed-down find-midpoint signals with origin $v$ and reflection vertex $v'$ are sent from $v$ in all directions away from $v'$. 

            And, for each vertex $v''$ of $\Graph$ that from the viewpoint of $v$ lies in a direction away from $v'$, which means that $v'' \neq v$ and there is a tree $\Tree \in \subtreesOf(v)$ such that $v' \notin \Tree$ and $v'' \in \Tree$, at time $\distanceOf(\general, v) + \distanceOf(v, v') + \distanceOf(v', v) + 3 \cdot \distanceOf(v, v'')$, a maybe-marked slowed-down find-midpoint signal with origin $v$ and reflection vertex $v'$ reaches $v''$, whereupon maybe-marked slowed-down find-midpoint signals with origin $v$ and reflection vertex $v'$ are sent from $v''$ in all directions away from $v$ (or, equivalently, away from $v'$). 

            In short, each maybe-marked reflected find-midpoint signal travels back to its origin and when it reaches its origin, it is slowed-down and spreads to all leaves away from its reflection vertex. Note that unmarked signals never become marked, but marked signals may become unmarked; precisely when the latter does or does not happen is answered in the present subsection.
    \end{aenumerate}
    What is said above is evident from the definition of the signal machine $\mathcal{S}$ (if it is carefully studied).
  \end{remark}

  When and why leaves send \emph{marked} reflected/slowed-down find-mid\-point signals is said in

  \begin{lemma}
  \label{lemma:when-do-leaves-sent-marked-signals}
    Let $v$ be a leaf of $\Graph$. At time $\distanceOf(\general, v)$, an initiate signal reaches $v$, for each vertex $w \in \Vertices_{\general, v} \smallsetminus \setOf{v}$, a find-midpoint signal with origin $w$ reaches $v$, and no other find-midpoint signal reaches $v$, whereupon
    \begin{aenumerate}
      \item for each vertex $w \in \Vertices_{\general, v} \smallsetminus \setOf{v}$, a marked reflected find-midpoint signal with origin $w$ and reflection vertex $v$ is sent from $v$ towards $w$, which is the only possible direction,
      \item a marked slowed-down find-midpoint signal with origin $v$ and reflection vertex $v$ is sent from $v$ towards $w$, and
      \item no other marked signal is sent from $v$.
    \end{aenumerate} 
    And before time $\distanceOf(\general, v)$ no signals reach and are sent from $v$, and after time $\distanceOf(\general, v)$ no marked signals are sent from $v$ (because after this time no initiate signal reaches $v$).
  \end{lemma}

  \begin{proof}
    This is a direct consequence of \cref{remark:how-initiate-and-find-midpoint-signals-spread} and the definition of $\mathcal{S}$.
  \end{proof}

  When does a signal that is sent from a vertex towards the general at a special time reach the next vertex is answered in

  \begin{lemma}
  \label{lemma:when-does-a-signal-that-is-sent-from-a-vertex-towards-the-general-reach-the-next-vertex}
    Let $v$ be a vertex of $\Graph$, let $\Tree$ be a tree of $\subtreesOf(v) \smallsetminus \maxSubtreesOf(v)$ such that either $v = \general$ or $\general \notin \Tree$, and let $v'$ be the one and only neighbour of $v$ in $\Tree$. Then,
    \begin{aenumerate}
      \item $\maxSubtreesOf(v')$ is a singleton set and its only element $\maxTree_{v'}$ contains $v$,
      \item $\displaystyle
              \maxVertices_v(\Tree) = \begin{dcases*}
                                        \setOf{v'}, &if $v'$ is a leaf,\\
                                        \bigDisjointUnionOf_{\Tree' \in \secondMaxSubtreesOf(v')} \maxVertices_{v'}(\Tree'), &otherwise,
                                      \end{dcases*}
            $
      \item $\radiusOf{\Tree}_v = \distanceOf(v, v') + \max_{\Tree' \in \secondMaxSubtreesOf(v')}\radiusOf{\Tree'}_{v'}$, and
      \item when a signal of speed $1$ is sent from $v'$ towards $v$ at time $\distanceOf(\general, v') + 2 \cdot \max_{\Tree' \in \secondMaxSubtreesOf(v')}\radiusOf{\Tree'}_{v'}$, it reaches $v$ at time $\distanceOf(\general, v) + 2 \cdot \radiusOf{\Tree}_v$, where in the case that $\secondMaxSubtreesOf(v')$ is empty, we define $\max_{\Tree' \in \secondMaxSubtreesOf(v')}\radiusOf{\Tree'}_{v'}$ as $0$. \qedhere
    \end{aenumerate}
  \end{lemma}

  \begin{proof}
    The first item is evident, the second and third follow from it, and the fourth follows from the third with $\distanceOf(\general, v') = \distanceOf(\general, v) + \distanceOf(v, v')$ and the fact that the signal needs the time span $\distanceOf(v, v')$ to traverse the edge from $v'$ to $v$.
  \end{proof}

  The set of all non-leaf vertices whose unique maximum-radius tree contains the general is named in

  \begin{definition} 
    The set of all non-leaf vertices $v$ of $\Graph$ such that $\maxSubtreesOf(v)$ is a singleton set and its only element $\maxTree_v$ contains the vertex $\general$, is denoted by $V_{\general}$\graffito{set $V_{\general}$ of vertices}\index[symbols]{Vgsubscript@$V_{\general}$}.
  \end{definition}

  \begin{remark}
    For each vertex $v \in V_{\general}$, each tree $\Tree \in \subtreesOf(v) \smallsetminus \setOf{\maxTree_v}$, and each vertex $v'$ of $\Tree$, we have $v' \in V_{\general}$. And, for each vertex $v \in V_{\general}$, the set $\secondMaxSubtreesOf(v)$ is non-empty and, for each tree $\Tree \in \secondMaxSubtreesOf(v)$, we have $\radiusOf{\Tree}_v = \max_{\Tree \in \secondMaxSubtreesOf(v)}\radiusOf{\Tree}_v$.
  \end{remark}

  When and why non-leaf vertices whose unique maximum-radius tree contains the general send \emph{marked} reflected/slowed-down find-midpoint signals is shown in

  \begin{lemma}
  \label{lemma:when-do-non-leaf-vertices-whose-unique-maximum-radius-tree-contains-the-general-sent-marked-signals}
    For each vertex $v \in V_{\general}$, at time $\distanceOf(\general, v) + 2 \cdot \max_{\Tree \in \secondMaxSubtreesOf(v)}\radiusOf{\Tree}_v$, 
    \begin{aenumerate}
      \item the count signal at $v$, before it is updated, has the memory $\directionOf(v) \smallsetminus \directionOf_v(\secondMaxSubtreesOf(v) \cup \setOf{\maxTree_v})$,
      \item for each tree $\Tree \in \secondMaxSubtreesOf(v)$, each leaf $\cev{v} \in \maxVertices_v(\Tree)$, and each vertex $w \in \Vertices_{\general, v}$, a marked reflected find-midpoint signal with origin $w$ and reflection vertex $\cev{v}$ reaches $v$ from direction $\directionOf_v(\Tree)$, and
      \item no other marked reflected find-midpoint signal reaches $v$, 
    \end{aenumerate}
    whereupon
    \begin{aenumerate}
      \item the count signal at $v$, after it is updated, has the memory $\directionOf(v) \smallsetminus \setOf{\directionOf_v(\maxTree_v)}$,
      \item for each tree $\Tree \in \secondMaxSubtreesOf(v)$, each leaf $\cev{v} \in \maxVertices_v(\Tree)$, and each vertex $w \in \Vertices_{\general, v} \smallsetminus \setOf{v}$, a marked reflected find-midpoint signal with origin $w$ and reflection vertex $\cev{v}$ is sent from $v$ in the direction $\directionOf_v(\maxTree_v)$ towards $w$, and no other marked reflected find-midpoint signal with origin $w$ and reflection vertex $\cev{v}$ is sent from $v$, and,
      \item for each tree $\Tree \in \secondMaxSubtreesOf(v)$, each leaf $\cev{v} \in \maxVertices_v(\Tree)$, and each direction $d \in \directionOf(v) \smallsetminus \setOf{\directionOf_v(\Tree)}$, a marked slowed-down find-midpoint signal with origin $v$ and reflection vertex $\cev{v}$ is sent from $v$ in direction $d$ (note that $\directionOf_v(\maxTree_v) \in \directionOf(v) \smallsetminus \setOf{\directionOf_v(\Tree)}$), and, no other marked slowed-down find-midpoint signal with origin $v$ and reflection vertex $\cev{v}$ is sent from $v$. 
    \end{aenumerate}
    And before time $\distanceOf(\general, v) + 2 \cdot \max_{\Tree \in \secondMaxSubtreesOf(v)}\radiusOf{\Tree}_v$, marked signals may reach $v$ but no marked signal is sent from $v$, and after that time, no marked signals as above are sent from $v$. 
  \end{lemma}

  \begin{proof}
    We prove this by induction on $n_v = \max_{(W, F) \in \subtreesOf(v) \smallsetminus \setOf{\maxTree_v}}\cardinalityOf{F}$, for $v \in V_{\general}$. For brevity though, we only treat the existence of signals and not their absence. 

    \proofPart{Base Case (see \cref{figure:when-does-count-signal-memorise-the-penultimate-direction:base-case})}
      Let $v \in V_{\general}$ such that $n_v = 1$. And, let $\Tree \in \subtreesOf(v) \smallsetminus \setOf{\maxTree_v}$. Then, $\Tree$ consists of one edge $e$ whose one end is $v$, whose other end is a leaf $v'$, and whose weight is $\radiusOf{\Tree}_v$. According to \cref{lemma:when-do-leaves-sent-marked-signals}, at time $\distanceOf(\general, v') = \distanceOf(\general, v) + \radiusOf{\Tree}_v$, for each vertex $w \in \Vertices_{\general, v'} \smallsetminus \setOf{v'} = \Vertices_{\general, v}$, a marked reflected find-midpoint signal with origin $w$ and reflection vertex $v'$ is sent from $v'$ towards $v$, in particular, one with origin $v$. These signals reach $v$ at time $\distanceOf(\general, v) + 2 \cdot \radiusOf{\Tree}_v$, whereupon
      \begin{aenumerate}
        \item the count signal at $v$ memorises the direction $\directionOf_v(\Tree)$ (because a marked reflected find-midpoint signal with origin $v$ reached $v$),
        \item for each vertex $w \in \Vertices_{\general, v} \smallsetminus \setOf{v}$, a maybe-marked reflected find-midpoint signal with origin $w$ and reflection vertex $v'$ is sent from $v$ in the direction $\directionOf_v(\maxTree_v)$ towards $w$, and
        \item for each direction $d \in \directionOf(v) \smallsetminus \setOf{\directionOf_v(\Tree)}$, a maybe-marked slowed-down find-midpoint signal with origin $v$ and reflection vertex $v'$ is sent from $v$ in direction $d$.
      \end{aenumerate}
      On the timeline, for the trees of $\subtreesOf(v) \smallsetminus \setOf{\maxTree_v}$ in non-decreasing order with respect to the radius and at the same time for trees with the same radius, the signals reach $v$ and are sent from $v$. For those trees whose radius is less than the second greatest radius among the trees of $\subtreesOf(v)$, which is $\max_{\Tree \in \secondMaxSubtreesOf(v)}\radiusOf{\Tree}_v$, the aforementioned maybe-marked signals that are sent from $v$ are unmarked (because the memory of the count signal, after it is updated, does neither contain the directions of $\directionOf_v(\secondMaxSubtreesOf(v))$ nor the direction $\directionOf_v(\maxTree_v)$). And, for the trees of $\secondMaxSubtreesOf(v)$, the aforementioned maybe-marked signals that are sent from $v$ are marked (because the count signal, after it is updated, has the memory $\directionOf(v) \smallsetminus \setOf{\directionOf_v(\maxTree_v)}$). In conclusion, at time $\distanceOf(\general, v) + 2 \cdot \max_{\Tree \in \secondMaxSubtreesOf(v)}\radiusOf{\Tree}_v$, what is to be proven holds.
      \begin{figure}
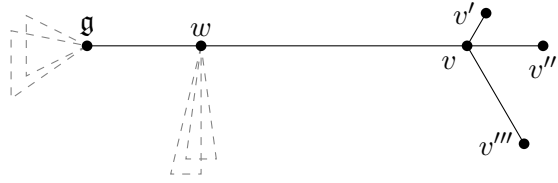

        \myfloatalign
        \figureWhenDoesCountSignalMemoriseThePenultimateDirectionBaseCase
        \caption{Schematic representation of the set-up of the base case of the proof of \cref{lemma:when-do-non-leaf-vertices-whose-unique-maximum-radius-tree-contains-the-general-sent-marked-signals}. The possible existence of subtrees of vertices is hinted at by dashed triangles.}
        \label{figure:when-does-count-signal-memorise-the-penultimate-direction:base-case}
      \end{figure}

    \proofPart{Inductive Step (see \cref{figure:when-does-count-signal-memorise-the-penultimate-direction:inductive-step})}
      Let $v \in V_{\general}$ such that $n_v \geq 2$ and such that what is to be proven holds for each vertex $v' \in V_{\general}$ with $n_{v'} < n_v$, which is the \emph{inductive hypothesis}. And, let $\Tree \in \subtreesOf(v) \smallsetminus \setOf{\maxTree_v}$, let $e$ be the edge of $\Tree$ whose one end is $v$, and let $v'$ be the other end of $e$. The vertex $v'$ is either a leaf or an element of $V_{\general}$ with $n_{v'} < n_v$.
      \begin{figure}
        \myfloatalign
        \figureWhenDoesCountSignalMemoriseThePenultimateDirectionInductiveStep
        \caption{Schematic representation of the set-up of the inductive step of the proof of \cref{lemma:when-do-non-leaf-vertices-whose-unique-maximum-radius-tree-contains-the-general-sent-marked-signals}. The possible existence of subtrees of vertices is hinted at by dashed triangles.}
        \label{figure:when-does-count-signal-memorise-the-penultimate-direction:inductive-step}
      \end{figure}

      In the first case, according to \cref{lemma:when-do-leaves-sent-marked-signals}, at time $\distanceOf(\general, v') = \distanceOf(\general, v) + \distanceOf(v, v')$, for each vertex $w \in \Vertices_{\general, v'} \smallsetminus \setOf{v'} = \Vertices_{\general, v}$, a marked reflected find-midpoint signal with origin $w$ and reflection vertex $v'$ is sent from $v'$ towards $v$; note that $v'$ is the one and only element of $\maxVertices_v(\Tree)$, the set $\secondMaxSubtreesOf(v')$ is empty, and we define $\max_{\Tree' \in \secondMaxSubtreesOf(v')}\radiusOf{\Tree'}_{v'}$ as $0$ (see \cref{lemma:when-does-a-signal-that-is-sent-from-a-vertex-towards-the-general-reach-the-next-vertex}). In the second case, according to the inductive hypothesis, at time $\distanceOf(\general, v') + 2 \cdot \max_{\Tree' \in \secondMaxSubtreesOf(v')}\radiusOf{\Tree'}_{v'}$, for each tree $\Tree' \in \secondMaxSubtreesOf(v')$, each leaf $\cev{v}' \in \maxVertices_{v'}(\Tree')$, and each vertex $w \in \Vertices_{\general, v'} \smallsetminus \setOf{v'} = \Vertices_{\general, v}$, a marked reflected find-midpoint signal with origin $w$ and reflection vertex $\cev{v}'$ is sent from $v'$ towards $v$; note that the set of all vertices $\cev{v}'$ is equal to $\maxVertices_v(\Tree)$ (see \cref{lemma:when-does-a-signal-that-is-sent-from-a-vertex-towards-the-general-reach-the-next-vertex}). 

      In both cases, the marked signals reach $v$ at time $\distanceOf(\general, v) + \distanceOf(v, v') + 2 \cdot \max_{\Tree' \in \secondMaxSubtreesOf(v')}\radiusOf{\Tree'}_{v'} + \distanceOf(v', v)$, which is equal to $\distanceOf(\general, v) + 2 \cdot \radiusOf{\Tree}_v$, whereupon 
      \begin{aenumerate}
        \item the count signal at $v$ memorises the direction $\directionOf_v(\Tree)$ (because at least one marked reflected find-midpoint signal with origin $v$ reached $v$),
        \item for each leaf $\cev{v} \in \maxVertices_v(\Tree)$ and each vertex $w \in \Vertices_{\general, v} \smallsetminus \setOf{v}$, a maybe-marked reflected find-midpoint signal with origin $w$ and reflection vertex $\cev{v}$ is sent from $v$ in the direction $\directionOf_v(\maxTree_v)$ towards $w$, and,
      \item for each leaf $\cev{v} \in \maxVertices_v(\Tree)$ and each direction $d \in \directionOf(v) \smallsetminus \setOf{\directionOf_v(\Tree)}$, a maybe-marked slowed-down find-midpoint signal with origin $v$ and reflection vertex $\cev{v}$ is sent from $v$ in direction $d$.
      \end{aenumerate}
      It follows verbatim as in the base case, that what is to be proven holds.
  \end{proof}

  %

  The set of all vertices whose unique maximum-radius tree does not contain the general is named in

  \begin{definition}
    The set of all vertices $v$ of $\Graph$ such that $\maxSubtreesOf(v)$ is a singleton set and its only element $\maxTree_v$ does \emph{not} contain the vertex $\general$, is denoted by $U_{\general}$\graffito{set $U_{\general}$ of vertices}\index[symbols]{Ugsubscript@$U_{\general}$}.
  \end{definition}

  \begin{remark}
    We have $\general \notin U_{\general}$. And, for each vertex $v \in U_{\general}$, we have $\Vertices_{\general, v} \smallsetminus \setOf{\general} \in U_{\general}$, the vertex $v$ is not a leaf and the set $\secondMaxSubtreesOf(v)$ is non-empty. 
  \end{remark}

  The maximal subtree of a non-general vertex that contains the general is named in

  \begin{definition}
    Let $v$ be a vertex of $\Graph$ such that $v \neq \general$. The tree of $\subtreesOf(v)$ that does contain the vertex $\general$ is denoted by $\Tree_v^{\general}$\graffito{tree $\Tree_v^{\general}$}\index[symbols]{Tvgcalligraphic@$\Tree_v^{\general}$}. And to avoid case differentiations, we define $\Tree_{\general}^{\general}$\graffito{$\Tree_{\general}^{\general} = 0$}\index[symbols]{Tggcalligraphic@$\Tree_{\general}^{\general}$} as the number $0$.
  \end{definition}

  The vertices of the maximal path from a vertex to the general whose vertices excluding its source have the subtree that contains the general as second-maximum-radius subtree is named in

  \begin{definition}
    Let $v$ be a vertex of $U_{\general}$, and let $p$ be the maximum-weight subpath of $p_{\general, v}$ such that $\targetOf(p) = v$ and, for each vertex $w$ on $p$ with $w \neq \sourceOf(p)$, we have $\Tree_w^{\general} \in \secondMaxSubtreesOf(w)$. The set of the vertices on $p$ is denoted by $W_{\general, v}$\graffito{set $W_{\general, v}$ of vertices}\index[symbols]{Wgvsubscript@$W_{\general, v}$}. See \cref{figure:W-general-v}.
    \begin{figure}
      \myfloatalign
      \figureWGeneralV
      \caption{Schematic representation of $\Graph$, where $v \in U_{\general}$, $W_{\general, v} = \setOf{w'', w', w, v}$, $W_{\general, w} = \setOf{w'', w', w}$, $W_{\general, w'} = \setOf{w'', w'}$, $W_{\general, w''} = \setOf{w''}$, and subtrees of vertices and their radii are depicted as dashed triangles of various sizes.}
      \label{figure:W-general-v}
    \end{figure}
  \end{definition}

  \begin{remark}
    We have $v \in W_{\general, v} \smallsetminus \setOf{\general} \subseteq U_{\general}$. And, for each $w \in W_{\general, v}$, we have $W_{\general, w} \subseteq W_{\general, v}$. And, if $\Tree_v^{\general} \in \secondMaxSubtreesOf(v)$, then $W_{\general, v} = W_{\general, v'} \cup \setOf{v}$, where $v'$ is the neighbour of $v$ on $p_{\general, v}$.
  \end{remark}

  When and why vertices whose unique maximum-radius tree does not contain the general send \emph{marked} reflected/slowed-down find-midpoint signals is shown in

  \begin{lemma}
  \label{lemma:when-do-vertices-whose-unique-maximum-radius-tree-does-not-contain-the-general-sent-marked-signals}
    For each vertex $v \in U_{\general}$, at time $\distanceOf(\general, v) + 2 \cdot \max_{\Tree \in \secondMaxSubtreesOf(v)}\radiusOf{\Tree}_v$,
    \begin{aenumerate}
      \item the count signal at $v$, before it is updated, has the memory $\directionOf(v) \smallsetminus \directionOf_v(\secondMaxSubtreesOf(v) \cup \setOf{\maxTree_v})$,
      \item for each tree $\Tree \in \secondMaxSubtreesOf(v)$ and each leaf $\cev{v} \in \maxVertices_v(\Tree)$, a maybe-marked reflected find-midpoint signal with origin $v$ and reflection vertex $\cev{v}$ reaches $v$ from direction $\directionOf_v(\Tree)$, where the reflected signal is marked if and only if $\Tree \neq \Tree_v^{\general}$, and,
      \item for each vertex $w \in W_{\general, v} \smallsetminus \setOf{v}$, each tree $\Tree \in \secondMaxSubtreesOf(w) \smallsetminus \setOf{\Tree_w^{\general}}$, and each leaf $\cev{w} \in \maxVertices_w(\Tree)$, a marked slowed-down find-midpoint signal with origin $w$ and reflection vertex $\cev{w}$ reaches $v$ from the direction towards $w$, which is the direction $\directionOf_v(\Tree_v^{\general})$,
    \end{aenumerate}
    whereupon
    \begin{aenumerate}
      \item the count signal at $v$, after it is updated, has the memory $\directionOf(v) \smallsetminus \setOf{\directionOf_v(\maxTree_v)}$,
      \item for each tree $\Tree \in \secondMaxSubtreesOf(v) \smallsetminus \setOf{\Tree_v^{\general}}$, each leaf $\cev{v} \in \maxVertices_v(\Tree)$, and each direction $d \in \directionOf(v) \smallsetminus \setOf{\directionOf_v(\Tree)}$, a marked slowed-down find-midpoint signal with origin $v$ and reflection vertex $\cev{v}$ is sent from $v$ in direction $d$ (note that $\setOf{\directionOf_v(\maxTree_v), \directionOf_v(\Tree_v^{\general})} \subseteq \directionOf(v) \smallsetminus \setOf{\directionOf_v(\Tree)}$), and,
      \item for each vertex $w \in W_{\general, v} \smallsetminus \setOf{v}$, each tree $\Tree \in \secondMaxSubtreesOf(w) \smallsetminus \setOf{\Tree_w^{\general}}$, each leaf $\cev{w} \in \maxVertices_w(\Tree)$, and each direction $d \in \directionOf(v) \smallsetminus \setOf{\directionOf_v(\Tree_v^{\general})}$, a marked slowed-down find-midpoint signal with origin $w$ and reflection vertex $\cev{w}$ is sent from $v$ in direction $d$ (note that $\directionOf_v(\maxTree_v) \in \directionOf(v) \smallsetminus \setOf{\directionOf_v(\Tree_v^{\general})}$ and that, if $w = \general$, then $\secondMaxSubtreesOf(w) \smallsetminus \setOf{\Tree_w^{\general}} = \secondMaxSubtreesOf(w)$). \qedhere
    \end{aenumerate} 
  \end{lemma}

  \begin{proof}
    We prove this by induction on $n_v = \cardinalityOf{\Vertices_{\general, v}}$, for $v \in U_{\general}$.

    \proofPart{Base Case (compare \cref{figure:marked-signals-inner-vertices:inductive-step})} 
      Let $v \in U_{\general}$ such that $n_v = 1$.

      First, for each tree $\Tree \in \subtreesOf(v) \smallsetminus \setOf{\Tree_v^{\general}, \maxTree_v}$, according to \cref{lemma:when-do-leaves-sent-marked-signals}, if the neighbour of $v$ in $\Tree$ is a leaf, or \cref{lemma:when-do-non-leaf-vertices-whose-unique-maximum-radius-tree-contains-the-general-sent-marked-signals}, otherwise, and \cref{lemma:when-does-a-signal-that-is-sent-from-a-vertex-towards-the-general-reach-the-next-vertex}, at time $\distanceOf(\general, v) + 2 \cdot \radiusOf{\Tree}_v$, for each leaf $\cev{v} \in \maxVertices_v(\Tree)$, a marked reflected find-midpoint signal with origin $v$ and reflection vertex $\cev{v}$ reaches $v$ from direction $\directionOf_v(\Tree)$, whereupon
      \begin{aenumerate}
        \item the count signal at $v$ memorises the direction $\directionOf_v(\Tree)$,
        \item for each leaf $\cev{v} \in \maxVertices_v(\Tree)$ and each direction $d \in \directionOf(v) \smallsetminus \setOf{\directionOf_v(\Tree)}$, a maybe-marked slowed-down find-midpoint signal with origin $v$ and reflection vertex $\cev{v}$ is sent from $v$ in direction $d$.
      \end{aenumerate}

      Secondly, according to \cref{remark:how-initiate-and-find-midpoint-signals-spread}, at time $\distanceOf(\general, v) + 2 \cdot \radiusOf{\Tree_v^{\general}}_v$, for each leaf $\cev{v} \in \maxVertices_v(\Tree_v^{\general})$, a maybe-marked reflected find-midpoint signal with origin $v$ and reflection vertex $\cev{v}$ reaches $v$ from direction $\directionOf_v(\Tree_v^{\general})$. This signal is actually unmarked, because it was reflected at $\cev{v}$ at time $\distanceOf(\general, v) + \radiusOf{\Tree_v^{\general}}_v = \distanceOf(\general, v) + \distanceOf(v, \cev{v})$, which, because $v \neq \general$ and $\cev{v} \in \Tree_v^{\general}$, is greater than the only time, namely $\distanceOf(\general, \cev{v})$, at which an initiate signal reaches $\cev{v}$.

      Thirdly, because $U_{\general}$ is non-empty, we have $\general \in V_{\general}$ and $v \in \maxTree_{\general}$. According to \cref{lemma:when-do-non-leaf-vertices-whose-unique-maximum-radius-tree-contains-the-general-sent-marked-signals}, at time $2 \cdot \max_{\Tree \in \secondMaxSubtreesOf(\general)}\radiusOf{\Tree}_{\general}$, for each tree $\Tree \in \secondMaxSubtreesOf(\general)$, and each leaf $\cev{\general} \in \maxVertices_{\general}(\Tree)$, a marked slowed-down find-midpoint signal with origin $\general$ and reflection vertex $\cev{\general}$ is sent from $\general$ towards $v$; note that the set of all vertices $\cev{\general}$ is equal to $\maxVertices_v(\Tree_v^{\general})$ (compare \cref{lemma:when-does-a-signal-that-is-sent-from-a-vertex-towards-the-general-reach-the-next-vertex}). The marked signals reach $v$ from the direction towards $\general$, which is the direction $\directionOf_v(\Tree_v^{\general})$, at time $2 \cdot \max_{\Tree \in \secondMaxSubtreesOf(\general)}\radiusOf{\Tree}_{\general} + 3 \cdot \distanceOf(\general, v) = \distanceOf(\general, v) + 2 \cdot \radiusOf{\Tree_v^{\general}}_v$, whereupon
      \begin{aenumerate}
        \item the count signal at $v$ memorises the direction $\directionOf_v(\Tree_v^{\general})$ and,
        \item for each leaf $\cev{\general} \in \maxVertices_v(\Tree_v^{\general})$ and each direction $d \in \directionOf(v) \smallsetminus \setOf{\directionOf_v(\Tree_v^{\general})}$, a maybe-marked slowed-down find-midpoint signal with origin $\general$ and reflection vertex $\cev{\general}$ is sent from $v$ in direction $d$.
      \end{aenumerate} 

      Altogether, on the timeline, for the trees of $\subtreesOf(v) \smallsetminus \setOf{\maxTree_v}$ in non-decreasing order with respect to the radius and at the same time for trees with the same radius, the signals reach $v$ and are sent from $v$. For those trees whose radius is less than the second greatest radius among the trees of $\subtreesOf(v)$, which is $\max_{\Tree \in \secondMaxSubtreesOf(v)}\radiusOf{\Tree}_v$, the aforementioned maybe-marked signals that are sent from $v$ are unmarked (because the memory of the count signal, after it is updated, does neither contain the directions of $\directionOf_v(\secondMaxSubtreesOf(v))$ nor the direction $\directionOf_v(\maxTree_v)$). And, for the trees of $\secondMaxSubtreesOf(v)$, the aforementioned maybe-marked signals that are sent from $v$ are marked (because the count signal, after it is updated, has the memory $\directionOf(v) \smallsetminus \setOf{\directionOf_v(\maxTree_v)}$). Note that, if $\Tree_v^{\general} \in \secondMaxSubtreesOf(v)$, then $W_{\general, v} \smallsetminus \setOf{v} = \setOf{\general}$, and otherwise, $W_{\general, v} \smallsetminus \setOf{v} = \emptyset$. In conclusion, at time $\distanceOf(\general, v) + 2 \cdot \max_{\Tree \in \secondMaxSubtreesOf(v)}\radiusOf{\Tree}_v$, what is to be proven holds. 

    \proofPart{Inductive Step (see \cref{figure:marked-signals-inner-vertices:inductive-step})}
      Let $v \in U_{\general}$ such that $n_v \geq 2$ and such that what is to be proven holds for each vertex $v' \in U_{\general}$ with $n_{v'} < n_v$, which is the \emph{inductive hypothesis}.
      \begin{figure}
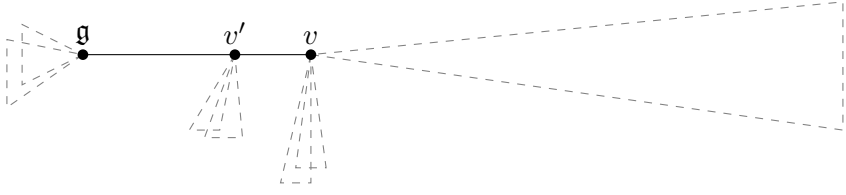

        \myfloatalign
        \figureMarkedSignalsInnerVerticesInducitveStep
        \caption{Schematic representation of the set-up of the inductive step of the proof of \cref{lemma:when-do-vertices-whose-unique-maximum-radius-tree-does-not-contain-the-general-sent-marked-signals}. The possible existence of subtrees of vertices and their radii is hinted at by dashed triangles of various sizes.}
        \label{figure:marked-signals-inner-vertices:inductive-step}
      \end{figure}

      First, for each tree $\Tree \in \subtreesOf(v) \smallsetminus \setOf{\Tree_v^{\general}, \maxTree_v}$, the same as in the base case happens.

      Secondly, as in the base case, at time $\distanceOf(\general, v) + 2 \cdot \radiusOf{\Tree_v^{\general}}_v$, for each leaf $\cev{v} \in \maxVertices_v(\Tree_v^{\general})$, an unmarked reflected find-midpoint signal with origin $v$ and reflection vertex $\cev{v}$ reaches $v$ from direction $\directionOf_v(\Tree_v^{\general})$.

      Thirdly, let $v'$ be the neighbour of $v$ in $\Tree_v^{\general}$. Then, $v' \in \Vertices_{\general, v} \smallsetminus \setOf{\general} \subseteq U_{\general}$ and $n_{v'} < n_v$. According to the inductive hypothesis, at time $\distanceOf(\general, v') + 2 \cdot \max_{\Tree' \in \secondMaxSubtreesOf(v')}\radiusOf{\Tree'}_{v'}$,
      \begin{aenumerate}
        \item for each tree $\Tree' \in \secondMaxSubtreesOf(v') \smallsetminus \setOf{\Tree_{v'}^{\general}}$ and each leaf $\cev{v}' \in \maxVertices_{v'}(\Tree')$, a marked slowed-down find-midpoint signal with origin $v'$ and reflection vertex $\cev{v}'$ is sent from $v'$ towards $v$, and,
        \item for each vertex $w' \in W_{\general, v'} \smallsetminus \setOf{v'}$, each tree $\Tree' \in \secondMaxSubtreesOf(w') \smallsetminus \setOf{\Tree_{w'}^{\general}}$, each leaf $\cev{w}' \in \maxVertices_{w'}(\Tree')$, a marked slowed-down find-midpoint signal with origin $w'$ and reflection vertex $\cev{w}'$ is sent from $v'$ towards $v$.
      \end{aenumerate}
      The marked signals reach $v$ from from direction $\directionOf_v(\Tree_v^{\general})$ at time $\distanceOf(\general, v') + 2 \cdot \max_{\Tree' \in \secondMaxSubtreesOf(v')}\radiusOf{\Tree'}_{v'} + 3 \cdot \distanceOf(v', v) = \distanceOf(\general, v) + 2 \cdot \radiusOf{\Tree_v^{\general}}_v$, whereupon
      \begin{aenumerate}
        \item the count signal at $v$ memorises the direction $\directionOf_v(\Tree_v^{\general})$ and,
        \item for each vertex $w \in W_{\general, v'}$, each tree $\Tree \in \secondMaxSubtreesOf(w) \smallsetminus \setOf{\Tree_w^{\general}}$, each leaf $\cev{w} \in \maxVertices_w(\Tree)$, and each direction $d \in \directionOf(v) \smallsetminus \setOf{\directionOf_v(\Tree_v^{\general})}$, a maybe-marked slowed-down find-midpoint signal with origin $w$ and reflection vertex $\cev{w}$ is sent from $v$ in direction $d$.
      \end{aenumerate}
      Note that, if $\Tree_v^{\general} \in \secondMaxSubtreesOf(v)$, then $W_{\general, v} \smallsetminus \setOf{v} = W_{\general, v'}$, and otherwise, $W_{\general, v} \smallsetminus \setOf{v} = \emptyset$. With that it follows verbatim as in the base case, that what is to be proven holds. 
  \end{proof}

  That midpoints of maximum-weight direction-preserving paths are recognised as such is shown in

  \begin{theorem} 
  \label{theorem:midpoints-of-maximum-weight-paths-are-recognised}
    Let $\hat{p}$ be a maximum-weight direction-preserving path in $\Graph$, let $\hat{v}$ be the vertex on $\hat{p}$ that is nearest to $\general$, let $\hat{v}_1$ and $\hat{v}_2$ be the two ends of $\hat{p}$ such that $\distanceOf(\hat{v}, \hat{v}_1) \leq \distanceOf(\hat{v}, \hat{v}_2)$, and let $\hat{\midpoint}$ be the midpoint of $\hat{p}$. At time $r + d/2$, at the midpoint $\hat{\midpoint}$,
    \begin{aenumerate}
      \item if $\distanceOf(\hat{v}, \hat{v}_1) = \distanceOf(\hat{v}, \hat{v}_2)$, two marked reflected find-midpoint signals with origin $\hat{v}$ and reflection vertices $\hat{v}_1$ and $\hat{v}_2$ collide, and
      \item otherwise, a marked slowed-down find-midpoint signal with origin $\hat{v}$ and reflection vertex $\hat{v}_1$ collides with a marked reflected find-midpoint signal with origin $\hat{v}$ and reflection vertex $\hat{v}_2$. \qedhere
    \end{aenumerate}
  \end{theorem}

  \begin{proof-sketch}
    From a broad perspective and ignoring boundary cases the following happens in the given order (see \cref{figure:midpoints-of-maximum-weight-paths-are-recognised:midpoint-is-not-nearest-to-general}). At time $0$, an initiate signal is sent from $\general$ towards $\hat{v}$. At time $\distanceOf(\general, \hat{v})$, this signal reaches $\hat{v}$, whereupon find-midpoint signals with origin $\hat{v}$ are sent from $\hat{v}$ in all directions. At time $\distanceOf(\general, \hat{v}) + \radiusOf{\Tree_{\hat{v}}^{\general}}_{\hat{v}}$, the slowest but unmarked reflected find-midpoint signal with origin $\hat{v}$ returns to $\hat{v}$ from the direction towards $\general$. At time $\distanceOf(\general, \hat{v}) + 2 \cdot \distanceOf(\hat{v}, \hat{v}_1)$, the slowest and marked reflected find-midpoint signal with origin $\hat{v}$ returns to $\hat{v}$ from the direction towards $\hat{v}_1$, whereupon a marked slowed-down find-midpoint signal with origin $\hat{v}$ is sent towards $\hat{v}_2$. At time $\distanceOf(\general, \hat{v}) + 2 \cdot \distanceOf(\hat{v}, \hat{v}_1) + 3 \cdot \distanceOf(\hat{v}, \hat{\midpoint}) = \distanceOf(\general, \hat{\midpoint}) + 2 \cdot \distanceOf(\hat{v}_1, \hat{\midpoint}) = r + d/2$, this signal reaches $\hat{\midpoint}$. And, at the same time, which is equal to $\distanceOf(\general, \hat{v}) + \distanceOf(\hat{v}, \hat{v}_2) + \distanceOf(\hat{v}_2, \hat{\midpoint})$, the slowest and marked reflected find-midpoint signal with origin $\hat{v}$ that is on its way to return to $\hat{v}$ from the direction towards $\hat{v}_2$ reaches $\hat{\midpoint}$. The two marked signals collide at $\hat{\midpoint}$ recognising it as the midpoint of a maximum-weight path.

    The slowest reflected find-midpoint signal with origin $\hat{v}$ that returns to $\hat{v}$ from the direction towards $\general$ is unmarked, because it reaches the leaf it is reflected at later than the initiate signal and is never marked in the first place. The slowest reflected find-midpoint signal with origin $\hat{v}$ that returns to $\hat{v}$ from the direction towards $\hat{v}_1$ is marked, because it reaches the leaf $\hat{v}_1$ alongside initiate signals and at each vertex on its way back, it is the penultimate marked signal to return (this is essentially \cref{lemma:when-do-non-leaf-vertices-whose-unique-maximum-radius-tree-contains-the-general-sent-marked-signals}). For the same reason, the signal from the direction towards $\hat{v}_2$ is marked when it reaches $\hat{\midpoint}$. And, the slowed-down find-midpoint signal with origin $\hat{v}$ that reaches $\hat{\midpoint}$ from the direction towards $\hat{v}$ is marked, because at each vertex it reaches on its way, it is the penultimate marked signal to do so (this is essentially \cref{lemma:when-do-vertices-whose-unique-maximum-radius-tree-does-not-contain-the-general-sent-marked-signals}). Note that, for each vertex $u$ on the path from $\general$ to $\hat{\midpoint}$ except for $\general$, the direction towards $\general$ is added to the memory of the count signal at $u$, because a marked slowed-down find-midpoint signal whose origin is \emph{not} $u$ reaches $u$ from that direction, and the other directions are added, because marked reflected find-midpoint signals with origin $u$ reach it from those directions; for all other vertices, the latter is the reason the directions are added.
  \end{proof-sketch}

  \begin{proof}
    Let $v_1$ and $v_2$ be the neighbours of $\hat{\midpoint}$ on $\hat{p}$ such that $\distanceOf(v_1, \hat{v}_1) < \distanceOf(v_2, \hat{v}_1)$ (or, equivalently, $\distanceOf(v_1, \hat{v}_2) > \distanceOf(v_2, \hat{v}_2)$).

    First, let $\distanceOf(\hat{v}, \hat{v}_1) = \distanceOf(\hat{v}, \hat{v}_2)$ (see \cref{figure:midpoints-of-maximum-weight-paths-are-recognised:general-lies-on-path-and-midpoint-is-nearest-to-general,figure:midpoints-of-maximum-weight-paths-are-recognised:midpoint-is-nearest-to-general}). Then, $\hat{v} = \hat{\midpoint}$, and, for each index $i \in \setOf{1, 2}$, we have $v_i \neq \general$ (because, if $\general$ lies on $\hat{p}$, then $\general = \hat{v} = \hat{\midpoint} \neq v_i$, and otherwise, $\general$ does not lie on $\hat{p}$ but $v_i$ does), and, $v_i$ is either a leaf or an element of $\Vertices_{\mathfrak{g}}$ (because $\maxSubtreesOf(v_i)$ consists of the tree of $\subtreesOf(v_i)$ that contains $\hat{v}_j$, where $j \in \setOf{1, 2} \smallsetminus \setOf{i}$, and this tree, namely $\hat{\Tree}_{v_i}$, contains $\general$). Hence, according to \cref{lemma:when-do-leaves-sent-marked-signals} or \cref{lemma:when-do-non-leaf-vertices-whose-unique-maximum-radius-tree-contains-the-general-sent-marked-signals}, for each index $i \in \setOf{1, 2}$, at time $\distanceOf(\general, v_i) + 2 \cdot \distanceOf(v_i, \hat{v}_i)$, a marked reflected find-midpoint signal with origin $\hat{\midpoint}$ and reflection vertex $\hat{v}_i$ is sent from $v_i$ towards $\hat{\midpoint}$, and at time $\distanceOf(\general, \hat{\midpoint}) + 2 \cdot \distanceOf(\hat{\midpoint}, \hat{v}_i) = r + d/2$, it reaches $\hat{\midpoint}$, where it collides with the other signal.
    \begin{figure}
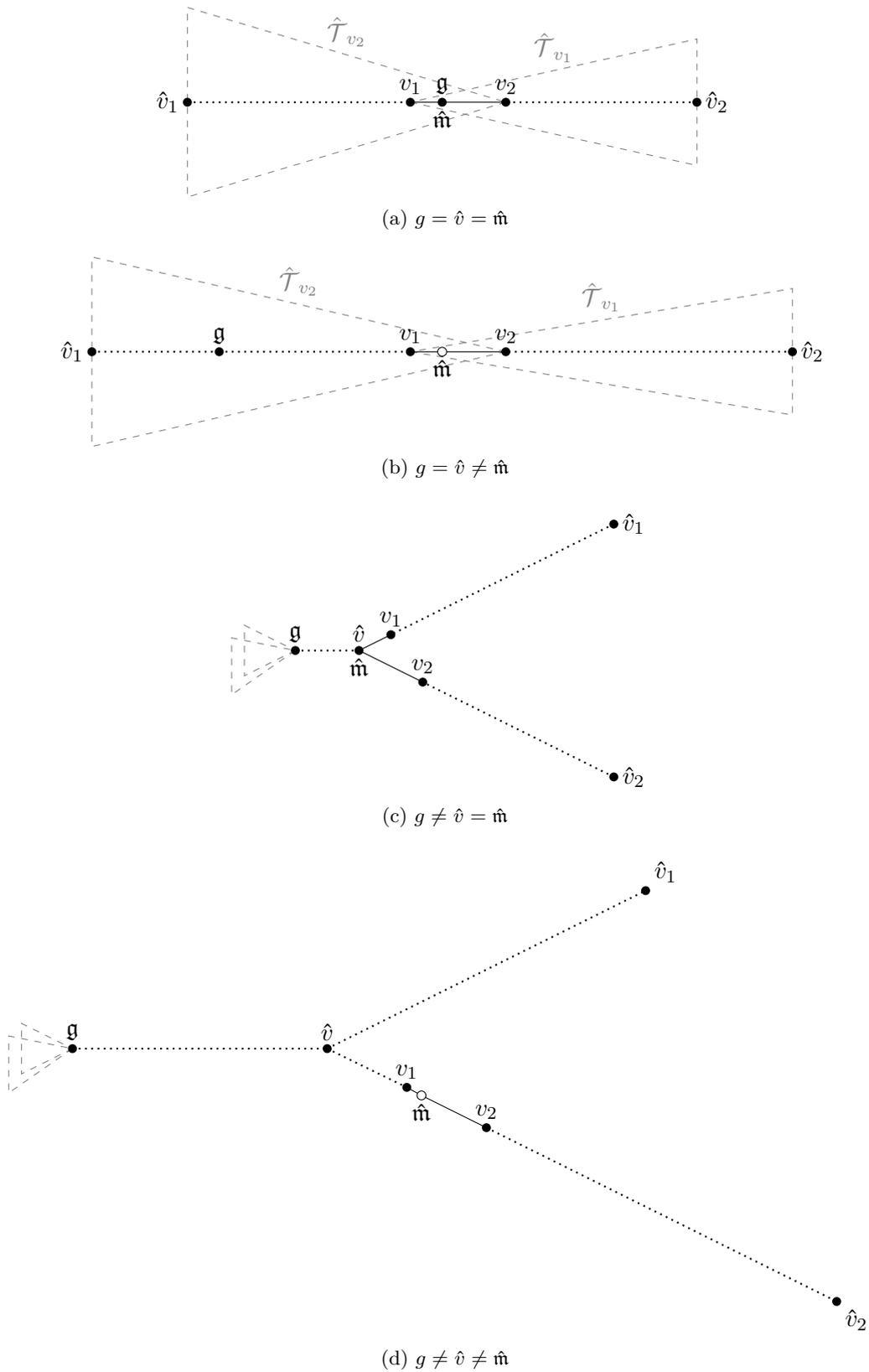

      \myfloatalign
      \begin{wide}
        \figureMidpointsOfMaximumWeightPathsAreRecognised
        \caption{Schematic representation of the set-up of the proof of \cref{theorem:midpoints-of-maximum-weight-paths-are-recognised}. Vertices are depicted as dots, points that may or may not be vertices are depicted as circles, paths that may consist of more than one edge are drawn dotted, and paths that consist of precisely one edge are drawn solid.}
        \label{figure:midpoints-of-maximum-weight-paths-are-recognised}
      \end{wide}
    \end{figure}

    Secondly, let $\distanceOf(\hat{v}, \hat{v}_1) \neq \distanceOf(\hat{v}, \hat{v}_2)$ (see \cref{figure:midpoints-of-maximum-weight-paths-are-recognised:general-lies-on-path-and-midpoint-is-not-nearest-to-general,figure:midpoints-of-maximum-weight-paths-are-recognised:midpoint-is-not-nearest-to-general}). Furthermore, let $X_{\general}$ be the union of the set of leaves and the set $V_{\general}$. Then, $\hat{v} \neq \hat{\midpoint}$, and either $v_1 = \general = \hat{v} \in X_{\general}$ or $v_1 \in U_{\general}$, and $v_2 \in X_{\general}$ (because otherwise, there would be a direction-preserving path with more weight than $\hat{p}$, which contradicts that $\hat{p}$ has maximum-weight).

    According to \cref{lemma:when-do-leaves-sent-marked-signals} or \cref{lemma:when-do-non-leaf-vertices-whose-unique-maximum-radius-tree-contains-the-general-sent-marked-signals}, at time $\distanceOf(\general, v_2) + 2 \cdot \distanceOf(v_2, \hat{v}_2)$, a marked reflected find-midpoint signal with origin $\hat{v}$ and reflection vertex $\hat{v}_2$ is sent from $v_2$ towards $\hat{\midpoint}$, and at time $\distanceOf(\general, v_2) + 2 \cdot \distanceOf(v_2, \hat{v}_2) + \distanceOf(v_2, \hat{\midpoint}) = r + d/2$ it reaches $\hat{\midpoint}$.

    If $v_1 = \general = \hat{v} \in X_{\general}$, then, according to \cref{lemma:when-do-leaves-sent-marked-signals} or \cref{lemma:when-do-non-leaf-vertices-whose-unique-maximum-radius-tree-contains-the-general-sent-marked-signals}, at time $2 \cdot \distanceOf(\hat{v}, \hat{v}_1)$ (which is $0$ in the case that $v_1$ is a leaf, because in this case $v_1 = \hat{v} = \hat{v}_1$), a marked slowed-down find-midpoint signal with origin $\hat{v}$ and reflection vertex $\hat{v}_1$ is sent from $\hat{v}$ towards $\hat{\midpoint}$, and at time $2 \cdot \distanceOf(\hat{v}, \hat{v}_1) + \distanceOf(\hat{v}, \hat{\midpoint}) = r + d/2$, it reaches $\hat{\midpoint}$.

    Otherwise, we have $\hat{v} \in W_{\general, v_1} \smallsetminus \setOf{v_1}$ (because, for each vertex $v \in \Vertices_{\hat{v}, v_1} \smallsetminus \setOf{\hat{v}}$, the set $\maxSubtreesOf(v)$ consists of the tree of $\subtreesOf(v)$ that contains $\hat{v}_2$, the tree $\Tree_v^{\general}$ contains $\hat{v}_1$, and due to the maximality of $\hat{p}$, except for $\hat{\Tree}_v$, there can be no tree of $\subtreesOf(v)$ with a greater radius than $\Tree_v^{\general}$), and hence, according to \cref{lemma:when-do-vertices-whose-unique-maximum-radius-tree-does-not-contain-the-general-sent-marked-signals}, at time $\distanceOf(\general, v_1) + 2 \cdot \distanceOf(v_1, \hat{v}_1)$, a marked slowed-down find-midpoint signal with origin $\hat{v}$ and reflection vertex $\hat{v}_1$ is sent from $v_1$ towards $\hat{\midpoint}$, and at time $\distanceOf(\general, v_1) + 2 \cdot \distanceOf(v_1, \hat{v}_1) + \distanceOf(v_1, \hat{\midpoint}) = r + d/2$, it reaches $\hat{\midpoint}$.

    In either case, at time $r + d/2$, at the midpoint $\hat{\midpoint}$, a marked slowed-down find-midpoint signal with origin $\hat{v}$ and reflection vertex $\hat{v}_1$ collides with a marked reflected find-midpoint signal with origin $\hat{v}$ and reflection vertex $\hat{v}_2$.
  \end{proof}

  \begin{remark} 
    On the other hand, for each direction-preserving path that does \emph{not} have maximum-weight, whenever its midpoint is found, one of the signals is unmarked and hence the midpoint is not falsely thought to be one of a maximum-weight path. The interested reader may prove that for herself.
  \end{remark}

  \subsection{Thaw Signals Traverse Midpoints and Thaw Synchronisation of Edges Just in Time}
  \label{subsection:thaw-signals-traverse-midpoints-and-thaw-synchronisation-of-edges-just-in-time}

  The inverse of a path is its traversal from target to source as given in

  \begin{definition}
    Let $p = (v_0, e_1, v_1, \dotsc, e_n, v_n)$ be a path in $\Graph$. The path $\invert(p) = (v_n, e_n, \dotsc, v_1, e_1, v_0)$ is called \define{inverse of $p$}\graffito{inverse $\invert(p)$ of $p$}\index[symbols]{invp@$\invert(p)$}. 
  \end{definition}

  \begin{remark}
    The inverse of an empty path is the empty path itself.
  \end{remark}

  \begin{remark}
    The weight and the midpoint of the inverse of a path is the same as the weight and the midpoint of the path itself.
  \end{remark}

  An undirected path does not know which of its ends is its source and which is its target and it can be represented as in

  \begin{definition}
    Let $\leftrightarrow$\graffito{equivalence relation $\leftrightarrow$ on $\Paths$}\index[symbols]{arrowleftright@$\leftrightarrow$} be the equivalence relation on $\Paths$ given by $p \leftrightarrow \invert(p)$. Each equivalence class $\equivalenceClassOf{p}_\leftrightarrow \in \Paths \modulo {\leftrightarrow}$ is called \graffito{undirected path $\equivalenceClassOf{p}_\leftrightarrow$}\define{undirected path}\index{path!undirected} and the non-negative real number $\weightOf(\equivalenceClassOf{p}_\leftrightarrow) = \weightOf(p)$ is called \define{weight of $\equivalenceClassOf{p}_\leftrightarrow$}\graffito{weight $\weightOf(\equivalenceClassOf{p}_\leftrightarrow)$ of $\equivalenceClassOf{p}_\leftrightarrow$}\index[symbols]{omegapequivalenceclass@$\weightOf(\equivalenceClassOf{p}_\leftrightarrow)$}.
  \end{definition}

  \begin{remark}
    The equivalence class of an empty path is the singleton set that consists of the empty path.
  \end{remark}

  For each path $p$, the set(s) of paths with the same source (or target), less weight, and midpoints on the continuum representation of $p$ are named in

  \begin{definition}
    For each direction-preserving path $p \in \Paths_{\directionPreserving}$, let
    \begin{equation*}
      \Sigma_p = \setOf{p' \in \Paths_{\directionPreserving} \suchThat
                        \sourceOf(p') = \sourceOf(p) \text{, }
                        \weightOf(p') < \weightOf(p) \text{, and }
                        \midpoint_{p'} \in \imageOf\continuumRepresentationOf{p}} \mathnote{set $\Sigma_p$}\index[symbols]{Sigmap@$\Sigma_p$}
    \end{equation*}
    and let
    \begin{equation*}
      T_p = \setOf{p' \in \Paths_{\directionPreserving} \suchThat
                   \targetOf(p') = \targetOf(p) \text{, }
                   \weightOf(p') < \weightOf(p) \text{, and }
                   \midpoint_{p'} \in \imageOf\continuumRepresentationOf{p}}. \mathnote{set $T_p$}\index[symbols]{Tp@$T_p$} \qedhere
    \end{equation*}
  \end{definition}

  \begin{remark}
    Note that $T_p = \invert(\Sigma_{\invert(p)})$.
  \end{remark}

  \begin{remark}
    Let $p$ be a direction-preserving path of $\Graph$. For each path $p' \in \Sigma_p \cup T_p$, we have $\distanceOf(\midpoint_p, \midpoint_{p'}) = \weightOf(p)/2 - \weightOf(p')/2$.
  \end{remark}

  The set of undirected direction-preserving paths can be turned into a graph by having an edge from an undirected direction-preserving path to each such path that has one end in common with the path, has less weight than the path, and has the greatest weight among the paths with the former two properties. The edges can be weighted by the distance of the midpoints of the undirected paths. The graph and its edge-weighting are introduced in

  \begin{definition}
    Let $\Vertices_{\thawKind} = \Paths_{\directionPreserving} \modulo {\leftrightarrow}$, let
    \begin{equation*}
      \Edges_{\thawKind} = \setOf{(\equivalenceClassOf{p}_\leftrightarrow, \equivalenceClassOf{p'}_\leftrightarrow) \in \Vertices_{\thawKind} \times \Vertices_{\thawKind} \suchThat p' \in \argmax_{p'' \in \Sigma_p \cup T_p} \weightOf(p'')},
    \end{equation*} 
    and let
    \begin{align*}
      \weightOf_{\thawKind} \from \Edges_{\thawKind} &\to \R_{\geq 0},\\
      (\equivalenceClassOf{p}_\leftrightarrow, \equivalenceClassOf{p'}_\leftrightarrow) &\mapsto \distanceOf(\midpoint_p, \midpoint_{p'}).
    \end{align*}
    The triple $\Graph_{\thawKind} = \ntuple{\Vertices_{\thawKind}, \Edges_{\thawKind}, \weightOf_{\thawKind}}$\graffito{edge-weighted directed acyclic graph $\Graph_{\thawKind}$}\index[symbols]{GTcalligraphictypewriter@$\Graph_{\thawKind}$} is an edge-weighted directed acyclic graph.
  \end{definition}

  \begin{remark}
    For each maximum-weight direction-preserving path $\hat{p}$ in $\Graph$, the in-degree of $\equivalenceClassOf{\hat{p}}_\leftrightarrow$ in $\Graph_{\thawKind}$ is $0$, because there are only edges to equivalence classes of less-weight paths.
  \end{remark}

  \begin{remark}
  \label{remark:weight-of-edge-is-half-weight-of-longer-path-minus-half-weight-of-shorter-path}
    For each edge $(\equivalenceClassOf{p}_\leftrightarrow, \equivalenceClassOf{p'}_\leftrightarrow)$ of $\Graph_{\thawKind}$, we have $\weightOf_{\thawKind}(\equivalenceClassOf{p}_\leftrightarrow,\allowbreak \equivalenceClassOf{p'}_\leftrightarrow) = \weightOf(p)/2 - \weightOf(p')/2$.
  \end{remark}

  \begin{remark}
    There is a bijection between the vertices of $\Vertices_{\thawKind}$ and the set of all midpoint signals that are created by the signal machine $\mathcal{S}$. The reason is that each midpoint signal memorises two words of directions in a set, one that leads from its position to one end of the path it designates the midpoint of and the other to the other end; because these words are stored in a set, the midpoint signal does not differentiate between source and target of its path.

    Each maximum-weight vertex of $\Graph_{\thawKind}$ is a maximum-weight undirected direction-preserving path in $\Graph$ and is, under suitable identifications and definitions, a longest undirected direction-preserving path in the continuum representation $M$ of $\Graph$. And, each minimum-weight vertex of $\Graph_{\thawKind}$ is one of weight $0$, is an undirected empty path in $\Graph$, and is, under suitable identifications, an empty path in $\Graph$ and in $M$, and a vertex in $\Graph$ and in $M$.

    For each path $p_{\thawKind}$ in $\Graph_{\thawKind}$ that starts in a maximum-weight vertex and ends in a minimum-weight vertex, in the time evolution of the signal machine $\mathcal{S}$, there is a thaw signal that traverses the midpoint signals represented by the vertices on the path in the order they occur on the path such that the time the thaw signal takes to get from the vertex $v_{\thawKind}$ on the path to the next vertex $v_{\thawKind}'$ on the path is precisely the weight of the edge $(v_{\thawKind}, v_{\thawKind}')$, in particular, the time the thaw signal takes to reach the end of its path is the weight of $p_{\thawKind}$.

    To show that the synchronisation of all edges is thawed and finishes at the same time, we have to show that, for each edge, there is a thaw signal that collides with the midpoint signal of the edge and reaches the end of its path at one end of the edge, and that all thaw signals reach the ends of their paths at the same time. The former is equivalent to showing that, for each edge $e \in \Edges$ with the two ends $v_0$ and $v_1$, there is a maximum-weight direction-preserving path $p$ in $\Graph$ such that there is a path in $\Graph_{\thawKind}$ from $\equivalenceClassOf{p}_{\leftrightarrow}$ to $\equivalenceClassOf{(v_0, e, v_1)}_{\leftrightarrow}$. And the latter is equivalent to showing that the weights of all paths from maximum-weight vertices to minimum-weight vertices in $\Graph_{\thawKind}$ are the same. See \cref{corollary:the-weight-of-paths-to-edges-and-vertices-in-the-thaw-graph}.

    For each non-empty direction-preserving path $p$ in $\Graph$, the midpoint of $p$ is found at time $\max\setOf{\distanceOf(\general, \sourceOf(p)),\allowbreak \distanceOf(\general, \targetOf(p))} + \weightOf(p)/2$ (see \cref{lemma:the-time-at-which-the-midpoint-of-a-path-is-found}) and, for each maximum-weight direction-preserving path $\hat{p}$ in $\Graph$, the midpoint of $\hat{p}$ is found at time $r + d/2$ (see \cref{lemma:the-time-at-which-the-midpoints-of-longest-paths-are-found}). If the thaw signals that emanate from the midpoints of maximum-weight direction-preserving paths at the time $r + d/2$ reach the ends of their paths after the time span $d/2$, then synchronisation finishes at the optimal time $r + d$. The condition is equivalent to showing that the weights of all paths from maximum-weight vertices to minimum-weight vertices in $\Graph_{\thawKind}$ are equal to $d/2$. See \cref{corollary:the-weight-of-paths-to-edges-and-vertices-in-the-thaw-graph}.
  \end{remark}

  The thaw signals that spread from the midpoint signal of a path $p$ eventually collide with the midpoint signals of all paths that have the same source or target as $p$, less weight, and whose midpoints lie on $p$. This is what is essentially shown in

  \begin{lemma}
  \label{lemma:existence-of-path-in-midpoint-signal-graph}
    Let $p$ be a path in $\Graph$ and let $p' \in \Sigma_p \cup T_p$. There is a path from $\equivalenceClassOf{p}_\leftrightarrow$ to $\equivalenceClassOf{p'}_\leftrightarrow$ in $\Graph_{\thawKind}$.
  \end{lemma}

  \begin{proof}
    \begin{description}
      \item[Case 1:] $p' \in \Sigma_p$. If $p' \in \argmax_{p'' \in \Sigma_p} \weightOf(p'')$, then $(\equivalenceClassOf{p}_\leftrightarrow,\allowbreak \equivalenceClassOf{p'}_\leftrightarrow) \in \Edges_{\thawKind}$ and the path consisting of this edge is one from $\equivalenceClassOf{p}_\leftrightarrow$ to $\equivalenceClassOf{p'}_\leftrightarrow$. Otherwise, there is a $q \in \argmax_{p'' \in \Sigma_p} \weightOf(p'')$ such that $\weightOf(p') < \weightOf(q)$. Then, $(\equivalenceClassOf{p}_\leftrightarrow, \equivalenceClassOf{q}_\leftrightarrow) \in \Edges_{\thawKind}$. And, because $\setOf{p', q} \subseteq \Sigma_p$, we have $\sourceOf(p') = \sourceOf(p) = \sourceOf(q)$. Thus, because $\continuumRepresentationOf{q}(\weightOf(q) / 2) = \midpoint_q \in \imageOf \continuumRepresentationOf{p}$, the paths $q$ and $p$ are direction-preserving, and the multigraph $\Graph$ is a tree, we have $\continuumRepresentationOf{q}\restrictedTo_{\closedInterval{0, \weightOf(q) / 2}} = \continuumRepresentationOf{p}\restrictedTo_{\closedInterval{0, \weightOf(q) / 2}}$ and, analogously, we have $\continuumRepresentationOf{p'}\restrictedTo_{\closedInterval{0, \weightOf(p') / 2}} = \continuumRepresentationOf{p}\restrictedTo_{\closedInterval{0, \weightOf(p') / 2}}$. Hence, because $\weightOf(p') < \weightOf(q)$,
            \begin{equation*}
              \midpoint_{p'}
              = \continuumRepresentationOf{p'}(\weightOf(p') / 2)
              = \continuumRepresentationOf{p}(\weightOf(p') / 2)
              = \continuumRepresentationOf{q}(\weightOf(p') / 2)
              \in \imageOf\continuumRepresentationOf{q}.
            \end{equation*}
            Therefore, $p' \in \Sigma_q$. Now, if $p' \in \argmax_{p'' \in \Sigma_q} \weightOf(p'')$, then $(\equivalenceClassOf{q}_\leftrightarrow,\allowbreak \equivalenceClassOf{p'}_\leftrightarrow) \in \Edges_{\thawKind}$, and the path consisting of the edges $(\equivalenceClassOf{p}_\leftrightarrow, \equivalenceClassOf{q}_\leftrightarrow)$ and $(\equivalenceClassOf{q}_\leftrightarrow, \equivalenceClassOf{p'}_\leftrightarrow)$ is a path from $\equivalenceClassOf{p}_\leftrightarrow$ to $\equivalenceClassOf{p'}_\leftrightarrow$. Otherwise, because $\Vertices_{\thawKind}$ is finite and $\weightOf(q) > \weightOf(p')$, it follows by induction that there is a path from $\equivalenceClassOf{q}_\leftrightarrow$ to $\equivalenceClassOf{p'}_\leftrightarrow$ and therefore one from $\equivalenceClassOf{p}_\leftrightarrow$ to $\equivalenceClassOf{p'}_\leftrightarrow$. 
      \item[Case 2:] $p' \in T_p$. Then, $\invert(p') \in \Sigma_{\invert(p)}$. Hence, according to the first case, there is a path from $\equivalenceClassOf{\invert(p)}_\leftrightarrow$ to $\equivalenceClassOf{\invert(p')}_\leftrightarrow$. Therefore, because $\equivalenceClassOf{p}_\leftrightarrow = \equivalenceClassOf{\invert(p)}_\leftrightarrow$ and $\equivalenceClassOf{p'}_\leftrightarrow = \equivalenceClassOf{\invert(p')}_\leftrightarrow$, there is a path from $\equivalenceClassOf{p}_\leftrightarrow$ to $\equivalenceClassOf{p'}_\leftrightarrow$. \qedhere
    \end{description}
  \end{proof}

  Each midpoint signal of a path eventually collides with a matching thaw signal that originated at the midpoint of a maximum-weight direction-preserving path. This is what is essentially shown in

  \begin{lemma}
  \label{lemma:existence-of-paths-in-the-thaw-graph}
    Let $p$ be a direction-preserving path in $\Graph$. There is a maximum-weight direction-preserving path $\hat{p}$ in $\Graph$ such that there is a path from $\equivalenceClassOf{\hat{p}}_\leftrightarrow$ to $\equivalenceClassOf{p}_\leftrightarrow$ in $\Graph_{\thawKind}$.
  \end{lemma}

  \begin{proof}
    If $p$ is a maximum-weight path, then the path $\hat{p} = p$ in $\Graph$ and the empty path $(\equivalenceClassOf{p}_\leftrightarrow)$ in $\Graph_{\thawKind}$ have the desired properties. From now on, let $p$ not be a maximum-weight path. Furthermore, let $\hat{p}$ be a maximum-weight path in $\Graph$ and let $p_d$ be the minimum-weight path in $\Graph$ such that $\sourceOf(p_d)$ lies on $p$ and $\targetOf(p_d)$ lies on $\hat{p}$. Then, there are paths $p_1$ and $p_2$ in $\Graph$ such that $p_1 \concat p_2 = p$ and $\targetOf(p_1) = \sourceOf(p_d)$ as well as $\targetOf(\invert(p_2)) = \sourceOf(p_d)$. And, there are paths $\hat{p}_1$ and $\hat{p}_2$ in $\Graph$ such that $\hat{p}_1 \concat \hat{p}_2 = \hat{p}$ and $\targetOf(p_d) = \sourceOf(\hat{p}_2)$ as well as $\targetOf(p_d) = \sourceOf(\invert(\hat{p}_1))$. Let
    \begin{equation*}
      q_1 = \begin{dcases*} 
              p_1, &if $\weightOf(p_1) \geq \weightOf(p_2)$,\\
              \invert(p_2), &otherwise,
            \end{dcases*}
    \end{equation*}
    let
    \begin{equation*}
      q_2 = \begin{dcases*} 
              p_2, &if $\weightOf(p_1) \geq \weightOf(p_2)$,\\
              \invert(p_1), &otherwise,
            \end{dcases*}
    \end{equation*}
    and let $q = q_1 \concat q_2$. Furthermore, let
    \begin{equation*}
      \hat{q}_1 = \begin{dcases*}
                    \hat{p}_1, &if $\weightOf(\hat{p}_2) \geq \weightOf(\hat{p}_1)$,\\
                    \invert(\hat{p}_2), &otherwise,
                  \end{dcases*}
    \end{equation*}
    let
    \begin{equation*}
      \hat{q}_2 = \begin{dcases*}
                    \hat{p}_2, &if $\weightOf(\hat{p}_2) \geq \weightOf(\hat{p}_1)$,\\
                    \invert(\hat{p}_1), &otherwise,
                  \end{dcases*}
    \end{equation*}
    and let $\hat{q} = \hat{q}_1 \concat \hat{q}_2$. Moreover, let $p' = q_1 \concat p_d \concat \hat{q}_2$. Then, $q \in \setOf{p, \invert(p)}$ as well as $\hat{q} \in \setOf{\hat{p}, \invert(\hat{p})}$. And, $\weightOf(q_1) \geq \weightOf(q_2)$ as well as $\weightOf(\hat{q}_1) \leq \weightOf(\hat{q}_2)$. And, $\targetOf(p') = \targetOf(\hat{q})$ as well as $\sourceOf(q) = \sourceOf(p')$. And, because $\weightOf(q_1) \geq \weightOf(q_2)$, we have $\midpoint_q \in \imageOf\continuumRepresentationOf{q_1} \subseteq \imageOf\continuumRepresentationOf{p'}$.
    \begin{description}
      \item[Case 1:] $p'$ is a maximum-weight path. Then, because $p$ is not a maximum-weight path, we have $\weightOf(q) = \weightOf(p) < \weightOf(p')$. Hence, because $\sourceOf(q) = \sourceOf(p')$ and $\midpoint_q \in \imageOf\continuumRepresentationOf{p'}$, we have $q \in \Sigma_{p'}$. In conclusion, according to \cref{lemma:existence-of-path-in-midpoint-signal-graph}, there is a path from $\equivalenceClassOf{p'}_\leftrightarrow$ to $\equivalenceClassOf{p}_\leftrightarrow = \equivalenceClassOf{q}_\leftrightarrow$.
      \item[Case 2:] $p'$ is not a maximum-weight path. Then, $\weightOf(q_1 \concat p_d \concat \hat{q}_2) = \weightOf(p') < \weightOf(\hat{q}) = \weightOf(\hat{q}_1 \concat \hat{q}_2)$ and thus $\weightOf(q_1 \concat p_d) < \weightOf(\hat{q}_1)$, in particular, $\weightOf(q_1) < \weightOf(\hat{q}_1)$. Hence, because $\weightOf(q_2) \leq \weightOf(q_1) < \weightOf(\hat{q}_1) \leq \weightOf(\hat{q}_2)$, we have $\weightOf(q) = \weightOf(q_1) + \weightOf(q_2) < \weightOf(q_1) + \weightOf(\hat{q}_2) \leq \weightOf(p')$. And, because $\weightOf(q_1 \concat q_d) < \weightOf(\hat{q}_1) \leq \weightOf(\hat{q}_2)$, we have $\midpoint_{p'} \in \imageOf\continuumRepresentationOf{\hat{q}_2} \subseteq \imageOf\continuumRepresentationOf{\hat{q}}$. And, recall that $\midpoint_{q} \in \imageOf\continuumRepresentationOf{p'}$. Altogether, because $\targetOf(p') = \targetOf(\hat{q})$ and $\sourceOf(q) = \sourceOf(p')$, we have $p' \in T_{\hat{q}}$ and $q \in \Sigma_{p'}$. Therefore, according to \cref{lemma:existence-of-path-in-midpoint-signal-graph}, there is a path from $\equivalenceClassOf{\hat{p}}_\leftrightarrow = \equivalenceClassOf{\hat{q}}_\leftrightarrow$ to $\equivalenceClassOf{p'}_\leftrightarrow$ and there is a path from $\equivalenceClassOf{p'}_\leftrightarrow$ to $\equivalenceClassOf{q}_\leftrightarrow = \equivalenceClassOf{p}_\leftrightarrow$. In conclusion, there is a path from $\equivalenceClassOf{\hat{p}}_\leftrightarrow$ to $\equivalenceClassOf{p}_\leftrightarrow$. \qedhere
    \end{description}
  \end{proof}

  The time it takes a thaw signal from the midpoint of a path to collide with the midpoint signal of a matching path is given by half the former path's length minus half the latter path's length. This is what is essentially shown in

  \begin{lemma}
  \label{lemma:the-weight-of-paths-in-the-thaw-graph}
    For each path $p_{\thawKind}$ in $\Graph_{\thawKind}$, the weight of $p_{\thawKind}$ is equal to $\weightOf(\sourceOf(p_{\thawKind}))/2 - \weightOf(\targetOf(p_{\thawKind}))/2$. 
  \end{lemma}

  \begin{proof}
    We prove the statement by induction. 

    \proofPart{Base Case}
      Each empty path $p_{\thawKind}$ in $\Graph_{\thawKind}$ has weight $0$, has the same source and target vertices, and $\weightOf(\sourceOf(p_{\thawKind}))/2 - \weightOf(\targetOf(p_{\thawKind}))/2$ is equal to $0$ as needed.

    \proofPart{Inductive Step}
      Let $p_{\thawKind} = (\equivalenceClassOf{p_0}_\leftrightarrow, \equivalenceClassOf{p_1}_\leftrightarrow, \dotsc, \equivalenceClassOf{p_n}_\leftrightarrow)$ be a non-empty path in $\Graph_{\thawKind}$ such that the weight of the subpath $(\equivalenceClassOf{p_1}_\leftrightarrow, \dotsc, \equivalenceClassOf{p_n}_\leftrightarrow)$ is equal to $\weightOf(\equivalenceClassOf{p_1}_\leftrightarrow)/2 - \weightOf(\equivalenceClassOf{p_n}_\leftrightarrow)/2$. Then, according to \cref{remark:weight-of-edge-is-half-weight-of-longer-path-minus-half-weight-of-shorter-path}, we have $\weightOf_{\thawKind}(\equivalenceClassOf{p_0}_\leftrightarrow, \equivalenceClassOf{p_1}_\leftrightarrow) = \distanceOf(\midpoint_{p_0}, \midpoint_{p_1}) = \weightOf(p_0)/2 - \weightOf(p_1)/2 = \weightOf(\equivalenceClassOf{p_0}_\leftrightarrow)/2 - \weightOf(\equivalenceClassOf{p_1}_\leftrightarrow)/2$. Hence, the weight of the path $p_{\thawKind}$ is equal to $\weightOf(\equivalenceClassOf{p_0}_\leftrightarrow)/2 - \weightOf(\equivalenceClassOf{p_1}_\leftrightarrow)/2 + \weightOf(\equivalenceClassOf{p_1}_\leftrightarrow)/2 - \weightOf(\equivalenceClassOf{p_n}_\leftrightarrow)/2 = \weightOf(\equivalenceClassOf{p_0}_\leftrightarrow)/2 - \weightOf(\equivalenceClassOf{p_n}_\leftrightarrow)/2$.
  \end{proof}

  The time it takes a thaw signal from a maximum-weight direction-preserving path to collide with the midpoint signal of a matching path is essentially given in

  \begin{corollary}
  \label{corollary:the-weight-of-paths-in-the-thaw-graph}
    For each maximum-weight direction-preserving path $\hat{p}$ in $\Graph$ such that there is a path from $\equivalenceClassOf{\hat{p}}_\leftrightarrow$ to $\equivalenceClassOf{p}_\leftrightarrow$ in $\Graph_{\thawKind}$, the weight of this path is $d/2 - \weightOf(p)/2$.
  \end{corollary}

  \begin{proof}
    This is a direct consequence of \cref{lemma:the-weight-of-paths-in-the-thaw-graph} with the fact that $\weightOf(\hat{p}) = d$.
  \end{proof}

  The time it takes a thaw signal from a maximum-weight direction-preserving path to collide with the midpoint signal of an edge it thaws and to reach the ends of the edge is essentially given in

  \begin{corollary}
  \label{corollary:the-weight-of-paths-to-edges-and-vertices-in-the-thaw-graph}
    Let $e$ be an edge of $\Graph$, and let $v_0$ and $v_1$ be the two ends of $e$. There is a maximum-weight direction-preserving path $\hat{p}$ such that there is a path from $\equivalenceClassOf{\hat{p}}_\leftrightarrow$ to $\equivalenceClassOf{(v_0, e, v_1)}_\leftrightarrow$ and all such paths have weight $d / 2 - \weightOf(e) / 2$, and there are also paths from $\equivalenceClassOf{\hat{p}}_\leftrightarrow$ to $\equivalenceClassOf{(v_0)}_\leftrightarrow$ and to $\equivalenceClassOf{(v_1)}_\leftrightarrow$ and all such paths have weight $d/2$.
  \end{corollary}

  \begin{proof}
    This is a direct consequence of \cref{lemma:existence-of-paths-in-the-thaw-graph} and corollary \ref{corollary:the-weight-of-paths-in-the-thaw-graph}. 
  \end{proof}

  \appendix


  \clearToOddPage
  \chapter{Zorn's Lemma}
  \label{chapter:zorns-lemma}

  \begin{definition} 
    Let $I$ be a set and let $\leq$ be a binary relation on $I$.
    The relation $\leq$ is called \define{preorder on $I$}\graffito{preorder $\leq$ on $I$}\index[symbols]{lessthanorequalto@$\leq$} and the tuple $(I, \leq)$ is called \define{preordered set}\graffito{preordered set $(I, \leq)$} if and only if the relation $\leq$ is reflexive and transitive. 
  \end{definition}

  \begin{definition}
    Let $\leq$ be a preorder on $I$. It is called \graffito{partial order $\leq$ on $I$}\define{partial order on $I$} and the preordered set $(I, \leq)$ is called \define{partially ordered set}\graffito{partially ordered set $(I, \leq)$} if and only if the relation $\leq$ is antisymmetric. 
  \end{definition}

  \begin{definition}
    Let $\leq$ be a partial order on $I$. It is called \graffito{total order $\leq$ on $I$}\define{total order on $I$} and the partially ordered set $(I, \leq)$ is called \graffito{totally ordered set $(I, \leq)$}\define{totally ordered set} if and only if the relation $\leq$ is total. 
  \end{definition}

  \begin{definition} 
    Let $\leq$ be a preorder on $I$ and let $i$ be an element of $I$. The element $i$ is called \define{maximal in $(I, \leq)$}\graffito{maximal in $(I, \leq)$} if and only if
    \begin{equation*}
      \ForEach i' \in I \Holds (i' \geq i \implies i' \leq i). \qedhere
    \end{equation*}
  \end{definition}

  \begin{definition}
    Let $\leq$ be a preorder on $I$ and let $J$ be a subset of $I$. The set $J$ is called \define{chain in $(I, \leq)$}\graffito{chain in $(I, \leq)$} if and only if the restriction of $\leq$ to $J$ is a total order on $J$. 
  \end{definition}

  \begin{lemma}[Zorn's Lemma; Kazimierz Kuratowski, 1922; Max August Zorn, 1935]
  \label{lemma:Zorns-lemma}
    Let $(I, \leq)$ be a preordered set such that each chain in $(I, \leq)$ has an upper bound. Then, there is a maximal element in $(I, \leq)$.
  \end{lemma}

  \begin{proof}
    See section~16 in \cite{halmos:1960} for the proof of the case that $\leq$ is a partial order. The general case follows from this case, because each preorder $\leq$ on $I$ induces a partial order on $I \modulo {\sim}$, where $\sim$ is the equivalence relation on $I$ given by $i \sim i'$ if and only if $i \leq i'$ and $i' \leq i$. 
  \end{proof}

  \clearToOddPage
  \chapter{Topologies}
  \label{chapter:topologies}

  Most of the theory of topologies as presented here may be found in more detail in Appendix~A in the monograph \enquote{\citetitle*{ceccherini-silberstein:coornaert:2010}}\cite{ceccherini-silberstein:coornaert:2010}.

  \section{Topologies}

  \begin{definition}
    Let $X$ be a set and let $\mathcal{T}$ be a set of subsets of $X$. The set $\mathcal{T}$ is called \define{topology on $X$}\graffito{topology $\mathcal{T}$ on $X$}\index[symbols]{Tcalligraphic@$\mathcal{T}$} if and only if
    \begin{aenumerate}
      \item $\setOf{\emptyset, X}$ is a subset of $\mathcal{T}$,
      \item for each family $\family{O_i}_{i \in I}$ of elements in $\mathcal{T}$, the union $\bigcup_{i \in I} O_i$ is an element of $\mathcal{T}$,
      \item for each finite family $\family{O_i}_{i \in I}$ of elements in $\mathcal{T}$, the intersection $\bigcap_{i \in I} O_i$ is an element of $\mathcal{T}$. \qedhere 
    \end{aenumerate}
  \end{definition}

  \begin{definition}
    Let $X$ be a set, and let $\mathcal{T}$ and $\mathcal{T}'$ be two topologies on $X$. The topology $\mathcal{T}$ is called
    \begin{aenumerate}
      \item \define{coarser than $\mathcal{T}'$}\graffito{coarser than $\mathcal{T}'$} if and only if $\mathcal{T} \subseteq \mathcal{T}'$;
      \item \define{finer than $\mathcal{T}'$}\graffito{finer than $\mathcal{T}'$} if and only if $\mathcal{T} \supseteq \mathcal{T}'$. \qedhere
    \end{aenumerate}
  \end{definition}

  \begin{definition} 
    Let $X$ be a set and let $\mathcal{T}$ be a topology on $X$. The tuple $(X, \mathcal{T})$ is called \define{topological space}\graffito{topological space $(X, \mathcal{T})$}, each subset $O$ of $X$ with $O \in \mathcal{T}$ is called \define{open in $X$}\graffito{open set $O$ in $X$}, each subset $A$ of $X$ with $X \smallsetminus A \in \mathcal{T}$ is called \define{closed in $X$}\graffito{closed set $A$ in $X$}, and each subset $U$ of $X$ that is both open and closed is called \define{clopen in $X$}\graffito{clopen set $U$ in $X$}.

    The set $X$ is said to be \define{equipped with $\mathcal{T}$}\graffito{equipped with $\mathcal{T}$} if and only if it shall be implicitly clear that $\mathcal{T}$ is the topology on $X$ being considered. The set $X$ is called \define{topological space}\graffito{topological space $X$} if and only if it is implicitly clear what topology on $X$ is being considered. 
  \end{definition}

  \begin{example}
    Let $X$ be a set. The set $\powerSetOf(X)$ is the finest topology on $X$. Itself as well as the topological space $(X, \powerSetOf(X))$ are called \defineX{discrete}{discrete!topology}\graffito{discrete topology}.
  \end{example}

  \begin{definition}
    Let $X$ be a set, let $(X', \mathcal{T}')$ be a topological space, and let $f$ be a map from $X$ to $X'$. The set
    \begin{equation*}
      \mathcal{T} = \setOf{f^{-1}(O') \suchThat O' \in \mathcal{T}'}
    \end{equation*}
    is a topology on $X$ and called \defineX{induced on $X$ by $f$}{induced on $X$ by $f$!topology}\graffito{induced topology on $X$ by $f$}. 
  \end{definition}

  \begin{definition}
    Let $(X, \mathcal{T})$ be a topological space, let $Y$ be a subset of $X$, and let $\iota$ be the canonical injection from $Y$ to $X$. The topology $\mathcal{S}$ on $Y$ induced by $\iota$ is called \defineX{subspace}{subspace!topology}\graffito{subspace topology} and the topological space $(Y, \mathcal{S})$ is called \define{subspace of $(X, \mathcal{T})$}\graffito{subspace of $(X, \mathcal{T})$}.
  \end{definition}

  \begin{remark}
    Let $(X, \mathcal{T})$ be a topological space and let $Y$ be a subset of $X$. The subspace topology on $Y$ is $\setOf{O \cap Y \suchThat O \in \mathcal{T}}$.
  \end{remark}

  \begin{definition}
    Let $(X, \mathcal{T})$ be a topological space and let $\mathcal{B}$ be a subset of $\mathcal{T}$.
    \begin{aenumerate}
      \item The set $\mathcal{B}$ is called \define{base of $\mathcal{T}$}\graffito{base $\mathcal{B}$ of $\mathcal{T}$}\index[symbols]{Bcalligraphic@$\mathcal{B}$} if and only if
            \begin{equation*}
              \ForEach O \in \mathcal{T} \Exists \family{B_i}_{i \in I} \text{ in } \mathcal{B} \SuchThat \bigcup_{i \in I} B_i = O.
            \end{equation*}
      \item The set $\mathcal{B}$ is called \define{subbase of $\mathcal{T}$}\graffito{subbase $\mathcal{B}$ of $\mathcal{T}$}\index{base of $\mathcal{T}$!sub-}\index[symbols]{Bcalligraphic@$\mathcal{B}$} if and only if
            \begin{equation*}
                \setOf{\bigcap_{i = 1}^n B_i \suchThat B_i \in \mathcal{B}, i \in \setOf{1,2,\dotsc,n}, n \in \N_+}
              \end{equation*}
              is a base of $\mathcal{T}$. \qedhere
    \end{aenumerate}
  \end{definition}

  \begin{definition} 
    Let $(X, \mathcal{T})$ be a topological space, let $x$ be a point of $X$, and let $N$ be a subset of $X$. The set $N$ is called \graffito{neighbourhood $N$ of $x$}\define{neighbourhood of $x$} if and only if there is an open subset $O$ of $X$ such that $x \in O$ and $O \subseteq N$.
  \end{definition}

  \begin{definition}
    Let $(X, \mathcal{T})$ be a topological space, let $x$ be a point of $X$, and let $\mathcal{B}_x$ be a set of neighbourhoods of $x$. The set $\mathcal{B}_x$ is called \define{neighbourhood base of $x$}\graffito{neighbourhood base $\mathcal{B}_x$ of $x$}\index[symbols]{Bxcalligraphic@$\mathcal{B}_x$} if and only if, for each neighbourhood $N$ of $x$, there is a neighbourhood $B_x \in \mathcal{B}_x$ such that $B_x \subseteq N$.
  \end{definition}

  \begin{definition}
    Let $(X, \mathcal{T})$ be a topological space and let $x$ be a point of $X$. The set of open neighbourhoods of $x$ is denoted by $\mathcal{T}_x$\graffito{set $\mathcal{T}_x$ of open neighbourhoods of $x$}\index[symbols]{Txcalligraphic@$\mathcal{T}_x$}.
  \end{definition}

  \section{Nets}

  \begin{definition}
    Let $\leq$ be a preorder on $I$. It is called \defineX{directed}{directed preorder on $I$}\graffito{directed preorder $\leq$ on $I$}\index{preorder on $I$!directed}\index[symbols]{lessthanorequalto@$\leq$} and the preordered set $(I, \leq)$ is called \define{directed set}\graffito{directed set $(I, \leq)$} if and only if
    \begin{equation*}
      \ForEach i \in I \ForEach i' \in I \Exists i'' \in I \SuchThat i \leq i'' \land i' \leq i''. \qedhere
    \end{equation*}
  \end{definition}

  \begin{definition}
    Let $\leq$ be a preorder on $I$, let $J$ be a subset of $I$, and let $i$ be an element of $I$. The element $i$ is called \graffito{upper bound of $J$ in $(I, \leq)$}\define{upper bound of $J$ in $(I, \leq)$} if and only if
    \begin{equation*}
      \ForEach i' \in J \Holds i' \leq i. \qedhere
    \end{equation*}
  \end{definition}

  \begin{remark}
    In the definition of directed sets, the element $i''$ is an upper bound of $\setOf{i, i'}$ in $(I, \leq)$.
  \end{remark}

  \begin{definition}
    Let $I$ be a set, let $\leq$ be a binary relation on $I$, and let $\family{m_i}_{i \in I}$ be a family of elements in $M$ indexed by $I$. The family $\family{m_i}_{i \in I}$ is called \define{net in $M$ indexed by $(I, \leq)$}\graffito{net $\family{m_i}_{i \in I}$ in $M$ indexed by $(I, \leq)$}\index[symbols]{miiinI@$\family{m_i}_{i \in I}$} if and only if the tuple $(I, \leq)$ is a directed set.
  \end{definition}

  \begin{definition} 
    Let $\net{m_i}_{i \in I}$ and $\net{m_j'}_{j \in J}$ be two nets in $M$. The net $\net{m_j'}_{j \in J}$ is called \define{subnet of $\net{m_i}_{i \in I}$}\graffito{subnet $\net{m_j'}_{j \in J}$ of $\net{m_i}_{i \in I}$} if and only if there is a map $f \from J \to I$ such that $\net{m_j'}_{j \in J} = \net{m_{f(j)}}_{j \in J}$ and
    \begin{equation*} 
      \ForEach i \in I \Exists j \in J \SuchThat \ForEach j' \in J \Holds \parens[\big]{j' \geq j \implies f(j') \geq i}. \qedhere
    \end{equation*}
  \end{definition}

  \begin{definition}
    Let $(X, \mathcal{T})$ be a topological space, let $\net{x_i}_{i \in I}$ be a net in $X$ indexed by $(I, \leq)$, and let $x$ be a point of $X$. The net $\net{x_i}_{i \in I}$ is said to \define{converge to $x$}\graffito{$\net{x_i}_{i \in I}$ converges to $x$} and $x$ is called \define{limit point of $\net{x_i}_{i \in I}$}\graffito{limit point $x$ of $\net{x_i}_{i \in I}$} if and only if
    \begin{equation*}
      \ForEach O \in \mathcal{T}_x \Exists i_0 \in I \SuchThat \ForEach i \in I \Holds (i \geq i_0 \implies x_i \in O). \qedhere
    \end{equation*}
  \end{definition}

  \begin{definition}
    Let $(X, \mathcal{T})$ be a topological space and let $\net{x_i}_{i \in I}$ be a net in $X$ indexed by $(I, \leq)$. The net $\net{x_i}_{i \in I}$ is called \define{convergent}\graffito{convergent net} if and only if there is a point $x \in X$ such that it converges to $x$.
  \end{definition}

  \begin{remark}
    Let $\net{m_i}_{i \in I}$ be a net that converges to $x$. Each subnet $\net{m_j'}_{j \in J}$ of $\net{m_i}_{i \in I}$ converges to $x$.
  \end{remark}

  \begin{definition}
    Let $(X, \mathcal{T})$ be a topological space and let $Y$ be a subset of $X$. The set
    \begin{equation*}
      \closureOf{Y} = \bigcap_{\substack{A \subseteq X \text{ closed}\\ Y \subseteq A}} A
      \mathnote{closure $\closureOf{Y}$ of $Y$ in $X$}
      \index[symbols]{Ybar@$\closureOf{Y}$}
    \end{equation*}
    is called \define{closure of $Y$ in $X$}.
  \end{definition}

  \begin{lemma} 
    Let $(X, \mathcal{T})$ be a topological space, let $Y$ be a subset of $X$, and let $x$ be an element of $X$. Then, $x \in \closureOf{Y}$ if and only if there is a net $\net{y_i}_{i \in I}$ in $Y$ that converges to $x$.
  \end{lemma}

  \begin{proof}
    See proposition~A.2.1 in \cite{ceccherini-silberstein:coornaert:2010}.
  \end{proof}

  \begin{definition}
     Let $(X, \mathcal{T})$ be a topological space. It is \define{Hausdorff}\graffito{Hausdorff} if and only if
     \begin{equation*}
       \ForEach x \in X \ForEach x' \in X \smallsetminus \setOf{x} \Exists O \in \mathcal{T}_x \Exists O' \in \mathcal{T}_{x'} \SuchThat O \cap O' = \emptyset. \qedhere
     \end{equation*}
  \end{definition}

  \begin{lemma} 
    Let $(X, \mathcal{T})$ be a topological space. It is Hausdorff if and only if each convergent net in $X$ has exactly one limit point. 
  \end{lemma}

  \begin{proof}
    See proposition~A.2.2 in \cite{ceccherini-silberstein:coornaert:2010}.
  \end{proof}

  \begin{definition}
    Let $(X, \mathcal{T})$ be a Hausdorff topological space, let $\net{x_i}_{i \in I}$ be a convergent net in $X$ indexed by $(I, \leq)$, and let $x$ be the limit point of $\net{x_i}_{i \in I}$. The point $x$ is denoted by $\lim_{i \in I} x_i$\graffito{the limit point $\lim_{i \in I} x_i$ of $\net{x_i}_{i \in I}$}\index[symbols]{limiinIxi@$\lim_{i \in I} x_i$} and we write $x_i \underset{i \in I}{\to} x$\graffito{$x_i \underset{i \in I}{\to} x$}\index[symbols]{arrow right limit@$x_i \underset{i \in I}{\to} x$}.
  \end{definition}

  \begin{definition}
    Let $(X, \mathcal{T})$ be a topological space, let $\net{x_i}_{i \in I}$ be a net in $X$ indexed by $(I, \leq)$, and let $x$ be an element of $X$. The point $x$ is called \define{cluster point of $\net{x_i}_{i \in I}$}\graffito{cluster point $x$ of $\net{x_i}_{i \in I}$} if and only if
    \begin{equation*}
      \ForEach O \in \mathcal{T}_x \ForEach i \in I \Exists i' \in I \SuchThat (i' \geq i \land x_{i'} \in O). \qedhere
    \end{equation*}
  \end{definition}

  \begin{lemma} 
    Let $(X, \mathcal{T})$ be a topological space, let $\net{x_i}_{i \in I}$ be a net in $X$ indexed by $(I, \leq)$, and let $x$ be an element of $X$. The point $x$ is a cluster point of $\net{x_i}_{i \in I}$ if and only if there is a subnet of $\net{x_i}_{i \in I}$ that converges to $x$.
  \end{lemma}

  \begin{proof}
    See proposition~A.2.3 in \cite{ceccherini-silberstein:coornaert:2010}.
  \end{proof}

  \begin{definition}
    Let $(X, \mathcal{T})$ and $(X', \mathcal{T}')$ be two topological spaces and let $f$ be a map from $X$ to $X'$. The map $f$ is called \define{continuous}\graffito{continuous map} if and only if
    \begin{equation*}
      \ForEach O' \in \mathcal{T}' \Holds f^{-1}(O') \in \mathcal{T}. \qedhere
    \end{equation*}
  \end{definition}

  \begin{lemma}
    Let $(X, \mathcal{T})$ and $(X', \mathcal{T}')$ be two topological spaces, let $f$ be a continuous map from $X$ to $X'$, let $\net{x_i}_{i \in I}$ be a net in $X$, and let $x$ be an element of $X$.
    \begin{aenumerate}
      \item If $x$ is a limit point of $\net{x_i}_{i \in I}$, then $f(x)$ is a limit point of $\net{f(x_i)}_{i \in I}$. 
      \item If $x$ is a cluster point of $\net{x_i}_{i \in I}$, then $f(x)$ is a cluster point of $\net{f(x_i)}_{i \in I}$. \qedhere
    \end{aenumerate}
  \end{lemma}

  \begin{proof}
    See the last paragraph of section~A.2 in \cite{ceccherini-silberstein:coornaert:2010}.
  \end{proof}

  \begin{lemma}
  \label{lemma:characterisation-of-continuity-by-continuity-at-each-point}
    Let $(X, \mathcal{T})$ and $(X', \mathcal{T}')$ be two topological spaces, and let $f$ be a map from $X$ to $X'$. The map $f$ is continuous if and only if, for each point $x \in X$ and each open neighbourhood $N'$ of $f(x)$, the preimage $f^{-1}(N')$ is an open neighbourhood of $x$.
  \end{lemma}

  \begin{proof}
    See proposition~2.4.2 of section~2.4 in \cite{waldmann:2014}.
  \end{proof}

  \begin{definition} 
    Let $\overline{\R} = \R \cup \setOf{-\infty,+\infty}$ be the affinely extended real numbers and let $\net{r_i}_{i \in I}$ be a net in $\overline{\R}$ indexed by $(I, \leq)$.
    \begin{aenumerate} 
      \item The limit of the net $\net{\inf_{i' \geq i} r_{i'}}_{i \in I}$ is called \graffito{limit inferior $\liminf_{i \in I} r_i$ of $\net{r_i}_{i \in I}$}\define{limit inferior of $\net{r_i}_{i \in I}$} and denoted by $\liminf_{i \in I} r_i$\index[symbols]{liminfiinIri@$\liminf_{i \in I} r_i$}.
      \item The limit of the net $\net{\sup_{i' \geq i} r_{i'}}_{i \in I}$ is called \graffito{limit superior $\limsup_{i \in I} r_i$ of $\net{r_i}_{i \in I}$}\define{limit superior of $\net{r_i}_{i \in I}$} and denoted by $\limsup_{i \in I} r_i$\index[symbols]{limsupiinIri@$\limsup_{i \in I} r_i$}. \qedhere
    \end{aenumerate}
  \end{definition}

  \section{Initial And Product Topologies}

  \begin{definition}
    Let $X$ be a set, let $I$ be a set, and, for each index $i \in I$, let $(Y_i, \mathcal{T}_i)$ be a topological space and let $f_i$ be a map from $X$ to $Y_i$. The coarsest topology on $X$ such that, for each index $i \in I$, the map $f_i$ is continuous, is called \defineX{initial with respect to $\family{f_i}_{i \in I}$}{initial!topology}\graffito{initial topology with respect to $\family{f_i}_{i \in I}$}.
  \end{definition}

  \begin{remark}
    The initial topology is explicitly constructed in section~A.3 in \cite{ceccherini-silberstein:coornaert:2010}.
  \end{remark}

  \begin{lemma}
  \label{lemma:map-to-initial-topology-continuous-if-and-only-if-gens-after-map-continuous}
    Let $(X, \mathcal{T})$ be a topological space, where $\mathcal{T}$ is the initial topology with respect to $\family{f_i \from X \to Y_i}_{i \in I}$, let $(Z, \mathcal{S})$ be a topological space, and let $g$ be a map from $Z$ to $X$. The map $g$ is continuous if and only if, for each index $i \in I$, the map $f_i \after g$ is continuous.
  \end{lemma}

  \begin{proof}
    See the last paragraph of section~A.3 in \cite{ceccherini-silberstein:coornaert:2010}. 
  \end{proof}

  \begin{lemma}
  \label{lemma:limit-and-cluster-points-in-the-initial-topology}
    Let $(X, \mathcal{T})$ be a topological space, where $\mathcal{T}$ is the initial topology with respect to $\family{f_i \from X \to Y_i}_{i \in I}$, let $\net{x_{i'}}_{i' \in I'}$ be a net in $X$, and let $x$ be a point in $X$. The point $x$ is a limit point or cluster point of $\net{x_{i'}}_{i' \in I'}$ if and only if, for each index $i \in I$, the point $f_i(x)$ is a limit point or cluster point of $\net{f_i(x_{i'})}_{i' \in I'}$.
  \end{lemma}

  \begin{proof}
    See the last paragraph of section~A.3 in \cite{ceccherini-silberstein:coornaert:2010}.
  \end{proof}

  \begin{definition}
    Let $\family{(X_i, \mathcal{T}_i)}_{i \in I}$ be a family of topological spaces, let $X$ be the set $\prod_{i \in I} X_i$, and, for each index $i \in I$, let $\pi_i$ be the canonical projection of $X$ onto $X_i$. The initial topology on $X$ with respect to $\family{\pi_i}_{i \in I}$ is called \defineX{product}{product!topology}\graffito{product topology} and is also known as \graffito{topology of pointwise convergence}\defineX{topology of pointwise convergence}{topology!of pointwise convergence}.
  \end{definition}

  \begin{remark}
    The product topology on $X$ has for a base the sets $\prod_{i \in I} O_i$, where, for each index $i \in I$, the set $O_i$ is an open subset of $X_i$, and the set $\setOf{i \in I \suchThat O_i \neq X_i}$ is finite.
  \end{remark}

  \begin{definition}
    Let $\family{(X_i, \mathcal{T}_i)}_{i \in I}$ be a family of discrete topological spaces and let $X$ be the set $\prod_{i \in I} X_i$. The product topology on $X$ is called \defineX{prodiscrete}{prodiscrete!topology}\graffito{prodiscrete topology}.
  \end{definition}

  \begin{lemma} 
  \label{lemma:product-of-Hausdorff-is-Hausdorff}
    Let $\family{(X_i, \mathcal{T}_i)}_{i \in I}$ be a family of Hausdorff topological spaces. The set $\prod_{i \in I} X_i$, equipped with the product topology, is Hausdorff.
  \end{lemma}

  \begin{proof}
    See proposition~A.4.1 in \cite{ceccherini-silberstein:coornaert:2010}.
  \end{proof}

  \begin{definition}
    Let $(X, \mathcal{T})$ be a topological space and let $D$ be a subset of $X$. The set $D$ is called \define{connected}\graffito{connected subset} if and only if
    \begin{equation*}
      \ForEach O \in \mathcal{T} \ForEach O' \in \mathcal{T} \Holds
      \begin{aligned}[t]
        \big(&O \cap D \neq \emptyset \land O' \cap D \neq \emptyset \land O \cap O' \cap D = \emptyset\\
             &{}\implies O \cup O' \nsupseteq D\big). \qedhere
      \end{aligned}
    \end{equation*}
  \end{definition}

  \begin{definition}
    Let $(X, \mathcal{T})$ be a topological space. It is called \define{totally disconnected}\graffito{totally disconnected topological space} if and only if, for each non-empty and connected subset $D$ of $X$, we have $\cardinalityOf{D} = 1$. 
  \end{definition}

  \begin{lemma} 
  \label{lemma:product-of-totally-disconnected-is-disconnected}
    Let $\family{(X_i, \mathcal{T}_i)}_{i \in I}$ be a family of totally disconnected topological spaces. The set $\prod_{i \in I} X_i$, equipped with the product topology, is totally disconnected.
  \end{lemma}

  \begin{proof}
    See proposition~A.4.2 in \cite{ceccherini-silberstein:coornaert:2010}.
  \end{proof}

  \begin{lemma}
    Let $\family{(X_i, \mathcal{T}_i)}_{i \in I}$ be a family of topological spaces and, for each index $i \in I$, let $A_i$ be a closed subset of $X_i$. The set $\prod_{i \in I} A_i$ is a closed subset of $\prod_{i \in I} X_i$, equipped with the product topology.
  \end{lemma}

  \begin{proof}
    See proposition~A.4.3 in \cite{ceccherini-silberstein:coornaert:2010}.
  \end{proof}

  \section{Compactness}

  \begin{definition} 
    Let $(X, \mathcal{T})$ be a topological space and let $\family{O_i}_{i \in I}$ be a family of elements of $\mathcal{T}$. The family $\family{O_i}_{i \in I}$ is called \graffito{open cover $\family{O_i}_{i \in I}$ of $X$}\define{open cover of $X$} if and only if $\bigcup_{i \in I} O_i = X$.
  \end{definition}

  \begin{definition}
    Let $(X, \mathcal{T})$ be a topological space. It is called \define{compact}\graffito{compact topological space} if and only if, for each open cover $\family{O_i}_{i \in I}$ of $X$, there is a finite subset $J$ of $I$ such that $\family{O_j}_{j \in J}$ is an open cover of $X$.
  \end{definition}

  \begin{definition}
    Let $M$ be a set and let $\family{D_i}_{i \in I}$ be a family of subsets of $M$. The family $\family{D_i}_{i \in I}$ is said to have the \graffito{finite intersection property of $\family{D_i}_{i \in I}$}\define{finite intersection property} if and only if, for each finite subset $J$ of $I$, we have $\bigcap_{j \in J} D_j \neq \emptyset$.
  \end{definition}

  \begin{lemma}
    Let $(X, \mathcal{T})$ be a topological space. It is compact if and only if, for each family $\family{A_i}_{i \in I}$ of closed subsets of $X$ that has the finite intersection property, we have $\bigcap_{i \in I} A_i \neq \emptyset$.
  \end{lemma}

  \begin{proof}
    See the first paragraph of section~A.5 in \cite{ceccherini-silberstein:coornaert:2010}.
  \end{proof}

  \begin{theorem} 
    Let $(X, \mathcal{T})$ be a topological space. The following three statements are equivalent:
    \begin{aenumerate}
      \item The space $(X, \mathcal{T})$ is compact.
      \item Each net in $X$ has a cluster point with respect to $\mathcal{T}$.
      \item Each net in $X$ has a convergent subnet with respect to $\mathcal{T}$. \qedhere
    \end{aenumerate}
  \end{theorem}

  \begin{proof}
    See theorem~A.5.1 in \cite{ceccherini-silberstein:coornaert:2010}.
  \end{proof}

  \begin{theorem}[Andrey Nikolayevich Tikhonov, 1935] 
  \label{theorem:Tychonoff}
    Let $\family{(X_i, \mathcal{T}_i)}_{i \in I}$ be a family of compact topological spaces. The set $\prod_{i \in I} X_i$, equipped with the product topology, is compact.
  \end{theorem}

  \begin{proof}
    See theorem~A.5.2 in \cite{ceccherini-silberstein:coornaert:2010}.
  \end{proof}

  \begin{corollary} 
  \label{corollary:tychonoff}
    Let $\family{(X_i, \mathcal{T}_i)}_{i \in I}$ be a family of finite topological spaces. The set $\prod_{i \in I} X_i$, equipped with the product topology, is compact.
  \end{corollary}

  \begin{proof}
    See the paragraph before corollary~A.5.3 in \cite{ceccherini-silberstein:coornaert:2010}.
  \end{proof}

  \clearToOddPage
  \chapter{Uniformities}
  \label{chapter:uniformities}

  Most of the theory of uniformities as presented here may be found in more detail in Appendix~B in the monograph \enquote{\citetitle*{ceccherini-silberstein:coornaert:2010}}\cite{ceccherini-silberstein:coornaert:2010}.

  \begin{definition}
    Let $X$ be a set. The set $\Delta_X = \setOf{(x,x) \suchThat x \in X}$ is called \define{diagonal in $X \times X$}\graffito{diagonal $\Delta_X$ in $X \times X$}\index[symbols]{DeltaX@$\Delta_X$}.
  \end{definition}

  \begin{definition}
    Let $X$ be a set, let $R$ be a subset of $X \times X$, and let $y$ be an element of $X$. The set $R[y] = \setOf{x \in X \suchThat (x,y) \in R}$ is called \defineX{$y$-th column in $R$}{column in $R$@$y$-th column in $R$}\graffito{$y$-th column $R[y]$ in $R$}\index[symbols]{Rybrackets@$R[y]$}.
  \end{definition}

  \begin{definition}
    Let $X$ be a set and let $R$ be a subset of $X \times X$. The set $R^{-1} = \setOf{(y, x) \suchThat (x, y) \in R}$ is called \define{inverse of $R$}\graffito{inverse $R^{-1}$ of $R$}\index[symbols]{Rminus1@$R^{-1}$}. The set $R$ is called \defineX{symmetric}{symmetric!binary relation}\graffito{symmetric subset $R$ of $X \times X$} if and only if $R^{-1} = R$. 
  \end{definition}

  \begin{definition}
    Let $X$ be a set, and let $R$ and $R'$ be two subsets of $X \times X$. The set $R' \after R = \setOf{(x, z) \suchThat \Exists y \in X \SuchThat (x, y) \in R \wedge (y, z) \in R'}$ is called \define{composition of $R$ and $R'$}\graffito{composition $R' \after R$ of $R$ and $R'$}\index[symbols]{RcircleRprime@$R' \after R$}.
  \end{definition}

  \begin{definition}
    Let $X$ be a set and let $\mathcal{U}$ be a set of subsets of $X \times X$. The set $\mathcal{U}$ is called \define{uniformity on $X$}\graffito{uniformity $\mathcal{U}$ on $X$}\index[symbols]{Ucalligraphic@$\mathcal{U}$} if and only if 
    \begin{gather*}
      \mathcal{U} \neq \emptyset,\\
      \ForEach E \in \mathcal{U} \Holds \Delta_X \subseteq E,\\
      \ForEach E \in \mathcal{U} \ForEach E' \subseteq X \times X \Holds (E \subseteq E' \implies E' \in \mathcal{U}),\\
      \ForEach E \in \mathcal{U} \ForEach E' \in \mathcal{U} \Holds E \cap E' \in \mathcal{U},\\
      \ForEach E \in \mathcal{U} \Holds E^{-1} \in \mathcal{U},\\
      \ForEach E \in \mathcal{U} \Exists E' \in \mathcal{U} \SuchThat E' \after E' \subseteq E. \qedhere
    \end{gather*}
  \end{definition}

  \begin{definition} 
    Let $X$ be a set and let $\mathcal{U}$ be a uniformity on $X$. The tuple $(X, \mathcal{U})$ is called \define{uniform space}\graffito{uniform space $(X, \mathcal{U})$} and each element $E \in \mathcal{U}$ is called \define{entourage of $X$}\graffito{entourage $E$ of $X$}.

    The set $X$ is said to be \define{equipped with $\mathcal{U}$}\graffito{equipped with $\mathcal{U}$} if and only if it shall be implicitly clear that $\mathcal{U}$ is the uniformity on $X$ being considered. The set $X$ is called \define{uniform space}\graffito{uniform space $X$} if and only if it is implicitly clear what uniformity on $X$ is being considered. 
  \end{definition}

  \begin{example}
    Let $X$ be a set and let $\mathcal{U}$ be the set of all subsets of $X \times X$ that contain $\Delta_X$. Then, $\mathcal{U}$ is the finest uniformity on $X$. Itself as well as the uniform space $(X, \mathcal{U})$ are called \defineX{discrete}{discrete!uniformity}\graffito{discrete uniformity}.
  \end{example}

  \begin{definition}
    Let $(X, \mathcal{U})$ be a uniform space. The set
    \begin{equation*}
      \setOf{O \subseteq X \suchThat \ForEach x \in O \Exists E \in \mathcal{U} \SuchThat E[x] = O}
    \end{equation*}
    is a topology on $X$ and called \define{induced by $\mathcal{U}$}\graffito{topology on $X$ induced by $\mathcal{U}$}.
  \end{definition}

  \begin{example}
    Let $(X, \mathcal{U})$ be a discrete uniform space. The topology induced by $\mathcal{U}$ is the discrete topology on $X$. However, there are non-discrete uniformities whose induced topologies are discrete. 
  \end{example}

  \begin{remark}
  \label{remark:characterisation-of-continuity-on-induced-topology}
    Let $(X, \mathcal{U})$ and $(X', \mathcal{U}')$ be two uniform spaces, let $\mathcal{T}$ and $\mathcal{T}'$ be the respective topologies on $X$ and $X'$ induced by $\mathcal{U}$ and $\mathcal{U}'$, and let $f$ be a map from $X$ to $X'$. The map $f$ is continuous if and only if, for each point $x \in X$ and each entourage $E' \in \mathcal{U}'$, there is an entourage $E \in \mathcal{U}$ such that $f(E[x]) \subseteq E'[f(x)]$.
  \end{remark}

  \begin{definition}
    Let $(X, \mathcal{U})$ be a uniform space and let $\mathcal{B}$ be a subset of $\mathcal{U}$.
    \begin{aenumerate}
      \item The set $\mathcal{B}$ is called \define{base of $\mathcal{U}$}\graffito{base $\mathcal{B}$ of $\mathcal{U}$}\index[symbols]{Bcalligraphic@$\mathcal{B}$} if and only if
            \begin{equation*}
              \ForEach E \in \mathcal{U} \Exists B \in \mathcal{B} \SuchThat B \subseteq E.
            \end{equation*}
      \item The set $\mathcal{B}$ is called \define{subbase of $\mathcal{U}$}\graffito{subbase $\mathcal{B}$ of $\mathcal{U}$}\index{base of $\mathcal{U}$!sub-}\index[symbols]{Bcalligraphic@$\mathcal{B}$} if and only if the set
            \begin{equation*}
              \setOf{\bigcap_{i = 1}^n B_i \suchThat B_i \in \mathcal{B}, i \in \setOf{1,2,\dotsc,n}, n \in \N_+}
            \end{equation*}
            is a base of $\mathcal{U}$. \qedhere
    \end{aenumerate}
  \end{definition}

  \begin{lemma}
    Let $X$ be a set and let $\mathcal{B}$ be a set of subsets of $X \times X$.
    \begin{aenumerate}
      \item The set $\mathcal{B}$ is a base of a uniformity on $X$ if and only if
            \begin{gather*}
              \mathcal{B} \neq \emptyset,\\
              \ForEach B \in \mathcal{B} \Holds \Delta_X \subseteq B,\\
              \ForEach B \in \mathcal{B} \ForEach B' \in \mathcal{B} \Exists B'' \in \mathcal{B} \SuchThat B'' \subseteq B \cap B',\\
              \ForEach B \in \mathcal{B} \Exists B' \in \mathcal{B} \SuchThat B' \subseteq B^{-1},\\
              \ForEach B \in \mathcal{B} \Exists B' \in \mathcal{B} \SuchThat B' \after B' \subseteq B.
            \end{gather*}
      \item The set $\mathcal{B}$ is a subbase of a uniformity on $X$ if and only if
            \begin{gather*}
              \mathcal{B} \neq \emptyset,\\
              \ForEach B \in \mathcal{B} \Holds \Delta_X \subseteq B,\\
              \ForEach B \in \mathcal{B} \Exists B' \in \mathcal{B} \SuchThat B' \subseteq B^{-1},\\
              \ForEach B \in \mathcal{B} \Exists B' \in \mathcal{B} \SuchThat B' \after B' \subseteq B. \qedhere
            \end{gather*}
    \end{aenumerate}
  \end{lemma}

  \begin{definition}
    Let $X$ be a set and let $\mathcal{B}$ be a base or subbase of a uniformity on $X$. The uniformity $\mathcal{U}$ that has $\mathcal{B}$ for a base or subbase respectively is uniquely determined and called \define{generated by $\mathcal{B}$}\graffito{uniformity generated by $\mathcal{B}$}.
  \end{definition}

  \begin{definition} 
    Let $I$ be a set and, for each index $i \in I$, let $f_i$ be a map from $X_i$ to $X_i'$. The map\index[symbols]{productfiiinI@$\prod_{i \in I} f_i$}
    \begin{align*}
      \prod_{i \in I} f_i \from \prod_{i \in I} X_i &\to \prod_{i \in I} X_i', \mathnote{product $\prod_{i \in I} f_i$ of $\family{f_i}_{i \in I}$}\\
      \family{x_i}_{i \in I} &\mapsto \family{f_i(x_i)}_{i \in I},
    \end{align*}
    is called \define{product of $\family{f_i}_{i \in I}$} and, if $I$ is the set $\setOf{1, 2, \dotsc, \cardinalityOf{I}}$, then it is also denoted by $f_1 \times f_2 \times \dotsb \times f_{\cardinalityOf{I}}$\graffito{$f_1 \times f_2 \times \dotsb \times f_{\cardinalityOf{I}}$}\index[symbols]{crossf1f2@$f_1 \times f_2 \times \dotsb \times f_{\cardinalityOf{I}}$}.
  \end{definition}

  \begin{definition}
    Let $X$ be a set, let $(X', \mathcal{U}')$ be a uniform space, and let $f$ be a map from $X$ to $X'$. The set
    \begin{equation*}
      \mathcal{B} = \setOf{(f \times f)^{-1}(E') \suchThat E' \in \mathcal{U}'}
    \end{equation*}
    is a base of a uniformity on $X$ and the uniformity on $X$ it generates is called \defineX{induced on $X$ by $f$}{induced on $X$ by $f$!uniformity}\graffito{uniformity on $X$ induced by $f$}. 
  \end{definition}

  \begin{definition}
    Let $(X, \mathcal{U})$ be a uniform space, let $Y$ be a subset of $X$, and let $\iota$ be the canonical injection from $Y$ to $X$. The uniformity $\mathcal{V}$ on $Y$ induced by $\iota$ is called \defineX{subspace}{subspace!uniformity}\graffito{subspace uniformity} and the uniform space $(Y, \mathcal{V})$ is called \define{subspace of $(X, \mathcal{U})$}\graffito{subspace of $(X, \mathcal{U})$}.
  \end{definition}

  \begin{remark}
    Let $(X, \mathcal{U})$ be a uniform space and let $Y$ be a subset of $X$. The subspace uniformity on $Y$ is $\setOf{E \cap (Y \times Y) \suchThat E \in \mathcal{U}}$.
  \end{remark}

  \begin{definition}
    Let $(X, \mathcal{U})$ and $(X', \mathcal{U}')$ be two uniform spaces and let $f$ be a map from $X$ to $X'$. The map $f$ is called \graffito{uniformly continuous map}\define{uniformly continuous} if and only if
    \begin{equation*}
      \ForEach E' \in \mathcal{U}' \Exists E \in \mathcal{U} \SuchThat (f \times f)(E) \subseteq E'. \qedhere
    \end{equation*}
  \end{definition}

  \begin{remark}
    The map $f$ is uniformly continuous if and only if
    \begin{equation*}
      \ForEach E' \in \mathcal{U}' \Holds (f \times f)^{-1}(E') \in \mathcal{U}. \qedhere
    \end{equation*}
  \end{remark}

  \begin{lemma}
  \label{lemma:characterisation-of-uniformly-continuous-by-base}
    Let $(X, \mathcal{U})$ and $(X', \mathcal{U}')$ be two uniform spaces, let $f$ be a map from $X$ to $X'$, and let $\mathcal{B}'$ be a base or subbase of $\mathcal{U}'$. The map $f$ is uniformly continuous if and only if
    \begin{equation*}
      \ForEach B' \in \mathcal{B}' \Holds (f \times f)^{-1}(B') \in \mathcal{U}. \qedhere
    \end{equation*}
  \end{lemma}

  \begin{proof}
    Compare the second paragraph of section~B.2 in \cite{ceccherini-silberstein:coornaert:2010}.
  \end{proof}

  \begin{lemma}
    Let $(X, \mathcal{U})$ and $(X', \mathcal{U}')$ be two uniform spaces, let $f$ be a map from $X$ to $X'$, and let $\mathcal{B}$ and $\mathcal{B}'$ be two bases or subbases of $\mathcal{U}$ and $\mathcal{U}'$ respectively. The map $f$ is uniformly continuous if and only if
    \begin{equation*}
      \ForEach B' \in \mathcal{B}' \Exists B \in \mathcal{B} \SuchThat (f \times f)(B) \subseteq B'. \qedhere
    \end{equation*}
  \end{lemma}

  \begin{proof}
    Compare the second paragraph of section~B.2 in \cite{ceccherini-silberstein:coornaert:2010}.
  \end{proof}

  \begin{lemma}
  \label{lemma:uniformly-continuous-map-is-continuous}
    Let $(X, \mathcal{U})$ and $(X', \mathcal{U}')$ be two uniform spaces and let $f$ be a uniformly continuous map from $X$ to $X'$. The map $f$ is continuous, where $X$ and $X'$ are equipped with the topologies induced by $\mathcal{U}$ and $\mathcal{U}'$ respectively.
  \end{lemma}

  \begin{proof}
    See proposition~B.2.2 in \cite{ceccherini-silberstein:coornaert:2010}.
  \end{proof}

  \begin{theorem}
  \label{theorem:continuous-map-on-compact-space-is-uniformly-continuous}
    Let $(X, \mathcal{U})$ and $(X', \mathcal{U}')$ be two uniform spaces such that $X$, equipped with the topology induced by $\mathcal{U}$, is compact. Each continuous map $f \from X \to X'$ is uniformly continuous.
  \end{theorem}

  \begin{proof}
    See theorem~B.2.3 in \cite{ceccherini-silberstein:coornaert:2010}.
  \end{proof}

  \begin{definition}
    Let $(X, \mathcal{U})$ and $(X', \mathcal{U}')$ be two uniform spaces, and let $f$ be a map from $X$ to $X'$. The map $f$ is called \graffito{uniform isomorphism}\define{uniform isomorphism} if and only if it is bijective, and both $f$ and $f^{-1}$ are uniformly continuous.
  \end{definition}

  \begin{definition}
    Let $(X, \mathcal{U})$ and $(X', \mathcal{U}')$ be two uniform spaces, and let $f$ be a map from $X$ to $X'$. The map $f$ is called \graffito{uniform embedding}\define{uniform embedding} if and only if it is injective and $f\restrictedTo_{X \to f(X)}$ is a uniform isomorphism.
  \end{definition}

  \begin{lemma} 
  \label{lemma:continuous-injective-map-from-compact-to-Hausdorff-is-uniform-embedding}
    Let $(X, \mathcal{U})$ and $(X', \mathcal{U}')$ be two uniform spaces such that $X$ is compact and $X'$ is Hausdorff. Furthermore, let $f$ be a continuous and injective map from $X$ to $X'$. The map $f$ is a uniform embedding.
  \end{lemma}

  \begin{proof}
    See proposition~B.2.5 in \cite{ceccherini-silberstein:coornaert:2010}.
  \end{proof}

  \begin{definition} 
    Let $X$ be a set, let $I$ be a set, and, for each index $i \in I$, let $(X_i, \mathcal{U}_i)$ be a uniform space and let $f_i$ be a map from $X$ to $X_i$. The coarsest uniformity on $X$ such that, for each index $i \in I$, the map $f_i$ is uniformly continuous, is called \defineX{initial with respect to $\family{f_i}_{i \in I}$}{initial!uniformity}\graffito{initial uniformity with respect to $\family{f_i}_{i \in I}$}.
  \end{definition}

  \begin{definition} 
    Let $\family{(X_i, \mathcal{U}_i)}_{i \in I}$ be a family of uniform spaces, let $X$ be the set $\prod_{i \in I} X_i$, and, for each index $i \in I$, let $\pi_i$ be the canonical projection from $X$ onto $X_i$. The initial uniformity on $X$ with respect to $\family{\pi_i}_{i \in I}$ is called \define{product}\graffito{product uniformity}.
  \end{definition}

  \begin{remark}
    The product uniformity on $X$ has for a base the sets $\prod_{i \in I} E_i$, where, for each index $i \in I$, the set $E_i$ is an entourage of $X_i$, and the set $\setOf{i \in I \suchThat E_i \neq X_i \times X_i}$ is finite. Note that
    \begin{equation*}
      \prod_{i \in I} E_i
      \subseteq \prod_{i \in I} X_i \times X_i
      = (\prod_{i \in I} X_i) \times (\prod_{i \in I} X_i)
      = X \times X. \qedhere
    \end{equation*}
  \end{remark}

  \begin{definition}
  \label{definition:prodiscrete-uniformity}
    Let $\family{(X_i, \mathcal{U}_i)}_{i \in I}$ be a family of discrete uniform spaces and let $X$ be the set $\prod_{i \in I} X_i$. The product uniformity on $X$ is called \defineX{prodiscrete}{prodiscrete!uniformity}\graffito{prodiscrete uniformity}.
  \end{definition}

  \begin{remark}
  \label{remark:base-of-prodiscrete-uniformity}
    Let $M$ and $Q$ be two sets. The prodiscrete uniformity on $Q^M = \prod_{m \in M} Q$ has for a base the sets
    \begin{equation*}
      \setOf{(c,c') \in Q^M \times Q^M \suchThat c\restrictedTo_F = c'\restrictedTo_F},
      \text{ for } F \subseteq M \text{ finite}. \qedhere
    \end{equation*}
  \end{remark}

  \clearToOddPage
  \chapter{Dual Spaces}
  \label{chapter:dual-spaces}

  Most of the theory of dual spaces as presented here may be found in more detail in Appendix~F in the monograph \enquote{\citetitle*{ceccherini-silberstein:coornaert:2010}}\cite{ceccherini-silberstein:coornaert:2010}.

  \begin{definition}
    Let $X$ be an $\R$-vector space and let $X$ be equipped with a topology. The $\R$-vector space $X$ is called \defineX{topological}{topological!vector space}\graffito{topological $\R$-vector space $X$} if and only if the maps $+ \from X \times X \to X$, $(x, x') \mapsto x + x'$, and $\cdot \from \R \times X \to X$, $(r, x) \mapsto r \cdot x$, are continuous, where $X \times X$ and $\R \times X$ are equipped with their respective product topology.
  \end{definition}

  \begin{definition}
    Let $X$ be an $\R$-vector space and let $Y$ be a subset of $X$. The set $Y$ is called
            \define{convex}\graffito{convex subset of $X$} if and only if
            \begin{equation*}
              \ForEach (y, y') \in Y \times Y \ForEach t \in \closedInterval{0, 1} \Holds t y + (1 - t) y' \in Y. \qedhere 
            \end{equation*}
  \end{definition}

  \begin{definition} 
    Let $X$ be a topological $\R$-vector space. It is called \define{locally convex}\graffito{locally convex topological $\R$-vector space}\index{convex!locally} if and only if the origin has a neighbourhood base of convex sets. 
  \end{definition}

  In the remainder of this chapter, let $(X, \normOf{\blank})$ be a normed $\R$-vector space. 

  \begin{definition} 
    The vector space
    \begin{equation*}
      X^* = \setOf{\psi \from X \to \R \suchThat \psi \text{ is linear and continuous}}
      \mathnote{topological dual space $X^*$ of $X$} 
      \index[symbols]{Xstar@$X^*$}
    \end{equation*}
    with pointwise addition and scalar multiplication is called \define{topological dual space of $X$} and each map $\psi \in X^*$ is called \graffito{continuous linear functional $\psi$}\define{continuous linear functional}.
  \end{definition}

  \begin{definition}
    The norm
    \begin{align*}
      \normOf{\blank}_{X^*} \from X^* &\to \R, \mathnote{operator norm $\normOf{\blank}_{X^*}$ on $X^*$}\index[symbols]{normXstar@$\normOf{\blank}_{X^*}$}\\
      \psi &\mapsto \sup_{x \in X \smallsetminus \setOf{0}} \frac{\absoluteValueOf{\psi(x)}}{\normOf{x}},
    \end{align*}
    is called \define{operator norm on $X^*$}. And, the topology on $X^*$ induced by $\normOf{\blank}_{X^*}$ is called \define{strong topology on $X^*$}\graffito{strong topology on $X^*$}.
  \end{definition}

  \begin{definition}
    Let $x$ be an element of $X$. The map
    \begin{align*}
      \evaluationMap_x \from X^* &\to \R, \mathnote{evaluation map $\evaluationMap_x$ at $x$}\index[symbols]{evx@$\evaluationMap_x$}\\
      \psi &\mapsto \psi(x),
    \end{align*}
    is called \define{evaluation map at $x$}.
  \end{definition}

  \begin{definition} 
    The initial topology on $X^*$ with respect to $\family{\evaluationMap_x}_{x \in X}$ is called \define{weak-$*$ topology on $X^*$}\graffito{weak-$*$ topology on $X^*$}.
  \end{definition}

  \begin{lemma}
    Let $\net{\psi_i}_{i \in I}$ be a net in $X^*$, let $\psi$ be an element of $X^*$, and let $X^*$ be equipped with the weak-$*$ topology. The net $\net{\psi_i}_{i \in I}$ converges to $\psi$ if and only if, for each element $x \in X$, the net $\net{\psi_i(x)}_{i \in I}$ converges to $\psi(x)$.
  \end{lemma}

  \begin{proof}
    This is a direct consequence of \cref{lemma:limit-and-cluster-points-in-the-initial-topology}. 
  \end{proof}

  \begin{lemma}
  \label{lemma:weak-star-coarser-than-strong-topology}
    The weak-$*$ topology on $X^*$ is coarser than the strong topology on $X^*$.
  \end{lemma}

  \begin{proof}
    This holds because the evaluation maps are continuous with respect to the strong topology on $X^*$.
  \end{proof}

  \begin{corollary}
    Let $\net{\psi_i}_{i \in I}$ be a net in $X^*$ that converges to $\psi$ with respect to the strong topology on $X^*$. The net $\net{\psi_i}_{i \in I}$ converges to $\psi$ with respect to the weak-$*$ topology on $X^*$.
  \end{corollary}

  \begin{proof}
    This is a direct consequence of \cref{lemma:weak-star-coarser-than-strong-topology}.
  \end{proof}

  \begin{lemma}
    Let $X^*$ be equipped with the weak-$*$ topology and let $\psi$ be an element of $X^*$. An open neighbourhood base of $\psi$ is given by the sets
      \begin{multline*}
        \operatorname{B}(\psi, F, \varepsilon) = \setOf{\psi' \in X^* \suchThat \ForEach x \in F \Holds \absoluteValueOf{\psi(x) - \psi'(x)} < \varepsilon},\\
        \text{ for } F \subseteq X \text{ finite and } \varepsilon \in \R_{> 0}.
        \mathnote{$\operatorname{B}(\psi, F, \varepsilon)$, for $\psi \in X^*$, $F \subseteq X$ finite, and $\varepsilon \in \R_{> 0}$}\index[symbols]{BphiFepsilon@$\operatorname{B}(\psi, F, \varepsilon)$} \qedhere
      \end{multline*}
  \end{lemma}

  \begin{proof}
    Compare the third paragraph of section~F.2 in \cite{ceccherini-silberstein:coornaert:2010}. 
  \end{proof}

  \begin{corollary}
    Let $X^*$ be equipped with the weak-$*$ topology. The space $X^*$ is locally convex.
  \end{corollary}

  \begin{proof}
    See the third paragraph of section~F.2 in \cite{ceccherini-silberstein:coornaert:2010}.
  \end{proof}

  \begin{lemma}
  \label{lemma:weak-star-topology-is-Hausdorff}
    The space $X^*$, equipped with the weak-$*$ topology, is Hausdorff.
  \end{lemma}

  \begin{proof}
    See the last paragraph of section~F.2 in \cite{ceccherini-silberstein:coornaert:2010}.
  \end{proof}

  \begin{theorem}[Stefan Banach, 1932; Leonidas Alaoglu, 1940]
  \label{theorem:banach-alaoglu}
    Let $X^*$ be equipped with the weak-$*$ topology. The unit ball $\setOf{\psi \in X^* \suchThat \normOf{\psi}_{X^*} \leq 1}$, equipped with the subspace topology, is compact.
  \end{theorem}

  \begin{proof}
    See theorem~F.3.1 in \cite{ceccherini-silberstein:coornaert:2010}.
  \end{proof}

  \clearToOddPage
  \chapter{Hall's Theorems} 
  \label{chapter:halls-theorems}

  Most of the theory concerning Hall's theorems as presented here may be found in more detail in Appendix~H in the monograph \enquote{\citetitle*{ceccherini-silberstein:coornaert:2010}}\cite{ceccherini-silberstein:coornaert:2010}.

  \section{Bipartite Graphs}

  \begin{definition}
    Let $X$ and $Y$ be two sets, and let $E$ be a subset of $X \times Y$. The triple $(X, Y, E)$ is called \define{bipartite graph}\graffito{bipartite graph $(X, Y, E)$}, each element $x$ of $X$ is called \define{left vertex}\graffito{left vertex $x$}\index{vertex!left}, each element $y$ of $Y$ is called \define{right vertex}\graffito{right vertex $y$}\index{vertex!right}, and each element $e$ of $E$ is called \define{edge}\graffito{edge $e$}.
  \end{definition}

  \begin{definition}
    Let $(X, Y, E)$ and $(X', Y', E')$ be two bipartite graphs. The graph $(X, Y, E)$ is called \define{bipartite subgraph of $(X', Y', E')$}\graffito{bipartite subgraph $(X, Y, E)$ of $(X', Y', E')$} if and only if $X \subseteq X'$, $Y \subseteq Y'$, and $E \subseteq E'$.
  \end{definition}

  \begin{definition}
    Let $(X, Y, E)$ be a bipartite graph, and let $(x, y)$ and $(x', y')$ be two elements of $E$. The edges $(x, y)$ and $(x', y')$ are called \define{adjacent}\graffito{adjacent edges} if and only if $x = x'$ or $y = y'$.
  \end{definition}

  \begin{definition}
    Let $(X, Y, E)$ be a bipartite graph.
    \begin{aenumerate}
      \item Let $x$ be an element of $X$. The set
            \begin{equation*}
              \mathcal{N}_r(x) = \setOf{y \in Y \suchThat (x, y) \in E}
              \mathnote{right neighbourhood $\mathcal{N}_r(x)$ of $x$}
              \index[symbols]{Nrxcalligraphic@$\mathcal{N}_r(x)$}
            \end{equation*}
            is called \define{right neighbourhood of $x$}\index{neighbourhood!right}.
      \item Let $A$ be a subset of $X$. The set $\mathcal{N}_r(A) = \bigcup_{a \in A} \mathcal{N}_r(a)$ is called \define{right neighbourhood of $A$}\graffito{right neighbourhood $\mathcal{N}_r(A)$ of $A$}\index[symbols]{NrAcalligraphic@$\mathcal{N}_r(A)$}. 
      \item Let $y$ be an element of $Y$. The set
            \begin{equation*}
              \mathcal{N}_l(y) = \setOf{x \in X \suchThat (x, y) \in E}
              \mathnote{left neighbourhood $\mathcal{N}_l(y)$ of $y$}
              \index[symbols]{Nlycalligraphic@$\mathcal{N}_l(y)$}
            \end{equation*}
            is called \define{left neighbourhood of $y$}\index{neighbourhood!left}.
      \item Let $B$ be a subset of $Y$. The set $\mathcal{N}_l(B) = \bigcup_{b \in B} \mathcal{N}_l(b)$ is called \define{left neighbourhood of $B$}\graffito{left neighbourhood $\mathcal{N}_l(B)$ of $B$}\index[symbols]{NlBcalligraphic@$\mathcal{N}_l(B)$}. \qedhere
    \end{aenumerate}
  \end{definition}

  \begin{definition}
    Let $(X, Y, E)$ be a bipartite graph. It is called
    \begin{aenumerate}
      \item \defineX{finite}{finite bipartite graph}\graffito{finite bipartite graph} if and only if the sets $X$ and $Y$ are finite.
      \item \defineX{locally finite}{locally finite bipartite graph}\graffito{locally finite bipartite graph}\index{finite!locally} if and only if, for each element $x \in X$, the set $\mathcal{N}_r(x)$ is finite, and for each element $y \in Y$, the set $\mathcal{N}_l(y)$ is finite. \qedhere
    \end{aenumerate}
  \end{definition}

  \begin{remark}
    Let $(X, Y, E)$ be a locally finite bipartite graph. Then, for each finite subset $A$ of $X$, the set $\mathcal{N}_r(A)$ is finite; and, for each finite subset $B$ of $Y$, the set $\mathcal{N}_l(B)$ is finite.
  \end{remark}

  \section{Matchings}

  \begin{definition}
    Let $(X, Y, E)$ be a bipartite graph and let $M$ be a subset of $E$. The set $M$ is called \define{matching}\graffito{matching $M$} if and only if, for each tuple $(e, e') \in M \times M$ with $e \neq e'$, the edges $e$ and $e'$ are not adjacent.
  \end{definition}

  \begin{remark}
    The set $M$ is a matching if and only if the maps $p \from M \to X$, $(x, y) \mapsto x$, and $q \from M \to Y$, $(x, y) \mapsto y$, are injective
  \end{remark}

  \begin{definition}
    Let $(X, Y, E)$ be a bipartite graph and let $M$ be a matching. The matching $M$ is called
    \begin{aenumerate}
      \item \define{left-perfect}\graffito{left-perfect matching} if and only if
            \begin{equation*}
              \ForEach x \in X \Exists y \in Y \SuchThat (x, y) \in M;
            \end{equation*}
      \item \define{right-perfect}\graffito{right-perfect matching} if and only if
            \begin{equation*}
              \ForEach y \in Y \Exists x \in X \SuchThat (x, y) \in M;
            \end{equation*}
      \item \define{perfect}\graffito{perfect matching} if and only if it is left-perfect and right-perfect. \qedhere
    \end{aenumerate}
  \end{definition}

  \begin{remark}
    The matching $M$ is left-perfect if and only if the map $p \from M \to X$, $(x, y) \mapsto x$, is surjective (and hence bijective); and it is right-perfect if and only if the map $q \from M \to Y$, $(x, y) \mapsto y$, is surjective (and hence bijective).
  \end{remark}

  \begin{remark} 
    Let $(X, Y, E)$ be a bipartite graph and let $M$ be a subset of $E$. The set $M$ is a
    \begin{aenumerate}
      \item left-perfect matching if and only if there is an injective map $\varphi \from X \to Y$ such that $M = \setOf{(x, \varphi(x)) \suchThat x \in X}$;
      \item right-perfect matching if and only if there is an injective map $\psi \from Y \to X$ such that $M = \setOf{(\psi(y), y) \suchThat y \in Y}$;
      \item perfect matching if and only if there is a bijective map $\varphi \from X \to Y$ such that $M = \setOf{(x, \varphi(x)) \suchThat x \in X}$. \qedhere
    \end{aenumerate}
  \end{remark}

  \section{Hall's Marriage Theorem}

  \begin{definition}
    Let $(X, Y, E)$ be a locally finite bipartite graph. It is said to satisfy the 
    \begin{aenumerate}
      \item \define{left Hall condition}\graffito{left Hall condition} if and only if
            \begin{equation*}
              \ForEach A \subseteq X \text{ finite} \Holds \cardinalityOf{\mathcal{N}_r(A)} \geq \cardinalityOf{A};
            \end{equation*}
      \item \define{right Hall condition}\graffito{right Hall condition} if and only if
            \begin{equation*}
              \ForEach B \subseteq Y \text{ finite} \Holds \cardinalityOf{\mathcal{N}_l(B)} \geq \cardinalityOf{B};
            \end{equation*}
      \item \define{Hall marriage conditions}\graffito{Hall marriage conditions} if and only if it satisfies the left and right Hall conditions. \qedhere
    \end{aenumerate}
  \end{definition}

  \begin{theorem}
    Let $(X, Y, E)$ be a locally finite bipartite graph. It satisfies the left or right Hall condition if and only if there is a left- or right-perfect matching respectively.
  \end{theorem}

  \begin{proof}
    See theorem~H.3.2 in \cite{ceccherini-silberstein:coornaert:2010}.
  \end{proof}

  \begin{theorem}
    Let $(X, Y, E)$ be a bipartite graph such that there is a left-perfect matching and there is a right-perfect matching. Then, there is a perfect matching.
  \end{theorem}

  \begin{proof}
    See theorem~H.3.4 in \cite{ceccherini-silberstein:coornaert:2010}.
  \end{proof}

  \begin{corollary}[Cantor–Schröder–Bernstein theorem; Georg Ferdinand Ludwig Philipp Cantor, Friedrich Wilhelm Karl Ernst Schröder, Felix Bernstein]
    Let $X$ and $Y$ be two sets such that there is an injective map $f$ from $X$ to $Y$ and there is an injective map $g$ from $Y$ to $X$. Then, there is a bijective map from $X$ to $Y$.
  \end{corollary}

  \begin{proof}
    See corollary~H.3.5 in \cite{ceccherini-silberstein:coornaert:2010}.
  \end{proof}

  \begin{theorem}[Hall's marriage theorem; Philip Hall, 1935; Marshall Hall, Jr., 1948]
    Let $(X, Y, E)$ be a locally finite bipartite graph. It satisfies the Hall marriage conditions if and only if there is a perfect matching.
  \end{theorem}

  \begin{proof}
    See theorem~H.3.6 in \cite{ceccherini-silberstein:coornaert:2010}.
  \end{proof}

  \section{Hall's Harem Theorem}

  \begin{definition}
    Let $X$ and $Y$ be two sets, let $f$ be a map from $X$ to $Y$, and let $k$ be a positive integer. The map $f$ is called \defineX{$k$-to-$1$ surjective}{surjectivekto1@$k$-to-$1$ surjective}\graffito{$k$-to-$1$ surjective map} if and only if
    \begin{equation*}
      \ForEach y \in Y \Holds \cardinalityOf{f^{-1}(y)} = k. \qedhere
    \end{equation*}
  \end{definition}

  \begin{definition}
    Let $(X, Y, E)$ be a bipartite graph, let $k$ be a positive integer, and let $M$ be a subset of $E$. The set $M$ is called \graffito{perfect $(1,k)$-matching}\define{perfect $(1,k)$-matching} if and only if
    \begin{equation*}
      \ForEach x \in X \Holds \cardinalityOf{\setOf{y \in Y \suchThat (x, y) \in E}} = k,
    \end{equation*}
    and
    \begin{equation*}
      \ForEach y \in Y \Holds \cardinalityOf{\setOf{x \in X \suchThat (x, y) \in E}} = 1. \qedhere
    \end{equation*}
  \end{definition}

  \begin{remark}
    The set $M$ is a perfect $(1,k)$-matching if and only if there is a $k$-to-$1$ surjective map $\psi \from Y \to X$ such that $\setOf{(\psi(y), y) \suchThat y \in Y} = M$.
  \end{remark}

  \begin{remark}
    The set $M$ is a perfect $(1,1)$-matching if and only if it is a perfect matching.
  \end{remark}

  \begin{definition}
    Let $(X, Y, E)$ be a locally finite bipartite graph and let $k$ be a positive integer. The graph $(X, Y, E)$ is said to satisfy the \define{Hall $k$-harem conditions}\graffito{Hall $k$-harem conditions} if and only if, for each finite subset $A$ of $X$, we have $\cardinalityOf{\mathcal{N}_r(A)} \geq k \cardinalityOf{A}$, and for each finite subset $B$ of $Y$, we have $\cardinalityOf{\mathcal{N}_l(B)} \geq k^{-1} \cardinalityOf{B}$.
  \end{definition}

  \begin{theorem}[Hall's harem theorem]
  \label{theorem:hall-harem}
    Let $(X, Y, E)$ be a locally finite bipartite graph and let $k$ be a positive integer. The graph $(X, Y, E)$ satisfies the Hall $k$-harem conditions if and only if there is a perfect $(1,k)$-matching.
  \end{theorem}

  \begin{proof}
    See theorem~H.4.2 in \cite{ceccherini-silberstein:coornaert:2010}.
  \end{proof}

  \clearToOddPage
  \manualmark
  \markboth{\spacedlowsmallcaps{\bibname}}{\spacedlowsmallcaps{\bibname}} 
  \refstepcounter{dummy}
  \addtocontents{toc}{\protect\vspace{\beforebibskip}} 
  \addcontentsline{toc}{chapter}{\tocEntry{\bibname}}
  \label{app:bibliography}
  \printbibliography

  \clearToOddPage
  \manualmark
  \markboth{\spacedlowsmallcaps{Symbols}}{\spacedlowsmallcaps{Symbols}} 
  \refstepcounter{dummy}
  \addtocontents{toc}{\protect\vspace{\beforebibskip}} 
  \addcontentsline{toc}{chapter}{\tocEntry{Symbols}} 
  \label{app:symbols}
  \printindex[symbols]

  \clearToOddPage
  \manualmark
  \markboth{\spacedlowsmallcaps{\indexname}}{\spacedlowsmallcaps{\indexname}} 
  \refstepcounter{dummy}
  \addtocontents{toc}{\protect\vspace{\beforebibskip}} 
  \addcontentsline{toc}{chapter}{\tocEntry{Index}}
  \label{app:index}
  \printindex

\end{document}